\documentclass{book}

\usepackage{amsmath}
\usepackage{accents}
\usepackage[]{algorithm2e}
\usepackage{amsthm}
\usepackage{amssymb}
\usepackage{array} 
\usepackage{bm}
\usepackage{booktabs}
\usepackage{enumerate}
\usepackage{fancyhdr}
\usepackage{kbordermatrix}
\usepackage{mathtools}
\usepackage{multirow}
\usepackage{nag}
\usepackage{natbib}
\usepackage{paralist}
\usepackage{pgfplots}
\usepackage{rotating}
\usepackage{subfigure}
\usepackage{supertabular}
\usepackage{tikz}
\usepackage{tkz-graph}
\usepackage{tkz-berge}
\usepackage{url}
\usepackage{verbatim}

\usetikzlibrary{arrows,arrows.meta}
\usetikzlibrary{matrix,calc,positioning}

\newtheorem{ass}{Assumption}[chapter]
\newtheorem{exm}{Example}[chapter]
\newtheorem{dfn}{Definition}[chapter]
\newtheorem{cor}{Corollary}[chapter]
\newtheorem{lem}{Lemma}[chapter]
\newtheorem{prp}{Proposition}[chapter]
\newtheorem{thm}{Theorem}[chapter]

\renewcommand{\kbldelim}{\{}
\renewcommand{\kbrdelim}{\}}

\newcommand{\mb}[1]{\mathbf{#1}} 
\newcommand{\mr}[1]{\mathrm{#1}} 
\newcommand{\mc}[1]{\mathcal{#1}} 
\newcommand{\mbb}[1]{\mathbb{#1}} 
\newcommand{\bbr}{\mbb{R}} 
\newcommand{\bbz}{\mbb{Z}} 
\newcommand{\unit}[1]{\ensuremath{\, \mathrm{#1}}} 

\newcommand{\ubar}[1]{\underaccent{\bar}{#1}}

\newcommand{\myp}[1]{\left( #1 \right)} 
\newcommand{\mys}[1]{\left[ #1 \right]} 
\newcommand{\myc}[1]{\left\{ #1 \right\}} 

\newcommand{\pdr}[2]{\frac{\partial #1}{\partial #2}} 
\newcommand{\pdrb}[2]{\frac{\partial^2 #1}{\partial #2^2}} 
\newcommand{\pdrc}[3]{\frac{\partial^2 #1}{\partial #2 \partial #3}} 
\newcommand{\ttt}{\texttt} 
\newcommand{\ds}{\displaystyle} 
\newcommand{\vect}[1]{\begin{bmatrix} #1 \end{bmatrix}} 
\newcommand{\minivect}[1]{\mys{\begin{matrix} #1 \end{matrix}}} 
\newcommand{\labeleqn}[2]{\begin{equation} #2 \label{eqn:#1} \end{equation}} 
\newcommand{\nolabeleqn}[1]{\begin{equation} #1 \end{equation}} 
\newcommand{\eqn}[1]{\eqref{eqn:#1}} 
\newcommand{\genfig}[4]{\begin{figure} \centering \includegraphics[#4]{#2} \caption{#3\label{fig:#1}} \end{figure}}
\newcommand{\stevefig}[3]{\begin{figure} \centering \includegraphics[width=#3]{#1.pdf} \caption{#2\label{fig:#1}} \end{figure}}
\newcommand{\citeapos}[1]{\citeauthor{#1}'s (\citeyear{#1})}

\newcommand{\solution}[1]{\textbf{Solution.} #1 $\blacksquare$}
\newcommand{\optimize}[1]{\[ \begin{array}{rrl} \ds #1 \end{array} \]} 
 
\newcommand{\optimizex}[4]{{ \everymath={\displaystyle} \[ 
                              \begin{array}{rrl}
                                #1_{#2} \,\, & \multicolumn{2}{l}{#3}  \\
                                \mathrm{s.t.} \,\, #4 \,.
                              \end{array}  \]}}

\newcommand{\diff}[1]{\emph{[#1]}} 

\makeatletter
\def\bbordermatrix#1{\begingroup \m@th
  \@tempdima 4.75\p@
  \setbox\z@\vbox{%
    \def\cr{\crcr\noalign{\kern2\p@\global\let\cr\endline}}%
    \ialign{$##$\hfil\kern2\p@\kern\@tempdima&\thinspace\hfil$##$\hfil
      &&\quad\hfil$##$\hfil\crcr
      \omit\strut\hfil\crcr\noalign{\kern-\baselineskip}%
      #1\crcr\omit\strut\cr}}%
  \setbox\tw@\vbox{\unvcopy\z@\global\setbox\@ne\lastbox}%
  \setbox\tw@\hbox{\unhbox\@ne\unskip\global\setbox\@ne\lastbox}%
  \setbox\tw@\hbox{$\kern\wd\@ne\kern-\@tempdima\left[\kern-\wd\@ne
    \global\setbox\@ne\vbox{\box\@ne\kern2\p@}%
    \vcenter{\kern-\ht\@ne\unvbox\z@\kern-\baselineskip}\,\right]$}%
  \null\;\vbox{\kern\ht\@ne\box\tw@}\endgroup}
\makeatother

\usepackage{appendix}
\usepackage[nottoc]{tocbibind}
\usepackage{setspace}
\usepackage{pdflscape}
\usepackage[hidelinks]{hyperref}
\usepackage{makeidx}

\newcommand{\HRule}{\rule{\linewidth}{0.5mm}}

\makeindex

\begin{document}

\frontmatter

\hypersetup{pageanchor=false}

\begin{titlepage}

\begin{center}

\HRule \\[0.4in]
{ \begin{onehalfspace}
 \huge \bfseries Transportation Network Analysis \\[0.3in]
 \large Volume I: Static and Dynamic Traffic Assignment
 \end{onehalfspace}
 }

\HRule \\[1in]
Stephen D. Boyles\\
Civil, Architectural, and Environmental Engineering\\
The University of Texas at Austin\\[0.4in]

Nicholas E. Lownes\\
Civil and Environmental Engineering\\
University of Connecticut\\[0.4in]

Avinash Unnikrishnan\\
Civil, Construction, and Environmental Engineering\\
The University of Alabama at Birmingham\\

\vfill

Version 1.0\\
January 6, 2025

\end{center}
\end{titlepage}

\hypersetup{pageanchor=true}

\section*{Preface}
This book is the product of more than fifteen years of teaching transportation network analysis, at the The University of Texas at Austin, the University of Washington, the University of Wyoming, the University of Connecticut, the University of Alabama at Birmingham, Portland State University, and West Virginia University.
The project began during a sabbatical visit by the second author to The University of Texas at Austin, and has continued since.
We are also developing a companion set of lecture slides, assignments, and solution sets which will be available upon request.
A second volume, covering transit, freight, and logistics, is also under preparation.

Any help you can offer to improve this text would be greatly appreciated, whether spotting typos, math or logic errors, inconsistent terminology, or any other suggestions about how the content can be better explained, better organized, or better presented. 
We will periodically release updated versions incorporating corrections and new material.

We gratefully acknowledge the support of the National Science Foundation under Grants 1069141/1157294, 1254921, 1562109/1562291, 1636154, 1739964, and 1826320.
Travis Waller (University of New South Wales), Chris Tamp\`{e}re (Katholieke Universiteit Leuven), and Xuegang Ban (University of Washington) hosted visits by the first author to their respective institutions, and provided wonderful work environments where much of this writing was done.

\section*{Target Audience}
\label{sec:prerequisites}

This book is primarily intended for first-year graduate students, but is also written with other potential audiences in mind.
The content should be fully accessible to highly-prepared undergraduate students, and certain specialized topics would be appropriate for advanced graduate students as well.
The text covers a large number of topics, likely more than would be covered in one or two semesters, and would also be useful for self-paced learners, or practitioners who may want in-depth learning on specific topics.
We have included some supplementary material in optional sections, marked with an asterisk, which we believe are interesting, but which can be skipped without loss of continuity.

The most important prerequisites for this book are an understanding of multivariate calculus, and the intellectual maturity to understand the tradeoffs involved in mathematical modeling.
Modeling does not involve one-size-fits-all approaches, and dogma about the absolute superiority of one model or algorithm over another is scarce.
Instead, the primary intent of this book is to present a survey of important approaches to modeling transportation network problems, as well as the context to determine when particular models or algorithms are appropriate for real-world problems.
Readers who can adopt this perspective will gain the most from the book.

Appendix~\ref{chp:mathbackground} covers the mathematics needed for the network models in this book.
Readers of this book will have significantly different mathematical backgrounds, and this appendix is meant to collect the necessary results and concepts in one place.
Depending on your background, some parts of it may need only a brief review, while other parts may be completely new.
  
While this book does not explicitly cover how to program these models and algorithms into a computer, if you have some facility in a programming language, it is highly instructive to try to implement them as you read.
Many network algorithms are tedious to apply by hand or with other manual tools (calculator, spreadsheet).
Computer programming will open the door to applying these models to interesting problems of large, realistic scale.

\section*{Difficulty Scale for Exercises}

Inspired by Donald Knuth's \emph{The Art of Computer Programming}, the exercises are marked with estimates of their difficulty.  The key reference points are:
\begin{description}
 \item[0:] A nearly trivial problem that you should be able to answer without any pencil-and-paper work.
 \item[20:] A straightforward problem that may require a few minutes of effort, but nothing too difficult if you have given the chapter a good read.
 \item[40:] A typical problem requiring some thought, or a few attempts at solution in different ways, but the answer should yield after some dedicated effort.
 \item[60:] A problem of above-average difficulty, where the correct approach is not obvious.  You may require a bit of scratch paper or computer work as you try out different approaches before settling on a solution.
 \item[80:] A highly challenging or involved problem, which may be appropriate as a course project or other long-term study.
 \item[100:] An open problem in the research literature, whose solution would be a substantial advance in the understanding of transportation networks.
\end{description}
\emph{The tens digit indicates the intellectual difficulty of the exercise, while the ones digit indicates the amount of calculation required.}  An exercise rated at 50 may require more cleverness and insight than one ranked at 49, but the ultimate solution is shorter.  Of course, each student will find different problems more challenging than others.

\tableofcontents

\mainmatter

\part{Preliminaries}

\chapter{Introduction to Transportation Networks}
\label{chp:introchapter}

This introductory chapter lays the groundwork for traffic assignment, providing some overall context for transportation planning in Section~\ref{sec:transportationnetworks}.
Some examples of networks in transportation are given in Section~\ref{sec:examplesofnetworks}.
The key idea in traffic assignment is the notion of equilibrium, which is presented in Section~\ref{sec:notionofequilibrium}.
The goals of traffic assignment are described in Section~\ref{sec:introductionassignment}.
Traffic assignment models can be broadly classified as static or dynamic.
Both types of models are described in this book, and Sections~\ref{sec:introductionstatic} and~\ref{sec:introductiondynamic} provide general perspective on these types of models.

\section{Transportation Networks}
\label{sec:transportationnetworks}

\index{planning|(}
\index{transportation planning|see{planning}}
Planning helps ensure that transportation spending and policies are as effective as possible.
As transportation engineers and researchers, we support this process by developing and running models which predict the impact of potential projects or policies --- for instance, what would be the impact on city traffic and emissions if an extra lane was added on a major freeway?
If the toll on a bridge was reduced?
If streetcar lines are installed downtown?
In this way, the benefits of projects can be compared with their costs, and funding and implementation priorities can be established.
Depending on the models used, a variety of measures of effectiveness can be considered, and one may want to know the impacts of a project on mobility, congestion, emissions, equity, toll revenue, transit ridership, infrastructure maintenance needs, or countless other metrics.

This is rather difficult.
To predict ridership on a new transit line with complete accuracy would require knowing how many trips every single possible rider makes, and the decision process each one of these potential riders uses when deciding whether or not to use transit.
Unlike trusses or beams, human beings can behave in ways that are impossible to predict, maddeningly inconsistent, and motivated by a variety of factors difficult to observe (many of which occur at a subconscious level).
Further, most transportation infrastructure lasts for decades, meaning that effective planning must also predict the impact of projects and policies decades into the future.
And so far, we've only considered the pure engineering dimension.
Introduce local, state, and federal politics into the mix, other stakeholders such as neighborhood associations and transit agencies, and a public with a variety of priorities (is it more important to reduce congestion, increase safety, or have livable communities?), and the picture only grows more complicated.
What to do?

\index{mathematical model|see {model, mathematical}}
Enter the mathematical model\index{model!mathematical}.
The purpose of a mathematical model is to translate a complicated, but important, real-world problem into precise, quantitative language that can be clearly and unambiguously analyzed.
By their nature, models cannot account for all of the possible factors influencing planning.
As the famous statistician George Box\index{Box, George} once quipped, ``All models are wrong, but some models are useful.''
A useful model is one which provides enough insight that good decisions can be made.
To do this, a model must capture the most important characteristics of the underlying system; must not require more input data than what is available for calibration; and must not require more time and memory than what available hardware permits.

Further, just because a model is useful does not mean it cannot be improved.
Indeed, this is the goal of transportation researchers around the world.
The usual pattern is to start with a model which is simple, transparent, insightful... and also wrong.
This simple model can then be improved in ways to make it more correct and useful, and this is the general pattern which will be seen in this book.
The first network models you will see are such gross simplifications of reality that you may question whether they can truly be of value in practice.
Perhaps they can, perhaps they can't; but in any case, they form a foundation for more advanced and realistic models which relax the assumptions made earlier on.
  
\textbf{For this reason, as a student of transportation planning, you should always be looking for the assumptions involved in everything you see.}
All models make assumptions which are not entirely correct.
The relevant questions are, how much does this assumption limit the applicability of the model, and how easy would it be to relax this assumption?
If you're looking for a research topic, finding an existing model and relaxing an assumption is often a good approach.
With this book, if you clearly understand all of the assumptions underlying each model, and how they differ from those made in other models introduced, you're 90\% of the way there.

\emph{Networks} are a type of mathematical model which are very frequently used in the study of transportation planning.
\index{network}
This introductory chapter gives a very brief overview of transportation networks, and provides a sketch for the remainder of the book.

\index{static model|see {network,static}}
\index{dynamic model|see {network,dynamic}}
This book covers both \emph{static} and \emph{dynamic} network models\index{network!static}\index{network!dynamic}.
Static models assume that network conditions are at steady-state, while dynamic models represent changes in congestion and demand patterns over the course of several hours or a day.
Static models were the first to be developed historically, and remain the most commonly-used in current transportation planning practice.
Dynamic models are more realistic in portraying congestion, but require more data for calibration and validation, and more computational resources to run.
Solving and interpreting the output of dynamic models is also more difficult.
As research progresses, however, more planners are using dynamic models, particularly for applications when travel conditions are changing rapidly during the analysis period.
This chapter will present a balanced perspective of the advantages and disadvantages of dynamic traffic assignment \emph{vis-\`{a}-vis} static assignment, but one of them is worth mentioning now: dynamic traffic assignment models are inherently mode-specific.

That is, unlike in static assignment (where it is relatively easy to build ``multimodal''\index{multimodal} networks mixing roadway, transit, air, and waterway infrastructure), the vast majority of dynamic traffic assignment models have been specifically tailored to modeling vehicle congestion on roadways.
In recent years, researchers have started integrating other modes into dynamic traffic assignment, and this area is likely to receive more attention in years to come.
However, the congestion model for each mode must be custom-built.
This is at once an advantage (in that congestion in different modes arises from fundamentally different sources, and perhaps ought to be modeled quite differently) and a disadvantage (a ``generic'' dynamic traffic assignment model can only be specified at a very high level).
For this reason, this book will focus specifically on automobile traffic on roadway networks.
This is not meant to suggest that dynamic traffic assignment cannot or should not be applied to other modes, but simply an admission that covering other modes would essentially require re-learning a new theory for each mode.
Developing such theories would make excellent research topics.
\index{planning|)}

\section{Examples of Networks}
\label{sec:examplesofnetworks}

\index{network|(}
Networks are fundamental to the study of large-scale transportation models representing an entire metropolitan area, a state, or multistate regions.
They can be applied in many contexts, including alternatives analysis, developing congestion pricing plans, identifying bottlenecks and critical infrastructure, shipping and freight logistics, multimodal planning, and disaster evacuation planning, to name only a few.
The reason network models are so useful, and so broadly applicable, is because a mathematical network is a simple, compact, and flexible way to represent a large, complicated system.

A network consists of \emph{links}\index{link} and \emph{nodes}\index{node}.
In transportation applications, a link usually represents a means of travel from one point to another: a road segment between two intersections, a bus route between two stops, and so on, as seen in Figure~\ref{fig:netexample}.
The nodes, in turn, are the endpoints of the links.
Links may correspond directly to physical infrastructure, as with the ``roadway links'' in the figure.
Other links, such as the ``bus route links'' simply represent an abstract connection between two points (here, bus stops), without being concerned with the specific physical route between them.
Quite often, nodes are adjacent to multiple links, so a node representing an intersection may adjoin multiple links representing road segments.
Nodes and links may also be more abstract; for instance, links in a multimodal network might represent a transfer from one transport mode to another.
The level of detail in a network varies from application to application.
For multistate freight models, major highways may be the only links, and major cities the only nodes.
For a city's planning model, all major and minor arterials may be included as well.
For a more detailed model, individual intersections may be ``exploded''\index{node!exploded} so that different links represent each turning movement (Figure~\ref{fig:explode}).
Other examples of transportation networks are shown in Table~\ref{tbl:networktypes}.
\index{network|)}
  
\stevefig{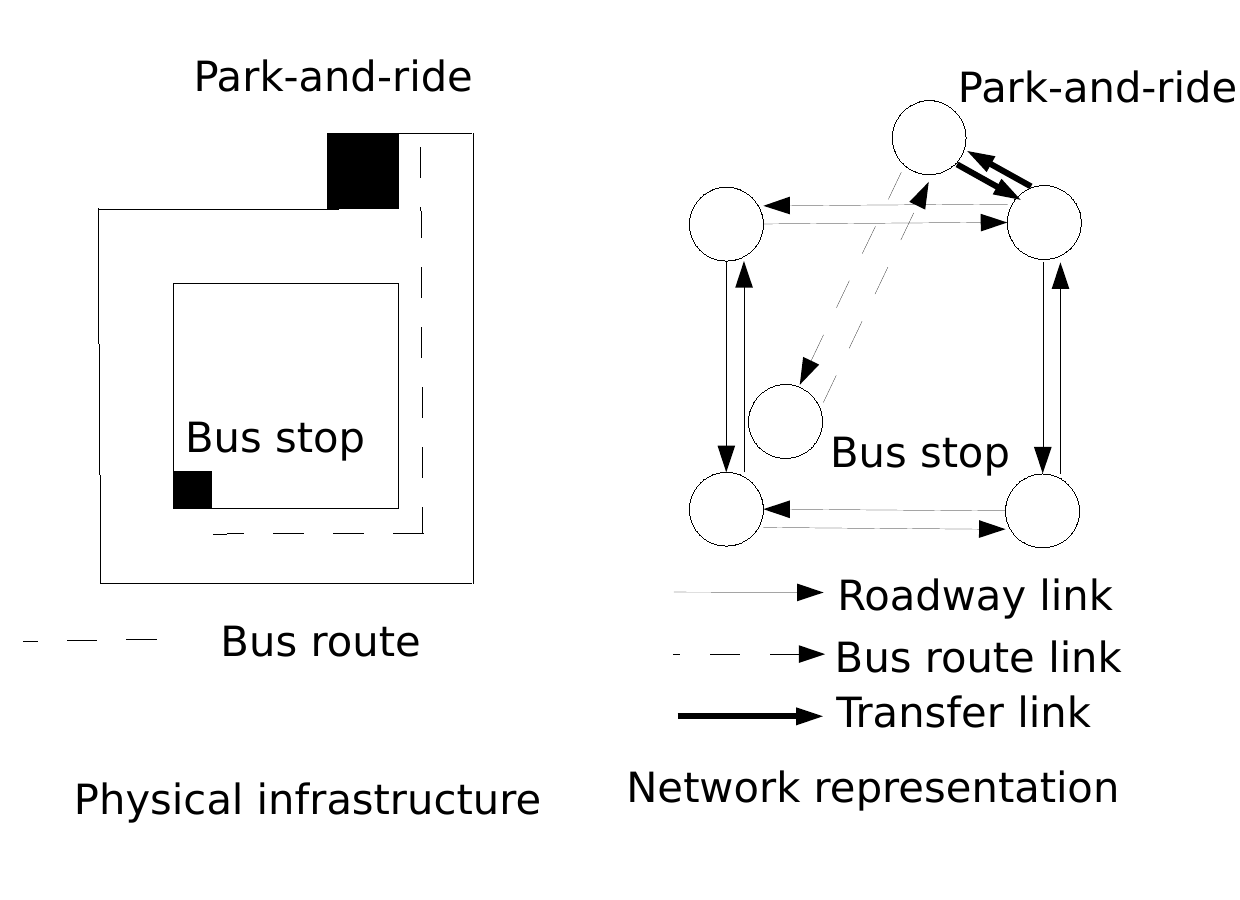}{Nodes and links in transportation networks.}{0.7\textwidth}

\stevefig{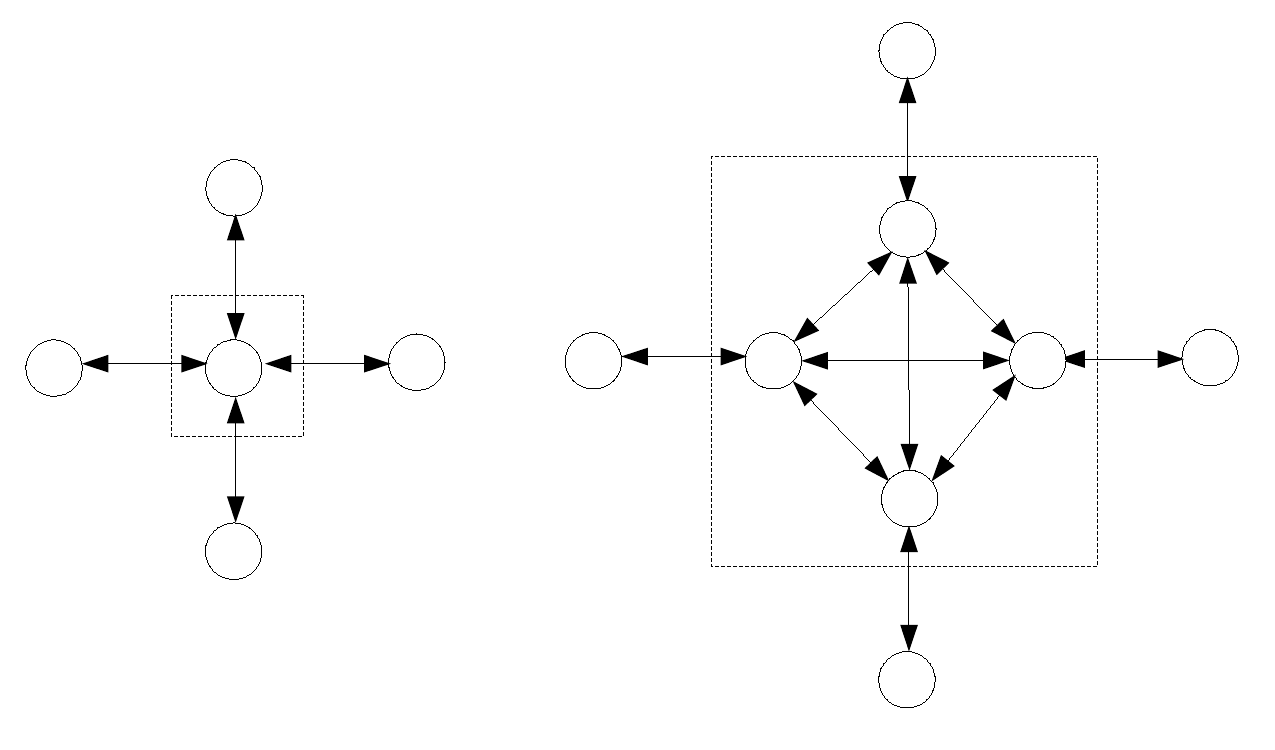}{Two representations of the same intersection.}{0.7\textwidth}

\begin{table}
	\centering
	\caption{Nodes and links in different kinds of transportation networks. \label{tbl:networktypes}}
	\begin{tabular}{lll}
	\hline
	Network type &  Nodes & Links \\
	\hline
	Roadway & Intersections & Street segments \\
	Public transit &  Bus or train stops  & Route segments  \\
	Freight & Factories, warehouses, retailers & Shipping options \\
	Air & Airports & Flights \\
	Maritime & Ports & Shipping channels \\
	\hline
	\end{tabular}
\end{table}

\section{The Notion of Equilibrium}
\label{sec:notionofequilibrium}

\index{equilibrium|(}
The nature of transportation systems is that of multiple interacting systems.
Congestion is determined by the choices travelers make: where, when, how often to travel, and by what mode.
At the same time, these choices depend on congestion: travelers may choose routes or departure times to avoid congestion.
These two ``systems'' (travel choices and system congestion) are thus interdependent and interrelated, with a circular or chicken-and-egg quality to their relationship.
This interdependency lies at the root of transportation analysis.
It is at once interesting, because of the complexity of transportation systems involving both humans and physical systems; challenging, because we must find a way to resolve this circular dependency; and frustrating, because obvious-looking policy interventions can actually be counterproductive.
Some examples of ``paradoxical'' effects will be seen in Chapters~\ref{chp:equilibrium} and~\ref{chp:dynamiceqm}.

\index{travel demand model|see {model, demand}}
\index{demand model|see {model, demand}}
\index{supply model|see {model, supply}}
The schematic in Figure~\ref{fig:taschematic} shows the dependency between the choices made by travelers (sometimes called the \emph{demand side}\index{model!demand side}), and the congestion and delay in the system (sometimes called the \emph{supply side}\index{model!supply side}) in the basic traffic assignment problem.
Each of these systems requires a distinct set of models.
Demand-side models should be behavioral in nature, identifying what factors influence travel choices, and how.
Supply-side models are often based in traffic flow theory, queueing theory, computer simulation, or empirical formulas describing how congestion will form.

\begin{figure}
\begin{center}
\begin{tikzpicture}
[squarednode/.style={rectangle,draw=black!60,fill=white,very thick,minimum size=5mm,
                    text centered,text width=10cm, node distance=2.5cm}]
\node[squarednode](tfm){Traffic flow model};
\node[squarednode](rcm)[below=of tfm]{Route choice model};
\path (rcm.north) -- (rcm.north west) coordinate[pos=0.33] (north-by-northwest);
\path (rcm.north) -- (rcm.north east) coordinate[pos=0.33] (north-by-northeast);
\path (tfm.south) -- (tfm.south west) coordinate[pos=0.33] (south-by-southwest);
\path (tfm.south) -- (tfm.south east) coordinate[pos=0.33] (south-by-southeast);
\draw[-{Latex[length=3mm]}] (north-by-northwest) -- node [left, midway] {Path flows} (south-by-southwest);
\draw[-{Latex[length=3mm]}] (south-by-southeast) -- node [right, midway] {Path travel times} (north-by-northeast);
\end{tikzpicture}
\caption{Traffic assignment as interacting systems.  \label{fig:taschematic}}
\end{center}
\end{figure}

It is not difficult to imagine other types of mutually-dependent transportation systems.
Figure~\ref{fig:complextaschematic} shows how one might model traffic assignment in a region with a privately-operated toll road.\index{tolls}
Now, there are three systems.
In addition to the demand side and supply side from before, the private toll operator can also influence the state of the network by choosing the toll in some way, such as maximizing toll revenue.
But this choice is not made in isolation: as the toll is increased, drivers will choose alternate routes, suggesting that driver choices are affected by tolls just as the toll revenue is determined by driver choices.
It is fruitful to think of other ways this type of system can be expanded.
For instance, a government agency might set regulations on the maximum and minimum toll values, but travelers can influence these policy decisions through the voting process.
  
\begin{figure}
\begin{center}
\begin{tikzpicture}
[squarednode/.style={rectangle,draw=black!60,fill=white,very thick,minimum size=5mm,
                    text centered,text width=10cm, node distance=2.5cm}]
\node[squarednode](tfm){Traffic flow model};
\node[squarednode](rcm)[below=of tfm]{Route choice model};
\node[squarednode](tsm)[below=of rcm]{Toll-setting model};
\path (rcm.north) -- (rcm.north west) coordinate[pos=0.33] (north-by-northwest);
\path (rcm.north) -- (rcm.north east) coordinate[pos=0.33] (north-by-northeast);
\path (tfm.south) -- (tfm.south west) coordinate[pos=0.33] (south-by-southwest);
\path (tfm.south) -- (tfm.south east) coordinate[pos=0.33] (south-by-southeast);
\path (tsm.north) -- (tsm.north west) coordinate[pos=0.33] (north-by-northwest-s);
\path (tsm.north) -- (tsm.north east) coordinate[pos=0.33] (north-by-northeast-s);
\path (rcm.south) -- (rcm.south west) coordinate[pos=0.33] (south-by-southwest-s);
\path (rcm.south) -- (rcm.south east) coordinate[pos=0.33] (south-by-southeast-s);
\draw[-{Latex[length=3mm]}] (north-by-northwest) -- node [left, midway] {Path flows} (south-by-southwest);
\draw[-{Latex[length=3mm]}] (south-by-southeast) -- node [right, midway] {Path travel times} (north-by-northeast);
\draw[-{Latex[length=3mm]}] (north-by-northwest-s) -- node [left, midway] {Selected toll} (south-by-southwest-s);
\draw[-{Latex[length=3mm]}] (south-by-southeast-s) -- node [right, midway] {Revenue} (north-by-northeast-s);
\end{tikzpicture}
\caption{A more complicated traffic assignment problem, with tolls.  \label{fig:complextaschematic}}
\end{center}
\end{figure}

The task of transportation planners is to somehow make useful predictions to assist with policy decision and alternatives analysis, despite the complexities which arise when mutually-dependent systems interact.
The key idea is that a good prediction is \emph{mutually consistent} in the sense that all of the systems should ``agree'' with the prediction.
\index{consistency}
As an example, in the basic traffic assignment problem (Figure~\ref{fig:taschematic}), a planning model will provide both a forecast of travel choices, and a forecast of system congestion.
\index{traffic assignment}
These should be consistent in the sense that inputting the forecasted travel choices into the supply-side model should give the forecast of system congestion, and inputting the forecasted system congestion into the demand-side model should give the forecast of travel choices.
Such a consistent solution is termed an \emph{equilibrium}.

The word equilibrium is meant to allude to the concept of economic equilibrium, as it is used in game theory.\index{game theory}
In game theory, several agents each choose a particular action, and depending on the choices of all of the agents, each receives a payoff (perhaps monetary, or simply in terms of happiness or satisfaction).
Each agent wants to maximize their payoff.
The objective is to find a ``consistent''\index{consistency} or equilibrium solution, in which all of the agents are choosing actions which maximize their payoff (keeping in mind that an agent cannot control another agent's decision).
A few examples are in order.

\index{equilibrium!existence}
Consider first a game with two players (call them Alice and Bob), who happen to live in a small town with only two bars (the Cactus Caf\'{e} and the Desert Drafthouse).
Alice and Bob have recently broken off their relationship, so they each want to go out to a bar.
If they attend different bars, both of them will be happy (signified by a payoff of $+1$), but if they attend the same bar an awkward situation will arise and they will regret having gone out at all (signified by a payoff of $-1$).
Table~\ref{tbl:alicebob} shows the four possible situations which can arise --- each cell in the table lists Alice's payoff first, followed by Bob's.
Two of these are boldfaced, indicating that they are equilibrium solutions: if Alice is at the Cactus Caf\'{e} and Bob at the Desert Drafthouse (or vice versa), they each receive a payoff of $+1$, which is the best they could hope to receive given what the other is doing.
The states where they attend the same bar are not equilibria; either of them would be better off switching to the other bar.
This is a game with two equilibria.\footnote{There is also a third equilibrium in which they each randomly choose a bar each weekend, but equilibria involving randomization are outside the scope of this book.}

\begin{table}
\begin{center}
\caption{Alice and Bob's game; Alice chooses the row and Bob the column. \label{tbl:alicebob}}
\begin{tikzpicture}
\matrix[matrix of math nodes, every node/.style={text width=3cm,align=center},row sep=0.
2cm,column sep=0.2cm] (m) {
$(-1,-1)$ & $\mathbf{(1, 1)}$ \\
$\mathbf{(1,1)}$ & $(-1,-1)$  \\
};
\draw (m.north east) rectangle (m.south west);
\draw (m.north) -- (m.south);
\draw (m.east) -- (m.west);

\coordinate (a) at ($(m.north west)!0.25!(m.north east)$);
\coordinate (b) at ($(m.north west)!0.75!(m.north east)$);
\node[above=5pt of a,anchor=base] {Cactus Caf\'{e}};
\node[above=5pt of b,anchor=base] {Desert Drafthouse};

\coordinate (c) at ($(m.north west)!0.25!(m.south west)$);
\coordinate (d) at ($(m.north west)!0.75!(m.south west)$);
\node[left=2pt of c,text width=3cm]  {Cactus Caf\'{e}};
\node[left=2pt of d,text width=3cm]  {Desert Drafthouse};

\node[left=4cm of m.west,align=center,anchor=center] {Alice};
\node[above=10pt of m.north] (firm b) {Bob};

\end{tikzpicture}
\end{center}
\end{table}

\index{equilibrium!efficiency}
\index{prisoners' dilemma|(}
A second game involves the tale of Erica and Fred, two criminals who have engaged in a decade-long spree of major art thefts.
They are finally apprehended by the police, but for a minor crime of shoplifting a candy bar from the grocery store.
The police suspect the pair of the more serious crimes, but have no hard evidence.
So, they place Erica and Fred in separate jail cells.
They approach Erica, offering her a reduced sentence in exchange for testifying against Fred for the art thefts, and separately approach Fred, offering him a reduced sentence if he would testify against Erica.
If they remain loyal to each other, they will be convicted only of shoplifting and will each spend a year in jail.
If Erica testifies against Fred, but Fred stays silent, then Fred goes to jail for 15 years while Erica gets off free.
(The same is true in reverse if Fred testifies against Erica.) 
If they both testify against each other, they will both be convicted of the major art theft, but will have a slightly reduced jail term of 14 years for being cooperative.
This game is diagrammed in Table~\ref{tbl:ericafred}, where the ``payoff'' is the negative of the number of years spent in jail, negative because more years in jail represents a worse outcome.
Surprisingly, the only equilibrium solution is for both of them to testify against each other.
From Erica's perspective, she is better off testifying against Fred \emph{no matter what Fred will do}.
If he is going to testify against her, she can reduce her sentence from 15 years to 14 by testifying against Fred.
If he is going to stay silent, she can reduce her sentence from one year to zero by testifying against him.
Fred's logic is exactly the same.
This seemingly-paradoxical result, known as the \emph{prisoner's dilemma}, shows that agents maximizing their own payoff can actually end up in a very bad situation when you look at their combined payoffs!  
\index{prisoners' dilemma|)}

\begin{table}
\begin{center}
\caption{Erica and Fred's game.\label{tbl:ericafred}}
\begin{tikzpicture}
\matrix[matrix of math nodes, every node/.style={text width=3cm,align=center},row sep=0.2cm,column sep=0.2cm] (m) {
$\mathbf{(-14,-14)}$ & $ (0,-15) $ \\
$(-15,0)$ & $(-1,-1)$   \\
};
\draw (m.north east) rectangle (m.south west);
\draw (m.north) -- (m.south);
\draw (m.east) -- (m.west);

\coordinate (a) at ($(m.north west)!0.25!(m.north east)$);
\coordinate (b) at ($(m.north west)!0.75!(m.north east)$);
\node[above=5pt of a,anchor=base] {Testify};
\node[above=5pt of b,anchor=base] {Remain silent};

\coordinate (c) at ($(m.north west)!0.25!(m.south west)$);
\coordinate (d) at ($(m.north west)!0.75!(m.south west)$);
\node[left=2pt of c,text width=2.2cm]  {Testify};
\node[left=2pt of d,text width=2.2cm]  {Remain silent};

\node[left=3cm of m.west,align=center,anchor=center] {Erica};
\node[above=10pt of m.north] (firm b) {Fred};

\end{tikzpicture}
\end{center}
\end{table}

\index{equilibrium!uniqueness}
A third game, far less dramatic than the first two, involves Ginger and Harold, who are retirees passing the time by playing a simple game.
Each of them has a penny, and on a count of three each of them chooses to reveal either the head or the tail of their penny.
If the pennies show the same (both heads or both tails), Ginger keeps them both.
If one penny shows heads and the other shows tails, Harold keeps them both.
(Table~\ref{tbl:gingerharold}).
In this case, there is \emph{no} equilibrium solution: if Ginger always shows heads, Harold will learn and always show tails; once Ginger realizes this, she will start showing tails, and so on \emph{ad infinitum}.

\begin{table}
\begin{center}
\caption{Ginger and Harold's game.  \label{tbl:gingerharold}}
\begin{tikzpicture}
\matrix[matrix of math nodes, every node/.style={text width=3cm,align=center},row sep=0.2cm,column sep=0.2cm] (m) {
$(+1,-1)$             & $(-1,+1)$         \\
$(-1,+1)$             & $(+1,-1)$         \\
};
\draw (m.north east) rectangle (m.south west);
\draw (m.north) -- (m.south);
\draw (m.east) -- (m.west);

\coordinate (a) at ($(m.north west)!0.25!(m.north east)$);
\coordinate (b) at ($(m.north west)!0.75!(m.north east)$);
\node[above=5pt of a,anchor=base] {Heads};
\node[above=5pt of b,anchor=base] {Tails};

\coordinate (c) at ($(m.north west)!0.25!(m.south west)$);
\coordinate (d) at ($(m.north west)!0.75!(m.south west)$);
\node[left=2pt of c,text width=1.5cm]  {Heads};
\node[left=2pt of d,text width=1.53cm]  {Tails};

\node[left=2.5cm of m.west,align=center,anchor=center] {Ginger};
\node[above=10pt of m.north] (firm b) {Harold};

\end{tikzpicture}
\end{center}
\end{table}

You may be wondering how these games are relevant to transportation problems.
In fact, the route choice decision can be seen as a game with a very large number of players.
Some drivers may choose to avoid the freeway, anticipating a certain level of congestion and trying to second-guess what others are doing --- but surely other drivers are engaging in the exact same process.\footnote{To borrow from Yogi Berra, nobody takes the freeway during rush hour anymore --- it's too congested.}
Each of these three games has bearing on the traffic assignment problem.
The game with Alice and Bob shows that some games have more than one equilibrium solution (an issue of equilibrium \emph{uniqueness}).
What does it mean for transportation planning if a model can give several seemingly valid predictions?  The game with Erica and Fred shows that agents individually doing what is best for themselves may lead to an outcome which is quite bad overall, an issue of equilibrium \emph{efficiency}.
As we will see later on, in transportation systems this opens the door for seemingly helpful projects (like capacity expansion on congested roads) to actually make things worse.
The game with Ginger and Harold is a case where there is no equilibrium at all (an issue of equilibrium  \emph{existence}).
If this could happen in a transportation planning model, then perhaps equilibrium is the wrong concept to use altogether.
These questions of uniqueness, efficiency, and existence are important, and will appear throughout the book.

The three example games described above can be analyzed directly, by enumerating all the possible outcomes.
However, transportation systems involve thousands or even millions of different ``players'' and an analysis by enumeration is hopeless.
The good news is that the number of players is so great that little is lost in assuming that the players can be treated as a continuum.\index{model!continuous}\footnote{This is analogous to solving structural design problems by assuming the usual stress-strain relationships, which assume a continuous material.
In reality, a beam or column is composed of many distinct atoms, not a homogeneous material --- but surely it is unnecessary to model each atom separately.
The continuum assumption works almost as well and is much, much easier to work with.}
This allows us to work with smooth functions, greatly simplifying the process of finding equilibria.
The remainder of this chapter introduces the basic traffic assignment problem in terms of the equilibrium concept and with a few motivating examples, but still in generally qualitative terms and restricted to small networks.
The following two chapters provide us with the mathematical vocabulary and network tools needed to formulate and solve equilibrium on realistic, large-scale systems.
\index{equilibrium|)}

\section{Traffic Assignment}
\label{sec:introductionassignment}

There are many possible measures of effectiveness\index{planning!measures of effectiveness} for evaluating the impacts of a roadway transportation project or policy.
However, many of these can be calculated if one can predict the number of drivers on each roadway segment.
These are called \emph{link flows}.\index{link!flow}
Predicting link flows allows a city or state government to evaluate different options.

If link flows are the output of a planning model, the main input is demographic data.
That is, \emph{given certain information about a population (number of people, income, amount of employment, etc.), we want to predict how many trips they will make, and how they will choose to travel.}
Census records form an invaluable resource for this, often supplemented with travel surveys.
Commonly, a medium-to-large random sample of the population is offered some money in exchange for keeping detailed diaries indicating all of the trips made within the next several weeks, including the time of day, reason for traveling, and other details.
  
\index{four-step model|see {planning, four-step model}}
\index{model!four-step|see {planning, four-step model}}
To get link flows from demographic data, most planners use the so-called \emph{four-step model}\index{planning!four-step model} (Figure~\ref{fig:fourstep}).
The first step is \emph{trip generation}\index{trip generation}: based on demographic data, how many trips will people make?
The second is \emph{trip distribution}\index{trip distribution}: once we know the total number of trips people make, what are the specific locations people will travel to?
The third is \emph{mode choice}\index{mode choice}: once we know the trip locations, will people choose to drive, take the bus, or use another mode?
The fourth and final step is \emph{route choice}\index{route choice|see {traffic assignment}}, also known as \emph{traffic assignment}: once we know the modes people will take to their trip destinations, what routes will they choose?
Thus, at the end of the four steps, the transition from demographic data to link flows has been accomplished.\footnote{In more sophisticated models, the four steps may be repeated again, to ensure that the end results are consistent with the input data.
There are also newer and arguably better alternatives to the four-step model.
}
\index{consistency}

\stevefig{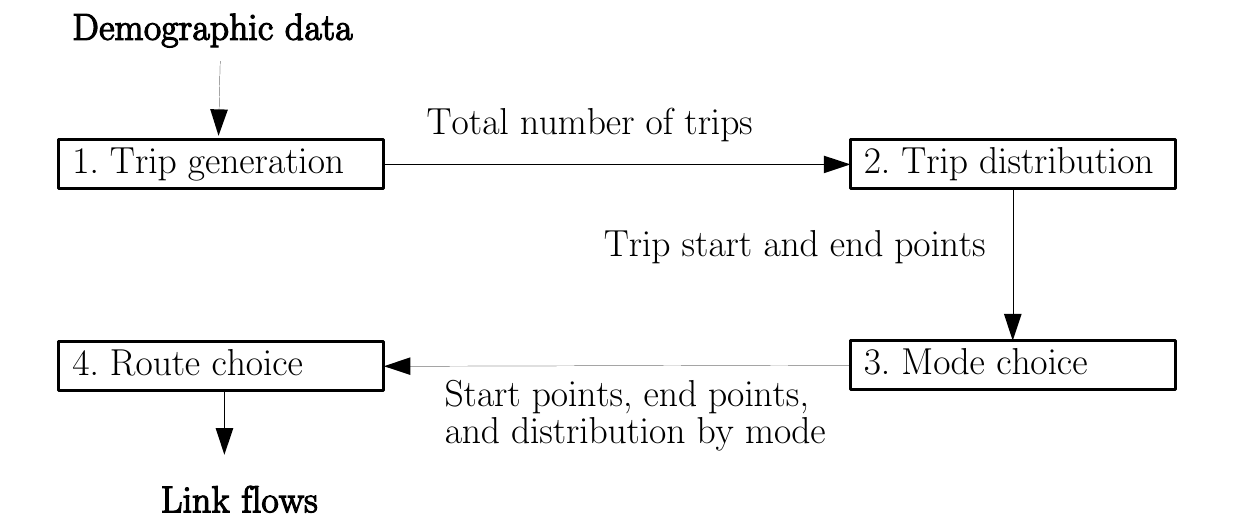}{Schematic of the four-step process.}{\textwidth}

Demographics are not uniform in a city; some areas are wealthier than others, some areas are residential while others are commercial, some parts are more crowded while other parts have a lower population density.
For this reason, planners divide a city into multiple \emph{zones}\index{zone}, and assume certain parameters within each zone are uniform.
Clearly this is only an approximation to reality, and the larger the number of zones, the more accurate the approximation.
(At the extreme, each household would be its own zone and the uniformity assumption becomes irrelevant.)
On the other hand, the more zones, the longer it takes to run each model, and at some point computational resources become limiting.
Typical networks used for large metropolitan areas have a few thousand zones.
Zones are often related to census tracts, to make it easy to get demographic information from census results.

The focus of this book is the last of the four steps, traffic assignment.
In the beginning, we assume that the first three steps have been completed, and we know the number of drivers traveling between each origin zone and destination zone.
From this, we want to know how many drivers are going to use each roadway segment, from which we can estimate congestion, emissions, toll revenue, or other measures of interest.\index{planning!measures of effectiveness}

We've already discussed several of the pieces of information we need in order to describe traffic assignment models precisely, including zones and travel demand.
The final piece of the puzzle is a representation of the transportation infrastructure itself: the transportation network described more in the next chapter.

\index{centroid|see {node, centroid}}
\index{centroid connector|see {link, centroid connector}}
It is usually convenient to use a node to represent each zone; such nodes are called \emph{centroids}\index{node!centroid}, and all trips are assumed to begin and end at centroids.
The set of centroids is thus a subset of the set of nodes, defined in the next chapter.
Centroids may coincide with physical nodes in the network.
Centroids may also represent artificial nodes which do not correspond to any one physical point, and are connected to the physical infrastructure with links called \emph{centroid connectors}\index{link!centroid connector} (dashed lines in Figure~\ref{fig:centroid}).

\stevefig{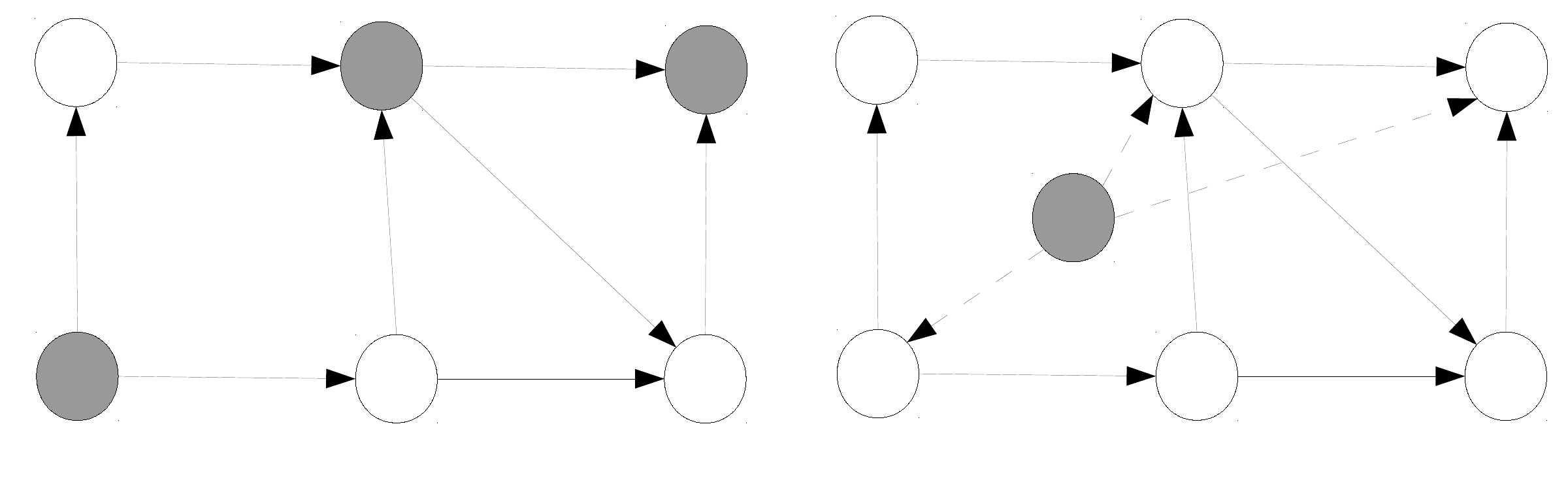}{Centroids (shaded) coinciding with existing infrastructure, and artificial centroids.  Dashed links on the right represent centroid connectors.}{\textwidth}

\section{Static Traffic Assignment}
\label{sec:introductionstatic}

\index{traffic assignment!static|see {static traffic assignment}}
\index{static traffic assignment|(}
Figure~\ref{fig:taschematic} is the template for all traffic assignment models, be they static or dynamic: the choices of travelers lead to congestion patterns in the network (as predicted by a traffic flow model), and these patterns in turn influence the choices travelers make.
The difference between static and dynamic traffic assignment lies in the traffic flow models used.
Historically, static assignment models were the first to be developed, and research into dynamic models arose from the need to improve earlier, static models.
Dynamic traffic assignment thus has a large number of parallels with static assignment; but where they differ, this difference is often intentional and important.
Understanding these distinctions is key to knowing when dynamic models are appropriate to use.

\subsection{Overview}

\index{BPR function|see {link performance function, BPR}}
In static assignment, the traffic flow model is based on \emph{link performance functions}\index{link performance function}, which map the flow on each link to the travel time on that link.
Mathematically, if the notation $(i,j)$\label{not:ij} is used to refer to a roadway link connecting two nodes $i$\label{not:i} and $j$\label{not:j}, then $x_{ij}$\label{not:xij} is the flow on link $(i,j)$ and the function $t_{ij}(x_{ij})$\label{not:tij} gives the travel time on link $(i,j)$ as a function of the flow on $(i,j)$.
These functions $t_{ij}(\cdot)$ are typically assumed to be nonnegative, nondecreasing, and convex, reflecting the idea that as more vehicles attempt to drive on a link, the greater the congestion and the higher the travel times will be.
A variety of link performance functions exist, but the most popular is the Bureau of Public Roads\index{BPR (Bureau of Public Roads)}\index{link performance function!BPR (Bureau of Public Roads)} (BPR) function, which takes the form
\labeleqn{bprdemo}{t_{ij}(x_{ij}) = t^0_{ij} \left(1 + \alpha \left( \frac{x_{ij}}{u_{ij}} \right) ^\beta \right)}
where $t^0_{ij}$\label{not:t0ij} and $u_{ij}$\label{not:uij} are the free-flow time and capacity of link $(i,j)$, respectively, and $\alpha$\label{not:alpha} and $\beta$\label{not:beta} are shape parameters which can be calibrated to data.
It is common to use $\alpha = 0.15$ and $\beta = 4$, but see Section~\ref{sec:tapdata} for more discussion on how to choose these parameters.

With such functions, the more travelers choose a path, the higher its travel time will be.
Since travelers seek to minimize their travel time, travelers will not choose a path with high travel time unless there is no other option available.
Indeed, if travelers only choose paths to minimize travel time, and if they have perfect knowledge of network conditions, then the network state can be described by the \emph{principle of user equilibrium}\index{user equilibrium}: all used paths between the same origin and destination have equal and minimal travel times, for if this were not the case travelers would switch from slower routes to faster ones, which would tend to equalize their travel times.
\index{equilibrium!user|see {user equilibrium}}

It is not difficult to show that this \emph{user equilibrium} state is not socially optimal, and that other configurations of traffic flow can reduce the average travel time (or even the travel time for all drivers) compared to the user equilibrium state.
In other words, individual drivers seeking to minimize their own travel times will not always minimize travel times throughout the network, and this latter \emph{system optimal}\index{system optimum} state can be contrasted with the user equilibrium one.

The prime advantage of using link performance functions like that in equation~\eqn{bprdemo} is that the user equilibrium and system optimum states can be found with relative ease, even in realistic networks involving tens of thousands of links.
Part~\ref{part:statictrafficassignment} of the book discusses this in detail, showing how the static assignment problem can be formulated using the mathematical tools of optimization, fixed point, and variational inequality problems.
These three representations of the equilibrium problem can be linked to powerful mathematical results which assure the existence and uniqueness of user equilibrium solutions under mild conditions.
 Efficient algorithms allow these states to be identified in a matter of minutes on large-scale networks.
  
For these reasons, static traffic assignment has been widely used in transportation planning practice for decades, and remains a powerful tool that can be used for performing alternatives analysis and generating forecasts of network conditions.

\subsection{Critique}

There are also a number of serious critiques of static assignment models, focused primarily on the link performance functions.
By definition, static models do not monitor how network conditions (either demand or congestion) change over time, and implicitly assume a steady-state condition.
This is clearly not the case.
There are additional, subtler and more fundamental problems with link performance functions.
This section describes a few of these problems.

First, not all vehicles on a link experience the same travel time.
Even if the demand for travel on a link exceeds the capacity, the first vehicles to travel on that link will not experience much delay at all, while vehicles which arrive later may experience a very high delay.
Together with the principle of user equilibrium, this means that the paths chosen by travelers will also depend on when they are departing.
Route choices during periods of high congestion will be different from route choices made while this congestion is forming or dissipating.
Furthermore, the travel time faced by a driver on a link depends crucially on the vehicles in front of them, and very little on the vehicles behind them.
(A driver must adjust their speed to avoid colliding with vehicles downstream; a driver has no obligation to change their speed based on vehicles behind them.)
This asymmetry is known as the \emph{anisotropic}\index{anisotropic property} property of traffic flow\index{traffic flow theory}, and it is violated by link performance functions --- an increase in flow on the link is assumed to increase the travel time of all vehicles, and directly using link performance functions in a dynamic model would lead to undesirable phenomena, such as vehicles entering a link immediately raising the travel time for all other vehicles on the link, even those at the downstream end.

\index{flow!in static traffic assignment|see {link, flow}}
\index{capacity|see {link, capacity}}
\index{volume!in static traffic assignment|see {link, flow}}
\index{link!volume|see {link, flow}}
Second, the use of the word ``flow'' in static assignment is problematic.
In traffic engineering, \emph{flow}\index{link!flow} is defined as the (time) rate at which vehicles pass a fixed point on the roadway, and \emph{capacity}\index{link!capacity} is the greatest possible value for flow.
By definition, flow cannot exceed capacity.
However, the BPR function~\eqn{bprdemo} imposes no such restriction, and it is common to see ``flow-to-capacity'' ratios much greater than one in static assignment.\footnote{There are several ways to add such a restriction, but these are less than satisfactory.
A link performance function which tends to $+\infty$ as capacity is reached introduces numerical issues in solving for equilibrium.
Explicitly adding link capacity constraints to the traffic assignment problem may make the problem infeasible, during peak periods there may be no way to assign all vehicles to the network without (temporarily) exceeding capacity.}
Instead, the $x_{ij}$ values in static assignment are better thought of as \emph{demand}\index{link!demand} rather than actual flow, since there is no harm in assuming that the demand for service exceeds the capacity, but it is impossible for the flow itself to exceed capacity.
And for purposes of calibration, demand is much harder to observe than actual flow.
These issues do not have clean resolutions.

Third, and related to the previous issue, link performance functions suggest that lower-capacity links have higher travel times under the same demand.
But consider what happens at a freeway bottleneck, such as the lane drop shown in Figure~\ref{fig:lanedrop}.
Congestion actually forms \emph{upstream} of a bottleneck, and downstream of the lane drop there is no reason for vehicles to flow at a reduced speed.
In reality, it is upstream links that suffer when the demand for traveling on a link exceeds the capacity, not the bottleneck link itself.

\stevefig{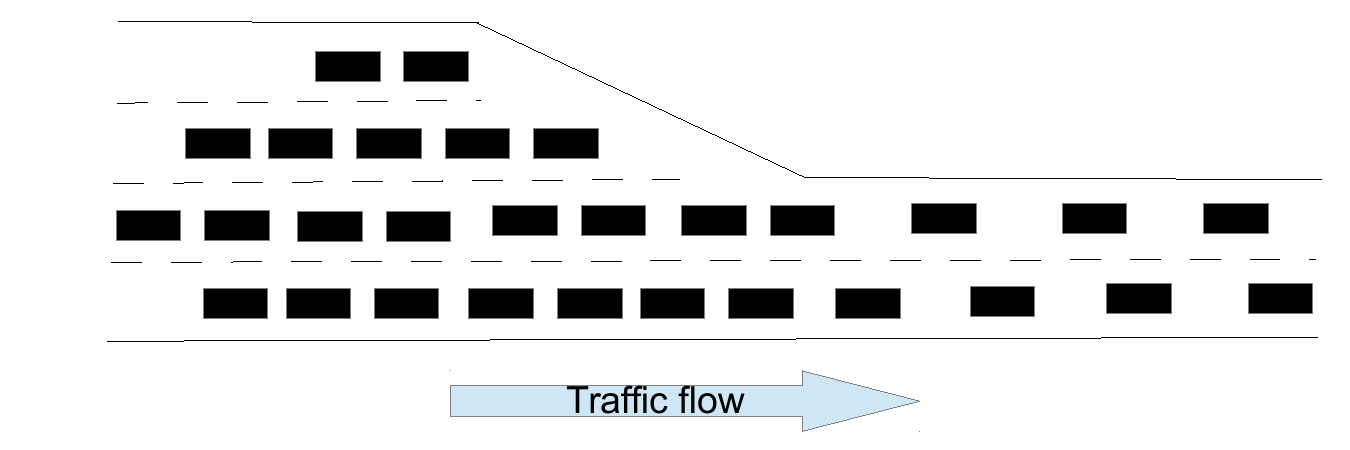}{Congestion arises upstream of a bottleneck, not on the link with reduced capacity.}{0.9\textwidth}

\index{spillback|see {queue spillback}}
Fourth, in congested urban systems it is very common for queues to fill the entire length of a link, causing congestion to spread to upstream links.\index{queue spillback}
This is observed on freeways (congested offramp queues) and in arterials (gridlock in central business districts) and is a major contributor to severe delay.
In addition to the capacity, which is a maximum flow rate, real roadways also have a \emph{jam density}\index{density!jam}, a maximum spatial concentration of vehicles.
If a link is at jam density, no more vehicles can enter, which will create queues on upstream links.
If these queues continue to grow, they will spread even further upstream to other links.
Capturing this phenomenon can greatly enhance the realism of traffic models.
\index{jam density|see {density, jam}}
 
For all of these reasons, dynamic traffic assignment models shift to an entirely different traffic flow model.
Rather than assuming simple, well-behaved link performance functions for each link we turn to traffic flow theory, to find more realistic ways to link traffic flow to congestion.
Some early research in dynamic traffic assignment attempted to retain the use of link performance functions --- for instance, modeling changes in demand by running a sequence of static assignment models over shorter time periods, with methods for linking together trips spanning multiple time periods.
While this addresses the obvious shortcoming of static models, that they cannot accommodate changes in demand or model changes in congestion over time, it does nothing to address the more serious and fundamental problems with link performance functions described above.
For this reason, this approach has largely been abandoned in favor of entirely different traffic flow models.
\index{static traffic assignment|)}

\section{Dynamic Traffic Assignment}
\label{sec:introductiondynamic}

\index{dynamic traffic assignment|(}
\index{traffic assignment!dynamic|see {dynamic traffic assignment}}

Dynamic traffic assignment arose from efforts to resolve the problems with static assignment noted in the previous section.
While it is being used more and more in practice, it has not completely supplanted static traffic assignment.
This is partially due to institutional inertia, but is also due to a few drawbacks associated with more realistic traffic flow models.
Both these advantages and drawbacks are discussed in this section.

\subsection{Overview}

\index{network!dynamic}
\index{traffic flow theory}
The biggest difference between static and dynamic traffic assignment is in the traffic flow models used.
A number of different traffic flow models are available, and a number of them are discussed in the following chapter.
At a minimum, a traffic flow model for dynamic traffic assignment must be able to track changes in congestion at a fairly fine resolution, on the order of a few seconds.
To do this, the locations of vehicles must be known at the same resolution.
Most of them also address some or all of the shortcomings of link performance functions identified above.

\index{dynamic user equilibrium|see {dynamic traffic assignment, equilibrium}}
\index{user equilibrium!dynamic|see {dynamic traffic assignment, equilibrium}}
The behavior of drivers is similar to that in static assignment in that drivers choose paths with minimum travel time.
However, since congestion (and therefore travel time) changes over time, the minimum-time paths also change with time.
Therefore, the principle of user equilibrium is replaced with a principle of \emph{dynamic user equilibrium}: All paths used by travelers departing the same origin, for the same destination, \emph{at the same time}, have equal and minimal travel times.
Travelers departing at different times may experience different travel times; all the travelers departing at the same time will experience the same travel time at equilibrium, regardless of the path they choose.
By virtue of representing demand changes over time, dynamic traffic assignment can also incorporate departure time choice of drivers, as well as route choice.
Some dynamic traffic assignment models include both of these choices, while others focus only on route or departure time choice.\index{departure time choice}
Which choices are appropriate to model depends on which is more significant for a particular application, as well as the data and computational resources available.
Most chapters of this book focus only on route choice, and as a general rule we will assume that departure times are fixed.
In a few places we show how departure time choice can be added in.

It is important to specify that this equilibrium is based on the travel times the drivers actually experience, not the ``instantaneous'' travel times\index{instantaneous travel time} at the moment they depart.
That is, we do not simply assume that drivers will pick the fastest route based on current conditions (as would be provided by most advanced traveler information services), but that they will anticipate changes in travel times which will occur during their journey.
This suggests that drivers are familiar enough with the network that they know how congestion changes with time.
This distinction is important --- dynamic traffic assignment equilibrates on experienced travel times,\index{experienced travel time} not instantaneous travel times.

\stevefig{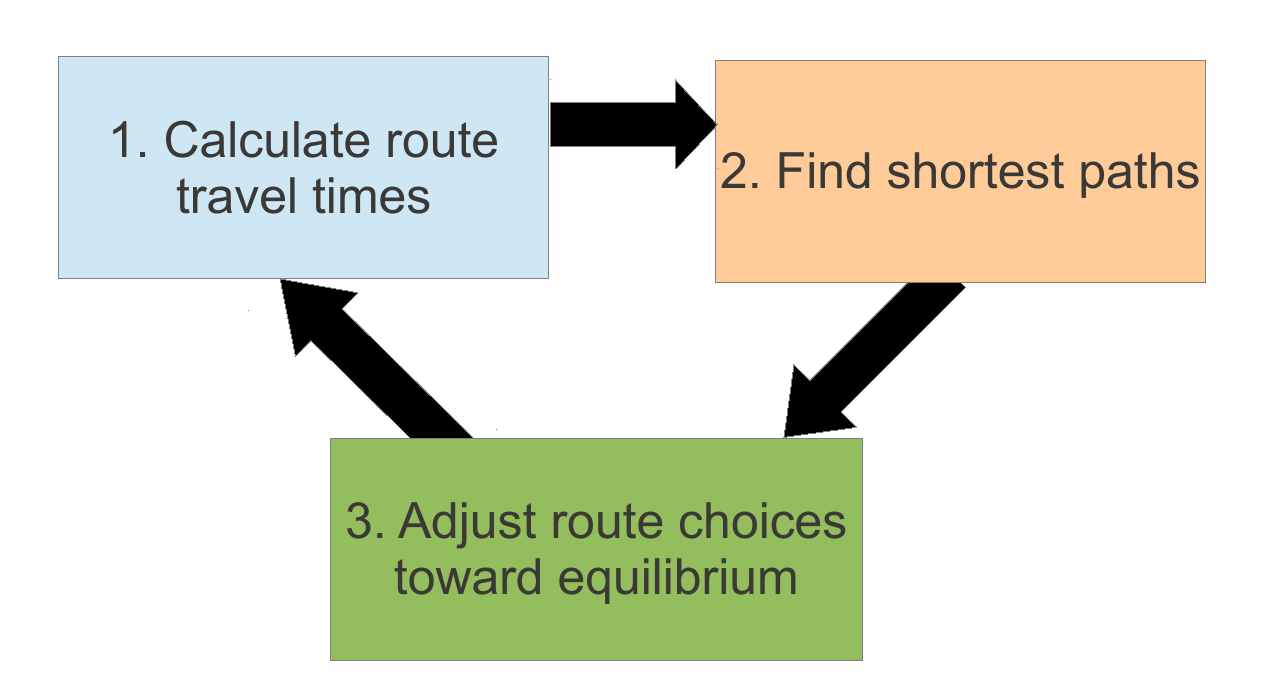}{Three-step iterative process for dynamic traffic assignment.}{0.7\textwidth}

This means that it is impossible to find the dynamic user equilibrium in one step.
Experienced travel times cannot be calculated at the moment of departure, but only after the vehicle has arrived at the destination.
Therefore, dynamic traffic assignment is an iterative process, where route choices are updated at each iteration until an (approximate) dynamic user equilibrium solution has been found.
This iterative process virtually always involves three steps, shown in Figure~\ref{fig:iterativeeqm}:
\begin{description}
\item[Network loading:\index{network loading} ] This is the process of using a traffic flow model to calculate the (time-dependent) travel times on each link, taking the routes and departure times of all drivers as inputs.
In static assignment, this step was quite simple, involving nothing more than evaluating the link performance functions for each network link.
In dynamic traffic assignment, this involves the use of a more sophisticated traffic flow model, or even the use of a traffic simulator.
Network loading is discussed in Chapter~\ref{chp:networkloading}.
\item[Path finding: ] Once network loading is complete, the travel time of each link, at each point in time, is obtained.
From this, we find the shortest path\index{shortest path} from each origin to each destination, \emph{at each departure time}.
Since we need experienced travel times, our shortest path finding must take into account that the travel time on each link depends upon the time a vehicle enters.
This requires \emph{time-dependent shortest path}\index{shortest path!time-dependent} algorithms, which are discussed in Chapter~\ref{chp:tdsp}.
\item[Route updating: ] Once time-dependent shortest paths have been found for all origins, destinations, and departure times, vehicles can be shifted from their current paths onto these new, shortest paths.
As in static assignment, this step requires care, because shifting vehicles will change the path travel times as well.
A few options for this are discussed in Chapter~\ref{chp:dynamiceqm}, along with other issues characterizing dynamic equilibrium.
Unfortunately, and in contrast with static assignment, dynamic user equilibrium need not always exist, and when it exists it need not be unique.
\end{description}

\subsection{Critique}

The primary advantage of dynamic traffic assignment, and a significant one, is that the underlying traffic flow models are much more realistic.
Link performance functions used in static assignment are fundamentally incapable of representing many important traffic phenomena.
However, dynamic traffic assignment is not without its drawbacks as well.
  
Dynamic assignment requires considerably more computational time and memory than static assignment.
As computers advance, this drawback is becoming less severe; but regardless of the computational resources available, one can run a large number of static traffic assignments in the time required for a single dynamic traffic assignment run.
It may be advantageous to examine a large number of scenarios with static assignment, rather than a single run with dynamic traffic assignment, particularly if there is a lot of uncertainty in the input data.

Dynamic assignment models also require more input data for calibration.
In addition to the usual link parameters such as capacity and free-flow time, dynamic traffic assignment models require a time-dependent origin-destination matrix, often known as a \emph{demand profile}.\index{demand profile}
Estimating even a static origin-destination matrix is difficult; estimating a dynamic demand profile can be even harder.

Furthermore, in addition to simply requiring more input data, dynamic traffic assignment also requires more accurate input data.\index{model!robustness}
Dynamic traffic assignment tends to be much more sensitive to the input data.
Unfortunately, these models are more sensitive precisely because they are more realistic --- features such as queue spillback mean that even a single erroneous capacity value can have ramifications throughout the entire network, not just on the link with the wrong value.
If there is great uncertainty in the inputs (for instance, when attempting to predict demand decades into the future), then using a dynamic traffic assignment model may actually be further away from the truth than a static model, despite its more ``realistic'' congestion model.

Separately, dynamic traffic assignment is a relatively young field relative to static assignment, and there is no consensus on a single formulation.
There are a large number of software packages (and an even larger number of academic models) which make differing assumptions and can lead to distinct results.
By contrast, the optimization, variational inequality, and fixed point formulations in static assignment are completely standard, and therefore these models are more transparent.
This book attempts to provide a high-level perspective, along with detailed discussions of a few of the most common modeling choices.

Finally, dynamic traffic assignment generally lacks neat, exact mathematical properties.
In many dynamic traffic assignment models, one can create examples where no dynamic user equilibrium exists, or where this equilibrium is not unique.
Proving convergence of iterative schemes is also more difficult.
  
All of these drawbacks must be traded off against the increased realism of dynamic traffic assignment models.
Both static and dynamic traffic assignment models have their place as distinct tools for transportation engineering, and you should learn to identify circumstances where one or the other is more appropriate.
As a general rule of thumb, dynamic models are most appropriate when the input data are known with high certainty (as in present-day traffic studies), and when queue lengths or other detailed congestion measures are required.
Static models, by contrast, perform best when there is considerable uncertainty in the input data, or when it is more important to run a large number of scenarios than to obtain detailed congestion information.
\index{dynamic traffic assignment|)}

\section{Historical Notes and Further Reading}

The first uses of networks were in solving mathematical puzzles like whether it is possible to walk through a city traversing each bridge exactly once \citep[the K\"{o}nigsberg bridge problem of][]{euler1736}\index{Euler, Leonhard}\index{K\"{o}nigsberg bridge problem}, or to visit each square on a chessboard exactly once using a knight~\citep{vandermonde1774}.\index{knight's tour}
The chemist Arthur Cayley\index{Cayley, Arthur} also studied network-like structures to count particular types of hydrocarbon compounds~\citep{cayley1857}.
These investigations led to the development of the mathematical field of graph theory.
There is no essential difference between what we call \emph{networks} in this book, and what mathematicians more commonly call \emph{graphs}.
The ``network'' terminology is more common in engineering and optimization, so that is what we will use here.
Our primary focus will be optimization problems defined over networks, and \cite{ahuja93} give a good treatment of applications, formulations, and algorithms for such problems.
If you are interested in the broader field of graph theory,\index{graph theory} \cite{diestel16} is a good reference. 
\index{graph|see {network}}

As alluded to in Section~\ref{sec:examplesofnetworks}, the network structure is very flexible and can represent more than just the ``obvious'' cases of physical transportation infrastructure with roads as links and intersections as nodes.
Networks have been used to represent systems as broad as ecosystems, social structures, communication systems, waterways, power systems, organizational structures, and so on.
Even within the transportation domain, phenomena such as elastic demand, destination choice, and mode choice can be accommodated within a network structure by adding links and nodes, and assigning them costs, in a special way~\citep[Chapter 9]{sheffi85}.

There is a long history of using equilibrium models in the field of economics.
We are using the word \emph{equilibrium} in the game theoretic sense, and the oligopoly pricing and production models of \cite{cournot1838}, \cite{bertrand1883}, and \cite{edgeworth1897} are the earliest ancestors of modern game theory.
The modern formulations of game theory were given by \cite{vonneumann44} and \cite{nash50}.
\cite{fudenbergtirole}, \cite{gibbons92}, \cite{rasmussen06}, and \cite{ritzberger02} provide overviews of this field and its historical development.
An early example of this type of analysis in the transportation field is \cite{pigou20}\index{Pigou, Arthur Cecil}, who gave an example showing that the equilibrium, ``user optimal'' solution need not be efficient, in terms of maximizing total utility across travelers.
(His example is described later in the book, in Section~\ref{sec:motivatingexamples}.)
Among traffic engineers, \cite{wardrop52}\index{Wardrop, John Glen} is credited with introducing two principles of route choice that correspond to the user optimal and system optimal states.

Traffic assignment is essentially solving an economic game for equilibrium on a network representing transportation infrastructure.
As discussed briefly in Section~\ref{sec:introductionassignment}, traffic assignment is part of the larger transportation planning process.
Transportation planning emerged in the aftermath of World War II,\index{World War II}\index{planning!history} due to the connected needs of suburbanization, population growth, and major investments in freeways and other infrastructure --- how should this new infrastructure be designed, how should alternatives be compared, and which projects deserve funding?
Comprehensive overviews of transportation planning can be found in \cite{domencich75} and \cite{meyer00}.

\index{planning!four-step model|(}
The ``four-step model'' is very commonly used in transportation planning, consisting of four steps in sequence: trip generation, trip distribution, mode choice, and traffic assignment (route choice).
Prior to executing these steps, a planner divides the region into zones in which trips originate and terminate.\index{trip generation|(}
Trip generation aims to determine the number of trips originating at zone $r$\label{not:r}, called the number of \emph{productions}\index{productions} $P_r$\label{not:Pr}, and the number of trips terminating at each zone $s$\label{not:s}, called the number of \emph{attractions}\index{attractions} $A_s$\label{not:As}.\footnote{The exact definition is slightly different for home-based and non-home-based trips, but this distinction is not important for our purposes because we only need the final demand matrix.  See Section 5.4 of \cite{meyer00} for more on this point.}
Productions and attractions are estimated based on land use, demographic data, and increasingly by observing cellular phone trajectories and other location data.
The \emph{Trip Generation Manual}, published by the Institute of Transportation Engineers~(\citeyear{itetripgeneration})\index{Institute of Transportation Engineers (ITE)}, is an example of a professional reference compiling trip generation data.
Other approaches to trip generation are described in \cite{kassoff69} and Chapter 4 of \cite{ortuzar11}. 
\index{trip generation|)}

\index{trip distribution|(}
Trip distribution converts zonal productions and attractions into interzonal trips $N_{rs}$\label{not:Nrs}, giving the number of trips which start at zone $r$ and end at zone $s$.
To match productions we must have $\sum_s N_{rs} = P_r$ for each origin $r$, and to match attractions we must have $\sum_r N_{rs} = A_s$ for each destination $s$, but these ``consistency'' constraints are not enough to specify the values $N_{rs}$.
The gravity model\index{gravity model} provides one way to find these values, with the formula
\labeleqn{gravitymodel}{
N_{rs} = \alpha_r P_r \beta_s A_s \phi_{rs}
\,,
}
where $\alpha_r$\label{not:alphar} and $\beta_s$\label{not:betas} are proportionality constants for each origin and destination, and $\phi_{rs}$\label{not:phirs} is a \emph{friction function}\index{friction function} indicating how onerous travel is between zones $r$ and $s$.
The constants $\alpha_r$ and $\beta_s$ are set to satisfy the consistency constraints.
Different forms of $\phi_{rs}$ are available, but the idea is that zones which are further away should have fewer trips between them, so $\phi_{rs}$ should shrink with distance, travel time, and/or accessibility.
With a little bit of imagination, one can view equation~\eqn{gravitymodel} as analogous to \citeapos{newton1687}\index{Newton, Isaac!gravitation} celebrated model of gravitation: factor out $\alpha_r \beta_s$ as the ``gravitational constant,'' interpret $P_r$ and $A_s$ as the masses of two objects, and take $\phi_{rs}$ to be the reciprocal of the squared distance between zones $r$ and $s$, and \emph{voil\`{a}}, we have Newton's formula.
This interpretation is a bit precious, and there is no real reason to believe that travelers behave in a way consistent with physical gravity.
However, by selecting a negative exponential form of the friction function $\phi_{rs} = \exp(-\gamma \kappa_{rs})$, where $\kappa_{rs}$\label{not:kappars} expresses the generalized cost of travel between zones $r$ and $s$, and $\gamma$\label{not:gamma} is a constant expressing the importance of this generalized cost in destination choice, we obtain a model with more solid foundations.
The gravity model with negative exponential friction can be interpreted using the language of utility maximization~\citep{mcfadden74td}, entropy maximization~\citep{wilson70}, and maximum likelihood~\citep{wilson81}, all of which are more plausible explanations of behavior than Newtonian gravity.
\index{trip distribution|)}

\index{mode choice|(}
Mode choice is commonly performed using a logit\index{logit} model from discrete choice,\index{discrete choice} following~\cite{mcfadden74}.
Assume that a traveler must make a choice among $n$ modes (such as driving, taking transit, or riding a bicycle).
Further assume that each traveler derives a utility associated from each mode choice.
Some of these factors are known to the planner (such as travel time and monetary cost by each mode, and demographic data such as income level and vehicle ownership), whereas others are not (such as personal tastes).
If we denote the former by $V_m$\label{not:Vm} for each mode $m$\label{not:m}, and the latter by a random variable $\epsilon$\label{not:epsilonutility} --- note here that we are assuming $\epsilon$ is independent and identically distributed for each mode --- then one can show that the probability that a given traveler selects mode $m$ is
\labeleqn{myfirstlogit}{
P_m = \left. \exp(V_m) \middle/ \sum_{m'=1}^n \exp(V_{m'}) \right.
}
if the random variables $\epsilon$ follow a Gumbel distribution.\footnote{The Gumbel distribution has a similar shape to the normal distribution, but with fatter tails.  The standard Gumbel distribution has cumulative distribution function $\exp(-\exp(-x))$.  The one used in equation~\eqn{myfirstlogit} has additional scaling and shifting parameters that are calibrated to data.}\index{Gumbel distribution}
These probabilities are used to split the total interzonal trips $N_{rs}$ calculated during trip distribution: the number of trips from zone $r$ to $s$ using mode $m$ is given by $N_{rs} P_m$, and it is these values that are used as the origin-destination matrix for traffic assignment in the remainder of the book.
The logit model is also discussed more in Section~\ref{sec:sue}.
\index{mode choice|)}

Although the four-step model has a long history of use, and remains the most common planning paradigm in practice, many researchers have suggested alternative techniques that can represent more sophisticated travel behavior.
For instance, considering travelers who make several stops as part of a trip (\emph{trip chaining}) in the four-step model is awkward.
Activity-based modeling is one alternative to the traditional practice; see \cite{bhat99} as a starting point for this growing field.
\index{planning!four-step model|)}
\index{trip chaining}
\index{activity-based model}

Static traffic assignment,\index{static traffic assignment} as described in this book, was first formulated by \cite{beckmann56}\index{Beckmann, Martin}, who also introduced mathematical representations which will be described further in Chapter~\ref{chp:trafficassignmentproblem}.
A great deal of subsequent research further developed and extended the basic static traffic assignment model.
The books by \cite{sheffi85}, \cite{patriksson94}, and \cite{bell97} discuss many of these further developments.

Dynamic traffic assignment\index{dynamic traffic assignment} is relatively newer.
The first dynamic traffic assignment models were developed by \cite{merchant78model, merchant78optimality}.
These models received increasing attention starting in the late 1990s, as more powerful computers became available, and \cite{peeta01} give a review of these advances.
Another perspective, tailored towards practicing engineers and planners, was provided in a primer~\citep{dtaprimer}.

\section{Exercises}
\label{exercises_transportationnetworks}

\begin{enumerate}
\item \diff{10} If all models are wrong, how can some of them be useful?
\item \diff{10} All of the links in Figure~\ref{fig:netexample} have a ``mirror'' connecting the same nodes, but in the opposite direction.
When might mirror links not exist?
\item \diff{23} A network is called \emph{planar}\index{network!planar} if it can be drawn in two dimensions without links crossing each other.
(The left network in Figure~\ref{fig:explode} is planar, but not the network on the right.)
Do we expect to see planar graphs in transportation network modeling?
Does it depend on the mode of transportation?
What other factors might influence whether a transportation network is planar?
\item \diff{35} Table~\ref{tbl:networktypes} shows how networks can represent five types of transportation infrastructure.
List at least five more systems (not necessarily in transportation) that can be modeled by networks, along with what nodes and links would represent.
\item \diff{45} What factors might determine how large of a geographic area is modeled in a transportation network (e.g., corridor, city, region, state, national)? 
Illustrate your answer by showing how they would apply to the hypothetical projects or policies at the start of Section~\ref{sec:transportationnetworks}.
\item \diff{45} What factors might determine the level of detail in a transportation network (e.g., freeways, major arterials, minor arterials, neighborhood streets)? 
Illustrate your answer by showing how they would apply to the hypothetical projects or policies at the start of Section~\ref{sec:transportationnetworks}.
\item \diff{21} Provide an intuitive explanation of the prisoner's dilemma, as described in the Erica-Fred game of Table~\ref{tbl:ericafred}.
Why does it happen?
Name other real-world situations which exhibit a similar phenomenon.
\item \diff{20} For each of the following games, list all equilibria (or state that none exist).
In which of these games do equilibria exist; in which are the equilibria unique; and in which are there some equilibria which are inefficient?  These games are all played by two players A and B: A chooses the row and B chooses the column, and each cell lists the payoffs to A and B, in that order.
  \begin{center}
  \setlength{\extrarowheight}{2pt}
  \begin{tabular}{c|c|c|}
     \multicolumn{1}{c}{(a)} & \multicolumn{1}{c}{$L$}  & \multicolumn{1}{c}{$R$} \\\cline{2-3}
     $U$ & $(8,13)$ & $(5,14)$ \\\cline{2-3}
     $D$ & $(10,10)$ & $(7,12)$ \\\cline{2-3}
  \end{tabular}
  \begin{tabular}{c|c|c|}
     \multicolumn{1}{c}{(b)} & \multicolumn{1}{c}{$L$}  & \multicolumn{1}{c}{$R$} \\\cline{2-3}
     $U$ & $(12,12)$ & $(2,10)$ \\\cline{2-3}
     $D$ & $(10,2)$ & $(4,5)$ \\\cline{2-3}
  \end{tabular}
  
  \begin{tabular}{c|c|c|}
     \multicolumn{1}{c}{(c)} & \multicolumn{1}{c}{$L$}  & \multicolumn{1}{c}{$R$} \\\cline{2-3}
     $U$ & $(8,6)$ & $(10,8)$ \\\cline{2-3}
     $D$ & $(2,3)$ & $(8,4)$ \\\cline{2-3}
  \end{tabular}
  \begin{tabular}{c|c|c|}
     \multicolumn{1}{c}{(d)} & \multicolumn{1}{c}{$L$}  & \multicolumn{1}{c}{$R$} \\\cline{2-3}
     $U$ & $(4,5)$ & $(5,6)$ \\\cline{2-3}
     $D$ & $(3,4)$ & $(6,3)$ \\\cline{2-3}
  \end{tabular}
  \end{center}
\item \diff{1} Explain the difference between the demand for travel on a link, and the flow on a link.
\item \diff{42} Explain why a less realistic model less sensitive to correct input data may be preferred to a more realistic model more sensitive to correct inputs, and in what circumstances.
Give specific examples.
\end{enumerate}

\chapter{Network Representations and Algorithms}
\label{chp:networkrepresentations}

This chapter introduces networks as they are used in the transportation field.
Section~\ref{sec:terminology} introduces the mathematical terminology and notation used to describe network elements.
Section~\ref{sec:treesacyclic} discusses two special kinds of networks which are used frequently in network analysis, acyclic networks and trees.
Section~\ref{sec:datastructures} then describes several ways of representing a network in a way computers can use when solving network problems.
This third section can be skipped if you do not plan to do any computer programming.
Section~\ref{sec:shortestpath} describes the shortest path problem, a classic network algorithm which plays a central role in traffic assignment.

\section{Terminology}
\label{sec:terminology}

\index{network|(}
Because of its relative youth, there are a variety of notational conventions in the transportation networks community.
A common notation is adopted in the book to present the methods and topics in a consistent manner, but you should be aware that other authors may use slightly different conventions and definitions.
Table~\ref{tbl:alternateterminology} gives an example of some terms which are often used synonymously (or nearly synonymously) with ours.
These differences are mainly a matter of style, not substance, but when reading other books or papers you should pay careful attention to the exact wording of their definitions.
\index{vertex|see {node}}
\index{arc|see {link}}
\index{edge|see {link}}
\index{arborescence|see {tree}}

\begin{table}
\begin{center}
\caption{A thesaurus of network terminology. \label{tbl:alternateterminology}}
\begin{tabular}{l|l}
Terminology in this book & Alternative terms \\
\hline
Network & Graph \\
Node    & Vertex \\
Link    & Arc, edge, line \\
Tree	  & Arborescence
\end{tabular}
\end{center}
\end{table}

The fundamental construct we will utilize is the \emph{network}.
A network is a collection of \emph{nodes}\index{node}, and a collection of \emph{links}\index{link} which connect the nodes.
A network is denoted $G = (N,A)$\label{not:G}\label{not:N}\label{not:A}, where $N$ is the set of nodes and $A$ is the set of links.
Figure~\ref{fig:graphnotation}(a) shows a simple network, with four nodes in the set $N = \{1, 2, 3, 4\}$ and five links in the set $A =\{(1,2), (1,3), (2,3), (2,4), (3,4)\}$.
Notice that the notation for each link contains the two nodes connected by the link: the upstream node is called the \emph{tail}\index{link!tail node} of the link, and the downstream node the \emph{head}.\index{link!head node}
We will often refer to the total number of nodes in a network as $n$\label{not:n} and the total number of links as $m$\label{not:mlinks}.
This book is solely concerned with the case where $n$ and $m$ are finite --- networks with infinitely many nodes and links are sometimes studied in theoretical mathematics, but rarely in transportation applications.
\index{node!tail|see {link, tail node}}
\index{node!head|see {link, head node}}
\index{tail|see {link, tail node}}
\index{head|see {link, head node}}

\begin{figure} 
\centering
   \begin{tikzpicture}[->,>=stealth',shorten >=1pt,auto,node distance=3cm,
   thick,main node/.style={circle,draw,font=\sffamily\Large\bfseries}]

      \node[main node] (1) {$1$};
      \node[main node] (2) [above right of=1] {$2$};
      \node[main node] (3) [below right of=1] {$3$};
      \node[main node] (4) [above right of=3] {$4$};

   \path[every node/.style={font=\sffamily\small}]
         (1) edge node[above left] {$(1,2)$} (2)
         edge node [below left] {$(1,3)$} (3)

         (2) edge node [right] {$(2,3)$} (3)
         edge node [above right] {$(2,4)$} (4)

         (3) edge node [below right] {$(3,4)$} (4);

   \end{tikzpicture}
\caption{Example network notation. \label{fig:graphnotation}}
\end{figure}

All of the networks in this book are \emph{directed}.\index{network!directed}
In a directed network, a link can only be traversed in one direction, specified by the ordering of the nodes.
Therefore, (1,2) and (2,1) represent different links: they connect the same nodes, but represent travel in opposite directions.
Unless stated otherwise, we assume that there are no ``self-loops'' $(i,i)$ from a node to itself, and no parallel links, so the notation $(i,j)$ is unambiguous.
This is not usually a restrictive assumption, since we can introduce artificial nodes to split up self-loops or parallel links, as shown in Figure~\ref{fig:networktransformations}.
The upper left panel in this figure shows a network with a self-loop (2,2).
By introducing a fourth node in the bottom left, we have divided the self-loop into two links (2,4) and (4,2).
In the upper right panel of the figure, there are two networks connecting nodes 2 and 3, so the notation (2,3) does not tell us which of these two links we are referring to.
By introducing a fourth node in the bottom right, we now have three links: (2,3), (2,4), and (4,3), which collectively represent both of the ways to travel between nodes 2 and 3 in the original network, but without parallel links.
These new nodes are \emph{artificial}\index{node!artificial} in the sense that they do not represent physical transportation infrastructure, but are inserted for modeling convenience.
Artificial nodes (and artificial links\index{link!artificial}) play an important role in simplifying certain network problems, as discussed throughout the book.

\stevefig{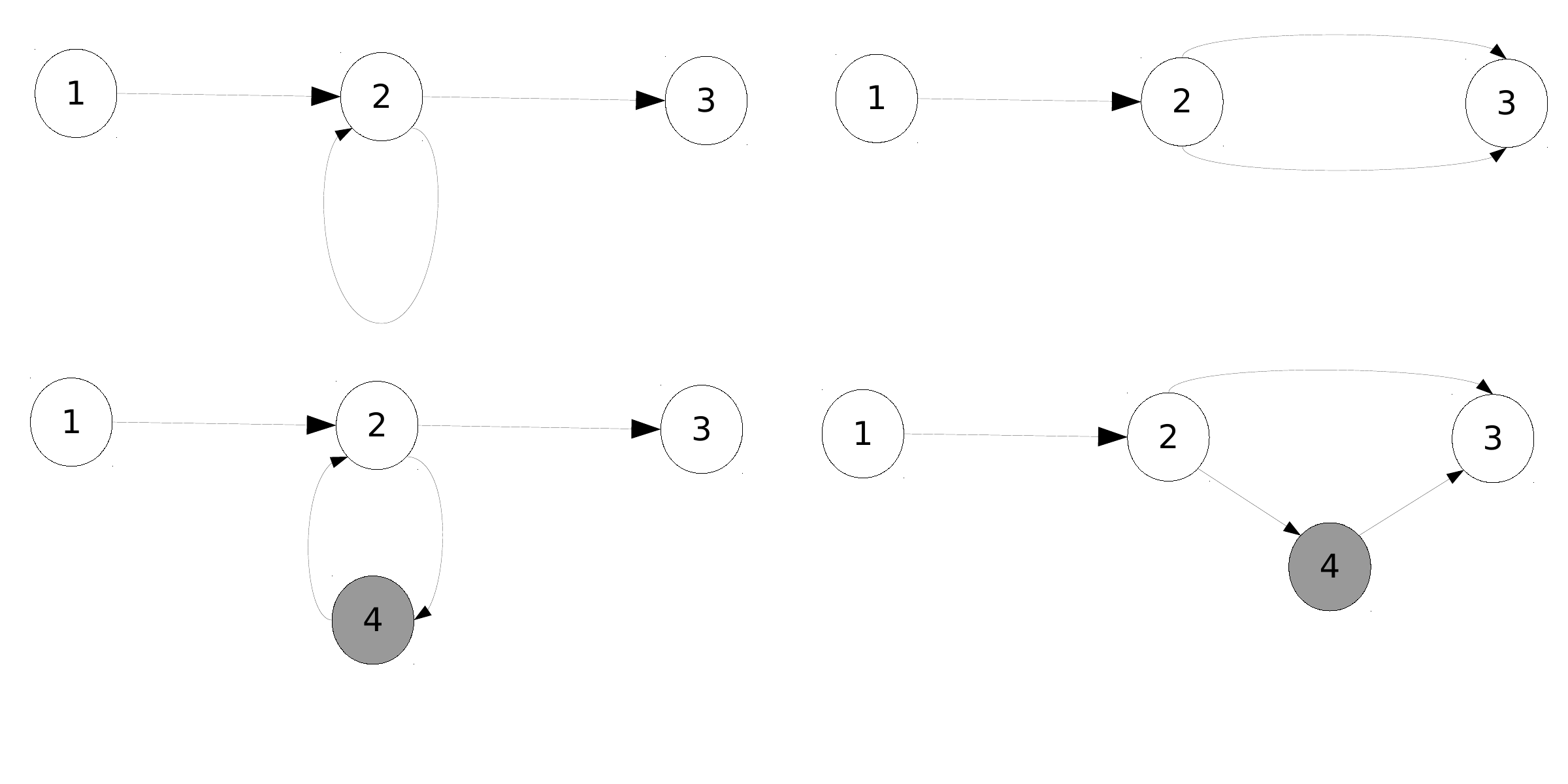}{Artificial links can be introduced to simulate self-loops and parallel links.}{0.8\textwidth} 

To introduce some more terminology, Figure~\ref{fig:genericgraph} shows a two-node network with nodes $i$ and $j$ connected by the link $(i, j)$.
In this figure, we say that link $(i, j)$ is \emph{incident}\index{link!incident} to both $i$ and $j$ because it is connected to both.
We would further say that $(i, j)$ is an \emph{outgoing}\index{link!outgoing} link of $i$ and an \emph{incoming}\index{link!incoming} link of $j$.
If $(i,j) \in A$, then we say that node $j$ is \emph{adjacent}\index{node!adjacent} to node $i$.
In this example, $j$ is adjacent to $i$, but $i$ is not adjacent to $j$.
Adjacency can also be applied to links; two links are adjacent\index{link!adjacent} if the head of the first link is the tail of the second link.
In Figure~\ref{fig:graphnotation}, both links (2,3) and (2,4) are adjacent to (1,2).

The \emph{forward star}\index{forward star} of a node $i$ is the set of all outgoing links, denoted $\Gamma(i)$\label{not:Gammai}; the \emph{reverse star}\index{reverse star}\label{not:Gamma-1i} is the set of all incoming links, denoted $\Gamma^{-1}(i)$.
In Figure~\ref{fig:graphnotation}, $\Gamma(2) = \{ (2,3), (2,4) \}$ and $\Gamma(3) = \{ (3,4) \}$, while $\Gamma^{-1}(2) = \{ (1,2) \}$ and $\Gamma^{-1}(3) = \{ (1,3), (2,3) \}$.

\begin{figure}
\centering
	\begin{tikzpicture}[->, >=stealth',shorten >=1pt,auto,node distance=3cm,
  	thick,main node/.style={circle,draw,font=\sffamily\Large\bfseries}]

		\node[main node] (i) {$i$};
 		\node[main node] (j) [right of=i] {$j$};

  	\path[every node/.style={font=\sffamily\small}]
    		(i) edge node[above] {$(i, j)$} (j);

	\end{tikzpicture}
\caption{Generic two-node network. \label{fig:genericgraph}}

\end{figure}

The \emph{degree}\index{node!degree} of a node is the total number of links incident to that node.
The node degree can be separated into the \emph{indegree}\index{node!indegree} and \emph{outdegree}\index{node!outdegree} of a node.
The indegree is the total number of incoming links at a node, while the outdegree is the number of outgoing links at a node.
The degree is then the sum of the indegree and outdegree.
Referring back to Figure~\ref{fig:graphnotation}, node 2 has an indegree of 1, an outdegree of 2 and a degree of 3.

A \emph{path}\index{path} $\pi$\label{not:pi} is a sequence of adjacent links connecting two nodes $i_0$ and $i_k$.
We can either write $\pi$ as an ordered set of links
\[\{ (i_0, i_1), (i_1, i_2), (i_2, i_3), \ldots, (i_{k-1}, i_k)\}\,,\]
or more compactly, by the nodes passed on the way with the notation
\[[i_0, i_1, i_2, i_3, \ldots, i_{k-1}, i_k]\,.\]
A path is a \emph{cycle}\index{cycle} if the starting and ending nodes are the same ($i_0 = i_k$).
Paths which contain a cycle are called \emph{cyclic}\index{path!cyclic}; paths which do not have a cycle as a component are called \emph{acyclic}\index{path!acyclic}.
Cyclic paths are uncommon in transportation problems, so unless stated otherwise, we only consider acyclic paths.
Let $\Pi^{rs}$\label{not:Pirs} denote the set of all acyclic paths connecting nodes $r$ and $s$, and let $\Pi$\label{not:Pi} denote the set of all acyclic paths in the entire network, that is, $\Pi = \cup_{(r,s) \in N^2} \Pi^{rs}$\label{not:union}.
A network is \emph{connected}\index{network!connected} if there is at least one path connecting any two nodes in the network, assuming that the links can be traversed in either direction (that is, ignoring the direction of the link); otherwise it is \emph{disconnected}.\index{network!disconnected}
A network is \emph{strongly connected}\index{network!strongly connected} if there is at least one path connecting any two nodes in the network, obeying the given directions of the links.
Figure~\ref{fig:connectedstrongly} shows networks which are strongly connected; connected, but not strongly connected; and disconnected.

\stevefig{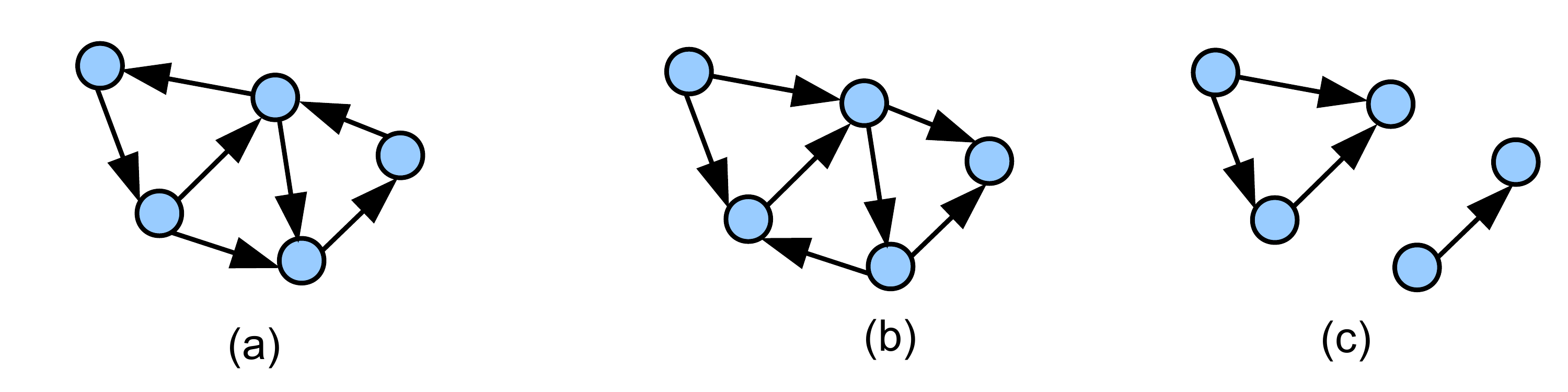}{Networks which are (a) strongly connected; (b) connected, but not strongly connected; (c) disconnected.}{\textwidth}

\section{Acyclic Networks and Trees}
\label{sec:treesacyclic}

Networks which do not have any cycles are called \emph{acyclic networks}.\index{network!acyclic}
Acyclic networks are extremely important, because the lack of cycles can greatly simplify analysis.
Even when networks actually have cycles in reality, many efficient transportation network algorithms temporarily divide the network into a set of acyclic subnetworks.

A defining characteristic of acyclic networks is that a \emph{topological order}\index{topological order} can be established.
A topological order is a labeling of each node with a number between 1 and $n$ ($n$ is the number of nodes), such that every link connects a node of lower topological order to a node of higher topological order.
Every path, therefore, consists of a sequence of nodes in increasing topological order.
As an example, in Figure~\ref{fig:graphnotation}, the nodes are labeled in a topological order: each path ($[1,3,4]$, $[1,2,4]$, and $[1,2,3,4]$) traverses nodes in increasing numerical order.
(The phrase ``a'' topological order is used because it may not be unique, and there may be several ways to label the nodes so that links always connect lower-numbered nodes to higher-numbered ones.)  In general networks, a topological order may not exist (see Figure~\ref{fig:topologicaldne}, where a cycle makes it impossible to label nodes so that the numbers increase when traversing any link).
The following theorem shows that the existence of a topological order is a \emph{defining} characteristic of acyclic networks.

\stevefig{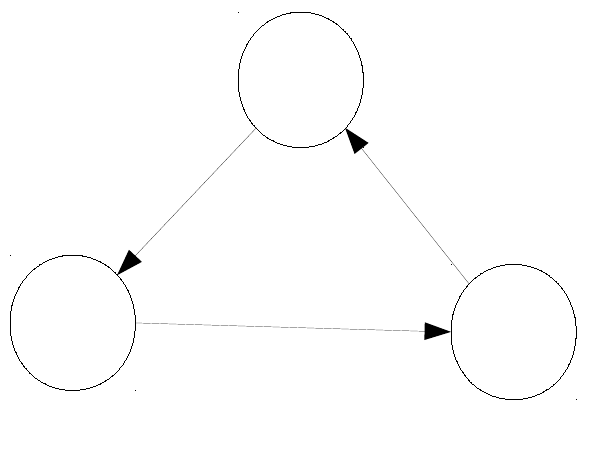}{Cyclic network with no topological order.}{0.3\textwidth}

\begin{thm}
	\label{thm:acyclictopological}
	A topological order exists on a network if and only if it is acyclic.
\end{thm}
\begin{proof}
	($\Rightarrow$)  Assume a topological order exists on $G$.
Then every path contains a sequence of nodes whose topological order is strictly increasing (if this were not so, then there is a link connecting a higher topological order node to a lower one, contradicting the definition of topological order).
Therefore, no path can repeat the same node more than once, and the network is acyclic.
	
	($\Leftarrow$)  Assume the network $G$ is acyclic.
Then we can prove existence of a topological order by construction.
Let $o(i)$\label{not:oi} be the topological order of node $i$.
We will describe a procedure to build such a topological order one step at a time.
At any point in time, a ``marked'' node is one which has a topological order assigned, and an ``unmarked'' node is one which does not yet have its topological order.
Because $G$ is acyclic, there is at least one node $r$ with no incoming links.
(You will be asked to prove this statement as an exercise.)  Let $o(r) = 1$.
Again, because $G$ is acyclic, there is at least one unmarked node with no incoming links from an unmarked node; mark this node with the topological order 2.
This process can be repeated until all nodes are marked; if at any point every unmarked node has an incoming link from another unmarked node, then a cycle exists, contradicting the assumption of acyclicity.
Otherwise, the topological order will have been constructed.
The topological order so constructed is valid because when a node is assigned an order, the only incoming links are from nodes which have already been marked, and thus have a lower topological order.
\end{proof}        

A \emph{tree}\index{tree} is a special type of acyclic network, which also shows up frequently in transportation network algorithms.
A tree is defined as a network in which there is a unique node $r$ (called the \emph{root}) which has the property that exactly one path exists between $r$ and every other node in the network.\footnote{In graph theory, this is usually called an \emph{arborescence}, and trees are defined for undirected graphs.
However, the transportation community generally uses the term ``tree'' in both cases, a convention followed in this book.}  Figure~\ref{fig:tree} shows a tree.
 Other useful properties of trees are:
\begin{itemize}
	\item A tree has at least two nodes with degree one.
	\item In a tree, there is at most one path between any two nodes.
	\item A tree would become disconnected if any of its links were deleted.
	\item Every node $i$ in a tree (except the root) has a unique incoming link; the tail node of that unique link is called the \emph{parent}\index{node!parent} of node $i$.
Similarly, the head nodes of the links in the forward star of $i$ are its \emph{children}.\index{node!child}
\end{itemize}

\begin{figure}[!ht]
\centering
	\begin{tikzpicture}[->,>=stealth',shorten >=1pt,auto,node distance=2cm,
  	thick,main node/.style={circle,draw,font=\sffamily\scriptsize\bfseries}]

  		\node[main node] (1) {$1$};
 		\node[main node] (2) [right of=1] {};
 		\node[main node] (3) [right of=2] {};
  		\node[main node] (4) [above right of=3] {};
		\node[main node] (5) [below right of=3] {};
		\node[main node] (6) [above right of=5] {};
		\node[main node] (7) [below right of=5] {};

  	\path[every node/.style={font=\tiny\small}]
    		(1) edge node[above left] {} (2)
			
    		(2) edge node [below right] {} (3)
			
		(3) edge node [above right] {} (4)
			edge node [above right] {} (5)

    		(5) edge node [below right] {} (6)
			edge node[below] {} (7);

	\end{tikzpicture}
\caption{A tree. \label{fig:tree}}
\end{figure}

\section{Data Structures}
\label{sec:datastructures}

In practice, the methods and techniques in this book are carried out by computer programs.
Thus, it is important not only to have a convenient way to theoretically represent a network (the network structure) but a way to store data so that it can be easily accessed and utilized by computer programs.
To demonstrate, Figure~\ref{fig:adjacencyexample} presents a network similar to those earlier, with additional information provided.
In this case, the extra information is the travel time of the link from $i$ to $j$, denoted $t_{ij}$.

\begin{figure}
\centering
\subfigure{
\begin{tikzpicture}[->, >=stealth',shorten >=1pt,auto,node distance=3cm,
  	thick,main node/.style={circle,draw,font=\sffamily\Large\bfseries}]

		\node[main node] (i) {$i$};
 		\node[main node] (j) [right of=i] {$j$};

  	\path[every node/.style={font=\sffamily\small}]
    		(i) edge node[above] {$t_{ij}$} (j);

	\end{tikzpicture}
}

\subfigure{
	\begin{tikzpicture}[->,>=stealth',shorten >=1pt,auto,node distance=3cm,
  	thick,main node/.style={circle,draw,font=\sffamily\Large\bfseries}]

  		\node[main node] (1) {$1$};
 		\node[main node] (2) [above right of=1] {$2$};
 		\node[main node] (3) [below right of=1] {$3$};
  		\node[main node] (4) [above right of=3] {$4$};

  	\path[every node/.style={font=\sffamily\small}]
    		(1) edge node[above left] {$2$} (2)
			edge node [below left] {$1$} (3)

    		(2) edge node [right] {$3$} (3)
			edge node [above right] {$4$} (4)

    		(3) edge node [below right] {$5$} (4);

	\end{tikzpicture}
}
\caption{Graph used for network representation examples. \label{fig:adjacencyexample}}
\end{figure}

The first data structure is the node-link incidence matrix,\index{node-link incidence matrix} shown in Figure~\ref{fig:nodelink}.
The columns are indexed by the network links, and the rows by the network nodes.
A ``$1$'' in the matrix indicates an outgoing link from that node, while a ``$-1$'' indicates an incoming link.
A ``$0$'' indicates that the link is not incident to that particular node.
This is a fairly inefficient data structure, especially for large networks, as in nearly any case there are going to be a large number of zeros in every row.
Note that the forward and reverse stars of node $i$ can respectively be identified with the columns in the matrix which have a $1$ or $-1$ in the $i$-th row.

\begin{figure}
\centering

$\bbordermatrix{~ & (1,2) & (1,3) & (2,3) & (2,4) & (3,4) \cr
		    1 & 1 & 1 & 0 & 0 & 0 \cr
		    2 & -1 & 0 & 1 & 1 & 0 \cr
		    3 & 0 & -1 & -1 & 0 & 1 \cr
		    4 & 0 & 0 & 0 & -1 & -1 \cr}$

\caption{Node-link incidence matrix for the network in Figure~\ref{fig:adjacencyexample}. \label{fig:nodelink}}
\end{figure}

The second data structure, a node-node adjacency matrix\index{node-node incidence matrix} is a simpler representation of a network but also suffers from inefficiencies.
Each row represents a node $i$ and each column a node $j$.
A ``$1$'' indicates the existence of a link $(i, j)$.
In this structure, a non-zero entry provides both existence and direction information.
An example of a node-node adjacency matrix is given in Figure~\ref{fig:nodenode}.
In very dense networks, with many more links than nodes ($m \gg n$), node-node adjacency matrices can be efficient.

\begin{figure}
\centering

$\bbordermatrix{~ & 1 & 2 & 3 & 4 \cr
		    1 & 0 & 1 & 1 & 0 \cr
		    2 & 0 & 0 & 1 & 1 \cr
		    3 & 0 & 0 & 0 & 1 \cr
		    4 & 0 & 0 & 0 & 0 \cr}$

\caption{Node-node adjacency matrix for the network in Figure~\ref{fig:adjacencyexample}. \label{fig:nodenode}}
\end{figure}

A shortcoming of the previous data structures is that they do not contain information beyond the existence and direction of the links in a network.
In transportation applications we are typically working with networks for which we want to store and use more information about each element.
For example, we may want to store information about the travel time, number of lanes (capacity), speed limit, signalization, etc. for each link.
You can maintain multiple matrices, one storing the travel times, another storing the capacity, and so forth, but this is inefficient and wastes storage space: each matrix has to store a `0' entry for any link which does not exist, and in transportation networks the number of links is much smaller than the number of links which could possibly exist (each node is only connected to three or four other nodes, on average).

Adjacency lists\index{adjacency list} give us an opportunity to store this information efficiently.
Many modern programming languages have built-in list or ``dictionary'' hash table data structures.
In such languages, you can create a custom data object or data structure to represent a single link.
You can then create lists or dictionaries that reference all of the links entering or leaving each node in the network, and use such lists for finding paths through a network or other computations.

A more basic way to store these lists is with an array; while the details of list and dictionary structures vary from language to language, virtually every programming language has an array data structure.
Storage is not wasted on a large number of zeros, as with the adjacency matrices.
A ``forward star'' representation\index{forward star representation|(} is presented below.
In the forward star representation, links are sorted in order of their tail node, as shown in Table~\ref{tbl:fsexample}.
(Links with the same tail node do not have to be further sorted in order of head node.)
A second array, called \ttt{point}, is used to indicate the position in this sorted list where the forward star of each node begins.
As shown in Table~\ref{tbl:fsexample}, \ttt{point(2)} = 3 because the link 3 is the first link adjacent to node 2; \ttt{point(3)} = 5 because link 5 is the first link adjacent to node 3, and so forth.

We adopt three conventions in the forward star representation to handle ``edge cases'' outside of the normal pattern:
\begin{enumerate}
	\item If there are no outgoing links, \ttt{point(i) = point (i+1)}
	\item The number of entries in the \ttt{point} array is one more than the number of nodes.
	\item We always set \ttt{point(n+1) = m + 1}.
\end{enumerate}

With these conventions, we can say that the forward star for any node $i$ consists of all links whose IDs are greater than or equal to \ttt{point(i)}, and strictly less than \ttt{point(i+1)}.
This representation is most useful when we frequently need to loop over all of the links leaving a particular node, which can be accomplished by programming a ``for'' loop between \ttt{point(i)} and \ttt{point(i+1) - 1}, and referring to the link arrays with these indices.
In the Python language, this is especially convenient: the links leaving node $i$ are exactly those given by the \ttt{range(point[i], point[i+1])} expression.

We can make this statement universally, because we defined the \ttt{point} array to have one more entry than the number of nodes.
If \ttt{point} only had $n$ entries, then referring to \ttt{point(n+1)} would be meaningless in the above statements.
Our conventions allow us to treat all nodes the same, without special tests to see if there are any links in the forward star, or if this node is the ``last'' node.
Another example of the forward star representation is seen in Figure~\ref{tbl:fsexample2}.
\index{forward star representation|)}

It is also possible to store a network in a ``reverse star'' representation along similar lines, which is most useful when we frequently need to loop over all the links entering a particular node.
The exercises explore this further, along with an array-based method for easily iterating over both the forward and reverse stars of a node.

\begin{table}
	\begin{center}
	\caption{Forward star representation of the network in Figure~\ref{fig:adjacencyexample}. \label{tbl:fsexample}}
	\begin{tabular}{ r c | c | c r c | c | c | c | c |}
		\multicolumn{1}{c}{Node} & 
		\multicolumn{1}{c}{} & 
		\multicolumn{1}{c}{Point} & 
		\multicolumn{1}{c}{} & 
		\multicolumn{1}{c}{Link} & 
		\multicolumn{1}{c}{} & 
		\multicolumn{1}{c}{tail} & 
		\multicolumn{1}{c}{head} & 
		\multicolumn{1}{c}{cost} & 
		\multicolumn{1}{c}{capacity} \\
		\cline{1-1}
		\cline{3-3}
		\cline{5-5}
		\cline{7-10}
		1 & & 1 & & 1 & & 1 & 3 & $\cdots$ & $\cdots$ \\ 
		\cline{3-3}
		\cline{7-10}
		2 & & 3 & & 2 & & 1 & 2 & $\cdots$ & $\cdots$ \\ 
		\cline{3-3}
		\cline{7-10}
		3 & & 5 & & 3 & & 2 & 3 & $\cdots$ & $\cdots$ \\ 
		\cline{3-3}
		\cline{7-10}
		4 & & 6 & & 4 & & 2 & 4 & $\cdots$ & $\cdots$ \\ 
		\cline{3-3}
		\cline{7-10}
		5 & & 6 & & 5 & & 3 & 4 & $\cdots$ & $\cdots$ \\ 
		\cline{3-3}
		\cline{7-10}
	\end{tabular}
	\end{center}
\end{table}

\begin{table}
	\begin{center}
	\caption{Forward star representation of a network with more nodes.\label{tbl:fsexample2}}
    \includegraphics[width=0.7\textwidth]{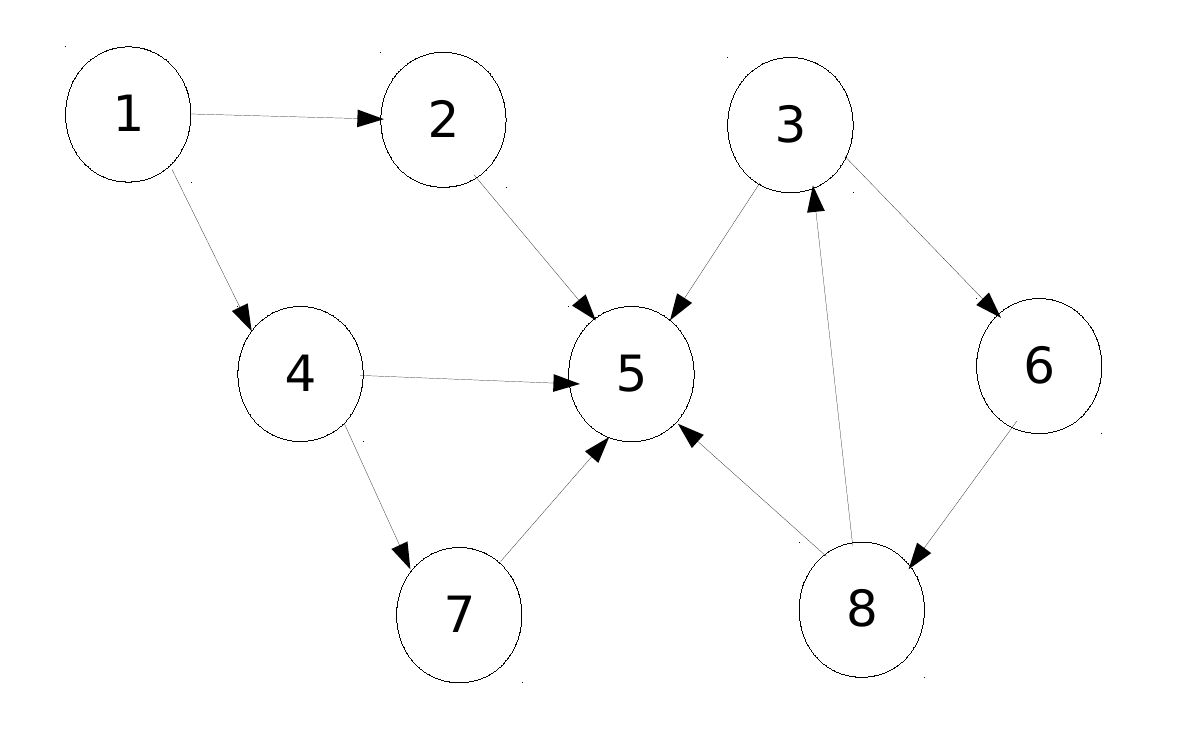}
	\begin{tabular}{ r c | c | c r c | c | c |}
		\multicolumn{1}{c}{Node} & 
		\multicolumn{1}{c}{} & 
		\multicolumn{1}{c}{Point} & 
		\multicolumn{1}{c}{} & 
		\multicolumn{1}{c}{Link} & 
		\multicolumn{1}{c}{} & 
		\multicolumn{1}{c}{tail} & 
		\multicolumn{1}{c}{head} \\ 
		\cline{1-1}
		\cline{3-3}
		\cline{5-5}
		\cline{7-8}
		1 & & 1 & & 1 & & 1 & 2 \\
		\cline{3-3}
		\cline{7-8}
		2 & & 3 & & 2 & & 1 & 4 \\
		\cline{3-3}
		\cline{7-8}
		3 & & 4 & & 3 & & 2 & 5 \\
		\cline{3-3}
		\cline{7-8}
		4 & & 6 & & 4 & & 3 & 5 \\
		\cline{3-3}
		\cline{7-8}
		5 & & 8 & & 5 & & 3 & 6 \\
		\cline{3-3}
		\cline{7-8}
		6 & & 8 & & 6 & & 4 & 5 \\
		\cline{3-3}
		\cline{7-8}
		7 & & 9 & & 7 & & 4 & 7 \\
		\cline{3-3}
		\cline{7-8}
		8 & & 10 & & 8 & & 6 & 8 \\
		\cline{3-3}
		\cline{7-8}
		9 & & 12 & & 9 & & 7 & 5 \\
		\cline{3-3}
		\cline{7-8}
		  &\multicolumn{1}{c}{}   & \multicolumn{1}{c}{}  & & 10 & & 8 & 3 \\
		\cline{7-8}
		  &\multicolumn{1}{c}{}   & \multicolumn{1}{c}{} & & 11 & & 8 & 5 \\
		\cline{7-8}
	\end{tabular}
	\end{center}
\end{table}

The purpose of this brief discussion of data structures is not to provide a comprehensive treatment of the topic.
Rather, it is intended to present some basic ideas related to network storage and data representation.
It should also prompt you to think about data structures carefully when working with networks.
Transportation network problems tend to be large in size and complicated in nature.
The difference between computer code that takes minutes to run, as opposed to hours, is often the way data is stored and passed.
\index{network|)}

\section{Shortest Paths}
\label{sec:shortestpath}

\index{shortest path|(}
As a first network algorithm, we'll discuss how to find the shortest (least cost) path between any two nodes in a network in an efficient and easily-automated way.
 In shortest path algorithms, we are given a network $G = (N, A)$; each link $(i,j) \in A$ has a \emph{fixed} cost $c_{ij}$\label{not:cij}.
The word ``cost''\index{link!cost} does not necessarily mean a monetary cost, and refers to whatever total quantity we are trying to minimize.
In traffic applications, we often use the link travel times $t_{ij}$ as the cost, to find the least travel-time path.
If there are tolls in a network, the cost of each link may be the sum of the toll on that link, and the travel time on that link multiplied by a ``value of time'' factor converting time into monetary units.
Costs may reflect still other factors, and some of these are explored in the exercises.
Because we can solve all of these problems in the same way, we may as well use a single name for all of these quantities we are trying to minimize, and ``cost'' has become the standard term.

The word ``fixed'' in the last paragraph is emphasized because in many transportation problems the cost is dependent on the total flow and vehicle route choices, and it is not reasonable to assume fixed costs.
However, even in such network problems the most practical solution methods involve solving several shortest path problems as part of an iterative framework, updating paths as travel times change.
The context for a shortest path algorithm is to find the least travel-time path between two nodes, at the current travel times --- with everybody's route choices held fixed.
By separating the process of identifying the shortest path from the process of shifting path flows toward the shortest path, we simplify the problem and end up with something which is relatively easy to solve.
This problem can also be phrased in the language of optimization, as shown in Appendix~\ref{sec:optimizationexamples}, by identifying an objective function, decision variables, and constraints.
Here, we develop specialized algorithms for the shortest path problem which are more efficient and more direct than general optimization techniques.

Although the shortest path problem most obviously fits into transportation networks, many other applications also exist in construction management, geometric design, operations research, and many other areas.
For instance, the fastest way to solve a Rubik's Cube\index{Rubik's Cube} from a given position can be solved using a shortest path algorithm, as can the fewest number of links needed to connect an actor to Kevin Bacon when playing Six Degrees of Separation.\index{Six Degrees of Separation}

A curious fact of shortest path algorithms is that finding the shortest path from a single origin to \emph{every other} node is only slightly harder than finding the shortest path from that origin to a single destination.
This also happens to be the reason why we can find shortest paths without having to list all of the zillions of possible paths from an origin to a destination, and adding up their paths.
This common reason is \emph{Bellman's principle}\index{Bellman's principle|(}, which states that \emph{any segment of a shortest path must itself be a shortest path between its endpoints}.
For instance, consider the network in Figure~\ref{fig:braesssp}, where the costs $c_{ij}$ are printed next to each link.
The shortest path from node 1 to node 4 is $[1,2,3,4]$.
Bellman's principle requires that $[1,2,3]$ also be a shortest path from nodes 1 to 3, and that $[2,3,4]$ be a shortest path from nodes 2 to 4.
It further requires that $[1,2]$ be a shortest path from node 1 to node 2, $[2,3]$ be a shortest path from node 2 to node 3, and $[3,4]$ be a shortest path from node 3 to node 4.
You should verify that this is true with the given link costs.

\stevefig{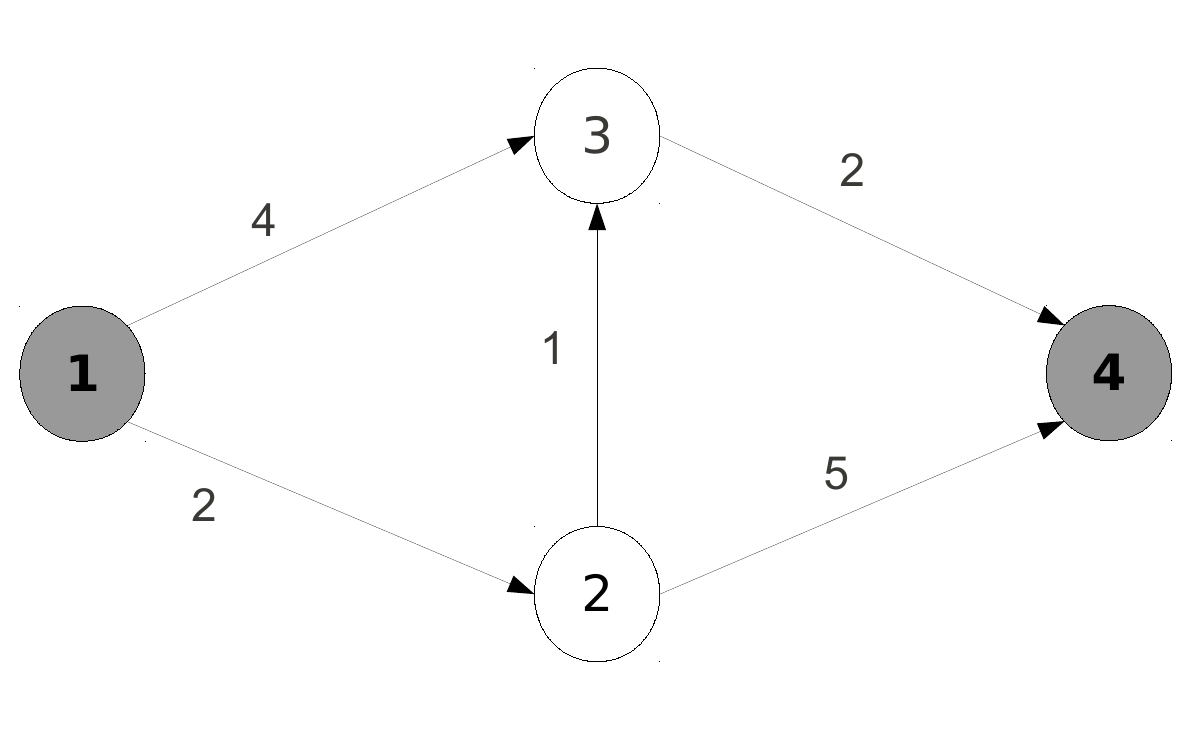}{Example network for shortest paths (link costs indicated).}{0.6\textwidth}

To see why this must be the case, assume that Bellman's principle was violated.
If the cost on link $(1,3)$ was reduced to 2, then $[1,2,3]$ is no longer a shortest path from node 1 to node 3 (that path has a cost of 3, while the single-link path $[1,3]$ has cost 2).
Bellman's principle then implies that $[1,2,3,4]$ is no longer the shortest path between nodes 1 and 4.
Why?  The first part of the path can be replaced by $[1,3]$ (the new shortest path between 1 and 3), reducing the cost of the path from 1 to 4: $[1,3,4]$ now has a cost of 4.
In general, if a segment of a path does \emph{not} form the shortest path between two nodes, we can replace it with the shortest path, and thus reduce the cost of the entire path.
Thus, the shortest path must satisfy Bellman's principle for all of its segments.
\index{Bellman's principle|)}

The implication of this is that \emph{we can construct shortest paths one node at a time}, proceeding inductively.
Let's say we want to find the shortest path from node $r$ to a node $i$, and furthermore let's assume that we've already found the shortest paths from $r$ to every node which is directly upstream of $i$ (nodes $f$, $g$, and $h$ in Figure~\ref{fig:recursion}).
The shortest path from $r$ to $i$ must pass through either $f$, $g$, or $h$; and according to Bellman's principle, the shortest path from $r$ to $i$ must be either (a) the shortest path from $r$ to $f$, plus link $(f,i)$; the shortest path from $r$ to $g$, plus link $(g,i)$; or the shortest path from $r$ to $h$, plus link $(h,i)$.
This is efficient because, rather than considering all of the possible paths from $r$ to $i$, we only have to consider three, which can be easily compared.
Furthermore, we can re-use the information we found when finding shortest paths to $f$, $g$, and $h$, and don't have to duplicate the same work when finding the shortest path to $i$.
This idea doesn't give a complete algorithm yet --- how did we find the shortest paths to $f$, $g$, and $h$, for instance? --- but gives the flavor of the shortest path algorithms presented next.

\stevefig{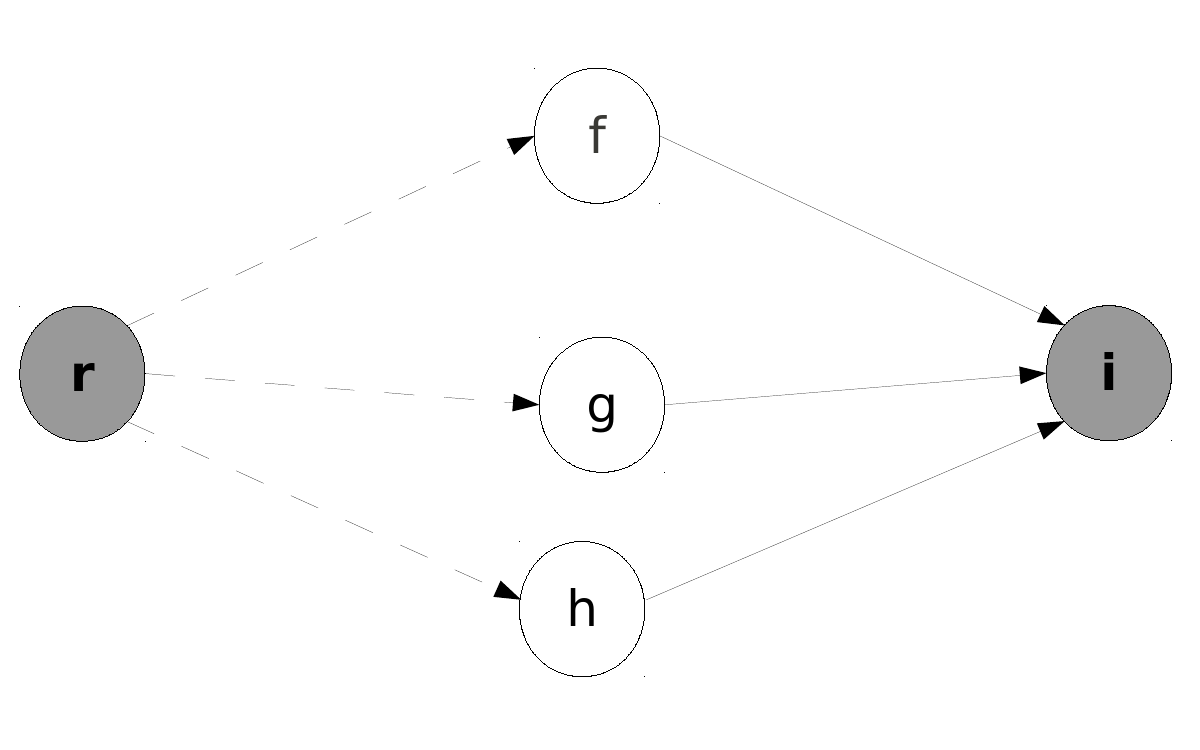}{Illustration of Bellman's principle (dashed lines indicated shortest paths between two nodes).}{0.6\textwidth}

Bellman's principle also gives us a compact way of expressing all of the shortest paths from an origin to every other node in the network: for each node, simply indicate the \emph{last} node in the shortest path from that origin.
This is called the backnode\index{backnode|see {node, backnode}}\index{node!backnode} vector $\mathbf{q}^r$\label{not:qr}, where each component $q^r_i\label{not:qir}$ is the node immediately preceding $i$ in the shortest path from $r$ to $i$.
If $i = r$, then $q^r_i$ is not well-defined ($i$ is the origin itself; what is the shortest path from the origin to itself, and if we can define it, what node immediately precedes the origin?) so we say $q^r_r \equiv -1$ by definition.
For the network in Figure~\ref{fig:braesssp} (with the original costs), we thus have $q^1_1 = -1$, $q^1_2 = 1$, $q^1_3 = 2$, and $q^1_4 = 3$, or, in vector notation, $\mathbf{q}^1 = \vect{-1 & 1 & 2 & 3}$. 

The backnode vector can be used as follows: say we want to look up the shortest path from node 1 to node 4.
Starting at the destination, the backnode of 4 is 3, which means ``the shortest path to node 4 is the shortest path to node 3, plus the link $(3,4)$.''  To find the shortest path to node 3, consult its backnode: ``the shortest path to node 3 is the shortest path to node 2, plus the link $(2,3)$.''  For the shortest path to node 2, its backnode says: ``the shortest path to node 2 is from node 1.''  This is the origin, so we've found the start of the path, and can reconstruct the original path to node 4: $[1,2,3,4]$.
More briefly, we can use the backnodes to trace the shortest path back to an origin, by starting from the destination, and reading back one node at a time.

We will also define $L^r_i$\label{not:Lir} to be the total cost on the shortest path from origin $r$ to node $i$ (the letter $L$ is used because these values are often referred to as node labels), with $L^r_r \equiv 0$, so in this example we would have $L^1_1 = 0$, $L^1_2 = 2$, $L^1_3 = 3$, and $L^1_4 = 5$.\index{node!cost}\index{cost label|see {node, cost}}

This chapter presents four shortest path algorithms.
The first, in Section~\ref{sec:acyclicsp}, only applies when the network is acyclic but is extremely fast and simple.
Section~\ref{sec:spgeneral} then presents a general method that can apply in any network, with or without cycles.
There are several ways to implement this general method; one of the most important is Dijkstra's algorithm.
All of these methods actually find the shortest path from the origin $r$ to \emph{every} other node, not just the destination $s$; this exploits Bellman's principle, because the shortest path from $r$ to $s$ must also contain the shortest path from $r$ to every node in that path.
In transportation network analysis, we usually have to find the shortest paths from all origins to all destinations, so this is a significant advantage of these methods.
However, there are occasions when we only need to find the shortest path from one origin to one destination.
In this case, there is a faster algorithm $A^*$ which is presented in Section~\ref{sec:astar}.

\subsection{Shortest paths on acyclic networks}
\label{sec:acyclicsp}

\index{network!acyclic|(}
Recall from Section~\ref{sec:treesacyclic} that an \emph{acyclic} network is one in which no cyclic paths exist; it is impossible to visit any node more than once on any path.
Conversely, a \emph{cyclic}\index{network!cyclic} network is one where a cycle does exist, and where (in theory) one could repeatedly drive around in a circle forever.
Transportation networks are usually cyclic, and in fact any network where you can reach every node from every other node must be cyclic.
(Why?)
However, we'll take a short diversion into acyclic networks for the purposes of finding shortest paths, for two reasons: (1) finding the shortest path on an acyclic network is much simpler and faster, and makes for a good first illustration; and (2) more advanced solution methods take advantage of the fact that people do \emph{not} use cyclic paths (because of the shortest path assumption), so we can only look at an acyclic portion of the network and thereby use the much faster shortest path algorithm for acyclic networks.
      
Once we have a topological order on a network, it becomes very easy to find the shortest path from any node $r$ to any other node $s$ (clearly $r$ has lower topological order than $s$), using the following algorithm which is based directly on Bellman's principle.\index{Bellman's principle}
The idea is that a shortest path to any node must consist of a shortest path to a node which is immediately upstream, extended by one link connecting that upstream node to the current node.
At any point in time, the $L$ and $q$ labels contain the cost and backnode corresponding to the shortest paths known so far.
The algorithm starts at node $r$, and \emph{scans}\index{scanning a node} subsequent nodes in topological order.
At each iteration we have already found the shortest path to whichever node is being scanned, and we consider all possible ways to extend this path by one link to downstream nodes. 
Whenever this produces a better path than the one previously known, we update the labels.
Once we reach node $s$, we are done, and can trace the backnode labels back from $s$ to $r$.

\begin{enumerate}
\item Initialize by setting $L^r_r \leftarrow 0$,\label{not:leftarrow} because the distance from the origin to itself is zero, $L^r_i \leftarrow \infty$ for $i \neq r$, and $q^r_i \leftarrow -1$ for all nodes $i \in N$ to indicate that we have not found any shortest paths yet.\footnote{We are intentional about using $\leftarrow$ to mean assignment, rather than $=$; see Appendix~\ref{sec:algorithms} for more on this distinction.}
Then set $L^r_r \leftarrow 0$, because the distance from the origin to itself is zero.
Finally set $i \leftarrow r$.
\item For every link $(i,j)$ whose tail is node $i$, compute the cost of extending the shortest path to node $i$ as $L^r_i + c_{ij}$.  
If $L^r_i + c_{ij} < L^r_j$, then this is a better path to $j$ than anything found so far, so update $L^r_j \leftarrow L^r_i + c_{ij}$ and $q^r_j \leftarrow i$.
\item Set $i$ to be the next node topologically.
If $i = s$, then stop:  we have found the shortest path from $r$ to $s$, which has cost $L^r_s$.
Otherwise, return to step 2.
\end{enumerate}
We can find the shortest paths in one pass over the nodes (in topological order), because there are no cycles which could make us loop back.

We demonstrate this algorithm on the network in Figure~\ref{fig:braesssp}, with node 1 as origin and node 4 as destination.
\begin{description}
\item[Step 1: ] Initialize: $\mathbf{L}^1 \leftarrow \vect{0 & \infty & \infty & \infty}$, $\mathbf{q}^1 \leftarrow \vect{-1 & -1 & -1 & -1}$, and $i \leftarrow 1$.
\item[Step 2: ] Scan node 1.
There are two links leaving this node: (1,2) and (1,3).
For (1,2), $L^1_1 + c_{12} = 0 + 2 = 2$ and $L^1_2 = \infty$, so we update the labels at node 2: $L^1_2 \leftarrow 2$ and $q^1_2 \leftarrow 1$.  
Likewise for (1,3), $L^1_1 + c_{13} = 0 + 4 = 4 < L^1_3 = \infty$, so update $L^1_3 \leftarrow 4$ and $q^1_3 \leftarrow 1$.
\item[Step 3: ] Move to the next node by setting $i \leftarrow 2$.
\item[Step 2: ] Scan node 2.
There are two links leaving this node: (2,3) and (2,4).
For (2,3), $L^1_2 + c_{23} = 2 + 1 = 3$, which is less than the current value of $L^1_3$ (4).
Therefore the path through node 2 has lower cost, and we update $L^1_3 \leftarrow 3$ and $q^1_2 \leftarrow 2$.
For (2,4), we update $L^1_4 \leftarrow 7$ and $q^1_4 \leftarrow 2$.
\item[Step 3: ] Move to the next node by setting $i \leftarrow 3$.
\item[Step 2: ] Scan node 3.
There is only one link leaving this node: (3,4).
Since $L^1_3 + c_{34} = 3 + 2 = 5$, which is less than the current value of $L^1_4$, so we update $L^1_4 \leftarrow 5$ and $q^1_4 \leftarrow 3$.
\item[Step 3: ] Move to the next node by setting $i \leftarrow 4$.
This is the destination, so stop.
\end{description}

\index{mathematical induction}
A simple induction proof shows that this algorithm must give the correct shortest paths to each node (except for those topologically before $r$, since such paths do not exist).
By ``give the correct shortest paths,'' we mean that the labels $\mb{L}$ give the lowest cost possible to each node when traveling from $r$, and that the backnodes $\mb{q}$ yield shortest paths when traced back to $r$.
Assume that the nodes are numbered in topological order.
Clearly $L^r_r$ and $q^r_r$ are set to the correct values when scanning $r$ (the shortest path from $r$ to itself is trivial), and are never changed again because the network proceeds in increasing topological order.
Now assume that the algorithm has just finished scanning node $k$, and is about to move to node $k + 1$.
We claim that the $L$ and $q$ values for node $k + 1$ must be correct. 
Let $(i, k+ 1)$ be the last link in a shortest path from $r$ to $k + 1$.
By Bellman's principle,\index{Bellman's principle} the first part of this path must be a shortest path from $r$ to $i$.
Since $i$ is topologically between $r$ and $k$, by the induction hypothesis the $L^r_i$ and $q^r_i$ labels at this node were correct when node $i$ was scanned, so $L^r_i + c_{i,k+1}$ is indeed the cost of the shortest path from $r$ to $k + 1$.
The same argument establishes that $L^r_h + c_{h,k+1}$, where $h$ is \emph{any} node immediately upstream of $k + 1$, is the cost of an actual path from $r$ to $k + 1$, and therefore at least as large as $L^r_i + c_{i,k+1}$.
Therefore, after the optimal backnode $i$ was scanned $L^r_{k+1}$ must have been equal to $L^r_i + c_{i,k+1}$, and it can never have been reduced further, and the backnode label must also have been set correctly.
\index{network!acyclic|)}

\subsection{Shortest paths on general networks}
\label{sec:spgeneral}

When there are cycles in the network, this previous approach can't be applied, because there is no clear sequence in which to examine nodes and apply Bellman's principle.
However, we can generalize the approach.
Rather than scanning nodes in a rigid order, we can fan out from the origin, keeping in mind that we may have to scan a node multiple times in case of a cycle.
To keep track of the nodes we need to examine, we define a scan eligible list $SEL$\label{not:SEL}, a set of nodes that we still need to examine before we have found all of the shortest paths.
This is a \emph{label correcting}\index{shortest path!label correcting} approach, because nodes may be scanned multiple times, and the labels updated.

\begin{enumerate}
\item Initialize the labels as $L^r_r \leftarrow 0$, $L^r_i \leftarrow \infty$ for $i \neq r$, and $\mathbf{q}^r \leftarrow \bm{-1}$.
Also initialize the scan eligible list\index{scan eligible list} to contain the origin only: $SEL \leftarrow \{ r \}$.
\item Choose a node $i \in SEL$ and delete it from that list.
\item Scan node $i$.
For every link $(i,j)$ whose tail is node $i$, compute the cost of extending the shortest path to node $i$ as $L^r_i + c_{ij}$.  
If $L^r_i + c_{ij} < L^r_j$, then this is a better path to $j$ than anything found so far, so update $L^r_j \leftarrow L^r_i + c_{ij}$, $q^r_j \leftarrow i$, and add $j$ to $SEL$.
\item If $SEL$ is empty, then terminate.
Otherwise, return to step 2.
\end{enumerate}

Repeating the same example, this algorithm works as follows:

\begin{description}
\item[Step 1: ]  We set $\mathbf{L}^1 \leftarrow \vect{0&\infty&\infty&\infty}$, $\mathbf{q}^1 \leftarrow \vect{-1&-1&-1&-1}$, and $SEL \leftarrow \myc{1}$.
\item[Step 2: ] The only choice for $i$ is node 1; at this point $SEL$ is empty.
\item[Step 3: ] Scan node 1 by considering the two links which leave this node: (1,2) and (1,3).
For (1,2), we compute $L^1_1 + c_{12} = 0 + 2 = 2$; this is lower than the current value of $L^1_2$ ($\infty$) so we update $L^1_2 \leftarrow 2$, $q^1_2 \leftarrow 1$, and add node 2 to $SEL$.
Similarly for (1,3), we update $L^1_3 \leftarrow 4$, $q^1_3 \leftarrow 1$, and add node 3 to $SEL$.
\item[Step 4: ] At this point $SEL = \myc{2, 3}$, so we return to step 2.
\item[Step 2: ] We can choose either node in $SEL$ for $i$; assume that we set $i \leftarrow 3$, so $SEL = \myc{2}$.
\item[Step 3: ] Scan node 3.
The only link leaving node 3 is (3,4); $L^1_3 + c_{34} = 4 + 2 = 6$, which is less than $L^1_4$ ($\infty)$, so update $L^1_4 \leftarrow 6$, $q^1_4 \leftarrow 3$, and add node 4 to $SEL$.
\item[Step 4: ] At this point $SEL = \myc{2, 4}$, so we return to step 2.
\item[Step 2: ] We can again choose either node in $SEL$; assume that we set $i \leftarrow 4$, so $SEL = \myc{2}$.
\item[Step 3: ] There are no links leaving node 4, so there is nothing to do in this step.
\item[Step 4: ] At this point $SEL = \myc{2}$, so we return to step 2.
\item[Step 2: ] There is only one choice for $i$; so we set $i \leftarrow 2$.
The scan list is now empty.
\item[Step 3: ] There are two links leaving node 2: (2,3) and (2,4).
For (2,3), $L^1_2 + c_{23} = 2 + 1 = 3$, which is better than the current value of $L^1_3$ (4).
So we update $L^1_3 \leftarrow 3$, $q^1_3 \leftarrow 2$, and add 3 to $SEL$.
For (2,4), $L^1_2 + c_{24} = 7$, which is worse than the current value of $L^1_4$ (6)
Therefore no changes are made to the labels at node 4.
\item[Step 4: ] At this point $SEL = \myc{3}$, so we return to step 3.
\item[Step 2: ] There is only one choice for $i$; so we set $i \leftarrow 3$, emptying the scan list.
\item[Step 3: ] Scan node 3; this process updates $L^1_4 \leftarrow 5$ and $q^1_4 \leftarrow 3$, and node 4 is added again to $SEL$.
\item[Step 4: ] At this point $SEL = \myc{4}$, so we return to step 3.
\item[Step 2: ] The only choice is $i \leftarrow 4$.
\item[Step 3: ] Scan node 4; there is nothing to do since there are no links in the forward star.
\item[Step 4: ] The scan list is empty, so terminate.
\end{description}

As you can see, this method required more steps than the algorithm for acyclic networks (because there is a possibility of revisiting nodes), but it does not rely on having a topological order and can work in any network.
One can show that the algorithm converges no matter how you choose the node from $SEL$, but it is easier to prove if you choose a systematic rule.
Here are a few of these rules that generally work well in practice; in these rules, whenever there is a tie you can break it arbitrarily:
\begin{itemize}
\item Choose the node which has been in $SEL$ \textbf{for the most number of iterations}, breaking ties arbitrarily.
This is commonly called the \emph{FIFO label-correcting method}, where the acronym FIFO indicates that the scan eligible list is managed as a ``first-in, first-out'' queue.
This method is simple to program and is reasonably efficient.  
Furthermore, it requires at most $mn$ iterations to compute all of the shortest paths, where $m$ and $n$ are the number of links and nodes in the network, respectively.
One way to prove this is to show that after $m$ iterations, you have certainly found the shortest paths from $r$ which consist of only a single link.
(You've probably found quite a few more shortest paths, but even in the worst case you'll have found at least these.)  After $2m$ iterations, you will have certainly found the shortest paths from $r$ which consist of one or two links only, and so on.
So, after $mn$ iterations, you will have found all of the shortest paths, since a shortest path cannot use more than $n$ links.
The exercises ask you to fill in the details of this proof sketch.

\item Follow the FIFO rule the first time you add each node to $SEL$.
However, \textbf{if a node is added to $SEL$ after it has already been scanned, scan that node next.}
In other words, nodes entering the scan list a second time (or a third time, etc.) are bumped to the front of the queue.
The intuition is that if you are re-scanning a node, it is likely to affect paths you have already found that pass through that node.
So it is better to update all of those costs immediately before proceeding with the regular order of nodes.
This is often called \emph{Pap\'{e}'s rule}\index{Pap\'{e}'s rule}, or the \emph{deque label-correcting method}\index{deque}, where ``deque'' means ``double-ended queue.''
Unlike a FIFO queue, nodes can enter from both the front and the back, although they always leave from the front.
In practice this works very well, and is usually faster than the FIFO rule.
However, there are pathological instances where this method is much, much slower.

\item Choose the node $i$ in $SEL$ \textbf{with the smallest $L^r_i$ label}.
This rule, which is an implementation of \emph{Dijkstra's algorithm}\index{Dijkstra's algorithm}, has several interesting properties.
If the link costs are nonnegative, it is possible to show that nodes are scanned once and only once.
For this reason this method is sometimes called a \emph{label setting}\index{shortest path!label setting} method, in contrast to the more general family of label correcting methods were nodes may need to be scanned several times.
This means that only $n$ iterations are required, the same as with acyclic networks!
(If some links have negative costs, nodes may need to be scanned multiple times.)
However, there is no such thing as a free lunch: it requires more effort at each iteration to identify the node in $SEL$ with minimal $L^r_i$ value.
A na\"{i}ve way to do this is to examine each node in the list one by one; this is sure to find the minimum, but if there are many nodes in $SEL$ this can take quite a bit of effort.
Much research has gone into identifying customized data structures that make this process faster.
Because Dijkstra's algorithm (and other label setting methods) only have to scan each node once, it may be possible to terminate even before $SEL$ is empty.
If you are only interested in finding shortest paths from the origin to only selected nodes as destinations (e.g., all of the zones in the network), you can stop as soon as all of the nodes of interest have been scanned.
This is \emph{not} true for the other rules for choosing $i$, because there is no guarantee that you've found the shortest path the first time you scan a node.
\end{itemize}

\subsection{Shortest paths from one origin to one destination}
\label{sec:astar}

\index{A$^*$|(}
The previous sections gave algorithms to find the shortest path from one origin to all destinations.
Slight modifications to these algorithms can find the shortest path from all origins to one destination (see Exercise~\ref{ex:reversesp}).
As discussed above, Bellman's principle lets us re-use information from one shortest path when finding another, and as a result even finding a shortest path from one origin to one destination can provide many other shortest paths ``for free.''  In traffic assignment, where we are modeling the flows of many drivers, it often makes sense to use an algorithm that finds many shortest paths simultaneously.

However, there are times when we are only concerned with a single origin and destination, and do not care about the ``free'' shortest paths to other destinations we get.
Examples are in-vehicle routing systems, or if the number of origin-destination pairs with positive demand is small compared to the possible number of these pairs (i.e., the demand matrix\index{OD matrix} is sparse), so there are relatively few paths we need to find.
In such cases we can use a more focused algorithm to find a single shortest path in less time than it takes to find shortest paths to all destinations.

The $A^*$ algorithm is a simple modification to Dijkstra's algorithm,\index{Dijkstra's algorithm} which in some cases can dramatically reduce the running time, in exchange for limiting the scope to one origin $r$ and one destination $s$.
This algorithm requires an additional value $g_i^s$\label{not:gis} for each node, representing an estimate of the cost of the shortest path from $i$ to $s$.
This estimate (often called the ``heuristic'') should be a \emph{lower bound} on the actual cost of the shortest path from $i$ to $s$.
Some examples of how these estimates are chosen are discussed below.

Once the estimates $g_i^s$ are chosen, the labeling algorithm from Section~\ref{sec:spgeneral} proceeds as before with a small modification to Dijkstra's rule for choosing nodes: rather than choosing a node $i$ from $SEL$ which minimizes $L_i^r$, we choose a node which minimizes $L_i^r + g_i^s$.
Everything else proceeds exactly as before.
As we noted in the previous section, when using Dijkstra's algorithm we can stop as soon as we scan all of the destination nodes.
Since $A^*$ is one-to-one, you can terminate when you are about to scan the destination (when $s$ is chosen as the node $i$ in Step 2.)

The intuition behind choosing a node with minimum $L_i^r + g_i^s$ is that we want to preferentially scan nodes which we think are closer to the destination.
Dijkstra's algorithm fans out in all directions from the origin (by simply looking at $L_i^r$), rather than directing the search towards a particular destination.
Note the $r$ superscript on $L$, and the $s$ superscript on $g$.
The $L$ values represent distances on paths \emph{from the origin}; the $g$ values represent estimates of the distance \emph{to the destination}.
This is why the $A^*$ method only applies to shortest path problems with a single origin and destination.\footnote{Of course, you can always repeat the algorithm with different origins and destinations on the same network.  You must then weigh the increased number of times you need to run the algorithm (for each OD pair, rather than just once per origin) against the faster run times for the single origin-destination problems.}

As an example, consider the network in Figure~\ref{fig:grid-astar}, where the origin is node 5, and the destination $s$ is node 3.
Before starting the algorithm, we need to choose the $g_i^s$ values in such a way that the shortest path distance from each node $i$ to the destination $s$ is \emph{at least} $g_i^s$.
In this network, we notice that each link has a cost of at least 2.
Therefore, the shortest path from each node to the destination has a cost at least twice the number of \emph{links} in the most direct path (the one with the least number of links).
So $g_1^3 = 4$, $g_2^3 = 2$, and so forth.
From the origin, $g_5^3 = 4$, because any path from the origin to the destination uses at least two links, each with a cost of 2.
From the opposite corner, $g_7^3 = 8$, because any path to the destination uses at least four links.

Table~\ref{tbl:astarclever} shows how the $A^*$ algorithm progresses.
Each row of the table represents an iteration, showing the node $i$ selected from $SEL$, the $L_i^r$ values for each node at the \emph{end} of that iteration, and the nodes in $SEL$, also at the end of the iteration.
The header of the table shows the $g_i^s$ values for each node.
For brevity, the $q$ labels are not shown.
As noted in the discussion of Dijkstra's algorithm at the end of the previous subsection, we can stop as soon as we scan the destination (node 3) because this is a label setting algorithm, where each node is only scanned once.
In this case, six iterations were needed before we scan the destination and stop.
Different choices of $g_i^s$ will lead to different performance, as demonstrated below.

\stevefig{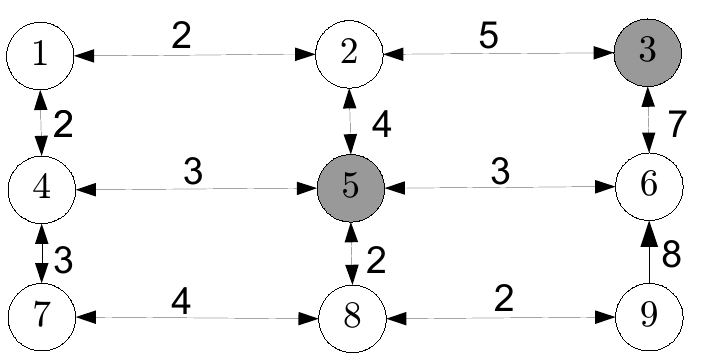}{Network for demonstrating the $A^*$ algorithm, with link costs shown.}{0.8\textwidth}

\begin{table}
\begin{center}
\caption{$A^*$ with lower bounds based on number of links to destination.
         \label{tbl:astarclever}}
\begin{tabular}{cc|ccccccccc|c}
&	$g_i^3$&	4&	2&	0&	6&	4&	2&	8&	6&	4&	\\
Iteration&	$i$&	$L_1$&	$L_2$&	$L_3$&	$L_4$&	$L_5$&	$L_6$&	$L_7$&	$L_8$&	$L_9$&	$SEL$\\
\hline
0&	---&	$\infty$&	$\infty$&	$\infty$&	$\infty$&	0&	$\infty$&	$\infty$&	$\infty$&	$\infty$&	$\{ 5 \}$\\
1&	5&	$\infty$&	4&	$\infty$&	3&	0&	3&	$\infty$&	2&	$\infty$&	$\{ 2, 4, 6, 8 \}$\\
2&	6&	$\infty$&	4&	10&	3&	0&	3&	$\infty$&	2&	11&	$\{ 2, 3, 4, 8, 9 \}$\\
3&	2&	6&	4&	9&	3&	0&	3&	$\infty$&	2&	11&	$\{ 1, 3, 4, 8, 9 \}$\\
4&	8&	6&	4&	9&	3&	0&	3&	6&	2&	4&	$\{ 1, 3, 4, 7, 9 \}$\\
5&	9&	6&	4&	9&	3&	0&	3&	6&	2&	4&	$\{ 1, 3, 4, 7 \}$\\
6&	3&	6&	4&	9&	3&	0&	3&	6&	2&	4&	$\{ 1, 4, 7 \}$\\
\end{tabular}
\end{center}
\end{table}

It can be shown that this algorithm will always yield the correct shortest path from $r$ to $s$ as long as the $g_i^s$ are lower bounds on actual shortest path costs from $i$ to $s$.
If this is not the case, $A^*$ is not guaranteed to find the shortest path from $r$ to $s$.
Some care must be taken in how these estimates are found.
Two extreme examples are:
\begin{itemize}
\item Choose $g_i^s = 0$ for all $i$.
This is certainly a valid lower bound on the shortest path costs (recall that label-setting methods assume nonnegative link costs), so $A^*$ will find the shortest path from $r$ to $s$. However, zero is a very poor estimate of the actual shortest path costs.
With this choice of $g_i^s$, $A^*$ will run exactly the same as Dijkstra's algorithm, and there is no time savings.
In the network of Figure~\ref{fig:grid-astar}, we would get the result shown in Table~\ref{tbl:astardumb}.
Nine iterations are needed to find the shortest path in this case.
\begin{table}
\begin{center}
\caption{Na\"{i}ve $A^*$ with the trivial lower bound (identical to Dijkstra's).
         \label{tbl:astardumb}}
\begin{tabular}{cc|ccccccccc|c}
         & $g_i^3$ & 0 & 0 & 0 & 0 & 0 & 0 & 0 & 0 & 0 & \\
Iteration&	$i$&	$L_1$&	$L_2$&	$L_3$&	$L_4$&	$L_5$&	$L_6$&	$L_7$&	$L_8$&	$L_9$&	$SEL$\\
\hline
0&	---&	$\infty$&	$\infty$&	$\infty$&	$\infty$&	0&	$\infty$&	$\infty$&	$\infty$&	$\infty$&	$\{ 5 \}$\\
1&	5&	$\infty$&	4&	$\infty$&	3&	0&	3&	$\infty$&	2&	$\infty$&	$\{ 2, 4, 6, 8 \}$\\
2&	4&	5&	4&	$\infty$&	3&	0&	3&	6&	2&	$\infty$&	$\{ 1, 2, 6, 7, 8 \}$\\
3&	8&	5&	4&	$\infty$&	3&	0&	3&	6&	2&	4&	$\{ 1, 2, 6, 7, 9 \}$\\
4&	6&	5&	4&	10&	3&	0&	3&	6&	2&	4&	$\{ 1, 2, 3, 7, 9 \}$\\
5&	2&	5&	4&	10&	3&	0&	3&	6&	2&	4&	$\{ 1, 3, 7, 9 \}$\\
6&	9&	5&	4&	10&	3&	0&	3&	6&	2&	4&	$\{ 1, 3, 7 \}$\\
7&	1&	5&	4&	9&	3&	0&	3&	6&	2&	4&	$\{ 3, 7 \}$\\
8&	7&	4&	4&	9&	3&	0&	3&	6&	2&	4&	$\{ 3 \}$\\
9&	3&	4&	4&	9&	3&	0&	3&	6&	2&	4&	$\emptyset$\\   
\end{tabular}
\end{center}
\end{table}
\item Choose $g_i^s$ to be the \emph{actual} shortest path cost from $i$ to $s$.
This is the tightest possible ``lower bound,'' and will make $A^*$ run extremely quickly --- in fact, it will only scan nodes along the shortest path, the best possible performance that can be achieved.
We see this in Table~\ref{tbl:astaromniscient} for the network in Figure~\ref{fig:grid-astar}: just three iterations are needed!
However, coming up with these ``estimates'' is just as hard as solving the original problem.
So in the end we aren't saving any effort; what we gain from $A^*$ is more than lost by the extra effort we need to compute $g_i^s$ in the first place.
\begin{table}
\begin{center}
\caption{Omniscient $A^*$ with perfect lower bounds.
         \label{tbl:astaromniscient}}
\begin{tabular}{cc|ccccccccc|c}
&	$g^3_i$&	7&	5&	0&	9&	9&	7&	12&	11&	13&		\\
Iteration&	$i$&	$L_1$&	$L_2$&	$L_3$&	$L_4$&	$L_5$&	$L_6$&	$L_7$&	$L_8$&	$L_9$&	$SEL$\\
\hline
0&	---&	$\infty$&	$\infty$&	$\infty$&	$\infty$&	0&	$\infty$&	$\infty$&	$\infty$&	$\infty$&	$\{ 5 \}$\\
1&	5&	$\infty$&	4&	$\infty$&	3&	0&	3&	$\infty$&	2&	$\infty$&	$\{ 2, 4, 6, 8 \}$\\
2&	2&	6&	4&	9&	3&	0&	3&	$\infty$&	2&	$\infty$&	$\{ 1, 3, 4, 6, 8 \}$\\
3&	3&	6&	4&	9&	3&	0&	3&	$\infty$&	2&	$\infty$&	$\{ 1, 4, 6, 8 \}$\\
\end{tabular}
\end{center}
\end{table}
\end{itemize}
So, there is a tradeoff between choosing tight bounds (the closer $g_i^s$ to the true costs, the faster $A^*$ will be) while not spending too long in computing the estimates (which might swamp the savings in $A^*$ itself).
Luckily, in transportation networks, there are several good bounds available which can be computed fairly quickly.
For instance:
\begin{itemize}
\item The Euclidean (``as the crow flies'') distance between $i$ and $s$, divided by the fastest travel speed in the network, is a lower bound on the travel time between $i$ and $s$.
\item Replace every link with a lower bound on its cost (say, free-flow travel time) and find shortest paths between all nodes and all destinations (repeatedly using one of the previous algorithms from this chapter).
This takes more time, but only needs to be done once and can be done as a preprocessing step.
As we will see in Chapter~\ref{chp:solutionalgorithms}, solving traffic assignment requires many shortest path computations.
The extra time spent finding these costs once might result in total time savings over many iterations.
\end{itemize}
You may find it instructive to think about other ways we can estimate $g_i^s$ values, and how they might be used in transportation settings.
\index{A$^*$|)}
\index{shortest path|)}

\section{Historical Notes and Further Reading}
\label{sec:graphreferences}

Networks arise frequently in transportation and optimization, particularly network flow problems such as shortest path, maximum flow,\index{maximum flow} and minimum cost flow.\index{minimum cost flow}
The text by \cite{ahuja93} discusses a number of ways to formulate and solve such optimization problems, along with a breadth of applications.
Other good surveys of network optimization are found in \cite{rockafellar98} and \cite{bertsekas_net}.
\cite{ahuja93} and \cite{tarjan83} also provide more detail on data structures used for working with networks in computer programs.
Readers interested in graph theory\index{graph theory} apart from optimization are referred to the texts by \cite{busaker65}, \cite{bondy76}, and \cite{diestel16}.
A detailed exposition of tree networks in particular is found in Section 2.3 of \cite{knuth1}.

Bellman's principle\index{Bellman's principle} was first identified in the context of dynamic programming~\index{dynamic programming}\citep{bellman57}, a broader class of optimization problems which includes the shortest path problem~\citep{bertsekas_dp1, bertsekas_dp2}.
Dynamic programming forms the foundation of Markov decision processes and reinforcement learning methods which are currently popular in artificial intelligence~\citep{powell07}.

The first label-correcting\index{shortest path!label correcting} shortest path algorithm was developed by~\cite{ford56}.
Subsequent experience has shown that the practical performance of this method depends significantly on the order in which nodes are retrieved from the scan list; good choices are ``first-in, first-out,'' (FIFO) in which the node selected from $SEL$ is always one of the \emph{oldest} there, and a double-ended queue implementation which provides an exception to FIFO in that when nodes are return to the scan list after being scanned earlier they are moved to the front of the queue.
These choices were first described in~\cite{bellman58} and~\cite{pape74}.
The quintessential label-setting shortest path method\index{shortest path!label setting} is that of \cite{dijkstra59}.
The $A^*$\index{A$^*$} method is essentially an extension to Dijkstra's algorithm, and was developed by \cite{hart68}.
For more thorough reviews of shortest path algorithms, the reader is referred to the survey paper by~\cite{deo84}, and to~\cite{ahuja93}.

Most recently, contraction hierarchies\index{shortest path!contraction hierarchy} have proved to be very efficient ways of solving shortest path problems for traffic assignment~\citep{schneck20}.
These algorithms are beyond the scope of this book, but the idea is to simplify the underlying network as a preprocessing step.
Once this preprocessing is done, shortest path problems can be solved very quickly~\citep{geisberger12}, even when the link costs change from one shortest path problem to the next~\citep{dibbelt16}.

In some applications, we may need to find multiple low-cost paths between an origin and destination.
The $k$-shortest path problem specifically aims to find $k$ distinct paths that have the lowest cost.
There are variations of this problem; for instance, sometimes the second shortest path may just be the shortest path with a small cycle added.
Sometimes this may be acceptable, other times you may need paths that differ more substantially.
\cite{yen71}, \cite{eppstein98}, and \cite{chen16} all discuss approaches for the $k$-shortest path problem.  \index{shortest path!$k$ shortest path}

There are many variants of the shortest path program, for instance, adding side constraints on the amount of some other resource consumed (such as money spent on tolls or battery charge from an electric vehicle).\index{shortest path!resource-constrained}
Algorithms for this type of problem are discussed in \cite{pugliese12}, \cite{lozano16}, and \cite{himmich20}, among others.
\cite{ziegelmann_diss} provides a review of approaches for this problem.
For electric vehicles which can charge \emph{en route}, the ``shortest path problem with relays'' is particularly appropriate, as discussed in \cite{laporte11} and \cite{baum19}.\index{shortest path!with relays}
Multicriteria optimization is another way to address issues related to multiple kinds of costs~\citep{chen13}.\index{shortest path!multicriteria}

\index{shortest path!stochastic|(}
Another common variant is aimed at addressing reliability and uncertainty in link costs, by having the link costs be drawn from some probability distribution rather than known exactly in advance.
If one aims to minimize the expected travel time, it is enough to replace each link cost with its expected value.
However, travelers may also care about reliability specifically.
This can be modeled by adding variance or standard deviation to the objective function~\citep{xing11,shahabi13,khani15,zhang16,zhang17,zhang19}, adding a constraint on variance~\citep{sivakumar94}, adopting a robust optimization perspective~\citep{yu98,montemanni04,shahabi15}, applying multiobjective optimization~\citep{sen01}, introducing a nonlinear utility function to represent arrival time preferences~\citep{loui83,eiger85,murthy96,boyles06,gao05}, or changing the objective entirely to maximize the probability of on-time arrival~\citep{fan05,nie09}.

Yet another way to address reliability and real-time information provision is to allow the traveler's path to change \emph{en route} based on information learned while traveling; this leads to a class of \emph{online}\index{shortest path!online} shortest path formulations.
For examples of such formulations, see \cite{andreatta88}, \cite{psaraftis93}, \cite{polychronopoulos96},  \cite{millerhooks01}, \cite{waller02}, \cite{provan03}, \cite{gao05}, \cite{boyles06}, \cite{boyles09diss}, and \cite{boyles16}.
This concept has also been applied in shortest path problems modeling route choice among public transit users, who may choose different bus or train routes depending on which arrives first.
These settings often refer to the ``shortest hyperpath'' problem to identify optimal strategies~\citep{nguyen88,decea93,wu94,khani19}.
\index{shortest path!stochastic|)}

It is often helpful to solve a shortest path problem where link costs vary over time.
This particular extension will be covered in Chapter~\ref{chp:tdsp}.

\section{Exercises}
\label{exercises_networkrepresentations}

\begin{enumerate}
\item \diff{9} In the network in Figure~\ref{fig:exercisenet}, list the indegree, outdegree, degree, forward star, and reverse star of each node. \label{ex:netex1}
\stevefig{exercisenet}{Network for Exercises~\ref{ex:netex1}--\ref{ex:netex3}}{0.5\textwidth}
\item \diff{3} In the network in Figure~\ref{fig:exercisenet}, list all of the paths between nodes 1 and 5. \label{ex:netex2}
\item \diff{4} State whether the network in Figure~\ref{fig:exercisenet} is or is not (a) cyclic; (b) a tree; (c) connected; (d) strongly connected. \label{ex:netex3}
\item \diff{15} For each of the following, either draw a network with the stated properties, or explain why no such network can exist: (a) connected, but not strongly connected; (b) strongly connected, but not connected; (c) cyclic, but not strongly connected.
\item \diff{25} If $m$ and $n$ are the number of links and nodes in a network, show that $m < n^2$.

\item \diff{25} If a network is connected, show that $m \geq n - 1$.
\item \diff{25} If a network is strongly connected, show that $m \geq n$.
\item \diff{10} Show that $\sum_{i \in N} |\Gamma(i)| = \sum_{i \in N} |\Gamma^{-1}(i)| = m$.
\item \diff{1} Why must transportation infrastructure networks be cyclic?
\item \diff{15} Although we do not usually expect drivers to use cyclic paths, there are some exceptions.
Name one.
\item \diff{15} Find a topological order for the network in Figure~\ref{fig:topologicalexercise}. \label{ex:topex}
\stevefig{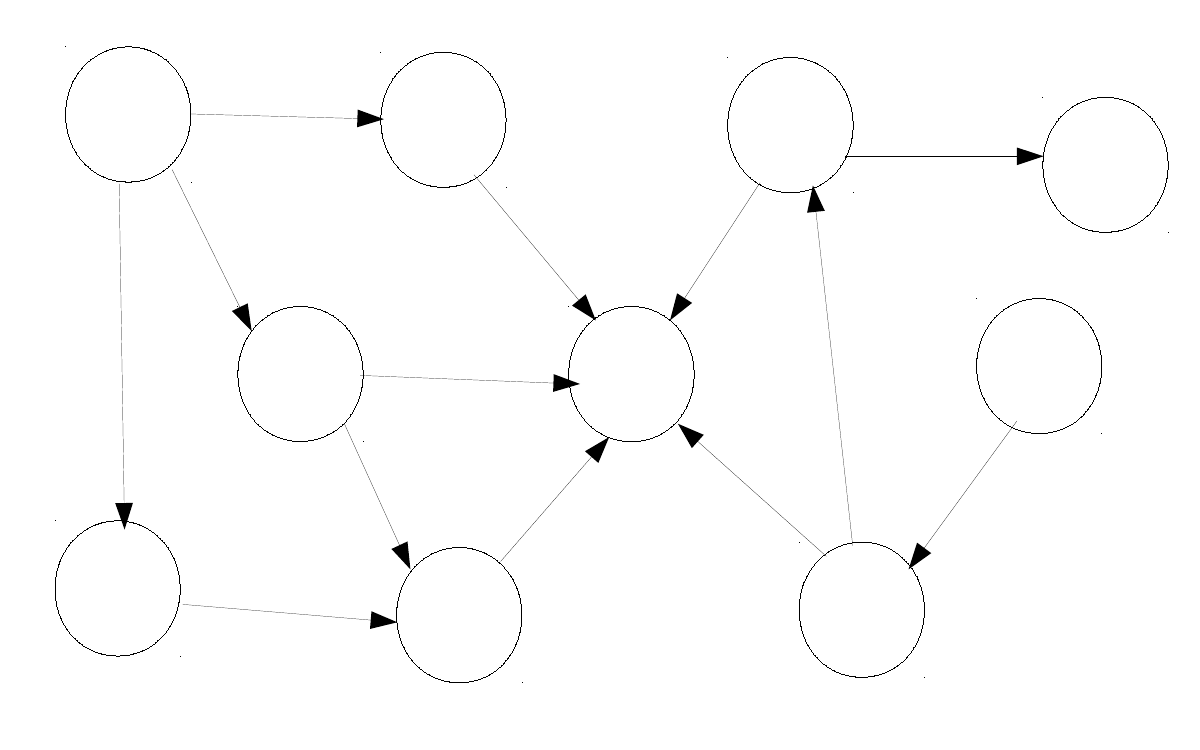}{Network for Exercise~\ref{ex:topex}}{0.5\textwidth}
\item \diff{11} Is the topological order for an acyclic network unique?  Either explain why it is, or provide a counterexample showing it is not.
\item \diff{42} Show that any acyclic network has at least one node with no incoming links.
This was the ``missing step'' in the proof of Theorem~\ref{thm:acyclictopological}.
(Since this is part of the proof, you can't use the result of this theorem in your answer.)
\item \diff{36} Let $i$ and $j$ be any two nodes in an acyclic network.
Give a procedure for calculating the number of paths between $i$ and $j$ which involves at most one calculation per node in the network.
\item \diff{57} Give a procedure for determining whether a given network is strongly connected.
Try to make your method require as few steps as possible.
\item \diff{58} Give a procedure for determining whether a given network is connected.
Try to make your method require as few steps as possible.
\item \diff{30} Show that there is at most one path between \emph{any} two nodes in a tree. 
\item \diff{30} Show that any node in a tree (excluding its root) has exactly one link in its reverse star (and thus one parent).
\item \diff{30} Show that any tree is an acyclic network.
\item \diff{30} Show that removing any link from a tree produces a disconnected network. 
\item \diff{50} Show that any tree has at least two nodes of degree one.
\item \diff{55} Consider a rectangular grid network of one-way streets, with $r$ rows of nodes and $c$ columns of nodes.
All links are directed northbound and eastbound.
How many paths exist between the lower-left node (southwest) and the upper-right (northeast) node?
\item \diff{14} Write down the node-node adjacency matrix of the network in Figure~\ref{fig:exercisenet}.
\item \diff{16} Consider the network defined by this node-node adjacency matrix:
\[\vect{ 0 & 0 & 0 & 1 & 1 & 0 & 0 
      \\ 1 & 0 & 1 & 0 & 0 & 0 & 0 
      \\ 0 & 0 & 0 & 1 & 0 & 0 & 0 
      \\ 0 & 0 & 0 & 0 & 1 & 0 & 0 
      \\ 0 & 0 & 0 & 0 & 0 & 0 & 0 
      \\ 1 & 1 & 0 & 0 & 0 & 0 & 0 
      \\ 0 & 1 & 0 & 0 & 0 & 1 & 0 }\]

   \begin{enumerate}[(a)]
   \item Draw the network.
   \item How many links does the network have?
   \item How many links enter node 1?  How many links leave node 1?
   \end{enumerate}
\item \diff{17} Is the network represented by the following node-node adjacency matrix strongly connected?
 \[\vect{   0 & 0 & 0 & 1 & 1 & 0 & 0 
         \\ 1 & 0 & 1 & 0 & 0 & 0 & 0 
         \\ 0 & 0 & 0 & 1 & 0 & 0 & 0 
         \\ 0 & 0 & 0 & 0 & 1 & 0 & 0 
         \\ 0 & 0 & 0 & 0 & 0 & 1 & 0 
         \\ 1 & 1 & 0 & 0 & 0 & 0 & 0 
         \\ 0 & 1 & 0 & 0 & 0 & 1 & 0 }\]
\item \diff{10} If the nodes in an acyclic network are numbered in a topological order, show that the node-node adjacency matrix is upper triangular.
\item \diff{37} Let $\mb{A}$ be the node-node adjacency matrix for a network.
What is the interpretation of the matrix product $\mb{A}^2$?
\item \diff{42} Let $\mb{A}$ be the node-node adjacency matrix for an acyclic network.
First show that $\sum_{n=1}^\infty \mb{A}^n$ exists, and give an interpretation of this sum.
\item \diff{65} A \emph{unimodular}\index{matrix!unimodular} matrix is a square matrix whose elements are integers and whose determinant is either $+1$ or $-1$.
A matrix is \emph{totally unimodular}\index{matrix!totally unimodular} if every nonsingular square submatrix is unimodular.
(Note that a totally unimodular matrix need not be square).
Show that every node-link incidence matrix is totally unimodular.
\item \diff{10} One disadvantage of the forward star representation is that it is time-consuming to identify the reverse star of a node --- one must search through the entire array to find every link with a particular head node.
Describe a ``reverse star'' data structure using arrays, where the reverse star can be easily identified.
\item \diff{52} By combining the forward star representation from the text and the reverse star representation from the previous exercise, we can quickly identify both the forward and reverse stars of every node.
However, a naive implementation will have two different sets of arrays, one sorted according to the forward star representation, and the other sorted according to the reverse star representation.
This duplication wastes space, especially if there are many attributes associated with each link (travel time, cost, capacity, etc.)  Identify a way to easily identify the forward and reverse stars of every node, with only \emph{one} set of arrays of link data, by adding an appropriate attribute to each link.
\item \diff{58} In the language of your choice, write computer code to do the following:
\begin{enumerate}
\item Produce the node-node adjacency matrix of a network when given the node-link incidence matrix.
\item Produce the node-link incidence matrix of a network when given the node-node adjacency matrix.
\item Produce the node-node adjacency matrix of a network when given the forward star representation of the network.
\item Produce the forward star representation of a network when given the node-node adjacency matrix.
\end{enumerate}
\item \diff{13} After solving a shortest path problem from node 3 to every other node, I obtain the backnode vector shown in Table~\ref{tbl:backnode}.
Write the shortest paths (a) from node 3 to node 5; (b) from node 3 to node 7; (c) from node 4 to node 8. \label{ex:backnode}
\begin{table}
\begin{center}    
\caption{Backnode vector for Exercise~\ref{ex:backnode}. \label{tbl:backnode}}
\begin{tabular}{|c|c|}
  \hline
  Node  & Backnode \\
  \hline
  1     &     4  \\
  2     &     6  \\
  3     &     $-1$  \\
  4     &     7  \\
  5     &     4  \\
  6     &     5  \\
  7     &     3  \\
  8     &     10 \\
  9     &     5  \\
  10    &     6  \\
  \hline
\end{tabular}
\end{center}
\end{table}
\item \diff{26} Find the shortest path from node 1 to every other node in the network shown in Figure~\ref{fig:lpfail}.
Report the final labels and backnodes ($L$ and $q$ values) for all nodes.
\label{ex:sppractice}
\stevefig{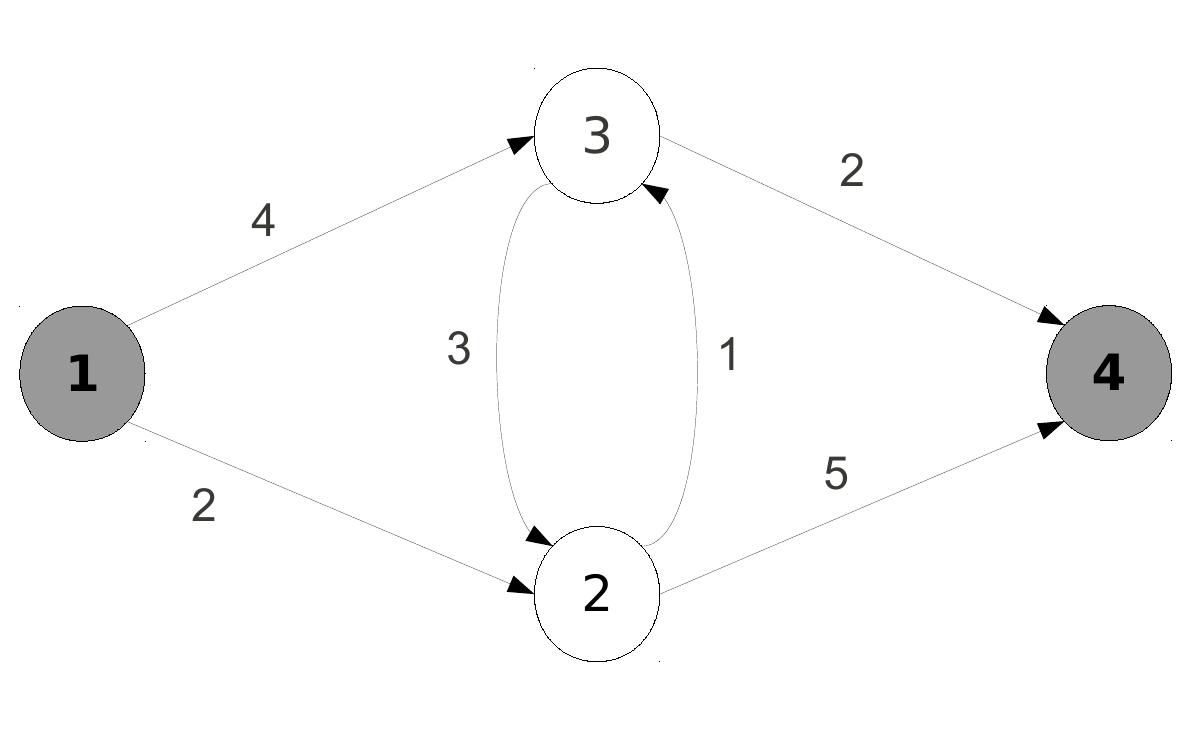}{Network for Exercises~\ref{ex:sppractice} and~\ref{ex:lpfail}.}{0.6\textwidth}
\item \diff{25} The network in Figure~\ref{fig:negcostsp} has a link with a negative cost.
Show that the label-correcting algorithm still produces the correct shortest paths in this network, while the label-setting algorithm does not. \label{ex:negativelabelcorrecting}
\stevefig{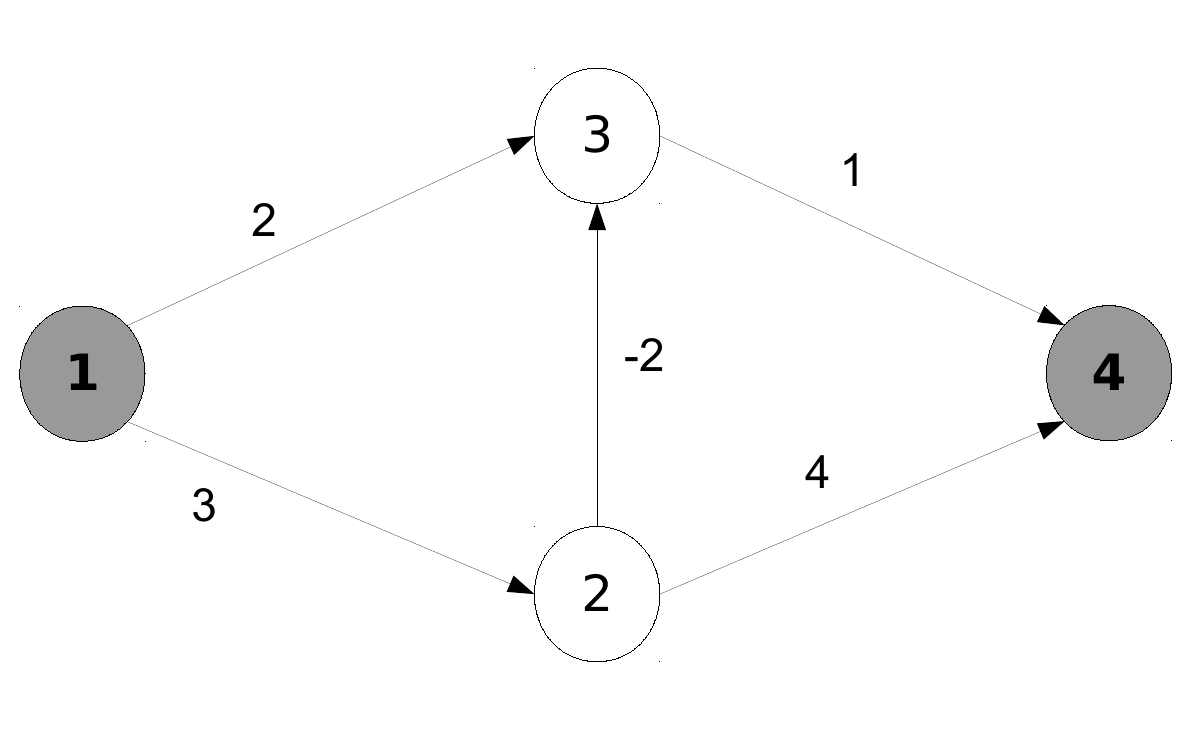}{Network for Exercise~\ref{ex:negativelabelcorrecting}.}{0.6\textwidth}
\item \diff{57} Prove or disprove the following statement: ``Any network with negative costs can be transformed into a network with nonnegative costs by adding a large enough constant to every link's cost.
We can then use the label-setting algorithm on this new network.
Therefore, the label-setting algorithm can find the shortest paths on any network, even if it has negative-cost links.''  (Proving this statement means showing that it is true in any network; to disprove it, it is sufficient to find a single counterexample.) 
\item \diff{37}.
Find the shortest paths on the network in Figure~\ref{fig:spexercisenet}, using both the label-correcting and label-setting algorithms.
Sketch the \emph{shortest path tree}\index{tree!shortest path} (that is, the network with only the links implied by the backnode vector) produced by each algorithm.
Which algorithm required fewer iterations? \label{ex:spexercisenet}
\stevefig{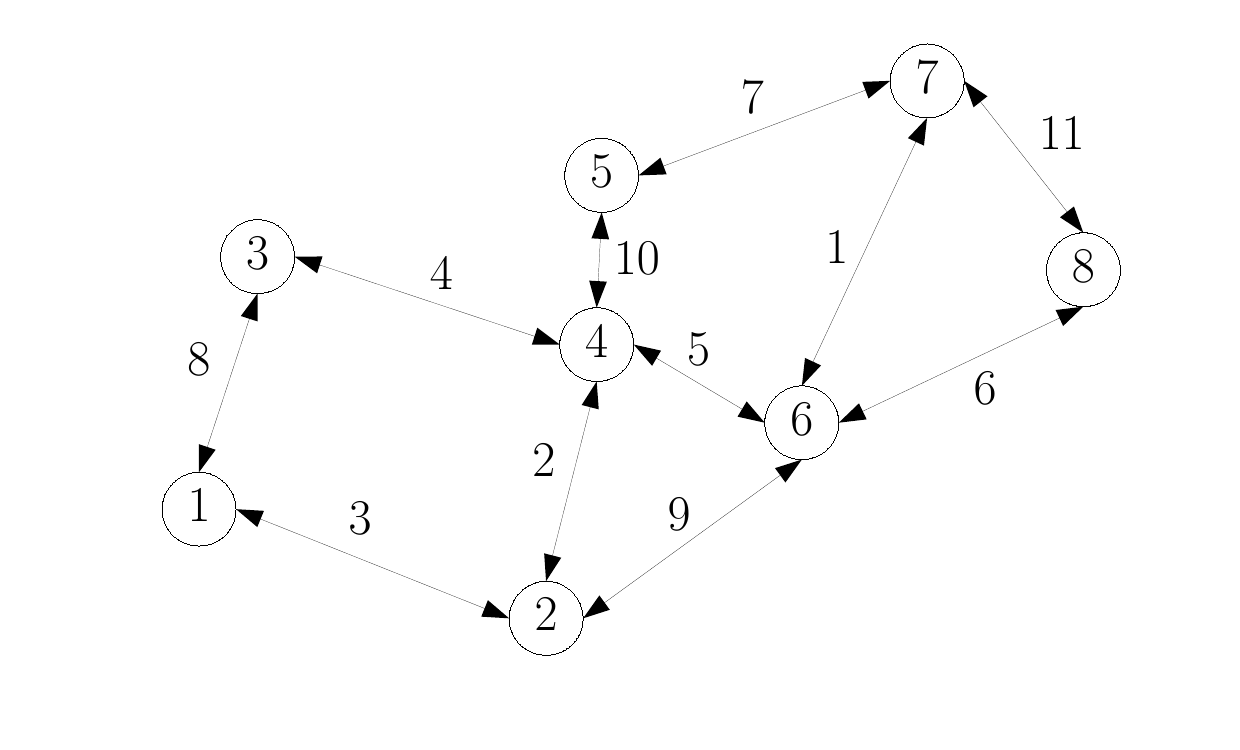}{Network for Exercise~\ref{ex:spexercisenet}, link labels are costs.}{0.7\textwidth}
\item \diff{33}.
You and your friends are camping at Yellowstone Park, when you suddenly realize you have to be at the Denver airport in ten hours to catch a flight.
You don't have Internet access, but you do have an atlas showing the distance and travel time between selected cities (Figure~\ref{fig:atlas}).
Assuming these travel times are accurate, can you make it to Denver in time for your flight? \label{ex:yellowstonedenver}
\genfig{atlas}{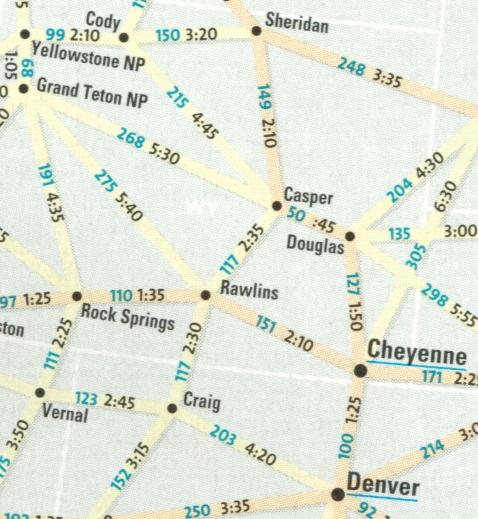}{Atlas page for Exercise~\ref{ex:yellowstonedenver}.}{width=0.5\textwidth}
\item \diff{43}.
Consider the following variation of the shortest path problem: instead of the costs $c_{ij}$ being fixed, instead assume that they are random, and the cost of link $(i,j)$ takes the values $\myc{c_{ij}^1, c_{ij}^2, \ldots, c_{ij}^s}$ with probabilities $\myc{p_{ij}^1, p_{ij}^2, \ldots, p_{ij}^s}$ independent of the cost of any other link.
How can you adapt the shortest path algorithms in this chapter to find the path with the least \emph{expected} cost?
\item \diff{53}.
Assume that each link $(i,j)$ in the network fails with probability $p_{ij}$, and that link failures are independent.
(Failure can represent some kind of damage or disruption, the probability of being detected in a military routing problem, etc.)  A path fails if any link in that path fails.
Explain how you can find the most reliable path (that is, the one with the least failure probability) with a shortest path algorithm.
\item \diff{23} Instead of trying to find the \emph{shortest} path between two nodes, let's try to find the \emph{longest} path between two nodes.
(As we will see later, there are actually cases when this is useful.) \label{ex:lpfail}
\begin{enumerate}[(a)]
\item Modify the algorithm presented in Section~\ref{sec:acyclicsp} to find the longest path in an acyclic network.
\item If you modify the label-correcting or label-setting algorithms for general networks in a similar way, will they find the longest paths successfully?  Try them on the network in Figure~\ref{fig:lpfail}. 
\end{enumerate} 

\item \diff{25}.
As a modification to the node-node adjacency matrix you might use to represent the network in a computer program (cf.\ Section~\ref{sec:datastructures}), you could store the cost of the links in this matrix, with $\infty$ where no link exists (as opposed to `1' where links exist and `0' where they do not).
Find the shortest path between nodes 7 and 4 using the modified adjacency matrix below, where the entry in row $i$ and column $j$ is the cost $c_{ij}$ if this link exists, and $\infty$ if it does not.
\[
\vect{ \infty & \infty & \infty & 5      & 3      & \infty & \infty \\
      4      & \infty & 7      & \infty & \infty & \infty & \infty \\
      \infty & \infty & \infty & 6      & \infty & \infty & \infty \\
      \infty & \infty & \infty & \infty & 2      & \infty & \infty \\
      \infty & \infty & \infty & \infty & \infty & 4      & \infty \\
      3      & 2      & \infty & \infty & \infty & \infty & \infty \\
      \infty & 1      & \infty & \infty & \infty & 4      & \infty }
\]
\item \diff{22}.
In the game ``Six Degrees of Kevin Bacon,''\index{Six Degrees of Separation} players are given the name of an actor or actress, and try to connect them to Kevin Bacon in as few steps as possible, winning if they can make the connection in six steps or less.
Two actors or actresses are ``connected'' if they were in the same film together.
For example, Alfred Hitchcock is connected to Kevin Bacon in three steps: Hitchcock was in \emph{Show Business at War} with Orson Welles, who was in \emph{A Safe Place} with Jack Nicholson, who was in \emph{A Few Good Men} with Kevin Bacon.
Beyonc\'{e} Knowles is connected to Kevin Bacon in two steps, since she was in \emph{Austin Powers: Goldmember}, where Tom Cruise had a cameo, and Cruise was in \emph{A Few Good Men} with Bacon.
Assuming that you have total, encyclopedic knowledge of celebrities and films, show how you can solve ``Six Degrees of Kevin Bacon'' as a shortest path problem.
Specify what nodes, links, costs, origins, and destinations represent in the network you construct.
\item \diff{14} In a network with no negative-cost cycles, show that $-nC$ is a lower bound on the shortest path cost between any two nodes in a network.
\label{ex:negcycletighterbound}
\item \diff{57} Show that if the label-correcting algorithm is performed, and that at each iteration you choose a node in $SEL$ which has been in the list the longest, at the end of $kn$ iterations the cost and backnode labels correctly reflect all shortest paths from $r$ which are no more than $k$ links long.
(Hint: try an induction proof.)
\item \diff{44} Assume that the label-correcting algorithm is terminated once the label for some node $i$ falls below $-mC$.
Show that the following the current backnode labels from $i$ will lead to a negative-cost cycle. \label{ex:negcycletraceback}
\item \diff{32} Modify three of the shortest path algorithms in this chapter so that they find shortest paths from all origins to one destination, rather than one origin to all destinations: \label{ex:reversesp}
\begin{enumerate}[(a)]
\item The acyclic shortest path algorithm from Section~\ref{sec:acyclicsp}.
\item The label-correcting algorithm from Section~\ref{sec:spgeneral}.
\end{enumerate} 
\item \diff{73} \label{ex:fasterspwithorigin} It is known that the label correcting algorithm will find the correct shortest paths as long as the initial labels $\mb{L}$ correspond to the distance of \emph{some} path from the origin (they do not necessarily need to be initialized to $+\infty$).
Assume that we are given a vector of backnode labels $\mb{q}$ which represents some tree (not necessarily the shortest paths) rooted at the origin $r$.
Develop a one-to-all shortest path algorithm that uses this vector to run more efficiently.
In particular, if the given backnode labels $\mb{q}$ \emph{do} correspond to a shortest path tree, your algorithm should recognize this fact and terminate in a number of steps linear in the number of network links.
\end{enumerate}

\chapter{Mathematical Techniques for Equilibrium}
\label{chp:mathematicalpreliminaries}

This chapter provides a survey of mathematical techniques used in network analysis. There is no attempt to be comprehensive; many indeed, books have been written about each of the sections in this chapter.
Rather, the intent is to cover topics which are used frequently in transportation network problems.

Several of the appendices may be useful at this point.
Appendix~\ref{chp:mathbackground} reviews definitions and facts related to vectors, matrices, sets, and functions which are needed, including the topics of convex functions, convex sets, and multivariable calculus (the gradient vector, and the Jacobian and Hessian matrices will play a particularly important role).
If these topics are new to you, you may wish to consult more extended treatments of these topics in other books or references.
Appendix~\ref{chp:basicoptimization} introduces basic concepts of optimization, focusing on what will ultimately be relevant for static and dynamic traffic assignment.
Optimization is a much richer and deeper field than what is used in traffic assignment, and Appendices~\ref{chp:fancyoptimization} and~\ref{chp:algorithmcomplexity} go into further detail here --- the material from these latter two appendices is not strictly necessary for this book, but we believe it to be of interest to many readers nevertheless.

This chapter focuses on mathematical material which is both specialized to the traffic assignment problem, and which is likely to be new to anticipated readers of the book.
This will involve discussion of three main techniques: the \emph{fixed point problem}, in Section~\ref{sec:fixedpoint}; the \emph{variational inequality}, in Section~\ref{sec:variationalinequality}, and \emph{convex optimization}, in Section~\ref{sec:convexoptimization}.
Chapter~\ref{chp:introchapter} characterized equilibrium as a ``consistent'' state in which no driver can improve his or her satisfaction by unilaterally changing routes.
Under the assumption of continuous flow variables, we can use calculus to greatly simplify the problem.
Each of the mathematical techniques discussed in this chapter formalizes this equilibrium principle in different ways.

Each approach has its own advantages from the standpoint of this book.
The fixed point formulation is perhaps the most intuitive and generally-applicable definition, but does not give much indication as to how one might actually find this equilibrium.
The variational inequality formulation lends itself to physical intuition and can also accommodate a number of variations on the equilibrium problem.
The convex optimization approach provides an intuitive interpretation of solution methods, provides an elegant proof of equilibrium uniqueness in link flows, and powers the best-known solution algorithms, but the connection between the equilibrium concept and optimization requires more mathematical explanation and is less obvious at first glance.

\section{Fixed Point Problems}
\label{sec:fixedpoint}

\index{fixed point problem|(}
A \textbf{fixed point} of a function $f$\label{not:f} is a value $x$\label{not:x} such that $f(x) = x$, that is, the value $x$ is unchanged by $f$.
As an example, the function $f_1(x) = 2x - 1$ has only one fixed point at $x = 1$, because $f_1(1) = 2 \times 1 - 1 = 1$, the function $f_2(x) = x^2$ has two fixed points at 0 and 1 ($0^2 = 0$ and $1^2 = 1$), while the function $f_3(x) = x^2 + 1$ has no fixed points at all.
A helpful visual illustration is that the fixed points of a function are the points of intersection between the function's graph and the 45-degree line $y = x$\label{not:y} (Figure~\ref{fig:fixedpointgraph}).

Equilibrium solutions can often be formulated as fixed points of a suitable function.
For instance, in the traffic assignment problem, a route choice model and a congestion model are mutually interdependent: drivers choose routes to avoid congestion, but congestion is determined by the routes drivers choose.
(Figure~\ref{fig:taschematic}).
An equilibrium solution is consistent\index{consistency} in the sense that drivers are satisfied with the travel times calculated by the paths they chose.\index{equilibrium!and fixed points}
Feeding the route choices into the congestion model, then feeding the resulting travel times into the route choice model, one can obtain the original route choices back.
This is reminiscent of fixed point problems: when you evaluate a function at its fixed point, after performing whatever calculations the function requires you obtain the fixed point again.
Fixed points thus arise naturally when dealing with these kinds of ``circular'' dependencies.

As a first example, consider the problem of trying to estimate the number of bus riders in a dense urban area.
The bus system is subject to congestion; when there are $x$ riders, the average delay to customers is given by the function $t = T(x)$\label{not:t}\label{not:Tminor} which is assumed continuous and increasing.
(Assuming that the fleet of buses and timetables are fixed, more riders mean longer boarding and offloading times, more frequent stops, and the possibility of denied boarding when buses are full.)  However, the number of bus riders depends on the congestion in the system --- as the buses become more crowded, some riders may switch to alternate modes of transportation or combine trips, so we can write $x = X(t)$\label{not:Xminor} for some function $X$ which is continuous and decreasing.
For a concrete example, assume that the system is designed such that $T(x) = 10 + x$, and that the ridership function is given by $X(t) = [20 - t]^+$\label{not:bracketplus}, when $x$ and $t$ are measured in appropriate units (say, thousands of passengers and minutes).\footnote{The notation $[\cdot]^+$ is used to mean the \emph{positive part} of the term in brackets, that is, $[y]^+ = \max\{y,0\}$.
If the term in brackets is negative, it is replaced by zero; otherwise it is unchanged.}

The goal is to find the ridership $x$; but $x = X(t)$ and $t = T(x)$, which means that we need to find some value of $x$ such that $x = X(T(x))$.
This is a fixed point problem!
Here the function $f$ is the composition of $X$ and $T$: 
\begin{align*}
   x &= X(T(x)) \\
     &= [20 - T(x)]^+ \\
     &= [20 - (10 + x)]^+ \\
     &= [10 - x]^+
\,.
\end{align*}
Assuming that $10 - x$ is nonnegative, we can replace the right-hand side of the last equation by $10 - x$.
Solving the resulting equation $x = 10 - x$, we obtain $x = 5$ as the fixed point.
Checking, when there are 5 riders the average travel time will be 15 minutes (based on $T$); and when the travel time is 15 minutes, there will indeed be 5 riders (based on $X$).
Finally, if $x = 5$, our assumption that $10 - x$ was nonnegative was true so this solution is valid.

For this example, it is relatively easy to find the fixed point by substituting the definitions of the functions and performing some algebra.
However, in more complex problems it will be difficult or impossible to calculate the fixed point directly.
Despite this, fixed points are important because of so-called ``fixed point theorems'' which guarantee the existence of a fixed point under certain conditions of the function $f$.
A fixed point theorem by Brouwer is provided below, and another by Kakutani is provided in Section~\ref{sec:eqmmultifunction}.
These fixed point theorems are \emph{non-constructive} because they give us no indication of how to find a fixed point, they simply guarantee that at least one exists.
Brouwer's theorem can be stated as:
\begin{thm} \emph{(Brouwer).}\label{thm:brouwersthm}\index{Brouwer's theorem}
Let $f$ be a continuous function from the set $K$\label{not:K} to itself, where $K$ is convex and compact.
Then there is at least one point $x \in K$ such that $f(x) = x$.
\index{fixed point problem!existence of solutions}
\index{continuous function!applications}
\index{convex set!applications}
\index{compact set!applications}
\end{thm}
The exercises ask you to show that each of the conditions in Brouwer's theorem is necessary; you might find it helpful to visualize these conditions geometrically similar to Figure~\ref{fig:fixedpointgraph}.

Notice the stipulation that $f$ be a function ``from the set $K$ to itself;'' this means that the range of its function must be contained in its domain.
Intuitively, this means that any ``output'' of the function $f$ must also be a valid ``input'' to that same function.
In other words, iteration is possible: starting from any value $x$ in its domain, you can apply the function over and over again to produce a sequence of values $x, f(x), f(f(x)), \ldots$.
If this condition does not hold, then Brouwer's theorem does not guarantee anything about a fixed point.

Let us apply Brouwer's theorem to the transit ridership example.
Both $X$ and $T$ are continuous functions, so their composition $X \circ T$ is continuous as well.
(Alternately, by substituting one function into the other we obtain $X(T(x)) = [10 - x]^+$, which is evidently continuous.)  What are the domain and range of $X(T(x))$?
Since $X(T(x))$ is the positive part of $20 - T(x)$, then $x = X(T(x)) \geq 0$.
Further note that because $x \geq 0$, $T(x) \geq 10$, so $x = X(T(x)) \leq 10$.
That is, we have shown that $x$ must lie between 0 and 10, so the function $X(T(x)) = [10 - x]^+$ can be defined from the set $[0, 10]$ to itself.
This set is convex and compact, so Brouwer's theorem would have told us that at least one fixed point must exist, even if we didn't know how to find it.\index{Brouwer's theorem!examples}\index{Brouwer's theorem!applications}\index{fixed point problem|)}

\stevefig{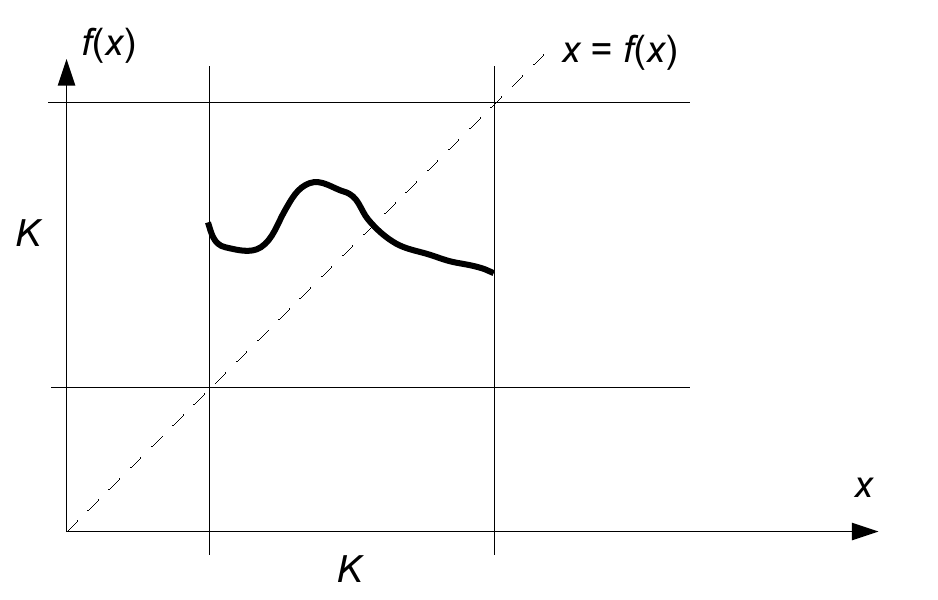}{Visualizing fixed point problems; the intersection of $f(x)$ and the 45-degree line is a point where $x = f(x)$.}{0.7\textwidth}

\section{Variational Inequalities}
\label{sec:variationalinequality}

\index{variational inequality|(}
Fixed point problems often lend themselves to elegant theorems like Brouwer's, which prove that a fixed point must exist.
However, such problems often lack easy solution methods.
The variational inequality can be more convenient to work with in this respect.
Variational inequalities can be motivated with a physical analogy.
Imagine an object, initially stationary, which is confined to move within some frictionless container (Figure~\ref{fig:container}) and cannot leave.
This object is acted on by a force whose magnitude and direction can be different at each point.
If the object is in the interior of the container, the object will begin to move in the same direction as the force at that point.
If the object starts at the edge of the container, it may not be able to move in the same direction as the force, but it might slide along the side of the object.
The problem is to determine where in the container (if anywhere) the object can be placed so that it will not move under the action of the force field.
Such a point is an equilibrium in the physical sense.
In the coming chapters, we will show how to connect this with the idea of traffic equilibrium introduced in Section~\ref{sec:notionofequilibrium}, by choosing the ``container'' and ``force field'' appropriately.
For now, it is enough to think about the problem in terms of physical forces and static equilibrium.

\stevefig{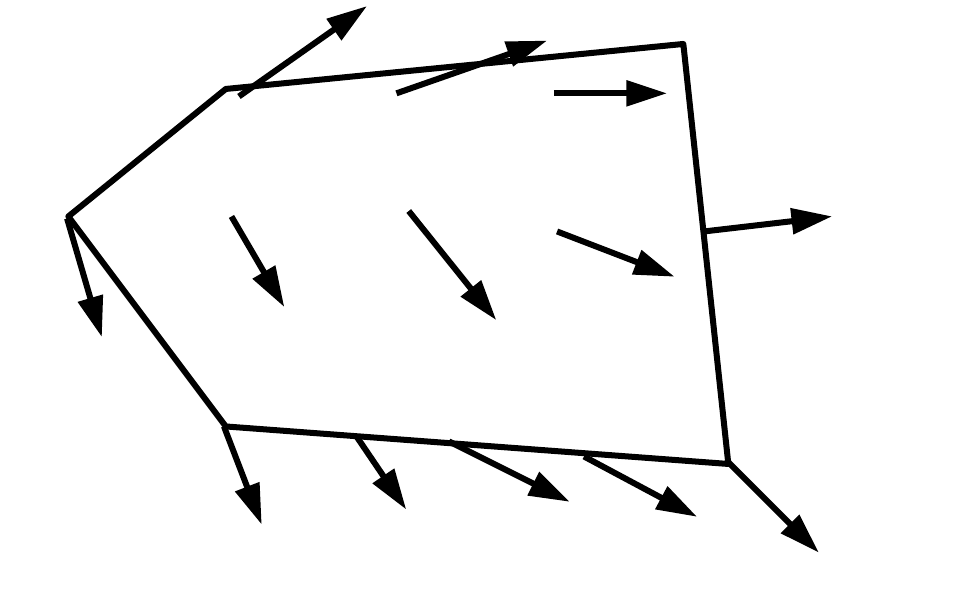}{A container with a force field.}{0.6\textwidth}

Figure~\ref{fig:viexamples} shows some examples.
In this figure, the direction of the force is drawn with black arrows, shown only at the points under consideration for clarity.
At three of the points (A, B, and C) the object will move: at A in the direction of the force, and at B and C sliding along the edge of the container in the general direction of the force.
At the other two points (D and E), the object will not move under the action of the force, being effectively resisted by the container wall.

\stevefig{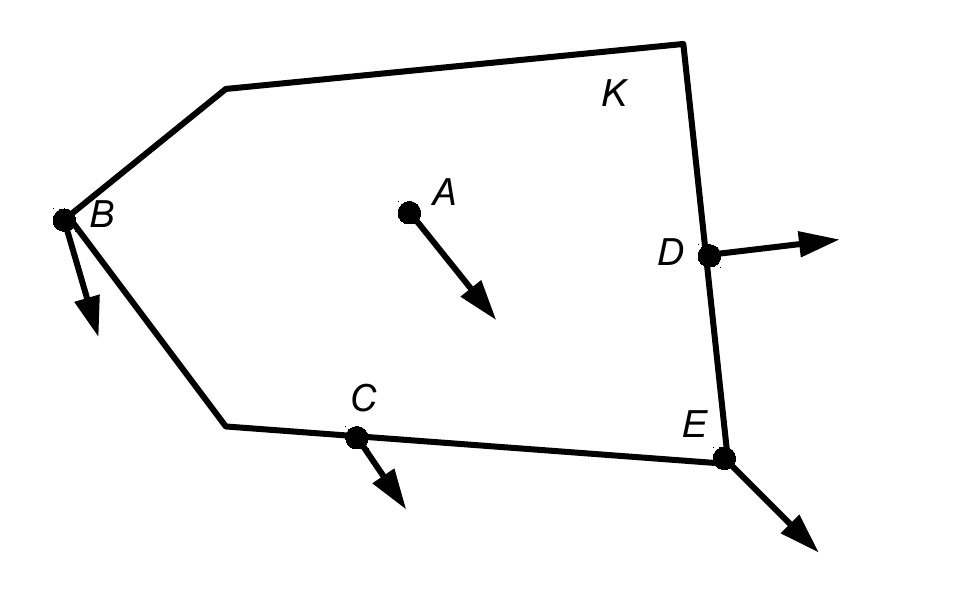}{Two solutions, and three non-solutions, to a variational inequality.}{0.6\textwidth}

How can we think about such problems in general?
A little thought should convince you that (1) if the object is on the boundary of the container, but not at a corner point, it will be unmoved if and only if the force is perpendicular to the boundary (point D in Figure~\ref{fig:viexamples}), and (2) if the object is at a corner of the container, it will be unmoved if and only if the force makes a right or obtuse angle with all of the boundary directions (point E).
These two cases can be combined together: a point on the boundary is an equilibrium if and only if the force makes a right or obtuse angle with all boundary directions.
In fact, if the force makes such an angle with all boundary directions, it will do so with any other direction pointing into the feasible set (Figure~\ref{fig:anglefig2}).
So, we see that \emph{a point is unmoved by the force if and only if the direction of the force at that point makes a right or obtuse angle with any possible direction the object could move in}.

\stevefig{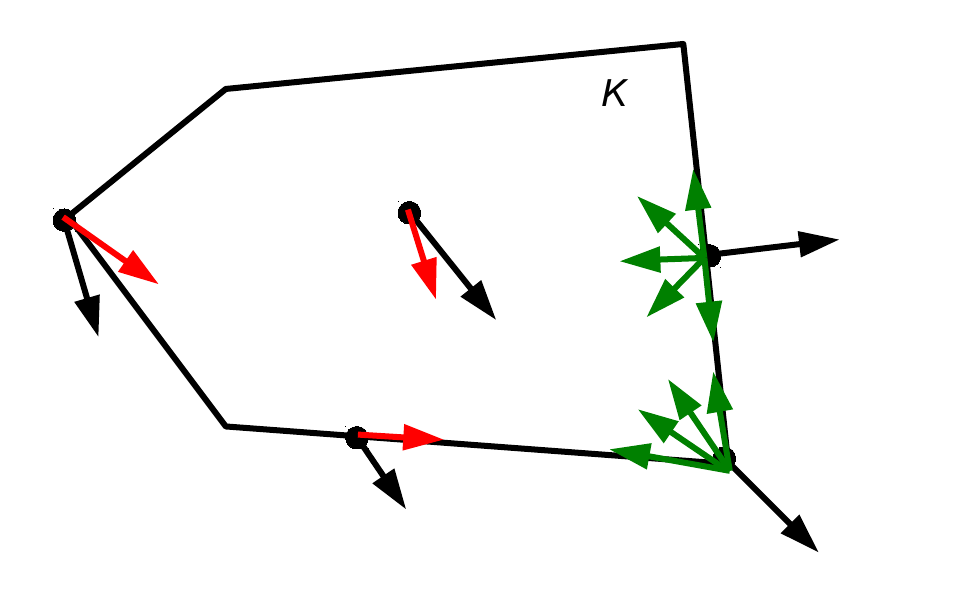}{Stable points (in green) make an obtuse or right angle with all feasible
directions; unstable points (in red) make an acute angle with some
feasible directions.}{0.6\textwidth}

The mathematical definition of a variational inequality is little more than translating the above physical problem into algebraic terminology.
The ``container'' is replaced by a set $K$ of $n$-dimensional vectors, which for our purposes can be assumed compact and convex (as in all of the figures so far).
The ``force field'' is replaced by a vector-valued function $\mb{f} : \bbr^n \rightarrow \bbr^n$\label{not:F}\label{not:bbr}\label{not:bbrn} which depends on $n$ variables and produces an $n$-dimensional vector as a result.
The geometric idea of a ``right or obtuse angle'' can be expressed using the dot product.
Recalling that the dot product of two vectors can be written as $\mb{a} \cdot \mb{b} = |a||b| \cos \theta$, with $\theta$ the angle between $\mb{a}$ and $\mb{b}$, saying that two vectors make a right or obtuse angle is equivalent to saying that their dot product is negative.
A ``solution'' to the variational inequality is a point which is unmoved by the force field.
(In the above example, we want to say that D and E are solutions to the variational inequality problem created by the container shape and force field, while A, B, and C are not.)
Therefore, rewriting the condition in the previous paragraph with this mathematical notation, we have the following definition:

\begin{dfn}
\label{dfn:videfinition}
\index{convex set!applications}
\index{compact set!applications}
Given a convex set $K \subseteq \bbr^n$\label{not:subseteq} and a function $\mb{f} : K \rightarrow \bbr^n$, we say that the vector $\mb{\hat{x}} \in K$ solves the $\mr{VI}(K,\mb{f})$ if, for all $\mb{x} \in K$, we have
\labeleqn{videfinition}{\mb{f}(\mb{\hat{x}}) \cdot (\mb{x} - \mb{\hat{x}}) \leq 0\,.}\label{not:vikf}\label{not:xhat}
\end{dfn}

In other words, as $\mb{x}$ ranges over all possible points in $K$, $\mb{x} - \mb{\hat{x}}$ represents all possible directions the object can move from $\mb{\hat{x}}$.
The point $\mb{\hat{x}}$ solves the variational inequality precisely when the dot product of the force and all possible directions is negative.

\index{variational inequality!and fixed points}
There is a relationship between fixed point problems and variational inequalities, which can be motivated again by the physical analogy of a force acting within a container.
The solutions to the variational inequality are quite literally ``fixed points'' in the sense that an object placed there will not move.
Consider some point $\mb{x}$ under the action of the force $\mb{f}$.
Assume furthermore that this force is constant and does not change magnitude or direction as this point moves.
Then, since $K$ is convex, the trajectory of the object can be identified with the curve $\mr{proj}_K (\mb{x} + \mu \mb{f}(\mb{x}))$\label{not:proj} where $\mu$\label{not:mu} ranges over all positive numbers and $\mr{proj}_K$ means projection onto the set $K$ as defined in Section~\ref{sec:sets}.
See Figure~\ref{fig:forceapplication2} for a few examples.\index{projection!applications}
If a point is a solution to $\mr{VI}(K,\mb{f})$, then the corresponding ``trajectory'' will simply be the same point no matter what $\mu$ is.
So, for the sake of convenience we arbitrarily choose $\mu = 1$ and look at the location of the point $\mr{proj}_K (\mb{x} + \mb{f}(\mb{x}))$.
If this is the same as the initial point $\mb{x}$, then $\mb{x}$ is a solution to the variational inequality.
So, finally, if we let 
\labeleqn{fixedpointvi}{\mb{f}(\mb{x}) = \mr{proj}_K (\mb{x} + \mb{f}(\mb{x}))\,,}
then the fixed points of $\mb{f}$ coincide exactly with solutions to $\mr{VI}(K,\mb{f})$.

\stevefig{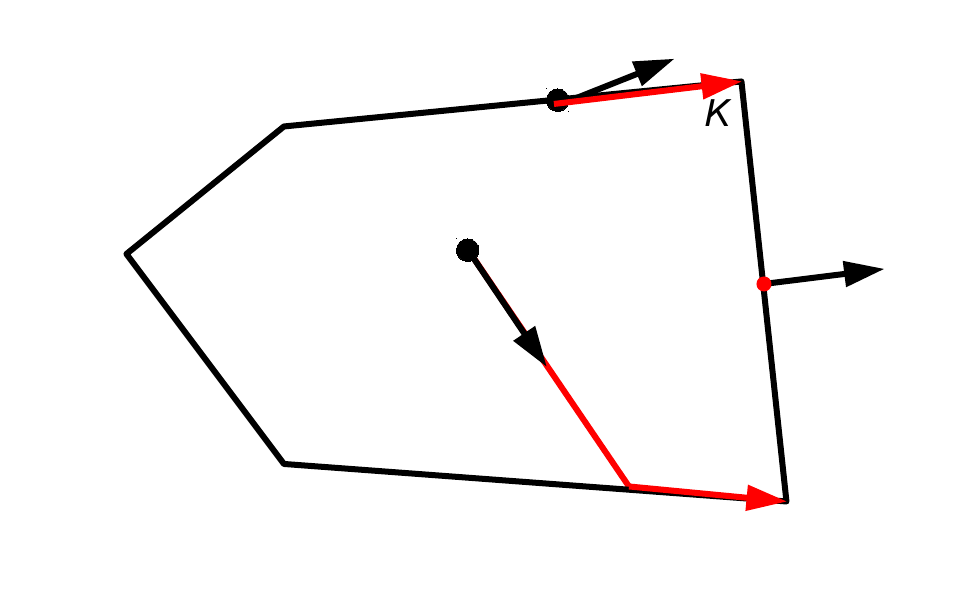}{Trajectories of different points (red) assuming the force (black) is constant.}{0.6\textwidth}

In many cases of practical interest, $\mb{f}$ will be a continuous function.
Furthermore, one can show that the projection mapping onto a convex set (such as $K$) is a well-defined (i.e., single-valued), continuous function.\index{convex set!applications}\index{projection!properties}\index{continuous function!applications}
Then by Proposition~\ref{prp:continuousdifferentiable}, the function $\mb{f}$ defined by equation~\eqn{fixedpointvi} is a continuous function.
So, if the set $K$ is compact in addition to being convex, then Brouwer's theorem shows that the variational inequality must have at least one solution:\index{Brouwer's theorem!applications}

\begin{thm}
\label{thm:viexistence}
\index{variational inequality!existence of solutions}
\index{continuous function!applications}
\index{convex set!applications}
\index{compact set!applications}
If $K \subseteq \bbr^n$ is a compact, convex set and $\mb{f} : \bbr^n \rightarrow \bbr^n$ is a continuous function, then the variational inequality $\mr{VI}(K,\mb{f})$ has at least one solution.
\end{thm}

You should convince yourself that all of these conditions are necessary: if the container $K$ is not bounded or not closed, or if the force field $\mb{f}$ is not continuous, then it is possible that an object placed at any point in the container will move under the action of the force field.
This result also says nothing about solution uniqueness.
Without further conditions, a variational inequality can have multiple solutions, even infinitely many.
\index{variational inequality|)}

\section{Convex Optimization}
\label{sec:convexoptimization}

\index{convex optimization|see {optimization, convex}}
\index{optimization!convex|(}
As we will see in the coming chapters, fixed points and variational inequalities are relatively intuitive ways to represent the idea of an equilibrium in transportation systems.
A fixed point captures the idea that ``at equilibrium, nobody can make a better choice than the one they are currently making; therefore the state of the system is the same from one day to the next.''
The transit example in Section~\ref{sec:fixedpoint} illustrates how such a fixed point problem might arise.
Variational inequalities can represent physical equilibrium problems, involving forces acting on a particle within a confined space.
We haven't yet drawn the connection between this kind of physical equilibrium problem, and the transportation behavior equilibrium problems we are studying in this book, but hopefully such a connection is plausible.
The key will be to define the ``forces'' and the ``container'' in a way that express our behavioral assumptions physically.

Both of these methods have disadvantages.
Fixed point theorems are ``non-constructive,'' which means they often lack methods guaranteed to find a fixed point, even if one exists.
Brouwer's theorem gives us conditions under which a fixed point exists, but tells us nothing about how to find it.
Sometimes applying $f$ repeatedly from a starting point will converge to a fixed point, but not always.
The force field analogy in variational inequalities suggests a natural algorithm (pick a starting point, and see where the force carries you), but again this is not always guaranteed to work.
It is also possible to have multiple solutions to a fixed point or variational inequality problem.
From the standpoint of transportation planning this is inconvenient --- how can you consistently rank alternatives if you have several different predictions for what might happen under each alternative?

Convex optimization is a more powerful tool in that we have uniqueness guarantees on solutions, and efficient algorithms that provably converge to an optimal point.
The downside is that it is not obvious how a user equilibrium problem can be formulated in terms of optimization.
Chapter~\ref{chp:trafficassignmentproblem} takes up this task; this section presents what you need to know about convex optimization for the derivation in that chapter to make sense.
If you have never encountered optimization problems before, please read Appendix~\ref{chp:basicoptimization} before proceeding further, to get familiar with the terminology and notation used in presenting optimization problems.
This subsection will focus on \emph{convex optimization}, a specific kind of optimization problem which is both easier to solve, and well-suited for solving transportation network problems.
If you are interested in other applications, Appendix~\ref{chp:fancyoptimization} discusses methods and properties of other kinds of optimization problems.

We will restrict ourselves to convex optimization here because it is simpler, more powerful, and sufficient to express traffic equilibrium.
Nonconvex optimization is much harder.
For instance, the function in Figure~\ref{fig:nonconvex} has many local minima and is unbounded below as $x \rightarrow -\infty$, both of which can cause serious problems if we're trying to minimize this function.\index{convex function!applications}
Usually, the best that a software program can do is find a local minimum.
If it finds one of the local minima for this function, it may not know if there is a better one somewhere else (or if there is, how to find it).
Or if it starts seeking $x$ values which are negative, we could run into the unbounded part of this function.

\genfig{nonconvex}{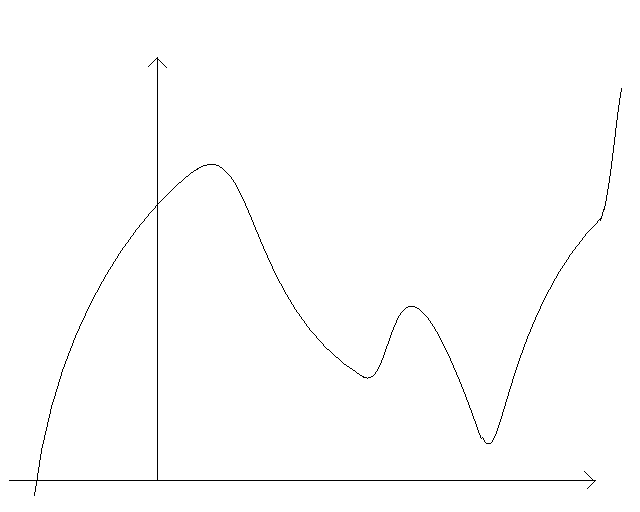}{A function which is not convex.}{width=0.5\textwidth}

On the other hand, some functions are very easy to minimize.
The function in Figure~\ref{fig:convex} only has one minimum point, is not unbounded below, and there are many algorithms which can find that minimum point efficiently.

\genfig{convex}{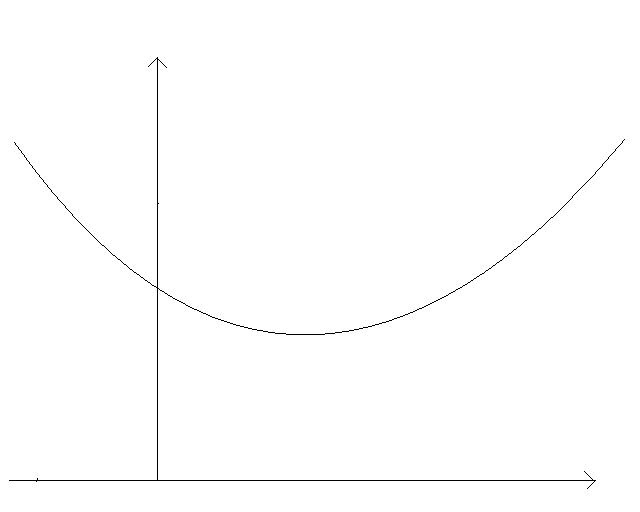}{A convex function.}{width=0.5\textwidth}

What distinguishes these is a property called \emph{convexity}, which is defined in Appendix~\ref{chp:mathbackground}.\index{convex function!and convex optimization}\index{convex set!applications}\index{compact set!applications}
If the feasible region is a convex set, and if the objective function is a convex function, then it is much easier to find the optimal solution.\index{convex set!and convex optimization}
Checking convexity of the objective function is not usually too difficult.
To check convexity of the feasible region, the following result is often useful.

\begin{thm}
Consider an optimization program with $n$ decision variables, whose constraints all take the form $g_i(\mb{x}) \leq 0$\label{not:gcon}\label{not:icon} or $h_j(\mb{x}) = 0$\label{not:hcon}\label{not:jcon}, where $i = 1, 2, \ldots, k$ indexes the inequality constraints and $j = 1, 2, \ldots, \ell$ indexes the equality constraints.
If each function $g_i(\mb{x})$ is convex, and if each function $h_j(\mb{x})$ is linear, then the feasible region $K = \{ \mathbf{x} \in \mathbb{R}^n : g_i(\mathbf{x}) \leq 0, h_j (\mathbf{x}) = 0, i \in \{1, \ldots, k\}, j \in \{1, \ldots, \ell\}\}$ is a convex set.\index{optimization!convex!showing convexity}
\end{thm}
\begin{proof}
Let $Y_i = \{ \mathbf{x} \in \mathbb{R}^n : g_i(\mathbf{x}) \leq 0 \}$\label{not:Yi} represent the values of $\mathbf{x}$ which satisfy the $i$-th inequality constraint, and let $Z_j = \{ \mathbf{x} \in \mathbb{R}^n : h_j(\mathbf{x}) = 0 \}$\label{not:Zj} be the values of $\mathbf{x}$ which satisfy the $j$-th equality constraint.
From Proposition~\ref{prp:levelset}, all of the sets $Y_i$ are convex.
From Example~\ref{exm:hyperplane}, all of the sets $Z_j$ are convex.
The feasible region $K$ is the set of vectors $\mathbf{x}$ which satisfy \emph{all} of the inequality and equality constraints, that is, the intersection of all of the sets $Y_i$ and $Z_j$.
By Proposition~\ref{prp:convexintersection}, therefore, $K$ is convex.
\end{proof}

This is a very common situation, where the functions representing the inequality constraints are convex, and the functions representing equality constraints are linear.
From this theorem, this means that the feasible region must be convex.

This subsection collects a few useful results on convex optimization, the problem of minimizing a convex function over a feasible region.
In convex optimization, every local minimum is a global minimum, every stationary point is a local minimum, and the set of global minima is a convex set.
Furthermore, if the function is strictly convex, the global minimum is unique.
Unlike in elementary calculus, you don't have to perform any ``second derivative tests'' on solutions to ensure they are truly minima, or distinguish between local and global minima.
The set of global minima being a convex set is useful because it means that all solutions are ``connected'' or ``adjacent'' in some sense --- there are no far-flung optimal solutions.

First, a few definitions; in everything that follows, we are trying to minimize a function $f(\mathbf{x})$ over a feasible region $K$.
(That is, $K$ is the set of all $\mathbf{x}$ which satisfy all of the constraints.)

\begin{prp} If $f$ is a convex function and $K$ is a convex set, then every local minimum of $f$ is also a global minimum.\index{optimization!convex!properties}\index{convex set!applications}\index{compact set!applications}
\end{prp}
\begin{proof}
By contradiction, assume that $\mathbf{x_1}$ is a local minimum of $f$, but not a global minimum.
Then there is some $\mathbf{x_2} \in K$ such that $f(\mathbf{x_2}) < f(\mathbf{x_1})$.
Because $f$ is a convex function, for all $\lambda \in (0, 1]$\label{not:lambda} we have
\[f((1 - \lambda) \mathbf{x_1} + \lambda \mathbf{x_2}) \leq (1 - \lambda) f(\mathbf{x_1}) + \lambda f(\mathbf{x_2}) = f(\mathbf{x_1}) + \lambda(f(\mathbf{x_2}) - f(\mathbf{x_1}))\,,\]
and furthermore all points $(1 - \lambda) \mathbf{x_1} + \lambda \mathbf{x_2}$ are feasible since $K$ is a convex set and $\mathbf{x_1}$ and $\mathbf{x_2}$ are feasible.
Since $f(\mathbf{x_2}) < f(\mathbf{x_1})$, this means that
\[f((1 - \lambda) \mathbf{x_1} + \lambda \mathbf{x_2}) < f(\mathbf{x_1})\]
even as $\lambda \rightarrow 0$, contradicting the assumption that $\mathbf{x_1}$ is a local minimum.
\end{proof}

\begin{prp}
If $f$ is a convex function and $K$ is a convex set, then the set of global minima is convex.\index{optimization!convex!properties}\index{convex set!applications}\index{compact set!applications}
\end{prp}
\begin{proof}
Let $\hat{K}$ be the set of global minima of $f$ over the feasible region $K$.
Choose any two global optima $\mathbf{\hat{x}_1} \in \hat{K}$ and $\mathbf{\hat{x}_2} \in \hat{K}$, and any $\lambda \in [0, 1]$.

Since $\mathbf{\hat{x}_1}$ and $\mathbf{\hat{x}_2}$ are global minima, $f(\mathbf{\hat{x}_1}) = f(\mathbf{\hat{x}_2})$; let $\hat{f}$ denote this common value.
Because $K$ is a convex set, the point $(1 - \lambda) \mathbf{\hat{x}_1} + \lambda \mathbf{\hat{x}_2}$ is also feasible.
Because $f$ is a convex function,
\begin{align*}
f((1 - \lambda) \mathbf{\hat{x}_1} + \lambda \mathbf{\hat{x}_2}) &\leq (1 - \lambda) f(\mathbf{\hat{x}_1}) + \lambda f(\mathbf{\hat{x}_2}) \\
                                                                 &= (1 - \lambda) \hat{f} + \lambda \hat{f} \\
                                                                 &= \hat{f}
\,.
\end{align*}
Therefore $f((1 - \lambda) \mathbf{\hat{x}_1} + \lambda \mathbf{\hat{x}_2}) \leq \hat{f}$
But at the same time, $f((1 - \lambda) \mathbf{\hat{x}_1} + \lambda \mathbf{\hat{x}_2}) \geq \hat{f}$ because $\hat{f}$ is the global minimum value of $f$.
So we must have $f((1 - \lambda) \mathbf{\hat{x}_1} + \lambda \mathbf{\hat{x}_2}) = \hat{f}$
which means that this point is also a global minimum and $(1 - \lambda) \mathbf{\hat{x}_1} + \lambda \mathbf{\hat{x}_2} \in \hat{K}$ as well, proving its convexity.
\end{proof}

\begin{prp}
\label{prp:strictconvexuniqueness}
(Uniqueness.)  If $f$ is a strictly convex function and $K$ is a convex set, the set of optimal solutions $\hat{K}$ has at most one element.\index{optimization!convex!uniqueness of solutions}\index{convex function!strict convexity!applications}\index{strictly convex function|see {convex function, strict convexity}}
\end{prp}
\begin{proof}
By contradiction, assume that $\hat{K}$ contains two distinct elements $\mathbf{\hat{x}_1}$ and $\mathbf{\hat{x}_2}$.
Repeating the proof of the previous proposition, because $f$ is strictly convex, the first inequality becomes strict and we must have $f((1 - \lambda) \mathbf{\hat{x}_1} + \lambda \mathbf{\hat{x}_2}) < \hat{f}$.
This contradicts the assumption that $\mathbf{\hat{x}_1}$ and $\mathbf{\hat{x}_2}$ are global minima.
Therefore, there is at most one optimal solution.
\end{proof}

Combining this with Weierstrass' theorem (Theorem~\ref{thm:weierstrass}), which gives conditions guaranteeing existence of an equilibrium solution, we have the following result:\index{Weierstrass' theorem!applications}

\begin{prp}
\label{prp:strictconvexexistunique}
\index{optimization!convex!uniqueness of solutions}
\index{optimization!convex!existence of solutions}
\index{continuous function!applications}
\index{convex set!applications}
\index{convex function!applications}
\index{convex function!strict convexity!applications}
\index{compact set!applications}
(Existence and uniqueness.)  If $f$ is a continuous, strictly convex function, and $K$ is a non-empty, compact, and convex set, there is exactly one optimal solution.
\end{prp}

These additional conditions are not onerous.
In practice, most objective functions are continuous; and in any case, one can show that a convex function \emph{must} be continuous, except possibly on its boundary.\index{continuous function!applications}
Requiring that the feasible region be non-empty and compact in addition to convex is not very limiting, either; in practice, your range of options is generally neither empty nor unbounded.\index{convex set!applications}\index{compact set!applications}

We next discuss techniques which are common in transportation network analysis.
We begin with simple problems: optimization problems with a single decision variable, and optimization problems with no constraints, before moving to more general problem classes.
Throughout this section, we will assume that the objective is a convex, differentiable function, and that the feasible region is a convex set.\index{convex set!applications}\index{continuous function!applications}
In each case, we will identify \emph{optimality conditions}: simple equations and inequalities that characterize optimality in the sense that any optimal solution must satisfy all of the optimality conditions, and that any solution satisfying all of the optimality conditions must be optimal.\index{optimality conditions|see {optimization, convex, optimality conditions}}\index{optimization!convex!optimality conditions|(}
For small problems, we can solve the optimality conditions to directly find an optimal solution.
For larger-scale problems, these conditions are more useful as a ``certificate'' of optimality, to know whether a solution we have found in another way is optimal, or close to optimal.

\subsection{Single-variable problems}
\label{sec:onedimension}

To start off, consider a simple minimization problem in one variable with no constraints:

\begin{equation*}
\begin{array}{rrl}
\displaystyle  \min_x \quad & f(x)     &     \,.   \\
\end{array}
\end{equation*}

Because $f$ is convex, any local minimum is also a global minimum.
So, all we need is to know when we've reached a local minimum.
For the unconstrained case, this is easy: we know from basic calculus that $\hat{x}$ is a local minimum if\label{not:fprime}
\[f'(\hat{x}) = 0\,.\]
We don't have to check whether $\hat{x}$ is a local minimum or a local maximum because $f$ is convex.

\begin{exm} Find the value of $x$ which minimizes $f(x) = x^2 - 3x + 5$. \label{exm:unconstrained1d} \end{exm}
\textbf{Solution}.
$f'(x) = 2x - 3$, which vanishes if $x = 3/2$.
Therefore $x = 3/2$ minimizes $f(x)$. $\blacksquare$

It's a little bit more complicated if we add constraints to the picture.
For instance, consider the function in Figure~\ref{fig:convex} (which could very well be the function from Example~\ref{exm:unconstrained1d}), but with the added constraint $x \geq 0$.
In this case, nothing is different, and the optimum still occurs where $f'(x)$ vanishes, that is, at $x = 3/2$.
But what if the constraint was $x \geq 2$?
In this case, $x = 3/2$ is infeasible, and $f'(x)$ is always strictly positive in the entire feasible region.
This means that $f$ is strictly increasing over the entire feasible region, so the minimum value is obtained at the smallest possible value of $x$, that is, $x = 2$.
So we see that sometimes the local minimum of a constrained optimization problem can be at a point where $f'(x)$ is nonzero.

To simplify things a little bit, assume that the constraint is of the form $x \geq 0$, that is, we are trying to solve
\begin{equation*}
\begin{array}{rrl}
\displaystyle  \min_x \quad & f(x)     &        \\
\mbox{s.t.}    \quad        & x  & \geq 0 \,. \\
\end{array}
\end{equation*}

\begin{figure}
   \centering
   \subfigure[Convex function where the constraint is not binding at the minimum.]{
   \includegraphics[height=0.4\textheight]{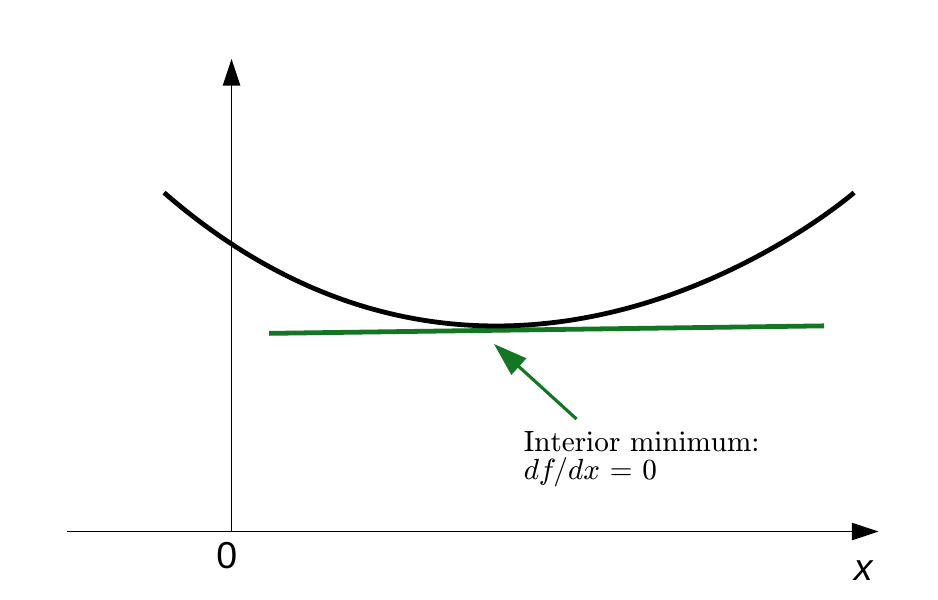}
   \label{fig:simpleconstraint_interior}
   }
   \subfigure[Convex function where the constraint is binding at the minimum.]{
   \includegraphics[height=0.4\textheight]{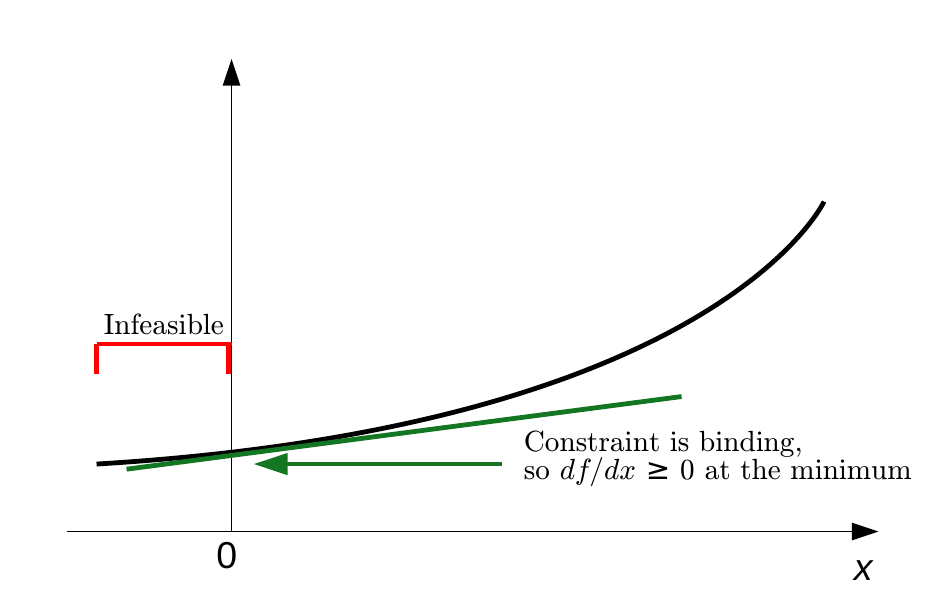}
   \label{fig:simpleconstraint_boundary}
   }
   \caption{Two possibilities for minimizing a convex function with a constraint. \label{fig:binding_nonbinding}}
\end{figure}

As Figure~\ref{fig:binding_nonbinding} shows, there are only two possibilities.
In the first case, the minimum occurs when $x$ is strictly positive.
We can call this an \emph{interior} minimum\index{interior optimum}, or we can say that the constraint $x \geq 0$ is \emph{nonbinding}\index{nonbinding constraint} at this point.
In this case, clearly $f'(x)$ must equal zero: otherwise, we could move slightly in one direction or the other, and reduce $f$ further.
The other alternative is that the minimum occurs for $x = 0$, as in Figure~\ref{fig:simpleconstraint_boundary}.
For this to be a minimum, we need $f'(0) \geq 0$ --- if $f'(0) < 0$, $f$ is decreasing at $x = 0$, so we could move to a slightly positive $x$, and thereby reduce $f$.

Let's try to draw some general conclusions.
For the interior case of Figure~\ref{fig:simpleconstraint_interior}, we needed $x \geq 0$ for feasibility, and $f'(x) = 0$ for optimality.
For the boundary case of Figure~\ref{fig:simpleconstraint_boundary}, we had $x = 0$ exactly, and $f'(x) \geq 0$.
So we see that in both cases, $x \geq 0$ and $f'(x) \geq 0$, and furthermore that at least one of these has to be exactly equal to zero.\index{extreme-point optimum}\index{binding constraint}
To express the fact that either $x$ or $f'(x)$ must be zero, we can write $x f'(x) = 0$.
So a solution $\hat{x}$ solves the minimization problem if and only if
\begin{eqnarray*}
   \hat{x} & \geq & 0 \\
   f'(\hat{x}) & \geq & 0 \\
   \hat{x} f'(\hat{x}) & = & 0 \,.
\end{eqnarray*}
Whenever we can find $\hat{x}$ that satisfies these three conditions, we know it is optimal.
These are often called \emph{first-order conditions}\index{first-order condition} because they are related to the first derivative of $f$.
The condition $\hat{x} f(\hat{x}) = 0$ is an example of a \emph{complementarity constraint}\index{complementarity constraint} because it forces either $\hat{x}$ or $f(\hat{x})$ to be zero.

\subsection{Bisection method}
\label{sec:bisection}

\index{optimization!line search!bisection|(}
\index{line search|see {optimization, line search}}
\index{bisection|see {optimization, line search, bisection}}
The bisection method allows us to solve one-dimensional problems over bounded feasible regions.
Consider the one-dimensional optimization problem\label{not:a}\label{not:b}
\begin{align*}
   \min_x      \quad & f(x) \\
   \mbox{s.t.} \quad & a \leq x \leq b \,, 
\end{align*}
where $f$ is continuously differentiable and convex.\index{continuous function!applications}\index{convex function!applications}
The bisection method works by constantly narrowing down the region where the optimal solution lies.
  After the $k$-th iteration\label{not:k}, the bisection method will tell you that the optimum solution lies in the interval $[a_k, b_k]$, with this interval shrinking over time (that is, $b_k - a_k < b_{k-1} - a_{k-1}$).
A natural termination criterion is to stop when the interval is sufficiently small, that is, when $b_k - a_k < \epsilon$, where $\epsilon$ is the precision you want for the final solution. 

The idea is that the sign of the derivative of the midpoint tells you where the optimum is.
If the derivative is negative at $a_k$, but positive at $b_k$, the optimum occurs somewhere in-between, at a point where it is zero.
So, if the derivative is positive at the midpoint, we know that the zero point has to happen somewhere to its left; if negative, somewhere to its right.
If we happen to get lucky, the derivative at the midpoint will be exactly zero, and you can stop --- but this is really rare.
Figure~\ref{fig:bisection} illustrates how bisection works.
(The exercises ask you to show that bisection converges to the optimum solution even if the derivative initially has the same sign at both endpoints; the case where $f'$ is positive at $a_k$ but negative at $b_k$ is impossible for a convex function.)

\genfig{bisection}{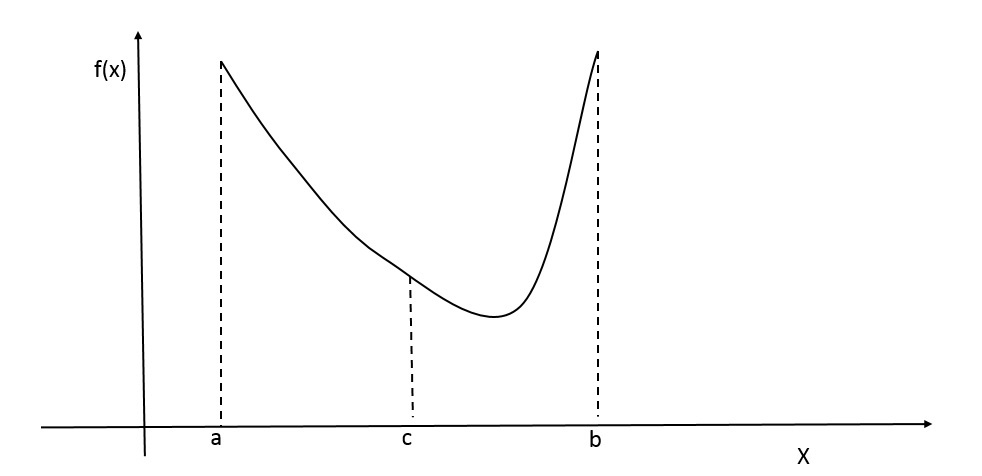}{Bisection method.}{width=\linewidth}

Here's how the algorithm works.
\begin{description}
   \item[Step 0: Initialize.]  Set the iteration counter $k = 0$, $a_0 = a$, $b_0 = b$.
   \item[Step 1: Evaluate midpoint.]  Identify the midpoint $c_k = (a_k + b_k)/2$ and calculate the derivative of $f$ there, $d_k = f'(c_k)$
   \item[Step 2: Bisect.]  If $d_k > 0$, set $a_{k+1} = a_k$, $b_{k+1} = c_k$.
Otherwise, set $a_{k+1} = c_k$, $b_{k+1} = b_k$.
   \item[Step 3: Iterate.]  Increase the counter $k$ by 1 and check the termination criterion.
If $b_k - a_k < \epsilon$, then terminate; otherwise, return to step 1.
\end{description}

\begin{exm}
Use the bisection algorithm to find the minimum of $f(x) = (x - 1)^2 + e^x$ on the interval $x \in [0, 2]$, to within a tolerance of $\epsilon = 0.01$.
\end{exm}
\solution{
You may find it useful to follow along in Table~\ref{tbl:bisection}.
We start off with $k = 0$, $a_0 = 0$, and $b_0 = 2$.
We calculate the derivative at the midpoint: $f'(x) = 2(x - 1) + e^x$, so $f'(1) = 2.71$, which is positive.
Since $f'$ is positive at $x = 1$, the minimum must occur to the \emph{left} of this point, that is, somewhere in the interval $[0, 1]$.
We set $a_1$ and $b_1$ equal to these new values, and repeat.
The new midpoint is $1/2$, and $f'(1/2) = 0.649$ is again positive, so the minimum must occur to the left of this point, in the interval $[0, 1/2]$.
The midpoint of $0$ and $1/2$ is $1/4$, and $f'(1/4) = -0.216$ is negative, so the minimum is to the \emph{right} of the midpoint.
Thus, the new interval is $[1/4, 1/2]$, and we repeat as shown in Table~\ref{tbl:bisection}.

During the eighth iteration, the width of the interval is $41/128 - 5/16 = 0.00781$, which is less than the tolerance of 0.01; therefore we stop, and return our best guess of the optimum as the midpoint of this interval: $\hat{x} \approx 81/256 = 0.31640625$.
The true minimum point occurs at $\hat{x} = 0.3149230578\ldots$; if we had chosen a smaller tolerance $\epsilon$, the algorithm would have narrowed the interval further, with both ends converging towards this point.
}

\begin{table}
\centering
\caption{Demonstration of the bisection algorithm with $f(x) = (x-1)^2 + e^x$, $x \in [0,2]$.
\label{tbl:bisection}}
\begin{tabular}{|c|cc|cc|}
  \hline
  $k$   &  $a_k$    &  $b_k$    & $c_k = (a_k + b_k)/2$ &  $d_k$       \\
  \hline
  0     &  0        &  2        &  1                    &  $2.71 > 0$     \\
  1     &  0        &  1        &  $1/2$                &  $0.649  > 0$   \\
  2     &  0        &  $1/2$    &  $1/4$                &  $-0.216 < 0$    \\
  3     &  $1/4$    &  $1/2$    &  $3/8$                &  $0.205  > 0$   \\
  4     &  $1/4$    &  $3/8$    &  $5/16$               &  $-0.00816 < 0$ \\
  5     &  $5/16$   &  $3/8$    &  $11/32$              &  $0.0978 > 0$   \\
  6     &  $5/16$   &  $11/32$  &  $21/64$              &  $0.0446 > 0$   \\
  7     &  $5/16$   &  $21/64$  &  $41/128$             &  $0.0182 > 0$   \\
  8     &  $5/16$   &  $41/128$ &  $81/256$             &  \\
  \hline
\end{tabular}
\end{table}

There are additional methods for solving one-dimensional convex optimization problems like these.
Appendix~\ref{chp:fancyoptimization} also describes the ``golden section'' method, which is useful when the objective function is not differentiable, and a method based on Newton's method, which is useful when the objective is twice differentiable.\index{optimization!line search!bisection|)}

\subsection{Multiple decision variables}
\label{sec:unconstrained}

Most interesting optimization problems have more than one decision variable.
However, we will keep the assumption that the only constraint on the decision variables is nonnegativity.
Using the vector $\mathbf{x}$ to refer to all of the decision variables, we solve the problem
\begin{equation*}
\begin{array}{rrl}
\displaystyle  \min_\mathbf{x} \quad & f(\mathbf{x})     &        \\
       \mbox{s.t.}    \quad        & \mathbf{x}  & \geq \bm{0} \,. \\     
\end{array}
\end{equation*}

Using the same logic as before, we can show that $\mathbf{\hat{x}}$ solves this problem if and only if the following conditions are satisfied for \emph{every} decision variable $x_i$:
\begin{eqnarray*}
\hat{x}_i & \geq & 0 \\
\pdr{f(\hat{x})}{x_i} & \geq & 0 \\
\hat{x}_i \pdr{f(\hat{x})}{x_i} & = & 0 \,. \\
\end{eqnarray*}
You should convince yourself that if these conditions are not met for each decision variable, then $\mathbf{\hat{x}}$ cannot be optimal: if the first condition is violated, the solution is infeasible; if the second is violated, the objective function can be reduced by increasing $x_i$; if the first two are satisfied but the third is violated, then $\hat{x}_i > 0$ and $\pdr{f(\hat{x})}{x_i} > 0$, and the objective function can be reduced by decreasing $x_i$.

This can be compactly written in vector form as
\labeleqn{kktcomplementarity}{\mathbf{0} \leq \mathbf{\hat{x}} \perp \nabla f(\mathbf{\hat{x}})  \geq \mathbf{0} \,,}
where the $\perp$\label{not:perp} symbol indicates orthogonality, i.e., that the dot product of $\mathbf{\hat{x}}$ and the gradient vector $\nabla f(\mathbf{\hat{x}})$\label{not:nabla} is zero\label{not:zero}.

Unfortunately, the bisection algorithm does not work nearly as well in higher dimensions.
It is difficult to formulate an extension that always works, and those that do are inefficient.
We'll approach solution methods for higher-dimensional problems somewhat indirectly, tackling a few other topics first: addressing constraints other than nonnegativity, and a few highlights of linear optimization.

\subsection{Constrained problems} 
\label{sec:constrained}

Constraints can take many forms.
For the purposes of this book, we can restrict attention to only two types of constraints: \emph{linear equality}\index{optimization!constraints!linear equality} constraints and \emph{nonnegativity} constraints.\index{optimization!constraints!nonnegativity}
The previous two sections showed you how to deal with nonnegativity constraints; this section discusses linear equality constraints.
A linear equality constraint is of the form
\labeleqn{linearequality}{\sum_i a_i x_i = b\,,}
where the $x_i$ are decision variables, and the $a_i$ and $b$ are constants.

\index{Lagrange multipliers|(}
We can handle these using the technique of Lagrange multipliers.
This technique is demonstrated in the following example for the case of a single linear equality constraint:
\begin{equation*}
\begin{array}{rrl}
  \displaystyle  \min_{x_1, x_2} \quad & x_1^2 + x_2^2     &        \\
           \mbox{s.t.}    \quad        & x_1 + x_2   & = 5 \,. \\    
\end{array}
\end{equation*}
(It is a useful exercise to verify that $x_1^2 + x_2^2$ is a strictly convex function, and that $\{ x_1, x_2 : x_1 + x_2 = 5\}$ is a convex set.)

The main idea behind Lagrange multipliers is that unconstrained problems are easier than constrained problems.
The technique is an ingenious way of nominally removing a constraint while still ensuring that it holds at optimality.
The equality constraint is ``brought into the objective function'' by multiplying the difference between the right- and left-hand sides by a new decision variable $\kappa$\label{not:kappa} (called the Lagrange multiplier), adding the original objective function.
This creates the Lagrangian function\index{Lagrangian function}\label{not:fancyL}
\labeleqn{lagrangian}{\mathcal{L}(x_1, x_2, \kappa) = x_1^2 + x_2^2 + \kappa(5 - x_1 - x_2) \,.}
It is possible to show that \emph{the optimal solutions of the original optimization problem correspond to stationary points of the Lagrangian function}, that is, to values of $x_1$, $x_2$, and $\kappa$ such that $\nabla \mathcal{L} (x_1, x_2, \kappa) = \mathbf{0}$.
To find this stationary point, take partial derivatives with respect to each variable and set them all equal to zero:
\begin{align}
   \pdr{\mathcal{L}}{x_1} = 2x_1 - \kappa &= 0 \label{eqn:lagrange1} \\
   \pdr{\mathcal{L}}{x_2} = 2x_2 - \kappa &= 0 \label{eqn:lagrange2} \\
   \pdr{\mathcal{L}}{\kappa} = 5 - x_1 - x_2 &= 0 \label{eqn:lagrange3} \,.
\end{align}
Notice that the third optimality condition~\eqn{lagrange3} \emph{is simply the original constraint}, so this stationary point must be feasible.
Equations~\eqn{lagrange1} and~\eqn{lagrange2} respectively tell us that $x_1 = \kappa / 2$ and $x_2 = \kappa / 2$; substituting these expressions into~\eqn{lagrange3} gives $\kappa = 5$, and therefore the optimal solution occurs for $x_1 = x_2 = 5/2$.

This technique generalizes perfectly well to the case of multiple linear equality constraints.
Consider the general optimization problem
\begin{equation*}
   \begin{array}{rrl}
      \displaystyle  \min_{x_1, \ldots, x_n} \quad & f(x_1, \ldots, x_n)      &        \\
               \mbox{s.t.}    \quad        & \sum_{i=1}^n a_{1i} x_i & = b_1 \\
                                           & \sum_{i=1}^n a_{2i} x_i & = b_2 \\
                                           & \vdots \qquad & \\
                                           & \sum_{i=1}^n a_{mi} x_i & = b_m \,,\\
   \end{array}
\end{equation*}
where $f$ is convex.\label{not:aij}\label{not:bj}
The corresponding Lagrangian is
\begin{multline*}\mathcal{L}(x_1, \ldots, x_n, \kappa_1, \ldots, \kappa_m)
= 
f(x_1, \ldots, x_n) + \kappa_1\myp{b_1 - \sum_{j=1}^n a_{1j} x_j} + \\ \kappa_2\myp{b_2 - \sum_{j=1}^n a_{2j} x_j} + \cdots + \kappa_m \myp{b_m - \sum_{j=1}^n a_{mj} x_j}
\,.
\end{multline*}

For an optimization problem that has \emph{both} linear equality constraints and nonnegativity constraints, we form the optimality conditions by combining the Lagrange multiplier technique with the complementarity technique from the previous section.
Thinking back to Section~\ref{sec:unconstrained}, in the same way that we replaced the condition $f'(\hat{x}) = 0$ for the unconstrained case with the three conditions $\hat{x} \geq 0$, $f'(\hat{x}) \geq 0$, and $x f'(\hat{x}) = 0$ when the nonnegativity constraint was added, we'll adapt the Lagrangian optimality conditions.
If the optimization problem has the form
\begin{equation*}
   \begin{array}{rrl}
      \displaystyle  \min_{x_1, \ldots, x_n} \quad & f(x_1, \ldots, x_n)      &        \\
               \mbox{s.t.}    \quad        & \sum_{i=1}^n a_{1i} x_i & = b_1 \\
                                           & \sum_{i=1}^n a_{2i} x_i & = b_2 \\
                                           & \vdots \qquad & \\
                                           & \sum_{i=1}^n a_{mi} x_i & = b_m \\
                                           & x_1, \ldots, x_n & \geq 0 \,, \\
   \end{array}
\end{equation*}
where $m \leq n$, then the Lagrangian is
\begin{multline*}\mathcal{L}(x_1, \ldots, x_n, \kappa_1, \ldots, \kappa_m)
= 
f(x_1, \ldots, x_n) + \kappa_1\myp{b_1 - \sum_{j=1}^n a_{1j} x_j} + \\ \kappa_2\myp{b_2 - \sum_{j=1}^n a_{2j} x_j} + \cdots + \kappa_m \myp{b_m - \sum_{j=1}^n a_{mj} x_j}
\,,
\end{multline*}
and the optimality conditions are
\begin{align*}
   \pdr{\mathcal{L}}{x_i} &\geq 0 & \forall i \in \{ 1, \ldots, n \} \\
   \pdr{\mathcal{L}}{\kappa_j} &= 0 & \forall j \in \{1, \ldots, m \} \\
   x_i & \geq 0 & \forall i \in \{1, \ldots, n \} \\
   x_i\pdr{\mathcal{L}}{x_i} & = 0 & \forall i \in \{1, \ldots, n \} \,.
\end{align*}
Be sure you understand what each of these formulas implies.
Each decision variable must be nonnegative, the partial derivative of $\mathcal{L}$ with respect to this variable must be nonnegative, and their product must equal zero (for the same reasons as discussed in Section~\ref{sec:unconstrained}).
For the Lagrange multipliers $(\kappa_1, \ldots, \kappa_m)$, the corresponding partial derivative of $\mathcal{L}$ must be zero.
Notice how this is a combination of the two techniques.

For small optimization problems, we can write down each of these conditions and solve for the optimal solution, as above.
However, for large-scale problems this process can be very inconvenient.
Later chapters in the book explain methods which work better for large problems.
As a final note, the full theory of Lagrange multipliers is more involved than what is discussed here.
However, it suffices for the case of a convex objective function and linear equality constraints.
Optimality conditions for some other cases are given in Appendix~\ref{chp:fancyoptimization}.\index{Lagrange multipliers|)}\index{optimization!convex!optimality conditions|)}\index{optimization!convex|)}

\section{Historical Notes and Further Reading}

\cite{pant02} present a history of fixed point problems and major results.
The main result on fixed point problems (Brouwer's theorem)\index{Brouwer's theorem} was presented in \cite{brouwer10}, although this and similar results were anticipated by Cauchy\index{Cauchy, Augustin-Louis} and Poincar\'{e}.\index{Poincar\'{e}, Henri}

Variational inequalities were first formulated in mechanics, to model the physics of elastic bodies deforming under their own weight.
Early summaries of this theory can be found in \cite{duvaut71} and \cite{glowinski76}.
A comprehensive and more contemporary treatment is given in the two volumes of \cite{facchinei03}.
For extended treatments of convex optimization, see textbooks of \cite{rockafellar97} and \cite{bertsekas03}.

Mathematical optimization has a long history, including important contributions by Fermat\index{Fermat, Pierre de}, Newton\index{Newton, Isaac}, the Bernoullis\index{Bernoulli, Johann}\index{Bernoulli, Jacob}, Lagrange\index{Lagrange, Joseph-Louis}, and Gauss.\index{Gauss, Carl Friedrich}
These earlier methods are largely based in calculus and analytical in nature (with Newton's method being a notable exception).
With the advent of computers in the early 20th century, and the logistics demands imposed by World War II,\index{World War II} the field of optimization took on additional foci centered on computation and solution of large-scale problems.
The seminal work of \cite{dantzig63}\index{Dantzig, George} in linear programming substantially expanded both optimization theory and the range of applications where optimization was used.

Optimization problems are often classified based on the structure of the objective function, decision variables, and constraints.
Contemporary treatments can be found in linear programming~\citep{bertsimas97}, nonlinear programming~\citep{bazaara_nlp,bertsekas_nlp}, integer programming~\citep{wolsey_ip}, stochastic optimization~\citep{birge97}, and network optimization~\citep{ahuja93}.
For further reading on convex optimization, \cite{rockafellar97} and \cite{bertsekas03} provide thorough reviews of properties of convex functions and sets, of their analysis, and of optimization problems defined on them.

\section{Exercises}
\label{sec:mathematicalpreliminaries_exercises}

\begin{enumerate}
   \item \diff{31} For each of the following functions, find all of its fixed points or state that none exist.
   \begin{enumerate}[(a)]
      \item $f(x) = x^2$, where $K = \bbr$.
      \item $f(x) = 1 - x^2$, where $K = [0, 1]$.
      \item $f(x) = e^x$, where $K = \bbr$.
      \item $f(x_1, x_2) = \vect{-x_1 \\ x_2}$ where $K = \{ (x_1, x_2) : -1 \leq x_1 \leq 1, 0 \leq x_2 \leq 1 \} $.
      \item $f(x,y) = \vect{-y \\ x} $ where $K = \{ (x,y) : x^2 + y^2 \leq 1 \}$ is the unit disc.
   \end{enumerate}
   \item \diff{54} Brouwer's theorem guarantees the existence of a fixed point for the function $f: K \rightarrow K$\label{not:fXY} if $f$ is continuous and $K$ is closed, bounded, and convex.
Show that each of these four conditions is necessary by creating examples of a function $f$ and set $K$ which satisfy only three of those conditions but do not have a fixed point.
Come up with such examples with each of the four conditions missing.
(The notation $f: K \rightarrow K$ means that the range of the function must be contained in its domain; every ``output'' from $f$ is also a valid ``input.'')
Hint: you will probably find it easiest to work with simple functions and sets whenever possible, e.g., something like $K = [0,1]$.
Visualizing fixed points as intersections with the diagonal line through the origin may help you as well.
   \item \diff{45} Find all of the solutions of each of the following variational inequalities $\mr{VI}(K, \mb{f})$.
   \begin{enumerate}[(a)]
      \item $F(x) = x + 1$, $K = [0, 1]$
      \item $F(x) = \frac{x}{2}$, $K = [-1, 1]$
      \item $F(x,y) = \vect{0 \\ -y}$, $K = \{(x,y) : x^2 + y^2 \leq 1\}$
   \end{enumerate}
   \item \diff{58} Theorem~\ref{thm:viexistence} guarantees the existence of a solution to the variational inequality $\mr{VI}(K,\mb{f})$ if $K$ is closed, bounded, and convex, and if $\mb{f}$ is continuous.
Show that each of these three conditions is necessary by creating counterexamples of functions $\mb{f}$ and sets $K$ which satisfy only three of these conditions, but for which $\mr{VI}(K,\mb{f})$ has no solution.
It may be helpful to include sketches.
\end{enumerate}

\part{Static Traffic Assignment}
\label{part:statictrafficassignment}

\chapter{Introduction to Static Assignment}
\label{chp:equilibrium}

\index{static traffic assignment|(}
This chapter introduces the static traffic assignment problem, laying the groundwork for more extended discussions in later chapters in this part.
It defines the main concepts and notation to be used throughout this part of the book (Section~\ref{sec:staticconcepts}) and explains how the equilibrium principle applies in static traffic assignment (Section~\ref{sec:principleofuserequilibrium}).
In this chapter we present the equilibrium idea in the simplest possible way.
Chapter~\ref{chp:trafficassignmentproblem} will formulate these ideas in terms of the mathematical concepts introduced in Chapter~\ref{chp:mathematicalpreliminaries}, an approach which will scale better for problems of a more realistic scale.
Section~\ref{sec:principleofuserequilibrium} also gives a basic way to solve for equilibrium using trial-and-error.
Chapter~\ref{chp:solutionalgorithms} will present approaches better-suited for real world problems.

As additional motivation, Section~\ref{sec:motivatingexamples} presents three ``paradoxes'' where the equilibrium solution behaves in counterintuitive ways.
Adding additional capacity to a network, or retiming signals with the aim of reducing delay, may actually \emph{increase} delay if proper care is not taken to anticipate how users of the system will react to these changes.
Finally, Section~\ref{sec:tapdata} discusses what data is needed to construct a static traffic assignment model.

\section{Notation and Concepts}
\label{sec:staticconcepts}

As discussed in Chapter~\ref{chp:introchapter}, there are many possible measures of effectiveness for evaluating the impacts of a transportation project or policy.
Recall that the traffic assignment problem assumes that the number of drivers traveling between each origin zone and destination zone is given, and we want to find the number of drivers using each roadway segment.
These link flow values provide information on congestion, emissions, toll revenue, or other measures of interest.
We are also given the underlying network, which has a set of links\index{link} $A$ representing the roadway infrastructure, and a set of nodes\index{node} $N$ representing junctions, and a set of centroids $Z$\label{not:Z} where trips are allowed to begin and end.
Every centroid\index{node!centroid} is represented by a node, so $Z$ is a subset of $N$.

Let $d^{rs}$\label{not:drs} denote the number of drivers whose trips start at centroid $r$ and end at centroid $s$.
It is often convenient to write these values in matrix form as an \emph{origin-destination (OD) matrix} where the number of trips between $r$ and $s$ is given by the value in the $r$-th row and $s$-th column.
Together, an origin $r$ and destination $s$ form an \emph{OD pair} $(r,s)$.\index{OD pair}\label{not:rs}

Associated with each link is its \emph{flow}\index{link!flow}, denoted $x_{ij}$, representing the total number of vehicles wanting to use link $(i,j)$ during the analysis period.
The flow is also known as the \emph{volume} or \emph{demand}.
For reasons described below, ``demand'' is actually the most accurate term, but ``flow'' and ``volume'' are the most common for reasons of tradition, and ``flow'' is used in this book.
The travel time on link $(i,j)$\index{link!travel time} is expressed as $t_{ij}$, and to represent congestion, we let this travel time be a function of $x_{ij}$ and write $t_{ij}(x_{ij})$.\index{link performance function|(}
Because of congestion effects, $t_{ij}$ is typically increasing and convex, that is, its first two derivatives are typically positive.
The function used to relate flow to travel time is called a \emph{link performance function}.
The most common used in practice is the Bureau of Public Roads (BPR) function,\index{link performance function!BPR (Bureau of Public Roads)} named after the agency which developed it:
\labeleqn{bpr}{t_{ij}(x_{ij}) = t^0_{ij} \left(1 + \alpha \left( \frac{x_{ij}}{u_{ij}} \right) ^\beta \right)\,,}
where $t^0_{ij}$ is the ``free-flow'' travel time (the travel time with no congestion), $u_{ij}$ is the practical capacity (typically the value of flow which results in a level of service of C or D), and $\alpha$ and $\beta$ are shape parameters which can be calibrated to data.
The values $\alpha = 0.15$ and $\beta = 4$ are commonly used if no calibration is done, but see Section~\ref{sec:tapdata} for more discussion on how they should be chosen.

Notice that this function is well-defined for any value of $x_{ij}$, even when flows exceed the stated ``capacity'' $u_{ij}$.
In the basic traffic assignment problem, there are no explicit upper bounds enforced on link flows.\index{link!capacity}
The interpretation of a ``flow'' greater than the capacity is actually that the \emph{demand} for travel on the link exceeds the capacity, and queues will form.
The delay induced by these queues is then reflected in the link performance function.
Alternately, one can choose a link performance function which asymptotically grows to infinity as $x_{ij} \rightarrow u_{ij}$, implicitly enforcing the capacity constraint.
However, this approach can introduce numerical issues in solution methods.
More discussion on this issue follows in Section~\ref{sec:commentary}, but the short answer is that properly addressing capacity constraints in traffic assignment requires a \emph{dynamic} traffic assignment model, which is the subject of Part~\ref{part:dynamictrafficassignment}.\index{link!demand}\index{link performance function|)}

\stevefig{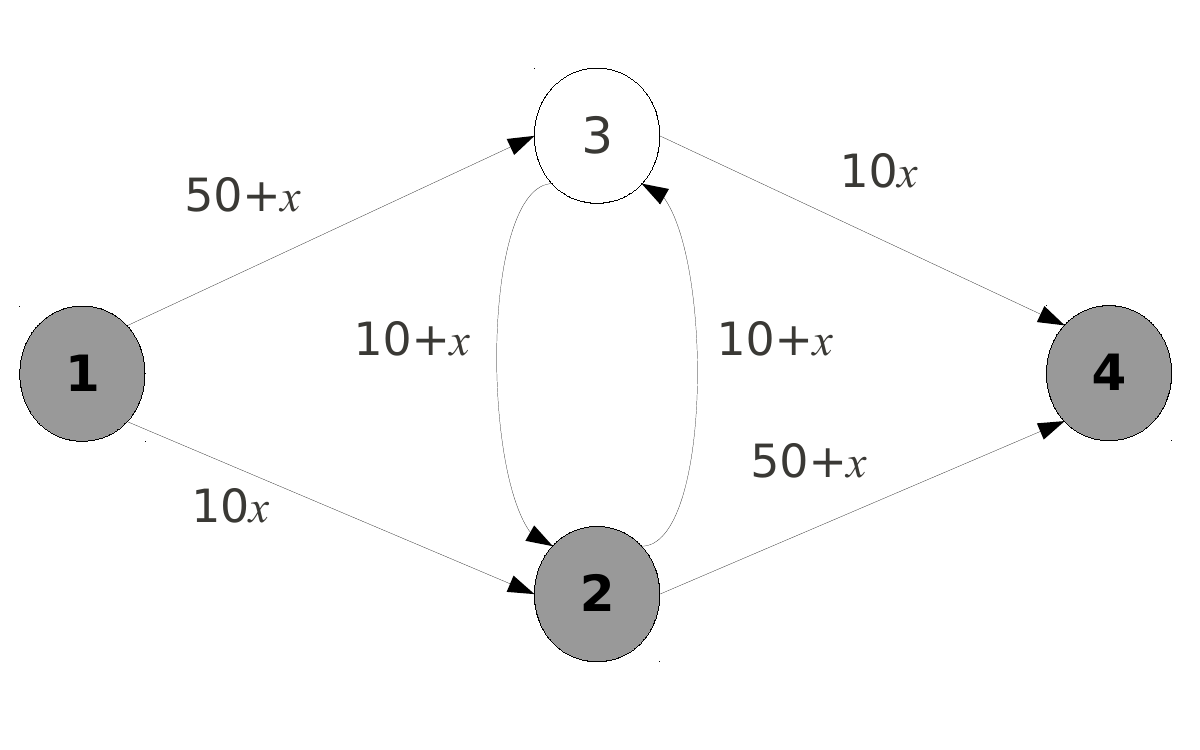}{Example network for demonstration, with link performance functions shown.}{0.8\textwidth}

From the modeler's perspective, the goal of traffic assignment is to determine the link flows in a network.
But from the standpoint of the travelers themselves, it is easier to think of them each choosing a path connecting their origin to their destination.
Let $d^{rs}$ denote the number of vehicles which will be departing origin $r$ to destination $s$, where $r$ and $s$ are both centroids.
Using $h^\pi$\label{not:hpi} to represent the number of vehicles who choose path $\pi$\index{path!flow}, a \emph{feasible assignment}\index{feasible assignment} is defined as a vector of path flows $\mb{h}$ satisfying the following conditions:
\begin{enumerate}
   \item $h^\pi \geq 0$ for all paths $\pi \in \Pi$.
That is, path flows must be nonnegative.
   \item $\sum_{\pi \in \Pi^{rs}}h^\pi = d^{rs}$ for all OD pairs $(r,s)$.
Together with the nonnegativity condition, this requires that every vehicle traveling from $r$ to $s$ is assigned to exactly one of the paths connecting these zones.
\end{enumerate}
Let $H$\label{not:H} denote the set of feasible assignments.

Link flows $x_{ij}$ and path flows $h^\pi$ are closely linked, and we can obtain the former from the latter.\index{link!flow!computing from path flow}
Let $\delta^\pi_{ij}$\label{not:deltaij} denote the number of times link $(i,j)$ is used by path $\pi$, so $\delta^\pi_{ij} = 0$ if path $\pi$ does not use link $(i,j)$, and $\delta^\pi_{ij} = 1$ if it does.
With this notation, we have
\labeleqn{linkpath}{x_{ij} = \sum_{r \in Z} \sum_{s \in Z} \sum_{\pi \in \Pi^{rs}} \delta^\pi_{ij} h^\pi \,.}
This can be more compactly written using matrix notation as
\labeleqn{linkpathmatrix}{\mathbf{x} = \bm{\Delta} \mathbf{h} \,,}
where $\mathbf{x}$ and $\mathbf{h}$ are the vectors of link and path flows, respectively, and $\bm{\Delta}$\label{not:Delta} is the \emph{link-path adjacency matrix}.\index{link-path adjacency matrix}
The number of rows in this matrix is equal to the number of links, and the number of columns is equal to the number of paths in the network, and the value in the row corresponding to link $(i,j)$ and the column corresponding to path $\pi$ is $\delta_{ij}^\pi$.

Given a feasible assignment, the corresponding link flows can be obtained by using equation~\eqn{linkpath} or~\eqn{linkpathmatrix}.
The set of \emph{feasible link assignments}\index{feasible link assignment} is the set $X$\label{not:X} of vectors $\mb{x}$ which satisfy~\eqn{linkpathmatrix} for some feasible assignment $\mb{h} \in H$.

Similarly, the path travel times $c^\pi$\index{path!travel time}\label{not:cpi} are directly related to the link travel times $t_{ij}$: the travel time of a path is simply the sum of the travel times of the links comprising that path.\index{path!travel time!computing from link travel time}
By the same logic, we have
\labeleqn{pathlinktimes}{c^\pi = \sum_{(i,j) \in A} \delta^\pi_{ij} t_{ij}}
or, in matrix notation,\label{not:transpose}
\labeleqn{pathlinkmatrix}{\mathbf{c} = \bm{\Delta}^T \mathbf{t} \,.}

A small example illustrates these ideas.
Consider the network in Figure~\ref{fig:netdefinitionsexample} with four nodes and six links.
The centroid nodes are shaded, and the link performance functions are as indicated.
For travelers from node 2 to node 4, there are two paths: $\{(2,3), (3,4)\}$ and $\{(2,4)\}$.
Using the compact notation for paths, we could also write these as $[2,3,4]$ and $[2,4]$.
Therefore, the set of acyclic paths between these nodes is $\Pi^{24} = \{ [2,3,4], [2,4] \}$.
You should verify for yourself that
\[ \Pi^{14} = \{ [1,2,4], [1,2,3,4], [1,3,4], [1,3,2,4] \}\,.\]
It will be useful to have specific indices for each path, i.e.
$\Pi^{24} = \{ \pi_1, \pi_2 \}$ so $\pi_1 = [2,3,4]$ and $\pi_2 = [2,4]$, and similarly $\Pi^{14} = \{ \pi_3, \pi_4, \pi_5, \pi_6 \}$ with $\pi_3 = [1,2,4]$, $\pi_4 = [1,2,3,4]$, $\pi_5 = [1,3,4]$, and $\pi_6 = [1,3,2,4]$.

Let's say that the demand for travel from centroid 1 to centroid 4 is 40 vehicles, and that the demand from 2 to 4 is 60 vehicles.
Then $d^{14} = 40$ and $d^{24} = 60$, and we must have 
\nolabeleqn{\sum_{\pi \in \Pi^{14}} h^\pi = h^3 + h^4 + h^5 + h^6 = d^{14} = 40}
and
\nolabeleqn{\sum_{\pi \in \Pi^{24}} h^\pi = h^1 + h^2 = d^{24} = 60 \,.}
Let's assume that the vehicles from each OD pair are divided evenly among all of the available paths, so $h^\pi = 10$ for each $\pi \in \Pi^{14}$ and $h^\pi = 30$ for each $\pi \in \Pi^{24}$.
We can now use equation~\eqn{linkpath} to calculate the link flows.
For instance the flow on link (1,2) is
\begin{multline}
\sum_{r \in Z} \left\{ \sum_{s \in Z} \left\{  \sum_{\pi \in \Pi^{rs}} \delta^\pi_{1,2} h^\pi \right\} \right\} = \\
= \left\{ \left\{ \delta^3_{(1,2)} h^3 + \delta^4_{(1,2)} h^4 + \delta^5_{(1,2)} h^5 + \delta^6_{(1,2)} h^6  \right\} \right\} + \left\{ \left\{ \delta^1_{(1,2)} h^1 + \delta^2_{(1,2)} h^2 \right\} \right\} = \\
= \left\{ \left\{ 1 \times 10 + 1 \times 10 + 0 \times 10 + 0 \times 10 \right\} \right\} + \left\{ \left\{ 0 \times 30 + 0 \times 30 \right\} \right \} = 20 \,.
\end{multline}               
where the braces show how the summations ``nest.''  Remember, this is just a fancy way of picking the paths which use link (1,2), and adding their flows.
The equation is for use in computer implementations or for large networks; when solving by hand, it's perfectly fine to just identify the paths using a particular link by inspection --- in this case, only paths $\pi^{1,4}_1$ and $\pi^{1,4}_2$ use link (1,2).
Repeating similar calculations, you should verify that $x_{13} = 20$, $x_{23} = 40$, $x_{32} = 10$, $x_{24} = 50$, and $x_{34} = 50$.

From these link flows we can get the link travel times by substituting the flows into the link performance functions, that is, $t_{12} = 10x_{12} = 200$, $t_{13} = 50 + x_{13} = 70$, $t_{23} = 10 + x_{23} = 50$, $t_{32} = 10 + x_{32} = 20$, $t_{24} = 50 + x_{24} = 100$, and $t_{34} = 10x_{34} = 500$.
Finally, the path travel times can be obtained by either adding the travel times of their constituent links, or by applying equation~\eqn{pathlinktimes}.
You should verify that  $c^1 = 550$, $c^2 = 100$, $c^3 = 300$, $c^4 = 750$, $c^5 = 570$, and $c^6 = 190$.

The role of traffic assignment is to choose one path flow vector $\mb{\hat{h}}$\label{not:hhat} for purposes of forecasting and ranking alternatives, out of all of the feasible assignments in the network.
An \emph{assignment rule}\index{assignment rule} is a principle used to determine this path flow vector.
The most common assignment rule in practice is that the path flow vector should place all vehicles on a path with minimum travel time between their origins and destinations, although other rules are possible as well and will be discussed later in the book.

\subsection{Commentary}
\label{sec:commentary}

\index{consistency|(}
The equations and concepts mentioned in the previous subsection can be related as follows.
Given a vector of path flows $\mathbf{h}$, we can obtain the vector of link flows $\mathbf{x}$ from equation \eqn{linkpathmatrix}; from these, we can obtain the vector of link travel times $\mathbf{t}$ by substituting each link's flow into its link performance function; from this, we can obtain the vector of path travel times $\mathbf{c}$ from equation \eqn{pathlinkmatrix}.
This process is shown schematically in Figure~\ref{fig:fundamentalloop}.
   
\stevefig{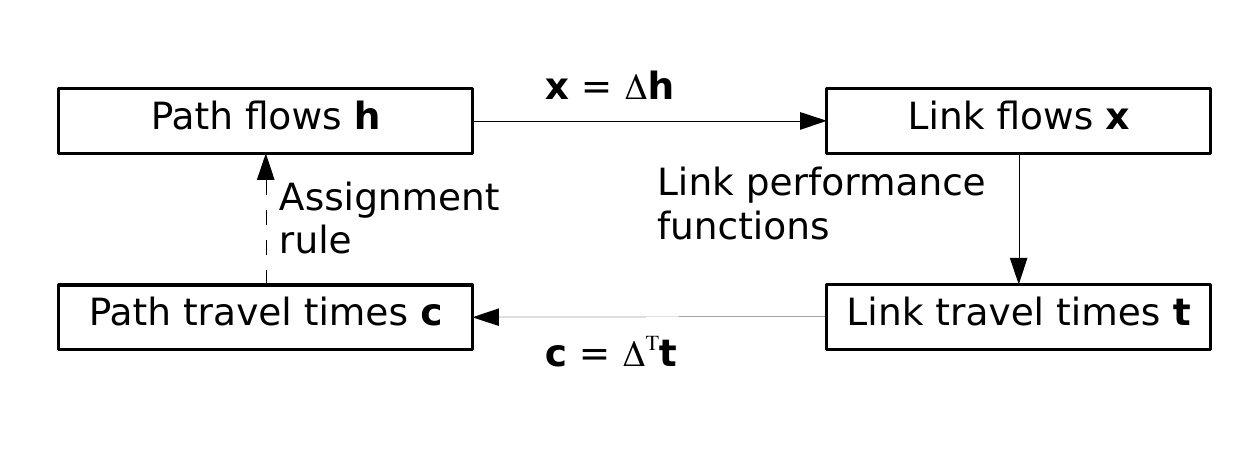}{Traffic assignment as an iterative process.}{\textwidth}   

The one component which does not have a simple representation is how to obtain path flows from path travel times using an assignment rule, ``completing the loop'' with the dashed line in the figure.
This is actually the most complicated step, and answering it will require most of the remainder of the chapter.
The main difficulty is that introducing some rule for relating path travel times to path flows creates a circular dependency: the path flows would depend on the path travel times, which depend on the link travel times, which depend on the link flows, which depend on the path flows, which depend on the path travel times and so on \emph{ad infinitum}.
We need to find a \emph{consistent} solution to this process, which in this case means a set of path flows which remain unchanged when we go around the circuit: the path flows must be consistent with the travel times we obtain from those same path flows.
Furthermore, the assignment rule must reflect the gaming behavior described in Section~\ref{sec:notionofequilibrium}.\index{consistency|)}

At this point, it is worthwhile to discuss the assumptions that go into the traffic assignment problem as stated here.
The first assumption concerns the time scale over which the network modeling is occurring.
This part of the book is focused entirely on what is called \emph{static} network modeling, in which we assume that the network is close to a ``steady state'' during whatever length of time we choose to model (whether a peak hour, a multi-hour peak period, or a 24-hour model), and the link flows and capacities are measured with respect to the entire time period.\index{link!capacity}\index{network!static}
That is, the capacity is the capacity \emph{over the entire analysis period}, so if a facility has a capacity of 2200 veh/hr and we are modeling a three-hour peak period, we would use a capacity of 6600.
Likewise, the link flows are the total flow over the three hours.\index{link!flow}

How long the analysis period should be is a matter of some balancing.\index{time horizon}
Obviously, the longer the analysis period, the less accurate the steady state assumption is likely to be.
In the extreme case of a 24-hour model, it is usually a (very big) stretch to assume that the congestion level will be the same over all 24 hours.
On the other hand, choosing too short an analysis period can be problematic as well.
In particular, if we are modeling a congested city or a large metropolitan area, trips from one end of the network to the other can easily take an hour or two, and it is good practice to have the analysis period be at least as long as most of the trips people are taking.

Properly resolving the issue of the ``steady state'' assumption requires moving to a \emph{dynamic} traffic assignment model.\index{dynamic traffic assignment}
Dynamic models have the potential to more accurately model traffic, but are harder to calibrate, are more sensitive to having correct input data, and require more computer time.
Further, dynamic models end up requiring very different formulations and approaches.
In short, while useful in some circumstances, they are not universally better than static models, and in any case they are surprisingly dissimilar.
As a result, a full discussion of dynamic models is deferred to the final part of this volume.

\index{model!continuous}
\index{continuum assumption}
Another assumption we make is that \emph{link and path flows can take any nonnegative real value}.
In particular, they are not required to be whole numbers, and there is no issue with saying that the flow on a link is 12.5 vehicles or that the number of vehicles choosing a path is, say, $\sqrt{2}$.
This is often called the \emph{continuum assumption}, because it treats vehicles as an infinitely-divisible fluid, rather than as a discrete number of ``packets'' which cannot be split.
The reason for this assumption is largely computational --- without the continuum assumption, traffic assignment problems become extremely difficult to solve, even on small networks.
Further, from a practical perspective most links of concern have volumes in the hundreds or thousands of vehicles per hour, where the difference between fractional and integer values is negligible.
Some also justify the continuum assumption by interpreting link and path flows to describe a long-term average of the flows, which may fluctuate to some degree from day to day.

\section{Principle of User Equilibrium}
\label{sec:principleofuserequilibrium}

\index{assignment rule|(}
\index{user equilibrium|(}
The previous section introduced the main ideas and components of static traffic assignment, but did not provide much explanation of assignment rules other than claiming that a reasonable assignment rule results in all used routes connecting an origin and destination to have equal and minimal travel time, and that the ``gaming'' equilibrium idea described in Section~\ref{sec:notionofequilibrium} is relevant to route choice.
This section explains this assignment rule in more detail.

Assignment rules are more complex than the other steps in traffic assignment shown in Figure~\ref{fig:fundamentalloop}, for two reasons.
The first is behavioral: unlike the other three steps, it does not follow from basic definitions, but instead represents more complicated human behavior.
The second is structural: adding the dashed link in the figure creates a cycle of dependency.
By linking path flows to path travel times, we now have four quantities which all depend on each other, and untangling this cycle requires a little bit of thought.\index{consistency}

To handle this, let's consider a simpler situation first.
Rather than trying to find path flows (which require stating \emph{every} driver's route choice), let's stick to a single driver.
Why do they pick the route that they do?  The number of potential paths in a network between two even slightly distant points in a network is enormous; in a typical city there are literally millions of possible routes one could theoretically take between an origin and destination.
The vast majority of these are ridiculous, say, involving extraneous trips to the outskirts of the city and then back to the destination even though the origin is nearby.

Why are such paths ridiculous on their face?  Table~\ref{tbl:criteria} gives a list of potential criteria which a desirable route would have, given in roughly decreasing order (at least according to my tastes for a trip to work or school).
The most important times of day to properly model are the morning and evening peak periods, when most of the trips made are work-related commutes.
For these types of trips, most people will choose the most direct route to their destination, where ``direct'' means some combination of travel time, cost, and distance (which are usually correlated).
For simplicity, we'll use ``low travel time'' as our starting point, specified as Assumption~\ref{ass:lowtime}.

\begin{table}
   \centering
   \caption{Potential criteria in choosing a route.
\label{tbl:criteria}}
   \begin{tabular}{l}
      \hline
      Low travel time                                 \\
      Reliable travel time                            \\
      Low out-of-pocket cost (tolls, fuel, etc.)      \\
      Short distance                                  \\
      Bias toward (or away from) freeways             \\
      Low accident risk                               \\
      Few potholes                                    \\
      Scenic view                                     \\
      \hline
   \end{tabular}
\end{table}

\begin{ass}
\label{ass:lowtime}
\index{shortest path assumption}
\emph{(Shortest path assumption.)}  Each driver wants to choose the path between their origin and their destination with the least travel time.
\end{ass}
Notice that this principle does not include the impact a driver has on other drivers in the system, and a pithy characterization of Assumption~\ref{ass:lowtime} is that ``people are greedy.''  If this seems unnecessarily pejorative, Section~\ref{sec:motivatingexamples} shows a few ways that this principle can lead to suboptimal flow distributions on networks.
Again, I emphasize that we adopt this assumption because we are modeling urban, peak period travel which is predominantly composed of work trips.
If we were modeling, say, traffic flows around a national park during summer weekends, quality of scenery may be considerably more important than being able to drive at free-flow speed, and a different assumption would be needed.

\index{link!cost}
Further, the basic model developed in the first few weeks could really function just as well with cost or some other criterion, as long as it is \emph{separable} and \emph{additive} (that is, you can get the total value of the criterion by adding up its value on each link --- the travel time of a route is the sum of the travel times on its component links, the total monetary cost of a path is the sum of the monetary costs of each link, but the total scenic quality of a path may not be the sum of the scenic quality of each link), and can be expressed as a function of the link flows $\mathbf{x}$.
So even though we will be speaking primarily of travel times, it is quite possible, and sometimes appropriate, to use other measures as well.
A second assumption follows from modeling commute trips:
\begin{ass}
\index{perfect knowledge assumption}
\label{ass:perfectknowledge}
Drivers have perfect knowledge of link travel times.
\end{ass}
In reality, drivers' knowledge is not perfect --- but commutes are typically habitual trips, so it is not unreasonable to assume that drivers are experienced and well-informed about congestion levels at different places in the network, at least along routes they might plausibly choose.
We will later relax this assumption, but for now we'll take it as it greatly simplifies matters and is not too far from the truth for commutes.
(Again, in a national park or other place with a lot of tourist traffic, this would be a poor assumption.)

Now, with a handle on individual behavior, we can try to scale up this assumption to large groups of travelers.
What will be the resulting state if there are a large number of travelers who all want to take the fastest route to their destinations?  For example, if you have to choose between two routes, one of which takes ten minutes and the other fifteen, you would always opt for the first one.
If you are the only one traveling, this is all well and good.
The situation becomes more complicated if there are others traveling.
If there are ten thousand people making the same choice, and all ten thousand pick the first route, congestion will form and the travel time will rapidly increase.
According to Assumption~\ref{ass:perfectknowledge}, drivers would become aware of this, and some people would switch from the first route to the second route, because the first would no longer be faster.
This process would continue: as long as the first route is slower, people would switch away to the second route.
If too many people switch to the second route, so the first becomes faster again, people would switch back.

With a little thought, it becomes clear that if there is \emph{any} difference in the travel times between the two routes, people will switch from the slower route to the faster one.
Note that Assumption~\ref{ass:lowtime} does not make any allowance for a driver being satisfied with a path which is a minute slower than the fastest path, or indeed a second slower, or even a nanosecond slower.
Relaxing Assumption~\ref{ass:lowtime} to say that people may be indifferent as long as the travel time is ``close enough'' to the fastest path leads to an interesting, but more complicated line of research based on the concept of ``bounded rationality,'' which is discussed later, in Section~\ref{sec:boundedrationality}.
It is much simpler to assume that Assumption~\ref{ass:lowtime} holds strictly, in which case there are only three possible stable states:
\begin{enumerate}
\item Route 1 is faster, even when all of the travelers are using it.
\item Route 2 is faster, even when all of the travelers are using it.
\item Most commonly, neither route dominates the others.
In this case, people use both Routes 1 and 2, and their travel times are \emph{exactly equal}.
\end{enumerate}

Because the third case is most common, this basic route choice model is called \emph{user equilibrium}: the two routes are in equilibrium with each other.
Why must the travel times be equal?  If the first route was faster than the second, people would switch from the second to the first.
This would decrease the travel time on the second route, and increase the travel time on the first, and people would continue switching until they were equal.
The reverse is true as well: if the second route were faster, people would switch from the first route to the second, decreasing the travel time on the first route and increasing the travel time on the second.
The \emph{only} outcome where nobody has any reason to change their decision, is if the travel times are equal on both routes.

This is important enough to state again formally:
\begin{cor}
\label{cor:ue}
\index{user equilibrium!principle of}
\emph{(Principle of user equilibrium.)}
Every used route connecting an origin and destination has equal and minimal travel time.
\end{cor}
Unused routes may of course have a higher travel time, and used routes connecting different origins and destinations may have different travel times, but any two used routes connecting the same origin and destination must have exactly the same travel times.
Notice that we call this principle a corollary rather than an assumption: the real assumptions are the shortest path\index{shortest path assumption} and full information assumptions.\index{perfect information assumption}
If you believe these are true, the principle of user equilibrium follows immediately and does not require you to assume anything more than you already have.
The next section describes how to solve for equilibrium, along with a small example.\index{user equilibrium|)}\index{assignment rule|)}

\subsection{A trial-and-error solution method}
\label{sec:solvingforequilibrium}

\index{static traffic assignment!algorithms!trial-and-error|(}
We can develop a simple method for solving for path flows using the principle of user equilibrium, using nothing but the definition itself.
For now, assume there is a single OD pair $(r,s)$.
The method is as follows:
\begin{enumerate}
\item Select a set of paths $\hat{\Pi}^{rs}$\label{not:Pihatrs} which you think will be used by travelers from this OD pair.\index{path!set of used paths}
\item Write equations for the travel times of each path in $\hat{\Pi}^{rs}$ as a function of the path flows.
\item Solve the system of equations enforcing equal travel times on all of these paths, together with the requirement that the total path flows must equal the total demand $d^{rs}$.
\item Verify that this set of paths is correct; if not, refine $\hat{\Pi}^{rs}$ and return to step 2.
\end{enumerate}
The first step serves to reduce a large potential set of paths to a smaller set of reasonable paths $\hat{\Pi}^{rs}$, which you believe to be the set of paths which will be used by travelers from $r$ to $s$.
Feel free to use engineering judgment here; if this is not the right set of paths, we'll discover this in step 4 and can adjust the set accordingly.
The second step involves applying equations~\eqn{pathlinktimes} and~\eqn{linkpath}.
Write the equations for link flows as a function of the flow on paths in $\hat{\Pi}^{rs}$ (assuming all other paths have zero flow).
Substitute these expressions for link flows into the link performance functions to get an expression for link travel times; then write the equations for path travel times as a function of link travel times.

If $\hat{\Pi}^{rs}$ truly is the set of used paths, all the travel times on its component paths will be equal.
So, we solve a system of equations requiring just that.
If there are $n$ paths, there are only $n - 1$ independent equations specifying equal travel times.
(If there are three paths, we have one equation stating paths one and two have equal travel time, and a second one stating paths two and three have equal travel time.
An equation stating that paths one and three have equal travel time is redundant, because it is implied by the other two, and so it doesn't help us solve anything.)  To solve for the $n$ unknowns (the flow on each path), we need one more equation: the requirement that the sum of the path flows must equal the total flow from $r$ to $s$ (the ``no vehicle left behind'' equation requiring every vehicle to be assigned to a path).
Solving this system of equations gives us path flows which provide equal travel times.

The last step is to verify that the set of paths is correct.
What does this mean?  There are two ways that $\hat{\Pi}^{rs}$ could be ``incorrect'':  either it contains paths that it shouldn't, or it omits a path that should be included.
In the former case, you will end up with an infeasible solution (e.g., a negative or imaginary flow on one path), and should eliminate the paths with infeasible flows, and go back to the second step with a new set $\hat{\Pi}^{rs}$.
In the latter case, you have a feasible solution and the travel times of paths in $\hat{\Pi}^{rs}$, but they are not minimal: you have missed a path which has a faster travel time than any of the ones in $\hat{\Pi}^{rs}$, so you need to include this path in $\hat{\Pi}^{rs}$ and again return to the second step.

Let's take a concrete example: Figure~\ref{fig:2link} shows 7000 travelers traveling from zone 1 to zone 2 during one hour, and choosing between the two routes mentioned above: route 1, with free-flow time 20 minutes and capacity 4400 veh/hr, and route 2, with free-flow time 10 minutes and capacity 2200 veh/hr.
That means we have
\begin{align}
t_1(x_1) &= 20 \left(1 + 0.15 \left( \frac{x_1}{4400} \right) ^4 \right)  \,,\\
t_2(x_2) &= 10 \left(1 + 0.15 \left( \frac{x_2}{2200} \right) ^4 \right) \,.
\end{align}

This example is small enough that there is no real distinction between paths and links (because each path consists of a single link), so link flows are simply path flows ($x_1 = h_1$ and $x_2 = h_2$), and path travel times are simply link travel times ($c_1 = t_1$ and $c_2 = t_2$).
As a starting assumption, we assume that both paths are used, so $\Pi^{12} = \{ \pi_1, \pi_2 \}$.
We then need to choose the path flows $h_1$ and $h_2$ so that $t_1(h_1) = t_2(h_2)$ (equilibrium) and $h_1 + h_2 = 7000$ (no vehicle left behind).
Substituting $h_2 = 7000 - h_1$ into the second delay function, the equilibrium equation becomes
\begin{equation}
20 \left(1 + 0.15 \left( \frac{h_1}{4400} \right) ^4 \right) = 10 \left(1 + 0.15 \left( \frac{7000 - h_1}{2200} \right) ^4 \right) \,.
\end{equation}
Using a numerical equation solver, we find that equilibrium occurs for $h_1 = 3376$, so $h_2 = 7000 - 3376 = 3624$, and $t_1(x_1) = t_2(x_2) = 21.0$ minutes.
Alternately, we can use a graphical approach.
Figure~\ref{fig:2linkbpr} plots the travel time on \emph{both} routes as a function of the flow on route 1 (because if we know the flow on route 1, we also know the flow on route 2).
The point where they intersect is the equilibrium: $h_1 = 3376$, $t_1 = t_2 = 21.0$.

\stevefig{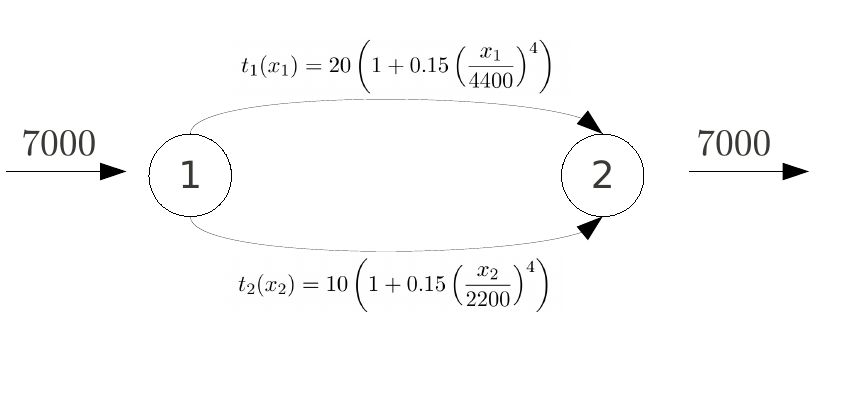}{Small example using the two links.}{\textwidth}
\stevefig{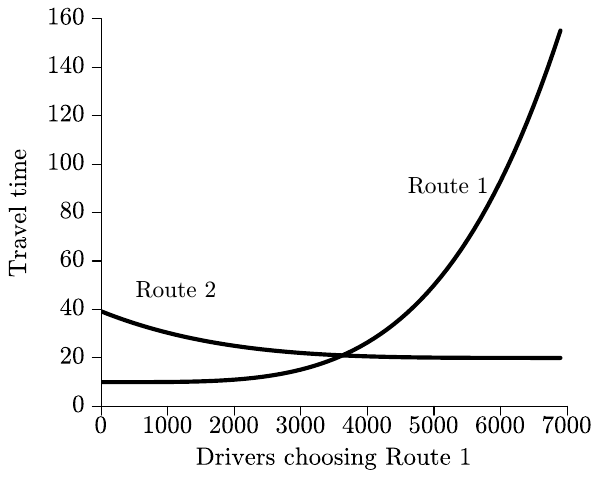}{The equilibrium point lies at the intersection of the delay functions.}{0.7\textwidth}

In the last step, we verify that the solution is reasonable (no paths have negative flow) and complete (there are no paths we missed which have a shorter travel time).
Both conditions are satisfied, so we have found the equilibrium flow and stop.

Let's modify the problem slightly, so the travel time on link 1 is now
\[ t_1(x_1) = 50 \left(1 + 0.15 \left( \frac{x_1}{4400} \right) ^4 \right)\,.\]
In this case, solving 
\[50 \left(1 + 0.15 \left( \frac{h_1}{4400} \right) ^4 \right) = 10 \left(1 + 0.15 \left( \frac{7000 - h_1}{2200} \right) ^4 \right)\]
leads to a nonsensical solution: none of the answers involve real numbers without an imaginary part.
The physical interpretation of this is that \emph{there is no way to assign 7000 vehicles to these two paths so they have equal travel times}.
Looking at a plot (Figure~\ref{fig:2linkbpr_modified}), we see that this happens because path 2 dominates path 1: even with all 7000 vehicles on path 2, it has a smaller travel time.

So, the solution to the modified problem is: $h_1 = 0$, $h_2 = 7000$, $t_1 = 50$, and $t_2 = 39.2$.
This still satisfies the principle of user equilibrium: path 1 is not \emph{used}, so it is fine for its travel time to not equal that of path 2.

To find this solution in the trial-and-error method, we would start by removing \emph{both} paths (since imaginary path flows are not meaningful), and then adding in the shortest path at free-flow (since there must be at least one used path).
When there is only one path in $\hat{\Pi}^{rs}$, the equilibrium state is easy to find, because the ``system of equations'' reduces trivially to placing all demand onto the single used path.

\stevefig{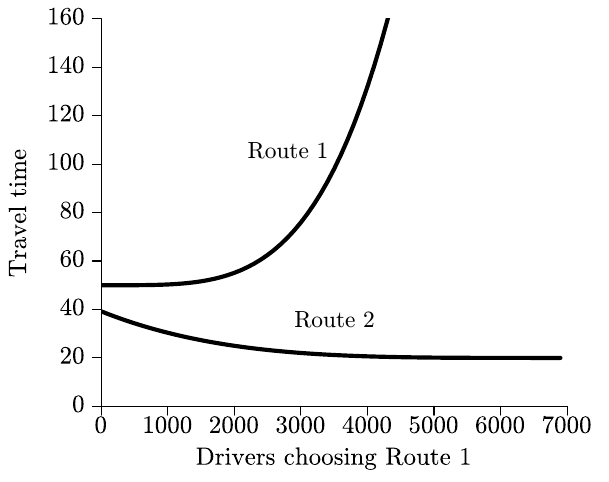}{There is no intersection point: path 1 is dominated by path 2.}{0.7\textwidth}

In more complicated networks, writing all of the equilibrium equations and solving them simultaneously is much too difficult.
Instead, a more systematic approach is taken, and will be described in Chapter~\ref{chp:solutionalgorithms}.\index{static traffic assignment!algorithms!trial-and-error|)}

\section{Three Motivating Examples}
\label{sec:motivatingexamples}

It is worth asking whether user equilibrium is the best possible condition.
``Best'' is an ambiguous term, but can be related towards our general goals as transportation engineers.
User equilibrium probably does not lead to the flow pattern with, say, the best emissions profile, simply because there's no reason to believe that the principle of user equilibrium has any connection whatsoever to emissions --- it is based on people trying to choose fastest paths.
But maybe it is related to congestion in some way, and you might find it plausible that each individual driver attempting to choose the best path for himself or herself would minimize congestion system-wide.
This section presents three examples which should shatter this innocent-seeming idea.
In this section, we will not concern ourselves with how a feasible assignment satisfying this assignment rule is found (the trial-and-error method from the previous section would suffice), but will focus instead on how equilibrium can be used to evaluate the performance of potential alternatives.

\index{paradoxes!Knight-Pigou-Downs|(}
In the first example, consider the network shown in Figure~\ref{fig:knightpigoudowns} where the demand between nodes 1 and 2 is $d^{12} = 30$ vehicles.
Using $h^\uparrow$ and $h^\downarrow$ to represent the flows on the top and bottom paths in this network, the set of feasible assignments are the two-dimensional vectors $\mb{h} = \vect{h^\uparrow & h^\downarrow}$ which satisfy the conditions $h^\uparrow \geq 0$, $h^\downarrow \geq 0$, and $h^\uparrow + h^\downarrow = 30$.
The figure shows the link performance functions for each of the two links in the network.
Notice that the travel time on the top link is constant irrespective of the flow, while the travel time on the bottom link increases with its flow.
However, at low values of flow, the bottom link is faster than the top one.
This corresponds to a scenario where one link is shorter, but more subject to congestion, while the other link is longer but free of congestion.

\stevefig{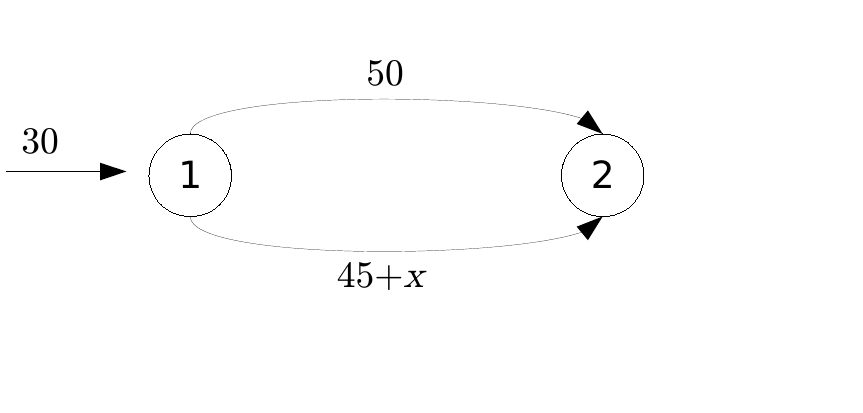}{A two-link network.}{0.8\textwidth}   

With our stated assignment rule, the solution to the traffic assignment problem is $h^\uparrow = 25$, $h^\downarrow = 5$, because this will result in link flows of 25 and 5 on the top and bottom links, respectively, giving travel times of 50 on both paths.
In this state, all vehicles in the network experience a travel time of 50.

Now, suppose that it is possible to improve one of the links in the network.
The intuitive choice is to improve the link which is subject to congestion, say, changing its link performance function from $45 + x^\downarrow$ to $40 + \frac{1}{2} x^\downarrow$, as shown in Figure~\ref{fig:knightpigoudownsimproved}.
(Can you see why this is called an improvement?)  However, in this case the path flow solution which corresponds to the assignment rule is $h^\uparrow = 10$, $h^\downarrow = 20$, because this results in equal travel times on both the top and bottom paths.
These travel times are still 50!  Even though the bottom link was improved, the effect was completely offset by vehicles switching paths away from the top link and onto the bottom link.

\stevefig{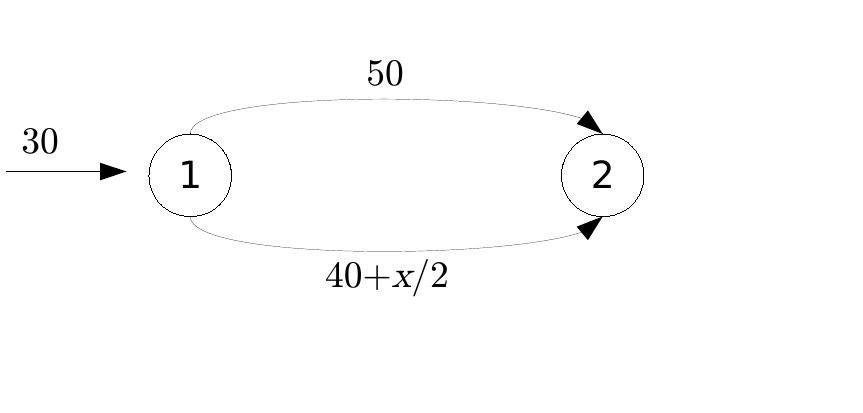}{The two-link network with an improvement on the bottom link.}{0.8\textwidth}   

In other words, improving a link (even the only congested link in a network) does not necessarily reduce travel times, because drivers can change their behavior in response to changes on the network.
(If we could somehow force travelers to stay on the same paths they were using before, then certainly some vehicles would experience a lower travel time after the improvement is made.)  This is called the \emph{Knight-Pigou-Downs}\index{paradoxes!Knight-Pigou-Downs|)}\index{Pigou, Arthur Cecil} paradox.

\index{paradoxes!Braess|(}
The second example was developed by Dietrich Braess\index{Braess, Dietrich}, and shows a case where building a new roadway link can actually worsen travel times for all travelers in the network, after the new equilibrium is established.
Assume we have the network shown in Figure~\ref{fig:braesssimple}a, with the link performance functions next to each link.
The user equilibrium solution can be found by symmetry: since the top and bottom paths are exactly identical, the demand of six vehicles will evenly split between them.
A flow of three vehicles on the two paths leads to flow of three vehicles on each of the four links; substituting into the link performance functions gives travel times of 53 on $(1,3)$ and $(2,4)$ and 30 on $(1,2)$ and $(3,4)$, so the travel time of both paths is 83, and the principle of user equilibrium is satisfied.

\stevefig{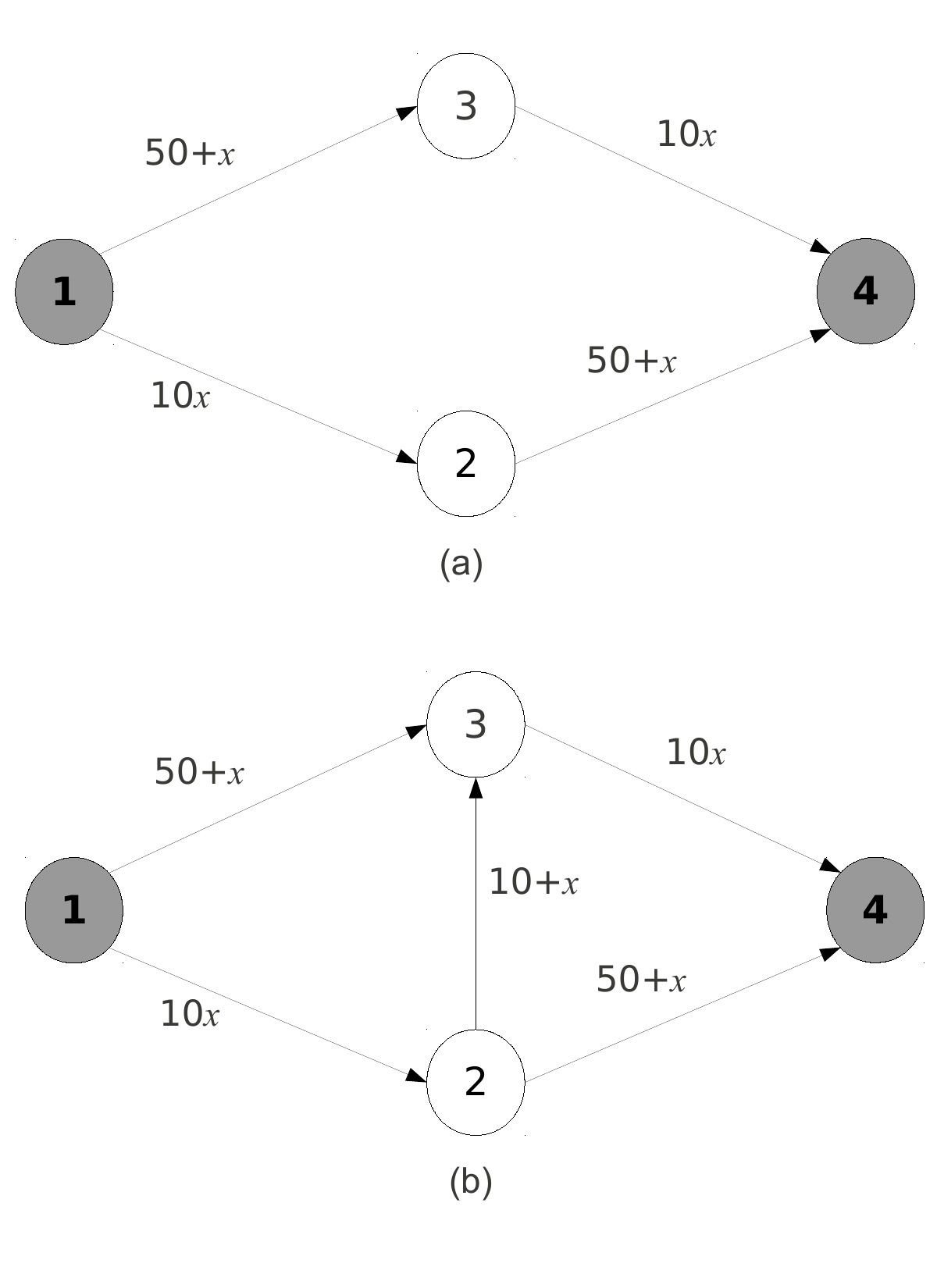}{Braess network, before and after construction of a new link.}{0.6\textwidth}   

Now, let's modify the network by adding a new link from node 2 to node 3, with link performance function $10 + x_{23}$ (Figure~\ref{fig:braesssimple}b).
We have added a new path to the network; let's label these as follows.
Path 1 is the top route $[1, 3, 4]$; path 2 is the middle route $[1, 2, 3, 4]$, and path 3 is the bottom route $[1, 2, 4]$.
Paths 1 and 3 each have a demand of three vehicles and a travel time of 83.
Path 2, on the other hand, has a flow of zero vehicles and a travel time of 70, so the principle of user equilibrium is violated: the travel time on the used paths is equal, but not minimal.

Assumption~\ref{ass:lowtime} suggests that someone will switch their path to take advantage of this lower travel time; let's say someone from path 1 switches to path 2, so we have $h_1 = 2$, $h_2 = 1$, and $h_3 = 3$.
From this we can predict new link flows: $x_{12} = 4$, $x_{13} = 2$, $x_{23} = 1$, $x_{24} = 3$, $x_{34} = 3$.
Substituting into link performance functions gives new link travel times: $t_{12} = 40$, $t_{13} = 52$, $t_{23} = 11$, $t_{24} = 53$, $t_{34} = 30$, and finally we can recover the new path travel times: $c_1 = 82$, $c_2 = 81$, and $c_3 = 93$.

This is still not an equilibrium; perhaps someone from path 3 will switch to path 2 in an effort to save 12 minutes of time.
So now $h_1 = h_2 = h_3 = 2$, the link flows are $x_{12} = x_{34} = 4$ and $x_{13} = x_{23} = x_{24} = 2$, the travel times are $t_{12} = t_{34} = 40$, $t_{13} = t_{24} = 52$ and $t_{23} = 12$, and the path travel times are $c_1 = c_2 = c_3 = 92$.
We have found the new equilibrium, and it is unconditionally worse than the old one.
Before adding the new link, each driver had a travel time of 83; now every driver has a travel time of 92.
Nobody is better off.
Everyone is worse off.
This is the famous \emph{Braess paradox}.
\index{paradoxes!Braess|)}

\index{paradoxes!Smith|(}
The third example, adapted from M. J. Smith,\index{Smith, Michael J.} shows how the traffic assignment problem can provide more insight into other traffic engineering problems such as signal timing.
Consider the network shown in Figure~\ref{fig:smithnetwork}.
Drivers choose one of two links, which join at a signalized intersection at the destination.
Assume that each of these links has the same free-flow time ($t^\uparrow_0 = t^\downarrow_0 = 1 \unit{min}$), but that the saturation flow of the bottom link is twice that of the top link ($s^\uparrow = 30 \unit{veh/min}$ and $s^\downarrow = 60 \unit{veh/min}$)\label{not:si}.
Vehicles enter the network at a demand level of $d = 35 \unit{veh/min}$.
The signal has a cycle length of $C = 1 \unit{min}$\label{not:C} and, for simplicity's sake, assume that there is no lost time so that the green times allocated to the top and bottom links equal the cycle length: $G^\uparrow + G^\downarrow = 1$.\label{not:Gi}

\stevefig{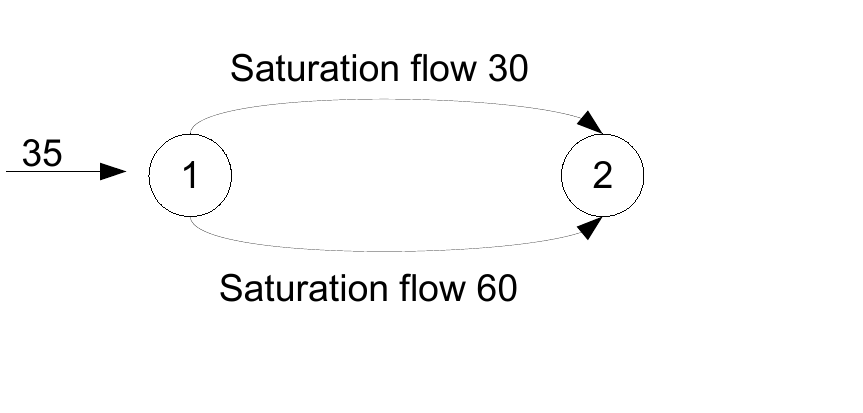}{Network for Smith's paradox.}{0.6\textwidth}

In traditional traffic signal analysis\index{traffic signal timing}, the \emph{capacity}\index{link!capacity} $u$ of an approach is the saturation flow scaled by the proportion of green time given to that approach, so $u_i = s_i (G_i / C)$ where $i$ can refer to either the top or bottom link.
The \emph{degree of saturation}\index{degree of saturation} $X$\label{not:Ximinor} for an approach is the ratio of the link flow to the capacity, so $X_i = x_i / u_i = (x_i C) / (s_i G_i)$.
With these quantities, the total travel time on each link (free-flow time plus signal delay) can be written as
\labeleqn{smithsignal}{t_i = 1 + \frac{9}{20} \mys{\frac{C(1 - G_i/C)^2}{1 - x_i / s_i} + \frac{X_i^2}{x_i (1 - X_i)}} \,.}

Assume that initially, the signal is timed such that $G^\uparrow = 40 \unit{sec}$ and $G^\downarrow = 20 \unit{sec}$.
Then the formulas for delay on the top and bottom links are functions of their flows alone, obtained from~\eqn{smithsignal} by substituting the corresponding green times into the expressions for $X_i$, and the equilibrium solution can be obtained by solving the equations $t^\uparrow = t^\downarrow$ and $x^\uparrow + x^\downarrow = 35$ simultaneously.
With these green times, the equilibrium solution occurs when $x^\uparrow = 23.6 \unit{veh/min}$ and $x^\downarrow =  11.4 \unit{veh/min}$, and the reader can confirm that the travel times on the two links are $t^\uparrow = t^\downarrow = 2.11 \unit{min}$.

So far, so good.
Now assume that it has been a long time since the signal was last retimed, and a traffic engineer decides to check on the signal and potentially change the timing.
A traditional rule in traffic signal timing is that the green time given to an approach should be proportional to the degree of saturation.
However, with the given solution, the degrees of saturation are $X^\uparrow = 0.982$ and $X^\downarrow = 0.953$, which are unequal --- the top link is given slightly less green time than the equisaturation rule suggests, and the bottom link slightly more.
Therefore, the engineer changes the green times to equalize these degrees of saturation, which occurs if $G^\uparrow = 48.3 \unit{sec}$ and $G^\downarrow = 11.7 \unit{sec}$, a smallish adjustment.
If drivers could be counted on to remain on their current routes, all would be well.
However, changing the signal timing changes the delay formulas~\eqn{smithsignal} on the two links.
Under the assumption that drivers always seek the shortest path, the equilibrium solution will change as drivers swap from the (now longer) path to the (now shorter) path.
Re-equating the travel time formulas, the flow rates on the top and bottom links are now $x^\uparrow = 23.8 \unit{veh/min}$ and $x^\downarrow = 11.2 \unit{veh/min}$, with equal delays $t^\uparrow = t^\downarrow = 2.26 \unit{min}$ on each link.
Delay has actually increased, because the signal re-timing (aimed at reducing delay) did not account for the changes in driver behavior after the fact.

A bit surprised by this result, our diligent traffic engineer notes that the degrees of saturation are still unequal, with $X^\uparrow = 0.984$ and $X^\downarrow = 0.959$, actually further apart than before the first adjustment.
Undeterred, the engineer changes the signal timings again to $G^\uparrow = 48.5 \unit{s}$ and $G^\downarrow = 11.5 \unit{s}$, which results in equal degrees of saturation under the new flows.
But by changing the green times, the delay equations have changed, and so drivers re-adjust to move toward the shorter path, leading to $x^\uparrow = 23.9 \unit{veh/min}$ and $x^\downarrow = 11.1 \unit{veh/min}$, and new delays of 2.43 minutes on each approach, even higher than before!

You can probably guess what happens from here, but Table~\ref{tbl:smithtable} tells the rest of the story.
As our valiant engineer stubbornly retimes the signals in a vain attempt to maintain equisaturation, flows always shift in response.
Furthermore, the delays grow faster and faster, asymptotically growing to infinity as more and more adjustments are made.
The moral of the story?  \emph{Changing the network will change the paths that drivers take}, and ``optimizing'' the network without accounting for how these paths will change is naive at best, and counterproductive at worst.
\index{paradoxes!Smith|)}

\begin{landscape}
\begin{table}
\begin{center}
\caption{Evolution of the network as signals are iteratively retimed.
\label{tbl:smithtable}}
\begin{tabular}{c|cc|cc|cc|cc}
Iteration & $G^\uparrow$ (s)  & $G^\downarrow$ & $x^\uparrow$ (v/min)   & $x^\downarrow$ & $t^\uparrow$ (min)  & $t^\downarrow$ & $X^\uparrow$ &   $X^\downarrow$ \\
\hline
0  &48   &12   &23.6 &11.4 &2.11 &2.11 &0.982   &0.953   \\
1  &48.3 &11.7 &23.8 &11.2 &2.26 &2.26 &0.984   &0.959   \\
2  &48.5 &11.5 &23.9 &11.1 &2.43 &2.43 &0.986   &0.965   \\
3  &48.7 &11.3 &24.1 &10.9 &2.63 &2.63 &0.988   &0.969   \\
4  &48.9 &11.1 &24.2 &10.8 &2.86 &2.86 &0.99 &0.974   \\
5  &49.1 &10.9 &24.3 &10.7 &3.12 &3.12 &0.991   &0.977   \\
10 &49.5 &10.5 &24.7 &10.3 &5.11 &5.11 &0.996   &0.989   \\
20 &49.88   &10.12   &24.91   &10.09   &16.58   &16.58   &0.9988  &0.9971  \\
50 &49.998  &10.002  &24.998  &10.002  &855.92  &855.93  &0.99998 &0.99995 \\
\hline
$\infty$ &50   &10   &25   &10   &$\infty$&$\infty$   &1 &1 \\
\end{tabular}
\end{center}
\end{table}
\end{landscape}

\section{Practical Considerations}
\label{sec:tapdata}

\index{planning|(}
Solving the traffic assignment problem on a network requires several pieces of data.
This section briefly describes some of the issues involved in collecting this data, and in calibrating such a model.
At a high level, three types of data are needed: (i) the network topology, (ii) the link performance functions, and (iii) the demand (OD) matrix.
Each of these has different considerations, and is discussed in turn.

Before discussing these specific data types, it is crucial to remember that \emph{traffic assignment models are used to predict future, hypothetical conditions}.
As sensing technologies become cheaper and more ubiquitous, it is becoming easier and easier to know the \emph{current} and \emph{historical} states of a traffic network, but that is not the domain of traffic assignment.
Rather, we are using these models to \emph{predict} or \emph{anticipate} how changes in the network (demand or supply) will affect congestion and other metrics of interest, and thereby provide policy and resource allocation recommendations.
Knowing the current and historical patterns on the network is important in building such a model, but they are not sufficient to make these decisions on their own.

\index{network!practical considerations|(}
Selecting a network topology means deciding what the nodes and links are in your network.
This involves questions of geographic scope and of detail.
In terms of geographic scope, how large of an area do you need to model?
The networks used by planning organizations commonly involve an entire metropolitan area: a major city, its surrounding suburbs, and some of the adjoining rural areas.
Such a network can be re-used for many applications in this metropolitan area.
Networks can also be purpose-built for a specific application, and in this case they often cover a smaller area.
It is important for the geographic scale to be large enough to include the main route alternatives being considered in your application.
If there is no route choice, the traffic assignment model will not tell you much of anything meaningful.

In terms of detail, you need to decide which links and nodes are important enough to model.
One option is to include \emph{all} of the roadway links within the geographic scope you have chosen.
This is the simplest choice, and there are tools designed to automatically generate networks this way, using open data sources.
However, in most cities the great majority of these links are neighborhood streets and other local roads.
These links are usually not of interest: they carry little traffic volume and are rarely congested.
Including all of these will often increase the size of your network by an order of magnitude, with little added value.
On the other hand, highways, freeways, and major arterials clearly need to be included.
Minor arterials are often included as well.\index{network!practical considerations|)}

Many of the node locations will be determined once the links are chosen, at the points where they intersect.
Additional nodes may be created to serve as demand centroids, a notion first introduced in Section~\ref{sec:introductionassignment}.
Considerations related to centroids\index{node!centroid} are discussed later in this section, along with the OD matrix.
If centroids are used, you will also need to add \emph{centroid connector}\index{link!centroid connector} links joining each centroid to the links representing the physical transportation infrastructure.
Care should be taken to place centroid connectors in a way that represents where travelers starting or ending trips in the neighborhood represented by the centroid will actually get on or off of the main roadway network.
This is an underappreciated point of network construction: the choice of connectors (how many and where) can play a surprisingly large role in assignment results.
This is even more true in dynamic traffic assignment, as will be discussed more in Part~\ref{part:dynamictrafficassignment}.

\index{link performance function|(}
Selecting link performance functions means choosing a function $t_{ij}(x_{ij})$ for each link $(i,j)$, giving its travel time in terms of the number of vehicles trying to use that link.
(As stated in Section~\ref{sec:introductionstatic}, although the word ``flow'' is commonly used to describe $x_{ij}$ in a static model, it is better to think of it as a ``demand'' which can possibly exceed the capacity of the link.)
The unit of $x_{ij}$ is in vehicles over the entire modeling horizon.
If your study is for one peak hour, then $x_{ij}$ is in vehicles per hour; if you are modeling a three-hour peak period, then $x_{ij}$ is in vehicles per three hours; and so forth.
This is important to keep in mind: many delay formulas, such as signal delay formulas in the Highway Capacity Manual\index{Highway Capacity Manual} or other professional references, make assumptions about the unit of flow, and failing to adapt these formulas to the units of $x_{ij}$ will cause major errors.

The gold standard is to estimate these functions on each link based on observed conditions on that specific link, over a long enough time period to know how demand relates to travel time on that link.
There are several complications in this process.
First, networks for a large metropolitan area commonly include tens of thousands of links, and you may not have data on all of them.
Freeways and major arterials are often equipped with permanent sensors recording flow, speed, and other traffic variables, but other links may not be.
Second, congested traffic behaves in a complex way.
Sensor data during \emph{uncongested} times usually shows a clear trend between flow and volume, but during congested times there is considerable scatter in the data.
This is because of phenomena like stop-and-go oscillations and hysteresis which are visible in traffic flow observations, but not simple to model, especially not with a link performance function (which is crude relative to state-of-the-art traffic flow theories).
Finally, traffic sensors can easily record speed (from which travel time can be inferred) and volume, but measuring \emph{demand} is quite a bit harder.
Observed volumes can never exceed the roadway capacity, by definition.
But demand --- what the $x_{ij}$ variables really represent --- certainly can.

For all of these reasons, the most common practice is to adopt a small number of functional forms that can be applied to specific links.
One example that we've already seen is the Bureau of Public Roads function\index{link performance function!BPR (Bureau of Public Roads)}
\labeleqn{bprvdf}{
t_{ij} = t^0_{ij} \myp{1 + \alpha \myp{\frac{x_{ij}}{u_{ij}}}^\beta}
\,,
}
where $t^0_{ij}$ and $u_{ij}$ are the free-flow travel time and capacity on link $(i,j)$, and $\alpha$ and $\beta$ are shape parameters.
The advantage of this function is that it can be applied with relatively little data on each link: free-flow travel time can be estimated as the link length divided by its speed limit, and capacity can be estimated by multiplying the number of lanes by a standard value (such as 2100 vehicles per hour per lane for a freeway).
Professional references contain more sophisticated ways to estimate free-flow speed\index{speed!free-flow speed} and capacity if you have information such as lane width, interchange density, and so on.
When available, sensor data can also provide estimates of $t^0_{ij}$ and $u_{ij}$.

The shape parameters are commonly chosen as $\alpha = 0.15$ and $\beta = 4$, although some researchers have suggested other values may be better choices.
In particular, the traditional values of 0.15 and 4 were estimated with $u_{ij}$ being the ``practical capacity'' of the link, which is an obsolete term used in a different way than ``capacity''\index{link!capacity} is used in contemporary transportation.
To modern transportation professionals, \emph{capacity} refers to the maximum flow rate that can be sustained on a link, roughly corresponding to a level of service E.
\emph{Practical capacity}, by contrast, describes a state when speeds are about 15\% lower than free-flow, roughly corresponding to a level of service C.
Practical capacity is about 80 percent of the true capacity of the link.
It is a mistake to use $\alpha = 0.15$ and $\beta = 4$ without also using ``practical capacity'' (by reducing capacity by 20\%).
\textbf{This mistake is surprisingly common in practice, and will lead to a systematic underestimation of congestion.}
If you are using the true capacity for $u_{ij}$, \cite{horowitz91} recommends $\alpha = 0.83$ and $\beta = 5.5$ for freeways, and \cite{dowling97} recommend $\alpha = 0.05$ and $\beta = 10$ for arterials, and $\alpha = 0.20$ and $\beta = 10$ for freeways.

Equation~\eqn{bprvdf} is not the only functional form for a link performance functions.
Other popular choices are the conical function\index{link performance function!conical}
\labeleqn{conical}{
t_{ij} = t^0_{ij} \mys{ 2
                        + \sqrt{ \beta^2 \myp{1 - \frac{x_{ij}}{u_{ij}}}^2
                                 + \gamma^2
                               }
                        - \beta \myp{1 - \frac{x_{ij}}{u_{ij}}}
                        - \gamma
                      }
\,,
}
where $\beta$ is the same shape parameter used in equation~\eqn{bprvdf} but $\gamma$\label{not:gammavdf} is fixed at $(2\beta - 1)/(2\beta - 2)$, and the Ak\c{c}elik function\index{link performance function!Ak\c{c}elik}
\labeleqn{akcelic}{
t_{ij} = t^0_{ij}
         + \frac{1}{4} T \mys{\myp{\frac{x_{ij}}{u_{ij}} - 1}
                             + \sqrt{ \myp{\frac{x_{ij}}{u_{ij}} - 1}^2
                                      + \frac{8Jx_{ij}}{Tu_{ij}^2}
                                    }
                            }
\,,
}
with $T$\label{not:T} the length of the analysis period and $J$\label{not:J} a parameter representing queueing dynamics on the link\footnote{\cite{akcelik91} recommends $J = 0.1$ for freeways; for arterials the situation is more complex and depends on the type and spatial density of junction controls along the link.}.
See Section~\ref{sec:sta_references} for additional references on different link performance functions.\index{link performance function|)}

\index{OD matrix|(}
\index{OD matrix!estimation|(}
The OD matrix is the most challenging input to prepare, for several reasons.
Unlike the network topology and physical characteristics of the roadway, the demand pattern is not directly observable.
Sensors on links can report volumes, speeds, and so forth, but cannot report vehicle origins or destinations.
There are emerging techniques for observing demand patterns directly from cell phone, GPS, or Bluetooth traces, but there are complications --- for instance, trip starts and ends are often ``fuzzed'' for privacy reasons, but this is exactly the data we need for an OD matrix.
Such data is also not fully representative of the traveling population.
There are many \emph{travel demand models}\index{model!demand} that estimate an OD matrix from publicly-available demographic data such as census records, often calibrated using travel surveys.
These methods are beyond the scope of this book, but Section~\ref{sec:sta_references} provides some references for interested readers.

Many have hoped for methods that can directly estimate an OD matrix from observed traffic volumes on links.
The main problem with this is \emph{overfitting}.\index{overfitting}
The OD matrix grows with the square of network size (an entry for every pair of zones), whereas the number of links that can provide observations grows linearly.
This means that estimating an OD matrix from link volumes is greatly underdetermined in practical networks, and in general there are infinitely many OD matrices which can match a given set of link volumes.
In other words, simply matching link volumes is not hard --- rather, the problem is that it's too easy.
Only one of these infinite number of OD matrices is actually the correct one which will respond in the right way when the network changes (remember, the goal of transportation network models is to evaluate counterfactual scenarios for guiding infrastructure decisions and policy).
Therefore, \emph{estimating an OD matrix solely from observed link flows cannot be recommended}.
You are likely to find an OD matrix which fits the observed data very well, but responds in a very different way when anything changes.

A better strategy is to use observed link flows or other field observations to improve on the OD matrix produced by another travel demand model.  
One way to do this is discussed in Chapter~\ref{chp:sensitivityanalysis}.
The idea is to use a travel demand model\index{model!demand} to produce a ``seed'' OD matrix.\index{OD matrix!seed}\index{optimization!bilevel}
We then solve for equilibrium using this matrix, and compare the predicted link flows to field observations.
The seed OD matrix can then be adjusted in a way that improves the fit with observed link flows, while not moving too far way from the seed point.
This retains much of the travel behavior which is embedded in travel demand models, and results in a better OD matrix to use for scenario comparisons.\index{OD matrix|)}\index{OD matrix!estimation|)}\index{planning|)}

\section{Historical Notes and Further Reading}
\label{sec:sta_references}

\index{user equilibrium!principle of!variations|(}
\cite{beckmann56}\index{Beckmann, Martin} were the first to pose the static traffic assignment assignment
problem in its full generality, although similar concepts were expressed earlier in \cite{pigou20}\index{Pigou, Arthur Cecil} and \cite{wardrop52}.\index{Wardrop, John Glen}\index{Wardrop, John Glen!principles}
The principle of user equilibrium is sometimes referred to as Wardrop's first principle.  
(His second principle corresponds to the system optimum state, discussed more in Section~\ref{sec:systemoptimal}.)

Researchers have explored subtle variations of these principles.
As used in this text, and as most commonly used in our field, our principle is that any path with positive flow must have minimal cost, which implies that if multiple paths connecting the same OD pair have positive flow, they must have equal cost.
\cite{dafermos69} and \cite{smith84_alternative} adopt a slightly different condition, where any path with positive flow must have a lower cost than any other path, \emph{accounting for the change in that other path's travel time if some flow on the original path were to switch to it}.
\cite{heydecker86} then extends this principle further, proposing an equilibrium principle where any path with positive flow must have a lower cost than any other path, \emph{accounting for the change in \emph{both} paths' travel times if some flow on the original path were to switch to it.}
These principles all agree in the standard case of continuous, increasing link performance functions which are separable in the sense that the travel time on a link depends only on its own flow, but the equilibria can differ if these assumptions are violated.\index{continuous function!applications}\index{link performance function!separability}
(See Section~\ref{sec:linkinteractions} for traffic assignment with nonseparable link performance functions.)

There are also variations of these principles which do not adopt the continuum assumption, but instead assume that there are a finite, unsplittable number of agents choosing routes in the network.\index{model!continuous}
The most direct extension is to consider the case of discrete (integer) vehicle flows, as in \cite{rosenthal73}.
In this case, there may be no assignment that exactly equalizes the costs on used paths, so the principle of user equilibrium must be modified (e.g., the costs on used paths may not be equal, but no driver can reduce their cost of travel by switching paths).
Another extension is to consider the case where there are a finite number of ``players'' controlling multiple vehicles that are assigned to minimizes each player's total cost. \citep{orda93}
Such models have been used to describe managing fleets of vehicles, or to describe telecommunications and other infrastructure networks.\index{telecommunication networks}
Game theorists have formulated the notions of \emph{congestion games} and \emph{potential games} to generalize all of the above variants, even outside the setting of a network. \citep{rosenthal73_congestion, monderer96_potential}\index{game theory!congestion game}
\index{game theory!potential game}See also \cite{bernstein90} for a comparison of discrete and continuous formulations of network equilibrium problems.\index{user equilibrium!principle of!variations|)}

The first of the three motivating examples is due to \cite{pigou20}; it is referred to as the ``Knight-Pigou-Downs'' paradox as the same effect was independently discovered by \cite{knight24} and \cite{downs62}.
\index{Pigou, Arthur Cecil}
The second example was first described by \cite{braess69}, and an up-to-date list of research and popular works further exploring this paradox is maintained at \url{https://homepage.ruhr-uni-bochum.de/Dietrich.Braess/#paradox}.
The third example is adapted from \cite{smith78}.

As an alternative to dividing a city into discrete zones and centroids, see \cite{daganzo80a, daganzo80b}.
For a review of link performance functions\index{link performance function}, including their history, theoretical basis, and practical applications, see \cite{pan_vdf}.
The BPR function was introduced in \cite{bpr64}, the conical function in \cite{spiess90}, and the Ak\c{c}elic function in \cite{akcelik91}.
For more on OD matrix estimation, see Sections~\ref{sec:odme} and~\ref{sec:sensitivity_litreview}.\index{link performance function!BPR (Bureau of Public Roads)}\index{link performance function!conical}\index{link performance function!Ak\c{c}elic}

\section{Exercises}
\label{exercises_tapconcepts}

\begin{enumerate}
\item \diff{15} Expand the diagram of Figure~\ref{fig:complextaschematic} to include a ``government model'' which reflects public policy and regulation regarding tolls.
What would this ``government'' agent influence, and how would it be influenced?
\item \diff{23} How realistic do you think link performance functions are?  Name at least two assumptions they make about traffic congestion, and comment on how reasonable you think they are.
\item \diff{23} How realistic do you think the principle of user equilibrium is?  Name at least two assumptions it makes (either about travelers or congestion), and comment on how reasonable you think they are.
\item \diff{12} Develop analogies to the equilibrium principle to represent the following phenomena: (a) better restaurants tend to be more crowded; (b) homes in better school districts tend to be more expensive; (c) technologies to improve road safety (such as antilock brakes) can increase reckless driving.
Specify what assumptions you are making in these analogies.
\item \diff{34} \emph{(All-or-nothing solutions\index{all-or-nothing assignment} are corner points.)}\index{convex set!corner point}  The point $\mb{x}$ is a \emph{corner point} of the convex set $X$ if it is impossible to write $\mb{x} = \lambda \mb{y} + (1 - \lambda) \mb{z}$ for $\mb{y} \in X$, $\mb{z} \in X$, and $\lambda \in (0, 1)$.
Show that if $H$ is the set of feasible path flows, $\mb{h}$ is a corner point if and only if it is an \emph{all-or-nothing assignment} (a path assignment where all of the flow from each OD pair is on a single path).
Repeat for the set of feasible link flows.
\item \diff{14} There are three routes available to travelers between a single OD pair.
For each case below, you are given the travel time $t$ on each route, as well as the number of travelers $x$ using each route.
Indicate whether each case satisfies the principle of user equilibrium.
If this principle is violated, explain why.

\begin{enumerate}[(a)]
\item $t_1 = 45$, $t_2 = 60$, and $t_3 = 30$ while $x_1 = 0$, $x_2 = 50$, and $x_3 = 0$.

\item $t_1 = 45$, $t_2 = 30$, and $t_3 = 30$ while $x_1 = 0$, $x_2 = 50$, and $x_3 = 25$.

\item $t_1 = 45$, $t_2 = 30$, and $t_3 = 30$ while $x_1 = 0$, $x_2 = 50$, and $x_3 = 0$.

\item $t_1 = 15$, $t_2 = 30$, and $t_3 = 45$ while $x_1 = 30$, $x_2 = 0$, and $x_3 = 0$.

\item $t_1 = 15$, $t_2 = 30$, and $t_3 = 45$ while $x_1 = 30$, $x_2 = 20$, and $x_3 = 10$.

\item $t_1 = 15$, $t_2 = 30$, and $t_3 = 30$ while $x_1 = 0$, $x_2 = 50$, and $x_3 = 30$.

\end{enumerate}
\item \diff{26} The network in Figure~\ref{fig:eqmverify} has 8 nodes, 12 links, and 4 zones.
The travel demand is $d^{13} = d^{14} = 100$, $d^{23} = 150$, and $d^{24} = 50$.
The dashed links have a constant travel time of 10 minutes regardless of the flow on those links; the solid links have the link performance function $15 + x/20$ where $x$ is the flow on that link.
\label{ex:eqmverify}
\begin{enumerate}[(a)]
\item For each of the four OD pairs with positive demand, list all acyclic paths connecting that OD pair.
In total, how many such paths are in the network?
\item Assume that the demand for each OD pair is divided evenly among all of the acyclic paths you found in part (a) for that OD pair.
What are the resulting link flow vector $\mathbf{x}$, travel time vector $\mathbf{t}$, and path travel time vector $\mathbf{C}$?
\item Does that solution satisfy the principle of user equilibrium?
\item What is the total system travel time?
\end{enumerate}
\stevefig{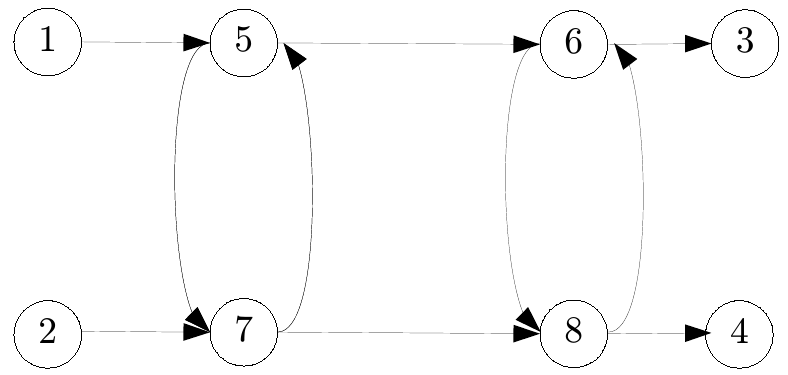}{Network for Exercise~\ref{ex:eqmverify}.}{0.5\textwidth}
\item \diff{32} Consider the network and OD matrix shown in Figure~\ref{fig:wyoeqm}.
The travel time on \emph{every} link is $10 + x/100$, where $x$ is the flow on that link.
Find the link flows and link travel times which satisfy the principle of user equilibrium.
\label{ex:wyoeqm}
\begin{figure}
\begin{center} 
\hfill
\includegraphics[width=0.5\textwidth]{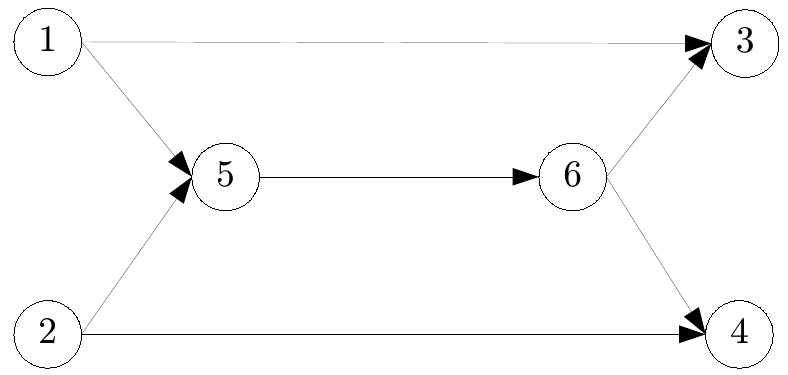}
\hfill
\raisebox{0.5in}{
\begin{tabular}{r|cc}
&  3  &  4 \\
\hline
1  &  5,000 & 0 \\
2  &  0     & 10,000 \\
\end{tabular}
}
\hfill
\caption{Network and OD matrix for Exercise~\ref{ex:wyoeqm}.
\label{fig:wyoeqm}}
\end{center}
\end{figure}
\item \diff{36} Find the equilibrium path flows, path travel times, link flows, and link travel times on the Braess network (Figure~\ref{fig:braesssimple}b) when the travel demand from node 1 to node 4 is (a) 2; (b) 7; and (c) 10.
\item \diff{52} In the Knight-Pigou-Downs network, the simple and seemingly reasonable heuristic to improve the congested link (i.e., changing its link performance function to lower the free-flow time) did not help.
Can you identify an equally simple heuristic that would lead you to improve the other link?  (Ideally, such a heuristic would be applicable in many kinds of networks, not just one with an uncongestible link.)
\item \diff{20} Give nontechnical explanations of why uniqueness, efficiency, and existence of equilibrium solutions have practical implications, not just theoretical ones.
Concrete examples may be helpful.
\item \diff{77} In the trial-and-error method, identify several different strategies for choosing the initial set of paths.
Test these strategies on networks of various size and complexity.
What conclusions can you draw about the effectiveness of these different strategies, and the amount of effort they involve? \label{ex:trialanderrortest}
\item \diff{55} The trial-and-error method generally involves the solution of a system of nonlinear equations.
\emph{Newton's method}\index{Newton's method} for solving a system of $n$ nonlinear equations is to move all quantities to one side of the equation, expressing the resulting equations in the form $\mb{f}(\mb{x}) = \mb{0}$ where $\mb{f}$ is a vector-valued function mapping the $n$-dimensional vector of unknowns $\mb{x}$ to another $n$-dimensional vector giving the value of each equation.
An initial guess is made for $\mb{x}$, which is then updated using the rule $\mb{x} \leftarrow \mb{x} - (J\mb{f}(\mb{x}))^{-1} \mb{f}(\mb{x})$, where $(J\mb{f}(\mb{x}))^{-1}$ is the inverse of the Jacobian matrix of $\mb{f}$\label{not:Jacobian}, evaluated at $\mb{x}$.
This process continues until (hopefully) $\mb{x}$ converges to a solution.
A \emph{quasi-Newton method} approximates the Jacobian with a diagonal matrix, which is equal to the Jacobian along the diagonal and zero elsewhere, which is much faster to calculate and convenient to work with.
Extend your experiments from Exercise~\ref{ex:trialanderrortest} to see whether Newton or quasi-Newton methods work better.\index{quasi-Newton method|see {Newton's method, quasi}}\index{Newton's method!quasi}
\item \diff{84} What conditions on the link performance functions are needed for Newton's method (defined in the previous exercise) to converge to a solution of the system of equations?  What about the quasi-Newton method?  Assume that the network consists of $n$ parallel links connecting a single origin to a single destination.
Can you guarantee that the solution to the system of equations only involves real numbers, and not complex numbers?  \emph{Hint: it may be useful to redefine the link performance functions for negative $x$ values.
This should not have any effect on the ultimate equilibrium solution, since negative link flows are infeasible, but cleverly redefining the link performance functions in this region may help show convergence.}
\item \diff{96} Repeat the previous exercise, but for a general network with any number of links and nodes, and where paths may overlap.
\end{enumerate}

\chapter{The Traffic Assignment Problem}
\label{chp:trafficassignmentproblem}

This chapter formalizes the user equilibrium traffic assignment problem defined in Chapter~\ref{chp:equilibrium}.
Using the mathematical language from Chapter~\ref{chp:mathematicalpreliminaries}, we are prepared to model and solve equilibrium problems even on networks of realistic size, with tens of thousands of links and nodes.
Section~\ref{sec:mathematicalformulations} begins with fixed point, variational inequality, and optimization formulations of the user equilibrium problem.
Section~\ref{sec:tapproperties} then introduces important existence and uniqueness properties of the user equilibrium assignment, as was first introduced with the small two-player games of Section~\ref{sec:notionofequilibrium}.
Section~\ref{sec:alternaterules} names several alternatives to the user equilibrium rule for assigning traffic to a network, including the system optimum principle, and the idea of bounded rationality.
This chapter concludes with Section~\ref{sec:ueso}, exploring the inefficiency of the user equilibrium rule and unraveling the mystery of the Braess paradox.

\section{Mathematical Formulations}
\label{sec:mathematicalformulations}

\index{variational inequality!and traffic assignment|(}
\index{static traffic assignment!formulation!as variational inequality|(}
The previous section introduced a trial-and-error method for solving the user equilibrium problem.
To solve equilibrium on larger networks, we need something more sophisticated than trial-and-error, but doing so requires a more convenient representation of the principle of user equilibrium in mathematical notation, using the tools defined in the preceding chapters.
This section does so by first formulating the equilibrium problem as the solution to a variational inequality; as the solution to a fixed point problem based on the variational inequality; and as the solution to a convex optimization problem.
The reason for presenting all three of these formulations (rather than just one) is that each is useful in different ways.
Just as true fluency in a language requires being able to explain the same concept in different ways, fluency in traffic assignment problems requires familiarity with all of these ways of understanding the principle of user equilibrium.
(It is also possible to formulate the principle of user equilibrium as a fixed point problem without referring to a variational inequality through the use of multifunctions, as described in the optional Subsection~\ref{sec:eqmmultifunction}.)

The variational inequality formulation proceeds as follows.
For each feasible path flow vector, we can associate a direction representing travelers' desire for lower travel time paths.
Figure~\ref{fig:feasibleset} shows the feasible set for the two-link equilibrium problem on the left, and a conceptual schematic of the feasible set for the general equilibrium problem on the right.
Notice that in all cases the feasible set is convex and has no ``gaps.'' For each feasible vector $\mb{h}$, we can calculate the path travel time vector $\mb{c}$ associated with these path flows.
Associate with each point a vector pointing in the \emph{opposite} direction as $\mb{c}$, that is, each point $\mb{h}$ is associated with the direction $-\mb{c}(\mb{h})$.
You can think about this as a force acting on the point $\mb{h}$.
The interpretation is that the path flow vector is ``pulled'' in the direction of decreasing travel times, and we will shortly show that an equilibrium is a point which is unmoved by this force.
Of course, this intuitive interpretation must be established mathematically as well, and we will show that a vector $\mb{h}$ satisfies the principle of user equilibrium if, and only if, the force criterion is satisfied.

\stevefig{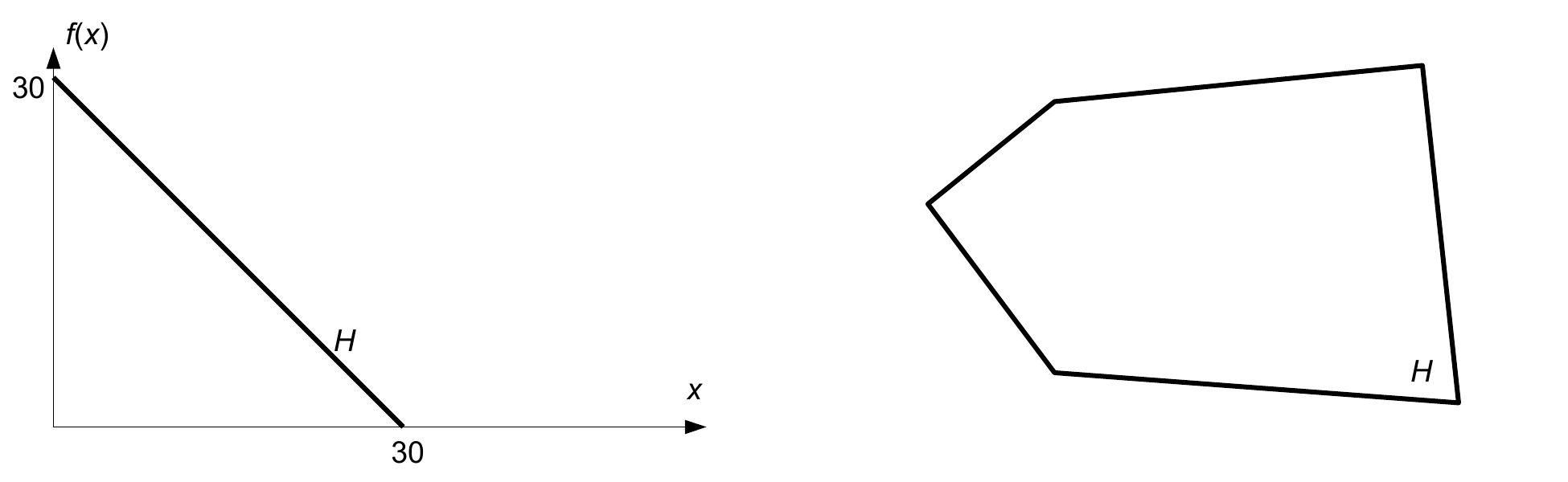}{Feasible path-flow sets $H$ for the simple two-link network (left) and a schematic representing $H$ for more complex problems (right). \index{feasible assignment!visualization}}{\textwidth}
\stevefig{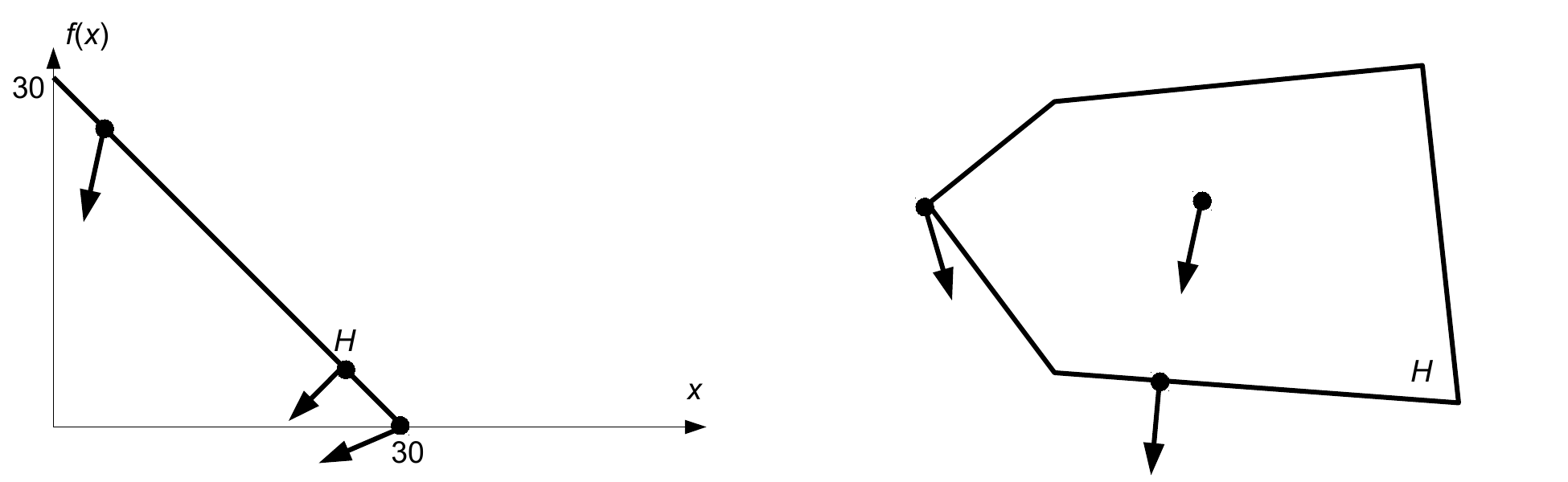}{``Force'' vectors $-\mb{c}$ for the simple two-link network (left) and a schematic representing $H$ for more complex problems (right).}{\textwidth}

A few points are shown in Figure~\ref{fig:forcevectors} as examples.
Recall that a non-corner point is unmoved by such a force only if the force is perpendicular to the boundary of the feasible set at that point.
At a corner point, the force must make a right or obtuse angle with all boundary directions.
So we can characterize stable path flow vectors $\mb{\hat{h}}$ as points where the force $-\mb{c}(\mb{\hat{h}})$ makes a right or obtuse angle with any feasible direction $\mb{h} - \mb{\hat{h}}$, where $\mb{h}$ is any other feasible path flow vector.
Recalling vector operations, saying that two vectors make a right or obtuse angle with each other is equivalent to saying that their dot product is nonpositive.
Thus, $\mb{\hat{h}}$ is a stable point if it satisfies 
\labeleqn{basicviunsimplified}{-\mb{c}(\mb{\hat{h}}) \cdot (\mb{h} - \mb{\hat{h}}) \leq 0 \qquad \forall \mb{h} \in H}
or, equivalently,
\labeleqn{basicvi}{\mb{c}(\mb{\hat{h}}) \cdot (\mb{\hat{h}} - \mb{h}) \leq 0 \qquad \forall \mb{h} \in H\,.}

This is a variational inequality (VI) in the form shown in Section~\ref{sec:variationalinequality}.
We now prove that solutions of this VI correspond to equilibria, as shown by the next result:  

\begin{thm}
\label{thm:stability}
A path flow vector $\mb{\hat{h}}$ solves the variational inequality~\eqn{basicvi} if and only if it satisfies the principle of user equilibrium.
\end{thm}
\begin{proof}
   The theorem can equivalently be written ``a path flow vector $\mb{\hat{h}}$ does not solve~\eqn{basicvi} if and only if it does not satisfy the principle of user equilibrium,'' which is easier to prove.
Assume $\mb{\hat{h}}$ does not solve~\eqn{basicvi}.
Then there exists some $\mb{h} \in H$ such that $\mb{c}(\mb{\hat{h}}) \cdot (\mb{\hat{h}} - \mb{h}) > 0$, or equivalently $\mb{c}(\mb{\hat{h}}) \cdot \mb{\hat{h}} > \mb{c}(\mb{\hat{h}}) \cdot \mb{h}$.
Now, $\mb{c}(\mb{\hat{h}}) \cdot \mb{\hat{h}}$ is the total system travel time when the path flows are $\mb{\hat{h}}$, and $\mb{c}(\mb{\hat{h}}) \cdot \mb{h}$ is the total system travel time \emph{if the travel times were held constant at $\mb{c}(\mb{\hat{h}})$ even when the flows changed to $\mb{h}$.}  For the latter to be strictly less than the former, switching from $\mb{\hat{h}}$ to $\mb{h}$ must have reduced at least one vehicle's travel time even though the path travel times did not change.
This can only happen if that vehicle was not on a minimum travel time path to begin with, meaning that $\mb{\hat{h}}$ does not satisfy user equilibrium.
   
Conversely, assume that $\mb{\hat{h}}$ is not a user equilibrium.
Then there is some OD pair $(r,s)$ and path $\pi$ such that $\hat{h}_\pi > 0$ even though $c^\pi(\mb{\hat{h}}) > \min_{\pi' \in \Pi^{rs}} c^{\pi'} (\mb{\hat{h}})$.
Let $\hat{\pi}$ be a minimum travel time path for this OD pair.
Create a new path flow vector $\mb{h}$ which is the same as $\mb{\hat{h}}$ except that $\hat{h}^\pi$ is reduced by some small positive amount $\epsilon$\label{not:epsilon} and $\hat{h}^{\hat{\pi}}$ is increased by $\epsilon$.
As long as $0 < \epsilon < \hat{h}^\pi$ the new point $\mb{h}$ remains feasible.
By definition all the components of $\mb{\hat{h}} - \mb{h}$ are equal to zero, except the component for $\pi$ is $\epsilon$ and the component for $\hat{\pi}$ is $-\epsilon$.
So $\mb{c}(\mb{\hat{h}}) \cdot (\mb{\hat{h}} - \mb{h}) = \epsilon( c^\pi(\mb{\hat{h}}) - c^{\hat{\pi}}(\mb{\hat{h}}) ) > 0$, so $\mb{\hat{h}}$ does not solve~\eqn{basicvi}.
\end{proof}

\index{fixed point problem!and traffic assignment|(}
\index{static traffic assignment!formulation!as fixed point|(}
This variational inequality formulation leads directly to a fixed point formulation, as was discussed in Section~\ref{sec:variationalinequality}.
Theorem~\ref{thm:stability} shows that the equilibrium path flow vectors are exactly the solutions of the VI~\eqn{basicvi}, that is, the stable points with respect to the fictitious force $-\mb{c}$.
Therefore, an equilibrium vector $\mb{\hat{h}}$ is a fixed point of the function $\mr{proj}_H (\mb{h} -\mb{c}(\mb{h}))$, that is, for all equilibrium solutions $\mb{\hat{h}}$ we have
\labeleqn{equilibriumfixedpoint}{\mb{\hat{h}} = \mr{proj}_H (\mb{\hat{h}} -\mb{c}(\mb{\hat{h}}))\,.}
\index{static traffic assignment!formulation!as variational inequality|)}
\index{variational inequality!and traffic assignment|)}
\index{fixed point problem!and traffic assignment|)}
\index{static traffic assignment!formulation!as fixed point|)}

\index{optimization!convex optimization!and traffic assignment|(}
\index{static traffic assignment!formulation!as convex optimization|(}
\index{Beckmann formulation!see {static traffic assignment, formulation, as convex optimization}}
Finally, user equilibrium can also be shown as the solution to a convex optimization problem.
The main decision variables are the path flows $\mb{h}$, and we need to choose one from the set of feasible assignments $H$; the trick is identifying an appropriate function $f(\mathbf{h})$ in terms of the path flows, which is minimized when the path flows satisfy the principle of user equilibrium.
So, the constraint set of our optimization problem consists of the equations defining $H$:\label{not:Z2}
\begin{align}
\sum_{\pi \in \Pi^{rs}} h^\pi &= d^{rs} 					&		\forall (r,s) \in Z^2 \label{eqn:tapdemand} \\
h^\pi &\geq 0													&		\forall \pi \in \Pi \label{eqn:tapnonneg}
\,.
\end{align}

\index{optimization!convex!optimality conditions|(}
Different assignment rules will lead to different objective functions.
It turns out that the objective function corresponding to the principle of user equilibrium has a somewhat unintuitive form.
Therefore, we will derive this function in the same way that it was originally derived, by working backwards.
\emph{Rather than writing down the optimization problem, and then deriving the optimality conditions, we start by writing down the optimality conditions we want, then determining what kind of optimization problem has that form.}

\index{Lagrange multipliers!application}
Writing the (unknown!) objective function as $f(\mathbf{h})$, using the procedures described in Chapter~\ref{chp:mathematicalpreliminaries} we can Lagrangianize the constraints~\eqn{tapdemand} introducing multipliers $\kappa^{rs}$ for each OD pair, providing the following optimality conditions:
\begin{align}
\pdr{f}{h^\pi} - \kappa^{rs} & \geq 0 & \forall (r,s) \in Z^2, \pi \in \Pi^{rs} \label{eqn:minimum} \\ 
h^\pi \myp{\pdr{f}{h^\pi} - \kappa^{rs}} & =0 & \forall (r,s) \in Z^2, \pi \in \Pi^{rs} \label{eqn:complementarity} \\ 
\sum_{\pi \in \Pi^{rs}} h^\pi &= d^{rs} 					&		\forall (r,s) \in Z^2  \\ 
h^\pi & \geq 0 & \forall \pi \in \Pi  
\,.
\end{align}
The last two of these are simply~\eqn{tapdemand} and~\eqn{tapnonneg}, requiring that the solution be feasible.
The condition~\eqn{complementarity} is the most interesting of these, requiring that the product of each path's flow and another term involving $f$ must be zero.
For this product to be zero, either $h^\pi$ must be zero, or $\pdr{f}{h^\pi} = \kappa^{rs}$ must be true; and by~\eqn{minimum}, for all paths $\pdr{f}{h^\pi} \geq \kappa^{rs}$.
That is, in a solution to this problem, whenever $h^\pi$ is positive we have $\pdr{f}{h^\pi} = \kappa^{rs}$; if a path is unused, then $\pdr{f}{h^\pi} \geq \kappa^{rs}$ --- in other words, at optimality all used paths connecting OD pair $(r,s)$ must have equal and minimal $\pdr{f}{h^\pi}$, which is equal to $\kappa^{rs}$.
According to the principle of user equilibrium, all used paths have equal and minimal travel time\dots~so if we can choose a function $f(\mathbf{h})$ such that $\pdr{f}{h^\pi} = c^\pi$, we are done!  

So, the objective function must involve some integral of the travel times.
A first guess might look something like
\labeleqn{firstguess}{f(\mb{h}) = \sum_{\pi \in \Pi} \int_?^? c^\pi~dh^\pi\,,}
where the bounds of integration and other details are yet to be determined.
The trouble is that $c^\pi$ is not a function of $h^\pi$ alone: the travel time on a path depends on the travel times on other paths as well, so the partial derivative of this function with respect to $h^\pi$ will not simply be $c^\pi$, but contain other terms as well.
However, these interactions are not arbitrary, but instead occur where paths overlap, that is, where they share common links.
In fact, if we try writing a similar function to our guess but in terms of \emph{link flows}, instead of \emph{path flows}, we will be done.

To be more precise, let $\mb{x}(\mb{h})$ be the link flows as a function of the path flows $\mb{h}$, as determined by equation~\eqn{linkpath}.
Then the function
\labeleqn{secondguess}{f(\mb{h}) = \sum_{(i,j) \in A} \int_0^{x_{ij}(\mb{h})} t_{ij}(x)~dx}
satisfies our purposes.
To show this, calculate the partial derivative of $f$ with respect to the flow on an arbitrary path $\pi$, using the fundamental theorem of calculus and the chain rule:
\labeleqn{pdrder}{\pdr{f}{h^\pi} = \sum_{(i,j) \in A} t_{ij}(x_{ij}(\mb{h})) \pdr{x_{ij}}{h^\pi} = \sum_{(i,j) \in A} \delta_{ij}^\pi t_{ij}(x_{ij}(\mb{h})) = c^\pi\,,}
where the last two equalities respectively follow from differentiating~\eqn{linkpath} and from~\eqn{pathlinktimes}.

Finally, we can clean up the notation a bit by simply introducing the link flows $\mathbf{x}$ as a new set of decision variables, adding equations~\eqn{linkpath} as constraints to ensure they are consistent with the path flows.
This gives the following optimization problem:
\begin{align}
\min_{\mathbf{x},\mathbf{h}} \quad & \sum_{(i,j) \in A} \int_0^{x_{ij}} t_{ij}(x) dx &  \label{eqn:beckmann} \\
\index{Beckmann function}
\mathrm{s.t.} \quad & x_{ij} = \sum_{\pi \in \Pi} h^\pi \delta_{ij}^{\pi} & \forall (i,j) \in A \label{eqn:map} \\
& d^{rs} = \sum_{\pi \in \Pi^{rs}} h^\pi & \forall (r,s) \in Z^2 \label{eqn:demand} \\
& h^\pi \geq 0 & \forall \pi \in \Pi \label{eqn:nonneg} 
\,.
\end{align}
The objective function~\eqn{beckmann} is often called the \emph{Beckmann function}, named for Martin Beckmann who first reported this objective in 1956 (see historical notes in Section~\ref{sec:tap_historical}).

Section~\ref{sec:tapproperties} shows that the objective function is convex, and the feasible region is a convex set, so this problem will be relatively easy to solve.
\index{optimization!convex optimization!and traffic assignment|)}
\index{static traffic assignment!formulation!as convex optimization|)}
\index{optimization!convex!optimality conditions|)}

\subsection{Interpretation of the Beckmann function}
\label{sec:beckmannmeaning}

\index{Beckmann function!interpretation|(}Many readers wonder about the objective function~\eqn{beckmann} in the optimization problem used to represent the traffic assignment problem.
What does it mean?
Unlike the optimization problems formulated in Appendices~\ref{chp:basicoptimization} and~\ref{chp:fancyoptimization}, which aimed to minimize a cost, maximize an area, and so forth, in this case there is no obvious physical interpretation for summing integrals of the link performance functions.
There is no mathematical reason why an optimization problem needs to have a physical interpretation --- in this case, all that matters is that the optimality conditions correspond to user equilibrium, so we just need to solve the optimization problem somehow and the answer will give us what we're looking for.
Many writers in our field leave the explanation at that, but we feel this kind of answer is frankly unsatisfying.

The best way to understand what the Beckmann function accomplishes is to draw an analogy to physics.
The Beckmann function can be interpreted as a \emph{potential energy}\index{potential energy} for a given set of route choices.
As taught in elementary physics, potential energy has several properties:
\begin{enumerate}
\item It is a property of the configuration of a system: the arrangement of masses for gravitational potential energy, the arrangement of charges for electric potential energy, the displacement of a spring for elastic potential energy, and so forth.
\item The change in potential energy when moving the system from one state to another equals the work done.
Work is defined in terms of a line integral, integrating the dot product of force and displacement between the initial and final states.
\item Reversing the second property, we can represent force as a derivative of potential energy.
Using the typical sign convention for work, we specifically can say that \emph{the force at a point is the negative gradient of the potential energy function}.
\item Equilibria (in the physical sense) correspond to local minima or maxima of potential energy.
These are points where the gradient vanishes, and therefore there is zero force.
At a strict local minimum, the equilibrium is stable; at a strict local maximum, it is unstable.
For gravitational potential energy, a ball located at the bottom of a valley, or at the top of a hill, are respective examples of these types of physical equilibria.
\item At a non-equilibrium point, a force is applied in the direction of decreasing potential energy.
The forces that push a ball downhill, that represent gravitational or electric attraction, or that pull a spring back towards its neutral position are all of this form.
\end{enumerate}

The Beckmann function can be understood in the same way, with a suitable interpretation of ``force'' and ``energy.''
We have already defined what a user equilibrium is (all travelers are on shortest paths, so all used paths between the same origin and destination have equal and minimal travel time), and want to form an interpretation that agrees with the definition of a physical equilibrium in points 4 and 5 in the above list.
The ``configuration of the system'' is the assignment of travelers to paths, as represented by the path flow vector $\mb{h}$ and the corresponding path travel times $\mb{c}$.
From the previous section, we know that the minimum points of Beckmann function correspond to user equilibrium, so already we can think of the Beckmann function as a potential energy which is minimized at equilibrium.
The negative gradient of the potential energy function represents a force.
For the Beckmann function, equation~\eqn{pdrder} gives the partial derivatives with respect to the flow on each path.
Therefore, if we treat the Beckmann function as a potential, the associated force associated with path $\pi$ is $-c^\pi$.

The behavioral interpretation is that travelers are trying to minimize their travel time, so each path exerts a repulsive force on them equal to its travel time.
In moving a traveler from one path A to another B, the repulsive force from path A is $-c^A$, and that from path B is $-c^B$.
If path A has a higher travel time, the repulsion from path A overcomes that of path B, and the ``net work'' done is $c^A - c^B$, decreasing the potential of the system and moving to a lower-energy (more stable) state.
On the other hand, if path B has the higher travel time, the repulsive force from path B is stronger than from path A, so moving the travelers from A to B requires us to ``work against'' the travelers' natural tendency, introducing more energy to the system and moving it to a less stable state.
Equilibrium occurs when the system cannot be moved to a lower-energy state: potential energy (the Beckmann function) is minimized, and no travelers can switch paths to reduce their individual travel time.

(This interpretation has entirely been about \emph{stable} equilibria.
It turns out that for the basic traffic assignment problem, as formulated in this chapter, the user equilibrium solution is always stable as long as the link performance functions are strictly increasing functions.
Later in the book, we discuss more complex situations with unstable user equilibria.
In some cases, we can view such equilibria as being local \emph{maxima} of a potential function.)

We emphasize that we are not using this interpretation to \emph{define} what a user equilibrium is: there is no ``real force'' pushing people to different paths.
Rather, the aim is to draw an analogy between the user equilibrium concept of route choice, and mechanical equilibrium of physical bodies, in hopes that it provides insight and intuition into the Beckmann function.
Like all such analogies, if it's not helpful to you as a reader, feel free to ignore it.
The critical feature of the Beckmann function is that its minimum points correspond to user equilibrium solutions, as shown by analyzing the optimality conditions.\index{Beckmann function!interpretation|)}

\subsection{Multifunctions and application to network equilibrium (*)}
\label{sec:eqmmultifunction}

\emph{(This section is optional and can be skipped.
However if you are curious about the ``other road'' to using fixed points to show equilibrium existence, read this section to learn about multifunctions and Kakutani's theorem.)}

\index{fixed point problem!and traffic assignment|(}
\index{static traffic assignment!formulation!as fixed point|(}
In the transit ridership example of Section~\ref{sec:fixedpoint}, we could formulate the solution to the problem as a fixed point problem directly, without first going through a variational inequality.
The traffic assignment problem is slightly more complex, because at equilibrium all used paths will have equal travel time.
So, if we try to apply the same technique as in the transit ridership problem, we will run into a difficulty with ``breaking ties.''  To see this concretely, consider the two-path network in Figure~\ref{fig:2linksimple}, and let $\mb{h}$ and $\mb{c}$ be the vectors of path flows and path travel times.
As Figure~\ref{fig:taschematic} suggests, we can try to define two functions: $\mb{H}(\mb{c})$\label{not:Hc} gives the vector of path flows representing path choices \emph{if the travel times were fixed at $\mb{c}$}, and $\mb{C}(\mb{h})$\label{not:Ch} gives the vector of path travel times when the path flows are $\mb{h}$.
The function $\mb{C}$ is clearly defined: simply calculate the cost on each path using the link performance functions.
However, if both paths are used at equilibrium, then $c_1 = c_2$, which means that \emph{any} path flow vector would be a valid choice for $\mb{H}(\mb{c})$.
Be sure you understand this difficulty: because $\mb{H}$ assumes that the travel times are fixed, if they are equal any path flow vector is consistent with our route choice assumptions.
If we relax the assumption that the link performance functions are fixed, then $\mb{H}$ is no simpler than solving the equilibrium problem in the first place.

\begin{figure}
\begin{center}
\begin{tikzpicture}[->,>=stealth',shorten >=1pt,auto,node distance=3cm,
thick,main node/.style={circle,,draw}]
\node[main node] (1) {$1$};
\node[main node] (2) [right of=1] {$2$};
\path (1) edge [bend left] node[above] {$50$} (2);
\path (1) edge [bend right] node[below] {$45 + x$} (2);
\end{tikzpicture}
$d^{12} = 30$
\end{center}
\caption{A simple two-link network for demonstration.\label{fig:2linksimple}}
\end{figure}

\index{correspondence|see {multifunction}}
\index{point-to-set map|see {multifunction}}
\index{set-valued map|see {multifunction}}
\index{set-valued function|see {multifunction}}
To resolve this, we introduce the concept of a \emph{multifunction}\index{multifunction}.\footnote{Also known by many other names, including \emph{correspondence}, \emph{point-to-set map}, \emph{set-valued map}, or \emph{set-valued function}.}  Recall that a regular function from $X$ to $Y$ associates each value in its domain $X$ with exactly one value in the set $Y$.
A multifunction, on the other hand, associates each value in the domain $X$ with some subset of $Y$.
If $F$\label{not:Fmulti} is such a multifunction the notation $F : X \rightrightarrows Y$\label{not:fXmY} can be used.
For example, consider the multifunction $\mb{H}(\mb{c}) : \bbr^2 \rightrightarrows H$ defined to reflect all path flows which are consistent with a given vector of path costs.
If $c_1 < c_2$ (and these travel times were fixed), then the only consistent path flows are to have everybody on path 1.
Likewise, if $c_1 > c_2$, then everyone would have to be on path 2.
Finally, if $c_1 = c_2$ then people could split in any proportion while satisfying the rule that drivers are using least-time paths.
That is, 
\labeleqn{multifn}{\mb{H}(c_1, c_2) = \begin{cases}
                                    \myc{ \minivect{30 & 0} }                    & c_1 < c_2 \\
                                    \myc{ \minivect{0 & 30} }                    & c_1 > c_2 \\
                                    \myc{ \minivect{h & 30 - h} : h \in [0,30]}  & c_1 = c_2
                                 \end{cases}\,. }
The notation in equation~\eqn{multifn} is chosen very carefully to reflect the fact that $\mb{H}$ is a multifunction, which means that its values are \emph{sets}, not a specific number.
In the first two cases, the set only consists of a single element, so this distinction may seem a bit pedantic; but in the latter case, the set contains a whole range of possible values.
In this way, multifunctions generalize the concept of a function by allowing $\mb{H}$ to take multiple ``output'' values for a single input.
This is graphically represented in Figure~\ref{fig:eqmmultifunction}, using the fact that the function can be parameterized in terms of a single variable $c_2 - c_1$.

\stevefig{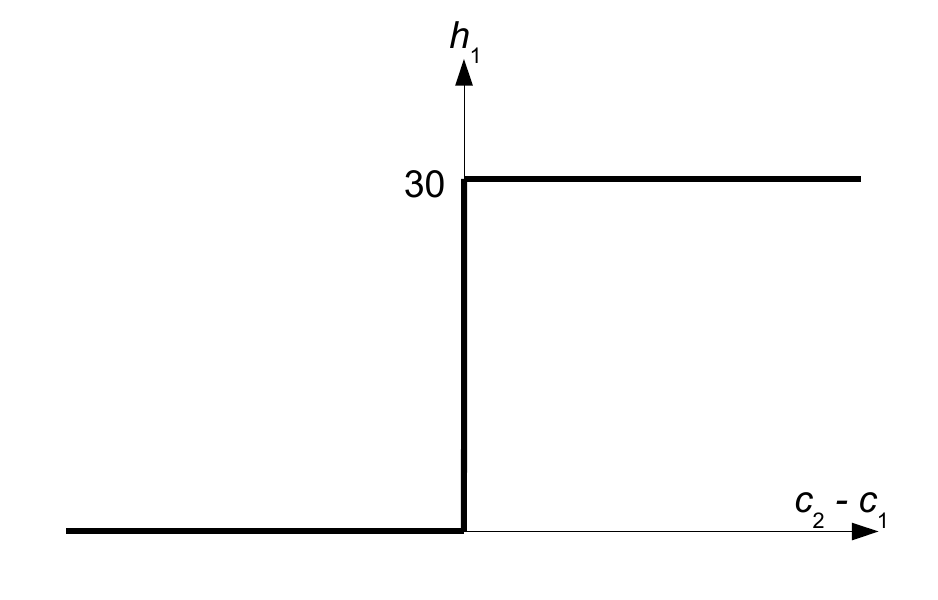}{Multifunction corresponding to the two-link network; notice multiple values when $c_1 = c_2$.}{0.5\textwidth}

\index{fixed point problem!for multifunctions}
Fixed point problems can be formulated for multifunctions as well as for functions.
If $F$ is an arbitrary multifunction defined on the set $K$, then a fixed point of $F$ is a point $x$ such that $x \in F(x)$.
Note the use of set inclusion $\in$ rather than equality $=$ because $F$ can associate multiple values with $x$.
Just as Brouwer's theorem guarantees existence of fixed points for functions, Kakutani's theorem guarantees existence of fixed points for multifunctions under general conditions:
\begin{thm} \emph{(Kakutani).}
\index{Kakutani's theorem}
\index{convex set!applications}
\index{compact set!applications}
Let $F: K \rightrightarrows K$ be a multifunction (from the set $K$ to itself), where $K$ is convex and compact.
If $F$ has a closed graph and $F(x)$ is nonempty and convex for all $x \in K$, then there is at least one point $x \in K$ such that $x \in F(x)$.
\end{thm}
The new terminology here is a \emph{closed graph}\index{closed graph}; we say that the multifunction $F$ has a closed graph if the set $\{ (\mb{x}, \mb{y}) : \mb{x} \in X, \mb{y} \in F(\mb{x}) \}$ is closed.
Again, each of these conditions is necessary; you might find it helpful to visualize these conditions geometrically similar to Figure~\ref{fig:fixedpointgraph}.

As a result of Kakutani's theorem, we know that an equilibrium solution must always exist in the traffic assignment problem.
Let $H$ denote the set of all feasible path flow vectors, that is, vectors $\mb{h}$ such that $d^{rs} = \sum_{\pi \in \Pi^{rs}} h^\pi$ for all OD pairs $(r,s)$ and $h^\pi \geq 0$ for all $\pi \in \Pi$.
Let the function $\mathbf{C}(\mathbf{h})$ represent the path travel times when the path flows are $\mb{h}$.\footnote{If you look back to Chapter~\ref{chp:equilibrium}, you will see that $\mathbf{C}(\mathbf{h}) = \bm{\Delta}^T \mathbf{t}(\bm{\Delta} \mathbf{h})$ in matrix notation.}  Let the multifunction  $\mb{H}(\mathbf{c})$ represent the set of path flows which could possibly occur if all travelers chose least cost paths given travel times $\mathbf{c}$.
Mathematically 
\labeleqn{hmapping}{
   \mb{H}(\mb{c}) = \myc{ \mb{h} \in H: h^\pi > 0 \mbox{ only if } c^\pi = \min_{\pi' \in \Pi^{rs}} c^{\pi'} \mbox{ for all } (r,s) \in Z^2, \pi \in \Pi^{rs}}
\,.
}
This is the generalization of equation~\eqn{multifn} when there are multiple OD pairs and multiple paths connecting each OD pair --- be sure that you can see how \eqn{multifn} is a special case of~\eqn{hmapping} when there is just one OD pair connected by two paths.

\begin{thm}
\label{thm:eqmkakutani}
\index{continuous function!applications}
\index{user equilibrium!existence in static}
\index{static traffic assignment!properties!existence}
If the link performance functions $t_{ij}$ are continuous, then at least one equilibrium solution exists to the traffic assignment problem.
\end{thm}
\begin{proof}
Let the multifunction $F$ be the composition of the multifunction $\mb{H}$ and the function $\mathbf{C}$ defined above, so $F(\mathbf{h}) = \mb{H}(\mathbf{C}(\mathbf{h}))$ is the set of path flow choices which are consistent with drivers choosing fastest paths when the travel times correspond to path flows $\mathbf{h}$.
An equilibrium path flow vector is a fixed point of $F$: if $\mb{h} \in F(\mb{h})$ then the path flow choices and travel times are consistent with each other.
The set $H$ of feasible $\mathbf{h}$ is convex and compact by Lemma~\ref{lem:compactconvexpath} in the next section.
Examining the definition of $\mb{H}$ in~\eqn{hmapping}, we see that for any vector $\mb{c}$, $\mb{H}(\mb{c})$ is nonempty (since there is at least one path of minimum cost), and compact by the same argument used in the proof of Lemma~\ref{lem:compactconvexpath}.
Finally, the graph of $F$ is closed, as you are asked to show in the exercises.
\end{proof}
\index{fixed point problem!and traffic assignment|)}
\index{static traffic assignment!formulation!as fixed point|)}

\section{Properties of User Equilibrium}
\label{sec:tapproperties}

This section uses the mathematical formulations from the previous section to explore basic properties of user equilibrium solutions.
In Section~\ref{sec:existenceuniqueness} we give fairly general conditions under which there is one (and only one) user equilibrium link flow solution in a network.
Although this section is mostly mathematical, these properties are actually of great practical importance.
We ultimately want to use traffic assignment as a tool to help evaluate and rank transportation projects, using the principle of user equilibrium to predict network conditions.
If it is possible that there is no feasible assignment which satisfies the principle of user equilibrium, then we would be at a loss as to how projects should be ranked, and we would need to find another assignment rule.
Or, if there could be multiple feasible assignments which satisfy the principle of user equilibrium, again it is unclear how projects should be ranked.

Interestingly, while the conditions for a unique \emph{link flow}\index{static traffic assignment!link flow solution} user equilibrium are fairly mild, there will almost certainly be multiple \emph{path flow}\index{static traffic assignment!path flow solution} solutions which satisfy the principle of equilibrium --- infinitely many, in fact.
As explained in Section~\ref{sec:maxentropy}, this suggests that the principle of user equilibrium is not strong enough to identify path flows, simply link flows.
If path flows are needed to evaluate a project, an alternative approach is needed, and this section explains the concepts of maximum entropy and proportionality which provide this alternative.

Lastly, when solving equilibrium problems on large networks both the link flow and path flow representations have significant limitations --- algorithms only using link flows tend to be slow, while algorithms only using path flows can require a large amount of memory and tend to return low-entropy path flow solutions.
A compromise is to use the link flows, but distinguish the flow on each link by its origin or destination.
While explored more in Chapter~\ref{chp:solutionalgorithms}, Section~\ref{sec:bushes} lays the groundwork by showing how flow can be decomposed this way, and that at equilibrium the links with positive flow from each origin or destination form an acyclic network.

\subsection{Existence and link flow uniqueness}
\label{sec:existenceuniqueness}

We would like to use the mathematical formulations from the previous section --- variational inequality, fixed point, and convex optimization --- to identify conditions under which one (and only one) user equilibrium solution exists.
This section relies heavily on mathematical formulations and results proved previously; interested readers may refer to those sections for more detail, and others can still read the proofs here to obtain the general idea.
Recall the following results from Chapter~\ref{chp:mathematicalpreliminaries}:
\begin{itemize}
\item \textbf{Brouwer's theorem~\ref{thm:brouwersthm}: } Let $K$ be a compact and convex set, and let $f : K \rightarrow K$ be a continuous function.
Then there is at least one fixed point $x \in K$, that is, a point such that $f(x) = x$.
\item \textbf{Proposition~\ref{prp:strictconvexuniqueness}: } Let $K$ be a convex set, and let $f : K \rightarrow \bbr$ be a strictly convex function.
Then there is a unique global minimum $\hat{x}$ of $f$ on $K$.
\end{itemize}

Since both of these results will rely on some properties of the sets of feasible path flows $H$ and feasible link flows $X$, we establish them first:

\begin{lem}
\label{lem:compactconvexpath}
\index{feasible assignment!properties}
\index{compact set!applications}
\index{convex set!applications}
The set of feasible path assignments $H$ and the set of feasible link assignments $X$ are both compact and convex.
\end{lem}
\begin{proof}
 To show that $H$ is compact, we must show that it is closed and bounded.
The set $H$ is defined by a combination of linear weak inequalities ($h^\pi \geq 0$) and linear equalities ($\sum h^\pi = d^{rs}$), so together Propositions~\ref{prp:setfacts1}b and~\ref{prp:setfacts2}a--b show that it is closed.
For boundedness, consider any $\mb{h} \in H$, and OD pair $(r,s)$, and any path $\pi \in \Pi^{rs}$.
We must have $h^\pi \leq d^{rs}$, since each $h^\pi$ is nonnegative and $d^{rs} = \sum_{\pi \in \Pi^{rs}} h^\pi$.
Therefore, if $D$\label{not:Dmax} is the largest entry in the OD matrix, we have $h^\pi \leq D$ for any path $\pi$.
Therefore $|\mb{h}| \leq \sqrt{|\Pi|} D$, so $H$ is contained in the ball $B_{ \sqrt{|\Pi|} D} (\mb{0})$\label{not:B} and $H$ is bounded.
 
Finally, for convexity, consider any $\mb{h_1} \in H$, any $\mb{h_2} \in H$, any $\lambda \in [0, 1]$, and the resulting vector $\mb{h} = \lambda \mb{h_1} + (1 - \lambda) \mb{h_2}$.
For any path $\pi$, $h_1^\pi \geq 0$ and $h_2^\pi \geq 0$, so $h^\pi = \lambda h_1^\pi + (1 - \lambda) h_2^\pi \geq 0$ as well.
Furthermore, for any OD pair $(r,s)$ we have
\begin{align*}
\sum_{(r,s) \in \Pi^{rs}} h^\pi &= \sum_{(r,s) \in \Pi^{rs}} [\lambda h_1^\pi + (1 - \lambda) h_2^\pi ] \\
                             &= \lambda \sum_{(r,s) \in \Pi^{rs}} h_1^\pi + (1 - \lambda) \sum_{(r,s) \in \Pi^{rs}} h_2^\pi \\
                             &= \lambda d^{rs} + (1 - \lambda) d^{rs} \\
                             &= d^{rs}
\,,
\end{align*}
so $\mb{h} \in H$ as well.
 
Every feasible link assignment $\mb{x} \in X$ is obtained from a linear transformation of some $\mb{h} \in H$ by~\eqn{map} so $X$ is also closed, bounded, and convex.
\end{proof}

To show existence of solutions, we use Brouwer's theorem based on the formulation of user equilibrium as a fixed point of the function $f(\mb{h}) = \mr{proj}_H (\mb{h} - \mb{c}(\mb{h}))$.

\begin{prp}
\label{prp:existence}
If the link performance functions $t_{ij}$ are continuous for each link $(i,j) \in A$, then there is a feasible assignment satisfying the principle of user equilibrium.
\index{static traffic assignment!properties!existence}
\index{link performance function!properties}
\index{continuous function!applications}
\index{user equilibrium!existence in static}
\end{prp}
\begin{proof}
\index{projection!applications}
Taking $H$ as the set of feasible assignments, the range of the function $f(\mb{h}) = \mr{proj}_H (\mb{h} - \mb{c}(\mb{h}))$ clearly lies in $H$ because of the projection.
By Lemma~\ref{lem:compactconvexpath}, $H$ is compact and convex, so if we show that $f$ is continuous, then Brouwer's Theorem guarantees existence of a fixed point.
The discussion in Section~\ref{sec:mathematicalformulations} showed that each fixed point is a user equilibrium solution, which will be enough to prove the result.

We use the result that the composition of continuous functions is continuous (Proposition~\ref{prp:continuousdifferentiable}).
The function $f$ is the composition of two other functions (call them $f_1$ and $f_2$), where $f_1(\mb{h}) = \mr{proj}_H(\mb{h})$ and $f_2(\mb{h}) = \mb{h} - \mb{c}(\mb{h})$.
By Proposition~\ref{prp:continuousprojection}, $f_1$ is continuous because $H$ is a convex set.
Furthermore, $f_2$ is continuous if $\mb{c}$ is a continuous function of $\mb{h}$, which is true by hypothesis.
Therefore, the conditions of Brouwer's Theorem are satisfied, and $f$ has at least one fixed point (which satisfies the principle of user equilibrium).
\end{proof}

The easiest way to show uniqueness of solutions is to make use of the convex optimization formulation, minimizing the function $f(\mb{x}) = \sum_{ij} \int_0^{x_{ij}} t_{ij}(x)~dx$ over the set $X$.

\begin{prp}
\label{prp:linkuniqueness}
If the link performance functions $t_{ij}$ are differentiable, and $t'_{ij}(x) > 0$ for all $x$ and links $(i,j)$, there is exactly one feasible link assignment satisfying the principle of user equilibrium.
\index{static traffic assignment!properties!uniqueness in link flows}
\index{link performance function!properties}
\index{monotone function!applications}
\index{user equilibrium!uniqueness in static}
\end{prp}
\begin{proof}
The previous section showed that user equilibrium link flow solutions correspond to minimum points of $f$ on $X$.
If we can show that $f$ is strictly convex, then there can only be one minimum point.
(Differentiability implies continuity, so we can assume that at least one minimum point exists by Proposition~\ref{prp:existence}).
Since $f$ is a function of multiple variables (each link's flow), to show that $f$ is convex we can write its Hessian matrix of second partial derivatives.
The first partial derivatives take the form $\pdr{f}{x_{ij}} = t_{ij}(x_{ij})$.
So, the diagonal entries of the Hessian\index{Hessian matrix!applications} take the form $\pdrb{f}{x_{ij}} = t'_{ij}(x_{ij})$, while the off-diagonal entries take the form $\pdrc{f}{x_{ij}}{x_{k\ell}} = 0$.\label{not:kl}
Since the Hessian is a diagonal matrix, and its diagonal entries are strictly positive by assumption, $f$ is a strictly convex function.
Therefore there is only one feasible link assignment satisfying the principle of user equilibrium.
\end{proof}
It is also possible to prove the same result under the slightly weaker condition that the link performance functions are continuous and strictly increasing; this is undertaken in the exercises.

To summarize the above discussion, if the link performance functions are continuous, then at least one user equilibrium solution exists; if in addition the link performance functions are strictly increasing, then \emph{exactly} one user equilibrium solution exists.
For many problems representing automobile traffic, these conditions seem fairly reasonable, if not universally so: adding one more vehicle to the road is likely to increase the delay slightly, but not dramatically.

\subsection{Path flow nonuniqueness, entropy, and proportionality}
\label{sec:maxentropy}

\index{static traffic assignment!properties!nonuniqueness in path flows|(}
The previous subsection gave a relatively mild condition for the equilibrium solution to be unique in terms of link flows.
However, a more subtle point is that the equilibrium solution need not be unique in terms of path flows.
Consider the example shown in Figure~\ref{fig:nonunique}.
With the paths as numbered in the figure, the reader can verify that both $\mathbf{h_1} = \vect{20 & 10 & 20 & 10}$ and $\mathbf{h_2} = \vect{30 & 0 & 10 & 20}$ produce the equilibrium link flows, and therefore satisfy the principle of user equilibrium (all paths having equal travel time).
Furthermore, any weighted average of $\mathbf{h_1}$ and $\mathbf{h_2}$ also produces the same link flows, so there are an infinite number of feasible path flow solutions which satisfy the principle of user equilibrium.

\stevefig{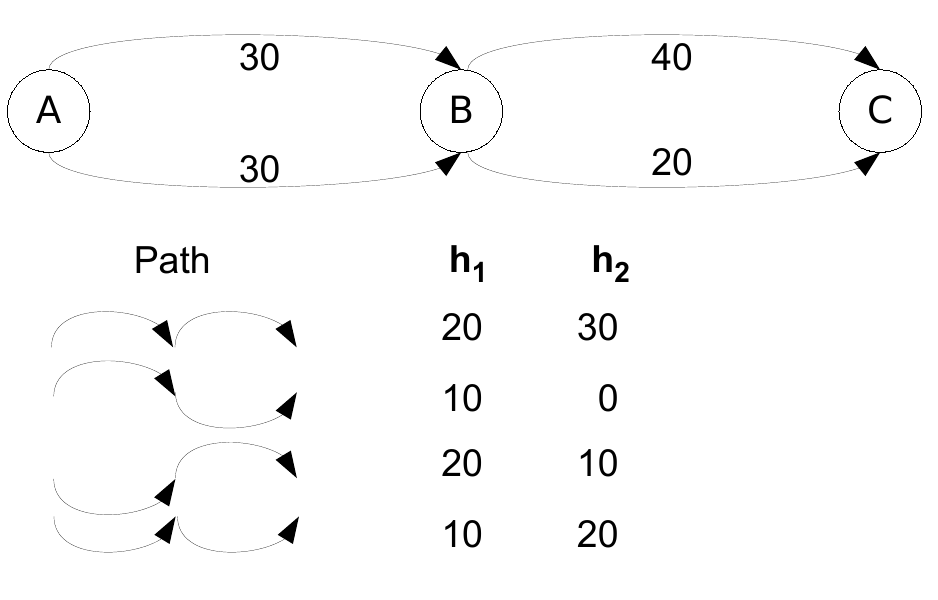}{Nonuniqueness of equilibrium in terms of path flows}{0.7\textwidth}

At first glance this may appear to contradict the proof of Proposition~\ref{prp:linkuniqueness}, which showed that the function $f$ was strictly convex\index{convex function!strict convexity!applications}.
Although $f$ is strictly convex \emph{as a function of $\mathbf{x}$}, it is \emph{not} strictly convex as a function of $\mathbf{h}$.
This is easy to see mathematically --- if there are two distinct path flow solutions $\mathbf{h_1}$ and $\mathbf{h_2}$ which satisfy the principle of user equilibrium, then $\lambda \mathbf{h_1} + (1 - \lambda) \mathbf{h_2}$ is an equilibrium as well for $\lambda \in (0, 1)$.
Since all three are equilibria, $f(\mathbf{h_1}) = f(\mathbf{h_2}) = f(\lambda \mathbf{h_1} + (1 - \lambda) \mathbf{h_2})$ even though strict convexity would require $f(\lambda \mathbf{h_1} + (1 - \lambda) \mathbf{h_2}) < \lambda f(\mathbf{h_1}) + (1 - \lambda) f(\mathbf{h_2}) = f(\mathbf{h_1})$.
From a more intuitive standpoint, in typical networks there are many more paths than links; therefore, it is very likely that multiple path flow solutions correspond to the same link flow solution.
Since path travel times only depend on the values of the link flows (because link flows determine link travel times, which are added up to get path travel times), and since the principle of user equilibrium is defined in terms of path travel times, it is reasonable to expect multiple path flow solutions to satisfy the equilibrium principle.

The existence of multiple path flow equilibria (despite a unique link flow equilibrium) is not simply a technical curiosity.
Rather, it plays an important role in using network models to evaluate and rank alternatives.
Thus far, we've been content to ask ourselves what the equilibrium link flows are, but in practice these numbers are usually used to generate other, more interesting measures of effectiveness such as the total system travel time and total vehicle-miles traveled.
If a neighborhood group is worried about increased traffic, link flows can support this type of analysis.
If we are concerned about safety, we can use link flows as one input in estimating crash frequency, and so forth.

Other important metrics, however, require more information than just the total number of vehicles on a link.
Instead, we must know the entire \emph{paths} used by drivers.
Examples include
\begin{description}
\item[Select link analysis:\index{select link analysis} ] In transportation planning parlance, a ``select link analysis'' goes beyond asking ``how many vehicles are on a link,'' to ``exactly which vehicles are using this link?'' This is used to produce nice visualizations, as well as to identify which origins and destinations are using a particular link.
This is needed to identify which neighborhoods and regions are affected by transportation improvements, both positively and negatively.

\item[Equity\index{equity} and environmental justice\index{environmental justice}: ] As a special case of the above, planners are focusing more and more on the notion of equity.
Most transportation projects involve both ``winners'' and ``losers.''  For instance, constructing a new freeway through an existing neighborhood may improve regional mobility, but cause significant harm to those whose neighborhood is disrupted in the process.
Historically, disadvantaged populations tend to bear the brunt of such projects.
To know who is benefiting or suffering from a project, we need to know the origins and destinations of the travelers on different links.

\item[Emissions\index{emissions}: ] Vehicles emit more pollutants toward the start of trips, when the engine is cold, than later on once the catalytic converter has warmed up.
Therefore, to use the output of a network model to predict emissions from vehicles on a link, ideally we need to distinguish between vehicles which have just started their trips from vehicles which have already traveled some distance.
This requires knowing path flows.

\item[Pavement loading\index{pavement management}: ] Heavy vehicles cause much more damage to pavement than lighter vehicles; a commonly-used relation in pavement engineering estimates that pavement damage is proportional to the \emph{fourth power} of vehicle weight.
If a network model is being used to forecast pavement deterioration, it is therefore important to know which links are being used by heavy vehicles.
Again, this requires more than just the total flow on a link; we must know what specific type of vehicle it is (which often depends heavily on the origin-destination pair).
\end{description}

In fact, software vendors often discuss such abilities as key features of their applications.\index{traffic assignment!software}
In this light, the fact that equilibrium path flows are nonunique should trouble you somewhat.
If there are multiple path flow equilibria (and therefore multiple values for emissions forecasts, equity analyses, etc.), which one is the ``right'' one?  Simply picking one possible path flow solution and hoping that it's the right one is not particularly rigorous, and inappropriate for serious planning applications.

So, the principle of user equilibrium is only strong enough to determine link flows, not path flows, even though having a path flow solution would be very helpful in analyzing engineering alternatives.
This has motivated the search for the ``most likely path flow'' solution which satisfies equilibrium.
Even though user equilibrium is not a strong enough concept to determine path flows uniquely, we can still ask ``of all the path flow vectors which satisfy the principle of user equilibrium, which do we think is most likely to occur?''  This section describes the concepts most commonly used to this end.
To my knowledge, there has been little to no research directly validating this principle, but nobody has come up with a better, generally agreed-upon solution to this issue.

\index{entropy|(}
A physical analogy is presented to motivate the idea.
Figure~\ref{fig:gas} shows the locations of $n$\label{not:nn} gas molecules within a box.
In the left panel, all $n$ molecules happen to be located in the top half of the box, while in the right panel $n/2$ molecules are located in the top half and $n/2$ in the bottom half.
Both of these situations are physically possible (that is, they satisfy the laws of mechanics, etc.) although you would certainly think that the scenario on the right is more likely or plausible than the scenario on the left.
In exactly the same way, although there are many path flow solutions which satisfy the ``law'' of user equilibrium, we are not forced to conclude that all such solutions are equally likely to occur.

\genfig{gas}{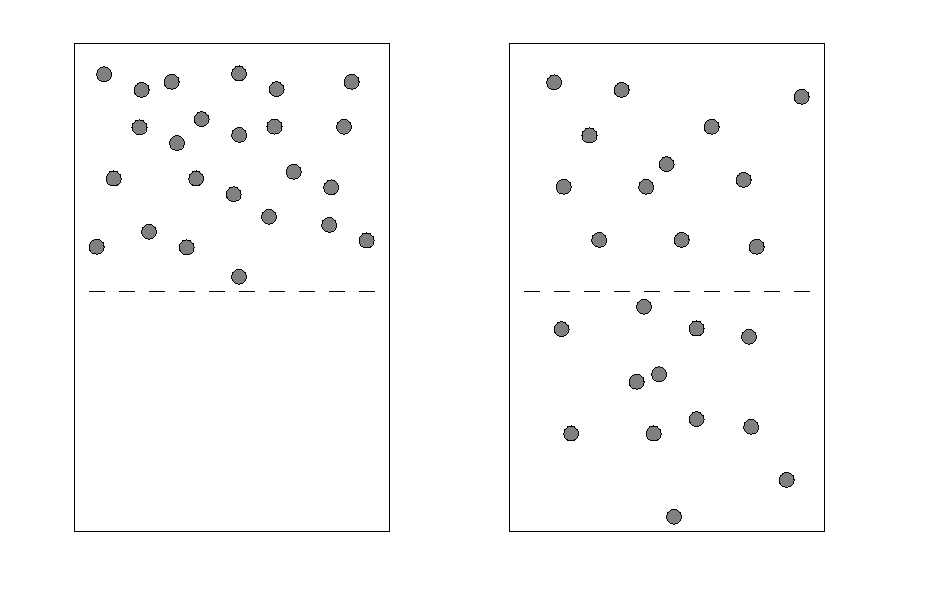}{Molecules of an ideal gas in a box. Which situation is more likely?}{width=0.5\textwidth}

Why does the scenario on the left seem so unlikely?  If the molecules form an ideal gas, the location of each molecule is independent of the location of every other molecule.
Therefore, the probability that any given molecule is in the top half of the box is $1/2$, and the probability that \emph{all} $n$ molecules are as in the left figure is\label{not:p}
\labeleqn{tophalf}{p^L = \myp{ \frac{1}{2}}^n \,.}
To find the probability that the distribution of molecules is as in the right figure ($p^R$), we use the binomial distribution: 
\labeleqn{binomial}{p^R = \binom{n}{n/2} \myp{\frac{1}{2}}^{n/2} \myp{\frac{1}{2}}^{n/2} = \frac{n!}{(n/2)! (n/2)!} \myp{\frac{1}{2}}^n \,.}
Since $n! \gg [(n/2)!]^2$\label{not:fact}\label{not:gg}, we have $p^R \gg p^L$, that is, the situation on the right is much, much likelier than that on the left even though both are physically possible.

We want to use the same principle for path flows.
The principle of user equilibrium assumes that each user is identical (in the sense that they choose paths only according to minimal travel time) and chooses their path independent of other drivers (except insofar as other drivers' choices affect travel times, drivers do not coordinate and strategize about their route choices with each other).
Therefore, we can apply the same logic as with the gas molecules.
For now, assume there is a single OD pair, with integer demand $d$, and let $\hat{\Pi} = \myc{\pi_1, \cdots, \pi_k} $ denote the set of minimum travel-time paths at the equilibrium solution.
For simplicity assume that these paths have constant travel times independent of flow.
Since drivers only care about travel time, each path is essentially identical, and the probability any particular driver chooses any particular path is $1/k$.
The probability that the path flow vector $\mathbf{h}$ takes a particular value $\vect{h_1 & \cdots & h_k}$ is then given by the multinomial distribution:
\labeleqn{multinomial}{p = \frac{d!}{h_1! h_2! \cdots h_k!} \myp{\frac{1}{k}}^d\,.}
and the most likely path flow vector is the one which maximizes this product.
Since $(1/k)^d$ is a constant, the most likely path flow vector simply maximizes
\labeleqn{entropyfactorial}{p = \frac{d!}{h_1! h_2! \cdots h_k!}\,.}
Now we introduce a common trick: since the logarithm function is strictly increasing, the path flows which maximize $p$ also maximize $\log p$,\label{not:log} which is
\labeleqn{entropylog}{\log p = \log d! - \sum_{\pi = 1}^k \log h_\pi!\,.}
To simplify further, we use Stirling's approximation\index{Stirling's approximation}, which states that $n! \approx n^n e^{-n} \sqrt{2 \pi n}$, or equivalently $\log n! \approx n \log n - n + (1/2) \log (2 \pi n)$.
This approximation is asymptotically exact in the sense that the ratio between $n!$ and $n^n e^{-n} \sqrt{2 \pi n}$ approaches 1 as $n$ grows large.
Further, when $n$ is large, $n$ is much larger than $\log n$, so the last term can be safely ignored and $\log n! \approx n \log n - n$.
Substituting into~\eqn{entropylog} we obtain
\labeleqn{entropystirling}{\log p \approx (d \log d - d) - \sum_\pi (h_\pi \log h_\pi - h_\pi)\,.}
Since $d = \sum_\pi h_\pi$, we can manipulate~\eqn{entropystirling} to obtain
\labeleqn{entropysingleod}{\log p \approx - \sum_\pi h_\pi \log (h_\pi / d)\,,}
and the path flows $h_1, \ldots, h_k$ maximizing this quantity (subject to the constraint $d = \sum_\pi h_\pi$ and nonnegativity) approximately maximize the probability that this particular path flow vector will occur in the field if travelers are choosing routes independent of each other.

To move towards the traffic assignment problem we're used to, we need to make the following changes:
\begin{enumerate}
\item In the traffic assignment problem, the demand $d$ and path flows $h_\pi$ do not need to be integers, but can instead take on real values.
The corresponding interpretation is that the demand is ``divided'' into smaller and smaller units, each of which is assumed to act independently of every other unit.
(This is how user equilibrium works with continuous flows.)  This doesn't cause any problems, and in fact helps us out: as we take the limit towards infinite divisibility, Stirling's approximation becomes exact --- so we can replace the $\approx$ in our formulas with exact equality.

\item There are multiple origins and destinations.
This doesn't change the basic idea, the formulas just become incrementally more complex as we sum equations of the form~\eqn{entropysingleod} for each OD pair.

\item Travel times are not constant, but are instead flow dependent.
Again, this doesn't change any basic ideas; we just have to add a constraint stating that the path flows correspond to an equilibrium solution.
Since the equilibrium link flows are unique, we simply have to state that the path flows correspond to the equilibrium link flows $\mathbf{\hat{x}}$.
\end{enumerate}
Putting all this together, we seek the path flows $\mathbf{h}$ which solve the optimization problem
\begin{align}
	\max_\mathbf{h} \qquad & -\sum_{(r,s) \in Z^2} \sum_{\pi \in \hat{\Pi}^{rs}} h_\pi \log (h_\pi / d^{rs}) & \label{eqn:entropy} \\
	\mathrm{s.t.} \qquad   & \sum_{\pi \in \Pi} \delta_{ij}^\pi h_\pi = \hat{x}_{ij}				& \forall (i,j) \in A \label{eqn:entropyeqm} \\
								  & \sum_{\pi \in \hat{\Pi}^{rs}} h_\pi = d^{rs}	& \forall (r,s) \in Z^2 \label{eqn:entropynvlb} \\
								  & h_\pi \geq 0												& \forall \pi \in \Pi \label{eqn:entropynonneg} 
\,.
\end{align}
The objective function~\eqn{entropy} is called the \emph{entropy} of a particular path flow solution, and the constraints~\eqn{entropynvlb},~\eqn{entropyeqm}, and~\eqn{entropynonneg} respectively require that the path flows are consistent with the OD matrix, equilibrium link flows, and nonnegativity.
The set $\hat{\Pi}^{rs}$ is the set of paths which are used by OD pair $(r,s)$.

The word \emph{entropy} is meant to be suggestive.
The connection with the thermodynamic concept of entropy may be apparent from the physical analogy at the start of this section.\footnote{Actually, the term here is more directly drawn from the fascinating field of information theory, which is unfortunately beyond the scope of this text.}  In physics, entropy can be interpreted as a measure of disorder in a system.
Both the left and right panels in Figure~\ref{fig:gas} are ``allowable'' in the sense that they obey the laws of physics.
However, the scenario in the right has much less structure and much higher entropy, and is therefore more likely to occur.

\index{proportionality!see {static traffic assignment, proportionality}}
It is not trivial to solve the optimization problem~\eqn{entropy}--\eqn{entropynonneg}.
It turns out that entropy maximization implies a much simpler condition, proportionality\index{static traffic assignment!proportionality}\index{entropy!and proportionality} among pairs of alternate segments.
Any two paths with a common origin and destination imply one or more pairs of alternate segments\index{pair of alternate segments} where the paths diverge, separated by links common to both paths.
In Figure~\ref{fig:nonunique}, there are two pairs of alternate segments: the top and bottom links between nodes A and B, and the top and bottom links between B and C.
It turns out that path flows $\mathbf{h_1}$ are the solution to the entropy maximization problem.

You might notice some regularity or patterns in this solution.
For instance, looking only at the top and bottom links between A and B, at equilibrium these links have equal flow (30 vehicles each).
Paths 1 and 3 are the same except for between nodes A and B, and these paths also have equal flow (20 vehicles each).
The same is true for paths 2 and 4 (10 vehicles each).
Or, more subtly, between B and C, the ratio of flows between the top and bottom link is 2:1, the same as the ratio of flows on paths 1 and 2 (which only differ between B and C), and the ratio of flows on paths 3 and 4.
This is no coincidence; in fact, we can show that \emph{entropy maximization implies proportionality}.\index{pair of alternate segments!and proportionality}

Before proving this fact, we show that it holds even for different OD pairs.
The network in Figure~\ref{fig:multiOD} has 40 vehicles traveling from origin A to destination F, and 120 vehicles from B to E, and the equilibrium link flows are shown in the figure.
Each OD pair has two paths available to it, one using the top link between C and D, and the other using the bottom link between C and D.
Using an upward-pointing arrow to denote the first type of path, and a downward-pointing arrow to describe the second, the four paths are $h^\uparrow_{AF}$, $h^\downarrow_{AF}$, $h^\uparrow_{BE}$, and $h^\downarrow_{BE}$.
Solving the optimization problem~\eqn{entropy}--\eqn{entropynonneg}, we find the most likely path flows are $h^\uparrow_{AF} = 10$, $h^\downarrow_{AF} = 30$, $h^\uparrow_{BE} = 30$, and $h^\downarrow_{BE} = 90$.
The equilibrium flows for the top and bottom links between C and D have a ratio of 1:3; you can see that this ratio also holds between $h^\uparrow_{AF}$ and $h^\downarrow_{AF}$, as well as between $h^\uparrow_{BE}$ and $h^\downarrow_{BE}$.

In particular, the obvious-looking solution $h^\uparrow_{AF} = 40$, $h^\downarrow_{AF} = 0$, $h^\uparrow_{BE} = 0$, $h^\downarrow_{BE} = 120$ has extremely low entropy, because it implies that for some reason all travelers from one OD pair are taking one path, and all travelers from the other are taking the other path, even though both paths are perceived identically by all travelers and even though travelers are making choices independent of each other.
This is exceptionally unlikely.

\begin{figure}
\begin{center}
\begin{tikzpicture}[->,>=stealth',shorten >=1pt,auto,node distance=3cm,
thick,main node/.style={circle,,draw}]
\node[main node] (C) {$C$};
\node[main node] (D) [right of=C] {$D$};
\node[main node] (A) [above left of=C] {$A$};
\node[main node] (B) [below left of=C] {$B$};
\node[main node] (E) [above right of=D] {$E$};
\node[main node] (F) [below right of=D] {$F$};
\path (C) edge [bend left] node[above] {$40$} (D);
\path (C) edge [bend right] node[below] {$120$} (D);
\path (A) edge node[above right] {$40$} (C);
\path (B) edge node[below right] {$120$} (C);
\path (D) edge node[above left] {$120$} (E);
\path (D) edge node[below left] {$40$} (F);
\end{tikzpicture}
\end{center}
\caption{Multiple OD pair example of proportionality.\label{fig:multiOD}}
\end{figure}

We now derive the proportionality condition by defining it a slightly more precise way:
\begin{thm}
\label{thm:proportionality}
	Let $\pi_1$ and $\pi_2$ be any two paths connecting the same OD pair.
If the path flows $\mathbf{h}$ solve the entropy-maximizing problem~\eqn{entropy}--\eqn{entropynonneg}, then the ratio of flows $h_1/h_2$ for paths $\pi_1$ and $\pi_2$ is identical regardless of the OD pair these paths connect, and only depends on the pairs of alternate segments distinguishing these two paths.

\end{thm}
The proof is a bit lengthy, and is deferred to the end of the chapter.
But even though the derivation is somewhat involved, the proportionality condition itself is fairly intuitive: when confronted with the same set of choices, travelers from different OD pairs should behave in the same way.
The proportionality condition implies nothing more than that the share of travelers choosing one alternative over another is the same across OD pairs.
Proportionality is also a relatively easy condition to track and enforce.

It would be especially nice if proportionality can be shown fully equivalent to entropy maximization.\index{static traffic assignment!proportionality!and entropy maximization}
Theorem~\ref{thm:proportionality} shows that entropy maximization implies proportionality, but is the converse true?  Unfortunately, the result is no, and examples can be created which satisfy proportionality without maximizing entropy.
Therefore, proportionality is a weaker condition than entropy maximization.
The good news is that proportionality ``gets us most of the way there.''  In the Chicago regional network, there are over 93 million equal travel-time paths at equilibrium; after accounting for the equilibrium and ``no vehicle left behind'' constraints~\eqn{entropyeqm} and~\eqn{entropynvlb}, there are still over 90 million degrees of freedom in how the path flows are chosen.
Accounting for proportionality reduces the number of degrees of freedom to 91 (a reduction of 99.9999\%!)  So, in practical terms, enforcing proportionality seems essentially equivalent to maximizing entropy.
\index{entropy|)}
\index{static traffic assignment!properties!nonuniqueness in path flows|)}

\subsection{Aggregation by origins or destinations}
\label{sec:bushes}

\index{static traffic assignment!origin-based solution|see {static traffic assignment, bush-based solution}}
\index{static traffic assignment!bush-based solution}
Thus far, we have looked at flow representations in one of two ways: path flows $\mb{h}$ or link flows $\mb{x}$.
In a way, these can be thought of as two extremes.
The path flow solution contains the most information, showing the exact path chosen by every traveler in the network.
The price paid for this level of detail is the amount of memory needed.
The number of paths in a typical network is very, very large and storing the entire vector of path flows is impractical.
The link flows, on the other hand, are much more compact and provide enough information for many common measures of effectiveness (such as vehicle-hours or vehicle-miles traveled, and volumes or delays on individual links).
The link flow vector can be thought of as an ``aggregated'' form of the path flow vector, where we lump together all the travelers using the same link (regardless of which specific path, they are using, or which specific OD pair they are from).
Naturally, some information is lost in this process of aggregation.
So, there is a tradeoff: path flow solutions contain a great deal of information, at the expense of requiring a large amount of memory; link flow solutions are more compact but contain less information.

So, it is natural to ask if there is an intermediate ``Goldilocks'' flow representation which contains more detail than link flow solutions, without invoking the combinatorial explosion in storage space associated with a path flow solution.
One such representation involves aggregating all travelers departing the same origin together.\footnote{Everything in this section would apply equally well to an aggregation by destination; the presentation from here on will be in terms of origins only to avoid repetition.}  On a link $(i,j)$, let $x_{ij}^r$\label{not:xijr} be the flow on this link who left origin $r$, that is,
\labeleqn{linkdisaggregate}{x_{ij}^r = \sum_{s \in Z} \sum_{\pi \in \Pi^{rs}} \delta_{ij}^\pi h^\pi\,,}
and let $\mb{x}^r$ be the vector of link flows associated with origin $r$.

Many of the most efficient algorithms for solving the traffic assignment problem use such a representation, often called an \emph{origin-based} or \emph{bush-based} representation.
The term ``origin-based'' describes how the aggregation occurs.
The term ``bush-based'' will make more sense once we get to Chapter~\ref{chp:solutionalgorithms}.
Let $\mb{\hat{h}}$ be a feasible assignment satisfying the principle of user equilibrium, and let $\mb{x}^r$ be the link flows associated with origin $r$ at this solution.

It turns out that, at equilibrium, the links with positive $x_{ij}^r$ values form an acyclic subnetwork.
As a result, they have a topological order, which allows network calculations to be done much more rapidly than in general networks.
This property is exploited heavily by bush-based traffic assignment algorithms, which are discussed in Section~\ref{sec:bushbased}.

The following result shows that this must be true when each link's travel time is strictly positive.
With slightly more effort, a similar result can be shown even when there are zero-time links.

\begin{prp}
\label{prp:eqmbushacyclic}
Let $\mb{x}^r$ be the link flows associated with origin $r$ at some feasible assignment satisfying the principle of user equilibrium.
If $t_{ij} > 0$ for all links $(i,j)$, then the subset of arcs $A^r = \{ (i,j) \in A : x^r_{ij} > 0 \}$\label{not:Ar} with positive flow from origin $r$ contains no cycle.
\index{static traffic assignment!bush-based solution!acyclicity}
\index{static traffic assignment!properties!acyclicity}
\end{prp}
\begin{proof}
By contradiction, assume that there is a cycle $[i_1, i_2, \ldots, i_k, i_1]$, where $(i_1, i_2), (i_2, i_3), \ldots, (i_k, i_1)$ are all in $A^r$.
Let $L_i^r$ represent the travel time on the shortest path from origin $r$ to node $i$.
Consider any link $(i,j) \in A^r$; this implies that $x_{ij}^r > 0$ and that some used path starting from origin $r$ includes $(i,j)$.
Because the principle of user equilibrium holds, this path must be a shortest path to some destination.
Clearly $L_i^r + t_{ij} \geq L_j^r$ by the definition of the labels $\mb{L}$; but if $L_i^r + t_{ij} > L_j^r$, then there is a shorter path to node $j$ than any of them passing through node $i$, and such a path cannot be used at equilibrium.
Therefore we must have $L_i^r + t_{ij} = L_j^r$ for all links in $A^r$; since all travel times are strictly positive this implies $L_i^r < L_j^r$ for all $(i,j) \in A^r$.
So, for the cycle under consideration, we have $L_1^r < L_2^r < \cdots < L_k^r < L_1^r$, which is impossible.
Therefore the links in $A^r$ cannot contain a cycle.
\end{proof}

\section{Alternative Assignment Rules}
\label{sec:alternaterules}

\index{assignment rule|(}
The previous section described the most important assignment rule, the principle of user equilibrium.
However, it is not the only possible assignment rule, and it is not hard to see that other factors may enter route choice beyond travel time, or that other assumptions of the user equilibrium model are not perfectly realistic.
Over time, researchers have introduced alternative assignment rules that can be used.
This text by necessity can only cover a few of these alternative assignment rules, and three representative rules were chosen:  system optimum assignment, assignment with perception errors, and bounded rationality.
System optimum assignment is the most direct contrast with user equilibrium, and imagines the ``best possible'' network state that could occur if each traveler's path could be chosen regardless of whether it is their individual shortest path.
Assignment with perception errors relaxes the assumption that drivers have perfect knowledge of travel times.
The interpretation is that every driver \emph{perceives} a particular travel time on each path (which may or may not equal its actual travel time), and chooses a path that they believe to be shortest based on these perceptions.
This section introduces the main ideas in this kind of model, which is termed \emph{stochastic user equilibrium}.
More details are provided in Section~\ref{sec:sue}.
Bounded rationality assumes that drivers are satisfied with a path which is ``close enough'' to shortest, and (for example) do not care if their path is longer than the shortest path by only a few seconds.

\subsection{System optimum assignment}
\label{sec:systemoptimal}

\index{system optimum|(}
Imagine for a moment that route choice was taken out of the hands of individual travelers, and instead placed in the hands of a dictator who could assign each traveler a path that they would be required to follow.
Suppose also that this dictator was benevolent, and wanted to act in a way to minimize average travel delay in the network.
The resulting assignment rule results in the \emph{system optimum} state.

While this scenario is a bit fanciful, the system optimum state is important for several reasons.
First, it provides a theoretical lower bound on the delay in a network, a benchmark for the best performance that could conceivably be achieved that more realistic assignment rules can be compared to.
Second, there are ways to actually achieve the system optimum state even without taking route choice out of the hands of travelers, if one has the ability to freely charge tolls or provide incentives on network links.
Third, there are some network problems where a single agent can exert control over the routes chosen by travelers, as in certain logistics problems where a dispatcher can assign specific routes to vehicles.

At first, it may not be obvious that this assignment rule is meaningfully different from the user equilibrium assignment rule.
After all, in user equilibrium each traveler is individually choosing routes, while in system optimum a single agent is choosing routes for everyone, but both are doing so to minimize travel times.
To see why they might be different, consider the network used for the Knight-Pigou-Downs paradox\index{Knight-Pigou-Downs paradox} in Section~\ref{sec:introductionstatic}: a two-link network (Figure~\ref{fig:knightpigoudowns}) with a constant travel time of 50 minutes on the top link, a link performance function $45 + x^\downarrow$ on the bottom link, and a total demand of 30 vehicles.
The user equilibrium solution was $x^\uparrow = 25$, $x^\downarrow = 5$, with equal travel times of 50 minutes on both top and bottom paths.

So, in the user equilibrium state the average travel time is 50 minutes for all vehicles.
Is it possible to do better?  We can write the average travel time as
\begin{multline}
\label{eqn:averagetime}
\frac{1}{50} \myp{50 x^\uparrow + (45 + x^\downarrow) x^\downarrow} = \frac{1}{50} \myp{50 (30 - x^\downarrow) + (45 + x^\downarrow) x^\downarrow} \\ = \frac{1}{50} \myp{(x^\downarrow)^2 - 5 x^\downarrow + 1500}\,.
\end{multline}
This is a quadratic function which obtains its minimum at $x^\downarrow = 2.5$.
At the solution $x^\uparrow = 27.5$, $x^\downarrow = 2.5$, the travel times are unequal ($t^\uparrow = 50$, $t^\downarrow = 47.5$), but the average travel time of 49.8 minutes is slightly less than the average travel time of 50 minutes at the user equilibrium solution.

Therefore, travelers individually choosing routes to minimize their own travel times may create more delay than would be obtained with central control.
This reveals an important, but subtle point about user equilibrium assignment, a point important enough that Section~\ref{sec:ueso} is devoted entirely to explaining this issue.
So, we will not belabor the point here, but it will be instructive to start thinking about why the user equilibrium and system optimum states need not coincide.

To be more precise mathematically, the system optimum assignment rule chooses the feasible assignment minimizing the average travel time.
Since the total number of travelers is a constant, we can just as well minimize the \emph{total system travel time}\index{total system travel time}, defined as 
\labeleqn{tstt}{TSTT = \sum_{\pi \in \Pi} h^\pi c^\pi = \sum_{(i,j) \in A} x_{ij} t_{ij}\,.}
It is not hard to show that calculating $TSTT$\label{not:TSTT} through path flows and link flows produces the same value.

The system optimum assignment problem can be written in the language of optimization (Section~\ref{sec:convexoptimization}) as
\begin{align}
\index{system optimum!optimization formulation}
\min_{\mathbf{x},\mathbf{h}} \quad & \sum_{(i,j) \in A} x_{ij} \cdot t_{ij} (x_{ij}) &  \label{eqn:soobj} \\
\mathrm{s.t.} \quad & x_{ij} = \sum_{\pi \in \Pi} h^\pi \delta_{ij}^{\pi} & \forall (i,j) \in A \label{eqn:somap} \\
& d^{rs} = \sum_{\pi \in \Pi^{rs}} h^\pi & \forall (r,s) \in Z^2 \label{eqn:sodemand} \\
& h^\pi \geq 0 & \forall \pi \in \Pi \label{eqn:sononneg} 
\,.
\end{align}
In the objective function~\eqn{soobj} we emphasize that $t_{ij}$ is a function of $x_{ij}$, and write out the multiplication with a dot to further clarify the equation --- in particular, $t_{ij}(x_{ij})$ means the function $t_{ij}$ evaluated at $x_{ij}$, \emph{not} $t_{ij}$ multiplied by $x_{ij}$.
We will use this convention consistently.
The constraints~\eqn{somap}--\eqn{sononneg} are identical to those of the user equilibrium problem, constraints \eqn{map}--\eqn{nonneg}.
The only difference is the objective function.
The objective function is convex if the link performance functions $t_{ij}$ are differentiable convex functions (as well as our usual assumption that they are nonnegative and nondecreasing), which is usually not a restrictive assumption.
\index{convex function!applications}
\index{system optimum|)}

\subsection{Perception errors}
\label{sec:perceptionerrors}

\index{stochastic user equilibrium|see {static traffic assignment, stochastic user equilibrium}}
\index{user equilibrium!stochastic|see {static traffic assignment, stochastic user equilibrium}}
\index{static traffic assignment!stochastic user equilibrium|(}
An alternative assignment rule tries to relax the assumption that travelers choose the path with exactly the shortest travel time.
After all, this assumption implicitly requires \emph{drivers to have perfect knowledge of the travel times on all paths in the network}.
In reality, we know this is not true: do you know the travel times on literally \emph{all} path between an origin and destination?  And can you accurately distinguish between a path with a travel time of 16 minutes, and one with a travel time of 15 minutes and 59 seconds?  And perhaps travelers take other factors into account anyway, travel time may be just one of several criteria used to choose paths.

Both of these factors can be modeled using discrete choice concepts.
In discrete choice models, an individual chooses one option from a set of alternatives to maximize his or her utility, consisting of both an observed portion and an unobserved portion.
Rather than requiring that each driver follow the true shortest path between each origin and destination, we assume that drivers follow the path they \emph{believe} to be shortest, but allow for some perception error between their belief and the actual travel times.
An alternative, mathematically equivalent, interpretation (explained below) is that drivers do in fact perceive travel times accurately, but care about factors other than travel time.
This leads to another important traffic assignment model commonly called \emph{stochastic user equilibrium} (SUE).
We discuss this model briefly here, and at length in Section~\ref{sec:sue}.

Consider a traveler leaving origin $r$ for destination $s$.
They must choose one of the paths $\pi$ connecting $r$ to $s$, that is, they must make a choice from the set $\Pi_{rs}$.
The most straightforward way to generalize the principle of user equilibrium to account for perception errors is to set the observed utility\index{utility} equal to the negative of path travel time, so\label{not:Upi}
\labeleqn{sueutility_brief}{U^\pi = -c^\pi + \epsilon^\pi\,,}
with the negative sign indicating that \emph{maximizing utility} for drivers means \emph{minimizing travel time}.
Assuming that the $\epsilon^\pi$ are independent, identically distributed Gumbel\index{Gumbel distribution} random variables, we can use the logit\index{logit} formula~\eqn{logitpath_brief} to express the probability that path $\pi$ is chosen:
\labeleqn{logitpath_brief}{p^\pi = \frac{\exp(-\theta c^\pi)}{\sum_{\pi' \in \Pi^{rs}} \exp(\theta C_{\pi'})}\,.}
In this formula, $\theta$\label{not:theta} is a parameter related to the variance of the Gumbel distribution (higher variance corresponds to lower $\theta$, and vice versa).
As $\theta$ approaches 0, drivers' perception errors are large relative to the path travel times, and each path is chosen with nearly equal probability: the errors are so large, the choice is essentially random.
As $\theta$ grows large, perception errors are small relative to path travel times, and the path with lowest travel time is chosen with higher and higher probability.
At any level of $\theta$, there is a strictly positive probability that each path will be taken.

For concreteness, the route choice discussion so far corresponds to the interpretation where the unobserved utility represents perception errors in the utility.
The other interpretation would mean that $\epsilon^\pi$ represents factors other than travel time which affect route choice (such as comfort, quality of scenery, etc.).
Either of these interpretations is mathematically consistent with the discussion here.

The fact that the denominator of~\eqn{logitpath_brief} includes a summation over \emph{all} paths connecting $r$ to $s$ is problematic, both theoretically and practically.
From a theoretical standpoint, it implies that drivers are considering literally every path between the origin and destination, even when this path is very circuitous or illogical (driving all around town when going to a store a half mile away).
Presumably no driver would ever choose such a path, no matter how ill-informed they are about travel conditions or how poorly they can estimate travel times --- yet the logit formula~\eqn{logitpath_brief} suggests that some travelers will indeed choose such paths.
From a practical standpoint, evaluating the formula~\eqn{logitpath_brief} first requires \emph{enumerating} all of these paths.\index{path!enumeration challenges}
Since the number of paths grows exponentially with the network size, any approach which requires an explicit listing of all the paths in a network will not scale to realistic-sized problems.

To address this fact, we can restrict the choice set somewhat.
Rather than using all paths $\Pi^{rs}$ connecting $r$ to $s$, we can restrict path choices to a subset of \emph{reasonable paths} (denoted $\hat{\Pi}^{rs}$).\index{path!set of used paths}
There are different ways to identify sets of reasonable paths, but they should address both the theoretical and practical difficulties in the previous paragraph.
That is, the paths in $\hat{\Pi}^{rs}$ correspond to a plausible-sounding behavioral principle (why might travelers only consider paths in this set) and lead to a formula which is efficiently computable even in large networks.
There may be some tension between these ideas, in that there may be very efficient formulas which do not correspond to realistic choices, or that the most realistic models of reasonable paths may not lead to an efficient formula.
Section~\ref{sec:sue} discusses these in more detail.

So, equation~\eqn{logitpath_brief} leads us an expression for path flows in terms of travel times, which is the assignment rule filling the place of the question mark in Figure~\ref{fig:fundamentalloop}.
To be explicit, the formula for path flows as a function of path travel times is 
\labeleqn{logitpathexplicit_brief}{h^\pi = d^{rs} \frac{\exp(-\theta c^\pi)}{\sum_{\pi' \in \Pi^{rs}} \exp(\theta C_{\pi'})}\,,}
where $(r,s)$ is the OD pair corresponding to path $\pi$.
The complete traffic assignment problem with this assignment rule can then be expressed as follows: find a feasible path flow vector $\hat{\mb{h}}$ such that $\mb{\hat{h} = H(C(\hat{h}))}$.
This is a standard fixed point problem.\index{static traffic assignment!stochastic user equilibrium!fixed point formulation}
Clearly $\mb{H}$ and $\mb{C}$ are continuous functions if the link performance functions are continuous, and the feasible path set is compact and convex, so Brouwer's theorem immediately gives existence of a solution to the SUE problem.
\index{static traffic assignment!stochastic user equilibrium!existence}
\index{Brouwer's theorem}
\index{continuous function!applications}
\index{compact set!applications}
\index{convex set!applications}

Notice that this was much easier than showing existence of an equilibrium solution to the original traffic assignment problem!  For that problem, there was no equivalent of~\eqn{logitpathexplicit_brief}.
Travelers were all using shortest paths, but if there were two or more shortest paths there was no rule for how those ties should be broken.
As a result, we had to reformulate the problem as a variational inequality and introduce an auxiliary function based on movement of a point under a force.
With this assignment rule, there is no need for such machinations, and we can write down the fixed point problem immediately.

That said, it \emph{is} also possible to formulate the assignment problem with this rule as the solution to a variational inequality, and as the solution to a convex minimization problem.
However, the objective function involves enumerating all of the used paths, which can pose computational difficulties in large networks.
The objective function, which will be derived in Section~\ref{sec:sue}, is
\labeleqn{sueobjective_brief}{\theta \sum_{(i,j) \in A}  \int_0^{x_{ij}} t_{ij}(x) dx + \sum_{\pi \in \hat{\Pi}} h^\pi \log h^\pi
\,.}

To demonstrate what this assignment rule looks like, refer again to the network in Figure~\ref{fig:2linksimple}.
Figures~\ref{fig:SUEflows} and~\ref{fig:SUEtimes} show how the flows and travel times in the network vary with $\theta$.
A few observations worth noting: first, the travel times on the two links are \emph{not} equal, due to the presence of perception errors.
However, as drivers' perceptions become more accurate ($\theta$ increasing), the travel times on the two paths become closer, and the solution asymptotically approaches the user equilibrium solution.
\index{static traffic assignment!stochastic user equilibrium|)}

\stevefig{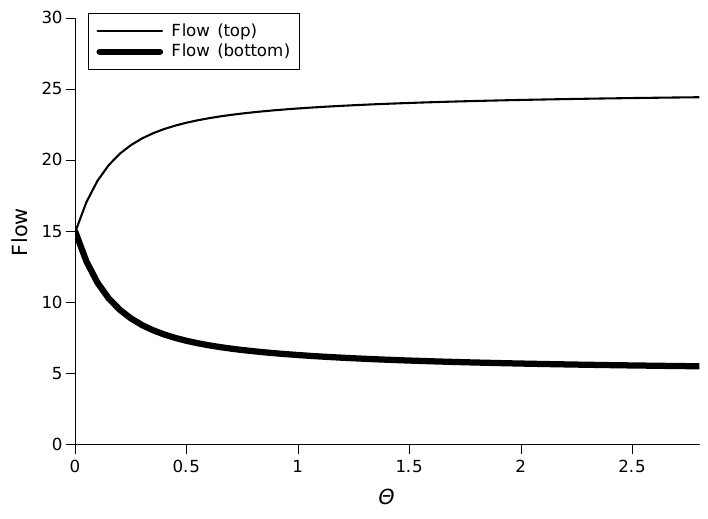}{Flows with different perception errors.}{0.8\textwidth}
\stevefig{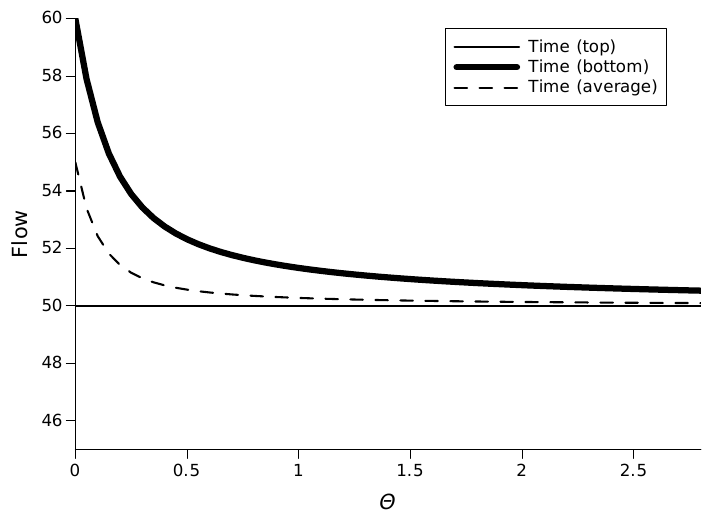}{Travel times with different perception errors.}{0.8\textwidth}

\subsection{Bounded rationality (*)}
\label{sec:boundedrationality}

\emph{(This optional section discusses an assignment rule which relaxes the assumption that all travelers are on the shortest path.)}

\index{bounded rationality|(}
An alternative assignment rule is based on another variation of shortest path-seeking behavior.
Under \emph{bounded rationality}, drivers will choose paths that are within some threshold of the true shortest path travel time.
The behavioral rationale is that a driver will switch paths from a longer path to a shorter one, but only if the new path is sufficiently faster than their current choice.
Perhaps a driver is unlikely to switch paths to save, say, a tenth of a second in travel time: the burden of learning a new path outweighs any travel time savings which are possible.

There is some evidence from behavioral economics that humans often make choices in a boundedly rational way.
Specifically in network assignment, bounded rationality can explain observed reactions to changes in a network in a way that the principle of user equilibrium cannot.
Of course, these advantages must be balanced against other concerns: there may be many distinct bounded rational assignments, and efficient methods for finding bounded rational assignments in large networks have not yet been developed.

The boundedly rational user equilibrium (BRUE) problem can be stated as follows.
Let $\epsilon$ represent some tolerance value in terms of travel time, and let $\kappa^{rs}$ be the travel time on the shortest path connecting origin $r$ to destination $s$.

\begin{dfn}
\emph{(Principle of boundedly rational user equilibrium [BRUE].)}  Every used path has travel time no more than $\epsilon$ in excess of the travel time on the shortest path connecting these nodes.
\index{user equilibrium!bounded rationality}
That is, the feasible assignment $\mb{h} \in H$ is a BRUE if, for each $(r,s) \in Z^2$ and $\pi \in \Pi^{rs}$, $h^\pi > 0$ implies $c^\pi \leq \kappa^{rs} + \epsilon$.
\end{dfn}

As an example, consider once again the network in Figure~\ref{fig:2linksimple}, but now assume we are seeking a BRUE where the tolerance is $\epsilon = 2$ minutes.
That is, travelers are willing to use any path as long as it is within one minute of being the shortest path.
The user equilibrium solution $x^\uparrow = 25$, $x^\downarrow = 5$ is clearly a BRUE as well.
However, there are others: if $x^\uparrow = 26$, $x^\downarrow = 4$, then the travel times on the top and bottom paths are 50 and 49 minutes, respectively.
Even though the majority of travelers (those on the top link) are not on the shortest path, the travel time on their path is close enough to the shortest path that the situation is acceptable, and this solution is BRUE as well.
However, the flows $x^\uparrow = 28$, $x^\downarrow = 2$ are not BRUE: the travel times are 50 and 47 minutes, and the difference between the (used) top path and the shortest path (3 minutes) exceeds the tolerance value (2 minutes).
You should convince yourself that $x^\uparrow \in [23,27]$ and $x^\downarrow = 30 - x^\uparrow$ characterize all of the BRUE solutions in this simple network.

To come up with a mathematical formulation of BRUE assignment, it is helpful to introduce an auxiliary variable $\rho^\pi$\label{not:rhopi} indicating the amount by which the travel time on path $\pi$ falls short of the maximum acceptable travel time on a path.
Specifically, define
\labeleqn{excesscost}{\rho^\pi = \mys{\kappa^{rs} + \epsilon - c^\pi}^+ \qquad \forall (r,s) \in Z^2, \pi \in \Pi^{rs}\,.}
So, in the example of the previous paragraph, if $\epsilon = 1$, $x^\uparrow = 26$ and $x^\downarrow = 4$, then $t^\uparrow = 50$, $t^\downarrow = 49$, and thus $\rho^\uparrow = [49 + 2 - 50]^+ = 1$ and $\rho^\downarrow = [49 + 2 - 49]^+ = 2$.
If the path $\pi$ is unacceptable ($c^\pi > \kappa^{rs} + \epsilon$), then $\rho^\pi = 0$ by~\eqn{excesscost}.
We can then define the \emph{effective travel time} $\hat{c}^\pi$\label{not:chatpi} as
\labeleqn{effectivetime}{\hat{c}^\pi = c^\pi + \rho^\pi\,.}

These new quantities are helpful because the principle of BRUE can now be formulated in a more familiar way:
\begin{prp}
\label{prp:brueequalminimal}
A feasible assignment $\mb{h} \in H$ is a BRUE if and only if every used path has equal and minimal effective travel time, that is, for each $(r,s) \in Z^2$ and $\pi \in \Pi^{rs}$, $h^\pi > 0$ implies $\hat{c}^\pi = \min_{\pi' \in \Pi^{rs}} \hat{c}^{\pi'}$.

\end{prp}
\begin{proof}
Let $\mb{h} \in H$ satisfy the principle of BRUE, and consider any OD pair $(r,s)$ and path $\pi \in \Pi^{rs}$ with $h^\pi > 0$.
Since $\mb{h}$ is BRUE, $c^\pi \leq \kappa^{rs} + \epsilon$, so $\rho^\pi = \kappa^{rs} + \epsilon - c^\pi$ and $\hat{c}^\pi = c^\pi + \kappa^{rs} + \epsilon - c^\pi = \kappa^{rs} + \epsilon$.
But for any $\pi' \in \Pi^{rs}$, $\hat{c}^{\pi'} = c^{\pi'} + \rho{\pi'} \geq c^{\pi'} + \kappa^{rs} + \epsilon - c^{\pi'} = \kappa^{rs} + \epsilon$, so $\hat{c}^\pi = \min_{\pi' \in \Pi^{rs}} \hat{c}^{pi'}$.

Contrariwise, let $\mb{h} \in H$ be such that $h^\pi > 0$ implies $\hat{c}^\pi = \min_{\pi' \in \Pi^{rs}} \hat{c}^{\pi'}$, and consider any OD pair $(r,s)$.
Let $\hat{\pi}$ be a shortest path connecting $r$ to $s$, so $\hat{c}^{\hat{\pi}} = \kappa^{rs} + [\kappa^{rs} + \epsilon - \kappa^{rs}]^+ = \kappa^{rs} + \epsilon$ and this must be minimal among all effective travel times on paths connecting $r$ to $s$.
For any path $\pi \in \Pi^{rs}$ with $h^\pi > 0$, we thus have $\hat{c}^\pi = \kappa^{rs} + \epsilon$, so $c^\pi = \hat{c}^\pi - \rho^\pi = \kappa^{rs} + \epsilon - \rho^\pi \leq \kappa^{rs} + \epsilon$ since $\rho^\pi \geq 0$.
Thus $\mb{h}$ satisfies the principle of BRUE as well.

\end{proof}

Proposition~\ref{prp:brueequalminimal} may give you hope that we can write a convex optimization problem whose solutions correspond to BRUE.
Unfortunately, this is not the case, for reasons that will be shown shortly.
However, we can write a variational inequality where both the path flows $\mb{h}$ and the auxiliary variables $\bm{\rho}$ are decision variables.
In what follows, let $\kappa^{rs}(\mb{h}) = \min_{\pi \in \Pi^{rs}} c^\pi(\mb{h})$ and, with slight abuse of notation, let $\kappa^\pi$ refer to the $\kappa^{rs}$ value corresponding to the path $\pi$, and $\bm{\kappa}$ be the vector of $\kappa^\pi$ values.

\begin{prp}
\label{prp:bruevi}
\index{bounded rationality!formulation!as variational inequality}
A vector of path flows $\mb{\hat{h}} \in H$ is a BRUE if there exists a vector of nonnegative auxiliary costs $\hat{\bm{\rho}} \in \bbr_+^{|\Pi|}$\label{not:bbrp} such that
\labeleqn{bruevi}{\vect{ \hat{\mb{c}}(\mb{\hat{h}},\hat{\bm{\rho}}) \\ \mb{c}(\mb{\hat{h}}) + \hat{\bm{\rho}} - \bm{\kappa}(\mb{\hat{h}}) - \epsilon} \cdot \vect{ \mb{\hat{h}} - \mb{h} \\ \hat{\bm{\rho}} - \bm{\rho} } \leq 0}
for all $\mb{h} \in H$, $\bm{\rho} \in \bbr_+^{|\Pi|}$.

\end{prp}
\begin{proof}
Assume that $(\mb{\hat{h}}, \hat{\bm{\rho}})$ solve the variational inequality~\eqn{bruevi}.
The first component of the VI $\hat{\mb{c}}(\mb{\hat{h}},\hat{\bm{\rho}}) \cdot (\mb{\hat{h}} - \mb{h}) \leq 0$ shows that all used paths have equal and minimal effective cost, and by Proposition~\ref{prp:brueequalminimal} satisfy the principle of BRUE, assuming that $\bm{\rho}$ is consistent with~\eqn{excesscost}.
To show this, consider the second component of the VI: $(\bm{\kappa}(\mb{\hat{h}}) + \epsilon - \hat{\bm{\rho}}) \cdot \hat{\bm{\rho}} - \bm{\rho} \leq 0$, or equivalently
\labeleqn{rhovi}{\sum_{\pi \in \Pi} (c^\pi(\mb{\hat{h}}) + \hat{\rho}_\pi - \kappa^\pi(\mb{\hat{h}}) - \epsilon) \hat{\rho}_\pi \leq \sum_{\pi \in \Pi} (c^\pi(\mb{\hat{h}}) + \hat{\rho}_\pi - \kappa^\pi(\mb{\hat{h}}) - \epsilon) \rho_\pi}
for all $\bm{\rho} \in \bbr_+^{|\Pi|}$.
If $\hat{\rho}_\pi > 0$ for any path, inequality~\eqn{rhovi} can be true only if $\rho^\pi = \kappa^\pi(\mb{\hat{h}}) + \epsilon - c^\pi(\mb{\hat{h}})$; or, if $\hat{\rho}_\pi = 0$, then inequality~\eqn{rhovi} can be true only if $\rho^\pi = 0 \leq \kappa^\pi(\mb{\hat{h}}) + \epsilon - c^\pi(\mb{\hat{h}})$.
In either case~\eqn{excesscost} is satisfied.

In the reverse direction, assume that $\mb{\hat{h}}$ is BRUE.
Choose $\hat{\bm{\rho}} \in \bbr_+^{|\Pi|}$ according to~\eqn{excesscost}.
By Proposition~\ref{prp:brueequalminimal} surely $\hat{\mb{c}}(\mb{\hat{h}},\hat{\bm{\rho}}) \cdot (\mb{\hat{h}} - \mb{h}) \leq 0$  for all $\mb{h} \in H$.
For any path where $\rho^\pi = \kappa^\pi(\mb{\hat{h}}) + \epsilon - c^\pi(\mb{\hat{h}})$, clearly $\sum_{\pi \in \Pi} (c^\pi(\mb{\hat{h}}) + \hat{\rho}_\pi - \kappa^\pi(\mb{\hat{h}}) - \epsilon) (\hat{\rho}_\pi - \rho^\pi) = 0$; and for any path where $\rho^\pi = 0$, $\kappa^\pi(\mb{\hat{h}}) + \epsilon - c^\pi(\mb{\hat{h}}) \geq 0$ and $\sum_{\pi \in \Pi} (c^\pi(\mb{\hat{h}}) + \hat{\rho}_\pi - \kappa^\pi(\mb{\hat{h}}) - \epsilon) (\hat{\rho}_\pi - \rho^\pi) \geq 0$.
In either case the second component of the VI is satisfied as well.

\end{proof}

This BRUE formulation has been helpful in explaining observed changes in network flows after a disruption.
For instance, assume that a link is removed from the network due to a disaster of some sort, and that flows adjust towards a new equilibrium in the network without the affected link.
When the link is restored, flows will adjust again.
If the principle of user equilibrium is true, the flows will move back to exactly the same values as before.
However, in practice there has been some ``stickiness'' observed, and not all drivers will return to the same routes they were initially on.
The BRUE framework provides a logical explanation for this: when the network is disrupted, certain drivers were forced to choose new paths.
When the network is restored, they will only switch back to their original paths if the travel time savings are sufficiently large.
Otherwise, they will remain on their new paths.

\index{bounded rationality!properties|(}
However, an implication of this finding is that the BRUE solution need not be unique, even in link flows.
(This was also seen in the small example at the start of this subsection.)  In general, the set of BRUE solutions is not convex either --- which immediately shows that there is no convex optimization problem whose solutions correspond to BRUE assignments, since the sets of minima to convex optimization problems are always convex.
However, there is always at least one BRUE solution by Proposition~\ref{prp:existence}, since any path flows satisfying the principle of user equilibrium also satisfy BRUE.

Still, the lack of uniqueness means that BRUE should be applied carefully.
While it is certainly more realistic than the principle of user equilibrium (since user equilibrium can always be obtained as a special case with $\epsilon = 0$), it is considerably more difficult to use BRUE to rank projects (which BRUE assignment should be used?) and at present it is not clear what value of $\epsilon$ best represents travel choices.
Whether the added realism offsets these disadvantages is application and network dependent.
\index{bounded rationality!properties|)}
\index{bounded rationality|)}

\section{User Equilibrium and System Optimum}
\label{sec:ueso}

\index{user equilibrium|(}
\index{system optimum|(}
The relationship between the user equilibrium and system optimum assignment rules is worth studying further, for a few reasons.
First, the fact that these rules lead to distinct solutions means that travelers cannot be relied upon to independently behave in a way that minimizes congestion\footnote{Perhaps not a very surprising observation to those used to driving in crowded cities.}, even though each driver is minimizing his or her own travel time.
This means that there is a role for transportation professionals to reduce congestion through planning and operational control.
Exploring why this happens leads to a discussion of economic externalities, which are illustrated in the Braess paradox discussed in Section~\ref{sec:motivatingexamples}.

Second, the mathematical structure of the system optimum assignment problem is actually very similar to that of the user equilibrium problem.
In fact, if one has a ``black box'' procedure for solving user equilibrium assignment, one can just as easily solve system optimum assignment using that same procedure, making some easy modifications to the network.
Likewise, if one has a procedure for solving system optimum assignment, one can just as easily solve user equilibrium by modifying the network.
So, when we discuss algorithms for solving the user equilibrium traffic assignment problem in Chapter~\ref{chp:solutionalgorithms}, they all apply equally well to the system optimum problem.
This relationship is explored in Section~\ref{sec:uesoequivalence}.

Finally, there is a more recent and interesting set of results which attempts to quantify how far apart the user equilibrium and system optimum assignments could possibly be, in terms of total system travel time.
At first glance it seems that these assignment rules would tend to give relatively similar total travel times, since both are based on travel time minimization (one individual and the other collective).
However, the Braess paradox shows that these can diverge in surprising ways.
Nevertheless, in many cases one can in fact bound how much inefficiency is caused by allowing individual drivers to pick their own routes, a set of results called the \emph{price of anarchy} and discussed in Section~\ref{sec:priceofanarchy}.

\subsection{Externalities}
\label{sec:externalities}

\index{externality|(}
Recall the three examples in Section~\ref{sec:motivatingexamples}, in which equilibrium traffic assignment exhibited counterintuitive results.
Perhaps the most striking of these was the Braess paradox\index{Braess paradox}, in which adding a new roadway link actually worsened travel times for all travelers.
This is counterintuitive, because if the travelers simply agreed to stay off of the newly-built road, their travel times would have been no worse than before.
This section describes exactly why this happens.
But before delving into the roots of the Braess paradox, let's discuss some of the implications:

\textbf{User equilibrium does \emph{not} minimize congestion.}  In principle, everyone could have stayed on their original paths and retained the lower travel time of 83; the only catch was that this was no longer an equilibrium solution with the new link, and the equilibrium solution with the link is worse.
Therefore, there may be traffic assignments with lower travel times (for everybody!) than the user equilibrium one.

\textbf{The ``invisible hand'' may not work in traffic networks.}  In economics, Adam Smith's\index{Smith, Adam} celebrated ``invisible hand''\index{invisible hand} suggests that each individual acting in his or her own self-interest also tends to maximize the interests of society as well.
The Braess paradox shows that this is not necessarily true in transportation networks.
The two drivers switched paths out of self-interest, to reduce their own travel times; but the end effect made things worse for everyone, including the drivers who switched paths.

\textbf{There may be room to ``improve'' route choices.} A corollary of the previous point, when the invisible hand does not function well, a case can be made for regulation or another coordination mechanism to make everyone better off.
In this case, this mechanism might take the form of certain transportation policies.
What might they be?

If you have studied economics, you might already have some idea of why the Braess paradox occurs.
The ``invisible hand'' requires certain things to be true if individual self-interest is to align with societal interest, and fails otherwise.
A common case where the invisible hand fails is in the case of \emph{externalities}: an externality is a cost or benefit resulting from a decision made by one person, and felt by another person who had no say in the first person's decision.
Examples of externalities include industrial pollution (we cannot rely on industries to self-regulate pollution, because all citizens near a factory suffer the effects of poor air quality, even though they had no say in the matter), driving without a seatbelt (if you get into an accident without a seatbelt, it is likely to be more serious, costing others who help pay your health care costs and delaying traffic longer while the ambulance arrives; those who help pay your costs and sit in traffic had no control over your decision to not wear a seatbelt), and education (educated people tend to commit less crime, which benefits all of society, even those who do not receive an education).
The first two examples involve costs, and so are called negative externalities; the latter involves a benefit, and is called a positive externality.

Another example of an externality is seen in the prisoner's dilemma\index{prisoners' dilemma} (the Erica-Fred game of Section~\ref{sec:notionofequilibrium}).
When Erica chooses to testify against Fred, she does so because it reduces her jail time by one year.
The fact that her choice also increases Fred's jail time by fourteen years was irrelevant (such is the nature of greedy decision-making).
When \emph{both} Erica and Fred behaved in this way, the net effect was to dramatically increase the jail time they experienced, even though each of them made the choices they did in an attempt to minimize their jail time.

The relevance of this economics tangent is that \emph{congestion} can be thought of as an externality.
Let's say I'm driving during the peak period, and I can choose an uncongested back road which is out of my way, or a congested freeway which is a direct route (and thus faster, even with the congestion).
If I choose the freeway, I end up delaying everyone with the bad luck of being behind me.
Perhaps not by much; maybe my presence increases each person's commute by a few seconds.
However, multiply these few seconds by a large number of drivers sitting in traffic, and we find that my presence has cost society as a whole a substantial amount of wasted time.
This is an externality, because those other drivers had no say over my decision to save a minute or two to take the freeway.
(This is why the user equilibrium assumption is said to be ``greedy.''  It assumes a driver will happily switch routes to save a minute even if it increases the total time people spend waiting in traffic by an hour.)  These concepts are made more precise in the next subsection.

As a society, we tend to adopt regulation or other mechanisms for minimizing the impact of negative externalities, such as emissions regulation or seatbelt laws (although you can probably think of many negative externalities which still exist).
Transfer payments are another option: if a factory had to directly compensate each citizen for the harm they suffer as a result of pollution, the externality no longer exists because the factory now has to account for the costs it imposed on society, and will only continue to operate if its profits exceed the cost of the pollution it causes.
What types of regulation or transfer payments might exist in the transportation world?
\index{externality|)}

\subsection{Mathematical equivalence}
\label{sec:uesoequivalence}

\index{user equilibrium!relation to system optimum|(}
\index{system optimum!relation to user equilibrium|(}
Recall that the user equilibrium traffic assignment was the solution to the following convex optimization problem:
\begin{align}
	\min_{\mathbf{x},\mathbf{h}}	\quad	& \sum_{(i,j) \in A} \int_0^{x_{ij}} t_{ij}(x) dx	&		\label{eqn:beckmann_65} \\
	\mathrm{s.t.}						\quad & x_{ij} = \sum_{\pi \in \Pi} h^\pi \delta_{ij}^{\pi}	& \forall (i,j) \in A \label{eqn:map_65} \\
													& d^{rs} = \sum_{\pi \in \Pi^{rs}} h^\pi 					&		\forall (r,s) \in Z^2 \label{eqn:demand_65} \\
													& h^\pi \geq 0													&		\forall \pi \in \Pi \label{eqn:nonneg_65}
\,,
\end{align}
while the system optimum traffic assignment solves this convex optimization problem:
\begin{align}
	\min_{\mathbf{x},\mathbf{h}}	\quad	& \sum_{(i,j) \in A} x_{ij} \cdot t_{ij}(x_{ij}) 	&		\label{eqn:sobeckmann_65} \\
	\mathrm{s.t.}						\quad & x_{ij} = \sum_{\pi \in \Pi} h^\pi \delta_{ij}^{\pi}	& \forall (i,j) \in A \label{eqn:somap_65} \\
													& d^{rs} = \sum_{\pi \in \Pi^{rs}} h^\pi 					&		\forall (r,s) \in Z^2 \label{eqn:sodemand_65} \\
													& h^\pi \geq 0													&		\forall \pi \in \Pi \label{eqn:sononneg_65} 
\,.
\end{align}

Inspecting these, the only difference is in their objective functions, and even these are quite similar: both are sums of terms involving each link in the network.
The problems share essentially the same structure, the only difference is the term corresponding to each link.
For the user equilibrium problem, the term for link $(i,j)$ is $\int_0^{x_{ij}} t_{ij}(x)~dx$, while for the system optimum problem it is $x_{ij} \cdot t_{ij}(x_{ij})$.

Suppose for a moment that we had some algorithm which can solve the user equilibrium traffic assignment problem for any input network and OD matrix.
(A number of these algorithms will be discussed in Chapter~\ref{chp:solutionalgorithms}.)
Now suppose that we had to solve the system optimum problem.
If we were to replace the link performance functions\index{link performance function!transformations} $t_{ij}(x_{ij})$ with modified functions $\hat{t}_{ij}(x_{ij})$\label{not:thatij} defined by
\labeleqn{marginalcost}{\hat{t}_{ij}(x_{ij}) = t_{ij}(x_{ij}) + x_{ij} \cdot t'_{ij}(x_{ij})\,,}
where $t'_{ij}$ is the derivative of the link performance functions, then
\labeleqn{equivalence1}{\int_0^{x_{ij}} \hat{t}_{ij}(x)~dx = \int_0^{x_{ij}} (t_{ij}(x) + x \cdot t'_{ij}(x))~dx\,.}
Integrating the second term by parts, this simplifies to
\labeleqn{equivalencefinal}{\int_0^{x_{ij}} \hat{t}_{ij}(x)~dx = x_{ij} \cdot t_{ij}(x_{ij})\,,}
which is the same term for link $(i,j)$ used in the system optimum function.

That is, solving user equilibrium with the modified functions $\hat{t}_{ij}$ produces the same link and path flow solution as solving system optimum with the original performance functions $t_{ij}$.
The interpretation of the formula~\eqn{marginalcost} is closely linked with the concept of externalities introduced in the previous subsection.
Equation~\eqn{marginalcost} consists of two parts: the actual travel time $t_{ij}(x_{ij})$, and an additional term $x_{ij} \cdot t'_{ij}(x_{ij})$.
If an additional vehicle were to travel on link $(i,j)$, it would experience a travel time of $t_{ij}(x_{ij})$.
At the margin, the travel time on link $(i,j)$ would increase by approximately $t'_{ij}(x_{ij})$, and this marginal increase in the travel time would be felt by the other $x_{ij}$ vehicles on the link.
So, the second term in~\eqn{marginalcost} expresses the additional, external increase in travel time caused by a vehicle using link $(i,j)$.
For this reason, the modified functions $\hat{t}_{ij}$ can be said to represent the \emph{marginal cost}\index{marginal cost}\index{link!marginal cost} on a link.

Since the system optimum solution is a ``user equilibrium'' with respect to the modified costs $\hat{t}_{ij}$, we can formulate a ``principle of system optimum'' similar to the principle of user equilibrium:
\begin{dfn}
  \label{dfn:so}
	\emph{(Principle of system optimum.)}
	Every used path connecting an origin and destination has equal and minimal marginal cost.
    \index{system optimum!principle of}
\end{dfn}

So, system optimum can be seen as a special case of user equilibrium, with modified link performance functions.
The converse is true as well.
Suppose we had an algorithm which would solve the system optimum problem for any input network and OD matrix.
If we replace the link performance functions $t_{ij}(x_{ij})$ with modified functions $\tilde{t}_{ij}(x_{ij})$\label{not:ttwiddleij} defined by\index{link performance function!transformations}
\labeleqn{averagecost}{\tilde{t}_{ij}(x_{ij}) = \frac{\int_0^{x_{ij}} t_{ij}(x)~dx}{x_{ij}}}
when $x_{ij} > 0$ and $\tilde{t}_{ij}(x_{ij}) = 0$ otherwise, a similar argument shows that the objective function for the system optimum problem with link performance functions $\tilde{t}_{ij}$ is identical to that for the user equilibrium problem with link performance functions $t_{ij}$.

As a result, despite very different physical interpretations, the user equilibrium and system optimum problems are essentially the same mathematically.
If you can solve one, you can solve the other by making some simple changes to the link performance functions.
So, there is no need to develop separate algorithms for the user equilibrium and system optimum problems.

\subsection{Price of anarchy}
\label{sec:priceofanarchy}

\index{price of anarchy|(}
The discussion of the user equilibrium and system optimum problems so far has gone in several directions.
On the one hand, it seems like the problems are fairly similar: both involve paths chosen according to some travel time minimization rule, and they are mathematically equivalent.
On the other hand, the real-world interpretations are very different (individual control vs.\ centralized control), and the Braess paradox shows that user equilibrium can behave in a far more counterintuitive manner than system optimum.
This subsection concludes the discussion by showing that, while the user equilibrium solution may be different than the system optimum solution, and worse in terms of higher total system travel time, in many cases there is a limit to how much worse it can be.

In particular, if the link performance functions are linear (that is, of the form $t_{ij}(x_{ij}) = a_{ij} + b_{ij} x_{ij}$\label{not:aij2}\label{not:bij}), the ratio between the total system travel time at user equilibrium, and the total system travel time at system optimum, can be no greater than $\frac{4}{3}$.
This is often called the \emph{price of anarchy}, since it represents the amount of additional delay which can potentially be caused by allowing individual drivers to choose their own routes, compared to a centrally-controlled solution.
This bound holds no matter what the network structure or the level of demand, as long as the link performance functions are affine.
This is an elegant and perhaps surprising result, which is actually not difficult to show using the variational inequality formulation of equilibrium.

\begin{thm}
Given a network with link performance functions of the form $t_{ij}(x_{ij}) = a_{ij} + b_{ij} x_{ij}$ with $a_{ij} \geq 0$ and $b_{ij} \geq 0$ for all links $(i,j)$, let feasible link flows $\mb{\hat{x}} \in X$ satisfy the principle of user equilibrium and $\mb{\hat{y}} \in X$ the system optimum principle.
Denoting the total system travel times corresponding to these solutions as $TSTT(\mb{\hat{x}})$ and $TSTT(\mb{\hat{y}})$, we have
\labeleqn{priceofanarchy}{\frac{TSTT(\mb{\hat{x}})}{TSTT(\mb{\hat{y}})} \leq \frac{4}{3}\,.}  
\end{thm}
\begin{proof}
Since $\mb{\hat{x}}$ satisfies the principle of user equilibrium, we have
\labeleqn{poavi}{\mb{t}(\mb{\hat{x}}) \cdot (\mb{\hat{x}} - \mb{\hat{y}}) \leq 0}
or, equivalently,
\begin{align}
   \mb{t}(\mb{\hat{x}}) \cdot \mb{\hat{x}} &\leq \mb{t}(\mb{\hat{x}}) \cdot \mb{\hat{y}} \\
                                   &= (\mb{t}(\mb{\hat{x}}) - \mb{t}(\mb{\hat{y}})) \cdot \mb{\hat{y}} + \mb{t}(\mb{\hat{y}}) \cdot \mb{\hat{y}}
\,,                             
\end{align}
which is the same as
\labeleqn{poatstt}{TSTT(\mb{\hat{x}}) \leq (\mb{t}(\mb{\hat{x}}) - \mb{t}(\mb{\hat{y}})) \cdot \mb{\hat{y}} + TSTT(\mb{\hat{y}})\,.}

If we can show that
\labeleqn{targetbound}{(\mb{t}(\mb{\hat{x}}) - \mb{t}(\mb{\hat{y}})) \cdot \mb{\hat{y}} \leq \frac{1}{4} TSTT(\mb{\hat{x}})}
then we have established the result, by substitution into~\eqn{poatstt}.
We do this by showing an even stronger result, namely
\labeleqn{strongerbound}{(t_{ij}(\hat{x}_{ij}) - t_{ij}(\hat{y}_{ij})) \hat{y}_{ij} \leq \frac{1}{4} t_{ij}(\hat{x}_{ij}) \hat{x}_{ij}}
for all links $(i,j)$ which is clearly sufficient.

If $\hat{x}_{ij} \leq \hat{y}_{ij}$, then the left-hand side of~\eqn{strongerbound} is nonpositive (the link performance functions are nondecreasing since $b_{ij} \geq 0$), while the right-hand side is nonnegative (since $a_{ij} \geq 0$), so the inequality is certainly true.
On the other hand, if $\hat{x}_{ij} > \hat{y}_{ij}$, then 
\labeleqn{simplifiedbound}{(t_{ij}(\hat{x}_{ij}) - t_{ij}(\hat{y}_{ij})) \hat{y}_{ij} = b_{ij}(\hat{x}_{ij} - \hat{y}_{ij}) \hat{y}_{ij}\,.}
Now, the function $b_{ij}(\hat{x}_{ij} - z)z$ is a concave quadratic function in $z$, which obtains its maximum when $z = \hat{x}_{ij} / 2$.
Therefore
\labeleqn{nextbound}{(t_{ij}(\hat{x}_{ij}) - t_{ij}(\hat{y}_{ij})) \hat{y}_{ij} \leq \frac{b_{ij} (\hat{x}_{ij})^2}{4} \leq \frac{1}{4} (a_{ij} + b_{ij} \hat{x}_{ij}) \hat{x}_{ij} = t_{ij}(\hat{x}_{ij}) \hat{x}_{ij}\,,}
proving the result.
\end{proof}

Furthermore, this bound is tight.
Consider the more extreme version of the Knight-Pigou-Downs network shown in Figure~\ref{fig:poatight} where the demand is 1 unit, the travel time on the upper link is 1 minute, and the travel time on the bottom link is $t^\downarrow = x^\downarrow$.
The user equilibrium solution is $x^\uparrow = 0$, $x^\downarrow = 1$, when both links have equal travel times of 1 minute, and the total system travel time is 1.
You can verify that the system optimum solution is $x^\uparrow = x^\downarrow = \frac{1}{2}$, when the total system travel time is $\frac{3}{4}$.
Thus the ratio between the user equilibrium and system optimum total system travel times is $\frac{4}{3}$.

\begin{figure}
\begin{center}
\begin{tikzpicture}[->,>=stealth',shorten >=1pt,auto,node distance=3cm,
thick,main node/.style={circle,,draw}]
\node[main node] (1) {$1$};
\node[main node] (2) [right of=1] {$2$};
\path (1) edge [bend left] node[above] {$1$} (2);
\path (1) edge [bend right] node[below] {$x$} (2);
\end{tikzpicture}
$d^{12} = 1$
\end{center}
\caption{The price of anarchy bound is tight.\label{fig:poatight}}
\end{figure}

You may be wondering if a price of anarchy can be found when we relax the assumption that the link performance functions are affine.
In many cases, yes; for instance, if the link performance functions are quadratic, then the price of anarchy is $\frac{3 \sqrt{3}}{3 \sqrt{3} - 2}$, if cubic, then the price of anarchy is $\frac{4 \sqrt[3]{4}}{4 \sqrt[3]{4} - 3}$, and so on.
In all of these cases the modified Knight-Pigou-Downs network can show that this bound is tight.
On the other hand, if the link performance functions have a vertical asymptote (e.g., $t_{ij} = (u_{ij} - x_{ij})^{-1}$), then the ratio between the total system travel times at user equilibrium and system optimum may be made arbitrarily large.
\index{price of anarchy|)}
\index{user equilibrium!relation to system optimum|)}
\index{system optimum!relation to user equilibrium|)}
\index{user equilibrium|)}
\index{system optimum|)}
\index{assignment rule|)}

\section{Historical Notes and Further Reading}
\label{sec:tap_historical}

The variational inequality formulation of the traffic assignment problem was presented by \cite{smith79}, \cite{dafermos80}, and \cite{smith83a}.
A convex optimization formulation of the traffic assignment problem was first presented in \cite{beckmann56}, although it did not take exactly the same form as in this chapter.
That reference actually presented an \emph{elastic demand} equilibrium formulation, a generalization which will be treated more in Section~\ref{sec:elasticdemand}.
The formulation given in this chapter is due to Stella Dafermos\index{Dafermos, Stella}, and can first be found in \cite{dafermos68} and \cite{dafermos69}.
The fixed point formulation described in the optional section relies on a more general fixed point theorem due to \cite{kakutani41}.
Although not described in this book, the traffic assignment problem can also be formulated as a nonlinear complementarity problem~\citep{aashtiani80}.\index{static traffic assignment!formulation!as nonlinear complementarity problem}

\index{user equilibrium!relation to system optimum|(}
\index{system optimum!relation to user equilibrium|(}
\cite{rossi89} proposed the use of entropy to distinguish a ``most likely'' path flow solution among all of those satisfying the principle of user equilibrium.
\cite{bargera99} and \cite{bargera06} discuss how the more intuitive (but slightly weaker) proportionality condition can be derived from entropy maximization.
\cite{bargera07} conducted a case study considering the practical effects of nonunique path flow solutions, looking at the range of path flows that correspond to user equilibrium.
In the Chicago network, they found wide variation among the path flow solutions that correspond to user equilibrium, concluding that finding high-entropy solutions is important for analyses that require path flows, and not just link flows.
\index{user equilibrium!relation to system optimum|)}
\index{system optimum!relation to user equilibrium|)}

The idea of aggregating solutions by origins or destination dates to \cite{dial_diss} and \cite{dial71}, in the context of a network loading with perception errors.\index{stochastic network loading}\index{static traffic assignment!bush-based solution}
That this idea is applicable to the traffic assignment problem was independently discovered by \cite{dial99_bobtail}, in an unpublished manuscript, and \cite{bargera_diss} in a doctoral thesis; see \cite{bargera02} and \cite{dial06a} for more concise journal publications reporting the main ideas.
Lemma 3 of \cite{bargera02} establishes the more general form of Proposition~\ref{prp:eqmbushacyclic} for the case of zero-cost links.
\index{static traffic assignment!properties!acyclicity}

The idea that user equilibrium and system optimum are generally different predates the formulations given in this chapter, having been identified by \cite{pigou20} and \cite{knight24}.
\cite{dafermos69} showed how they have the same mathematical structure, as demonstrated in Section~\ref{sec:uesoequivalence}.
In practice, the relative difference between the solutions in terms of total travel time is relatively small, typically a few percent; see, for instance, analyses on networks representing the cities of Istanbul~\citep{gunay96} and Chicago~\citep{boycexiong04}.
In absolute terms, the savings are more significant.
The solutions can also differ more substantially in terms of flows on specific links.
Nevertheless, the solutions are often similar enough that relative rankings of network design strategies are roughly consistent regardless of whether the user equilibrium or system optimum assignment rule is used, and sometimes the system optimum assignment is easier to solve in such contexts~\citep{leblanc84}.
The paradoxes from Chapter~\ref{chp:equilibrium} also illustrate this distinction.

The work on the price of anarchy grew out of an attempt to see just how different user equilibrium and system optimum could possibly be --- not just in typical networks in practice, or in the paradoxes chosen to illustrate a point --- but over \emph{all} instances of a certain type.
The first results on the price of anarchy in traffic assignment are due to \cite{roughgarden02}, and the specific proof given in this book is from \cite{correa04}.
The former reference contains a number of related results, including bounds for other classes of link performance functions, and performance bounds for other optimization problems that have traffic assignment as a subproblem.
The assignment rule involving perception errors was first developed by \cite{vonfalkenhausen66} for constant link travel times, and by \cite{daganzo77} for flow-dependent link travel times; more discussion on this model and its historical context are found in Sections~\ref{sec:sue} and~\ref{sec:staticextensions_references}
The bounded rationality formulations and results presented in this chapter are drawn from \cite{mahmassani87}, \cite{lou10}, and \cite{di13}.

The traffic assignment problem presented in this chapter has been extended in many directions.
Chapter~\ref{chp:staticextensions} discusses several of these extensions, and the historical notes in Section~\ref{sec:staticextensions_references} mention even more.

\section{Exercises}
\label{exercises_tap}

\begin{enumerate}
\item \diff{34} Consider a network with two parallel links connecting a single origin to a single destination; the link performance function on each link is $2 + x$ and the total demand is $d = 10$.
\begin{enumerate}
   \item Write the equations and inequalities defining the set of feasible path flows $H$, and draw a sketch.
   \item What are the vectors $-\mathbf{C}(\mathbf{h})$ for the following path flow vectors?  (1) $\mathbf{h} = [0,10]$  (2) $\mathbf{h} = [5,5]$ (3) $\mathbf{h} = [10,0]$  Draw these vectors at these three points on your sketch.
   \item For each of the vectors from part (b), identify the point $\mathrm{proj}_H (\mathbf{h} - \mathbf{C}(\mathbf{h}))$ and include these on your sketch.
\end{enumerate}
\item \diff{45} Consider a network with two parallel links connecting a single origin to a single destination; the link performance function on the first link is $10 + x$, and the link performance function on the second link is $x^2 / 20$.
The total demand is $d = 30$.
\begin{enumerate}
	\item Demonstrate that the solution $\mb{\hat{h}} = \vect{20 & 10}$ is not an equilibrium, and then provide a vector $\mb{h}$ showing that the variational inequality~\eqn{basicvi} is not satisfied.
	\item Demonstrate that the solution $\mb{\hat{h}} = \vect{10 & 20}$ is an equilibrium, and prove that the variational inequality~\eqn{basicvi} is satisfied no matter what $\mb{h} \in H$ is.
\end{enumerate}
\item \diff{43} For the Braess network of Figure~\ref{fig:braesssimple}(b), write down the optimality conditions corresponding to Beckmann's formulation~\eqn{beckmann}--\eqn{nonneg}.
\item \diff{27} Which of the following multifunctions $F$ are closed?  Each of these multifunctions maps $[0,1]$ to subsets of $[0,1]$.
Draw sketches of each of these multifunctions.
\label{ex:multifunctionpractice}
\begin{enumerate}
	\item $F(x) = \myc{y : 0 \leq y \leq x}$
	\item $F(x) = \myc{y : 0 < y < 1}$
	\item $F(x) = \myc{0} \cup \myc{1 - x} \cup \myc{1}$
	\item $F(x) = \myc{x^2}$
\end{enumerate}
\item \diff{62} Show that if $f(x)$ is a continuous function with a compact domain, the single-valued ``multifunction'' $F(x) = \myc{f(x)}$ is closed.
\item \diff{23} For each of the multifunctions in Exercise~\ref{ex:multifunctionpractice}, identify all of the fixed points.
\item \diff{62} Specify the multifunction $\mb{H}(\mb{c})$ for the Braess network (Figure~\ref{fig:braesssimple}b), and identify all of its fixed points.
\item \diff{42} Complete the proof of Theorem~\ref{thm:eqmkakutani} by showing that the graph of $F$ is indeed closed.
\item \diff{51} Create a simple network where one or more link performance functions are \emph{not} continuous, and where no user equilibrium solution exists.
(Don't worry about creating ``realistic'' functions for this problem.)
\item \diff{51} Create a simple network where one or more link performance functions are \emph{not} strictly increasing, and where the user equilibrium link flows are not unique.
\item \diff{11} Show that the BPR link delay function~\eqn{bpr} does not satisfy the condition $t'_{ij}(x_{ij}) > 0$ for all $x_{ij}$ if $\beta > 1$.
\label{ex:bprfail}
\item \diff{46} In spite of Exercise~\ref{ex:bprfail}, the BPR link delay functions are strictly increasing, and the resulting link flow equilibrium solution is still unique.
Generalize the proof of Proposition~\ref{prp:linkuniqueness} to handle the case when the link delay functions are differentiable and strictly increasing, by showing that the Beckmann function is still strictly convex even if $t'_{ij}(x_{ij})$ is not always strictly positive.
\item \diff{51} Show that the user equilibrium and system optimum link flows are the same if there is no congestion, that is, if $t_a(x_a) = \tau_a$ for some constant $\tau_a$.
\item \diff{10} Consider the network of Figure~\ref{fig:decomposition}, ignoring the labels next to the link in the figure.
The demand in this network is given by $d^{19} = 400$.
Does the following path-flow solution satisfy the principle of user equilibrium?  100 vehicles on path [1,2,3,6,9], 100 vehicles on path [1,2,3,6,9], 100 vehicles on path [1,4,5,6,9], and 100 vehicles on path [1,4,5,2,3,6,9].
You can answer this without any calculation.
\label{ex:decomposition1}
\item \diff{45} In the network in Figure~\ref{fig:decomposition}, the demand is given by $d^{13} = d^{19} = d^{17} = d^{91} = d^{93} = d^{97} = 100$.
The flow on each link is shown in the figure.
\label{ex:decomposition2}
\begin{enumerate}
   \item Find a path flow vector $\mb{h}$ which corresponds to these link flows.
   \item Show the resulting origin-based link flows $\mb{x^1}$ and $\mb{x^9}$ from the two origins in the network.
   \item Show that the links used by these two origins form an acyclic subnetwork by finding topological orders for each subnetwork.
   \item Determine whether these link flows satisfy the principle of user equilibrium.
\end{enumerate}
\stevefig{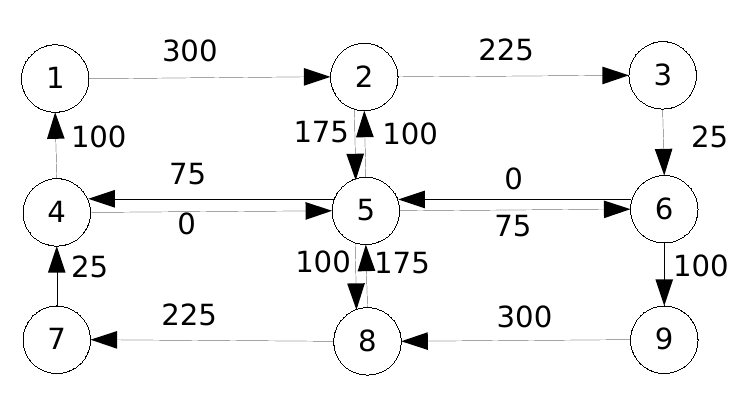}{Network for Exercises~\ref{ex:decomposition1} and~\ref{ex:decomposition2}.
Each link has link performance function $10 + x_{ij}$}{0.6\textwidth}
\item \diff{35} Consider the network in Figure~\ref{fig:prop1}, along with the given (equilibrium) link flows.
There is only one OD pair, from node 1 to node 3.
Identify three values of path flows which are consistent with these link flows, \emph{in addition} to the most likely (entropy-maximizing) path flows.
\label{ex:prop1}
\stevefig{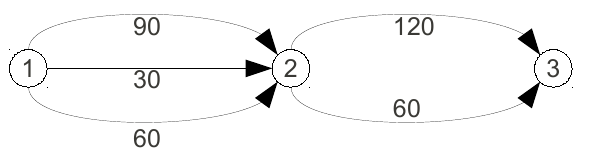}{Network for Exercise~\ref{ex:prop1}.}{0.6\textwidth}
\item \diff{37}  In the network shown in Figure~\ref{fig:prop2}, 320 vehicles travel from A to C, 640 vehicles travel from A to D, 160 vehicles travel from B to C, and 320 vehicles travel from B to D.
The equilibrium link flows are shown.
\label{ex:prop2}
\begin{enumerate}[(a)]
   \item Give a path flow solution which satisfies proportionality and produces the equilibrium link flows.
   \item At the proportional solution, what fraction of flow on the top link connecting 3 and 4 is from origin A?
   \item Is there a path flow solution which produces the equilibrium link flows, yet has \emph{no} vehicles from origin A on the top link connecting 3 and 4?  If so, list it.
If not, explain why.
\end{enumerate}
\stevefig{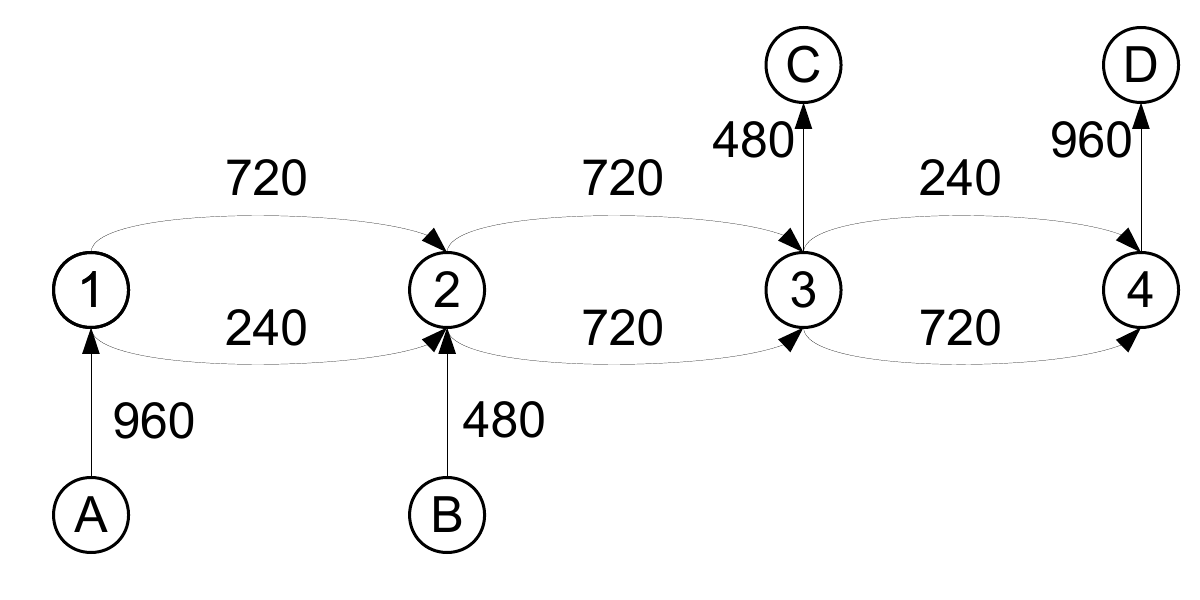}{Network for Exercise~\ref{ex:prop2}.}{0.6\textwidth}
\item \diff{38} Consider the network in Figure~\ref{fig:hwprop}, where 2 vehicles travel from 1 to 4 (with a value of time of \$20/hr), and 4 vehicles travel from 2 to 4 (with a value of time of \$8/hr).
The equilibrium volumes are 3 vehicles on Link 1, and 3 vehicles on Link 2.
\label{ex:hwprop}
\begin{enumerate}
   \item Assuming that the vehicles in this network are discrete and cannot be split into fractions, identify every combination of path flows which give the equilibrium link volumes (there should be 20).
Assuming each combination is equally likely, show that the proportional division of flows has the highest probability of being realized.
   \item What is the average value of travel time on Link 1 at the most likely path flows?  What are the upper and lower limits on the average value of travel time on this link?
\end{enumerate}
\stevefig{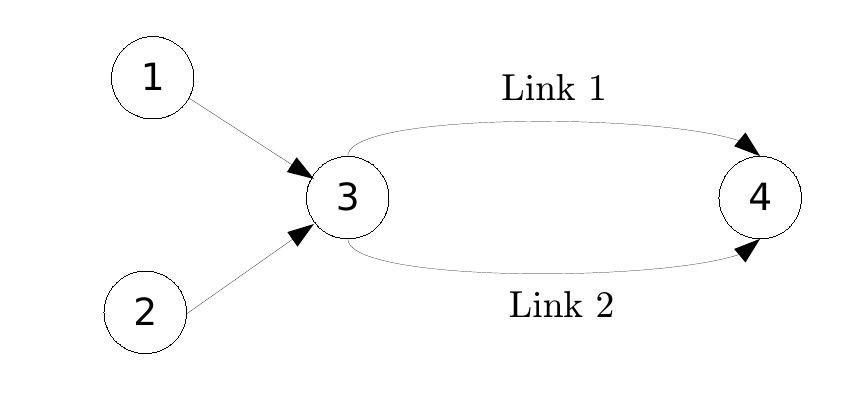}{Network for Exercise~\ref{ex:hwprop}.}{0.6\textwidth}
\item \diff{65} Derive the optimality conditions for the system optimum assignment, and provide an interpretation of these conditions which intuitively relates them to the concept of system optimality.
\item \diff{33} Find the system optimum assignment in the Braess network (Figure~\ref{fig:braesssimple}b), assuming a demand of 6 vehicles from node 1 to node 4.
\item \diff{51} Consider a network where every link performance function is linear, of the form $t_{ij} = a_{ij} x_{ij}$.
Show that the user equilibrium and system optimum solutions are the same.
\item \diff{22} Prove that the system optimum objective function is convex if the link performance functions are all nonnegative, nondecreasing, differentiable, and convex. 
\item \diff{34} Calculate the ratio of the total system travel time between the user equilibrium and system optimum solutions in the Braess network (Figure~\ref{fig:braesssimple}b).
\item \diff{46} Find the set of boundedly rational assignments in the Braess network (Figure~\ref{fig:braesssimple}b).
\end{enumerate}

\textbf{Proof of Theorem~\ref{thm:proportionality}:}
\label{sec:proportionalityproof}

\index{entropy!and proportionality|(}
Begin by forming the Lagrangian for the entropy-maximizing problem:\label{not:betaij}\label{not:gammars}
\begin{multline}
\label{eqn:entropylagrangian}
\mathcal{L}(\mathbf{h}, \bm{\beta}, \bm{\gamma}) = -\sum_{(r,s) \in Z^2} \sum_{\pi \in \hat{\Pi}^{rs}} h_\pi \log \myp{\frac{h_\pi}{d^{rs}}} + \sum_{(i,j) \in A} \beta_{ij} \myp{\hat{x}_{ij} - \sum_{\pi \in \Pi} \delta_{ij}^\pi h_\pi} \\ + \sum_{(r,s) \in Z^2} \gamma_{rs} \myp{d^{rs} - \sum_{\pi \in \hat{\Pi}^{rs}} h_\pi} 
\,,
\end{multline}
using $\beta$ and $\gamma$ to denote the Lagrange multipliers.
Note that the nonnegativity constraint can be effectively disregarded since the objective function is only defined for strictly positive $\mathbf{h}$.
Therefore at the optimum solution the partial derivative of $\mathcal{L}$ with respect to any path flow must vanish:
\labeleqn{pdrLh}{\pdr{\mathcal{L}}{h_\pi} = -1 - \log\myp{\frac{h_\pi}{d^{rs}}} - \sum_{(i,j) \in A} \delta_{ij}^\pi \beta_{ij} - \gamma_{rs} = 0\,, }
where $(r,s)$ is the OD pair connected by path $\pi$.
Solving for $h^\pi$ we obtain
\labeleqn{Hfromlagrange}{h^\pi = d^{rs} \exp\myp{-1 - \sum_{(i,j) \in A} \delta_{ij}^\pi \beta_{ij} - \gamma_{rs}}\,.}
Likewise we have
\labeleqn{pdrLgamma}{\pdr{\mathcal{L}}{\gamma_{rs}} = d^{rs} - \sum_{\pi' \in \hat{\Pi}^{rs}} h_{\pi'} = 0\,,}
so
\labeleqn{gammafromlagrange}{d^{rs} = \sum_{\pi' \in \hat{\Pi}^{rs}} h_{\pi'} = d^{rs} \exp\myp{-1 - \gamma_{rs}} \sum_{\pi \in \hat{\Pi}^{rs}} \exp\myp{-\sum_{(i,j) \in A} \delta_{ij}^\pi \beta_{ij} }\,,}
substituting the result from~\eqn{Hfromlagrange}.
Therefore we can solve for $\gamma_{rs}$:
\labeleqn{gammafromsubstitution}{\gamma_{rs} = -1 + \log \myp{\sum_{\pi' \in \hat{\Pi}^{rs}} \exp\myp{-\sum_{(i,j) \in A} \delta_{ij}^{\pi'} \beta_{ij} }}\,.}
Finally, substituting~\eqn{gammafromsubstitution} into~\eqn{Hfromlagrange} and simplifying, we obtain
\labeleqn{finalh}{h^\pi = \frac{d^{rs}}{\sum_{\pi' \in \hat{\Pi}^{rs}} \exp\myp{-\sum_{(i,j) \in A} \delta_{ij}^{\pi'} \beta_{ij} }} \exp \myp{- \sum_{(i,j) \in A} \delta_{ij}^\pi \beta_{ij}}  \,.}
Noting that the fraction in~\eqn{finalh} only depends on the OD pair $(r,s)$, we can simply write
\labeleqn{finalhnoOD}{h^\pi = K_{rs} \exp \myp{ - \sum_{(i,j) \in A} \delta_{ij}^\pi \beta_{ij}} \,,  }
where $K_{rs}$\label{not:Krs} is a constant associated with OD pair $(r,s)$.

Let $A_1$ be the set of links not in either $\pi_1$ or $\pi_2$, $A_2$ the set of links common to both paths, $A_3$ the links in $\pi_1$ but not $\pi_2$, and $A_4$ the links in $\pi_2$ but not $\pi_1$.
Then the ratio $h_1/h_2$ can be written
\begin{align}
	\frac{h_1}{h_2} &= \frac{K_{rs} \exp \myp{ - \sum_{(i,j) \in A} \delta_{ij}^{\pi_1} \beta_{ij}} }{K_{rs} \exp \myp{ -\sum_{(i,j) \in A} \delta_{ij}^{\pi_2} \beta_{ij}}} \\
							&= \frac{K_{rs} \exp \myp{ \begin{array}{l} - \sum_{(i,j) \in A_1} \delta_{ij}^{\pi_1} \beta_{ij} - \sum_{(i,j) \in A_2} \delta_{ij}^{\pi_1} \beta_{ij} \\ \mbox{\hspace{0.25in}} - \sum_{(i,j) \in A_3} \delta_{ij}^{\pi_1} \beta_{ij} - \sum_{(i,j) \in A_4} \delta_{ij}^{\pi_1} \beta_{ij} \end{array} }}{K_{rs} \exp \myp{ \begin{array}{l} -\sum_{(i,j) \in A_1} \delta_{ij}^{\pi_2} \beta_{ij} -\sum_{(i,j) \in A_2} \delta_{ij}^{\pi_2} \beta_{ij} \\ \mbox{\hspace{0.25in}} -\sum_{(i,j) \in A_3} \delta_{ij}^{\pi_2} \beta_{ij} -\sum_{(i,j) \in A_4} \delta_{ij}^{\pi_2} \beta_{ij} \end{array} }}  \\
							&= \frac{K_{rs} \exp \myp{ - \sum_{(i,j) \in A_2} \delta_{ij}^{\pi_1} \beta_{ij} - \sum_{(i,j) \in A_3} \delta_{ij}^{\pi_1} \beta_{ij} } }{K_{rs} \exp \myp{ -\sum_{(i,j) \in A_2} \delta_{ij}^{\pi_2} \beta_{ij} -\sum_{(i,j) \in A_4} \delta_{ij}^{\pi_2} \beta_{ij} }}  \\
							&= \frac{K_{rs} \exp \myp{ - \sum_{(i,j) \in A_2} \delta_{ij}^{\pi_1} \beta_{ij}} \exp \myp{- \sum_{(i,j) \in A_3} \delta_{ij}^{\pi_1} \beta_{ij} } }{K_{rs} \exp \myp{ -\sum_{(i,j) \in A_2} \delta_{ij}^{\pi_2} \beta_{ij}} \exp \myp{ - \sum_{(i,j) \in A_4} \delta_{ij}^{\pi_2} \beta_{ij} }}   \\
							&= \frac{\exp \myp{- \sum_{(i,j) \in A_3} \delta_{ij}^{\pi_1} \beta_{ij} }}{\exp \myp{- \sum_{(i,j) \in A_4} \delta_{ij}^{\pi_1} \beta_{ij} }}
\,.
\end{align}
where the steps of the derivation respectively involve substituting~\eqn{finalh}, expanding $A$ into $A_1 \cup A_2 \cup A_3 \cup A_4$, using the definitions of the $A_i$ sets to identify sums where $\delta_{ij}^\pi = 0$, splitting exponential terms, and canceling common factors.

Thus in the end we have
\labeleqn{hratio}{\frac{h_1}{h_2} = \frac{\exp \myp{- \sum_{(i,j) \in A_3} \delta_{ij}^{\pi_1} \beta_{ij} }}{\exp \myp{- \sum_{(i,j) \in A_4} \delta_{ij}^{\pi_1} \beta_{ij}}}\,,}
regardless of the OD pair $h_1$ and $h_2$ connect, and this ratio only depends on $A_3$ and $A_4$, that is, the pairs of alternate segments distinguishing $\pi_1$ and $\pi_2$.
\index{entropy!and proportionality|)}

\chapter{Algorithms for Traffic Assignment}
\label{chp:solutionalgorithms}

\index{static traffic assignment!algorithms|(}
This chapter presents algorithms for solving the basic traffic assignment problem (TAP), which was defined in Chapter~\ref{chp:trafficassignmentproblem} as the solution $\mb{\hat{x}}$ to the variational inequality
\labeleqn{tapvi}{\mb{t}(\mb{\hat{x}}) \cdot (\mb{\hat{x}} - \mb{x}) \leq 0 \qquad \forall \mb{x} \in X\,,}
which can also be expressed in terms of path flows $\mb{\hat{h}}$ as
\labeleqn{tapvipaths}{\mb{c}(\mb{\hat{h}}) \cdot (\mb{\hat{h}} - \mb{h}) \leq 0 \qquad \forall \mb{h} \in H\,,}
or equivalently as the solution $\mb{\hat{x}}$ to the optimization problem
\begin{align}
\min_{\mathbf{x},\mathbf{h}}  \quad & \sum_{(i,j) \in A} \int_0^{x_{ij}} t_{ij}(x) dx     & \label{eqn:tapoptstart} \\
\mathrm{s.t.}                 \quad & x_{ij} = \sum_{\pi \in \Pi} h^\pi \delta_{ij}^{\pi} & \forall (i,j) \in A  \\
                                    & \sum_{\pi \in \Pi^{rs}} h^\pi = d^{rs}              & \forall (r,s) \in Z^2  \\
                                    & h^\pi \geq 0                                        & \forall \pi \in \Pi \label{eqn:tapoptend} 
\end{align}

While there are general algorithms for variational inequalities or nonlinear optimization problems, TAP involves tens of thousands or even millions of variables and constraints for practical problems.
So, these general algorithms are outperformed by specialized algorithms which are designed to exploit some specific features of TAP.
The most significant features to exploit are the network structure embedded in the constraints, and the fact that the objective function is separable by link.

This chapter presents four types of algorithms for TAP.
The first three are aimed at finding an equilibrium solution (either link flows $\mb{\hat{x}}$ or path flows $\mb{\hat{h}}$).
Broadly speaking, these algorithms can be divided into link-based, path-based, and bush-based algorithms, according to the way the solution is represented.
As discussed in Section~\ref{sec:maxentropy}, this is not trivial even though it plays a major role in certain types of applications, such as select link analysis.

The chapter focuses entirely on the basic TAP, and not the alternative assignment rules from Section~\ref{sec:alternaterules}.
System optimal assignment can be transformed mathematically into a user equilibrium problem with modified cost functions, so all of these algorithms can be easily adapted for the system optimal problem, but most of these algorithms cannot be directly applied to the variations with bounded rationality or perception errors.

Section~\ref{sec:algointroduction} provides an introduction, justifying the need for efficient algorithms for TAP and discussing issues such as convergence criteria which are relevant to all algorithms.
The next three sections respectively present link-based, path-based, and bush-based algorithms for finding user equilibrium solutions.

\section{Introduction to Assignment Algorithms}
\label{sec:algointroduction}

This introductory section touches on three topics: first, the advantages and disadvantages of link-based, path-based, and bush-based algorithms; second, the general framework for solving equilibrium problems; and third, the question of convergence criteria, which is pertinent to all of these.

\subsection{Why do different algorithms matter?}

It may not be clear why we need to present so many different algorithms for the same problem.
Why not simply present the one ``best'' algorithm for solving TAP?
Link-based, path-based, and bush-based algorithms all exhibit advantages and disadvantages relative to each other.
In this introductory section, these algorithms are compared qualitatively.
As you read through the following sections, which include details and specifications of each algorithm, keep these concepts in mind.

\index{static traffic assignment!algorithms!link-based}
\textbf{Link-based algorithms} only keep track of the aggregate link flows $\mb{x}$, not the path flows $\mb{h}$ which lead to these link flows.
Since the number of links in a large network is much less than the number of paths, link-based algorithms are very economical in terms of computer memory.
For this reason, link-based algorithms were the first to be used in practice decades ago, when computer memory was very expensive.
Link-based algorithms also tend to be easier to parallelize.
With modern desktop machines containing multiple cores, this parallelization can reduce computation time significantly.
Finally, link-based algorithms tend to be easy to code and implement (even in a spreadsheet).
A significant drawback is that they tend to be slow, particularly when high precision is demanded.
The first few iterations of link-based algorithms make good progress, but they stall quickly, and final convergence to the equilibrium solution can be extremely slow.

\index{static traffic assignment!algorithms!path-based}
\textbf{Path-based algorithms}, by contrast, keep track of the path flow vector $\mb{h}$, which can be used to generate the link flows $\mb{x}$ whenever needed.
Although this representation requires more memory than a link-based solution, it also retains a great deal of information which is lost when one aggregates path flows to link flows.
This information can be exploited to converge much faster than link-based algorithms, especially when a very precise solution is needed.
Because the number of paths in a network is so large, much of the effort involved in coding path-based algorithms is associated with clever data structures and algorithm schemes which are aimed at minimizing the number of paths which need to be stored.
This increases the complexity of the coding of these algorithms, even if the algorithmic concepts are not difficult.

\index{static traffic assignment!algorithms!bush-based}
\textbf{Bush-based algorithms} are the most recent developed, and aim to offer speed comparable to path-based algorithms, while requiring less memory (although still significantly more than link-based algorithms).
Bush-based algorithms do this by selectively aggregating the path flows: as shown in Section~\ref{sec:bushes}, at equilibrium the set of paths used by all travelers from the same origin (or traveling to the same destination) forms an acyclic subnetwork (a bush).
Bush-based algorithms aim to identify these bushes for each origin and destination, exploiting the fact that calculations in acyclic networks are very fast.
Bush-based algorithms also tend to return higher-entropy solutions than path-based algorithms, which is important when interpreting the path flow solution.
The downside to these algorithms is that their design involves more complexity, and the success of these algorithms depends highly on implementation and programming skill.

\subsection{Framework}

\index{static traffic assignment!algorithms!framework|(}
It turns out that none of these solution methods get to the right answer immediately, or even after a finite number of steps.
There is no ``step one, step two, step three, and then we're done'' recipe for solving large-scale equilibrium problems.\index{consistency}
Instead, an iterative approach is used where we start with some feasible assignment (link or path), and move closer and closer to the equilibrium solution as you repeat a certain set of steps over and over, until you're ``close enough'' to quit and call it good.
One iterative algorithm you probably saw in calculus was Newton's method for finding zeros of a function.
In this method, one repeats the same step over and over until the function is sufficiently close to zero.

Broadly speaking, all equilibrium solution algorithms repeat the following three steps:
\begin{enumerate}
\item Find the shortest (least travel time) path between each origin and each destination.
\item Shift travelers from slower paths to faster ones.
\item Recalculate link flows and travel times after the shift, and return to step one unless we're close enough to equilibrium.
\end{enumerate}

The shortest path computation can be done quickly and efficiently even in large networks, as was described in Section~\ref{sec:shortestpath}.
The third step is even more straightforward, and is nothing more than re-evaluating the link performance functions on each link with the new volumes.
The second step requires the most care; the danger here is shifting either too few travelers onto faster paths, or shifting too many.
If we shift too few, then it will take a long time to get to the equilibrium solution.
On the other hand, systematically shifting too many can be even more dangerous, because it creates the possibility of ``infinite cycling'' and never finding the true equilibrium.

\begin{figure}
\begin{center}
\begin{tikzpicture}[->,>=stealth',shorten >=1pt,auto,node distance=3cm,
thick,main node/.style={circle,,draw}]

\node[main node] (1) {$1$};
\node[main node] (2) [right of=1] {$2$};

\path (1) edge [bend left] node[above] {$10 + x_1$} (2);
\path (1) edge [bend right] node[below] {$20 + x_2$} (2);

\end{tikzpicture}

50 vehicles travel from node 1 to node 2
\end{center}
\caption{Two-link example for demonstration. \label{fig:twolink}}
\end{figure}

In the simple example in Figure~\ref{fig:twolink}, by inspection the equilibrium is for thirty travelers to choose the top route, and twenty to choose the bottom route, with an equal travel time of 40 minutes on both paths.
Solving this example using the above process, initially (i.e., with nobody on the network) the fastest path is the top one (step one), so let's assign all 50 travelers onto the top path (step two).
Performing the third step, we recalculate the travel times as 60 minutes on the top link, and 20 on the bottom.
This is not at all an equilibrium, so we go back to the first step, and see that the bottom path is now faster, so we have to shift some people from the top to the bottom.
If we wanted, we could shift travelers one at a time, that is, assigning 49 to the top route and 1 to the bottom, seeing that we still haven't found equilibrium, so trying 48 and 2, then 47 and 3, and so forth, until finally reaching the equilibrium with 30 and 20.
Clearly this is not efficient, and is an example of shifting too few travelers at a time.

At the other extreme, let's say we shift \emph{everybody} onto the fastest path in the second step.
That is, we go from assigning 50 to the top route and 0 to the bottom, to assigning 0 to the top and 50 to the bottom.
Recalculating link travel times, the top route now has a travel time of 10 minutes, and the bottom a travel time of 70.
This is even worse!\footnote{By ``worse'' we mean farther from equilibrium.}
Repeating the process, we try to fix this by shifting everybody back (50 on top, 0 on bottom), but now we're just back in the original situation.
If we kept up this process, we'd keep bouncing back and forth between these solutions.
This is clearly worse than shifting too few, because we never reach the equilibrium no matter how long we work!
You might think it's obvious to detect if something like this is happening.
With this small example, it might be.
Trying to train a computer to detect this, or trying to detect cycles with over millions OD pairs (as is common in practice), is much harder.
The key step in all of the algorithms for finding equilibria is determining how much flow to shift.\index{static traffic assignment!algorithms!framework|)}

\subsection{Convergence criteria}

\index{gap function|(}
A general issue is how one chooses to stop the iterative process, that is, how one knows when a solution is ``good enough'' or close enough to equilibrium.
This is called a \emph{convergence criterion}.
Many convergence criteria have been proposed over the years; perhaps the most common in practice is the \emph{relative gap}, which is defined first.
Unfortunately, the relative gap has been defined in several different ways, so it is important to be familiar with all of the definitions.

\index{gap function!relative gap|(}
Remembering that the multiplier $\kappa^{rs}$ represents the time spent on the fastest path between origin $r$ and destination $s$, one definition of the relative gap $\gamma_1$\label{not:gammagap} is
\labeleqn{relgap}{
\gamma_1 = \frac{\sum_{(i,j) \in A} t_{ij} x_{ij} }{\sum_{(r,s) \in Z^2} \kappa^{rs} d^{rs} } - 1 = \frac{\mathbf{t} \cdot \mathbf{x}}{\bm{\kappa} \cdot \mathbf{d}} - 1
\,.}
The numerator of the fraction is the total system travel time\index{total system travel time} (TSTT).
The denominator is called the \emph{shortest path travel time}\index{shortest path travel time} (SPTT), and reflects what the total system travel time would theoretically be if all travelers could be shifted to the current shortest paths \emph{without changing the travel times}.
The relative gap is always nonnegative, and it is equal to zero if and only if the flows $\mb{x}$ satisfy the principle of user equilibrium.
It is these properties which make the relative gap a useful convergence criterion: once it is close enough to zero, our solution is ``close enough'' to equilibrium.
For most practical purposes, a relative gap of $10^{-4}$--$10^{-6}$ is small enough.

A second definition of the relative gap $\gamma_2$ is based on the Beckmann function itself.
Let $f$ denote the Beckmann function, and $\hat{f}$\label{not:fhat} its value at equilibrium (which is a global minimum).
In many algorithms, given a current solution $\mb{x}$, it is not difficult to generate upper and lower bounds on $\hat{f}$ based on the current solution, respectively denoted $\bar{f}(\mb{x})$ and $\ubar{f}(\mb{x})$.\label{not:fbar}
A trivial upper bound is its value at the current solution: $\bar{f} = f(\mb{x})$, since clearly $\hat{f} \leq f(\mb{x})$ for any feasible link assignment $\mb{x}$.
Sometimes, a corresponding lower bound can be identified as well.
Assuming that these bounds can become tighter over time, and that in the limit both $\bar{f} \rightarrow \hat{f}$ and $\ubar{f} \rightarrow \hat{f}$, the difference or gap $\bar{f}(\mb{x}) - \ubar{f}(\mb{x})$ can be used as a convergence criterion.
These values are typically normalized, leading to one definition of the relative gap:
\labeleqn{relgap2}{
\gamma_2 = \frac{\bar{f} - \ubar{f}}{\ubar{f}}
\,,}
or a slightly modified version
\labeleqn{relgap3}{
\gamma_3 = \frac{\bar{f} - \max \ubar{f}}{\max \ubar{f}}
\,,}
where $\max \ubar{f}$ is the greatest lower bound found to date, in case the sequence of $\ubar{f}$ values is not monotone over iterations.
A disadvantage of these definitions of the relative gap is that different algorithms calculate these upper and lower bounds differently.
While they are suitable as termination criteria in an algorithm, it is not possible (or at least not easy) to directly compare the relative gap calculated by one algorithm to that produced by another to assess which of two solutions is closer to equilibrium.

\index{relative gap|see {gap function, relative gap}}
One drawback of the relative gap (in all of its forms) is that it is unitless and does not have an intuitive meaning.
Furthermore, it can be somewhat confusing to have several slightly different definitions of the relative gap, even though they all have the same flavor.\index{gap function!relative gap|)}
A more recently proposed metric is the \emph{average excess cost}\index{gap function!average excess cost|(}\label{not:aec}, defined as
\labeleqn{aec}{
AEC =
\frac{\mathbf{t} \cdot \mathbf{x} - \bm{\kappa} \cdot \mathbf{d}}{\mathbf{d} \cdot \mathbf{1}}
= \frac{TSTT - SPTT}{\mathbf{d} \cdot \mathbf{1}}
\,.}
This quantity represents the average difference between the travel time on each traveler's \emph{actual} path, and the travel time on the \emph{shortest} path available to him or her.
Unlike the relative gap, $AEC$ has units of time, and is thus easier to interpret.\index{gap function!average excess cost|)}\index{average excess cost|see {gap function, average excess cost}}

Another convergence measure with time units is the \emph{maximum excess cost}\index{gap function!maximum excess cost}\index{maximum excess cost|see {gap function, maximum excess cost}}\label{not:mec}, which relates directly to the principle of user equilibrium.
The maximum excess cost is defined as the largest amount by which a used path's travel time exceeds the shortest path travel time available to that traveler:
\labeleqn{mec}{
MEC = \max_{(r,s) \in Z^2} \myc{\max_{\pi \in \Pi^{rs} : h^\pi > 0} \myc{ c^\pi - \kappa^{rs}}}
\,.}
This is often a few orders of magnitude higher than the average excess cost.
One disadvantage of the maximum excess cost is that it is only applicable when the path flow solution is known.
This is easy in path-based or bush-based algorithms.
However, since many path-flow solutions correspond to the same link-flow solution (cf.\ Section~\ref{sec:maxentropy}), $MEC$ is not well suited for link-based algorithms.\index{gap function|)}

Finally, this section concludes with two convergence criteria which are inferior to those discussed thus far.
The first is to simply use the Beckmann function itself; when it is sufficiently close to the global optimal value $\hat{f}$, terminate.
A moment's thought should convince you that this criterion is not practical: there is no way to know the value of $\hat{f}$ until the problem has already been solved.
(Upper and lower bounds are possible to calculate, though, as with $\gamma_2$.)  A more subtle version is to terminate when the Beckmann function stops decreasing, or (in a more common form) to terminate the algorithm when the link or path flows stabilize from one iteration to the next.
The trouble with these convergence criteria is that they cannot distinguish between a situation when the flows stabilize because they are close to the equilibrium solution, and when they stabilize because the algorithm ``gets stuck'' and cannot improve further due to a flaw in its design or a bug in the programming.
For this reason, it is always preferable to base the termination criteria on the equilibrium principle itself.

\section{Link-Based Algorithms}

\index{static traffic assignment!algorithms!link-based|(}
Link-based algorithms for traffic assignment are the simplest to understand and implement, and require the least amount of computer memory.
Given a current set of link flows $\mb{x}$, a link-based algorithm attempts to move closer to equilibrium by performing two steps.
First, a target point\index{all-or-nothing assignment} $\mb{x^*}$\label{not:xstar} is identified; moving in the direction of $\mb{x^*}$ from $\mb{x}$ should lead towards an equilibrium solution.
Then, a step of size $\lambda \in [0, 1]$ is taken in this direction, updating the link flows to $\lambda \mb{x^*} + (1 - \lambda) \mb{x}$.
This is a specific way of implementing the first two steps of the framework specified in the previous section: $\lambda$ represents the fraction of flow which is shifted from the paths at the current solution $\mb{x}$ to the paths at the target solution $\mb{x^*}$.
If $\lambda = 1$, all of the flow has shifted to the target solution; if $\lambda = 0$, no flow has shifted at all.
Since the set of feasible link flow solutions $X$ is convex, as long as $\mb{x}$ and $\mb{x^*}$ are feasible, we can be assured that the convex combination $\lambda \mb{x^*} + (1 - \lambda) \mb{x}$ is feasible as well.

Link-based algorithms differ in two primary ways: first, how the target $\mb{x^*}$ is chosen; and second, how the step size $\lambda$ is chosen.
If $\lambda$ is too small, convergence to equilibrium will be very slow, but if $\lambda$ is too large, the solution may never converge at all --- the example from Figure~\ref{fig:twolink} in the previous section is an example of what can happen if $\lambda = 1$ for all iterations.
This section presents three link-based algorithms, in increasing order of sophistication (but in increasing order of convergence speed.)  The first is the method of successive averages, which is perhaps the simplest equilibrium algorithm.
The second is the Frank-Wolfe algorithm, which can be thought of as a version of the method of successive averages with a more intelligent choice of step size $\lambda$.
(Frank-Wolfe was the most common used in practice for several decades.)  The third is the conjugate Frank-Wolfe algorithm, which can be thought of as a version of Frank-Wolfe with a more intelligent choice of target $\mb{x^*}$.

All of these algorithms use the following framework; the only difference is how $\mb{x^*}$ and $\mb{\lambda}$ are calculated.

\begin{enumerate}
\item Generate an initial solution $\mb{x} \in X$.
\item Generate a target solution $\mb{x^*} \in X$.
\item Update the current solution: $\mb{x} \leftarrow \lambda \mb{x^*} + (1 - \lambda) \mb{x}$ for some $\lambda \in [0,1]$.
\item Calculate the new link travel times.
\item If the convergence criterion is satisfied, stop; otherwise return to step 2.
\end{enumerate}

\subsection{Method of successive averages}
\label{sec:convexmsa}

\index{static traffic assignment!algorithms!method of successive averages|(}
Although the method of successive averages is not competitive with other equilibrium solution algorithms, its simplicity and clarity in applying the three-step iterative process make it an ideal starting place.
To specify the method of successive averages, we need to specify how the target solution $\mb{x^*}$ is chosen, and how the step size $\lambda$ is chosen.

The target $\mb{x^*}$ is an all-or-nothing assignment.\index{all-or-nothing assignment}
That is, \emph{assuming that the current travel times are fixed}, identify the shortest path between each origin and destination, and load all of the demand for that OD pair onto that path.
Thus, $\mb{x^*}$ is the state which would occur if literally every driver was to switch paths onto what is currently the shortest path.
Of course, if we were to switch everybody onto these paths, which would occur if we choose $\lambda = 1$, those paths would almost certainly not be ``shortest'' anymore.
But $\mb{x^*}$ can still be thought of as a target, or a direction in which travelers would feel pressure to move.

So, what should $\lambda$ be?
As discussed above, there are problems if you shift too few travelers, and potentially even bigger problems if you shift too many.
The method of successive averages adopts a reasonable middle ground: initially, we shift a lot of travelers, but as the algorithm progresses, we shift fewer and fewer until we settle down on the average.
The hope is that this avoids both the problems of shifting too few (at first, we're taking big steps, so hopefully we get somewhere close to equilibrium quickly) and of shifting too many (eventually, we'll only be moving small amounts of flow so there is no worry of infinite cycling).

Specifically, on the $i$-th iteration, the method of successive averages uses $\lambda = 1/(i + 1)$.
So, the first time through, half of the travelers are shifted to the current shortest paths.
The second time through, a third of the people shift to the current shortest paths (and two thirds stay on their current path).
On the third iteration, a fourth of the people shift to new paths, and so on.
(The method of successive averages can also be applied with different step size rules; see Exercise~\ref{ex:msastepsize}.)

At this point, it's worth using the Beckmann formulation\index{static traffic assignment!formulation!as convex optimization} to show that the choice of an all-or-nothing assignment for $\mb{x^*}$ has mathematical justification, in addition to the intuitive interpretation of shifting towards shortest path.
Let $\mb{x}$ be the current set of link flows, $\mb{x^*}$ an all-or-nothing assignment.
As a result of a shift of size $\lambda$, the Beckmann function will change as well, and we want to show that it's possible to choose $\lambda$ in some way to guarantee that it will decrease.
That is, we want to show that we can reduce the Beckmann function (and thus move closer to the equilibrium solution) by taking a (correctly-sized) step in the direction $\mathbf{x^* - x}$.
Define $f(\mathbf{x}(\lambda)) = f((1 - \lambda) \mathbf{x} + \lambda \mathbf{x^*})$ to be the Beckmann function after taking a step of size $\lambda$.
Using the multivariate chain rule, the derivative of $f(\lambda)$ is\footnote{In this and all similar equations, we are \emph{evaluating} the function $t_{ij}$ at the specific value $(1 - \lambda) x_{ij} + \lambda x^*_{ij}$, and then \emph{multiplying} the result by $(x^*_{ij} - x_{ij})$.
Unfortunately, ``parentheses'' have multiple uses in mathematics.
In places where this is most likely to be confusing, we have used a dot to emphasize multiplication.
}
\[
\frac{df}{d\lambda} = \sum_{(i,j) \in A} \frac{\partial f}{\partial x_{ij}} \frac{dx_{ij}}{d\lambda} = \sum_{(i,j) \in A} t_{ij} ((1 - \lambda) x_{ij} + \lambda x^*_{ij}) \cdot (x^*_{ij} - x_{ij})
\,.\]
Evaluating this derivative at $\lambda = 0$ gives
\[
\frac{df}{d\lambda}(0) = \sum_{(i,j) \in A} t_{ij} (x_{ij}) \cdot (x^*_{ij} - x_{ij})
\,.\]
Now, $\mathbf{x^*}$ was specifically chosen to put all vehicles on the shortest paths at travel times $\mathbf{t(x)}$, and so $\sum_{(i,j) \in A} x^*_{ij} t_{ij}(x_{ij}) \leq \sum_{(i,j) \in A} x_{ij} t_{ij}(x_{ij})$, and therefore $\frac{df}{d\lambda}(0) \leq 0$.
Furthermore, if we are not at the equilibrium solution already, $\sum_{(i,j) \in A} t_{ij}(x_{ij}) \cdot x^*_{ij} $ is strictly less than $\sum_{(i,j) \in A} t_{ij}(x_{ij}) \cdot  x_{ij}$.
This implies $\frac{df}{d\lambda}(0) < 0$ or, equivalently, we can decrease the Beckmann function if we take a small enough step in the direction $\mathbf{x^* - x}$, by shifting people from longer paths onto shorter ones.

Two examples of the method of successive averages are shown below.
A proof of convergence is sketched in Exercise~\ref{ex:msaproof}.

\paragraph{Small network example}

Here we solve the small example of Figure~\ref{fig:twolink} by the method of successive averages, using the relative gap $\gamma_1$ to measure how close we are to equilibrium.

\begin{description}
\item[Initialization.] Find the shortest paths: with no travelers on the network, the top link has a travel time of 10, and the bottom link has a travel time of 20.
Therefore the top link is the shortest path, so $\mathbf{x^*} = \vect{50 & 0}$.
We take this to be the initial solution $\mathbf{x} \leftarrow \mathbf{x^*} = \vect{50 & 0}$.
Recalculating the travel times, we have $t_1 = 10 + x_1 = 60$ and $t_2 = 20 + x_2 = 20$ (or, in vector form, $\mathbf{t} = \vect{60 & 20}$).
\item[Iteration 1.] With the new travel times, the shortest path is now the bottom link, so $\kappa = 20$ and the relative gap is 
\[
\gamma_1 = \frac{\mathbf{t} \cdot \mathbf{x}}{\bm{\kappa} \cdot \mb{d}} - 1 = \frac{50 \times 60 + 0 \times 20}{20 \times 50} - 1 = 2
\,.\]
This is far too big, so we continue with the second iteration.
If everyone were to take the new shortest path, the flows would be $\mathbf{x^*} = \vect{0 & 50}$.
Because this is the first iteration, we shift 1/2 of the 
travelers onto this path, so $\mathbf{x} \leftarrow (1/2) \mathbf{x^*} + (1/2) \mathbf{x} = \vect{0 & 25} + \vect{25 & 0} = \vect{25 & 25}$.
The new travel times are thus $\mathbf{t} = \vect{35 & 45}$.
\item[Iteration 2.] With the new travel times, the shortest path is now the top link, so $\kappa = 35$ and the relative gap is 
\[
\gamma_1 = \frac{\mathbf{t} \cdot \mathbf{x}}{\bm{\kappa} \cdot \mb{d}} - 1 = \frac{25 \times 35 + 25 \times 45}{35 \times 50} - 1 = 0.143
\,.\]
If everyone were to take the new shortest path, the flows would be $\mathbf{x^*} = \vect{50 & 0}$.
Because this is the second iteration, we shift 1/3 of the 
travelers onto this path, so $\mathbf{x} \leftarrow (1/3) \mathbf{x^*} + (2/3) \mathbf{x} = \vect{50/3 & 0} + \vect{50/3 & 50/3} = \vect{100/3 & 50/3}$.
The new travel times are thus $\mathbf{t} = \vect{43.33 & 36.67}$.
\item[Iteration 3.]  With the new travel times, the shortest path is now the bottom link, so $\kappa = 36.67$ and the relative gap is $\gamma_1 = 0.121$.
A bit better, but still too big, so we carry on.
Here $\mathbf{x^*} = \vect{0 & 50}$, $\mathbf{x} \leftarrow (1/4) \mathbf{x^*} + (3/4) \mathbf{x} = \vect{0 & 50/4} + \vect{25 & 50/4} = \vect{25 & 25}$.
The new travel times are $\mathbf{t} = \vect{35 & 45}$.
Note that we have returned to the same solution found in Iteration 1.
\emph{Don't despair; this just means the last shift was too big.
Next time we'll shift fewer vehicles (because $\lambda$ is smaller with each iteration).}
\item[Iteration 4.]  With the new travel times, the shortest path is now the top link, so $\kappa = 35$ and the relative gap is $\gamma_1 = 0.143$.
The new target is $\mathbf{x^*} = \vect{50 & 0}$, $\mathbf{x} \leftarrow (1/5) \mathbf{x^*} + (4/5) \mathbf{x} = \vect{30 & 20}$.
The new travel times are $\mathbf{t} = \vect{40 & 40}$.
With the new travel times, the shortest path is the top link, so $\kappa = 40$ and the relative gap is $\gamma_1 = 0$, so we stop.
\emph{In fact, either path could have been chosen for the shortest path.
Whenever there is a tie between shortest paths, you are free to choose among them.}
\end{description}

\paragraph{Larger network example}

\begin{figure} 
\hfill
\includegraphics[width=0.6\textwidth]{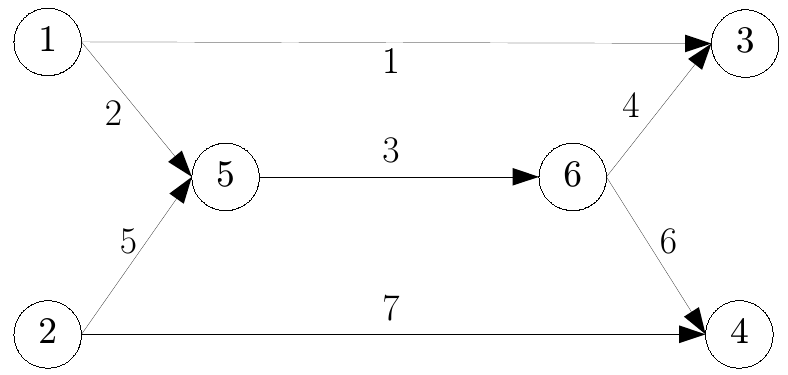}
\hfill
\raisebox{1in}{
	\begin{tabular}{r|cc}
			&	3	&	4 \\
		\hline
		1	&	5,000 & 0 \\
		2	&	0		& 10,000 \\
	\end{tabular}
}
\hfill
\caption{Larger example with two OD pairs.  (Link numbers shown.)  \label{fig:hw1net}}
\end{figure}

Here we apply the method of successive averages to a slightly larger network with two OD pairs, shown in Figure~\ref{fig:hw1net}, where each link has the link performance function $t(x) = 10 + x/100$.

There are four paths in this network; for OD pair (1,3) these are denoted $[1,3]$ and $[1,5,6,3]$ according to their link numbers, and for OD pair (2,4) these are $[2,5,6,4]$ and $[2,4]$.
In this example, we'll calculate the average excess cost, rather than the relative gap.

\begin{description}
\item[Initialization.] Find the shortest paths: with no travelers on the network, paths $[1,3]$, $[1,5,6,3]$, $[2,5,6,4]$, and $[2,4]$ respectively have travel times of 10, 30, 30, and 10.
Therefore $[1,3]$ is shortest for OD pair (1,3), and $[2,4]$ is shortest for OD pair (2,4), so $\mathbf{x^*} = \vect{5000 & 0 & 0 & 0 & 0 & 0 & 10000}$.\footnote{For each OD pair, we add the total demand from the OD matrix onto each link in the shortest path.}  Since is the first iteration, we simply set
\[
\mathbf{x} \leftarrow \mathbf{x^*} = \vect{5000 & 0 & 0 & 0 & 0 & 0 & 10000}
\,.\]
Recalculating the travel times, we have
\[
\mathbf{t} = \vect{60 & 10 & 10 & 10 & 10 & 10 & 110}
\,.\]
\item[Iteration 1.] With the new travel times, the shortest path for (1,3) is now $[1,5,6,3]$, with a travel time of 30, so $\kappa^{13} = 30$.
Likewise, the new shortest path for (2,4) is $[2,5,6,4]$, so $\kappa^{24} = 30$ and the average excess cost is 
\begin{multline*}
AEC = \frac{\mathbf{t} \cdot \mathbf{x} - \bm{\kappa} \cdot \mb{d}}{\mathbf{d} \cdot \mathbf{1}} \\ = \frac{5000 \times 60 + 10000 \times 110 - 30 \times 5000 - 30 \times 10000}{5000 + 10000} = 63.33
\,.
\end{multline*}
This is far too big and suggests that the average trip is 63 minutes slower than the shortest paths available!
If everyone were to take the new shortest paths, the flows would be
\[
\mathbf{x^*} = \vect{0 & 5000 & 15000 & 5000 & 10000 & 10000 & 0}
\,.\]
(Be sure you understand how we calculated this.)
Because this is iteration 1, we shift 1/2 of the travelers onto this path, so
\[
\mathbf{x} \leftarrow (1/2) \mathbf{x^*} + (1/2) \mathbf{x} = \vect{2500 & 2500 & 7500 & 2500 & 5000 & 5000 & 5000}
\,.\]
The new travel times are thus
\[
\mathbf{t} = \vect{35 & 35 & 85 & 35 & 60 & 60 & 60}
\,.\]
\item[Iteration 2.] With the new travel times, the shortest path for (1,3) is now $[1,3]$, with $\kappa^{13} = 35$.
The new shortest path for (2,4) is $[2,4]$, so $\kappa^{24} = 60$ and the average excess cost is 
\begin{multline*}
AEC = \\ \frac{\myp{\begin{array}{l} 2500 \times 35 + 2500 \times 35 + 7500 \times 85 + 2500 \times 35 + 5000 \times 60 \\ \mbox{\hspace{0.5in}} + 5000 \times 60 + 5000 \times 60 - 35 \times 5000 - 60 \times 10000  \end{array}}}{15000} \\ = 68.33
\,.
\end{multline*}
This is still big (and in fact worse), but we persistently continue with the second iteration.
If everyone were to take the new shortest paths, the flows would be
\[
\mathbf{x^*} = \vect{5000 & 0 & 0 & 0 & 0 & 0 & 10000}
\,,\]
so
\[
\mathbf{x} \leftarrow (1/3) \mathbf{x^*} + (2/3) \mathbf{x} = \vect{3333 & 1667 & 5000 & 1667 & 3333 & 3333 & 6667}
\,.\]
The new travel times are 
\[
\mathbf{t} = \vect{43.3 & 26.7 & 60 & 26.7 & 43.3 & 43.3 & 76.7}
\,.\]
\item[Iteration 3.]  With the new travel times, the shortest path for (1,3) is still $[1,3]$, with $\kappa^{13} = 43.3$, and the shortest path for (2,4) is still $[2,4]$ with $\kappa^{24} = 76.7$ and the average excess cost is $AEC = 23.6$.
Continuing the fourth iteration, as before 
\[\mathbf{x^*} = \vect{5000 & 0 & 0 & 0 & 0 & 0 & 10000}
\,,\]
so
\[
\mathbf{x} \leftarrow (1/4) \mathbf{x^*} + (3/4) \mathbf{x} = \vect{3750 & 1250 & 3750 & 1250 & 2500 & 2500 & 7500}
\,.\]
The new travel times are
\[
\mathbf{t} = \vect{47.5 & 22.5 & 47.5 & 22.5 & 35 & 35 & 85}
\,.\]
\item[Iteration 4.]  With the new travel times, the shortest path for (1,3) is still $[1,3]$, with $\kappa^{13} = 47.5$, and the shortest path for (2,4) is still $[2,4]$ with $\kappa^{24} = 85$ and the average excess cost is $AEC = 9.42$.
Continuing the fifth iteration, as before 
\[
\mathbf{x^*} = \vect{5000 & 0 & 0 & 0 & 0 & 0 & 10000}
\,,\]
so
\[
\mathbf{x} \leftarrow (1/5) \mathbf{x^*} + (4/5) \mathbf{x} = \vect{4000 & 1000 & 3000 & 1000 & 2000 & 2000 & 8000}
\,.\]
The new travel times are
\[
\mathbf{t} = \vect{50 & 20 & 40 & 20 & 30 & 30 & 90}
\,.\]
\item[Iteration 5.]  With the new travel times, the shortest path for (1,3) is still $[1,3]$, with $\kappa^{13} = 50$, and the shortest path for (2,4) is still $[2,4]$ with $\kappa^{24} = 90$ and the average excess cost is $AEC = 3.07$.
\emph{Note that the shortest paths have stayed the same over the last three iterations.
This means that we really could have shifted more flow than we actually did.
The Frank-Wolfe algorithm, described in the next section, fixes this problem.} 
We have
\[
\mathbf{x^*} = \vect{5000 & 0 & 0 & 0 & 0 & 0 & 10000}
\,,\]
so
\[
\mathbf{x} \leftarrow (1/6) \mathbf{x^*} + (5/6) \mathbf{x} = \vect{4167 & 833 & 2500 & 833 & 1667 & 1667 & 8333}
\,.\]
The new travel times are
\[
\mathbf{t} = \vect{51.7 & 18.3 & 35 & 18.3 & 26.7 & 26.7 & 93.3}
\,.\]
\item[Iteration 6.]  With the new travel times, the shortest path for (1,3) is still $[1,3]$, with $\kappa^{13} = 51.7$, but the shortest path for (2,4) is now $[2,5,6,4]$ with $\kappa^{24} = 88.3$.
The average excess cost is $AEC = 3.80$.
\emph{Note that the OD pairs are no longer behaving ``symmetrically,'' the shortest path for (1,3) stayed the same, but the shortest path for (2,4) has changed.} 
We have
\[
\mathbf{x^*} = \vect{5000 & 0 & 10000 & 0 & 10000 & 10000 & 0}
\,.\]
so
\[
\mathbf{x} \leftarrow (1/7) \mathbf{x^*} + (6/7) \mathbf{x} = \vect{4286 & 714 & 3571 & 714 & 2857 & 2857 & 7142}
\,.\]
The new travel times are
\[
\mathbf{t} = \vect{52.9 & 17.1 & 45.7 & 17.1 & 38.6 & 38.6 & 81.4}
\,.\]
\end{description}

This process continues over and over until the average excess cost is sufficiently small.
Even with such a small network, the method of successive averages requires a very long time to converge.
An average excess cost of 1 is obtained after eleven iterations, 0.1 after sixty-three iterations, 0.01 after three hundred thirty-two, and the rate of convergence only slows down from there.
\index{static traffic assignment!algorithms!method of successive averages|)}

\subsection{Frank-Wolfe}
\label{sec:frankwolfe}

\index{static traffic assignment!algorithms!Frank-Wolfe|(}
One of the biggest drawbacks with the method of successive averages is that it has a fixed step size.
Iteration $i$ moves exactly $1/(i+1)$ of the travelers onto the new shortest paths, no matter how close or far away we are from the equilibrium.
Essentially, the method of successive averages decides its course of action before it even gets started, then sticks stubbornly to the plan of moving $1/(i+1)$ travelers each iteration.
The Frank-Wolfe algorithm fixes this problem by using an \emph{adaptive} step size.
At each iteration, Frank-Wolfe calculates exactly the right amount of flow to shift to get as close to equilibrium as possible.

We might try to do this by picking $\lambda$ to minimize the relative gap or average excess cost, but this turns out to be harder to compute.
Instead, we pick $\lambda$ to solve a ``restricted'' variational inequality where the feasible set is the line segment connecting $\mb{x}$ and $\mb{x^*}$.
It turns out that this is the same as choosing $\lambda$ to minimize the Beckmann function~\eqn{tapoptstart} along this line segment.
Both approaches for deriving the step size $\lambda$ are discussed below.

Define $X'$\label{not:Xprime}\label{not:xprime} to be the link flows lying on the line segment between $\mb{x}$ and $\mb{x^*}$.
That is, $X' = \{ \mb{x'} : \mb{x} = \lambda \mb{x^*} + (1 - \lambda) \mb{x} \mbox{ for some } \lambda \in [0,1] \}$.
The restricted variational inequality\index{variational inequality!restricted} is: find $\mb{\hat{x}'} \in X'$\label{not:xprimehat} such that $\mathbf{t}(\mathbf{\hat{x}'}) \cdot (\mathbf{\hat{x}'} - \mathbf{x'}) \leq 0$ for all $\mathbf{x'} \in X'$.

This variational inequality is simple enough to be solved as a single equation.
The set $X$ has two endpoints ($\mb{x}$ and $\mb{x^*}$, corresponding to $\lambda = 0$ and $\lambda = 1$, respectively).
For now, assume that the solution $\mb{\hat{x}'}$ to the variational inequality is not at one of these endpoints.\footnote{Exercise~\ref{ex:cornersolution} asks you to show that the solution methods provided below will still give the right answer even in these cases.}  In this case, the force vector $-\mb{t(\hat{x}')}$ is perpendicular to the direction $\mb{x^* - x}$.
(Figure~\ref{fig:fwvifig}), so $\mb{-t(\hat{x}') \cdot (x^* - x)} = 0$.
Writing this equation out in terms of individual components, we need to solve
\labeleqn{fwvizero}{\sum_{ij} t_{ij}(\hat{x}'_{ij}) \left( x^*_{ij} - x_{ij} \right) = 0}
or equivalently
\labeleqn{fwvibalance}{\sum_{ij} t_{ij}(\hat{x}'_{ij}) x^*_{ij} = \sum_{ij} t_{ij}(\hat{x}'_{ij}) x_{ij}}

The same equation can be derived based on the Beckmann function.\index{Beckmann function|(}
Recall the discussion above, where we wrote the function $f(\mb{x}(\lambda)) = f((1 - \lambda) \mathbf{x} + \lambda \mathbf{x^*})$ to be the value of the Beckmann function after taking a step of size $\lambda$, and furthermore found the derivative of $f(\mb{x}(\lambda))$ to be 
\labeleqn{zetader}{
\frac{df}{d\lambda} = \sum_{(i,j) \in A} t_{ij} ((1 - \lambda) x_{ij} + \lambda x^*_{ij}) (x^*_{ij} - x_{ij})
\,.}
It is not difficult to show that $f(\mb{x}(\lambda))$ is a convex function of $\lambda$, so we can find its minimum by setting the derivative equal to zero, which occurs if the condition
\labeleqn{balancing}{
\sum_{(i,j) \in A} x^*_{ij} t_{ij} ((1 - \lambda) x_{ij} + \lambda x^*_{ij})  = \sum_{(i,j) \in A} x_{ij} t_{ij} ((1 - \lambda) x_{ij} + \lambda x^*_{ij})
\,.}
is satisfied, which is the same as~\eqn{fwvibalance}.\footnote{We again clarify that in equation~\eqn{zetader} we are \emph{evaluating} the function $t_{ij}$ at the value $(1 - \lambda) x_{ij} + \lambda x^*_{ij}$ and then \emph{multiplying} the result by $(x^*_{ij} - x_{ij}$; and that in equation~\eqn{balancing} we are multiplying $x^*_{ij}$ (or $x_{ij})$ by the value of the function $t_{ij}$ evaluated at $(1 - \lambda) x_{ij} + \lambda x^*_{ij}$.}

Study these equations carefully: the coefficients $x_{ij}$ and $x^*_{ij}$ are constants and do not change with $\lambda$; the only part of this condition which is affected by $\lambda$ are the travel times.
You can interpret this equation as trying to find a balance between $\mathbf{x}$ and $\mathbf{x^*}$ in the following sense: different values of $\lambda$ correspond to shifting a different number of travelers from their current paths to shortest paths, which will result in different travel times on all the links.
You want to pick $\lambda$ so that, \emph{after} you make the switch, both the old paths $\mathbf{x}$ and the old shortest paths $\mathbf{x^*}$ are equally attractive in terms of their travel times.\index{Beckmann function|)}

\stevefig{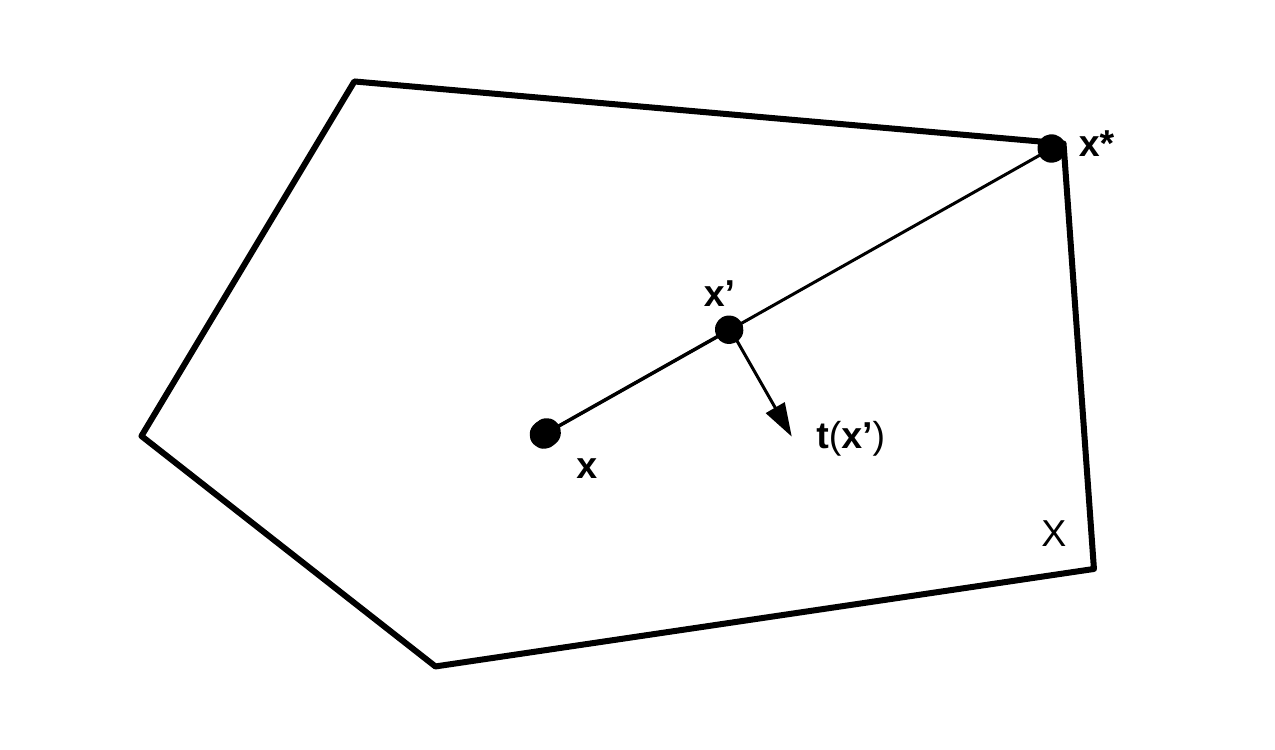}{A solution to the restricted variational inequality in the Frank-Wolfe method.}{0.8\textwidth}

It is convenient to write equation~\eqn{fwvizero} or~\eqn{balancing} as a function in terms of $\lambda$ as 
\labeleqn{balancingfinal}{
\zeta'(\lambda) = \sum_{ij} t_{ij}(\lambda x^*_{ij} + (1 - \lambda) x_{ij} ) \left( x^*_{ij} - x_{ij} \right) = 0
\,,}
which we need to solve for $\lambda \in [0, 1]$ and $\zeta(\lambda)$\label{not:zeta} is used as a shorthand for $f(\mb{x}(\lambda))$.
Since the link performance functions are typically nonlinear, we cannot expect to be able to solve this equation analytically to get an explicit formula for $\lambda$.
General techniques such as Newton's method or an equation solver can be used; but it's not too difficult to use an enlightened trial-and-error method such as a binary search or bisection, as was described in Section~\ref{sec:bisection}.
The examples in this section will use the bisection method.

To summarize: in Frank-Wolfe, $\mb{x^*}$ is an all-or-nothing assignment\index{all-or-nothing assignment} (just as with the method of successive averages).
The difference is that $\lambda$ is chosen to solve~\eqn{balancingfinal}.
So, there is a little bit more work at each iteration (we have to solve an equation for $\lambda$ instead of using a pre-computed formula as in the method of successive averages), but the reward is much faster convergence to the equilibrium solution.
Two examples of Frank-Wolfe now follow; you are asked to provide a proof of correctness in Exercise~\ref{ex:fwdescent} and~\ref{ex:fwproof}.

\paragraph{Small network example}

Here we solve the small example of Figure~\ref{fig:twolink} by Frank-Wolfe.
Some steps are similar to the method of successive averages, and therefore omitted.
Here, when we do the bisection method, we do five interval reductions, so we are within $1/2^5 = 1/32$ of the correct $\lambda^*$ value.
When solving by computer, you would usually perform more steps than this, because the bisection calculations are very fast.

\begin{description}
\item[Initialization.]  As before, we load everybody on the initial shortest path, so $\mathbf{x} = \mathbf{x^*} = \vect{50 & 0}$ and $\mathbf{t} = \vect{60 & 20}$.
\item[Iteration 1.]  As before, the relative gap is $\gamma_1= 2$.
With the new shortest paths, $\mathbf{x^*} = \vect{0 & 50}$.
Begin the bisection method.
\begin{description}
    \item[Bisection Iteration 1.] Initially $\lambda^* \in [0, 1]$.
Calculate $\zeta'(1/2) = (0 - 50) \times (10 + 25) + (50 - 0) \times (20 + 25) = 500 > 0$ so we discard the upper half.
    \item[Bisection Iteration 2.] Now we know $\lambda^* \in [0, 1/2]$.
Calculate \[\zeta'(1/4) = (0 - 50) \times (10 + 37.5) + (50 - 0) \times (20 + 12.5) = -750 < 0\] so we discard the lower half.
    \item[Bisection Iteration 3.] Now we know $\lambda^* \in [1/4, 1/2]$.
Calculate \[\zeta'(3/8) = (0 - 50) \times (10 + 31.25) + (50 - 0) \times (20 + 18.75) = -125 < 0\] so we discard the lower half.
    \item[Bisection Iteration 4.] Now we know $\lambda^* \in [3/8, 1/2]$.
Calculate \[\zeta'(7/16) = (0 - 50) \times (10 + 28.125) + (50 - 0) \times (20 + 21.875) = 187.5\] so we discard the upper half.
    \item[Bisection Iteration 5.] Now we know $\lambda^* \in [3/8, 7/16]$.
Calculate \[\zeta'(13/32) = (0 - 50) \times (10 + 29.6875) + (50 - 0) \times (20 + 20.3125) = 31.25\] so we discard the upper half.
\end{description}
From here we take the midpoint of the last interval $[3/8,13/32]$ to estimate $\lambda^* \approx 25/64 = 0.390625$, so $\mathbf{x} = (25/64) \mathbf{x^*} + (39/64) \mathbf{x} = \vect{30.47 & 19.53}$ and $\mathbf{t} = \vect{40.47 & 39.53}$.
\item[Iteration 2.]  The relative gap is calculated as $\gamma_1 = 0.014$.
(This is an order of magnitude smaller than the relative gap the method of successive averages found by this point.)
The shortest paths are still $\mathbf{x^*} = \vect{0 & 50}$, and we begin bisection.
\begin{description}
	\item[Bisection Iteration 1.] Initially $\lambda^* \in [0, 1]$.
Calculate \[\zeta'(1/2) = 900 > 0\] so we discard the upper half.
	\item[Bisection Iteration 2.] Now we know $\lambda^* \in [0, 1/2]$.
Calculate \[\zeta'(1/4) = 435 > 0\] so we discard the upper half.
	\item[Bisection Iteration 3.] Now we know $\lambda^* \in [0, 1/4]$.
Calculate \[\zeta'(1/8) = 203 > 0\] so we discard the upper half.
	\item[Bisection Iteration 4.] Now we know $\lambda^* \in [0, 1/8]$.
Calculate \[\zeta'(1/16) = 87 > 0\] so we discard the upper half.
	\item[Bisection Iteration 5.] Now we know $\lambda^* \in [0, 1/16]$.
Calculate \[\zeta'(1/32) = 29 > 0\] so we discard the upper half.
\end{description}
The midpoint of the final interval is $\lambda^* \approx 1/64$, so $\mathbf{x} = (1/64) \mathbf{x^*} + (63/64) \mathbf{x} = \vect{29.99 & 20.01}$ and $\mathbf{t} = \vect{39.99 & 40.01}$.
\item[Iteration 3.]  The relative gap is now $\gamma_1 = 0.00014$, so we quit and claim we have found flows that are ``good enough'' (the difference in travel times between the routes is less than a second).
\end{description}

Alternately, using calculus, we could have identified $\lambda^*$ during the second iteration as \emph{exactly} 0.40, which would have found the exact equilibrium after only one step.

\paragraph{Large network example}

Here we apply Frank-Wolfe to the network shown in Figure~\ref{fig:hw1net}, using the same notation as in the method of successive averages example.

\begin{description}
\item[Initialization.] Path $[1,3]$ is shortest for OD pair (1,3), and path $[2,4]$ is shortest for OD pair (2,4), so
\[
\mathbf{x^*} = \vect{5000 & 0 & 0 & 0 & 0 & 0 & 10000}
\]
and
\[
\mathbf{x} = \mathbf{x^*} = \vect{5000 & 0 & 0 & 0 & 0 & 0 & 10000}
\,.\]
Recalculating the travel times, we have
\[
\mathbf{t} = \vect{60 & 10 & 10 & 10 & 10 & 10 & 110}
\,.\]
\item[Iteration 1.] With the new travel times, the shortest path for (1,3) is now $[1,5,6,3]$, and the new shortest path for (2,4) is $[2,5,6,4]$, so $AEC = 63.33$
If everyone were to take the new shortest paths, the flows would be
\[
\mathbf{x^*} = \vect{0 & 5000 & 15000 & 5000 & 10000 & 10000 & 0}
\,.\]
Begin the bisection method to find the right combination of $\mathbf{x^*}$ and $\mathbf{x}$.
\begin{description}
	\item[Bisection Iteration 1.] Initially $\lambda^* \in [0, 1]$.
Calculate \begin{multline*} \zeta'(1/2) = (0 - 5000) \times (10 + 2500/100) + \ldots \\ + (0 - 10000) \times (10 + 5000/100) = 2050000 > 0\end{multline*} so we discard the upper half.
	\item[Bisection Iteration 2.] Now we know $\lambda^* \in [0, 1/2]$.
Calculate \begin{multline*}\zeta'(1/4) = (0 - 5000) \times (10 + 3750/100) + \ldots \\ + (0 - 10000) \times (10 + 7500/100) = 550000 > 0\end{multline*} so we discard the upper half.
	\item[Bisection Iteration 3.] Now we know $\lambda^* \in [0, 1/4]$.
Calculate \begin{multline*}\zeta'(1/8) = (0 - 5000) \times (10 + 4375/100) + \ldots \\ + (0 - 10000) \times (10 + 8750/100) = -200000 < 0\end{multline*} so we discard the lower half.
	\item[Bisection Iteration 4.] Now we know $\lambda^* \in [1/8, 1/4]$.
Calculate \begin{multline*}\zeta'(3/16) = (0 - 5000) \times (10 + 4062/100) + \ldots \\ + (0 - 10000) \times (10 + 8125/100) = 175000 > 0\end{multline*} so we discard the upper half.
	\item[Bisection Iteration 5.] Now we know $\lambda^* \in [1/8, 3/16]$.
Calculate \begin{multline*}\zeta'(5/32) = (0 - 5000) \times (10 + 4219/100) + \ldots \\ + (0 - 10000) \times (10 + 8437/100) = -12500 < 0\end{multline*} so we discard the lower half.
\end{description}
The final interval is $\lambda* \in [5/32, 3/16]$, so the estimate is $\lambda^* = 11/64$ and
\[
\mathbf{x} = \frac{11}{64} \mathbf{x^*} + \frac{53}{64} \mathbf{x} = \vect{4141 & 859 & 2578 & 859 & 1719 & 1719 & 8281}
\,.\]
The new travel times are thus
\[
\mathbf{t} = \vect{51.4 & 18.6 & 35.8 & 18.6 & 27.2 & 27.2 & 92.8}
\,.\]
\item[Iteration 2.] With the new travel times, the shortest path for (1,3) is now $[1,3]$, but the shortest path for (2,4) is still $[2,5,6,4]$.
The relative gap is $AEC = 2.67$ (roughly 30 times smaller than the corresponding point in the method of successive averages!)
We have
\[
\mathbf{x^*} = \vect{5000 & 0 & 10000 & 0 & 10000 & 10000 & 0}
\,.\]
We begin the bisection method to find the right combination of $\mathbf{x^*}$ and $\mathbf{x}$.
\begin{description}
	\item[Bisection Iteration 1.] Initially $\lambda^* \in [0, 1]$.
Calculate \[\zeta'(1/2) = 637329 > 0\] so we discard the upper half.
	\item[Bisection Iteration 2.] Now we know $\lambda^* \in [0, 1/2]$.
Calculate \[\zeta'(1/4) = 154266 > 0\] so we discard the upper half.
	\item[Bisection Iteration 3.] Now we know $\lambda^* \in [0, 1/4]$.
Calculate \[\zeta'(1/8) = 36063 > 0\] so we discard the upper half.
	\item[Bisection Iteration 4.] Now we know $\lambda^* \in [0, 1/8]$.
Calculate \[\zeta'(1/16) = 7741 > 0\] so we discard the upper half.
	\item[Bisection Iteration 5.] Now we know $\lambda^* \in [0, 1/16]$.
Calculate \[\zeta'(1/32) = 1302 > 0\] so we discard the upper half.
\end{description}
The final interval is $\lambda* \in [0, 1/32]$, so the estimate is $\lambda^* = 1/64$ and
\[
\mathbf{x} = \frac{1}{64} \mathbf{x^*} + \frac{63}{64} \mathbf{x} = \vect{4154 & 845 & 2694 & 845 & 1848 & 1848 & 8152}
\,.\]
The new travel times are thus
\[
\mathbf{t} = \vect{51.5 & 18.5 & 36.9 & 18.5 & 28.5 & 28.5 & 91.5}
\,.\]
\end{description}

At this point, the average excess cost is around 1.56 min; note that Frank-Wolfe is able to decrease the relative gap much faster than the method of successive averages.
However, we're still quite far from equilibrium if you compute the actual path travel times.
In this case, even though we're allowing the step size to vary for each iteration, we are forcing travelers from all OD pairs to shift in the same proportion.
In reality, OD pairs farther from equilibrium should see bigger flow shifts, and OD pairs closer to equilibrium should see smaller ones.
This can be remedied by more advanced algorithms.
\index{static traffic assignment!algorithms!Frank-Wolfe|)}

\subsection{Conjugate Frank-Wolfe}

\index{static traffic assignment!algorithms!conjugate Frank-Wolfe|(}
If the distinction between Frank-Wolfe and the method of successive averages is that Frank-Wolfe chooses the step size $\lambda$ in a more clever way, the distinction between conjugate Frank-Wolfe and plain Frank-Wolfe is that conjugate Frank-Wolfe chooses the target $\mb{x^*}$ in a more clever way.
To understand why conjugate Frank-Wolfe is more clever, we first need to understand why using the all-or-nothing assignment as target can be problematic.

\index{static traffic assignment!algorithms!Frank-Wolfe|(}
Viewed in terms of the set of feasible link assignments $X$, the all-or-nothing assignments\index{all-or-nothing assignment} correspond to corner points\index{convex set!corner point} of $X$.
That is, Frank-Wolfe must limit itself to the corner points of the feasible region when determining where to move.
In Figure~\ref{fig:fwfail}, Frank-Wolfe is constrained to follow the trajectory shown by the thin lines, and is unable to take a direct step like that indicated by the thick arrow.
While these directions are effective in the early iterations, as the algorithm approaches the equilibrium point its converge slows down dramatically, and ``zigzagging'' behavior is observed.
\index{static traffic assignment!algorithms!Frank-Wolfe|)}

\stevefig{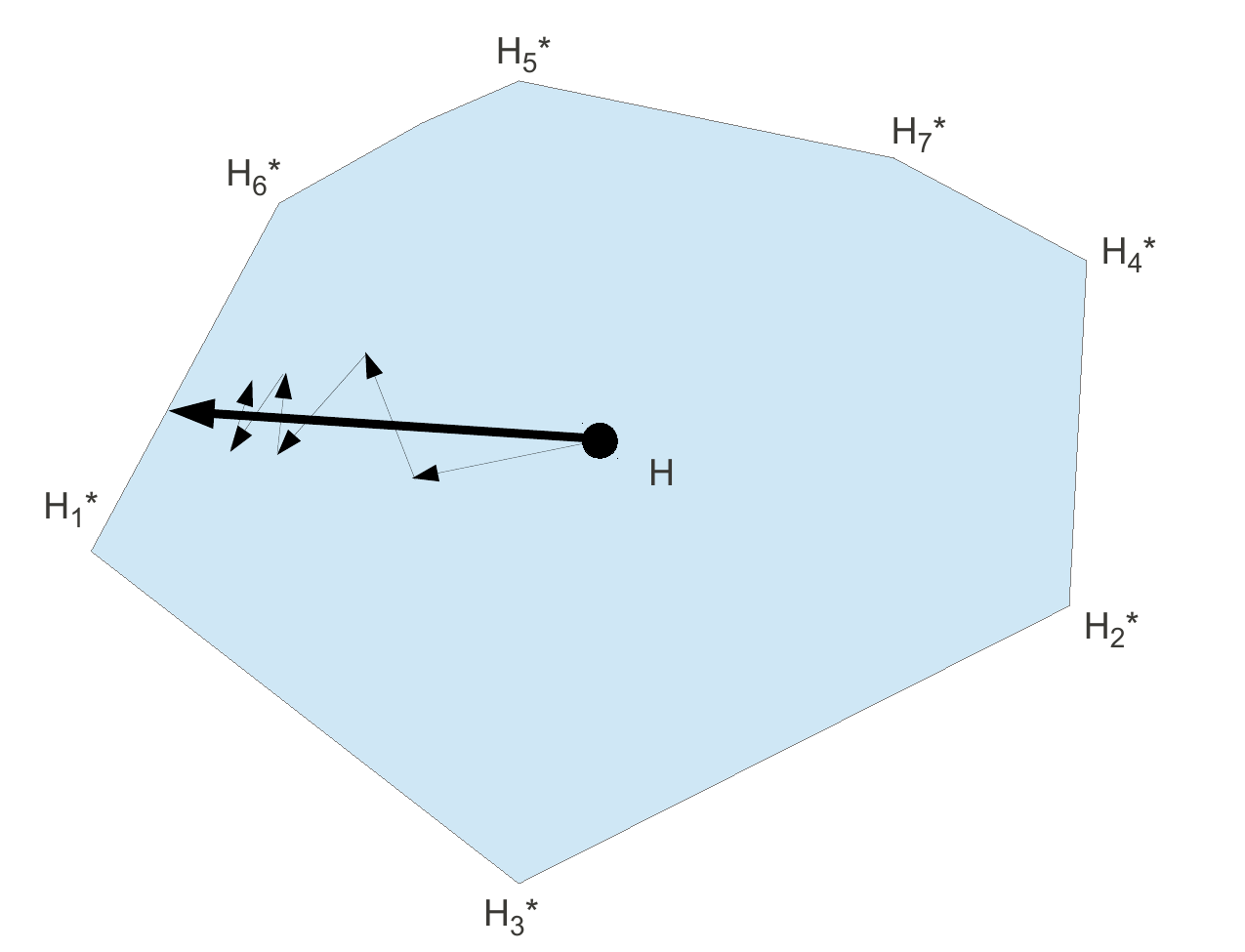}{Frank-Wolfe can only move towards extreme points of the feasible region.}{0.8\textwidth}

So, how can we choose the target solution in a smarter way, so that steps in the direction of the target still move toward equilibrium, while granting more flexibility than the use of an all-or-nothing assignment?
Conjugate Frank-Wolfe provides one approach towards doing so.

\begin{figure}
\begin{center}
$
\begin{array}{c}
\includegraphics[height=0.19\textheight]{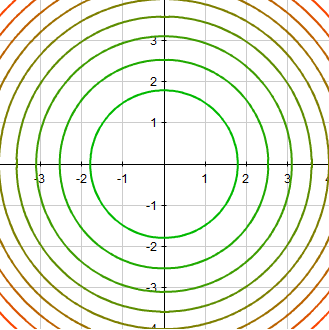} 
\end{array}
$  \hspace{1in} $\mb{Q} = \vect{1 & 0 \\ 0 & 1} \qquad \mb{b} = \vect{0 \\ 0}$ \\
(a) $f(x_1, x_2) = x_1^2 + x_2^2$ \\
\vspace{0.1in}
$
\begin{array}{c}
\includegraphics[height=0.19\textheight]{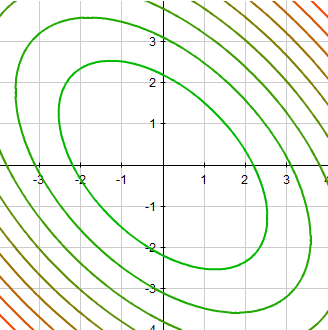} 
\end{array}
$  \hspace{1in} $\mb{Q} = \vect{1 & 1/2 \\ 1/2 & 1} \qquad \mb{b} = \vect{0 \\ 0}$ \\
(b) $f(x_1, x_2) = x_1^2 + x_2^2 + x_1 x_2$ \\
\vspace{0.1in}
$
\begin{array}{c}
\includegraphics[height=0.19\textheight]{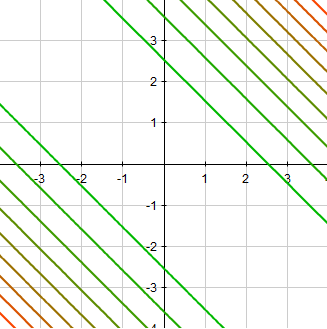} 
\end{array}
$  \hspace{1in} $\mb{Q} = \vect{1 & 1 \\ 1 & 1} \qquad \mb{b} = \vect{0 \\ 0}$ \\
(c) $f(x_1, x_2) = x_1^2 + x_2^2 + 2x_1 x_2$ \\
\vspace{0.1in}
$
\begin{array}{c}
\includegraphics[height=0.19\textheight]{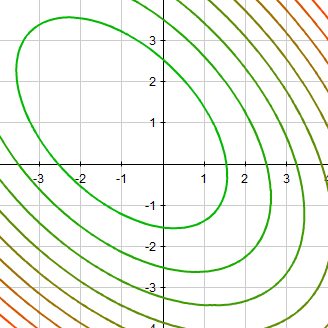} 
\end{array}
$  \hspace{1in} $\mb{Q} = \vect{1 & 1/2 \\ 1/2 & 1} \qquad \mb{b} = \vect{1 \\ -1}$ \\
(d) $f(x_1, x_2) = x_1^2 + x_2^2 + x_1 x_2 + x_1 - x_2$ \\
\caption{Four examples of convex quadratic functions of the form $f(\mb{x}) = \mb{x}^T Q \mb{x} + \mb{b}^T \mb{x}$. \label{fig:quadraticcontours}}
\end{center}
\end{figure}

\index{method of successive averages|see {static traffic assignment, algorithms}}
\index{Frank-Wolfe|see {static traffic assignment, algorithms}}
\index{conjugate Frank-Wolfe|see {static traffic assignment, algorithms}}
\index{conjugacy|(}
Understanding conjugate Frank-Wolfe requires introducing the concept of conjugacy, which is done here.
Temporarily ignoring the context of the traffic assignment problem, assume that we are trying to find the minimum point of a convex quadratic function $f(x_1, x_2)$ of two variables, when there are no constraints.
Figure~\ref{fig:quadraticcontours} shows a few examples of these functions.
Any quadratic function of two variables can be written in the form\label{not:Qbold}\label{not:bbold}
\labeleqn{genquadratic}{f(\mb{x}) = \mb{x}^T \mb{Q} \mb{x} + \mb{b}^T\mb{x}\,,}
where $\mb{x} = \vect{x_1 & x_2}$ and $\mb{b}$ are two-dimensional vectors, and $\mb{Q}$ is a $2 \times 2$ matrix.
Figure~\ref{fig:quadraticcontours} shows the $\mb{Q}$ matrix and $\mb{b}$ vector corresponding to each example.

\begin{figure}
\begin{center}
$
\begin{array}{c}
\includegraphics[height=0.2\textheight]{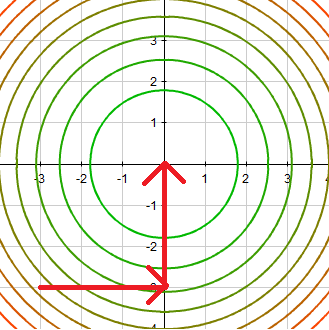} 
\end{array}
$  \hspace{1in} $\mb{Q} = \vect{1 & 0 \\ 0 & 1} \qquad \mb{b} = \vect{0 \\ 0}$ \\
(a) $f(x_1, x_2) = x_1^2 + x_2^2$ \\
\vspace{0.1in}
$
\begin{array}{c}
\includegraphics[height=0.2\textheight]{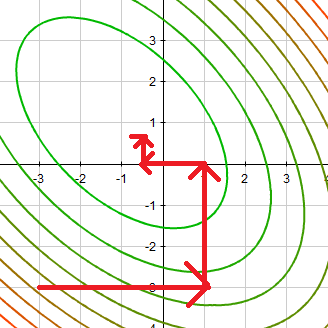} 
\end{array}
$  \hspace{1in} $\mb{Q} = \vect{1 & 1/2 \\ 1/2 & 1} \qquad \mb{b} = \vect{1 \\ -1}$ \\
(b) $f(x_1, x_2) = x_1^2 + x_2^2 + x_1 x_2 + x_1 - x_2$ \\
\caption{Searching in orthogonal directions finds the optimum in two steps only if $Q = I$.
\label{fig:quadraticorthogonal}}
\end{center}
\end{figure}

How would you go about finding the minimum of such a function?
Given some initial solution $(x_1, x_2)$, one idea is to fix $x_1$ as a constant, and find the value of $x_2$ which minimizes $f$.
Then, we can fix $x_2$, and find the value of $x_1$ which minimizes $f$, and so on.
This process will converge to the minimum, as shown in Figure~\ref{fig:quadraticorthogonal}, but in general this convergence is only asymptotic, and the process will never actually reach the minimum.
The exception is when $\mb{Q}$ is the identity matrix, as in Figure~\ref{fig:quadraticorthogonal}(a).
In this case, the exact optimum is reached in only two steps.

In fact, it is possible to reach the exact optimum in only two steps even when $\mb{Q}$ is not the identity matrix, by changing the search directions.
The process described above (alternately fixing $x_1$, and then $x_2$) can be thought of as alternating between searching in the direction $\vect{0 & 1}$, then searching in the direction $\vect{1 & 0}$.
As shown in Figure~\ref{fig:quadraticconjugate}, by making a different choice for the two search directions, the minimum can always be obtained in exactly two steps.

\genfig{quadraticconjugate}{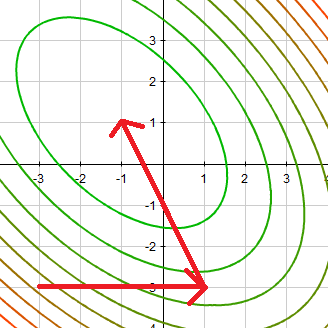}{Searching in conjugate directions always leads to the optimum in two steps.}{height=0.2\textheight}

This happens if the two directions $\mb{d_1}$\label{not:dbold} and $\mb{d_2}$ are \emph{conjugate}, that is, if 
\labeleqn{conjugate}{\mb{d_1}^T \mb{Q} \mb{d_2} = 0}
Conjugacy generalizes the concept of orthogonality (or perpendicularity).
If $\mb{Q}$ is the identity matrix, equation~\eqn{conjugate} reduces to $\mb{d_1} \cdot \mb{d_2} = 0$, the definition of perpendicular vectors.

Now, returning to the traffic assignment problem, we want to use the concept of conjugacy to choose a more intelligent search direction.
In particular, we want the target $\mb{x^*}$ to be chosen so that the search direction $\mb{x^*} - \mb{x}$ is conjugate to the previous search direction.
Before the derivation, there are a few differences between traffic assignment and the unconstrained quadratic program used to introduce conjugacy which should be addressed.
\index{conjugacy|)}

\begin{itemize}
\item The Beckmann function\index{Beckmann function|(} is not a function of two variables.
This is not a huge problem.
Finding the unconstrained minimum of a quadratic function of $n$ variables requires $n$ conjugate steps (so, in the above examples with two variables, two steps sufficed), due to a result known as the expanding subspace theorem.
Of course...
\item The Beckmann function is in general not a quadratic function.
(Can you think of a case when it is?)
Instead of the matrix $\mb{Q}$, we will instead use the Hessian of the Beckmann function $Hf$.\label{not:Hessian}
Therefore, we cannot hope for exact convergence in a finite number of iterations.
However, when the solution gets closer and closer to equilibrium, the Beckmann function can be better and better approximated by a quadratic by taking the first two terms of its Taylor series.
This is good news, because zigzagging in plain Frank-Wolfe becomes worse and worse as the equilibrium solution is approached.\index{Beckmann function|)}
\item The optimization problem~\eqn{tapoptstart}--\eqn{tapoptend} has constraints.
Again, this is not a major problem, since the feasible region is convex.
As long as we ensure that the target point $\mb{x^*}$ is feasible, any choice of $\lambda \in [0,1]$ will retain feasibility.
\end{itemize}

So, how can we make sure that the new target vector $\mb{x^*}$ is chosen so that the search direction is conjugate to the previous direction, and that $\mb{x^*}$ is feasible?
Since $X$ is a convex set, feasibility can be assured by choosing $\mb{x^*}$ to be a convex combination of the old target vector $(\mb{x^*_{old}})$ and the all-or-nothing assignment $\mb{x^{AON}}$:\index{convex set!applications}
\labeleqn{newtarget}{
\mb{x^*} = \alpha \mb{x^*_{old}} + (1 - \alpha) \mb{x^{AON}}
\,,}
for some $\alpha \in [0, 1]$.\label{not:alphacfw}
Choosing $\alpha = 0$ would make the all-or-nothing assignment the target (as in plain Frank-Wolfe), while choosing $\alpha = 1$ would make the target in this iteration the same as in the last.
In fact, $\alpha$ should be chosen so that the new direction is conjugate to the last, that is,
\labeleqn{conjugacy}{
(\mb{x^*_{old}} - \mb{x})^T \mb{H} (\mb{x^*} - \mb{x}) = 0
\,,}
where $\mb{H}$ is the Hessian\index{Hessian matrix!applications} of the Beckmann function evaluated at the current solution $\mb{x}$.
Substituting~\eqn{newtarget} into~\eqn{conjugacy}, we can solve for $\alpha$ as follows:
\begin{align}
(\mb{x^*_{old}} - \mb{x})^T \mb{H} (\mb{x^*} - \mb{x}) &= 0 \\
\iff (\mb{x^*_{old}} - \mb{x})^T \mb{H} (\alpha \mb{x^*_{old}} + (1 - \alpha) \mb{x^{AON}} - \alpha \mb{x} - (1 - \alpha) \mb{x}) &= 0 \\
\iff (\mb{x^*_{old}} - \mb{x})^T \mb{H} (\alpha (\mb{x^*_{old}} - \mb{x}) + (1 - \alpha) (\mb{x^{AON}} - \mb{x})) &= 0 
\end{align}
or, after rearrangement,
\begin{align}
\alpha [(\mb{x^*_{old}} - \mb{x})^T \mb{H} (\mb{x^*_{old}} - \mb{x^{AON}})] &= -(\mb{x^*_{old}} - \mb{x})^T \mb{H} (\mb{x^{AON}} - \mb{x}) \\
\iff \alpha &= \frac{(\mb{x^*_{old}} - \mb{x})^T \mb{H} (\mb{x^{AON}} - \mb{x})}{(\mb{x^*_{old}} - \mb{x})^T \mb{H} (\mb{x^*_{AON}} - \mb{x^{old}})} \label{eqn:alphaone}
\end{align}

Now, for the traffic assignment problem, the Hessian takes a specific form.
Since the Beckmann function is
\labeleqn{beckmannfunction}{
f(\mb{x}) = \sum_{(i,j) \in A} \int_0^{x_{ij}} t_{ij}(x)~dx
\,,}
its gradient\index{Beckmann function!gradient} is simply the vector of travel times at the current flows
\labeleqn{beckmanngradient}{
\nabla f (\mb{x}) = \mr{vect} \{ t_{ij} (x_{ij}) \}
\,,}
and its Hessian\index{Beckmann function!Hessian} is the diagonal matrix of travel time derivatives at the current flows
\labeleqn{beckmannhessian}{
Hf(\mb{x}) = \mr{diag} \{ t'_{ij} (x_{ij}) \}
\,.}
So, the matrix products in equation~\eqn{alphaone} can be written out explicitly, giving
\labeleqn{alphatemp}{
\alpha = \frac{\sum_{(i,j) \in A} ((x^*_{old})_{ij} - x_{ij})(x^{AON}_{ij} - x_{ij}) t'_{ij} }{\sum_{(i,j) \in A} ((x^*_{old})_{ij} - x_{ij})(x^{AON}_{ij} - (x^*_{old})_{ij}) t'_{ij} }
\,,}
where the derivatives $t'_{ij}$ are evaluated at the current link flows $x_{ij}$.

Almost there!
A careful reader may have some doubts about the formula in~\eqn{alphatemp}.
First, it is possible that the denominator can be zero, and division by zero is undefined.
Second, to ensure feasibility of $\mb{x^*}$, we need $\alpha \in [0, 1]$, even though it is not obvious that this formula always lies in this range (and in fact, it need not do so).
Furthermore, $\alpha = 1$ is undesirable, because then the current target point is the same as the target point in the last iteration.
If the previous line search was exact, there will be no further improvement and the algorithm will be stuck in an infinite loop.
Finally, what should you do for the first iteration, when there is no ``old'' target $\mb{x^*}$?

To address the first issue, the easiest approach is to simply set $\alpha = 0$ if the denominator of~\eqn{alphatemp} is zero (i.e., if the formula is undefined, simply take a plain Frank-Wolfe step by using the all-or-nothing solution as the target).
As for the second and third issues, if the denominator is nonzero we can project the right-hand side of~\eqn{alphatemp} onto the interval $[0, 1 - \epsilon]$ where $\epsilon > 0$ is some small tolerance value.
That is, if equation~\eqn{alphatemp} would give a value greater than $1 - \epsilon$, set $\alpha = 1 - \epsilon$; if it would give a negative value, use zero.
Finally, for the first iteration, simply use the all-or-nothing solution as the target: $\mb{x^*} = \mb{x^{AON}}$.

So, to summarize the discussion, choose $\alpha$ in the following way.
If it is the first iteration or the denominator of~\eqn{alphatemp} is zero, set $\alpha = 0$.
Otherwise set
\labeleqn{alphafinal}{\alpha = \mr{proj}_{[0,1 - \epsilon]} \myp{\frac{\sum_{(i,j) \in A} t'_{ij} ((x^*_{old})_{ij} - x_{ij})(x^{AON}_{ij} - x_{ij})}{\sum_{(i,j) \in A} t'_{ij} ((x^*_{old})_{ij} - x_{ij})(x^{AON}_{ij} - (x^*_{old})_{ij})}} }
Then the target solution $\mb{x^*}$ is calculated using~\eqn{newtarget}.
The value of the step size $\lambda$ is chosen in the same way as in Frank-Wolfe, by performing a line search (e.g., using bisection or Newton's method) to solve~\eqn{balancingfinal}.

\paragraph{Large network example}

Here we apply conjugate Frank-Wolfe to the network shown in Figure~\ref{fig:hw1net}, using the same notation as in the Frank-Wolfe and method of successive averages examples.
The tolerance $\epsilon$ is chosen to be a small positive constant, 0.01 in the following example.

\begin{description}
\item[Initialization.] Generate the initial solution by solving an all-or-nothing assignment.
Path $[1,3]$ is shortest for OD pair (1,3), and path $[2,4]$ is shortest for OD pair (2,4), so
\[
\mathbf{x^{AON}} = \vect{5000 & 0 & 0 & 0 & 0 & 0 & 10000}
\]
and
\[
\mathbf{x} = \mathbf{x^{AON}} = \vect{5000 & 0 & 0 & 0 & 0 & 0 & 10000}
\,.\]
Recalculating the travel times, we have
\[
\mathbf{t} = \vect{60 & 10 & 10 & 10 & 10 & 10 & 110}
\,.\]
\item[Iteration 1.] Proceeding in the same way as in the large network example for Frank-Wolfe, the all-or-nothing assignment in this case is
\[
\mathbf{x^{AON}} = \vect{0 & 5000 & 0 & 10000 & 15000 & 5000 & 10000}
\,.\]
which is used as the target $\mb{x^*}$ since this is the first iteration of conjugate Frank-Wolfe.
Repeating the same line search process, the optimal value of $\lambda$ is 19/120, producing the new solution
\[
\mathbf{x} = \vect{4208 & 792 & 8417 & 1583 & 2375 & 792 & 1583}
\,.\]
\item[Iteration 2.] Again as with regular Frank-Wolfe, based on the updated travel times the all-or-nothing assignment is
\[
\mathbf{x^{AON}} = \vect{5000 & 0 & 0 & 10000 & 10000 & 0 & 10000}
\,.\]
From here Frank-Wolfe and conjugate Frank-Wolfe take different paths.
Rather than using $\mb{x^{AON}}$ as the target vector, conjugate Frank-Wolfe generates a conjugate search direction.
First calculate the right-hand side of~\eqn{alphatemp}.
Since for this problem $t_{ij} = 1/100$ for all links (regardless of the flow), the formula is especially easy to compute using $\mathbf{x}$ and $\mathbf{x^*}$ from the previous iteration and $\mathbf{x^{AON}}$ as just now computed.
The denominator of~\eqn{alphatemp} is nonzero, and the formula gives $-2.37$; projecting onto the set $[0, 1 - \epsilon]$ thus gives $\alpha = 0$.
So, calculating the target $\mb{x^*}$ from equation~\eqn{newtarget} with $\alpha = 0$ we have
\[
\mathbf{x^{*}} = \vect{5000 & 0 & 0 & 10000 & 10000 & 0 & 10000}
\]
and, using a line search between $\mb{x}$ and $\mb{x^*}$, find that $\lambda = 0.0321$ is best, resulting in
\[
\mathbf{x} = \vect{4233 & 766 & 8146 & 1853 & 2620 & 766 & 1854}
\,.\]
\item[Iteration 3.] With the new flows $\mathbf{x}$, the travel times are now
\[
\mb{t} = \vect{52.3 & 17.7 & 91.5 & 28.5 & 36.2 & 17.7 & 28.5}
\]
and the all-or-nothing assignment is 
\[
\mathbf{x^{AON}} = \vect{5000 & 0 & 10000 & 0 & 0 & 0 & 0}
\,.\]
Calculating the right-hand side of~\eqn{alphatemp}, we see that the denominator is nonzero, and the formula gives 0.198, which can be used as is since it lies in $[0, 1 - \epsilon]$.
So, calculating the target $\mb{x^*}$ from equation~\eqn{newtarget} with $\alpha = 0.198$ we have
\[
\mathbf{x^{*}} = \vect{5000 & 0 & 8024  & 1976 & 1976 & 0 & 1976}
\,.\]
Note that unlike any of the other link-based algorithms in this section, the target flows are not an all-or-nothing assignment (i.e., not an extreme point of $X$).
Performing a line search between $\mb{x}$ and $\mb{x^*}$, find that $\lambda = 0.652$ is best, resulting in
\[
\mathbf{x} = \vect{4733 & 267 & 8067 & 1933 & 2200 & 267 & 1933}
\,.\]
which solves the equilibrium problem exactly, so we terminate.
\end{description}

In this example, conjugate Frank-Wolfe found the exact equilibrium solution in three iterations.
This type of performance is not typical (even though it is generally faster than regular Frank-Wolfe or the method of successive averages).
In this example, the link performance functions are linear, so the Beckmann function is quadratic.
Iterative line searches with conjugate directions lead to the exact solution of quadratic programs in a finite number of iterations, as suggested by the above discussion.
This performance cannot be assured with other types of link performance functions.

An even faster algorithm known as \emph{biconjugate} Frank-Wolfe chooses its target so that the search direction is conjugate to both of the previous two search directions.\index{static traffic assignment!algorithms!biconjugate Frank-Wolfe}
This method converges faster than conjugate Frank-Wolfe, but is not explained here because the details are a little more complicated even though the idea is the same.
Exercise~\ref{ex:biconjugate} provides formulas for the target solution and asks you to show that they satisfy the necessary conditions.
\index{static traffic assignment!algorithms!conjugate Frank-Wolfe|)}
\index{static traffic assignment!algorithms!link-based|)}

\section{Path-Based Algorithms}
\label{sec:pathbasedalgorithms}

\index{static traffic assignment!algorithms!path-based|(}
This section introduces equilibrium algorithms which work in the space of path flows $H$, rather than the space of link flows $X$.
These tend to be faster, especially when high precision solutions are needed, but they require more computer memory.
Furthermore, achieving the full potential of these algorithms for rapid convergence requires considerably more programming skill than for link-based algorithms.

\index{static traffic assignment!algorithms!link-based method drawbacks|(}
Before explaining path-based algorithms, it is worth explaining why link-based algorithms are slow to converge.
The following criticisms are specifically aimed at the Frank-Wolfe algorithm, but apply to other link-based algorithms as well.

\begin{description}
\item[It treats all OD pairs equally.]  If an OD pair is close to equilibrium, only a small flow shift among its paths is needed, while if an OD pair is far from equilibrium, a larger flow shift is needed.
The Frank-Wolfe method uses the same $\lambda$ shift for all OD pairs regardless of how close or far away each one is from equilibrium.
\item[It uses a restricted set of search directions.]  If one imagines the space of feasible traffic assignments, the ``all-or-nothing'' $\mathbf{x^*}$ solutions generated by the Frank-Wolfe algorithm represent extreme points or corners of this region.
In other words, the Frank-Wolfe algorithm is only capable of moving towards a corner.
Initially, this is fine, but as one nears the optimal solution, this results in extensive zig-zagging when a much more direct path exists  (Figure~\ref{fig:fwfail}).
Conjugate Frank-Wolfe is a bit better in this regard, but still uses a fairly restrictive set of target points.
\item[It is unable to erase cyclic flows.]  Consider the network in Figure~\ref{fig:fwcycle}, with the flows as shown.
Such a flow might easily arise if $[1,2,3,4]$ is the shortest path during the first iteration of Frank-Wolfe, $[1,3,2,4]$ is the shortest path during the second, and $\lambda = 0.3$.
With only one OD pair, it is impossible for both links $(2,3)$ and $(3,2)$ to be used at equilibrium, as discussed in Section~\ref{sec:bushes}.\index{static traffic assignment!properties!acyclicity}
However, the Frank-Wolfe method will always leave some flow on both links unless $\lambda = 1$ at any iteration (which is exceedingly rare, especially in later iterations when $\lambda$ is typically very close to zero).
\end{description}
\index{static traffic assignment!algorithms!link-based method drawbacks|)}

\stevefig{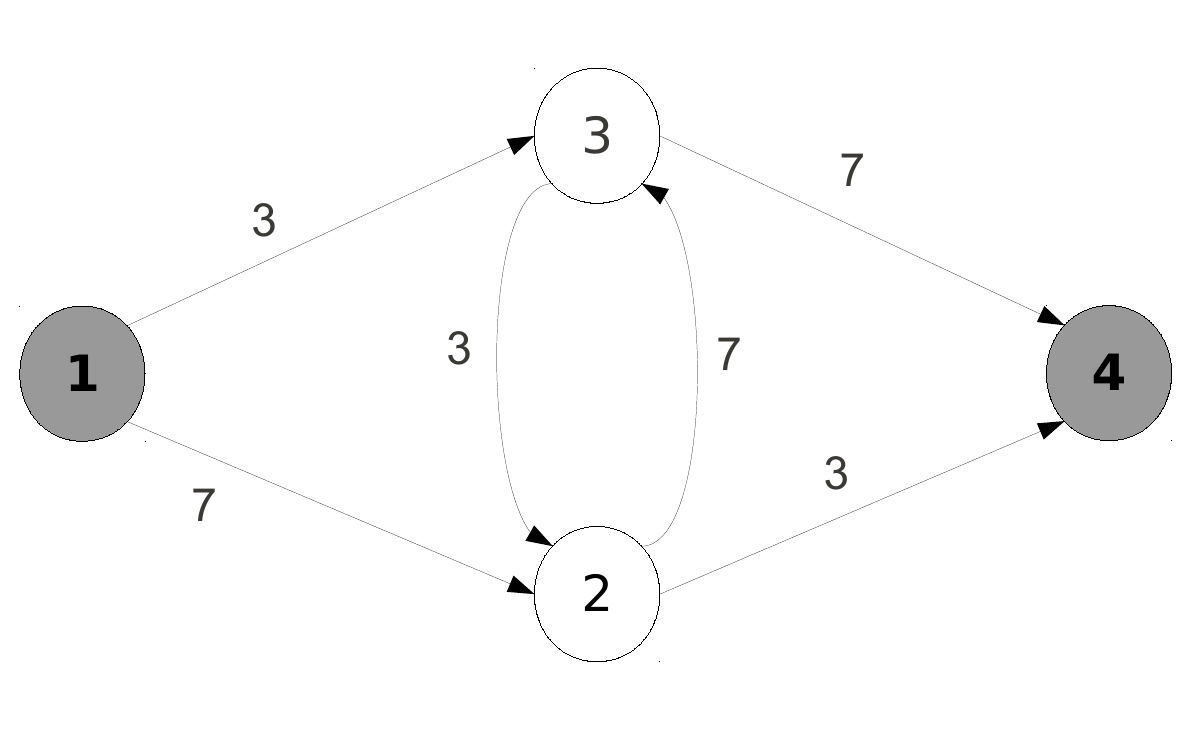}{A persistent cycle in Frank-Wolfe. (Link flows shown.)}{0.6\textwidth}

These difficulties can all be avoided by tracking the path flows $\mb{h}$, rather than the link flows $\mb{x}$.
The path flows contain much more information, tracking flow by origin and destination, as opposed to link flows which are aggregated together.
On balance, the number of elements in the path flow vector is many orders of magnitude larger than that of the link flows, easily numbering in the millions for realistic networks.
Algorithms which require us to first list off all paths in the network are not tractable.
Instead, path-based algorithms only track the paths which an OD pair actually uses, that is the set $\hat{\Pi}^{rs} = \{ \pi \in \mathcal{\pi}^{rs} : h^\pi > 0 \}$.\footnote{Those familiar with other types of optimization problems might recognize this as a column generation\index{column generation} scheme.}  This is often referred to as the set of \emph{working paths}\index{working path!see {path,set of used paths}}\index{path!set of used paths} for each OD pair.\footnote{You may recognize similarities with the ``trial-and-error'' method from Chapter~\ref{chp:equilibrium}.
Path-based algorithms are essentially a more clever form of this method.}  A rough description of path-based algorithms can then be described as
\begin{enumerate}
\item Initialize $\hat{\Pi}^{rs} \leftarrow \emptyset$\label{not:emptyset} for all OD pairs.
\item Repeat the following steps for each OD pair $(r,s)$:
\begin{enumerate}[(a)]
    \item Find the shortest path $\hat{\pi}_{rs}$.
Add it to $\hat{\Pi}^{rs}$ if it's not already used.
    \item Shift travelers among paths to get closer to equilibrium.
    \item Update travel times.
\end{enumerate}
\item Drop paths from $\hat{\Pi}^{rs}$ if they are no longer used; return to step 2 unless a convergence criterion is satisfied.
\end{enumerate}

On the surface, this scheme looks quite similar to the link-based methods presented earlier.
Why might it converge faster?
Recall the three factors described above.
First, each OD pair is now being treated independently.
With the method of successive averages and Frank-Wolfe, the same step-size $\lambda$ was applied across all links (and therefore, across all OD pairs).
If one OD pair is very close to equilibrium, while another is far away, we should probably make a finer adjustment to the first OD pair, and a larger adjustment to the second one.
Link-based methods allow no such finesse, and instead bluntly apply the same $\lambda$ to all origins.
In practice, this means that $\lambda$ becomes very small after a few iterations: we can't move very far without disturbing an OD pair which is already close to equilibrium.
As a result of this, it takes a really long time to solve OD pairs which are far from equilibrium.
Equivalently, the set of search directions is broader in that we can vary the step size by OD pairs; and lastly, it is quite possible to erase cyclic flows in a path-based context, because the extra precision allows us to take larger steps.

\index{projected gradient|see {static traffic assignment, algorithms, manifold suboptimization}}
Step 2b is where path-based algorithms differ.
This section describes two path-based algorithms: \emph{gradient projection} and \emph{manifold suboptimization}.
The latter is sometimes called \emph{projected gradient}, because both algorithms use the same two ingredients: exploiting the fact that the \emph{gradient} is the direction of steepest ascent (and therefore, in a minimization problem, we should move in the opposite direction to descend as quickly as possible); and having to consider the constraints in the problem by using a \emph{projection} operation to stay in the feasible set.
(Recall from Chapter~\ref{chp:mathematicalpreliminaries} that projection involves finding the point within a set which lies closest to another point.)
Where they differ is in the order these steps are applied.

In gradient projection, we first take a step in the opposite direction of the gradient, which will typically result in an infeasible point.
The projection is done \emph{after} the flow shift, not before: we do the projection after we make use of the gradient, so to speak.
We then apply a projection operation to return to the feasible set.
In manifold suboptimization, we first calculate the gradient, then project the gradient onto the feasible set so that our search direction respects the demand constraint.

\subsection{Gradient projection}
\label{sec:gradientprojection}

\index{static traffic assignment!algorithms!gradient projection|(}
The gradient projection method identifies the direction of steepest descent for the Beckmann function, and calculates a new point in this direction.
In case this point is infeasible, a projection\index{projection!applications} operation is applied to return to the closest feasible point.
Writing the Beckmann function in terms of path flows\index{Beckmann function!in path flows} $\mb{h}$, rather than link flows as is customary, we have
\labeleqn{beckmannpath}{
f(\mb{h}) = \sum_{(i,j) \in A} \int_0^{\sum_{\pi \in \Pi} \delta_{ij}^\pi h^\pi} t_{ij}(x)~dx
\,,}
and its partial derivative with respect to any path flow variable is
\labeleqn{beckmannpathder}{
\pdr{f}{h^\pi} = \sum_{(i,j) \in A} \delta_{ij}^\pi t_{ij} \myp{ \sum_{\pi' \in \Pi} \delta_{ij}^{\pi'} h^{\pi'} } = \sum_{(i,j) \in A} \delta_{ij}^\pi t_{ij}(x_{ij}) = c^\pi
\,,}
so the direction of steepest descent $\mb{s}$ is the negative gradient:\label{not:vect}
\labeleqn{neggrad}{
-\nabla_\mb{h} f = -\mr{vect} (c^\pi)
\,.}

So, a first attempt at gradient projection would be to update $\mb{h} \leftarrow \mr{proj}_H (\mb{h} - \nabla_\mb{h} f)$.
Unfortunately, projecting onto the set $H$ is not particularly easy.
If we apply a suitable change of variables, though, the projection can be made much easier.

Define the \emph{basic path}\index{path!basic path} for OD pair $(r,s)$ to be a path $\hat{\pi}_{rs}$\label{not:pirshat} with minimum travel time.
All of the other paths are called \emph{nonbasic paths}.\index{path!nonbasic path}
We can eliminate the basic path flow variable by expressing it in terms of the nonbasic path flows:
\labeleqn{nonbasic}{
h^{\hat{\pi}_{rs}} = d^{rs} - \sum_{\pi \in \hat{\Pi}^{rs} : \pi \neq \hat{\pi}_{rs}} h^\pi
\,.}
Substituting this into the path-based Beckmann function~\eqn{beckmannpath}, the partial derivative with respect to one of the nonbasic path flows is now
\labeleqn{nonbasicgrad}{
\pdr{\hat{f}}{h^\pi} = c^\pi - c^{\hat{\pi}_{rs}} \qquad \forall (r,s) \in Z^2, \pi \in \hat{\Pi}^{rs} - \myc{\hat{\pi}_{rs}}
\,,}
denoting the Beckmann function as $\hat{f}$ instead of $f$ because the function has been modified by using~\eqn{nonbasic} to eliminate some of the path flow variables.

So, the change in path flows will be the negative of the gradient.
Since the gradient is given by~\eqn{nonbasicgrad} and $c^\pi \geq c^{\hat{\pi}_{rs}}$ (because the basic path is by definition the shortest one), moving in this direction means that every nonbasic path flow will \emph{decrease}.

Since the transformation~\eqn{nonbasic} eliminated the demand satisfaction constraint, the only remaining constraint is that the nonbasic path flow variables be nonnegative.
Projecting onto this set is trivial: if any of the nonbasic path flow variables is negative after taking a step, simply set it to zero.
At this point, the basic path flow can be calculated through equation~\eqn{nonbasic}.\index{projection!applications}

Furthermore, the larger the difference in travel times, the larger the corresponding element of the gradient will be.
This suggests that more flow be shifted away from paths with higher travel times.
We can go a step further, and estimate directly how much flow should be shifted from a nonbasic path to a basic path to equalize the travel times, using Newton's method.\index{Newton's method}

Let $\Delta h$\label{not:Deltah} denote the amount of flow we shift \emph{away} from nonbasic path $\pi$ and \emph{onto} the basic path $\hat{\pi}$, and let $c_\pi(\Delta h)$ and $c_{\hat{\pi}} (\Delta h)$ denote the travel times on path $\pi$ and $\hat{\pi}$ after we make such a shift.
We want to choose $\Delta h$ so these costs are equal, that is, so 
\labeleqn{gdefinition}{
g(\Delta h) = c_\pi(\Delta h) - c_{\hat{\pi}}(\Delta h) = 0
\,;}
that is, $g$ is simply the difference in travel times between the two paths.\label{not:g}
To apply Newton's method, we need to find the derivative of $g$ with respect to $\Delta h$.

Using the relationships between link travel times and path travel times, we have
\[
g(\Delta h) = \sum_{(i,j)  \in \mathcal{A}} (\delta_{ij}^\pi - \delta_{ij}^{\hat{\pi}}) t_{ij}(x_{ij}(\Delta h))
\,,\]
so
\[
g'(\Delta h) = \sum_{(i,j) \in \mathcal{A}} (\delta_{ij}^\pi - \delta_{ij}^{\hat{\pi}}) \frac{dt_{ij}}{dx_{ij}} \frac{dx_{ij}}{d \Delta h}
\]
by the chain rule.
For each arc, there are four possible cases:
\begin{description}
\item [Case I:] $\delta_{ij}^\pi = \delta_{ij}^{\hat{\pi}} = 0$, that is, neither path $\pi$ nor path $\hat{\pi}$ uses link $(i,j)$.
Then $(\delta_{ij}^\pi - \delta_{ij}^{\hat{\pi}}) \frac{dt_{ij}}{dx_{ij}} \frac{dx_{ij}}{d \Delta h} = 0$ and this link does not contribute to the derivative.
\item [Case II:] $\delta_{ij}^\pi = \delta_{ij}^{\hat{\pi}} = 1$, that is, both paths $\pi$ and path $\hat{\pi}$ use link $(i,j)$.
Then $(\delta_{ij}^\pi - \delta_{ij}^{\hat{\pi}}) \frac{dt_{ij}}{dx_{ij}} \frac{dx_{ij}}{d \Delta h} = 0$ and this link again does not contribute to the derivative.
(Another way to think of it: since both paths use this arc, its total flow will not change if we shift travelers from one path to another.)
\item [Case III:] $\delta_{ij}^\pi = 1$ and $\delta_{ij}^{\hat{\pi}} = 0$, that is, path $\pi$ uses link $(i,j)$, but path $\hat{\pi}$ does not.
Then $(\delta_{ij}^\pi - \delta_{ij}^{\hat{\pi}}) \frac{dt_{ij}}{dx_{ij}} \frac{dx_{ij}}{\Delta h} =    \frac{dt_{ij}}{dx_{ij}} \frac{dx_{ij}}{d \Delta h} = -\frac{dt_{ij}}{dx_{ij}}$ since $\frac{dx_{ij}}{d \Delta h} = -1$.
\item [Case IV:] $\delta_{ij}^\pi = 0$ and $\delta_{ij}^{\hat{\pi}} = 1$, that is, path $\hat{\pi}$ uses link $(i,j)$, but path $\pi$ does not.
Then $(\delta_{ij}^\pi - \delta_{ij}^{\hat{\pi}}) \frac{dt_{ij}}{dx_{ij}} \frac{dx_{ij}}{\Delta h} = -\frac{dt_{ij}}{dx_{ij}} \frac{dx_{ij}}{d \Delta h} = -\frac{dt_{ij}}{dx_{ij}}$ since $\frac{dx_{ij}}{d \Delta h} = 1$.
\end{description}
Putting it all together, the only terms which contribute to the derivative are the links which are used by either $\pi$ or $\hat{\pi}$, but not both.
Let $A_1$, $A_2$, $A_3$, and $A_4$ denote the sets of links falling into the four cases listed above.
Then 
\[
g'(\Delta h) = -\sum_{(i,j) \in A_3 \cup A_4} \frac{dt_{ij}}{dx_{ij}}
\,.\]
which is simply the negative sum of the derivatives of these links, evaluated at the current link flows.

Then, starting with an initial guess of $\Delta h = 0$, one step of Newton's method gives an improved guess of
\[
\Delta h = 0 - g(0)/g'(0) = \frac{c_\pi - c_{\hat{\pi}}}{\sum_{a \in A_3 \cup A_4} \frac{dt_{ij}}{dx_{ij}}}
\,.\]
That is, the recommended Newton shift is given by the difference in path costs, divided by the sum of the derivatives of the link performance functions for links used by one path or the other, but not both.
Therefore, the updated nonbasic and basic path flows are given by
\[
h_{\hat{\pi}} \leftarrow h_{\hat{\pi}} + \frac{c_\pi - c_{\hat{\pi}}}{\sum_{a \in A_3 \cup A_4} \frac{dt_{ij}}{dx_{ij}}}
\]
and
\[
h_\pi \leftarrow h_\pi - \frac{c_\pi - c_{\hat{\pi}}}{\sum_{a \in A_3 \cup A_4} \frac{dt_{ij}}{dx_{ij}}}
\,.\]
This process is repeated for every nonbasic path.

This is demonstrated on the example in Figure~\ref{fig:hw1net} as follows.

\begin{description}
\item[Iteration 1, Step 1.] Initially $\hat{\Pi}^{13} = \hat{\Pi}^{24} = \emptyset$.
\item[Iteration 1, Step 2a, OD pair (1,3).] Find the shortest path for (1,3): with no travelers on the network, the top link has a travel time of 10.
This is not in the set of used paths, so include it: $\hat{\Pi}^{13} = \{[1,3]\}$.
\item[Iteration 1, Step 2b, OD pair (1,3).] Since there is only one used path, we simply have $h_{[1,3]}^{13} = 5000$.
\item[Iteration 1, Step 2c, OD pair (1,3).] Update travel times: $$\mathbf{t} = \vect{60 & 10 & 10 & 10 & 10 & 10 & 10 }$$
\item[Iteration 1, Step 2a, OD pair (2,4).] Find the shortest path for (2,4): with no travelers on the network, link 7 has a travel time of 10.
This is not in the set of used paths, so include it: $\hat{\Pi}^{24} = \{[2,4]\}$.
\item[Iteration 1, Step 2b.] Since there is only one used path, $h_{[2,4]}^{24} = 10000$.
\item[Iteration 1, Step 2c.] Update travel times: \[ \mathbf{t} = \vect{60 & 10 & 10 & 10 & 10 & 10 & 110 } \]
\item[Iteration 1, Step 3.] All paths are used, so return to step 2.
The relative gap is $\gamma_1= 2.11$.
\item[Iteration 2, Step 2a, OD pair (1,3).] With the new travel times, the shortest path is [1,5,6,3].
This is not part of the set of used paths, so we add it: $\hat{\Pi}^{13} = \{[1,3], [1,5,6,3]\}$.
\item[Iteration 2, Step 2b, OD pair (1,3).] The difference in travel times between the paths is 30 minutes; and the sum of the derivatives of links 1, 2, 3, and 4 is 0.04.
So we shift 30/0.04 = 750 vehicles from [1,3] to [1,5,6,3], $h_{[1,3]}^{13} = 4250$ and $h_{[1,5,6,3]}^{13} = 750$.
\item[Iteration 2, Step 2c, OD pair (1,3).] Update travel times: \[ \mathbf{t} = \vect{52.5 & 17.5 & 17.5 & 17.5 & 10 & 10 & 110 } \]  \emph{Note that the two paths have \emph{exactly} the same cost after only one step!
This is because Newton's method is exact for linear functions.}
\item[Iteration 2, Step 2a, OD pair (2,4).] With the new travel times, the shortest path is now [2,5,6,4].
This is not part of the set of used paths, so we add it: $\hat{\Pi}^{24} = \{[2,4], [2,5,6,4]\}$.
\item[Iteration 2, Step 2b, OD pair (2,4).] The difference in travel times between the paths is 72.5 minutes; and the sum of the derivatives of links 5, 3, 6, and 7 is 0.04.
So we shift 72.5/0.04 = 1812.5 vehicles from [2,4] to [2,5,6,4], $h_{[2,4]}^{24} = 8187.5$ and $h_{[2,5,6,4]}^{24} = 1812.5$.
\item[Iteration 2, Step 2c, OD pair (2,4).] Update travel times: \[ \mathbf{t} = \vect{52.5 & 17.5 & 35.625 & 17.5 & 28.125 & 28.125 & 91.875 } \] \emph{Note that the two paths again have exactly the same cost.
However, the equilibrium for the first OD pair has been disturbed.}
\item[Iteration 2, Step 3.] All paths are used, so return to step 2.
The relative gap is $\gamma_1= 0.0115$.
\item[Iteration 3, Step 2a, OD pair (1,3).] The shortest path is now [1,3], which is already in the set of used paths, so nothing to do here.
\item[Iteration 3, Step 2b, OD pair (1,3).] The difference in travel times between the paths is 18.125  minutes; and the sum of the derivatives of links 1, 2, 3, and 4 is 0.04.
So we shift 18.125/0.04 = 453 vehicles from [1,5,6,3] to [1,3], $h_{[1,3]}^{13} = 4703$ and $h_{[1,5,6,3]}^{13} = 297$.
\item[Iteration 3, Step 2c, OD pair (1,3).] Update travel times: \[ \mathbf{t} = \vect{57.0 & 13.0 & 31.1 & 13.0 & 28.1 & 28.1 & 91.9 } \] \emph{Again the first OD pair is at equilibrium, up to rounding error.}
\item[Iteration 3, Step 2a, OD pair (2,4).] The shortest path is again [2,5,6,4].
This is already in the set of used paths, so nothing to do here.
\item[Iteration 3, Step 2b, OD pair (2,4).] The difference in travel times between the paths is 4.6 minutes; and the sum of the derivatives of links 5, 3, 6, and 7 is 0.04.
So we shift 4.6/0.04 = 114 vehicles from [2,4] to [2,5,6,4], $h_{[2,4]}^{24} = 8073$ and $h_{[2,5,6,4]}^{24} = 1927$.
\item[Iteration 3, Step 2c, OD pair (2,4).] Update travel times: \[ \mathbf{t} = \vect{57.0 & 13.0 & 32.2 & 13.0 & 29.3 & 29.3 & 90.7 } \]  
\item[Iteration 3, Step 3.] All paths are used, so return to step 2.
The relative gap is $\gamma_1= 0.00028$.
\end{description}
Note that after three iterations of gradient projection, the relative gap is two orders of magnitude smaller than that from the Frank-Wolfe algorithm.
Although not demonstrated here for reasons of space, the performance of gradient projection relative to Frank-Wolfe actually \emph{improves} from here on out.
Frank-Wolfe usually does most of its work in the first few iterations, and then converges \emph{very} slowly after that.\footnote{It's been said that Frank-Wolfe converges, but just barely.}  On the other hand, gradient projection maintains a steady rate of progress throughout, with a nearly constant proportionate decrease in gap from one iteration to the next.
\index{static traffic assignment!algorithms!gradient projection|)}

\subsection{Manifold suboptimization}
\label{sec:manifoldsuboptimization}

\index{manifold suboptimization|see {static traffic assignment, algorithms}}
\index{static traffic assignment!algorithms!manifold suboptimization|(}
The manifold suboptimization algorithm is also based on the Beckmann formulation, applying an algorithm of Rosen from nonlinear optimization.
As with gradient projection, moving in the direction of the negative gradient will travel along the direction of steepest descent.
Such a direction will tend towards the equilibrium solution, which minimizes the Beckmann function.
However, we need to be careful not to leave the region of feasible path flows $H$, making sure that our search direction still satisfies the demand constraint and retains nonnegative path flows.
Unlike gradient projection, which fixes an infeasible move by projecting back to the feasible region, manifold suboptimization avoids ever leaving the feasible region in the first place by modifying the steepest descent direction.
This algorithm goes by several names in the literature; sometimes called projected gradient and at other times even called gradient projection (even though it is different from the gradient projection method of Section~\ref{sec:gradientprojection}).
To avoid confusion, we refer to it by manifold suboptimization, a convention also followed by~\cite{bertsekas_nlp}. 

Recall from equation~\eqn{neggrad} that the gradient of the Beckmann function\index{Beckmann function!gradient} is
\labeleqn{beckmanngrad}{
\nabla_\mb{h} f = \mr{vect} (c^\pi)
\,.}
Assuming that the current path flow solution $\mb{h}$ is feasible, we must move in a direction that does not violate any of the constraints, only using the working paths $\hat{\Pi}^{rs}$.
For instance, if the demand constraint $\sum_{\pi \in \hat{\Pi}^{rs}} h^\pi = d^{rs}$ is satisfied for all OD pairs $(r,s)$, it must remain so after taking a step in the direction $\Delta \mb{h}$:
\labeleqn{feasiblespace}{
\sum_{\pi \in \hat{\Pi}^{rs}} (h^\pi + \Delta  h^\pi) = d^{rs}
\,.}
This in turn implies
\begin{align}
\sum_{\pi \in \hat{\Pi}^{rs}} \Delta  h^\pi &= d^{rs} - \sum_{\pi \in \hat{\Pi}^{rs}} h^\pi  \\
\sum_{\pi \in \hat{\Pi}^{rs}} \Delta  h^\pi &= 0
\end{align}

So, we need to find the projection\index{projection!applications} of the steepest descent direction $\mb{s}$ onto the space $\Delta H \equiv \{ \Delta \mb{h} : \sum_{\pi \in \hat{\Pi}i^{rs}} \Delta  h^\pi = 0 \}$.\label{not:DeltaH}
It turns out that this projection has a remarkably simple closed form expression.

\begin{prp}
Using $\bar{c}^{rs} = (1/|\hat{\Pi}^{rs}|) \sum_{\pi \in \hat{\Pi}^{rs}} c^\pi$ to denote the average travel time of the working paths for OD pair $(r,s)$, the direction $\mb{s'}$ with components $s'_\pi = c^\pi - \bar{c}^{rs}$ is the projection of the steepest descent direction $\mb{s}$ onto the set $\Delta H$.
\end{prp}
\begin{proof}
To show that $\mb{s'}$ is the projection of $\mb{s}$ onto $\Delta H$, we must show that $\mb{s'} \in \Delta H$, and that $\mb{s} - \mb{s'}$ is orthogonal to $\Delta H$.
Regarding the first part, we have
\labeleqn{projmember}{
\sum_{\pi \in \hat{\Pi}^{rs}} s'_\pi = \sum_{\pi \in \hat{\Pi}^{rs}} (c^\pi - \bar{c}^{rs}) = \sum_{\pi \in \hat{\Pi}^{rs}} c^\pi - |\hat{\Pi}^{rs}| \bar{c}^{rs} = 0
\,,}
so $\mb{s'} \in \Delta H$.

Regarding the second part, let $\Delta \mb{h}$ be any vector in $\Delta H$.
We now show that $(\mb{s} - \mb{s'}) \cdot \Delta \mb{h} = 0$.
We have
\begin{align}
(\mb{s} - \mb{s'}) \cdot \Delta \mb{h} &= \sum_{\pi \in \hat{\Pi}^{rs}} [c^\pi - (c^\pi - \bar{c}^{rs}) ] \Delta h^\pi \\
&= -\bar{c}^{rs} \sum_{\pi \in \hat{\Pi}^P{rs}} \Delta h^\pi \\
&= 0
\end{align}
since $\Delta \mb{h} \in \Delta H$.
\end{proof}

So, given current path flows $\mb{h}^{rs}$ for OD pair $(r,s)$, we use $\Delta \mb{h}^{rs} = \mr{vect} (c^\pi - \bar{c}^{rs})$ as the search direction.
To update the path flows $\mb{h}^{rs} \leftarrow \mb{h}^{rs} + \mu \Delta \mb{h}^{rs}$, we need an expression for the step size.
For any path $\pi \in \hat{\Pi}^{rs}$ for which $\Delta h^\pi < 0$, the new path flows would be infeasible if $\mu > h^\pi / \Delta h^\pi$.
Therefore, the largest possible step size\label{not:mubar} is
\labeleqn{maxstep}{
\bar{\mu} = \min_{\pi \in \hat{\Pi}^{rs} : \Delta h^\pi < 0} \frac{h^\pi}{\Delta h^\pi}
\,.}
The actual step size $\mu \in [0, \bar{\mu}]$ should be chosen to minimize the Beckmann function.
This can be done either through bisection or one or more iterations of Newton's method.
\index{optimization!line search!bisection}
\index{Newton's method}
\index{static traffic assignment!algorithms!manifold suboptimization|)}
\index{static traffic assignment!algorithms!path-based|)}

\section{Bush-Based Algorithms}
\label{sec:bushbased}

\index{static traffic assignment!algorithms!bush-based|(}
\index{static traffic assignment!bush-based solution|(}
As seen in the previous section, path-based algorithms converge faster than link-based algorithms, and allow us find much more accurate equilibrium solutions in a fraction of the time.
The prime disadvantage is a huge memory requirement, with potentially millions of paths available for use in large networks.
Furthermore, many of these paths are ``redundant'' in some way, as they overlap: the same sequence of links might be used by many different paths.

Bush-based algorithms (also known as origin-based algorithms) try to address these limitations.
Rather than treating each OD pair separately (as link-based algorithms do), bush-based algorithms simultaneously consider every destination associated with a single origin.
Instead of considering every possible used path (as path-based algorithms do, and there are potentially very many of these), they maintain a set of links called a \emph{bush},\index{bush} which are the only links which travelers from that origin are permitted to use.
One can think of a bush as the set of links obtained by superimposing all of the paths used by travelers starting from an origin.
In particular, a bush must be:
\begin{itemize}
\item \emph{Connected}; that is, using only links in the bush, it is possible to reach every node which was reachable in the original network.
\item \emph{Acyclic}; that is, no path using only bush links can pass the same node more than once.
This is not restrictive, because travelers trying to minimize their own travel time would never cycle back to the same node, and greatly speeds up the algorithm, because acyclic networks are much simpler and admit much faster methods for finding shortest paths and other quantities of interest.
\end{itemize}
There is no particular reason why bushes have to be ``origin-based'' rather than ``destination-based,'' and all of the results in this section can be derived in a parallel way for bushes terminating at a common destination, rather than starting at a common origin.

An example of a bush is shown in panels (a) and (b) of Figure~\ref{fig:bushexample}, where the thickly-shaded links are part of the bush, and the lightly-shaded links are not.
The bush in panel (a) is a special type of bush known as a \emph{tree}, which has exactly one path from the origin to every node.
The thickly-shaded links in panel (c) do not form a bush, because it is not connected; there is no way to reach the nodes at the bottom of the network only using bush links.
Likewise, the thick links in panel (d) do not form a bush either, because a cycle exists and it would be possible to revisit some nodes multiple times using the bush links (find them!).

\stevefig{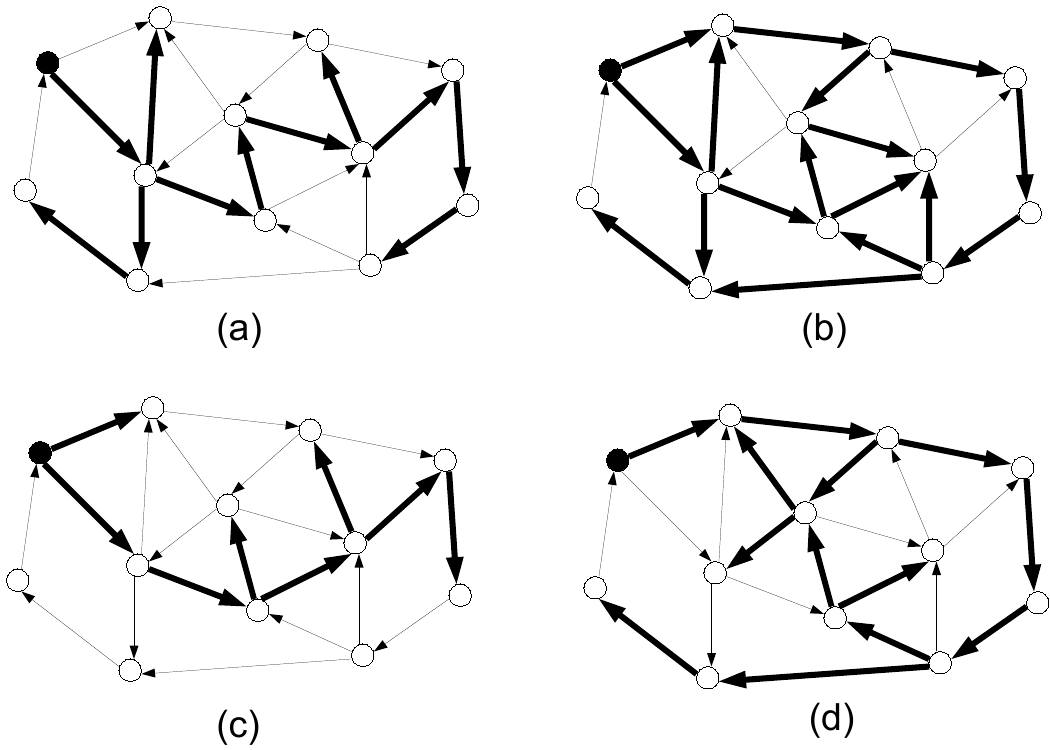}{Examples of bushes (panels (a) and (b)) and non-bushes (panels (c) and (d)).}{\textwidth}

Notice that because a bush is acyclic, it can never include both directions of a two-way link.
This implies that at equilibrium, \emph{on every street travelers from the same origin must all be traveling in the same direction.}  (This follows from Proposition~\ref{prp:eqmbushacyclic} in Section~\ref{sec:bushes}).
Interestingly, link-based and path-based algorithms cannot enforce this requirement easily; and this is yet another reason that bush-based algorithms are a good option for solving equilibrium.
There is one bush for every origin; this means that if there are $z$ origins and $m$ links, we need to keep track of at most $zm$ values.
By contrast, a link-based approach (such as Frank-Wolfe) requires storage of only $m$ values to represent a solution, while a path-based approach could conceivably require $z^2 2^m$ values.\footnote{These values are very approximate, but give you an idea of the scale.}  

The first well-known origin-based algorithm was developed by Hillel Bar-Gera in his dissertation (circa 2000), and was simply called origin-based assignment (OBA).
Bob Dial developed another, simpler method that was published in 2006 as ``Algorithm B,'' which also seems to work faster than OBA.
Yu (Marco) Nie compared both algorithms and developed additional variations by combining features of both, and contributing some ideas of his own.
Most recently, Guido Gentile has developed the LUCE algorithm, and Hillel Bar-Gera has provided a new algorithm called TAPAS which simultaneously solves for equilibrium and proportional link flows (approximating entropy maximization).

All bush-based algorithms operate according to the same general scheme:
\begin{enumerate}
\item Start with initial bushes for each origin (the shortest path tree\index{tree!shortest path} with free-flow times is often used as a starting point).
\item Shift flows within each bush to bring each origin closer to equilibrium.
\item Improve the bushes by adding links which can reduce travel times, and by removing unused links.
Return to step 2.
\end{enumerate}

Step 1 is fairly self-explanatory: collecting the shortest paths from an origin to every destination into one bush results in a connected and acyclic set of links, and we can simply load the travel demand to each destination on the unique path in the bush.

Step 2 is where most bush-based algorithms differ.\footnote{If this is starting to sound familiar, good! How to shift flows to move closer to equilibrium is also where link-based and path-based algorithms differ most.}
This section explains three different methods for shifting flows within a bush: Algorithm B, origin-based assignment (OBA), and linear user cost equilibrium (LUCE).
Each of these algorithms is based on labels calculated for bush nodes and links, and Section~\ref{sec:bushfundamentals} defines each of these.
The three algorithms themselves are presented in Section~\ref{sec:bushshift}.

Step 3 requires a little bit of thought to determine how to adjust the bush links themselves to allow movement towards equilibrium when Step 2 is performed on the new bushes.
Section~\ref{sec:bushadjust} shows how this can be done.

\subsection{Bush labels}
\label{sec:bushfundamentals}

\index{bush!labels|(}
Algorithm B, OBA, and LUCE all make use of ``bush labels'' associated with each node and link in a bush.
All of these labels can be calculated in a straightforward, efficient manner using the topological order\index{topological order} within each bush.
As you read this section, pay attention to which of these labels are calculated in ascending topological order (i.e., the formulas only involve nodes with lower topological order) and which are calculated in descending order (based only on nodes with higher topological order.)  In the discussion below, assume that a bush is given and fixed, with origin node $r$.
The set of bush links is denoted $\mc{B}$\label{not:Bush} and, for the purposes of this section, assume that the travel times $t_{ij}$ and travel time derivatives $t'_{ij}$ of all bush links are given and constant.
All of these labels are only defined for bush links, and the formulas only involve bush links.
From the standpoint of Step 2 of bush-based algorithms, non-bush links are completely irrelevant.
Table~\ref{tbl:bushchart} shows all of the labels defined in this section, and which labels are used in which algorithms, and in the bush-updating steps described in Section~\ref{sec:bushadjust}.

We start with two different ways to represent the travel patterns on each bush, starting with $x$ labels.\index{bush!labels!link flow}
The label $x^r_{ij}$ associated with each link indicates the number of travelers starting at node $r$ (the root of bush $\mc{B}$) and traveling on bush link $(i,j)$.
The superscript $\mc{B}$ indicates that we are only referring to the flow on this link associated with the bush $\mc{B}$, and the \emph{total} link flow $x_{ij}$ is the sum of $x_{ij}^r$ across all bushes $\mc{B}$.
However, using these superscripts tends to clutter formulas, and often times it is clear that we are only referring to flows within the context of a specific bush.
In this case, we can simply write $x_{ij}$ with it being understood that this label refers to the flow on a bush link.
Within this section, we are only concerned with a single bush and the superscript will be omitted for brevity.
\index{static traffic assignment!bush-based solution|)}

The network in Figure~\ref{fig:bushgrid}(a) will be used as to demonstrate the labels introduced in this section.
The thick links comprise the bush, and the link performance functions for all links in the network are shown.
The origin in this case is node 7, and the demand is 10 vehicles from node 7 to node 3.
You can verify that a topological ordering of the nodes on this bush is 7, 8, 9, 6, 4, 1, 5, 2, 3.

\stevefig{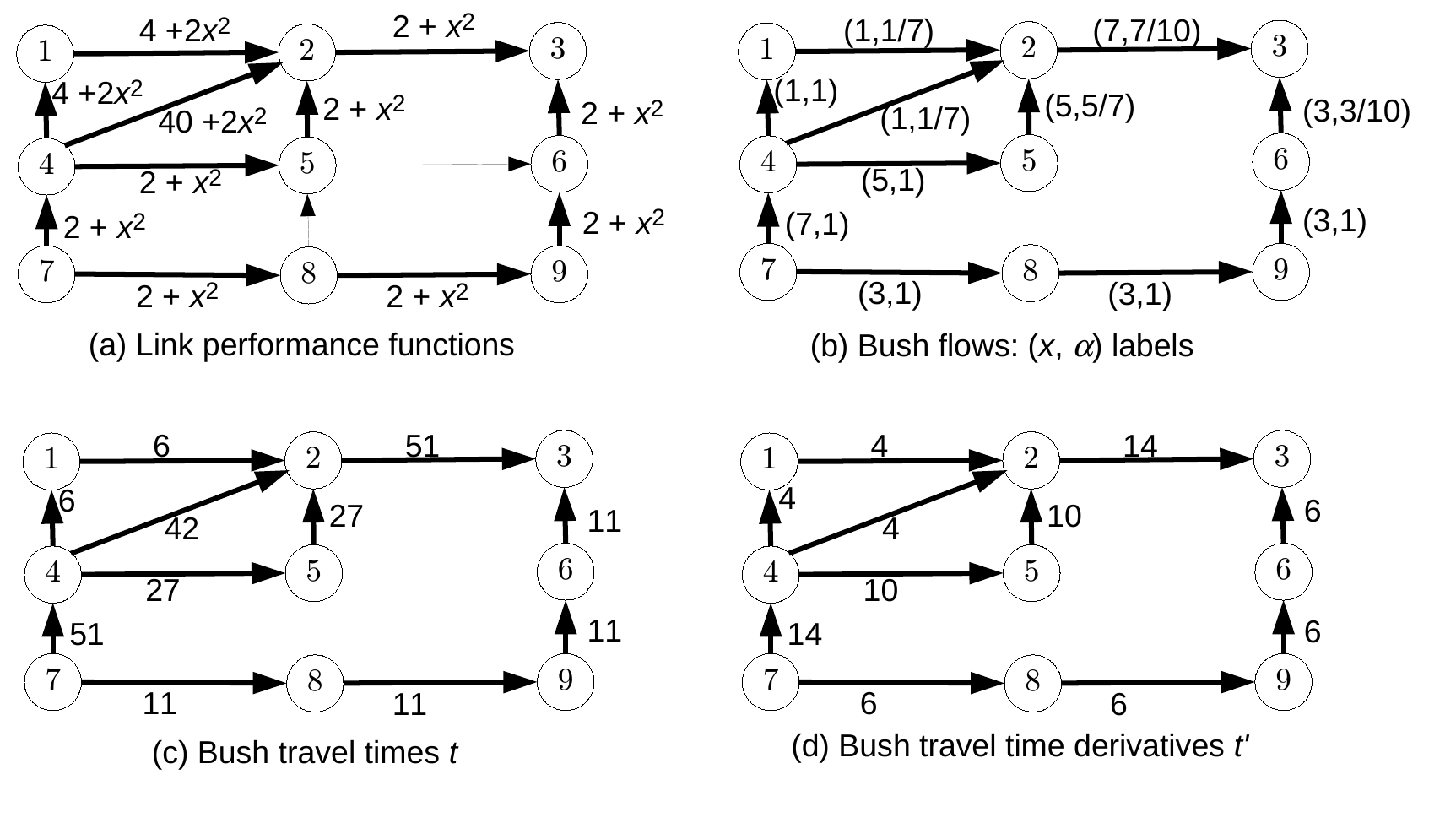}{Example bush used to demonstrate label calculation.}{\textwidth}

The corresponding label $x_i$ associated with each node\index{bush!labels!node flow} indicates the total number of vehicles using node $i$ on the bush $\mc{B}$, including flow which is terminating at node $i$.
The flow conservation equations relating $x_{ij}$ and $x_i$ labels are as follows:\label{not:hi}
\[
x_i = \sum_{(h,i) \in \mc{B}} x_{hi} = d^{ri} + \sum_{(i,j) \in \mc{B}} x_{ij}   \qquad \forall i \in N  \label{eqn:linknodeflow}
\,,\]
where the first expression defines the node flow in terms of \emph{incoming} link flows, and the second in terms of \emph{outgoing} link flows.
The two definitions are equivalent.

It is sometimes convenient to refer to the \emph{fraction}\index{bush!labels!link fraction} of the flow at a node coming from a particular link.
For any node $i$ with positive node flow ($x_i > 0$), define $\alpha_{hi}$\label{not:alphahi} to be the proportion of the node flow contributed by the incoming link $(h,i)$, that is
\labeleqn{alpha}{
\alpha_{hi} = x_{hi} / x_i
\,.}
Clearly each $\alpha_{hi}$ is nonnegative, and, by flow conservation, the sum of the $\alpha_{hi}$ values entering each node $i$ is one.
The definition of $\alpha_{hi}$ is slightly trickier when $x_i = 0$, because the formula~\eqn{alpha} then involves a division by zero.
To accommodate this case, we adopt this rule: \emph{when $x_i = 0$, the proportions $\alpha_{hi}$ may take any values whatsoever, as long as they are nonnegative and sum to one.}  It is important to be able to define $\alpha_{hi}$ values even in this case, because the flow-shifting algorithms may cause $x_i$ to become positive, and in this case we need to know how to distribute this new flow among the incoming links.

If we are given $\alpha$ labels for each bush link, it is possible to calculate the resulting node and link flows $x_i$ and $x_{ij}$, using this recursion:
\begin{align}
x_i &= d^{ri} + \sum_{(i,j) \in \mc{B}} x_{ij}  & \forall i \in N               \label{eqn:alphaxone} \\
x_{ij} &= \alpha_{ij} x_j                       & \forall (i,j) \in \mc{B}      \label{eqn:alphaxtwo}
\end{align}
The sum in~\eqn{alphaxone} is empty for the node with the highest topological order.
So we can start there, and then proceed with the calculations in backward topological order.

Figure~\ref{fig:bushgrid}(b) shows the $x$ and $\alpha$ labels for the example bush.
You should verify that the formulas~\eqn{alphaxone} and~\eqn{alphaxtwo} are consistent with these labels.
The link travel times and link travel time derivatives are shown in panels (c) and (d) of this figure.

As has been used earlier in the text for shortest path algorithms, $L$ is used to denote travel times on shortest paths.\index{bush!labels!shortest path}
The superscript $\mc{B}$ can be appended to these labels when it is necessary to indicate that these labels are for shortest paths \emph{specifically on the bush $\mc{B}$}, although it will usually be clear from context which bush is meant.
This section will omit such a superscript to avoid cluttering the formulas, and it should be understood that $L$ means the shortest path \emph{only on the bush under consideration}.
The same convention will apply to the other labels in this section.
There are $L$ labels associated with each node $i$, and with each link $(i,j) \in \mc{B}$: $L_i$ denotes the distance on the shortest path from $r$ to $i$ using only bush links, and $L_{ij}$\label{not:Lij} is the travel time which would result if you follow the shortest path to node $i$, then take link $(i,j)$.
These labels are calculated using these equations:
\begin{align}
L_r   & = 0                      &                          \\
L_{ij} & = L_i + t_{ij}           & \forall (i,j) \in \mc{B} \\
L_i   & = \min_{(h,i) \in \mc{B}} \{ L_{hi} \}  & \forall i \neq r  \label{eqn:spkey}
\end{align}
The $U$ labels are used to denote travel times on longest paths within the bush, and are calculated in a similar way.\index{bush!labels!longest path}
Like the $L$ labels, $U$\label{not:Uij} labels are calculated for both nodes and bush links, using the formulas:
\begin{align}
U_r   & = 0                      &                          \\
U_{ij} & = U_i + t_{ij}           & \forall (i,j) \in \mc{B} \\
U_i   & = \max_{(h,i) \in \mc{B}} \{ U_{hi} \}  & \forall i \neq r  \label{eqn:lpkey}
\end{align}
The $M$\label{not:Mij} labels are used to denote the \emph{average} travel times\index{bush!labels!mean cost} within the bush, recognizing that some travelers will be on longer paths and other travelers will be on shorter ones.
The node label $M_i$ represents the average travel time between origin $r$ and node $i$ across all bush paths connecting these nodes, weighted by the number of travelers using each of these paths.
The label $M_{ij}$ indicates the average travel time of vehicles after they finish traveling on link $(i,j)$, again averaging across all of the bush paths starting at the origin $r$ and ending with link $(i,j)$.
These can be calculated as follows:
\begin{align}
M_r   & = 0                      &                          \\   
M_{ij} & = M_i + t_{ij}           & \forall (i,j) \in \mc{B} \\
M_i   & = \sum_{(h,i) \in \mc{B}}  \alpha_{hi} M_{hi}   & \forall i \neq r 
\end{align}

Figure~\ref{fig:bushgrid2} shows the $L$, $U$, and $M$ labels corresponding to the example bush in panels (a) and (b).
Panel (a) shows the labels associated with links ($L_{ij}$, $M_{ij}$, and $U_{ij}$), while panel (b) shows the labels associated with nodes ($L_i$, $M_i$, and $U_i$).

\stevefig{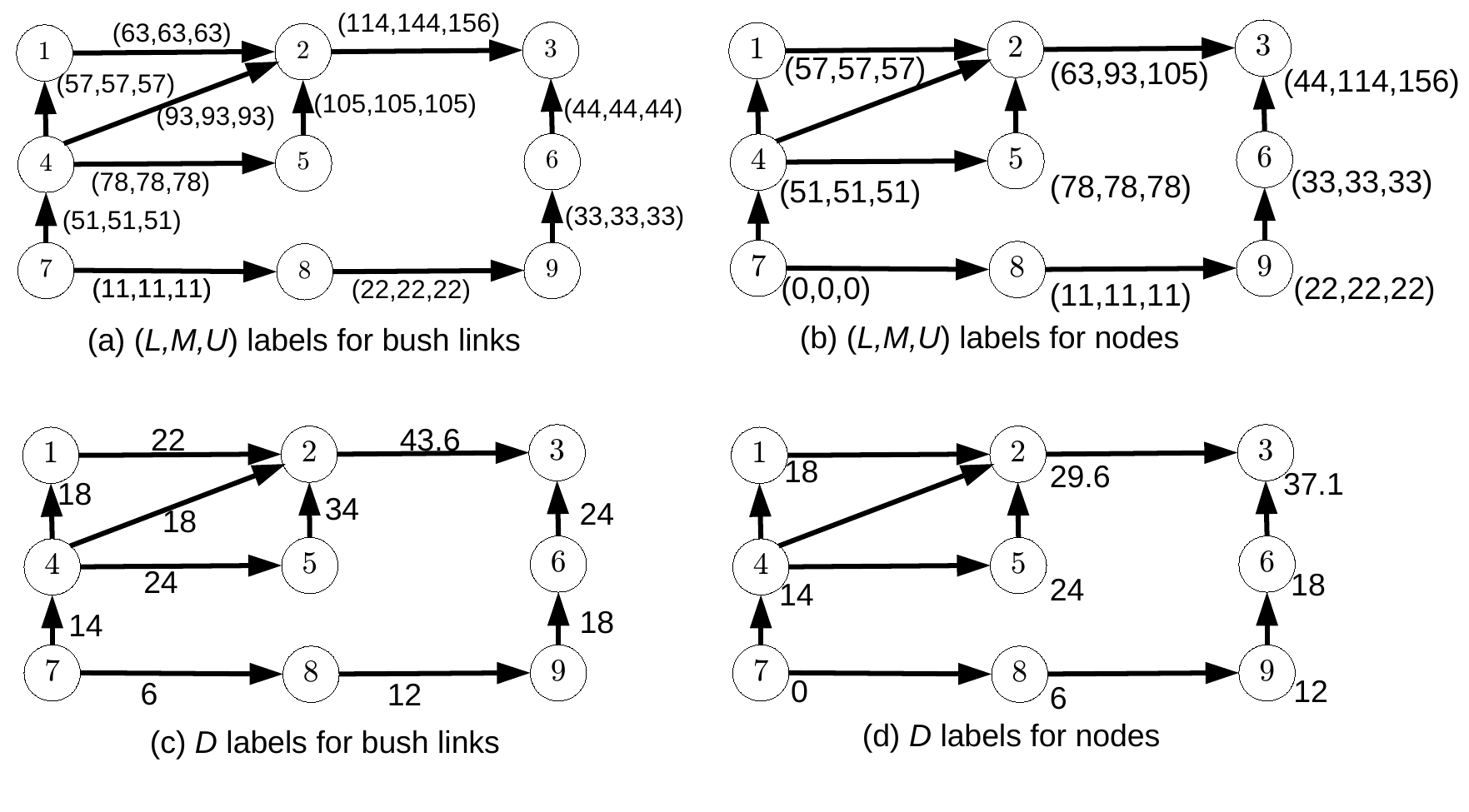}{Continuation of label calculation demonstration.}{\textwidth}

It is also useful to know how the average travel times (the $M$ labels) will change with a marginal increase in flow on a particular link or through a particular node.\index{bush!labels!cost derivative}
The set of $D$\label{not:Dij} labels are used to represent this.
As we will see shortly, these play the role of the travel time derivatives\footnote{For a mnemonic, you can associate the labels $L$, $U$, and $M$ with the \emph{lower}, \emph{upper}, and \emph{mean} travel times to nodes, while $D$ refers to the \emph{derivative}.}, which are needed in flow shifting rules based on Newton's method (like those described for path-based algorithms in the previous section.)  Ideally, the $D_{ij}$ labels would represent the derivative of $M_j$ with respect to $\alpha_{ij}$.
However, there is no known method for calculating these derivatives in an efficient way (``efficient'' meaning they can be calculated in a single topological pass).
The following formulas for $D_{ij}$ can be used as an approximation that can be calculated in a single topological pass.

\begin{align}
D_r    &= 0                      &           \\
D_{ij} &= D_i + t'_{ij}          & \forall (i,j) \in \mc{B} \\
D_i    &= \sum_{(h,i) \in \mc{B}} \alpha^2_{hi} D_{hi} + \sum_{(h,i) \in \mc{B}} \sum_{\substack{(g,i) \in \mc{B} \\ (g,i) \neq (h,i) }} \alpha_{hi} \alpha_{gi} \sqrt{D_{hi} D_{gi}} & \forall i \neq r
\end{align}

Figure~\ref{fig:bushgrid2} shows the $D$ labels associated with the example bush in panels (c) and (d).
Panel (c) shows the $D_{ij}$ labels associated with links, and panel (d) shows the $D_i$ labels associated with nodes.
\index{bush!labels|)}

\begin{table}
\centering
\caption{Bush labels used in different algorithms.
\label{tbl:bushchart}}
\begin{tabular}{|c|cccc|}
\hline
Label    & B & OBA & LUCE & Updating \\
\hline
$x$      & X &     &      &   \\
$\alpha$ &   &  X  &  X   &   \\
$L$      & X &     &      & X \\
$U$      & X &     &      & X \\
$M$      &   &  X  &  X   &   \\
$D$      &   &  X  &  X   &   \\
\hline
\end{tabular}
\end{table}

\subsection{Shifting flows on a bush}
\label{sec:bushshift}

With the labels defined in the previous subsection, we are now ready to state the Algorithm B, OBA, and LUCE flow shifting procedures.
Any of these can be used for Step 2 of the generic bush-based algorithm presented at the start of the section.
You will notice that all of these procedures follow the same general form: calculate bush labels (different labels for different algorithms) in \emph{forward} topological order, then scan each node in turn in \emph{reverse} topological order.
When scanning a node, use the labels to identify vehicles entering the node from higher-cost approaches, and shift them to paths using lower-cost approaches.
Update the $x$ and/or $\alpha$ labels accordingly, then proceed to the previous node topologically until the origin has been reached.

All three of these algorithms also make use of \emph{divergence nodes}\index{node!divergence node} (also called \emph{last common nodes}\index{last common node|see {node, divergence}}\index{node!last common|see {node, divergence}} or \emph{pseudo-origins}\index{pseudo-origin|see {node, divergence}}\index{origin!pseudo|see {node, divergence}} in the literature) as a way to limit the scope of these updates.
While the definition of divergence nodes is slightly different in these algorithms, the key idea is to find the ``closest'' node which is common to all of the paths travelers are being shifted among.
The rest of this subsection details these definitions and the role they play.

\subsubsection{Algorithm B}
\label{sec:algorithmb}

\index{static traffic assignment!algorithms!Algorithm B|(}
Algorithm B identifies the longest and shortest paths to reach a node, and shifts flows between them to equalize their travel times.
That is, when scanning a node $i$, only two paths (the longest and shortest) are considered.
It is easy to determine these paths using the $L$ and $U$ labels, tracing back the shortest and longest paths by identifying the links used for the minimum or maximum in equations~\eqn{spkey} and~\eqn{lpkey}.
Once these paths are identified, the divergence node $a$ is the last node common to both of these paths.

The shortest and longest path segments between nodes $a$ and $i$ form a \emph{pair of alternate segments}\index{pair of alternate segments}; let $\sigma_L$\label{not:sigma} and $\sigma_U$ denote these path segments.
Within Algorithm B, flow is shifted from the longest path to the shortest path, using Newton's method\index{Newton's method} to determine the amount of flow to shift:
\labeleqn{bshift}{\Delta h = \frac{(U_i - U_a) - (L_i - L_a)}{\sum_{(g,h) \in \sigma_L \cup \sigma_U} t'_{gh}} }
There is also the constraint $\Delta h \leq \min_{(i,j) \in \sigma_U} x_{ij}$ which must be imposed to ensure that all links retain nonnegative flow after the shift.
This flow is subtracted from each of the links in the longest path segment $\sigma_U$, and added to each of the links in the shortest path segment $\sigma_L$.
The derivation of this formula parallels that of gradient projection in Section~\ref{sec:gradientprojection}.
By shifting flow from the longer path to the shorter path, we will either (i) equalize the travel times on the two paths or (ii) shift all the flow onto the shorter path, and have it still be faster than the longer path.
In either case, we move closer to satisfying the equilibrium condition.

The steps of Algorithm B are as follows:
\begin{enumerate}
\item Calculate the $L$ and $U$ labels in forward topological order.
\item Let $i$ be the topologically last node in the bush.
\item Scan $i$ by performing the following steps:
\begin{enumerate}
    \item Use the $L$ and $U$ labels to determine the divergence node $a$ and the pair of alternate segments $\sigma_L$ and $\sigma_U$.
    \item Calculate $\Delta h$ using equation~\eqn{bshift} (capping $\Delta h$ at $\min_{(i,j) \in \sigma_U} x_{ij}$ if needed).
    \item Subtract $\Delta h$ from the $x$ label on each link in $\sigma_U$, and add $\Delta h$ to the $x$ label on each link in $\sigma_L$.
\end{enumerate}
\item If $i = r$, go to the next step.
Otherwise, let $i$ be the previous node topologically and return to step 3.
\item Update all travel times $t_{ij}$ and derivatives $t'_{ij}$ using the new flows $x$ (remembering to add flows from other bushes.)
\end{enumerate}

Demonstrating on the example in Figures~\ref{fig:bushgrid} and~\ref{fig:bushgrid2}, we start with the $L$ and $U$ labels as shown in Figure~\ref{fig:bushgrid2}(a) and (b), and start by letting $i = 3$, the last node topologically.
The longest path in the bush from the origin ($r = 7$) to node $i = 3$ is [7,4,5,2,3], and the shortest path is [7,8,9,6,3], as can be easily found from the $L$ and $U$ labels.
The divergence node is the last node common to both of these paths, which in this case is the origin, so $a = 7$.
Equation~\eqn{bshift} gives $\Delta h = 1.56$, so we shift this many vehicles away from the longest path and onto the shortest path, giving the flows in Figure~\ref{fig:algbstep}(a).

The second-to-last node topologically is node 2, and we repeat this process.
The longest and shortest paths from the origin to node 2 in the bush are [7,4,5,2] and [7,4,1,2], respectively.\footnote{Notice that we have not yet updated the travel time labels based on the shift at node 7.
You may do so if you wish, but it is not required for Algorithm B to work, and in all of the examples in this section travel times are not updated until all flow shifts are complete for the bush.}  The last node common to both of these paths is 4, so the divergence node is $a = 4$ and we shift flow between the pair of alternate segments [4,5,2] and [4,1,2].
Using equation~\eqn{bshift} gives $\Delta h = 1.50$.
Shifting this many vehicles from the longer segment to the shorter one gives the flows in Figure~\ref{fig:algbstep}(b).

The previous node topologically is node 5.
Since there is only one incoming bush link to node 5, there is nothing for Algorithm B to do.
To see why, notice that the longest and shortest bush paths are [7,4,5] and [7,4,5].
The divergence node would be $a = 5$, which is the same as $i$, and the ``pair of alternate segments'' is the empty paths [5] and [5].
Intuitively, since there is only one way to approach node 5, there are no ``alternate routes'' to divert incoming flow.
In fact, the same is true for all of the previous nodes topologically (1, 4, 6, 9, 8, 7 in the reverse of the order given above.), so there are no more flow shifts on the bush.
\index{static traffic assignment!algorithms!Algorithm B|)}

\stevefig{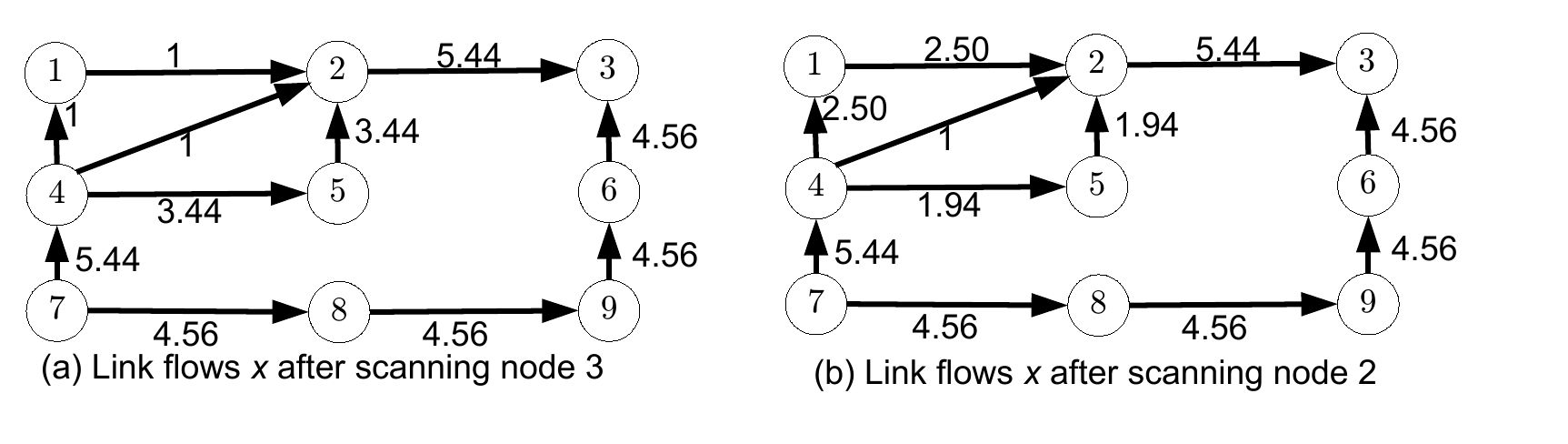}{Algorithm B flow shifts.}{\textwidth}

\subsubsection{Origin-based assignment (OBA)}

\index{static traffic assignment!algorithms!origin-based assignment|(}
Rather than considering just two paths at a time, as Algorithm B did, OBA shifts flow among many paths simultaneously.
This can be both a blessing (in that its moves affect more paths simultaneously) and a curse (in that its moves are not as sharp as Algorithm B, which can focus on the two paths with the greatest travel time difference).
Because OBA shifts flow among many paths, it makes use of the average cost labels $M$ and the derivative labels $D$, rather than the shortest and longest path costs (and the direct link travel times derivatives) used by Algorithm B.
OBA also makes use of a Newton-type shift, dividing a difference in travel times by an approximation of this difference's derivative.

When scanning a node $i$ in OBA, first identify a least-travel time approach, that is, a link $(\hat{h},i)$ such that $M_{\hat{h} i} \leq M_{hi}$ for all other approaches $(h,i)$.
This link is called the \emph{basic approach}\index{basic approach} in analogy to the basic path concept used in path-based algorithms.
, for each nonbasic approach $(h,i)$, the following amount of flow is shifted from $x_{hi}$ to $x_{\hat{h} i}$:
\labeleqn{obashift}{\Delta x_{hi} = \frac{M_{hi} - M_{\hat{h} i}}{D_{hi} + D_{\hat{h}i} - 2D_a} \,,}
with the constraint that $\Delta x_{hi} \leq x_{hi}$ to prevent negative flows.
It can be shown that this formula would equalize the mean travel times on the two approaches if they were linear functions.
In reality, they are not, but the formula is still used as an approximation.

Here $a$ is the divergence node, defined for OBA as the node with the highest topological order which is common to \emph{all} paths in the bush from $r$ to $i$, excluding $i$ itself.
This is a reasonable place to truncate the search, since shifting flow among path segments between $a$ and $i$ will not affect the flows on any earlier links in the bush.
This is a generalization of the definition used for Algorithm B, needed since there are more than two paths subject to the flow shift.

After applying the shift $\Delta x_{hi}$ to each of the links entering node $i$, we have to update flows on the other links between $a$ and $i$ to maintain flow conservation.
This is done by assuming that the $\alpha$ values stay the same elsewhere, meaning that any increase or decrease in the flow passing through a node propagates backward to its incoming links in proportion to the contribution each incoming link provides to that total flow.
Thus, the $x$ labels can be recalculated using equations~\eqn{alphaxone} and~\eqn{alphaxtwo} to links and nodes topologically between $a$ and $i$.

The steps of OBA are as follows:
\begin{enumerate}
\item Calculate the $M$ and $D$ labels in forward topological order.
\item Let $i$ be the topologically last node in the bush.
\item Scan $i$ by performing the following steps:
\begin{enumerate}
    \item Determine the divergence node $a$ corresponding to node $i$, and the basic approach $(\hat{h}, i)$.
    \item For each nonbasic approach $(h,i)$ calculate $\Delta x_{hi}$ using equation~\eqn{obashift} (with $\Delta x_{hi} \leq x_{hi}$), subtract $\Delta x_{hi}$ from $x_{hi}$ and add it to $x_{\hat{h} i}$.
    \item Update $\alpha_{hi}$ for every bush link terminating at $i$.
    \item Apply equations~\eqn{alphaxone} and~\eqn{alphaxtwo} to all nodes and links topologically between $a$ and $i$, updating their $x_{ij}$ and $x_i$ values.
\end{enumerate}
\item If $i = r$, go to the next step.
Otherwise, let $i$ be the previous node topologically and return to step 3.
\item Update all travel times $t_{ij}$ and derivatives $t'_{ij}$ using the new flows $x$ (remembering to add flows from other bushes.)
\end{enumerate}

Again demonstrating on the example in Figure~\ref{fig:bushgrid}, we use the $M$ and $D$ labels shown in Figure~\ref{fig:bushgrid2}.
As with Algorithm B, we scan nodes in reverse topological order.
Starting with node 3, we compare the labels $M_{23}$ and $M_{63}$ to determine the basic approach.
Since $M_{63}$ is lower, this is the least-cost approach and named basic.
For OBA, the divergence node corresponding to node 3 is node 7, since it is the only node common to all paths between nodes 7 and 3, except for node 3 itself.
Formula~\eqn{obashift} is then applied to determine the amount of flow $\Delta x_{23}$ that should be shifted from approach $(2,3)$ to the basic approach $(6,3)$, producing $\Delta x_{23} = 1.48$.
Shifting this flow updates the $\alpha$ values as shown in Figure~\ref{fig:obastep}(b).
To maintain flow conservation, we also need to adjust the flow on other links in the bush.
Holding the $\alpha$ proportions fixed at all other links, these changes at node 3 are propagated back proportional to the original flows, producing the link flows in Figure~\ref{fig:obastep}(a).

Next scanning node 2, we examine $M_{12}$, $M_{42}$, and $M_{52}$ to determine the basic approach: since $M_{12}$ is smallest, it is deemed the basic approach.
The divergence node $a$ corresponding to node 2 is node 4, since all paths from 7 to 2 pass through node 4 (and no other node with higher topological order).
We now apply two shifts, calculating $\Delta x_{42}$ and $\Delta x_{52}$ with equation~\eqn{obashift}.
This equation gives $\Delta x_{42} = 2.5$, but since $x_{42} = 0.79$ it is only possible to move 0.79 units of flow from (4,2) to (1,2).
For approach (5,2), $\Delta x_{52} = 1.50$, so we move 1.50 units of flow from (5,2) to (1,2).
Assuming that the $\alpha$ values at all other nodes remains constant, propagating this change back produces the $x$ and $\alpha$ labels in Figure~\ref{fig:obastep}(c) and (d).
Notice that we only had to update flows between the divergence node ($a = 4$) and the node being scanned $(i=2)$, since no other links would be affected by shifting flow among approaches to node 2.

As with Algorithm B, no changes are made for the remaining nodes scanned, since they only have one incoming link.
This incoming link is trivially the ``basic'' approach, and there are no nonbasic approaches to shift flow from.
\index{static traffic assignment!algorithms!origin-based assignment|)}

\stevefig{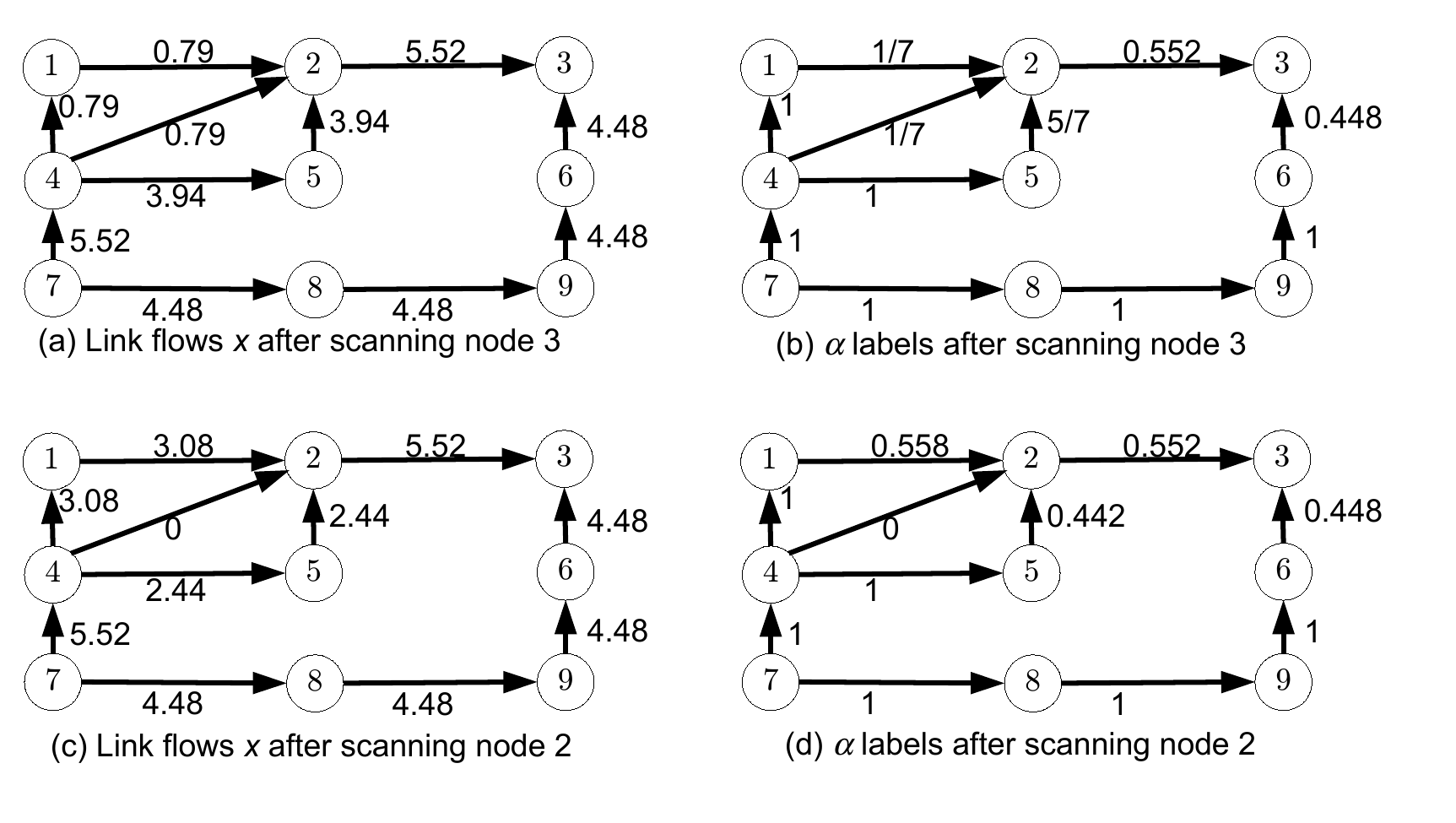}{OBA flow shifts.}{\textwidth}

\subsubsection{Linear user cost equilibrium (LUCE)}

\index{static traffic assignment!algorithms!linear user cost equilibrium (LUCE)|(}
The LUCE algorithm is similar in structure to OBA, using the same definition of a divergence node and the same sets of labels ($M$ and $D$).
Where it differs is in how nodes are scanned.
When LUCE scans a node, it attempts to solve a ``local'' user equilibrium problem, based on linear approximations to the approach travel times.
To be concrete, when node $i$ is scanned, the labels $M_{hi}$ and $D_{hi}$ respectively indicate the average travel time for travelers arriving at $i$ via $(h,i)$, and the derivative of this average travel time as additional vehicles arrive via this link.
Both of these are calculated based on the current bush flows, so if the flow on approach $(h,i)$ changes by $\Delta x_{hi}$, then the new average travel time $M'_{hi}$ can be approximated by
\labeleqn{luceapprox}{M'_{hi} \approx M_{hi} + \Delta x_{hi} D_{hi} \,.}
The LUCE algorithm chooses $\Delta x_{hi}$ values for each approach according to the following principles: 
\begin{enumerate}[(a)]
    \item Flow conservation must be obeyed, that is, $\sum_{(h,i) \in \mc{B}} \Delta x_{hi} = 0$.
    \item No link flow can be made negative, that is, $x_{hi} + \Delta x_{hi} \geq 0$ for all $(h,i) \in \mc{B}$.
    \item The $M'_{hi}$ values should be equal and minimal for any approach with positive flow after the shift, that is, if $x_{hi} + \Delta x_{hi} > 0$, then $M'_{hi}$ must be less than or equal to the $M'$ label for any other approach to $i$.
\end{enumerate}
The $\Delta x_{hi}$ values satisfying this principles can be found using a ``trial and error'' algorithm, like that introduced in Section~\ref{sec:solvingforequilibrium}.
Choose a set of approaches (call it $\mc{A}$), and set $\Delta x_{hi} = -x_{hi}$ for all $(h,i)$ not in this set $\mc{A}$.
For the remaining approaches, solve the linear system of equations which set $M'_{hi}$ equal for all $(h,i) \in \mc{A}$ and which have $\sum_{(h,i) \in \mc{A}} \Delta x_{hi} = 0$.
The number of equations will equal the number of approaches in $\mc{A}$.
After obtaining such a solution, you can verify whether the three principles are satisfied.
Principle (a) will always be satisfied.
If principle (b) is violated, approaches with negative $x_{hi} + \Delta x_{hi}$ should be removed from $\mc{A}$.
If principle (c) is violated, then some approach not in $\mc{A}$ has a lower $M'$ value, and that approach should be added to $\mc{A}$.
In either of the latter two cases, the entire process should be repeated with the new $\mc{A}$ set.

The steps of LUCE are as follows:
\begin{enumerate}
\item Calculate the $M$ and $D$ labels in forward topological order.
\item Let $i$ be the topologically last node in the bush.
\item Scan $i$ by performing the following steps:
\begin{enumerate}
    \item Determine the divergence node $a$ corresponding to node $i$, and the basic approach $(\hat{h}, i)$.
    \item Use the process described above to find $\Delta x_{hi}$ values satisfying the local equilibrium principles (a)--(c) above, and add $\Delta x_{hi}$ to $x_{hi}$ for each approach to $i$.
    \item Update $\alpha_{hi}$ for every bush link terminating at $i$.
    \item Use equations~\eqn{alphaxone} and ~\eqn{alphaxtwo} to all nodes and links topologically between $a$ and $i$, updating their $x_{ij}$ and $x_i$ values.
\end{enumerate}
\item If $i = r$, go to the next step.
Otherwise, let $i$ be the previous node topologically and return to step 3.
\item Update all travel times $t_{ij}$ and derivatives $t'_{ij}$ using the new flows $x$ (remembering to add flows from other bushes.)
\end{enumerate}
Applying LUCE to the same example, we again start by scanning node 3.
Assuming that both approaches (2,3) and (6,3) will continue to be used, we solve the following equations simultaneously for $\Delta x_{23}$ and $\Delta x_{63}$:
\begin{align*}
M_{23} + D_{23} \Delta x_{23} &= M_{63} + D_{63} \Delta x_{63} \\
\Delta x_{23} + \Delta x_{63} &= 0
\end{align*}
Substituting the $M$ and $D$ labels from Figure~\ref{fig:bushgrid2} and solving, we obtain $\Delta x_{23} = -1.48$ and $\Delta x_{63} = +1.48$.
Updating flows as in OBA (using the divergence node $a = 7$ and assuming all $\alpha$ values at nodes other than 3 are fixed) gives the $x$ and $\alpha$ labels shown in Figure~\ref{fig:lucestep}(a) and (b).
Notice that this step is exactly the same as the first step taken by OBA.
The interpretation of LUCE and solving a local, linearized equilibrium problem provides insight into how the OBA flow shift formula~\eqn{obashift} was derived.

The shift is slightly different when there are three approaches to a node, as happens when we proceed to scan node 2.
First assuming that all three approaches will be used, we solve these three equations simultaneously, enforcing flow conservation and that the travel times on the three approaches should be the same:
\begin{align*}
M_{12} + D_{12} \Delta x_{12} &= M_{42} + D_{42} \Delta x_{42} \\
M_{52} + D_{52} \Delta x_{52} &= M_{42} + D_{42} \Delta x_{42} \\
\Delta x_{12} + \Delta x_{42} + \Delta x_{52} &= 0
\end{align*}
Substituting values from Figure~\ref{fig:bushgrid2} and solving these equations simultaneously gives $\Delta x_{12} = 2.82$, $\Delta x_{42} = -1.85$, and $\Delta x_{52} = -0.97$.
This is problematic, since $x_{42} = 0.79$ and it is impossible to reduce its flow further by 1.85.
This means that approach (4,2) should not be used, so we fix $\Delta x_{42} = -0.79$ as a constant and re-solve the system of equations for $\Delta x_{12}$ and $\Delta x_{52}$:
\begin{align*}
M_{12} + D_{12} \Delta x_{12} &= M_{42} - 0.79 D_{42} \\
M_{52} + D_{52} \Delta x_{52} &= M_{42} - 0.79 D_{42} \\
\Delta x_{12} -0.79 + \Delta x_{52} &= 0
\end{align*}
This produces $\Delta x_{12} = 2.06$, and $\Delta x_{52} = -1.27$, alongside the fixed value $\Delta x_{42} = -0.79$.
Updating flows on other links as in OBA produces the $x$ and $\alpha$ labels shown in Figure~\ref{fig:lucestep}(b) and (c).
No other shifts occur at lower topologically-ordered nodes, because there is only one incoming link, and flow conservation demands that no flow increase or decrease take place.
\index{static traffic assignment!algorithms!linear user cost equilibrium (LUCE)|)}

\stevefig{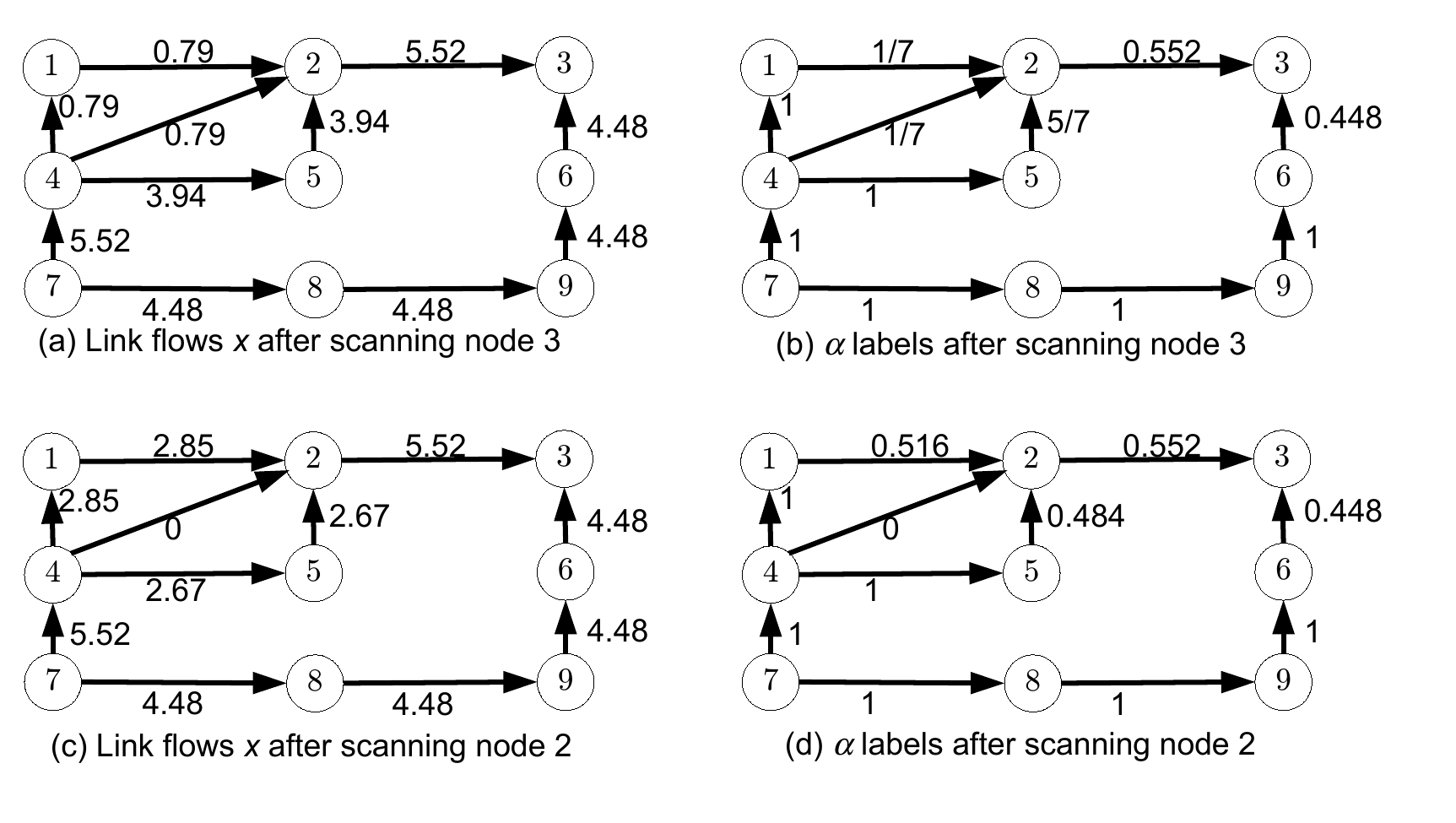}{LUCE flow shifts.}{\textwidth}

\subsection{Improving bushes}
\label{sec:bushadjust}

\index{bush!changing links|(}
After shifting flow for the current bushes as described in the previous section, the next step is to determine whether the bushes themselves need to change (adding or dropping links) to allow us to move closer to equilibrium on the entire network.

To begin, any link which is unused can be dropped from a bush without disturbing the solution (unless it is needed for connectivity); and it can always be added back later if we need to.
As far as which links to add, one approach is to add any ``shortcut'' links, defined based on the bush travel times.
If we solve each bush precisely to equilibrium in Step 2, then the distance from the origin $r$ to each node $i$ is the same on any used path; call this $L_i$.
A ``shortcut'' link $(i,j)$ is any link for which $L_i + t_{ij} < L_j$, where the $L$ labels have been re-calculated after eliminating the zero-flow bush links.
Such a link can justifiably be called a shortcut, because it provides a faster path to node $j$ than currently exists in the bush.

If this is our rule for adding links to the bush, it is easy to see that no cycles are created.
Reasoning by contradiction, assume that a cycle \emph{is} created, say, $[i, j, k, \ldots, i]$.
If link travel times are strictly positive, it is easy to see that $L_j > L_i$, $L_k > L_j$, and so forth for any successive pair of nodes in the cycle: if $(i,j)$ was in the previous bush, then $L_j = L_i + t_{ij}$ by equilibrium; of $(i,j)$ was just added, in which case $L_i + t_{ij} < L_j$.
Applying this identically cyclically, we have $L_i < L_j < L_k < \ldots < L_i$, a contradiction since we cannot have $L_{i_1} < L_{i_1}$.

The requirement that we solve each bush exactly to equilibrium in step 2 is vital to ensuring no cycles are created.
If this is not the case (and in practice, we can only solve equilibrium approximately), another criterion is needed.
To see why, the argument used to show that no cycles would be created relied critically on the assumption the difference in minimum cost labels for adjacent nodes in the cycle was exactly the travel time of the connecting link (except possibly for a link just added), which is only possible if all used paths have equal and minimal travel time.
Another one, which is almost as easy to implement, is to add links to the bush based on the \emph{maximum} travel time to a node using bush links, as opposed to the \emph{minimum} travel time.
Although calculating longest paths in general networks is difficult, in acyclic networks it can be done just as efficiently as finding shortest paths (Section~\ref{sec:acyclicsp}).

Assume that we have re-calculated the longest path labels $U_i$ to each bush node, after eliminating unused bush links as described above.
We now define a ``modified shortcut'' link as any link for which $U_i + t_{ij} < U_j$ (this is the same definition as a shortcut link but with maximum costs labels used instead of minimum cost labels).
Using a similar argument, we can show that adding modified shortcut links cannot create cycles, even if the bush is not at equilibrium.
Arguing again by contradiction, assume that the cycle $[i, j, k, \ldots, i]$ is created by adding modified shortcut links, and consider in turn each adjacent pair of $U$ labels.
By the definition of maximum cost, if the link $(i,j)$ was in the bush before, we have $U_i + t_{ij} \leq U_j$; and if the link $(i,j)$ was just added $U_i + t_{ij} < U_j$.
Applying this identity cyclically, we have $U_{i} \leq U_{j} \leq U_{k} \leq \ldots \leq U_{i}$.
Furthermore, since there were no cycles in the bush during the previous iteration, at least one of these links must be new, and for this link the inequality must be strict.

As an example, consider the bush updates which occur after performing the LUCE example in the previous section.
Figure~\ref{fig:luceupdate} shows the updated link travel times in panel (a), and the remaining bush links after zero-flow links are removed in panel (b).
Re-calculating $L$ and $U$ labels with the new travel times and bush topology gives the values in Figure~\ref{fig:luceupdate}(c).
At this point the three unused bush links (4,2), (5,6), and (8,5), are examined to determine whether $U_i + t_{ij} < U_j$ for any of them.
This is true for (5,6) and (8,5), since $41.6 + 2 < 66.3$ and $22.1 + 2 < 41.6$, but false for (4,2) since $32.5 + 42 \geq 72.9$.
So, (5,6) and (8,5) are added to the bush, as shown in Figure~\ref{fig:luceupdate}(d).
From here, one can return to the flow shifting algorithm to update flows further.
\index{bush!changing links|)}
\index{static traffic assignment!algorithms!bush-based|)}

\stevefig{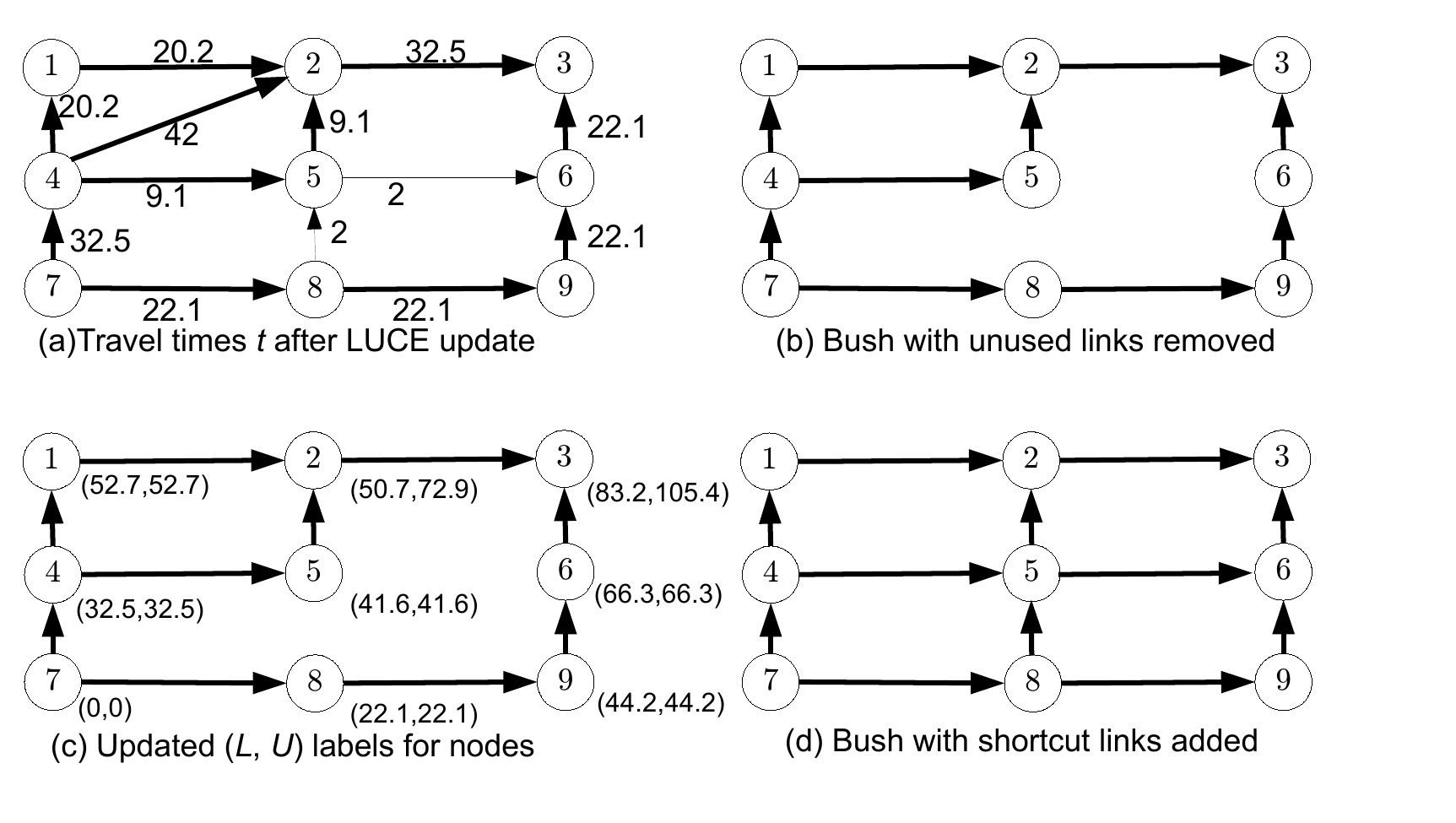}{Updating a bush by removing unused links and adding shortcuts.}{\textwidth}

\section{Likely Path Flow Algorithms}
\label{sec:likelypathflowalgorithms}

\index{static traffic assignment!properties!nonuniqueness in path flows|(}
\index{entropy}
Recall from Section~\ref{sec:maxentropy} that there is generally an infinite number of path flow vectors $\mb{h}$ which satisfy the principle of user equilibrium.
This contrasts with the fact that the equilibrium link flow vector $\mb{\hat{x}}$ is unique as long as link performance functions are increasing.
That section introduced the principle of entropy maximization as a way to select a equilibrium path flow believed to be the most likely to occur in practice.
The entropy-maximizing link flows solve the optimization problem
\begin{align}
\max_\mathbf{h} \qquad & -\sum_{(r,s) \in Z^2} \sum_{\pi \in \hat{\Pi}^{rs}} h_\pi \log (h_\pi / d^{rs}) & \label{eqn:entropyalgo} \\
\mathrm{s.t.} \qquad   & \sum_{\pi \in \Pi} \delta_{ij}^\pi h_\pi = \hat{x}_{ij}           & \forall (i,j) \in A \label{eqn:entropyeqmalgo} \\
                       & \sum_{\pi \in \hat{\Pi}^{rs}} h_\pi = d^{rs}   & \forall (r,s) \in Z^2 \label{eqn:entropynvlbalgo} \\
                       & h_\pi \geq 0                                   & \forall \pi \in \Pi \label{eqn:entropynonnegalgo} 
\end{align}

That section also introduced the proportionality condition, showing in Theorem~\ref{thm:proportionality} that the ratio of the flows on any two paths connecting the same OD pair depends \emph{only} on the pairs of alternate segments where the paths differ, not what that OD pair happens to be, or on any links where the paths coincide.
The proportionality condition is easier to verify, and entropy-maximizing solutions must satisfy proportionality.
Proportionality does not imply entropy maximization, but in practical terms they seem nearly equivalent.
Therefore, this section focuses on finding proportional solutions, hence the term ``likely path flow algorithms'' rather than ``\emph{most} likely.''
\index{static traffic assignment!proportionality}

The proof that entropy maximization implies proportionality on page~\pageref{sec:proportionalityproof} was based on the Lagrangian\index{Lagrange multipliers!applications} of the entropy-maximization problem:
\begin{multline}
\label{eqn:entropylagrangianalgo}
\mathcal{L}(\mathbf{h}, \bm{\beta}, \bm{\gamma}) = -\sum_{(r,s) \in Z^2} \sum_{\pi \in \hat{\Pi}^{rs}} h_\pi \log \myp{\frac{h_\pi}{d^{rs}}} + \sum_{(i,j) \in A} \beta_{ij} \myp{\hat{x}_{ij} - \sum_{\pi \in \Pi} \delta_{ij}^\pi h_\pi} \\ + \sum_{(r,s) \in Z^2} \gamma_{rs} \myp{d^{rs} - \sum_{\pi \in \hat{\Pi}^{rs}} h_\pi}
\,.\end{multline}
After some algebraic manipulation, we obtained the formula
\labeleqn{mlpfprimaldual}{
h^\pi = K_{rs} \exp \myp{ -\sum_{(i,j) \in A} \delta_{ij}^\pi \beta_{ij}}
\,,}
where $K_{rs}$ is a proportionality constant associated with the OD pair $(r,s)$ corresponding to path $\pi$.
The value of $K_{rs}$ can be found from the constraint that total path flows must equal demand, $\sum_{\pi \in \hat{\Pi}_{rs}} h^\pi = d^{rs}$.
The rest of the equation shows that the entropy-maximizing path flows are determined solely by the Lagrange multipliers $\beta_{ij}$ associated with each link, and therefore that the ratio between two path flows for the same OD pair only depends on the links where they differ.

Here we describe three ways to calculate likely path flows.
The first method is ``primal,'' and operates directly on the path flow vector $\mb{h}$ itself.
The second method is ``dual,'' operating on the Lagrange multipliers $\bm{\beta}$ in the entropy maximization problem, which can then be used to determine the path flows.
Both of these methods presume that an equilibrium link flow solution $\mathbf{\hat{x}}$ has already been found.
The third method, ``traffic assignment by paired alternate segments'' (TAPAS), is an algorithm which simultaneously solves for equilibrium and likely paths.

Primal and dual methods for likely path flows have two major issues in common:
\begin{enumerate}
\item The sets of equilibrium paths $\hat{\Pi}_{rs}$ must be determined in some way.
If $\mathbf{\hat{x}}$ is exactly a user equilibrium, we can identify the shortest path for each OD pair, and then set $\hat{\Pi}_{rs}$\index{path!set of used paths} to be all the paths with that same cost between $r$ and $s$.
In practice, though, we cannot solve for the equilibrium link flows exactly, but only to a finite precision.
So we cannot guarantee that all of the used paths for an OD pair have the same travel time, and we cannot guarantee that the unused paths are actually unused at the true equilibrium solution.
\item The system of constraints~\eqn{entropyeqmalgo} and~\eqn{entropynvlbalgo} has many more variables than equations; the number of equilibrium paths is typically much larger than the number of OD pairs and network links.
For primal methods, we need to find and specify these ``degrees of freedom'' so we know how to adjust path flows while keeping the solution feasible (matching equilibrium link flows and the OD matrix).
For dual methods, it means that many of the link flow constraints~\eqn{entropyeqmalgo} are redundant, and their corresponding $\beta_{ij}$ Lagrange multipliers are not needed.
If these redundant constraints and Lagrange multipliers are included, then there are infinitely many $\bm{\beta}$ values that correspond to entropy maximization.
\end{enumerate}

To address the first issue, we can include all paths in $\hat{\Pi}$ that are within some threshold $\epsilon$ of the shortest-path travel time for each OD pair.
Some care must be taken to ensure that the resulting path sets are ``consistent,'' in that travelers from different origins and destinations consider the same sets of alternatives, when their path sets overlap.\index{path!set of used!consistency}\index{consistent set of used paths|see {path, set of used, consistency}}
Figure~\ref{fig:consistency} illustrates the difficulties in obtaining a consistent solution.
The top panel of the figure shows the travel times at the true equilibrium solution, and the bottom  panel shows the travel times at an approximate equilibrium solution, such as the one obtained after stopping one of the algorithms in this chapter after a finite number of iterations.

\begin{figure}
\centering
\begin{tikzpicture}[->,>=stealth',shorten >=1pt,auto,node distance=2cm, thick,main node/.style={circle,draw}]
\node[main node] (1) {1};
\node[main node] (2) [above of=1] {2};
\node[main node] (3) [right of=1] {3};
\node[main node] (4) [above of=3] {4};
\node[main node] (5) [right of=3] {5};
\node[main node] (6) [above of=5] {6};
\node[main node] (7) [right of=5] {7};
\node[main node] (8) [above of=7] {8};

\path
(1) edge node[left] {50} (2)
    edge node[above] {5} (3)
(3) edge node[above] {5} (5)
(4) edge node[above] {5} (2)
(5) edge node[left] {30} (6)
    edge node[above] {10} (7)
(6) edge node[above] {5} (4)
(7) edge node[left] {10} (8)
(8) edge node[above] {10} (6);
\end{tikzpicture} \\
(a) Equilibrium link travel times \\
\begin{tikzpicture}[->,>=stealth',shorten >=1pt,auto,node distance=2cm, thick,main node/.style={circle,draw}]
\node[main node] (1) {1};
\node[main node] (2) [above of=1] {2};
\node[main node] (3) [right of=1] {3};
\node[main node] (4) [above of=3] {4};
\node[main node] (5) [right of=3] {5};
\node[main node] (6) [above of=5] {6};
\node[main node] (7) [right of=5] {7};
\node[main node] (8) [above of=7] {8};

\path
(1) edge node[left] {49} (2)
    edge node[above] {5} (3)
(3) edge node[above] {4} (5)
(4) edge node[above] {6} (2)
(5) edge node[left] {31} (6)
    edge node[above] {11} (7)
(6) edge node[above] {5} (4)
(7) edge node[left] {10} (8)
(8) edge node[above] {9} (6);
\end{tikzpicture} \\
(b) Link travel times at an approximate equilibrium solution
\caption{Link travel times at the true equilibrium solution (top) and at an approximate equilibrium solution (bottom).
\label{fig:consistency}}
\end{figure}

In this figure, assume there are two OD pairs, one from node 1 to node 2, and another from node 3 to node 4.
At the true equilibrium solution, there are three equal-cost paths between nodes 1 and 2 --- $[1,2]$, $[1,3,5,6,4,2]$, and $[1,3,5,7,8,6,4,2]$ --- and two equal-cost paths between nodes 3 and 4 --- $[3,5,6,4]$ and $[3,5,7,8,6,4]$.
Travelers whose paths cross nodes 5 and 6 consider two possible alternative routes between these nodes ($[3,5,6,4]$ and $[3,7,8,6,4]$), regardless of which OD pair they came from.
So defining $\hat{\Pi}_{12}$ and $\hat{\Pi}_{34}$ to include the equal-cost paths listed above is consistent.

Now consider the bottom panel of the figure.
There is now a unique shortest path for each OD pair --- $[1,2]$ and $[3,5,6,4]$ --- but this does not reflect the true equilibrium, simply the imprecision in an approximate solution.
The threshold $\epsilon$ can be set to reflect this.
Assume that $\epsilon = 1$, so that any path within 1 minute of the shortest path travel time is included in $\hat{\Pi}$.
This choice gives the path sets $\hat{\Pi}_{12} = \myc{ [1,2], [1,3,5,6,4,2]}$ and $\hat{\Pi}_{34} = \myc{ [3,5,6,4], [3,5,7,8,6,4] }$.
This choice is \emph{not} consistent, because it assumes travelers passing between nodes 5 and 6 consider different choices depending on their OD pair: travelers starting at node 1 only consider the segment $[5,6]$, while travelers starting at node 3 consider both $[5,6]$ and $[5,7,8,6]$ as options.
If travelers are choosing routes to minimize cost, it should not matter what their origin or destination is.
Increasing $\epsilon$ to a larger value would address this problem, but in a larger network would run the risk of including paths which are not used at the true equilibrium solution.

So some care must be taken in how $\epsilon$ is chosen.
In the approximate equilibrium solution, we can define the \emph{acceptance gap}\index{acceptance gap} $g_a$\label{not:ga} to be the greatest difference between a used path's travel time and the shortest path travel time for its OD pair, and the \emph{rejection gap}\index{rejection gap} $g_r$ to be the smallest difference between the travel time of an unused path, and the shortest path for an OD pair.
Figure~\ref{fig:rejectiongap}, the OD pairs are from 1 to 6 and 4 to 9, and the links are labeled with their travel times.
The thick lines show the links used by these OD pairs, and the thin lines show unused links.
For travelers between nodes 1 and 6, the used paths have travel times of 21 (shortest) and 22 minutes, and the unused path has a travel time of 26 minutes.
For travelers between nodes 4 and 9, the used paths have travel times 18 (which is shortest) and 20 minutes, while the unused path has a travel time of 21 minutes.
The acceptance gap $g_a$ for this solution is 2 minutes (difference between 20 and 18), and the rejection gap is 5 (difference between 26 and 21).

\begin{figure}
\centering
\begin{tikzpicture}[->,>=stealth',shorten >=1pt,auto,node distance=2cm, thick,main node/.style={circle,draw}]
\node[main node] (1) {1};
\node[main node] (2) [right of=1] {2};
\node[main node] (3) [right of=2] {3};
\node[main node] (4) [above of=1] {4};
\node[main node] (5) [right of=4] {5};
\node[main node] (6) [right of=5] {6};
\node[main node] (7) [above of=4] {7};
\node[main node] (8) [right of=7] {8};
\node[main node] (9) [right of=8] {9};

\path
(2) edge node[above] {10} (3)
(3) edge node[right] {10} (6)
(4) edge node[right] {7} (7)
(7) edge node[above] {10} (8);

\path[line width=1mm]
(1) edge node[above] {8} (2)
    edge node[right] {7} (4)
(2) edge node[right] {9} (5)
(4) edge node[above] {9} (5)
(5) edge node[above] {5} (6)
    edge node[right] {5} (8)
(6) edge node[right] {6} (9)
(8) edge node[above] {4} (9);
\end{tikzpicture}
\caption{Example demonstrating rejection and acceptance gaps.
Thick links are parts of used paths at the equilibrium solution.
Link labels are travel times.
\label{fig:rejectiongap}}
\end{figure}

It is possible to show that if $\epsilon$ is at least equal to the acceptance gap, but less than half of the rejection gap, then the resulting sets of paths $\hat{\Pi}_{rs}$ are consistent for the proportionality condition.
That is, we need
\labeleqn{twoconsistency}{g_a \leq \epsilon \leq \frac{g_r}{2}\,.}
In the network of Figure~\ref{fig:rejectiongap}, any $\epsilon$ choice between 2 and 2.5 will thus lead to a consistent set of paths.

Exercise~\ref{ex:twoconsistency} gives a more formal definition of consistency and asks you to prove this statement.
It is not possible to choose such an $\epsilon$ unless the equilibrium problem is solved with enough precision for the acceptance gap to be less than half of the rejection gap.
To fully maximize entropy, rather than just satisfying proportionality, demands a more stringent level of precision in the equilibrium solution.

The second issue involves redundancies in the set of equations enforcing the OD matrix and equilibrium link flow constraints.
Properly resolving this issue requires using linear algebra to analyze the structure of the set of equations (and this in fact is the key to bridging the gap between proportionality and entropy maximization), but for proportionality a simpler approach is possible.

\begin{figure}
\centering
\begin{tikzpicture}[->,>=stealth',shorten >=1pt,auto,node distance=3cm, thick,main node/.style={circle,draw}]
\node[main node] (1) {E};
\node[main node] (A) [above left of=1] {A};
\node[main node] (C) [below left of=1] {C};
\node[main node] (2) [right of=1] {F};
\node[main node] (3) [right of=2] {G};
\node[main node] (B) [above right of=3] {B};
\node[main node] (D) [below right of=3] {D};

\path
(A) edge node[above right] {15} node[below left] {5} (1)
(C) edge node[above left] {45} node[below right] {6} (1)
(1) edge [bend left] node[above] {40} node[below] {1} (2)
(1) edge [bend right] node[above] {20} node[below] {2} (2)
(2) edge [bend left] node[above] {30} node[below] {3} (3)
(2) edge [bend right] node[above] {30} node[below] {4} (3)
(3) edge node[above left] {15} node[below right] {7} (B)
(3) edge node[above right] {45} node[below left] {8} (D);
\end{tikzpicture}
\begin{tabular}{cc|cc}
\multicolumn{2}{c|}{OD pair $(A,B)$} & \multicolumn{2}{c}{OD pair $(C,D)$} \\
Path ID & Links & Path ID & Links \\
\hline
1 & 5, 1, 3, 7  & 5       & 6, 1, 3, 8 \\
2 & 5, 1, 4, 7  & 6       & 6, 1, 4, 8 \\
3 & 5, 2, 3, 7  & 7       & 6, 2, 3, 8 \\
4 & 5, 2, 4, 7  & 8       & 6, 2, 4, 8
\end{tabular} 
\caption{Example demonstrating redundancies and ``degrees of freedom'' when choosing a path flow solution.
Links are labeled with equilibrium flows (above) and link IDs (below).
\label{fig:dof}}
\end{figure}

A redundancy in a system of equations can be interpreted as a ``degree of freedom,'' a dimension along which a solution can be adjusted.
Consider the network in Figure~\ref{fig:dof}, which has two OD pairs (A to B, and C to D).
The equilibrium link flows are shown in the figure, along with the link IDs and an indexing of the eight paths.
The OD matrix and link flow constraints are reflected in the six equations
\begin{align}
h_1 + h_2 + h_3 + h_4 &= 15 \label{eqn:propexab} \\
h_5 + h_6 + h_7 + h_8 &= 45 \label{eqn:propexcd} \\
h_1 + h_2 + h_5 + h_6 &= 40 \label{eqn:propex1} \\
h_3 + h_4 + h_7 + h_8 &= 20 \label{eqn:propex2} \\
h_1 + h_3 + h_5 + h_7 &= 30 \label{eqn:propex3} \\
h_2 + h_4 + h_6 + h_8 &= 30 \label{eqn:propex4}
\end{align}
Equations~\eqn{propexab} and~\eqn{propexcd} reflect the constraints that the total demand among all paths from A to B, and from C to D, must equal the respective values in the OD matrix.
Equations~\eqn{propex1}--\eqn{propex4} reflect the constraints that the flow on links 1--4 must match their equilibrium values.
Similar equations for links 5--8 are omitted, since they are identical to~\eqn{propexab} and~\eqn{propexcd}, as can easily be verified.

This system has eight variables but only six equations, so there must be at least two independent variables --- in fact, there are four, since some of the six equations are redundant.
For example, adding equations~\eqn{propexab} and~\eqn{propexcd} gives you the same result as adding equations~\eqn{propex1} and~\eqn{propex2}, so one of them --- say,~\eqn{propex2} can be eliminated.
Likewise, equation~\eqn{propex4} can be eliminated, since adding~\eqn{propexab} and~\eqn{propexcd} is the same as adding~\eqn{propex3} and~\eqn{propex4}.
These choices are not unique, and there are other equivalent ways of expressing the same redundancies.

Each of these redundancies corresponds to an independent way to adjust the path flows without affecting either total path flows between OD pairs, or total link flows.
A primal method uses these redundancies to adjust the path flows directly, increasing the entropy of the solution without sacrificing feasibility of the original path flow solution.
A dual method uses these redundancies to eliminate unnecessary link flow constraints --- for instance, equations~\eqn{propex2} and~\eqn{propex4} in the example above --- so that the entropy-maximizing $\beta_{ij}$ values are unique.

\subsection{Primal method}
\label{sec:mlpfprimal}

\index{static traffic assignment!algorithms!primal method for proportionality|(}
To find a proportional solution with a primal method, we need (1) the equilibrium path set $\hat{\Pi}$; (2) the equilibrium link flows $\mathbf{\hat{x}}$; (3) an initial path flow solution $\mathbf{h^0}$, and (4) a list of ``redundancies'' in the link flow constraints, each of which corresponds to a way to change path flows while maintaining feasibility.
The algorithm then adjusts the $\mathbf{h}$ values, increasing the entropy at each iteration, until termination.\footnote{There are alternative solution representations that can make this algorithm much faster, but would make the explanations more complicated.
You are encouraged to think about how to implement this algorithm efficiently, without having to list all paths explicitly.}

The equilibrium path set and link flows were already described above.
Depending on the algorithm used to solve for the equilibrium link flows, a path flow solution may already be available.
If you solved for equilibrium with a path-based algorithm, its solution already contains the flows on each path.
If you used a bush-based algorithm, a corresponding path flow can be found using the ``within origins'' adjustment technique described below.
If you used a link-based algorithm, it is not as easy to directly identify a path flow solution from the final output --- but it may be possible to track a path flow solution as the algorithm progresses (for instance, in Frank-Wolfe, each ``target'' all-or-nothing solution can be clearly identified with a path flow solution, and this can be averaged with previous path flow solutions using the same $\lambda$ values).

To achieve proportionality, it is not necessary to identify \emph{all} of the redundancies in the link flow constraints (it would be needed to fully maximize entropy).
It is enough to identify \emph{pairs of alternate segments}\index{pair of alternate segments} between two nodes which are used by multiple OD pairs.
In Figure~\ref{fig:dof}, links 1 and 2 are alternate segments between nodes E and F, and links 3 and 4 are alternate segments between nodes F and G.

The primal algorithm alternates between two steps: (1) adjusting $\mathbf{h}$ to achieve proportionality for each origin $r$; and (2) for each pair of alternate segments, adjusting $\mathbf{h}$ to achieve proportionality between those nodes.
These steps are described below; you can verify that each one of these steps preserves feasibility of the solution, and increases entropy.
If the solution remains unchanged after performing both of these steps, then proportionality has been achieved and we terminate.
Other stopping criteria can be introduced, by measuring the deviation from proportionality and terminating once this deviation is sufficiently small.

\subsubsection{Proportionality within origins}

Section~\ref{sec:bushes} described how a user equilibrium solution can be described by the total flow from each origin $r$ on each link $(i,j)$, denoted $x^r_{ij}$, and how at equilibrium the links with positive $x^r_{ij}$ values must form an acyclic subnetwork, a bush.\index{static traffic assignment!bush-based solution}
It is easy to obtain such a solution if we have a path flow solution $\mathbf{h}$:
\labeleqn{originbasedfrompath}{
x^r_{ij} = \sum_{s \in Z} \sum_{\pi \in \hat{\Pi}_{rs}} \delta_{ij}^\pi h^\pi
\,.}
With these values, we can calculate the fraction of the flow from origin $r$ approaching any node $j$ from one specific link entering that node $(i,j)$:
\labeleqn{originbasedproportion}{
\alpha^r_{ij} = \left.
x^r_{ij} \middle/ \sum_{(h,j) \in \Gamma^{-1}(j)} x_{ij}^r \right.
\,}
with $\alpha^r_{ij}$ defined arbitrarily if the denominator is zero.
(This use of $\alpha$ is the same as in the bush-based algorithms described in Section~\ref{sec:bushbased}).

We can ensure that the proportionality condition holds between \emph{all} paths associated with an origin $r$ by updating the path flows according to the formula
\labeleqn{withinoriginpathupdate}{
h^\pi \leftarrow d^{rs} \prod_{(i,j) \in \pi} \alpha_{ij}^r \qquad \forall s \in Z, \pi \in \hat{\Pi}_{rs}
\,,}
that is, by applying the \emph{aggregate} approach proportions across all paths from this origin to each individual path.
One can show updating the path flows with this formula will not change either the total OD flows or link flows, maintaining feasibility, and will also increase entropy if there is any change in $\mathbf{h}$.

This process is repeated for each origin $r$.

\subsubsection{Proportionality between origins}

To achieve proportionality between different origins, we consider pairs of alternate segments between two nodes, used by multiple OD pairs.
In Figure~\ref{fig:dof}, an example of a pair of alternate segments is links 1 and 2.
In a larger network, these alternate segments can contain multiple links.
These links connect the same two nodes (E and F), and are parts of equal-cost paths for both OD pairs A-B and C-D.
If we move some flow from link 1 to link 2 from OD pair A-B, and move the equivalent amount of flow from link 2 to link 1 from OD pair C-D, we will not disturb the total link flows or OD flows, but will change $\mathbf{h}$ and allow us to increase entropy.

Let $\zeta = \{ \sigma^\zeta_1, \sigma^\zeta_2 \}$\label{not:zetapas}\label{not:sigmapas} be a pair of alternate segments, and let $Z^2(\zeta)$\label{not:Z2zeta} be the set of OD pairs which have paths using both of the alternate segments in $\zeta$.
For each segment $\sigma_i$, the segment flow $g(\sigma_i)$ is defined as the sum of flows on all paths using that segment:
\labeleqn{segmentflowformula}{
g(\sigma^\zeta_i) = \sum_{(r,s) \in Z^2(\zeta)} \sum_{\pi \in \hat{\Pi}_{rs} : \sigma^\zeta_i \subseteq \pi} h^\pi
\,,}
where the notation $\sigma_i \subseteq \pi$ means that all of the links in the segment $\sigma_i$ are in the path $\pi$.

If the equilibrium path set $\hat{\Pi}$ is consistent in the sense of~\eqn{twoconsistency}, then any path $\pi$ which uses one segment of the pair has a ``companion'' path which is identical, except it uses the other segment of the pair, denoted $\pi^c(\sigma)$\label{not:pic}.
For example, in Figure~\ref{fig:dof}, if $\sigma$ is the pair of alternate segments between nodes E and F, the companion of path 1 is path 3, and the companion of path 6 is path 8.

To achieve proportionality between origins for the alternate segments in $\sigma$, we calculate the ratios between $g(\sigma_i)$ values and apply the same ratios to the path flows for each OD pair using this set of alternate segments:\label{not:superseteq}
\labeleqn{betweenoriginpathupdate}{h^\pi \leftarrow (h^\pi + h^{\pi^c(\sigma)}) \frac{g(\sigma^\zeta_i)}{g(\sigma^\zeta_1) + g(\sigma^\zeta_2)} \quad \forall (r,s) \in Z^2(\zeta), i \in \{ 1, 2 \}, \pi \in \hat{\Pi}_{rs}, \pi \supseteq \sigma^\zeta_i\,.}
It is again possible to show that updating path flows with this formula leaves total OD flows and link flows fixed, and can only increase entropy.

\subsubsection{Example}

This section shows how the primal algorithm can solve the example in Figure~\ref{fig:dof}.
Assume that the initial path flow solution is $h_1 = 15$, $h_5 = 15$, $h_6 = 10$, $h_8 = 20$, and all other path flows zero.
As shown in the first column of Table~\ref{tbl:primalalgo}, this solution satisfies the OD matrix and the resulting link flows match the equilibrium link flows, so it is feasible.
The entropy of this solution, calculated using~\eqn{entropyalgo}, is 47.7.\footnote{When computing this formula, $0 \log 0$ is taken to be zero, since $\lim_{x \rightarrow 0^+} x \log x = 0$.}  Two pairs of alternate segments are identified: links 1 and 2 between nodes E and F, and links 3 and 4 between nodes F and G.

This table summarizes the progress of the algorithm in successive columns; you may find it helpful to refer to this table when reading this section.
The bottom section of the table shows the origin-based link flows corresponding to the path flow solution, calculated using~\eqn{originbasedfrompath}.

The first iteration applies the within-origin formula~\eqn{withinoriginpathupdate} to origin A, and to origin C.
To apply the formula to origin A, the origin-based proportions are first calculated with equation~\eqn{originbasedproportion}: $\alpha^A_1 = 1$, $\alpha^A_2 = 0$, $\alpha^A_3 = 1$, and $\alpha^A_4 = 0$, and thus $h_1 \leftarrow 15$, $h_2 \leftarrow 0$, $h_3 \leftarrow 0$, and $h_4 \leftarrow 0$.
(There is no change.)  For origin C, we have $\alpha^C_1 = 5/9$, $\alpha^C_2 = 4/9$, $\alpha^C_3 = 1/3$, and $\alpha^C_4 = 2/3$, and thus $h_5 \leftarrow 8 \frac{1}{3}$, $h_6 \leftarrow 16 \frac{2}{3}$, $h_7 \leftarrow 6\frac{2}{3}$, and $h_8 \leftarrow 13\frac{1}{3}$.
The entropy of this new solution has increased to 59.6.

Next, we apply the between-origin formula to the pair of alternate segments between nodes F and G.
Paths 1, 3, 5, and 7 use link 3; and their companion paths (2, 4, 6, and 8, respectively) use link 4.
The segment flow for link 3 is the sum of the path flows that use it (30), and similarly the segment flow for link 4 is also 30.
Thus formula~\eqn{betweenoriginpathupdate} requires paths 1 and 2 to have equal flow, paths 3 and 4 to have equal flow, and so on for each path and its companion.
Redistributing the flow between paths and companions in this way gives the result in the third column of Table~\ref{tbl:primalalgo}, and the entropy has increased to 72.5.

The between-origin formula is applied a second time to the other pair of alternate segments, between nodes E and F.
Paths 1, 2, 5, and 6 use link 1, and the companion paths using link 2 are 3, 4, 7, and 8, respectively.
The segment flows for links 1 and 2 are 40 and 20, so for each path and its companion, the path using link 1 should have twice the flow of the path using link 2.
The fourth column of Table~\ref{tbl:primalalgo} shows the results, and the entropy has increased again to 79.8.

This solution achieves proportionality (and in fact maximizes entropy).
The proportionality conditions can either be checked directly, or noticed when running the algorithm a second time does not change any path flows.
For larger networks with a more complicated structure, the algorithm generally requires multiple iterations, and only converges to proportionality in the limit.
\index{static traffic assignment!algorithms!primal method for proportionality|)}

\begin{table}
\centering
\caption{Demonstration of primal proportionality algorithm.
\label{tbl:primalalgo}}
\begin{tabular}{|cc|cccc|}
\hline
                 &         & Initial          & Within-origin & F/G PAS & E/F PAS \\
\hline 
Entropy          &         & 47.7             & 59.6          & 72.5    & 79.8 \\
\hline
Path flows       & $h_1$   & 15               & 15            & 7.5     & 5 \\
                 & $h_2$   & 0                & 0             & 7.5     & 5 \\
                 & $h_3$   & 0                & 0             & 0       & 2.5 \\
                 & $h_4$   & 0                & 0             & 0       & 2.5 \\
                 & $h_5$   & 15               & 8.33          & 12.5    & 15 \\
                 & $h_6$   & 10               & 16.67         & 12.5    & 15 \\
                 & $h_7$   & 0                & 6.67          & 10      & 7.5 \\
                 & $h_8$   & 20               & 13.33         & 10      & 7.5 \\
\hline 
OD flows         & $(A,B)$ & 15               & 15            & 15      & 15 \\
                 & $(C,D)$ & 45               & 45            & 45      & 45 \\
\hline 
Link flows       & $x_1$   & 40               & 40            & 40      & 40 \\
(total)          & $x_2$   & 20               & 20            & 20      & 20 \\
                 & $x_3$   & 30               & 30            & 30      & 30 \\
                 & $x_4$   & 30               & 30            & 30      & 30 \\
\hline 
Link flows       & $x_1^A$ & 15               & 15            & 15      & 10 \\
(origin A)       & $x_2^A$ & 0                & 0             & 0       & 5 \\
                 & $x_3^A$ & 15               & 15            & 7.5     & 7.5\\
                 & $x_4^A$ & 0                & 0             & 7.5     & 7.5 \\
\hline 
Link flows       & $x_1^B$ & 25               & 25            & 25      & 30 \\
(origin B)       & $x_2^B$ & 20               & 20            & 20      & 15 \\
                 & $x_3^B$ & 15               & 15            & 22.5    & 22.5 \\
                 & $x_4^B$ & 30               & 30            & 22.5    & 22.5 \\
\hline
\end{tabular}
\end{table}

\subsection{Dual method}

\index{static traffic assignment!algorithms!dual method for proportionality|(}
An alternative approach involves the optimality conditions directly.
Recall from equation~\eqn{mlpfprimaldual} that entropy-maximizing (and thus proportional) flow on each path is $K_{rs} \exp \myp{ -\sum_{(i,j) \in A} \delta_{ij}^\pi \beta_{ij}}$, where $K_{rs}$ is an OD-specific constant chosen so that the path flows sum to the total demand $d_{rs}$.
We can adopt this condition as a \emph{formula} for $\mb{h}$, and no matter what values are chosen for $\beta_{ij}$, the path flows we calculate satisfy proportionality.
The difficulty is that they will generally not be feasible, unless the link flows $\mb{x}$ corresponding to $\mb{h}$ happen to equal their equilibrium values $\mb{\hat{x}}$.
A \emph{dual} algorithm tries to adjust these $\beta_{ij}$ values until $\mb{x} = \mb{\hat{x}}$, at which point we terminate with proportional (and in fact entropy-maximizing) path flows.

This contrasts with the \emph{primal} approach in the previous section, which always maintained feasibility (the link flows in Table~\ref{tbl:primalalgo} never changed from the equilibrium values) and worked toward optimality, expressed in equation~\eqn{mlpfprimaldual}.
The algorithm in this section always maintains optimality, and works toward feasibility.
Unlike the primal algorithm, the path flows in the dual algorithm are not feasible until termination.

The idea behind the algorithm is simple enough: start with initial values for $\beta_{ij}$ on each link; calculate $\mb{h}$ from equation~\eqn{mlpfprimaldual}; calculate $\mb{x}$ from $\mb{h}$; and see which links have too much flow or too little flow.
Adjust the $\beta_{ij}$ values accordingly, and iterate until the flow on every link is approximately equal to its equilibrium value.
There are a few details to take care of: how to find an initial solution, when to terminate, how to adjust the $\beta_{ij}$ values, and dealing with redundancies in the system of constraints.
The first two details are fairly simple: any initial solution will do; $\bm{\beta} = \mb{0}$ is simplest.
Terminate when $\mb{x}$ is ``close enough'' to $\mb{\hat{x}}$ according to some measure.

For adjusting the $\beta_{ij}$ values, notice from equation~\eqn{mlpfprimaldual} that increasing $\beta_{ij}$ will decrease $x_{ij}$, and vice versa.
So a natural update rule is
\labeleqn{dualsearch}{
\beta_{ij} \leftarrow \beta_{ij} + \alpha(x_{ij} - \hat{x}_{ij})
\,.}
where $\alpha$ is a step size, and $x_{ij} - \hat{x}_{ij}$ is the difference between the link flows currently implied by $\bm{\beta}$, and the equilibrium values.
In addition to this intuitive interpretation, this search direction is also proportional to the gradient of the least-squares function\label{not:phi}
\labeleqn{dualleastsquares}{
\phi(\mb{x}) = \sum_{(i,j) \in A} (x_{ij} - \hat{x}_{ij})^2
\,.}
This function is zero only at a feasible solution, and~\eqn{dualsearch} is a steepest descent direction in terms of $\mb{x}$.\footnote{To be precise, it is not a steepest descent direction in terms of $\bm{\beta}$.
An alternative derivation of equation~\eqn{dualsearch} is explored in Exercise~\ref{ex:dualdirectionderivation}.}  

Equation~\eqn{dualleastsquares} can also be used to set the step size $\alpha$.
One can select a trial sequence of $\alpha$ values (say, $1, 1/2, 1/4, 1/8, \ldots$), evaluating each $\alpha$ value in turn and stopping once the new $\mb{x}$ values reduce~\eqn{dualleastsquares}.
A more sophisticated step size rule chooses $\alpha$ using Newton's method, to approximately maximize entropy.
Newton's method also has the advantage of scaling the step size based on the effect changes in $\beta$ have on link flows.
Exercise~\ref{ex:dualstepsize} develops this approach in more detail.

A last technical detail concerns redundancies in the system of link flow equations, as discussed at the end of Section~\ref{sec:likelypathflowalgorithms}.
Redundancies in the link flow equations mean that the $\beta_{ij}$ values maximizing entropy may not be unique.
To resolve this issue, redundant link flow constraints can be removed, and their $\beta_{ij}$ values left fixed at zero.
Practical experience shows that this can significantly speed convergence.

\subsubsection{Example}

The dual method is now demonstrated on the same example as the primal algorithm; see Figure~\ref{fig:dof}.
You may find it helpful to refer to Table~\ref{tbl:dualalgo} when reading this section to track the progress of the algorithm.
The format is similar to Table~\ref{tbl:primalalgo}, except for additional rows showing the $\beta_{ij}$ values.

This table summarizes the progress of the algorithm in successive columns; you may find it helpful to refer to this table when reading this section.
The bottom section of the table shows the origin-based link flows corresponding to the path flow solution, calculated using~\eqn{originbasedfrompath}.

To begin, as discussed at the end of Section~\ref{sec:likelypathflowalgorithms}, two of the link flow constraints are redundant, and their $\beta_{ij}$ values are fixed at zero.
Assume that links 2 and 4 are chosen for this purpose, so $\beta_2 = \beta_4 = 0$ throughout the algorithm.
(The algorithm would perform similarly for other choices of the two redundant links; note that we are also continuing to ignore the link flow constraints associated with links 5--8, since these are identical to the OD matrix constraints~\eqn{propexab} and~\eqn{propexcd}.)

Initially, $\beta_1 = \beta_3 = 0$.
This means that $\exp \myp{ -\sum_{(i,j) \in A} \delta_{ij}^\pi \beta_{ij}} = 1$ for all paths, so $h_1 = h_2 = h_3 = h_4 = K_{AB}$ and $h_5 = h_6 = h_7 = h_8 = K_{CD}$.
To ensure that the sum of each OD pairs' path flows equals the total demand, we need $K_{AB} = 3.75$ and $K_{CD} = 11.25$, and equation~\eqn{mlpfprimaldual} gives the flows on each path, as shown in the Iteration 0 column of Table~\ref{tbl:dualalgo}.

The table also shows the link flows corresponding to this solution: links 1--4 all have 30 vehicles, whereas the equilibrium solution has $x_1 = 40$ and $x_2 = 20$.
The least-squares function~\eqn{dualleastsquares} has the value $(40 - 30)^2 + (20 - 30)^2 = 200$.
Trying an initial step size of $\alpha = 1$ would give $\beta_1 = 0 + 1 \times (30 - 40) = -10$.
The other $\beta_{ij}$ values are unchanged: $\beta_3$ remains at zero because it has the correct link flow, while $\beta_2$ and $\beta_4$ are permanently fixed at zero because their link flow constraints were redundant.
Re-applying equation~\eqn{mlpfprimaldual} with this new value of $\beta_1$ (and recalculating $K_{AB}$ and $K_{CD}$ to satisfy the OD matrix) would give $x_1 = 60$, $x_2 = 0$, and $x_3 = x_4 = 30$.
This has a larger least-squares function than before (800 vs.\ 200), so we try again with $\alpha = 1/2$.
This is slightly better (the least-squares function is 768), but still worse than the current solution.
After two more trials, we reach $\alpha = 1/8$, which produces a lower mismatch (88).

This step size is accepted, and we proceed to the next iteration.
The path and link flows are shown in the Iteration 1 column of Table~\ref{tbl:dualalgo}.
The flows on links 1 and 2 are closer to their equilibrium values than before.
Continuing as before, we find that $\alpha = 1/8$ is again the acceptable step size with the new link flow values, so $\beta_1 \leftarrow -1.25 + \frac{1}{8}(36.2 - 40) = -0.42$, producing the values in the Iteration 2 column.
Over additional iterations, the algorithm converges to the final values shown in the rightmost column.

It is instructive to compare the dual algorithm in Table~\ref{tbl:dualalgo} with the primal algorithm in Table~\ref{tbl:primalalgo}.
Notice how the dual algorithm always maintains proportionality, and the link flows gradually converge to their equilibrium algorithms.
By contrast, the primal algorithm maintains the link flows at their equilibrium values, and gradually converges to proportionality.
The entropy also does not change monotonically, and at times it is higher than the maximum entropy value --- this can only happen for an infeasible solution.
\index{static traffic assignment!algorithms!dual method for proportionality|)}

\begin{table}
\centering
\caption{Demonstration of dual proportionality algorithm.
\label{tbl:dualalgo}}
\begin{tabular}{|cc|ccccc|}
\hline
                 &           & Iteration 0 & 1       & 2       & $\cdots$ & $\infty$ \\
\hline
Entropy          &           & 83.2        & 73.4    &81.9     &          & 79.8   \\
\hline
                 & $\beta_1$ & 0           & $-1.25$ & $-0.42$ &          & $-0.692$ \\
                 & $\beta_3$ & 0           & 0       & 0       &          & 0 \\
\hline
Path flows       &$h_1$      & 3.75        & 5.83    & 4.53    &          & 5 \\
                 &$h_2$      & 3.75        & 5.83    & 4.53    &          & 5 \\
                 &$h_3$      & 3.75        & 1.67    & 2.97    &          & 2.5 \\
                 &$h_4$      & 3.75        & 1.67    & 2.97    &          & 2.5 \\
                 &$h_5$      & 11.25       & 17.49   & 13.58   &          & 15 \\
                 &$h_6$      & 11.25       & 17.49   & 13.58   &          & 15 \\
                 &$h_7$      & 5.01        & 8.92    & 6.53    &          & 7.5 \\
                 &$h_8$      & 5.01        & 8.92    & 6.53    &          & 7.5 \\
\hline   
OD flows         &$(A,B)$    & 15          &15       & 15      &          & 15   \\
                 &$(C,D)$    & 45          &45       & 45      &          & 45   \\
\hline   
Total link flows &$x_1$      & 30          & 46.6    & 36.2    &          & 40   \\
                 &$x_2$      & 30          & 13.4    & 23.8    &          & 20   \\
                 &$x_3$      & 30          &30       & 30      &          & 30   \\
                 &$x_4$      & 30          &30       & 30      &          & 30   \\
\hline   
From origin A    &$x_1^A$    & 7.5         & 11.7    & 9.1     &          & 10 \\
                 &$x_2^A$    & 7.5         & 3.3     & 5.9     &          & 5 \\
                 &$x_3^A$    & 7.5         & 7.5     & 7.5     &          & 7.5\\
                 &$x_4^A$    & 7.5         & 7.5     & 7.5     &          & 7.5   \\
\hline   
From origin B    &$x_1^B$    & 22.5        & 35.0    & 27.1    &          & 30 \\
                 &$x_2^B$    & 22.5        & 10.0    & 17.8    &          & 15 \\
                 &$x_3^B$    & 22.5        &22.5     & 22.5    &          & 22.5   \\
                 &$x_4^B$    & 22.5        &22.5     & 22.5    &          & 22.5   \\
\hline
\end{tabular}
\end{table}

\subsection{Traffic assignment by paired alternate segments}

\index{static traffic assignment!algorithms!traffic assignment by paired alternate segments (TAPAS)|(}
\index{pair of alternate segments|(}
The primal and dual methods described in the preceding sections assumed that user equilibrium link flows were already available, and then found likely path flows as a post-processing step.
Traffic assignment by paired alternate segments (TAPAS) is an algorithm which finds the equilibrium solution and proportional path flows simultaneously.
Interestingly, accomplishing both tasks at once does not slow down the algorithm.
TAPAS is in fact among the fastest of the traffic assignment algorithms currently known.

Recall from Section~\ref{sec:bushbased} that path- and bush-based algorithms find the equilibrium solution by shifting flow from longer paths to shorter ones, and that these paths often differ on a relatively small set of links (in gradient projection, we denoted these by the set $A_3 \cup A_4$; in bush-based algorithms, by the concept of a divergence node).
The main insights of TAPAS are that these algorithms tend to shift flow repeatedly between the same sets of links, and that these links are common to paths used by between different origins and destinations.
As a result, it makes sense to store these path segments from one iteration to the next, rather than having to expend effort finding them again and again.
Furthermore, since these links are common to multiple origins, we can apply proportionality concepts at the same time to find a high-entropy path flow solution.

\stevefig{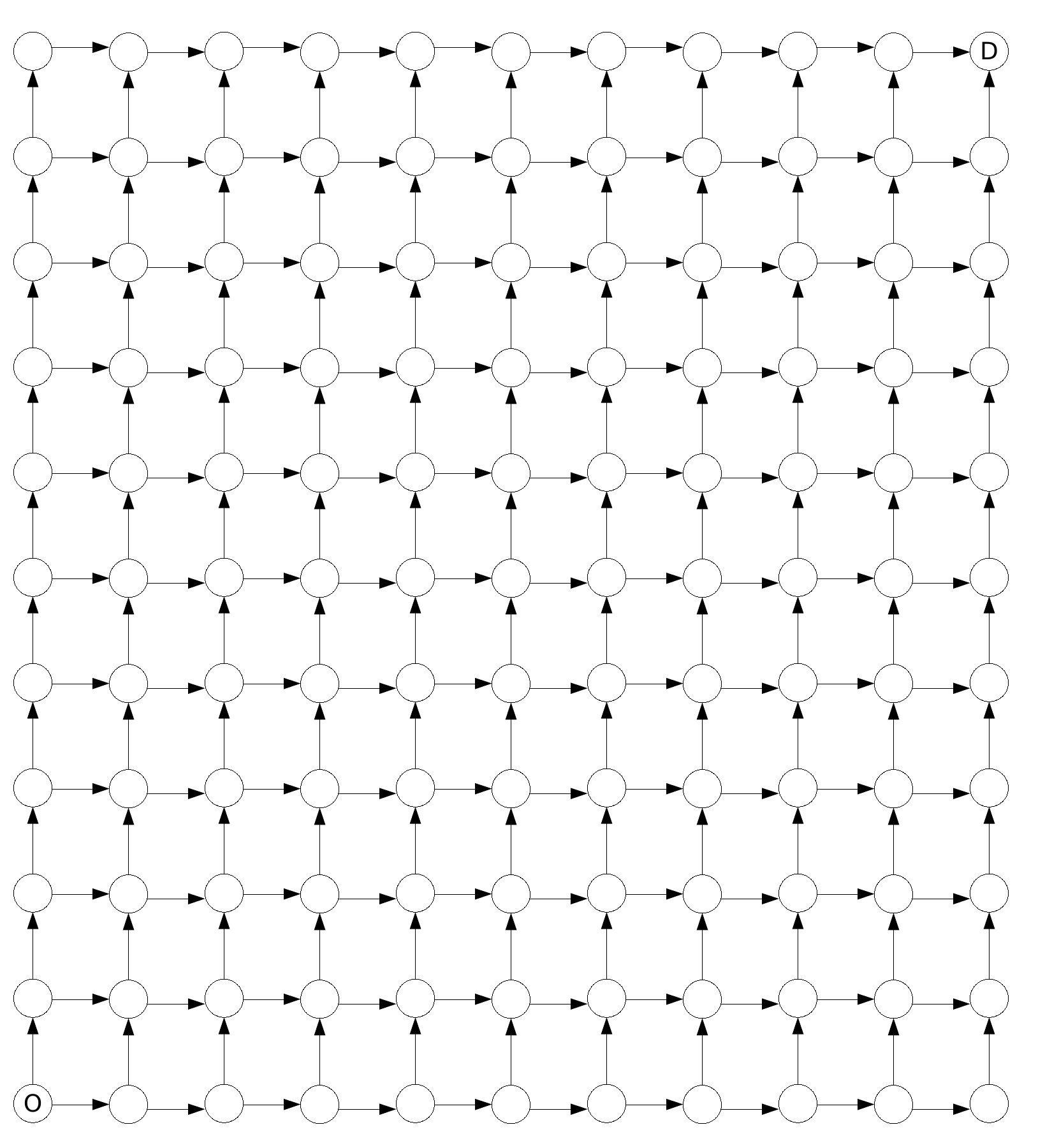}{A grid network consisting of 100 city blocks.}{0.8\textwidth}

Pairs of alternate segments\index{pair of alternate segments} can also form a concise representation of the equilibrium conditions.
In the grid network of Figure~\ref{fig:tapasintuition}, the number of paths between the origin in the lower-left and the destination in the upper-right is rather large (in fact there are $184{,}756$) even though the network is a relatively modestly-sized grid of ten rows and columns.
If all paths are used at equilibrium, expressing the equilibrium condition by requiring the travel times on all paths be equal requires $184{,}755$ equations.

The network can also be seen as a hundred ``city blocks,'' each of which can either be traversed in the clockwise direction (north, then east) or in the counterclockwise direction (east, then north).
If travel times on all paths in the network are identical, then the travel time around each block must be the same in the clockwise and counterclockwise directions.
In fact, the converse is true as well: if the travel time is the same around every block in both orientations, then the travel times on all paths in the network are the same as well.
We can thus express equality of all network paths with only 100 equations!

Furthermore, other nodes in the network may serve as origins and destinations, not just nodes at two corners.
Expressing equilibrium in terms of path travel time equality requires additional equations for each new OD pair, but the same 100 equations expressing equality of travel times around each block are sufficient no matter how many nodes serve as origins or destinations.

This discussion implies that most of the $184{,}755$ path travel-time equations are redundant.
It is not trivial to identify these linear dependencies, but methods based on paired alternate segments are a way to do so.
The two ways to travel around each block can be seen as a pair of alternate segments between their southwestern and northeastern nodes, and the set of all such pairs of alternate segments is ``spanning'' in the sense that one can shift flow between any two paths with the same origin and destination by shifting flow between a sequence of pairs of alternate segments.
There are some subtleties involving the equivalence of equilibrium conditions on pairs of alternate segments and on paths (see Exercise~\ref{ex:paspathsubtleties}),
but this example shows how they can often simplify the search for a user equilibrium solution.

The TAPAS algorithm represents network flows aggregated by origin, with $x_{ij}^r$ denoting the flow on link $(i,j)$ which started at origin $r$, following equation~\eqn{linkdisaggregate}.
The algorithm also involves a set of pairs of alternate segments (PASs).
Each PAS $\zeta$ is defined by two path segments $\sigma_1^\zeta$ and $\sigma_2^\zeta$ starting and ending at the same nodes, and by a list of \emph{relevant origins} $Z_\zeta$\label{not:Zzeta} indicating a subset of zones which have flow on both path segments at the current solution.
The steps  of the algorithm involve maintaining a set of PASs (creating new ones, updating relevant origins, and optionally discarding inactive ones), and adjusting the link flows $x_{ij}^r$ to move towards user equilibrium and proportionality.
Notice that TAPAS does not store the path flows $\mb{h}$ explicitly, for computational reasons.
Instead, the path flows are represented \emph{implicitly}, obtained from the $x_{ij}^r$ values using the proportional split formulas previously introduced in Section~\ref{sec:mlpfprimal}:
\labeleqn{tapasimplicit}{
\alpha^r_{ij} = \left.x^r_{ij} \middle/ \sum_{(h,j) \in \Gamma^{-1}(j)} x_{ij}^r \right.
\,}
where
\labeleqn{tapasimplicitalpha}{
x^r_{ij} = \sum_{s \in Z} \sum_{\pi \in \hat{\Pi}_{rs}} \delta_{ij}^\pi h^\pi
\,.}

An example of a PAS is shown in Figure~\ref{fig:initialtapasexample}.
For each link, the upper and lower labels give the flows on that link from Origin 1 and Origin 2, respectively.
There are two path segments: $[2,4]$ and $[2,3,4]$, and there is one relevant origin (Origin 2).
You might expect that this PAS is also relevant to Origin 1; and indeed at the ultimate equilibrium solution this will be true.
However, in large networks it is not immediately obvious which PASs are relevant to which origins, and the TAPAS algorithm must discover this during its steps.

\stevefig{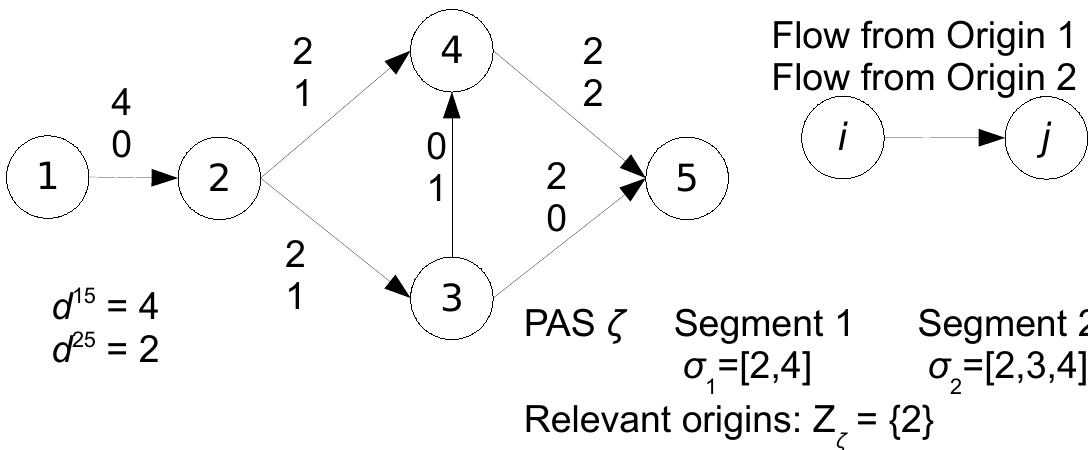}{Example for demonstrating pairs of alternate segments.}{\textwidth}

There are three main components to the algorithm: PAS management, flow shifts, and proportionality adjustments.
PAS management involves identifying new PAS, updating the lists of relevant origins, and removing inactive ones.
Flow shifts move the solution closer to user equilibrium, by shifting vehicles from longer paths to shorter ones.
Proportionality adjustments maintain the total link flows at their current values, but adjust the origin-specific link flows to increase the entropy of the path flow solution implied by~\eqn{tapasimplicit}.
One possible way to perform these steps is as follows; the rest of this subsection fills out the details of each step.

\begin{enumerate}
\item Find an initial origin-disaggregated solution, and initialize the set of PASs to be empty.
\item Update the set of PASs by determining whether new ones should be created, or whether existing ones are relevant to more origins.
\item Perform flow shifts within existing PASs.
\item Perform proportionality adjustments within existing PASs.
\item Check for convergence, and return to step 2 unless done.
\end{enumerate}

This algorithmic description may appear vague.
Like many of the fastest algorithms currently available, the performance of the algorithm depends on successfully balancing these three components of the algorithm.
The right amount of time to spend on each component is network- and problem-specific, and implementations that make such decisions adaptively, based on the progress of the algorithm, can work well.

\subsubsection{Updating the set of pairs of alternate segments}

In its second step, the TAPAS algorithm must update the set of PASs.
Since flow shifts mainly occur within PASs, finding the user equilibrium solution relies on being able to identify new alternative routes which are shorter than the ones currently being used.
Given the origin-disaggregated solution $\mb{x}^r$ for each origin $r$, we can search for routes in the following way.

Solving a shortest path algorithm over the entire network produces a tree rooted at an origin $r$, containing paths to every node.
(Note that already having an origin-disaggregated solution can greatly accelerate the process of finding shortest paths; see Exercise~\ref{ex:fasterspwithorigin} from Chapter~\ref{chp:networkrepresentations}.)  At equilibrium, essentially all of these links should be used.\footnote{The only exceptions would be to links leading to nodes not used by this origin, or if multiple paths are tied for being shortest with one of them having zero flow.}
By comparing this tree to the links which are used in the disaggregate solution $\mb{x}^r$, we can identify any links in the shortest path tree not currently being used by an origin $r$ even though there is flow to their head nodes.

More specifically, let $\hat{A}^r$\label{not:Ahatr} denote the set of these links.
A link $(i,j)$ is in $\hat{A}^r$ if $(i,j)$ is in the shortest path tree rooted at $r$, if $x_{ij}^r = 0$ (it is currently unused by origin $r$), and if there is some other link $(h,j)$ for which $x_{hj}^r$ (flow from origin $r$ is reaching the head node $j$ in another way).
We then look for a PAS $\zeta$ whose two segments $\sigma_1^\zeta$ and $\sigma_2^\zeta$ end with the links $(h,j)$ and $(i,j)$, which will allow us to shift flow onto the shortest path segment.

Two possibilities exist: either such a PAS already exists (in which case we add origin $r$ to the relevant set $Z_\zeta$ if it is not already listed), or we create a new one.
To create a new PAS, we must choose two path segments $\sigma_1^\zeta$ and $\sigma_2^\zeta$ which start and end at the same node.
The first segment $\sigma_1^\zeta$ should consist of links with positive flow from origin $r$, and the second one $\sigma_2^\zeta$ should consist of links from the shortest path tree.
We also know they must both end at node $j$, but must choose an appropriate node for them to start (a divergence node).
As discussed in Section~\ref{sec:bushshift}, there are several ways to choose divergence nodes.
For TAPAS, an ideal divergence node results in short path segments $\sigma_1^\zeta$ and $\sigma_2^\zeta$.
Short segments both result in faster computation, and intuitively are more likely to be relevant to more origins.

Putting these concepts together, we can search for a divergence node $a$ for which there is a segment $\sigma_1^\zeta$ starting at $a$, ending with link $(h,j)$, and only using links with positive flow from origin $r$; and a segment $\sigma_2^\zeta$ starting at $a$, ending with link $(i,j)$, and only using links in the current shortest path tree.
Among all such divergence nodes and segments, we want one for which the two path segments are short.

As an example, consider again the network from Figure~\ref{fig:initialtapasexample}.
The link performance functions and current travel times are shown in Figure~\ref{fig:tapasnewpas}.
The bold links show the shortest path tree rooted at Origin 1, and see that there are two links used by this origin which are not part of this tree: links $(2,4)$ and $(3,5)$.
For link $(2,4)$, we see that the two segments of the PAS $a$ from the previous example include a segment of links used by this origin $[2,4]$, and a segment of links from the shortest path tree $[2,3,4]$ that have a common divergence node 2.
At this point, we declare Origin 1 relevant to this PAS by adding it to the set $Z_a$.

\stevefig{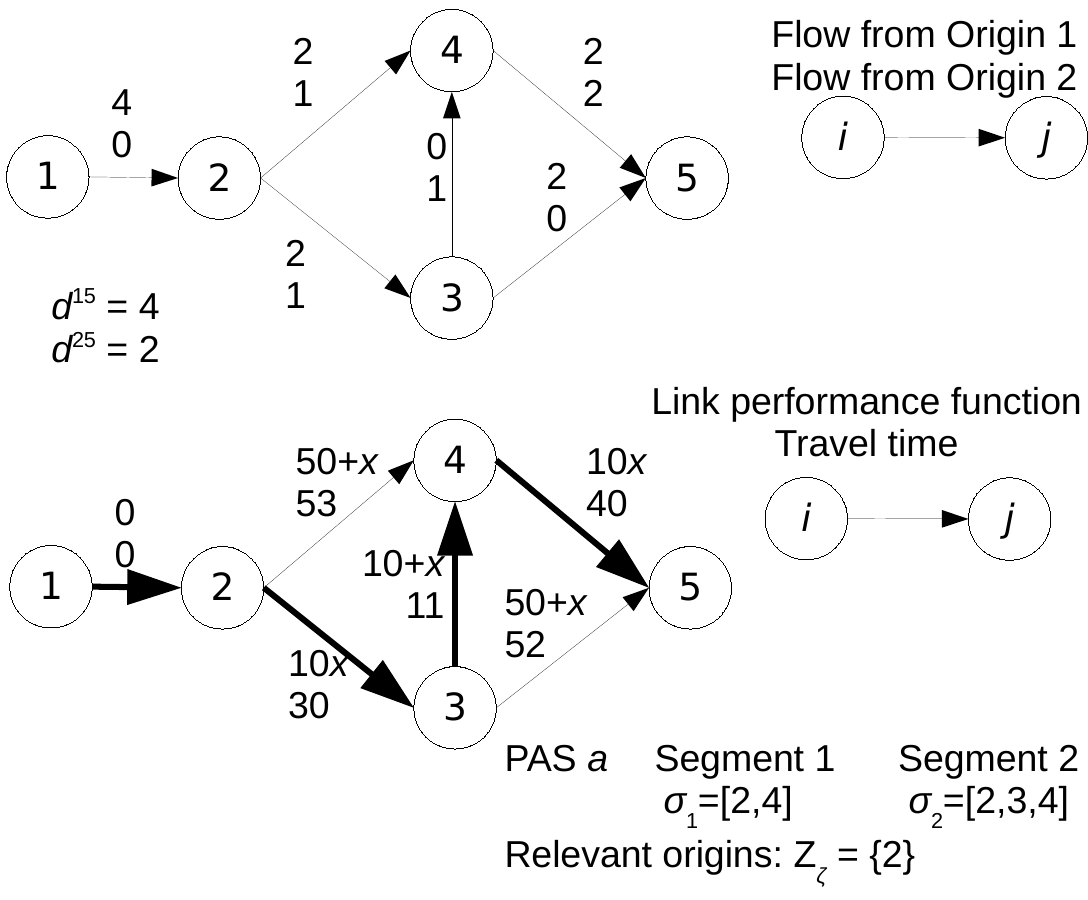}{Creating a new PAS. The top panel shows the origin-specific flows, the bottom panel the link performance functions and current times. Bold links are the shortest path tree for Origin 1.}{0.8\textwidth}

For link $(3,5)$, we need to create a new PAS.
There are two possibilities for choosing segments: one choice is $[3,5]$ and $[3,4,5]$; and the other choice is $[2,3,5]$ and $[2,3,4,5]$ (both involve a segment of used links and a segment from the shortest path tree, starting at a common divergence node).
The first one is preferred, because it has fewer links --- and in fact the common link $(2,3)$ in the second PAS is irrelevant, since shifting flow between segments will not change flow on such a link at all).
By being shorter, there are potentially more relevant origins.
If node 3 were also an origin, it could be relevant to the first choice of segments, but not the second.
Therefore we create a new PAS $b$, and set $\sigma_b^1 = [3,5]$, $\sigma_b^2 = [3,4,5]$, and $Z_b = \{ 1 \}$.
(The choice of which segment is the first and second one is arbitrary.)  

Repeating the same process with the shortest path tree from Origin 2, we verify that it is relevant to PAS $a$ (which it already is), and add it as relevant to PAS $b$, so $Z_a = Z_b = \{ 1, 2\}$.
(Both origins are now relevant to both PASs).

\subsubsection{Flow shifts}

TAPAS uses flow shifts to find an equilibrium solution.
There are two types of flow shifts: the most common involves shifting flow between the two segments on an existing PAS.
The second involves identifying and eliminating cycles of used links for particular origins.

For the first type, assume we are given a PAS $\zeta$, and without loss of generality assume that the current travel time on the first segment $\sigma_1^\zeta$ is greater than that on the second $\sigma_2^\zeta$.
We wish to shift flow from the first segment to the second one to either equalize their travel times, or to shift all the flow to the second path if it is still shorter.
The total amount of flow we need to shift to equalize the travel times is approximately given by Newton's method:\index{Newton's method}
\labeleqn{tapasdesiredshift}{
\Delta h = \frac{\sum_{(i,j) \in \sigma_1^\zeta} t_{ij} - \sum_{(i,j) \in \sigma_2^\zeta} t_{ij}} {\sum_{(i,j) \in \sigma_1^\zeta} t'_{ij} + \sum_{(i,j) \in \sigma_2^\zeta} t'_{ij} }
\,.}
We must also determine whether such a shift is feasible (would shifting this much flow force an $x_{ij}^r$ value to become negative?) and, unlike the algorithms earlier in this chapter, how much of this flow shift comes from each of the relevant origins in $Z_\zeta$.

To preserve feasibility, for any relevant origin $r$, we must subtract the same amount from $x_{ij}^r$ for each link in the longer segment, and add the same amount to each link in the shorter segment.
Call this amount $\Delta h^r$.
The non-negativity constraints require $\Delta h^r \leq \min_{(i,j) \in \pi_1^\zeta} \myc{ x_{ij}^r }$; let $\overline{\Delta h}^r$\label{not:Deltahath} denote the right-hand side of this inequality, which must hold for every relevant origin.

If we have
\labeleqn{tapasfeasibility}{
\Delta h \leq \sum_{r \in Z_\zeta} \overline{\Delta h}^r
\,,}
then the desired shift is feasible, and we choose the origin-specific shifts $\Delta h^r$ to be proportional to their maximum values $\overline{\Delta h}^r$ to help maintain proportionality.
If this shift is not feasible, then we shift as much as we can by setting $\Delta h^r = \overline{\Delta h}^r$ for each relevant origin.

\stevefig{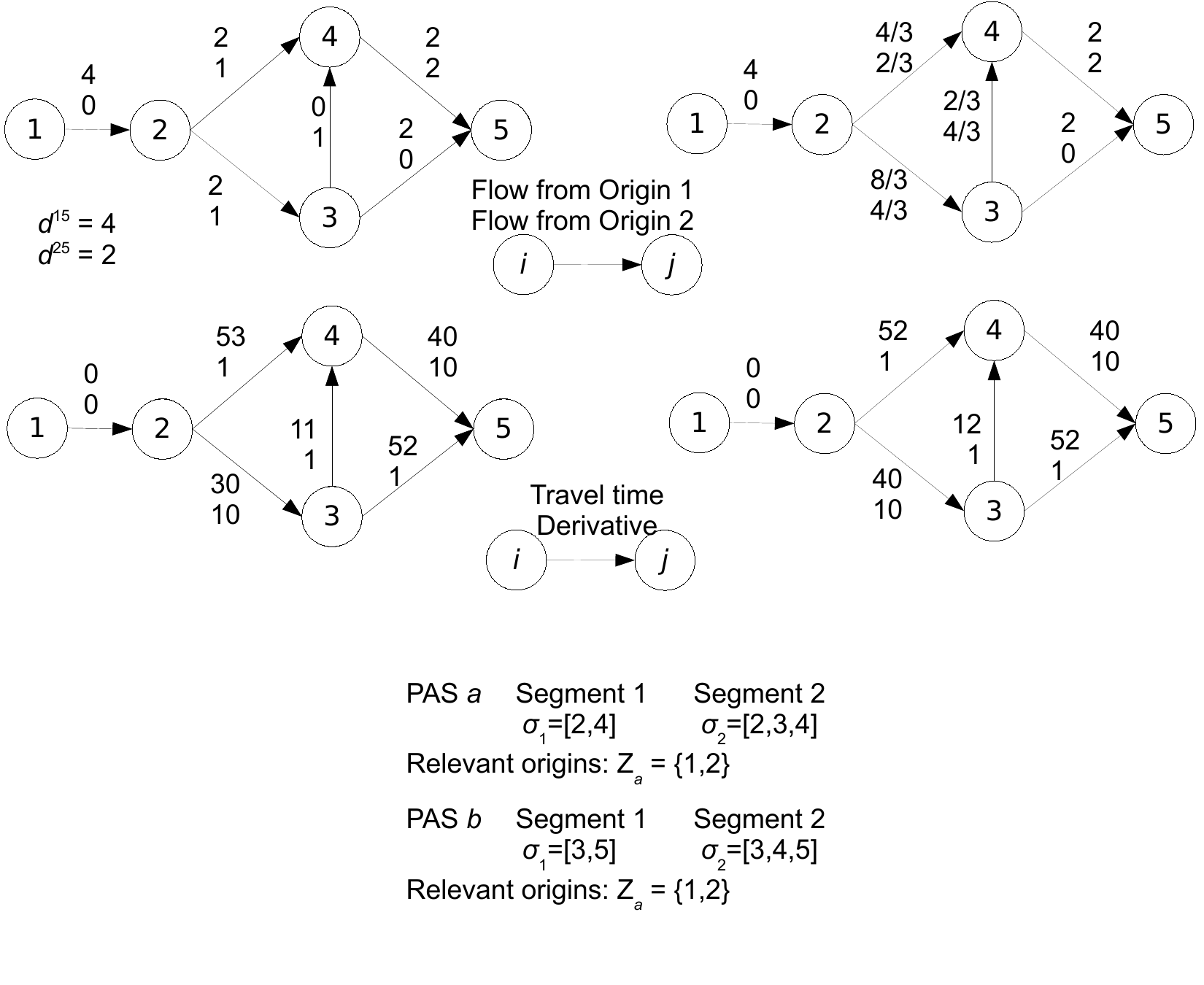}{Example of flow shifts using TAPAS; left panel shows initial flows and times, right panel shows flows and times after a shift for PAS $a$.}{\textwidth}

An example of such a shift is shown in Figure~\ref{fig:tapasflowshiftexample}, continuing the example from before.
The left side of the figure shows the state of the network prior to the flow shift.
The top panel shows the current origin-specific link flows; the middle panel the current travel times and travel time derivatives; and the bottom panel shows the structure of both PASs.
Starting with PAS $a$, we first calculate the desired total shift from equation~\eqn{tapasdesiredshift}:
\labeleqn{tapasdesiredshiftexample}{
\Delta h = \frac{53 - (30 + 11)}{2 + (10 + 1)} = 1
\,.}
For origin 1, we can subtract at most 2 units of flow from segment 1, and for origin 2, we can subtract at most 1.
Removing any more would result in negative origin flows on link $(2,4)$.
We thus split $\Delta h$ in proportion to these maximum allowable values, yielding
\labeleqn{tapasshiftdisaggregate}{
\Delta h^1 = 2/3 \qquad \Delta h^2 = 1/3
}
and producing the solution shown in the right half of Figure~\ref{fig:tapasflowshiftexample}.
Since the link performance functions are linear, Newton's method is exact, and travel times are equal on the two segments of PAS $a$.

Moving to the second PAS, we see that it is at equilibrium as well, and no flow shift is done: the numerator of equation~\eqn{tapasdesiredshift} is zero.
In fact, the entire network is now at equilibrium, but the origin-based link flows do not represent a proportional solution.
Proportionality adjustments are discussed below.

The second kind of flow shift involves removing flow from cycles.
In TAPAS, there may be occasions where cycles are found among the links with positive flow ($x_{ij}^r > 0$).
Such cycles can be detected using the topological ordering algorithm described in Section~\ref{sec:treesacyclic}.

In such cases, we can subtract flow from every link in the cycle, maintaining feasibility and reducing the value of the Beckmann function.
(See Exercise~\ref{ex:removecycle}).
Let $\overline{X}$ denote the minimum value of $x_{ij}^r$ among the links in such a cycle.
After subtracting this amount from every $x_{ij}^r$ in the cycle, the solution is closer to equilibrium and the cycle of positive-flow links no longer exists.

Figure~\ref{fig:tapascycle} shows an example of how this might happen.
This network has only a single origin, and two PASs.
Applying the flow shift formula, we move 1 vehicle from segment $[1,2]$ to $[1,3,2]$, and 1 vehicle from segment $[2,4]$ to segment $[2,3,4]$.
This produces the flow solution in the lower-left of the figure, which contains a cycle of flow involving links $(2,3)$ and $(3,2)$.
If we subtract 1 unit of flow from both of those links, we have the flow solution in the lower-right.
This solution is feasible and has a lower value of the Beckmann function, as you can verify.

\stevefig{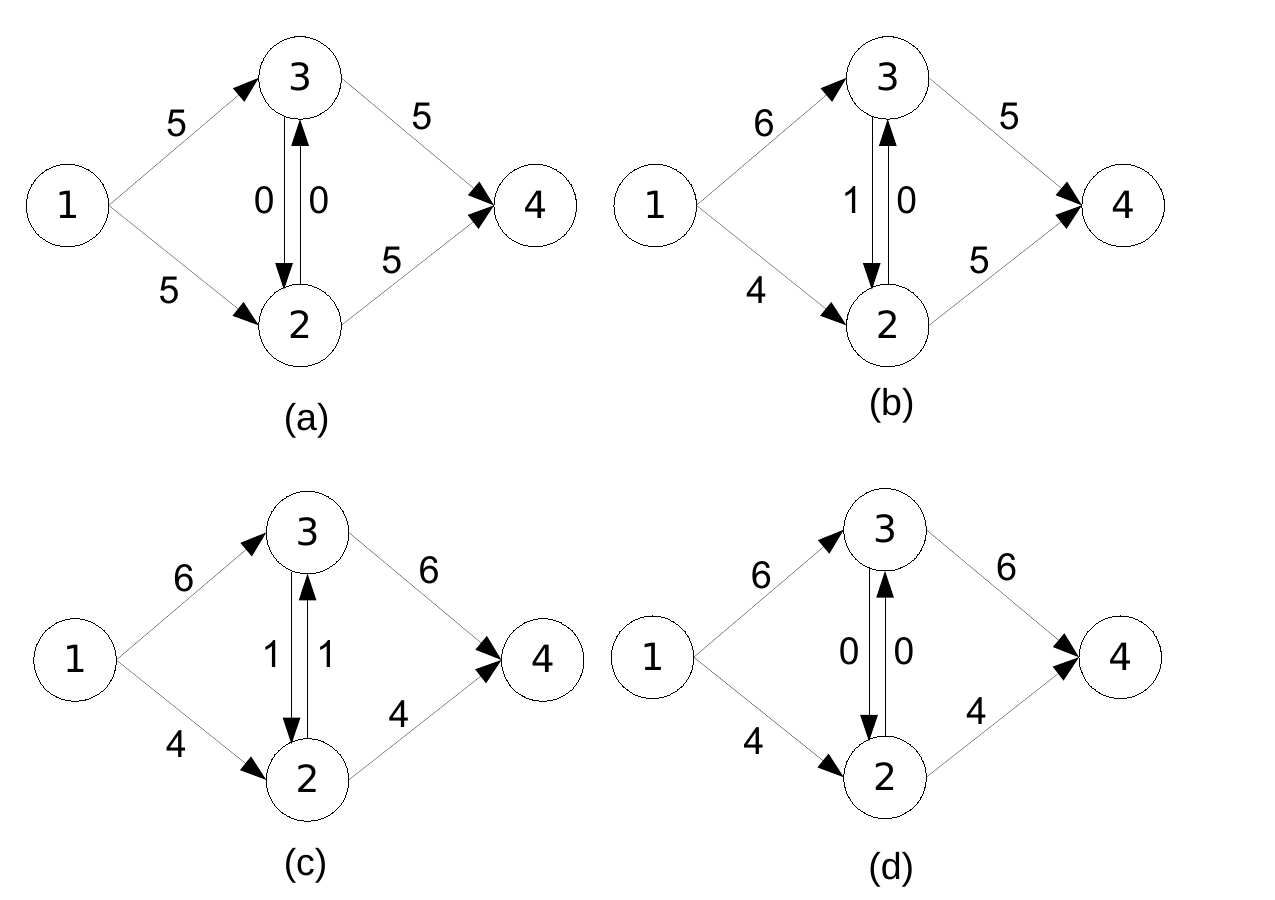}{(a) Initial flow solution; (b) After a flow shift at a PAS ending at 3; (c) After a flow shift at a PAS ending at 4; (d) After removing a cycle of flow.}{\textwidth}

\subsubsection{Proportionality adjustments}

\index{static traffic assignment!proportionality|(}
TAPAS uses proportionality adjustments to increase the entropy of the path flow solution.
Note that the path flow solution is not explicitly stored, since the number of used paths can grow exponentially with network size.
Rather, a path flow solution is implied by the bush, using the procedure described in Section~\ref{sec:mlpfprimal}, and the definitions of $\alpha$ and $h$ in equations~\eqn{tapasimplicit} and~\eqn{tapasimplicitalpha}.
That is, we calculate approach proportions $\alpha_{ij}^r$ using the formula

In a proportionality adjustment, we shift flows between segments in the PAS \emph{without} changing the total flow on each link, so some origins will shift flow from the links in $\sigma_1^\zeta$ to those in $\sigma_2^\zeta$, while other origins will shift flow from $\sigma_2^\zeta$ to $\sigma_1^\zeta$.
Following the previous section, we will use $\Delta h^r$ to denote the amount of flow shifted from each link in $\sigma_1^\zeta$ to each link in $\sigma_2^\zeta$, using negative numbers to indicate flow shifting from $\sigma_2^\zeta$ to $\sigma_1^\zeta$.
To maintain total link flows at their current levels, we will require
\labeleqn{tapasproportionalfeasible}{
\sum_{r \in Z_\zeta} \Delta h^r = 0
\,.}

To describe the problem more formally, let $b$ denote the node at the downstream end of the PAS, and let $x_b$ denote the total flow through this node as in Equation~\eqn{linknodeflow} .
We can calculate the flow on the segments $\sigma_1^\zeta$ and $\sigma_2^\zeta$ for each relevant origin $r$ with the formulas
\begin{align}
\label{eqn:segment1flow} g^r(\sigma_1^\zeta) &= x^r_b \prod_{(i,j) \in \sigma_1^\zeta} \frac{x_{ij}^r}{x_j^r} \\
\label{eqn:segment2flow} g^r(\sigma_2^\zeta) &= x^r_b \prod_{(i,j) \in \sigma_2^\zeta} \frac{x_{ij}^r}{x_j^r} \,, \\
\end{align}
assuming positive flow through all nodes in the segment ($x_j^r > 0$).\footnote{What would happen if this were not true?}
If proportionality were satisfied, we would have
\labeleqn{tapasproportionality}{
\frac{g^r(\sigma_1^\zeta)}{g^r(\sigma_1^\zeta) + g^r(\sigma_2^\zeta)} = \frac{\sum_{r' \in Z_\zeta} g^{r'}(\sigma_1^\zeta) }{\sum_{r' \in Z_\zeta} \myp{g^{r'}(\sigma_1^\zeta) + g^{r'}(\sigma_2^\zeta)}}
}
for all relevant origins.

After applying segment shifts of size $\Delta h^r$, the new segment flows will be given by
\begin{align}
\label{eqn:segment1flow2} g^r(\sigma_1^\zeta) &= x^r_b \prod_{(i,j) \in \pi_1} \frac{x_{ij}^r - \Delta h^r}{x_j^r - \Delta h^r [r \neq j]} \\
\label{eqn:segment2flow2} g^r(\sigma_2^\zeta) &= x^r_b \prod_{(i,j) \in \pi_2} \frac{x_{ij}^r + \Delta h^r}{x_j^r + \Delta h^r [r \neq j]} \,, \\
\end{align}
using brackets for an indicator function.\label{not:indicator}
We aim to find $\Delta h^r$ values satisfying constraint~\eqn{tapasproportionalfeasible} and~\eqn{tapasproportionality}, where the segment flows are computed with equations~\eqn{segment1flow2} and~\eqn{segment2flow2}.

Solving this optimization problem exactly is a bit difficult because equations~\eqn{tapasimplicit} and~\eqn{tapasimplicitalpha} are nonlinear.
A good approximation method is developed in Exercise~\ref{ex:tapasproportionalapprox}.
A simpler heuristic is to adapt the ``proportionality between origins'' technique from Section~\ref{sec:mlpfprimal} and approximate the (nonlinear) formulas~\eqn{segment1flow2} and~\eqn{segment2flow2} by the (linear) formulas
\begin{align}
\label{eqn:segment1heuristic} g^r(\sigma_1^\zeta) &\approx g^r_0(\sigma_1^\zeta) - \Delta h^r \\
\label{eqn:segment2heuristic} g^r(\sigma_2^\zeta) &\approx g^r_0(\sigma_2^\zeta) + \Delta h^r 
\end{align}
where $g^r_0(\sigma_1^\zeta)$ and $g^r_0(\sigma_2^\zeta)$ are the \emph{current} segment flows (with zero shift).

Substituting equations~\eqn{segment1heuristic} and~\eqn{segment2heuristic} into~\eqn{tapasproportionality} and simplifying, we obtain
\labeleqn{tapasheuristic}{
\Delta h^r = g^r_0(\sigma_1^\zeta) - (g^r_0(\sigma_1^\zeta) + g^r_0(\sigma_2^\zeta)) \frac{\sum_{r' \in Z_\zeta} g^{r'}_0(\sigma_1^\zeta)} {\sum_{r' \in Z_\zeta} \myp{g^{r'}_0(\sigma_1^\zeta) + g^{r'}_0(\sigma_2^\zeta)}}
}
as a flow shift to heuristically move toward proportionality.

This heuristic is in fact exact if the PAS is ``isolated'' in the sense that flow does not enter or leave the segments in the middle, so $\alpha_{ij}^r = 1$ for all links in $\sigma_1^\zeta$ and $\sigma_2^\zeta$ except for the last links of each segment.

To illustrate how this procedure works, consider the example in Figure~\ref{fig:tapasproportionshift}.
In this example, there are multiple destinations in addition to multiple origins.
The figure shows the path flow solution implied by the link flow solution at the left.
We emphasize that TAPAS does not maintain the path flow solution explicitly, and the path flows are constructed from the link flows using equation~\eqn{tapasimplicit}.
For instance, the flow on path $[1,3,4]$ from Origin 1 is calculated as $3 \times \frac{1}{3} \times \frac{1}{1} = 1$.
This solution does not satisfy proportionality.
All the vehicles from Origin 1 passing between nodes 2 and 4 use segment [2,4], while all of those from Origin 2 passing between these nodes use segment [2,3,4].
Since the total flow on the two segments are equal (two vehicles on each), flow from both origins should split equally between the two segments.

\stevefig{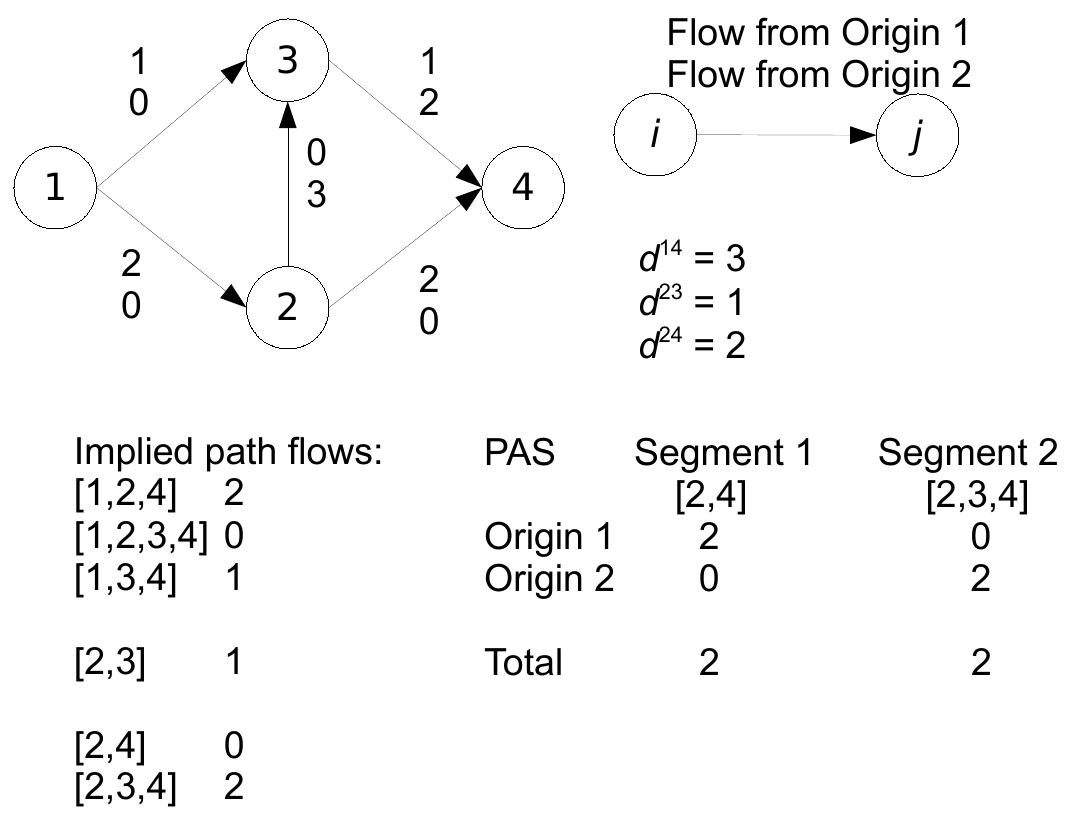}{Example of a proportionality adjustment in TAPAS.}{0.8\textwidth}

Applying equation~\eqn{tapasheuristic} gives the shifts shown in Figure~\ref{fig:tapasproportionalitypasa}.
This figure also shows the new origin-specific link flows and implied path flows.
In this case, the link flows on the segments of PAS $a$ now satisfy proportionality.
This is a case where the PAS is isolated, because no vehicles entered and left the segments in the middle, and the heuristic formula~\eqn{tapasheuristic} is exact.

\stevefig{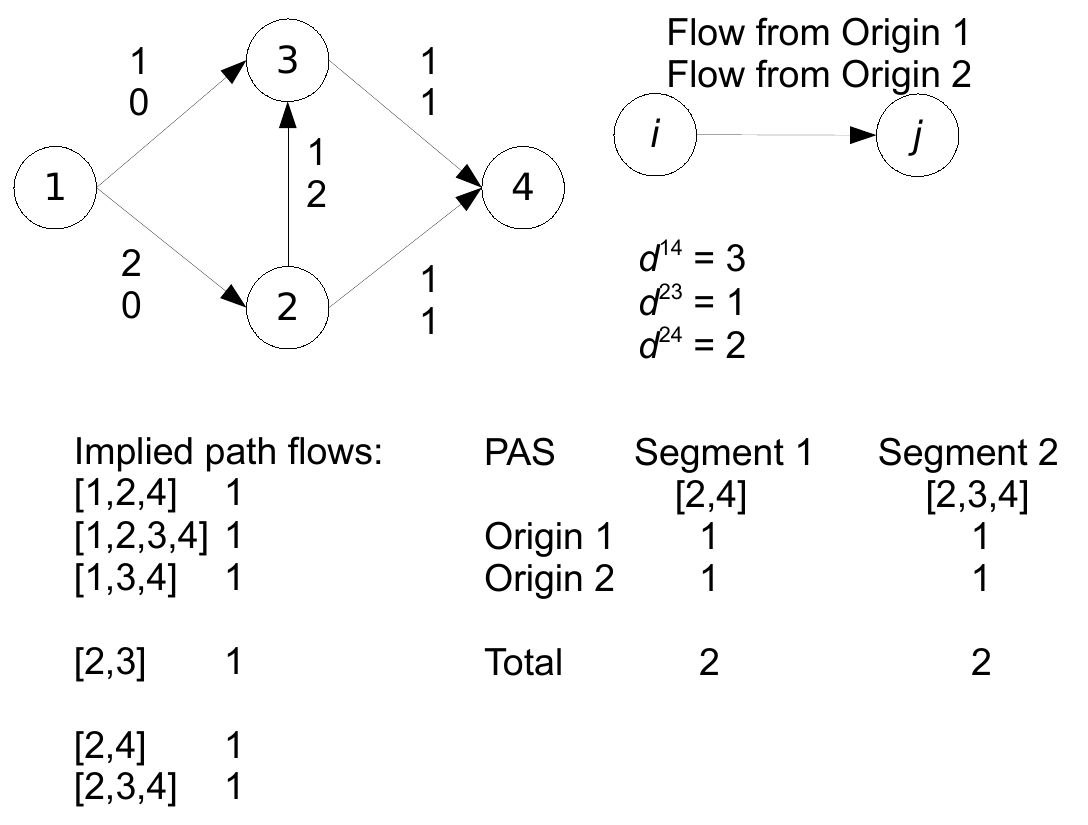}{Updated flows after a proportionality adjustment for PAS $a$.}{0.8\textwidth}

To show how the formula is inexact for a non-isolated PAS, consider the modification of this example shown in Figure~\ref{fig:tapasnonisolatedshift}.
The only change is that node 3 is now a destination for Origin 1, with a demand of 1, and that as a result the flow from Origin 1 on link (1,3) is increased by one vehicle.
Repeating the same process as above, and applying~\eqn{tapasheuristic}, we again one swap one vehicle between each pair of segments.
This results in the situation in Figure~\ref{fig:tapasnonisolatedshift2}.
The origin-specific link flows are the same as in Figure~\ref{fig:tapasproportionalitypasa}, except for the additional vehicle from Origin 1 on link (1,3).
But the implied path flows are quite different!
This is because the additional vehicle shifted onto link (3,4) was ``split'' between incoming links (1,3) and (2,3), according to equation~\eqn{tapasimplicit}, rather than allocated solely to (2,3).
(Again, TAPAS does not store the flows on individual paths, and must calculate them implicitly using this formula.)  As a result, proportionality is still not satisfied: between segments [2,4] and [2,3,4], Origin 1 splits in the ratio of 3:2, whereas Origin 2 splits in the ratio 1:1.
This is closer to proportionality from before, but not exact.
Repeated applications of the heuristic shift formula will converge to a proportional solution, in this case.
\index{static traffic assignment!proportionality|)}
\index{static traffic assignment!algorithms!traffic assignment by paired alternate segments (TAPAS)|)}
\index{pair of alternate segments|)}
\index{static traffic assignment!properties!nonuniqueness in path flows|)}

\stevefig{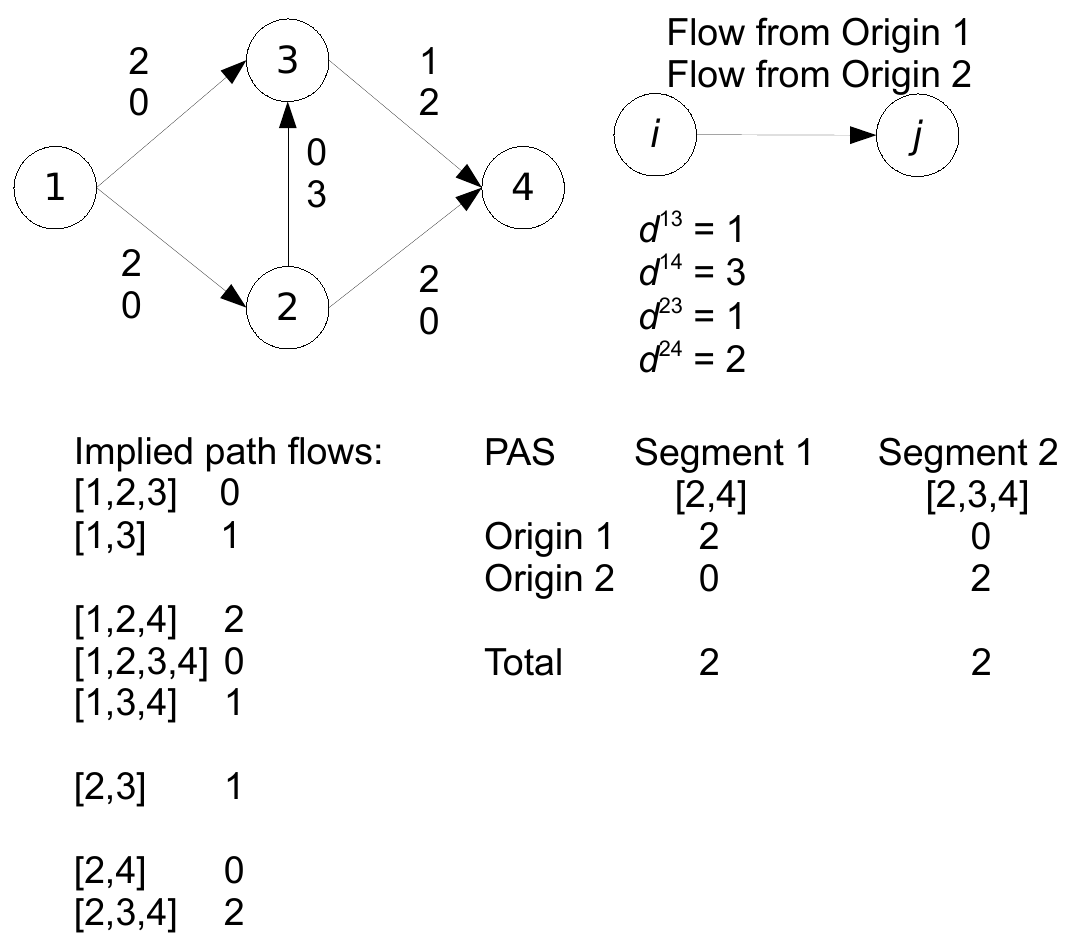}{Example of a proportionality adjustment for a non-isolated PAS.}{0.8\textwidth}

\stevefig{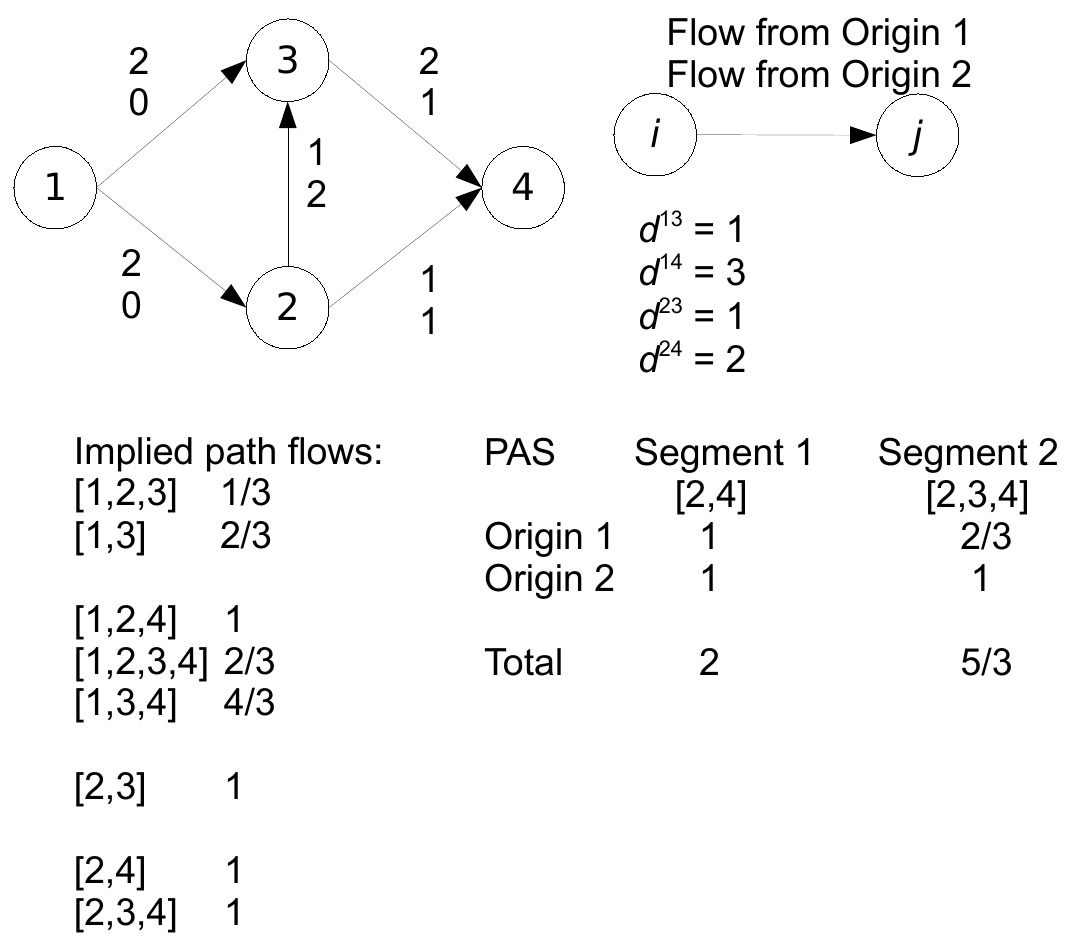}{The heuristic formula is not always exact for a non-isolated PAS.}{0.8\textwidth}

\section{Historical Notes and Further Reading}
\label{sec:algosreference}

Extensive reviews of link-based and path-based algorithms for the traffic assignment problem are found in \cite{patriksson94} and \cite{florian95}.
Bush-based algorithms are not described in these references, having been first developed independently in \cite{dial99_bobtail} and \cite{bargera_diss}.
The classification of algorithms as link-based,\index{static traffic assignment!algorithms!link-based} path-based,\index{static traffic assignment!algorithms!path-based} or bush-based\index{static traffic assignment!algorithms!bush-based} was first proposed by \cite{nie10}.
Many gap measures have been proposed as convergence criteria.
Some authors refer to path-based methods as \emph{route-based},\index{route|see {path}} and some authors refer to bush-based methods as \emph{origin-based}.

For the relative gap\index{gap function!relative gap} variants $\gamma_1$, $\gamma_2$, and $\gamma_3$ defined in the text, see \cite{dtaprimer}, \cite{patriksson94}, and \cite{boyce04} for more detail.
The average excess cost\index{gap function!average excess cost} and maximum excess cost\index{gap function!maximum excess cost} were proposed in \cite{bargera02}.
In studying the Philadelphia network, \cite{boyce04} found that freeway link flows stabilized once a relative gap of $10^{-4}$ was reached.
This specific recommendation was generalized by \cite{patil21} by considering other networks and other outputs from traffic assignment.
Aggregate values, such as total system travel time and vehicle-miles traveled, stabilize around a relative gap of $10^{-4}$; typical link flows (not just freeways) around $10^{-5}$; and (entropy-maximizing) path flows around $10^{-6}$.

For the specific algorithms discussed in this text, the method of successive averages\index{static traffic assignment!algorithms!method of successive averages} and Frank-Wolfe\index{static traffic assignment!algorithms!Frank-Wolfe} are both instances of the more general ``convex combinations'' method, and can in fact be applied to any convex optimization problem~\citep{bertsekas_nlp}.
The Frank-Wolfe method itself was proposed in \cite{frank56}, and the conjugate version\index{static traffic assignment!algorithms!conjugate Frank-Wolfe} in \cite{mitradjieva13}.
Notable link-based algorithms not presented in this chapter are the simplicial decomposition methods\index{static traffic assignment!algorithms!simplicial decomposition} of \cite{smith83a} and \cite{Law1984}; but see Chapter~\ref{chp:staticextensions} for discussion of this method in the setting of equilibrium with link interactions.

The gradient projection\index{static traffic assignment!algorithms!gradient projection} and manifold suboptimization\index{static traffic assignment!algorithms!manifold suboptimization} algorithms were presented in \cite{jayakrishnan94} and \cite{florian09}, respectively.
Another notable path-based algorithm not described here is the disaggregate simplicial decomposition\index{static traffic assignment!algorithms!simplicial decomposition} method of \cite{larsson92}.

\cite{nie10} described a general framework for bush-based algorithms, uniting earlier work on origin-based assignment~\citep{bargera02}, Algorithm B~\citep{dial06a}, and local user cost equilibrium~\citep{gentile14}.
See \cite{xie13} for a discussion about the close relationships between origin-based assignment and local user cost equilibrium.
Interestingly, many of the concepts in bush-based algorithms were anticipated in the study of routing in telecommunications networks;\index{telecommunication networks} see \cite{gallager77a} and \cite{bertsekas84} for examples of such work.

The primal method\index{static traffic assignment!algorithms!primal method for proportionality} for maximizing path flow entropy is described at greater length in \cite{bargera06}.
The dual method\index{static traffic assignment!algorithms!dual method for proportionality} is the conjugate gradient method of \cite{larsson01}; another dual method not described here is iterative balancing; see \cite{bell97}.
Traffic assignment by paired alternate segments was presented in \cite{bargeratapas}.\index{static traffic assignment!algorithms!traffic assignment by paired alternate segments (TAPAS)}\index{static traffic assignment!algorithms|)}

\index{simplicial decomposition|see {static traffic assignment, algorithms, simplicial decomposition}}
\index{gradient projection|see {static traffic assignment, algorithms, gradient projection; dynamic traffic assignment, algorithms, gradient projection}}
\index{Algorithm B|see {static traffic assignment, algorithms, Algorithm B}}
\index{traffic assignment by paired alternate segments|see {static traffic assignment, algorithms, traffic assignment by paired alternate segments}}
\index{origin-based assignment|see {static traffic assignment, algorithms, origin-based assignment}}
\index{manifold suboptimization|see {static traffic assignment, algorithms, manifold suboptimization}}

\section{Exercises}
\label{exercises_solutionalgorithms}

\begin{enumerate}
\item \diff{32} One critique of the BPR link performance function\index{link performance function!BPR (Bureau of Public Roads)} is that it allows link flows to exceed capacity.
An alternative link performance ``function'' is $t_{ij} = t^0_{ij} / (u_{ij} - x_{ij})$ if $x_{ij} < u_{ij}$, and $\infty$ otherwise, where $t^0_{ij}$ and $u_{ij}$ are the free-flow time and capacity of link $(i,j)$.
First show that $t_{ij} \rightarrow \infty$ as $x_{ij} \rightarrow u_{ij}$.
How would using this kind of link performance function affect the solution algorithms discussed in this chapter?
\item \diff{33} Show that the relative gap $\gamma_1$ and average excess cost are always nonnegative, and equal to zero if and only if the link or path flows satisfy the principle of user equilibrium.
\item \diff{61} Some relative gap definitions require a lower bound $\underline{f}$ on the value of the Beckmann function at optimality.
Let $\mb{x}$ denote the current solution, $f(\mb{x})$ the value of the Beckmann function at the current solution, and $TSTT(\mb{x})$ and $SPTT(\mb{x})$ the total system travel time and shortest path travel time at the current solution, respectively.
Show that $f(\mb{x}) + SPTT(\mb{x}) - TSTT(\mb{x})$ is a lower bound on the Beckmann function at user equilibrium.
\item \diff{10} What is the value of the lower bound $\underline{f} = f(\mb{x}) + SPTT(\mb{x}) - TSTT(\mb{x})$ if $\mb{x}$ satisfies the principle of user equilibrium?
\item \diff{23} Let $(\mb{h}, \mb{x})$ and $(\mb{g}, \mb{y})$ be two feasible solutions to the Beckmann formulation~\eqn{tapoptstart}--\eqn{tapoptend}, and let $\lambda \in [0, 1]$.
Show that $(\lambda \mb{h} + (1 - \lambda) \mb{g}, \lambda \mb{x} + (1 - \lambda) \mb{y})$ is also feasible, directly from the constraints (without appealing to convexity.)
\item \diff{35}  In the network in Figure~\ref{fig:aonex}, all trips originate at node A.
The links are labeled with the current travel times, and the nodes are labeled with the number of trips whose destination is that node.
\label{ex:aonex}
\begin{enumerate}
\item Find the shortest paths from node A to all other nodes, and report the cost and backnode labels upon termination.
\item What would be the target link flow solution $\mathbf{\hat{x}}$ in the method of successive averages or the Frank-Wolfe algorithm?
\end{enumerate}
\stevefig{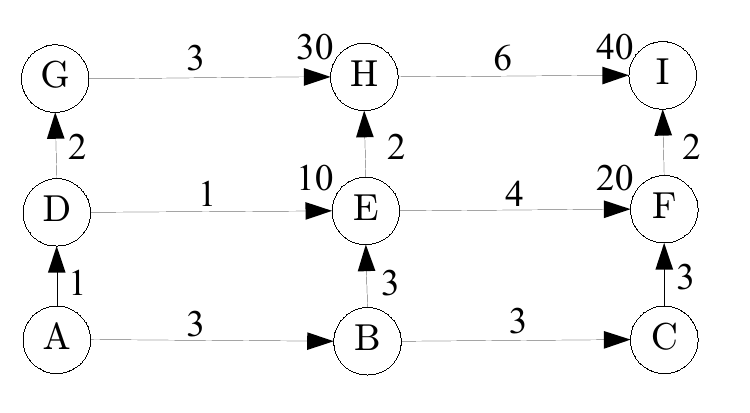}{Network for Exercise~\ref{ex:aonex}.}{0.6\textwidth}
\item \diff{73}  All-or-nothing assignments\index{all-or-nothing assignment} $\mb{\hat{x}}$ play a major role in link-based algorithms.
A na\"{i}ve way to calculate these is to start with zero flows on each link; then find the shortest path from each origin $r$ to each destination $s$; then add $d^{rs}$ to each link in this path.
This method may require adding up to $|Z|^2$ terms for each link, in case every shortest path uses the same link.
Formulate a more efficient algorithm which requires solving one shortest path problem per origin, and which requires adding no more than $|Z|$ terms for each link.
\emph{(Hint: Do not wait until the end to calculate $\mb{\hat{x}}$ and find a way to build $\mb{x}$ as you go.)}
\item \diff{42} Consider the network in Figure~\ref{fig:xstarpractice}, with a single origin and two destinations.
Each link has the link performance function $10 + x^2$, and the boldface links indicate the links used in the previous $\mb{\hat{x}}$ target.
Report the new target link flows $\mb{\hat{x}}$, the step size $\lambda$, and the new resulting link flows, according to (a) Frank-Wolfe and (b) conjugate Frank-Wolfe.
\label{ex:xstarpractice}
\stevefig{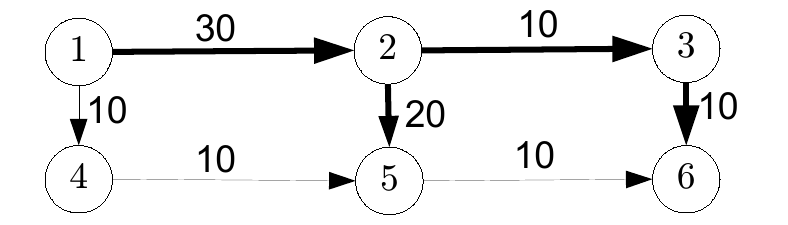}{Network for Exercise~\ref{ex:xstarpractice}, boldface links indicate previous $\mb{\hat{x}}$ target.}{0.6\textwidth}
\item \diff{47} Consider the network in Figure~\ref{fig:msafwnetex}, where 8 vehicles travel from node 1 to node 4.
Each link is labeled with its delay function.
For each of the algorithms listed below, report the resulting link flows, average excess cost, and value of the Beckmann function.
\label{ex:msafwnetex}
\begin{enumerate}
\item Perform three iterations of the method of successive averages.
\item Perform three iterations of the Frank-Wolfe algorithm.
\item Perform three iterations of conjugate Frank-Wolfe.
\item Perform three iterations of gradient projection.
\item Perform three iterations of manifold suboptimization.
\item Perform three iterations of Algorithm B (for each iteration, do one flow update and one bush update)
\item Perform three iterations of origin-based assignment (for each iteration, do one flow update and one bush update)
\item Perform three iterations of linear user cost equilibrium (for each iteration, do one flow update and one bush update)
\item Compare and discuss the performance of these algorithms.
\end{enumerate}
\stevefig{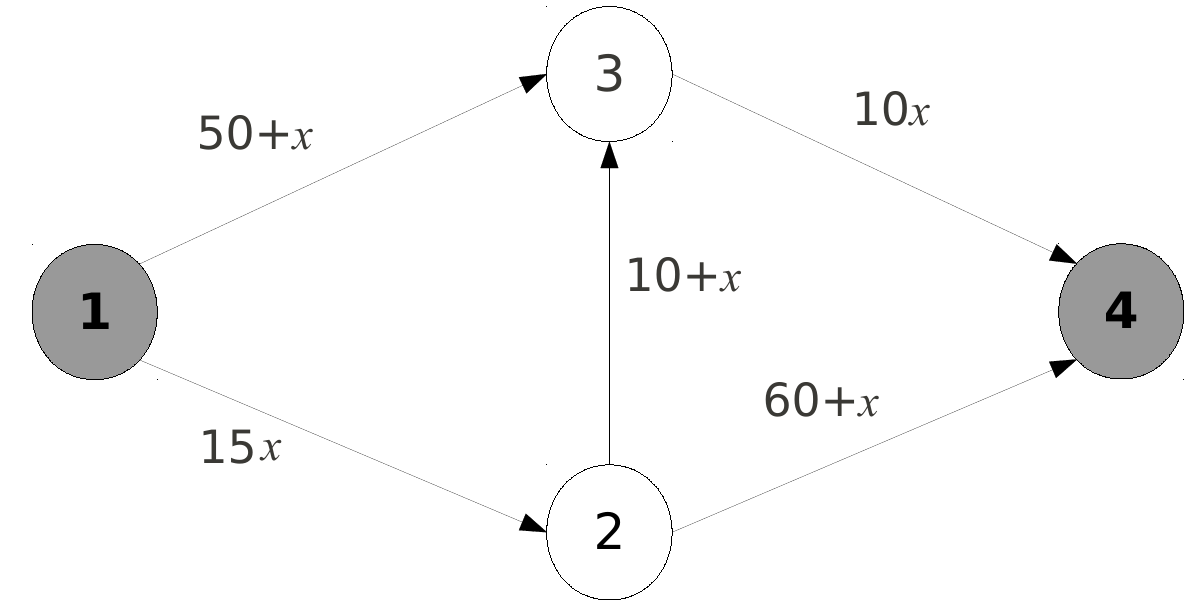}{Network for Exercise~\ref{ex:msafwnetex}.}{0.6\textwidth}
\item \diff{48} Consider the network in Figure~\ref{fig:exgridnet}.
The cost function on the light links is $3 + (x_a/200)^2$, and the delay function on the thick links is $5 + (x_a/100)^2$.
1000 vehicles are traveling from node 1 to 9, and 1000 vehicles from node 4 to node 9.
For each of the algorithms listed below, report the resulting link flows, average excess cost, and value of the Beckmann function.\label{ex:gridnetfw}
\begin{enumerate}
\item Perform three iterations of the method of successive averages.
\item Perform three iterations of the Frank-Wolfe algorithm.
\item Perform three iterations of conjugate Frank-Wolfe.
\item Perform three iterations of projected gradient.
\item Perform three iterations of gradient projection.
\item Perform three iterations of Algorithm B (for each iteration, do one flow update and one bush update)
\item Perform three iterations of origin-based assignment (for each iteration, do one flow update and one bush update)
\item Perform three iterations of linear user cost equilibrium (for each iteration, do one flow update and one bush update)
\item Compare and discuss the performance of these algorithms.
\end{enumerate}
\stevefig{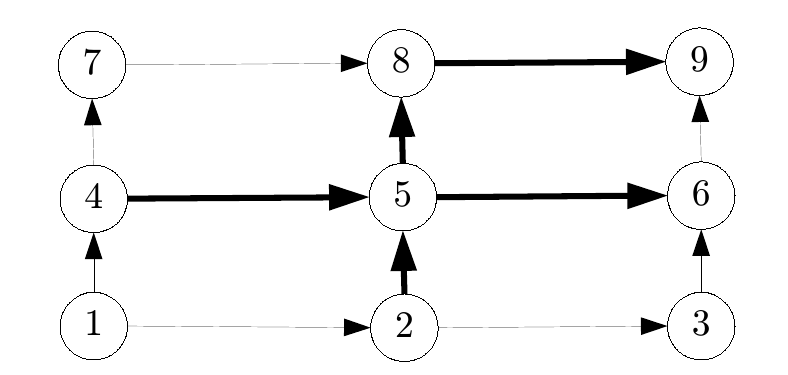}{Network for Exercise~\ref{ex:gridnetfw}}{0.6\textwidth}
\item \diff{43} The method of successive averages, as presented in the text, uses the step size $\lambda_i = 1/(i+1)$ at iteration $i$.
Other choices of step size can be used, and Exercise~\ref{ex:msaproof} shows that the algorithm converges whenever $\lambda_i \in [0, 1]$, $\sum \lambda_i = \infty$ and $\sum \lambda_i^2 < \infty$.
Which of the following step size choices guarantee convergence? \label{ex:msastepsize}
\begin{enumerate}
\item $\lambda_i = 1/(i+2)$
\item $\lambda_i = 4/(i+2)$
\item $\lambda_i = 1/i^2$
\item $\lambda_i = 1/(\log i)$
\item $\lambda_i = 1/\sqrt{i}$
\item $\lambda_i = 1/i^{2/3}$
\end{enumerate}
\item \diff{65} \emph{(Proof of convergence for the method of successive averages.)\index{static traffic assignment!algorithms!method of successive averages}}  Consider the method of successive averages applied to the vector of link flows.
This produces a sequence of link flow vectors $\mb{x^1}, \mb{x^2}, \mb{x^3}, \ldots$ where $\mb{x^i}$ is the vector of link flows at iteration $i$.
We can also write down the sequence $f^1, f^2, f^3, \ldots$ of the values taken by the Beckmann function for $\mb{x^1}, \mb{x^2}$, etc.
To show that this algorithm converges to the optimal solution, we have to show that either $\mb{x^i} \rightarrow \mb{\hat{x}}$ or $f^i \rightarrow \hat{f}$ as $i \rightarrow \infty$, where $\mb{\hat{x}}$ is the user equilibrium solution and $\hat{f}$ the associated value of the Beckmann function.
This exercise walks through one proof of this fact, for any version of the method of successive averages for which  $\lambda_i \in [0, 1]$, $\sum \lambda_i = \infty$ and $\sum \lambda_i^2 < \infty$.
\label{ex:msaproof}
\begin{enumerate}
\item Assuming that the link performance functions are differentiable, show that for any feasible $\mb{x}$ and $\mb{y}$ there exists $\theta \in [0,1]$ such that
\begin{multline} \label{eqn:msataylor} f(\mb{y}) = f(\mb{x}) + \sum_{(i,j) \in A} t_{ij}(x_{ij}) (y_{ij} - x_{ij}) + \\ \frac{1}{2} \sum_{(i,j) \in A} t'_{ij}((1 - \theta) x_{ij} + \theta y_{ij}) (y_{ij} - x_{ij})^2 \,.
\end{multline}
\item Setting $\mb{x} = \mb{x^i}$ and $\mb{y} = \mb{x^{i+1}}$, recast equation~\eqn{msataylor} into an expression for the difference in the values of the Beckmann function between two consecutive iterations of the method of successive averages, in terms of $\lambda_i$ and $\mb{x^*_i}$.
\item Sum the resulting equation over an infinite number of iterations to obtain a formula for the limiting value $\hat{f}$ of the sequence $f^1, f^2, \ldots$.
\item Use the facts that $\sum \lambda_i = \infty$, $\sum \lambda_i^2 < \infty$, and that $\hat{f}$ has a finite value to show that the limiting values of $SPTT(\mb{x^i})$ and $TSTT(\mb{x^i})$ must be equal, implying that the limit point is a user equilibrium.
\end{enumerate}
\item \diff{34} The derivation leading to~\eqn{balancingfinal} assumed that the solution to the restricted VI was not at the endpoints $\lambda = 0$ or $\lambda = 1$.
Show that if you are solving~\eqn{balancingfinal} using either the bisection method from Section~\ref{sec:bisection}, or Newton's method (with a ``projection'' step ensuring $\lambda \in [0, 1]$), you will obtain the correct solution to the restricted VI even if it is at an endpoint.
\label{ex:cornersolution}
\item \diff{55} \emph{(Linking Frank-Wolfe\index{static traffic assignment!algorithms!Frank-Wolfe} to optimization.)}  At some point in the Frank-Wolfe algorithm, assume that the current link flows are $\textbf{x}$ and the target link flows $\mathbf{x^*}$ have just be found, and we need to find new flows $\mathbf{x'}(\lambda) = \lambda \mathbf{x^*} + (1 - \lambda) \mathbf{x}$ for some $\lambda \in [0, 1]$.
Let $z(\lambda)$ be the value of the Beckmann function at $\mathbf{x'}(\lambda)$.
\label{ex:fwdescent}
\begin{enumerate}
\item Using the multi-variable chain rule, we can show that $z$ is differentiable and $z'(\lambda)$ is the dot product of the gradient of the Beckmann function evaluated at $\mathbf{x'}(\lambda)$ and the direction $\mathbf{x^* - x}$.
Calculate the gradient of the Beckmann function and use this to write out a formula for $z'(\lambda)$.
\item Is $z$ a convex function of $\lambda$?
\item Show that $z'(0) = 0$ only if $\mathbf{x}$ is an equilibrium, and that otherwise $z'(0) < 0$.
\item Assume that the solution of the restricted variational inequality in the Frank-Wolfe algorithm is for an ``interior'' point $\lambda^* \in (0, 1)$.
Show that $z'(\lambda^*) = 0$.
\item Combine the previous answers to show that the Beckmann function never increases after an iteration of the Frank-Wolfe algorithm (and always decreases strictly if not at an equilibrium).
\end{enumerate}
\item \diff{74} \emph{(Proof of convergence for Frank-Wolfe.)\index{static traffic assignment!algorithms!Frank-Wolfe}}  Exercise~\ref{ex:fwdescent} shows that the sequence of Beckmann function values $f^1, f^2, \ldots$ from subsequent iterations of Frank-Wolfe is nonincreasing.
Starting from this point, show that this sequence has a limit, and that the resulting limit corresponds to the global minimum of the Beckmann function (demonstrating convergence to equilibrium.)  Your solution may require knowledge of real analysis.
\label{ex:fwproof}
\item \diff{11} When is the Beckmann function quadratic?
\item \diff{33} Identify conjugate directions for the following quadratic programs:
\begin{enumerate}
\item $f(x_1, x_2) = x_1^2 + x_2^2$
\item $f(x_1, x_2) = x_1^2 + x_2^2 + \frac{2}{3} x_1 x_2$
\item $f(x_1, x_2) = 2x_1^2 + 3x_2^2 + \frac{1}{9} x_1 x_2$
\end{enumerate}
\item \diff{38} The \emph{biconjugate} Frank-Wolfe method\index{static traffic assignment!algorithms!biconjugate Frank-Wolfe} chooses a target vector $\mb{x^*}$ so that the search direction $\mb{x^*}$ is conjugate to the last two search directions, rather than just the last one.
Let $\mb{x^{AON}}$ reflect the all-or-nothing solution at the current point, $\mb{x^*_{-1}}$ the target vector used at the last iteration, and $\mb{x^*_{-2}}$ the target vector used two iterations ago.
Also let $\lambda_{-1}$ be the step size used for the last iteration, and define
\labeleqn{bfwmu}{
\mu = -\frac{\sum_{(i,j) \in A} (\lambda_{-1} (x^*_{-1})_{ij} + (1 - \lambda_{-1}) (x^*_{-2})_{ij} - x_{ij}) (x^{AON}_{ij} - x_{ij}) t'_{ij}}
            {\sum_{(i,j) \in A} (\lambda_{-1} (x^*_{-1})_{ij} + (1 - \lambda_{-1}) (x^*_{-2})_{ij} - x_{ij}) ((x^*_{-2})_{ij} - (x^*_{-1})) t'_{ij}}
}
and
\labeleqn{bfwnu}{
\nu = \frac{\mu \lambda_{-1}}{1 - \lambda_{-1}} -
            \frac{\sum_{(i,j) \in A} ((x^*_{-1})_{ij} - x_{ij}) (x^{AON}_{ij} - x_{ij}) t'_{ij}}
                 {\sum_{(i,j) \in A}  ((x^*_{-1})_{ij} - x_{ij})^2 t'_{ij}}
\,.}
Show that the formula
\labeleqn{bfwtarget}{
\mb{x^*} = \frac{1}{1 + \mu + \nu} \myp{ \mb{x^{AON}} + \nu \mb{x^*_{-1}} + \mu \mb{x^*_{-2}}}
}
gives both a feasible target point $\mb{x^*}$, and one conjugate to the two previous search directions based on the Hessian of the Beckmann function at the current solution.
\label{ex:biconjugate}
\item \diff{41} Section~\ref{sec:algorithmb} includes an example for Algorithm B, and Figure~\ref{fig:algbstep} shows the bush link flows at the end of a flow shifting operation.
Perform a second iteration of flow shifting on the same bush, recalculating $L$ and $U$ labels, and scanning all nodes in reverse topological order.
Report the new link flows and travel times.
Also report the new bush after eliminating unused links and adding shortcuts.
Does your answer depend on whether you use the $L$ or $U$ labels to define ``shortcuts''?
\item \diff{43} Figure~\ref{fig:bushpractice} shows a bush with the current link flows labeled.
Each link has delay function $5 + 2x^2$.
Calculate the $L$, $U$, $M$, $D$, and $\alpha$ labels for all links and nodes in the bush, and identify the divergence node for each bush node.
\label{ex:bushpractice}
\begin{enumerate}[(a)]
\item Perform one step of flow shifting using Algorithm B.
\item Perform one step of flow shifting using OBA (starting from the original flows).
\item Perform one step of flow shifting using LUCE (starting from the original flows).
\item Report the maximum excess cost before and after each of these flow shifts.
Which algorithm reduced this gap measure by the most?
\end{enumerate}
\stevefig{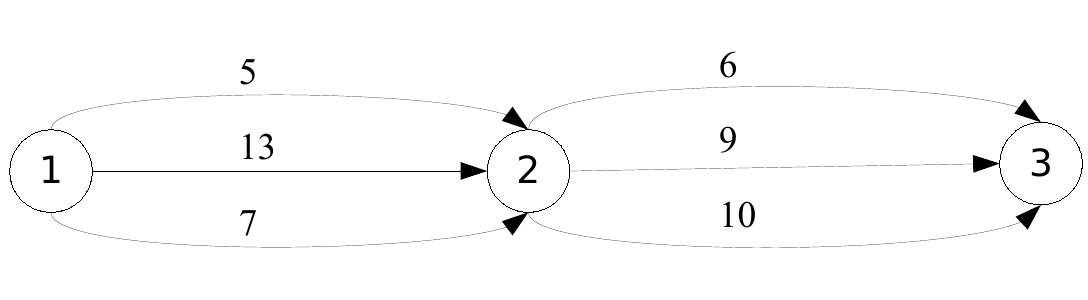}{Bush for Exercise~\ref{ex:bushpractice} with current flows.}{0.8\textwidth}
\item \diff{73}.\index{path!set of used paths!consistency}
This exercises walks through a proof of the formula~\eqn{twoconsistency} for choosing the threshold $\epsilon$ for finding proportional path flows.
A set of paths $\hat{\Pi} = \cup \hat{\Pi}_{rs}$ is \emph{2-consistent} if there are no paths $\pi_1, \pi_2, \pi_3, \pi_4$ satisfying the following conditions: (1) $\pi_1$ and $\pi_3$ are in $\hat{\Pi}$; (2) at least one of $\pi_2$ and $\pi_4$ is not in $\hat{\Pi}$; (3) $\pi_1$ and $\pi_2$ connect the same OD pair; (4) $\pi_3$ and $\pi_4$ connect the same OD pair; and (5) $\pi_1$ and $\pi_3$ use the same links as $\pi_2$ and $\pi_4$, exactly the same number of times.
For instance, in Figure~\ref{fig:consistency}, the path set $\hat{\Pi} = \myc{ [1,2], [1,3,5,6,4,2], [3,5,6,4], [3,5,7,8,6,4] }$ is not 2-consistent, because we can choose $\pi_1 = [1,3,5,6,4,2]$, $\pi_2 = [1,3,5,7,8,6,4,2]$, $\pi_3 = [3,5,6,4]$, and $\pi_4 = [3,5,7,8,6,4]$.
Conditions (1)--(4) are clearly satisfied.
For condition (5), look at any link in the network, and count the number of times that link is used in $\pi_1$ and $\pi_3$, and the number of times it is used in $\pi_2$ and $\pi_4$; every link is either unused in both pairs, used in exactly one path in both pairs, or used in both paths in each pair.
Now assume that $g_a \leq \epsilon \leq g_r / 2$, as in~\eqn{twoconsistency}, and choose $\hat{\Pi}$ to be all paths whose travel time is within $\epsilon$ of the shortest for its OD pair.
Show that there are no four paths $\pi_1$, $\ldots$, $\pi_4$ satisfying all of the conditions in the previous paragraph.
(Hint: argue by contradiction, and apply each of the conditions, using the assumptions about the acceptance and rejection gaps to bound each path's travel time relative to the shortest path travel time.)
\label{ex:twoconsistency}
\item \diff{38}.
\label{ex:dualdirectionderivation} The proof of Theorem~\ref{thm:proportionality} started from the entropy-maximizing Lagrangian~\eqn{entropylagrangian}, which Lagrangianized both the link flow constraints (with multipliers $\beta_{ij}$) and the OD matrix constraints (with multipliers $\gamma_{rs}$).
Alternatively, we can Lagrangianize only the link flow constraints, and replace $(h_\pi/d^{rs})$ with $h_\pi$ (why?), giving the equation 
\labeleqn{altlagrangian}{ \hat{\mc{L}}(\mb{x}, \bm{\beta}) = \sum_{\pi \in \Pi} h_\pi \log h_\pi + \sum_{(i,j) \in A} \beta_{ij} \myp{ \hat{x}_{ij} - \sum_{\pi \in \Pi} \delta_{ij}^\pi h_\pi } \,.}
Show that the gradient of this alternative Lagrangian with respect to $\bm{\beta}$ has components given by~\eqn{dualsearch}.
That is, show that
\[
\pdr{\hat{\mc{L}}}{\beta_{ij}} = x_{ij} - \hat{x}_{ij}
\,\]
where $x_{ij}$ is computed from the current path flows $h_\pi$.
\item \diff{59}.
\label{ex:dualstepsize} The dual algorithm step~\eqn{dualsearch} can be compactly written as $\bm{\beta} \leftarrow \bm{\beta} + \alpha \Delta \mb{\beta}$, where $\Delta \beta_{ij} = x_{ij} - \hat{x}_{ij}$.
Let $f(\alpha)$ denote the value of the alternative Lagrangian~\eqn{altlagrangian} after a step of size $\alpha$ is taken, and the new $\bm{\beta}$ and $\mb{h}$ values are calculated.
Newton's method can be used to find an $\alpha$ value which approximately minimizes $f(\alpha)$, maximizing entropy in the direction $\Delta \mb{\beta}$.
The Newton step is $\alpha = -f'(0)/f''(0)$.
\begin{enumerate}[(a)]
\item Show that $f'(0) = || \nabla_{\bm{\beta}} {\hat{\mc{L}}} ||^2$.
(See Exercise~\ref{ex:dualdirectionderivation}.)
\item Show that $f''(0) = \nabla^T_{\bm{\beta}} {\hat{\mc{L}}} H_{\bm{\beta}} {\hat{\mc{L}}} \nabla_{\bm{\beta}} {\hat{\mc{L}}}$, where $H_{\bm{\beta}}$ is the Hessian of the alternate Lagrangian with respect to $\bm{\beta}$.
\item Show that 
\[ \frac{f'(0)}{f''(0)} =
    \frac{ \ds \sum_{(i,j) \in A} (x_{ij} - \hat{x}_{ij})^2}
    {\ds \sum_{(i,j) \in A} \sum_{(k,\ell) \in A}
        (x_{ij} - \hat{x}_{ij}) (x_{kl} - \hat{x}_{k\ell})
        \sum_{\pi \in \Pi}
            h^\pi \delta_{ij}^\pi \delta_{k \ell}^\pi
    }
\,.\]
\end{enumerate}
\item \diff{30}.
\label{ex:paspathsubtleties} In the discussion surrounding Figure~\ref{fig:tapasintuition}, we argued that satisfying the equilibrium conditions around a ``spanning'' set of PASs (one for each block) was sufficient for establishing equilibrium on the entire network.
Consider the network in Figure~\ref{fig:paspathsubtleties}, where the demand from origin 1 to destination 4 is 100 vehicles, and there are two PASs: one between segments $[1,4]$ and $[1,2,3,4]$, and another between segments $[2,4]$ and $[2,3,4]$.
These are spanning, in the sense that by shifting flows between these two PASs we can obtain any feasible path flow solution from any other.
They also satisfy the equilibrium conditions: for the first PAS, because there is no flow on either segment\footnote{There is flow on link (1,2), but not on the entire segment $[1,2,3,4]$.}; for the second, because the travel times are equal on the two segments.
Yet the network is not at equilibrium, since $[1,4]$ is the only shortest path and it is unused.
Explain this apparent inconsistency.
\stevefig{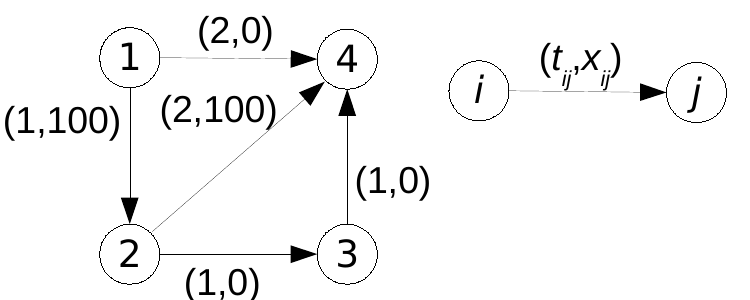}{Network for Exercise~\ref{ex:paspathsubtleties}.}{0.5\textwidth}
\item \diff{62}.
\label{ex:tapasbfs} Develop one or more algorithms to find ``short'' segments when generating a new PAS.
These methods should require a number of steps that grows at most linearly with network size.
\item \diff{22}.
\label{ex:removecycle} Show that the cycle-removing procedure described in the TAPAS algorithm maintains feasibility of the solution (flow conservation at each node, and non-negativity of link flows), and that the Beckmann function decreases strictly (assuming link performance functions are positive).
\item \diff{68}.
\label{ex:tapasproportionalapprox}  This exercise develops a technique for approximately solving equations~\eqn{tapasproportionalfeasible} and~\eqn{tapasproportionality}, better than the heuristic given in the text.
\begin{enumerate}
\item Define $\Delta h^r(\rho)$ to be the amount of flow that must be shifted from origin $r$'s flows on $\sigma_1$ to $\sigma_2$, to adjust the proportion $g^r(\sigma_1)^\zeta / (g^r(\sigma_1)^\zeta + g^r(\sigma_2)^\zeta)$ to be exactly $\rho$.
A negative value of this function indicates shifting flow in the reverse direction, from $\sigma_2$ to $\sigma_1$.
Show that this function is defined over $[0,1]$, and its range is $[\underline{\Delta},\overline{\Delta}]$, where $\underline{\Delta} = -\min_{(i,j) \in \sigma_2} x^r_{ij}$ and $\overline{\Delta} = \min_{(i,j) \in \sigma_1} x^r_{ij}$.
Furthermore show that $\Delta h^r(0) = \underline{\Delta}$ and $\Delta h^r(1) = \overline{\Delta}$.
\item Show that equation~\eqn{tapasproportionality} is satisfied if the flow shift $\Delta h^r(\rho)$ is applied to all relevant origins $r \in Z_\zeta$.
\item Therefore, it is enough to find a value of $\rho$ for which
$$U(\rho) \equiv \sum_{r \in Z_\zeta} \Delta h^r(\rho) = 0\,,$$
in order to satisfy~\eqn{tapasproportionalfeasible}.
The function $U$ is continuous and defined on the interval $[0,1]$.
Show that $U(0) \leq 0$ and $U(1) \geq 0$, ensuring that a zero exists in this interval.
\item Develop a quadratic approximation for $\Delta h^r(\rho)$ based on three known points: $\Delta h^r(0) = \underline{\Delta}$, $\Delta h^r(1) = \overline{\Delta}$, and $\Delta h^r(\rho_r) = 0$, where $\rho_r$ is the current proportion $g^r(\sigma_1)^\zeta / (g^r(\sigma_1)^\zeta + g^r(\sigma_2)^\zeta)$.
\item By summing these, develop a quadratic approximation for $U(\rho)$, and give an explicit formula for its root.
\end{enumerate}
\end{enumerate}

\chapter{Sensitivity Analysis and Applications}
\label{chp:sensitivityanalysis}

\index{static traffic assignment!sensitivity analysis|(}
This chapter shows how a sensitivity analysis can be conducted for the traffic assignment problem (TAP), identifying how the equilibrium assignment will change if the problem parameters (such as the OD matrix or link performance functions) are changed.
This type of analysis is useful in many ways: it can be used to determine the extent to which errors or uncertainty in the input data create errors in the output data.
It can be used as a component in so-called ``bilevel'' optimization problems, where we seek to optimize some objective function while enforcing that the traffic flows remain at equilibrium.
This occurs most often in the network design problem, where one must determine how to improve network links to reduce total costs, and in the OD matrix estimation problem, where one attempts to infer the OD matrix from link flows, or improve upon an existing estimate of the OD matrix.

After exploring the sensitivity analysis problem using the familiar Braess network, the first objective in the chapter is calculating derivatives of the equilibrium link flows with respect to elements in the OD matrix.
It turns out that this essentially amounts to solving another, easier, traffic assignment problem with different link performance functions and constraints.
The remainder of the chapter shows how these derivatives can be used in the network design and OD matrix estimation problems, which are classic transportation examples of bilevel programs.

\section{Sensitivity Analysis Preliminaries}

\index{static traffic assignment!sensitivity analysis!demand change|(}
Figure~\ref{fig:braessddnet} shows the Braess network.
When this network was first introduced, the demand between node 1 and node 4 was $d^{14} = 6$, and the equilibrium solution was found to be $x_{13} = x_{23} = x_{24} = 2$ and $x_{12} = x_{34} = 4$, with a travel time of 92 minutes on all three paths.
What if, instead, the demand $d^{14}$ took another value?  Figure~\ref{fig:braessdd} presents four plots showing how the equilibrium solution varies according to the demand level.
Panel (a) shows the flows $x_{12}$ and $x_{34}$, panel (b) shows the flows $x_{13}$ and $x_{24}$, panel (c) shows the flow $x_{23}$, and panel (d) shows the shortest path travel time between nodes 1 and 4, at the corresponding equilibrium solution.
You can check that when $d^{14} = 6$, the original equilibrium solution is shown in this figure.
\index{static traffic assignment!sensitivity analysis!demand change|)}

\stevefig{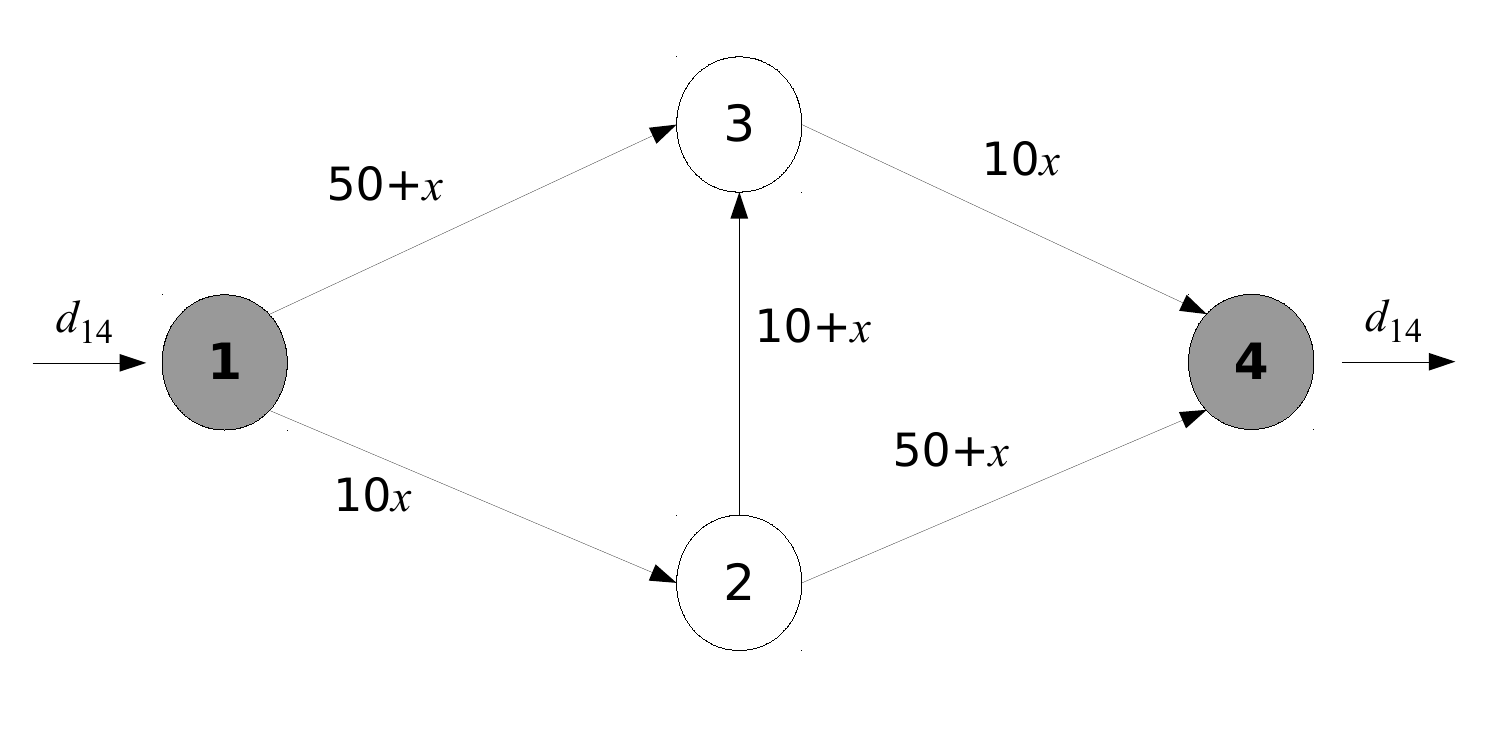}{Braess network with varying demand from 1 to 4.}{0.7\textwidth}
\stevefig{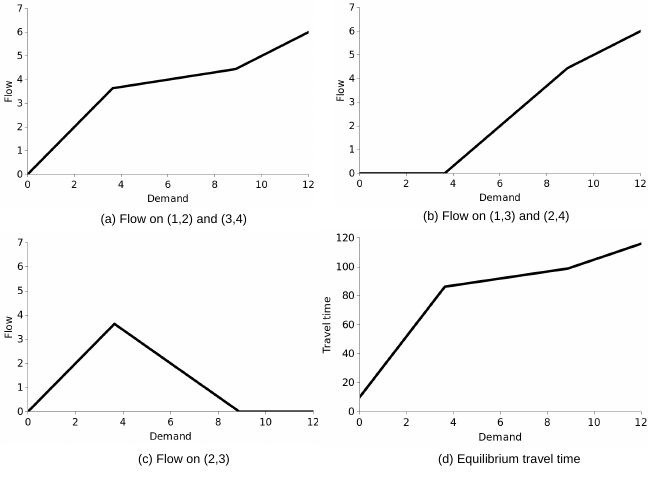}{Sensitivity analysis of the Braess network to $d^{14}$.}{\textwidth}

\index{static traffic assignment!sensitivity analysis!link change|(}
Instead of the OD matrix, we also could have changed the link performance functions in the network.
Now assume that $d_{14}$ is fixed at its original value of 6, but that the link performance function on link (2,3) can vary.
Let $t_{23}(x_{23}) = y + x_{23}$, where $y$\label{not:yndp} is the free-flow time, resulting in the network shown in Figure~\ref{fig:braessdxnet}.
In the base solution $y = 10$, but conceivably the ``free-flow time'' could be changed.
If the speed limit were increased, $y$ would be lower; if traffic calming were implemented, $y$ would be higher.
In an extreme case, if the link were closed entirely you could imagine $y$ takes an extremely large value, large enough that no traveler would use the path.
One can also effectively decrease $y$ by providing incentives for traveling on this link (a direct monetary payment, a discount at an affiliated retailer, etc.), and conceivably this incentive could be so large that $y$ is negative.
The resulting sensitivity analysis is provided in Figure~\ref{fig:braessdx}
\index{static traffic assignment!sensitivity analysis!link change|)}

\stevefig{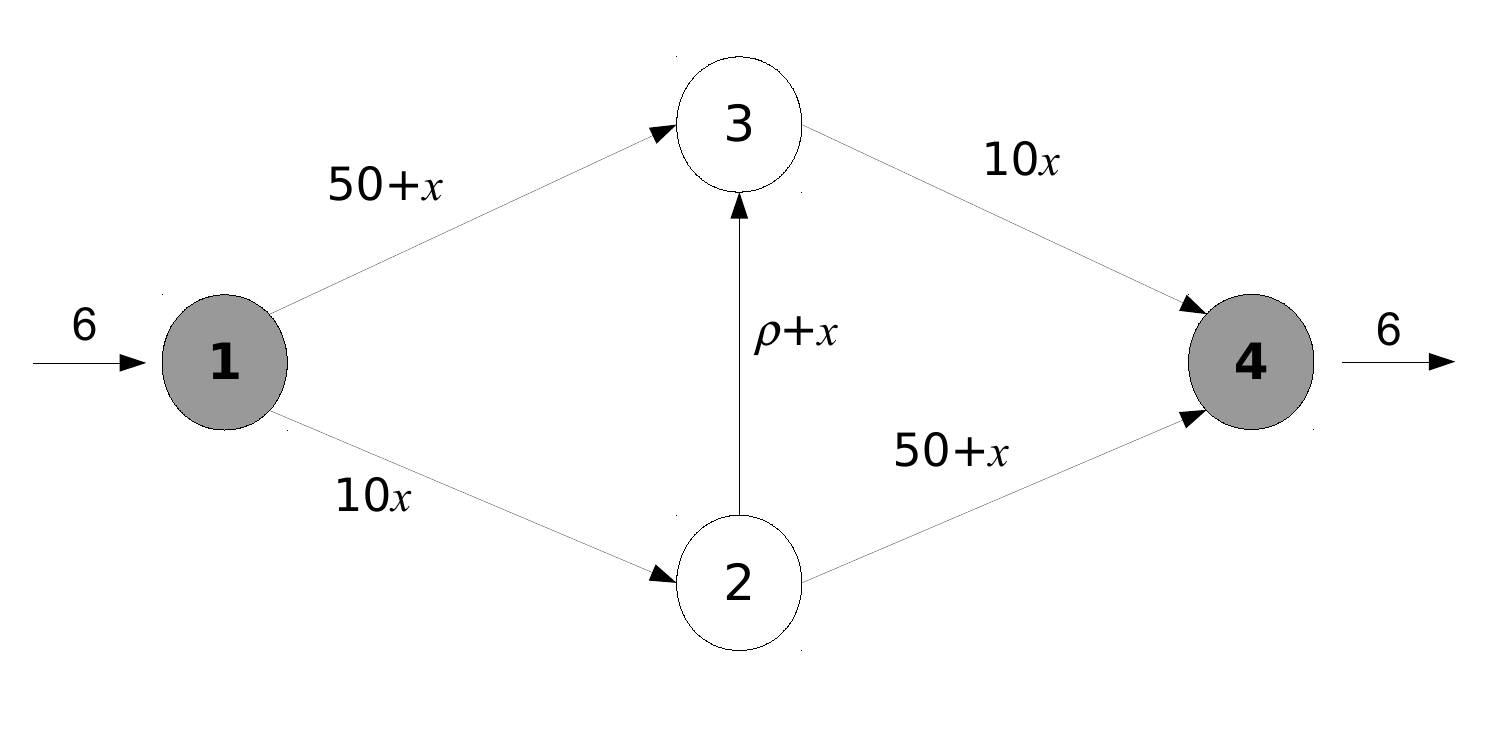}{Braess network with varying free-flow time on (2,3).}{0.7\textwidth}
\stevefig{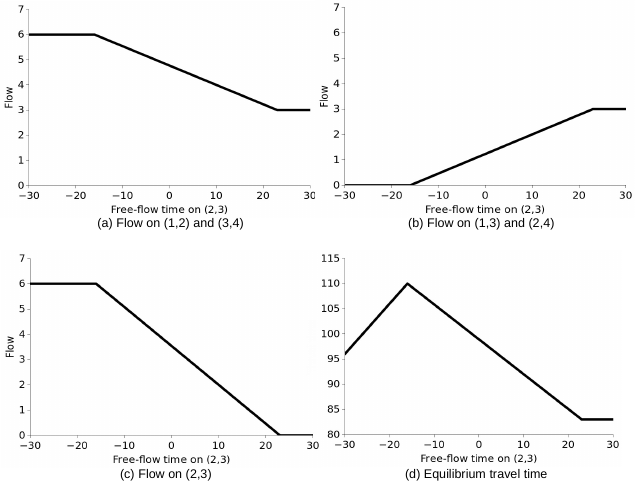}{Sensitivity analysis of the Braess network to $y$.}{\textwidth}

Examining the plots in Figure~\ref{fig:braessdd} and~\ref{fig:braessdx}, we see that the relationships between the equilibrium solution (link flows and travel times) and the demand or free-flow time are all piecewise linear.
Each ``piece'' of these piecewise linear functions corresponds to a particular subset of the paths being used --- for instance, in Figure~\ref{fig:braessdd}, when the demand is lowest, only the middle path is used.
When the demand is highest, only the two outer paths are used.
When the demand is at a moderate level, all three paths are used.
Within each of these regions, the relationship between the demand and the equilibrium solution is linear.
These pieces meet at so-called \emph{degenerate}\index{degenerate solution} solutions, where the equilibrium solution does not use all of the minimum travel-time paths.
(For instance, when $d_{14} = 40/11$ the equilibrium solution requires all drivers to be assigned to the middle path, even though all three have equal travel times.)

In general networks involving nonlinear link performance functions, these relationships cannot be expected to stay linear.
However, they are still defined by piecewise functions, with each piece corresponding to a certain set of paths being used, and with the pieces meeting at degenerate solutions.
The goal of the sensitivity analyses in this chapter is to identify derivatives of the equilibrium solution (link flows and travel times) at a given point.
For these derivatives to be well-defined, we therefore assume that \emph{the point at which our sensitivity analysis occurs is not degenerate.}  That is, all minimum-travel time paths have positive flow.\index{degenerate solution!impact on sensitivity analysis}
This assumption is not too restrictive, because there are only a finite number of degenerate points; for instance, if we pick the demand value at random, the probability of ending up at a degenerate point is zero.

This sensitivity analysis is still local, because the information provided by a derivative grows smaller as we move farther away from the point where the derivative is taken.
For a piecewise function, the derivative provides no information whatsoever for pieces other than the one where the derivative was taken.

In this chapter, we show how this kind of sensitivity analysis can be used in two different ways.
In the network design problem, this type of sensitivity analysis can be used to determine where network investments are most valuable.
In Figure~\ref{fig:braessdx}, the fact that the equilibrium travel time increases when $y$ decreases (around the base solution $y = 10$) highlights the Braess paradox: investing money to improve this link will actually increase travel times throughout the network.
If we were to conduct a similar analysis for other links in the network, we would see that the equilibrium travel time would decrease with improvements to the link.
In the OD matrix estimation problem, we can use this sensitivity analysis to help calibrate an OD matrix to given conditions.

\section{Calculating Sensitivities}

This section derives sensitivity formulas showing the derivative of the equilibrium link flows and travel times with respect to two parameters: (1) a change in an entry of the OD matrix $d^{rs}$, and (2) a change to a parameter in the link performance functions (such as the free-flow time or capacity in a BPR function).
In this section, assume that we are given some initial OD matrix or link performance functions, and the corresponding equilibrium solution.
For our purposes, it will be most convenient if this equilibrium solution is expressed in bush-based form,\index{static traffic assignment!bush-based solution} that is, with vectors $\mb{\hat{x}}^r$\label{not:xhatijr} showing the flow on each link corresponding to each origin $r$.\index{degenerate solution!impact on sensitivity analysis|(}
In this case, the non-degenerate condition requires that unused links are not part of the equilibrium bushes, that is, if $\hat{x}_{ij}^r = 0$ for any link $(i,j)$ and any origin $r$, then link $(i,j)$ does not correspond to any shortest path starting from node $r$ --- in terms of the distance labels $L_i^r$, we have $\hat{x}_{ij}^r > 0$ if and only if $L_j^r = L_i^r + t_{ij}$.
Let $\mc{B}^r = \{ (i,j) \in A : \hat{x}_{ij}^r > 0 \}$ denote the equilibrium bush for origin $r$.\index{bush!equilibrium bush}

The non-degeneracy assumption is important, because one can show that \emph{if the change to the OD matrix or link performance functions is small and the original equilibrium solution is non-degenerate, all of the equilibrium bushes remain unchanged.}  Equivalently, even after drivers shift flows to find the new equilibrium, the set of used paths will remain the same as it was before.
Furthermore, one can show that the equilibrium solution is differentiable, and the derivatives of the equilibrium link flows or travel times with respect to values in the OD matrix or link performance function parameters can be interpreted as the sensitivities of the equilibrium solution.\index{degenerate solution!impact on sensitivity analysis|)}

There are several ways to calculate the values of these derivatives: historically, the first researchers used matrix-based formulas, and subsequent researchers generalized these formulas using results from the theory of variational inequalities.
We adopt a different approach, using the bush-based solution representation, because it leads to an easy solution method and is fairly straightforward.
This approach is based on the fact that the equilibrium solution (travel times $t_{ij}$ and bush flows $x_{ij}^r$) must satisfy the following equations for each origin $r$:
\begin{align}
\label{eqn:bushreducedcost}  L^r_j - L^r_i - t_{ij}(\hat{x}_{ij}) &= 0 & \forall (i,j) \in \mc{B}^r \\
\label{eqn:bushorigin} L^r_r &= 0 & \\
\label{eqn:bushflowc} \sum_{(h,i) \in \Gamma^{-1}(i)} \hat{x}^r_{hi} - \sum_{(i,j) \in \Gamma(i)} \hat{x}^r_{ij} &= d^{ri} & \forall i \in N \backslash \{r\} \\
\label{eqn:bushflowcorigin} \sum_{(h,r) \in \Gamma^{-1}(r)} \hat{x}^r_{hr} - \sum_{(r,j) \in \Gamma(r)} \hat{x}^r_{rj} &= -\sum_{s \in Z} d^{rs} &  \\
\label{eqn:bushflownnonbush} \hat{x}^r_{ij} &= 0 & \forall (i,j) \notin \mc{B}^r
\end{align}
Equations~\eqn{bushreducedcost}--\eqn{bushorigin} reflect the equilibrium condition, and equations~\eqn{bushflowc}--\eqn{bushflowcorigin} represent flow conservation.
The number of equations for each origin is no more than the sum of the number of links and nodes in the network.

Furthermore, these conditions must remain true even as the problem data (OD matrix and link performance functions) are perturbed.
Since derivatives of the equilibrium solution exist under the non-degeneracy assumption, we can differentiate equations~\eqn{bushreducedcost}--\eqn{bushflowcorigin} to identify the relationships which must hold true among these derivatives.
For brevity, in this chapter we use $\xi^r_{ij}$\label{not:xirij} to denote the derivative of $x^r_{ij}$, and $\Lambda^r_i$\label{not:Lambdari} to denote the derivative of $L^r_i$.
These derivatives are taken with respect to either an OD matrix entry or a link performance function parameters, as described separately below.

\subsection{Changes to the OD matrix}
\label{sec:odsa}

\index{static traffic assignment!sensitivity analysis!demand change|(}
Assume first that we change a single entry in the OD matrix corresponding to origin $\hat{r}$\label{not:rhat} and destination $\hat{s}$,\label{not:shat} so $\xi^r_{ij} = dx^r_{ij}/dd^{\hat{r}\hat{s}}$ and $\Lambda_i = dL_i/dd^{\hat{r}\hat{s}}$.
Then differentiating each of equations~\eqn{bushreducedcost}--\eqn{bushflowcorigin} with respect to $d^{\hat{r}\hat{s}}$ gives the following equations for each $r \in Z$:
\begin{align}
 \label{eqn:ddbushreducedcost}  \Lambda^r_j - \Lambda^r_i - t'_{ij} \sum_{r' \in Z} \xi^{r'}_{ij} &= 0 & \forall (i,j) \in \mc{B}^{r} \\
 \label{eqn:ddbushroot} \Lambda^r_r &= 0 \\
 \label{eqn:ddbushflowc} \sum_{(h,i) \in \Gamma^{-1}(i)} \xi^r_{hi} - \sum_{(i,j) \in \Gamma(i)} \xi^r_{ij} &= \begin{cases}
                                                                                                                1 & \mbox{ if } r = \hat{r} \mbox{ and } i = \hat{s} \\
                                                                                                                0 & \mbox{ otherwise}
                                                                                                               \end{cases}
                                                                                                               & \forall i \in N, i \neq r\\
 \label{eqn:ddbushflowcorigin} \sum_{(h,r) \in \Gamma^{-1}(r)} \xi^r_{hr} - \sum_{(r,j) \in \Gamma(r)} \xi^r_{rj} &= \begin{cases}
                                                                                                                -1 & \mbox{ if } r = \hat{r} \\
                                                                                                                0 & \mbox{ otherwise}
                                                                                                               \end{cases} \\
 \label{eqn:ddbushflownnonbush} \xi^r_{ij} &= 0 & \forall (i,j) \notin \mc{B}^r
\end{align}
where $t'_{ij}$ is the derivative of the link performance function, evaluated at the current equilibrium solution $\mb{\hat{x}}$ (and thus treated as a constant in these equations).
Equations~\eqn{ddbushreducedcost} enforce the fact that the equilibrium bushes must remain the same.
That is, the shortest path labels $L_i$ and travel times $t_{ij}$ must change in such a way that every link on the bush is part of a minimum travel time path to its head node.
Equations~\eqn{ddbushflowc} and~\eqn{ddbushflowcorigin} enforce flow conservation.
For all bushes except for $\hat{r}$, the total flow from the origin to each destination is the same, so flow is allowed to redistribute among the bush links, but the flows starting or ending at a node cannot change.
For the bush corresponding to $\hat{r}$, a unit increase in demand from $\hat{r}$ to $\hat{s}$ must be reflected by an additional vehicle leaving $\hat{r}$ and an additional vehicle arriving at $\hat{s}$.

All together, the system of equations~\eqn{ddbushreducedcost}--\eqn{ddbushflowcorigin} involves variables $\Lambda^r_i$ for each origin $r$ and node $i$, and $\xi^r_{ij}$ for each origin $r$ and link $(i,j)$.
Furthermore, for each origin, it contains an equation for each link and each node.
Therefore, this linear system of equations can be solved to obtain the sensitivity values.\footnote{A careful reader will note that one of the flow conservation equations for each origin is redundant, but this is of no consequence to what follows.}

However, there is an easier way to solve for $\Lambda^r_i$ and $\xi^r_{ij}$.
Using the techniques in Section~\ref{sec:convexoptimization}, you can show that the equations~\eqn{ddbushreducedcost}--\eqn{ddbushflowcorigin} are exactly the optimality conditions\index{optimization!convex!optimality conditions} to the following minimization problem:
\begin{multline}
 \label{eqn:ddsaobjective} \min_{\bm{\xi}^r, \bm{\Lambda}^r} \quad \frac{1}{2} \sum_{(i,j) \in A} t'_{ij} \myp{ \sum_{r \in Z} \xi^r_{ij}}^2 \\+ \sum_{r \in Z} \sum_{i \in N} \Lambda_i^r \myp{ \sum_{(h,i) \in \Gamma^{-1}(i)} \xi^r_{hi} - \sum_{(i,j) \in \Gamma(i)} \xi^r_{ij} - \Delta_{ri}} 
 \end{multline}
\begin{align}
 \mr{s.t.} \quad  & \Lambda^r_r = 0 \qquad \forall r \in Z \\
                  & \xi_{ij}^r = 0 \qquad \forall (i,j) \notin \mc{B}^r  \label{eqn:ddsafc}
\end{align}
where $\Delta_{ri}$\label{not:Deltari} represents the right-hand side of equation~\eqn{ddbushflowc} or~\eqn{ddbushflowcorigin}, that is, $\Delta_{\hat{r} \hat{s}} = 1$, $\Delta_{\hat{r} \hat{r}} = -1$, and $\Delta_{ri} = 0$ otherwise.

This optimization problem can be put in a more convenient form by interpreting $\Lambda^r_i$ as the Lagrange multiplier\index{Lagrange multipliers!applications} for the flow conservation equation corresponding to node $i$ in bush $r$, and the second term in~\eqn{ddsaobjective} as the Lagrangianization of this equation.
Furthermore, defining $\xi_{ij} = \sum_{r \in Z} \xi^r_{ij}$, the optimization problem can be recast in the following equivalent form:
\begin{align}
  \label{eqn:ddfinalobjective} \min_{\bm{\xi^r}} \quad & \int_0^{\xi_{ij}} t'_{ij} \xi~d\xi & \\
  \mr{s.t.} \qquad & \sum_{(h,i) \in \Gamma^{-1}(i)} \xi^r_{hi} - \sum_{(i,j) \in \Gamma(i)} \xi^r_{ij} = \Delta_{ri}  & \forall r \in Z, i \in N \\
                   & \xi_{ij}^r = 0 & \forall r \in Z, (i,j) \notin \mc{B}^r 
\end{align}
This is essentially a traffic assignment problem (TAP) in bush-based form (see Chapter~\ref{chp:solutionalgorithms}), with the following changes:
\begin{itemize}
\item The original link performance functions have been replaced by \emph{linear} link performance functions with slope equal to the derivative of the original link performance function at the original equilibrium solution.
(Remember that $t'_{ij}$ is a constant, the value of the link performance function derivative at the equilibrium solution.)
\item The equilibrium bushes for each origin are \emph{fixed} at the bushes for the original equilibrium solution.
\item The only entry in the OD matrix is one unit of demand from $\hat{r}$ to $\hat{s}$.
\item There are no non-negativity conditions.
This is because the solution variables $\xi_{ij}$ represent changes in the original link flows, and it is possible for these changes to be negative as well as positive (cf.\ Figure~\ref{fig:braessdd}).
\end{itemize}
If you have access to an implementation of a bush-based algorithm\index{static traffic assignment!algorithms!bush-based} for solving TAP, it is easy to modify the program to take account of these distinctions, and to find the link flow sensitivities $\xi^r_{ij}$.
From here, the values of $\Lambda^r_i$ can be found by solving a shortest path problem with link travel times $\frac{dt_{ij}}{dx_{ij}} \xi_{ij}$.

As a demonstration, we use the Braess network of Figure~\ref{fig:braessddnet}, working around the base demand $d_{14} = 6$ and base equilibrium solution $\hat{x}_{12} = \hat{x}_{34} = 4$, $\hat{x}_{13} = \hat{x}_{23} = \hat{x}_{24} = 2$.
At this level of demand, all paths are used and the equilibrium bush contains all of the links in the original network.
Furthermore, at this equilibrium solution the derivatives of the link performance functions are $t'_{12} = t'_{34} = 10$ and $t'_{13} = t'_{23} = t'_{24} = 1$.
The linear link performance functions based on these derivatives are shown in Figure~\ref{fig:braessddsensitivitynet}.
The OD matrix is replaced by a single unit of flow traveling from 1 (our $\hat{r}$) to 4 (our $\hat{s})$.

Solving the problem without non-negativity constraints produces the $\xi_{ij}$ values shown in Figure~\ref{fig:braessddsensitivitydxsoln}.
Substituting these into the link performance functions in Figure~\ref{fig:braessddsensitivitynet} shows that the equilibrium is satisfied: all paths have equal travel times of $31/13$.
This is exactly the slope of the piece of the equilibrium travel time (Figure~\ref{fig:braessdd}d) around $d_{14} = 6$, that is, the equilibrium travel time in the sensitivity problem gives the derivative of the equilibrium travel time in the original network.

\stevefig{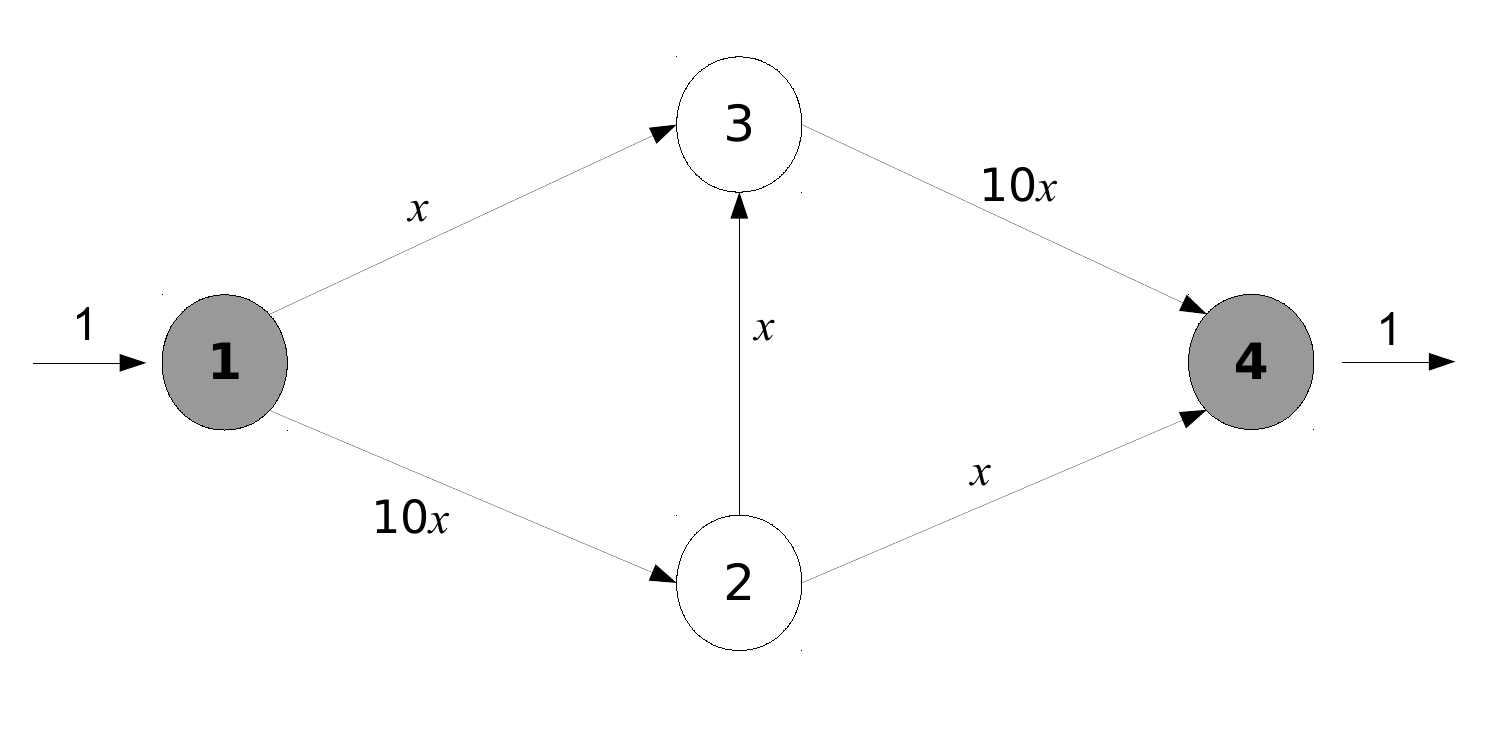}{Modified traffic assignment problem for sensitivity to $d_{14}$}{0.75\textwidth}
\stevefig{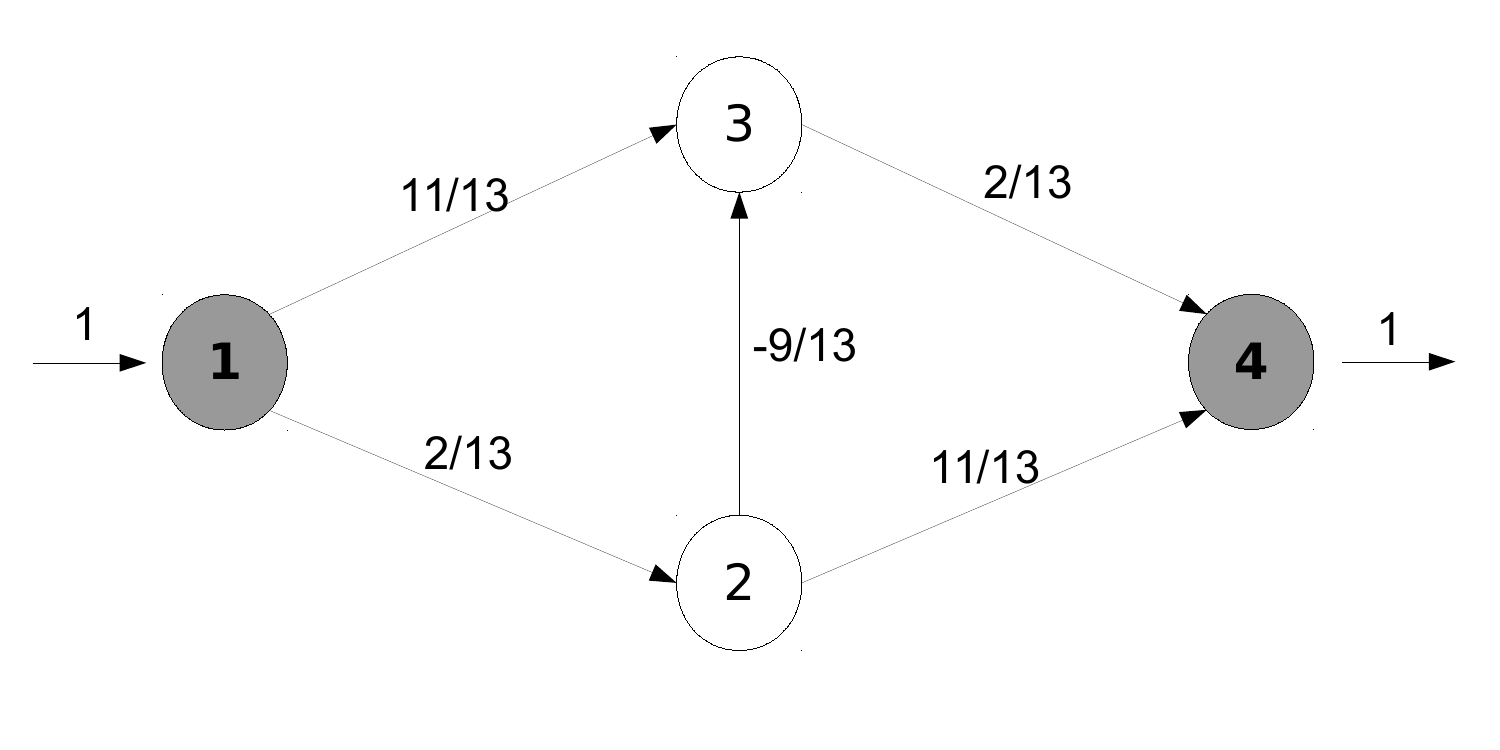}{Solution to traffic assignment problem for sensitivity to $d_{14}$}{0.75\textwidth}
\index{static traffic assignment!sensitivity analysis!demand change|)}

\subsection{Changes to a link performance function}
\label{sec:linksa}

\index{static traffic assignment!sensitivity analysis!link change|(}
Now assume that we change a parameter $y$ (which may represent the free-flow time, capacity, or any other parameter) in the link performance function corresponding to link $(i,j)$, so $\xi^r_{ij} = dx^r_{ij}/dy$ and $\Lambda^r_i = dL^r_i/dy$.
This means that the link performance function $t_{ij}$ now depends on both its flow $x_{ij}$, and the parameter $y$, so we write $t_{ij}(x,y)$.
We will write $t'_{ij,x}$\label{not:tijx} to mean the partial derivative with respect to link flow\footnote{This is what we wrote as $t'_{ij}$ in the previous section, when the link performance function only depended on $x_{ij}$.}, and $t'_{ij,y}$\label{not:tijy} to mean the partial derivative with respect to the improvement parameter $y$.
Then differentiating each of equations~\eqn{bushreducedcost}--\eqn{bushflowcorigin} with respect to $y$ gives the following equations for each $r \in Z$:
\begin{align}
\label{eqn:cdbushreducedcost}  \Lambda^r_j - \Lambda^r_i - t'_{ij,x} \sum_{r' \in Z} \xi^{r'}_{ij} - t'_{ij,y} &= 0 & \forall (i,j) \in \mc{B}^{r} \\
\label{eqn:cdbushroot} \Lambda^r_r &= 0 \\
\label{eqn:cdbushflowc} \sum_{(h,i) \in \Gamma^{-1}(i)} \xi^r_{hi} - \sum_{(i,j) \in \Gamma(i)} \xi^r_{ij} &= 0
                                                                                                               & \forall i \in N, i \neq r \\
\label{eqn:cdbushflowcorigin} \sum_{(h,r) \in \Gamma^{-1}(r)} \xi^r_{hr} - \sum_{(r,j) \in \Gamma(r)} \xi^r_{rj} &= 0 \\
\label{eqn:cdbushflownnonbush} \xi^r_{ij} &= 0 & \forall (i,j) \notin \mc{B}^r                                                                                                               
\end{align}
where as before, $t'_{ij,x}$ and $t'_{ij,y}$ are evaluated at the current, equilibrium solution $\mb{\hat{x}}$.
Equations~\eqn{cdbushreducedcost} enforce the fact that the equilibrium bushes must remain the same, taking into account both the change in travel time on $(i,j)$ due to the change in its link performance function as well as changes in all links' travel times from travelers shifting paths.
Equations~\eqn{cdbushflowc} and~\eqn{cdbushflowcorigin} enforce flow conservation.
These equations are simpler than for the case of a change to the OD matrix, because the total number of vehicles on the network remains the same, and these vehicles can only shift amongst the paths in the bush.
There is no change in the flow originating or terminating at any node, and the $\xi$ variables must form a circulation.

As with a change in an OD matrix entry, the system of equations~\eqn{cdbushreducedcost}--\eqn{cdbushflowcorigin} is a linear system involving, for each origin, variables for each node and link.
Repeating the same steps as before, this system of equations can be seen as the optimality conditions for the following optimization problem: 
\begin{align}
\label{eqn:cdfinalobjective} \min_{\bm{\xi^r}} \quad & \int_0^{\xi_{ij}} \myp{ t'_{ij,x} \xi + t'_{ij,y} } ~d\xi & \\
  \mr{s.t.} \qquad & \sum_{(h,i) \in \Gamma^{-1}(i)} \xi^r_{hi} - \sum_{(i,j) \in \Gamma(i)} \xi^r_{ij} = 0 & \forall r \in Z, i \in N \\
                   & \xi_{ij}^r = 0 & \forall r \in Z, (i,j) \notin \mc{B}^r 
\end{align}
This is essentially a traffic assignment problem in bush-based form, with the following changes:
\begin{itemize}
\item The original link performance functions have been replaced by \emph{affine} link performance functions with slope equal to the derivative of the original link performance function at the original equilibrium solution, and intercept equal to the derivative of the link performance function with respect to the parameter $y$.
\item The equilibrium bushes for each origin are \emph{fixed} at the bushes for the original equilibrium solution.
\item All entries in the OD matrix are zero.
\item There are no non-negativity conditions.
\end{itemize}

As a demonstration, we use the Braess network of Figure~\ref{fig:braessdxnet}, working around the base free-flow time $y = 10$ and base equilibrium solution $\hat{x}_{12} = \hat{x}_{34} = 4$, $\hat{x}_{13} = \hat{x}_{23} = \hat{x}_{24} = 2$.
That is, we replace the link performance function $t_{23}(x_{23}) = 10 + x_{23}$ with the function $t_{23}(x_{23},y) = y + x_{23}$ and see what happens when $y$ varies.
At the base value $y = 10$, all paths are used and the equilibrium bush contains all of the links in the original network.
Furthermore, at this equilibrium solution the derivatives of the link performance functions with respect to flows are $t'_{12,x} = t'_{34,x} = 10$ and $t'_{13,x} = t'_{23,x} = t'_{24,x} = 1$.
For link (2,3), we add the constant term $t'_{23,y} = 1$.
The link performance functions based on these derivatives are shown in Figure~\ref{fig:braessdxsensitivitynet}.
The OD matrix is set equal to zero, since there is no change in the total demand through the network.

Solving the problem without non-negativity constraints produces the $\xi_{ij}$ values shown in Figure~\ref{fig:braessdxsensitivitydxsoln}.
Substituting these into the link performance functions in Figure~\ref{fig:braessdxsensitivitynet} shows that the equilibrium is satisfied: all paths have equal travel times of $-9/13$.
This is the slope of the piece of the equilibrium travel time (Figure~\ref{fig:braessdx}d) around $y = 10$, that is, the equilibrium travel time in the sensitivity problem gives the derivative of the equilibrium travel time in the original network.

\stevefig{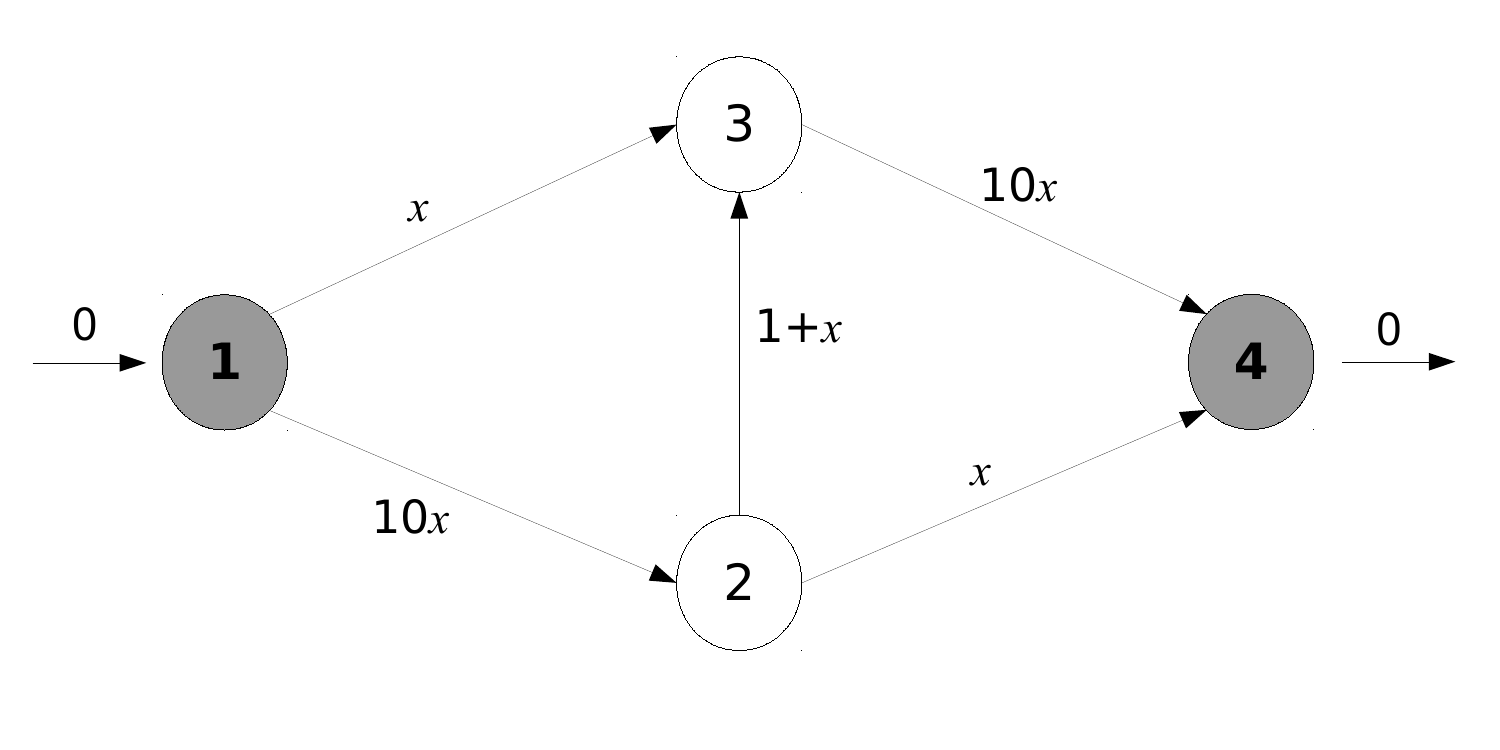}{Modified traffic assignment problem for sensitivity to $y_{23}$}{0.75\textwidth}
\stevefig{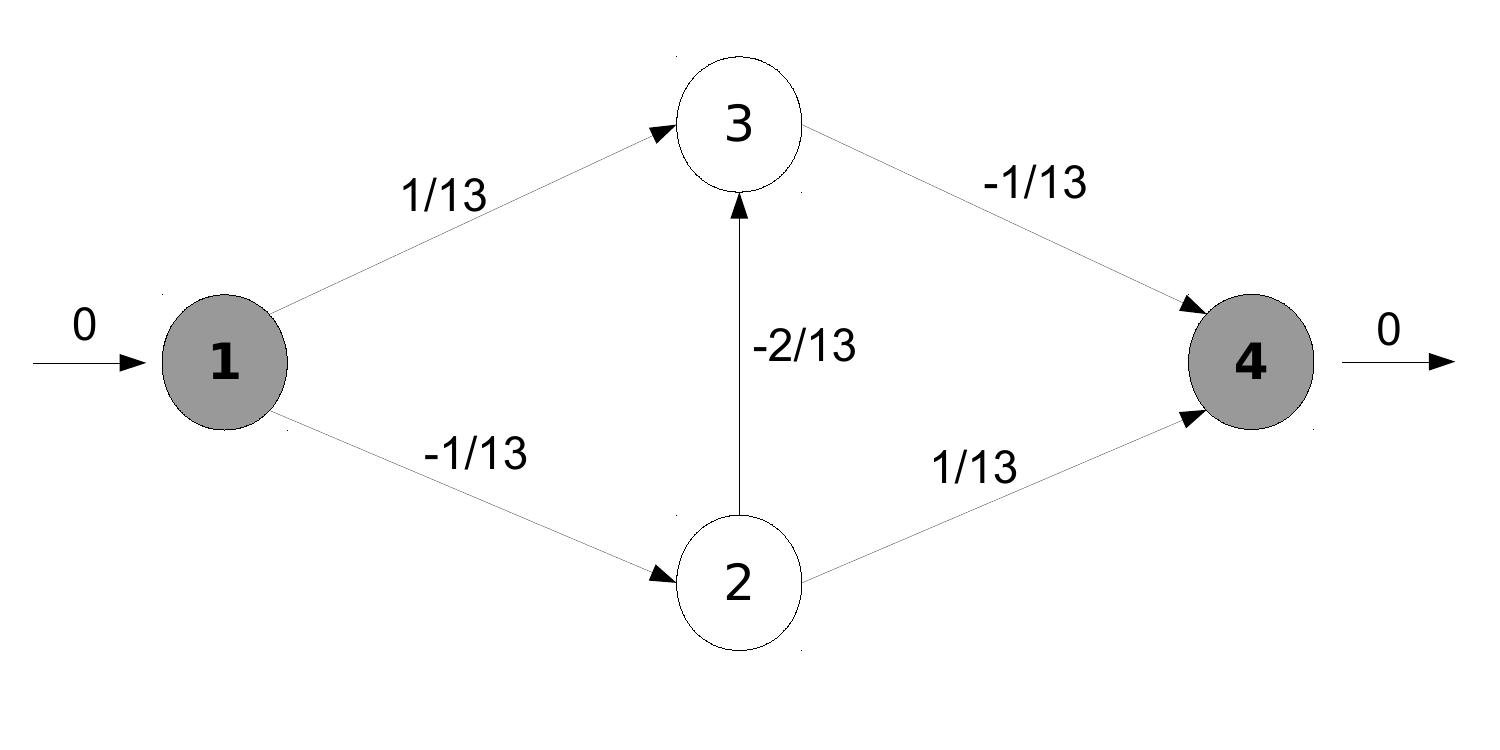}{Solution to traffic assignment problem for sensitivity to $y_{23}$}{0.75\textwidth}
\index{static traffic assignment!sensitivity analysis!link change|)}

\section{Network Design Problem}
\label{sec:ndp}

\index{optimization!bilevel|(}
\index{mathematical program with equilibrium constraints|see {optimization, bilevel}}
\index{network design|(}
In the network design problem, one must determine how best to spend funds on improving links in the transportation network.
This is a challenging optimization problem, and in cases of any practical interest one cannot hope to identify a globally optimal investment policy.
This is mainly because \emph{in transportation systems, the planner cannot compel travelers to choose a particular path.}  Instead, after any improvement is made to the links, flows will redistribute according to the principle of user equilibrium.\index{user equilibrium}
The goal is to find the best investment policy, knowing and anticipating how travelers will respond once the network has been changed.

Specifically, assume that the link performance functions for each link $(i,j)$ now depend on the amount of money $y_{ij}$ invested in that link (perhaps increasing its capacity through widening, or decreasing its free-flow time) as well as on the flow $x_{ij}$.
One example of such a link performance function is
\labeleqn{improvedbpr}{t_{ij}(x_{ij}, y_{ij}) = t^0_{ij} \myp{ 1 + \alpha \myp{ \frac{x_{ij}}{u_{ij} + K_{ij} y_{ij}} }^\beta } }
where $K_{ij}$\label{not:Kij} represents the capacity improvement if a single unit of money is invested on it.

A planning agency may have many different objectives and constraints when determining how to improve a network.
This section develops and explores one specific variation of the network design problem, but there are many other variations which have been proposed in the literature.
The version presented here is a fairly standard one, which can be extended in a number of different ways.
In this variation, the objective of the planning agency is to minimize the total cost, given by the sum of total system travel time $TSTT$\index{total system travel time} (converted to monetary units by a conversion factor $\Theta$,\label{not:Theta} which also reflects duration of the analysis horizon and discounting) and the costs of the network improvements themselves.
The optimization problem is
\begin{align}
\label{eqn:ndptstt}  \min_{\mb{x}, \mb{y}}   & f(\mb{x}, \mb{y}) = \Theta \sum_{(i,j) \in A} x_{ij} t_{ij}(x_{ij}, y_{ij}) + \sum_{(i,j) \in A} y_{ij} & \\
\label{eqn:ndpeqm}   \mr{s.t.} \qquad        & \mb{x} \in \arg \min_{\mb{x} \in X} \sum_{(i,j) \in A} \int_0^{x_{ij}} t_{ij}(x, y_{ij})~dx & \\
\label{eqn:ndpyfeas}                         & y_{ij} \geq 0 & \forall (i,j) \in A
\end{align}
Most of this problem is familiar: the objective~\eqn{ndptstt} is to minimize the sum of total system travel time (converted to units of money) and construction cost, and constraint~\eqn{ndpyfeas} requires that money can only be spent on links (not ``recovered'' from them with a negative $y_{ij}$ value).
Also note that since $\mb{y} = \mb{0}$ is a feasible solution (corresponding to the ``do-nothing'' alternative), in the optimal solution to this problem the cost savings (in the form of reduced $TSTT$) must at least be equal to the construction costs, guaranteeing that the optimal investment policy has greater benefit than cost.
The key equation here is~\eqn{ndpeqm}, which requires that the link flows $\mb{x}$ satisfy the principle of user equilibrium by minimizing the Beckmann function.
In other words, \emph{one of the constraints of the network design problem is itself an optimization problem}.
This is why the network design problem is called a bilevel program.
This type of problem is also known as a mathematical program with equilibrium constraints.
This class of problems is extremely difficult to solve, because the feasible region is typically nonconvex.

To see why, consider two feasible solutions $(\mb{x^1}, \mb{y^1})$ and $(\mb{x^2}, \mb{y^2})$ to the network design problem.
The link flows $\mb{x^1}$ are the equilibrium link flows under investment policy $\mb{y^1}$, and link flows $\mb{x^2}$ are the equilibrium link flows under investment policy $\mb{y^2}$.
If the feasible region were a convex set, then any weighted average of these two solutions would themselves be feasible.
Investment policy $\frac{1}{2} \mb{y^1} + \frac{1}{2} \mb{y^2}$ still satisfies all the constraints on $\mb{y}$ (all link investments are nonnegative).
However, the equilibrium link flows under this policy cannot be expected to be the average of $\mb{x^1}$ and $\mb{x^2}$, because the influence of $y_{ij}$ on $t_{ij}$ can be nonlinear and the sets of paths which are used in $\mb{x^1}$ and $\mb{x^2}$ can be completely different.
In other words, the equilibrium link flows after averaging two investment policies need not be the average of the equilibrium link flows under those two policies separately.

Unfortunately, solving optimization problems with nonconvex feasible regions is a very difficult task.
Therefore, solution methods for the network design problem are almost entirely heuristic in nature.
These heuristics can take many forms; one popular approach is to adapt a metaheuristic method, such as those discussed in Appendix~\ref{sec:metaheuristics}.
\index{metaheuristic}
\index{optimization!metaheuristic|see {metaheuristic}}

Another approach is to develop a more tailored heuristic based on specific insights about the network design problem.
This approach, being more educational, is adopted here.
Specifically, we can use the sensitivity analysis from the previous sections to identify derivatives of the objective function $f$ with respect to each link investment $y_{ij}$, and use this to move in a direction which reduces total cost.

Specifically, notice that constraint~\eqn{ndpeqm} actually makes $\mb{x}$ a function of $\mb{y}$, since the solution to the user equilibrium problem is unique in link flows.
That is, the investment policy $\mb{y}$ determines the equilibrium link flows $\mb{x}$ exactly.
So, the objective function can be made a function of $\mb{y}$ alone, written $f(\mb{x}(\mb{y}), \mb{y})$.
The derivative of this function with respect to an improvement on any link is then
\labeleqn{ndpder0}{\pdr{f}{y_{ij}} = \Theta \sum_{(k,\ell) \in A} \pdr{f}{x_{k\ell}} \pdr{x_{k\ell}}{y_{ij}} + 1 \,,}
or, substituting the derivative of~\eqn{ndptstt} with respect to each link flow, 
\labeleqn{ndpder}{\pdr{f}{y_{ij}} = \Theta \myc{ \sum_{(k,\ell) \in A} \pdr{x_{k\ell}}{y_{ij}} \myp{ t_{k \ell}(x_{k\ell}, y_{k\ell}) + x_{k\ell} \pdr{t_{k \ell}}{x_{k\ell}}(x_{k\ell}, y_{k\ell})  } + x_{ij} \pdr{t_{ij}}{y_{ij}} } + 1\,.}
In turn, the partial derivatives $\pdr{x_{k\ell}}{y_{ij}}$ can be identified using the technique of Section~\ref{sec:linksa} as the marginal changes in link flows throughout the network when the link performance function of $(i,j)$ is perturbed.

The vector of all the derivatives~\eqn{ndpder} forms the gradient of $f$ with respect to $\mb{y}$.
This gradient is the direction of steepest \emph{ascent}, that is, the direction in which $f$ is increasing fastest.
Since we are solving a minimization problem, we should move in the opposite direction.
Taking such a step, and ensuring feasibility, gives the updating equation
\labeleqn{ndpupdate}{\mb{y} \leftarrow \mys{\mb{y} - \mu \nabla_\mb{y} f}^+ \,,}
where $\mu$ is a step size to be determined, and the $[\cdot]^+$ operation is applied to each component of the vector.
This suggests the following algorithm:
\begin{enumerate}
\item Initialize $\mb{y} \leftarrow \mb{0}$.
\item Calculate the link flows $\mb{x}(\mb{y})$ by solving the traffic assignment problem with link performance functions $\mb{t}(\mb{x}, \mb{y})$.
\item For each link $(i,j)$ determine $\pdr{f}{y_{ij}}$ by solving the sensitivity problem described in Section~\ref{sec:linksa} and using~\eqn{ndpder}.
\item Update $\mb{y}$ using~\eqn{ndpupdate} for a suitable step size $\mu$.
\item Test for convergence, and return to step 2 if not converged.
\end{enumerate}
Two questions are how $\mu$ should be chosen in step 4, and how convergence should be tested in step 5.
The difficulty in step 4 is that the derivatives provided by a sensitivity analysis are only local, and in particular if $\mu$ is large enough that~\eqn{ndpupdate} changes the set of used paths, the derivative information is meaningless.
However, if $\mu$ is small enough one will see a decrease in the objective function if at all possible.
So, one could start by testing a sequence of $\mu$ values (say, $1, 1/2, 1/4, \ldots$), evaluating the resulting $\mb{y}$ values, $\mb{x}$ values, and $f$, stopping as soon as $f$ decreases from its current value.
(Note that this is a fairly computationally intensive process, since the traffic assignment problem must be solved for each $\mu$ to get the appropriate $\mb{x}$ value.)  Other options include using a stricter stopping criterion such as the Armijo rule\index{Armijo rule} (Appendix~\ref{sec:unconstrainedstepsize}), which would ensure that the decrease is ``sufficiently large''; or using bisection\index{optimization!line search!bisection} to try to choose the value of $\mu$ which minimizes $f$ in analogy to Frank-Wolfe.\index{static traffic assignment!algorithms!Frank-Wolfe}
All of these methods require solving multiple traffic assignment problems at each iteration.

Regarding convergence in step 5, one can either compare the progress made in decreasing $f$ over the last few iterations, or the changes in the investments $\mb{y}$.
This choice of stopping criterion is different in nature than those used in Chapter~\ref{chp:solutionalgorithms} for solution methods to TAP.
In that case, we can prove theoretical convergence to the equilibrium solution, and we can design our stopping criteria directly on the equilibrium condition.
The method described above, by contrast, is \emph{not} proven to converge to the global optimum, and can get stuck in solutions which are locally optimal but not globally so.
(This often happens in nonconvex optimization problems.)  By terminating the algorithm when no more progress is made, we are checking that we have found a local optimum, but cannot guarantee that we have found a global one.

As stated above, the network design problem does not have a convex feasible region, and even if one were to eliminate $\mb{x}$ by writing $\mb{x}$ as a function of $\mb{y}$, the resulting function can be shown to be nonconvex and have multiple local optimal solutions.
Therefore, the method described above is not guaranteed to converge to a global optimum solution.
This is perhaps a bit disappointing; but, as stated above, at present there is no way to ensure global optimality within a reasonable amount of computation time.
Therefore, alternative heuristics can only be compared by the quality of solutions obtained for a given quantity of computational effort.

For an example, consider the Braess network shown in Figure~\ref{fig:braess_ndp}, where the link performance functions are shown and $\Theta = \frac{1}{20}$.
Notice now that the coefficients of $x_{ij}$ in the link performance functions now depend on the amount of money $y_{ij}$ invested as well, with decreasing marginal returns as more money is spent.
When $\mb{y} = \mb{0}$, the solution to the user equilibrium problem $\mb{x}(\mb{y})$ is the now-familiar solution to the Braess network, which divides flow evenly among all three paths in the network.
The total system travel time is 552 vehicle-minutes, so the value of the objective~\eqn{ndptstt} function is $\frac{1}{20} (552) + 0 = 27.6$.
This completes the first two steps of the algorithm.

\stevefig{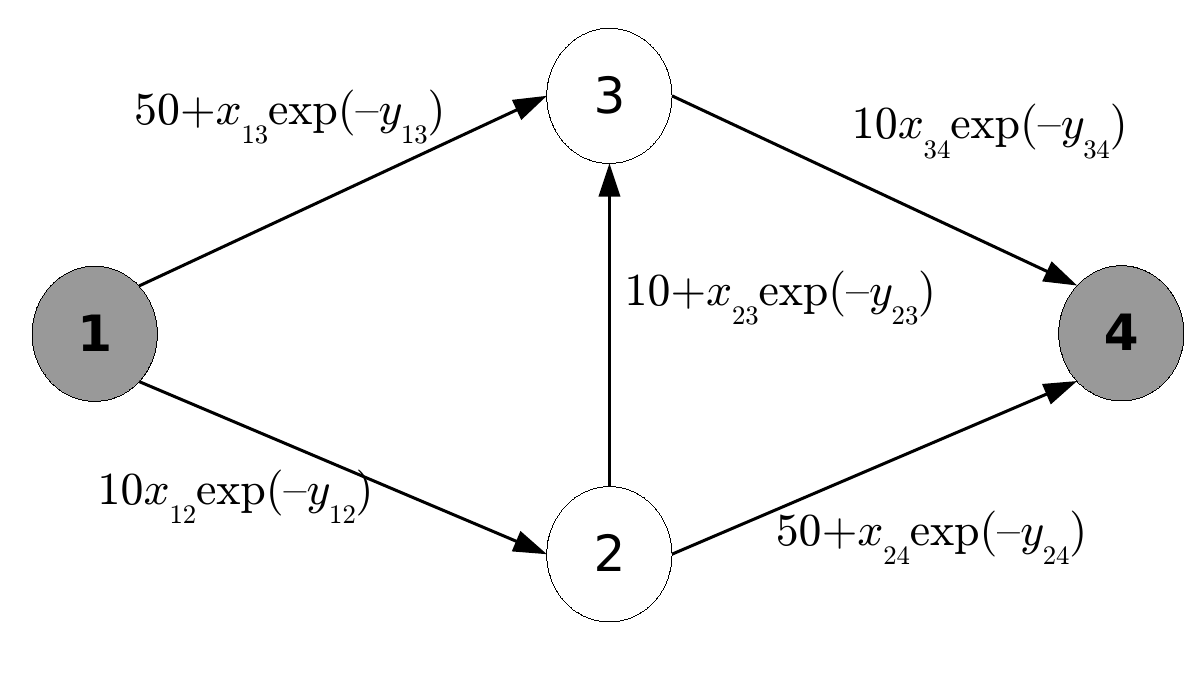}{Network design problem example.}{0.75\textwidth}

For the third step, we must solve five sensitivity problems, one for each link, to determine $\pdr{f}{y_{ij}}$.
These five problems are shown in Figure~\ref{fig:braessndpsensitivityproblems}, where we have substituted the initial values $\mb{y} = \mb{0}$ and the equilibrium link flows $\mb{x}$.
In all of these, notice that the demand is zero, and the link performance functions are affine.
Within each problem, the link being improved has a slightly different link performance function, accounting for the effect of the link improvement.
The other links' performance functions only reflect their change due to shifting flows.
For instance, in the problem in the upper left, link (1,3) is not improved, so its link performance function is simply $\xi_{13} \pdr{t_{13}}{x_{13}} = 10 \xi_{12} \exp(-y_{13}) = 10 \xi_{12}$ since $y_{13} = 0$.
Since link (1,2) is being improved, in addition to the term $\xi_{12} \pdr{t_{12}}{x_{12}} = \xi_{12} \exp(-y_{12}) = \xi_{14}$, we add the constant term $\pdr{t_{12}}{y_{12}} = -x_{14} \exp(-y_{12}) = -4$ since $x_{12} = 4$ and $y_{12} = 0$.

\stevefig{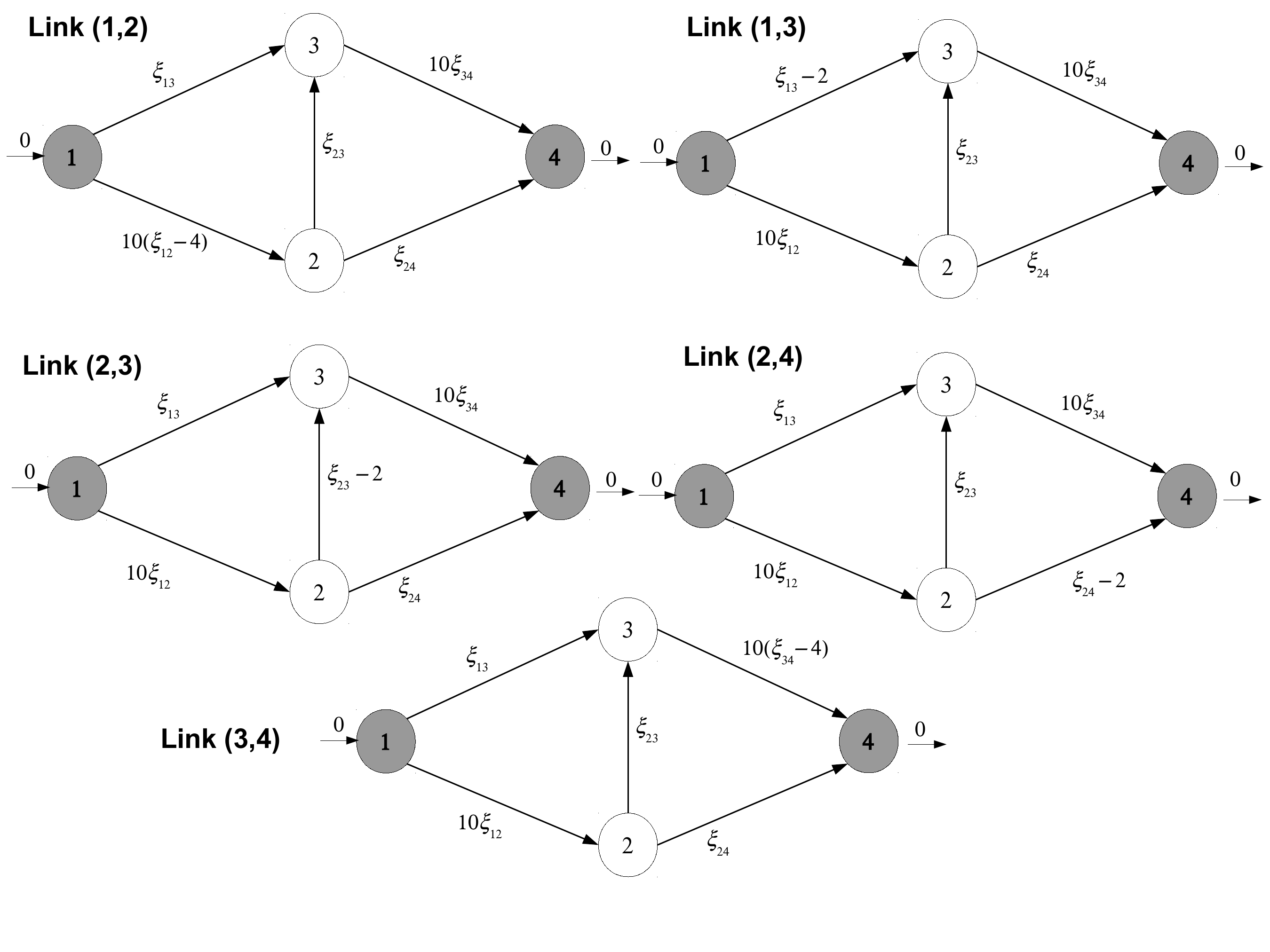}{Five sensitivity problems for network design.}{\textwidth}

The solutions (in terms of $\xi_{ij}$) to the five sensitivity problems are shown in Table~\ref{tbl:braessndpsensitivitysolutions}, as can be verified by substituting these $\xi$ values into the networks in Figure~\ref{fig:braessndpsensitivityproblems}.
Substituting these $\xi$ values into equation~\eqn{ndpder}, along with the current values of the travel times and link flows, gives the gradient
\labeleqn{ndpgradient}{\nabla_\mb{y} f = \vect{ y_{12} \\ y_{13} \\ y_{23} \\ y_{24} \\ y_{34} } = \vect{ -0.85 \\ 0.49 \\ 1.42 \\ 0.49 \\ -0.85} \,. }

\begin{table}
\begin{center}
\caption{Solutions to Braess sensitivity problems.
\label{tbl:braessndpsensitivitysolutions}}
\begin{tabular}{c|ccccc}
      & \multicolumn{5}{c}{Sensitivity problem} \\
Link  & (1,2) & (1,3) & (2,3) & (2,4) & (3,4) \\
\hline
(1,2) & 3.36 & $-0.17$ & 0.15 & 0.014 & $-0.28$ \\
(1,3) & $-3.36$ & 0.17 & $-0.15$ & $-0.014$ & 0.28 \\
(2,3) & 3.08 & $-0.15$ & 0.31 & $-0.15$ & 3.08 \\
(2,4) & 0.28 & $-0.014$ & $-0.15$ & 0.17 & $-3.36$ \\
(3,4) & $-0.28$ & 0.014 & 0.15 & $-0.17$ & 3.36 \\
\end{tabular}
\end{center}
\end{table}

Each component of the gradient shows how the objective of the network design problem will change if a unit of money is spent improving a particular link.
In the derivative formula~\eqn{ndpder}, the term in parentheses represents the marginal change in total system travel time, and the addition of unity at the end of the formula represents the increase in total expenditures.
If the derivative~\eqn{ndpder} is negative, then the reduction in total system travel time from a marginal improvement in the link will outweigh the investment cost.
If it lies between zero and one for a link, then a marginal investment in the link will reduce total system travel time, but the cost of the improvement will outweigh the value of the travel time savings.
If it is greater than one, then total system travel time would actually increase if the link is improved, as in the Braess paradox.\index{Braess paradox}
So, in this example, the gradient~\eqn{ndpgradient} shows that improvements on links (1,2) and (3,4) will be worthwhile; improvements on links (1,3) and (2,4) would reduce $TSTT$ but not by enough to outweigh construction cost; and an improvement on link (2,3) would actually be counterproductive and worsen congestion.

So, proceeding to step 4 of the network design algorithm, we choose a trial value $\mu = 1$ and apply~\eqn{ndpupdate} to obtain a candidate solution where $y_{12} = y_{34} = 0.85$ and all other $y_{ij}$ values remain at zero.
Re-solving equilibrium with the new link performance functions, the new equilibrium solution is to load all vehicles on the middle path, that is, $x_{12} = x_{23} = x_{34} = 6$ and $x_{13} = x_{24} = 0$.
The total system travel time is now 405 vehicle-minutes, so the objective is $\frac{1}{20}(405) + 2(0.85) = 22.0$.
This reduces the objective from its current value of 27.6, so we accept the step size $\mu = 1$, and return to step 2.

Notice that solving the network design problem requires solving a very large number of traffic assignment subproblems: once for each iteration to determine $\mb{x}$; modified sensitivity problems for each link to calculate derivatives; and again multiple times per iteration to identify $\mu$.
Solving practical problems can easily require solution of thousands or even millions of traffic assignment problems.
In bilevel programs such as network design, having an efficient method for solving traffic assignment problems is critical.
Path-based and bush-based algorithms can be efficiently ``warm-started,''\index{static traffic assignment!algorithms!warm starting} making them good choices for this application.
\index{network design|)}

\section{OD Matrix Estimation}
\label{sec:odme}

\index{OD matrix!estimation|(}
As described in Section~\ref{sec:tapdata}, practical application of the traffic assignment problem requires two main inputs: information about the physical network (roadway topology, link performance functions) as well as information about travel demand patterns (the OD matrix).
The former is much easier to obtain and validate.
Standard link performance functions (such as the BPR function) can be used only knowing the free-flow time and capacity of links, which can be easily estimated; and in principle, this information can be directly inferred from field measurements and traffic sensors.
The OD matrix is quite a bit harder to estimate, for several reasons.
First, there is no practical way to directly observe the complete OD matrix; at best we can work with a sample of travelers who consent to reveal their travel patterns.
Second, at least until recently, it was difficult to obtain this information without travelers explicitly reporting their origins and destinations --- traffic sensors tell you what is happening on a specific link, but not where the people are coming from or where they are going.
Recently, GPS and Bluetooth technologies have become more commonplace and can in principle provide origins and destinations automatically (a rough estimate could even be obtained from cellular phone traces) --- but there are major privacy issues associated with using such data, as well as a number of data inference issues involved in translating these data into an OD matrix.
Third, the number of entries in an OD matrix is much larger than the number of links in the network: the number of OD pairs is generally proportional to the \emph{square} of the number of nodes, while the number of links is generally proportional to the number of nodes themselves.
It is typical to see practical networks which have tens of thousands of links, but millions of OD pairs.

Therefore, it is natural to find ways to determine an OD matrix which should be used for traffic assignment.
One approach, used in the field of travel demand modeling,\index{model!demand} is based on developing behavioral models of how households make travel choices.
Using demographic and other features, there are models to estimate the total number of trips made by households over some period of time, their destinations, their mode choices, and so forth.
The first three steps of ``four-step'' model\index{planning!four-step model} described in Chapter~\ref{chp:introchapter} are one way to do this.
Another approach is to attempt to infer an OD matrix from traffic counts on specified links in the network.
Both methods have advantages and disadvantages: travel demand models provide more insights about underlying behavior (critical when developing long-range forecasts of future demand), and can directly lead to an OD matrix.
Unfortunately, the data needed to calibrate demand models is more expensive and cumbersome to work with, often involving recruiting survey participants.

Traffic sensors, on the other hand, collect data automatically, inexpensively, and without privacy issues, but there is a major dimensionality issue.
Since the number of OD pairs is much larger than the number of links, one cannot hope to uniquely determine the OD matrix solely from link counts; the problem is massively underdetermined.
In a network where all nodes are zones, it is trivial to find an OD matrix which perfectly matches link counts, by creating an OD matrix where all trips travel from one node to an adjacent node.
Even though this matrix is entirely unrealistic, there is no deviation whatsoever from the counts.
Even worse, due to unavoidable errors in traffic count records (Figure~\ref{fig:unrealisticod}), sometimes this trivial matrix will match counts much better than a more ``realistic'' matrix!  In this figure, it is likely that the total flow on the freeway is approximately 1500 vehicles from node 1 to node 3; but this does not match the counts as well as 1510 vehicles from 1 to 2, and 1490 from 2 to 3.
Simply matching counts does not provide the behavioral insight needed to identify what a ``realistic'' matrix is.\index{overfitting}

\stevefig{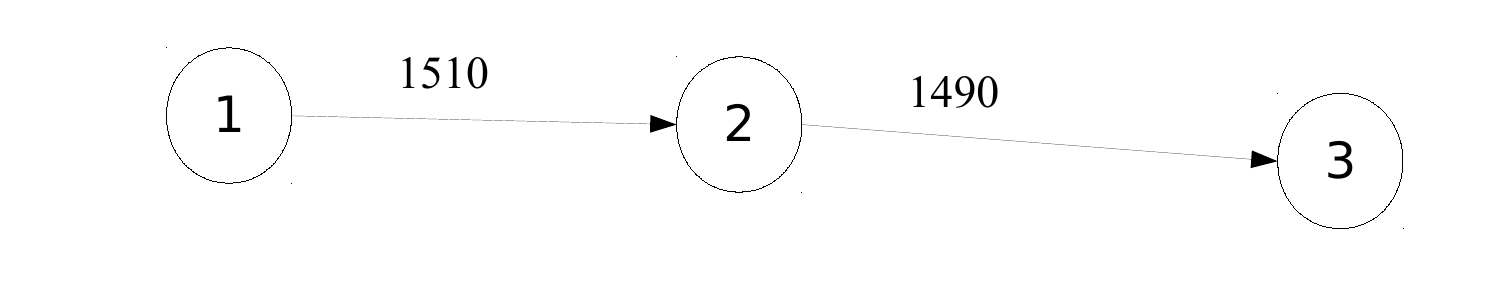}{Observed link volumes from traffic sensors along a freeway.}{0.75\textwidth}

This section provides a method to reconcile both approaches.
The idea is that an initial OD matrix $\mb{d^*}$\index{seed matrix} is already available from a travel demand model.
While this matrix is hopefully close to the true value, it also contains sampling and model estimation errors and can never be fully accurate.
However, there are sensors on a subset of links $\bar{A} \subset A$, and there is a vector $\mb{x^*}$ of traffic volume counts on these links.
The intent is to use these traffic counts to try to improve the initial OD matrix $\mb{d^*}$.
The following optimization problem expresses this:
\begin{align}
\label{eqn:odmeobj}  \min_{\mb{d}, \mb{x}}   & f(\mb{d}, \mb{x}) = \Theta \sum_{(r,s) \in Z^2} \myp{ d_{rs} - d^*_{rs} }^2  + (1 - \Theta) \sum_{(i,j) \in \bar{A}} \myp{x_{ij} - x^*_{ij}}^2 & \\
\label{eqn:odmeeqm}   \mr{s.t.} \qquad        & \mb{x} \in \arg \min_{\mb{x} \in X(\mb{d})} \sum_{(i,j) \in A} \int_0^{x_{ij}} t_{ij}(x)~dx & \\
\label{eqn:odmefeas}                         & d_{rs} \geq 0 \qquad \qquad \forall (r,s) \in Z^2 \,,&
\end{align}
where $\Theta$ is a parameter ranging from zero to one and $X(\mb{d})$ is the set of feasible link flows when the OD matrix is $\mb{d}$.

The objective function~\eqn{odmeobj} is of the least-squares\index{least squares} type, and attempts to minimize both the deviation between the final OD matrix $\mb{d}$ and the initial estimate $\mb{d^*}$, and the deviation between the equilibrium link flows $\mb{x}$ associated with the OD matrix $\mb{d}$, and the actual observations $\mb{x^*}$ on the links with sensors.
The hope is to match traffic counts reasonably well, while not wandering into completely unrealistic OD matrices.
The factor $\Theta$ is used to weight the importance of matching the initial OD matrix estimate, and the link flows.
It can reflect the relative degree of confidence in $\mb{d^*}$ and $\mb{x^*}$; if the travel demand was obtained from high-quality data and a large sample size, whereas the traffic count data is old and error-prone, a $\Theta$ value close to one is appropriate.
Conversely, if the travel demand model is less trustworthy but the traffic count data is highly reliable, a lower $\Theta$ value is appropriate.
In practice, a variety of $\Theta$ values can be chosen, and the resulting tradeoffs between matching the initial estimate and link flows can be seen.

As indicated by the constraint~\eqn{odmeeqm}, this optimization problem is also a bilevel program, because the mapping from an OD matrix $\mb{d}$ to the resulting link flows $\mb{x}$ involves the equilibrium process.
The fact that the optimization problem is bilevel also means that we cannot expect to find the global optimum OD matrix, and that heuristics should be applied.
A sensitivity-based heuristic, like the one used for the network design problem, would determine how the link flows would shift if the OD matrix is perturbed, and use this information to find a ``descent direction'' which would reduce the value of the objective.

Following the same technique as in Section~\ref{sec:ndp}, the constraint~\eqn{odmeeqm} defines $\mb{x}$ uniquely as a function of $\mb{d}$ due to the uniqueness of the link flow solution to the traffic assignment problem.
(Also note that the dependence on $\mb{d}$ appears through the feasible region, requiring that the minimization take place over $X(\mb{d})$.)  So, we can rewrite the objective function as a function of $\mb{d}$ alone, by defining $\mb{x}(\mb{d})$ as the equilibrium link flows in terms of the OD matrix and transforming the objective to $f(\mb{d}, \mb{x}(\mb{d}))$.
Then, the partial derivative of the objective with respect to any entry $d^{rs}$ in the OD matrix is given by
\labeleqn{odmeder}{\pdr{f}{d^{rs}} = 2\Theta (d_{rs} - d^*_{rs}) + 2(1 - \Theta) \sum_{(i,j) \in \bar{A}} (x_{rs} - x^*_{rs}) \pdr{x_{ij}}{d^{rs}}\,,}
where the partial derivatives $\pdr{x_{ij}}{d^{rs}}$ are found from the sensitivity formulas in Section~\ref{sec:odsa}.

The vector of all the derivatives~\eqn{odmeder} forms the gradient of $f$ with respect to $\mb{d}$.
So, taking a step in the opposite direction, and ensuring that the values in the OD matrix remain non-negative provide the following update rule:
\labeleqn{odmeupdate}{\mb{d} \leftarrow \mys{\mb{d} - \mu \nabla_\mb{d} f}^+\,,}
where $\mu$ is a step size to be determined, and the $[\cdot]^+$ operation is applied to each component of the vector.
This leads to the following algorithm for OD matrix estimation:
\begin{enumerate}
\item Initialize $\mb{d} \leftarrow \mb{d^*}$.
\item Calculate the link flows $\mb{x}(\mb{d})$ by solving the traffic assignment problem with OD matrix $\mb{d}$.
\item For each OD pair $(r,s)$ determine $\pdr{f}{d^{rs}}$ by solving the sensitivity problem in Section~\ref{sec:odsa} and using~\eqn{odmeder}.
\item Update $\mb{d}$ using~\eqn{odmeupdate} for a suitable step size $\mu$.
\item Test for convergence, and return to step 2 if not converged.
\end{enumerate}
The comments about the step size $\mu$ and convergence criteria from the network design problem (Section~\ref{sec:ndp}) apply equally as well here: $\mu$ can be determined using an iterative line search procedure, choosing smaller values until the objective function decreases, and the algorithm can be terminated when it fails to make additional substantial progress in reducing the objective.

\begin{figure} 
\hfill
\includegraphics[width=0.6\textwidth]{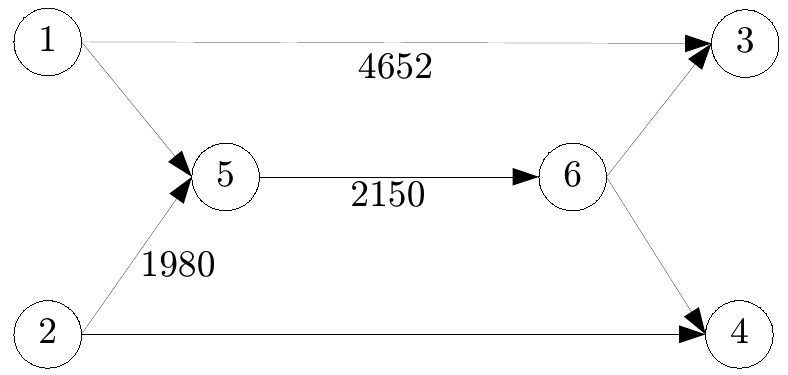}
\hfill
\raisebox{1in}{
\begin{tabular}{r|cc}
& 3       & 4 \\
\hline
1 & 5,000 & --- \\
2 & ---   & 10,000 \\
\end{tabular}
}
\hfill
\caption{Example network for OD matrix estimation.
(Observed link volumes and initial OD matrix shown.)  \label{fig:odmenet}}
\end{figure}

This procedure is demonstrated using the network in Figure~\ref{fig:odmenet}.
In this network, traffic counts are available on three of the links in the network, and an initial OD matrix is available based on a travel demand model.
The link performance function on every link is $t_{ij} = 10 + x_{ij} / 100$, and $\Theta$ is given as $\frac{1}{10}$.

We begin by initializing the OD matrix to the initial matrix $\mb{d^*}$, and solving a traffic assignment problem.

As the reader can verify, the equilibrium link flows on the three links with available counts are $x_{13} = 4733 \frac{1}{3}$, $x_{25} = 1933 \frac{1}{3}$, and $x_{56} = 2200$.
The first sum in the objective~\eqn{odmeobj} gives the fit of the OD matrix to the initial matrix (zero) and the the second sum gives the fit of the equilibrium link flows to the observed link flows (11292).
When weighted with $\Theta = 0.1$, this gives an initial objective value of 10164.
Next, we solve two sensitivity problems, one for $d_{13}$ and one for $d_{24}$.
(If there was reason to believe there were trips between 1 and 4, and between 2 and 3, we would solve sensitivity problems for those OD pairs as well.)  These are shown in Figure~\ref{fig:odmesensitivityproblems}.
The link performance functions are the same in both problems; the only difference is in the demand.
The solutions to these sensitivity problems are shown in Table~\ref{tbl:odmesensitivitysolutions}.

\stevefig{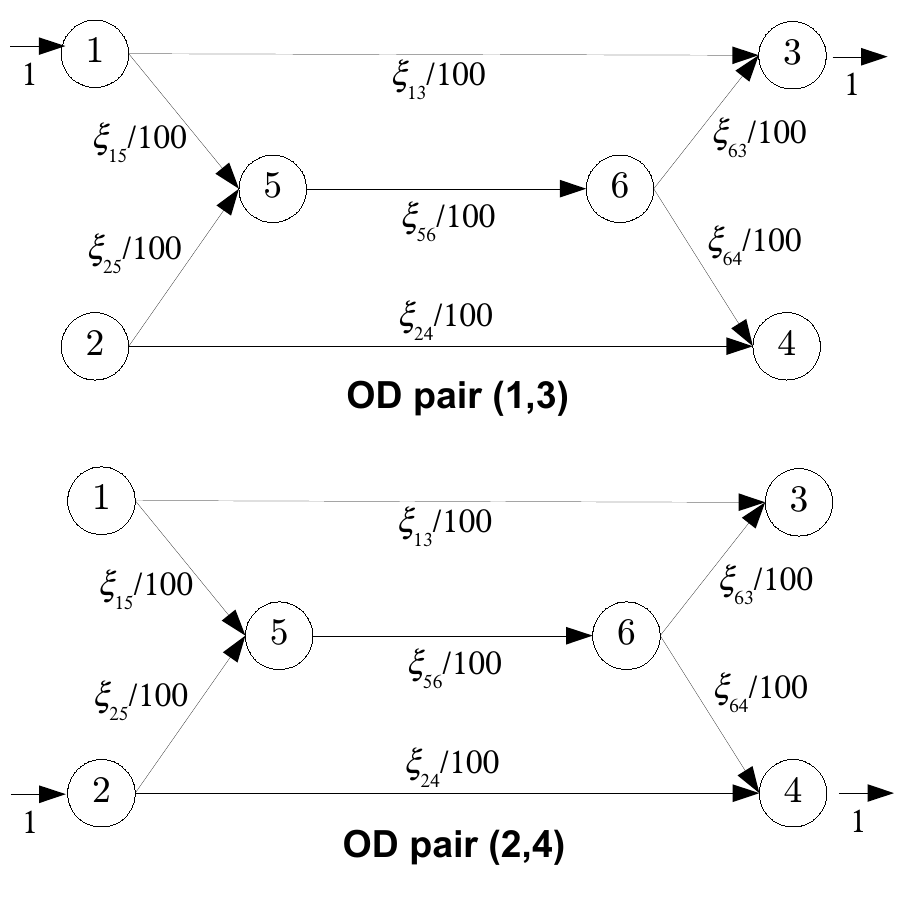}{Two sensitivity problems for OD matrix estimation.}{0.8\textwidth}

\begin{table}
\begin{center}
\caption{Solutions to OD matrix estimation sensitivity problems.
\label{tbl:odmesensitivitysolutions}}
\begin{tabular}{c|cc}
      & \multicolumn{2}{c}{Sensitivity problem} \\
Link  & $d_{13}$ & $d_{24}$ \\
\hline
(1,3) & 0.733 & 0.067 \\
(1,5) & 0.267 & $-0.067$ \\
(2,4) & 0.067 & 0.733 \\
(2,5) & $-0.067$ & 0.267 \\
(5,6) & 0.200 & 0.200 \\
(6,3) & 0.267 & $-0.067$ \\
(6,4) & $-0.067$ & 0.267 
\end{tabular}
\end{center}
\end{table}

We now have all the information needed to calculate the gradient of the objective function, using equation~\eqn{odmeder}.
Substituting the values of $x_{ij}$, $\xi_{ij}$, the observed link flows, and the initial OD matrix, we calculate
\labeleqn{odmegradient}{\nabla_\mb{d} f = \vect{ \pdr{f}{d_{13}} \\ \pdr{f}{d_{24}} } = \vect{ 131 \\ 5.36 }\,. }
Taking a trial step in this direction with $\mu = 1$, the updating rule~\eqn{odmeupdate} gives the candidate OD matrix $d_{13} = 4869$, $d_{24} = 9995$.
The resulting equilibrium link flows include $x_{13} = 4637$, $x_{25} = 1941$, and $x_{56} = 2150$, so the ``fit'' of the equilibrium link flows and traffic counts has improved from 11293 to 2293.
The ``fit'' of the OD matrix has worsened from 0 to 17179, but with the weight $\Theta = 0.1$, the overall objective still decreases from 10164 to 3782.
Therefore, we accept the step size $\mu = 1$, and return to step 3 to continue updating the OD matrix.
\index{OD matrix!estimation|)}
\index{optimization!bilevel|)}

\section{Historical Notes and Further Reading}
\label{sec:sensitivity_litreview}

The derivation of the equilibrium sensitivity analysis in this chapter follows
that in \cite{boyles_partbcontraction} and \cite{jafari16}.
There are alternative ways to derive the same results, using the implicit function theorem\index{implicit function theorem}~\citep{tobin88,cho00,yang05} or results from sensitivity of variational inequalities\index{variational inequality!sensitivity analysis}~\citep{patriksson04,lu08}.
In particular, for application of the latter approach to the network design problem, see \cite{josefsson07}.

Sensitivity analysis is not the only way solve bilevel optimization problems such as network design or OD matrix estimation.
The network design problem has wide application outside of transportation planning, and is used for designing infrastructure and logistics systems in many domains.
\cite{magnanti84} provide an overview of applications and algorithms based on the structure of the optimization problem.
\cite{farahani13} describe more recent progress and applications for network design.
OD matrix estimation has been part of travel demand forecasting for decades.
\cite{cascetta88}, and Chapter 7 of~\cite{bell97}, present reviews of methods that have been applied to this end.
A key difficulty is dealing with the dimensionality of the problem: OD matrices have many more degrees of freedom than the number of links where counts are available.
There may be advantages in adjusting behavioral parameters earlier in the travel demand process, as in trip generation or distribution, where the number of parameters to be estimated is of the same order of magnitude as the number of observations available~\citep{alexander21}.

An alternative approach to sensitivity analysis is to approximate the lower-level equilibrium problem by a ``nicer'' function, reducing the bilevel optimization problem to a single level problem, with a constraint showing how the link flows will (approximately) change given changes in the upper-level decision variables.
Such approaches go by several names, including metamodeling\index{metamodel}~\citep{osorio13,osorio19} and response surface methods.\index{response surface}
More recently, the lower-level problem has been approximated using neural networks\index{neural network} or other machine learning\index{machine learning} techniques~\citep{bagloee18}.
\cite{dempe20} provide an overview of current approaches in bilevel optimization.

Owing to the difficulties in solving bilevel optimization problems, metaheuristic methods are common.
Two such heuristics (simulated annealing and genetic algorithms) are discussed in Appendix~\ref{chp:fancyoptimization}.
These approaches are relatively generic and can be applied to a wide variety of optimization problems.
They lack useful performance guarantees and do not exploit specific problem structure in the ways that the methods presented in other chapters of this book do.
These drawbacks are less impactful for bilevel optimization, given their seeming intractability.
Many bilevel programs, including some variants of network design, are $NP$-hard,\index{NP-hard} a class of problems which has been studied extensively.
Efficient and exact algorithms remain elusive for such problems.
Appendix~\ref{chp:algorithmcomplexity} describes $NP$-hardness and other complexity classes in greater detail.
\index{static traffic assignment!sensitivity analysis|)}

\section{Exercises}
\label{exercises_sensitivityanalysis}

\begin{enumerate}
\item \diff{63} Given a nondegenerate equilibrium solution in a network with a single origin and continuous link performance functions, show that the equilibrium bushes remain unchanged in a small neighborhood of the current OD matrix.
\label{ex:bushfixed}
\item \diff{65} Given a nondegenerate equilibrium solution in a network with a single origin and differentiable link performance functions, show that the derivatives $dx_{ij}/dd^{rs}$ exist at the current equilibrium solution.
\label{ex:bushdifferentiable}
\item \diff{56} Extend the results in Exercises~\ref{ex:bushfixed} and~\ref{ex:bushdifferentiable} to networks with multiple origins.
\item \diff{84} Show that the entropy-maximizing path flow solution is a continuous function of the OD matrix.
\item \diff{26} Verify that the optimality conditions for~\eqn{ddsaobjective}--\eqn{ddsafc} include the conditions~\eqn{ddbushreducedcost}--\eqn{ddbushflownnonbush}.
\item \diff{33} In the modified Braess network of Figure~\ref{fig:sabraessmod}, find the sensitivity of each link's flow to the demand $d^{14}$.
\label{ex:ddpractice}
\item \diff{35} In the modified Braess network of Figure~\ref{fig:sabraessmod}, find the sensitivity of each link's flow to the link performance function parameter $y^{23}$.
What value of this parameter minimizes the equilibrium travel time?  Suggest what sort of real-world action would correspond to adjusting this parameter to this optimal value.
\label{ex:dxpractice}
\stevefig{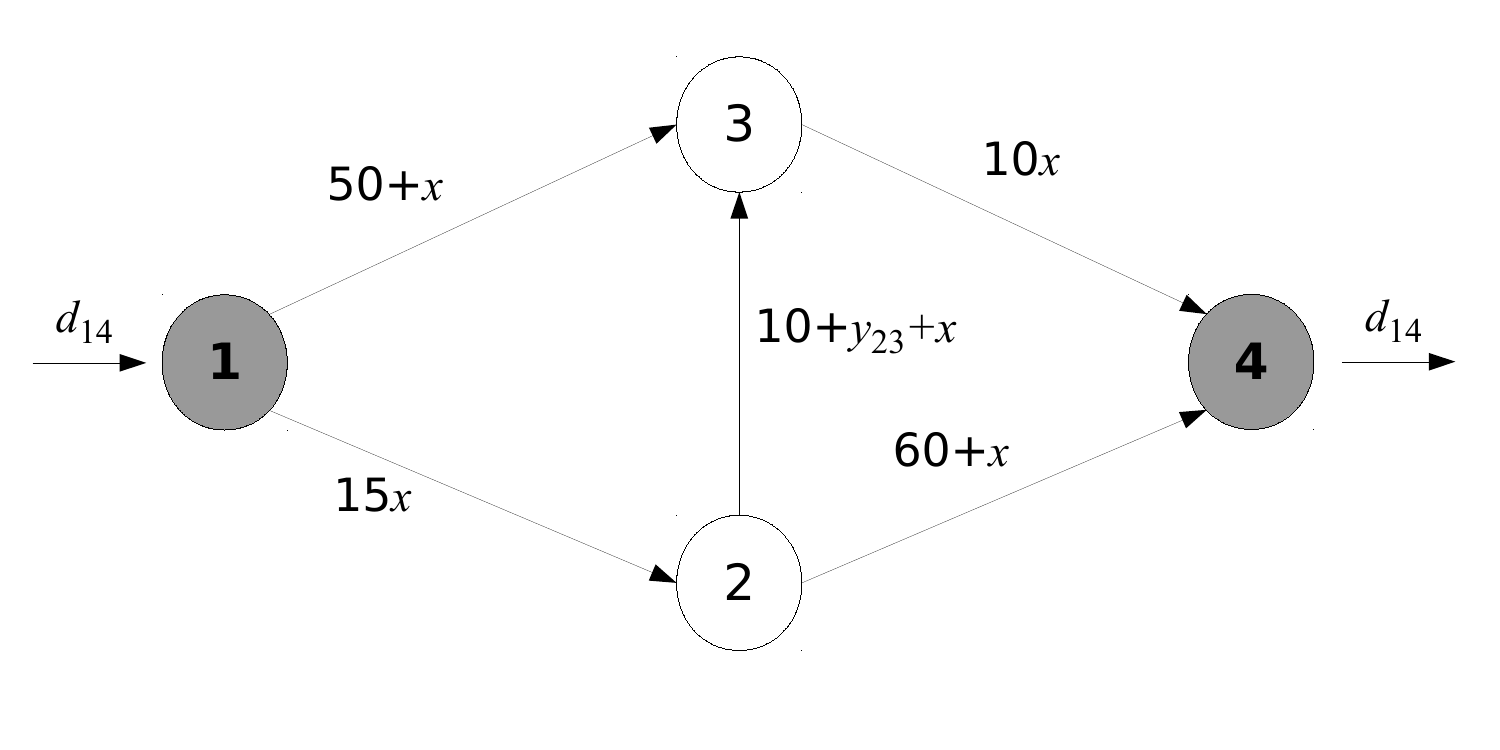}{Network for use in Exercises~\ref{ex:ddpractice} and~\ref{ex:dxpractice}.}{0.6\textwidth}
\item \diff{44} In the network design problem for Figure~\ref{fig:braess_ndp_mod}, give the gradient of the objective function at the initial solution $\mb{y} = \mb{0}$.
Assume that $\Theta = \frac{1}{20}$.
\label{ex:ndpstart}
\stevefig{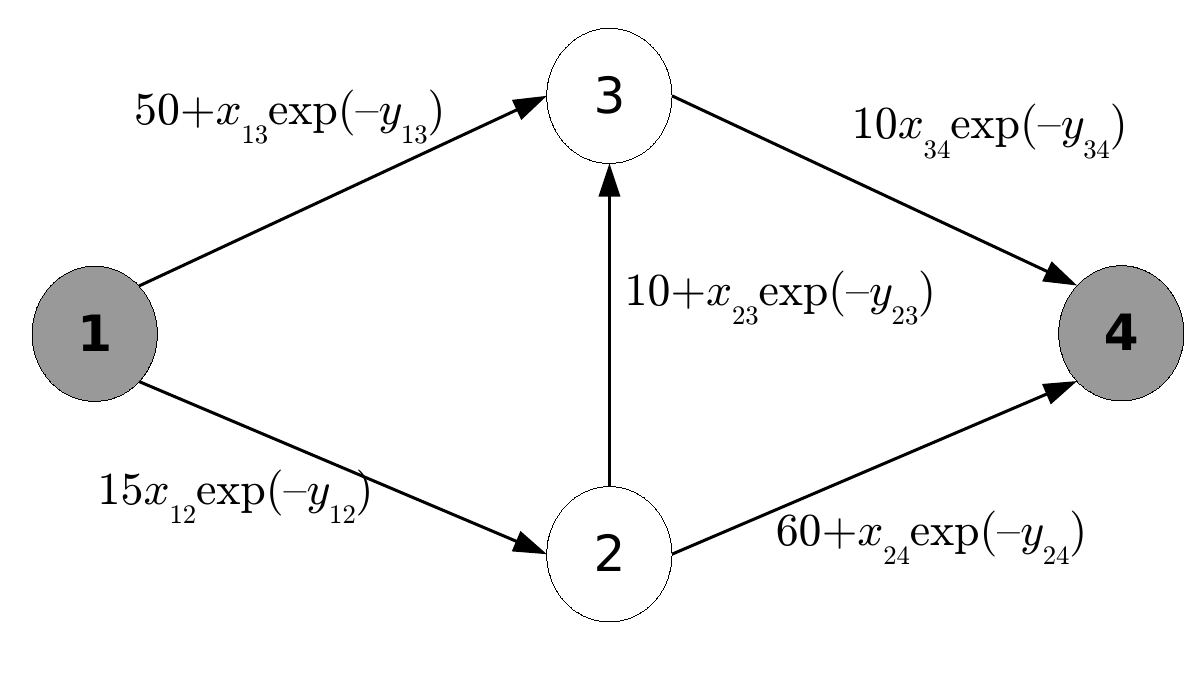}{Network for use in Exercises~\ref{ex:ndpstart} and~\ref{ex:ndpend}.}{0.6\textwidth}
\item \diff{48} Continue Exercise~\ref{ex:ndpstart} by performing three iterations of the algorithm given in the text.
What is the resulting total system travel time and construction cost? \label{ex:ndpend}
\item \diff{25} Write out the network design optimization problem for the network in Figure~\ref{fig:ndp2link}, with $\Theta = 1$.
Show that the feasible region for this problem is not convex by constructing a counterexample.
\label{ex:ndp2link}
\stevefig{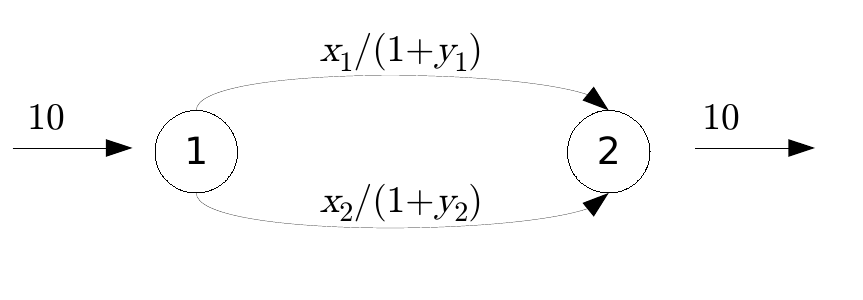}{Network for use in Exercise~\ref{ex:ndp2link}.}{0.6\textwidth}
\item \diff{79} Design a heuristic for the network design problem, based on simulated annealing as discussed in Section~\ref{sec:simulatedannealing}.
Compare the performance of this heuristic with the algorithm given in the text for several networks.
\item \diff{79} Design a heuristic for the network design problem, based on genetic algorithms as discussed in Section~\ref{sec:geneticalgorithm}.
Compare the performance of this heuristic with the algorithm given in the text for several networks.
\item \diff{47} In the OD matrix estimation problem of Figure~\ref{fig:odmebraess}, give the gradient of the objective function at the initial solution $d^*_{14} = 6$, $d^*_{24} = 4$.
What is the value of the objective function if $\Theta = \frac{1}{2}$? \label{ex:odmestart}
\stevefig{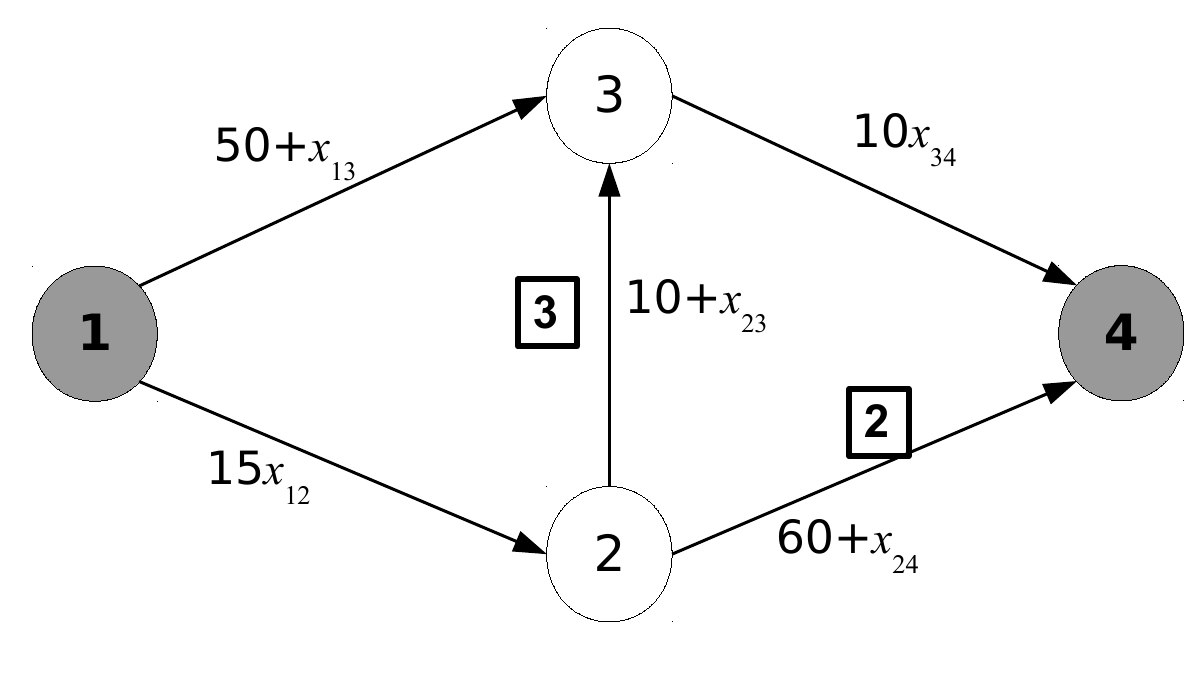}{Network and observed link flows (in boxes) for Exercise~\ref{ex:odmestart}.}{0.6\textwidth}
\item \diff{49} Continue Exercise~\ref{ex:odmestart} by performing three iterations of the algorithm given in the text.
What is the resulting OD matrix and objective function value? \label{ex:odmecontinued}
\item \diff{59} Generate five additional OD matrices corresponding to the network in Figure~\ref{fig:odmebraess}, with different values of $\Theta$.
Which of these OD matrices seems most reasonable to you, and why?
\item \diff{79} Design a heuristic for the OD matrix estimation problem, based on simulated annealing as discussed in Section~\ref{sec:simulatedannealing}.
Compare the performance of this heuristic with the algorithm given in the text for several networks.
\item \diff{79} Design a heuristic for the OD matrix estimation problem, based on genetic algorithms as discussed in Section~\ref{sec:geneticalgorithm}.
Compare the performance of this heuristic with the algorithm given in the text for several networks.
\item \diff{76} Perform the following ``validation'' exercise: create a small network with a given OD matrix, and find the equilibrium solution.
Then, given the equilibrium link flows, try to compute your original OD matrix using the algorithm given in the text.
Do you get your original OD matrix back?  
\end{enumerate}

\chapter{Extensions of Static Assignment}
\label{chp:staticextensions}

The basic traffic assignment problem (TAP) was defined in Chapter~\ref{chp:trafficassignmentproblem} as follows: we are given a network $G = (N,A)$, link performance functions $t_{ij}(x_{ij})$, and the demand values $d^{rs}$ between each origin and destination.
The objective is to find a feasible vector of path flows (or link flows) which satisfy the principle of user equilibrium, that is, that every path with positive flow has the least travel time among all paths connecting that origin and destination.
We formulated this as a VI (find $\mathbf{\hat{h}} \in H$ such that $\mathbf{c}(\mathbf{\hat{h}}) \cdot (\mathbf{\hat{h}} - \mathbf{h}) \leq 0$ for all $\mathbf{h} \in H$) and as the solution to the following convex optimization problem:
\begin{align}
   \min_{\mathbf{x},\mathbf{h}}  \quad & \sum_{(i,j) \in A} \int_0^{x_{ij}} t_{ij}(x) dx  &     \label{eqn:beckmann_ch9} \\
   \mathrm{s.t.}                 \quad & x_{ij} = \sum_{\pi \in \Pi} h^\pi \delta_{ij}^{\pi} & \forall (i,j) \in A \label{eqn:pathlinkmap} \\
                                       & \sum_{\pi \in \Pi^{rs}} h^\pi = d^{rs}              &     \forall (r,s) \in Z^2 \label{eqn:demand_chpext} \\
                                       & h^\pi \geq 0                                     &     \forall \pi \in \Pi \label{eqn:nonneg_chpext} 
\end{align}

This formulation remains the most commonly used version of traffic assignment in practice today.
However, it is not difficult to see how some of the assumptions may not be reasonable.
This chapter shows extensions of the basic TAP which relax these assumptions.
This is the typical course of research: the first models developed make a number of simplifying assumptions, in order to capture the basic underlying behavior.
Then, once the basic behavior is understood, researchers develop progressively more sophisticated and realistic models which relax these assumptions.

This chapter details three such extensions.
Section~\ref{sec:elasticdemand} relaxes the assumption that the OD matrix is known and fixed, leading to an \emph{elastic demand} formulation.
Section~\ref{sec:linkinteractions} relaxes the assumption tat the travel time on a link depends only on the flow on that link (and not on any other link flows, even at intersections).
Section~\ref{sec:sue} relaxes the assumption that travelers have accurate knowledge and perception of all travel times in a network, leading to the important class of \emph{stochastic user equilibrium} models.

For simplicity, all of these variations are treated independently of each other.
That is, the OD matrix is assumed known and fixed in all sections \emph{except} Section~\ref{sec:elasticdemand}, and so forth.
This is done primarily to keep the focus on the relevant concept of each section, but also to guard the reader against the temptation to assume that a model which relaxes \emph{all} of these assumptions simultaneously is necessarily better than one which does not.
While realism is an important characteristic of a model, it is not the only relevant factor when choosing a mathematical model to describe an engineering problem.
Other important concerns are computation speed, the existence of enough high-quality data to calibrate and validate the model, transparency, making sure the sensitivity of the model is appropriate to the level of error in input data, ease of explanation to decision makers, and so on.
All of these factors should be taken into account when choosing a model, and you can actually do worse off by choosing a more ``realistic'' model when you don't have adequate data for calibration --- the result may even give the impression of ``false precision'' when in reality your conclusions cannot be justified.

\section{Elastic Demand}
\label{sec:elasticdemand}

\index{elastic demand|see {static traffic assignment, elastic demand}}
\index{OD matrix!elastic|see {static traffic assignment, elastic demand}}
\index{static traffic assignment!elastic demand|(}
The assumption that the OD matrix is known and fixed can strain credibility, particularly when considering long time horizons (20--30 years) or when projects are major enough to influence travel decisions at all levels, not just route choice.
For instance, consider the ``induced demand'' phenomenon where major expansion of roadway capacity ends up increasing the amount of demand.
This is partly due to changes in route choice (which the basic TAP accounts for), but is also due to changes in other kinds of travel choices, such as departure time, mode, destination, or trip frequency.

Therefore, it is desirable to develop a model which can relax the assumption of an exogenous OD matrix known \emph{a priori}.
This section describes how TAP can be extended to accommodate this relaxation.
This is called the \emph{elastic demand} formulation of TAP.
The elastic demand model is at once more and less useful than basic TAP: more useful because it can provide a more accurate view of the impacts of transportation projects; less useful because it is harder to calibrate.

\subsection{Demand functions}

\index{demand function|(}
The new idea in the traffic assignment problem with elastic demand is the \emph{demand function}, which relates the demand for travel between an origin and destination to the travel time between these zones.
Specifically, let $D^{rs}$\label{not:Drs} be a function relating the demand between $r$ and $s$ to the travel time $\kappa^{rs}$ on the shortest path between these zones.\footnote{An alternative is to have $D^{rs}$ be a function of the \emph{average} travel time on the used paths between $r$ and $s$, not the shortest.
At equilibrium it doesn't matter because all used paths have the same travel time as the shortest, but in the process of finding an equilibrium the alternative definition of the demand function can be helpful.
For our purposes, though, definition in terms of the shortest path time is more useful because it facilitates a link-based formulation.}  Generally, $D^{rs}$ is a nonincreasing function --- as the travel time between $r$ and $s$ increases, the demand for travel between these nodes is lower.
It will also be highly useful to assume that $D^{rs}$ is invertible as well, which will require it to be strictly decreasing.
The inverse demand function will give the travel time between $r$ and $s$ corresponding to a given demand level.

As an example, let $D^{rs}(\kappa^{rs}) = 200 \exp(-\kappa^{rs} / 100)$.
If the travel time between $r$ and $s$ is $\kappa^{rs} = 10$ by the shortest path, then the demand is $d^{rs} = 181$.
If the travel time was 20 minutes between these zones, the demand will be 163.7, which is lower as fewer drivers choose to travel between $r$ and $s$.
The inverse demand function is $D^{-1}_{rs}(d^{rs}) =  100 \log (200 / d^{rs})$, and can be used to calculate the $\kappa^{rs}$ value corresponding to a given $d^{rs}$: when the demand is 181, the shortest path travel time is $\kappa^{rs} = D^{-1}_{rs} (181) = 10$, and so on.\footnote{In the text, we typically indicate OD pairs with a superscript, as in $d^{rs}$, and link variables with a subscript, as in $x_{ij}$.
In elastic demand, we will often need to refer to inverse functions, and writing $(D^{rs})^{-1}$ is clumsy.
For this reason, OD pairs may also be denoted with a subscript, as in $D^{-1}_{rs}$.
This is purely for notational convenience and carries no significance.}  Or, if $D^{rs} = 50 - \frac{1}{2}\kappa^{rs}$, then $D^{-1}_{rs} = 100 - 2 d^{rs}$ and a travel time of 10 minutes corresponds to a demand of 45 vehicles as can be seen by substituting either number into the corresponding equation.

An attentive reader may have noticed a potential issue with the demand function $D^{rs} = 50 - \frac{1}{2} \kappa^{rs}$, namely that the demand would be negative if $\kappa^{rs} > 100$.
In reality, the demand would simply equal zero if the travel time exceeded 100.
We could patch this by redefining $D^{rs}$ as $[50 - \frac{1}{2} \kappa^{rs}]^+$, but then $D^{rs}$ is no longer invertible.
Instead, we can allow $D^{rs}$ to take negative values, but replace the relation $d^{rs} = D^{rs}(\kappa^{rs})$ with $d^{rs} = [D^{rs}(\kappa^{rs})]^+$.
This allows us to have $D^{rs}$ be strictly decreasing (and thus invertible), but still allow the travel demand to be zero when costs are sufficiently high.
While this ``trick'' may seem a bit trivial (or at least not very useful), it will eventually allow us to formulate the elastic demand equilibrium problem as a variational inequality and convex program, as shown below.

The demand function can be used to define a \emph{consumer surplus}\index{consumer surplus} $CS$\label{not:CS} representing the benefits of mobility in a region, defined by
\labeleqn{consumersurplus}{CS = \sum_{(r,s) \in Z^2} \myp{\int_0^{d^{rs}} D^{-1}_{rs} (\omega)~d\omega - d^{rs} \kappa^{rs}}}
(We are using $\omega$\label{not:omega} for the dummy variable of integration instead of $d$ because $\int D^{-1}(d)~dd$ is notationally awkward.)   The interpretation of this formula is as follows.
Each driver has a certain travel time threshold: if the travel time is greater than this threshold, the trip will not be made, and if the travel time is less than this threshold, the trip will be made.
Different drivers have different thresholds, and the demand function represents the aggregation of these thresholds: when the travel time is $\kappa$,\label{not:kappagen} $D(\kappa)$ represents the number of travelers whose threshold is $\kappa$ or higher.
If my threshold is, say, 15 minutes and the travel time is 10 minutes, the difference (5 minutes) can be thought of as the ``benefit'' of travel to me: the trip is worth 15 minutes of my time, but I was able to travel for only 10.
Adding this up for all travelers provides the total benefits of travel, which is what $CS$ represents.
Figure~\ref{fig:consumersurplus} shows the connection between this concept and equation~\eqn{consumersurplus}: assume that drivers are numbered in decreasing order of their threshold values.
Then $D^{-1}(1)$ gives the threshold value for the first driver, $D^{-1}(2)$ gives the threshold value for the second driver, and so forth.
At equilibrium all drivers experience a travel time of $\kappa$, so the benefit to the first driver is $D^{-1}(1) - \kappa$, the benefit to the second driver is $D^{-1}(2) - \kappa$, and so forth.
Adding over all drivers gives equation~\eqn{consumersurplus}.
\index{demand function|)}

\stevefig{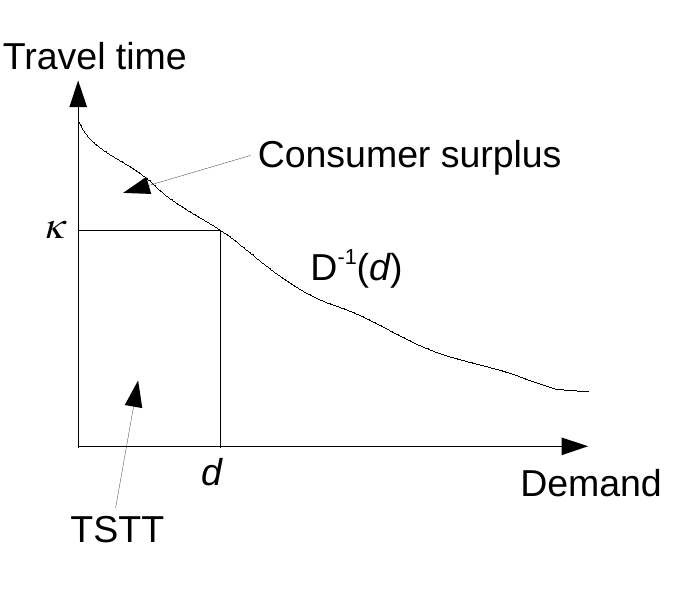}{Relationship between inverse demand function, consumer surplus, and total system travel time.}{0.5\textwidth}

\subsection{Gartner's transformation}
\label{sec:gartner}

\index{Gartner's transformation|(}
\index{static traffic assignment!elastic demand!network transformation|see {Gartner's transformation}}
Before moving to variational inequality and optimization formulations of the elastic demand problem, we'll take a short digression and show how the elastic demand problem can be cleverly transformed into a traditional equilibrium problem with fixed demand.
This transformation  works if each demand function $D^{rs}$ is bounded above.
Repeat the following for each OD pair $(r,s)$.
Let $\overline{d}^{rs}$\label{not:dbarrs} be such an upper bound for OD pair $(r,s)$.
Create a new link $(r,s)$ directly connecting origin $r$ to destination $s$, and make its link performance function $t_{rs}(x_{rs}) = D^{-1}_{rs} (\overline{d}^{rs} - x_{rs})$ where $x_{rs}$ is the flow on that new link.
If $D^{rs}$ is decreasing, then $D^{-1}_{rs} (\overline{d}^{rs} - x_{rs})$ is increasing in $x_{rs}$ so this is a valid link performance function.

Now, solve a \emph{fixed} demand problem where the demand from each origin to each destination is $\overline{d}^{rs}$.
An equilibrium on this network corresponds to an elastic demand equilibrium on the original network as follows: the flows on the links common to both networks represent flows on the actual traffic network; the flow on the new direct connection links represent drivers who choose not to travel due to excess congestion.
Think of $\overline{d}^{rs}$ as the total number of people who might possibly travel from $r$ to $s$; those that actually complete their trips travel on the original links and those who choose not to travel choose the direct connection link.
At equilibrium, all used paths connecting $r$ to $s$ (including the direct connection link) have the same travel time $\kappa_{rs}$; therefore, the flow on the direct connection link $x_{rs}$ must be such that $D^{-1}_{rs} (\overline{d}^{rs} - x_{rs}) = \kappa^{rs}$, or equivalently $x_{rs} = \overline{d}^{rs} - D^{rs}(\kappa^{rs})$, which is exactly the number of drivers who choose not to travel when the equilibrium times are $\kappa^{rs}$.
That is, the demand $d^{rs} = \overline{d}^{rs} - x_{rs}$.

The downside of this approach is that it requires creating a large number of new links.
In a typical transportation network, the number of links is proportional to the number of nodes and of the same order of magnitude (so a network of 1,000 nodes may have 3--4,000 links).
However, the number of OD pairs is roughly proportional to the \emph{square} of the number of nodes, since every node could potentially be both an origin and a destination.
So, a network with 1,000 nodes could have roughly 1,000,000 OD pairs.
Implementing the Gartner transformation requires creating a new link for every one of these OD pairs, which would result in 99.9\% of the network links being the artificial arcs for the transformation!  
\index{Gartner's transformation|)}

\subsection{Variational inequality formulation}

\index{static traffic assignment!elastic demand!variational inequality formulation|(}
\index{variational inequality!and traffic assignment|(}
Because the elastic demand problem can be expressed as a version of the regular traffic assignment problem through the Gartner transformation, we immediately have a variational inequality formulation of the elastic demand equilibrium problem.
Partition the vectors of link flows and travel times into regular and direct-connect (Gartner transformation)\index{Gartner's transformation} links, using $\mathbf{x}$ to represent regular link flows, $\mathbf{x}^\rightarrow$\label{not:xarrow} flows on direct-connect links, and $\mb{t}$ and $\mb{t}^\rightarrow$\label{not:tarrow} similarly.
Then the variational inequality is
\labeleqn{baseelasticdemandvi}{\vect{\mb{t(\hat{x})} \\ \mb{{t}^\rightarrow(\hat{x}^\rightarrow)}} \cdot \vect{\mb{\hat{x} - x} \\ \mb{\hat{x}^\rightarrow - x^\rightarrow} } \leq 0 }
Since $x_{rs} = \overline{d}^{rs} - d^{rs}$ and $t_{rs}(x_{rs}) = D^{-1}_{rs}(\overline{d}^{rs} - x_{rs})$, the variational inequality can be written in terms of the link flows and OD demands as
\labeleqn{elasticdemandvi}{\vect{\mb{t(\hat{x})} \\ \mb{-D^{-1}(\hat{d})}} \cdot \vect{\mb{\hat{x} - x} \\ \mb{\hat{d} - d} } \leq 0  }
which is the customary form.\label{not:hatdrs}
\index{static traffic assignment!elastic demand!variational inequality formulation|)}
\index{variational inequality!and traffic assignment|)}

\subsection{Optimization formulation}

\index{static traffic assignment!elastic demand!convex optimization formulation|(}
The Gartner transformation can also lead directly to a convex programming formulation of the elastic demand problem, in a similar way as the variational inequality was derived in the previous subsection.
However, it is also instructive to derive the convex programming formulation from first principles.

As discussed above, the demand is related to the demand function by $d^{rs} = [D^{rs}(\kappa^{rs})]^+$.
Put another way, \emph{the demand must always be at least as much as the demand function; further, if the demand is greater than zero then it must equal the demand function}.\index{complementarity constraint}
Thinking laterally, you might notice this is similar to the principle of user equilibrium: the travel time on any path must always be at least as large as the shortest path travel time; further, if the demand is positive then the path travel time must equal the shortest path travel time.
When deriving the Beckmann function, we showed that the latter statements could be expressed by $c^\pi \geq \kappa^{rs}$ and $h^\pi (c^\pi - \kappa^{rs}) = 0$ (together with the nonnegativity condition $h^\pi \geq 0$).
The same ``trick'' applies for the relationship between demand and the demand function: $d^{rs} \geq D^{rs}(\kappa^{rs})$, $d^{rs} ( d^{rs} - D^{rs}(\kappa^{rs})) = 0$, and the nonnegativity condition $d^{rs} \geq 0$.

It will turn out to be easier to express the latter conditions in terms of the inverse demand functions $D^{-1}$, rather than the ``forward'' functions $D$, because the convex objective function we will derive will be based on the Beckmann function.
The Beckmann function\index{Beckmann function} involves link performance functions (with units of time).
Since the inverse demand functions also are measured in units of time, it will be easier to combine them with the link performance functions than the regular demand functions (which have units of vehicles).
Expressed in terms of the inverse demand functions, the conditions above become $\kappa^{rs} \geq D^{-1}_{rs}(d^{rs})$, $d^{rs} (D^{-1}_{rs}(d^{rs}) - \kappa^{rs}) = 0$, and $d^{rs} \geq 0$.

So, this is the question before us.
What optimization problem has the following as its optimality conditions?
\begin{align}
\label{eqn:pathdual} c^\pi & \geq \kappa^{rs} & \forall (r,s) \in Z^2, \pi \in \Pi^{rs} \\
\label{eqn:pathcs}   h^\pi (c^\pi - \kappa^{rs}) &= 0 & \forall (r,s) \in Z^2, \pi \in \Pi^{rs} \\
\label{eqn:dmnddual} \kappa^{rs} & \geq D^{-1}_{rs}(d^{rs}) & \forall (r,s) \in Z^2 \\
\label{eqn:dmndcs}   d^{rs} (D^{-1}_{rs}(d^{rs}) - \kappa^{rs}) &= 0 & \forall (r,s) \in Z^2 \\
\label{eqn:nvlb}     \sum_{\pi \in \Pi^{rs}} h^\pi &= d^{rs} & \forall (r,s) \in Z^2 \\
\label{eqn:pathnng}  h^\pi & \geq 0 & \forall \pi \in \Pi \\
\label{eqn:dmndnng}  d^{rs} & \geq 0 & \forall (r,s) \in Z^2
\end{align}
The Beckmann formulation is a good place to start, since it already includes \eqn{pathdual}, \eqn{pathcs}, \eqn{nvlb}, and \eqn{pathnng}.
So let's start by conjecturing that the Lagrangian takes the form
\labeleqn{lagrangesketch}{\mc{L}(\mathbf{h}, \mathbf{d}, \bm{\kappa}) = \sum_{(i,j) \in A} \int_0^{\sum_{\pi \in \Pi} \delta_{ij}^\pi h^\pi} t_{ij}(x)~dx + \sum_{(r,s) \in Z^2} \kappa^{rs} \myp{d^{rs} - \sum_{\pi \in \Pi^{rs}} h^\pi} + F(\mathbf{d})}
where $f(\mathbf{d})$ is some function involving the OD matrix.
(Note also that $\mc{L}$ is now a function of $\mathbf{d}$ in addition to $\mathbf{h}$ and $\bm{\kappa}$, since the demand is a decision variable.)  You can check that the optimality conditions related to $\bm{\kappa}$ and $\mb{h}$ are already included in the list of optimality conditions above.
Assuming that there is a nonnegativity constraint on the demand,~\eqn{dmndnng} follows immediately as well.
What's left is to show that the conditions $\partial \mc{L} / \partial d^{rs} \geq 0$ and $d^{rs} (\partial \mc{L} / \partial d^{rs}) = 0$ correspond to~\eqn{dmnddual} and~\eqn{dmndcs}.

Calculating from~\eqn{lagrangesketch}, we have \[ \pdr{\mc{L}}{d^{rs}} = \kappa^{rs} + \pdr{f}{d^{rs}}\,, \] so if $\pdr{f}{d^{rs}} = -D^{-1}_{rs}(d^{rs})$, we are done (both equations will be true).
Integrating, $f(\mathbf{d}) = -\sum_{(r,s) \in Z^2} \int_0^{d^{rs}} D^{-1}_{rs}(\omega)~d\omega$ gives us what we need.
De-Lagrangianizing the ``no vehicle left behind'' constraint, we obtain the optimization problem associated with the elastic demand problem:
\index{Lagrange multipliers!applications}
\begin{align}
\label{eqn:elasticfirst}
\min_{\mathbf{x},\mathbf{h},\mathbf{d}}   \quad & \sum_{(i,j) \in A} \int_0^{x_{ij}} t_{ij}(x)~dx -\sum_{(r,s) \in Z^2} \int_0^{d^{rs}} D^{-1}_{rs}(\omega)~d\omega    &     \\
\mathrm{s.t.}                 \quad & x_{ij} = \sum_{\pi \in \Pi} h^\pi \delta_{ij}^{\pi} & \forall (i,j) \in A \\
                                   & \sum_{\pi \in \Pi^{rs}} h^\pi = d^{rs}              &     \forall (r,s) \in Z^2 \\
                                   & h^\pi \geq 0                                     &     \forall \pi \in \Pi \\
                                   & d^{rs} \geq 0                                    & \forall (r,s) \in Z^2
\label{eqn:elasticlast}
\end{align}
\index{static traffic assignment!elastic demand!convex optimization formulation|)}

From this convex optimization formulation, we immediately know that an elastic demand equilibrium solution exists as long as the link performance functions $t_{ij}$ and demand functions $D^{rs}$ are continuous.
If they are additionally strictly monotone (increasing for $t_{ij}$, decreasing for $D^{rs}$), then the objective function is strictly convex, and this elastic demand equilibrium is unique.
\index{static traffic assignment!elastic demand!existence of equilibrium}
\index{static traffic assignment!elastic demand!uniqueness of equilibrium}
\index{continuous function!applications}

\subsection{Solution method}
\label{sec:fwed}

\index{static traffic assignment!elastic demand!Frank-Wolfe|(}
\index{static traffic assignment!algorithms!Frank-Wolfe|(}
This section shows how the Frank-Wolfe algorithm can be used to solve the optimization problem~\eqn{elasticfirst}--\eqn{elasticlast}.
This is certainly not the only choice, and it is worthwhile for you to think about how other algorithms from Chapter~\ref{chp:solutionalgorithms} could also be used instead of Frank-Wolfe.
The implementation of this algorithm is quite similar to how Frank-Wolfe works for the basic traffic assignment problem, with three changes.
First, since the OD matrix is a decision variable along with the link flows, we must keep track of both $\mb{d}$ as well as $\mb{x}$; therefore, in addition to the target link flows $\mb{x^*}$ we will have a target OD matrix $\mb{d^*}$, and in addition when we update the link flows we must update the OD matrix as well.
Each of these pairs of flows and OD matrices should be consistent with each other, in that the link flows $\mb{x}$ must be a feasible network loading when the demand is $\mb{d}$ (both before and after updating), and similarly $\mb{x^*}$ must correspond to $\mb{d^*}$.
Luckily, this will not be difficult.

The second change to the Frank-Wolfe algorithm is how the restricted variational inequality\index{variational inequality!restricted} is solved.
Instead of solving $\mb{t(\bar{x}') \cdot (x^* - x)} \leq 0$ for $\mb{\bar{x}'} = \lambda \mb{x^*} + (1 - \lambda) \mb{x}$, $\lambda \in [0, 1]$, we must solve the variational inequality~\eqn{elasticdemandvi} in $\mb{\bar{x}'}$ and $\mb{\bar{d}'}$.
As before, the usual solution involves an ``interior'' $\lambda \in (0,1)$, in which case we must solve the equation
\begin{multline}
\label{eqn:fwequation}
\sum_{(i,j) \in A} t_{ij}(\lambda x^*_{ij} + (1 - \lambda) x_{ij}) \cdot (x^*_{ij} - x_{ij}) \\ - \sum_{(r,s) \in Z^2} D^{-1}_{rs}(\lambda d^*_{rs} + (1 - \lambda) d^{rs}) \cdot (d^*_{rs} - d^{rs}) = 0
\end{multline}
in $\lambda$.
This is simply the variational inequality~\eqn{elasticdemandvi} written out in terms of its components, substituting $\lambda \mb{x^*} + (1 - \lambda) \mb{x}$ for $\mb{\bar{x}'}$ and $\lambda \mb{d^*} + (1 - \lambda) \mb{d}$ for $\mb{\bar{d}'}$.

Third, the stopping criterion (relative gap or average excess cost) needs to be augmented with a measure of how well the OD matrix matches the values from the demand functions.
A simple measure is the \emph{total misplaced flow}\index{total misplaced flow|see {gap function, total misplaced flow}}\index{gap function!total misplaced flow} defined as $TMF = \sum_{(r,s) \in Z^2} |d^{rs} - [D^{rs}(\kappa^{rs})]^+|$.\label{not:TMF}
The total misplaced flow is always nonnegative, and is zero only if all of the entries in the OD matrix are equal to the values given by the demand function (or zero if the demand function is negative).
We should keep track of both total misplaced flow and one of the equilibrium convergence measures (relative gap or average excess cost), and only terminate the algorithm when both of these are sufficiently small.

Implementing these changes, the Frank-Wolfe algorithm for elastic demand is as follows:
\begin{enumerate}
\setcounter{enumi}{0}
\item Choose some initial OD matrix $\mb{d}$ and initial link flows $\mb{x}$ corresponding to that OD matrix.
\item Find the shortest path between each origin and destination, and calculate convergence measures (total misplaced flow, and either relative gap or average excess cost).
If both are sufficiently small, stop.
\item Improve the solution:
\begin{enumerate}
  \item Calculate a target OD matrix $\mb{d^*}$ using the demand functions: $d^*_{rs} = [D^{rs}(\kappa^{rs})]^+$ for all OD pairs $(r,s)$.
  \item Using the target matrix $\mb{d^*}$, find the link flows if everybody were traveling on the shortest paths found in step 2, store these in $\mathbf{x^*}$.
  \item Solve the restricted variational inequality by finding $\lambda$ such that~\eqn{fwequation} is true.
  \item Update the OD matrix and link flows: replace $\mb{x}$ with $\lambda \mb{x^*} + (1 - \lambda) \mb{x}$ and replace $\mb{d}$ with $\lambda \mb{d^*} + (1 - \lambda) \mb{d}$.
\end{enumerate}
\item Return to step 1.
\end{enumerate}

This algorithm can also be linked to the convex programming formulation described above.
Given a current solution $(\mb{x,d})$, it can be shown that the derivative of the objective function in the direction towards $(\mb{x^*, d^*})$ is nonpositive (and strictly negative if the current solution does not solve the elastic demand problem), and that the solution of the restricted variational inequality~\eqn{fwequation} minimizes the objective function along the line joining $(\mb{x,d})$ to $(\mb{x^*,d^*})$.
The algebra is a bit tedious and is left as an exercise at the end of the chapter.

\paragraph{Example}

Here we solve the small example of Figure~\ref{fig:2link_elastic} with the Frank-Wolfe algorithm, using the average excess cost to measure how close we are to equilibrium.
The demand function is $D(\kappa) = 50 - \kappa$, so its inverse function is $D^{-1}(d) = 50 - d$.

\begin{description}
\item[Initialization.]  Arbitrarily set $d = 50$, then arbitrarily load all 50 vehicles onto the two links; say $\mathbf{x} = \vect{50 & 0}$.
\item[Iteration 1.]  The link travel times are now $\mathbf{t} = \vect{60 & 20}$, so the shortest path travel time $\kappa = 20$ and the demand function indicates that the demand should be $D = 50 - 20 = 30$.
The average excess cost is $(60 \times 50 - 20 \times 50) / 50 = 40$, and the total misplaced flow is $|50 - 30|$ = 20.
The target demand is what the demand function indicates $d^* = 30$, and this flow should all be loaded on the bottom path, so $\mathbf{x^*} = \vect{0 & 30}$.
Solving the equation
\[ (10 + 50(1 - \lambda))(-50) + (20 + 30\lambda)(30) - (50 - (30\lambda + 50(1 - \lambda)))(-20) = 0 \]
we obtain $\lambda = 12/19 \approx 0.632$ so the new demand and flows are $d = 37.36$ and $\mathbf{x} = \vect{18.42 & 18.95}$.
\item[Iteration 2.]  The link travel times are now $\mathbf{t} = \vect{28.42 & 38.95}$, so $\kappa = 28.42$, $D = 21.58$, $AEC = 5.19$, and $TMF = 15.78$.
Both convergence measures have decreased from the first iteration, particularly the average excess cost.
Thus, $d^* = 21.58$, $\mathbf{x^*} = \vect{21.58 & 0}$, and we solve
\begin{multline*}
(10 + 21.58 \lambda + 18.42 (1 - \lambda))(21.58 - 18.42) + (20 + 18.95(1 - \lambda))(0 - 18.95) \\ - (50 - (21.58\lambda + 37.36(1 - \lambda)))(21.58 - 37.36) = 0 
\end{multline*}
so $\lambda = 0.726$ and the new demand and flows are $d = 25.9$ and $\mathbf{x} = \vect{20.71 & 5.19}$
\item[Iteration 3].
The link travel times are now $\mathbf{t} = \vect{30.71 & 25.19}$, so $\kappa = 25.19$, $D = 24.81$, $AEC = 4.42$, and $TMF = 1.09$.
Assuming that these are small enough to terminate, we are done.
\end{description}
\index{static traffic assignment!elastic demand|)}
\index{static traffic assignment!elastic demand!Frank-Wolfe|)}
\index{static traffic assignment!algorithms!Frank-Wolfe|)}

\stevefig{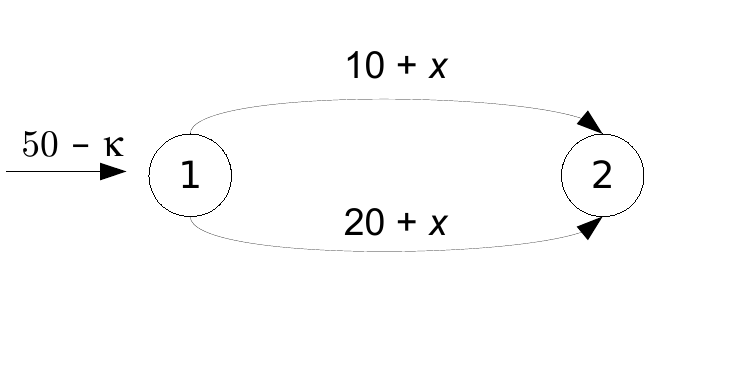}{Example network for elastic demand equilibrium.}{0.7\textwidth}

\section{Link Interactions}
\label{sec:linkinteractions}

\index{static traffic assignment!link interactions|(}
\index{link performance function!link interactions|(}
This section develops another extension of TAP, in which we relax the assumption that the travel time on a link depends only on the flow on that link.
There are a few reasons why this kind of extension may be useful:
\begin{description}
\item[Junction interactions:]  At an intersection, the delay on a particular approach often depends on flows from competing approaches.
As an example, consider a freeway onramp which has to yield to mainline traffic at a merge.
Because merging traffic must find an acceptably large gap in the main lanes, the travel time on the onramp depends on the flow on the main lanes as well as the flow on the onramp.
Similar arguments hold at arterial junctions controlled by two-way or four-way stops, at signalized intersections with permissive phases (e.g., left-turning traffic yielding to gaps in oncoming flow), or at actuated intersections where the green times are determined in real-time based on available flow.
The basic TAP cannot model the link interactions characterizing these types of links.
\item[Overtaking traffic:]  On rural highways, overtaking slow-moving vehicles often requires finding a (fairly large) gap in oncoming flow.
If the oncoming flow is small, the effect of slow-moving vehicles on average travel time is negligible.
However, as the oncoming flow becomes larger and larger, the ability to overtake is diminished and traffic speeds will tend to be determined by the slowest-moving vehicle on the highway.
Since traffic moving in different directions on the same highway is modeled with different links, a link interaction model is needed to capture this effect.
\item[Multiclass flow:]  Consider a network model where there are two types of vehicles (say, passenger cars and semi trucks, or passenger cars and buses).
Presumably these vehicles may choose routes differently or even have a different roadway network available to them; heavy vehicles are prohibited from some streets, and buses must drive along a fixed route.
This type of situation can be modeled by creating a ``two-layer'' network, with the two layers representing the links available to each class.
However, where these links represent the same physical roadway, the link performance functions should be connected to each other (truck volume influences passenger car speed and vice versa) even if they are not identical (truck speed need not be the same as passenger car speed).
Link interaction models therefore allow us to model multiclass flow as well.
\end{description}

However, there are a few twists to the story, some of which are explored below.
Section~\ref{sec:formulation} presents a mathematical formulation of the link interactions model, but shows that a convex programming formulation is not possible except in some rather unlikely cases.
Section~\ref{sec:properties} explores the properties of the link interactions model, in particular addressing the issue of uniqueness --- even though link flow solutions to TAP are unique under relatively mild assumptions, this is not generally true when there are link interactions.
Section~\ref{sec:interactionalgorithms} gives us two solution methods for the link interactions model, the diagonalization method and simplicial decomposition.
Diagonalization is easier to implement, but simplicial decomposition is generally more powerful.

\subsection{Formulation}
\label{sec:formulation}

In the basic TAP, the link performance function for link $(i,j)$ was a function of $x_{ij}$ alone, that is, we could write $t_{ij}(x_{ij})$.
Now, $t_{ij}$ may depend on the flow on multiple links.
For full generality, our notation will allow $t_{ij}$ to depend on the flows on \emph{any or all} other links in the network: the travel time is given by the function $t_{ij}(x_1, x_2, \cdots, x_{ij}, \cdots, x_m)$ or, more compactly, $t_{ij}(\mb{x})$ using vector notation.
We assume these are given to us.
Everything else is the same as in vanilla TAP: origin-destination demand is fixed, and we seek an equilibrium solution where all used paths have equal and minimal travel time.

Now, how to formulate the equilibrium principle?  It's not hard to see that the variational inequality for TAP works equally well here:\index{static traffic assignment!link interactions!variational inequality}
\labeleqn{vi}{\mb{c(\hat{h}) \cdot (\hat{h} - h)} \leq 0 \qquad \forall \mathbf{h} \in H}
where the only difference is that the link performance functions used to calculate path travel times $C$ are now of the form $t_{ij}(\mathbf{x})$ rather than $t_{ij}(x_{ij})$.
But this is of no consequence.
Path flows $\mb{\hat{h}}$ solve the variational inequality if and only if 
\labeleqn{vi_interpretation}{\mb{c(\hat{h}) \cdot \hat{h} \leq c(\hat{h}) \cdot h}}
for any other feasible path flows $\mb{h}$.
That is, if the travel times were fixed at their current values, then it is impossible to reduce the total system travel time by changing any drivers' route choices.
This is only possible if all used paths have equal and minimal travel time.
Similarly, the link-flow variational inequality
\labeleqn{vi_link}{\mb{t(\hat{x}) \cdot (\hat{x} - x)} \leq 0 \,,}
where $\mb{x}$ is any feasible link flow, also represents the equilibrium problem with link interactions.

\index{static traffic assignment!link interactions!convex optimization|(}
The ease of translating the variational inequality formulation for the case of link interactions may give us hope that a convex programming formulation exists as well.
The feasible region is the same, all we need is to find an appropriate objective function.
Unfortunately, this turns out to be a dead end.
For example, the obvious approach is to amend the Beckmann function in some way, for instance, changing $\ds \sum_{(i,j) \in A} \int_0^{x_{ij}} t_{ij}(x)~dx$ to
\labeleqn{beckmanntry}{\sum_{(i,j) \in A} \int_\mb{0}^\mb{x} t_{ij}(\mb{y})~d\mb{y}}
where the simple integral in the Beckmann function is replaced with a line integral between the origin $\mb{0}$ and the current flows $\mb{x}$.
Unfortunately, this line integral is in general not well-defined, since its value depends on the path taken between the origin and $\mb{x}$.

The one exception is if the vector of travel times $\mb{t}(\mb{x})$ is a gradient map (that is, it is a conservative vector field).
In this case, the fundamental theorem of line integrals implies that the value of this integral is \emph{independent} of the path taken between $\mb{0}$ and $\mb{x}$.
For $\mb{t(x)}$ to be a gradient map, its Jacobian must be symmetric.\index{Jacobian matrix!applications}
That is, for every pair of links $(i,j)$ and $(k,\ell)$, we need the following condition to be true:\index{link performance function!symmetric interactions}
\labeleqn{symmetry}{\pdr{t_{ij}}{x_{k\ell}} = \pdr{t_{k\ell}}{x_{ij}} }
That is, \emph{regardless of the current flow vector $\mb{x}$}, the marginal impact of another vehicle added to link $(i,j)$ on the travel time of $(k,\ell)$ must equal the marginal impact of another vehicle added to link $(k,\ell)$ on the travel time of $(i,j)$.
This condition is very strong.
Comparing with the motivating examples used to justify studying link interactions, the symmetry condition is not usually satisfied: the impact of an additional unit of flow on the mainline on the onramp travel time is much greater than the impact of an additional unit of onramp flow on mainline travel time.
The impact of semi truck flow on passenger car travel time is probably greater than the impact of passenger car flow on truck travel time at the margin.
Symmetry may perhaps hold in the case of overtaking on a rural highway, but even then it is far from clear.
So, when modeling link interactions we cannot hope for condition~\eqn{symmetry} to hold.
If it does so, consider it a happy accident: the function~\eqn{beckmanntry} is then an appropriate convex optimization problem.
\index{static traffic assignment!link interactions!convex optimization|)}

\subsection{Properties}
\label{sec:properties}

This section explores the properties of the link interaction equilibrium problem defined by the variational inequality~\eqn{vi}.
The first question concerns existence of an equilibrium.
Because~\eqn{vi} is essentially the same variational inequality derived for TAP, the arguments used to derive existence of an equilibrium (based on Brouwer's theorem)\index{Brouwer's theorem!applications} carry over directly and we have the same result:
\begin{prp}
If each link performance function $t_{ij}(\mb{x})$ is continuous in the vector of link flows $\mb{x}$, then at least one solution exists satisfying the principle of user equilibrium.
\index{static traffic assignment!link interactions!existence of equilibrium}
\end{prp}

\index{static traffic assignment!link interactions!uniqueness of equilibrium|(}
However, uniqueness turns out to be trickier.
Consider a network of two parallel links where the demand is 6 vehicles, and the link performance functions are $t_1 = x_1 + 2x_2$ and $t_2 = 2x_1 + x_2$.
Setting the two links' travel times equal to each other and using $x_1 + x_2 = 6$, it is easy to see that one equilibrium is the solution $x_1 = x_2 = 3$ when the travel time on both links is 9.
However, this is not the only equilibrium solution: if all of the drivers were to choose link 1, then $x_1 = 6$, $x_2 = 0$, $t_1 = 6$, and $t_2 = 12$.
The top link is the only used path, but it has the least travel time so this solution also satisfies the principle of user equilibrium.
Likewise, if $x_1 = 0$ and $x_2 = 6$, then $t_1 = 12$ and $t_2 = 6$ and again the only used path has the least travel time.
Therefore, this network has three equilibrium solutions; compare with Figure~\ref{fig:asymmetric1}.

\stevefig{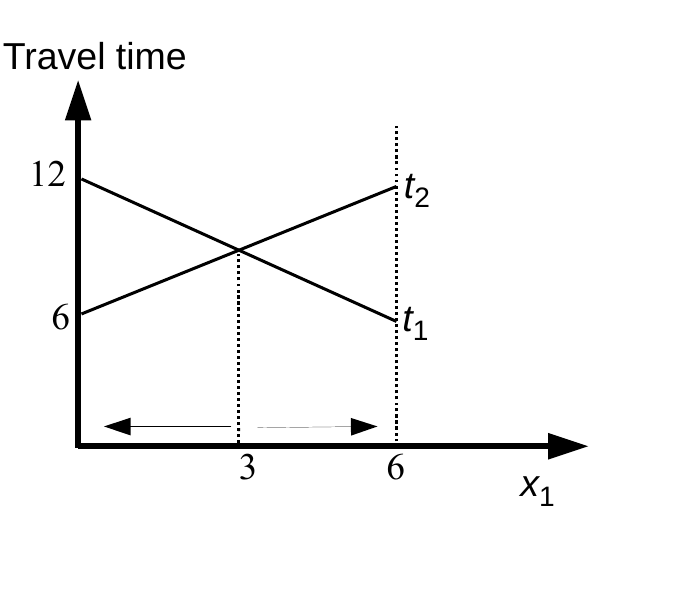}{Change in path travel times as $x_1$ varies, ``artificial'' two-link network.}{0.6\textwidth}

To make this situation less artificial, we can change the link performance functions to represent a more realistic scenario.
Assume that the rate of demand is 1800 vehicles per hour, and that link 1 has a constant travel time of 300 seconds independent of the flow on either link.
Link 2 is shorter with a free-flow time of 120 seconds, but must yield to link 1 using gap acceptance principles.
In traffic operations, gap acceptance\index{traffic flow theory!gap acceptance} is often modeled with two parameters: the \emph{critical gap}\index{traffic flow theory!critical gap} $t_c$,\label{not:tc} and the \emph{follow-up gap}\index{traffic flow theory!follow-up gap} $t_f$.\label{not:tf}
The critical gap is the smallest headway required in the main stream for a vehicle to enter.
Given that the gap is large enough for one vehicle to enter the stream, the follow-up gap is the incremental amount of time needed for each additional vehicle to enter.
For this example, let $t_c$ be 4 seconds and $t_f$ be 2 seconds.
Then, assuming that flows on both links 1 and 2 can be modeled as Poisson arrivals,\index{Poisson process} the travel time on link 2 can be derived as
\labeleqn{secondarydelay}
{t_2(x_1, x_2) = \frac{1}{u} + \frac{\Lambda}{4}
                    \left[
                       \frac{x_{2}}{u} - 1 +
                       \sqrt{\left(
                                \frac{x_{2}}{u} - 1
                             \right)^2 + \frac{8x_2}{u^2 \Lambda}
                            }
                    \right]}
where $\Lambda$\label{not:Lambda} is the length of the analysis period and $u$ is the capacity of link 2 defined by
\labeleqn{secondarycapacity}{u = \frac{x_{2} \exp(-x_{1} t_c)}{1 - \exp(-x_{1} t_f)}}
Figure~\ref{fig:asymmetric2} shows the travel times on the two paths as $x_1$ varies.
Again, there are three equilibria: (1) $x_1 = 0$, $x_2 = 1800$, where $t_1 = 300$ and $t_2 = 182$; (2) $x_1 = 362$, $x_2 = 1438$, where $t_1 = t_2 = 300$; and (3) $x_1 = 892$, $x_2 = 908$, where $t_1 = t_2 = 300$.

\genfig{asymmetric2}{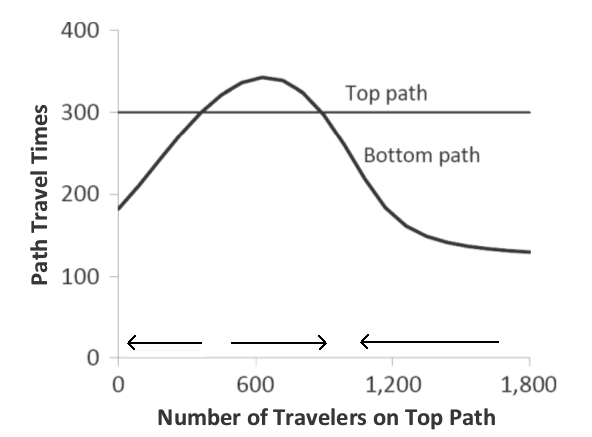}{Change in path travel times as $x_1$ varies, ``realistic'' two-link network.
Arrows indicate ``pressure'' for drivers to move to reduce travel time.
}{width=0.6\textwidth}

So, even in realistic examples we cannot expect equilibrium to be unique when there are link interactions.
The practical significance is that it raises doubt about which equilibrium solution should be used for project evaluation or ranking.
For instance, consider a candidate project which would improve the free-flow time on link 1 from 300 to a smaller value; this would correspond to lowering the horizontal line in Figure~\ref{fig:asymmetric2}.
If we are at one of the equilibria where the travel times are equal, such a change will indeed reduce the travel times experienced by drivers.
However, if we are at the equilibrium where the top path is unused, such a change will have no impact whatsoever.

While a complete study of the methods used to distinguish among multiple equilibria is beyond the scope of this section, a simple stability criterion is explained here: an equilibrium solution is \emph{stable}\index{equilibrium!stability} if small perturbations to the solution  would incentivize drivers to move back towards that initial equilibrium --- that is, if we reassign a few drivers to different paths, the path travel times will change in such a way that those drivers would want to move back to their original paths.
By contrast, an \emph{unstable} equilibrium does not have this property: if a few drivers are assigned to different paths, the path travel times will change in such a way that even more drivers would want to switch paths, and so on until another equilibrium is found.

In the simple two-link network we've been looking at, stability can be identified using graphs such as those in Figures~\ref{fig:asymmetric1} and~\ref{fig:asymmetric2}.
The arrows on the bottom axis indicate the direction of the path switching which would occur for a given value of $x_1$.
When the arrow is pointing to the left, $t_1 > t_2$ so travelers want to switch away from path 1 to path 2, resulting in a decrease in $x_1$ (a move further to the left on the graph).
When $t_2 > t_1$, travelers want to switch away from path 2 to path 1, resulting in an increase in $x_1$, indicated by an arrow pointing to the right.
At the equilibrium solutions, there is no pressure to move in any feasible direction.
So, for the first example, the only stable equilibria are the ``extreme'' solutions with all travelers on either the top or bottom link.
The equilibrium with both paths used is unstable in the sense that any shift from one path to another amplifies the difference in travel times and encourages even more travelers to shift in that direction.
In the second example, the first and third equilibria are stable, but the second is unstable.

Based on these two examples, an intuitive explanation for the presence of stability with link interactions can be provided.
For the regular traffic assignment problem with increasing link performance functions, shifting flow away from a path $\pi$ and onto another path $\pi'$ always decreases the travel time on $\pi$ and increases the travel time on $\pi'$.
Therefore, if the paths have different travel times, flow will shift in a way that always tends to equalize the travel times on the two paths.
Even where there are link interactions, the same will hold true \emph{if the travel time on a path is predominantly determined by the flow on that path}.
However, when the link interactions are very strong, the travel time on a path may depend more strongly by the flow on a different path.
In the first example, notice that each link's travel time is influenced more by the \emph{other} link's flow than its own.
In the second example, for certain ranges of flow the travel time on the merge path is influenced more by the flow on the priority path.
In such cases, there is no guarantee that moving flow from a higher-cost path to a lower-cost path will tend to equalize their travel times.
In the first example, we have an extreme case where moving flow to a path \emph{decreases} its travel time while \emph{increasing} the travel time of the path the flow moved away from!

To make this idea more precise, the following section introduces the mathematical concept of strict monotonicity.

\subsubsection{Strict Monotonicity}

\index{strictly monotone function|(}
Let $\mb{f}(\mb{x})$ be a vector-valued function whose domain and range are vectors of the same dimension.
For instance, $\mb{t}(\mb{x)}$ maps the vector of link flows to the vector of link travel times; the dimension of both of these is the number of links in the network.
We say that $\mb{f}$ is \emph{strictly monotone} if for any two distinct vectors $\mb{x}$ and $\mb{y}$ in its domain, the dot product of $\mb{f(x) - f(y)}$ and $\mb{x - y}$ is strictly positive.

For example, let $\mb{f}(\mb{x})$ be defined by $f_1(x_1, x_2) = 2x_1$ and $f_2(x_1, x_2) = 2x_2$.
Then for any distinct vectors $\mb{x}$ and $\mb{y}$, we have
\begin{multline*}
 \mb{(f(x) - f(y)) \cdot (x - y)} = (\vect{2x_1 & 2x_2} - \vect{2y_1 & 2y_2}) \cdot (\vect{x_1 & x_2} - \vect{y_1 & y_2}) \\ = 2((x_1 - y_1)^2 + (x_2 - y_2)^2)  
\end{multline*}

Since $\mb{x} \neq \mb{y}$, the right-hand side is always greater than zero, so $f$ is monotone.
As another example, let the function $\mb{g}(\mb{x})$ be defined by $g_1(x_1, x_2) = 2x_2$ and $g_2(x_1, x_2) = 2x_1$.
If we choose $\mb{x} = \vect{0 & 1}$ and $\mb{y} = \vect{1 & 0}$, then
\[ \mb{(g(x) - g(y)) \cdot (x - y)} = (\vect{2 & -2}) \cdot (\vect{-1 & 1}) = -4 \]
so $\mb{g}$ is not strictly monotone.
Note that proving strict monotonicity requires a general argument valid for \emph{any} distinct vectors $\mb{x}$ and $\mb{y}$; showing that a function is not strictly monotone only requires a single counterexample.

\textbf{Warning!}  It is very common for students to think that the link performance functions are strictly monotone if they are strictly increasing functions of the flow on each link.
This is not true: in the first example in this section, all link performance functions are increasing in each flow variable but if we compare $\mb{x} = \vect{0 & 6}$ and $\mb{y} = \vect{5 & 1}$, we have
\[ (\mb{t(x) - t(y)}) \cdot (\mb{x - y}) = (\vect{12 & 6} - \vect{7 & 11}) \cdot (\vect{0 & 6} - \vect{5 & 1}) = -50 \]
so these link performance functions are not strictly monotone.
Roughly speaking, strict monotonicity requires the diagonal terms of the Jacobian of $\mb{t}$ to be large compared to the off-diagonal terms.
The precise version of this ``roughly speaking'' fact is the following:

\begin{prp}
\label{prp:strictmonotoneiffpd}
If $\mb{f}$ is a continuously differentiable function whose domain is convex, then $\mb{f}$ is strictly monotone if and only if its Jacobian is positive definite at all points in the domain.
\index{Jacobian matrix!applications}
\index{positive definite matrix!applications}
\index{continuous function!applications}
\index{convex set!applications}
\end{prp}

With this definition of monotonicity in hand, we can provide the uniqueness result we've been searching for:

\begin{prp}
Consider an instance of the traffic assignment problem with link interactions.
If the link performance functions $\mb{t(x)}$ are continuous and strictly monotone, then there is exactly one user equilibrium solution.
\index{continuous function!applications}
\end{prp}
\begin{proof}
Since $\mb{t(x)}$ is continuous, we are guaranteed existence of at least one equilibrium solution from Brouwer's theorem; let $\mb{\hat{x}}$ be such an equilibrium and let $\mb{\tilde{x}}$\label{not:xtwiddle} be any other feasible link flow solution.\index{Brouwer's theorem!applications}
We need to show that $\mb{\tilde{x}}$ cannot be an equilibrium.
Arguing by contradiction, assume that $\mb{\tilde{x}}$ is in fact an equilibrium.
Then it would solve the variational inequality~\eqn{vi_link}, so\index{variational inequality!applications}
\[ \mb{t(\tilde{x}) \cdot (\tilde{x} - \hat{x})}\leq 0 \]
Adding a clever form of zero to the left hand side, this would imply
\labeleqn{uniqueproof}{ \mb{(t(\tilde{x}) - t(\hat{x})) \cdot (\tilde{x} - \hat{x}) + t(\hat{x}) \cdot (\tilde{x} - \hat{x})} \leq 0 }
But since the link performance functions are strictly monotone, the first term on the left-hand side is strictly positive.
Furthermore, since $\mb{\hat{x}}$ is an equilibrium the variational inequality~\eqn{vi_link} is true, so $\mb{t(\hat{x}) \cdot (\hat{x} - \tilde{x})} \leq 0$, which implies that the second term on the right-hand side is nonnegative.
Therefore, the left-hand side of~\eqn{uniqueproof} is strictly positive, which is a contradiction.
Therefore $\mb{\tilde{x}}$ cannot satisfy the principle of user equilibrium.
\end{proof}
\index{strictly monotone function|)}
\index{static traffic assignment!link interactions!uniqueness of equilibrium|)}

\subsection{Algorithms}
\label{sec:interactionalgorithms}

This section presents two algorithms for the traffic assignment problem with link interactions.
If the link performance functions are strictly monotone, it can be shown that both of these algorithms converge to the unique equilibrium solution.
Otherwise, it is possible that these algorithms may not converge, although they will typically do so if they start sufficiently close to an equilibrium.
In any case, these algorithms may be acceptable heuristics even when strict monotonicity does not hold.

\subsubsection{Diagonalization}

\index{static traffic assignment!link interactions!diagonalization algorithm|(}
The diagonalization method is a variation of Frank-Wolfe, which differs only in how the step size $\lambda$ is found.
Recall that the Frank-Wolfe step size is found by solving the equation
\labeleqn{fwstepsize}{\sum_{(i,j) \in A} t_{ij}(\lambda x^*_{ij} + (1 - \lambda) x_{ij}) (x^*_{ij} - x_{ij}) = 0}
since this minimizes the Beckmann function along the line segment connecting $\mb{x}$ to $\mb{x^*}$.
Since there is no corresponding objective function when there are asymmetric link interactions, it is not clear that a similar approach will necessarily work.
(And in any case, $t_{ij}$ is no longer a function of $x_{ij}$ alone, so the formula as stated will not work.)

To make this formula logical, construct a temporary link performance function $\tilde{t}_{ij}(x_{ij})$ which only depends on its own flow.
This is done by assuming that the flow on all other links is constant: $\tilde{t}_{ij}(x_{ij}) = t_{ij}(x_1, \ldots, x_{ij}, \ldots, x_m)$.
For example, if $t_1(x_1, x_2, x_3) = x_1 + x_2^2 + x_3^3$ and the current solution is $x_1 = 1$, $x_2 = 2$, and $x_3 = 3$, then $\tilde{t}_1(x_1) = 31 + x_1$, since this is what we would get if $x_2$ and $x_3$ were set to constants at their current values of 2 and 3, respectively.

The step size $\lambda$ is then found by adapting the Frank-Wolfe formula, using $\tilde{t}_{ij}(x_{ij})$ in place of $t_{ij}(\mathbf{x})$.\index{static traffic assignment!algorithms!Frank-Wolfe}
That is, in the diagonalization method $\lambda$ solves
\labeleqn{diagstepsize}{\sum_{(i,j) \in A} \tilde{t}_{ij}(\lambda x^*_{ij} + (1 - \lambda) x_{ij}) (x^*_{ij} - x_{ij}) = 0}

At each iteration, new $\tilde{t}$ functions are calculated based on the current solution.
The complete algorithm is as follows:
\begin{enumerate}
\item Find the shortest path between each origin and destination, and calculate the relative gap (unless it is the first iteration).
If the relative gap is sufficiently small, stop.
\item Shift travelers onto shortest paths:
\begin{enumerate}[(a)]
  \item Find the link flows if everybody were traveling on the shortest paths found in step 1, store these in $\mathbf{x^*}$.
  \item If this is the first iteration, set $\mathbf{x} \leftarrow \mathbf{x^*}$ and move to step 3.
Otherwise, continue with step c.
  \item Using the current solution $\mathbf{x}$, form the diagonalized link performance functions $\tilde{t}_{ij}(x_{ij})$ for each link.
  \item Find $\lambda$ which solves equation~\eqn{diagstepsize}.
  \item Update $\mathbf{x} \leftarrow \lambda \mathbf{x^*} + (1 - \lambda) \mathbf{x}$.
\end{enumerate}
\item Calculate the new link travel times and the relative gap.
Increase the iteration counter $i$ by one and return to step 1.
\end{enumerate}

As an example, consider the modified Braess network shown in Figure~\ref{fig:asymmetricbraess}.
At each merge node, the travel time on each of the incoming links depends on the flow on \emph{both} links which merge together.
The link flow vectors are indexed $\mb{x} = \vect{x_{12} & x_{13} & x_{23} & x_{24} & x_{34}}$.

\stevefig{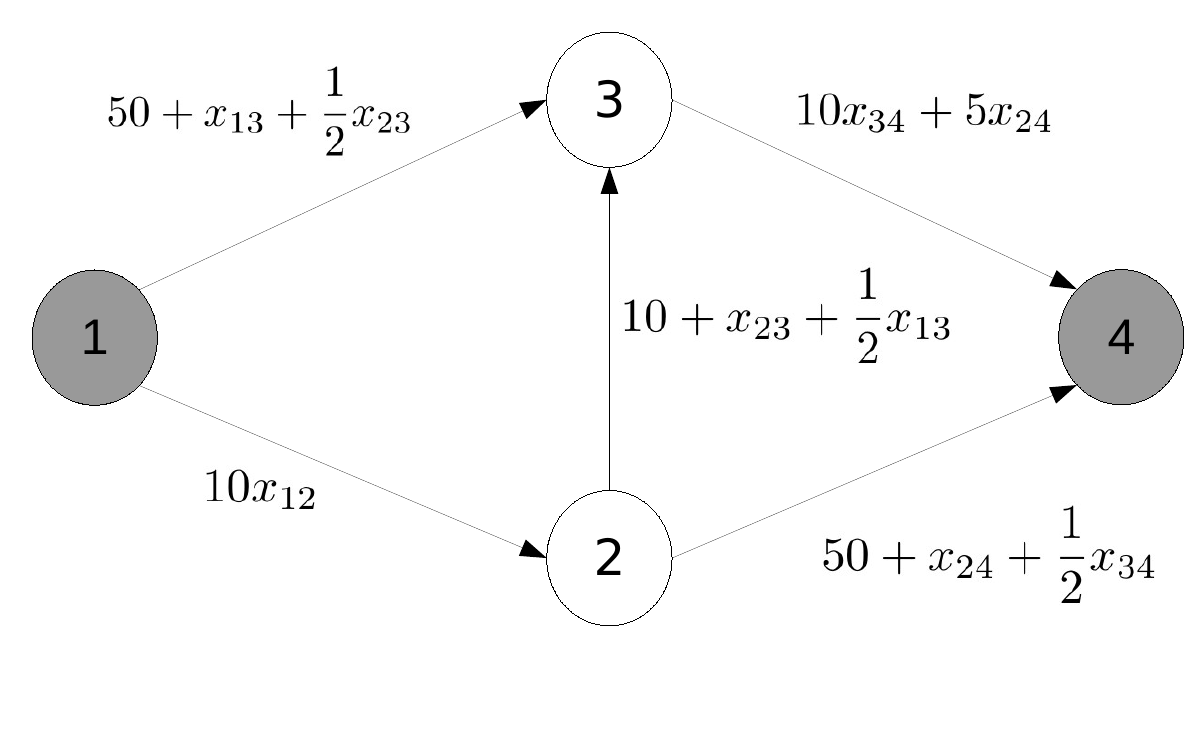}{Braess network with link interactions.}{0.7\textwidth}

\begin{description}
\item[Iteration 1.]  In the first iteration, load all travelers onto shortest paths, so $\mb{x} = \mb{x^*} = \vect{6 & 0 & 6 & 0 & 6}$, $\mb{t} = \vect{60 & 53 & 16 & 53 & 60}$ and the average excess cost is 23.
\item[Iteration 2.]  The all-or-nothing loading on shortest paths given $\mb{t}$ is $\mb{x^*} = \vect{0 & 6 & 0 & 0 & 6}$.
The diagonalized link performance functions are obtained by assuming the flows on all other links are constant at $\mb{x}$: $\hat{t}_{12} = 10x_{12}$, $\hat{t}_{13} = 53 + x_{13}$, $\hat{t}_{23} = 10 + x_{23}$, $\hat{t}_{24} = 53 + x_{24}$, and $\hat{t}_{34} = 10x_{34}$.
So we solve the equation~\eqn{diagstepsize} for $\lambda$:
\[ (53 + 6 \lambda) 6 + 10(6(1 - \lambda))(-6) + (10 + 6(1 - \lambda))(-6) = 0\]
omitting terms where $x_{ij} = x^*_{ij}$ because they are zero in~\eqn{diagstepsize}.
The solution is $\lambda = 23/72$, so we update $\mb{x} = \vect{4 \frac{1}{12} & 1 \frac{11}{12} & 4 \frac{1}{12} & 0 & 6}$ and $\mb{t} = \vect{40 \frac{5}{6} & 53 \frac{23}{24} & 15 \frac{1}{24} & 53 & 60}$ (using the regular cost functions, not the diagonalized ones.)  The average excess cost is now 21.43.
\item[Iteration 3.] The all-or-nothing assignment is $\mb{x^*} = \vect{6 & 0 & 0 & 6 & 0}$ and the diagonalized link performance functions are $\hat{t}_{12} = 10x_{12}$, $\hat{t}_{13} = 52\frac{1}{24} + x_{13}$, $\hat{t}_{23} = 10\frac{23}{24} + x_{23}$, $\hat{t}_{24} = 53 + x_{24}$, and $\hat{t}_{34} = 10x_{34}$.
The solution to~\eqn{diagstepsize} is $\lambda = 0.284$, which gives $\mb{x} = \vect{4.63 & 1.37 & 2.92 & 1.70 & 4.30}$, $\mb{t} = \vect{46.3 & 52.8 & 13.6 & 53.9 & 51.5}$, so the average excess cost is 6.44.
\end{description}
and so on until convergence is reached.
\index{static traffic assignment!link interactions!diagonalization algorithm|)}

\subsubsection{Simplicial decomposition}
\label{sec:staticsimplicial}

\index{static traffic assignment!algorithms!simplicial decomposition|(}
\index{static traffic assignment!link interactions!simplicial decomposition algorithm|(}
An alternative to diagonalization is the simplicial decomposition algorithm.
This algorithm is introduced at this point (rather than in Chapter~\ref{chp:solutionalgorithms}) for several reasons.
First, it was historically the first provably convergent algorithm for the equilibrium problem with link interactions.
Second, although it is an improvement on Frank-Wolfe, for the basic TAP it is outperformed by the path-based and bush-based algorithms presented in that chapter.
However, like those algorithms, it overcomes the ``zig-zagging'' difficulty that Frank-Wolfe runs into (cf.\ Figure~\ref{fig:fwfail}).

The price of this additional flexibility is that more computer memory is needed.
Frank-Wolfe and the method of successive averages are exceptionally economical in that they only require two vectors to be stored: the current link flows $\mathbf{x}$ and the target link flows $\mathbf{x^*}$.
In simplicial decomposition, we will ``remember'' all of the target link flows found in earlier iterations, and exploit this longer-term memory by allowing ``combination'' moves towards several of these previous targets simultaneously.
In the algorithm, the set $\mc{X}$\label{not:fancyX} is used to store all target link flows found thus far.

A second notion in simplicial decomposition is that of a ``restricted equilibrium.''  Given a set $\mc{X} = \myc{\mb{x^*_1}, \mb{x^*_2}, \cdots, \mb{x^*_k}}$ and a current link flow solution $\mathbf{x}$, we say that $\mathbf{x}$ is a restricted equilibrium\index{equilibrium!restricted}\index{restricted equilibrium} if it solves the variational inequality
\labeleqn{restrictedvi}{\mb{t(x) \cdot (x - x')} \leq 0 \qquad \forall \mb{x'} \in X(\mb{x}, \mathcal{X})}
where $X(\mb{x}, \mathcal{X})$\label{not:XxX} means the set of link flow vectors which are obtained by a convex combination of $\mb{x}$ and any of the target vectors in $\mathcal{X}$.\footnote{A flow vector $\mb{x'}$ is a convex combination\index{convex combination} of $\mb{x}$ and the target vectors in $\mathcal{X}$ if there exist nonnegative constants $\lambda_0, \lambda_1, \cdots, \lambda_k$ such that $\mb{x'} = \lambda_0 \mb{x} + \lambda_1 \mb{x^*_1} + \cdots \lambda_k \mb{x^*_k}$. 
}
This is a ``restricted variational inequality''\index{variational inequality!restricted} like that used for Frank-Wolfe in Section~\ref{sec:frankwolfe}, but where the feasible set now consists of combinations of \emph{all} of the vectors in $\mc{X}$, rather than just the line segment connecting $\mb{x}$ and $\mb{x^*}$.

Equivalently, $\mathbf{x}$ is a restricted equilibrium if none of the targets in $\mc{X}$ lead to improving directions in the sense that the total system travel time would be reduced by moving to some $\mb{x^*_i} \in \mc{X}$ while fixing the travel times at their current values.
That is,
\labeleqn{extremepoint}{\mb{t(x) \cdot (x - x^*_i)} \leq 0 \qquad \forall \mb{x^*_i} \in \mc{X}}

At a high level, simplicial decomposition works by iterating between adding new target vectors to $\mc{X}$, and then finding a restricted equilibrium using the current vectors in $\mc{X}$.
Frank-Wolfe can be seen as a special case of simplicial decomposition, where $\mc{X}$ only consists of the current target vector (forgetting any from past iterations).

In practice, it is too expensive to exactly find a restricted equilibrium at each iteration.
Instead, several ``inner iteration'' steps are taken to move towards a restricted equilibrium with the current set $\mc{X}$ before looking to add another target.
In each inner iteration, the current solution $\mb{x}$ is adjusted to $\mb{x} + \mu \Delta \mb{x}$, where $\mu$ is a step size and $\Delta \mb{x}$ is a direction which moves toward restricted equilibrium.
Smith (1984) shows that one good choice for this direction is
\labeleqn{smith}{\Delta \mathbf{x} = \frac{\sum_{\mathbf{x^*_i} \in \mathcal{X}} \left[ \mathbf{t(x)} \cdot \mathbf{(x - x^*_i)} \right]^+ \mathbf{(x^*_i - x)}}{\sum_{\mathbf{x^*_i} \in \mathcal{X}} \left[ \mathbf{t(x)} \cdot \mathbf{(x - x^*_i)} \right]^+ }}
This rather intimidating-looking formula is actually quite simple.
It is nothing more than a weighted average of the directions $\mathbf{x^*_i - x}$ (potential moves toward each target in $\mc{X}$), where the weight for each potential direction is the extent to which it improves upon the current solution: $\left[ \mathbf{t(x)} \cdot \mathbf{(x - x^*_i)} \right]$ is the reduction in total system travel time obtained by moving from $\mb{x}$ to $\mb{x^*_i}$ while holding travel times constant.
If this term is negative, there is no need to move in that direction, so the weight is simply set to zero.
The denominator is simply the sum of the weights, which serves as a normalizing factor.

The step size $\mu$ is chosen through trial-and-error.
One potential strategy is to iteratively test $\mu$ values in some sequence (say, $1, 1/2, 1/4, \ldots$) until we have found a solution acceptably closer to restricted equilibrium than $\mb{x}$.\footnote{This is similar to the Armijo\index{Armijo rule} rule described in Appendix~\ref{sec:unconstrainedstepsize}.
}
``Acceptably closer'' can be calculated using the \emph{Smith gap}\index{gap function!Smith gap|(}\index{Smith gap|see {gap function, Smith gap}}\label{not:gammaS}
\labeleqn{raec}{\gamma_S = \sum_{\mathbf{x^*_i} \in \mc{X}} ( [\mb{t(x) \cdot (x - x^*_i)} ]^+)^2}
which is similar to the gap measures described in Chapter~\ref{chp:solutionalgorithms} in that it is zero if and only if $\mb{x}$ is a restricted equilibrium, and positive otherwise.
It reflects a restricted equilibrium because the summation is only over the vectors in $\mc{X}$: if any of them reflect an ``improvement'' over the current solution, in that total travel time would be reduced (if link travel times were held fixed), the corresponding term in square brackets is positive, and $\gamma_S > 0$.
It can be shown that taking a small enough step in the direction~\eqn{smith} will reduce the Smith gap, and a line search or trial-and-error method can be used to determine what an acceptable step size is.
Squaring the term in brackets ensures that $\gamma_S$ is differentiable, which plays a role in deriving the convergence rate of this method.

The Smith gap could thus be converted into a stopping criterion for the (unrestricted) traffic assignment problem, by extending the sum to include every possible all-or-nothing assignment, not just the ones in $\mc{X}$, but this idea will not be explored further in this book.
\index{gap function!Smith gap|)}

Putting all of this together, the simplicial decomposition algorithm can be stated as:
\begin{enumerate}
\item Initialize the set $\mathcal{X} \leftarrow \emptyset$
\item Find shortest paths for all OD pairs.
\item Form the all-or-nothing assignment $\mathbf{x^*}$ based on shortest paths.
\item If $\mathbf{x^*}$ is already in $\mathcal{X}$, stop.
\item Add $\mathbf{x^*}$ to $\mathcal{X}$.
\item \textbf{Subproblem:} Find a restricted equilibrium $\mathbf{x}$ using only the vectors in $\mathcal{X}$.
\begin{enumerate}[(a)]
 \item Find the improvement direction $\Delta \mathbf{x}$ using equation~\eqn{smith}.
 \item Update $\mathbf{x} \leftarrow \mathbf{x} + \mu \Delta \mathbf{x}$, with $\mu$ sufficiently small (to reduce $\gamma_S$).
 \item Update travel times.
 \item Return to step 1 of subproblem unless $\gamma_S$ is small enough.
\end{enumerate}
\item Return to step 2.
\end{enumerate}

Below we apply this algorithm to the example in Figure~\ref{fig:asymmetricbraess}, choosing $\mu$ via trial and error from the sequence $1/2, 1/4, \ldots$ and stopping at the first value that reduces $\gamma_S$.
Other ways of choosing $\mu$ are also possible.

\begin{description}
\item[Iteration 1.]  We set $$\mb{x^*_1} = \vect{6 & 0 & 6 & 0 & 6}$$ and $\mc{X} = \mb\{x^*_1\}$.
For the subproblem, the only possible solution is $\mb{x} = \mb{x^*_1}$, which has $\gamma_S = 0$ (it is trivially a restricted equilibrium) and travel times are $$\mb{t} = \vect{60 & 53 & 16 & 53 & 60}\,.$$
\item[Iteration 2.]  The new all-or-nothing assignment is $$\mb{x^*_2} = \vect{0 & 6 & 0 & 0 & 6}\,,$$ and $\mc{X} = \myc{\mb{x^*_1, x^*_2}}$.
For the first iteration of the subproblem, notice that $\mb{t} \cdot (\mb{x - x^*_1}) = 0$ and $\mb{t} \cdot (\mb{x - x^*_2}) = 138$, so Smith's formula~\eqn{smith} reduces to $$\Delta \mb{x} = \frac{0}{138} (\mb{x^*_1 - x}) + \frac{138}{138} (\mb{x^*_2 - x}) = \vect{-6 & 6 & -6 & 0 & 0}\,.$$  Taking a step of size $\mu = 1/2$ gives us $$\mb{x} = \vect{3 & 3 & 3 & 0 & 6}\,.$$  The new travel times are $$\mb{t} = \vect{30 & 54\frac{1}{2} & 14\frac{1}{2} & 53 & 60}\,,$$ so $\mb{t} \cdot \mb{x} = 657$, $\mb{t} \cdot \mb{x^*_1} = 627$, and $\mb{t} \cdot \mb{x^*_2} = 687$, so
\[ \gamma_S = ([657-627]^+)^2 + ([657-687]^+)^2 = 30^2 + 0^2 = 900 \,.\]
Assume this is ``small enough'' to complete the subproblem.
\item[Iteration 3.]  The new all-or-nothing assignment is $$\mb{x^*_3} = \vect{6 & 0 & 0 & 6 & 0}$$ and $\mc{X} = \myc{\mb{x^*_1, x^*_2, x^*_3}}$.
For the first iteration of the subproblem, calculate $\mb{t} \cdot (\mb{x - x^*_1}) = 30$, $\mb{t} \cdot (\mb{x - x^*_2}) = -30$, and $\mb{t} \cdot (\mb{x - x^*_3}) = 159$.
So~\eqn{smith} gives
\begin{multline*}
\Delta \mb{x} = \frac{30}{189} (\mb{x^*_1 - x}) + \frac{0}{189} (\mb{x^*_2 - x}) + \frac{159}{189} (\mb{x^*_3 - x}) \\ = \vect{3 & -3 & -2.05 & 5.05 & -5.05}\,.
\end{multline*}
Taking a step of size $\mu = 1/2$ would give $$\mb{x} = \vect{4.5 & 1.5 & 1.98 & 2.52 & 3.48}$$ and $$\mb{t} = \vect{45 & 52.5 & 12.7 & 54.3 & 47.4}$$ which has $\gamma_S = 1308$.
Assume that this is no longer ``small enough'' to return to the master problem, so we begin a second subproblem iteration.
Smith's formula~\eqn{smith} now gives $\Delta \mb{x} = 0.2 (\mb{x - x^*_2}) + 0.8 (\mb{x - x^*_3})$, and the trial solution $\mb{x} + \frac{1}{2} \Delta \mb{x}$ has $\gamma_S = 1195$, which is an improvement.
In this case, choosing a smaller $\mu$ would work even better; for instance $\mu = 1/4$ would reduce the Smith gap to 557.
There is thus a tradeoff between spending more time on finding the ``best'' value of $\mu$, or spending more time on finding new search directions and vectors for $\mc{X}$.
Balancing these is an important question for implementation.
\end{description}

The algorithm can continue from this point or terminate if this average excess cost is small enough.
\index{static traffic assignment!algorithms!simplicial decomposition|)}
\index{static traffic assignment!link interactions!simplicial decomposition algorithm|)}
\index{static traffic assignment!link interactions|)}
\index{link performance function!link interactions|)}

\section{Stochastic User Equilibrium}
\label{sec:sue}

\index{static traffic assignment!stochastic user equilibrium|(}
This section describes another extension to the basic TAP.
To this point in the text, we have been using the principle of user equilibrium to determine link and path flows, requiring all used paths between the same origin and destination to have equal and minimal travel time.
We derived this principle by assuming that all drivers choose the least-travel time path between their origin and destination.
However, this assumption implicitly requires drivers to have perfect knowledge of the travel times on all routes in the network.
In reality, we know this is not true: do you know the travel times on literally \emph{all} routes between an origin and destination?  And can you accurately distinguish between a route with a travel time of 16 minutes, and one with a travel time of 15 minutes and 59 seconds?  Relaxing these assumptions leads us to the \emph{stochastic user equilibrium} (SUE) model.

In SUE, rather than requiring that each driver follow the true shortest path between each origin and destination, we assume that drivers follow the path they \emph{believe} to be shortest, but allow for some perception error between their belief and the actual travel times.
An alternative, mathematically equivalent, interpretation (explained below) is that drivers do in fact perceive travel times accurately, but care about factors other than travel time.
This section explains the development of the SUE model.
The mathematical foundation for the SUE model is in discrete choice concepts, which are briefly reviewed in Section~\ref{sec:discretechoice}.
The specific application of discrete choice to the route choice decision is taken up in Section~\ref{sec:logitchoice}.
These sections address the ``individual'' perspective of logit route choice.

The next steps to creating the SUE model are an efficient network loading model (a way to find the choices of \emph{all} drivers efficiently), and then finally the equilibrium model which combines the network loading with updates to travel times, to account for the mutual dependence between travel times and route choices.\index{consistency}
These are undertaken in Sections~\ref{sec:stochasticloading} and~\ref{sec:sueformulation}, respectively.
For the most part, this discussion assumes a relatively simplistic model for perception errors in travel times; Section~\ref{sec:generalchoice} briefly discusses how more general situations can be handled.

\subsection{Discrete choice modeling}
\label{sec:discretechoice}

\index{discrete choice|(}
This section provides a brief overview of discrete choice concepts.
The application to route choice is in the following section.
Discrete choice is a large area of scholarly inquiry in and of itself, and so the discussion is restricted to what is needed in the chapter.

Consider an individual who must make a choice from a set of options.
For instance, when purchasing groceries, you must choose one store from a set of alternatives (the grocery stores in your city).
When dressing in the morning, you must choose one set of clothes among all of the clothing you own.
And, more relevant to transportation, when choosing to drive from one point to another, you must choose one route among all of the possible routes connecting your origin to your destination.

Mathematically, let $I$\label{not:I} be a finite set of alternatives.
Each alternative $i \in I$ is associated with a \emph{utility}\index{utility} value $U_i$ representing the amount of happiness or satisfaction you would have if you were to choose option $i$.
We assume that you would choose an alternative $i^* \in \arg \max_i \myc{U_i}$ which maximizes the utility you receive.
Now, the utility $U_i$ consists of two parts: an \emph{observable utility} $V_i$, and an \emph{unobservable utility} denoted by the random variable $\epsilon$:
\labeleqn{discretechoice}{U_i = V_i + \epsilon_i}

\index{utility!observable and unobservable|(}
The difference between observable and unobservable utility can be explained in different ways.
One interpretation is that $V_i$ represents the portion of the utility that is due to objective factors visible to the modeler; when choosing a grocery store, that might include the distance from your home, the price, the variety of items stocked, etc.
The unobserved utility consists of subjective factors that the modeler cannot see (or chooses not to include in the model), even though they are real insofar as they affect your choice.
In the grocery store example, this might include your opinion on the taste of the store brands, the cleanliness of the store, and so on.
Then, by modeling the unobserved utility as a random variable $\epsilon$, we can express choices in terms of probabilities.
(The modeler does not know all of the factors affecting your choice, so they can only speak of probabilities of choosing different options based on what is observable.)

A second interpretation is that the observed utility $V_i$ actually represents all of the factors that you care about.
However, for various reasons you are incapable of knowing all of these reasons with complete accuracy.
(You probably have a general sense of the prices of items at a grocery store, but very few know the exact price of every item in a store's inventory.)  Then the random variable $\epsilon_i$ represents the \emph{error} between the true utility ($V_i$) and what you believe the utility to be ($U_i$).
Either interpretation leads to the same mathematical model.
\index{utility!observable and unobservable|)}

\index{logit|(}
Depending on the distribution we choose for the random variables $\epsilon_i$, different discrete models are obtained.
A classic is the \emph{logit} model, which is obtained when the unobserved utilities $\epsilon_i$ are assumed to be independent across alternatives, and to have Gumbel distributions\index{Gumbel distribution} with zero mean and identical variance.
Under this assumption, the probability of choosing alternative $i$\label{not:p_i} is given by
\labeleqn{logit}{p_i = \frac{\exp(\theta V_i)}{\sum_{j \in I} \exp(\theta V_j)}}
where $\theta$ is a nonnegative parameter reflecting the relative magnitude of the observed utility relative to the unobserved utility.
Notice what happens in this formula as $\theta$ takes extreme values: if $\theta = 0$, then all terms in the numerator and denominator are 1, and the probability of choosing any alternative is exactly the same.
(Interpretation: the unobserved utility $\epsilon$ is much more important than the observed utility, so the observed utility has no impact on the decision made.
Since the unobserved utility has the same distribution for every alternative, each is equally likely.)  Or, if $\theta$ grows large, then the denominator of~\eqn{logit} will be dominated by whichever terms have the largest observed utility $V_i$.
If there is some alternative $i^*$ for which $V_{i^*} > V_i$ for all $i \neq i^*$ (the observed utility for $i^*$ is strictly greater than any other alternative), then as $\theta \rightarrow \infty$, the probability of choosing $i^*$ approaches 1 and the probability of choosing any other alternative approaches 0.
Another important consequence of~\eqn{logit} is that the probability of choosing any alternative is always strictly positive; there is \emph{some} chance that the unobserved utility will be large enough that any option could be chosen.

While the logit model is nice in that we have a closed-form expression~\eqn{logit} for the probabilities of choosing any alternative, the logit assumptions are very strong --- particularly the assumption that the $\epsilon_i$ are independent and identically distributed.
The exercises at the end of the chapter explore some examples demonstrating how these assumptions can lead to unreasonable results.
Another common assumption is that the $\epsilon$ are drawn from a multivariate normal distribution, which allows for correlation among the unobserved utilities for different alternatives.
This leads to the \emph{probit}\index{probit} choice model, which is more flexible and arguably realistic.
However, unlike the logit model, the probit model does not have a closed-form expression for probabilities like~\eqn{logit}.
Instead, Monte Carlo sampling methods are used to estimate choices.

The majority of this section is focused on logit-based models.
While probit models are more general and arguably more realistic, the logit model has two major advantages from the perspective of a book like this.
First, computations in logit models can often be done analytically, simplifying explanations and making it possible to give examples you can easily verify.
This helps you better understand the main ideas in stochastic user equilibrium and build intuition.
Second, logit models admit faster solution algorithms, algorithms which scale relatively well with network size.
This is an important practical advantage for logit models.
Nevertheless, Section~\ref{sec:generalchoice} provides some discussion on probit and other models and what needs to change from the logit discussion below.
\index{logit|)}
\index{discrete choice|)}

\subsection{Logit route choice}
\label{sec:logitchoice}

This section specializes the discrete choice framework from the previous section to route choice in networks.
Consider a traveler leaving origin $r$ for destination $s$.
They must choose one of the paths $\pi$ connecting $r$ to $s$, that is, they must make a choice from the set $\Pi_{rs}$.
The most straightforward way to generalize the principle of user equilibrium to account for perception errors is to set the observed utility equal to the negative of path travel time, so
\labeleqn{sueutility}{U^\pi = -c^\pi + \epsilon^\pi}
with the negative sign indicating that \emph{maximizing utility} for drivers means \emph{minimizing travel time}.
Assuming that the $\epsilon^\pi$ are independent, identically distributed Gumbel random variables, we can use the logit formula~\eqn{logit} to express the probability that path $\pi$ is chosen:
\labeleqn{logitpath}{p^\pi = \frac{\exp(-\theta c^\pi)}{\sum_{\pi' \in \Pi^{rs}} \exp(-\theta c_{\pi'})}}
The comments in the previous section apply to the interpretation of this formula.
As $\theta$ approaches 0, drivers' perception errors are large relative to the path travel times, and each path is chosen with nearly equal probability.
(The errors are so large, the choice is essentially random.)  As $\theta$ grows large, perception errors are small relative to path travel times, and the path with lowest travel time is chosen with higher and higher probability.
At any level of $\theta$, there is a strictly positive probability that each path will be taken.

For concreteness, the route choice discussion so far corresponds to the second interpretation of SUE, where the unobserved utility represents perception errors in the utility.
The first interpretation would mean that $\epsilon^\pi$ represents factors other than travel time which affect route choice (such as comfort, quality of scenery, etc.).
Either of these interpretations is mathematically consistent with the discussion in Section~\ref{sec:discretechoice}.

The fact that the denominator of~\eqn{logitpath} includes a summation over \emph{all} paths connecting $r$ to $s$ is problematic, from a computation standpoint.
The number of paths can grow exponentially with network size.
Any use of stochastic user equilibrium in a practical setting, therefore, requires a way to compute link flows without explicitly calculating the sum in~\eqn{logitpath}.

\index{reasonable paths|(}
This is done by carefully defining which paths are in the choice set for travelers.
The notation $\Pi^{rs}$ in~\eqn{logitpath} in this book means the set of all acyclic paths connecting origin $r$ to destination $s$.
Following Chapter~\ref{chp:solutionalgorithms}, we will use the notation $\hat{\Pi}^{rs}$ to define the set of paths being considered by travelers in SUE; these sets are sometimes called sets of \emph{reasonable paths}.
With a suitable definition of this set, using $\hat{\Pi}^{rs}$ in place of $\Pi^{rs}$ in equation~\eqn{logitpath} leads to tractable computation schemes.

Two possibilities are common: selecting an acyclic subset of links, and choosing $\hat{\Pi}^{rs}$ to contain the paths using these links only; or setting $\hat{\Pi}^{rs}$ to consist of literally \emph{all} paths (even cyclic ones) connecting origin $r$ to destination $s$.
Both of these are discussed next.
The key to both of these definitions of $\hat{\Pi}^{rs}$ is that we can determine how many of the travelers passing through a given node came from each of the available incoming links, \emph{without needing to know the specific path they are on}.
This is known as the Markov property,\index{Markov property} and is discussed at more length in the optional Section~\ref{sec:markovproperty}.

\subsubsection{Totally acyclic paths}
\label{sec:totallyacyclic}

For a particular origin $r$ and destination $s$, choose a set of allowable links $\hat{A}^{rs}$,\label{not:Ahatrs} and let $\hat{\Pi}^{rs}$ consist of all paths starting at $r$, ending at $s$, and only containing links from $\hat{A}^{rs}$.
We require two conditions on the set of allowable links:
\begin{enumerate}
\item There is at least one path from $r$ to $s$ using allowable links; this ensures that $\hat{\Pi}^{rs}$ is nonempty.
\item The set of allowable links contains no cycle; this ensures that all paths in $\hat{\Pi}^{rs}$ are acyclic.
\end{enumerate}
We say a set of paths $\hat{\Pi}^{rs}$ is \emph{totally acyclic}\index{totally acyclic path set} if it can be generated from an allowable link set satisfying these conditions.

This definition is closely related to the idea of a bush from Section~\ref{sec:bushes}.
If $\hat{A}^{rs}$ contains a path from $r$ to every destination, it is also a bush;\index{bush} and if furthermore the sets $\hat{A}^{rs}$ are the same for all destinations $s$, we can do the network loading for all the travelers leaving origin $r$ simultaneously, rather than separately for each destination.

Note that there are collections of acyclic paths which are not totally acyclic.
In the network in Figure~\ref{fig:braesstotallyacyclic}, if we choose $\hat{\Pi}^{14} = \{ [1,2,3,4], [1,3,2,4] \}$, both paths in this set are acyclic, but the set of allowable links needs to include every link in the network:
\[ \hat{A}^{14} = \myc{ (1,2), (1,3), (2,3), (2,4), (3,2), (3,4) } \]
This set contains the cycle $[2,3,2]$, so it is not possible to generate a reasonable path set containing $[1,2,3,4]$ and $[1,3,2,4]$ from an acyclic set of allowable links.

\stevefig{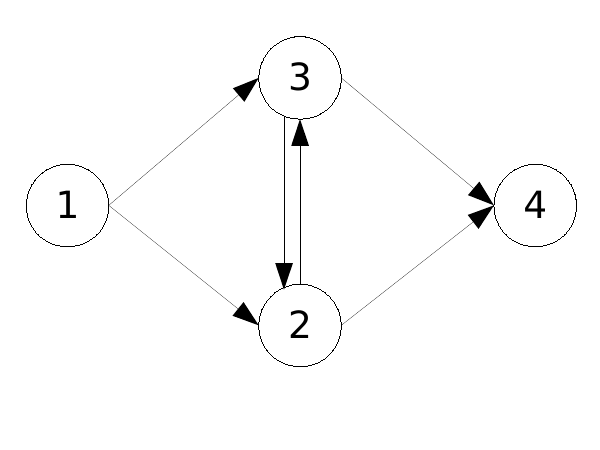}{A set of acyclic paths need not be totally acyclic.}{0.4\textwidth}

The advantage of totally acyclic path sets is that we can define a topological order on the nodes (see Section~\ref{sec:treesacyclic}).
With this topological order,\index{topological order} we can efficiently make computations using the logit formula~\eqn{logitpath} without having to enumerate all the paths.
This procedure is described in the next section.

We next describe two ways to form sets of totally acyclic paths.
For each link, define a positive value $c^0_{ij}$\label{not:c0ij} for each link which is constant and independent of flow --- examples include the free-flow travel time or distance on the link.
For each origin $r$ and node $i$, let $L_i^r$ denote the length of the shortest path from $r$ to $i$, using the quantities $c^0_{ij}$ as the link costs.
Likewise, for each destination $s$ and node $i$, let $\ell_i^s$\label{not:lis} denote the length of the shortest path from $i$ to $s$, again using the quantities $c^0_{ij}$ as the link costs.

Consider the following sets of paths:
\begin{enumerate}
\item The set of all paths, for which the head node of each link is further away from the origin than the tail node, based on the quantities $c^0_{ij}$.
That is, the sets $\hat{\Pi}^r$ containing all paths starting at $r$ and satisfying $L_j^r > L_i^r$ for each link $(i,j)$ in the path.
\item The set of all paths, for which the head node of each link is closer to the destination than the tail node, based on the quantities $c^0_{ij}$.
That is, the sets $\hat{\Pi}^s$ containing all paths ending at $s$ satisfying $\ell_j^s < \ell_i^s$ for each link $(i,j)$ in the path.
(This is like the first one, but oriented toward the destination, rather than the origin.)
\item The set of paths which satisfy both of the above criteria: $\hat{\Pi}^{rs} = \hat{\Pi}^r \cap \hat{\Pi}^s$.
\end{enumerate}
For instance, if $c^0_{ij}$ reflects the physical length on each link, then for a given origin, $\hat{\Pi}^r$ would consist of all of the paths which start at that origin and always move away from it, never ``doubling back.''  Likewise, $\hat{\Pi}^s$ would consist of all paths which always move closer to their destination node $s$, without any diversions that lead it away.
The third option has paths which both continually move away from their origin and toward their destination.
Exercise~\ref{ex:acyclic} asks you to show that all three possibilities for $\hat{\Pi}$ are totally acyclic.

Figure~\ref{fig:acyclicdefinitions} illustrates these three definitions.
The top of the figure shows a network with node A as origin and node F as destination, and the links are labeled with their $c^0_{ij}$ values.
The nodes are labeled with their $L_i$ values (above each node) and $\ell_i$ values below.
The bottom of the figure shows the links satisfying each of the three criteria ($L_j > L_i$; $\ell_j < \ell_i$; and both of these simultaneously).
The paths in these networks are the allowable paths in the original network.
Notice that in all cases, there are no cycles in these links (even though the original network had the cycle [2,5,2]).

\stevefig{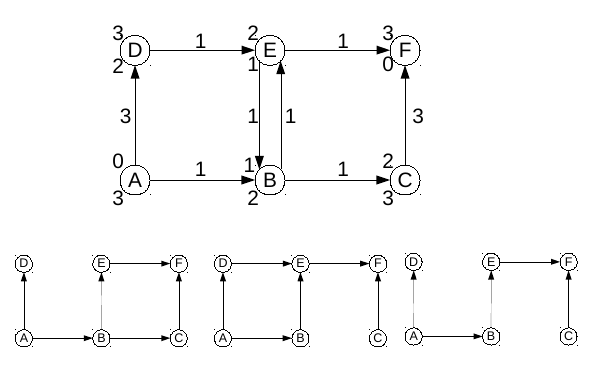}{Demonstration of three alternative definitions of totally acyclic paths.
Top panel shows $c^0_{ij}$ values on links, shortest path distances $L_i$ from A above nodes, and shortest path distances $\ell_i$ to F below nodes.}{\textwidth}

\stevefig{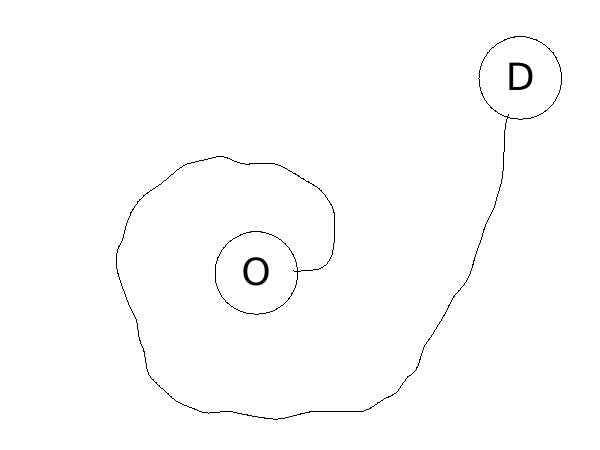}{Is this path reasonable?}{0.4\textwidth}

Of these principles, the third imposes stricter conditions on which paths are in the reasonable set.
The first and second are weaker, and includes some paths which may not seem reasonable to you.
For instance, the spiral path in Figure~\ref{fig:reasonablespiral} satisfies the first condition, since the distance from the origin is always increasing.
However, it does not satisfy the third condition, since at times the distance to the destination increases as well.

However, a major advantage of the first two principles is that we can aggregate travelers by origin or destination.
With the first principle, the destination of travelers can be ignored for routing purposes --- if a path is reasonable for a travel from an origin $r$ to a node $i$, that path segment is reasonable for travel to any node beyond $i$ as well.
This allows us to aggregate travelers by origin (as in Section~\ref{sec:bushes}) and calculate a ``one-to-all'' path set for each origin, rather than having separate path sets for each OD pair.
A similar destination-based aggregation is possible with the second principle.

\subsubsection{Full cyclic path set}

\index{static traffic assignment!stochastic user equilibrium!cyclic paths|(}
Instead of restricting the path set to create a totally acyclic collection, an alternative is to have $\hat{\Pi}^{rs}$ consist of literally all paths from origin $r$ to destination $s$, \emph{even including cycles}.
This will often mean these sets are infinite.
For example, consider the network in Figure~\ref{fig:braesstotallyacyclic}.
Under this definition, the (cyclic) path $[1,2,3,2,4]$ is part of $\hat{\Pi}^{rs}$, as is $[1,2,3,2,3,2,4]$, and so on.

Including cyclic paths, especially paths with arbitrarily many repetitions of cycles, may seem counterintuitive.
There are several reasons why this definition of $\hat{\Pi}^{rs}$ is nevertheless useful.
One reason is that requiring total acyclicity is in fact quite a strong condition.
In Figure~\ref{fig:braesstotallyacyclic}, there is no totally acyclic path set that includes both $[1,2,3,4]$ and $[1,3,2,4]$ as paths --- if we want to allow one path as reasonable, then by symmetry the other should be reasonable as well.
But any path set including both of those includes both links $(2,3)$ and $(3,2)$, which form a cycle.

So, it is desirable to have an alternative to total acyclicity that still does not require path enumeration.
As shown in the following section, it is possible to compute the link flows resulting from~\eqn{logitpath} without having to list all the paths, if all cyclic paths are included.
The intuition is that all travelers at a given node can be treated identically in terms of which link they move to next: in Figure~\ref{fig:braesstotallyacyclic}, we can split the vehicles arriving at node 2 between links $(2,3)$ and $(2,4)$ without having to distinguish whether they came via link (1,2), as if on the path [1,2,4] or [1,2,3,4], or whether they came via link (3,2), as if on the path [1,3,2,4], or even [1,2,3,2,4].

Furthermore, some modelers are philosophically uncomfortable with including restrictions like those in the previous section, without evidence that those rules really represent traveler behavior.
Determining which sets of paths travelers actually consider (and why) is complicated, and still not well-understood.\footnote{Emerging data sources, such as Bluetooth readers, are providing more complete information on observed vehicle trajectories.
This may provide more insight on this subject.}  One school of thought is that it is therefore better to impose no restrictions at all, essentially taking an ``agnostic'' position with respect to the sets $\hat{\Pi}^{rs}$, rather than imposing restrictions which may not actually represent real behavior.

As an example, assume that every link in Figure~\ref{fig:braesstotallyacyclic} has the same travel time of 1 unit, and that $\theta = \log 2$.
Then there are two paths of length 2 ([1,2,4] and [1,3,4]), two paths of length 3 ([1,2,3,4] and [1,3,2,4]), two paths of length 4 ([1,2,3,2,4] and [1,3,2,3,4]), and so on.
Therefore the denominator in the logit formula is
\labeleqn{logitdenomexample}{
\sum_{\pi \in \hat{\Pi}} \exp(-\theta c^\pi) = \frac{1}{2} + \frac{1}{2} + \frac{1}{4} + \frac{1}{4} + \frac{1}{8} + \frac{1}{8} + \cdots = 2
}
and the probability of choosing one of the length-2 paths is $(1/2) / 2 = 1/4$, the probability of choosing one of the length-3 paths is 1/8, and so on.
The flows on each link can be calculated by multiplying these path flows by the number of times that path uses a link.
For instance, to calculate the flow on link (1,2), observe that it is used by paths [1,2,4], [1,2,3,4], [1,2,3,2,4], [1,2,3,2,3,4], and so on, with respective probabilities $1/4$, $1/8$, $1/16$, $1/32$, etc.
Thus the total flow on this link is the sum of these, or $1/2$.
The flows on links (1,3), (2,4), and (3,4) are also found to be 1/2 by the same technique.
Calculating the flow on links (2,3) and (3,2) is trickier, because some paths use these multiple times.
For example, path [1,2,3,2,3,4] uses link (2,3) twice, so even though the probability of selecting this path is $1/2^5 = 1/32$, it actually contributes twice this ($1/16$) to the flow.
It is possible to show that the flow on these links is also 1/2, giving the final flows in Figure~\ref{fig:braesstotallycyclicfinalflows}.

\stevefig{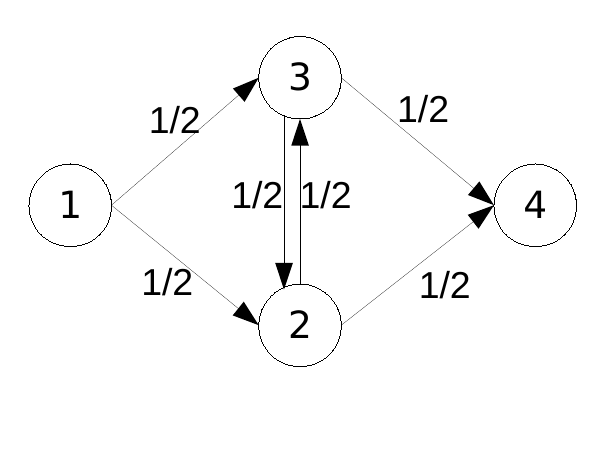}{Flows from stochastic network loading when all links have unit cost, $\theta = \log 2$, and the full cyclic path set is allowed.}{0.5\textwidth}

As this example shows, direct calculations involving this path set usually involve summing infinite series.
As will show in Section~\ref{sec:stochasticloading}, there is an alternative method that allows us to make these computations without explicitly calculating such sums.
\index{static traffic assignment!stochastic user equilibrium!cyclic paths|)}
\index{reasonable paths|)}

\subsection{The Markov property and the logit formula (*)}
\label{sec:markovproperty}

\emph{(This optional section gives mathematical reasons why totally acyclic path sets and complete path sets both allow for efficient computation of the logit formula.)}

\index{Markov property|(}
Both path set definitions above (``totally acyclic paths'' and the ``full cyclic path set'') allow for the computation of the logit formula to be disaggregated by node and by link, without having to enumerate the paths in the network.
The key to this is the \emph{Markov property}.
An informal statement of this property is that if we randomly select a traveler passing through a node, and want to know the probability that they leave that node by a particular link, there is no information provided by knowing which link they used to arrive to that node.

For example, consider the network in Figure~\ref{fig:triforce}, where the demand is $d_{13} = 40$ vehicles and $\theta = \log 2$.
All links have unit cost.
Assume first that all paths in this network are allowed.
Then the right panel of Figure~\ref{fig:triforcemarkov} shows the flow on each path, and the left panel shows the flow on each link.

\stevefig{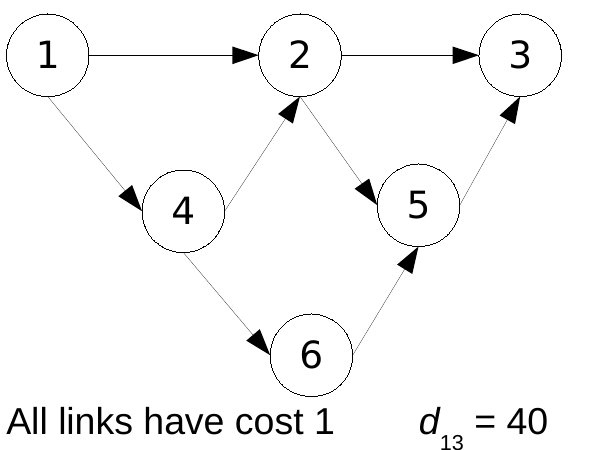}{Network for demonstrating the Markov property.}{0.4\textwidth}

In this network 18 vehicles pass through node 2.
Suppose we pick one of them at random, and want to know the probability that the next link in this vehicle's path is (2,3), as opposed to (2,5).
From examining the path flows in Figure~\ref{fig:triforcemarkov}, we can see that this probability is $(8 + 4) / (8 + 4 + 4 + 2) = 2/3$ (ignoring the flow on path [1,4,6,5,3] in the denominator, since these trips do not pass through node 2).
Now, suppose that these vehicles also reported the segment of their path that led them to node 2 --- that is, they also report whether they came via segment [1,2] or segment [1,4,2].
Does this change our answers in any way?  

\stevefig{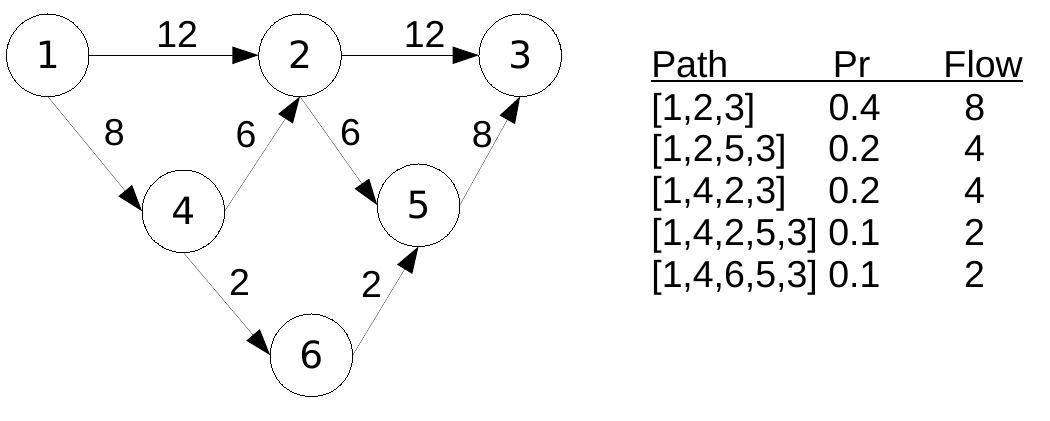}{Link and path flows when all paths are allowed.}{\textwidth}

If we know they came from segment [1,2], then they are either on path [1,2,3] or [1,2,5,3], and the probability that they continue on (2,3) is $8 / (8 + 4) = 2/3$.
If we know they came from segment [1,4,2], then they are either on path [1,4,2,3] or path [1,4,2,5,3], and the probability that they continue on (2,3) is $4 / (4 + 2)$, which is still $2/3$.
So knowing the first segment of their trip does not provide any additional information as to the remaining segment.

The situation changes if we modify the allowable path set to only include three paths, [1,2,3], [1,2,5,3], and [1,4,2,3].\footnote{A natural way this path set might arise is to include paths that are only within a small threshold of the shortest path cost, thus including paths with cost 2 or 3 but excluding paths with cost 4.}  Here Figure~\ref{fig:triforcenonmarkov} shows the corresponding path flows and link flows.

\stevefig{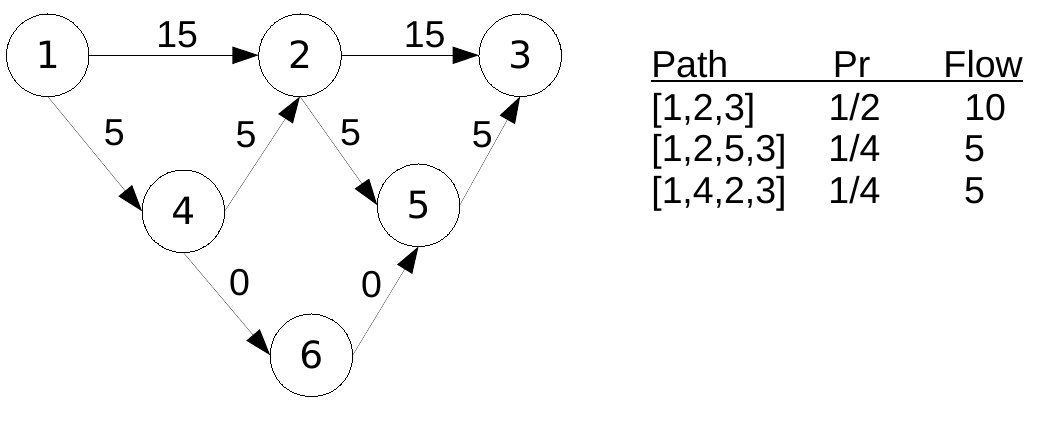}{Link and path flows when only a subset of paths is allowed.}{\textwidth}

Let us ask the same question of the travelers passing through node 2.
Without knowing anything further, the probability that they continue on link (2,3) is $(10 + 5)/20 = 3/4$.
However, if we know they came from [1,2], then the probability that they continue on (2,3) is $10/15 = 2/3$.
If we know they came from [1,4,2], then the probability that they continue on (2,3) is 1, and there is no other option!  So in this case, knowing the first segment of the path \emph{does} give us additional information about the rest of their journey.

It turns out that the Markov property will be very useful, and will allow us to efficiently evaluate the logit formula without enumerating paths.
Informally, we can do computations using just the link flows (the left panels in Figures~\ref{fig:triforcemarkov} and~\ref{fig:triforcenonmarkov}) without having to use the path flows (the right panels of these figures) --- we can get the path flows on the right from the link flows on the left.
In large networks, the link-based representation is much more compact and efficient.

A more formal statement of this property is as follows.
To keep the formulas clean, assume that there is a single origin $r$ and destination $s$; in a general network, we can apply the same logic separately to each OD pair.
In logit assignment, the path set $\hat{\Pi}$ is said to satisfy the Markov property if there exist values $P_{ij}$\label{not:Pij} for each link such that
\labeleqn{markov}{
p^\pi = \frac{\exp(-\theta c^\pi)}
               {\sum_{\pi' \in \hat{\Pi}} \exp(-\theta c^{\pi'})}
       = \prod_{(i,j) \in \pi} P_{ij}^{\delta_{ij}^\pi}       
\,
}
where $\delta_{ij}^\pi$ is the number of times path $\pi$ uses link $(i,j)$.

That is, that the probability of a traveler selecting any path can be computed by multiplying $P_{ij}$ values across its links.
The $P_{ij}$ values can be interpreted as conditional probabilities: given that a traveler is passing through $i$, they express the probability that their path leaves that node through link $(i,j)$.

\subsubsection{The segment substitution property}

\index{segment}
\index{reasonable paths!segment substitution property}
To make the connection between totally acyclic and complete path sets and the Markov property, we first show that both path set choices satisfy the \emph{segment substitution property}.
Given any path $\pi = [r, i_1, i_2, \ldots, s]$, a \emph{segment} $\sigma$\label{not:sigmaseg} is any contiguous subset of one or more of its nodes.
For example, the path $[1,2,3,2,4]$ contains $[1,2,3]$, $[2,3]$, $[2]$, and $[2,3,2,4]$.
It does not contain the segment $[1,4]$; even though both of those nodes are in the path, they do not appear consecutively.
Note that a segment can consist of a single link, such as $[2,3]$, or even a single node, such as $[2]$.
We use the notation $\oplus$\label{not:oplus} to indicate joining segments, so $[1,2,3,2,4] = [1,2,3] \oplus [3,2,4]$.
If two segments are being joined, the end node of the first must match the starting node of the second.

The set of reasonable paths $\hat{\Pi}$ satisfies the segment substitution property if, for any pair of reasonable paths which pass through the same two nodes, the paths formed by exchanging the segments between those nodes are also reasonable.
That is, if $\pi = \sigma_1 \oplus \sigma_2 \oplus \sigma_3$ is reasonable, and if there is another reasonable path $\pi' = \sigma'_1 \oplus \sigma'_2 \oplus \sigma'_3$ with $\sigma_2$ and $\sigma'_2$ starting and ending at the same node, then the paths $\sigma_1 \oplus \sigma'_2 \oplus \sigma_3$ and $\sigma'_1 \oplus \sigma_2 \oplus \sigma'_3$ are also reasonable.

With the allowable path set in Figure~\ref{fig:triforcenonmarkov}, the segment substitution property is not satisfied.
There are paths [1,2,5,3] and [1,4,2,3], both passing through nodes 1 and 2, which can be decomposed in the following way: 
\begin{align}
[1,2,5,3] &= [1] \oplus [1,2] \oplus [2,5,3] \label{eqn:segsubfail1} \\
[1,4,2,3] &= [1] \oplus [1,4,2] \oplus [2,3] \label{eqn:segsubfail2}
\end{align}
Notice that the corresponding pairs of segments on the right-hand sides all start and end at the same nodes.
We can generate two new paths by ``crossing'' the middle segments of~\eqn{segsubfail1} and~\eqn{segsubfail2}: [1,2,3] and [1,4,2,5,3].
The first of these is allowable, but the second is not.

By contrast, you can verify that segment substitution is satisfied for the allowable path set in Figure~\ref{fig:triforcemarkov}.
No matter which pairs of paths you choose, swapping the segments results in another allowable path.

A reasonable path set cannot have the Markov property unless it satisfies the
segment substitution property, as shown in the above example.
Without the segment substitution property, we could not know how the vehicles at node $i$ would split without knowing the specific paths they were on.
Segment substitution ensures that all travelers passing through node $i$ are considering the same set of outgoing links (and, indeed, the same set of path segments continuing on to the destination).

It is fairly easy to show that both path set definitions considered above --- sets of totally acyclic paths, and the set of all paths (even cyclic ones) --- satisfy the segment substitution property; see Exercise~\ref{ex:segmentsubstitution}.

\subsubsection{Decomposing the logit formula}

Given a reasonable path set $\hat{\Pi}$ satisfying the segment substitution property, and any two nodes $a$ and $b$, let $\Sigma_{ab}$\label{not:Sigmaab} denote the set of segments which start and end at these nodes, and are part of a reasonable path.
Define the quantity\label{not:Vab}
\labeleqn{logitv}{
V_{ab} = \sum_{\sigma \in \Sigma_{ab}} \exp(-\theta c^\sigma)
}
where $c^\sigma = \sum_{(i,j) \in \sigma} c_{ij}$\label{not:csigma} is the travel time of a segment.
In particular, $V_{rs}$ is a sum over all of the reasonable paths from origin to destination, and corresponds to the denominator of the logit formula.
Thus,
the probability that a traveler chooses a particular path $\pi$ is simply $\exp(-\theta c^\pi)/V_{rs}$, and the flow on this path is $d^{rs} \exp(-\theta c^\pi) / V_{rs}$.
Also note the special case $V_{ii}$ where the start and end nodes are the same.
In this case $\Sigma_{ii}$ consists only of the single-node segment $[i]$ with zero cost, so $V_{ii} = 1$.

Furthermore, note that the numerator of the logit formula can be factored by segment, so that if $\pi = \sigma_1 \oplus \sigma_2 \oplus \sigma_3$, we have
\labeleqn{logitfactor}{
\exp(-\theta c^\pi) = \exp(-\theta (c^{\sigma_1} + c^{\sigma_2} + c^{\sigma_3}))
= \exp(-\theta c^{\sigma_1}) \exp(-\theta c^{\sigma_2}) \exp(-\theta c^{\sigma_3})
\,.}

To calculate the flow on a link $x_{ij}$, we need to add the flow from all of the reasonable paths which use this link.
Every reasonable path using link $(i,j)$ takes the form $\sigma_1 \oplus [i,j] \oplus \sigma_2$, where $\sigma_1$ goes from the origin $r$ to node $i$ and $\sigma_2$ goes from node $j$ to the destination $s$.\footnote{If the link starts at the origin or ends at the destination, we may have $i = r$ or $j = s$, in which case $\sigma_1$ or $\sigma_2$ will consist of a single node, $[r]$ or $[s]$.}
Therefore 
\begin{align}
x_{ij} &= \sum_{\pi \in \hat{\Pi}^{rs}} \delta_{ij}^\pi h^\pi \\
       &= \frac{d^{rs}}{V_{rs}} \sum_{\pi \in \hat{\Pi}^{rs}} \delta_{ij}^\pi \exp(-\theta c^\pi) \\
       &= \frac{d^{rs}}{V_{rs}} \sum_{\sigma_1 \in \Sigma_{ri}} \sum_{\sigma_2 \in \Sigma_{js}} \exp(-\theta(c^{\sigma_1} + t_{ij} + c^{\sigma_2})) \\
       &= \frac{d^{rs}}{V_{rs}} \sum_{\sigma_1 \in \Sigma_{ri}} \sum_{\sigma_2 \in \Sigma_{js}} \exp(-\theta c^{\sigma_1}) \exp(-\theta t_{ij}) \exp(-\theta c^{\sigma_2})\\
       &= \frac{d^{rs}}{V_{rs}} \myp{\sum_{\sigma_1 \in \Sigma_{ri}} \exp(-\theta c^{\sigma_1})} \exp(-\theta t_{ij}) \myp{\sum_{\sigma_2 \in \Sigma_{js}} \exp(-\theta c^{\sigma_2})} \\
       &= d^{rs} \frac{V_{ri} \exp(-\theta t_{ij}) V_{js}}{V_{rs}} \label{eqn:finallogitdecomposition}
\end{align}
The third equality groups the sum over paths according to the starting and ending segment.
This equation for $x_{ij}$ will be used extensively in the rest of this section.

Similarly, to find the number of vehicles passing through a particular node $i$ (call this $x_i$), observe that every path through $i$ can be divided into a segment from $r$ to $i$, and a segment from $i$ to $s$.
Grouping the paths according to these segments and repeating the algebraic manipulations above gives the formula
\labeleqn{finallogitnodedecomposition}{
    x_i = d^{rs} \frac{V_{ri} V_{is}}{V_{rs}}
}

With equations~\eqn{finallogitdecomposition} and~\eqn{finallogitnodedecomposition} in hand, it is easy to show the Markov property holds in any reasonable path set with the segment substitution property.
Define
\labeleqn{logitmarkov}{
P_{ij} = \frac{x_{ij}}{x_i} = \frac{\exp(-\theta t_{ij}) V_{js}}{V_{is}}
    \,
}
and then multiply these values together for the links in a path, say, $\pi =
[r, i_1, i_2, \cdots, i_k, s]$:\index{Markov property!definition}
\begin{align}
\prod_{(i,j) \in \pi} P_{ij}
 &= \prod \frac{\exp(-\theta t_{ij}) V_{js}}{V_{is}} \\
 &= \exp \myp{
        -\theta \sum_{(i,j) \in \pi} t_{ij}
     }
     \frac{V_{i_1 s} V_{i_2 s} \cdots V_{i_k s} V_{ss}}
          {V_{rs} V_{i_1 s} V_{i_2 s} \cdots V_{i_k s}} \\
 &= \exp(-\theta c^\pi) \frac{V_{ss}}{V_{rs}}
\end{align}
But $V_{ss} = 1$ and $V_{rs}$ is simply the denominator in the logit formula,
so this product is exactly $p^\pi$.
Therefore the logit path flow assignment in
the set of reasonable paths satisfies the Markov property.
\index{Markov property|)}

\subsection{Stochastic network loading}
\label{sec:stochasticloading}

\index{static traffic assignment!stochastic network loading|(}
The stochastic network loading problem is to determine the flows on each link $x_{ij}$, given their travel times $t_{ij}$, according to a particular discrete choice model.
This section describes how to do so for the logit model, using formula~\eqn{logitpath}.
For stochastic network loading, we assume that these travel times are fixed and constant, and therefore unaffected by the path and link flows we calculate.

In the larger stochastic user equilibrium problem (where travel times can depend on flows), stochastic network loading plays the role of a subproblem in an iterative algorithm.
This is analogous to how shortest paths and all-or-nothing loadings are often used as a subproblem in the classical traffic assignment problem: although travel times do depend on link flows, to find an equilibrium we can solve a number of shortest path problems, temporarily fixing the link costs at particular values.

Throughout this section, we are assuming a single origin $r$ and destination $s$ to simplify the notation.
If there are many origins and destinations, these procedures should be repeated for each, and the link flows added to obtain the total link flows.
(There is no harm in doing so, since we are assuming link travel times are constant, and therefore the different OD pairs do not interact with each other.)  Depending on the choice of path set, it may be possible to aggregate all travelers from the same origin or destination, and load them at once.
This is more efficient than doing a separate loading for each OD pair, but requires a more limited definition of the reasonable path sets.

In principle, the stochastic network loading procedure is straightforward.
Given the link travel times $t_{ij}$, we can calculate the path travel times $c^\pi$ as in Section~\ref{sec:staticconcepts}.
The path flows $h^\pi$ can then be calculated from the logit formula~\eqn{logitpath}, from which the link flows can be calculated by addition, again as in Section~\ref{sec:staticconcepts}.

While conceptually straightforward, this procedure faces the practical difficulty that the logit formula requires summations over the set of all reasonable paths $\hat{\Pi}^{rs}$, which can grow exponentially with network size.
In a realistic-sized network, this renders the above procedure unusable, or at least computationally taxing.\index{path!enumeration challenges}

This section describes how the link flows can be calculated \emph{without} explicitly enumerating paths or using the logit formula.
This is possible for both of the $\hat{\Pi}$ definitions discussed in the previous section: a totally acyclic set of paths, or including the full path set, even including all cycles.

These two procedures have several features in common.
As was described in Section~\ref{sec:markovproperty}, both of these path set definitions have the Markov property.
This section derived an important formula, repeated here:
\labeleqn{repeatmarkov}{
x_{ij} = d^{rs} \frac{V_{ri} \exp(-\theta t_{ij}) V_{js}}{V_{rs}}
\,,
}
where $V_{ab}$ is the sum of $\exp(-\theta c^\sigma)$ for all reasonable path segments $\sigma$ starting at node $a$ and ending at node $b$.
This formula
is important because it allows us to calculate the flow on each link without
having to enumerate all of the paths that use that link.

It is also possible to show (see Exercise~\ref{ex:rescale}) that we can replace
$t_{ij}$ with $L_i + t_{ij} - L_j$, where $\mb{L}$ is a vector of node-specific constants in this formula; this reduces numerical errors in computations.
It is common to use shortest path distances at free-flow (hence the use of the notation $L$).
Doing this for the link and segment costs helps avoid numerical issues involved with calculating exponentials of large values.
With this re-scaling, we define the \emph{link likelihood}\label{not:Lijlik} as
\labeleqn{linklikelihood}{L_{ij} = \exp(\theta(L_j - L_i - t_{ij}))}
and thus
\labeleqn{repeatmarkov2}{x_{ij} = d^{rs} \frac{V_{ri} L_{ij} V_{js}}{V_{rs}}\,.} 
\index{static traffic assignment!stochastic network loading!link likelihood}

It remains to describe how the $V_{ab}$ values can be efficiently calculated.
This section shows how this can be done for both totally acyclic path sets, and the complete set of all paths (even cyclic ones).

In the case of a totally acyclic path set, the relevant $V$ values can be calculated in a single pass over the network in topological order, and then the link flows calculated in a second pass over the network in reverse topological order.
In the case of the full cyclic path set, the cycles create dependencies in the link weights and in the link flow formulas which prevent them from being directly evaluated.
But we can still calculate them explicitly using matrix techniques.

\subsubsection{Totally acyclic paths}

\index{static traffic assignment!stochastic network loading!totally acyclic paths|(}
\index{totally acyclic paths|(}
\index{Dial's method|(}
The defining feature of a totally acyclic path set is that the collection of links used by reasonable paths has no cycles.
This allows us to define a topological order on the nodes, so that each allowable link connects a lower-numbered node to a higher-numbered one.
In acyclic networks, it is often easy to perform calculations recursively, in increasing or decreasing topological order.
Stochastic network loading is one of these cases, using Dial's method.

Define $W_i$ and $W_{ij}$\label{not:Wi} to be the \emph{weight} of a node and link, respectively.
The node weight $W_i$ is a shorthand for $V_{ri}$, and $W_{ij}$ is a shorthand for $V_{ri} L_{ij}$.
These can be calculated recursively, using the following procedure:
\begin{enumerate}
\item Calculate the link likelihoods $L_{ij}$ using equation~\eqn{linklikelihood} for all allowable links; set $L_{ij} \leftarrow 0$ for any link not in the allowable set.
\item Set the current node $i$ to be the origin $r$, and initialize its weight: $W_r \leftarrow 1$.
\item For all links $(i,j)$ leaving node $i$, set $W_{ij} \leftarrow W_i L_{ij}$.
\item If $i$ is the destination, stop.
Otherwise, set $i$ to be the next node in topological order.
\item Calculate $W_i \leftarrow \sum_{(h,i) \in \Gamma^{-1}(i)} W_{hi}$ by summing the weights of the incoming links.
\item Return to step 3.
\end{enumerate}

With the node and link weights in hand, we proceed to calculate the flows on each link $x_{ij}$ and the flow through each node $x_i$.
Using the link weights, we can rewrite equation~\eqn{finallogitnodedecomposition} for the flow through each node as
\labeleqn{noderewrite}{
x_i = d^{rs} \frac{V_{ri} V_{is}}{V_{rs}} = d^{rs} \frac{W_i V_{is}}{W_s}
\,.
}
Combining with equation~\eqn{repeatmarkov2} for link flow, we have
\labeleqn{linkrewrite}{
x_{hi} = d^{rs} \frac{V_{rh} L_{hi} V_{is}}{V_{rs}} 
       = x_i \frac{W_{hi}}{W_i}
\,,
}
so if we know the flows to some node $i$, we can calculate the flows to its
incoming links $(h,i)$.
So, the link flows can be calculated in reverse topological
order:
\begin{enumerate}
\item Initialize all flows to zero: $x_i = 0$ for all nodes, and $x_{ij} = 0$ for all links.
\item Set the current node $i$ to be the destination $s$, and initialize its flow: $x_i = d^{rs}$, since all vehicles must reach the destination.
\item For all links $(h,i)$ entering node $i$, set $x_{hi} = x_i W_{hi} / W_i$.
\item If $i$ is the origin, stop.
Otherwise, set $i$ to be the previous node in topological order.
\item Compute $x_i = \sum_{(i,j) \in \Gamma(i)} x_{ij}$ as the sum of the flows on outgoing links.
\item Return to step 3.
\end{enumerate}

As an example, Dial's method is demonstrated on the network shown in Figure~\ref{fig:braessdial}, where the demand is 2368 vehicles from node 1 to node 5 and $\theta = 1$.
For convenience, the results of the calculations in this example are shown in Tables~\ref{tbl:nodedial} and~\ref{tbl:linkdial}.
As a preliminary step, we calculate the shortest path from node 1 to all other nodes at free-flow conditions (assuming $\mb{x= 0}$) using standard techniques.
Details are omitted since the technique is familiar, and the resulting shortest path labels $L_i$ are shown in Table~\ref{tbl:nodedial}.

Next, we must identify the allowable links.
This example will adopt the first principle from the previous section, where the allowable links $(i,j)$ are those for which $L_i < L_j$.
Every link is thus allowable except for (2,3), because this link connects node 2 (with shortest path label $L_2 = 2$) to node 3 ($L_3 = 1$).
Therefore, link (2,3) will be excluded from the remaining steps of Dial's method, since no vehicles will use this link.

The next step is to calculate the link likelihoods $L_{ij}$ for the allowable links $(i,j) \in A_B$.
So, $L_{12} = \exp(2 - 0 - 2) = 1$, $L_{35} = \exp(3 - 1 - 3) = e^{-1} \approx 0.368$, and so forth.
Weights are now calculated in forward topological order.
The topological ordering for the reasonable bush is 1, 3, 2, 4, and 5 in that order.
Notice that the original network has a cycle $[2,3,2]$ and so no topological order can exist on the full network.
But when restricted to the allowable set, the cycle disappears and node 3 must come before node 2 topologically.
For the origin, $W_1 = 1$ by definition.
The weights on links leaving node 1 can now be calculated: $W_{12} = W_1 L_{12} = 1$, and $W_{13} = 1 \times 1 = 1$.
Proceeding to node 3 (the next in topological order), $W_3$ is calculated as the sum of the weights on incoming links: $W_3 = W_{13} = 1$ and thus $W_{32} = 1 \times 1 = 1$ and $W_{35} = 1 \times e^{-1} \approx 0.368$.
Node 2 is next in topological order, and its weight is the sum of the weights on its incoming links: $W_2 = W_{12} + W_{32} = 2$, so $W_{24} = 2$.
Node 4 is next, and we have $W_4 = 2$ and $W_{45} = 2$.
Finally, the weight of node 5 is $W_5 = W_{35} + W_{45} = 2.368$.

Link and node flows are now calculated in \emph{reverse} topological order, starting with node 5 and then proceeding to nodes 4, 2, 3, and 1.
For node 5, the node flow $X_5$ is simply the demand destined to this node (since there are no outgoing links), so $X_5 = 2368$.
This flow is now distributed among the two incoming links $(3,5)$ and $(4,5)$ in proportion to their weights.
So, $x_{35} = X_5 W_{35} / (W_{35} + W_{45}) = 2368 \times 0.368 / (2 + 0.368) = 368$ and $x_{45} = 2000$.
Proceeding upstream, the node flow at 4 is simply the flow on link $(4,5)$, the only outgoing link (since no vehicles have node 4 as their destination), and $X_4 = 2000$.
There is only one incoming link $(2,4)$, so $x_{24} = 2000$.
(This follows trivially from the formula in step 3c since $W_{4} = W_{24}$.)
Thus $X_2 = 2000$.
Since there are two incoming reasonable links to node 2 with equal weight, they receive equal flow (again, following from the formula in 3c), setting $x_{32} = x_{12} = 1000$.
Then $X_3 = x_{32} + x_{35} = 1368$, and $x_{13} = X_3 = 1368$.
Finally, the node flow at the origin 1 is $X_1 = 1368 + 1000 = 2368$, as it should be.

Dial's method is now complete, having calculated the flows on each link.
This method is completely consistent with~\eqn{logitpath}, if we were to restrict attention to the reasonable paths in the network.
If we were to use this formula directly, we would first enumerate the three reasonable paths $[1,3,5]$, $[1,3,2,4,5]$, and $[1,2,4,5]$ and calculate their costs: $C^{[1,3,5]} = 4$, $C^{[1,3,2,4,5]} = C^{[1,2,4,5]} = 3$.
Equation~\eqn{logitpath} then gives
\[
 p^{[1,3,5]} = \frac{e^{-4}}{e^{-4} + e^{-3} + e^{-3}} = 0.155 \quad
   p^{[1,3,2,4,5]} = \frac{e^{-3}}{e^{-4} + e^{-3} + e^{-3}} = 0.422 \]
\[   p^{[1,2,4,5]} = \frac{e^{-3}}{e^{-4} + e^{-3} + e^{-3}} = 0.422 \]
as the path choice proportions.
Multiplying each of these by the total demand (2368) gives path flows
\[ h^{[1,3,5]} = 368 \qquad
   h^{[1,3,2,4,5]} = 1000 \qquad
   h^{[1,2,4,5]} = 1000
\,.   
\]
As you can verify, these path flows correspond to the same link flows shown in Table~\ref{tbl:linkdial}.
\index{static traffic assignment!stochastic network loading!totally acyclic paths|)}
\index{totally acyclic paths|)}
\index{Dial's method|)}

\stevefig{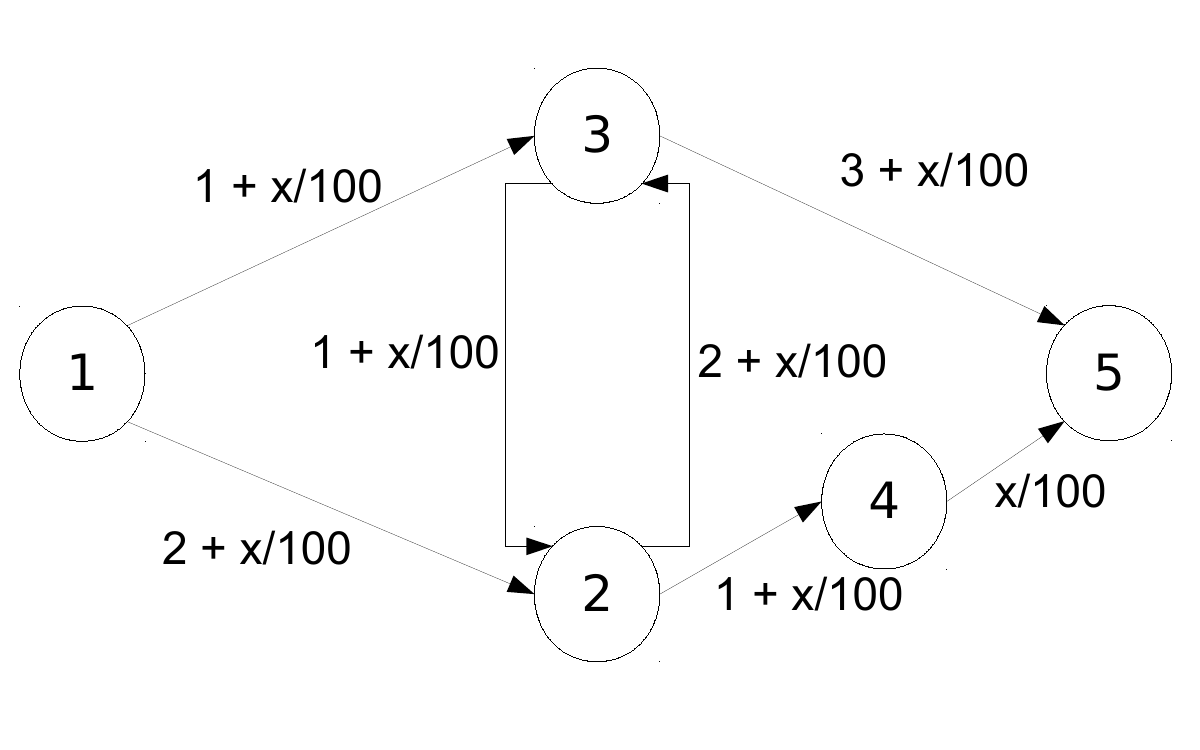}{Network for demonstration of Dial's method.}{0.6\textwidth}

\begin{table}
\begin{center}
\caption{Node calculations in Dial's method example.
\label{tbl:nodedial}}
\begin{tabular}{c|ccc}
Node $i$ & $L_i$ & $W_i$ & $X_i$ \\
\hline
1        & 0      &  1              & 2368 \\
2        & 2      &  2              & 1368 \\
3        & 1      &  1              & 2000 \\
4        & 3      &  2              & 2000 \\
5        & 3      & $2 + e^{-1}$    & 2368 \\
\end{tabular}
\end{center}
\end{table}

\begin{table}
\begin{center}
\caption{Link calculations in Dial's method example.
\label{tbl:linkdial}}
\begin{tabular}{c|ccc}
Link $(i,j)$ & $L_{ij}$ & $W_{ij}$ & $x_{ij}$ \\
\hline
(1,2) & 1        & 1       & 1000 \\
(1,3) & 1        & 1       & 1368 \\
(2,3) & 0     & 0      & 0 \\
(2,4) & 1        & 2        & 2000 \\
(3,2) & 1        & 1        & 1000 \\
(3,5) & $e^{-1}$ & $e^{-1}$ & 368 \\
(4,5) & 1        & 2        & 2000 \\
\end{tabular}
\end{center}
\end{table}

\subsubsection{Full cyclic path set}

\index{static traffic assignment!stochastic network loading!cyclic paths|(}
We can also efficiently calculate the link flows from logit network loading, if the set of reasonable paths contains all paths, even cyclic ones.
Again we will use formula~\eqn{repeatmarkov2}.
If we have an efficient way to compute $V_{ab}$ for all pairs of nodes $a$ and $b$, we can substitute them into this formula to directly obtain the link flows.

Since the set of reasonable paths contains cycles, we cannot calculate these $V$ values inductively on topological order, as was done above.
Rather, a different approach is needed.

Let $\mb{V}$ be the $n \times n$ matrix whose components are $V_{ab}$.
Recall that $V_{ab}$ is defined as the sum of $\exp(-\theta c^\sigma)$ for all segments $\sigma$ starting at node $a$ and ending at node $b$.
We will calculate this sum by dividing the sum into segments of the same length.

We begin by calculating the part of $V_{ab}$ which corresponds to the segments of length one (that is, the segments consisting of a single link.)  If there are no parallel links in the network, then this is simply $L_{ab}$: there is at most one such segment, which must be $\sigma = [a,b]$, and if it exists $\exp(-\theta c^\sigma) = \exp(-\theta t_{ab}) = L_{ab}$.
If it does not exist, then $L_{ab} = 0$, which is again the part of $V_{ab}$ corresponding to segments of length one (which is empty if no such segment exists.)  We can proceed similarly if there are parallel links; see Exercise~\ref{ex:parallellinks}.

Now, let $\mb{L}$ be the $n \times n$ matrix whose components are $L_{ab}$, and form the matrix product $\mb{L}^2$.
Its components are
\labeleqn{l2components}{(L^2)_{ab} = \sum_{c \in N} L_{ac} L_{cb}}
by the definition of matrix multiplication.\footnote{The parentheses are intentional: $(L^2)_{ab}$ is the component of matrix $\mb{L^2}$ in row $a$ and column $b$.
This is \emph{not} the same as $(L_{ab})^2$, the square of the component of matrix $\mb{L}$ in row $a$ and column $b$.}

For a given node $c$, the product $L_{ac} L_{cb}$ is zero unless there is a link both from $a$ to $c$, and from $c$ to $b$.
So, we can restrict the sum in~\eqn{l2components} to be over nodes $c$ for which there are links $(a,c)$ and $(c,b)$ --- which is precisely the nodes $c$ for which there is a segment $[a,c,b]$ of length two.
Furthermore, for such segments,
\labeleqn{l2decompose}{
\exp(-\theta c^{[a,c,b]}) = \exp(-\theta [t_{ac} + t_{cb}]) = \exp(-\theta t_{ac}) \exp(-\theta t_{cb}) = L_{ac} L_{cb} \,.}
So, \emph{the components of $\mb{L^2}$ are exactly the portion of the sums defining $V_{ab}$ for segments of length two!}

We demonstrate this using the example from Figure~\ref{fig:braesstotallyacyclic}, recalling that $d^{14} = 1$, $\theta = \log 2$, and that all links have unit cost.
In this network, the matrices $\mb{L}$ and $\mb{L^2}$ take the form
\labeleqn{ll2example}{
\mb{L} = \vect{ 0 & 1/2 & 1/2 & 0 \\
                0 & 0   & 1/2 & 1/2 \\
                0 & 1/2 & 0   & 1/2 \\
                0 & 0   & 0   & 0 }
\qquad
\mb{L^2} = \vect{ 0 & 1/4 & 1/4 & 1/2 \\
                0 & 1/4   & 0 & 1/4 \\
                0 & 0  & 1/4   & 1/4 \\
                0 & 0   & 0   & 0 }
\,.}
Looking at the first row of $\mb{L^2}$, we see that $(L^2)_{12} = 1/4$, because there is one segment which starts at node 1, ends at node 2, and contains two links ([1,3,2]) and $\exp(-\theta c^{[1,3,2]}) = 1/4$.
Similarly $(L^2)_{13} = 1/4$, and $(L^2)_{14} = 1/2$ because there are two segments starting at node 1, ending at node 2, and containing two links ([1,2,4] and [1,3,4]).
The sum $\exp(-\theta c^{[1,2,4]}) + \exp(-\theta c^{[1,3,4]})$ is indeed 1/2.
In this matrix, $(L^2)_{23} = 0$, because there are no segments of length two starting at node 2 and ending at node 3.

Proceeding a step further, we have
\labeleqn{l3components}{(L^3)_{ab} = \sum_{c \in N} (L^2)_{ac} L_{cb}\,.}
By the same logic, we see that this sum expresses the component of $V_{ab}$ corresponding to segments of length three.
Every segment of length three connecting $a$ to $b$ consists of a segment of length two connecting $a$ to some node $c$, followed by a link $(c,b)$.
Group the sum of all segments of length three by this final link, and note that $(L^2)_{ac}$ already contains the relevant portion of the sum for the first segment.

Thus, by induction,\index{mathematical induction} the components of matrix $\mb{L}^n$ contains the portion of the sum defining $\mb{V}$ corresponding to segments of length $n$.
Therefore
\labeleqn{vseries}{\mb{V} = \mb{L^0} + \mb{L^1} + \mb{L^2} + \mb{L^3} + \cdots\,,}
where the infinite sum is needed since we allow paths with an arbitrary number of cycles.
Assuming that this sum exists, we can calculate it as follows:\label{not:Identity}
\begin{align}
\mb{V} &= \mb{I} + \mb{L} + \mb{L^2} + \mb{L^3} + \cdots \\
       &= \mb{I} + \mb{L} (\mb{I} + \mb{L} + \mb{L^2} + \cdots) \\
       &= \mb{I} + \mb{LV} \,.
\end{align}
Therefore $\mb{V} - \mb{LV} = \mb{I}$, or
\labeleqn{finalv}{\mb{V} = (\mb{I} - \mb{L})^{-1}\,.}
After calculating the matrix $\mb{V}$ with this formula, we can directly read off its components $V_{ab}$ and use them to calculate the link flows using~\eqn{repeatmarkov2}.

To complete the example, we have
\labeleqn{vexample}{
\mb{V} = \vect{1 & 1 & 1 & 1 \\
               0 & 4/3 & 2/3 & 1 \\
               0 & 2/3 & 4/3 & 1 \\
               0 & 0 & 0 & 1}
}
and, for example, the flow on link (2,3) is given by
\labeleqn{xample}{
x_{23} = d^{14} \frac{V_{12} L_{23} V_{34}}{V_{14}} = 1 \frac{1 \cdot \frac{1}{2} \cdot 1}{1} = \frac{1}{2}\,.
}
Repeating this process will give the flow on every link --- and unlike in the previous section, does not involve summing an infinite series term-by-term.
The formula~\eqn{repeatmarkov2} handles all of the paths.
\index{static traffic assignment!stochastic network loading!cyclic paths|)}
\index{static traffic assignment!stochastic network loading|)}

\subsection{Stochastic user equilibrium}
\label{sec:sueformulation}

Stochastic network loading methods, described in the previous section, are the analogue of the all-or-nothing assignments we used to find $\mb{x^*}$ in the method of successive averages, or Frank-Wolfe, in Chapter~\ref{chp:solutionalgorithms}.
Recall that in those methods, we identified $\mb{x^*}$ by finding the link flows that would be observed if all the travel times were held constant at the current values.
In TAP, all drivers aim to take the shortest path, so $\mb{x^*}$ was obtained by loading all flow onto the shortest paths at the current travel times.
In logit SUE, drivers do not always take the shortest path, but instead choose paths by~\eqn{logitpath}.
The methods in the previous section thus calculate an ``$\mb{x^*}$'' in the sense that it reflects the link flows which would arise if travel times were held constant.

The full development of the SUE model requires relaxing the assumption of constant travel times, in the same way as finding a traditional user equilibrium requires more than a single shortest path computation.
SUE is very easy to define as a fixed point problem.
Let $\mb{C(h)}$ denote the vector of path travel times as a function of the vector of path flows (this is the same as before).
However, the logit formula directly gives us a complementary function $\mb{H(c)}$ giving path flows as a function of path travel times, with components
\labeleqn{logitpathexplicit}{h^\pi = d^{rs} \frac{\exp(-\theta c^\pi)}{\sum_{\pi' \in \Pi^{rs}} \exp(\theta c_{\pi'})}}
where $(r,s)$ is the OD pair corresponding to path $\pi$.
Therefore, the SUE problem can be expressed as follows: find a feasible path flow vector $\mb{h^*}$ such that $\mb{h^* = H(C(h^*))}$.
This is a standard fixed-point problem.
Clearly $\mb{H}$ and $\mb{C}$ are continuous functions if the link performance functions are continuous, and the feasible path set is compact and convex, so Brouwer's theorem immediately gives existence of a solution to the SUE problem.
\index{static traffic assignment!stochastic user equilibrium!fixed point formulation}
\index{Brouwer's theorem!applications}
\index{continuous function!applications}
\index{compact set!applications}
\index{convex set!applications}
\index{fixed point problem!and stochastic user equilibrium}

Notice that this was much easier than showing existence of an equilibrium solution to the original traffic assignment problem!  For that problem, there was no equivalent of~\eqn{logitpathexplicit}.
Travelers were all using shortest paths, but if there were two or more shortest paths there was no rule for how those ties should be broken.
As a result, we had to reformulate the problem as a variational inequality and introduce an auxiliary function based on movement of a point under a force.
For the SUE problem, there is no need for such machinations, and we can write down the fixed point problem immediately.

However, fixed-point theorems do not offer much help in terms of actually finding the SUE solution.
It turns out that convex programming and variational inequality formulations exist as well, and we will get to them shortly.
But first, as a practical note, we mention that our definition of the reasonable path set $\hat{\Pi}$ should be specified \emph{without reference to the final travel times.}  The reason is that until we have found the equilibrium solution, we do not know what the link and path travel times will be.
If the sets of reasonable paths vary from iteration to iteration, as the flows and travel times change, there may be problems with convergence or solution consistency.
This is why the methods described above for generating totally acyclic path sets relied on constants $c^0_{ij}$ which were independent of flow: constants such as physical length or free-flow times.
If using the full cyclic path set, there is no concern; at all iterations every path is reasonable.

\subsubsection{Path-flow formulation and the method of successive averages}

\index{static traffic assignment!stochastic user equilibrium!convex optimization formulation|(}
Consider the following optimization problem, developed by Caroline Fisk:\index{Fisk, Caroline}
\begin{align}
\min_{\mathbf{x},\mathbf{h}} \quad & \theta \sum_{(i,j) \in A}  \int_0^{x_{ij}} t_{ij}(x) dx + \sum_{\pi \in \hat{\Pi}} h^\pi \log h^\pi  & \label{eqn:fisk} \\
\mathrm{s.t.}                \quad & x_{ij} = \sum_{\pi \in \hat{\Pi}} h^\pi \delta_{ij}^{\pi}	                             & \forall (i,j) \in A \label{eqn:suemap} \\
                                   & \sum_{\pi \in \hat{\Pi}^{rs}} h^\pi = d^{rs}                                       & \forall (r,s) \in Z^2 \label{eqn:suedemand} \\
                                   & h^\pi \geq 0                                                                   & \forall \pi \in \hat{\Pi} \label{eqn:suenonneg} 
\end{align}
Comparing with Beckmann's formulation from Chapter~\ref{chp:trafficassignmentproblem}, we see that it is identical except that the sums on path variables are now over the set of reasonable paths, rather than all paths, and that the objective is first scaled by $\theta$, and then an additional term is added.
We show below that these modifications ensure that the optimal path flows satisfy the logit formula~\eqn{logitpath}.
This objective function is strictly convex in the \emph{path flows} (see Exercise~\ref{ex:fiskunique}), so the SUE path flow solution is unique.
This is a different situation than the classical traffic assignment problem, which has a unique equilibrium solution in link flows, but not in path flows.

First, we show that the non-negativity constraint~\eqn{suenonneg} can be ignored without any problem, and in fact that at optimality the flow on each path is \emph{strictly} positive.
The $\log h^\pi$ term in the objective is not defined at all if a path flow is negative; if $h^\pi = 0$ we can define $h^\pi \log h^\pi = 0$, since this is the limiting value by l'H\^{o}pital's rule.
But the derivative of $h \log h$ becomes infinitely steep as $h$ approaches zero: $\frac{d}{dh} (h \log h) = \log h + 1 \rightarrow -\infty$ as $h \rightarrow 0$.
Therefore the objective function cannot be minimized at $h = 0$; no matter what the change in the other terms in the objective would be, it is better to have a very slightly positive $h$ value than a zero one.
This is a very useful property, since the optimality conditions are much simpler if there are no non-negativity constraints.

Next, as with Beckmann's formulation, we Lagrangianize the demand constraint~\eqn{suedemand}, and substitute~\eqn{suemap} in place of $\mb{x}$, to obtain a Lagrangian in terms of path flow variables only:\index{Lagrange multipliers!applications}
\begin{multline}
\label{eqn:fisklagrangian}
\mc{L}(\mb{h}, \bm{\kappa}) = \theta \sum_{(i,j) \in A} \int_0^{\sum_{\pi \in \hat{\Pi}} h^\pi \delta_{ij}^{\pi}} t_{ij}(x) dx + \sum_{\pi \in \hat{\Pi}} h^\pi \log h^\pi + \\ \sum_{(r,s) \in Z^2} \kappa^{rs} \myp{d^{rs} - \sum_{\pi \in \hat{\Pi}^{rs}} h^\pi}
\end{multline}
Since there are no non-negativity conditions, it is enough to find the stationary points of the Lagrangian, that is, the points where $\pdr{\mc{L}}{h} = 0$.
Therefore, we require
\labeleqn{fisklagrangeder}{
\pdr{\mc{L}}{h^\pi} = \theta c^\pi + 1 + \log h^\pi - \kappa^{rs} = 0 \,,
}
or, solving for $h^\pi$,
\labeleqn{fiskpathflow}{
h^\pi = \exp(\kappa^{rs} - 1) \exp(-\theta c^\pi)\,.
}
To find the value of the Lagrange multiplier $\kappa^{rs}$, we substitute~\eqn{fiskpathflow} into the demand constraint~\eqn{suedemand}, and find that $\kappa^{rs}$ must be chosen such that
\labeleqn{logitpathflow}{
\exp(\kappa^{rs} - 1) = \frac{d^{rs}}{\sum_{\pi' \in \hat{\Pi}^{rs}} \exp(-\theta c^{\pi'})}\,.
}
Substituting into~\eqn{fiskpathflow} gives the logit formula, and therefore the path flows solving Fisk's convex optimization problem are those solving logit stochastic user equilibrium.
\index{static traffic assignment!stochastic user equilibrium!convex optimization formulation|)}

\index{static traffic assignment!stochastic user equilibrium!method of successive averages|(}
This convex optimization problem can be solved using the method of successive averages.
In this way, SUE can be solved quite simply, using a familiar method from the basic TAP (Section~\ref{sec:convexmsa}).
The major change is that $\mathbf{x^*}$ is calculated using a method from the previous section, rather than by finding an all-or-nothing assignment loading all flow onto shortest paths.
Two other minor changes are discussed below.
The algorithm is as follows:
\begin{enumerate}
\item Choose an initial feasible link assignment $\mathbf{x}$.
\item Update link travel times based on $\mathbf{x}$.
\item Calculate target flows $\mathbf{x^*}$:
\begin{enumerate}[(a)]
 \item For each OD pair $(r,s)$, use a method from the previous section to calculate OD-specific flows $\mathbf{x^*_{rs}}$.
 \item Calculate $\mathbf{x^*} \leftarrow \sum_{(r,s) \in Z^2} \mathbf{x^*_{rs}}$
\end{enumerate}
\item Update $\mathbf{x} \leftarrow \lambda \mathbf{x^*} + (1 - \lambda) \mathbf{x}$ for some $\lambda \in [0, 1]$.
\item If $\mathbf{x}$ and $\mathbf{x^*}$ are sufficiently close, terminate.
Otherwise, return to step 2.
\end{enumerate}
(It is often possible to do the calculations in step 3 per origin or per destination, rather than separately for each OD pair.)

Notice that the termination criteria is slightly different than before.
With the classical traffic assignment problem, we argued that it was a \emph{bad} idea to compare the current solution to the previous solution and terminate when they are sufficiently small.
Here, we are not doing that, but are doing something which looks similar: comparing the current solution $\mathbf{x}$ with the ``target'' solution $\mathbf{x^*}$.
This works because, unlike in the regular traffic assignment problem, the mapping $\mathbf{H(C)}$ is always continuous and well-defined.
This means that $\mathbf{x^*}$ varies slowly with $\mathbf{x}$: small changes in link flows mean small changes in link and path travel times, which means small changes in path flows from~\eqn{logitpathexplicit}.
With classic traffic assignment, a small change in link flows could mean a shift in the shortest path, which would result in a dramatic change in $\mathbf{x^*}$, moving all flow from an OD pair to that new shortest path.
Furthermore, when there are ties there are multiple possible $\mathbf{x^*}$ values.
So, in this case we had to introduce auxiliary convergence measures like the relative gap or average excess cost.
SUE is simpler in that we can simply compare the current and target solutions --- in fact, as defined earlier, the relative gap and average excess cost do not make sense, since the equilibrium principle is no longer defined by all travelers being on the shortest path.

There is another significant consequence of the fact that the distance between $\mathbf{x}$ and $\mathbf{x^*}$ shrinks as we approach the equilibrium solution.
In Section~\ref{sec:convexmsa}, we had to use a decreasing step size, with $\lambda \rightarrow 0$ over iterations to ensure convergence to a single point.
In SUE this is not necessary, and convergence is in fact achievable with a \emph{constant} step size, as long as it is not too large.
In practice, you can start with $\lambda = 1$, and continue using it as long as the distance between $\mathbf{x}$ and $\mathbf{x^*}$ is shrinking across iterations.
Whenever you are unable to make further progress, reduce $\lambda$ by half and continue using that new value as a step size for as long as it is effective.
Doing so will provide faster convergence than shrinking $\lambda$ at each iteration.

For all of these reasons, the method of successive averages works much better for SUE than for deterministic assignment.

\stevefig{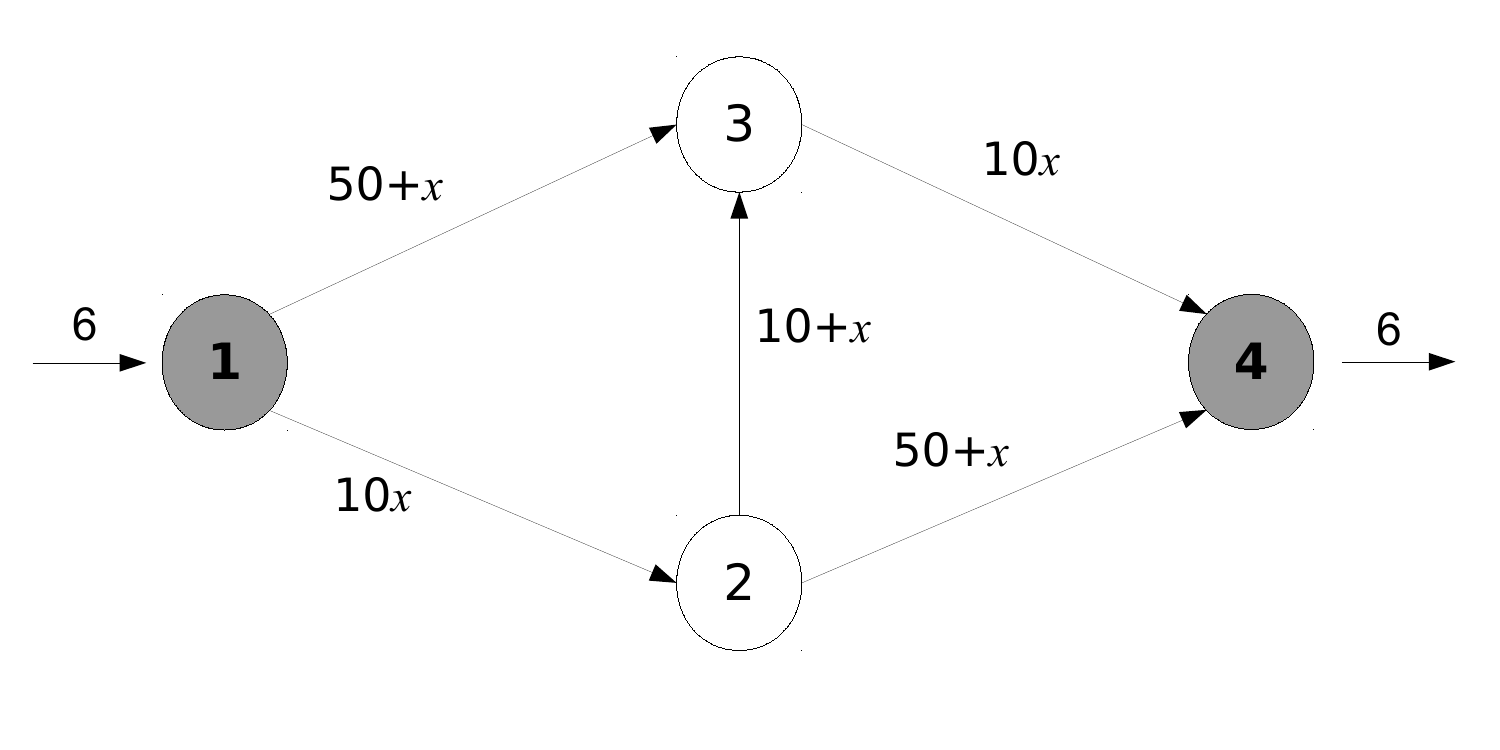}{Network for demonstrating method of successive averages in stochastic assignment.}{0.5\textwidth}

To demonstrate this algorithm, consider the Braess paradox network shown in Figure~\ref{fig:braesssue}, where the demand from node 1 to node 4 is 6 vehicles, $\theta = 0.01$, and all paths are allowable.
Assume that we choose the initial solution $\mb{x}$ by performing Dial's method on this network, using the \emph{free-flow} travel times.
Table~\ref{tbl:suetable} shows the calculations; in this table, we first calculate the $L_i$ values for nodes (based on shortest paths), then the link likelihoods $L_{ij}$; then the node and link weights $W_i$ and $W_{ij}$, and finally the node and link flows $x_i$ and $x_{ij}$.
This gives initial link flows of $\mb{x} = \vect{4.282 & 1.718 & 2.563 & 1.718 & 4.282}$.
Recalculating the link performance functions with these flows gives new link travel times, which we put into Dial's method again, giving the result $\mb{x^*} = \vect{3.976 & 2.024 & 1.951 & 2.024 & 3.976}$.
\emph{Notice that $\mb{x^*}$ and $\mb{x}$ are quite close to each other!}  This is a very different situation than when the method of successive averages is applied to the classical traffic assignment problem (compare with the examples in Section~\ref{sec:convexmsa}), where $\mb{x^*}$ was an extreme-point solution quite far away from $\mb{x}$.

So the link flows $\mb{x}$ are updated by averaging $\mb{x^*}$ into the old $\mb{x}$ values, using a weight of $\lambda = 1/2$, producing the results in the rightmost column of Figure~\ref{fig:braesssue}.
The process is repeated until convergence.
\index{static traffic assignment!stochastic user equilibrium!method of successive averages|)}

\begin{table}
\begin{center}
\caption{Method of successive averages for stochastic user equilibrium. \label{tbl:suetable}}
\begin{tabular}{|c|ccc|cccc|}
\hline
      & \multicolumn{3}{c|}{Initialization} & \multicolumn{4}{c|}{Iteration 1}           \\
\hline     
Nodes & $L_i$    & $W_i$    & $x_i$          & $L_i$    & $W_i$    & $x_i$      &          \\
\hline
1     & 0        & 1        & 6              & 0        & 1        & 6          &          \\
2     & 0        & 1        & 4.282          & 42.8     & 1        & 3.976      &          \\
3     & 10       & 1.670    & 4.282          & 51.7     & 1.964    & 3.976      &          \\
4     & 10       & 2.341    & 6              & 94.5     & 2.964    & 6          &          \\
\hline
Links & $L_{ij}$ & $W_{ij}$ & $x_{ij}$       & $L_{ij}$ & $W_{ij}$ & $x^*_{ij}$ & $x_{ij}$ \\
\hline
(1,2) & 1        & 1        & 4.282          & 1        & 1        & 3.976      & 4.129    \\
(1,3) & 0.670    & 0.670    & 1.718          & 1        & 1        & 2.024      & 1.871    \\
(2,3) & 1        & 1        & 2.563          & 0.964    & 0.964    & 1.951      & 2.257    \\
(2,4) & 0.670    & 0.670    & 1.718          & 1        & 1        & 2.024      & 1.871    \\
(3,4) & 1        & 1.670    & 4.282          & 1        & 1.964    & 3.976      & 4.129    \\
\hline
\end{tabular}
\end{center}
\end{table}

\subsubsection{Disaggregate link-flow formulation and Frank-Wolfe}

\index{static traffic assignment!stochastic user equilibrium!convex optimization formulation|(}
For the classical TAP, the Frank-Wolfe algorithm was much faster than the method of successive averages, because it chose $\lambda$ adaptively, to maximize the reduction in the objective function at each iteration.
In theory, it is possible to do the same thing with the convex program described above.
In practice, it is harder because this program makes use of the path-flow variables $\mb{h}$.
With the classical traffic assignment problem, we could express solutions and the objective solely in terms of the link flows $\mb{x}$.
The addition of the $h \log h$ terms to the SUE objective function renders this impossible.
Furthermore, the number of paths grows very quickly with network size (and if we are choosing the reasonable path set to include fully cyclic paths, the number of paths is usually infinite).

The method of successive averages avoids this problem by not actually referring to the objective function at any point --- if you review the steps above, you see that the objective function is never calculated.
Its role is implicit, guaranteeing that the algorithm will eventually converge, since the direction $\mb{x^*} - \mb{x}$ is one along which the objective is decreasing.
If we can find an efficient way to evaluate the objective function, then we can develop an analogue to the Frank-Wolfe method for classical assignment.

It turns out that we can reformulate the objective function in terms of the destination-aggregated flows on each link $x_{ij}^s$ (see Section~\ref{sec:bushes}),\index{static traffic assignment!bush-based solution} using the Markov property of the logit loading.
There is an equivalent disaggregation by origin, if we reverse our interpretation of the Markov property; see Exercise~\ref{ex:markovreverse}.
The following results were proved in Section~\ref{sec:markovproperty}:\footnote{Technically, this section only proved them for the case of a single origin-destination pair, but they hold for each destination as well.}
\begin{itemize}
\item There exist values $P_{ij}^s$ for each link $(i,j)$ and destination $s$, such that
\labeleqn{markovgeneral}{h^\pi = d^{rs} \prod_{(i,j) \in \pi} (P_{ij}^s)^{\delta_{ij}^\pi} }
for any path $\pi$ connecting origin $r$ to destination $s$.
\item The $P_{ij}^s$ values can be interpreted as the conditional probability that a vehicle arriving at node $i$ and destined for node $s$ will have link $(i,j)$ as the next link in its path.
Therefore,
\labeleqn{markovdecompose}{P_{ij}^s = \frac{x_{ij}^s}{x_i^s} }
for each destination $s$, node $i \neq s$, and outgoing link $(i,j)$.
\end{itemize}

Using these properties, we derive an equivalent formula for the additional term to the Beckmann function.

\begin{prp}Let $\mb{h}$ be a feasible path flow vector satisfying the Markov property, and let $x_{ij}^s$ and $x_i^s$ be the corresponding destination-aggregated link and node flows.
Then
\labeleqn{akamatsudecomposition}{
\sum_{\pi \in \hat{\Pi}} h^\pi \log \myp{ \frac{h^\pi}{d^{rs}}} = \sum_{s \in Z} \myp{
    \sum_{(i,j) \in A} x_{ij}^s \log x_{ij}^s - \sum_{\substack{i \in N \\ i \neq s}} x_i^s \log x_i^s
   }
}
\end{prp}
\begin{proof}
We treat each destination $s$ separately; summing over all destinations gives the result.

Using~\eqn{markovgeneral} and~\eqn{markovdecompose}, we have
\begin{align}
\sum_{\pi \in \hat{\Pi}^s} h^\pi \myp{ \frac{h^\pi}{d^{rs}}} &=
 \sum_{\pi \in \hat{\Pi}^s} h^\pi \log \myp{
    \prod_{(i,j) \in A}   \myp{ \frac{x^s_{ij}}{x^s_i} }^{\delta_{ij}^\pi}
    } \\
&= \sum_{\pi \in \hat{\Pi}^s} h^\pi \sum_{(i,j) \in A} \delta_{ij}^\pi \log \frac{x^s_{ij}}{x^s_i} \\
&= \sum_{(i,j) \in A} \myp{ \sum_{\pi \in \Pi^s} h^\pi \delta_{ij}^\pi } \log \frac{x^s_{ij}}{x^s_i} \\
&= \sum_{(i,j) \in A} x_{ij}^s \log \frac{x^s_{ij}}{x^s_i} \\
&= \sum_{(i,j) \in A} x_{ij}^s \log x_{ij}^s - \sum_{(i,j) \in A} x_{ij}^s \log {x^s_i} \\
&= \sum_{(i,j) \in A} x_{ij}^s \log x_{ij}^s - \sum_{\substack{i \in N \\ i \neq s} } x_{i}^s \log {x^s_i}
\end{align}
by grouping terms in the last sum by the tail node.
\end{proof}

Furthermore, the left-hand side of equation~\eqn{akamatsudecomposition} is very similar to the second term in the objective function~\eqn{fisk}.
In fact,
\labeleqn{simplifyentropy}{
\sum_{\pi \in \hat{\Pi}^s} h^\pi \log \myp{ \frac{h^\pi}{d^{rs}}}
=
\sum_{\pi \in \hat{\Pi}^s} h^\pi \log h^\pi - \sum_{r \in Z} d^{rs} \log d^{rs}
}
using properties of logarithms and the fact that $\sum_{\pi \in \Pi^{rs}} = d^{rs}$.
They differ only by $\sum_{r \in Z} d^{rs} \log d^{rs}$, which is a constant --- and recall from Proposition~\ref{prp:optimizationconstants} that adding a constant does not affect the optimal solution to an optimization problem.

This means that we can replace the objective function~\eqn{fisk} with
\labeleqn{akamatsu}{\sum_{(i,j) \in A} \int_0^{x_{ij}} t_{ij}(x) dx + \frac{1}{\theta}  \sum_{s \in Z} \myp{
    \sum_{(i,j) \in A} x_{ij}^s \log x_{ij}^s - \sum_{\substack{i \in N \\ i \neq s}} x_i^s \log x_i^s
} }
which does not require path enumeration.
\index{static traffic assignment!stochastic user equilibrium!Frank-Wolfe|(}We can thus develop the following analogue to the Frank-Wolfe algorithm:
\begin{enumerate}
\item Choose an initial, feasible destination-aggregated link assignment $\mathbf{x^s}$, and the corresponding aggregated link flows $\mb{x}$.
\item Update link travel times based on $\mathbf{x}$.
\item Calculate target flows:
\begin{enumerate}[(a)]
 \item For each destination $s$, use a method from the previous section to calculate OD-specific flows $\mathbf{x_s^*}$.
 \item Calculate $\mathbf{x^*} \leftarrow \sum_{s \in Z} \mathbf{x_s^*}$
\end{enumerate}
\item Find the value of $\lambda$ minimizing~\eqn{akamatsu} along the line $(1 - \lambda) \vect{\mb{x} & \mb{x^s}} + \lambda \vect{\mb{x^*} & \mb{x^*_s}}$.
\item Update $\mathbf{x} \leftarrow \lambda \mathbf{x^*} + (1 - \lambda) \mathbf{x}$ .
\item If $\mathbf{x}$ and $\mathbf{x^*}$ are sufficiently close, terminate.
Otherwise, return to step 2.
\end{enumerate}
\index{static traffic assignment!stochastic user equilibrium!Frank-Wolfe|)}
\index{static traffic assignment!stochastic user equilibrium!convex optimization formulation|)}

\subsection{Logit loading and most likely path flows (*)}

\emph{(This optional section draws a connection between the logit-based stochastic equilibrium model and the concept of entropy maximization in most likely path flows.)}

\index{static traffic assignment!stochastic user equilibrium!entropy|(}
\index{entropy!and logit loading|(}
The use of the $h \log h$ terms to represent logit assignment in the previous section may have reminded you of the most likely path flows problem in deterministic user equilibrium, discussed in Section~\ref{sec:maxentropy}.
In fact, these two concepts are related to each other quite closely!  This optional section explores this relationship, and describes another algorithm for most likely based path flows.

Recall that the most likely path flows were defined as those maximizing entropy, and solving the following optimization problem:
\begin{align}
\max_\mathbf{h} \qquad & -\sum_{\pi \in \hat{\Pi}} h_\pi \log (h_\pi / d)    & \label{eqn:logitentropy} \\
  \mathrm{s.t.} \qquad & \sum_{\pi \in \Pi} \delta_{ij}^\pi h_\pi = x^*_{ij} & \forall (i,j) \in A \label{eqn:logitentropyeqm} \\
                       & \sum_{\pi \in \hat{\Pi}} h_\pi = d	& \\
                       & h_\pi \geq 0   & \forall \pi \in \Pi \label{eqn:entropylogitnonneg} 
\end{align}
where we have assumed there is a single origin-destination pair for simplicity.

We can transform this into an equivalent optimization that more closely resembles the stochastic user equilibrium optimization problems.
First, we replace maximization by minimization by changing the sign of the objective.
Second, we can remove $d$ from the objective, because
\begin{align}
\sum_{\pi \in \hat{\Pi}} h_\pi \log (h_\pi / d) &= \sum_{\pi \in \hat{\Pi}} h_\pi \log h_\pi - \sum_{\pi \in \hat{\Pi}} h_\pi \log d \\
&= \sum_{\pi \in \hat{\Pi}} h_\pi \log h_\pi - d \log d 
\end{align}
and the constant $d \log d$ can be ignored.
This gives
\begin{align}
\min_\mathbf{h} \qquad & \sum_{\pi \in \hat{\Pi}} h_\pi \log h_\pi \label{eqn:logitentropy2} \\
\mathrm{s.t.} \qquad   & \sum_{\pi \in \Pi} \delta_{ij}^\pi h_\pi = x^*_{ij} & \forall (i,j) \in A \label{eqn:logitentropyeqm2} \\
& \sum_{\pi \in \hat{\Pi}} h_\pi = d & \label{eqn:entropylogitnvlb} \\
& h_\pi \geq 0   & \forall \pi \in \Pi \label{eqn:entropylogitnonneg2} 
\end{align}

We make one more change: rather than insisting that the path flows match the equilibrium link flows exactly, we change constraint~\eqn{logitentropyeqm2} to require that the average travel time across all travelers be $c^*$.
That is, we change the constraint to
\labeleqn{averagepathconstraint}{\sum_{\pi \in \hat{\Pi}} h^\pi c^\pi = c^* d\,.}
Arguing as in the previous section, we can ignore the non-negativity constraint on $h^\pi$, since it will not be binding at optimality.

We then Lagrangianize the other two constraints, introducing multipliers $\theta$ and $\kappa$ for~\eqn{averagepathconstraint} and~\eqn{entropylogitnonneg2}, respectively.\index{Lagrange multipliers!applications}
The resulting Lagrangian function is
\labeleqn{mlpfsuelagrangian}{\mathcal{L}(\mb{h},\theta,\kappa) = \sum_{\pi \in \hat{\Pi}} h^\pi \log h^\pi + \theta \myp{\sum_{\pi \in \hat{\Pi}} h^\pi c^\pi - c^* d} + \kappa \myp{d - \sum_{\pi \in \hat{\Pi}} h^\pi}}
without non-negativity constraints.
The optimality conditions are thus
\begin{align}
\pdr{\mc{L}}{h^\pi} &= 1 + \log h^\pi + \theta c^\pi - \kappa = 0 & \forall \pi \in \hat{\Pi} \\
\pdr{\mc{L}}{\theta} &= \sum_{\pi \in \hat{\Pi}} h^\pi c^\pi - c^* d = 0 & \\
\pdr{\mc{L}}{\kappa} &= \sum_{\pi \in \hat{\Pi}} h^\pi - d = 0 & \\
\end{align}
The last two of these are simply the constraints~\eqn{logitentropyeqm2} and~\eqn{entropylogitnvlb}.
The first can be solved for the path flows, giving
\labeleqn{mlpfsuepathflow}{h^\pi = \exp(\kappa - 1 - \theta c^\pi)\,.}
To satisfy the demand constraint~\eqn{entropylogitnvlb}, $\kappa$ must be chosen so that the sum of~\eqn{mlpfsuepathflow} over all paths gives the total demand $d$.

Omitting the algebra, we must have
\labeleqn{mlpfkappa}{\kappa = \log d (1 - \log \sum_{\pi' \in \Pi} \exp(-\theta c^\pi))}
or
\labeleqn{mlpflogit}{h^\pi = d \frac{\exp(-\theta c^\pi)}{\sum_{\pi' \in \Pi} \exp(-\theta c^\pi)}\,.}
But this is just the logit formula!  The Lagrange multiplier $\theta$ must be chosen to satisfy the remaining constraint on the average path cost.

This derivation shows that the most likely path flows and stochastic user equilibrium problems have a similar underlying structure.
If we relax the requirement that all travelers be on shortest paths, and simply constrain the average cost of travel, the most likely path flows coincide with a logit loading, where the parameter $\theta$ is the Lagrange multiplier for this constraint.
As $\theta$ approaches infinity, the average cost of travel approaches its value at the deterministic user equilibrium solution, and the stochastic user equilibrium path flows approach the most likely path flows in the corresponding deterministic problem.
This provides another algorithmic approach for solving for most likely path flows, in addition to those discussed in Section~\ref{sec:likelypathflowalgorithms}.
In practice, this algorithm is difficult to implement, because of numerical issues that arise as $\theta$ grows large.
\index{static traffic assignment!stochastic user equilibrium!entropy|)}
\index{entropy!and logit loading|)}

\subsection{Alternatives to logit loading}
\label{sec:generalchoice}

\index{logit!limitations}
The majority of this section has focused on the logit model for stochastic route choice, because it demonstrates the main ideas simply, and leads to computationally-efficient solution techniques.
There are several serious criticisms of logit assignment.
For instance, the assumption that the error terms $\epsilon^\pi$ are independent across paths is hard to defend if paths overlap significantly.
Common methods for creating totally acyclic route sets can also create unreasonable artifacts, and allowing all cyclic paths can also be unreasonable if the network topology creates many such paths with low travel times.
See Exercise~\ref{ex:logitproblems} for some concrete examples.
The main alternative to the logit model is the \emph{probit}\index{probit} model, in which the $\epsilon$ terms have a multivariate normal distribution, with a (possibly nondiagonal) covariance matrix to allow for correlation between these terms.

The framework of stochastic user equilibrium can be generalized to other distributions of the error terms, including correlation.
The full development of this general framework is beyond the scope of this book, but we provide an overview and summary.
The objective function~\eqn{fisk} must be replaced with\index{static traffic assignment!stochastic user equilibrium!non-logit formulation}
\labeleqn{sueobjective}{f(\mathbf{x}) = -\sum_{(r,s) \in Z^2} d^{rs} E\mys{\min_{\pi \in \Pi^{rs}} (c^\pi + \epsilon)} + \sum_{(i,j) \in A} x_{ij} t_{ij} - \sum_{(i,j) \in A} \int_0^{x_{ij}} t_{ij}(x)~dx } 
where the expectation\label{not:E} is taken with respect to the ``unobserved'' random variables $\epsilon$, and $t_{ij}$ in the last two terms are understood to be functions of $x_{ij}$.
It can be shown that this function is convex, and therefore that the SUE solution is unique.
However, evaluating this function is harder.
The first term in~\eqn{sueobjective} involves an expectation over all paths connecting an OD pair.
In discrete choice, this is known as the \emph{satisfaction function},\index{satisfaction function} and expresses the expected perceived travel time on a path chosen by a traveler.
In the case of the logit model, this expectation can be computed in closed form; for most distributions it cannot, and must be evaluated through Monte Carlo sampling or another approximation.

The method of successive averages can still be used, even without evaluating the objective function.
Step 4a needs to be replaced with a stochastic network loading, using the current travel times and whatever distribution of $\epsilon$ is chosen.
This often requires Monte Carlo sampling\index{Monte Carlo sampling} as well: for (multiple) samples of the $\epsilon$ terms, the shortest paths can be found using one of the standard algorithms, and the resulting flows averaged together to form an estimate of $\mb{x^*}$.

Because we are using an estimate of $\mb{x^*}$, it is possible that the target direction is not exactly right, and that it is not in a direction in which $f(\mathbf{x})$ is decreasing.
Nevertheless, as long as it is correct ``on average'' (i.e., the sampling is done in an unbiased manner), one can show that the method of successive averages will still converge to the stochastic user equilibrium solution.
\index{static traffic assignment!stochastic user equilibrium|)}

\section{Historical Notes and Further Reading}
\label{sec:staticextensions_references}

\index{static traffic assignment!elastic demand}
The original formulation of the user equilibrium traffic assignment problem in \cite{beckmann56} actually modeled demand as elastic.
(In this book we chose to focus on the fixed-demand problem in previous chapters, which is a special case.)  The transformation of the elastic demand problem to a fixed-demand problem with artificial links was reported in \cite{gartner80}, based on earlier work by \cite{murchland70}.
While the original Frank-Wolfe algorithm can be applied to the elastic demand equilibrium problem, the modified version presented in this chapter is faster.
It is essentially a specialized version of the double-stage algorithm of \cite{evans76} developed for a combined trip distribution and assignment model.

\index{static traffic assignment!link interactions}
The variational inequality formulation of the traffic assignment problem with link interactions is due to \cite{smith79}, \cite{dafermos80}, and \cite{aashtiani80}.
The diagonalization method was demonstrated in \cite{fisk82}, and its convergence proved by \cite{dafermos82_relax}.
The simplicial decomposition method is adapted from \cite{smith83a}.
The example with multiple equilibria in the simple merge network is taken from \cite{boylesvortices}.
Other examples of multiple equilibria are given in \cite{netter72} and \cite{marcotte04b}.

\index{static traffic assignment!stochastic user equilibrium}
Stochastic network loading using the logit model and totally acyclic path sets was described in \cite{dial71}.
Logit loading with the full cyclic path set was developed by \cite{bell95}, \cite{akamatsu96}, and \cite{akamatsu97}.
The stochastic user equilibrium model was first proposed in \cite{daganzo77}.
\cite{powell82} proved the convergence of the method of successive averages for this problem, for all well-behaved distributions of the unobserved utility term.
The convex programming case for the logit model was given in \cite{fisk80}, and the optimization formulation for more general distributions was given in \cite{daganzo82}.

\index{static traffic assignment!and trip distribution}
\index{static traffic assignment!and mode choice}
There are other extensions to the basic traffic assignment problem which are not treated in this chapter.
One major class of these extensions integrates other steps of the planning process with user equilibrium route choice.
\cite{bruynooghe68}, \cite{florian75}, and \cite{evans76} describe models and algorithms that integrate trip distribution with route choice.
\cite{abdulaal79_combine} and \cite{dafermos82_relax} integrated mode choice with equilibrium traffic assignment.
\cite{florian78} present a model integrating all three of the latter steps in the four-step model (trip distribution, mode choice, and route choice); \cite{sheffi78,sheffi80} showed how network transformations can encode these choices as a standard traffic assignment problem, with appropriately chosen artificial links and link performance function.

\index{static traffic assignment!multiclass}
Traffic assignment can also be solved as a \emph{multiclass} problem, in which different groups of travelers have different route choice behaviors~\cite{dafermos72}.
For instance, one group of ``selfish'' travelers may route according to user equilibrium principles, and another group of ``selfless'' travelers may route according to system optimum~\citep{roughgarden02,yang07,sharon18}.
More commonly, all travelers aim to minimize their travel times, but perhaps have different levels of perception error (different $\theta$ values in stochastic user equilibrium), as in \cite{huang07}, or have different values of time.
The latter is important in networks with tolls or other monetary costs on links, as different travelers will trade off travel time and monetary cost differently.
If this is the only ways vehicle classes differ, it is not too hard to solve the case with a finite number of values of time~\citep{nagurney00}.
Researchers have also investigated the case where the value of time is continuously distributed among the population~\citep{dial96,dial97,marcotte98_t2}.
The general multiclass equilibrium problem is hard to solve without relatively strong assumptions on how the classes interact with each other~\citep{hammond84,florian82,toint96,marcotte04b}.

\index{static traffic assignment!side constraints}
Some researchers have extended the traffic assignment problem to include \emph{side constraints}, such as enforcing link capacities\index{link!capacity} strictly ($x_{ij} \leq u_{ij}$ on each link).
These constraints require changes to the definition of user equilibrium; used paths may have different travel times if the shorter of them is at capacity and cannot accept additional flow.
Corresponding algorithmic changes are also needed.
See \cite{larsson95}, \cite{larsson99}, \cite{larsson04}, \cite{prashker04}, and \cite{feng20} for examples of formulations and algorithms that handle this variant.

\index{static traffic assignment!travel time uncertainty}
There are also variants of traffic assignment which aim to represent uncertainty in the model inputs, both on the demand side~\citep[e.g.,][]{clark05}, and on the supply side (to reflect incidents or other disruptions).
One instance of the latter is the \emph{user equilibrium with recourse} model, in which links exist in discrete states according to a given probability distribution (e.g., normal operating conditions, mild incident, severe incident), and each state is associated with a different link performance function.
Travelers receive information on link states as they travel, and can update their paths \emph{en route} based on the information they receive.
Such models are useful for identifying locations of variable message signs or other traveler information devices~\citep{boyles_crp}.
\cite{unnikrishnan_uer} first presented the user equilibrium with recourse model, and \cite{rambha18} the system optimal version.
The solution algorithm in \cite{unnikrishnan_uer} is link-based; a faster method (analogous to bush-based methods for the basic problem) is presented in \cite{rambha17}.
\index{user equilibrium!with recourse}

\section{Exercises}
\label{exercises_staticextensions}

\begin{enumerate}
\item \diff{25} Verify that each of the demand functions $D$ below is strictly decreasing and bounded above (for $\kappa \geq 0$), then find the inverse functions $D^{-1}(d)$.
\begin{enumerate}
\item $D(\kappa) = 50 - \kappa$
\item $D(\kappa) = 1000 / (\kappa + 1)$
\item $D(\kappa) = 50 - \kappa^2$
\item $D(\kappa) = (\kappa + 2) / (\kappa + 1)$
\end{enumerate}
\item \diff{13} The total misplaced flow reflects consistency of a solution with (in this case, that the OD matrix should be given by the demand functions).
Suggest another measure for how close a particular OD matrix and traffic assignment are to satisfying this consistency condition.
Your proposed measure should be a nonnegative and continuous function of values related to the solution (e.g., $d^{rs}$, $x_{ij}$, $\kappa^{rs}$, etc.), which is zero if and only if the OD matrix is completely consistent with the demand functions.
Compare your new measure with the total misplaced flow, and comment on any notable differences.
\item \diff{31} Verify that the elastic demand objective function~\eqn{elasticfirst} is convex, given the assumptions made on the demand functions.
\item \diff{47} Using the network in Figure~\ref{fig:braessed}, solve the elastic demand equilibrium problem with demand given by $d^{14} = 15 - \kappa^{14}/30$.
Perform three iterations of the Frank-Wolfe method, and report the average excess cost and total misplaced flow for the solution.
\item \diff{48} Using the network in Figure~\ref{fig:braessed}, solve the elastic demand equilibrium problem with demand given by $d^{14} = 10 - \kappa^{14}/30$.
\label{ex:braessed}
\begin{enumerate}
\item Perform three iterations of the Frank-Wolfe method designed for elastic demand (Section~\ref{sec:fwed}).
\item Transform the problem to an equivalent fixed-demand problem using the Gartner transformation from Section~\ref{sec:gartner}, and perform three iterations of the original Frank-Wolfe method.
\item Compare the performance of these two methods: after three iterations, which is closer to satisfying the equilibrium and demand conditions?
\end{enumerate}
\stevefig{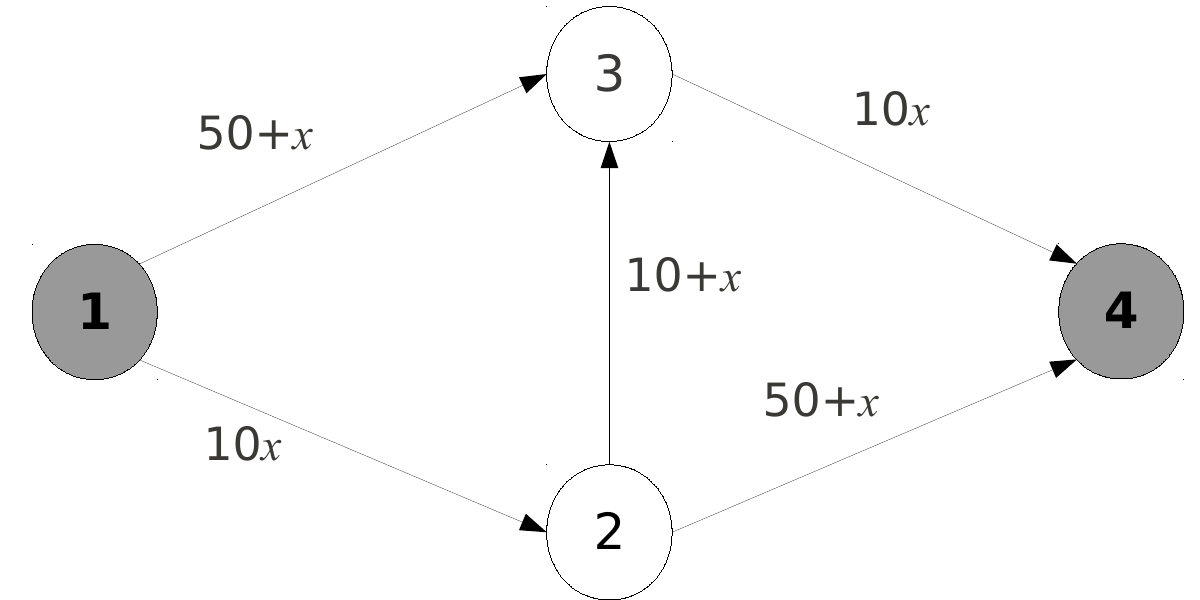}{Network for Exercise~\ref{ex:braessed}.}{0.6\textwidth}
\item \diff{49} Using the network in Figure~\ref{fig:exgridneted}, solve the elastic demand equilibrium problem with demand functions $q^{19} = 1000 - 50(u^{19} - 50)$ and $q^{49} = 1000 - 75(u^{49} - 50)$.
The cost function on the light links is $3 + (x_a/200)^2$, and the cost function on the thick links is $5 + (x_a/100)^2$.
1000 vehicles are traveling from node 1 to 9, and 1000 vehicles from node 4 to node 9.
Perform three iterations of the Frank-Wolfe method and report the link flows and OD matrix.
\label{ex:exgridneted}
\stevefig{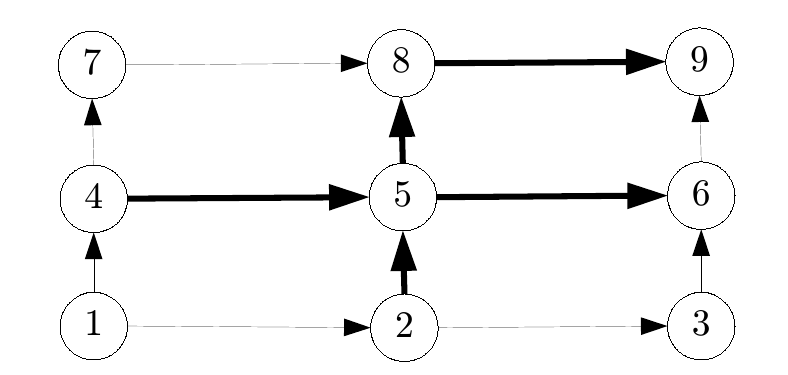}{Network for Exercise~\ref{ex:exgridneted}.}{0.6\textwidth}
\item \diff{13} Assume that $(\mb{x}, \mb{d})$ and $(\mb{x^*}, \mb{d^*})$ are both feasible solutions to an elastic demand equilibrium problem.
Show that $(\lambda \mb{x} + (1 - \lambda) \mb{x^*}, \lambda \mb{d} + (1 - \lambda) \mb{d^*})$ is also feasible if $\lambda \in [0, 1]$.
This ensures that the Frank-Wolfe solutions are always feasible, assuming we start with $\mb{x}$ and $\mb{d}$ values which are consistent with each other, and always choose targets in a consistent way.
\item \diff{57} (Calculating derivative formulas.)  Let $\mathbf{x}$ and $\mathbf{d}$ be the current, feasible, link flows and OD matrix, and let $\mathbf{x^*}$ and $\mathbf{d^*}$ be any other feasible link flows and OD matrix.
Let $\mathbf{x'} = \lambda \mathbf{x^*} + (1 - \lambda) \mathbf{x}$ and $\mathbf{d'} = \lambda \mathbf{d^*} + (1 - \lambda) \mathbf{d}$.
\begin{enumerate}
    \item Let $f(\mathbf{x'}, \mathbf{d'})$ be the objective function for the elastic demand equilibrium problem.
    Recognizing that $\mathbf{x'}$ and $\mathbf{d'}$ are functions of $\lambda$, calculate $df/d\lambda$.
    \item For the elastic demand problem, show that $\frac{df}{d\lambda} \big|_{\lambda = 0} \leq 0$ if $\mathbf{d^*}$ and $\mathbf{x^*}$ are chosen in the way given in the text.
    That is, the objective function is nonincreasing in the direction of the ``target.'' (You can assume that the demand function values are strictly positive if that would simplify your proof.)
\end{enumerate}
\item \diff{24} Consider a two-link network.
For each pair of link performance functions shown below, determine whether or not the symmetry condition~\eqn{symmetry} is satisfied.
\begin{enumerate}
    \item $t_1 = 4 + 3x_1 + x_2$, $t_2 = 2 + x_1 + 4x_2$
    \item $t_1 = 7 + 3x_1 + 4x_2$, $t_2 = 12 + 2x_1 + 4x_2$
    \item $t_1 = 4 + x_1 + 3x_2$, $t_2 = 2 + 3x_1 + 2x_2$	
    \item $t_1 = 3x^2_1 + x_2$, $t_2 = 4 + x_1 + 4x^3_2$
    \item $t_1 = 3x^2_1 + 2x^2_2$, $t_2 = 2x^2_1 + 3x^3_2$		
    \item $t_1 = 50 + x_1$, $t_2 = 10x_2$
\end{enumerate}
\item \diff{34} Determine which of the pairs of link performance functions in the previous exercise are strictly monotone.
\item \diff{49} Consider the network in Figure~\ref{fig:exgridneted} with a fixed demand of 1000 vehicles from 1 to 9 and 1000 vehicles from 4 to 9.
The link performance function on every link $a$ arriving at a ``merge node'' (that is, nodes 5, 6, 8, and 9) is $4 + (x_a / 150)^2 + (x_{a'} / 300)^2$ where $a'$ is the other link arriving at the merge.
Verify that the link interactions are symmetric, and perform five iterations of the Frank-Wolfe method.
Report the link flows and travel times.
\item \diff{49} Consider the network in Figure~\ref{fig:exgridneted} with a fixed demand of 1000 vehicles from 1 to 9 and 1000 vehicles from 4 to 9.
The link performance function on every link $a$ arriving at a ``merge node'' (that is, nodes 5, 6, 8, and 9) is  $3 + (x_a / 200)^2 + (x_{a'} / 400)^2$ if the link is light, and $5 + (x_a / 100)^2 + (x_{a'} / 200)^2$ if the link is thick, where $a'$ is the other link arriving at the merge.
Show that not all link interactions are symmetric, and then perform five iterations of the diagonalization method and report the link flows, travel times, and relative gap.

\item \diff{49} See the network in Figure~\ref{fig:interactions2}, where the numbers indicate the label for each link.
The link performance functions are:
\begin{align*}
    t_1(\mathbf{x}) & = 1 + 4x_1 + 2x_2 \\
    t_2(\mathbf{x}) & = 2 + x_1 + 2x_2 \\
    t_3(\mathbf{x}) & = 3 + 2x_3 + x_4 \\
    t_4(\mathbf{x}) & = 2 + 3x_4 \\
\end{align*}
and the demand from node A to node C is 10 vehicles.
For both methods below, start with an initial solution loading all flow on links 1 and 4.
\label{ex:interactions2}
\begin{enumerate}
    \item Use three iterations of the diagonalization method to try to find an equilibrium solution, and report the link flows and average excess cost.
(As before, this means finding three $\mathbf{x^*}$ vectors \emph{after} your initial solution.)
    \item Use three iterations of simplicial decomposition to try to find an equilibrium solution, and report the link flows and average excess cost.
Three iterations means that $\mathcal{X}$ should have three vectors in it when the algorithm terminates (unless it terminates early because the $\mathbf{x^*}$ you find is already in $\mathcal{X}$).
For each subproblem, make the number of improvement steps one less than the size of $\mathcal{X}$ (so when $\mathcal{X}$ has 1 vector, perform 0 steps; when it has 2 vectors, perform 1 step, and so on).
For each of these steps, try the sequence of $\mu$ values $1/2, 1/4, 1/8, \ldots$, choosing the first that reduces the restricted average excess cost.
\end{enumerate}
\stevefig{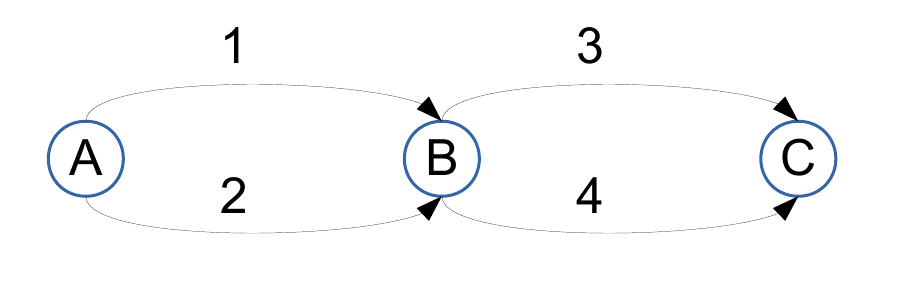}{Network for Exercise~\ref{ex:interactions2}.}{0.6\textwidth}
\item \diff{64} Prove Proposition~\ref{prp:strictmonotoneiffpd}.
\item \diff{42} Show that the three methods for generating paths described in Section~\ref{sec:totallyacyclic} indeed yield totally acyclic path sets, and that they always include the shortest path from the origin to the destination.
\label{ex:acyclic}
\item \diff{32} Show that any set of totally acyclic paths satisfies the segment substitution property.
Then show that the set of all paths (including all cyclic paths) also satisfies this property.
\label{ex:segmentsubstitution}
\item \diff{53} Complete the example following equation~\eqn{logitdenomexample} by forming the infinite sum for the flow on links (2,3) and (3,2), and showing that they are equal to $1/2$.
\item \diff{62} Reformulate the Markov property, and the results in Section~\ref{sec:markovproperty} in terms of conditional probabilities for the link a vehicle used to \emph{arrive} at a given node, rather than the link a vehicle will choose to depart a given node.
\label{ex:markovreverse}
\item \diff{34} In the logit formula~\eqn{logitpath}, show that the same path choice probabilities are obtained if each link travel time $t_{ij}$ is replaced with $L_i + t_{ij} - L_j$, where $\mb{L}$ is a vector of node-specific constants.
(In practice, these are usually the shortest path distances from the origin.)  \label{ex:rescale}
\item \diff{51} What would have to change in the ``full cyclic path set'' stochastic network loading procedure, if there were multiple links with the same tail and head nodes? \label{ex:parallellinks}
\item \diff{35} Consider the network in Figure~\ref{fig:dialgrid} with 1000 vehicles traveling from node A to node I, where the link labels are the travel times.
Use the first criterion in Section~\ref{sec:totallyacyclic} to define the path set, and identify the link flows corresponding to these travel times.
Assume $\theta = 1$.
\label{ex:dialgrid}
\item \diff{35} Repeat Exercise~\ref{ex:dialgrid} for the third definition of a path set.
Assume $\theta = 1$.
\label{ex:dialgrid2}
\stevefig{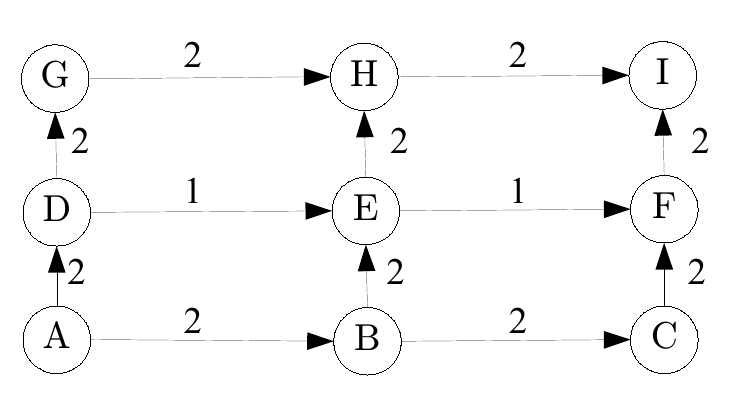}{Network for Exercises~\ref{ex:dialgrid} and~\ref{ex:dialgrid2}.}{0.6\textwidth}
\item \diff{34} Show that the objective function~\eqn{fisk} is strictly convex.
\label{ex:fiskunique}
\item \diff{49} Consider the network in Figure~\ref{fig:exgridneted} with a fixed demand of 1000 vehicles from 1 to 9 and 1000 vehicles from 4 to 9.
The light links have delay function $3 + (x_{ij}/200)^2$, and the dark links have delay function $5 + (x_{ij}/100)^2$.
Assume that drivers choose paths according to the stochastic user equilibrium principle, with $\theta = 1/5$.
Perform three iterations of the method of successive averages, and report the link flows.
Assume all paths are allowed.
\item \diff{24} (Limitations of logit assignment).\index{logit!limitations|(}
This problem showcases three ``problem instances'' for the stochastic network loading models described in this chapter (the logit model, and one proposed definition of allowable paths).
Throughout, assume that the first definition in Section~\ref{sec:totallyacyclic} is used to define the allowable path set.
These instances motivated the development of probit and other, more sophisticated, stochastic equilibrium models.
\label{ex:logitproblems}
\begin{enumerate}
\item Consider the network in Figure~\ref{fig:logitproblems}(a), where the numbers by each link represents the travel time and $z \in (0, 1)$.
Let $p_\uparrow$ represent the proportion of vehicles choosing the top path.
In a typical probit model, we have $p^\uparrow_{PROBIT} = 1 - \Phi\left( \sqrt{\frac{z}{2\pi - z}}\right)$, where $\Phi(\cdot)$\label{not:Phi} is the standard cumulative normal distribution function.
Calculate $p^\uparrow_{LOGIT}$ for the \emph{logit} model as a function of $z$, and plot $p^\uparrow_{PROBIT}$ and $p^\uparrow_{LOGIT}$ as $z$ varies from 0 to 1.
Which do you think is more realistic, and why?
\item Consider networks of the type shown in Figure~\ref{fig:logitproblems}(b), where there is a top path consisting of a single link, and a number of bottom paths.
The network is defined by an integer $m$; there are $m - 1$ intermediate nodes in the bottom paths, and each consecutive pair of intermediate nodes is connected by two parallel links with travel time $8 / m$.
In a common probit model, $p^\uparrow_{PROBIT} = \Phi(-0.435\sqrt{m})$.
What is $p^\uparrow_{LOGIT}$ as a function of $m$?  Again create plots for small values of $m$ (say 1 to 8), indicate which you think is more realistic, and explain why.
\item In the networks in Figure~\ref{fig:logitproblems}(c), identify the proportion of travelers choosing each link if $\theta = 1$.
The left and right panels show a network before and after construction of a new link.
Again identify the proportion of travelers choosing each link if $\theta = 1$.
Do your findings seem reasonable?
\end{enumerate}
\index{logit!limitations|)}
\stevefig{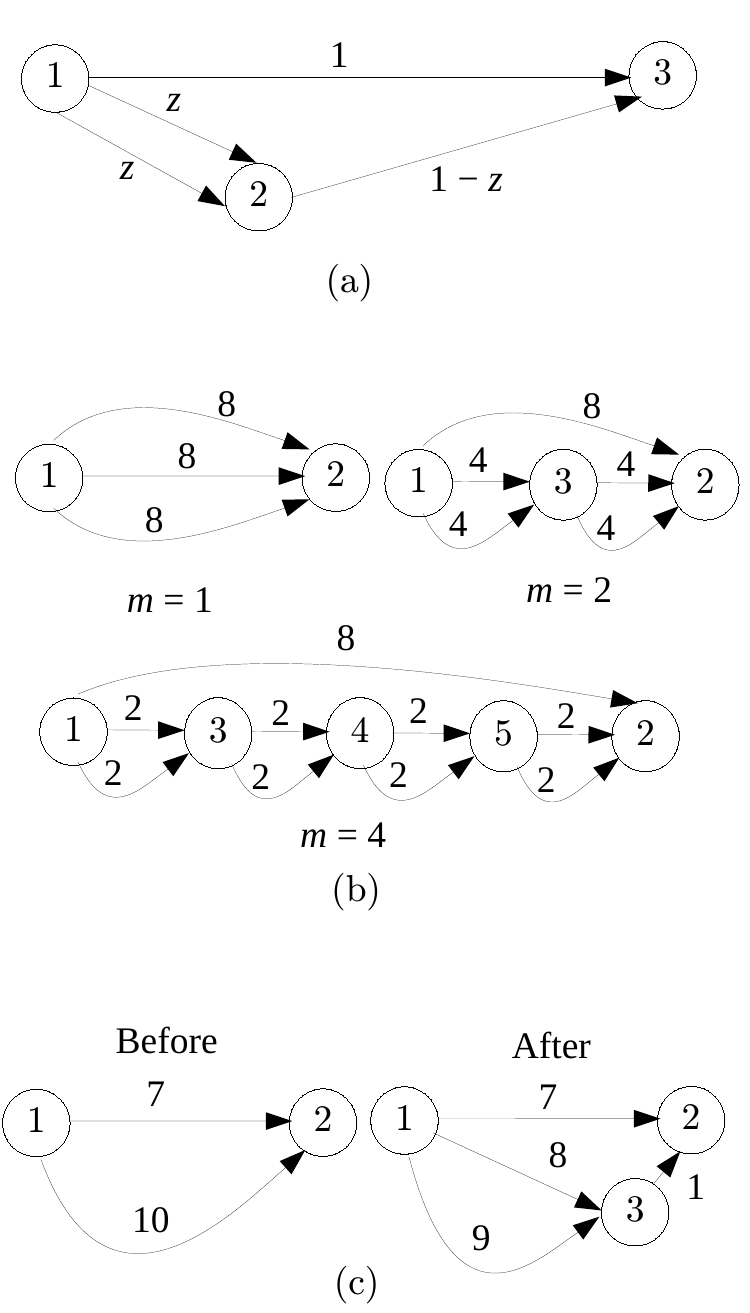}{Networks for Exercise~\ref{ex:logitproblems}.
The label on each link is its \emph{constant} travel time.}{0.6\textwidth}
\end{enumerate}
\index{static traffic assignment|)}

\part{Dynamic Traffic Assignment}
\label{part:dynamictrafficassignment}

\chapter{Network Loading}
\label{chp:networkloading}

\index{dynamic traffic assignment|(}
\index{network loading|(}
\index{dynamic traffic assignment!network loading|see {network loading}}
\index{dynamic traffic assignment!link model|see {link model}}
\index{dynamic traffic assignment!node model|see {node model}}
This chapter discusses \emph{network loading}, the process of modeling the state of traffic on a network, given the route and departure time of every vehicle.
In static traffic assignment, this is a straightforward process based on evaluating link performance functions.
In dynamic traffic assignment, however, network loading becomes much more complicated due to the additional detail in the traffic flow models --- but it is exactly this complexity which makes dynamic traffic assignment more realistic than static traffic assignment.
Rather than link performance functions, dynamic network loading models generally rely on some concepts of traffic flow theory.
There are a great many theories, and an equally great number of dynamic network loading models, so this chapter will focus on those most commonly used.

To give the general flavor of network loading, we start with two simple \emph{link models}, the point queue and spatial queue (Section~\ref{sec:simplelinks}), which describe traffic flow on a single link.
We next present three simple \emph{node models} describing how traffic streams behave at junctions (Section~\ref{sec:nodemodels}).
With these building blocks we can perform network loading on simple networks, showing how link and node models interact to represent traffic flow in a modular way (Section~\ref{sec:combiningnodelink}).

However, the point queue and spatial queue models have significant limitations in representing traffic flow.
The most common network loading models for dynamic traffic assignment are based on the hydrodynamic model of traffic flow, reviewed in Section~\ref{sec:trafficflowtheory}.
This theory is based in fluid mechanics, and assumes that the traffic stream can be modeled as the motion of a fluid, but it can be derived from certain car-following models as well, which have a more behavioral basis.
The cell transmission model and link transmission model are link models based on this theory, and both of these are discussed in Section~\ref{sec:fancylinks}.
Section~\ref{sec:fancynodes} concludes the chapter with a discussion of more sophisticated node models that can represent general intersections.

This chapter aims to present several alternative network loading schemes as part of the general dynamic traffic assignment framework in Figure~\ref{fig:iterativeeqm10}, so that any of them can be combined with the other steps in a flexible way.

\stevefig{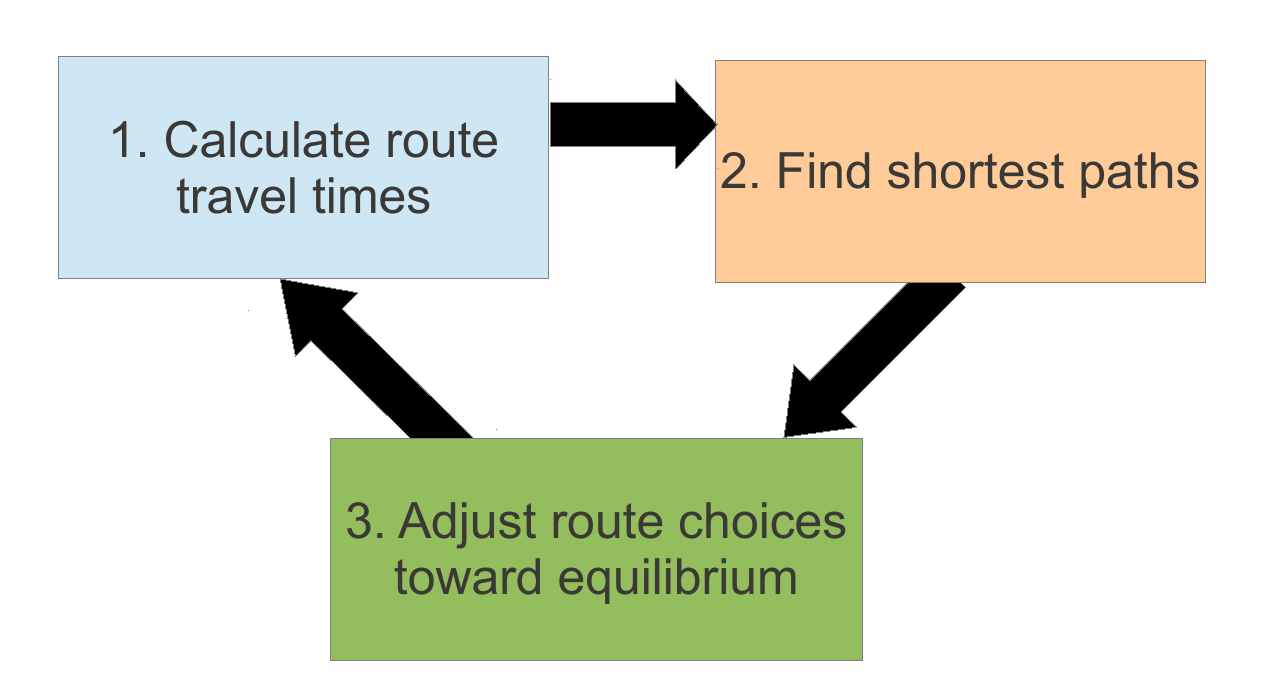}{Three-step iterative process for dynamic traffic assignment.}{0.7\textwidth}

\section{Link Model Concepts}
\label{sec:simplelinks}

\index{link model|(}
Link models are efficient ways to represent congestion and traffic flow.
Efficiency is key, since link models must be implemented on each link in a large network over the entire analysis period; and furthermore, due to the iterative nature of solving for equilibrium, the network loading must be repeated many times with different path flow values as input.
For both of these reasons, a link model involving complicated computations will limit the size of the network you can model.
This section presents the basic concepts of link models, and demonstrates them in the point queue and spatial queue models.
More sophisticated link models will be presented later in this chapter, in Section~\ref{sec:fancylinks}.

In this section, we are solely concerned with a single link.
Call the length of this link $L$,\label{not:L} so $x = 0$\label{not:xcoord} corresponds to the upstream end of the link, and $x = L$ corresponds to the downstream end.
\index{dynamic traffic assignment!time discretization|(}All of the link models we discuss in this book operate in \emph{discrete time}.
That is, we divide the analysis period into small time intervals of length $\Delta t$, and assume uniform conditions within each time interval.
These time intervals are very small relative to the analysis period; at a minimum, no vehicle should be able to traverse more than one link in a single time interval, so $\Delta t$\label{not:Deltat} can be no greater than the shortest free-flow travel time on a link.
In practice, $\Delta t$ values on the order of 5--10 seconds are a reasonable choice.
For convenience in this section, we assume that the unit of time is chosen such that $\Delta t = 1$, so we will only track the state of the network for times $t \in \myc{0, 1, 2, \ldots, T}$,\label{not:tindex} where $T$ is the length of the analysis period (in units of $\Delta t$).

In any discrete time model, it is important to clarify what exactly we mean when we index a variable with a time interval, such as $N(t)$ or $y(t)$.
Does this refer to the value of $x$ at the start of the $t$-th time interval, the end of the interval, the average of a continuous value through the interval, or something else?  This chapter will consistently use the following convention: when referring to a quantity \emph{measured at a single instant in time}, such as the number of vehicles on a link, or the instantaneous speed of a vehicle, we will take such measurements at the \emph{start} of a time interval.
When referring to a quantity \emph{measured over time}, such as the number of vehicles passing a fixed point on the link, we take such measurements \emph{over} the time interval.
(Figure~\ref{fig:discretekey}) This distinction can also be expressed by referring, say, to the number of vehicles on a link \emph{at} time $t$ (meaning the start of the $t$-th interval), or the number of vehicles that exit a link \emph{during} time $t$.
\index{dynamic traffic assignment!time discretization|)}

\stevefig{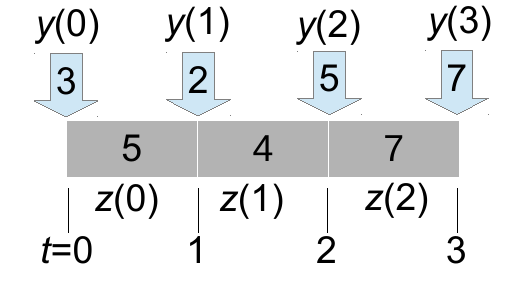}{Convention for indexing discrete time intervals for two hypothetical variables $y$ (measured at an instant) and $z$ (measured over time).}{0.4\textwidth}

We load the network in increasing order of time.
That is, we start with the network state at $t = 0$ (usually assuming an empty condition).
Then, with these values known, we compute the network state at $t = 1$, then at $t = 2$, and so forth.
That is, for any point in time, any values at an earlier point in time can be treated as ``known'' values.

The definition of the ``network state'' at a timestep is intentionally left a bit vague at this point.
It contains all information necessary for modeling traffic flow (e.g., the locations of vehicles, traffic signal indications).
The collection of network states at all time steps must also have enough information to check consistency with the equilibrium principle (cf.\ Section~\ref{sec:notionofequilibrium}) after the network loading is complete, and for adjusting vehicle flows if equilibrium is not satisfied.
The exact components of the network state vary from one model to the next (for instance, signals may be modeled in detail, approximately, or not at all), and as you read about the link and node models presented in this chapter, you should think about what information you would need to store in order to perform the computations for each link and node model.

At a minimum, it is common to record the cumulative number of vehicles\index{cumulative count}\index{link!cumulative count|see {cumulative count}} which have entered and left each link at each timestep, since the start of the modeling period.
These values are denoted by $N^\uparrow(t)$\label{not:Nup} and $N^\downarrow(t)$,\label{not:Ndown} respectively; the arrows are meant as a mnemonic for ``upstream'' and ``downstream'', since they can be thought of as counters at the ends of the links.
So, for instance, $N^\uparrow(3)$ is the number of vehicles which have crossed the upstream end of the link by the start of the 3rd time interval (cumulative entries), and $N^\downarrow(5)$ is the number of vehicles which have crossed the downstream end of the link by the start of the 5th time interval (cumulative exits).
We will assume that the network is empty when $t = 0$ (no vehicles anywhere), and as a result $N^\uparrow(t) - N^\downarrow(t)$ gives the number of vehicles currently on the link at any time $t$.
The values $N^\uparrow$ and $N^\downarrow$ are defined at the start of each time step, for integer values of $t$.
Some formulas may call for the value of a discrete variable at a non-integer point in time, such as $N^\uparrow(3.4)$, in which case a linear interpolation is used between the neighboring values $N^\uparrow(3)$ and $N^\uparrow(4)$.
If possible, the time step should be chosen to minimize or eliminate these interpolation steps, which are time-consuming and can introduce numerical errors.

It is possible to use different time interval lengths for different links and nodes, and this can potentially reduce the necessary computation time.
There are also \emph{continuous time} dynamic network loading models, where flows and other traffic variables are assumed to be functions defined for any real time value, not just a finite number of points.
These are not discussed here to keep the focus on the basic network loading ideas, and to avoid technical details associated with infinitesimal calculations.

\subsection{Sending and receiving flow}

\index{link model!demand|see {link model, sending flow}}
\index{link model!supply|see {link model, receiving flow}}
The main outputs of a link model are the \emph{sending flow} and \emph{receiving flow}, calculated for each discrete time interval.
Some authors use the term \emph{demand} in place of ``sending flow,'' and \emph{supply} in place of ``receiving flow.''
In this part of the book, we reserve the term demand to refer to values in an origin-destination matrix, but if you read additional work in the field you should be aware of both terminologies.

\index{link model!sending flow|(}
The sending flow at time $t$, denoted $S(t)$,\label{not:S} is the number of vehicles which would leave the link during the $t$-th time interval (that is, between times $t$ and $t + \Delta t$) if there was no obstruction from downstream links or nodes (you can imagine that the link is connected to a wide, empty link downstream).
You can also think of this as the flow that is ready to leave the link during this time interval.
The sending flow is calculated at the downstream end of a link.

To visualize sending flow, Figure~\ref{fig:sendingflowexample} shows two examples of links.
In the left panel, there is a queue at the downstream end of the link.
If there is no restriction from downstream, the queue would discharge at the full capacity of the link, and the sending flow would be equal to the capacity of the link multiplied by $\Delta t$.
In the right panel, the link is uncongested and vehicles are traveling at free-flow speed.
The sending flow will be less than the capacity, because relatively few vehicles are close enough to the downstream end of the link to exit within the next time interval.
The vertical line in the figure indicates the distance a vehicle would travel at free-flow speed during one time step, so in this case the sending flow would be 3 vehicles.
Note that the \emph{actual} number of vehicles which can leave the link during the next time step may be less, depending on downstream conditions: perhaps the next link is congested, or perhaps there is a traffic signal at the downstream end.
These considerations are irrelevant for calculating sending flow, which only depends on the link itself.
Node models, introduced in Section~\ref{sec:nodemodels}, will account for constraints from other links in the network.
\index{link model!sending flow|)}

\stevefig{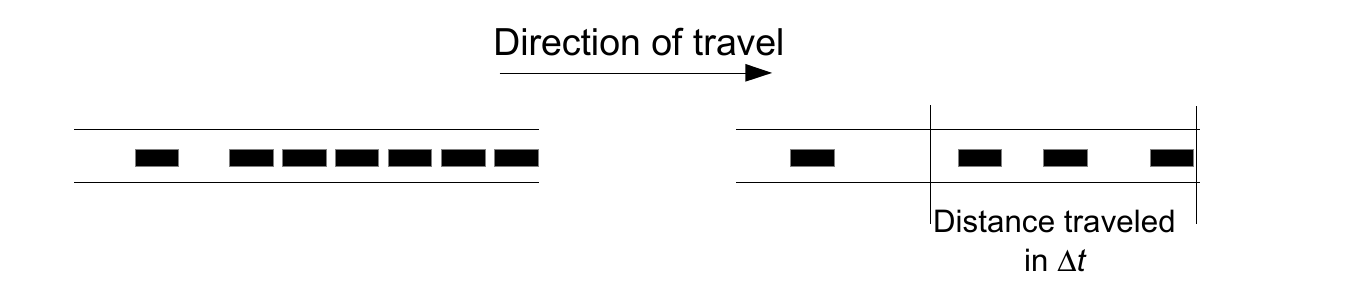}{Calculating the sending flow when there is a queue on the link (left) and when the link is at free-flow (right).}{\textwidth}

\index{link model!receiving flow|(}
The receiving flow during time $t$, denoted $R(t)$,\label{not:R} is the number of vehicles which would enter the link during the $t$-th time interval if the upstream link could supply a very large (even infinite) number of vehicles: you can imagine that the upstream link is wide, and completely full at jam density.
You can also think of this as the maximum amount of flow which can enter the link during this interval, taking into account the available free space.
The receiving flow is calculated at the upstream end of the link.

To visualize receiving flow, Figure~\ref{fig:receivingflowexample} shows two examples of links.
In the left panel, the upstream end of the link is empty.
This means that vehicles could enter the link at its full capacity, and the receiving flow would equal the link's capacity multiplied by $\Delta t$.
The actual number of vehicles which will enter the link may be less than this, if there are few vehicles upstream --- like the sending flow, the receiving flow is simply an upper bound indicating how many vehicles could potentially enter the link.
In the right panel, there is a stopped queue which nearly fills the entire link.
Here the vertical line indicates how far into the link a vehicle would travel at free-flow speed during one time step.
Assuming that the stopped vehicles remain stopped throughout the $t$-th time interval, the receiving flow is the number of vehicles which can physically fit into the link, in this case 2.
\index{link model!receiving flow|)}

\stevefig{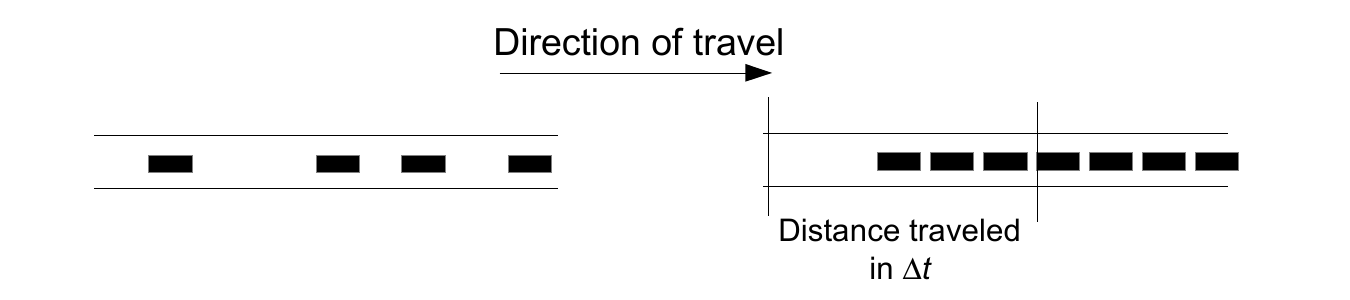}{Calculating the receiving flow when the link is at free-flow (left) and when there is a queue on the link (right).}{\textwidth}

Each link model has a slightly different way of calculating the sending and receiving flows, which correspond to different assumptions on traffic behavior with the link, or to different calculation methods.
The next two subsections present simple link models.
Notice how the different traffic flow assumptions in these models lead to different formulas for sending and receiving flow.

\subsection{Point queue}
\label{sec:pointqueue}

\index{link model!point queue|(}
Point queue models divide each link into two sections:
\begin{itemize}
\item A \textbf{physical section} which spans the length of the link and is assumed uncongestible: vehicles will always travel over this section at free-flow speed.
\item A \textbf{point queue} at the downstream end of the link which occupies no physical space, but conceptually holds vehicles back to represent any congestion delay on the link.
\end{itemize}
These are shown in Figure~\ref{fig:pqonlink}.

\stevefig{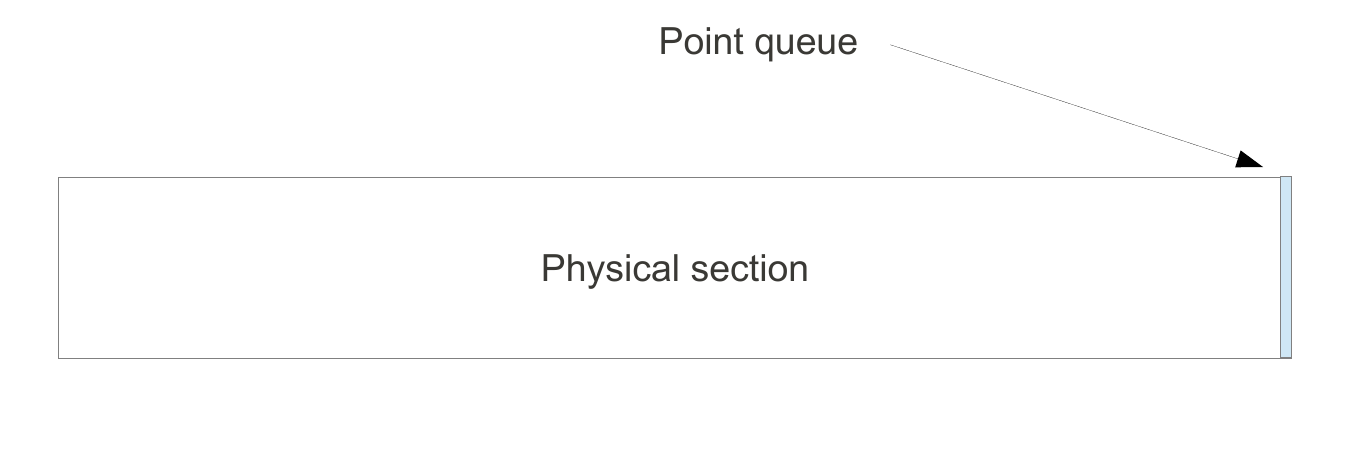}{Division of a link in point queue models.}{0.8\textwidth}

Point queues can be thought of in several ways.
One can imagine a wide link which necks down at its downstream end --- the physical section reflects the wide portion of the link, and the point queue represents the vehicles which are delayed as the capacity is reduced at the downstream bottleneck.
One can also imagine a traffic signal at the downstream end, and magical technology (flying cars?) which allows vehicles to ``stack'' vertically at the signal --- there can be no congestion upstream of this ``stack,'' since vehicles can always fly to the top.
(Figure~\ref{fig:stackedvehiclepq}) One may even imagine that there is no physical meaning to either of these, and that the physical section and point queue merely represent the delays incurred from traveling the link at free-flow, and the additional travel time due to congestion.

\stevefig{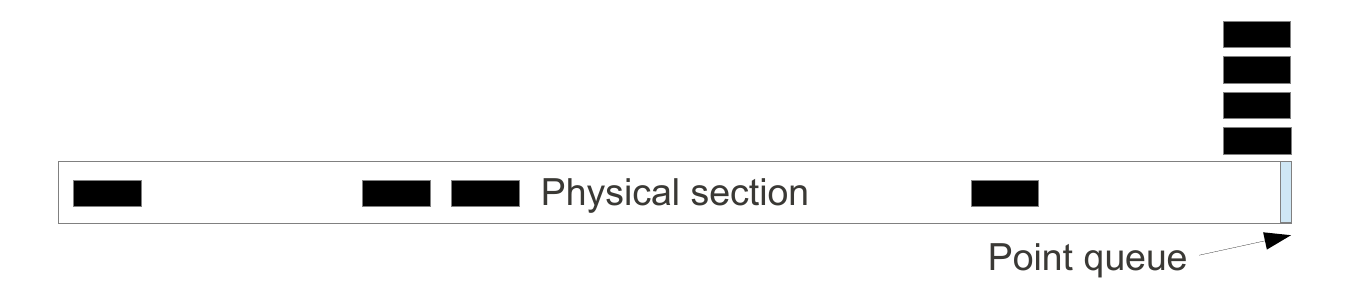}{Envisioning a point queue as vehicles stacking vertically.}{0.75\textwidth}

The point queue discharges vehicles at a maximum rate of $q_{max}^\downarrow$\label{not:qarrow} (measured in vehicles per unit time), called the \emph{capacity}.\index{link!capacity}
The capacity imposes an upper limit on the sending flow, so we always have
\labeleqn{pqsendingflowcap}{S(t) \leq q^\downarrow_{max} \Delta t\,.}
However, if the queue is empty, or if only a few vehicles are in the queue, the discharge rate may be less than this.
Once the queue empties, the only vehicles which can exit the link are ones reaching the downstream end from the uncongested physical section.
Since we assume that all vehicles in this section travel at the free-flow speed\index{speed!free-flow speed} (which we will denote $u_f$), this means that only the vehicles that are closer than $u_f \Delta t$ to the downstream end can possibly leave.

We can use the cumulative counts\index{cumulative count} $N^\uparrow$ and $N^\downarrow$ to count the number of vehicles which are close enough to the downstream end to exit in the next time step.
Since the entire physical section is traversed at the free-flow speed $u_f$,\label{not:uf} a vehicle whose distance from the downstream end of the link is exactly $u_f \Delta t$ distance units must have passed the upstream end of the link exactly $(L - u_f \Delta t) / u_f$ time units ago.
We call this a ``threshold'' vehicle, since any vehicle entering the link after this one has not yet traveled far enough, while any vehicle entering the link before this one is close enough to the downstream end to exit.
The number of vehicles between the threshold vehicle and the downstream end of the link can thus be given by
\labeleqn{thresholdnum}{N^\uparrow\myp{t - \frac{L - u_f \Delta t }{u_f}} - N^\downarrow(t) = N^\uparrow\myp{t + \Delta t - \frac{L}{u_f}} - N^\downarrow(t)\,.}

The sending flow is the smaller of the number of vehicles which are close enough to the downstream end to exit, given by equation~\eqn{thresholdnum}, and the capacity of the queue.
Thus\index{link model!point queue!sending flow}
\labeleqn{pqsendingflow}{S(t) = \min\{ N^\uparrow(t + \Delta t - L/u_f) - N^\downarrow(t), q^\downarrow_{max} \Delta t\}\,.}

The receiving flow for the point queue model is easy to calculate.
Since the physical section is uncongestible, the link capacity is the only limitation on the rate at which vehicles can enter.
The capacity of the upstream end of the link may be different than the capacity of the downstream end of the link (perhaps due to a lane drop, or a stop sign at the end of the link), so we denote the capacity of the upstream end by $q^\uparrow_{max}$.
The receiving flow is given by\index{link model!point queue!receiving flow}
\labeleqn{pqreceivingflow}{R(t) = q^\uparrow_{max} \Delta t\,.}
In real traffic networks, queues occupy physical space and cannot be confined to a single point.
The spatial queue model in the next subsection shows one way to reflect this.

Table~\ref{tbl:pqexample} shows how the point queue model operates, depicting the state of a link over ten time steps.
The $N^\uparrow$ and $N^\downarrow$ columns express the number of vehicles which have entered and left the link at each time step, as well as the sending and receiving flows during each timestep.
The difference between $N^\uparrow$ and $N^\downarrow$ represents the number of vehicles on the link at any point in time.
In this example, we assume that the unit of time is chosen so that $\Delta t = 1$, the free-flow speed is $u_f = L / 3$ (so a vehicle takes 3 time steps to traverse the link under free-flow conditions), the upstream capacity is $q^\uparrow_{max} = 10$, and the downstream capacity is $q^\downarrow_{max} = 5$.
Initially, the sending flow is zero, because no vehicles have reached the downstream end of the link.
The sending flow then increases as flow exits, but eventually reaches the downstream capacity.
At this point, a queue forms and vehicles exit at the downstream capacity rate.
Eventually, the queue clears, the link is empty, and the sending flow returns to zero.
The receiving flow never changes from the upstream capacity, even when a queue is present.
In this example, notice that $N^\downarrow(t+1) = N^\downarrow(t) + S(t)$.
This happens because we are temporarily ignoring what might be happening from downstream.
Depending on downstream congestion, $N^\downarrow(t+1)$ could be less than $N^\downarrow(t) + S(t)$; but it could never be greater, because the sending flow is always a limit on the number of vehicles that can exit.
Also notice that for all time steps, $N^\uparrow(t + 1) \leq N^\uparrow(t) + R(t)$, because the receiving flow is a limit on the number of vehicles that can enter the link.

\begin{table}
\begin{center}
\caption{Point queue example, with $\Delta t = 1$, $u_f = L / 3$, $q^\uparrow_{max} = 10$, and $q^\downarrow_{max} = 5$.
\label{tbl:pqexample}}
\begin{tabular}{c|cc|cc}
$t$   & $N^\uparrow$ & $N^\downarrow$  & $R$ & $S$ \\
\hline
0     & 0            & 0               & 10  & 0   \\
1     & 1            & 0               & 10  & 0   \\
2     & 5            & 0               & 10  & 0   \\
3     & 10           & 0               & 10  & 1   \\
4     & 17           & 1               & 10  & 4   \\
5     & 27           & 5               & 10  & 5   \\
6     & 30           & 10              & 10  & 5   \\
7     & 30           & 15              & 10  & 5   \\
8     & 30           & 20              & 10  & 5   \\
9     & 30           & 25              & 10  & 5   \\
10    & 30           & 30              & 10  & 0   \\
\end{tabular}
\end{center}
\end{table}

In practice, the upstream and downstream capacities are often assumed the same, in which case we just use the notation $q_{max}$ to refer to capacities at both ends.
The network features which would make the capacities different upstream and downstream (such as stop signs or signals) are usually better represented with node models, discussed in the next section.
\index{link model!point queue|)}

\subsection{Spatial queue}
\label{sec:spatialqueue}

\index{link model!spatial queue|(}
The spatial queue model is similar to the point queue model, except that a maximum queue length is now enforced.
When the queue reaches this maximum length, no further vehicles are allowed to enter the link.
That is, the queue now occupies physical space, and because the link is finite in length, the link can become completely blocked.
This will result in \emph{queue spillback}, as vehicles on upstream links will be unable to enter the blocked link.
Like the point queue model, we assume that the queue is always at the downstream end of the link: there is at most one uncongested physical section at the upstream end of the link, and at most one stopped queue at the downstream end, in that order.
This is still a simplification of real traffic (where links can have multiple congested and uncongested sections), but by allowing the length of the physical section to shrink as the queue grows, one can model the queue spillback phenomenon which is common in congested networks.

The sending flow for the spatial queue model is calculated in exactly the same way as for the point queue model: the smaller of the number of vehicles close enough to the downstream end of the link to exit in the next time step, and the capacity\index{link!capacity} of the link:\index{link model!spatial queue!sending flow}
\labeleqn{sqsendingflow}{S(t) = \min\{ N^\uparrow(t + \Delta t - L/u_f) - N^\downarrow(t), q^\downarrow_{max} \Delta t\}\,.}

The receiving flow includes an additional term to reflect the finite space on the link for the queue, alongside the link capacity.
Since vehicles in the queue are stopped, the space they occupy is given by the jam density\index{density!jam} $k_j$,\label{not:kj} expressed in vehicles per unit length.
The maximum number of vehicles the link can hold is $k_j L$, while the number of vehicles currently on the link is $N^\uparrow(t) - N^\downarrow(t)$.
The receiving flow cannot exceed the difference between these:\index{link model!spatial queue!receiving flow}
\labeleqn{spreceivingflow}{R(t) = \min\{ k_j L - (N^\uparrow(t) - N^\downarrow(t)), q^\uparrow_{max} \Delta t\}\,.}

By assuming that the queue is always at the downstream end of the link, the spatial queue model essentially assumes that all vehicles in a queue move together.
In reality, there is some delay between when the head of the queue starts moving, and when the vehicle at the tail of the queue starts moving --- when a traffic light turns green, vehicles start moving one at a time, with a slight delay between when a vehicle starts moving and when the vehicle behind it starts moving.
These delays cannot be captured in a spatial queue model.
To represent this behavior, we will need a better understanding of traffic flow theory.
Section~\ref{sec:trafficflowtheory} will present this information, and we will ultimately build more realistic link models.

Table~\ref{tbl:sqexample} shows how the spatial queue model operates.
We are representing the same scenario in the point queue example (Table~\ref{tbl:pqexample}), where 30 vehicles enter the link.
The difference between the $N^\uparrow$ values in that earlier table shows the demand for vehicles to enter the link (1, 4, 5, 7, 10, and 3 vehicles in the initial time intervals), but we will now limit the number of vehicles which can fit on the link by introducing a jam density.
In cases where there is more demand to enter than receiving flow on the link, the unserved demand will remain in a queue and will enter the link when space becomes available.

In particular, assume that the maximum number of vehicles which can fit on the link is $k_j L = 20$.
As a result, it takes longer for the 30 vehicles to enter the link, and in time intervals 4, 5, and 6 there is a queue waiting to enter the upstream end of the link.
The queue is still able to completely discharge by the end of the ten time steps.
As before, the difference between $N^\downarrow(t)$ and $N^\downarrow(t+1)$ is never more than $S(t)$ (in this example, exactly equal because we are ignoring downstream conditions), and the difference between $N^\uparrow(t)$ and $N^\uparrow(t+1)$ is never more than $R(t)$.
\index{link model!spatial queue|)}
\index{link model|)}

\begin{table}
\begin{center}
\caption{Spatial queue example, with $\Delta t = 1$, $u_f = L / 3$, $q^\uparrow_{max} = 10$, $q^\downarrow_{max} = 5$, and $k_j L = 20$.
\label{tbl:sqexample}}
\begin{tabular}{c|cc|cc}
$t$   & $N^\uparrow$ & $N^\downarrow$  & $R$ & $S$ \\
\hline
0     & 0            & 0               & 10  & 0   \\
1     & 1            & 0               & 10  & 0   \\
2     & 5            & 0               & 10  & 0   \\
3     & 10           & 0               & 10  & 1   \\
4     & 17           & 1               & 4   & 4   \\
5     & 21           & 5               & 4   & 5   \\
6     & 25           & 10              & 5   & 5   \\
7     & 30           & 15              & 5   & 5   \\
8     & 30           & 20              & 10  & 5   \\
9     & 30           & 25              & 10  & 5   \\
10    & 30           & 30              & 10  & 0   \\
\end{tabular}
\end{center}
\end{table}

\section{Node Model Concepts}
\label{sec:nodemodels}

\index{node model|(}
Node models complement the link models discussed in the previous section, by representing how flows between different links interact with each other.
Each link computes its own sending and receiving flow, and node models combine this information to determine the actual rate of flow from one link to the next.
Separating these computations into link and node models both simplifies the process of network loading, and makes it more computationally efficient, because calculations of sending and receiving flows can be done independently, in parallel, for each link.
It is even possible to use different link models on different links (say, one link with a spatial queue model and another with a point queue).

\stevefig{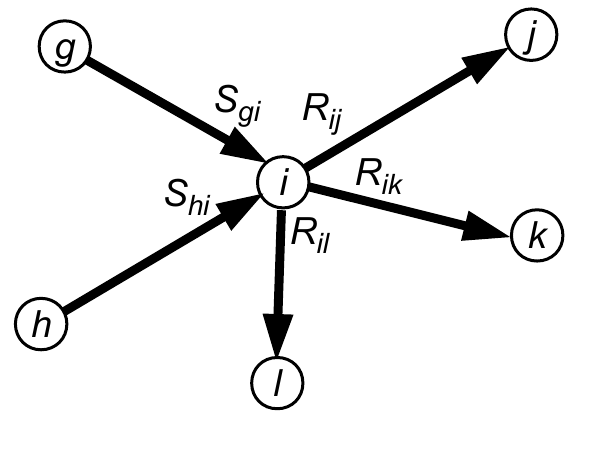}{Sending and receiving flows interact at nodes to produce turning movement flows.}{0.5\textwidth}

All that a node model needs is the sending flow of each of the links which enter this node, and the receiving flow for each link leaving this node.
(Figure~\ref{fig:genjunction})  That is, when processing node $i$ at time $t$, we assume that $S_{hi}(t)$ and $R_{ij}(t)$ have respectively been calculated for each incoming link $(h,i) \in \Gamma^{-1}(i)$ and each outgoing link $(i,j) \in \Gamma(i)$.
The task is to determine how many vehicles move from each incoming link to each outgoing link during the $t$-th time interval; denote this value by $y_{hij} (t)$.\label{not:yhij}
This is called the \emph{turning movement flow}\index{turning movement flow} from $(h,i)$ to $(i,j)$, or more compactly, the turning movement $[h,i,j]$.\label{not:hij}
Let $\Xi(i)$\label{not:Xii} denote the set of allowable turning movements at node $i$; this provides a natural way to model turn prohibitions, U-turn prohibitions, and so forth.

As an example of this notation, consider node $i$ in Figure~\ref{fig:genjunction}.
If all of the turning movements are allowed, then $\Xi(i)$ has six elements: $[g,i,j]$, $[g,i,k]$, $[g,i,l]$, $[h,i,j]$, $[h,i,k]$, and $[h,i,l]$.
This set might not include all six elements --- for instance, if the left turn from approach $(g,i)$ to $(i,j)$ is prohibited, then $\Xi(i)$ would exclude $[g,i,j]$.

Many different node models can be used.
This section discusses three simple node models --- links in series, diverges, and merges --- to illustrate the general concepts.
Later in this chapter, we discuss node models that can be used for signalized intersections, all-way stops, and other types of intersections.
There are many variations of all of these, but they share some common principles, which are discussed here.

Any node model must satisfy a number of constraints and principles, for all time intervals $t$:
\begin{enumerate}
\item Vehicles will not voluntarily hold themselves back.
That is, it should be impossible to increase any turning movement flow $y_{hij}(t)$ without violating one of the constraints listed below, or one of the other constraints imposed by a specific node model.
\item Each turning movement flow must be nonnegative, that is, $y_{hij}(t) \geq 0$ for each $[h,i,j] \in \Xi(i)$.
\item Any turning movement which is not allowable has zero flow, that is, $y_{hij}(t) = 0$ whenever $[h,i,j] \notin \Xi(i)$.
\item For each incoming link, the sum of the turning movement flows out of this link cannot exceed the sending flow, since by definition those are the only vehicles which could possibly leave that link in the next time step:
\labeleqn{sendingflowrestriction}{\sum_{(i,j) \in \Gamma(i)} y_{hij}(t) \leq S_{hi}(t) \qquad \forall (h,i) \in \Gamma^{-1}(i) \,.}
The sum on the left may be less than $S_{hi}(t)$, because it is possible that some vehicles cannot leave $(h,i)$ due to obstructions from a downstream link or from the node itself (such as a red signal).
\item For each outgoing link, the sum of the turning movement flows into this link cannot exceed the receiving flow, since that is the maximum number of vehicles that link can accommodate:
\labeleqn{receivingflowrestriction}{\sum_{(h,i) \in \Gamma^{-1}(i)} y_{hij}(t) \leq R_{ij}(t) \qquad \forall (i,j) \in \Gamma(i) \,.}
The sum on the left may be less than $R_{ij}(t)$, because there may not be enough vehicles from upstream links to fill all of the available space in the link.
(Recall from Chapter~\ref{chp:networkrepresentations} that $\Gamma(i)$ is the set of nodes immediately downstream of node $i$, and $\Gamma^{-1}(i)$ the set of nodes immediately upstream.)
\item Route choices must be respected.
That is, the values $y_{hij}(t)$ must be compatible with the directions travelers wish to go based on their chosen paths; we cannot reassign them on the fly in order to increase a $y_{hij}$ value.
\item The \emph{first-in, first-out (FIFO) principle\index{first-in, first-out (FIFO)}} must be respected.
This is closely related to the previous property.
We cannot allow vehicles to ``jump ahead'' in queue to increase a $y_{hij}$ value, unless there is a separate turn lane or other roadway geometry which can separate vehicles on different routes.
(Figure~\ref{fig:fifoqueuespillback}).
\item \index{invariance principle|(}The \emph{invariance principle} must be respected.
If the outflow from a link is less than its sending flow ($\sum_h y_{hij}(t) < S_{hi}(t)$), then recalculating $y_{hij}$ using a \emph{larger} value of $S_{hi}$ should not change the result.
Likewise, if the inflow to a link is less than its receiving flow ($\sum_j y_{hij}(t) < R_{ij}(t)$), then recalculating $y_{hij}$ with a \emph{larger} value of $R_{ij}$ should not change the result.
In other words, if the sending (or receiving) flow is not ``binding,'' then its specific value cannot matter for the actual flows.
\end{enumerate}

\stevefig{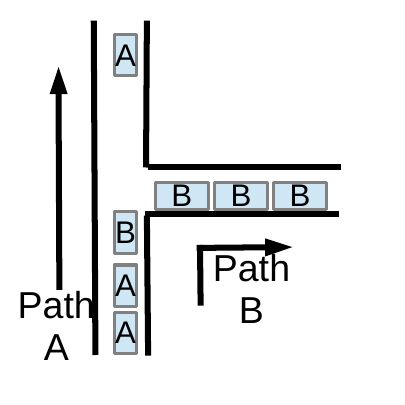}{The first-in, first-out principle implies that vehicles waiting to turn will obstruct vehicles further upstream.}{0.6\textwidth}

The invariance principle warrants additional explanation.
Consider the first scenario, where $y_{hij} < S_{hi}$, and consider what might happen if $S_{hi}$ were to increase.
Note that this requires that $S_{hi}$ is currently less than the capacity of link $(h,i)$ (or else we can't increase it).
Since $y_{hij} < S_{hi}$, not all of the vehicles that want to turn from $(h,i)$ onto $(i,j)$ are able to do so, perhaps because of a restriction due to the receiving flow of $R_{ij}$.
Since $S_{hi} < q^{hi}_{max}$, the sending flow is determined by the number of vehicles close to the downstream end of $(h,i)$, not the capacity of the link.
You can envision this scenario as a link $(h,i)$ flowing below capacity, connected to a link $(i,j)$ with a long queue that limits the number of vehicles that enter the link.

In the current time interval, the sending flow is less than the link capacity.
However, since $y_{hij}(t) < S_{hi}(t)$, not all of the vehicles in this sending flow are able to leave $(h,i)$.
They will remain in the next time interval, and form a queue.
When there is a queue on the link, $S_{hi}(t)$ will equal the full capacity of the link, because there are vehicles already waiting to leave --- we have moved from the situation in the right panel of Figure~\ref{fig:sendingflowexample} to the one on the left.
If the flow $y_{hij}(t+\Delta t)$ in the next time interval were to increase in response, then the queue might be entirely cleared in the next time interval, moving us back to the situation in the right panel of Figure~\ref{fig:sendingflowexample}.
This would reduce the sending flow, which might reduce $y_{hij}$, which would result in the queue forming again, and so on, with the flows oscillating between time intervals.
This is unrealistic, and is an artifact introduced by choosing a particular discretization, not a traffic phenomenon one would expect in the field.

The other case is similar, but involving cases where $y_{hij} < R_{ij}$ and $R_{ij} < q^{max}_{ij}$.
This kind of oscillating behavior between time steps does not reflect traffic physics, and should not be present in a model.
The node models presented in this section satisfy all of these principles.
The exercises present a node model which does not satisfy the invariance principle, and asks you to compare it to the node models presented in the main text.
\index{invariance principle|)}

\subsection{Links in series}
\label{sec:linksinseries}

\index{node model!links in series|(}
The simplest node to model is one with exactly one incoming link, and exactly one outgoing link.
(Figure~\ref{fig:serieslink}).
This may seem like a trivial node, since there is no real ``intersection'' here.
However, they are often introduced to reflect changes within a link.
For instance, if a freeway reduces from three lanes in a direction to two lanes, this reduction in capacity and jam density can be modeled by introducing a node at the point where the lane drops.
In this way, each link can have a homogeneous capacity and jam density, and we can simplify the notation --- in the link models above, we distinguished between $q^\uparrow_{max}$ and $q^\downarrow_{max}$ at the two ends of the link.
Now we can just use $q_{max}$\label{not:qmax} for the entire link, including both ends.

\stevefig{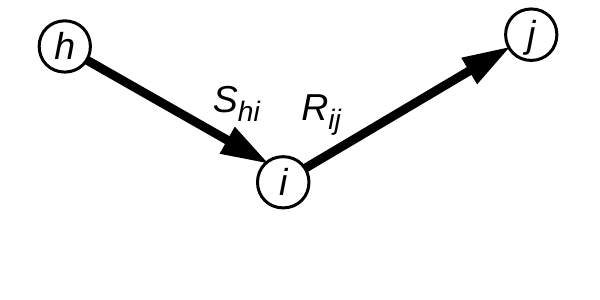}{Two links in series.}{0.5\textwidth}

For concreteness, let the incoming link be $(h,i)$ and the outgoing link be $(i,j)$.
In this case, there is only one turning movement, so $\Xi(i) = \{ [h,i,j] \}$, and the only turning movement flow we need concern ourselves with is $y_{hij}(t)$.
In the case of two links in series, the formula is simple:\index{node model!links in series!turning movement flow}
\labeleqn{seriesflow}{y_{hij}(t) = \min \{ S_{hi}(t), R_{ij}(t) \} \,,}
that is, the number of vehicles moving from link $(h,i)$ to $(i,j)$ during the $t$-th time interval is the lesser of the sending flow from the upstream link in that time interval, and the receiving flow of the downstream link.
For instance, if the upstream link sending flow is 10 vehicles, while the downstream link receiving flow is 5 vehicles, a total of 5 vehicles will successfully move from the upstream link to the downstream one, because that is all there is space for.
If the upstream link sending flow is 3 vehicles and the downstream link receiving flow is 10 vehicles, 3 vehicles will move from the upstream link to the downstream one, because that is all the vehicles that are available.

We can check that this node model satisfies all of the desiderata from Section~\ref{sec:nodemodels}.
Going through each of these conditions in turn:
\begin{enumerate}
\item The flow $y_{hij}(t)$ is chosen to be the minimum of the upstream sending flow, and the downstream receiving flow; any value larger than this would violate either the sending flow constraint or the receiving flow constraint.
\item The sending and receiving flows should both be nonnegative, regardless of the link model, so $y_{hij}(t)$ is as well.
\item There is only one turning movement, so this constraint is trivially satisfied.
\item The formula for $y_{hij}(t)$ ensures it cannot be greater than the upstream sending flow.
(This condition simplifies since there is only one incoming and outgoing link, so the summation and ``for all'' quantifier can be disregarded.)
\item The formula for $y_{hij}(t)$ ensures it cannot be greater than the downstream receiving flow.
(This condition simplifies in the same way.)
\item Route choice is irrelevant when two links meet in series, since all incoming vehicles must exit by the same link.
\item FIFO is also irrelevant, since all vehicles entering the node behave in the same way.
(This would not be the case if there was more than one exiting link, and vehicles were on different paths.)  So we don't have to worry about FIFO when calculating the turning movement flow.
\item To see that the formula satisfies the invariance principle, we have to check two conditions.
If the outflow from $(h,i)$ is less than the sending flow, this means that $y_{hij}(t) = R_{ij}(t)$, and the second term in the minimum of equation~\eqn{seriesflow} is binding.
Increasing the sending flow (the first term in the minimum) further would not affect its value.
Similarly, if the inflow to $(i,j)$ is less than its receiving flow, this means that $y_{hij}(t) = S_{hi}(t)$, and the first term in the minimum is binding.
Increasing the receiving flow (the second term) would not affect its value either.
\end{enumerate}
\index{node model!links in series|)}

\subsection{Merges}
\label{sec:merge}

\index{node model!merge|(}
A merge node has only one outgoing link $(i,j)$, but more than one incoming link, here labeled $(g,i)$ and $(h,i)$, as in Figure~\ref{fig:merge}.
This section only concerns itself with the case of only two upstream links, and generalizing to the case of additional upstream links is left as an exercise.
Here $\Xi(i) = \{ [g,i,j], [h,i,j] \}$, and we want to calculate the rate of flow from the upstream links to the downstream one, that is, the flow rates $y_{gij}(t)$ and $y_{hij}(t)$.
As you might expect, the main quantities of interest are the upstream sending flows $S_{gi}(t)$ and $S_{hi}(t)$, and the downstream receiving flow $R_{ij}(t)$.
We assume that these values have already been computed by applying a link model.

For brevity, we will omit the time index in the rest of this section --- it is implicit that all calculations are done with the sending and receiving flows at the current time step.

\stevefig{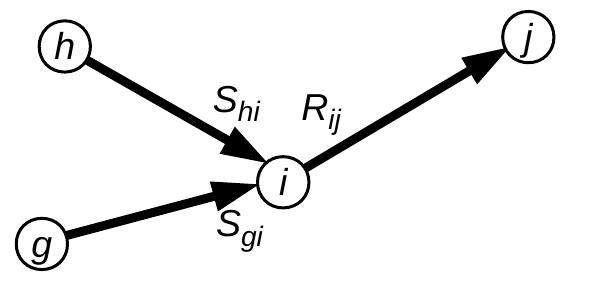}{Prototype merge node.}{0.5\textwidth}

There are three possibilities, one corresponding to free flow conditions at the merge, one corresponding to congestion with queues growing on both upstream links, and one corresponding to congestion on only one upstream link.
For the merge to be freely flowing, both upstream links must be able to transmit all of the flow which seeks to leave them, and the downstream link must be able to accommodate all of this flow.
Mathematically, we need $S_{gi} + S_{hi} \leq R_{ij}$, and if this is true then we simply set $y_{gij} = S_{gi}$, and $y_{hij} = S_{hi}$.

In the second case, there is congestion (so $S_{gi} + S_{hi} > R_{ij}$), and furthermore, flow is arriving fast enough on both upstream links for a queue to form at each of them.
Empirically, in such cases the flow rate from the upstream links is approximately proportional to the capacity on these links, that is,
\labeleqn{proportionality}{\frac{y_{gij}}{y_{hij}} = \frac{q^{gi}_{max}}{q^{hi}_{max}}}
A little thought should convince you that this relationship is plausible.
Furthermore, in the congested case, all of the available downstream capacity will be used, so
\labeleqn{flowconservation}{y_{gij} + y_{hij} = R_{ij}}
Substituting~\eqn{proportionality} into~\eqn{flowconservation} and solving, we obtain
\labeleqn{mergecase2}{y_{gij} = \frac{q_{max}^{gi}}{q_{max}^{gi} + q_{max}^{hi}}R_{ij}}
with a symmetric expression for $y_{hij}$.

The third case is perhaps a bit unusual.
The merge is congested ($S_{gi} + S_{hi} > R_{ij}$), but a queue is only forming on one of the upstream links.
This may happen if the flow on one of the upstream links is much less than the flow on the other.
In this case, the proportionality rule allows all of the sending flow from one link to enter the downstream link, with room to spare.
This ``spare capacity'' can then be consumed by the other approach.
If link $(g,i)$ is the link which cannot send enough flow to meet the proportionality condition, so that
\labeleqn{mergecase3condition}{S_{gi} < \frac{q_{max}^{gi}}{q_{max}^{gi} + q_{max}^{hi}}R_{ij}\,,}
then the two flow rates are $y_{gij} = S_{gi}$ and $y_{hij} = R_{ij} - S_{gi}$: one link sends all of the flow it can, and the other link consumes the remaining capacity.
The formulas are reversed if it is link $(h,i)$ that cannot send enough flow to meet its proportionality condition.

Exercise~\ref{ex:mergemedianformula} asks you to show that the second and third cases can be handled by the single equation
\labeleqn{mergemedianformula}{y_{gij} = \operatorname{med} \myc{S_{gi}, R_{ij} - S_{hi}, \frac{q_{max}^{gi}}{q_{max}^{gi} + q_{max}^{hi}}R_{ij}}\,,}
where $\operatorname{med}$ refers to the median of a set of numbers.
This formula applies whenever $S_{gi} + S_{hi} > R_{ij}$.
An analogous formula holds for the other approach by swapping the $g$ and $h$ indices.
Exercise~\ref{ex:mergedesiderata} asks you to verify the desiderata of Section~\ref{sec:nodemodels} are satisfied by this equation.

For certain merges, it may not be appropriate to assign flow proportional to the capacity of the incoming links.
Rules of the road, signage, or signalization might allocate the capacity of the downstream link differently.
In general, the share of the downstream receiving flow that is allocated to approaches $(g,i)$ and $(h,i)$ can be written as $\beta_{gi}$\label{not:betaconflict} and $\beta_{hi}$, respectively, with $\beta_{gi} + \beta_{hi} = 1$ and both of them nonnegative.
A more general form of the merge equation can then be written as\label{not:med}
\labeleqn{genericmergemedianformula}{y_{gij} = \operatorname{med} \myc{S_{gi}, R_{ij} - S_{hi}, \beta_{gi} R_{ij}}\,,}
for the second and third cases, with a similar formula for $(h,i)$.
(The first case is unchanged, because the shares $\beta$ are only relevant when there is insufficient receiving flow for the arriving vehicles.)

Finally, rather than listing the merge model into cases explicitly, we can give an iterative algorithm which will internally determine which case we are in, and the values of $y_{gij}$ and $y_{hij}$.
You may wonder why we need this algorithm, since the merge model is not complicated and there are only a few cases.
The advantages of the algorithm are twofold: it is easier to extend the algorithm to the case of more than 2 incoming links, than to write out all the cases explicitly; and the ideas in this algorithm will be used when we study more complicated nodes in Section~\ref{sec:fancynodes}.
This algorithm assumes $\beta_{gi}$ and $\beta_{hi}$ are strictly positive; Exercise~\ref{ex:zeromerge} asks you to generalize this to the case where one of them is zero.
\begin{enumerate}
\item \index{node model!merge!algorithm}Initialize the set of active turning movements $\Omega \leftarrow \{ [g,i,j], [h,i,j] \}$;\label{not:Omega} the remaining sending flows\label{not:Stwiddle} $\tilde{S}_{gi} \leftarrow S_{gi}$ and $\tilde{S}_{hi} \leftarrow S_{hi}$; the remaining receiving flow\label{not:Rtwiddle} $\tilde{R}_{ij} \leftarrow R_{ij}$; and the transition flows $y_{gij} \leftarrow 0$ and $y_{hij} \leftarrow 0$.
\item Set $\alpha_{gij} \leftarrow \beta_{gi}$ if $[g,i,j] \in A$\label{not:alphahij} and 0 if not; likewise set $\alpha_{hij} \leftarrow \beta_{hi}$ if $[h,i,j] \in A$ and 0 if not; and set $\alpha_{ij} \leftarrow \alpha_{gij} + \alpha_{hij}$.
\item Compute $$\theta = \min \myc{\frac{\tilde{S}_{gi}}{\alpha_{gij}}, \frac{\tilde{S}_{hi}}{\alpha_{hij}}, \frac{\tilde{R}_{ij}}{\alpha_{ij}}}\,,$$\label{not:thetastep} treating division by zero as $+\infty$.
Whichever quantity is limiting in this minimum will determine which is exhausted first: the sending flow from $(g,i)$, the sending flow from $(h,i)$, or the receiving flow from $(i,j)$.
\item Increase flows for active turning movements, and update unallocated sending and receiving flows: if $[g,i,j] \in A$, update $y_{gij} \leftarrow y_{gij} + \theta \alpha_{gij}$, $\tilde{S}_{gi} \leftarrow \tilde{S}_{gi} - \theta \alpha_{gij}$, and $\tilde{R}_{ij} \leftarrow \tilde{R}_{ij} - \theta \alpha_{gij}$.
If $[h,i,j] \in A$, update $y_{hij} \leftarrow y_{hij} + \theta \alpha_{hij}$; $\tilde{S}_{hi} \leftarrow \tilde{S}_{hi} - \theta \alpha_{hij}$, and $\tilde{R}_{ij} \leftarrow \tilde{R}_{ij} - \theta \alpha_{hij}$.
\item If $\tilde{S}_{gi} = 0$, remove $[g,i,j]$ from $A$.
If $\tilde{S}_{hi} = 0$, remove $[h,i,j]$ from $A$.
If $\tilde{R}_{ij} = 0$, set $\Omega \leftarrow \emptyset$.
\item Terminate if $\Omega = \leftarrow \emptyset$; otherwise return to step 2.
\end{enumerate}
The idea behind the algorithm is to initially assume that both approaches will use the merge in proportion to their $\beta$ values, but if one of them exhausts their sending flow first, the other approach is able to send additional flow until either its sending flow is exhausted, or the downstream link's receiving flow is exhausted.
These correspond to the three cases above.
It's worth coming up with a few scenarios for sending and receiving flow, and confirming that the algorithm gives the same answers as the explicit formulas earlier in the section.
\index{node model!merge|)}

\subsection{Diverges}
\label{sec:diverge}

\index{node model!diverge|(}
A diverge node is one with only one incoming link $(h,i)$, but more than one outgoing link, as in Figure~\ref{fig:diverge}.
This section concerns itself with the case of only two downstream links.
The exercises ask you to generalize to the case of three downstream links, using the same concepts.
Let these two links be called $(i,j)$ and $(i,k)$, so $\Xi(i) = \{ [h,i,j], [h,i,k] \}$.
Our interest is calculating the rate of flow from the upstream link to the downstream ones, that is, the flow rates $y_{hij}$ and $y_{hik}$.
We assume that the sending flow $S_{hi}$ and the receiving flows $R_{ij}$ and $R_{ik}$ have already been calculated.
Unlike links in series or merges, we also need to represent some model of route choice, since some drivers may choose link $(i,j)$, and others link $(i,k)$.
Let $p_{hij}$\label{not:phij} and $p_{hik}$ be the ``splitting'' proportions\index{splitting proportion} of drivers choosing these two turning movements during the $t$-th time interval, respectively.
Naturally, $p_{hij}$ and $p_{hik}$ are nonnegative, and $p_{hij} + p_{hik} = 1$.
Like the sending and receiving flows, these values can change with time, but to avoid cluttering formulas we will leave the time indices off of $p$ values unless it is unclear which time step we are referring to.

\stevefig{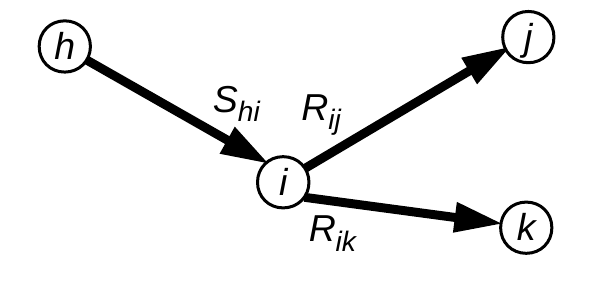}{Prototype diverge node.}{0.5\textwidth}

There are two possibilities, one corresponding to free flow conditions at the diverge, and the other corresponding to congestion and not all vehicles in the sending flow being able to leave the upstream link.
What does ``free flow'' mean?  For the diverge to be freely flowing, \emph{both} of the downstream links must be able to accommodate the flow which seeks to enter them.
The rates at which vehicles want to enter the two links are $p_{hij} S_{hi}$ and $p_{hik} S_{hi}$.
We will call these quantities the \emph{oriented sending flows}\index{sending flow!oriented}\index{oriented sending flow|see {sending flow, oriented}}
\labeleqn{orientedpreview}{S_{hij} = S_{hi} p_{hij} \qquad S_{hik} = S_{hi} p_{hik}}
to reflect the portions of the sending flow intended for the turning movements $[h,i,j]$ and $[h,i,k]$.
If both downstream links can accommodate the oriented sending flows, we need $S_{hij} \leq R_{ij}$ and $S_{hik} \leq R_{ik}$.
In this case we simply have $y_{hij} = S_{hij} = p_{hij} S_{hi}$ and $y_{hik} = S_{hik} = p_{hik} S_{hi}$: all of the flow which wants to leave the diverge can.

The case of congestion is slightly more interesting, and requires making assumptions about how drivers will behave.
One common assumption is that flow waiting to enter one link at a diverge will obstruct every other vehicle on the link (regardless of which link it is destined for).
This most obviously represents the case where the upstream link has only a single lane, so any vehicle which has to wait will block any vehicle behind it; but this model is commonly used even in other cases.\footnote{For instance, this can represent drivers attempting to ``queue jump'' by cutting into the turn lane at the last moment.
Or, if the turn lane is long, it may be appropriate to treat the diverge at the point where the turn lane begins, as opposed to at the physical diverge.}  When there is congestion, only some fraction $\phi$\label{not:phidiverge} of the upstream sending flow can move.
The assumption that any vehicle waiting blocks every vehicle upstream implies that this same fraction applies to both of the downstream links, so $y_{hij} = \phi S_{hij}$ and $y_{hik} = \phi S_{hik}$.

So, how to calculate $\phi$?  The inflow rate to a link cannot exceeds its receiving flow, so $y_{hij} = \phi S_{hij} \leq R_{ij}$ and $y_{hik} = \phi  S_{hik} \leq R_{ik}$, or equivalently $\phi \leq R_{ij} / S_{hij}$ and $\phi \leq R_{ik} / S_{hik}$.
Every vehicle which can move will, so
\labeleqn{basicdiverge}{\phi = \min \left\{ \frac{R_{ij}}{S_{hij}}, \frac{R_{ik}}{S_{hik}} \right\}}
Furthermore, we can introduce the uncongested case into this equation as well, and state
\labeleqn{generaldiverge}{\phi = \min \left\{ \frac{R_{ij}}{S_{hij}}, \frac{R_{ik}}{S_{hik}}, 1 \right\}}
regardless of whether there is congestion at the diverge or not.
Why?  If the diverge is at free flow, then $\phi = 1$, but $R_{ij}/ S_{hij} \geq 1$ and $R_{ik} / S_{hik} \geq 1$.
Introducing 1 into the minimum therefore gives the correct answer for free flow.
Furthermore, if the diverge is not at free flow, then either $R_{ij} / S_{hij} < 1$ or $R_{ik} / S_{hik} < 1$, so adding 1 does not affect the minimum value.
Therefore, this formula is still correct even in the congested case.
Exercise~\ref{ex:divergedesiderata} asks you to verify the desiderata of Section~\ref{sec:nodemodels} are satisfied by this equation.

We emphasize that for the purposes of network loading, the splitting proportions $p_{hij}(t)$ and $p_{hik}(t)$ that determine the oriented sending flows at each time step are taken as inputs and assume fixed.
The performance of the diverge node depends heavily on these values, and in the larger dynamic traffic assignment process the $p$ values will change from iteration to iteration to reflect changes in the driver route choice.
Within a single network loading, however, we treat them as given constants.
\index{node model!diverge|)}
\index{node model|)}

\section{Combining Node and Link Models}
\label{sec:combiningnodelink}

\index{link model!combining with node model|(}
\index{node model!combining with link model|(}
This section describes how node and link models are combined, to complete the network loading process.
We must also describe what happens at origins and destination (zone) nodes.\index{link!centroid connector}
For the algorithm in this section, we assume that the only links connected to zones are special links called \emph{centroid connectors}, which do not represent a specific roadway so much as a collection of small local streets used by travelers entering or leaving a specific neighborhood.
It is common to give centroid connectors between origins and ordinary nodes a very large (even infinite) jam density, and centroid connectors between ordinary nodes and destinations a very large (even infinite) capacity.
In both cases the free flow time should be small.
These considerations reflect the ideas that centroid connectors should not experience significant congestion (or else they should be modeled as proper links in the network), and simply convey flow from origin nodes and to destination nodes with as little interference as possible.
It is common to forbid travelers from using centroid connectors except to start and end their trips, that is, to exclude the use of centroid connectors as ``shortcuts.''  This can be done either by transforming the underlying network, having origins and destinations be distinct nodes adjacent to ``one-way'' centroid connectors, or by excluding such paths when finding paths for travelers, as discussed in Chapter~\ref{chp:tdsp}.

If centroid connectors are set up in this way, then flow entering the network can simply be added to the upstream ends of their centroid connectors, and flow leaving the network at destinations can simply vanish, without any constraints in either case.
The network loading algorithm can then be stated as follows:
\begin{enumerate}
\item Initialize all counts and the time index: $N^\uparrow_{ij}(0) \leftarrow 0$ and  $N^\downarrow_{ij}(0) \leftarrow 0$ for all links $(i,j)$, $t \leftarrow 0$.
\item Use a link model to calculate sending and receiving flows $S_{ij}$ and $R_{ij}$ for all links.
\item Use a node model to calculate transition flows $y_{ijk}$ for all nodes $j$ except for zones.
\item Update cumulative counts: for each non-zone node $i$, perform
\labeleqn{downstreamupdate}{N^\downarrow_{hi} (t + \Delta t) \leftarrow N^\downarrow_{hi}(t) + \sum_{(i,j) \in \Gamma(i)} y_{hij}}
for each upstream link $(h,i)$, and 
\labeleqn{upstreamupdate}{N^\uparrow_{ij} (t + \Delta t) \leftarrow N^\uparrow_{ij}(t) + \sum_{(h,i) \in \Gamma^{-1}(i)} y_{hij}} 
for each downstream link $(i,j)$.
\item Load trips: for all origins $r$, let $D_r(t) = \sum_{s \in Z} d^{rs}(t)$\label{not:Drt} be the total demand starting at this node, and for each centroid connector $(r,i)$, let $p_{ri}$ be the fraction of demand beginning their trips on that connector.
Set $N^\uparrow_{ri}(t + \Delta t) \leftarrow N^\uparrow_{ri}(t) + p_{ri} D_r (t) $ for each connector $(r,i)$.
\item Terminate trips at each destination $s$: for each centroid connector $(i,s)$, set $N^\downarrow_{is}(t + \Delta t) \leftarrow N^\downarrow_{is}(t) + S_{is}$.
\item Increment the time: $t \leftarrow t + \Delta t$.
If $t$ equals the time horizon $T$, then stop.
Otherwise, return to step 2.
\end{enumerate}

Figure~\ref{fig:netloadingexample} and Table~\ref{tbl:netloadingexample} illustrate this process on a network with three links in series, using the spatial queue link model.
The link parameters are shown in Figure~\ref{fig:netloadingexample}; notice that link $(i,j)$ has a larger capacity at the upstream end than at the downstream end.
The node model at $i$ is the ``links in series'' model discussed in Section~\ref{sec:linksinseries}.
The node model at $j$ is a modified version of the ``links in series'' node model: when $5 \leq t < 10 $, $y_{ijs} = 0$ regardless of the sending and receiving flows of $(i,j)$ and $(j,s)$; at all other times the ``links in series'' model is used.
This might reflect a red traffic signal, or a closed drawbridge between these time intervals, since no flow can move through the node.
In this example, the flow downstream is first interrupted when moving from the centroid connector $(r,i)$ to the link $(i,j)$, whose capacity is lower than the rate at which vehicles are being loaded at origin $r$.
This can be seen by examining the difference between the $N^\uparrow$ and $N^\downarrow$ values for link $(r,i)$ during the initial timesteps.
The difference between these values gives the total number of vehicles on the link at that instance in time.
Since $N^\uparrow$ is increasing at a faster rate than $N^\downarrow$, vehicles are accumulating on the link, in a queue at the downstream end.
There is no queue at node $j$, because no bottleneck exists there.
Compare $N^\uparrow_{ij}(t)$ and $N^\downarrow_{ij}(t)$ when $3 \leq t \leq 5$.
Both the upstream and downstream count values increase at the same rate, which means there is no net accumulation of vehicles.

At $t = 5$, the flow through node $j$ drops to zero, which introduces a further bottleneck.
The impacts of the bottleneck are first seen at $t = 6$: 40 vehicles are now on link $(i,j)$, up from 30.
As a result, the receiving flow $R_{ij}$ drops to zero, because the link is full.
Therefore, the node model at $i$ restricts any additional inflow to link $(i,j)$, and the queue on $(r,i)$ grows at an even faster rate than before, even though the number of new vehicles loaded onto the network has dropped.
At $t = 10$, the bottleneck at node $j$ is released, and vehicles begin to move again.
At $t = 25$, the upstream and downstream counts are equal on all links, which means that the network is empty.
All vehicles have reached their destination.
\index{link model!combining with node model|)}
\index{node model!combining with link model|)}

\stevefig{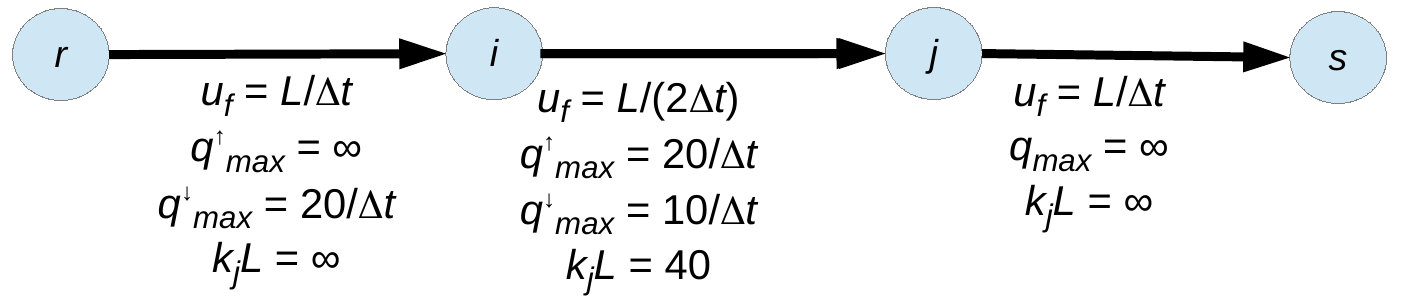}{Example for network loading algorithm.}{\textwidth}

\begin{sidewaystable}
\begin{center}
\caption{Example of network loading algorithm.
\label{tbl:netloadingexample}}
\begin{tabular}{c|c|cccc|c|cccc|c|cccc}
  & & \multicolumn{4}{c|}{Link $(r,i)$} & Node $i$ & \multicolumn{4}{|c|}{Link $(i,j)$} & Node $j$ & \multicolumn{4}{c}{Link $(j,s$)} \\
$t$&$d_{rs}$&$R_{ri}$&$N^\uparrow_{ri}$&$N^\downarrow_{ri}$&$S_{ri}$&$y_{rij}$&$R_{ij}$&$N^\uparrow_{ij}$&$N^\downarrow_{ij}$&$S_{ij}$&$y_{ijs}$&$R_{js}$&$N^\uparrow_{js}$&$N^\downarrow_{js}$&$S_{js}$ \\
\hline
0  &15   &$\infty$ &0 &0 &0 &0 &20   &0 &0 &0 &0 &$\infty$ &0 &0 &0\\
1  &15   &$\infty$ &15   &0 &15   &15   &20   &0 &0 &0 &0 &$\infty$ &0 &0 &0\\
2  &15   &$\infty$ &30   &15   &15   &15   &20   &15   &0 &0 &0 &$\infty$ &0 &0 &0\\
3  &15   &$\infty$ &45   &30   &15   &10   &10   &30   &0 &10   &10   &$\infty$ &0 &0 &0\\
4  &15   &$\infty$ &60   &40   &20   &10   &10   &40   &10   &10   &10   &$\infty$ &10   &0 &10\\
5  &15   &$\infty$ &75   &50   &20   &10   &10   &50   &20   &10   &0 &$\infty$ &20   &10   &10\\
6  &10   &$\infty$ &90   &60   &20   &0 &0 &60   &20   &10   &0 &$\infty$ &20   &20   &0\\
7  &10   &$\infty$ &100  &60   &20   &0 &0 &60   &20   &10   &0 &$\infty$ &20   &20   &0\\
8  &10   &$\infty$ &110  &60   &20   &0 &0 &60   &20   &10   &0 &$\infty$ &20   &20   &0\\
9  &10   &$\infty$ &120  &60   &20   &0 &0 &60   &20   &10   &0 &$\infty$ &20   &20   &0\\
10 &10   &$\infty$ &130  &60   &20   &0 &0 &60   &20   &10   &10   &$\infty$ &20   &20   &0\\
11 &10   &$\infty$ &140  &60   &20   &10   &10   &60   &30   &10   &10   &$\infty$ &30   &20   &10\\
12 &10   &$\infty$ &150  &70   &20   &10   &10   &70   &40   &10   &10   &$\infty$ &40   &30   &10\\
13 &0 &$\infty$ &160  &80   &20   &10   &10   &80   &50   &10   &10   &$\infty$ &50   &40   &10\\
14 &0 &$\infty$ &160  &90   &20   &10   &10   &90   &60   &10   &10   &$\infty$ &60   &50   &10\\
15 &0 &$\infty$ &160  &100  &20   &10   &10   &100  &70   &10   &10   &$\infty$ &70   &60   &10\\
16 &0 &$\infty$ &160  &110  &20   &10   &10   &110  &80   &10   &10   &$\infty$ &80   &70   &10\\
17 &0 &$\infty$ &160  &120  &20   &10   &10   &120  &90   &10   &10   &$\infty$ &90   &80   &10\\
18 &0 &$\infty$ &160  &130  &20   &10   &10   &130  &100  &10   &10   &$\infty$ &100  &90   &10\\
19 &0 &$\infty$ &160  &140  &20   &10   &10   &140  &110  &10   &10   &$\infty$ &110  &100  &10\\
20 &0 &$\infty$ &160  &150  &10   &10   &10   &150  &120  &10   &10   &$\infty$ &120  &110  &10\\
21 &0 &$\infty$ &160  &160  &0 &0 &10   &160  &130  &10   &10   &$\infty$ &130  &120  &10\\
22 &0 &$\infty$ &160  &160  &0 &0 &20   &160  &140  &10   &10   &$\infty$ &140  &130  &10\\
23 &0 &$\infty$ &160  &160  &0 &0 &20   &160  &150  &10   &10   &$\infty$ &150  &140  &10\\
24 &0 &$\infty$ &160  &160  &0 &0 &20   &160  &160  &0 &0 &$\infty$ &160  &150  &10\\
25 &0 &$\infty$ &160  &160  &0 &0 &20   &160  &160  &0 &0 &$\infty$ &160  &160  &0\\
\end{tabular}
\end{center}
\end{sidewaystable}

\section{Elementary Traffic Flow Theory}
\label{sec:trafficflowtheory}

\index{traffic flow theory|(}
This section provides an introduction to the hydrodynamic theory of traffic flow, the basis of several widely-used link models which are more realistic than the point or spatial queue models.
In this theory, traffic is modeled as a compressible fluid.
Its primary advantage is its simplicity: it can capture many important congestion phenomena, while remaining tractable for large networks.
Of course, vehicles are not actually molecules of a fluid, but are controlled by drivers with heterogeneous behavior, who drive vehicles of heterogeneous size, power, and so on.
By treating vehicles as identical particles, we are ignoring such distinctions; as an example of these limitations, we assume that no overtaking occurs within a link (because the particles are identical, none of them has any reason to move faster than another).
But seen as a first-order approximation, the hydrodynamic theory provides a simple and tractable way to model network traffic.
It is also possible to derive certain aspects of the hydrodynamic theory from more behavioral models representing car-following.

In contrast to the link models we will ultimately use in dynamic traffic assignment, fluid-based traffic models are formulated in \emph{continuous} space and time, so it makes sense to talk about the state of traffic at any point in time (not just at the integer timesteps $0, 1, \ldots, T$).
Link models which are based on fluid models will convert these continuous quantities to discrete ones.

\subsection{Traffic state variables}

We start by modeling a roadway link as a one-dimensional object, using $x$ to index the distance from the upstream end of the link, and using $t$ to index the current time.
At any point $x$, and at any time $t$, the state of traffic can be described by three fundamental quantities: the \emph{flow}\index{flow} $q$,\label{not:q} the \emph{density}\index{density} $k$,\label{not:kdensity} and the \emph{speed}\index{speed} $u$.\label{not:u}
Each of these can vary over space and time, so $q(x,t)$, $k(x,t)$, and $u(x,t)$ can be seen as functions defined over all $x$ values on the link, and all $t$ values in the analysis period.
However, when there is no ambiguity (e.g., only looking at a single point at a single time), we can simply write $q$, $k$, and $u$ without providing the space and time coordinates.
\index{volume|see {flow}}

Both $x$ and $t$ are treated as continuous variables,\index{continuum model} as are the vehicles themselves, allowing derivatives of $q$, $k$, and $u$ to be meaningfully defined.
This assumption is imported from fluid mechanics, where the molecules are so small and numerous that there is essentially no error in approximating the fluid as a continuum.
For traffic flow, this assumption is not so trivial, and is one of the drawbacks of hydrodynamic models.

Flow is defined as the rate at which vehicles pass a stationary point, and commonly has units of vehicles per hour.
In the field, flow can be measured using point detectors (such as inductive loops) which record the passage of each vehicle at a fixed location.
Flow can be thought of as a temporal concentration of vehicles.
Density, on the other hand, is defined to be the spatial concentration of vehicles at a given time, and is measured in vehicles per unit length (commonly vehicles per kilometer or vehicles per mile).
Density can be obtained from taking a photograph of a link, noting the concentration of vehicles at different locations at a single instant in time.
The speed is the instantaneous rate at which the vehicles themselves are traveling, and is measured in units such as kilometers per hour or miles per hour.
Speed can be directly measured from radar detectors.

These three quantities are not independent of each other.
As a start, there is the basic relationship\index{flow!relationship with speed and density}\index{speed!relationship with flow and density}\index{density!relationship with flow and speed}
\labeleqn{fundamentalrelationship}{q = uk}
which must hold at each point and time.
(You should check the dimensions of the quantities in this formula to verify their compatibility.)
If this equation is not evident to you from the definitions of flow, density, and speed, imagine that we want to know the number of vehicles $\Delta N$\label{not:DeltaN} which will pass a fixed point over a small, finite time interval $\Delta t$.
If the speed of vehicles is $u$, any upstream vehicle within a distance of $u \Delta t$ from the fixed point will pass during the next $\Delta t$ time units, and the number of such vehicles is
\labeleqn{fundamentaldifference}{\Delta N = k u \Delta t\,.}
The flow rate is approximately $\Delta N/\Delta t$, and the formula~\eqn{fundamentalrelationship} then follows from taking limits as the time increment $\Delta t$ shrinks to zero.

\stevefig{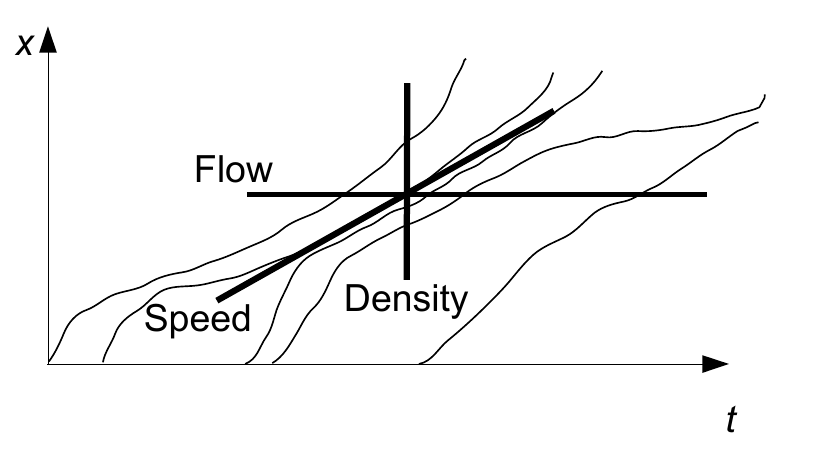}{Trajectory diagram illustrating speed, flow, and density.}{0.7\textwidth}

Figure~\ref{fig:trajectorydiagram} is a \emph{trajectory diagram}\index{trajectory diagram}\index{space-time diagram|see {trajectory diagram}} showing the locations of vehicles on the link over time --- the horizontal axis denotes time, and the vertical axis denotes space, with the upstream end of the link at the bottom and the downstream end at the top.
Speed, flow, and density can all be interpreted in terms of these trajectories.
The speed of a vehicle at any point in time corresponds to the slope of its trajectory there.
Flow is the rate at which vehicles pass a fixed point: on a trajectory diagram, a fixed point in space is represented by a horizontal line.
Time intervals when more trajectories cross this horizontal line have higher flow, and when fewer trajectories cross this line, the flow is lower.
Density is the spatial concentration of vehicles at a particular instant in time: on a trajectory diagram, a specific instant is represented by a vertical line.
Where more trajectories cross this vertical line, the density is higher, and where fewer trajectories cross, the density is lower.

However, this equation by itself is not enough to describe anything of real interest in traffic flow.
Another equation, based on vehicle conservation principles, is described in the next subsection.
The Lighthill-Whitham-Richards model, described at the end of this section, makes a further assumption about the relationships of three state variables.
These relationships are enough to specify the network loading problem as the solution to a well-defined system of partial differential equations.

\subsection{Cumulative counts and conservation}

The hydrodynamic theory can be simplified by introducing a fourth variable $N$, again defined at each location $x$ and at each time $t$.
This new variable $N$ is referred to as the \emph{cumulative count},\index{cumulative count} and has a similar interpretation to the $N^\uparrow$ and $N^\downarrow$ counts defined at the start of Section~\ref{sec:simplelinks}, but now applied at any point on the link, not just the upstream and downstream ends.
Imagine that each vehicle is labeled with a number; for instance, the vehicle at the upstream end of the link at $t = 0$ may be labeled as zero.
The next vehicle that enters the link is then labeled as one, the next vehicle as two, and so forth.
The numbering does not have to start at zero.
What is important that each entering vehicle be given consecutive numbers.
Then $N(x,t)$\label{not:Nxt} gives the number of the vehicle at location $x$ at time $t$ (keeping in mind that we are modeling vehicles as a continuous fluid, so $N(x,t)$ need not be an integer).
The contours of $N(x,t)$ then give the trajectories of individual vehicles in the traffic stream.
(Figure~\ref{fig:cumulativetrajectory})  

\stevefig{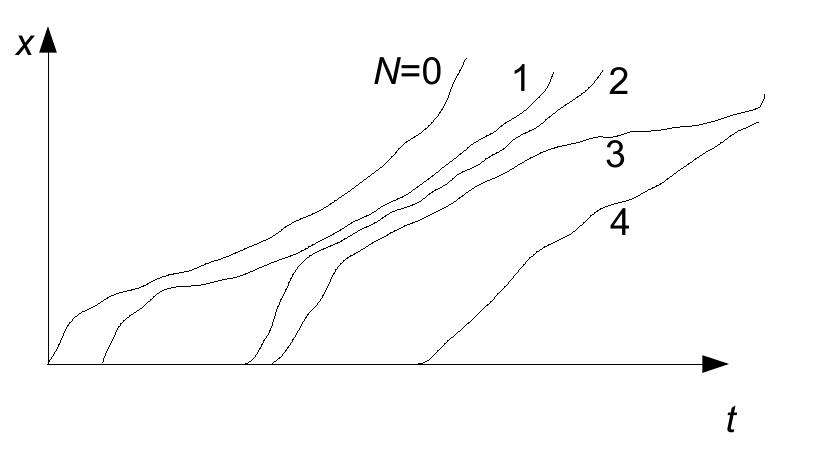}{Vehicle trajectories are contours of $N(x,t)$}{0.7\textwidth}

Two potential points of confusion: in a trajectory diagram like Figure~\ref{fig:cumulativetrajectory}, space (indexed by $x$) is conventionally denoted on the vertical axis, and time (indexed by $t$) on the horizontal axis.
This has the unfortunate side-effect of making the ``$x$-axis'' the vertical one.
It is conventional to list the spatial component before the time component, so a point like $(1,2)$ refers to $x = 1$ and $t = 2$.
This means that on a trajectory diagram, the \emph{vertical} component is listed first, as opposed to typical Cartesian coordinates where the \emph{horizontal} component is given first.
Unfortunately, both of these conventions are so well-established in the transportation engineering literature that it is best to simply highlight them and become comfortable using them.

The quantity $N(x,t)$ is called a cumulative count for the following reason: at the upstream end of the link $(x = 0)$, the quantity $N(0, t)$ gives the cumulative number of vehicles which have entered the link up to time $t$.
Furthermore, at any location $x$, the difference $N(x, t') - N(x,t)$ gives the number of vehicles which passed point $x$ between times $t$ and $t'$.
Taking the limit as $t'$ approaches $t$ provides the relationship\index{cumulative count!relationship with flow}
\labeleqn{cumulativeflow}{\pdr{N}{t} = q\,,}
which holds for every $x$ and $t$ where $N$ is differentiable.
That is, at any fixed spatial location, the rate at which $N$ increases in time is exactly the flow rate $q$.

The cumulative counts can also be related to density with a similar argument.
At any point in time $t$, the difference $N(x', t) - N(x, t)$ gives the number of vehicles which lie between locations $x$ and $x'$ at time $t$.
However, we have to be careful regarding the sign convention regarding $N$ described above: since vehicles are numbered in the order they enter the link, in any platoon of vehicles $N$ decreases as we move from the following vehicles to the lead vehicle.
So, if $x' > x$, then $N(x', t) \leq N(x, t)$.
Thus, taking the limit as $x'$ approaches $x$, we must have\index{cumulative count!relationship with density}
\labeleqn{cumulativedensity}{\pdr{N}{x} = -k\,,}
which again applies wherever $N$ is differentiable, and where the negative sign in the formula results from our sign convention.

In this respect, the cumulative counts $N$ can be seen as the most basic description of traffic flow: if we are given $N(x,t)$ at all points $x$ and times $t$, we can calculate $q(x,t)$ and $k(x,t)$ by using equations~\eqn{cumulativeflow} and~\eqn{cumulativedensity}, and therefore $u(x,t)$ everywhere by using~\eqn{fundamentalrelationship}.

Furthermore, wherever the flow and density are themselves continuously differentiable functions, Clairaut's theorem states that the mixed second partial derivatives of $N$ must be equal, that is,
\labeleqn{clairaut}{\pdrc{N}{x}{t} = \pdrc{N}{t}{x}\,,}
so substituting the relationships~\eqn{cumulativeflow} and~\eqn{cumulativedensity} and rearranging, we have\index{cumulative count!vehicle conservation}\index{conservation equation}\index{flow!relationships with density}\index{density!relationships with flow}
\labeleqn{vehicleconservation}{\pdr{q}{x} + \pdr{k}{t} = 0\,.}
This is an expression of vehicle conservation, that is, vehicles do not appear or disappear at any point.
This equation must hold everywhere that these derivatives exist.\footnote{Of course, vehicle conservation must hold even when these derivatives do not exist, it is just that the formula~\eqn{vehicleconservation} is meaningless there. We have to enforce flow conservation in a different way at such points.}
Equations~\eqn{cumulativeflow} and~\eqn{cumulativedensity} are useful in another way.
If $(x_1, t_1)$ and $(x_2, t_2)$ are any two points in space and time, the difference in cumulative count number between these points is given by the line integral
\labeleqn{curvecount}{N(x_2, t_2) - N(x_1, t_1) = \int_C q~dt - k~dx\,,}
where $C$\label{not:Ccurve} is any curve connecting $(x_1, t_1)$ and $(x_2, t_2)$.
Because vehicles are conserved, this line integral does not depend on the specific path taken.
This is helpful, because we can choose a path which is easy to integrate along.
In what follows, we often choose the straight line connecting these two points as the integration path.

Two issues concerning cumulative counts are often confusing, and are worth further explanation.
First, the cumulative counts can only be meaningfully compared within the same link.
Each link maintains its own counts, and the number associated with a vehicle may change when traveling between links.
Instead of thinking of the cumulative count as being a label permanently associated with a vehicle, it is better to think about it as a label given to the vehicle \emph{by a specific link}, and each link maintains its labels independently of all of the other links.
Each link simply gives successive numbers to each new vehicle entering the link, without having to coordinate its numbering scheme with other links.
This is needed to ensure that link models can function independently, and because in complex networks it is usually impossible to assign ``permanent'' numbers to vehicles such that any two vehicles entering a link consecutively have consecutive numbers.
Second, within any given link, only the \emph{difference} in cumulative counts is meaningful, the absolute numbers do not have specific meaning.
For instance, it may be relevant that 10 vehicles entered the link between times 5 and 6, but the specific numbers of these vehicles are not important.
A common convention is to have the first vehicle entering the link be assigned the number 0, the next vehicle the number 1, and so forth, but this is not required.
In cases where there are already vehicles on the link at the start of the modeling period, the first vehicle entering the link may be assigned a higher number (because the vehicles already on the link must have lower numbers, and a modeler's aesthetics may prefer nonnegative vehicle counts).
But there would be nothing wrong with a negative cumulative count either.
In this regard, different choices of the ``zero point'' are analogous to the different zero points in the Fahrenheit and Celsius temperature scales: either one will give you correct answers as long as you are consistent, and a negative number is not necessarily cause for alarm.

For example, consider a link where $N(x,t) = 1000t - 100x$, with $t$ measured in hours and $x$ measured in miles.
The flow at any point and time is $\partial N/\partial t$, which is a constant of 1000 vehicles per hour.
Likewise, the density is $-\partial N / \partial x = 100$ vehicles per mile.
Using the basic relationship~\eqn{fundamentalrelationship}, the speed must be 10 miles per hour uniformly on the link.

As a more involved example, consider a link which is 1 mile long, with all times measured in minutes and distances in miles.
If we are given that
\labeleqn{givenmap}{N(x,t) = 60t - 120 x + \frac{60x^2}{t+1}\,,}
we can calculate the densities and flows everywhere:
\begin{align}
\label{eqn:calculateddensity} k(x,t) &= -\pdr{N}{x} = 120\myp{1 - \frac{x}{t+1}} \\
\label{eqn:calculatedflow} q(x,t) &= \pdr{N}{t} = 60\myp{1 - \myp{\frac{x}{t+1}}^2} \,.
\end{align}
Exercise~\ref{ex:completelwrexample} asks you to verify that the conservation relationship~\eqn{vehicleconservation} is satisfied by explicit computation.

Ordinarily, the $N(x,t)$ map is not given --- indeed, the goal of network loading is to calculate it.
For if we know $N(x,t)$ everywhere, we can calculate flow, density, and speed everywhere, using equations~\eqn{cumulativeflow},~\eqn{cumulativedensity}, and~\eqn{fundamentalrelationship}.
So let's assume that we \emph{only} know the density and flow maps~\eqn{calculateddensity} and~\eqn{calculatedflow}, and try to recover information about the cumulative counts.
For the given $N$ map, the vehicle at $x = 1/2$ at $t = 0$ has the number 0.
(As discussed above, we do not necessarily have to set the zero point at $N(0,0)$.)  To calculate the number of the vehicle at $x = 1$ and $t = 1$ (the downstream end of the link, one minute later), we can use equation~\eqn{curvecount}.

As this equation involves a line integral, we must choose a path between $(x,t) = (1/2, 0)$ and $(1, 1)$.
Because of the conservation relationship~\eqn{vehicleconservation}, we can choose any path we wish.
For the purposes of an example, we will calculate this integral along three different paths, and verify that they give the same answer.
Figure~\ref{fig:integralpaths} shows the three paths of integration.

\begin{description}
\item[Path A:] This path consists of the line segment from $(1/2, 0)$ to $(1,0)$, followed by the segment from $(1,0)$ to $(1,1)$.
Because these line segments are parallel to the axes, this reduces the line integral to two integrals, one over $x$ alone, and the other over $t$ alone.
We thus have
\begin{align*}
N(1, 1) &= \int_A q~dt - k~dx = -\int_{1/2}^1 k(x,0)~dx + \int_0^1 q(1,t)~dt \\
        &= -\int_{1/2}^1 120\myp{1 - x}~dx + \int_0^1 60\myp{1 - \myp{\frac{1}{t+1}}^2} \\
        &= -15 + 30 = 15\,,
\end{align*}
and the vehicle at the downstream end of the link at $t = 1$ has the number 15.
\item[Path B:] This path consists of the line segment from $(1/2, 0)$ to $(1/2, 1)$, followed by the segment from $(1/2, 1)$ to $(1,1)$.
As before, we have
\begin{align*}
N(1,1) &= \int_B q~dt - k~dx = \int_0^1 q(1/2, t)~dt - \int_{1/2}^1 k(x,1)~dx \\
       &= \int_0^1 60\myp{1 - \myp{\frac{1/2}{t+1}}^2}~dt - \int_{1/2}^1 120\myp{1 - \frac{x}{2}}~dx \\
       &= 52.5 - 37.5 = 15\,.
\end{align*}
\item[Path C:] This path is the line segment directly connecting $(1/2, 0)$ to $(1, 1)$.
Although this line is not parallel to either axis, the integral actually ends up being the easiest to evaluate,
because $x/(t+1)$ is constant along this line, and equal to $1/2$.
Therefore $k(x,t) = 60$ at all points along this line, and $q(x,t) = 45$.
Since $dx = (1/2)dt$ on this line segment, we have
\begin{align*}
N(1,1) &= \int_C q~dt - k~dx = \int_0^1 (45dt - 60(1/2)dt) = \int_0^1 15~dt \\
       &= 15\,.
\end{align*}
\end{description}
All three integrals gave the same answer (as they must), which we can verify by checking $N(1,1)$ with equation~\eqn{givenmap}.
So, we can choose whichever integration path is easiest.
In this example, the integrals in Path B involved the most work.
The integral in Path C required a bit more setup, but the actual integral ended up being very easy, since $q$ and $k$ were constants along the integration path.
Such a path is called a \emph{characteristic},\index{characteristic} and will be described in more detail later in this chapter.

\stevefig{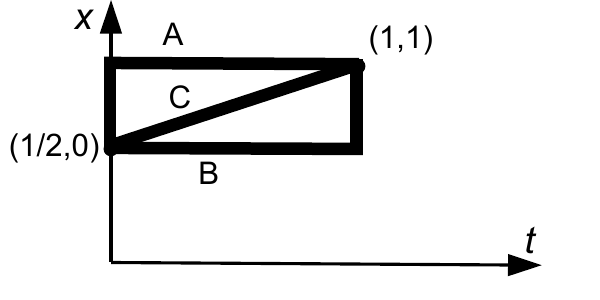}{Three possible paths for the line integral between $(1/2, 0)$ and $(1,1)$.}{0.6\textwidth}

\subsection{The Lighthill-Whitham-Richards model}
\label{sec:lwr}

\index{Lighthill-Whitham-Richards model|see {traffic flow theory, Lighthill-Whitham-Richards model}}
\index{traffic flow theory!Lighthill-Whitham-Richards (LWR) model|(}
At this point, we have two relationships between the flow, density, and speed variables: the basic relationship~\eqn{fundamentalrelationship} and the conservation relationship~\eqn{vehicleconservation}.
These first two relationships can be derived directly from the definitions of these variables, and can describe a wide range of fluid phenomena --- at this point nothing yet has been specific to vehicle flow.
The Lighthill-Whitham-Richards (LWR) model provides a third relationship, completing the hydrodynamic theory.\footnote{Lighthill and Whitham published this model in 1955, as the sequel to a paper on flow in rivers. Richards independently proposed an equivalent model in 1956. All three are now given credit for this model.}

Specifically, the LWR model postulates that the flow at any point is a function of the density at that point, that is,\index{flow!relationships with density}\index{density!relationships with flow}
\labeleqn{fundamentaldiagram}{q(x,t) = Q(k(x,t))}
for some function $Q$.\label{not:Q}
Equivalently, by the relationship~\eqn{fundamentalrelationship}, we can assume that the speed at any point depends only on the density at that point.
It must be emphasized that this relationship, unlike~\eqn{fundamentalrelationship} and~\eqn{vehicleconservation}, is an assumed behavior and does not follow from basic principles.
In dynamic network loading, we typically assume that this function $Q$ is uniform over space and time on a link, an assumption we adopt in this section for simplicity.
It is possible to generalize the results in this section when $Q$ varies over space and time.

The function $Q$ is commonly called the \emph{fundamental diagram};\index{fundamental diagram} an example of such a diagram is shown in Figure~\ref{fig:fundamentaldiagram}.
Fundamental diagrams are concave functions\footnote{A function $f$ is concave if $-f$ is convex.}, and typically assumed to be continuous and piecewise differentiable.
They have two zeros: one at $k = 0$ (zero density means zero flow, because no vehicles are present), and another at the \emph{jam density}\index{density!jam} $k_j$, corresponding to a maximum density where there is no flow because all vehicles are stopped.
At intermediate values of density, the flow is positive, although for a given flow value $q$, there can be two possible density values corresponding to this flow, one corresponding to uncongested conditions and the other to congested conditions.

\stevefig{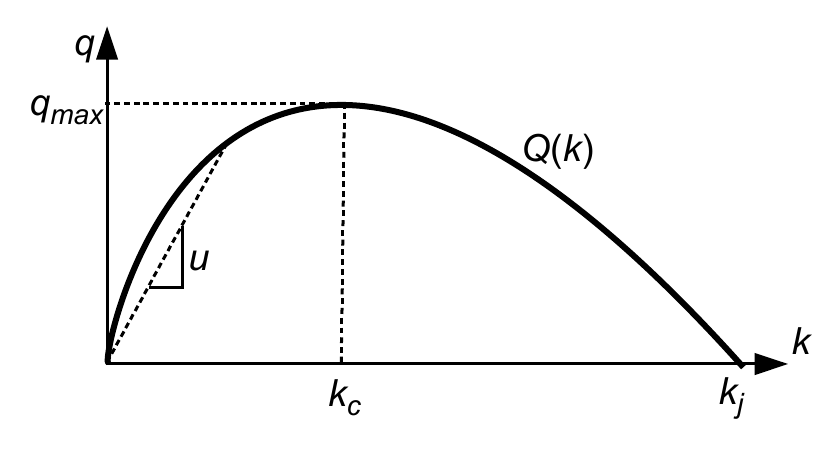}{Fundamental diagram for the LWR model.}{0.7\textwidth}

If $Q$ is concave, and has two zeros at $0$ and $k_j$, then there is an intermediate point where $Q(k)$ is maximal.
This maximal value of $Q$ is called the \emph{capacity}\index{flow!capacity} of the link, denoted $q_{max}$, and the \emph{critical density}\index{density!critical} $k_c$\label{not:kc} is defined to be a value such that $Q(k_c) = q_{max}$.
The $k$ values for which $Q(k) < q_{max}$ and $k < k_c$ are referred to as \emph{subcritical}, and reflect uncongested traffic flow; the $k$ values for which $Q(k) < q_{max}$ and $k > k_c$ are \emph{supercritical}, and reflect congested flow.\footnote{The definitions of subcritical and supercritical in the transportation field are exactly opposite to how these terms are used in fluid mechanics.
This is a bit annoying.
The difference in convention reflects differences in the ``default state'' of flow.
Traffic engineers view the default state of traffic flow as being uncongested, when vehicles can move freely and (as we show later) shockwaves only travel downstream.
Most fluids have to be moving rather quickly for the same state to occur, and the default state of rivers and many other fluid systems is a slower rate of travel, where waves can move both upstream and downstream --- what traffic engineers would describe as a ``congested'' state.}

Using $q = uk$, the speed at any point can be seen as the slope of the secant line connecting the origin to the point on the fundamental diagram corresponding to the density at that point.
That is, in the LWR model, the density $k$ at a point completely determines the traffic state.
The flow $q$ at that point is obtained from the fundamental diagram $Q$, and the speed $u$ can then be obtained from equation~\eqn{fundamentalrelationship1}.

To summarize, the three equations relating flow, density, and speed are:
\begin{align}
 \label{eqn:fundamentalrelationship1} q(x,t) - k(x,t) u(x,t) &= 0\\
 \label{eqn:fundamentaldiagram2} q(x,t) - Q(k(x,t)) &= 0 \\
 \label{eqn:conservation1} \pdr{q}{x} + \pdr{k}{t} &= 0 \,,
\end{align}
and these equations must hold everywhere (with the exception that~\eqn{conservation1} may not be defined if $q$ or $k$ is not differentiable at a point).

Together with initial conditions (such as the values of $k$ along the link at $t = 0$) and boundary conditions (such as the ``inflow rates'' $q$ at the upstream end $x = 0$ throughout the analysis period, or restrictions on $q$ at the downstream end from a traffic signal), this system of equations can in principle be solved to yield $k(x,t)$ everywhere.
Exercise~\ref{ex:completelwrexample} asks you to verify that the $N(x,t)$ map used in the example in the previous section is consistent with the fundamental diagram $Q(k) = \frac{1}{4} k(240 - k)$.

\index{shockwaves|(}
The points where $k$ is not differentiable are known as \emph{shockwaves}, and often correspond to abrupt changes in the density.
Figure~\ref{fig:shockwave} shows an example of several shockwaves associated with the changing of a traffic light.
Notice that in region A, the density is subcritical (uncongested); in region B, traffic is at jam density; and in region C, traffic is at critical density and flow is at capacity.
The speed of a shockwave can still be determined from conservation principles, even though the conservation equation~\eqn{conservation1} does not apply because the density and flow derivatives do not exist at a shock.

\stevefig{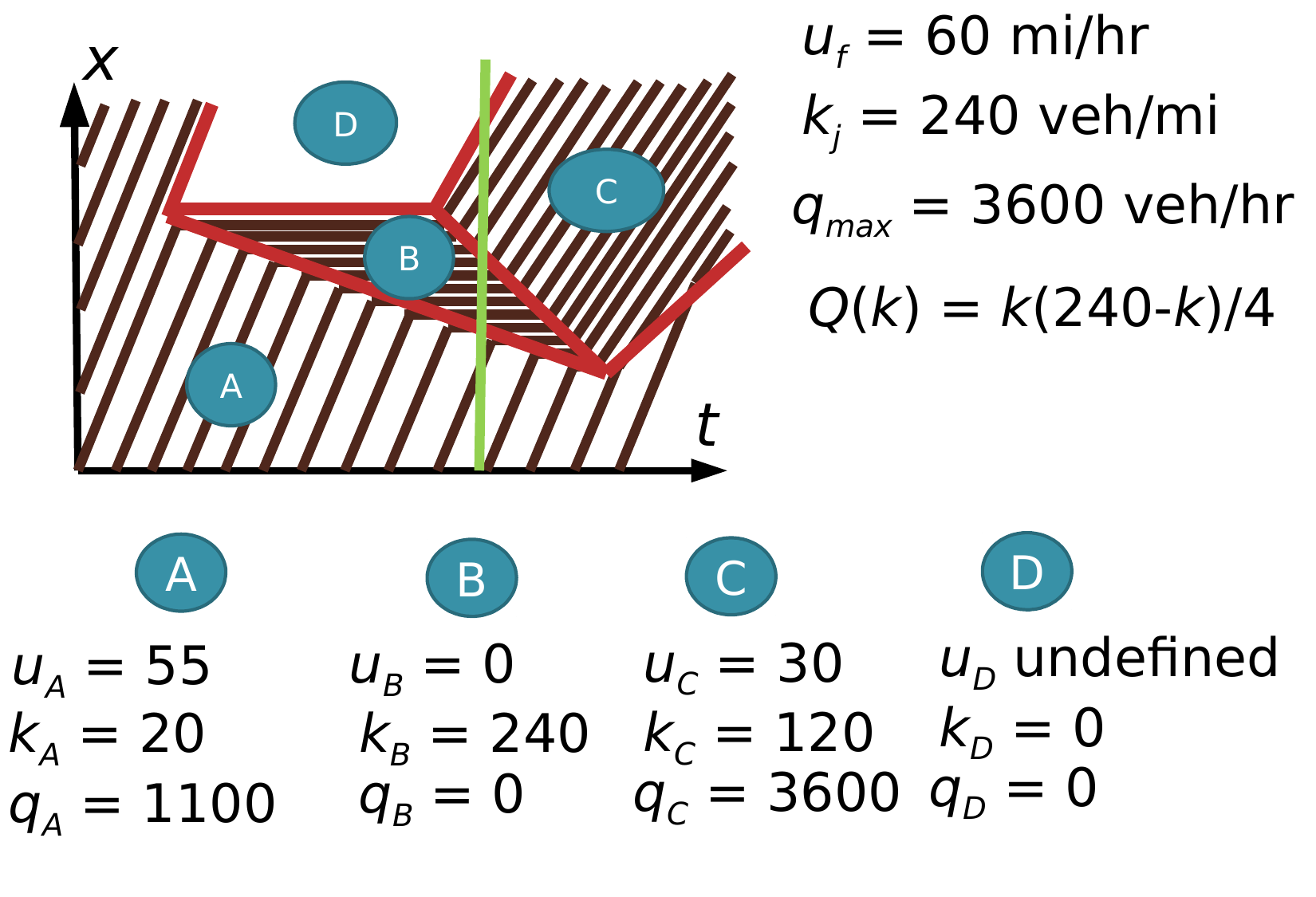}{Shockwaves associated with a traffic signal; vehicle trajectories indicated in brown and shockwaves in red.}{\textwidth}

Assume that $k_A$ and $k_B$ are the densities immediately upstream and immediately downstream of the shockwave (Figure~\ref{fig:shockwavederivation}).
The corresponding flow rates $q_A = Q(k_A)$ and $q_B = Q(k_B)$ can be calculated from the fundamental diagram, and finally the speeds are obtained as $u_A = q_A / k_A$ and $u_B = q_B / k_B$.
\stevefig{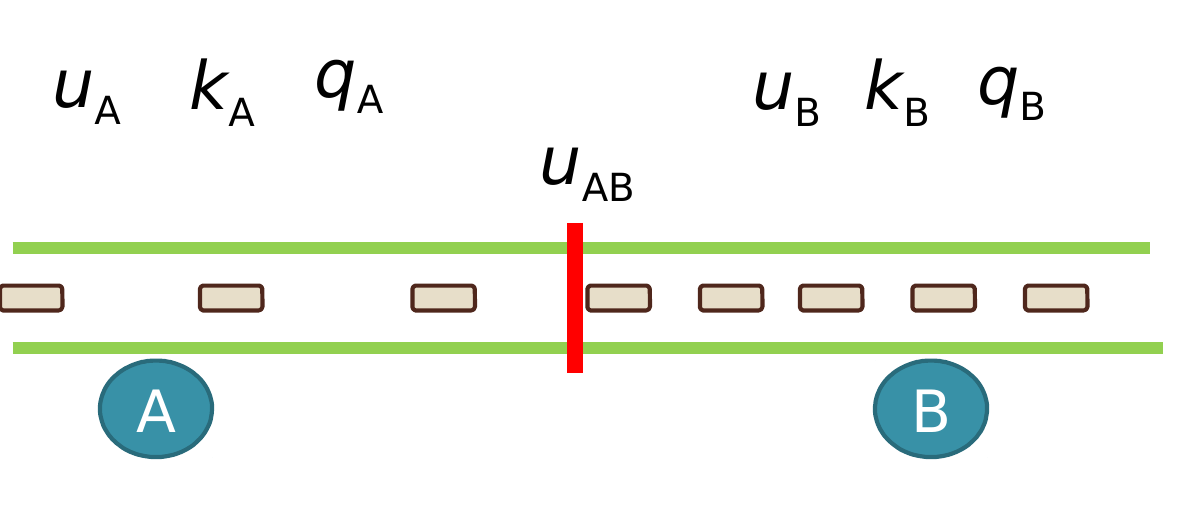}{Flow conservation produces an equation for shockwave speed.}{0.8\textwidth}
Furthermore, let $u_{AB}$\label{not:uAB} denote the speed of the shockwave.
Then the speed of vehicles in region A \emph{relative to the shockwave} is $u_A - u_{AB}$, and the rate at which vehicles cross the shockwave from region A is $(u_A - u_{AB}) k_A$; this is nothing more than equation~\eqn{fundamentalrelationship} as viewed from the perspective of an observer moving with the shockwave.

Likewise, the relative speed of the vehicles in region B is $u_B - u_{AB}$, and the rate at which vehicles cross the shockwave and enter region B is $(u_B - u_{AB}) k_B$.
Obviously these two quantities must be equal, since vehicles do not appear or disappear at the shock.
Equating these flow rates from the left and right sides of the shockwave, we can solve for the shockwave speed:
\labeleqn{shockwavespeed}{u_{AB} = \frac{q_A - q_B}{k_A - k_B}\,.}
Notice that this calculated speed is the same regardless of whether A is the upstream region and B the downstream region, or vice versa.
This equation also has a nice geometric interpretation: the speed of the shockwave is the slope of the line connecting regions A and B on the fundamental diagram (Figure~\ref{fig:shockwavefundamentaldiagram}).

For instance, in Figure~\ref{fig:shockwave}, in region A the flow and density are 1100 vehicles per hour and 20 vehicles per mile, and in region B the flow and density are 0 vehicles per hour and 240 vehicles per mile.
Therefore, using~\eqn{shockwavespeed}, the shockwave between regions A and B has a speed of $(1100 - 0) / (20 - 240) = -5$ miles per hour.
The negative sign indicates that the shockwave is moving \emph{upstream}.
Since A represents uncongested traffic, and region B represents the stopped queue at the traffic signal, the interpretation is that the queue is growing at 5 miles per hour.
Tracing the derivation of equation~\eqn{shockwavespeed}, the rate at which vehicles enter the shockwave is $(55 + 5) \times 20$ from the perspective of region A, or 1200 vehicles per hour.
You should check that the same figure is obtained from the perspective of region B.
Vehicles are entering the queue faster than the upstream flow rate (1200 vs.\ 1100 vph) because the queue is growing upstream, moving to meet vehicles as they arrive.
\index{shockwaves|)}

\stevefig{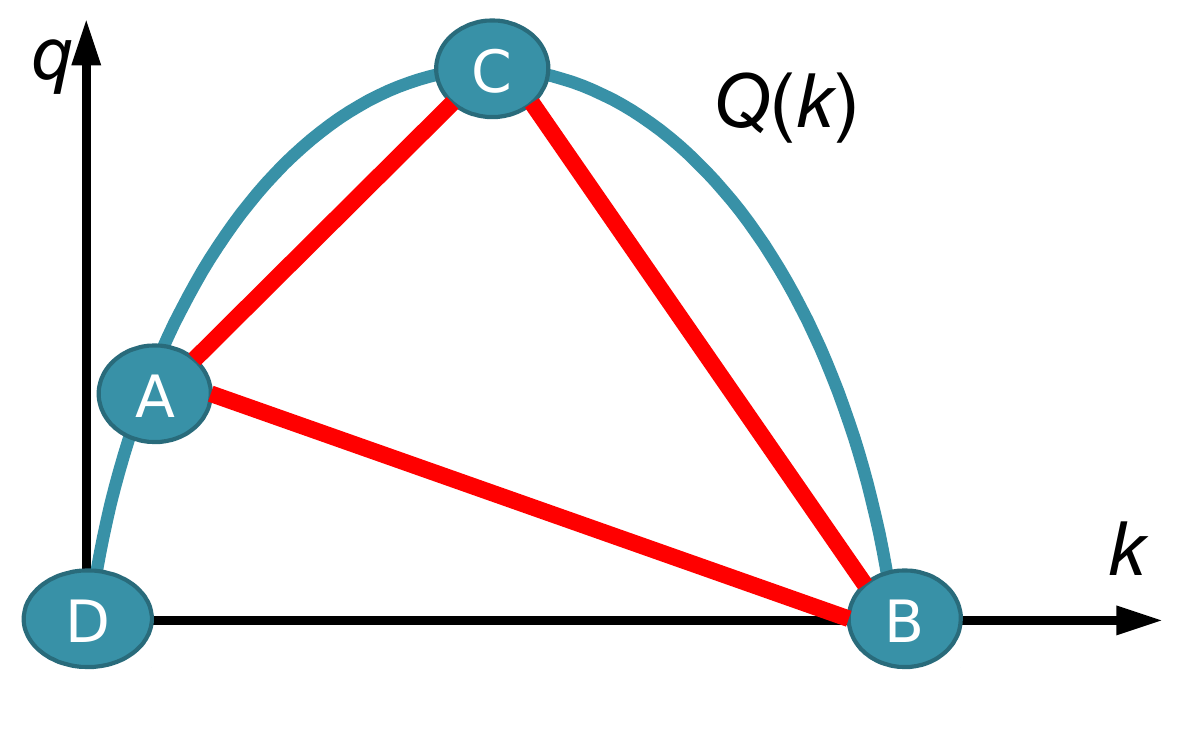}{Identifying shockwave speed on the fundamental diagram.}{0.7\textwidth}

The system of equations~\eqn{fundamentalrelationship1}--\eqn{conservation1} can then be solved, introducing shockwaves as necessary to accommodate the initial and boundary conditions, using~\eqn{shockwavespeed} to determine their speed.
This is the LWR model.
The theory presented above does not immediately suggest a technique for actually solving the system of partial differential equations, which is the topic of the next subsection.

Notice also that the fundamental diagram determines the maximum speed at which a shockwave can move.
Because the fundamental diagram is concave, its slope at any point can never be greater than the free-flow speed $u_f = Q'(0)$, nor can it be less than the slope at jam density $Q'(k_j)$.
The absolute values of these slopes give the fastest speeds shockwaves can move in the downstream and upstream directions, respectively, because shockwave speeds are the slopes of lines connecting points on the fundamental diagram.
This leads to an important notion, the \emph{domain of dependence}.\index{domain of dependence}
Consider a point $(x,t)$ in space and time.
Through this point, draw lines with slopes $Q'(0)$ and $Q'(k_j)$.
The area between these lines to the left of $(x,t)$ represents all the points which can potentially influence the traffic state at $(x,t)$ (labeled as region A in Figure~\ref{fig:lightcone}), and the area between the lines to the right of $(x,t)$ represents all of the points the traffic state at $(x,t)$ can potentially influence (labeled as region B).
In the LWR model, points outside of these regions (the regions labeled C) are independent of what happens at $(x,t)$ (say, a signal turning red or green): in the past, they are either too recent or too distant to affect what is happening at $(x,t)$.
In the future, they are too soon or too distant to be affected by an event at $(x,t)$.
\emph{This is a crucial fact for dynamic network loading.}  If you want to know what is happening at $(x,t)$, it is sufficient to know what has happened in the past, in the domain of dependence.
We do not need to know what is happening simultaneously at other points in the network, and as a result we can perform network loading in an decentralized fashion, performing calculations in parallel since they do not depend on each other.\footnote{If you have studied relativistic physics, there are many similarities with the notion of light cones and causality.}

\stevefig{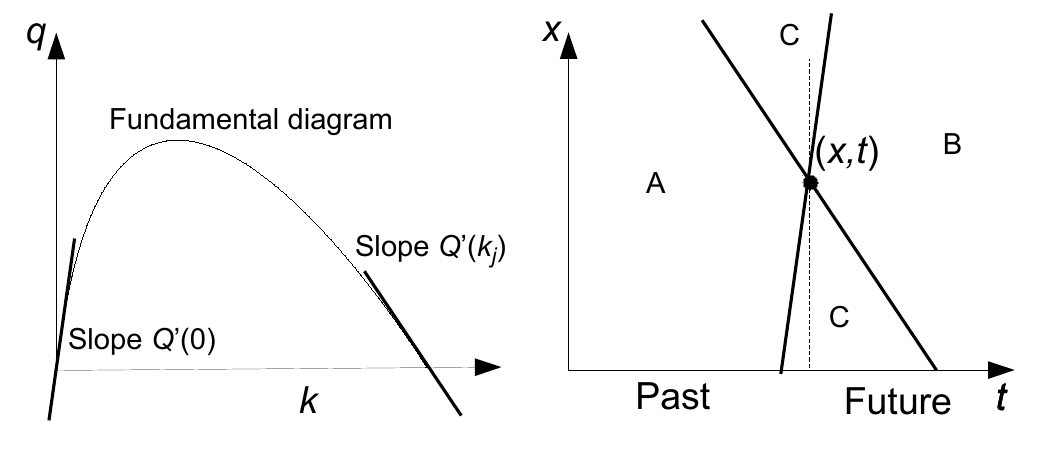}{The fundamental diagram (left) and the domain of dependence (right).}{\textwidth}

\subsection{Characteristics and the Newell-Daganzo method}

\index{characteristics|(}
By substituting equation~\eqn{fundamentaldiagram2} into~\eqn{conservation1}, and setting aside~\eqn{fundamentalrelationship1} for the moment, we can obtain a partial differential equation in $k$ alone:
\labeleqn{densitypde}{\frac{dQ}{dk} \pdr{k}{x} + \pdr{k}{t} = 0 \,. }
It is in this form that the LWR model is most frequently solved, obtaining $k(x,t)$ values everywhere.
This type of partial differential equation can be approached using the \emph{method of characteristics}, which is briefly described below.

A \emph{characteristic} is a curve in $(x,t)$ space, along which the density $k(x,t)$ varies in a predictable way.
If the fundamental diagram $Q$ does not vary in space and time, then straight lines form characteristics of~\eqn{densitypde}.
To demonstrate, assume that we know the density $k(x_1, t_1)$ at some point, and consider a straight line $C$ through $(x_1, t_1)$ with some slope $v$.
Then the directional derivative of $k$ along this line is given by\label{not:ftotder}
\labeleqn{kderiv}{k' = \nabla k \cdot \vect{ v & 1 } = \pdr{k}{x} v + \pdr{k}{t}\,.}
Substituting the conservation equation~\eqn{densitypde} and rearranging, we have
\labeleqn{kderivconserve}{k' = \myp{v - \frac{dQ}{dk}} \pdr{k}{x}\,.}
This gives us a formula for the change in density along the line.
But what is really useful is noticing that if $v = \frac{dQ}{dk}$, then~\eqn{kderivconserve} vanishes and \emph{the density is constant along the line}.
Put another way, the density is constant along any line whose slope is equal to the slope of the \emph{tangent} of the fundamental diagram at that density value.\footnote{Note that the speed of the characteristic is different from the speed of the vehicles themselves, or the speed of shockwaves (which are given by slopes of secant lines on the fundamental diagram).}

To see this, return to the example we have been using where $Q(k) = \frac{1}{4} k (240 - k)$ and $N(x,t) = 60t - 120 x + \frac{60x^2}{t+1}$.
We have calculated the density and flow at all points and times in equations~\eqn{calculateddensity} and~\eqn{calculatedflow}.
In our calculation, we observed that the line integral along path C was the simplest to evaluate, because $k(x,t)$ was a constant 60 vehicles per mile along the line $x/(t+1) = 1/2$.
This path is in fact a characteristic: the derivative of the fundamental diagram at this point is $Q'(60) = 60 - k / 2 = 30$ miles per hour.
Rearranging the equation of the line, we have $x = t/2 + 1/2$; in other words, this is a line whose location moves half a mile in one minute: 30 miles per hour, the same as the derivative of the fundamental diagram.

These slopes will be different at points where the density is different, which means that characteristic lines can potentially intersect.
This indicates the presence of a shockwave separating regions of different density.
One interpretation of these characteristic lines is as ``directions of influence,'' since the density at a point will determine the density at all later points along this line.
In uncongested regions, where $k$ is subcritical, $dQ/dk$ is positive, meaning that the characteristics have positive slope: uncongested states will propagate downstream.
In congested regions with supercritical density, $dQ/dk$ is negative, and the characteristics have negative slope: congested states will propagate upstream.

\index{Newell-Daganzo method|(}
Furthermore, by combining knowledge of characteristics with equation~\eqn{curvecount}, we can determine the cumulative count and the density at any point, given sufficient initial and boundary data.
Suppose we wish to calculate the density at a point $(x,t)$.
If we knew the density at this point, then we would know the slope of the characteristic through this point.
This characteristic could be traced back until it intersected a point $(x_0, t_0)$ where the cumulative count $N$ was known, either because it corresponds to an initial or boundary point, or a point where $N$ has already been calculated.
Then, applying~\eqn{curvecount}, we would have
\labeleqn{newellcountlineintegral}{N(x,t) = N(x_0, t_0) + \int_C q~dt - k~dx\,,}
where $C$ is the straight line between $(x_0, t_0)$ and $(x,t)$.
Since this line is a characteristic, we have $dx/dt = dQ/dk$, and so
\labeleqn{newellsimplified}{N(x,t) = N(x_0, t_0) + \int_C \myp{q - k\frac{dQ}{dk}}dt\,.}
Furthermore, as a characteristic, $k$ (and therefore $q$) are constant.
So the integral is easy to evaluate, and we have
\labeleqn{newellcount}{N(x,t) = N(x_0, t_0) + \myp{q - k\frac{dQ}{dk}} (t - t_0)\,.}
The only trouble is that we do not actually know the density at $(x,t)$.
Each possible value of density corresponds to a slightly different cumulative count, based on equation~\eqn{newellcount}.
The insight of the Newell-Daganzo method is that \emph{the correct value of the cumulative count is the lowest possible value}.
That is, imagine that the density at $(x,t)$ is $k$, and let $(x_k, t_k)$ be the known point corresponding to the characteristic slope of density $k$.
Then\label{not:inf}
\labeleqn{newellfinal}{N(x,t) = \inf_{k \in [0, k_j]} \myc {N(x_k, t_k) + \myp{q - k\frac{dQ}{dk}} (t - t_k)}
\,.\footnote{In this formula, `inf' means the \emph{infimum}, or greatest lower bound, of the set in question.
On a finite set, the infimum and minimum are the same thing.
Some infinite sets do not have a minimum, but still have an infimum; for instance there is no smallest positive real number, but the infimum of that set is zero.
}
}
Rigorously validating this insight requires knowledge of the calculus of variations, which is beyond the scope of this text.
An intuitive justification is that~\eqn{newellcount} represents an upper bound on the cumulative count imposed by a boundary or initial point --- we know the number of the vehicle passing that boundary or initial point, and~\eqn{newellcount} expresses one possible value of the cumulative count at the point in question, if the density took a particular value and no shockwave intervened.
However, due to the possible presence of shockwaves, there may be another, more restrictive constraint imposed on the cumulative count.
The equation~\eqn{newellfinal} thus finds the ``most restrictive'' boundary or initial condition, which gives the correct cumulative count.

\stevefig{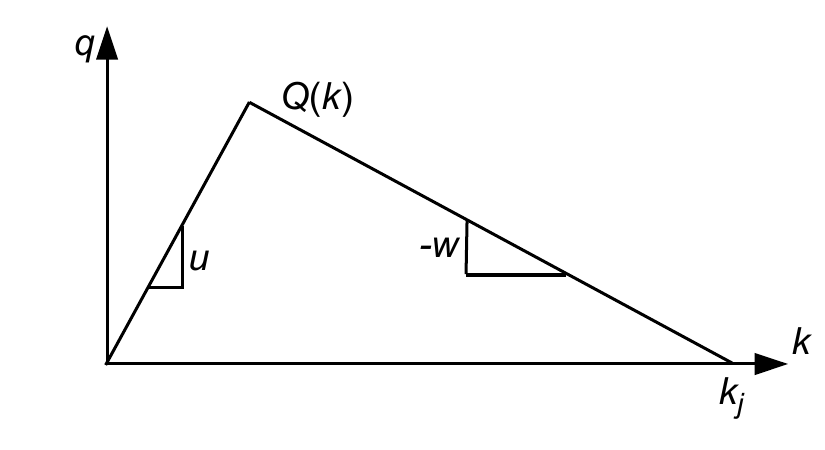}{A triangular fundamental diagram, defined by three parameters $u_f$, $w$, and $k_j$.}{0.7\textwidth}

Better yet, this method becomes exceptionally easy if we assume that the fundamental diagram takes a simple form, such as a triangular shape (Figure~\ref{fig:triangularfd}), where the equation is given by
\labeleqn{triangularfd}{Q(k) = \min \myc{u_f k, w(k_j - k)}\,.}
In this case, there are only two possible characteristic speeds: one ($+u_f$) corresponding to uncongested conditions, and the other ($-w$) corresponding to congested conditions.
The uncongested speed $u_f$ is the free-flow speed, and $w$ is known as the \emph{backward wave speed}.\index{speed!backward wave speed}
In the case where the characteristic speed is $u_f$, the vehicle speed $u$ equals the characteristic speed $u_f$, and both are equal to $q/k$.
Therefore, the line integral along the characteristic in~\eqn{newellcount} is
\labeleqn{uncongestedcharacteristic}{\myp{q - k u_f} (t - t_0) = \myp{q - uk}(t - t_0) = 0}
because $q = uk$.

In other words, \emph{the vehicle number is constant along characteristics at free-flow speed.}
For the congested characteristic with slope $-w$, we have
\labeleqn{congestedcharacteristic}{\myp{q - k (-w)} (t - t_0) = w \myp{\frac{q}{w} + k}(t - t_0) = k_j w (t - t_0)}
since $k + \frac{q}{w} = k_j$, as can be seen in Figure~\ref{fig:triangularfd}.
This expression can also be written as $k_j (x_0 - x)$.
Since these are the only two characteristics which can prevail at any point, equation~\eqn{newellfinal} gives
\labeleqn{newelltriangular}{N(x,t) = \min \myc{ N(x_U, t_U) , N(x_C, t_C) + k_j (x_C - x) }}
where $(x_U, t_U)$ is the known point intersected by the uncongested characteristic, and $(x_C, t_C)$ is the known point intersected by the congested characteristic.

\stevefig{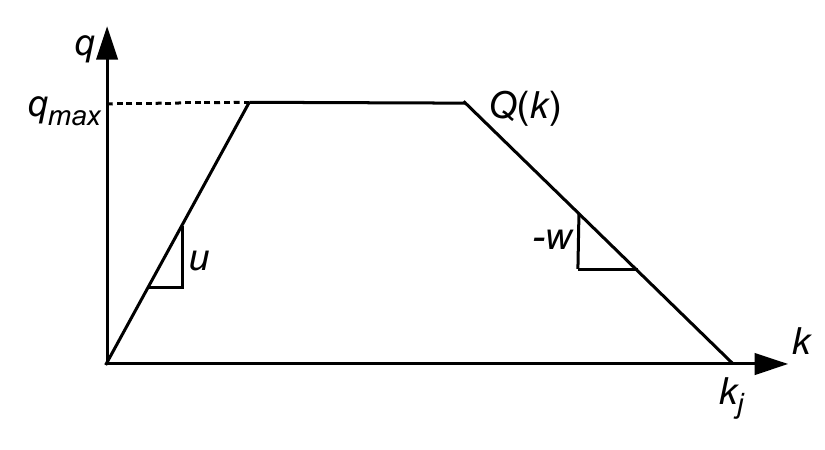}{A trapezoidal fundamental diagram, defined by four parameters $u_f$, $q_{max}$, $w$, and $k_j$.}{0.7\textwidth}

Trapezoidal fundamental diagrams, such as that in Figure~\ref{fig:trapezoidalfd}, are also commonly used, given by the equation
\labeleqn{trapezoidalfd}{Q(k) = \min \myc{ u_f k, q_{max}, w(k_j - k)}\,.}
In this case, there is a third possible characteristic speed, corresponding to the flat region where flow is at capacity.
Since the derivative of the fundamental diagram at this point is zero, this characteristic has zero speed, represented by a horizontal line on a space-time diagram.
Tracing this characteristic back to a known point $(x, t_R)$, the change in vehicle number between this known point and the unknown point $(x,t)$ is just
\labeleqn{capacitycharacteristic}{ q_{max} (t - t_R) }
since $dx = 0$ in the line integral~\eqn{curvecount}.
(Note that the location of this third known point is the same as the location of the point we are solving for, since the characteristic is stationary.)  This adds a third term to the minimum in~\eqn{newelltriangular}, giving
\labeleqn{newelltrapezoidal}{N(x,t) = \min \myc{ N(x_U, t_U) , N(x, t_R) + q_{max}(t - t_R), N(x_C, t_C) + k_j (x_C - x) }}
for trapezoidal fundamental diagrams.

The trapezoidal fundamental diagram requires four parameters to calibrate: the free-flow speed $u_f$, the capacity $q_{max}$, the jam density $k_j$, and the backward wave speed $-w$.
The first three of these are fairly straightforward to estimate from traffic engineering principles.
The backward wave speed is a bit trickier; empirically it is often a third to a half of the free-flow speed.

As a demonstration of this method, consider a link which is 1 mile long.
The fundamental diagram is shown in the left panel of Figure~\ref{fig:newellexample}, and has the equation
\labeleqn{newellexamplefd}{Q(k) = \min\myc{ k, 120 - \frac{1}{2}k}}
when flow is measured in vehicles per minute, and density in vehicles per mile.
Initially, vehicles on the link flow at an uncongested 48 veh/min --- this state has existed for a long time in the past, and vehicles continue entering the link at this rate.
However, at the downstream end of the link there is an obstruction which prevents any vehicles from passing, such as a red light or an incident blocking all lanes.
This will cause a queue of stopped vehicles to form.
Assume that we want to know how many vehicles lie between the obstruction and three given points: half a mile upstream of the obstruction, 30 seconds after it begins; an eighth of a mile upstream of the obstruction, 30 seconds after it begins; and an eighth of a mile upstream of the obstruction, one minute after it begins.
We do not know, \emph{a priori}, whether these points lie within the queue or in the portion of the link which is still uncongested.

The first step is to establish a coordinate system.
As always, we set $x = 0$ at the upstream end of the link.
For this problem, it will be convenient to set $t = 0$ at the time when the obstruction begins, and to count vehicles starting from the first vehicle stopped at the obstruction, that is, $N(1,0) = 0$.
This way, $N(x,t)$ will immediately give the number of vehicles between point $x$ and the obstruction at time $t$.
Next, we use the given data from the problem to construct initial and boundary conditions where we already know $N(x,t)$.
Since no vehicles can pass the obstruction, we know that $N(1,t) = 0$ for all $t$.
Since the link is initially uncongested at a flow rate of 48 veh/min, the fundamental diagram~\eqn{newellexamplefd} gives the initial density to be 48 veh/mi.
Therefore, the initial condition is $N(x,0) = 48 - 48x$, and the vehicle number at the origin of the coordinate system is 48.
Since vehicles continue to enter the link at a rate of 48 veh/min, we have $N(0,t) = 48 + 48t$ along the upstream boundary of the link.
The three points where we must calculate $N(x,t)$ are labeled as A, B, and C in the right panel of Figure~\ref{fig:newellexample}.

We start with point A, half a mile upstream of the obstruction and 30 seconds after it begins.
There are two possible characteristics at this point, one with slope $+1$ (corresponding to uncongested conditions) and one with slope $-\frac{1}{2}$ (corresponding to congested conditions).
We can trace back these characteristics until they reach a point where $N(x,t)$ is known, in this case an initial or boundary condition.
These points of intersection are labeled D and E in Figure~\ref{fig:newellexample}.
From the initial condition, we know that $N(D) = 48$ and $N(E) = 12$.
Along the uncongested characteristic, there is no change in the cumulative count, while along the congested characteristic the cumulative count increases at a rate of $k_j = 240$ veh/mi for each mile traveled.
Therefore, equation~\eqn{uncongestedcharacteristic} tells us that $N(A) = 48 + 0 = 48$ if point A is uncongested, while equation~\eqn{congestedcharacteristic} tells us that $N(A) = 12 + \frac{1}{4} 240 = 72$ if point A is congested.
The correct value is the smaller of the two: $N(A) = 48$, and this point is uncongested (the queue has not yet reached this point).

\stevefig{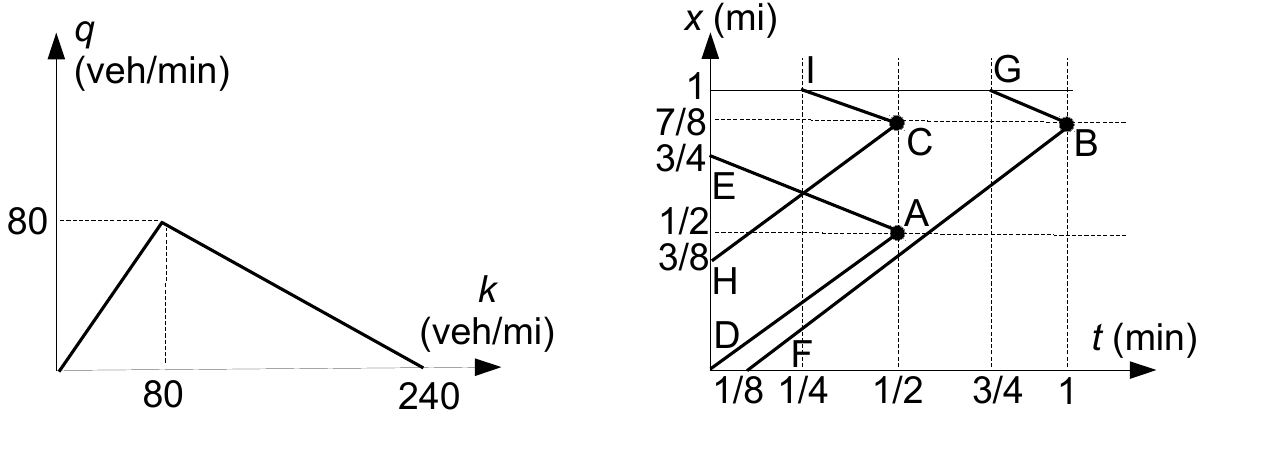}{Example of the Newell-Daganzo method.}{\textwidth}

We next move to point B, an eighth of a mile upstream of the obstruction and 60 seconds after it begins.
Tracing back the two possible characteristics from point B leads us to the points labeled F and G, and from the boundary conditions we know $N(F) = 54$ and $N(G) = 0$.
So $N(B)$ is the lesser of $N(F) + 0 = 54$ and $N(G) + \frac{1}{8} 240 = 30$.
This means that there are 30 vehicles between point B and the obstruction, and since the congested characteristic produced the lower value, point B lies within the queue.

At point C, the two characteristics lead to the points labeled H and I, where $N(H) = 30$ from an initial condition and $N(I) = 0$ from a boundary condition.
So $N(C) = \min \{ 30 + 0, 0 + \frac{1}{8} 240 \} = 30$.
In this case, both characteristics led to the \emph{same} value of $N(C)$.
This indicates that the shockwave passes exactly through point C.

An alternative approach for this problem would be to explicitly calculate the location of the shockwave, determine the densities in each region, and apply~\eqn{curvecount} directly.
For this problem, that approach would be simpler than the Newell-Daganzo method.
However, if there were multiple shockwaves introduced into the problem (say, from multiple red/green cycles), it would become very tedious to track the locations of all of the shockwaves and determine which region the points lie in.
The Newell-Daganzo method can be applied just as easily in such a case, once the boundary conditions are determined.
\index{Newell-Daganzo method|)}
\index{characteristics|)}
\index{traffic flow theory!Lighthill-Whitham-Richards (LWR) model|)}
\index{traffic flow theory|)}

\section{LWR-based Link Models}
\label{sec:fancylinks}

\index{link model|(}
The hydrodynamic traffic flow model developed by Lighthill, Whitham, and Richards forms the basis for several popular link models.
The LWR model is simple enough to be usable in large-scale dynamic network loading, while capturing enough key properties of traffic flow for its results to be meaningful.
Through shockwaves, we can capture how congestion grows and shrinks over time.
These shockwaves allow us to model queue spillback (when a congestion shockwave reaches the upstream end of a link) and to account for delays in queue startup (unlike the spatial queue model).
This section describes two link models based on the LWR model --- the cell transmission model, which is essentially an explicit solution scheme for the LWR system of partial differential equations, and the link transmission model, which uses the Newell-Daganzo method to directly calculate sending and receiving flows.
Lastly, we show how the point and spatial queue models can be seen as special cases of LWR-based link models, with a suitable choice of the fundamental diagram.

\subsection{Cell transmission model}

\index{link model!cell transmission model|(}
In the cell transmission model, in addition to discretizing time into intervals of length $\Delta t$, we also discretize space, dividing the link into cells of length $\Delta x$.\label{not:Deltax}
These two discretizations are not chosen independently.
Rather, they are related by
\labeleqn{ctmdiscretization}{\Delta x = u_f \Delta t \,,}
that is, the length of each cell is the distance a vehicle would travel in time $\Delta t$ at free flow.
The reasons for this choice are discussed at the end of this section.

With this discretization in mind, we will use the notation $n(x,t)$\label{not:nxt} to describe the number of vehicles in cell $x$ at time $t$, where $x$ and $t$ are both integers --- we must convert the continuous LWR variables $k$ and $q$ into discrete variables for dynamic network loading, which we will call $n$ and $y$.

If the cell size is small, we can make the approximation that
\labeleqn{ctmdensityapprox}{n(x,t) \approx k(x,t) \Delta x\,,}
essentially assuming that the density within the cell is constant.
Further define $y(x,t)$\label{not:yxt} to be the number of vehicles which \emph{enter} cell $x$ during the $t$-th time interval.
Making a similar assumption, we can make the approximation
\labeleqn{ctmflowapprox}{y(x,t) \approx q(x,t) \Delta t\,.} 
In a space-time diagram showing vehicle trajectories, such as Figure~\ref{fig:spacetimediagramctm}, $n(x,t)$ and $y(x,t)$ respectively correspond to the number of trajectories crossing the vertical line at $t$ between locations $x$ and $x + \Delta x$, and the number of trajectories crossing the horizontal line at $x$ between times $t$ and $t + \Delta t$.

\stevefig{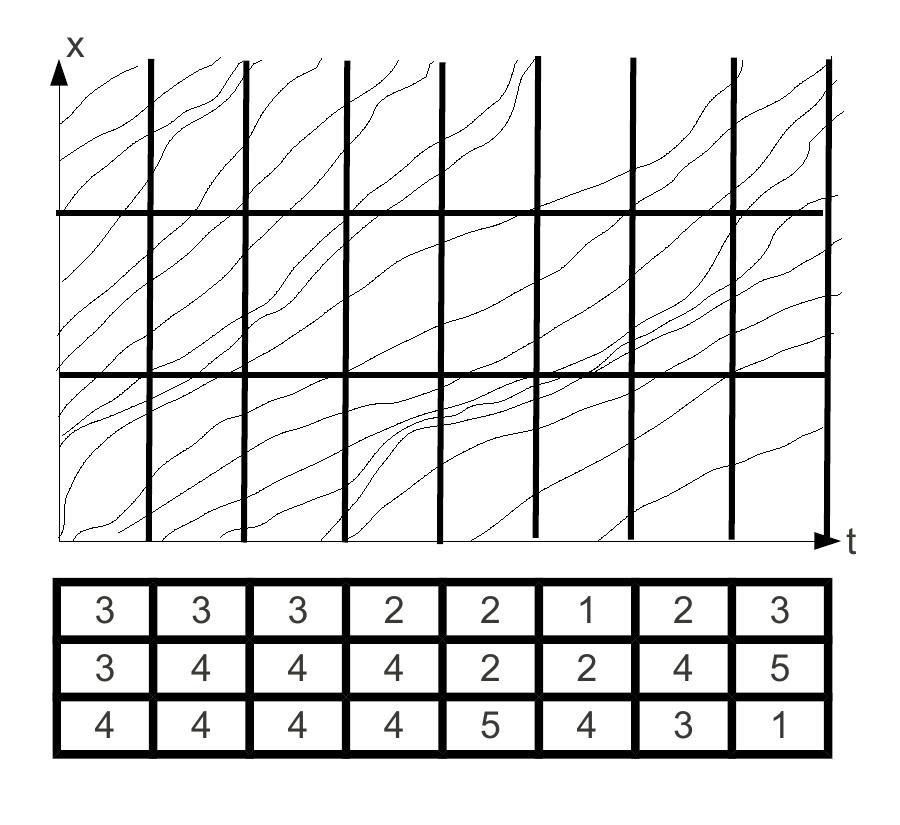}{Discretizing space and time into cells, $n(x,t)$ values below.}{0.7\textwidth}

The cell transmission model provides methods for solving for $n(x,t)$ for all integer values of $x$ and $t$; these can then be converted to density values $k(x,t)$ through~\eqn{ctmdensityapprox}.
These density values can then be used to calculate flows $q(x,t)$ through the fundamental diagram, and $u(x,t)$ values through equation~\eqn{fundamentalrelationship}, finally providing an approximate solution to the system of partial differential equations ~\eqn{fundamentalrelationship1}--\eqn{conservation1}.
Recall that the goal of a link model is to determine sending and receiving flow at each time step.
For this reason, we will be content with determining how the $n(x,t+\Delta t)$ values can be calculated, given the $n(x,t)$ values (which are already known), and the $y(x,t)$ values, which must be calculated.

Since $q(x,t) = Q(k(x,t))$, substitution into equations~\eqn{ctmdensityapprox} and~\eqn{ctmflowapprox} give
\labeleqn{ctmdiscrete}{y(x,t) \approx Q \myp{ \frac{n(x,t)}{\Delta x}} \Delta t\,.}
Substituting the particular form of the fundamental diagram gives an equation for $y(x,t)$ in terms of $n(x,t)$.
Using the trapezoidal diagram of Figure~\ref{fig:trapezoidalfd}, we have 
\labeleqn{ctmtransition1}{y(x,t) \approx  \min \myc{ u_f n(x,t) \frac{\Delta t}{\Delta x}, q_{max} \Delta t, w \Delta t \myp{ k_j - \frac{n(x,t)}{\Delta x} }}\,.}
Using the fact that $\Delta x / \Delta t = u_f$, the first term in the minimum is simply $n(x,t)$.
The third term can be simplified by defining $\bar{n}(x) = k_j \Delta x$ to be the maximum number of vehicles which can fit into a cell and $\delta = w / u_f$ to be the ratio between the backward wave speed\index{speed!backward wave speed} and free-flow speed.\index{speed!free-flow speed}
Then, factoring out $1 / \Delta x$ from the term in parentheses and again using $\Delta x / \Delta t = u_f$, the third term simplifies and we finally obtain
\labeleqn{ctmtransitionfinal}{y(x,t) \approx \min \myc{ n(x,t), q_{max} \Delta t, \delta(\bar{n}(x) - n(x,t))}\,.}

\emph{There is one more point which is subtle, yet incredibly important.}  Being a ``flow'' variable, $y(x,t)$ is calculated at a single point (over a time interval), while $n(x,t)$ is calculated at a single time (over a longer spatial interval).
As shown in Figure~\ref{fig:ctmdiscrete}, the $x$ in $y(x,t)$ refers to a single location, while the $x$ in $n(x,t)$ refers to an entire cell.
So, when we are calculating the flow across the (single) point $x$, which is the boundary between two cells, do we look at the adjacent cell upstream $n(x-1,t)$, or the adjacent cell downstream $n(x,t)$?  

\stevefig{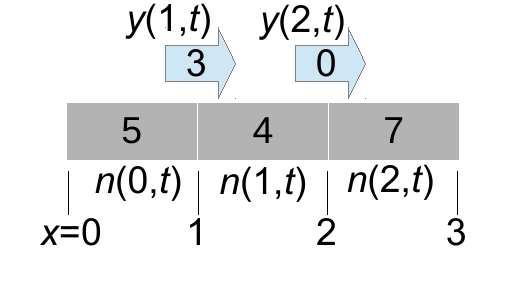}{Where discrete values are calculated in the cell transmission model.}{0.4\textwidth}

The correct answer depends on the fundamental diagram, and the meaning of characteristics.\index{characteristic}
In uncongested conditions, corresponding to the increasing part of the fundamental diagram and the first term in the minimum, the traffic state moves from upstream to downstream (because the characteristic has positive speed).
In congested conditions, corresponding to the decreasing part of the fundamental diagram and the third term in the minimum, the traffic state moves from downstream to upstream (because the characteristic has negative speed.)  So, if traffic is uncongested at the (single) point $x$, we need to refer to the upstream cell, while if traffic is congested at the (single) point $x$, we must refer to the downstream cell.
So, the final expression for the cell transmission model flows is
\labeleqn{ctmtransitionflows}{y(x,t) = \min \{ n(x-\Delta x,t), q_{max} \Delta t, \delta \myp{ \bar{n}(x) - n(x,t)} \}\,.}
This expression also has a nice intuitive interpretation.
The number of vehicles moving from cell $x - 1$ to cell $x$ is limited either by the number of vehicles in the upstream cell (the first term), the capacity of the roadway (the second term), or by the available space in the downstream cell (the third term).

Expression~\eqn{ctmtransitionflows} is the discrete equivalent of the differential equation~\eqn{fundamentaldiagram2}.
We now derive the discrete form of the partial differential equation giving the conservation law~\eqn{conservation1}.
The derivative $\pdr{q}{x}$ can be approximated as
\labeleqn{flowderapprox}{\pdr{q}{x}(x,t) \approx \frac{1}{\Delta t \Delta x} \myp{ y(x + \Delta x, t) - y(x, t)}}
and the derivative $\pdr{k}{t}$ can be approximated as
\labeleqn{densityderapprox}{\pdr{k}{t}(x,t) \approx \frac{1}{\Delta t \Delta x} \myp{ n(x, t + \Delta t) - n(x, t) }\,.}
Substituting into~\eqn{conservation1}, we have
\labeleqn{conservationtemp}{\frac{1}{\Delta t \Delta x} \myp{y(x+\Delta x,t) - y(x,t) + n(x, t+\Delta t) - n(x,t)} = 0 }
or, in a more convenient form,
\labeleqn{conservationctm}{n(x,t + \Delta t) = n(x,t) + y(x,t) - y(x+\Delta x,t)\,.}
This also has a simple intuitive interpretation: the number of vehicles in cell $x$ at time $t + \Delta t$ is simply the number of vehicles in cell $x$ at the previous time $t$, plus the number of vehicles which flowed into the cell during the $t$-th time interval, minus the number of vehicles which left.

Together, the equations~\eqn{ctmtransitionflows} and~\eqn{conservationctm} define the cell transmission model for trapezoidal fundamental diagrams.
There are only two pieces of ``missing'' information, at the boundaries of the link.
Refer again to Figure~\ref{fig:ctmdiscrete}.
How should $y(0,t)$ and $y(L, t)$ be calculated?  For $y(0,t)$, the first term in formula~\eqn{ctmtransitionflows} involves $n(-\Delta x, t)$, while for $y(L, t)$, the third term in the formula involves $n(L + \Delta x, t)$, and both of these cells are ``out of range.''  The answer is that these boundary flows are used to calculate the sending and receiving flows for the link, and a node model will then give the actual values of the link inflows $y(0,t)$ and link outflows $y(L, t)$.

Remember that the sending flow is the maximum number of vehicles which could leave the link if there was no obstruction from downstream.
In terms of~\eqn{ctmtransitionflows}, this means that the third term in the minimum (which corresponds to downstream congestion) is ignored.
Then, the first two terms in the minimum (which only refer to cells on the link) are the possible restrictions on the flow leaving the link, so, using $C$ to refer to the downstream-most cell on the link, we have\index{link model!cell transmission model!sending flow}
\labeleqn{ctmsending}{S(t) = \min \{ n(C, t), q_{max} \Delta t \}\,.
}
Likewise, the receiving flow is the maximum number of vehicles which could enter the link if there were a large number of vehicles wanting to enter from upstream.
In terms of~\eqn{ctmtransitionflows}, this means that the first term in the minimum (which corresponds to the number of vehicles wanting to enter) is ignored.
The second two terms in the minimum refer to cells on the link, and\index{link model!cell transmission model!receiving flow}
\labeleqn{ctmreceiving}{R(t) = \min \{ q_{max} \Delta t, \delta \myp{\bar{n}(0) - n(0, t)} \}\,.}

It remains to explain the choice of the discretization~\eqn{ctmdiscretization}.
An intuitive explanation for linking the cell length and time discretizations in this way is that this choice limits vehicles to moving at most one cell between time steps.
In fact, in uncongested conditions, \emph{all} vehicles in a cell will move to the next cell downstream in the next time interval, simplifying calculations --- for instance, the simplification in~\eqn{ctmtransition1} only works because of the choice made in~\eqn{ctmdiscretization}.
The underlying mathematical reason has to do with the speed of characteristics,\index{characteristics} which for the trapezoidal fundamental diagram lie between $-w$ and $+u_f$.
Empirically, $w < u_f$, so the fastest moving characteristic (in either direction) is one with speed $u_f$.
More important than vehicles moving at most one cell between time steps is that \emph{characteristics cannot move more than one cell between time steps}, in either the upstream or downstream directions.
This condition is needed for stability of finite-difference approximations to partial differential equations, and corresponds to the \emph{Courant-Friedrich-Lewy condition}\index{Courant-Friedrich-Lewy condition} for explicit solution methods such as the cell transmission model.

Table~\ref{tbl:ctmexample} provides a demonstration of how the cell transmission model works.
In this example, a link is divided into three cells, implying that it takes three time steps for a vehicle to traverse the link at free-flow.
The fundamental diagram is such that at most 10 vehicles can move between cells in one time step, at most 30 vehicles can fit into one cell at jam density, and the ratio between the backward and forward characteristic speeds is $\delta = 2/3$.
There is a red light at the downstream end of the link which turns green at $t = 10$, and remains green thereafter.
In this table, $d(t)$ represents the number of vehicles that wish to enter the link during the $t$-th timestep (perhaps the sending flow from an upstream link), while $R(t)$ is the receiving flow for the link, calculated from~\eqn{ctmreceiving}.
The middle columns of the table show the main cell transmission model calculations: the number of vehicles in each cell at the start of each timestep, $N(x,t)$, and the number of vehicles moving into each cell during the $t$-th time step.
These values are calculated from~\eqn{conservationctm} and~\eqn{ctmtransitionflows}, respectively, along with the initial condition that the link is empty, that is, $N(x,0) = 0$ for all $x$.
The rightmost columns of the table show the link's sending flow, calculated from~\eqn{ctmsending}, and the actual flow which leaves the link, denoted $y(3,t)$.
This latter value is constrained to be zero as long as the light at the downstream end of the link is red.

Notice that the table has non-integer values: we do not need to round cell occupancies and flows to whole values, since the LWR model assumes vehicles are a continuously-divisible fluid.
Preserving non-integer values also ensures that the cell transmission model remains accurate no matter how small the timestep $\Delta t$ is (in fact, its accuracy should increase as this happens).
Insisting that flows and occupancies be rounded to whole numbers can introduce significant error if the timestep is small, unless one is careful with implementation.

Table~\ref{tbl:ctmcolorexample} shows only the cell occupancies at each timestep, color-coded according to the density in the cells.
In this example, the link is initially at free-flow, until the first vehicles encounter the red light and must stop.
A queue forms, and a shockwave begins moving backward.
When this shockwave reaches the cell at the upstream end of the link, the receiving flow of the link decreases, and the inflow to the link is limited.
When the light turns green, a second shockwave begins moving backward as the queue clears.
Once this shockwave overtakes the first, vehicles can begin entering the link again.
For a few time steps, the inflow is greater than $d(t)$, representing demand which was blocked when the receiving flow was restricted by the queue and which was itself queued on an upstream link (the ``queue spillback'' phenomenon).
Unlike the point queue and spatial queue link models, the cell transmission model tells us what is happening in the interior of a link, not just at the endpoints.
This is both a blessing and a curse: sometimes this additional information is helpful, while other times we may not be concerned with such details.
The link transmission model, described next, can simplify computations if we do not need information on the internal state of a link.
\index{link model!cell transmission model|)}

\begin{table}
\begin{center}
\caption{Cell transmission model example, on a link with three cells, $q_{max} \Delta t = 10$, $\bar{n} = 30$, and $\delta = \frac{2}{3}$. \label{tbl:ctmexample}}
\setlength{\tabcolsep}{2pt}
\begin{tabular}{|c|lc|cccccc|cc|}
\hline
  &  &  &  &Cell 0  &  &Cell 1  &  &Cell 2  &  & \\
$t$ & $d(t)$ {\footnotesize (*)} &$R(t)$  & $y(0,t)$  & $n(0,t)$  & $y(1,t)$  & $n(1,t)$  & $y(2,t)$  & $n(2,t)$  & $S(t)$ & $y(3,t)$  \\
\hline
0  &10   &10   &10   &0 &0 &0 &0 &0 &0 &0 \\
1  &10   &10   &10   &10   &10   &0 &0 &0 &0 &0 \\
2  &10   &10   &10   &10   &10   &10   &10   &0 &0 &0 \\
3  &10   &10   &10   &10   &10   &10   &10   &10   &10   &0 \\
4  &10   &10   &10   &10   &10   &10   &6.7  &20   &10   &0 \\
5  &9 &10   &9 &10   &10   &13.3 &2.2  &26.7 &10   &0 \\
6  &8 &10   &8 &9 &5.9  &21.1 &0.7  &28.9 &10   &0 \\
7  &7 &10   &7 &11.1 &2.5  &26.3 &0.2  &29.6 &10   &0 \\
8  &6 &9.6  &6 &15.6 &1 &28.5 &0.1  &29.9 &10   &0 \\
9  &5 &6.3  &5 &20.6 &0.4  &29.4 &0 &30   &10   &0 \\
10 &4 &3.2  &3.2  &25.2 &0.1  &29.8 &0 &30   &10   &10   \\
11 &3 (+0.8) &1.2  &1.2  &28.3 &0.1  &29.9 &6.7  &20   &10   &10   \\
12 &2 (+2.6) &0.4  &0.4  &29.4 &4.5  &23.3 &8.9  &16.7 &10   &10   \\
13 &1 (+4.2) &3.1  &3.1 &25.3 &7.4  &18.9 &9.6  &15.6 &10   &10   \\
14 &0 (+2.1) &6    &2.1 &21.0 &8.9  &16.7 &9.9  &15.2 &10   &10   \\
15 &0 &10   &0 &14.2   &9.5  &15.7 &10   &15.1 &10   &10   \\
16 &0 &10   &0 &4.7  &4.7  &15.3 &10   &15   &10   &10   \\
17 &0 &10   &0 &0 &0 &10  &10  &15   &10   &10   \\
18 &0 &10   &0 &0 &0 &0 &0 &15 &10   &10   \\
19 &0 &10   &0 &0 &0 &0 &0 &5  &5  &5  \\
20 &0 &10   &0 &0 &0 &0 &0 &0 &0 &0 \\   
\hline
\end{tabular}

{\footnotesize 
(*) Numbers in parentheses indicate unserved earlier demand waiting in queue.
}
\end{center}
\end{table}

\begin{table}
\begin{center}
\caption{Cell occupancies from the example in Table~\ref{tbl:ctmexample}, green is zero density and red is jam density.
\label{tbl:ctmcolorexample}}
\includegraphics[width=0.5\textwidth]{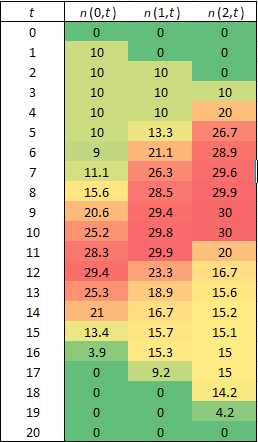}
\end{center}
\end{table}

\subsection{Link transmission model}

\index{link model!link transmission model|(}
The link transmission model allows us to calculate sending and receiving flows for links with any trapezoidal fundamental diagram.\footnote{In fact, it can be generalized to any piecewise linear fundamental diagram without too much difficulty; see Exercise~\ref{ex:ltmextension}.}  In contrast to the cell transmission model, it only involves calculations at the ends of the links --- details of what happen in the middle of the link are ignored.
As a consequence, the link transmission model does not require us to keep track of information within the link.
However, it does require us to keep track of information on the past state of the link, whereas the cell transmission model calculations only involve quantities at the current time step, and the previous time step.
The link transmission model can also overcome the ``shock spreading'' phenomenon, where backward-moving shockwaves in the cell transmission model can diffuse across multiple cells, even though they are crisp in the LWR model.
(See Exercise~\ref{ex:shockspreading}.)

Assume that the trapezoidal fundamental diagram is parameterized as in Figure~\ref{fig:trapezoidalfd}.
There are three characteristic speeds, $+u_f$ at free-flow, 0 at capacity flow, and $-w$ at congested flow.
As with the point queue and spatial queue models, we apply the Newell-Daganzo method to calculate sending and receiving flows.
We will only need to refer to cumulative counts $N$ at the upstream and downstream ends of the link, that is, at $N(0, \cdot)$ and $N(L, \cdot)$, respectively.
In keeping with the notation introduced in Section~\ref{sec:simplelinks}, we will refer to these as $N^\uparrow$ and $N^\downarrow$.
If there is no obstruction from a downstream link or node, then the end of the link will be uncongested, and the characteristic at this end will either have slope $+u_f$ or slope zero.
The Newell-Daganzo method thus gives
\labeleqn{ltmnewellsending}{N^\downarrow(t+\Delta t) = \min \{ N^\uparrow(t + \Delta t - L / u_f), N^\downarrow(t) + q_{max} \Delta t \}\, }
and the equation for the sending flow is obtained as the difference between $N^\downarrow(t+\Delta t)$ and $N^\downarrow(t)$:
\labeleqn{ltmsending}{S(t) = \min \{ N^\uparrow(t + \Delta t - L / u_f) - N^\downarrow(t), q_{max} \Delta t \} \,.}
\index{link model!link transmission model!sending flow}

For the receiving flow, we have to take into account the two relevant characteristic speeds of 0 and $-w$, since the receiving flow is calculated assuming an inflow large enough that the upstream end of the link is congested (or at least at capacity).
The stationary characteristic corresponds to the known point $N(0, t)$, while the backward-moving characteristic corresponds to the known point $N(L, t - L/w)$.
Thus, applying the last two terms of equation~\eqn{trapezoidalfd} would give
\labeleqn{ltmnewellreceiving}{N^\uparrow(t+\Delta t) = \min \{ N^\uparrow(t) + q_{max} \Delta t, N^\downarrow(t + \Delta t - L/w) + k_j L \} }
and\index{link model!link transmission model!receiving flow}
\labeleqn{ltmreceiving}{R(t) = \min \{ q_{max} \Delta t , N^\downarrow(t + \Delta t - L/w) + k_j L - N^\uparrow(t) \}  \,.
}
It is possible to show that equations~\eqn{ltmsending} and~\eqn{ltmreceiving} ensure that the number of vehicles on the link is always nonnegative, and less than $k_j L$.

The link transmission model is demonstrated on an example similar to the one used for the cell transmission model; the only difference is that the ratio of backward-to-forward characteristics has been adjusted from 2/3 to 3/4.
In particular, $L / u_f = 3 \Delta t$, and $L / w = 4 \Delta t$, so forward-moving characteristics require three time steps to cross the link, and backward-moving characteristics require four time steps.
The total number of vehicles which can fit on the link is $k_j L = 90$.
Otherwise, the example is the same: the demand profile is identical, and a red light prevents outflow from the link until $t = 10$.
The results of the calculations are shown in Table~\ref{tbl:ltmexample}.
The rightmost column shows the number of vehicles on the link, which is the difference between the upstream and downstream cumulative counts at any point in time.
Notice that inflow to the link is completely blocked during the 12th and 13th time intervals, even though the number of vehicles on the link is less than the jam density of 90.
This happens because the queue has started to clear at the downstream end, but the clearing shockwave has not yet reached the upstream end of the link.
The vehicles at the upstream end are still stopped, and no more vehicles can enter.
In contrast, the spatial queue model would allow vehicles to start entering the link as soon as they began to leave.
\index{link model!link transmission model|)}

\begin{table}
\begin{center}
\caption{Link transmission model example, with $L/u_f = 3 \Delta t$, $L/w = 4 \Delta t$, and $k_j L = 90$.
\label{tbl:ltmexample}}
\setlength{\tabcolsep}{2pt}
\begin{tabular}{|c|ccc|cc|ccc|}
\hline
$t$ & $d(t)$ &$R(t)$ & Inflow & $N^\uparrow(t)$ & $N^\downarrow(t)$ & $S(t)$ & Outflow & Vehicles on link \\
\hline
0  &10 &10   &10 &0 &0 &0 &0 &0\\
1  &10 &10   &10 &10   &0 &0 &0 &10\\
2  &10 &10   &10 &20   &0 &0 &0 &20\\
3  &10 &10   &10 &30   &0 &10   &0 &30\\
4  &10 &10   &10 &40   &0 &10   &0 &40\\
5  &9  &10   &9  &50   &0 &10   &0 &50\\
6  &8  &10   &8  &59   &0 &10   &0 &59\\
7  &7  &10   &7  &67   &0 &10   &0 &67\\
8  &6  &10   &6  &74   &0 &10   &0 &74\\
9  &5  &10   &5  &80   &0 &10   &0 &80\\
10 &4  &5 &4  &85   &0 &10   &10   &85\\
11 &3  &1 &1  &89   &10   &10   &10   &79\\
12 &2  &0 &0  &90   &20   &10   &10   &70\\
13 &1  &0 &0  &90   &30   &10   &10   &60\\
14 &0  &10   &5  &90   &40   &10   &10   &50\\
15 &0  &10   &0  &95   &50   &10   &10   &45\\
16 &0  &10   &0  &95   &60   &10   &10   &35\\
17 &0  &10   &0  &95   &70   &10   &10   &25\\
18 &0  &10   &0  &95   &80   &10   &10   &15\\
19 &0  &10   &0  &95   &90   &5 &5 &5\\
20 &0  &10   &0  &95   &95   &0 &0 &0\\
\hline
\end{tabular}
\end{center}
\end{table}

\subsection{Point and spatial queues, and the LWR model (*)}
\label{sec:pqsqlwr}

\emph{(This optional section shows how the previously-introduced point and spatial queue models can be seen as special cases of the LWR model.)}

\index{link model!point queue!relationship to LWR}
\index{link model!spatial queue!relationship to LWR}
The first link models introduced in this chapter were the point and spatial queue models, in Section~\ref{sec:simplelinks}.
These were presented as simple link models to illustrate concepts like the sending and receiving flow, rather than realistic depictions of traffic flow.
Nevertheless, it is possible to view the point and spatial queue models as special cases of the LWR model by making an appropriate choice of the fundamental diagram, as shown in this section.
Applying the Newell-Daganzo method with this fundamental diagram gives us a second way to derive the expressions for sending and receiving flow for these models.

The point queue model is equivalent to assuming that the flow-density relationship is as shown in Figure~\ref{fig:pqfundamentaldiagram}.
This diagram is unlike others we've seen, because there is no jam density.
This represents the idea that the point queue occupies no physical space: no matter how many vehicles are in queue, there is nothing to prevent additional vehicles from entering the link and joining the queue.
It is also the simplest diagram which we have seen so far, and is defined by only two parameters: the free-flow speed $u_f$ and the capacity $q_{max}$.
(Even the triangular fundamental diagram in Figure~\ref{fig:triangularfd} required a third parameter, either $-w$ or $k_j$.)

\stevefig{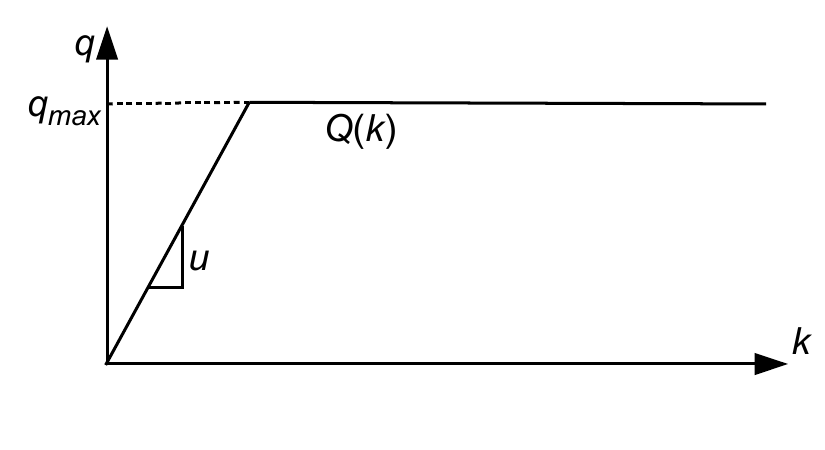}{Fundamental diagram in the point queue model.}{0.6\textwidth}

The Newell-Daganzo method leads to a simple expression for the sending and receiving flows in a point queue model.
To calculate the sending flow $S(t)$, we need to examine the downstream end of the link, so $x = L$.
Since we are solving in increasing order of time, we already know $N^\downarrow(0)$, $N^\downarrow(\Delta t), \ldots, N^\downarrow(t)$ (the number of vehicles which have left the link at each time interval).
Likewise, we know how many vehicles have entered the link at earlier points in time, so we know $N^\uparrow(0)$, $N^\uparrow(\Delta t), \ldots, N^\uparrow(t)$.
For the sending flow, we are assuming that there are no obstructions from downstream.
If this were the case, then we can calculate $N^\downarrow(t+ \Delta t)$ using the Newell-Daganzo method, and
\labeleqn{sendingnewell}{S(t) = N^\downarrow(t+\Delta t) - N^\downarrow(t)\,.}

\stevefig{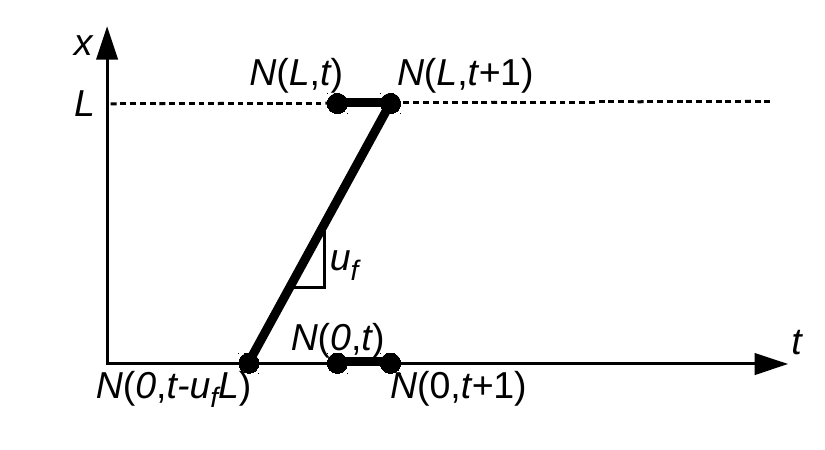}{Point queue model characteristics for sending and receiving flow.}{0.8\textwidth}

In the point queue model, there are two possible wave speeds, $u_f$ (corresponding to free-flow conditions) and zero (corresponding to flow moving at capacity).
Figure~\ref{fig:pqcharacteristic} shows how these characteristics can be traced back to known data, either at the upstream end $x = 0$, or at the downstream end $x = L$.
Therefore, $t_R = t$, $x_U = 0$, and $t_U = (t + \Delta t) - L / u_f$.
We can think of this as a special case of a ``trapezoidal'' diagram where the jam density $k_j$ is infinite.
Applying equation~\eqn{newelltrapezoidal}, we see that if $k_j = \infty$ then the minimum must occur in one of the first two terms, so
\labeleqn{pqnewellsending}{N^\downarrow(t+\Delta t) = \min \{ N^\uparrow(t + \Delta t - L / u_f), N^\downarrow(t) + q_{max} \Delta t\} }
and
\labeleqn{pqsending}{S(t) = \min \{ N^\uparrow(t + \Delta t - L / u_f) - N^\downarrow(t), q_{max} \Delta t\} \,.}
The two terms in the minimum in~\eqn{pqsending} correspond to the case when the queue is empty, and when there are vehicles in queue.
In the first term, since there is no queue, we just need to know how many vehicles will finish traversing the physical section of the link between $t$ and $t + \Delta t$; this is exactly the difference between the total number of vehicles which have entered by time $t + \Delta t - L / u_f$ and the total number that have left by time $t$.
When there is a queue, the vehicles exit the link at the full capacity rate.

In these expressions, it is possible that $(t + \Delta t) - L / u_f$ is not an integer, that is, it does not line up with one of the discretization points exactly.
In this case the most accurate choice is to interpolate between the known time points on either side (remember that we chose $\Delta t$ so that $L / u_f \geq 1$).
If you are willing to sacrifice some accuracy for efficiency, you can choose to round to the nearest integer, or to adjust the length of the link so that $L / u_f$ is an integer.

For the receiving flow, we look at the upstream end of the link.
We can treat the same points as known --- $N^\uparrow(0)$, $N^\uparrow(\Delta t), \ldots, N^\uparrow(t)$ and $N^\downarrow(0)$, $N^\downarrow(\Delta t)$, $\ldots$, $N^\downarrow(t)$.
\emph{Since the fundamental diagram for the point queue model has no decreasing portions, the known data at the downstream end can never be relevant.}  (A line connecting one of these points to the unknown point $N^\uparrow(t+ \Delta t)$ must have negative slope, see Figure~\ref{fig:pqcharacteristic}.)   Furthermore, for the receiving flow, the characteristic with positive slope $+u_f$, corresponding to upstream conditions, is irrelevant because we are assuming an unlimited number of vehicles are available to move from upstream --- and therefore its term in~\eqn{newelltrapezoidal} will never be the minimum.
We are only left with the middle term, corresponding to capacity, so 
\labeleqn{pqnewellreceiving}{N^\uparrow(t+\Delta t) = N^\uparrow(t) + q_{max} \Delta t}
and
\labeleqn{pqreceiving}{R(t) = q_{max} \Delta t \,.
}

\stevefig{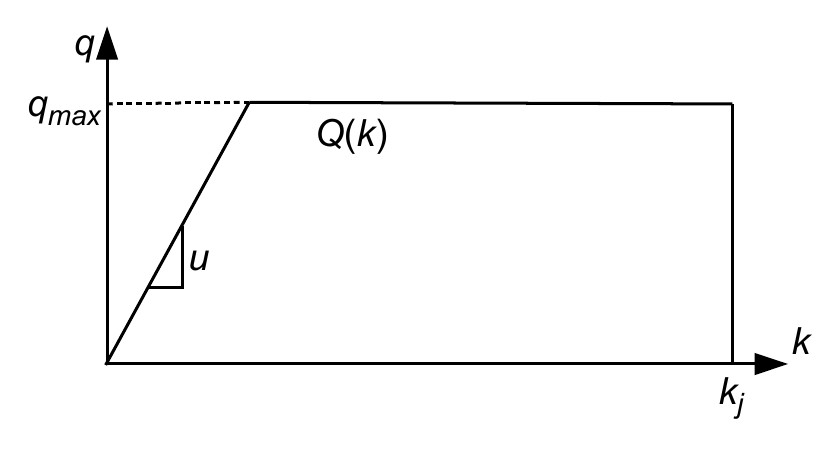}{Fundamental diagram in the spatial queue model.}{0.6\textwidth}

In terms of the fundamental diagram, the spatial queue model takes the form in Figure~\ref{fig:sqfundamentaldiagram}.
This diagram requires three parameters to calibrate: the free-flow speed $u_f$, the capacity $q_{max}$, and the jam density $k_j$.
Notice, however, that the fundamental diagram is discontinuous, and immediately drops from $q_{max}$ to zero once jam density is reached.
This implies that backward-moving shockwaves can have infinite speed in the spatial queue model --- a physical interpretation is that when vehicles at the front of the queue begin moving, vehicles at the rear of the queue immediately start moving as well.
In reality, there is a delay before vehicles at the rear of the queue begin moving, and this can be treated as an artifact arising from simplifying assumptions made in the spatial queue model.\footnote{Alternatively, connected and autonomous vehicles may be able to exhibit such behavior if an entire platoon of vehicles coordinates its acceleration.}

There are thus three possible characteristic speeds: $+u_f$ at free-flow, 0 at capacity flow, and $-\infty$ when the queue reaches jam density.
The Newell-Daganzo method is applied in much the same way as was done for the point queue model.
In particular, the sending flow expression is exactly the same, because the two characteristics with nonnegative velocity are the same.
We thus have
\labeleqn{sqnewellsending}{N(L,t+\Delta t) = \min \{ N^\uparrow(t + \Delta t - L / u_f), N^\downarrow(t) + q_{max} \Delta t \} }
and
\labeleqn{sqsending}{S(t) = \min \{ N^\uparrow(t + \Delta t - L / u_f) - N^\downarrow(t), q_{max} \Delta t \} \,.}

For the receiving flow, we have to take into account the new shockwave speed.
Dealing with an infinite speed can be tricky, since, taken literally, would mean that the upstream cumulative count $N^\uparrow(t + \Delta t)$ could depend on the downstream cumulative count at the same time $N^\downarrow(t + \Delta t)$.
Since we are solving the model in forward order of time, however, we do not know the value $N^\downarrow(t + \Delta t)$ when calculating $N^\uparrow(t + \Delta t)$.
In an acyclic network, we could simply do the calculations such that $N^\downarrow(t + \Delta t)$ is calculated first before $N^\uparrow(t + \Delta t)$, using the concept of a topological order.
In networks with cycles --- virtually all realistic traffic networks --- this will not work.
Instead, what is best is to approximate the ``infinite'' backward wave speed with one which is as large as possible, basing the calculation on the most recent known point $N^\downarrow(t)$ (Figure~\ref{fig:sqcharacteristic}).
Effectively, this replaces the infinite backward wave speed with one of speed $L / \Delta t$.
Equation~\eqn{newelltrapezoidal} thus gives
\labeleqn{sqnewellreceiving}{N^\uparrow(t+\Delta t) = \min \{ N^\uparrow(t) + q_{max} \Delta t, N^\downarrow(t) + k_j L \} }
and
\labeleqn{sqreceiving}{R(t) = \min \{ q_{max}\Delta t , (N^\downarrow(t) + k_j L) - N^\uparrow(t) \}  \,.
}
Equation~\eqn{sqreceiving} will ensure that the number of vehicles on the link will never exceed $k_j L$, assuming that this is true at time zero, as you are asked to show in Exercise~\ref{ex:sqcompliance}.

\stevefig{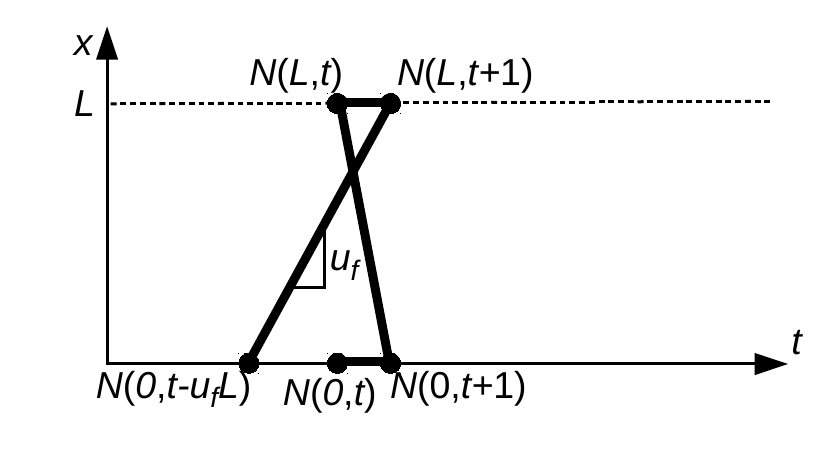}{Spatial queue model characteristics for sending and receiving flow; the upstream-moving wave is an approximation of the vertical component of the fundamental diagram.}{0.8\textwidth}

\subsection{Discussion}

\index{link model!point queue|(}
\index{link model!spatial queue|(}
\index{link model!cell transmission model|(}
\index{link model!link transmission model|(}
\index{link model!comparison|(}
This chapter has presented four different link models: point queues, spatial queues, the cell transmission model, and the link transmission model.
Although not initially presented this way, all four can be seen as special cases of the Lighthill-Whitham-Richards model.
The point and spatial queue models can be derived from particularly simple forms of the fundamental diagram (as well as from physical first principles, as in Section~\ref{sec:simplelinks}), while the cell transmission model and link transmission model are more general methods which can handle more sophisticated fundamental diagrams (typically triangular or trapezoidal in practice).
The cell transmission model directly solves the LWR system of partial differential equations by discretizing in space and time, and applying a finite-difference approximation.
The link transmission model is based on the Newell-Daganzo method.
The primary distinction between these methods is that the Newell-Daganzo method only requires tracking the cumulative counts $N$ at the upstream and downstream ends of each link in time, while the cell transmission model also requires tracking the number of vehicles $n$ at intermediate cells within the link.
However, the cell transmission model does not require storing any values from previous time steps, and can function entirely using the number of vehicles in each cell at the current time.
The Newell-Daganzo method requires that some past cumulative counts be stored, for the amount of time needed for a wave to travel from one end of the link to the other.
Which is more desirable depends on implementation details, and on the specific application context --- at times it may be useful to know the distribution of vehicles within a link (as the cell transmission model gives), while for other applications this may be an irrelevant detail.
One final advantage of the Newell-Daganzo method is that the values it gives are exact.
In the cell transmission model, backward-moving shockwaves will tend to ``spread out'' as a numerical issue involved in the discretization; this will not happen when applying the Newell-Daganzo method.
The exercises explore this issue in more detail.
On the other hand, the cell transmission model is easier to explain to decision-makers, and its equations have intuitive explanations in terms of vehicles moving within a link and the amount of available space.
The Newell-Daganzo method is a ``deeper'' method requiring knowledge of partial differential equations, and seems more difficult to convey to nontechnical audiences.

The point queue and spatial queue models have their places as well, despite their strict assumptions.
The major flaw in the point queue model, from the standpoint of realism, is its inability to model queue spillbacks which occur when links are full.
On the other hand, by ignoring this phenomenon, the point queue model is much more tractable, and is amenable even to closed-form expressions of delay and sensitivity to flows.
It is also more robust to errors in input data, because queue spillback can introduce discontinuities in the network loading.
There are cases where this simplicity and robustness may outweigh the (significant) loss in realism induced by ignoring spillbacks.
The spatial queue model can represent spillbacks, but will tend to underestimate its effect due to its assumption of infinitely-fast backward moving shockwaves.
Nevertheless, it can also lead to simpler analyses than the link transmission model.
\index{link model!point queue|)}
\index{link model!spatial queue|)}
\index{link model!cell transmission model|)}
\index{link model!link transmission model|)}
\index{link model!comparison|)}
\index{link model|)}
\index{point queue|see {link model, point queue}}
\index{spatial queue|see {link model, spatial queue}}
\index{cell transmission model|see {link model, cell transmission model}}
\index{link transmission model|see {link model, link transmission model}}
\index{merge model|see {node model, merge}}
\index{diverge model|see {node model, diverge}}

\section{Fancier Node Models} 
\label{sec:fancynodes}

\index{node model|(}
Section~\ref{sec:nodemodels} introduced node modeling concepts, and how the sending and receiving flows from adjacent links are mapped to transition flows showing how many vehicles move from each incoming link to each outgoing link in a single time step.
In analogy with Section~\ref{sec:fancylinks}, we now expand our discussion of node models beyond the simple intersection types presented thus far.
A general intersection can have any number of incoming and outgoing links.
General intersections can represent signal-controlled intersections, stop-controlled intersections, roundabouts, and so forth.
Node models for general intersections are less standardized, and less well-understood, than node models for merges and diverges, and representations of these can vary widely in dynamic traffic assignment implementations.
This section presents four alternative models for general intersections: the first for simple traffic signals, the second for intersections where drivers strictly take turns or have equal priority for all turning movements (such as an all-way stop), and the third where approaches have different priority levels and crossing conflicts must be considered, such as a two-way stop or signal with permitted phasing.
A fourth node model is presented which allows one to use arbitrary models for turning movement capacity and intersection delay.
As you read through this section, it would be helpful to think about how these models can be adapted or extended to represent more sophisticated signal types, roundabouts, and other types of junctions.
Many of the ideas in these models are based on the ideas described for merges and diverges.

\subsection{Basic signals}

\index{node model!basic signal|(}
For the purposes of this section, a ``basic signal'' is one which (1) only has protected phases (no permitted turns which must yield to oncoming traffic) and (2) all turning movements corresponding to the same approach move simultaneously (for instance, there are no turning lanes with separate phases from the through movement).
See Figure~\ref{fig:basicsignal} for an illustration.
In this case, the node can be modeled as a diverge intersection, where the ``upstream'' link varies over time, depending on which approach has the green indication.
Flows from other approaches (which have red indications) are set to zero.
Following the notation for diverges, we use $p_{hij}$ to reflect the proportion of the sending flow from approach $(h,i)$ which wishes to leave via link $(i,j)$: these values must be nonnegative, and $\sum_{(i,j):(h,i,j) \in \Xi(i)} p_{hij} = 1$ for all approaches $(h,i)$.
The algorithm is as follows:

\begin{enumerate}
\item Let $(h^*, i) \in \Gamma^{-1}(i)$ be the approach which has the green indication at the current time.
\item Calculate the fraction of flow which can move:
\labeleqn{basicfraction}{\phi = \min_{(i,j) : (h^*, i, j) \in \Xi(i)} \myc{ \frac{R_{ij}}{p_{h^*ij}S_{h^*i}} , 1 }\,.}
\item Calculate the transition flows for each turning movement:
\labeleqn{basictransition}{y_{hij} = \begin{cases}    \phi S_{hi} p_{hij} & \mbox{if } h = h^* \mbox{ and } (h, i, j) \in \Xi(i) \\
                                                    0                    & \mbox{otherwise}
                                   \end{cases}
}
\end{enumerate}

\stevefig{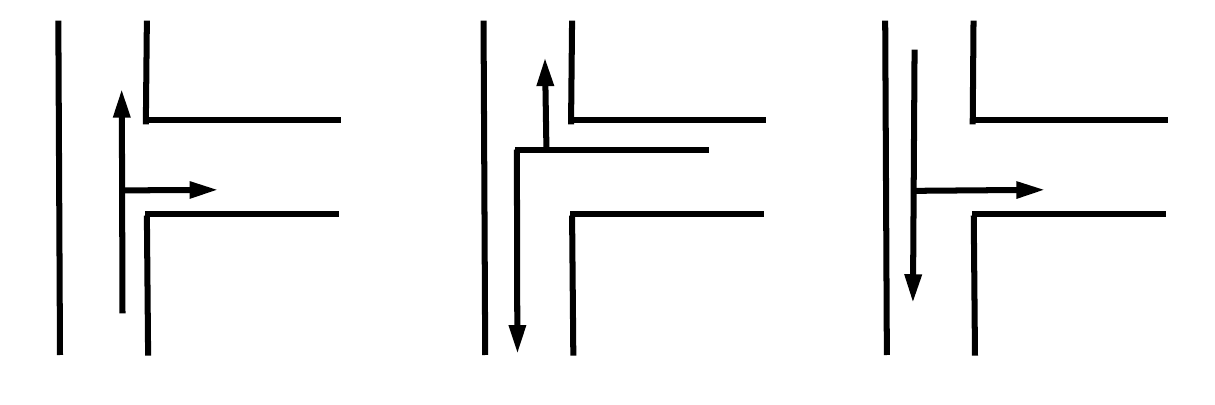}{Phasing plan for a basic signal.}{\textwidth}

In this implementation, one must be a little bit careful if the green times in the signal are not multiples of the time step $\Delta t$.
It is possible to round the green times so that they are multiples of $\Delta t$, but this approach can introduce considerable error over the analysis period: for instance, assume that $\Delta t$ is equal to six seconds, and a two-phase intersection has green times of 10 seconds and 14 seconds, respectively.
Rounding to multiples of the time step would give both phases twelve seconds each, which seems reasonable enough; but over a three-hour analysis period, the phases would receive 75 and 105 minutes of green time in reality, as compared to 90 minutes each in simulation.
In highly congested situations, this can introduce considerable error.
This issue can be avoided if, instead of rounding, one gives the green indication to the approach which would have green in reality at that time.
In the example above, the intersection has a cycle length of 24 seconds.
So, when $t = 60$ seconds, we are 12 seconds into the third cycle; and at this point the green indication should be given to the second phase.
In this way, there is no systematic bias introduced into the total green time each approach receives.
\index{node model!basic signal|)}

\subsection{Equal priority movements}

\index{node model!general, equal priority|(}
An all-way stop intersection is characterized by turn-taking: that is, vehicles have the opportunity to depart the intersection in the order in which they arrive.
No turning movement has priority over any other, but the turning movements from different approaches interact with each other and may compete for space on the same outgoing link.
This is different from the basic signal model, where the phasing scheme ensures that at most one approach is attempting to use an outgoing link at any given point of time.

Intersections with equal priority movements have characteristics of both diverges and merges.
Like a diverge, if a vehicle is unable to turn into a downstream link because of an obstruction, we assume that the vehicle obstructs all other vehicles from the same approach, respecting the FIFO principle.
This means that the outflows for all of the turning movements corresponding to any approach must follow the same proportions as the number of drivers \emph{wishing} to use all of these movements.
Similar to a merge, we assume that if there are high sending flows from all the approaches, the fraction of the receiving flow allocated to each approach is divided up proportionally.
However, instead of allocating the receiving flow $R_{ij}$ to approach $(h,i)$ based on the full capacity $q_{max}^{hi}$, we instead divide up the receiving flow based on the \emph{oriented capacity}\index{capacity!oriented}\index{oriented capacity|see {capacity, oriented}}\label{not:qmaxhij}
\labeleqn{orientedcapacity}{q_{max}^{hij} = q_{max}^{hi} p_{hij}\,,}
where $p_{hij}$ is the proportion of the flow from approach $(h,i)$ which wishes to exit on link $(i,j)$.
(If turning movement $[h,i,j]$ is not in the allowable set $\Xi(i)$, then $p_{hij} = 0$.)

Multiplying the capacity by this proportion reflects the fact that an upstream approach can only make use of an available space on a downstream link if there is a vehicle wishing to turn.
In Figure~\ref{fig:orientedcapacity}(a), each incoming link uses a unique exiting link, and thus can claim its full capacity.
In Figure~\ref{fig:orientedcapacity}(b), half the vehicles on link (4,2) want to turn right and half wish to go straight, whereas all the vehicles on link (3,2) wish to turn right.
Link (3,2) therefore has twice as many opportunities to fill available space on link (2,1), and thus its rightful share is twice that of link (4,2).
For any two approaches $[h,i,j]$ and $[h',i,j]$ using the same outgoing link, we thus require that
\labeleqn{orientedreceivingflow}{\frac{y_{hij}}{y_{h'ij}} = \frac{q_{max}^{hij}}{q_{max}^{h'ij}}\,.}
\emph{assuming that both approaches are fully competing for the link $(i,j)$}.
If an approach has a small sending flow, it may use less of its assigned receiving flow than equation~\eqn{orientedreceivingflow} allocates, and this unused receiving flow may be used by other approaches.

\stevefig{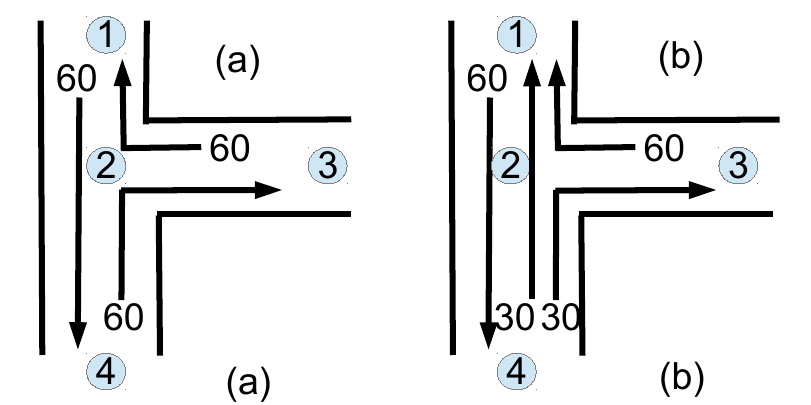}{Oriented capacities for a three-leg intersection, where all approaches have $q_{max} = 60$. (a) All approaches map to a unique outgoing link. (b) The flow on approach (4,2) is split between two outgoing links.}{0.6\textwidth}

As in Section~\ref{sec:diverge}, we use the oriented sending flows\index{sending flow!oriented}
\labeleqn{orientedsendingflow}{S_{hij} = S_{hi} p_{hij}}
to reflect the number of vehicles that wish to use the turning movement $[h,i,j]$.
Following the same principles as the merge model, if the oriented sending flow from an approach is less than its proportionate share of a downstream link's receiving flow, its unused share will be divided among the other approaches with unserved sending flow still remaining, in proportion to their oriented capacities.
If the oriented sending flow for a turning movement is greater than the oriented receiving flow for that movement, then by the FIFO principle applied to diverges, it will restrict flow to all other downstream links by the same proportion, and for any two turning movements $(h,i,j)$ and $(h,i,j')$ from the same approach we must have
\labeleqn{proportionalflow}{\frac{y_{hij}}{y_{hij'}} = \frac{S_{hij}}{S_{hij'}} = \frac{p_{hij}}{p_{hij'}}\,.
}
We can rearrange this equation to show that the ratio $y_{hij} / S_{hij}$ (the ratio of actual flow and desired flow for any turning movement) is uniform for all the turning movements approaching from link $(h,i)$ --- this ratio plays the same role as $\phi$ in a diverge.

The presence of multiple incoming and outgoing links causes another complication, in that the flows between approaches are all linked together.
If an approach is restricted by the receiving flow of a downstream link, flow from that approach is restricted to \emph{all} other downstream links.
This means that the approach may not fully consume its ``rightful share'' of another downstream link, thereby freeing up additional capacity for a different approach.
Therefore, we cannot treat the approaches or downstream links separately or even sequentially in a fixed order, because we do not know \emph{a priori} how these will be linked together.

However, there is an algorithm which generates a consistent solution despite these mutual dependencies.
In this algorithm, each approach link can be \emph{sending-constrained}, or \emph{receiving-constrained by a downstream link}.
If an approach is sending-constrained, its oriented sending flow to all downstream links is less than its rightful share, and therefore all of the sending flow can move.
If an approach is receiving-constrained by link $(i,j)$, then the approach is unable to move all of its sending flow, and the fraction which can move is dictated by link $(i,j)$.
(That is, receiving flow on $(i,j)$ is the most restrictive constraint for the approach).
The algorithm must determine which links are sending-constrained, and which are receiving-constrained by a downstream link.

To find such a solution, we define two sets of auxiliary variables, $\tilde{S}_{hij}$ to reflect the amount of unallocated sending flow for movement $[h,i,j]$, and $\tilde{R}_{ij}$ to reflect the amount of unallocated receiving flow for outgoing link $(i,j)$.
These are initialized to the oriented sending flows and link receiving flows, and reduced iteratively as flows are assigned and the available sending and receiving flows are used up.
The algorithm also uses the notion of \emph{active} turning movements; these are turning movements whose flows can still be increased.
A turning movement $[h,i,j]$ becomes inactive either when $\tilde{S}_{hij}$ drops to zero (all vehicles that wish to turn have been assigned), or when $\tilde{R}_{ij'}$ drops to zero for \emph{any} outgoing link $(i,j')$ that approach $(h,i)$ is using (that is, for which $p_{hij'} > 0$).
Allocating all of the receiving flow for one outgoing link can thus impact flow on turning movements which use other outgoing links, because of the principle that vehicles wishing to turn will block others, as expressed in equation~\eqn{proportionalflow}.
The set of active turning movements will be denoted by $\Omega$; a turning movement remains active until we have determined whether it is sending-constrained or receiving-constrained.

At each stage of the algorithm, we will increase the flows for all active turning movements.
We must increase these flows in a way which is consistent both with the turning fractions~\eqn{proportionalflow}, and with the division of receiving flow for outgoing links given by equation~\eqn{orientedreceivingflow}, and we will use $\alpha_{hij}$ to reflect the rate of increase for turning movement $(h,i,j)$.
The absolute values of these $\alpha$ values do not matter, only their proportions, so you can scale them in whatever way is most convenient to you.
Often it is easiest to pick one turning movement $(h,i,j)$ and fix its $\alpha$ value either to one, or to its oriented sending flow.
The turning proportions from $(h,i)$ then fix the $\alpha$ values for all other turning movements from the same approach.
You can then use equation~\eqn{orientedreceivingflow} to determine $\alpha$ values for turning movements competing for the same outgoing link, then use the turning fractions for the upstream link on those turning movements, and so on until a consistent set of $\alpha$ values has been determined.

The algorithm then increases the active turning movement flows in these proportions until some movement becomes inactive, because its sending flow or the receiving flow on its outgoing link becomes exhausted.
The process is then repeated with the smaller set of turning movements which remain active, and continues until all possible flows have been assigned.
This algorithm is a bit more involved than the node models seen thus far, and you may find it helpful to follow the example below as you read through the algorithm steps.

\begin{enumerate}
\item Initialize by calculating oriented capacities and sending flows using equations~\eqn{orientedcapacity} and~\eqn{orientedsendingflow}; by setting $y_{hij} \leftarrow 0$ for all $[h,i,j] \in \Xi(i)$; by setting $\tilde{S}_{hij} \leftarrow S_{hij}$ for all $[h,i,j] \in \Xi(i)$ and $\tilde{R}_{ij} \leftarrow R_{ij}$ for all $(i,j) \in \Gamma(i)$; and by declaring active all turning movements with positive sending flow: $\Omega \leftarrow \myc{ [h,i,j] \in \Xi(i) : S_{hij} > 0 }$.
\item Identify a set of $\alpha_{hij}$ values for all active turning movements which is consistent with the turning fractions ($\alpha_{hij} / \alpha_{hij'} = p_{hij} / p_{hij'}$ for all of the turning movements from the same incoming link) and oriented capacities ($\alpha_{hij} / \alpha_{h'ij} = q_{max}^{hij} / q_{max}^{h'ij}$ for all of the turning movements to the same outgoing link).
\item For each outgoing link $(i,j) \in \Gamma(i)$, identify the rate $\alpha_{ij}$ at which its receiving flow will be reduced, by adding $\alpha_{hij}$ for all active turning movements whose outgoing link is $(i,j)$: $\alpha_{ij} \leftarrow \sum_{[h,i,j] \in \Omega} \alpha_{hij}$.
\item Determine the point at which some turning movement will become inactive, by calculating the largest possible step size
\labeleqn{inmstepsize}{\theta = \min \myc{ \min_{[h,i,j] \in \Omega} \myc{\frac{\tilde{S}_{hij}}{\alpha_{hij}}}, \min_{(i,j) \in \Gamma(i) : \alpha_{ij} > 0} \myc{\frac{\tilde{R}_{ij}}{\alpha_{ij}} }}\,.}
\item Increase flows for active turning movements, and update unallocated sending and receiving flows: for all $[h,i,j] \in \Omega$ update $y_{hij} \leftarrow y_{hij} + \theta \alpha_{hij}$, $\tilde{S}_{hij} \leftarrow \tilde{S}_{hij} - \theta \alpha_{hij}$, and $\tilde{R}_{ij} \leftarrow \tilde{R}_{ij} - \theta \alpha_{hij}$.
\item Update the set of active turning movements, by removing from $\Omega$ any turning movement $[h,i,j]$ for which $\tilde{S}_{hij} = 0$ or for which $\tilde{R}_{ij'} = 0$ for \emph{any} $(i,j') \in \Gamma(i)$ which is being used ($p_{hij'} > 0$).
\item If there are any turning movements which are still active ($\Omega \neq \emptyset$), return to step 3.
Otherwise, stop.
\end{enumerate}

\stevefig{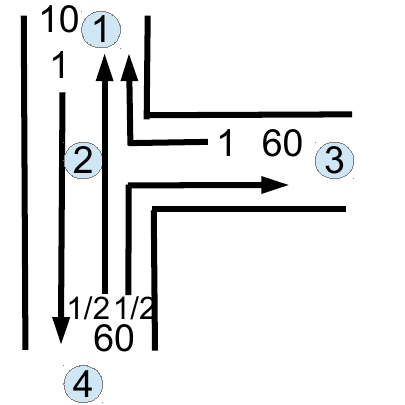}{Example of equal priorities algorithm, showing sending flows and turning proportions.
All links have capacity of 60.}{0.3\textwidth}

As a demonstration, consider the intersection in Figure~\ref{fig:allwaystopexample}, where all links (incoming and outgoing) have the same capacity of 60 vehicles per time step, and the sending flows $S_{hi}$ and turning proportions $p_{hij}$ are shown.
None of the downstream links is congested, so their receiving flows are equal to the capacity.
(In the figure, two sets of numbers are shown for each approach; the ``upstream'' number is the sending flow and the ``downstream'' number(s) are the proportions.)  The oriented capacities can be seen in Figure~\ref{fig:orientedcapacity}(b).
For this example, the six turning movements will be indexed in the following order:
\labeleqn{turnmovementorder}{\Xi(2) = \myc{[1,2,3], [1,2,4], [3,2,1], [3,2,4], [4,2,1], [4,2,3]}}
All vectors referring to turning movements will use this ordering for their components.

Step 1 of the algorithm initializes the oriented sending flows using equation~\eqn{orientedsendingflow},
\labeleqn{exampleorientedsendingflow}{\mb{S} = \vect{ S_{hij}} = \vect{0 & 10 & 60 & 0 & 30 & 30}\,,}
and the oriented capacities using equation~\eqn{orientedcapacity}
\labeleqn{exampleorientedcap}{\mb{q_{max}} = \vect{ q_{max}^{hij}} = \vect{ 0 & 60 & 60 & 0 & 30 & 30 }\,.}
The step also initializes the turning movement flows and auxiliary variables:
\labeleqn{turnmovementinit}{\mb{y} \leftarrow \vect{ 0 & 0 & 0 & 0 & 0 & 0 }, } 
\labeleqn{auxsendingflow}{\mb{\tilde{S}} \leftarrow \vect{{S}_{hij}} = \vect{0 & 10 & 60 & 0 & 30 & 30}\,,}
and
\labeleqn{auxreceivingflow}{\mb{\tilde{R}} = \vect{\tilde{R}_{21} & \tilde{R}_{23} & \tilde{R}_{24}} = \vect{60 & 60 & 60}\,.}
The set of active turning movements is $\Omega = \myc{[1,2,4], [3,2,1], [4,2,1], [4,2,3]}$.

Step 2 of the algorithm involves calculation of a set of consistent $\alpha_{hij}$ values.
One way of doing this is to start by setting $\alpha_{423} \leftarrow S_{423} = 30$.
The turning fractions from $(4,2)$ then require that $\alpha_{421} = \alpha_{423} = 30$.
The allocation rule for outgoing link $(2,1)$ then forces $\alpha_{321} = 60$: the oriented capacity for $[3,2,1]$ is twice that of $[4,2,1]$, and the $\alpha$ values must follow the same proportion.
Turning movement $[1,2,4]$ is independent of all of the other turning movements considered thus far, so we can choose its $\alpha$ value arbitrarily; say, $\alpha_{124} \leftarrow S_{124} = 10$.
(You should experiment around with different ways of calculating these $\alpha$ values, and convince yourself that the final flows are the same as long as the proportions of $\alpha$ values for interdependent turning movements are the same.)  We thus have
\labeleqn{flowincrements}{\bm{\alpha} = \vect{ \alpha_{hij}} = \vect{0 & 10 & 60 & 0 & 30 & 30}\,.}
The $\alpha$ values for inactive turning movements have been set to zero for clarity; their actual value is irrelevant because they will not be used in any of the steps that follow.

With these flow increments, flow on outgoing links $(2,1)$, $(2,3)$, and $(2,4)$ will be $\alpha_{21} = 90$, $\alpha_{23} = 30$, and $\alpha_{24} = 10$, as dictated by Step 3.

In Step 4, we determine how much we can increase the flow at the rates given by $\bm{\alpha}$ until some movement becomes inactive.
We have
\labeleqn{exampletheta}{\theta = \min \kbordermatrix { & \tilde{S}_{124}&\tilde{S}_{321}&\tilde{S}_{421}&\tilde{S}_{423}&\tilde{R}_{21} & \tilde{R}_{23} & \tilde{R}_{24} \\ 
                                                       & \frac{10}{10}, & \frac{60}{60},& \frac{30}{30},& \frac{30}{30},& \frac{60}{90},& \frac{60}{30},& \frac{60}{10} }}
or $\theta = 2/3$.

We can now adjust the flows, as in Step 5.
We increase the flow on each active turning movement by $\frac{2}{3} \alpha_{hij}$, giving
\labeleqn{newflows}{\mb{y} \leftarrow \vect{ 0 & 6 \frac{2}{3} & 40 & 0 & 20 & 20 }\,.}
We subtract these flow increments from the auxiliary sending and receiving flows, giving
\labeleqn{newauxsending}{\mb{\tilde{S}} \leftarrow \vect{ 0 & 3 \frac{1}{3} & 20 & 0 & 10 & 10}}
and
\labeleqn{newauxreceiving}{\mb{\tilde{R}} \leftarrow \vect{ 0 & 40 & 53 \frac{1}{3}}\,.}

Step 6 updates the set of active turning movements.
With the new $\mb{\tilde{S}}$ and $\mb{\tilde{R}}$ values, we see that $[3,2,1]$ and $[4,2,1]$ have become inactive, since there is no remaining receiving flow on link $(2,1)$.
Furthermore, this inactivates movement $[4,2,3]$: even though there are still travelers that wish to turn in this direction, and space on the downstream link ($\tilde{S}_{423}$ and $\tilde{R}_{23}$ are still positive), they are blocked by travelers waiting to use movement $[4,2,1]$.
So, there is only one active movement remaining, $\Omega \leftarrow \myc{[1,2,4]}$, and we must return to step 3.

In Step 3, we must recalculate the $\alpha_{ij}$ values because some of the turning movements are inactive.
With the new set $\Omega$, we have $\alpha_{21} = \alpha_{23} = 0$ and $\alpha_{24} = 10$.
The new step size is 
\labeleqn{newtheta}{\theta = \min \kbordermatrix{ &\tilde{S}_{124}& \tilde{R}_{24} \\
                                                  & \frac{3 \frac{1}{3}}{10},& \frac{53 \frac{1}{3}}{10}} = \frac{1}{3}\,.}
We then increase the flows, increasing $y_{124}$ by $10 \times \frac{1}{3}$ to 10, decreasing $\tilde{S}_{124}$ to zero, and decreasing $\tilde{R}_{24}$ to 50.
This change inactivates movement $[1,2,4]$.
Since there are no more active turning movements, the algorithm terminates, and the final vector of flows is
\labeleqn{turnmovementfinal}{\mb{y} = \vect{0 & 10 & 40 & 0 & 20 & 20}\,.}
\index{node model!general, equal priority|)}

\subsection{Intersections with priority}

\index{node model!general, unequal priority|(}
Intersections which allow crossing conflicts are more complex to model than the intersection types described above.
These include intersections with stop control only on some of the approaches, but not all, or signalized intersections with permitted movements that must yield to another traffic stream.
One way to model this type of intersection is to introduce a set $C$ of conflict points.
Like outgoing links, we associate a receiving flow $R_c$ indicating the maximum number of vehicles which can pass through this conflict point during a single time step $\Delta t$.
For each conflict point $c \in C$, let $\Xi(c)$ denote the turning movements which make use of this conflict point.
Then, for each turning movement$[h,i,j] \in \Xi(c)$, we define a (strictly positive) priority parameter $\beta^c_{hij}$.
These priority parameters are interpreted through their ratios: the share of the receiving flow allocated to turning movement $[h,i,j]$ relative to that allocated to turning movement $[h',i,j']$ is in the same proportion as the ratio
\labeleqn{priorityratio}{\beta^c_{hij} p_{hij} /\beta^c_{h'ij'} p_{h'ij'} \,.}

These play an analogous role to the oriented capacities defined in the previous subsection, in suggesting how the conflict point receiving flow should be divided among the competing approaches.
They are used more generally, however, to reflect priority rules.
For instance, a through movement often has priority over a turning movement which crosses it.
Even if the capacity and sending flows of the through lane and turn lane are the same, the through lane should have access to a greater share of the crossing point's receiving flow.
If, at full saturation, ten through vehicles move for every turning vehicle, then the $\beta$ value for the through movement should be ten times the $\beta$ value for the turning movement.

The node model is simplified if we assume that the $\beta$ values are all strictly positive.
You might find this unrealistic: in the example above, if the turning movement must yield to the through movement, then at full saturation perhaps \emph{no} turning vehicles could move.
In practice, however, priority rules are not strictly obeyed as traffic flows near saturation.
Polite through drivers may stop to let turning drivers move, or aggressive turning drivers might force their way into the through stream.
(Think of what would happen at a congested freeway if vehicles merging from an onramp took the ``yield'' sign literally!)  The requirement of strictly positive $\beta$ values thus has some practical merits, as well as mathematical ones.
The exercises explore ways to generalize this node model, including strict priority, and cases where different turning movements may consume different amounts of the receiving flow (for instance, if they are moving at different speeds).

We can now adapt the algorithm for equal priority for the case of intersections with different priorities.
The algorithm is augmented by adding receiving flows $R_c$ and auxiliary receiving flows $\tilde{R}_c$ for each conflict point, and we extend some of the computations to include the set of conflict points.
First, we require that the ratio of $\alpha$ values for two turning movements using the same crossing point follow the ratio of the $\beta$ values, assuming saturated conditions:
\labeleqn{orientedreceivingflowpartial}{\frac{\alpha_{hij}}{\alpha_{h'ij'}} = \frac{\beta^c_{hij} p_{hij}}{\beta^c_{h'ij'} p_{h'ij'}} \qquad \forall c \in C; [h,i,j], [h',i,j'] \in \Xi(c)\,.}
Note that the $\beta$ values are multiplied by the relevant oriented sending flows for the movements.
As with the oriented capacity, this respects the fact that the more flow is attempting to turn in a particular direction, the more opportunities or gaps will be available for it to claim.

Second, we must calculate the inflow rates to conflict points given the $\alpha$ values from active turning movements:
\labeleqn{conflictpointflow}{\alpha_c = \sum_{[h,i,j] \in \Xi(c) \cap \Omega} \alpha_{hij}\,.}
Third, the calculation of the step size $\theta$ must now include obstructions from conflict points:
\labeleqn{inmstepsize2}{\theta = \min \myc{ \min_{[h,i,j] \in \Omega} \myc{\frac{\tilde{S}_{hij}}{\alpha_{hij}}}, \min_{(i,j) \in \Gamma(i) : \alpha_{ij} > 0} \myc{\frac{\tilde{R}_{ij}}{\alpha_{ij}}}, \min_{c \in C : \alpha_{c} > 0} \myc{\frac{\tilde{R}_{c}}{\alpha_{c}} }}\,.}

With these modifications, the algorithm proceeds in the same way as before.
Specifically,
\begin{enumerate}
\item Initialize by calculating oriented capacities and sending flows using equations~\eqn{orientedcapacity} and~\eqn{orientedsendingflow}; by setting $y_{hij} \leftarrow 0$ for all $[h,i,j] \in \Xi(i)$; by setting $\tilde{S}_{hij} \leftarrow S_{hij}$ for all $[h,i,j] \in \Xi(i)$, $\tilde{R}_{ij} \leftarrow R_{ij}$ for all $(i,j) \in \Gamma(i)$, and $\tilde{R}_c \leftarrow R_c$ for each $c \in C$; and by declaring active all turning movements with positive sending flow: $\Omega \leftarrow \myc{ [h,i,j] \in \Xi(i) : S_{hij} > 0 }$.
\item Identify a set of $\alpha_{hij}$ values for all active turning movements which is consistent with the turning fractions ($\alpha_{hij} / \alpha_{hij'} = p_{hij} / p_{hij'}$ for all turning movements from the same incoming link), oriented capacities ($\alpha_{hij} / \alpha_{h'ij} = q_{max}^{hij} / q_{max}^{h'ij}$ for all turning movements to the same outgoing link), and conflict points based on equation~\eqn{orientedreceivingflowpartial}.
\item For each outgoing link $(i,j) \in \Gamma(i)$ and conflict point $c \in C$, identify the rate $\alpha_{ij}$ at which its receiving flow will be reduced, by adding $\alpha_{hij}$ for all active turning movements whose outgoing link is $(i,j)$: $\alpha_{ij} \leftarrow \sum_{[h,i,j] \in \Omega} \alpha_{hij}$ for links, and equation~\eqn{conflictpointflow} for conflict points.
\item Determine the point at which some turning movement will become inactive, by calculating $\theta$ using equation~\eqn{inmstepsize2}.
\item Increase flows for active turning movements, and update unallocated sending and receiving flows: for all $[h,i,j] \in \Omega$ update $y_{hij} \leftarrow y_{hij} + \theta \alpha_{hij}$, $\tilde{S}_{hij} \leftarrow \tilde{S}_{hij} - \theta \alpha_{hij}$; for all $(i,j) \in \Gamma(i)$ update $\tilde{R}_{ij} \leftarrow \tilde{R}_{ij} - \theta \alpha_{ij}$, and for all $c \in C$ update $\tilde{R}_c \leftarrow \tilde{R}_c - \theta \alpha_c$.
\item Update the set of active turning movements, by removing from $\Omega$ any turning movement $[h,i,j]$ for which $\tilde{S}_{hij} = 0$, for which $\tilde{R}_{ij'} = 0$ for \emph{any} $(i,j') \in \Gamma(i)$ which is being used ($p_{hij'} > 0$), or for which $\tilde{R}_c = 0$ for \emph{any} movement $[h,i,j'] \in \Xi(c)$ and conflict point $c$.
\item If there are any turning movements which are still active ($\Omega \neq \emptyset$), return to step 3.
Otherwise, stop.
\end{enumerate}

\stevefig{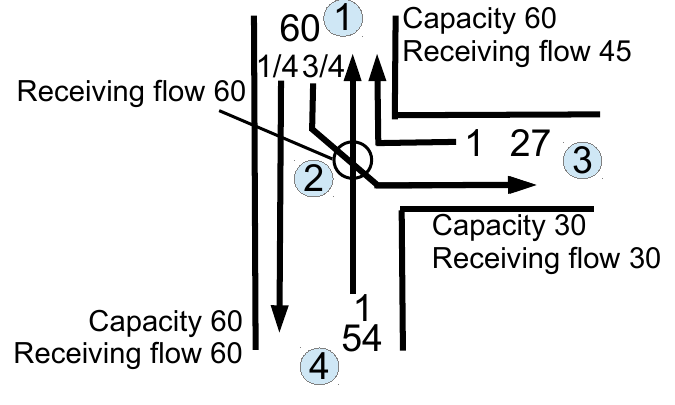}{Example of intersection priority algorithm, showing sending flows and turning proportions. Link capacity and receiving flows shown.}{0.6\textwidth}

As a demonstration, consider the intersection in Figure~\ref{fig:twowaystopexample}, where all links (incoming and outgoing) have the same capacity of 60 vehicles per time step, and the sending flows $S_{hi}$ and turning proportions $p_{hij}$ are shown.
(In the figure, two sets of numbers are shown for each approach; the ``upstream'' number is the sending flow and the ``downstream'' number(s) are the proportions.)    Links $(2,4)$ and $(2,3)$ have receiving flows of 60 vehicles, while link $(2,1)$ has a receiving flow of only 45 vehicles.
There is one conflict point, indexed $c$, which is marked with a circle in Figure~\ref{fig:twowaystopexample}).
Conflict point $c$ has a receiving flow of 60, and the turning movement $[1,2,3]$ must yield to the through movement $[4,2,1]$, as reflected by the ratio $\beta_{421}/\beta_{123} = 5$.
For this example, the six turning movements will be indexed in the following order:
\labeleqn{turnmovementorder2}{\Xi(2) = \myc{[1,2,3], [1,2,4], [3,2,1], [3,2,4], [4,2,1], [4,2,3]}\,.}
All vectors referring to turning movements will use this ordering for their components.

Step 1 of the algorithm initializes the oriented sending flows using equation~\eqn{orientedsendingflow},
\labeleqn{orientedsendingflow2}{\mb{S} = \vect{ S_{hij}} = \vect{45 & 15 & 27 & 0 & 54 & 0}\,,}
and the oriented capacities using equation~\eqn{orientedcapacity}
\labeleqn{exampleorientedcap2}{\mb{q_{max}} = \vect{ q_{max}^{hij}} = \vect{ 30 & 30 & 30 & 0 & 60 & 0 }\,.}
The step also initializes the turning movement flows and auxiliary variables:
\labeleqn{turnmovementinit2}{\mb{y} \leftarrow \vect{ 0 & 0 & 0 & 0 & 0 & 0 }, } 
\labeleqn{auxsendingflow2}{\mb{\tilde{S}} \leftarrow \vect{{S}_{hij}} = \vect{45 & 15 & 27 & 0 & 54 & 0}\,,}
and
\labeleqn{auxreceivingflow2}{\mb{\tilde{R}} = \vect{\tilde{R}_{21} & \tilde{R}_{23} & \tilde{R}_{24} & \tilde{R}_c } = \vect{45 & 30 & 60 & 60}\,.}
The set of active turning movements is $\Omega = \myc{[1,2,3], [1,2,4], [3,2,1], [4,2,1]}$.

Step 2 of the algorithm involves calculation of a set of consistent $\alpha_{hij}$ values.
One way of doing this is to start by setting $\alpha_{124} \leftarrow 1$ (again, this choice is arbitrary, and any positive number would work).
The turning fractions from $(1,2)$ then require that $\alpha_{123} = 3\alpha_{124} = 3$.
The allocation rule~\eqn{orientedreceivingflowpartial} for conflict point $c$ forces $\alpha_{421} = 20$: the $\beta$ value for $[4,2,1]$ is five times that of the $\beta$ value for $[1,2,3]$, and the turning fraction for $[4,2,1]$ is a third higher.
Since $\alpha_{421} = 20$, we must have $\alpha_{321} = 10$ because the oriented capacity of $[3,2,1]$ is half that of $[4,2,1]$.
This gives the flow increments
\labeleqn{flow2increments}{\bm{\alpha} = \vect{ \alpha_{hij}} = \vect{3 & 1 & 10 & 0 & 20  & 0}\,.}

With these flow increments, flow on outgoing links $(2,1)$, $(2,3)$, $(2,4)$, and the conflict point $c$ will be $\alpha_{21} = 30$, $\alpha_{23} = 3$, $\alpha_{24} = 1$, and $\alpha_c = 23$, as dictated by Step 3.

In Step 4, we determine how much we can increase the flow at the rates given by $\bm{\alpha}$ until some movement becomes inactive.
We have
\labeleqn{example2theta}{\theta = \min \kbordermatrix { & \tilde{S}_{123}&\tilde{S}_{124}&\tilde{S}_{321}&\tilde{S}_{421}&\tilde{R}_{21} & \tilde{R}_{23} & \tilde{R}_{24} & \tilde{R}_c\\ 
& \frac{15}{1}, & \frac{45}{3},& \frac{27}{10},& \frac{54}{20},& \frac{45}{30},& \frac{60}{3},& \frac{60}{1},& \frac{60}{23}}}
or $\theta = \frac{3}{2}$.

We can now adjust the flows, as in Step 5.
We increase the flow on each active turning movement by $\frac{3}{2} \alpha_{hij}$, giving
\labeleqn{newflows2}{\mb{y} \leftarrow \vect{ 4 \frac{1}{2} & 1 \frac{1}{2} & 15 & 0 & 30 & 0 }\,.}
We subtract these flow increments from the auxiliary sending and receiving flows, giving
\labeleqn{newauxsending2}{\mb{\tilde{S}} \leftarrow \vect{ 40 \frac{1}{2} & 13 \frac{1}{2} & 12 & 0 & 24 & 0}}
and
\labeleqn{newauxreceiving2}{\mb{\tilde{R}} \leftarrow \vect{ 0 & 25 \frac{1}{2} & 58 \frac{1}{2} & 25 \frac{1}{2} }\,.}

Step 6 updates the set of active turning movements.
With the new $\mb{\tilde{S}}$ and $\mb{\tilde{R}}$ values, we see that $[3,2,1]$ and $[4,2,1]$ have become inactive, since there is no remaining receiving flow on link $(2,1)$.
So, there are two active movements remaining, $\Omega \leftarrow \myc{[1,2,3], [1,2,4]}$, and we must return to step 3.

In Step 3, we must recalculate the $\alpha_{hij}$ and $\alpha_{ij}$ values because some of the turning movements are inactive.
With the new set $\Omega$, and starting with $\alpha_{124} = 1$, we compute $\alpha_{123} = 3$, and therefore  $\alpha_{21} = 0$, $\alpha_{23} = 3$, $\alpha_{24} = 1$, and $\alpha_c = 3$.
The new step size is 
\labeleqn{newtheta2}{\theta = \min \kbordermatrix{ &\tilde{S}_{123} &\tilde{S}_{124} & \tilde{R}_{23} & \tilde{R}_{24} & \tilde{R}_c \\
                                                  &\frac{40 \frac{1}{2}}{3},&\frac{13 \frac{1}{2}}{1},&\frac{25 \frac{1}{2}}{3},&\frac{58 \frac{1}{2}}{1},&\frac{25 \frac{1}{2}}{3}} = \frac{17}{2}\,.}
We then increase the flows, increasing $y_{123}$ by $3 \times \frac{17}{2}$ to 30 and $y_{124}$ by $1 \times \frac{17}{2}$ to 10, decreasing $\tilde{S}_{123}$ to 15, $\tilde{S}_{124}$ to 5, $\tilde{R}_{23}$ to 0, $\tilde{R}_{24}$ to 50, and $\tilde{R}_c$ to 0.
This change inactivates movement $[1,2,3]$; and movement $[1,2,4]$ is then inactivated because these movement's flows are blocked by vehicles waiting to take $[1,2,3]$.
Since there are no more active turning movements, the algorithm terminates, and the final vector of flows is
\labeleqn{turnmovementfinal2}{\mb{y} = \vect{30 & 10 & 15 & 0 & 30 & 0}\,.}
\index{node model!general, unequal priority|)}

\subsection{Smoothing movement delays and capacities}

The node models described above aim to explicitly model the intersection dynamics in some manner or another.
This is most evident in the ``basic signal'' model, where the movements with positive flow at any time step correspond to the movements which have a green indication.
An alternative to this type of explicit model is to instead propagate flows based on average conditions: for a certain traffic loading, one can calculate the maximum and actual long-term flow rates for each movement, as well as the average delay vehicles would encounter.
The \emph{Highway Capacity Manual}, for instance, contains detailed procedures for estimating turning movement capacities and delays for all sorts of intersections: unsignalized, signalized, with permitted or protected turns, signal progression, under nonuniform arrivals, and so forth.

As a specific example, for a signalized intersection the capacity of a turning movement $[h,i,j]$ is defined as\label{not:Nhij}
\labeleqn{hcmcapacity}{
q^{max}_{hij} = N_{hij} s_{hij} \frac{G_{hij}}{C}
\,,
}
where $N_{hij}$ is the number of lanes in the turning movement, $s_{hij}$ is the per-lane saturation flow for the turning movement, $g_{hij}$ is the length of the green interval for the movement, and $C$ is the cycle length for the signal.
The manual includes detailed procedures for computing $s_{hij}$; a typical ``base rate'' is 1900 vehicles per hour per lane in urban areas, and 1750 vehicles per hour per lane in rural areas.
This base rate can be modified to account for lane width, the presence of heavy vehicles, roadway grade, parking, bus stops, and many other factors based on empirical data.

The manual also includes a procedure to estimate the signal delay $d_{hij}$\label{not:dhij} associated with the turning movement $[h,i,j]$.
A simple formula, assuming vehicles arrive at a uniform rate, is
\labeleqn{hcmdelay}{
d_{hij} = \frac{1}{2} \frac{C (1 - G_{hij}/C)^2}{1 - (\min \{ 1, S_{hij} / q^{max}_{hij} \} G_{hij} / C)}
\,.
}
More sophisticated formulas account for fluctuations in the arrival rate, signal progression along a corridor, and so forth.
The signal delay in equation~\eqn{hcmdelay} is \emph{not} associated with a queue that forms due to inflows exceeding capacity.
Even if inflow is below capacity, some vehicles will arrive on a red indication and experience delay.
We will use $d_{hij}$ to denote this kind of delay for an arbitrary kind of intersection.

There are arguments to be made for using these types of formulas or procedures in node modeling, rather than explicitly trying to model details of traffic signal timing, gap acceptance for permitted movements, and so forth.
First, these formulas allow us to build on the considerable amount of research which has been done in this area, including extensive field data and simulation data, rather than trying to derive a completely  new set of formulas.
Second, these formulas allow more nuances of intersection and signal configurations to be taken into account; changes in capacity or delay due to, say, street parking or heavy vehicles can be handled rather seamlessly.
Finally, it avoids the risk of overcalibration --- typical time steps for dynamic network loading are on the order of seconds, but route choice and other travel decisions are made at a more granular level (someone might choose their departure time within a few minutes, but likely not within seconds), so one might argue there is little sense in modeling traffic conditions at a finer level than travelers make decisions.

Of course, there are also decent arguments to be made for more explicit intersection models (or else this chapter would have skipped immediately to this section!).
Most obviously, they more closely model the underlying physical process, rather than invoking another model or formula which may have been derived under different assumptions than the link models.
Mathematically, the network loading equations are nonlinear, and in nonlinear systems replacing the inputs with an average value does not generally produce the average values of the outputs.
That is to say, modeling traffic flow assuming an average capacity and delay will not fully replicate the expected traffic conditions.

Creating a node model based on average movement delays and capacities is conceptually very similar to the all-way stop model described above, with two key differences.
First, we need some procedure for calculating the turning movement capacities $q_{max}^{hij}$ and the turning movement delays $d_{hij}$, based on the sending flows and turning proportions from incoming links.
This book will not give formulas for calculating these; you are encouraged to consult the \emph{Highway Capacity Manual} or other references for equations, or you may wish to come up with your own formulas based on simulation results or field data.
Second, we need to impose the delays $d_{hij}$ on vehicles at the node, possibly holding back vehicles until they have waited for at least $d_{hij}$ time units.
Unlike the previous node models, where delays arose naturally (say, during a red indication no vehicles would be allowed to move), in this case we must hold vehicles back manually, relying on the formula used to calculate $d_{hij}$ to produce a reasonable value.

To implement the second step, our link models must be able to track histories of upstream and downstream cumulative counts, and compute a \emph{delayed sending flow}\index{sending flow!delayed} $\hat{S}_{hi}$\label{not:Shat} which will only allow vehicles to leave the link $(h,i)$ if they have waited the appropriate amount of time.
The link model still produces $S_{hi}$ as usual, using whatever formula is appropriate to the link (point queue, cell transmission model, etc.); the delayed sending flow $\hat{S}_{hi}$ is an additional value calculated in equation~\eqn{delayedsendingflow} below.
This section presents a method which assumes that the delays $d_{hij}$ are identical for all of the turning movements corresponding to the same approach $(h,i)$ (which we will write as $d_{hi}$).
In general the delays for turning movements can be different (vehicles waiting to make a permitted turn might have to wait longer than those traveling straight through).
It is possible to generalize this node model to handle this case with concepts from Chapter~\ref{chp:dynamiceqm} that allow us to track cumulative counts along different paths (not just on links).

The algorithm involves the following steps:
\begin{enumerate}
   \item Calculate the turning movement capacities $q_{max}^{hij}$ and the approach delays $d_{hi}$.
   \item Obtain the delayed sending flows from the incoming link models:
      \labeleqn{delayedsendingflow}{\hat{S}_{hi}(t) = \min \myc{  N^\uparrow_{hi} \myp{t + \Delta t - \frac{L_{hi}}{u_f} - d_{hi}} - N^\downarrow_{hi}(t), S_{hi}(t)} \,.} 
   \item Perform the ``all-way stop'' algorithm, using the delayed sending flows $\hat{S}_{hi}$ in place of $S_{hi}$.
\end{enumerate}
\index{node model|)}

\section{Historical Notes and Further Reading}

The idea of dividing networks into link models\index{link model} and node models\index{node model} (which operate independently of each other, possibly with different models) was suggested by \cite{yperman_diss} and \cite{nie08}.
The sending flow and receiving flow concepts date to \cite{daganzo95kinematic} and \cite{lebacque96} (who used the terms ``demand'' and ``supply''), although the presentation in this chapter more closely follows that of \cite{yperman_diss}.
The point queue and spatial queue models were described in \cite{vickrey69} and \cite{zhang06}.

The list of node model desiderata in this chapter is that of \cite{tampere11}; the invariance principle\index{invariance principle} specifically is discussed more in \cite{lebacque05}.
The diverge and merge equations are that of \cite{daganzo95}, and the use of capacities to determine the ratio of sending flows at a congested merge is from \cite{ni05_merge}.
For an alternative merge model (which does \emph{not} satisfy the invariance principle), see \cite{jin03}.

The hydrodynamic traffic flow theory described in Section~\ref{sec:lwr} was independently developed by \cite{lighthill1955kinematic} and \cite{richards56}.\index{traffic flow theory!Lighthill-Whitham-Richards (LWR) model}
\cite{newell93a}, \cite{newell93b}, and \cite{newell93c} recognized that the cumulative vehicle counts $N$, and the analysis of characteristics resulting from a triangular fundamental diagram, greatly simplify the solution of the model, a theory completed by \cite{daganzo05a} and \cite{daganzo05b}.
Interestingly, an equivalent model (used for soil erosion) was separately developed by \cite{luke72}.
Cumulative counts have been used for other transportation engineering problems for some time; see \cite{moskowitz63} and \cite{makigami71} for an example from the research literature.
Karl Moskowitz,\index{Moskowitz, Karl} the lead author on the first of these publications, was a highly influential traffic engineer working at the California Division of Highways (present-day Caltrans).
For this reason, the cumulative counts $N(x,t)$ are sometimes referred to as the \emph{Moskowitz function}.
\index{Moskowitz function|see {cumulative count}}

There are alternative means of solving the LWR model not presented in this book, through recognizing it as a Hamilton-Jacobi system of partial differential equations \citep{leveque92,evans98} which can be solved using the Lax-Hopf formula\index{Lax-Hopf formula} or viability theory\index{viability theory} \citep{lax57,hopf70,claudel10}; or for the purposes of a link model, representing sending and receiving flows using a ``double queue,'' one at each end of the link~\citep{osorio09a,osorio11}, a representation suitable for a stochastic version of the LWR model.
The hydrodynamic model can also be derived from car-following principles, or from cellular automata models of traffic flow~\citep{daganzo06}.

The cell transmission model was reported in \cite{daganzo94} and \cite{daganzo95}, essentially a Godunov scheme\index{Godunov scheme} for solving the LWR system \citep{godunov59, lebacque96}.
The link transmission model was developed by \cite{yperman_diss}; see also \cite{gentile10}.

The more sophisticated node models reported later in the chapter are adapted from \cite{tampere11}, \cite{flotterod11}, and \cite{corthout12}.
For more examples of the ``smoothed'' node models, see \cite{durlin05,durlin08}, and Chapters 5--6 of \cite{yperman_diss}.
\cite{han14} and \cite{han15} explore the relationship between using a ``smoothed'' node model and explicitly tracking signal phases, and identify conditions under which the representations are more and less equivalent.

Finally, network loading can be accomplished by entirely different means than that reported in this chapter, without the use of explicit link and node models discretized in space and time.
For instance, the discretization can be done in the space of \emph{vehicle trajectories} \citep{bargera_dta06}.
In mathematical terms, this involves converting from Eulerian coordinates\index{Eulerian coordinates} ($x$ and $t$) to Lagrangian coordinates\index{Lagrangian coordinates} (with the cumulative count $N$ in place of either $x$ or $t$).
For more on this alternative, and reformulations of the LWR model with this change of variables, see \cite{laval13}.

Another common alternative is to use traffic simulation\index{traffic simulation} to perform the network loading.
Examples include the software packages VISSIM~\citep{fellendorf94}, AIMSUN~\citep{barcelo98}, DynaMIT~\citep{ben1998dynamit}, VISTA~\citep{vista}, DYNASMART~\citep{mahmassani00}, Dynameq~\citep{mahut03}, and DynusT.
\index{network loading|)}

\section{Exercises}
\label{sec:networkloading_exercises}

\begin{enumerate}
\item \diff{14} Table~\ref{tbl:pqsqex} shows cumulative inflows and outflows to a link with a capacity of 10 vehicles per time step, and a free-flow time of 2 time steps.
Use the point queue model to calculate the sending and receiving flow for each time step.
\label{ex:pqex}
\begin{table}
\begin{center}
\caption{Upstream and downstream counts for Exercises~\ref{ex:pqex} and~\ref{ex:sqex}.
\label{tbl:pqsqex}}
\begin{tabular}{|c|cc|}
\hline
$t$ & $N^\uparrow(t)$ & $N^\downarrow(t)$ \\
\hline
0 & 0 & 0 \\
1 & 5 & 0 \\
2 & 10 & 0 \\
3 & 15 & 2 \\
4 & 20 & 4 \\
5 & 20 & 6 \\
6 & 20 & 10 \\
7 & 20 & 15 \\
8 & 20 & 20 \\
9 & 20 & 20 \\
\hline
\end{tabular}
\end{center}
\end{table}
\item \diff{14} Repeat Exercise~\ref{ex:pqex} with the spatial queue model.
Assume that the jam density is such that at most 20 vehicles can fit on the link simultaneously.
\label{ex:sqex}
\item \diff{33} In the point queue model, if the inflow and outflow rates $q_{in}$ and $q_{out}$ are constants with $q_{in} \geq q_{out}$, show that the travel time experienced by the $n$-th vehicle is $L/u_f + \frac{n}{q_{in}} \myp{\frac{q_{in}}{q_{out}} - 1}$.
The same result holds for the spatial queue model, if there is no spillback.
\label{ex:pqtime}
\item \diff{24} In the spatial queue model, show that if the total number of vehicles on a link is at most $k_j L$ at time $t$, then the total number of vehicles on the link at time $t + \Delta t$ is also at most $k_j L$.
\label{ex:sqcompliance}
\item \diff{13} In a merge, approach 1 has a sending flow of 50 and a capacity of 100.
Approach 2 has a sending flow of 100 and a capacity of 100.
Report the number of vehicles moving from each approach if (a) the outgoing link has a receiving flow of 150; (b) the outgoing link has a receiving flow of 120; and (c) the outgoing link has a receiving flow of 100.
\item \diff{32} Extend the merge model of Section~\ref{sec:merge} to a merge node with three incoming links.
\item \diff{43} Show that the formula~\eqn{mergemedianformula} indeed captures both case II and case III of the merge model presented in Section~\ref{sec:merge}.
\label{ex:mergemedianformula}
\item \diff{32} Show that the merge model of Section~\ref{sec:merge} satisfies all the desiderata of Section~\ref{sec:nodemodels}.
\label{ex:mergedesiderata}
\item \diff{41} Consider an alternative merge model for the congested case, which allocates the receiving flow proportional to the \emph{sending flows} of the incoming links, rather than proportional to the \emph{capacities} of the incoming links as was done in Section~\ref{sec:merge}.
Show that this model does \emph{not} satisfy the invariance principle.
\item \diff{21} Modify the algorithm at the end of Section~\ref{sec:merge} to handle the case when $\beta_{gi} = 0$.
\label{ex:zeromerge}
\item \diff{13} In a diverge, the incoming link has a sending flow of 120, 25\% of the vehicles want to turn onto outgoing link 1, and the remainder want to turn onto outgoing link 2.
Report the number of vehicles moving to outgoing links 1 and 2 if their respective receiving flows are (a) 80 and 100; (b) 80 and 60; (c) 10 and 40.
\item \diff{10} Extend the diverge model of Section~\ref{sec:diverge} to a diverge node with three outgoing links.
\item \diff{21} Develop a model for a diverge node with two outgoing links, in which flows waiting to enter one link do \emph{not} block flows entering the other link.
When might this model be more appropriate?
\item \diff{25} Show that the diverge model of Section~\ref{sec:diverge} satisfies all the desiderata of Section~\ref{sec:nodemodels}.
\label{ex:divergedesiderata}
\item \diff{10} In the network loading procedure in Section~\ref{sec:combiningnodelink}, we specified that centroid connectors starting at origins should have high jam density, and those ending at destinations should have high capacity.
Would anything go wrong if centroid connectors starting at origins also had high capacity?  What if centroid connectors ending at destinations had high jam density?
\item \diff{12} Draw trajectory diagrams which reflect the following situations: (a) steady-state traffic flow, no vehicles speeding up or slowing down; (b) vehicles approaching a stop sign, then continuing; (c) a slow semi truck merges onto the roadway at some point in time, then exits at a later point in time.
Draw at least five vehicle trajectories for each scenario.
\item \diff{33} The relationship $q = uk$ can be used to transform the fundamental diagram (which relates density and flow) into a relationship between density and speed, and vice versa.
For each of these three speed-density relationships, derive the corresponding fundamental diagram by writing an expression for $q$ in terms of $k$, and produce a plot.
\label{ex:greenshieldsgreenberg}
\begin{enumerate}[(a)]
 \item The Greenshields (linear) model: $u = u_f(1 - k/k_j)$.
 \item The Greenberg (logarithmic) model: $u = C \log (k_j / k)$ where $C$ is a constant.
 \item The Underwood (exponential) model: $u = u_f \exp (-k/k_c)$.
 \item The Pipes model: $u = u_f(1 - (k/k_j)^n)$ where $n$ is a constant.
 Plot the fundamental diagram for the case $n = 2$.
 (The Greenshields model is a special case when $n = 1$.)
\item What features of the Greenberg and Underwood models make them less suitable for dynamic network loading?  (Hint: Draw plots of these speed-density diagrams.) 
\end{enumerate}
\item \diff{51} Which of the following statements are true with the LWR model and a concave fundamental diagram?
\begin{enumerate}[(a)]
 \item Speed uniquely defines the values of flow and density.
 \item It is possible for higher density to be associated with higher speed.
 \item No shockwave can move downstream faster than the free-flow speed.
 \item Flow uniquely defines the values of speed and density.
 \item With a triangular fundamental diagram, traffic speed is constant for subcritical densities.
 \item Density uniquely defines the values of speed and flow.
\end{enumerate}
\item \diff{44} The relationship $q = uk$ can transform the fundamental diagram (which relates density and flow) into a relationship between speed and flow.
The speed-flow relationship is traditionally plotted with the speed on the vertical axis and flow on the horizontal axis.
\begin{enumerate}[(a)]
\item Show that for any concave fundamental diagram and any flow value $q$ less than the capacity, there are exactly two possible speeds $u_1$ and $u_2$ producing the flow $q$, one corresponding to subcritical (uncongested) conditions and the other corresponding to supercritical (congested) conditions.
\item Derive and plot the speed-flow relationship for the Greenshields model of Exercise~\ref{ex:greenshieldsgreenberg}.
Express this relationship with two functions $u_1(q)$ and $u_2(q)$ corresponding to uncongested and congested conditions, respectively; these functions should have a domain of $[0, q_{max}]$ and intersect at capacity.
\end{enumerate}
\item \diff{22} For each of these fundamental diagrams, derive the speed-density function (that is, the travel speed for any given density value), and provide a sketch.
\begin{enumerate}[(a)]
\item $Q(k) =  C k (k_j - k)$, where $C$ is a constant and $k_j$ is the jam density.
\item $Q(k) = \min \{ u_f k, w(k_j - k) \}$
\item $Q(k) = \min \{ u_f k, q_{max}, w(k_j - k) \}$
\end{enumerate}
\item \diff{44} Consider a long, uninterrupted freeway with a capacity of 4400 vehicles per hour, a jam density of 200 vehicles per mile, and a free-flow speed of 75 miles per hour.
Initially, freeway conditions are uniform and steady with a subcritical flow of 2000 vehicles per hour.
An accident reduces the roadway capacity to 1000 veh/hr for thirty minutes.
Draw a shockwave diagram to show the effects of this accident, reporting
the space-mean speed, volume, and density in each region of your diagram, and the speed and direction of each shockwave.
Assume that the fundamental diagram takes the shape of the Greenshields model (Exercise~\ref{ex:greenshieldsgreenberg}), and that a stopped queue discharges at capacity.
\item \diff{46} Consider a roadway with a linear-speed density relationship (cf.\ Exercise~\ref{ex:greenshieldsgreenberg})  whose capacity is 2000 veh/hr and free-flow speed is 40 mi/hr.
Initially, the flow is 1000 veh/hr and uncongested.
A traffic signal is red for 45 seconds, causing several shockwaves.
When the light turns green, the queue discharges at capacity.
\begin{enumerate}[(a)]
\item Sketch a time-space diagram, indicating all of the shockwaves which are formed.
\item Calculate the speed and direction of each shockwave from your diagram.
\item What is the minimum green time needed to ensure that no vehicle has to stop more than once before passing the intersection? (Neglect any yellow time, reaction time, etc.
Assume that when the signal is green, people move immediately, and that when it is red, people stop immediately.)
\end{enumerate}
\item \diff{56} Consider a single-lane roadway with a triangular fundamental diagram, a free-flow speed of 60 mi/hr, a backward wave speed of 30 mi/hr, and a jam density of 200 veh/mi.
Initially, traffic flow is uncongested, and the volume is half of capacity.
A slow-moving truck enters the roadway at time $t = 0$, and travels at 20 mi/hr.
This vehicle turns off of the roadway one mile later.
\begin{enumerate}[(a)]
\item What is the capacity of the roadway?
\item At time $t = 2$ minutes, you are a quarter of a mile behind the truck.
Use the Newell-Daganzo method to determine how many vehicles are between you and the truck.
\item In total, how many shockwaves are generated by the slow-moving truck?  Sketch them on a trajectory diagram.
\end{enumerate}
\item \diff{11} Show that a shockwave connecting two uncongested (subcritical) traffic states always moves downstream, while a shockwave connecting two congested (supercritical) traffic states always moves upstream.
This is related to the observation in the chapter that ``uncongested states propagate downstream, and congested states propagate upstream.''
\item \diff{36} This exercise asks you to fill in some details of the example in Section~\ref{sec:trafficflowtheory} where the fundamental diagram was $Q(k) = \frac{1}{240} k (240 - k)$ and the cumulative count map was $N(x,t) = 60t - 120 x + \frac{60x^2}{t+1}$.
Times are measured in minutes, and distances in miles.
\label{ex:completelwrexample}
\begin{enumerate}[(a)]
\item Calculate the capacity, jam density, and free-flow speed associated with this fundamental diagram.
\item Verify that the conservation relationship~\eqn{vehicleconservation} is satisfied by the flow and density maps $q(x,t)$ and $k(x,t)$.
\item Verify that the density and flow maps are consistent with the given fundamental diagram.
\item Calculate the speed $u(x,t)$ at each point and time.
Are vehicles accelerating, decelerating, or maintaining a constant speed?
\end{enumerate}
\item \diff{45} Consider the network in Figure~\ref{fig:pqnetworkloadingex}, where each link has a free-flow time of 5 minutes and a
capacity shown on the figure, and vehicles split equally at each diverge (that is, $p_{ABC} = p_{ABD} = p_{BCD} = p_{BCE} = 1/2$ at all times).
Vehicles enter the network at a rate of 80 veh/min for 20 minutes, and then the inflow rate drops to zero.
Perform dynamic network loading, using point queues for the link models.
For each link in the network, plot the cumulative counts $N^\uparrow$ and $N^\downarrow$  over time, as well
as the sending flow and receiving flow over time.
At what time does the last vehicle leave the network? \label{ex:pqnetworkloadingex}
\stevefig{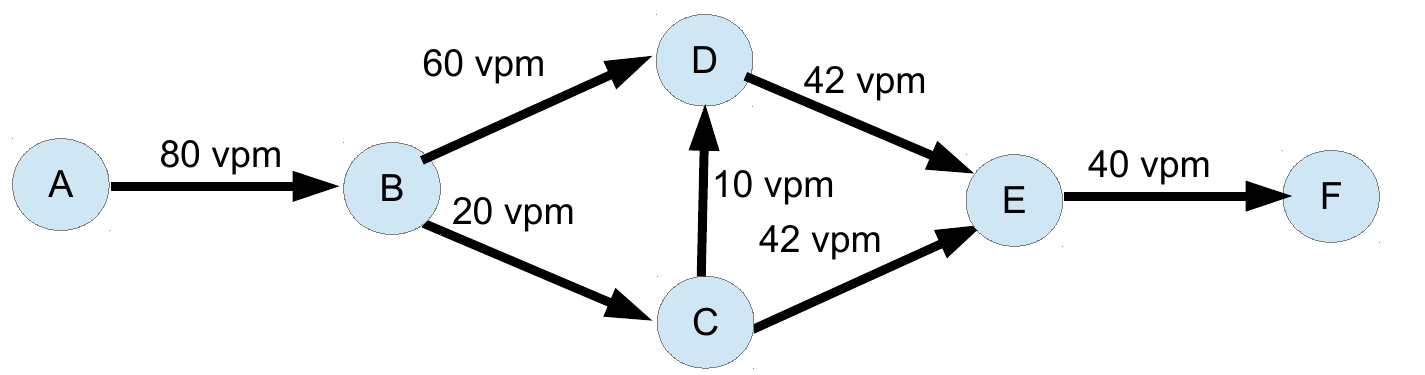}{Network for Exercise~\ref{ex:pqnetworkloadingex}}{\textwidth}
\item \diff{25} Write the formula for the fundamental diagram $Q(k)$ in the cell transmission model example depicted in Table~\ref{tbl:ctmexample}.
\item \diff{13} A link is divided into four cells; on this link the capacity is 10 vehicles per time step, each cell can hold at most 40 vehicles, and the ratio of backward wave speed to free-flow speed is 0.5.
Currently, the number of vehicles in each cell is as in Table~\ref{tbl:ctmex} (Cell 1 is at the upstream end of the link, Cell 4 at the downstream end.)   Calculate the number of vehicles that will move between each pair of cells in the current time interval (that is, the $y_{12}$, $y_{23}$, and $y{34}$ values.), and the number of vehicles in each cell at the start of the next time interval.
Assume no vehicles enter or exit the link.
\label{ex:ctmex}.
\begin{table}
\begin{center}
\caption{Current cell occupancies for Exercise~\ref{ex:ctmex}.
\label{tbl:ctmex}}
\begin{tabular}{|cccc|}
\hline
Cell 1 & Cell 2 & Cell 3 & Cell 4 \\
\hline
8 & 10 & 30 & 5 \\
\hline
\end{tabular}
\end{center}
\end{table}
\item \diff{23} Table~\ref{tbl:ltmex} shows cumulative inflows and outflows to a link with a capacity of 10 vehicles per time step, a free-flow time of 2 time steps, and a backward wave time of 4 time steps.
At jam density, there are 20 vehicles on the link.
Use the link transmission model to calculate the sending flow $S(7)$ and the receiving flow $R(7)$.
\label{ex:ltmex}
\begin{table}
\begin{center}
\caption{Upstream and downstream counts for Exercise~\ref{ex:ltmex}.
\label{tbl:ltmex}}
\begin{tabular}{|c|cc|}
\hline
$t$ & $N^\uparrow(t)$ & $N^\downarrow(t)$ \\
\hline
0 & 0 & 0 \\
1 & 5 & 0 \\
2 & 10 & 0 \\
3 & 15 & 2 \\
4 & 16 & 4 \\
5 & 17 & 6 \\
6 & 20 & 10 \\
7 & 20 & 15 \\
\hline
\end{tabular}
\end{center}
\end{table}
\item \diff{44} \emph{(Exploring shock spreading.)}  A link is seven cells long; at most 15 vehicles can fit into each cell, the capacity is 5 vehicles per timestep, and $w / u_f = 1/2$.
Each time step, 2 vehicles wish to enter the link, and will do so if the receiving flow can accommodate.
There is a traffic signal at the downstream end of the link.
During time steps 0--9, and from time step 50 onward, the light is green and all of the link's sending flow can leave.
For the other time steps, the light is red, and the sending flow of the link is zero.
\label{ex:shockspreading}
\begin{enumerate}[(a)]
\item Use the cell transmission model to propagate flow for 80 time steps, portraying the resulting cell occupancies in a time-space diagram (time on the horizontal axis, space on the vertical axis).
At what time interval does the receiving flow first begin to drop; at what point does it reach its minimum value; and what is that minimum value?  Is there any point at which the entire link is at jam density?  
\item Repeat, but with $w/u_f = 1$.
\item Repeat, but instead use the link transmission model (with the same time step) to determine how much flow can enter or leave the link.
\end{enumerate}
\item \diff{68} Consider the network in Figure~\ref{fig:netloadingex}.
The figure shows each link's length, capacity, jam density, free-flow speed, and backward wave speed.
The inflow rate at node A is 4320 veh/hr for $0 \leq t < 30$ ($t$ measured in seconds), 8640 veh/hr for $25 \leq t < 120$, and 0 veh/hr thereafter.
The splitting proportion towards node $C$ is $\frac{2}{3}$ for $0 \leq t < 50$, $\frac{1}{6}$ for $50 \leq t < 90$, and $\frac{1}{2}$ for $t \geq 90$.
\label{ex:netloadingex}
\begin{enumerate}[(a)]
\item Use a point queue model to propagate the vehicle flow with the time step $\Delta t = 5$ s.
Plot the turning movement flows $y_{ABC}$, $y_{ABD}$, $y_{BDE}$, and $y_{BCE}$ from $t = 0$ until the last vehicle has left the network.
($q_{12}$ is the rate at which flow leaves the downstream end of link 1 to enter the upstream end of link 2).
\item Use the cell transmission model to propagate the vehicle flow with the time step $\Delta t = 5$ s.
Plot the same flow rates as in the previous part.
\item Use the link transmission model to propagate the vehicle flow with the time step $\Delta t = 10$ s.
Plot the same flow rates as in the previous part.
\item Comment on any differences you see in these plots for the three flow models.
\end{enumerate}
\stevefig{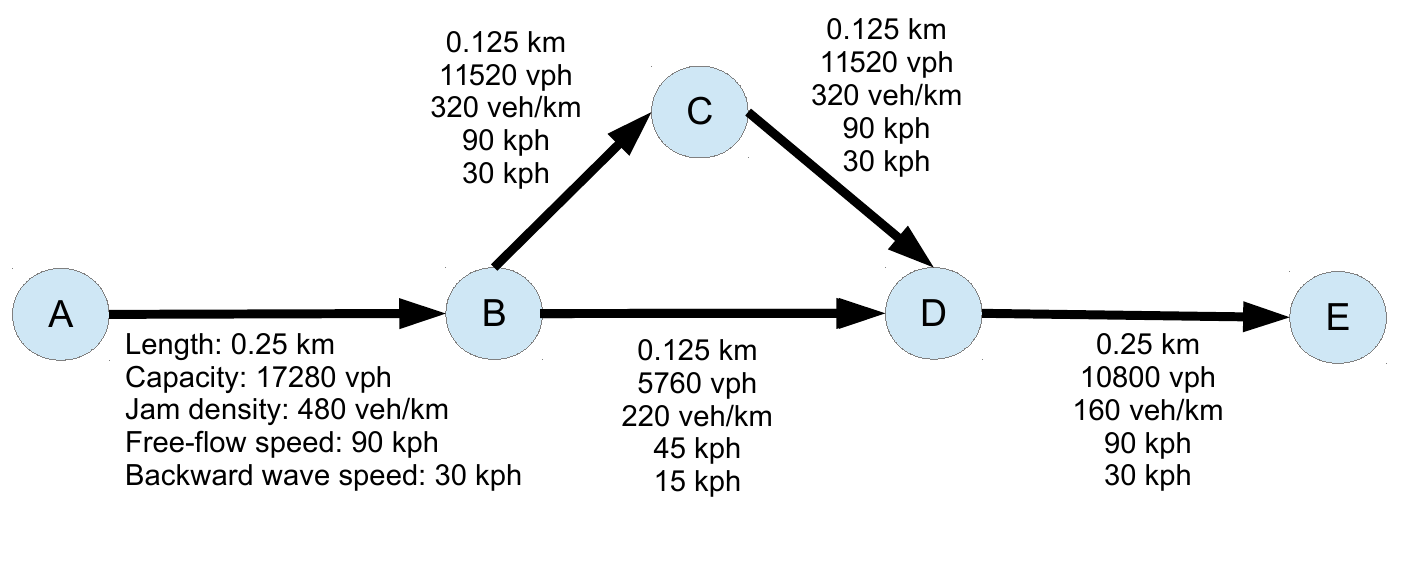}{Network for Exercise~\ref{ex:netloadingex}}{\textwidth}
\item \diff{21} Assuming that a cell initially has between 0 and $\bar{n}$ vehicles, show the cell transmission model formula~\eqn{ctmtransitionflows} ensures that it will have between 0 and $\bar{n}$ vehicles at all future time steps, regardless of upstream or downstream conditions.
\item \diff{21} On a link, we must have $N^\uparrow(t) \geq N^\downarrow(t)$ at all time steps.
Assuming this is true for all time steps before $t$, show that the link transmission model formulas~\eqn{ltmsending} and~\eqn{ltmreceiving} ensure this condition holds at $t$ as well.
\begin{figure}
\begin{center}
\begin{tikzpicture}[
declare function={q(\x) = (
(\x <= 40) * (50 * \x) 
+ and(\x >= 40, \x <= 140) * (10*(\x - 40) + 2000)
+ and(\x >= 140, \x <= 150) * (3000)
+ and(\x >= 150, \x <= 250) * (10*(150 - \x) + 3000)
+ and(\x >= 250, \x <= 300) * 40 * (300 - \x)
);
}]
\begin{axis}[ticks=none, axis x line=middle, axis y line=middle,
ymin = 0, ymax = 4000, ytick = {}, ylabel = {$q$ (veh/hr)}, xmin = 0, xmax = 450, xlabel = {$k$ (veh/mi)}]
\addplot[samples at ={0, 39.99, 40.01, 139.99, 140.01, 149.99, 150.01, 249.99, 250.01, 299.99}]{q(x)};
\filldraw (axis cs: 40,2000) circle (1pt) node [right,font=\small] {$(40,2000)$};
\filldraw (axis cs: 140,3000) circle (1pt) node [left,font=\small] {$(140,3000)$};
\filldraw (axis cs: 150,3000) circle (1pt) node [right,font=\small] {$(150,3000)$};
\filldraw (axis cs: 250,2000) circle (1pt) node [right,font=\small] {$(250,2000)$};
\filldraw (axis cs: 300,0) circle (1pt) node [above left,font=\small] {$300$};
\end{axis}
\end{tikzpicture}

\caption{Fundamental diagram for Exercises~\ref{ex:hexafundactm} and~\ref{ex:hexafundaltm}. \label{fig:hexafunda}}
\end{center}
\end{figure}
\item \diff{33} Write cell transmission model formulas for sending and receiving flow when the fundamental diagram is given by the piecewise-linear curve in Figure~\ref{fig:hexafunda}.
\label{ex:hexafundactm}
\item \diff{43} Write link transmission model formulas for sending and receiving flow when the fundamental diagram is given by the piecewise-linear curve in Figure~\ref{fig:hexafunda}.
\label{ex:hexafundaltm}
\item \diff{42} Generalize the cell transmission model formula~\eqn{ctmtransitionflows} to handle an arbitrary piecewise-linear fundamental diagram (not necessarily triangular or trapezoidal).
\item \diff{65} Generalize the link transmission model formulas~\eqn{ltmsending} and~\eqn{ltmreceiving} to handle an arbitrary piecewise-linear fundamental diagram.
\label{ex:ltmextension}
\item \diff{35} Figure~\ref{fig:lamarguadalupe} represents the intersection of Lamar and Guadalupe, showing sending and receiving flows, saturation flows, turning movement proportions, and the signal timing plan (assume no lost time due to clearance intervals or startup delay).
Note that the receiving flow on northbound Lamar is quite low, because of congestion spilling back from a nearby signal just downstream.
No U-turns are allowed, and drivers may not turn left from Guadalupe onto southbound Lamar.
\label{ex:lamarguadalupe}
\begin{enumerate}[(a)]
\item Find the transition flows $y_{ijk}$ for all five turning movements at the current time step, using the ``smoothed signal'' node model.
\item The southbound receiving flow on Guadalupe is now reduced to 50 due to congestion further downstream.
Find the updated transition flow rates for all turning movements.
\end{enumerate}
\stevefig{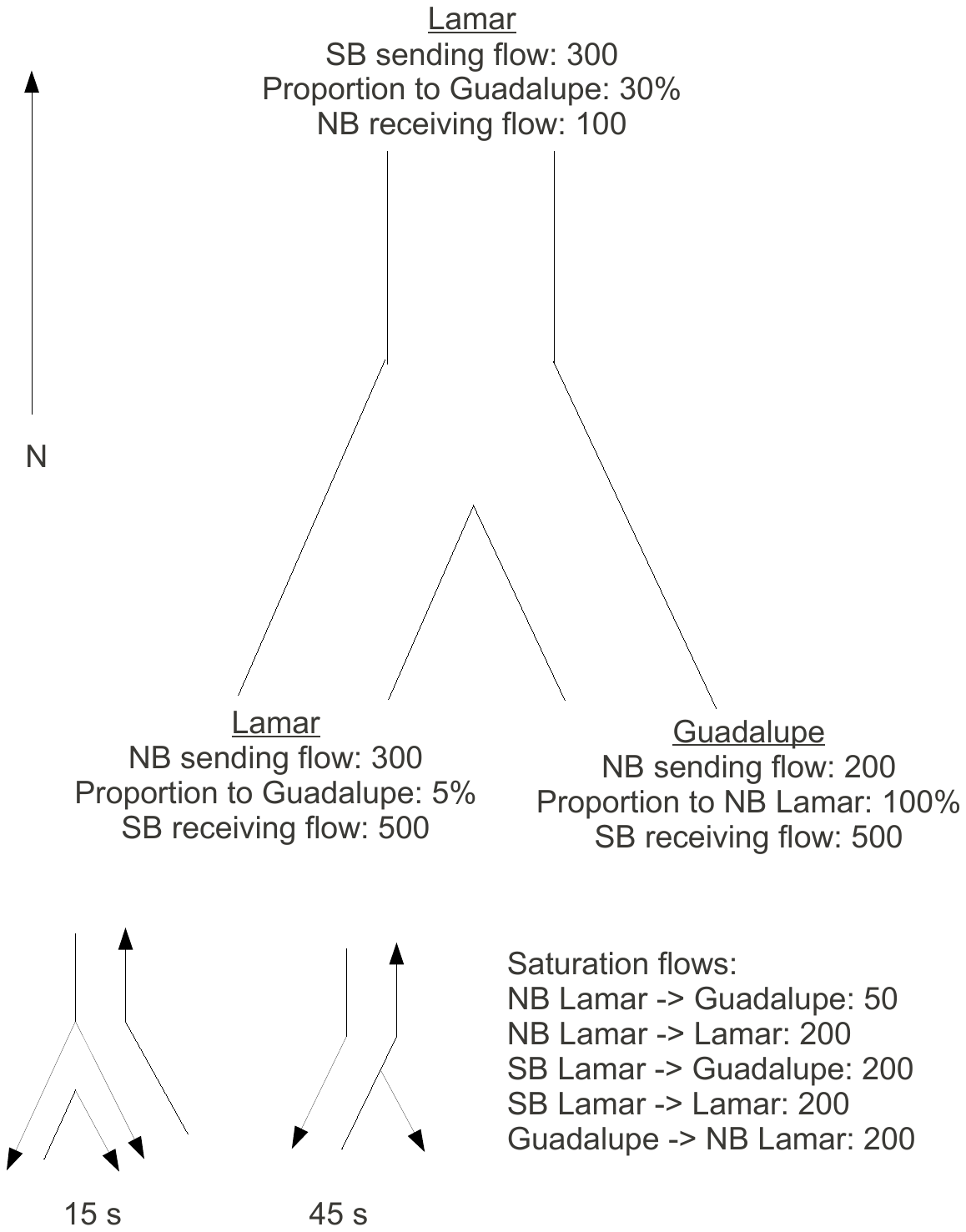}{Intersection for Exercise~\ref{ex:lamarguadalupe}}{0.6\textwidth}
\item \diff{51} Extend the ``basic signal'' node model to account for turns on red, where a vehicle facing a red indication may make a turn in the direction nearest to them (usually right-on-red in countries that drive on the right, left-on-red in countries that drive on the left).
Vehicles turning on red must yield to traffic which has a green indication.
\item \diff{53} Extend the ``basic signal'' node model to the case where there are turn lanes, and not all turn lanes from an approach have the same green time.
\item \diff{61} Extend the ``basic signal'' node model to account for permitted turns on green (a turning movement which has a green indication, but must yield to oncoming traffic).
\item \diff{36} Consider the four-legged intersection shown in Figure~\ref{fig:canvascross}.
The set of turning movements is $\Xi = \{ [1,2,3], [1,2,4], [3,2,1], [3,2,4], [5,2,4] \}$.
For the current time step, the sending and receiving flow values are $S_{12} = 6$, $S_{32} = 6$, $S_{52} = 1$, $R_{12} = 14$, $R_{23} = 14$, and $R_{24} = 14$.
Half of the drivers approaching from link (1,2) want to turn left, and half want to go straight.
Half of the drivers approaching from link (3,2) want to turn right, and half want to go straight.
Finally, the intersection geometry and signal timing are such that the capacities of the incoming links are $q_{12}^{max} = 200$, $q_{32}^{max} = 2000$, and $q_{52}^{max} = 20$.

Apply the ``equal priorities'' algorithm, and report the transition flows ($y$ values) for each turning movement.
\label{ex:canvascross}
\genfig{canvascross}{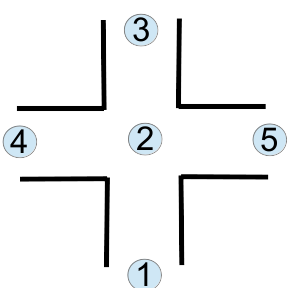}{Intersection for Exercise~\ref{ex:canvascross}.}{width=0.3\textwidth}
\item \diff{36} Repeat Exercise~\ref{ex:canvascross}, but with these values of the sending and receiving flows: $S_{12} = 18$, $S_{32} = 18$, $S_{52} = 3$, $R_{12} = 42$, $R_{23} = 42$, and $R_{24} = 12$.
\item \diff{68} Modify the ``partial stop control'' node model to allow cases of absolute priority (if movement $[h,i,j]$ has absolute priority over $[h',i,j']$ at conflict point $c$, then $\alpha^c_{h'ij'}/\alpha^c_{hij}$ would be zero.)  You will need to decide how priority will be granted for any combination of sending flows wishing to use a conflict point (including flows of equal absolute priority, and when not all sending flows are present), and your formulas can never divide by zero.
\item \diff{89} Show that all of the node models in Section~\ref{sec:fancynodes} satisfy all of the desiderata from Section~\ref{sec:nodemodels}.

\end{enumerate}

\chapter{Time-Dependent Shortest Paths}
\label{chp:tdsp}

\index{shortest path!time-dependent|(}
This chapter discusses how travelers make choices when traveling in networks whose state varies over time.
Two specific choices are discussed: how drivers choose a route when link costs are time-varying (Sections~\ref{sec:tdspconcept} and~\ref{sec:tdspalgos}), and how drivers choose a departure time (Section~\ref{sec:departuretimechoice}).
This chapter is the complement of the previous one.
In network loading, we assumed that the travelers' choices were known, and we then determined the (time-varying) flow rates and congestion pattern throughout the network.
In this chapter, we take this congestion pattern as known, and predict the choices that travelers would make given this congestion pattern.
In particular, by taking the congestion level as fixed, we can focus the question on an \emph{individual} traveler and do not need to worry about changes in congestion based on these choices just yet.

Chapter~\ref{chp:networkrepresentations} presented the shortest path problem in networks with constant travel times.
This chapter addresses shortest path problems where the link travel times can vary based on when the link is traversed.
Accounting for these time-varying link travel times is essential in dynamic traffic assignment.
The presentation in this chapter is self-contained, so it is not necessary to have read Chapter~\ref{chp:networkrepresentations} before reading this chapter.
Nevertheless, it is instructive to compare how static and dynamic shortest path algorithms are similar and different, and if you want to pursue further studies in network modeling it would be very helpful to read these two chapters together for comparison.

\section{Time-Dependent Shortest Path Concepts}
\label{sec:tdspconcept}

The \emph{time-dependent shortest path} problem involves finding a path through a network of minimum cost, when the cost of links varies with time.
As with the static shortest path problem, the ``cost'' of a link can include travel time, monetary costs, a combination of these, or any other disutility which can be added across links; and by a ``shortest'' path we mean one with least cost.
If a traveler enters link $(i,j)$ at time $t$, they will experience a cost of $c_{ij}(t)$.\label{not:cijt}
Unlike the static shortest path problem, however, we must always keep track of the travel time on a link even if the cost refers to a separate quantity: because the network state is dynamic, we must always know the time at which a traveler enters a link.
The travel time experienced by a traveler entering link $(i,j)$ at time $t$ is denoted by $\tau_{ij}(t)$.\label{not:tauijt}
This problem can be formulated either in discrete time (where $t$ and $\tau_{ij}(t)$ are limited to being integer multiples of the timestep $\Delta t$) or in continuous time (where $t$ can take any real value within a stated range).

There are several different variations of the time-dependent shortest path problem, all stemming from dynamic link costs and travel times.
One can consider time-dependent shortest paths where \emph{waiting}\index{shortest path!time-dependent!waiting} at intermediate nodes is allowed, or one can forbid it.
In public transit networks, waiting at intermediate nodes is logical, but in road networks, one would not expect drivers to voluntarily stop and wait in the middle of a route.
One can also restrict attention to problems where the link travel times satisfy the \emph{first-in, first-out}\index{first-in, first-out (FIFO)} (FIFO) property, where it is impossible for one to leave a link earlier by entering later, that is, for any link $(i,j)$ and any distinct $t_1 < t_2$, we have
\labeleqn{fifo}{t_1 + \tau_{ij}(t_1) \leq t_2 + \tau_{ij}(t_2)\,.}
If this is true, more efficient algorithms can be developed.
As an example, if the cost of a link is its travel time ($c_{ij}(t) = \tau_{ij}(t)$), there is no benefit to waiting in a FIFO network.
In discrete time, equation~\eqn{fifo} is equivalent to requiring
\labeleqn{fifodiscrete}{\tau_{ij}(t) - \tau_{ij}(t + \Delta t) \leq \Delta t\,,}
that is, a link's travel time cannot decrease by more than $\Delta t$ in one time step.
In continuous time, if $\tau_{ij}(t)$ is a continuous, piecewise differentiable function, we must have
\labeleqn{fifocontinuous}{\frac{d}{dt} \tau_{ij}(t) \geq -1\,,}
everywhere that $\tau_{ij}(t)$ has a derivative, which expresses the same idea.

One can also distinguish time-dependent shortest path problems by whether the departure time is fixed, whether the arrival time is fixed, or neither.
In the first case, the driver has already decided when they will leave, and want to find the route to the destination with minimum cost when leaving at that time, regardless of the arrival time at the destination.
In the second case, the arrival time at the destination is known (perhaps the start of work), and the traveler wants to find the route with minimum cost arriving at the destination at that time (regardless of departure time).
In the third case, both the departure and arrival times are flexible (as with many shopping trips), and the traveler wants to find a minimum cost route without regard to when they leave or arrive.
In this way, we can model the departure time choice decision simultaneously with the route choice decision.
Section~\ref{sec:departuretimechoice} develops this approach further, showing how we can incorporate penalty ``costs'' associated with different departure and arrival times.
Strictly speaking, one can imagine a fourth variant where \emph{both} departure and arrival time are fixed, but this problem is not always well-posed; there may be no path from the origin to the destination with those exact departure and arrival times.
In such cases, we can allow free departure and/or arrival times, but severely penalize departure/arrival times that differ greatly from the desired times.

In comparison with the static shortest path problem, the fixed departure time variant is like the one origin-to-all destinations shortest path problem, and the fixed arrival time variant is like the all origins-to-one  destination shortest path problem.

To be specific, in this chapter, we will first focus on time-dependent shortest paths with fixed departure times and free arrival times in Section~\ref{sec:tdspalgos}, but all of these approaches can be adapted to solve the fixed arrival time/free departure time version without much difficulty.
The origin $r$ and departure time $t_0$\label{not:t0} will be specified, and we will find the paths to all other nodes along which the sum of the costs $c_{ij}$ is minimal, given the changes in the costs which will occur during travel.
We will study two specific variants of the problem.
In the first, time is continuous, but the network is assumed to follow the FIFO principle and the link costs must be the travel times.
In the other variants, we do not require the FIFO assumption, and the link costs need not be the same as the travel times; but in exchange we will restrict ourselves to discrete time.
Section~\ref{sec:departuretimechoice} will then treat the case of free departure and arrival times, where only the origin $r$ will be specified.

\subsection{Time-expanded networks}

\index{network!time-expanded|(}
\index{time-expanded network|see {network, time-expanded}}
For discrete-time shortest path problems, we can form what is known as the \emph{time-expanded network}.
This technique transforms a time-dependent shortest path problem into a static shortest path problem, that can be solved using algorithms for this simpler problem (see Section~\ref{sec:shortestpath}).
However, the time-expanded network contains many more links and nodes than the original network, and can impose computational burdens.

If the time step is $\Delta t$, then the possible arrival times at any node are $0, \Delta t, 2\Delta t, \ldots, T \Delta t$ where $T$ is the time horizon under consideration.\footnote{In transportation systems, it is often reasonable to assume that after a large enough amount of time, all congestion will dissipate and travel times can be treated as constants equal to free-flow time. All routing after this point can be done with a static shortest path algorithm.}
In the time-expanded network, we create $T + 1$ copies of each node, one for each possible arrival time.
So, if the original network has $n$ nodes, the time-expanded network has $nT$ nodes.
Each node in the time-expanded network is written in the form $i:t$,\label{not:it} denoting the physical node $i$ at the time interval $t \Delta t$.
This labeling is shown in Figure~\ref{fig:timeexpandednet2}.
The original network is called the \emph{physical network}\index{network!physical}\index{physical network|see {network, physical}} when we need to distinguish it from the corresponding time-expanded network.

\stevefig{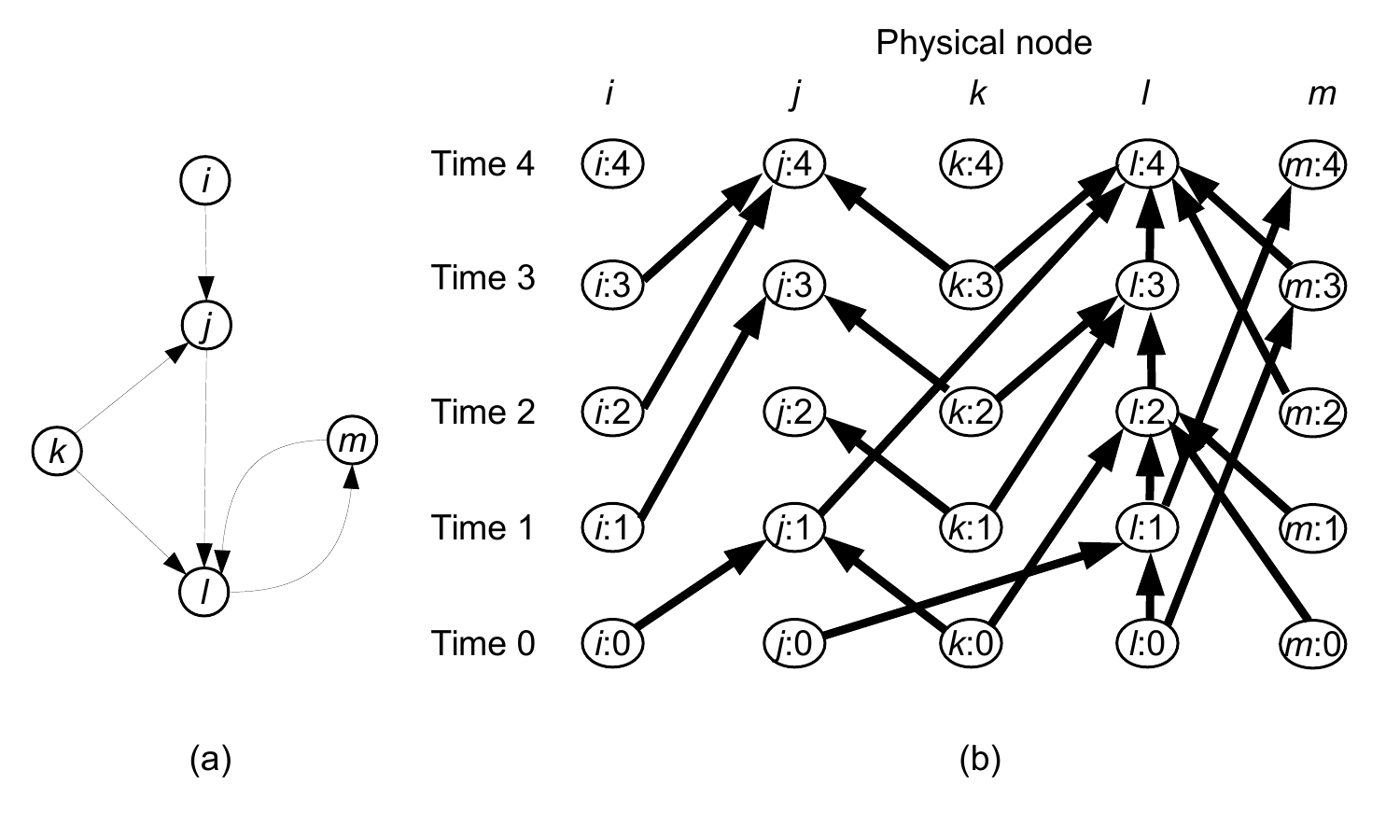}{(a) Original network (b) Corresponding time-expanded network with five physical nodes and five time intervals.}{\textwidth}

Links in the time-expanded network represent both the physical connection between nodes, as well as the time required to traverse the link.
In particular, for each link $(i,j)$ in the original network, and for each time interval $t \Delta t$, if the travel time is $\tau_{ij}(t \Delta t) = k \Delta t$, we create a link $(i:t, j:(t + k))$ in the time expanded network, assuming that $t + k \leq T$, and set the cost of this link equal to $c_{ij}(t)$.
A little bit of care must be taken when links would arrive at a node later than the time horizon.
If the time horizon is large enough, this should not be a significant issue.
In this chapter, we will not create a time-expanded link if its head node falls outside of the time horizon.
Exercises~\ref{ex:boundary1} and~\ref{ex:boundary2} introduce two other ways to treat boundary issues associated with the time horizon.

The advantage of the time-expanded network is that it reduces the time-dependent shortest path problem to the static one: solving the one-to-all static shortest path problem from node $i$ at time $t_0$ corresponds exactly to solving the fixed-departure time-dependent shortest path problem.
Furthermore, if all link travel times are positive, the time-expanded network is acyclic, with the time labels forming a natural topological order.
Shortest paths on acyclic networks can be solved rather quickly, so this is a significant advantage.
Even if some link travel times are zero (as may occur with some artificial links or centroid connectors), the time-expanded network remains acyclic unless there is a \emph{cycle} of zero-travel time links in the network; and if that is the case, it is often possible to collapse the zero-travel time cycle into a single node.
Time-expanded networks can also unify some of the variants of the time-dependent shortest path problem.
If waiting is allowed at a node $i$, we can represent that with links of the form $(i:t, i:(t+\Delta t))$ connecting the same physical node to itself, a time step later.
The FIFO principle in a time-dependent network means that two links connecting the same physical nodes will never cross (although they may terminate at the same node).

A disadvantage is that the number of time intervals $T$ may be quite large.
In dynamic traffic assignment, network loading often uses a time step on the order of a few seconds.
If this same time step is used for time-dependent shortest paths, a typical planning period of a few hours means that $T$ is in the thousands.
Given a physical network of thousands of nodes, the time-expanded network can easily exceed a million nodes.
Even though the time-expanded network is acyclic, this is a substantial increase in the underlying network size, increasing both computation time and amount of computer memory required.
With clever implementations, it is often possible to avoid explicitly generating the entire time-expanded network, and only generate links and nodes as needed.

\subsection{Bellman's principle in time-dependent networks}

The number of paths between any two points in a network can be very large, growing exponentially with the network size.
Therefore, any approach based on enumerating all paths and comparing their cost will not scale well to realistic networks.
So, efficient algorithms for finding shortest paths, whether time-dependent or not, must be more clever.
The key insight, \emph{Bellman's principle},\index{Bellman's principle} was introduced in Section~\ref{sec:shortestpath} for the static shortest path problem.
Briefly reviewing, Bellman's principle states that any segment of a shortest path between two nodes must itself be a shortest path between the endpoints of that segment.
Otherwise, the shortest path between the endpoints of the segment could be spliced into the original path, reducing its cost (Figure~\ref{fig:bellman}).

\stevefig{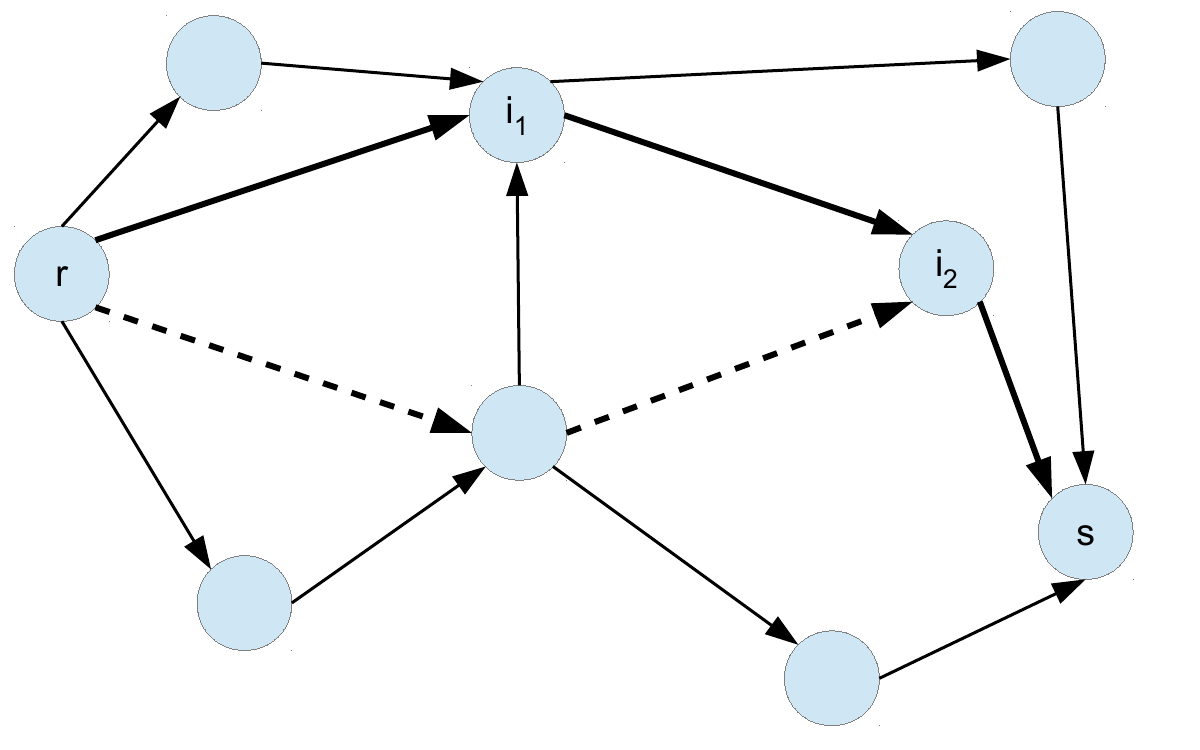}{Bellman's principle illustrated.
If the dashed path were a shorter route from $s$ to $i_2$, then it would also form a ``shortcut'' for the route from $s$ to $d$.}{0.6\textwidth}

In a time-dependent shortest path problem, the same general idea applies but we must be slightly more careful about how the principle is defined.
In a FIFO network where link costs are equal to the travel time, the principle holds identically: since it is always better to arrive at nodes as soon as possible, any ``shortcut'' between two nodes in a path can be spliced into the full path, thereby reducing its total travel time.
When the FIFO principle does not hold, or if link costs are different than travel time, this may not be true.
Figure~\ref{fig:bellmancounterexample} shows two counterexamples.
In the first, because the FIFO principle is violated and waiting is not allowed, the path segment $[1,2,3]$ has a shorter travel time between nodes 1 and 3 than path $[1,3]$.
However, when the link $(3,4)$ is added to these paths, the path $[1,2,3,4]$ has a higher travel time than path $[1,3,4]$ because link $(3,4)$ is entered at different times, and so its travel time is different.
In the second, path $[1,3]$ has a lower cost than path $[1,2,3]$, even though path $[1,3,4]$ has higher cost than path $[1,2,3,4]$.
In the second case, this is true even though arriving at any node later can only increase cost, and the same result holds even if waiting is allowed.

\stevefig{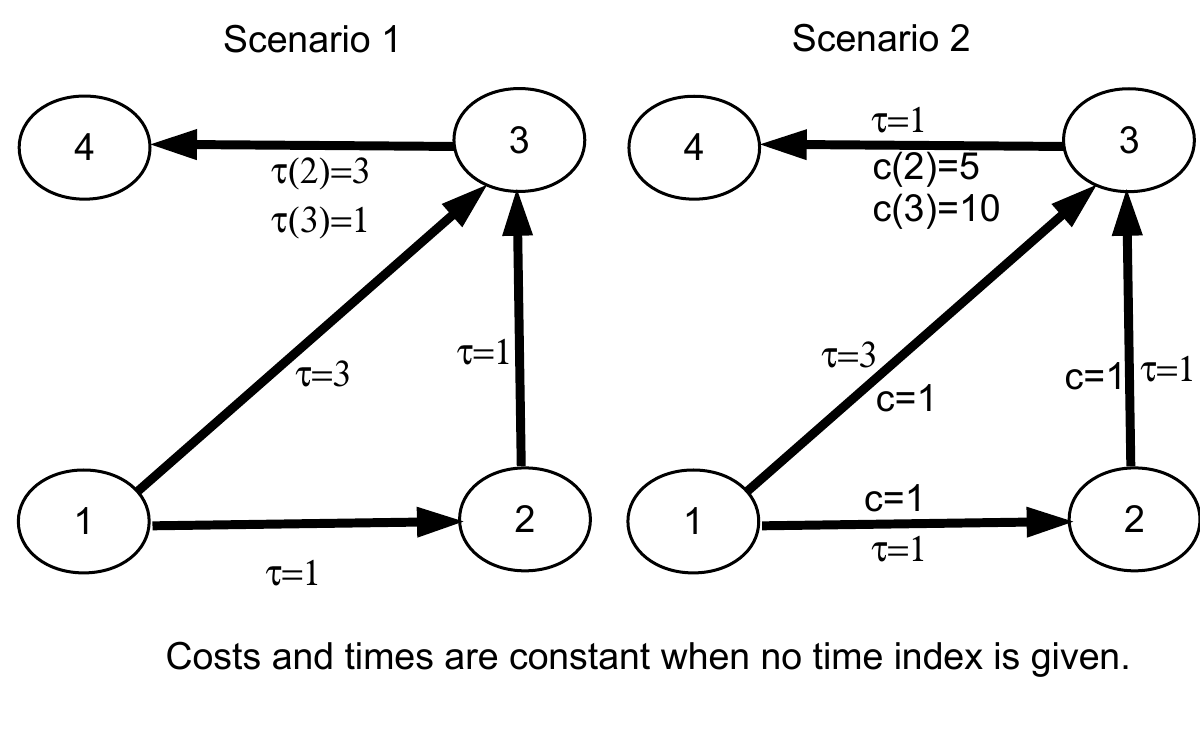}{Two counterexamples to the na\"{i}ve application of Bellman's principle in time-expanded networks.}{\textwidth}

In these cases, the correct approach is to apply Bellman's principle to the time-expanded network, an approach which also works in FIFO networks.
Bellman's principle \emph{does} apply to the time-expanded network, corresponding to the following principle in the physical network: \textbf{Let $\pi$ be a shortest path between nodes $r$ and $s$ when departing at time $t_0$.
If $i$ and $h$ are two nodes in this path, and if the arrival times at these nodes are $t_i$ and $t_h$, respectively, then the segment of $\pi$ between $i$ and $h$ must be a shortest path between these nodes when departing at $t_i$ and arriving at $t_h$.}  All of the algorithms in this chapter implicitly use this principle, by allowing us to construct shortest paths a single link at a time: if we already know the shortest path $\pi_i$ from the origin at $t_0$ to some other node $i$ at time $t_i$, then any other shortest path passing through node $i$ at time $t_i$ can be assumed to start with $\pi_i$.
\index{Bellman's principle!time-dependent}
\index{network!time-expanded|)}

\section{Time-Dependent Shortest Path Algorithms}
\label{sec:tdspalgos}

This section provides three algorithms for the time-dependent shortest path problem with fixed departure times.
In the first, time can be modeled as either discrete or continuous, but the network must satisfy the FIFO principle and the cost of a link must be its travel time.
In the second and third, time must be discrete, but FIFO need not hold and the cost of a link may take any value.
The second algorithm finds the time-dependent shortest path from a single origin and departure time to all destinations, while the third algorithm finds the time-dependent shortest paths from a single origin and \emph{all} departure times to a single destination.
There are many other possible variants of time-dependent shortest path algorithms, some of which are explored in the exercises --- in particular, Exercise~\ref{ex:alldeparturetimes} asks you to develop a time-dependent shortest path algorithm for all origins and departure times simultaneously, which often arises in dynamic traffic assignment software.
Nevertheless, the algorithms here should give the general flavor of how they function.
Which one of these algorithms is best, or whether another variant is better, depends on the particular dynamic traffic assignment implementation.
All of these algorithms use labels with similar (or even identical) names, but the meanings of these labels are slightly different in each.

\subsection{FIFO networks}
\label{sec:fifonets}

\index{shortest path!time-dependent!in FIFO networks|(}
For this algorithm to apply, assume that the FIFO principle holds, in either its discrete or continuous form.
Also assume that the cost of a link is simply its travel time, so $c_{ij}(t) = \tau_{ij}(t)$ for all links $(i,j)$ and times $t$.
We are also given the origin $r$ and departure time $t_0$.
Because the FIFO principle holds, waiting at an intermediate node is never beneficial, and so it suffices to find the earliest time we can reach a node and ignore all later times.
These properties allow us to adapt Dijkstra's algorithm (Section~\ref{sec:spgeneral}) for the static shortest path problem, a label-setting approach.
Two labels are defined for each node: $L_i$ gives the earliest possible arrival time to node $i$ found so far, when departing origin $r$ at time $t_0$.
The backnode label $q_i$ gives the previous node in a path corresponding to this earliest known arrival time.
By the dynamic version of Bellman's principle in FIFO networks, this suffices for reconstructing the shortest path from $r$ to $i$.
For nodes $i$ where we have not yet calculated the earliest possible arrival time, we will set $L_i$ to $\infty$, and for nodes $i$ where the backnode is meaningless (either because we have not yet found a path there, or because $i$ is the origin), we set $q_i$ to $-1$.
We also maintain a \emph{scan eligible list} $SEL$ of nodes.
The scan eligible list contains nodes we still have to examine before we can ensure we have found the time-dependent shortest paths from $i$.
In general, nodes can enter and leave $SEL$ multiple times.
For the specific case of FIFO networks, the modified version of Dijkstra's algorithm given below can guarantee that nodes enter and leave $SEL$ at most once, and that once a node has left $SEL$ its labels will never change again.
As a result, if you are only interested in finding a shortest path to one node, or to a subset of nodes, you can terminate as soon as every node you are interested in has left $SEL$.

The algorithm functions as follows:
\begin{enumerate}
\item Initialize by setting $L_r \leftarrow t_0$, $L_i \leftarrow \infty$ for $i \neq r$, and $\mb{q} \leftarrow \mb{-1}$.  Also initialize the scan eligible list to contain the origin only: $SEL \leftarrow \myc{r}$.
\item Choose a node $i \in SEL$ with minimum $L_i$ value, and delete it from the list.
\item Scan node $i$.  For every link $(i,j)$ whose tail is node $i$, compute the time at which you would arrive at node $j$ if you followed the shortest path to $i$, and then used link $(i,j)$.  This arrival time is $L_i + \tau_{ij}(L_i)$.  If this arrival time is within the time horizon, update $L_j \leftarrow L_i + \tau_{ij}(L_i)$, $q_j \leftarrow i$, and add $j$ to $SEL$.
\item If $SEL$ is empty, then terminate.  Otherwise, return to step 2.
\end{enumerate}

As an example of this algorithm, consider the network in Figure~\ref{fig:fifoalgodemo}.
The time-dependent travel times are shown in this figure.
The FIFO assumption is satisfied: for links $(1,3)$ and $(2,4)$ the travel times are constant; for links $(1,2)$ and $(3,4)$ the travel times are increasing (so arriving earlier always means leaving earlier); and for link $(2,3)$ the travel time is decreasing, but at a slow enough rate that you cannot leave earlier by arriving later, cf.\ equation~\eqn{fifocontinuous}.
Assume that the initial departure time from node 1 is at $t_0 = 2$, and that the time horizon is large enough that $L_i + \tau_{ij}(L_i)$ is always within the time horizon whenever step 3 is encountered.
The steps of the algorithm are explained below, and summarized in Table~\ref{tbl:fifoalgodemo}.
The first row of the table (iteration zero) shows the state of the algorithm \emph{just after} step 1 is performed.
The remaining rows of the table show the state of the algorithm \emph{just before} step 4 is performed, increasing the iteration number each time we reach this step.

Initially, all cost labels are initialized to $\infty$ (except for the origin, which is assigned 2, the departure time), all backnode labels are initialized to $-1$, and the only node in the scan list is the origin (node 1).
The node in $SEL$ with the least $L$ value is node 1 (indeed it is the only node in the list), which is selected as the node $i$ to scan.
At the current time of 2, link $(1,2)$ has a travel time of 6, and link $(1,3)$ has a travel time of 10.
Following these links would result in arrival at nodes 2 and 3 at times 8 and 12, respectively.
Each of these is less than their current values of $\infty$, so the cost and backnode labels are adjusted accordingly.
At the next iteration, node 2 is the node in $SEL$ with the least $L$ value, so $i = 2$.
At time 8, link $(2,3)$ has a travel time of 1, and link $(2,4)$ has a travel time of 5.
Following these links, one would arrive at nodes 3 and 4 at times 9 and 13, respectively.
Both of these values are less than the current $L$ values for these nodes, so their $L$ and $q$ labels are changed.
At the next iteration, node 3 is the node in $SEL$ with the least $L$ value, so $i = 3$.
At this time, link $(3,4)$ would have a travel time of $4 \frac{1}{2}$, and choosing it means arriving at node 4 at time $13 \frac{1}{2}$.
This is greater than the current value ($L_4 = 13$), so no labels are adjusted.
Finally, node 4 is chosen as the only node in $SEL$.
Since it has no outgoing links ($\Gamma(4)$ is empty), there is nothing to do in step 4, and since all nodes are finalized the algorithm terminates.

At this point, we can trace back the shortest paths using the backnode labels: the shortest paths to nodes 2, 3, and 4 are $[1, 2]$, $[1,2,3]$, and $[1,2,4]$, respectively; and following these paths one arrives at the nodes at times 8, 9, and 13.\index{shortest path!time-dependent!in FIFO networks|)}

\begin{figure}
\begin{center}
\begin{tikzpicture}[->,>=stealth',shorten >=1pt,auto,node distance=3cm,
thick,main node/.style={circle,,draw}]

\node[main node] (1L) {$1$};
\node[main node] (2L) [below right=1.5cm and 3cm of 1L] {$2$};
\node[main node] (3L) [above right=1.5cm and 3cm of 1L] {$3$};
\node[main node] (4L) [below right=1.5cm and 3cm of 3L] {$4$};      

\path (1L) edge node[above left] {$10$} (3L);
\path (1L) edge node[below left] {$4 + t$} (2L);
\path (2L) edge node[right] {$\max \myc{5 - \frac{t}{2}, 1 }$} (3L);      
\path (2L) edge node[below right] {$5$} (4L);
\path (3L) edge node[above right] {$\frac{t}{2}$} (4L);      

\end{tikzpicture}
\caption{Network and travel time functions for demonstrating the FIFO time-dependent shortest path algorithm.
\label{fig:fifoalgodemo}}
\end{center}
\end{figure}

\begin{table}
\begin{center}
\caption{FIFO time-dependent shortest path algorithm for the network in Figure~\ref{fig:fifoalgodemo}, departing node 1 at $t_0 = 2$.
\label{tbl:fifoalgodemo}}
\begin{tabular}{|c|c|cccc|cccc|c|}
\hline
Iteration & $i$ & $L_1$ & $L_2$ & $L_3$ & $L_4$ & $q_1$ & $q_2$ & $q_3$ & $q_4$ & $SEL$\\
\hline
0         & --- & 2 & $\infty$ & $\infty$ & $\infty$ & $-1$ & $-1$ & $-1$ & $-1$ & $\{ 1 \}$ \\
1         & 1   & 2 & 8 & 12 & $\infty$ & $-1$ & 1 & 1 & $-1$ & $\{2, 3 \}$ \\
2         & 2 & 2 & 8 & 9 & 13 & $-1$ & 1 & 2 & 2 & $\{ 3, 4 \}$ \\
3         &  3 & 2 & 8 & 9 & 13 & $-1$ & 1 & 2 & 2 & $\{ 4 \}$ \\
4         & 4 & 2 & 8 & 9 & 13 & $-1$ & 1 & 2 & 2 & $\emptyset$ \\
\hline
\end{tabular}
\end{center}
\end{table}

\subsection{Discrete-time networks, one departure time}
\label{sec:iotalgo}

\index{shortest path!time-dependent!in general networks|(}
This algorithm applies in any discrete-time network, regardless of whether or not the FIFO principle holds, and regardless of whether $c_{ij}(t) = \tau_{ij}(t)$ or not.
We are given the origin $r$ and departure time $t_0$, and work in the time-expanded network.
Each node $i:t$ in the time-expanded network is associated with two labels.
The label $L_i^t$ denotes the cost of the shortest path from $r$ to $i$ known so far, when departing at time $t_0$ and arriving at time $t$.
The backnode label $q_i^t$ provides the previous node on a path corresponding to cost $L_i^t$.
As before, $L_i^t = \infty$ signifies that no path to $i$ at time $t$ is yet known, and $q_i^t = -1$ signifies that the backnode is meaningless.
Since the FIFO principle may not hold, it is not always advantageous to arrive at a node as early as possible.
To reflect this, we work in the time-expanded network, where we can naturally identify the best time to arrive at nodes.
We presume that all links have strictly positive travel time, so that the time-expanded graph is acyclic and each link connects a node of earlier time to a node of later time.
If there are zero-travel time links in the physical network, but no cycles of zero-travel time links, then step 4 is assumed to proceed in topological order by physical node.
In this case, we can adapt the algorithm used to find static shortest paths in acyclic networks from Section~\ref{sec:acyclicsp}, applying it to the time-expanded network in the following way:
\begin{enumerate}
\item Initialize $L_r^{t_0} \leftarrow 0$, $L_i^t \leftarrow \infty$ for all node-time combinations except for $r : t_0$, and $\mb{q} \leftarrow \mb{-1}$.  
\item Initialize the current time to the departure time, $t \leftarrow t_0$
\item For each time-expanded node $i:t$ for which $L_i^t < \infty$, and for each time-expanded link $(i:t, j:t')$, perform the following steps:
\begin{enumerate}
  \item Set $L_j^{t'} \leftarrow \min \myc{ L_j^{t'} , L_i^t + c_{ij}(t)}$.
  \item If $L_j^{t'}$ changed in the previous step, update $q_j^{t'} \leftarrow i:t$.
\end{enumerate}
\item If $t = T$, then terminate.
Otherwise, move to the next time step ($t \leftarrow t + 1$) and return to step 4.
\end{enumerate}

At the conclusion of this algorithm, we have the least-cost paths for each possible arrival time at each destination.
To find the least-cost path to a particular destination $s$ (at any arrival time), you can consult the $L_s^t$ labels at all times $t$, and trace back the path for the arrival time $t$ with the least $L_s^t$ value.

This algorithm is demonstrated on the network in the \emph{right} panel of Figure~\ref{fig:bellmancounterexample}, and its progress is summarized in Table~\ref{tbl:alg2demo}.
The table shows the state of the algorithm \emph{just before} step 4 is executed.
For brevity, this table only reports $L_i^t$ and $q_i^t$ labels for nodes and arrival times which are reachable in the network (that is, $i$ and $t$ values for which $L_i^t < \infty$ at the end of the algorithm).
All other cost and backnode labels are at $\infty$ and $-1$ throughout the entire duration of the algorithm.

Initially, all cost labels are set to $\infty$ and all backnode labels to $-1$, except for the origin and departure time: $L_1^0 = 0$.
The algorithm then sets $t = 0$, and scans over all physical nodes which are reachable at this time.\footnote{Reachability is expressed by the condition $L_i^t < \infty$; a node which cannot be reached at this time will still have its label set to the initial value.
Any node which \emph{can} be reached at this time will have a finite $L_i^t$ value, since step 4a will always reduce an infinite value.}  Only node 1 can be reached at this time, and the possible links are $(1:0,2:1)$ and $(1:0,3:3)$.
Following either link incurs a cost of 1, which is lower than the (infinite) values of $L_2^1$ and $L_3^3$, so the cost and backnode labels are updated.

The algorithm then sets $t = 1$.
Only node 2 is reachable at this time, and the only link is $(2:1, 3:2)$.
Following this link incurs a cost of 1; in addition to the cost of 1 already involved in reaching node 2, this gives a cost of 2 for arriving at node 3 at time 2.
The cost and time labels for $3:2$ are updated.
Since the costs and times are different, notice that arriving at node 3 at a later time (3 vs.\ 2) incurs a lower cost (1 vs.\ 2).
This is why we need to track labels for different arrival times, unlike the algorithm in the previous section.

The next time step is $t = 2$.
Only node 3 is reachable at this time (from the path $[1,2,3]$), and the only link is $(3:2,4:3)$.
Following this link incurs a cost of 5, resulting in a total cost of 7, and the labels for $4:3$ are updated.
Time $t = 3$ is next, and again only node 3 is reachable at this time --- but from the path $[1,3]$.
The only link is $(3:3,4:4)$, and following this link incurs a cost of 10, for a total cost of 11.
Labels are updated for $4:4$.
There are no further label changes in the algorithm (all nodes have already been scanned at all reachable times), and it terminates as soon as $t$ is increased to the time horizon.

After termination, the $L_3^t$ labels show that we can reach node 3 either with a cost of 1 (arriving at time 3) or a cost of 2 (arriving at time 2).
The least-cost path thus arrives at time 3, and it is $[1,3]$.
The $L_4^t$ labels show that we can reach node 4 either with a cost of 7 (arriving at time 3) or a cost of 11 (arriving at time 4).
The least-cost path arrives at time 3, and it is $[1,2,3,4]$.
(The least-cost path to node 2 is $[1,2]$, since there is only one possible arrival time there).
Notice that the na\"{i}ve form of Bellman's principle is not satisfied: the least-cost path to node 2 is \emph{not} a subset of the least-cost path to node 3.
This was why we needed to keep track of different possible arrival times to nodes --- the least-cost path to node 3 \emph{arriving at time 2} is indeed a subset of the least-cost path to node 4 \emph{arriving at time 3}.
\index{shortest path!time-dependent!in general networks|)}

\begin{table}
\begin{center}
\caption{Discrete time-dependent shortest path algorithm for the network in the right panel of Figure~\ref{fig:bellmancounterexample}, departing node 1 at $t_0 = 0$, showing labels $L_i^t$ and $q_i^t$.
\label{tbl:alg2demo}}
\begin{tabular}{|c|cccccc|cccccc|}
\hline
Iteration & $L_1^0$ & $L_2^1$ & $L_3^2$ & $L_3^3$ & $L_4^3$ & $L_4^4$ & $q_1^0$ & $q_2^1$ & $q_3^2$ & $q_3^3$ & $q_4^3$ & $q_4^4$ \\
\hline
 0 & 0 & $\infty$ & $\infty$ & $\infty$ & $\infty$ & $\infty$ & $-1$ & $-1$ & $-1$ & $-1$ & $-1$ & $-1$\\
 1 & 0 & 1 & $\infty$ & 1 &  $\infty$ & $\infty$  & $-1$ & 1:0 & $-1$ & 1:0 & $-1$ & $-1$ \\
 2 & 0 & 1 & 2 & 1 &  $\infty$ & $\infty$ & $-1$ & 1:0 & 2:1 & 1:0 & $-1$ & $-1$ \\
 3 & 0 & 1 & 2 & 1 & 7 & $\infty$ & $-1$ & 1:0 & 2:1 & 1:0 & 3:2 & $-1$ \\
 4 & 0 & 1 & 2 & 1 & 7 & 11 & $-1$ & 1:0 & 2:1 & 1:0 & 3:2 & 3:3 \\
\hline
\end{tabular}
\end{center}
\end{table}

\section{Departure Time Choice}
\label{sec:departuretimechoice}

\index{departure time choice|(}
Travelers often have some flexibility when choosing their departure or arrival times, and may choose to leave earlier or later in order to minimize cost.
For example, commuters with flexible work hours may want to time their commutes to avoid congestion, or when dynamic congestion charges are lower.
Both flexible departures and arrivals can be incorporated into the shortest path algorithms described in the previous section --- in fact, you may have already noticed that the algorithms in Sections~\ref{sec:fifonets} and~\ref{sec:iotalgo} allow flexible arrival times.
This section explores these choices more systematically, showing how the departure time can also be made flexible, and ways to represent different departure time behaviors.

\subsection{Artificial origins and destinations}
\label{sec:artificialod}

\index{node!super-origin}
\index{node!super-destination}
In the time-expanded network, departure and arrival time choice can be modeled by adding artificial ``super-origin'' and ``super-destination'' nodes which reflect the start and end of a trip without regard to the time.
A super-origin is connected to time-expanded nodes which correspond to allowable departure times, and a super-destination is connected to time-expanded nodes which correspond to allowable arrival times.
Initially, we will assign these links a cost of zero, which means the traveler is indifferent among any of these departure or arrival times.
Section~\ref{sec:penalties} will describe how arrival and departure time preferences can be modeled.

Figure~\ref{fig:supersourcesink} shows how super-origins and super-destinations can be added to the physical network of Figure~\ref{fig:timeexpandednet2}(a), assuming that nodes $i$ and $k$ are origins, and nodes $j$ and $l$ are destinations.
The figure is drawn as if departures were allowed only for times 0, 1, and 2, but arrivals are allowed at any time.
The network is similar to the time-expanded network in Figure~\ref{fig:timeexpandednet2}(b), but now includes four artificial nodes, and artificial links corresponding to allowable departure and arrival times.

\stevefig{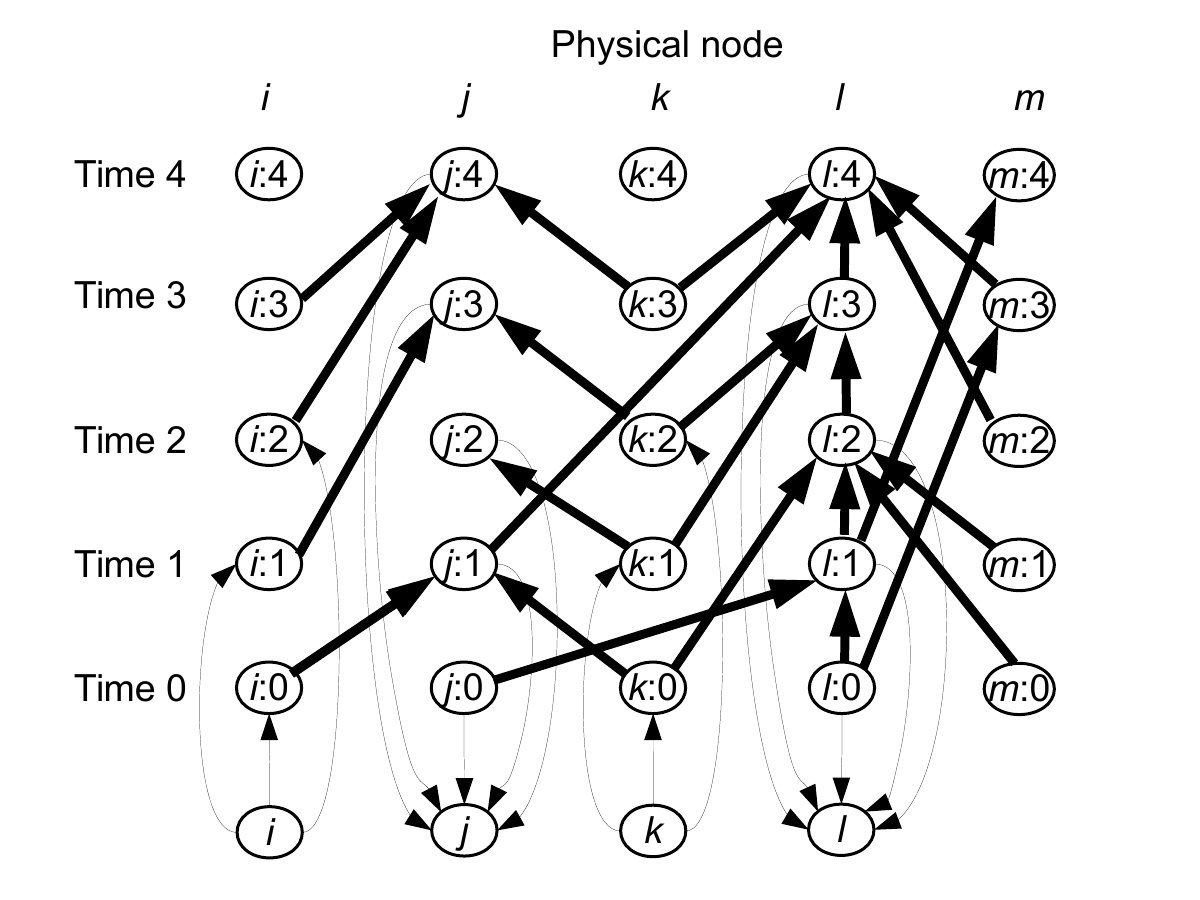}{Time-expanded network with artificial origins and destinations.
Artificial links shown with light arrows for clarity.}{0.8\textwidth}

\index{shortest path!time-dependent!departure time choice|(}
Once these artificial links and nodes are added to the network, the algorithm from Section~\ref{sec:iotalgo} can be applied directly, with only very minor changes.
In what follows, $r$ is the origin and the algorithm finds least-cost paths to all destinations, for any allowable departure and arrival time:
\begin{enumerate}
\item Initialize $\mb{q} \leftarrow \mb{-1}$ and $\mb{L} \leftarrow \infty$ for all nodes in the time-expanded network.
\item Set $L_r \leftarrow 0$ for the super-origin $r$, and for each artificial link $(r, r : t)$ set $L_r^t \leftarrow 0$ and $q_r^t \leftarrow r$.
\item Initialize the current time $t$ to the earliest possible departure time (the lowest $t$ index for which an artificial link $(r, r:t)$ exists).
\item For each time-expanded node $i:t$ for which $L_i^t < \infty$, and for each time-expanded link $(i:t, j:t')$, perform the following steps:
\begin{enumerate}
  \item Set $L_j^{t'} \leftarrow \min \myc{ L_j^{t'} , L_i^t + c_{ij}(t)}$.
  \item If $L_j^{t'}$ changed in the previous step, update $q_j^{t'} \leftarrow i:t$.
\end{enumerate}
Likewise, for each artificial link $(i:t, i)$ reaching a super-destination node $i$, perform the following steps:
\begin{enumerate}
  \item Set $L_i \leftarrow \min \myc{ L_i , L_i^t(t)}$.
  \item If $L_i$ changed in the previous step, update $q_i \leftarrow i:t$.
\end{enumerate}   
\item If $t = T$, then terminate.
Otherwise, move to the next time step ($t \leftarrow t + 1$) and return to step 4.
\end{enumerate}
The algorithm initializes labels differently in step 2; step 3 starts at the earliest possible departure time rather than the fixed time $t_0$; and step 4 is expanded to update labels both at adjacent time-expanded nodes and super-destinations.
All other steps work in the same way.

To demonstrate this algorithm, consider the network in Figure~\ref{fig:alltimesalgodemo}, where the time horizon is $T = 20$ and the destination is node 4, and where the cost of a link is equal to its travel time.
Assume that waiting at intermediate nodes is not allowed.
Table~\ref{tbl:alltimesalgodemo} shows the cost and backnode labels at the conclusion of the algorithm.
Each iteration of the algorithm generates one row of this table, starting with $T = 1$ and working up to $T = 20$.
Whenever $L_i^t = \infty$ is seen in Table~\ref{tbl:alltimesalgodemo}, there is no way to arrive at node $i$ at time $t$ given the time discretization and travel time functions.

Table~\ref{tbl:alltimesalgodemo} also shows the labels for the super-origin and super-destination, below the labels for the time-expanded nodes.
The backnode label for the super-destination tells us that the least-cost path arrives at node 4 at time 6; the $q_4^6$ label tells us the least-cost path there comes through node 3 at time 1; $q_3^1$ tells us the least-cost path there comes through node 1 at time 0, and $q_1^0$ brings us to the super-origin.
Therefore, we should depart the origin at time 0, and follow the path $[1,3,4]$ to arrive at the destination at time 6, with a total cost of $c_4 = 6$.
\index{shortest path!time-dependent!departure time choice|)}

\begin{figure}
\begin{center}
\begin{tikzpicture}[->,>=stealth',shorten >=1pt,auto,node distance=3cm,
thick,main node/.style={circle,,draw}]

\node[main node] (1L) {$1$};
\node[main node] (2L) [below right=1.5cm and 3cm of 1L] {$2$};
\node[main node] (3L) [above right=1.5cm and 3cm of 1L] {$3$};
\node[main node] (4L) [below right=1.5cm and 3cm of 3L] {$4$};      

\path (1L) edge node[above left] {$1 + 2t$} (3L);
\path (1L) edge node[below left] {$5$} (2L);
\path (2L) edge node[right] {$\max \myc{10 - t, 0 }$} (3L);      
\path (2L) edge node[below right] {$1 + t$} (4L);
\path (3L) edge node[above right] {$5$} (4L);      

\end{tikzpicture}
\caption{Network and travel time functions for demonstrating the all departure times time-dependent shortest path algorithm.
\label{fig:alltimesalgodemo}}
\end{center}
\end{figure}

\begin{table}
\begin{center}
\caption{Labels at conclusion of time-dependent shortest path algorithm with departure time choice, to node 4 in Figure~\ref{fig:alltimesalgodemo}.
\label{tbl:alltimesalgodemo}}
\begin{tabular}{|c|cccc|cccc|}
\hline
$t$ & $L_1^t$ & $L_2^t$ & $L_3^t$ & $L_4^t$ & $q_1^t$ & $q_2^t$ & $q_3^t$ & $q_4^t$ \\ 
\hline
20 & $0$ & $5$ & $5$ & $10$  & $1$  & $1:15$  & $2:20$  & $3:15$ \\
19 & $0$ & $5$ & $5$ & $10$  & $1$  & $1:14$  & $2:19$  & $3:14$ \\
18 & $0$ & $5$ & $5$ & $10$  & $1$  & $1:13$  & $2:18$  & $3:13$ \\
17 & $0$ & $5$ & $5$ & $10$  & $1$  & $1:12$  & $2:17$  & $3:12$ \\
16 & $0$ & $5$ & $5$ & $10$  & $1$  & $1:11$  & $2:16$  & $3:11$ \\
15 & $0$ & $5$ & $5$ & $10$  & $1$  & $1:10$  & $2:15$  & $3:10$ \\
14 & $0$ & $5$ & $5$ & $\infty$  & $1$  & $1:9$  & $2:14$  & $-1$ \\
13 & $0$ & $5$ & $5$ & $12$  & $1$  & $1:8$  & $2:13$  & $2:6$ \\
12 & $0$ & $5$ & $5$ & $10$  & $1$  & $1:7$  & $2:12$  & $3:7$ \\
11 & $0$ & $5$ & $5$ & $11$  & $1$  & $1:6$  & $2:11$  & $2:5$ \\
10 & $0$ & $5$ & $5$ & $\infty$  & $1$  & $1:5$  & $2:10$  & $-1$ \\
9 & $0$ & $5$ & $\infty$ & $8$  & $1$  & $1:4$  & $-1$  & $3:4$ \\
8 & $0$ & $5$ & $\infty$ & $\infty$  & $1$  & $1:3$  & $-1$  & $-1$ \\
7 & $0$ & $5$ & $5$ & $\infty$  & $1$  & $1:2$  & $1:2$  & $-1$ \\
6 & $0$ & $5$ & $\infty$ & $6$  & $1$  & $1:1$  & $-1$  & $3:1$ \\
5 & $0$ & $5$ & $\infty$ & $\infty$  & $1$  & $1:0$  & $-1$  & $-1$ \\
4 & $0$ & $\infty$ & $3$ & $\infty$  & $1$  & $-1$  & $1:1$  & $-1$ \\
3 & $0$ & $\infty$ & $\infty$ & $\infty$  & $1$  & $-1$  & $-1$  & $-1$ \\
2 & $0$ & $\infty$ & $\infty$ & $\infty$  & $1$  & $-1$  & $-1$  & $-1$ \\
1 & $0$ & $\infty$ & $1$ & $\infty$  & $1$  & $-1$  & $1:0$  & $-1$ \\
0 & $0$ & $\infty$ & $\infty$ & $\infty$  & $1$  & $-1$  & $-1$  & $-1$ \\
\hline
\end{tabular}
\end{center}
$(L,q)$ labels at super-origin 1: $(0, -1)$ \\
$(L,q)$ labels at super-destination 4: $(6, 4:6)$
\end{table}

\subsection{Arrival and departure time preferences}
\label{sec:penalties}

While there may be flexibility in departure or arrival time, travelers are usually not completely indifferent about when they depart or arrive.
For instance, there may be a well-defined arrival deadline (start of work, or check-in time before a flight), and a strong desire to arrive before this deadline.
Departing extremely early would ensure arriving before the deadline, but carries opportunity costs (by departing later, the traveler would have more time to do other things).
Both situations can be modeled by attaching a cost to when a traveler departs the origin, and to when they arrive at the destination.

A common way to model arrival costs is with the \emph{schedule delay}\index{schedule delay} concept.
In this model, travelers have a preferred arrival time $t^*$\label{not:tstar} at the destination, and arriving either earlier or later than $t^*$ is undesirable.
The most general formulation involves a function $f(t)$\label{not:fdelay} denoting the cost (or disutility) of arriving at the destination at time $t$.
This function is typically convex and has a minimum at $t^*$.
One such function is
\labeleqn{vickreydelay}{f(t) = \beta [t^* - t]^+ + \gamma [t - t^*]^+\,,}
where $[\cdot]^+ = \max \myc{ \cdot, 0 }$ expresses the positive part of the quantity in brackets.
In equation~\eqn{vickreydelay}, the first bracketed term then represents the amount by which the traveler arrived early, compared to the preferred time, and the second bracketed term represents the amount by which the traveler arrived late.
The coefficients $\beta$\label{not:betasd} and $\gamma$\label{not:gammasd} then weight these terms and convert them to cost units; generally $\beta < \gamma$ to reflect the fact that arriving early by a certain amount of time, while undesirable, is usually not as bad as arriving late by that same amount of time.

Another possible function is nonlinear, taking a form such as
\labeleqn{deviance}{f(t) = \beta ([t^* - t]^+)^2 + \gamma ([t - t^*]^+)^2\,.}
In this function, the penalty associated with early or late arrival grows faster and faster the farther the arrival time from the target.
This may occur if, for instance, being ten minutes late is more than ten times as bad as being one minute late.
Special cases of these functions arise when $\beta = \gamma$ (the function becomes symmetric), or when $\beta = 0$ (there is no penalty for early arrival, but only for late arrival).
This function is also differentiable everywhere, in contrast to equation~\eqn{vickreydelay} which is not differentiable at $t^*$.
For certain algorithms this may be advantageous.

To each of these schedule delay functions $f(t)$, one can add the cost of the path arriving at time $t$ (the sum of the link costs $c_{ij}$ along the way) to yield the total cost of travel.
We assume that travelers will choose both the departure time and the path to minimize this sum.
The algorithm from Section~\ref{sec:artificialod} can be used to find both this ideal departure time and the path.
The only change is that artificial links connecting time-expanded destination nodes $(s:t)$ to super-destinations $s$ now have a cost of $f(t)$, rather than zero.
At termination, the departure time $t^*_0$ minimizing $L_r^{t^*_0}$ corresponds to the least total cost, and the backnode labels trace out the path.\index{shortest path!time-dependent!departure time choice|(}
\begin{enumerate}
\item Initialize $\mb{q} \leftarrow \mb{-1}$ and $\mb{L} \leftarrow \infty$ for all nodes in the time-expanded network.
\item Set $L_r \leftarrow 0$ for the super-origin $r$, and for each artificial link $(r, r : t)$ set $L_r^t \leftarrow 0$ and $q_r^t \leftarrow r$.
\item Initialize the current time $t$ to the earliest possible departure time (the lowest $t$ index for which an artificial link $(r, r:t)$ exists).
\item For each time-expanded node $i:t$ for which $L_i^t < \infty$, and for each time-expanded link $(i:t, j:t')$, perform the following steps:
\begin{enumerate}
  \item Set $L_j^{t'} \leftarrow \min \myc{ L_j^{t'} , L_i^t + c_{ij}(t)}$.
  \item If $L_j^{t'}$ changed in the previous step, update $q_j^{t'} \leftarrow i:t$.
\end{enumerate}
Likewise, for each artificial link $(i:t, i)$ reaching a super-destination node $i$, perform the following steps:
\begin{enumerate}
  \item Set $L_i \leftarrow \min \myc{ L_i , L_i^t + f(t)}$.
  \item If $L_i$ changed in the previous step, update $q_i \leftarrow i:t$.
\end{enumerate}   
\item If $t = T$, then terminate.
Otherwise, move to the next time step ($t \leftarrow t + 1$) and return to step 4.
\end{enumerate}

To demonstrate this algorithm, again consider the network in Figure~\ref{fig:alltimesalgodemo}, but with the arrival time penalty function
\labeleqn{scheduledelayalgodemo}{f(t) = [10 - t]^+ + 2[t - 10]^+\,,}
which suggests that the traveler wishes to arrive at time 10, and that late arrival is twice as costly as early arrival.

As before, the time horizon is $T = 20$ and the destination is node 4, and the cost of each link is equal to its travel time.
Table~\ref{tbl:scheduledelayalgodemo} shows the cost and backnode labels at the conclusion of the algorithm, which runs in exactly the same way as before except that the artificial destination links $(4:t,4)$ have a cost equal to $f(t)$.
At the conclusion of the algorithm, we can identify the total cost on the shortest path from node 1 to node 4, now including the penalty for arriving early or late at the destination.

Table~\ref{tbl:scheduledelayalgodemo} also shows the labels for the super-origin and super-destination, below the labels for the time-expanded nodes.
The backnode label for the super-destination tells us that the least-cost path arrives at node 4 at time 9; the $q_4^9$ label tells us the least-cost path there comes through node 3 at time 4; $q_3^4$ tells us the least-cost path there comes through node 1 at time 1, and $q_1^1$ brings us to the super-origin.
Therefore, we should depart the origin at time 1, and follow the path $[1,3,4]$ to arrive at the destination at time 9, with a total cost of $c_4 = 9$.
Of this cost, 8 units are due to travel time (difference between arrival and departure times), and 1 unit is due to the arrival time penalty from equation~\eqn{scheduledelayalgodemo} with $t = 9$.

Departing earlier, at $t = 0$, on the same path would reduce the travel time to 6, but increase the early arrival penalty to 4.
The total cost of leaving at $t = 0$ is thus higher than departing one time step later.
Leaving at $t = 2$ and following the same path increases the travel time cost to 10.
Since this means arriving at time 12, there is a late penalty cost of 4 added, resulting in a total travel cost of 14.
By comparing all possible paths and departure times, you can verify that it is impossible to have a total cost less than 9.
\index{shortest path!time-dependent!departure time choice|)}

\begin{table}
\begin{center}
\caption{Labels at conclusion of arrival-penalty time-dependent shortest path, to node 4 in Figure~\ref{fig:alltimesalgodemo}.
\label{tbl:scheduledelayalgodemo}}
\begin{tabular}{|c|cccc|cccc|}
\hline
$t$ & $L_1^t$ & $L_2^t$ & $L_3^t$ & $L_4^t$ & $q_1^t$ & $q_2^t$ & $q_3^t$ & $q_4^t$ \\   
\hline                                                                                               
20 & $0$ & $5$ & $5$ & $10$  & $1$  & $1:15$  & $2:20$  & $3:15$ \\                                  
19 & $0$ & $5$ & $5$ & $10$  & $1$  & $1:14$  & $2:19$  & $3:14$ \\                                  
18 & $0$ & $5$ & $5$ & $10$  & $1$  & $1:13$  & $2:18$  & $3:13$ \\                                  
17 & $0$ & $5$ & $5$ & $10$  & $1$  & $1:12$  & $2:17$  & $3:12$ \\                                  
16 & $0$ & $5$ & $5$ & $10$  & $1$  & $1:11$  & $2:16$  & $3:11$ \\                                  
15 & $0$ & $5$ & $5$ & $10$  & $1$  & $1:10$  & $2:15$  & $3:10$ \\                                  
14 & $0$ & $5$ & $5$ & $\infty$  & $1$  & $1:9$  & $2:14$  & $-1$ \\                                 
13 & $0$ & $5$ & $5$ & $12$  & $1$  & $1:8$  & $2:13$  & $2:6$ \\                                    
12 & $0$ & $5$ & $5$ & $10$  & $1$  & $1:7$  & $2:12$  & $3:7$ \\                                    
11 & $0$ & $5$ & $5$ & $11$  & $1$  & $1:6$  & $2:11$  & $2:5$ \\                                    
10 & $0$ & $5$ & $5$ & $\infty$  & $1$  & $1:5$  & $2:10$  & $-1$ \\                                 
9 & $0$ & $5$ & $\infty$ & $8$  & $1$  & $1:4$  & $-1$  & $3:4$ \\                                   
8 & $0$ & $5$ & $\infty$ & $\infty$  & $1$  & $1:3$  & $-1$  & $-1$ \\                               
7 & $0$ & $5$ & $5$ & $\infty$  & $1$  & $1:2$  & $1:2$  & $-1$ \\                                   
6 & $0$ & $5$ & $\infty$ & $6$  & $1$  & $1:1$  & $-1$  & $3:1$ \\                                   
5 & $0$ & $5$ & $\infty$ & $\infty$  & $1$  & $1:0$  & $-1$  & $-1$ \\                               
4 & $0$ & $\infty$ & $3$ & $\infty$  & $1$  & $-1$  & $1:1$  & $-1$ \\                               
3 & $0$ & $\infty$ & $\infty$ & $\infty$  & $1$  & $-1$  & $-1$  & $-1$ \\                           
2 & $0$ & $\infty$ & $\infty$ & $\infty$  & $1$  & $-1$  & $-1$  & $-1$ \\                           
1 & $0$ & $\infty$ & $1$ & $\infty$  & $1$  & $-1$  & $1:0$  & $-1$ \\                               
0 & $0$ & $\infty$ & $\infty$ & $\infty$  & $1$  & $-1$  & $-1$  & $-1$ \\                           
\hline  
\end{tabular}
\end{center}
$(L,q)$ labels at super-origin 1: $(0, -1)$ \\
$(L,q)$ labels at super-destination 4: $(9, 4:9)$
\end{table}

Considering costs associated with departure time, rather than arrival time, is done in essentially the same way, by assigning a nonzero cost to the artificial links connecting the super-origin $r$ to the time-expanded nodes $r(t)$.
It is thus possible to have only departure time penalties, only arrival time penalties, both, or neither, depending on whether the artificial origin links and artificial destination links have nonzero costs.
\index{departure time choice|)}

\section{Dynamic $A^*$}

\index{A$^*$!dynamic|(}
Because the time-expanded network transforms the time-dependent shortest problem into the classical, static shortest path problem from Section~\ref{sec:shortestpath}, any of the algorithms described there can be applied.
The time-expanded network is acyclic, so it is fastest to use the algorithm from Section~\ref{sec:acyclicsp} --- and indeed this is all that the algorithm in Section~\ref{sec:iotalgo} is, using the time labels as a topological order.

Section~\ref{sec:astar} also presented the $A^*$ algorithm, which provides a single path from one origin to one destination, rather than all shortest paths from one origin to all destinations, or all origins to one destination.
By focusing on a single origin and destination, $A^*$ can often find a shortest path much faster than a one origin-to-all destinations algorithm.
The tradeoff is that the algorithm has to repeated many times, once for every OD pair, rather than once for every origin or destination.
In static assignment, one-to-all or all-to-one algorithms are preferred because, because in many cases, running $A$* for each OD pair takes more time than running a one-to-all algorithm for each origin.

In dynamic traffic assignment with fixed departure times, however, the number of ``origins'' in the time-expanded network is multiplied by the number of departure times.
If one were to write a full time-dependent OD matrix, the number of entries in this matrix is very large: in a network with 1000 centroids and 1000 time steps, there are 1 \emph{billion} entries, one for every origin, every destination, and departure time.
This is much larger than the number of vehicles that will be assigned, so almost every entry in this matrix will be zero.
In such cases, $A^*$ can work much better, only being applied to origins, destinations, and departure times with a positive entry in the matrix.

As discussed in Section~\ref{sec:astar}, an effective estimate for $A^*$ in traffic assignment problems is to use the free-flow travel costs.
As a preprocessing step at the start of traffic assignment, you can use an all-to-one static shortest path algorithm to find the least-cost travel cost $g^s_j$ from every node $j$ to every destination $s$ at free flow.
For the remainder of the traffic assignment algorithm, you can then use $g^s_j$ as the estimates for $A^*$.
This is quite effective in practical networks.
\index{A$^*$!dynamic|)}

\section{Vickrey's bottleneck (*)}
\label{sec:singlebottleneck}

\index{Vickrey bottleneck|(}
\index{departure time choice!equilibrium|(}
\emph{(The Vickrey bottleneck is a classic transportation science model describing departure-time choice equilibrium on a single link.
It is not strictly related to the rest of the book's content on dynamic traffic assignment in larger networks, but it is worthy of study in its own right.)}

Two kinds of travel choices are characteristic of dynamic traffic assignment: the route choice (a spatial decision), and the departure time choice (a temporal decision).
Static assignment models only involve the spatial decision and have no concept of time.
The Vickrey bottleneck is a transportation model which only involves the \emph{temporal} decision, and there is no notion of route choice.
Imagine a single link represented by the point queue model of Section~\ref{sec:pointqueue}, with a downstream capacity of $s$\label{not:svickrey} and an infinite upstream capacity.
A total of $N$\label{not:Nvickrey} travelers will travel on this link; they choose the time at which they enter the link, and the point queue model determines the time at which they exit and arrive at the destination.\footnote{Unfortunately, the words \emph{departure} and \emph{arrival} on their own are ambiguous.
A traveler departs their origin when they enter the link, and arrives at their destination when they exit the link.
But from the perspective of the link, one can also say that they arrive on the link when they enter, and depart the link when they exit.
To avoid confusion, we will use the words ``departure'' and ``arrival'' in the former sense, from the perspective of the traveler and not the link.}
We will assume that the physical section has \emph{zero} length, so that vehicles entering the link immediately reach the point queue at the downstream end.
We can do this, because the time spent traveling on the physical section is constant and unaffected by any of the choices travelers make, and because this is the only link in the network.
If the physical section did have a positive length, this would only amount to adding or subtracting a constant in different places in the derivation below.
Neglecting the length of the physical section simplifies the analysis, to highlight more salient features of the model.

Travelers are all identical, and wish to minimize the sum of their travel time on the link, and a schedule delay penalty based on their arrival time relative to their preferred arrival time.
This penalty takes the same form as in equation~\eqn{vickreydelay} above,
\labeleqn{scheduledelay}{f(t) = \beta [t^* - t]^+ + \gamma [t - t^*]^+\,.}
We will assume $0 < \beta < 1$ and $\gamma > 0$ (typically $\gamma > 1$ for behavioral reasons, but this is not strictly required).
We will model the traveler choices and queueing behavior in continuous time, and aim to characterize the equilibrium solution at which no traveler can reduce their total cost (travel time plus schedule delay) by departing at a different time.

It is easier to study this equilibrium in terms of \emph{arrival times} first, and use the arrival time representation to derive the departure time choices.
We will neglect time spent traveling on the physical section of the link, because this time is common to all travelers and is not affected by anyone's departure time choice.
We can thereby focus on what happens at the queue.
Let $Q(t)$\label{not:Qvickrey} denote the time spent waiting in the point queue for a traveler \emph{arriving} at the destination at time $t$.
The total cost experienced by this traveler is\label{not:Cvickrey}
\labeleqn{totalcost}{C(t) = Q(t) + f(t) = Q(t) + \beta[t^* - t]^+ + \gamma[t - t^*]^+\,.}
At equilibrium, all travelers experience the same total cost $c$.
This means that travelers who arrive closer to the preferred arrival time (smaller schedule delay) spent longer in queue, and that drivers who spend less time in queue arrive further from the preferred arrival time (larger schedule delay).
Figure~\ref{fig:waterfill} shows how the components of the equilibrium cost $c$ vary based on the arrival time.

\stevefig{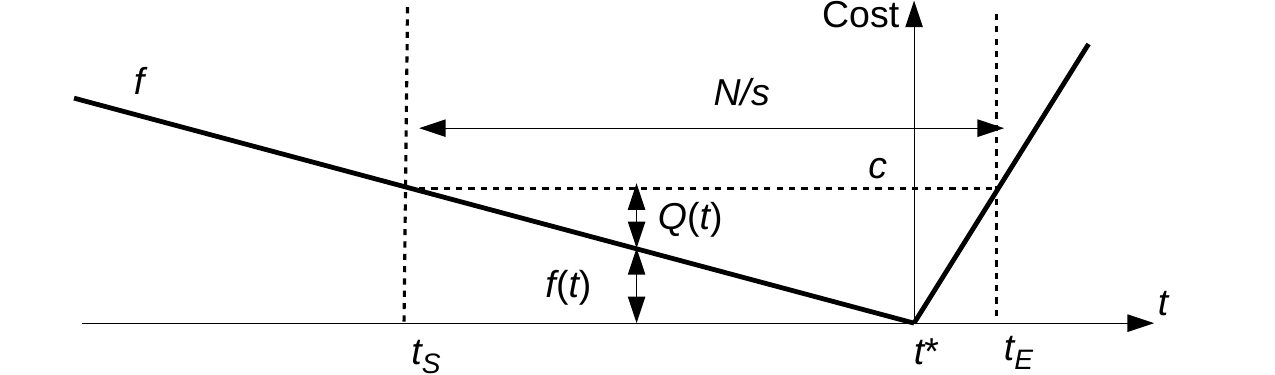}{Equilibrium costs consist of queueing and schedule delay.}{0.8\textwidth}

Let $t_S$\label{not:tS} and $t_E$\label{not:tE} denote the times at which the first and last vehicles arrive at the destination.
(The mnemonic is `s' for start of arrivals, and `e' for end.)
The downstream capacity is $s$, so it takes a time $N/s$ for all of the vehicles to exit the link; therefore $t_E = t_S + N/s$.
We also know that the first and last travelers to arrive experience the same total cost.
They experience no queueing delay (the first traveler faces no queue because nobody is in front of them; the last traveler faces no queue because they can time their departure so they arrive right as the queue dissipates), so we know $f(t_E) = f(t_S)$, or $\beta(t^* - t_S) = \gamma(t_E - t^*)$.
Solving these two equations gives the time window for arrivals:
\labeleqn{arrivals}{
t_S = t^* - \frac{\gamma}{\beta + \gamma} \frac{N}{s} \qquad t_E = t^* + \frac{\beta}{\beta + \gamma} \frac{N}{s}
\,.
}
The equilibrium cost experienced by \emph{all} of the travelers is thus
\labeleqn{eqmcost}{
c = f(t_S) = f(t_E) = \frac{\beta \gamma}{\beta + \gamma} \frac{N}{s}
\,.
}
and therefore the traveler arriving at time $t$ experiences a queueing delay
\labeleqn{eqmqueueingdelay}{
Q(t) = c - f(t) = \frac{\beta \gamma}{\beta + \gamma} \frac{N}{s} - \beta[t^* - t]^+ - \gamma[t - t^*]^+
\,.
}
The \emph{total system cost} $TSC$\label{not:TSC} experienced by all travelers is simply the product of the individual cost and the number of travelers.
\labeleqn{vickreytstt}{
TSC = c \frac{N}{s} = \frac{\beta \gamma}{\beta + \gamma} \frac{N^2} {s}
\,.
}
We can also interpret this graphically as the area of the rectangle in Figure~\ref{fig:waterfill}, multiplied by $s$ (because the horizontal axis measures time, and $s$ vehicles arrive per unit of time.)

This is a full description of the bottleneck equilibrium in terms of arrivals.
It's more natural to describe traveler choices instead of \emph{departures}, so we now transform these above results to identify the departure rates that produce the queueing profile~\eqn{eqmqueueingdelay}.
Figure~\ref{fig:recoverprofile} gives a graphical intuition of how these are related.
If we know the arrival time $t$ and queueing delay $Q(t)$, then the traveler's departure time must have been $t' = t - Q(t)$.
The departure time for a specified traveler can be obtained by drawing a 45-degree line upwards and to the left from their queueing delay experienced; this line intersects the horizontal equilibrium cost line at their departure time.
In Figure~\ref{fig:recoverprofile}, these lines have been drawn at a uniform spacing in terms of arrival time.
We see that the departure profile is uneven: vehicles that arrive early depart at a faster rate, and vehicles that arrive late depart at a slower rate.

\stevefig{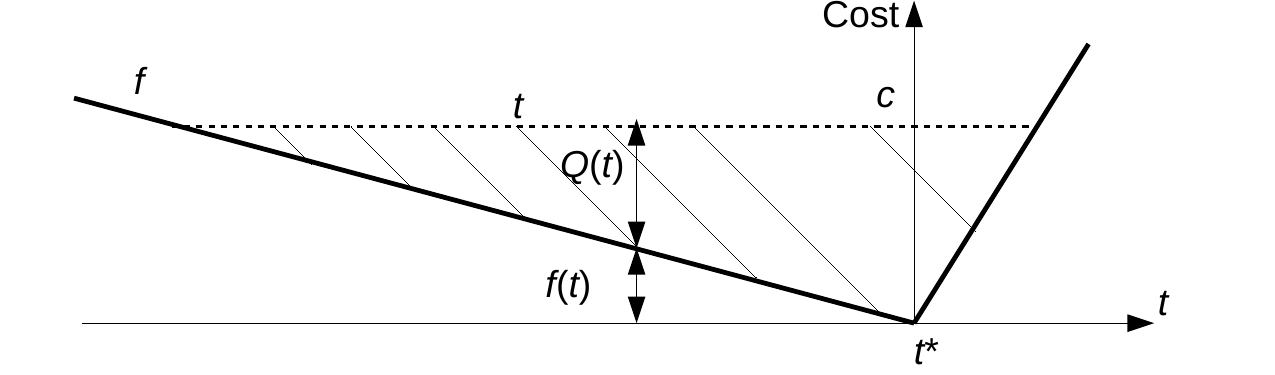}{Relating arrival queueing delay to departure rates.}{0.8\textwidth}

We now derive the exact expressions for the departure profile.
For this purpose, let $Q'(t')$ denote the queueing delay experienced by a traveler departing at time $t'$.
Essentially, $Q'(t')$ gives the queueing delay that a traveler leaving at $t'$ \emph{will} experience, while $Q(t)$ gives the queueing delay that a traveler arriving at $t$ \emph{has} experienced.
Of course, for a specific traveler, these two values must agree for their specific departure and arrival times $t'$ and $t$.
The traveler leaving at $t'$ arrives at $t' + Q'(t')$, so
\labeleqn{queuetransform}{
Q'(t') = Q(t' + Q'(t'))
\,.
}
Taking derivatives of both sides, the chain rule gives
\labeleqn{queuetransformderstep}{
\frac{dQ'}{dt'} = \frac{dQ}{dt} \myp{ 1 - \frac{dQ'}{dt'}}
}
or
\labeleqn{departureder}{
\frac{dQ'}{dt'} = \frac{dQ/dt}{1 + dQ/dt}
\,.
}
(In these expressions, $dQ'/dt'$ is evaluated at the departure time $t'$, and $dQ/dt$ is evaluated at the arrival time $t$.)
During the ``early arrival'' period, equation~\eqn{eqmqueueingdelay} tells us that $dQ/dt = \beta$, and during the ``late arrival'' period $dQ/dt = -\gamma$.
Therefore, during the early period we have $dQ'/dt' = \beta/(1 - \beta)$, and during the late period $dQ'/dt' = -\gamma(1 + \gamma)$.

The last step is to relate the change in departure queueing delay to the actual departure rates, given the point queue dynamics.
In Section~\ref{sec:pointqueue}, we used $N^\uparrow(t)$ and $N^\downarrow(t)$ to respectively denote the cumulative number of travelers who have entered the link up through time $t$, and the cumulative number of travelers who have exited the link up through time $t$.
We are assuming the physical section has zero length, so $N^\uparrow(t)$ also gives the number of vehicles that have entered the point queue up through time $t$.
Recall that the difference between these gives the number of vehicles waiting in the queue at time $t$, and the time spent waiting in the point queue is the quotient of the number of vehicles in the queue and the downstream capacity.
Therefore, a traveler departing at time $t'$ must experience a queueing delay of
\labeleqn{pqdelay}{
Q'(t') = \frac{N^\uparrow(t') - N^\downarrow(t')}{s}
\,,
}
so that if we denote the inflow rate by $q(t') = dN^\uparrow/dt'$, we have
\labeleqn{pqdelayder}{
\frac{dQ'}{dt'} = \frac{q(t') - s}{s}
\,.
}
Using the previously-computed values of $dQ'/dt'$, we obtain the departure rates as $s/(1 - \beta)$ during the early arrival period, and $s / (1 + \gamma)$ during the late arrival period.
The transition from early to late arrivals happens for the traveler departing at $t^* - Q(t^*) = t^* - c$.
The queue entry and exit profiles are plotted in Figure~\ref{fig:vickreyprofile}.

\stevefig{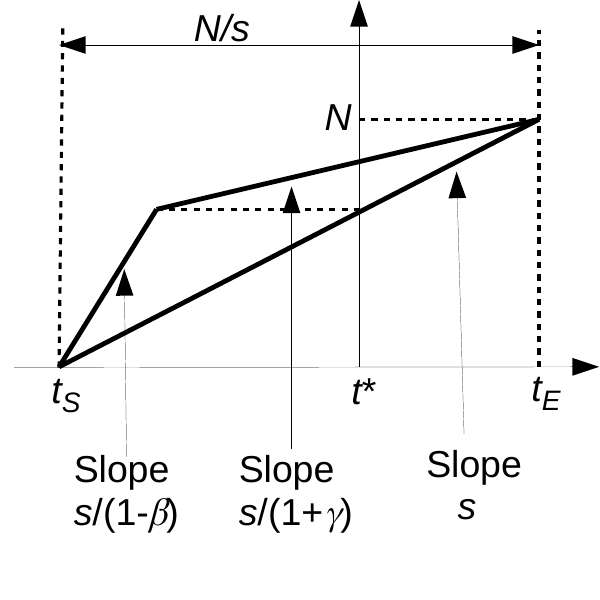}{Equilibrium departures and arrivals in Vickrey's bottleneck.}{0.5\textwidth}

\index{user equilibrium!departure-time choice}
\index{system optimal!departure-time choice}
This user equilibrium solution is not socially optimal.
To see this, imagine an alternate departure profile where travelers depart over the same interval $[t_S, t_E]$, but at a uniform rate of $s$.
In this case, there is no queueing delay, because travelers enter the point queue at the same rate at which they exit.
There is still schedule delay, and in fact the schedule delay is identical to the equilibrium profile because \emph{each traveler arrives at exactly the same time as in the user equilibrium}.
The only difference is that we have now spaced out their departures to be uniform, rather than heavier in the early arrival period, and lighter in the late arrival period.
The total cost in this case is simply the area under $f(t)$ in Figure~\ref{fig:waterfill} multiplied by $s$ as a scale factor; from this figure we immediately see that the total cost is half that of the user equilibrium.
(In the language of Section~\ref{sec:priceofanarchy}, the ``price of anarchy'' in Vickrey's bottleneck is 2.)
This departure profile has a much lower cost, but it is not an equilibrium because different travelers experience different cost.
The first and last travelers are envious of the traveler who departs and arrives at the preferred arrival time, with no queueing delay whatsoever.

How might we obtain this socially optimal solution?
Classical economists would propose replacing the queueing delay with a \emph{dynamic toll} exactly equal to the equilibrium $Q(t)$ (converting to units of money).
After doing so, the socially optimal departure profile is now an equilibrium, since all travelers face equal cost; travelers arriving closer to the preferred arrival time pay a monetary toll, rather than through waiting in queue.
Why is this any better?
Time spent waiting in queue is completely lost and unrecoverable.
Tolls, on the other hand, do not disappear.
The money simply changes hands to another part of the economy, and can either be used to fund infrastructure, or redistributed back to the travelers in the form of reduced taxation elsewhere, or a direct rebate.
We could even take the toll revenues collected and divide them equally among all the travelers in the system.
As long as the method of redistribution is not affected by a traveler's choice, their individual optimum will agree with the system optimum behavior, eliminating all queueing delay.
This is the classic justification for dynamic congestion pricing.
\index{tolling}

Of course, real systems are much more complex than the Vickrey bottleneck.
Networks have thousands of links and nodes, and there is not just as single bottleneck.
Travelers are not all uniform: they do not all have the same preferred arrival time, early and late penalty factors, or value of time.
(The toll example given above cannot represent any kind of equity issues associated with tolls, because it assumes a homogeneous population.)
Nevertheless, as a model, it provides high-level insights on how certain types of congested systems behaves, and provides a basis for more specific research that can relax these assumptions and provide guidance on a particular real-world context.
Literally hundreds of research papers have used or extended Vickrey's bottleneck in different ways, and Section~\ref{sec:tdsplitreview} describes a few of them.
\index{Vickrey bottleneck|)}
\index{departure time choice!equilibrium|)}

\section{Historical Notes and Further Reading}
\label{sec:tdsplitreview}

Several authors have discussed how the time-dependent shortest path problem differs from the static case; both label-setting~\citep{dreyfus69} and label-correcting~\citep{cooke66,ziliaskopoulos94} methods exist.
\cite{chabini99} showed how the time index can be used as a topological order in the time-expanded network.
The use of free-flow times as estimates in time-dependent $A^*$ is from \cite{boyles_wydot}.

Many researchers have studied equilibrium with departure time choice, using schedule delay concepts.
One such study is the famous ``single bottleneck'' model of \cite{vickrey69}, where there is no route choice, but only departure time choice.
As of 2019, over two hundred research works have been based on this model~\citep{li20}.
A number of review papers highlight such work, including \cite{small92}, \cite{arnott98}, \cite{lindsey01}, and \cite{small15}.
As just a few examples of how the bottleneck model has been extended, researchers have investigated ways to integrate the route choice decision \citep[e.g.,][]{arnott92}, mode choice decision \citep{tabuchi93}, traveler heterogeneity \citep{vickrey73,lindsey04,ramadurai10}, and uncertainty in bottleneck capacity and/or travel demand \citep{arnott99eer}.
Vickrey's bottleneck is commonly used for research in dynamic congestion pricing \citep{depalma86}.
Several authors have attempted to address the equity issues associated with the simple dynamic pricing scheme in Section~\ref{sec:singlebottleneck}, including \cite{liu10} and \cite{helsel17}.
Departure time choice has been combined with route choice in dynamic traffic assignment models; \cite{levin15departuretime} show one way to do this.
Activity-based modeling provides another way to model departure time choices.
A full discussion of activity-based modeling is beyond the scope of this book; see \cite{bhat99} for additional discussion.

Section~\ref{sec:graphreferences} described how the shortest path problem can be generalized to include stochastic link costs, where the link costs are drawn from some probability distribution and the path with minimum expected cost is sought.
You may recall that this problem was not too difficult to address in the static case.
However, if travel times are \emph{both} time-dependent and stochastic, more care is needed~\citep{hall86,fu98}, because Bellman's principle\index{Bellman's principle} need not hold.
Examples of algorithms to handle this issue are given in \cite{hall86} and \cite{millerhooks00}.
\index{shortest path!time-dependent|)}

\section{Exercises}
\label{sec:tdsp_exercises}

\begin{enumerate}
\item \diff{13} Table~\ref{tbl:fifoex} shows time-dependent costs on five links, for different entry times.
Which links have costs satisfying the FIFO principle? \label{ex:fifoex}
\begin{table}
\begin{center}
\caption{Link costs for Exercises~\ref{ex:fifoex} and~\ref{ex:canvasctd}.
\label{tbl:fifoex}}
\begin{tabular}{|c|ccccc|}
\hline
Entry time & Link 1 & Link 2 & Link 3 & Link 4 & Link 5 \\
\hline
0 & 5 & 1 & 3 & 1 & 1 \\
1 & 6 & 1 & 3 & 6 & 1 \\
2 & 4 & 2 & 3 & 5 & 2 \\
3 & 5 & 3 & 3 & 4 & 4 \\
4 & 7 & 5 & 3 & 3 & 2 \\
5 & 8 & 8 & 3 & 2 & 1 \\
\hline
\end{tabular}
\end{center}
\end{table}
\item \diff{23} Prove that waiting is never beneficial in a FIFO network where link costs are equal to travel time.
\item \diff{34} Prove that there is an acyclic time-dependent shortest path in a FIFO network, if link costs are equal to travel time.
\item \diff{53} In this chapter, we assumed it is impossible to travel beyond the time horizon, by not creating time-expanded links if they arrive at a downstream node after $\bar{T}$.
Another alternative is to assume that travel beyond $\bar{T}$ is permitted, but that travel times and costs stop changing after that point and take constant values.
(Perhaps free-flow times after the peak period is over.)   \label{ex:boundary1}
\begin{enumerate}
  \item Modify the time-expanded network concept to handle this assumption.
(The network should remain finite.)
  \item Modify the algorithm in Section~\ref{sec:iotalgo} to work in this setting.
\end{enumerate}
\item \diff{63} Another way to handle the time horizon is to assume that the link travel times and costs are \emph{periodic} with length $\bar{T}$.
(For instance, $\bar{T}$ may be 24 hours, and so entering a link at $t = 25$ hours would be the same as at $t = 1$ hour.)  First repeat Exercise~\ref{ex:boundary1} with this assumption.
Then prove that the modified algorithm you create will converge to the correct time-dependent shortest paths for any possible departure time.
\label{ex:boundary2}
\item \diff{31} Modify the algorithm in Section~\ref{sec:artificialod} to become an all-to-one algorithm, that finds least-cost paths from all origins and all departure times to one destination (arrival at any time is permitted).
\label{ex:alldeparturetimes} 
\item \diff{34} Show that the classical form of Bellman's principle (Section~\ref{sec:shortestpath}) holds in a FIFO, time-dependent network where the link costs are equal to the link travel time.
\item \diff{32} Find the time-dependent shortest path when departing node 1 in the network shown in the \emph{left} panel of Figure~\ref{fig:bellmancounterexample}, departing at time 0.
\label{ex:bellmanfail1}
\item \diff{23} Tables~\ref{tbl:canvasbacknode} and~\ref{tbl:canvascosts} show the backnode and cost labels for a time-dependent shortest path problem, where the destination is node 4, and the time horizon is 6.
(Figure~\ref{fig:canvasnet} shows the network topology.)  What is the shortest path from node 1 to node 4, when departing at time 1? \label{ex:canvasex}
\begin{table}
\begin{center}
\caption{Backnode labels for Exercise~\ref{ex:canvasex}.
\label{tbl:canvasbacknode}}
\begin{tabular}{|c|cccc|}
\hline
Entry time & Node 1 & Node 2 & Node 3 & Node 4 \\
\hline
6 & $-1$ & $-1$ & $1:1$ & $2:1$ \\
5 & $-1$ & $-1$ & $-1$ & $3:5$ \\
4 & $-1$ & $-1$ & $2:2$ & $-1$ \\
3 & $-1$ & $-1$ & $-1$ & $-1$ \\
2 & $-1$ & $1:1$ & $-1$ & $-1$ \\
1 & $-1$ & $-1$ & $-1$ & $-1$ \\
\hline
\end{tabular}
\end{center}
\end{table}
\begin{table}
\begin{center}
\caption{Cost labels for Exercise~\ref{ex:canvasex}.
\label{tbl:canvascosts}}
\begin{tabular}{|c|cccc|}
\hline
Entry time & Node 1 & Node 2 & Node 3 & Node 4 \\
\hline
6 & $\infty$ & $\infty$ & 6 & 6 \\
5 & $\infty$ & $\infty$ & $\infty$ & 5 \\
4 & $\infty$ & $\infty$ & 4 & $\infty$ \\
3 & $\infty$ & $\infty$ & $\infty$ & $\infty$ \\
2 & $\infty$ & 1 & $\infty$ & $\infty$ \\
1 & 0 & $\infty$ & $\infty$ & $\infty$ \\
\hline
\end{tabular}
\end{center}
\end{table} 
\genfig{canvasnet}{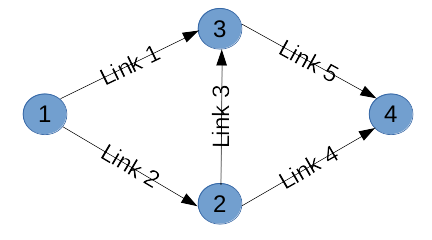}{Network for Exercises~\ref{ex:canvasex} and~\ref{ex:canvasctd}.}{width=0.5\textwidth}
\item \diff{22} In the same network as Exercise~\ref{ex:canvasex}, fill in the backnode and cost labels for time 0.
The link costs are as in Table~\ref{tbl:fifoex}, and waiting is not allowed at nodes.
\label{ex:canvasctd}
\item \diff{45} Consider the network in Figure~\ref{fig:tdspexercisenet}.
\label{ex:tdspexercisenet}
\begin{enumerate}[(a)]
\item Verify that the travel times satisfy the FIFO principle.
\item Find the shortest paths between nodes 1 and 4 when departing at $t = 0$, $t = 2$, and $t = 10$.
\item For what departure times would the travel times on paths $[1,2,4]$ and $[1,3,4]$ be equal?
\end{enumerate}
\stevefig{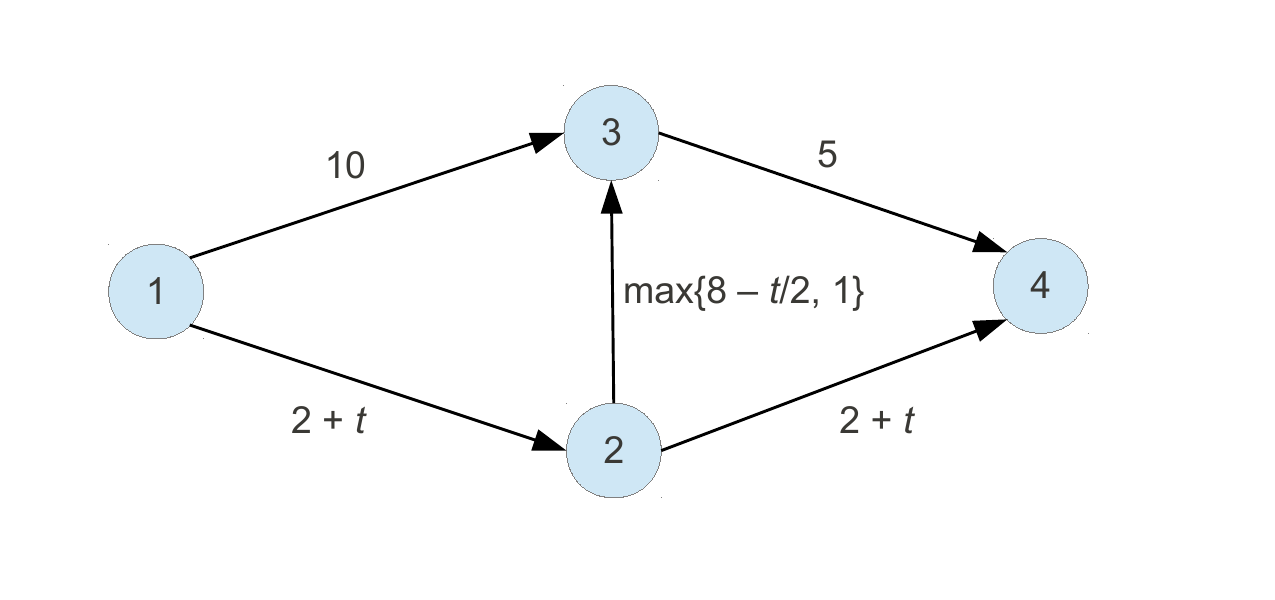}{Network for Exercise~\ref{ex:tdspexercisenet}.}{0.8\textwidth}

\label{ex:canvasexmod}
\begin{table}
\begin{center}
\caption{Backnode labels for Exercise~\ref{ex:canvasexmod}.
\label{tbl:canvasexmod}}
\begin{tabular}{|c|cccc|}
\hline
Entry time & Node 1 & Node 2 & Node 3 & Node 4 \\
\hline
6 & $-1$ & $-1$ & $-1$ & $3:3$ \\
5 & $-1$ & $-1$ & $-1$ & $3:4$ \\
4 & $-1$ & $-1$ & $2:3$ & $2:3$ \\
3 & $-1$ & $1:1$ & $1:1$ & $-1$ \\
2 & $-1$ & $-1$ & $-1$ & $-1$ \\
1 & $-1$ & $-1$ & $-1$ & $-1$ \\
\hline
\end{tabular}
\end{center}
\end{table}
\item \diff{42} Table~\ref{tbl:canvasexmod} shows the backnode labels after completing the algorithm in Section~\ref{sec:iotalgo}.
\begin{enumerate}
 \item Draw the network corresponding to these backnode labels.
 \item What is the shortest path from node 1 when departing at time 1?
 \item Does the links in this network obey the FIFO principle?
\end{enumerate}
\item \diff{23} Solve the time-dependent shortest path problem on the network in Figure~\ref{fig:alltimesalgodemo} using the algorithm in Section~\ref{sec:fifonets}.
Comment on the amount of work needed to solve the problem using this algorithm, compared to the way it was solved in the text.
\item \diff{51} In the schedule delay equation~\eqn{vickreydelay}, we typically assume $0 \leq \alpha \leq 1 \leq \beta$ for peak-hour commute trips.
What counterintuitive behavior would occur if any of these three inequalities were violated?
\item \diff{12} Find the optimal departure times and paths from nodes 2 and 3 in the network in Figure~\ref{fig:alltimesalgodemo}, with equation~\eqn{scheduledelayalgodemo} as the arrival time penalty.
(You can answer this question directly from Table~\ref{tbl:alltimesalgodemo}.)
\item \diff{36} Find the optimal departure times and paths from nodes 1, 2, and 3 in the network in Figure~\ref{fig:alltimesalgodemo}, if the arrival time penalty function is changed to $f(t) = ([12-t]^+)^2 + 2([t - 12]^+)^2$.
\label{ex:opttime}
\item \diff{0} According to the penalty function in Exercise~\ref{ex:opttime}, what is the desired arrival time?
\item \diff{88} Implement all of the algorithms described in this chapter, and test them on transportation networks with different characteristics (number of origins, destinations, nodes, links, time intervals, ratio of links to nodes, with or without waiting, FIFO or non-FIFO, etc.).
What algorithms perform best under what circumstances?  For circumstances most resembling real-world transportation networks, what performs best?
\end{enumerate}

\chapter{Dynamic User Equilibrium}
\label{chp:dynamiceqm}

\index{dynamic traffic assignment!equilibrium|(}
The previous two chapters presented separate perspectives on dynamic traffic modeling.
Chapter~\ref{chp:networkloading} described the network loading problem, in which driver behavior (as represented by path choices and possibly departure times) was known, and where we sought to represent the resulting traffic and congestion patterns on the network.
Chapter~\ref{chp:tdsp} described the time-dependent shortest path and departure time choice problems, in which the network state (as represented by travel times and costs) was known, and where we sought to identify how drivers would behave.
This chapter synthesizes these two perspectives through the concept of dynamic user equilibrium, defined as driver choices and network traffic which are mutually consistent, given the network loading and driver behavior assumptions.
The dynamic user equilibrium principle, and how the network loading and driver behavior models can be connected, are the subjects of Section~\ref{sec:dueintro}.
Section~\ref{sec:duealgos} then presents several algorithms that can be used to solve for dynamic equilibria.
Such algorithms are almost always heuristics, since realistic network loading models are not amenable to exact analysis, and indeed Section~\ref{sec:dueproperties} shows that dynamic user equilibrium need not exist; and if it exists, it need not be unique.
This section also provides examples to show that the dynamic user equilibrium solution need not minimize total travel time in a network, and that providing additional capacity to a network can increase total travel time, a dynamic analogue of the Braess paradox.
The chapter concludes with a discussion on implementing a dynamic traffic assignment model, including data collection, network construction, and validation.

\section{Towards Dynamic User Equilibrium}
\label{sec:dueintro}

This section has two major goals, aimed at reconciling the network loading problem of Chapter~\ref{chp:networkloading} and the time-dependent shortest path algorithms of Chapter~\ref{chp:tdsp}.
In particular, we will need to solve these problems sequentially, and repeatedly, as shown in Figure~\ref{fig:loadingtdsp}.
The first goal of this section is to provide the mechanics to link the output of each problem to the input needed by the other.
Following these, we will be equipped to define dynamic user equilibrium formally.
Fixed point and variational inequality formulations are given as well.

\subsection{Flow representation}
\label{sec:flowrepresentation}

\index{dynamic traffic assignment!flow representation|(}
The first question to address is how to represent the choices of all the travelers on the network.
In describing network loading, Chapter~\ref{chp:networkloading} took an aggregate approach.
All of the link models in that chapter can be seen as fluid models, where vehicles are infinitely divisible rather than discrete entities.
However, the algorithms in Chapter~\ref{chp:tdsp} took a disaggregate approach, and found paths and departure times for a single vehicle.
There are several ways to reconcile these perspectives.

\index{dynamic traffic assignment!flow representation!matrix}
The method adopted in this chapter is to store the number of vehicles departing on each path during each time interval.
Denote by $h_t^\pi$\label{not:hpit} the number of vehicles which start traveling on path $\pi$ during the $t$-th interval.
These values do not need to be integers; all that matters is that they are nonnegative, and that all paths connecting the same origin and destination sum to the total demand between these zones for each time interval.
The $h^\pi_t$ values can be collected into a single matrix $H$,\label{not:Hdta} whose dimensions are equal to the number of paths and time intervals.
We will use $d_t^{rs}$\label{not:drst} to reflect the total demand between origin $r$ and destination $s$ departing at time $t$, and as before, use $T$ and $\Pi$ to respectively denote the time horizon, and the set of network paths.

The main advantage of this notation is that the behavior is clear: by tracking some auxiliary variables in the network loading (as described in Section~\ref{sec:flowtracking}, at any point in time we can see exactly which vehicles are on which links, and can trace these to the paths the vehicles must follow).
These paths can be connected directly with the paths found in a time-dependent shortest path algorithm.
A disadvantage of this approach is that the number of paths grows exponentially with the network size.
In practice, ``column generation'' schemes are popular, in which paths $\pi$ are only identified when found by a shortest path algorithm, and $h^\pi_t$ values only need be calculated and stored for paths to which travelers have been assigned.

\index{dynamic traffic assignment!flow representation!splitting proportion}
\index{splitting proportion!flow representation}
This is not the only possible way to represent travel choices.
A link-based representation requires fewer variables.
For each turning movement $[h,i,j]$, each destination $s$,  and each time interval $t$, let $\alpha_{hij,s}^t$\label{not:alphahijst} denote the proportion of travelers arriving at node $i$, via approach $(h,i)$, during the $t$-th time interval, who will exit onto link $(i,j)$ \emph{en route} to destination $s$.
This approach mimics the flow splitting rules often found in traffic microsimulation software.
The number of variables required by this representation grows with the network size, but not at the exponential rate required by a path-based approach.
One can show that the path-based and link-based representation of route choice are equivalent in the sense that there exist $\alpha$ values which represent the same network state as any feasible set of $H$ values, regardless of the network loading model, and that one can identify the $H$ values corresponding to a given set of $\alpha$ values (see Exercise~\ref{ex:pathlinkequivalence}).
The primary disadvantage of this representation is that the behavior is less clear: one cannot trace the path of any vehicle throughout the network deterministically (although one can do so stochastically, making a turn at each node by treating the $\alpha$ values as probabilities, identifying the arrival times at successive nodes using the procedure in Section~\ref{sec:timetracing}).

\index{dynamic traffic assignment!flow representation!discrete vehicles}
Both of these methods adopt the aggregate, continuous-flow perspective of Chapter~\ref{chp:networkloading}.
Yet another way to represent travel choices is to adopt the individual perspective of the algorithms in Chapter~\ref{chp:tdsp}, and explicitly model each vehicle as a discrete agent with one specific path.
This is the most behavioral way to represent choices: each individual vehicle is assigned one path, with no divisibility or fractional flows.
Downsides of this approach are scalability --- if the demand doubles but the network topology is unchanged, the individual vehicle approach would require twice as many variables, whereas the continuous methods would not require any more variables --- and increased difficulty in implementing the network loading algorithms.
It is certainly possible to implement discrete versions of merges, diverges, the cell transmission model, and so forth, but one must be careful about rounding.
For instance, if the time step is chosen so that the capacity of a link is less than half a vehicle per timestep, consistently rounding numbers to the nearest integer would mean no vehicles can ever exit.
Stochastic rounding, or accumulating a continuous flow value which can then be rounded, can address these difficulties.
\index{dynamic traffic assignment!flow representation|)}

\stevefig{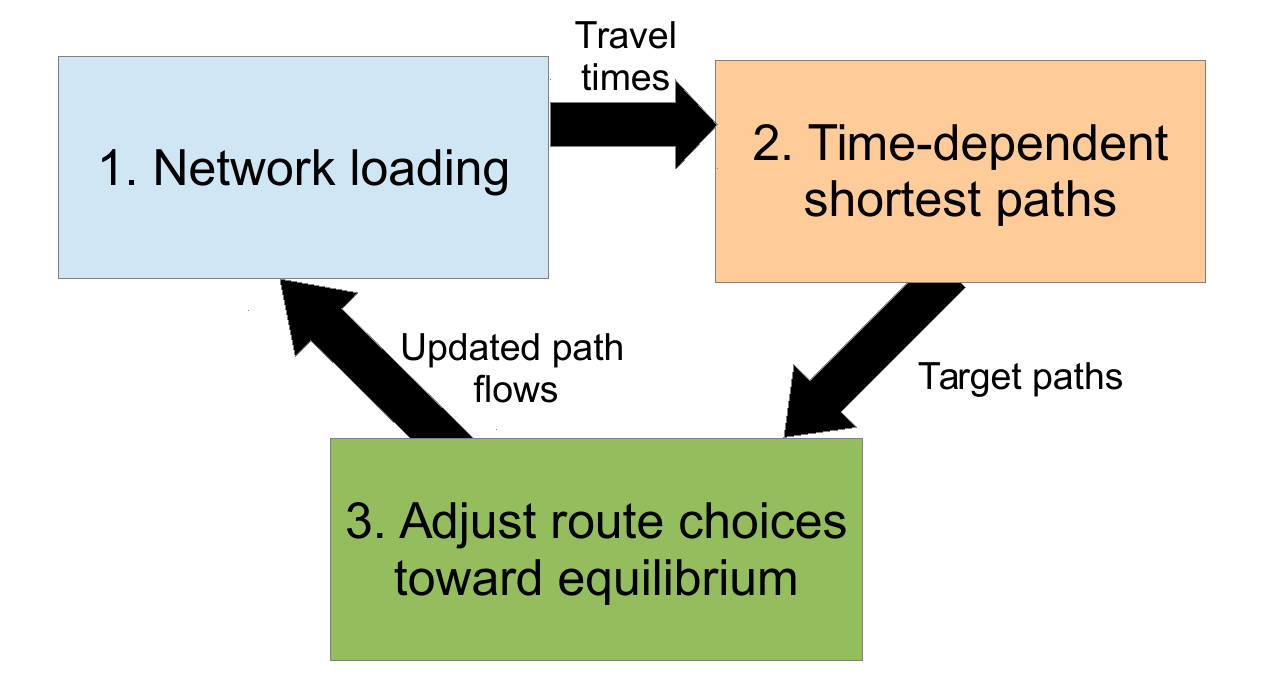}{Iterative dynamic traffic assignment framework, with details from previous chapters.}{0.6\textwidth}
\subsection{Travel time calculation}
\label{sec:timetracing}

\index{network loading!computing travel times|(}
Once network loading is complete, a time-dependent shortest path algorithm can be applied to determine the cost-minimizing routes for travelers.
To do so, we need the travel times $\tau_{ij}(t)$ on each link $(i,j)$ for travelers entering at each time $t$.
The network loading models do not provide this directly, but rather give the cumulative counts $N$ at the upstream and downstream ends of each link, at each time interval.
However, this information will suffice for calculating the travel times.
For link $(i,j)$ at time $t$, the upstream count $N^\uparrow_{ij}(t)$ gives the total number of vehicles which have entered the link by time $t$, while the downstream count $N^\downarrow_{ij}(t)$ gives the total number of vehicles which left the link by time $t$.

Assume for a moment that $N^\uparrow_{ij}(t)$ and $N^\downarrow_{ij}(t)$ are strictly increasing, continuous functions of $t$.
If this is the case, we can define inverse functions $T^\uparrow_{ij}(n)$\label{not:Tupijn} and $T^\downarrow_{ij}(n)$,\label{not:Tdownijn} respectively giving the times when the $n$-th\label{not:nveh} vehicle entered the link and left the link.
The travel time for the $n$-th vehicle is then the difference between these: $T^\downarrow_{ij}(n) - T^\uparrow_{ij}(n)$.
Graphically, this can be seen as the horizontal difference between the upstream and downstream $N$ curves.
(Figure~\ref{fig:cumulativecounttraveltime}).
Then, to find the travel time for a vehicle entering the link at time $t$, we simply evaluate this difference for $n = N^\uparrow_{ij}(t)$:
\labeleqn{continuoustt}{\tau_{ij}(t) = T^\downarrow_{ij}(N^\uparrow_{ij}(t)) - T^\uparrow_{ij}(N^\uparrow_{ij}(t)) = T^\downarrow_{ij}(N^\uparrow_{ij}(t)) - t\,,}
since $T^\uparrow_{ij}$ and $N^\uparrow_{ij}$ are inverse functions.

\stevefig{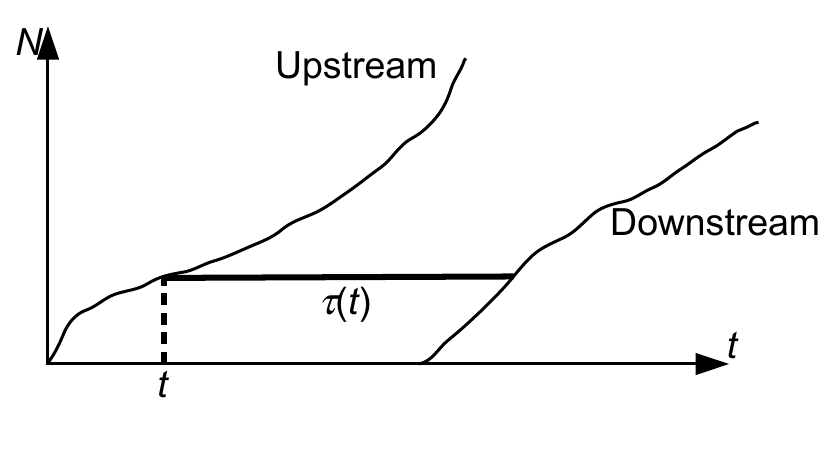}{Obtaining travel times from plots of the cumulative counts $N^\uparrow$ and $N^\downarrow$.}{0.6\textwidth}

A little bit of care must be taken because $N^\uparrow_{ij}(t)$ and $N^\downarrow_{ij}(t)$ are not strictly increasing functions of time, unless there is always a positive inflow and outflow rate for link $(i,j)$.
Furthermore, we often introduce a time discretization.
For both of these reasons, the inverse functions $T^\uparrow_{ij}(n)$ and $T^\downarrow_{ij}(n)$ may not be well-defined.
We thus need to modify equation~\eqn{continuoustt} in a few ways:
\begin{itemize}
\item We must ensure that $\tau_{ij}(t)$ is always at least equal to the free-flow travel time on the link.
The danger is illustrated in Figure~\ref{fig:noflow}, where no vehicles enter or leave the link for an extended period of time.
The horizontal distance between the $N^\uparrow$ and $N^\downarrow$ curves at their closest point is small (the dashed line in the figure), but this does not reflect the actual travel time of any vehicle.
In reality, a vehicle entering the link when it is completely empty, and when  there is no downstream bottleneck, would experience free-flow conditions on the link.
\item If the link has no outflow for an interval of time, then $N^\downarrow_{ij}(t)$ will be constant over that interval.
This frequently happens with traffic signals.
In this case, there are multiple values of time where $N^\downarrow_{ij}(t) = n$.
The correct way to resolve this is to define $T^\downarrow_{ij}(n)$ to be the \emph{earliest} time for which $N^\downarrow_{ij}(t) = n$:
\labeleqn{inverseproper}{T^\downarrow_{ij}(n) = \min_t \myc{ t : N^\downarrow_{ij}(t) = n}\,.
}
\item In discrete time, the time at which the $n$-th vehicle departs may not line up with a multiple of $\Delta t$, so there may be no known point where $N^\downarrow_{ij}(t)$ is exactly equal to $n$.
In this case, it is appropriate to interpolate between the last time point where $N^\downarrow_{ij}(t) < n$, and the first time point where $N^\downarrow_{ij}(t) \geq n$.
\end{itemize}
With these modifications to how $T^\downarrow_{ij}$ is calculated, the formula~\eqn{continuoustt} can be used to calculate the travel times on each link and arrival time.

\stevefig{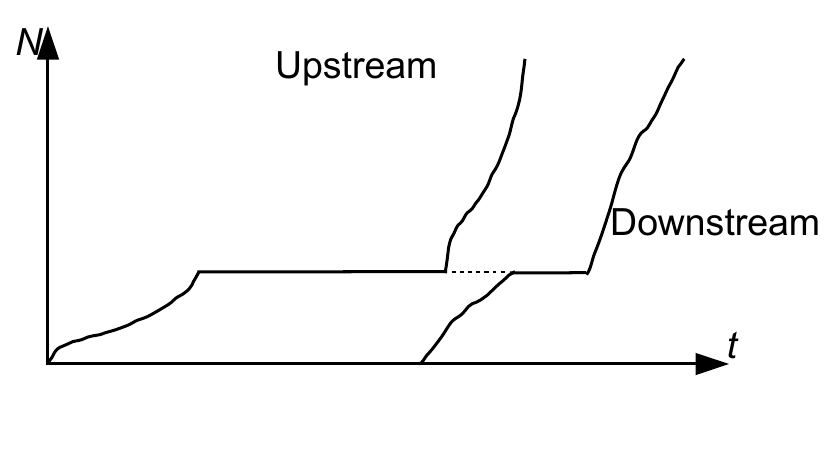}{Cumulative counts are flat when there is no inflow or outflow.}{0.6\textwidth}

The travel time on a path $\pi$ for a traveler departing at time $t$, denoted $C^\pi(t)$, can then be calculated sequentially.
If the path is $\pi = \mys{ r, i_1, i_2, \ldots, s }$, then the traveler departs origin $r$ at time $t$, and arrives at node $i_1$ at time $t + \tau_{ri_1}(t)$.
The travel time on link $(i_1, i_2)$ is then $\tau_{i_1 i_2} (t + \tau_{ri_1}(t))$, so the traveler arrives at $i_2$ at time $t + \tau_{ri_1}(t) + \tau_{i_1 i_2} (t + \tau_{ri_1}(t))$, and so forth.
Writing out this formula can be a bit cumbersome, but calculating it in practice is quite simple: it is nothing more than accumulating the travel times of the links in the path, keeping track of the time at which each link is entered.
\index{network loading!computing travel times|)}

\subsection{Determining splitting proportions}
\label{sec:flowtracking}

\index{splitting proportions!computing from loading|(}
The network loading requires knowledge of when and where vehicles enter the network, and the splitting proportions $p$ at diverges and general intersections.
In dynamic traffic assignment, this is reflected indirectly, through the variables $h^\pi(t)$, denoting the number of travelers departing path $\pi$ at time $t$.
We do not specifically give the proportions $p$ as a function of time, because we do not know when any traveler will reach any node before actually performing the network loading.
Path choice is a behavioral parameter associated with travelers, so it is easier to reconcile path choices with behavior if they are expressed this way, rather than by simply specifying the turning fractions at each node.

That said, the network loading does in fact need $p_{hij}(t)$ values for each turning movement $[h,i,j]$ at a diverge or general intersection, for each time period $t$.
These are obtained by examining the vehicles comprising the sending flow $S_{hi}(t)$, and calculating the fraction of these vehicles whose path includes link $(i,j)$ as the next link.
One way to do this is to disaggregate the cumulative count values $N^\uparrow(t)$ and $N^\downarrow(t)$ calculated at the upstream and downstream ends of each link.
For each path $\pi$ in the network, and for every link $(h,i)$ and time $t$, define $N^\uparrow_{hi,\pi}(t)$\label{not:Huphipit} and $N^\downarrow_{hi,\pi}(t)$\label{not:Hdownhipit} to be the total number of vehicles using path $\pi$ which have respectively entered and left link $(h,i)$ by time $t$.
Clearly we have
\labeleqn{pathdisaggregate}{N^\uparrow_{hi}(t) = \sum_{\pi \in \Pi} N^\uparrow_{hi,\pi}(t) \,; \qquad N^\downarrow_{hi}(t) = \sum_{\pi \in \Pi} N^\downarrow_{hi,\pi}(t)}
for all times $t$ and all links $(h,i)$.

Then, the sending flow for each link and time interval can be disaggregated in the same way, with $S_{hi,\pi}(t)$ defined as the number of vehicles in the sending flow which are using path $\pi$.
The procedure for doing this is conceptually simple, but writing out the equations and executing the steps requires some care.
We will explain the idea first with an example, and then give the general principle.
Consider the link $(h,i)$ in Figure~\ref{fig:flowdisagg_simple}, where there are two paths corresponding to different branches at the diverge node $i$.
Each circle in the figure represents a vehicle: empty circles reflect vehicles whose path involves a left turn at the diverge, and filled circles represent vehicles whose path involves a right turn at the diverge.
We have just queried the link model for $(h,i)$ for the sending flow at the current time $t$.
As indicated in the figure, the sending flow consists of four vehicles, so $S_{hi}(t) = 4$.
Evidently, three of them will turn right and one will turn left, so we should have $p_{hij}(t) = 1/4$ and $p_{hik}(t) = 3/4$.

\stevefig{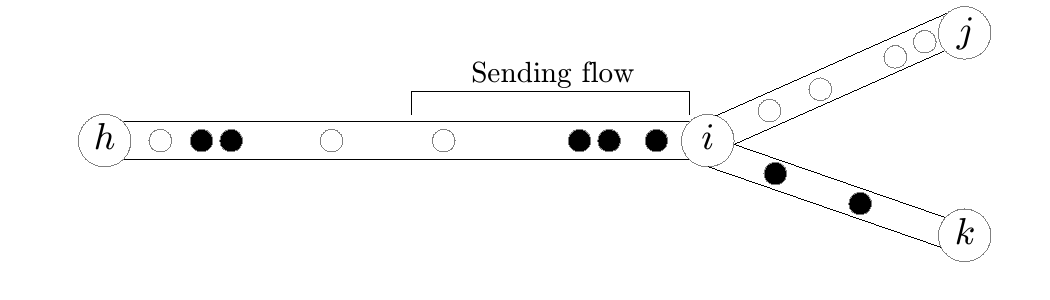}{What proportion of the sending flow is turning in each direction?}{0.8\textwidth}

How can we compute these proportions?
If we are tracking vehicles as discrete entities, doing so is trivial (simply look at the four vehicles that have been on the link the longest).
With our representation of path flows using a continuous matrix $H$, some more effort is required.
We know that vehicles must leave the link in the same order in which they enter.
If we keep track of the paths vehicles are on when they \emph{enter} the link, we can use this information to determine the proportions as they \emph{leave}.
To do this, we will use the disaggregate cumulative counts $N^\uparrow_{hi,\pi}(t)$ computed at the entry point to the link.
In Figure~\ref{fig:flowdisagg}(a), we have labeled each vehicle to show its $N_{hi}$ value, as well as the disaggregate value corresponding to its specific path.
In all the panels of this figure, please keep two points in mind: first, that discrete vehicles are shown for clarity, even though our models assume vehicles are a continuous fluid; and second, that $N$ values are defined at all points on the link, not just where the individual vehicles are shown.
On the first point, you can imagine that $N$ is roughly a step function which increases sharply by 1 where a vehicle is drawn, and nearly constant in between.
On the second point, recall that $N_{hi}(x,t)$ indicates the total number of vehicles that have crossed position $x$ on link $(h,i)$ up through time $t$, and likewise $N_{hi,\circ}(x,t)$ and $N_{hi,\bullet}$ respectively indicate the total number of vehicles on paths $\circ$ and $\bullet$ that have crossed position $x$ up through time $t$.
Therefore $N_{hi}(x,t) = N_{hi,\circ}(x,t) + N_{hi,\bullet}(x,t)$ for all positions $x$ and times $t$ on the link, as you should verify.

\stevefig{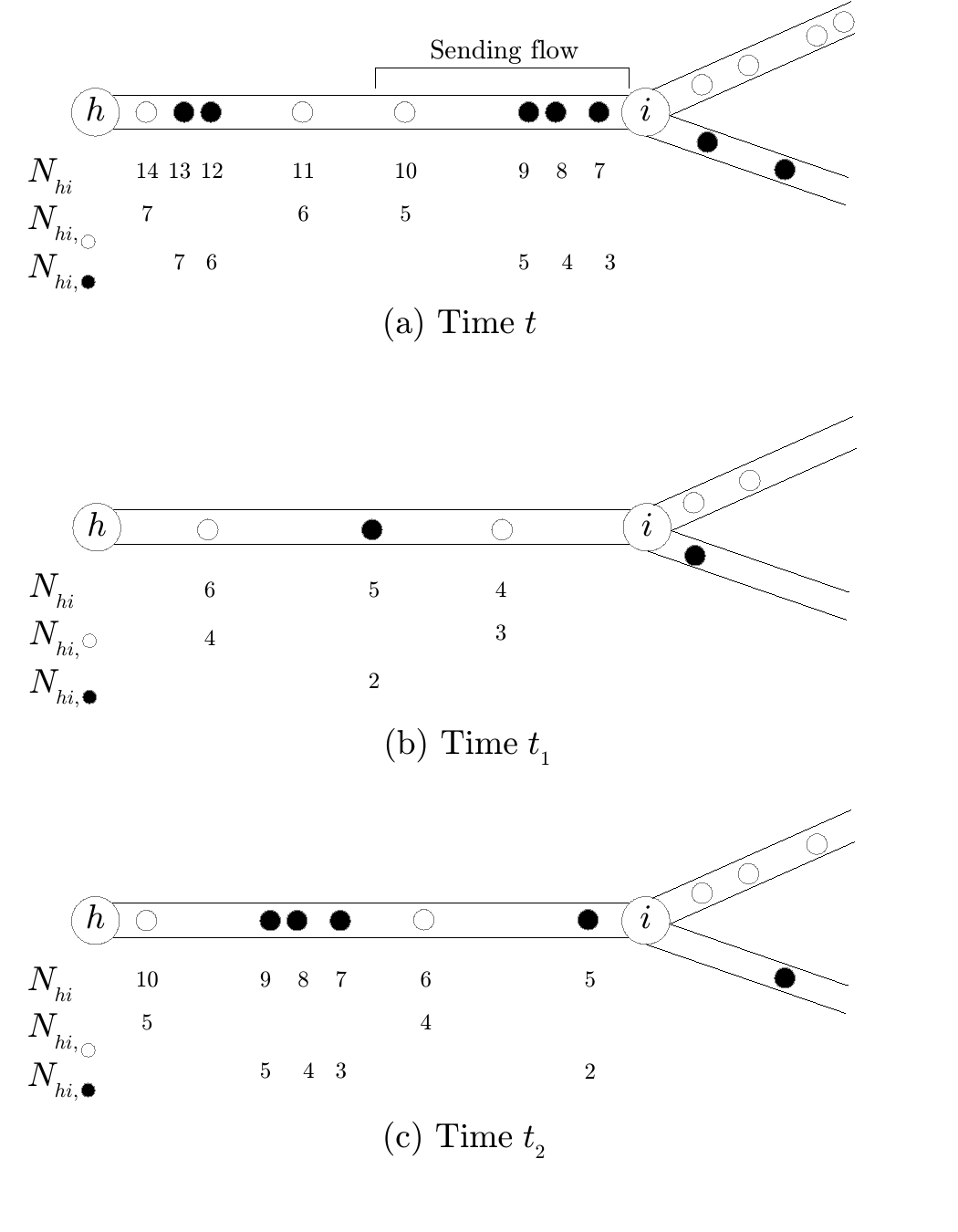}{Determining sending flow proportions based on vehicle entry times to the link.}{0.8\textwidth}

Figure~\ref{fig:flowdisagg}(a) is exactly the same situation as Figure~\ref{fig:flowdisagg_simple}, just with vehicle labels added.
Panels (b) and (c) in this figure show the same link at earlier points in time $t_1$ and $t_2$.
The time $t_1$ is chosen to be the time at which the first vehicle in the sending flow (vehicle 7) is about to enter the link.
The time $t_2$ is chosen to be the time at which the first vehicle \emph{not} in the sending flow (vehicle 11) is about to enter the link.
By comparing the state of the link at times $t_1$ and $t_2$, we can see which paths were on the vehicles that entered the link between these times.
These are exactly the vehicles which are in the sending flow, and so they will give us the diverge proportions.

Comparing the aggregate counts at the upstream end of the link, we have $N^\uparrow_{hi}(t_1) = 6$ and $N^\uparrow_{hi}(t_2) = 10$.
The difference of these is 4, as it must in order to match the value of $S_{hi}(t)$.
Looking at the disaggregate counts, we have $N^\uparrow_{hi,\circ}(t_1) = 4$ and $N^\uparrow_{hi,\circ}(t_2) = 5$, as well as $N^\uparrow_{hi,\bullet}(t_1) = 2$ and $N^\uparrow_{hi,\bullet}(t_1) = 5$.
Therefore, there is $5 - 4 = 1$ vehicle turning left, and $5 - 2 = 3$ vehicles turning right, and so our diverge proportions are $1/4$ and $3/4$.

We now express the process in a general way.
Assume that we are calculating sending flow for link $(h,i)$ at time $t$, and have determined $S_{hi}(t)$.
At this point in time, the total number of vehicles which has left this link is $N^\downarrow_{hi}(t)$,  Therefore, the vehicles in the sending flow are numbered in the range $N^\downarrow_{hi}(t)$ to $N^\downarrow_{hi}(t) + S_{hi}(t)$.
Using the inverse functions~\eqn{continuoustt}, the times at which these vehicles entered link $(h,i)$ are in the range $T^\uparrow(N^\downarrow_{hi}(t))$ to $T^\uparrow(N^\downarrow_{hi}(t) + S_{hi}(t))$.
Denote these two times by $t_1$ and $t_2$, respectively.
Then the disaggregate sending flow $S_{hi,\pi}(t)$\label{not:Shipit} is the number of vehicles on path $\pi$ which entered link $(h,i)$ between $t_1$ and $t_2$, that is,
\labeleqn{sendingdisaggregate}{S_{hi,\pi}(t) = N^\uparrow_{hi,\pi} (t_2) - N^\uparrow_{hi,\pi} (t_1)\,.}

The proportion of travelers in the sending flow $S_{hi}$ wishing to turn onto link $(i,j)$ can now be calculated.
Let $\Pi_{[h,i,j]}$\label{not:Pihij} be the set of paths which arrive at node $(h,i)$ and immediately continue onto link $(i,j)$.
Then 
\labeleqn{proportiontravelers}{p_{hij}(t) = \frac{\sum_{\pi \in \Pi_{[h,i,j]}} S_{hi,\pi}(t)}{S_{hi}(t)}\,,}
for use in diverge or general intersection node models.
Then, after a node model produces the actual flows $y_{hij}(t)$ between links, the disaggregated cumulative counts must be updated as well.

Let $y_{hi}(t) = \sum_{(i,j) \in \Gamma(i)} y_{hij}(t)$\label{not:yhi} be the total flow which leaves link $(h,i)$ during time $t$.
Then all vehicles which entered the link between $t_1$ and $T^\uparrow_{hi}(N^\downarrow(t)_{hi} + y_{hi})$ can leave the link (use $t_3$ to denote this latter time), so we update the downstream aggregate counts to
\labeleqn{downstreamdisaggregate}{N^\downarrow_{hi,\pi}(t + \Delta t) = N^\uparrow_{hi,\pi}(t_3)}
and the upstream counts to
\labeleqn{upstreamdisaggregate}{N^\uparrow_{ij,\pi}(t + \Delta t) = N^\uparrow_{ij,\pi}(t) + (N^\downarrow_{hi,\pi}(t + \Delta t) - N^\downarrow_{hi,\pi}(t)) \,.}
where $(h,i)$ is the link immediately upstream of $(i,j)$ on path $\pi$.

\stevefig{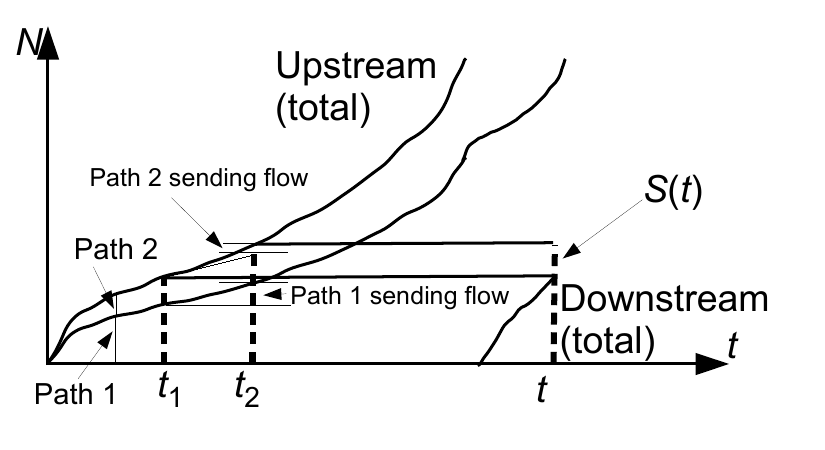}{Disaggregating the sending flow by path.}{0.8\textwidth}

Figure~\ref{fig:disaggregatesendingflowillustration} demonstrates these equations in a different way.
In this diagram, the two upstream counts are a ``stacked line chart'' showing how much of the upstream cumulative count $N^\uparrow$ consisted of flow on path 1, and how much consisted of flow on path 2.
We do not need to track this information at the downstream end, where we only show the aggregate cumulative count $N^\downarrow$.
The sending flow $S(t)$ has already been computed by a link model, but we need to see how much of this consists of flow on path 1, compared to path 2.
We identify the times $t_1$ and $t_2$ when the first and last vehicles in the sending flow entered the link, and compare the path-disaggregated upstream counts at these times to determine the proportion $p$.

These formulas are an approximation, since they assume that the vehicles on different paths in the sending flow are uniformly distributed.
In reality, vehicles towards the start of the sending flow may have a different mix of paths than vehicles towards the end of the sending flow.
This means that the $p$ value (calculated based on the entire sending flow) and the actual vehicles which are moved (which are from the start of the sending flow) may not be entirely consistent.
In Figure~\ref{fig:flowdisagg}, imagine that only three vehicles are able to leave link $(h,i)$ after we perform the node model computation at the diverge.
In this case, all of the vehicles should be turning right, even though the diverge model was computed assuming that a quarter of the vehicles would turn left.
Effectively, the procedure described in this section assumes that all of the vehicles in the sending flow are ``uniformly mixed'' regarding their path, rather than treated in a strict first-in, first-out manner.
In practice, this violation is typically small, and is ignored.
If you wish to enforce first-in, first-out\index{first-in, first-out (FIFO)} discipline exactly, an iterative process can be used to ensure that the proportions $p$ are consistent with the proportion of vehicles between $t_1$ and $t_3$.
\index{splitting proportions!computing from loading|)}

\subsection{Principle of dynamic user equilibrium}
\label{sec:dueprinciple}

With the procedures in the previous two subsections, we can now iteratively perform network loading and find time-dependent shortest paths.
The dynamic traffic assignment problem is to find a mutually consistent solution between these two models.
Let $H$ be a matrix of path flows indicating the number of vehicles departing on each path at each time, and let $C$\label{not:Ctimes} be a matrix of path travel times, indicating the travel time for each path and departure time.
Then the network loading can be concisely expressed by
\labeleqn{concisednl}{C = \mathcal{N} (H)}
where $\mathcal{N}$\label{not:fancyN} is a function encompassing whatever network loading is performed, mapping the path flows to path travel times.
Likewise, let the assumed behavioral rule be denoted by $\mathcal{B}$,\label{not:fancyB} which indicates the allowable path flow matrices $H$ when the travel times are $C$.
For instance, one may assume that departure times are fixed, but only minimum-travel time paths must be used; that minimum-travel time paths must be used but the departure times must minimize total cost (travel time plus schedule delay); or some other principle.
In general, there may be multiple $H$ matrices which satisfy this principle, as when multiple paths have the same, minimal travel time, in which case $\mathcal{B}$ actually denotes a \emph{set} of $H$ matrices.
So, we can express the behavioral consistency rule by
\labeleqn{concisebeh}{H \in \mathcal{B}(C)\,.}

We are now in a position to define the solution to the dynamic traffic assignment problem as a fixed point problem:\index{fixed point problem!applications}\index{dynamic traffic assignment!equilibrium!as fixed point} \emph{the path flow matrix $H$ is a dynamic user equilibrium if}
\labeleqn{concisedta}{H \in \mathcal{B}(\mathcal{N}(H))\,.}
That is, the path flow matrix must be consistent with the driver behavior assumptions, when the travel times are obtained from that same path flow matrix.
The most common behavioral rule is that departure times are fixed, but only minimum-travel time paths will be used.
In this case the principle can be stated more intuitively as \emph{all used paths connecting the same origin to the same destination at the same departure time must have equal and minimal travel time.}  If departure times can be varied to minimize total cost, then the principle can be stated as \emph{all used paths connecting the same origin to the same destination have the same total cost, regardless of departure time.}

\index{dynamic traffic assignment!equilibrium!as variational inequality}
The dynamic user equilibrium solution can also be stated as a solution to a variational inequality.
For the case of fixed departure times, the set of feasible path flow matrices is given by
\labeleqn{feasiblefixed}{\mc{H} = \myc{ H \in \bbr^{T \times |\Pi|}_+ : \sum_{\pi \in \Pi_{rs}} h^\pi_t = d^{rs}_t \quad \forall (r,s) \in Z^2, t \in \myc{1, \ldots, T} }\,,}
that is, matrices with nonnegative entries where the sum of all path flows connecting the origin $r$ to destination $s$ at departure time $t$ is the corresponding value in the time-dependent OD matrix $d^{rs}_t$.
Then, the dynamic user equilibrium solution $\hat{H}$\label{not:Hhat} satisfies
\labeleqn{duevi}{\mathcal{N}(\hat{H}) \cdot (\hat{H} - H) \leq 0 \qquad \forall H \in \mc{H}\,, }
where the product $\cdot$ is the Frobenius product,\index{matrix!Frobenius product} computed by treating the matrices as vectors and calculating their dot product.

In the case of departure time choice, we can define $C = \mathcal{S} (H)$ to be the matrix of total costs for each path and departure time, where $\mathcal{S}$\label{not:fancyS} is obtained by composing the schedule delay function~\eqn{vickreydelay} with the network loading mapping $\mathcal{N}$.
In this case, the set of feasible path flow matrices is given by
\labeleqn{feasibledtc}{\mc{H} = \myc{ H \in \bbr^{T \times |\Pi|}_+ : \sum_{\pi \in \Pi_{rs}} \sum_t h^\pi_t = d^{rs} \quad \forall (r,s) \in Z^2 }\,,}
where $d^{rs}$ is the OD matrix giving total flows between each origin and destination throughout the analysis period.
The variational inequality in this case is
\labeleqn{dtcvi}{\mathcal{S}(\hat{H}) \cdot (\hat{H} - H) \leq 0 \qquad \forall H \in \mc{H}\,.}

Although these fixed point and variational inequality formulations encode the equilibrium principle, and concretely specify the problem at hand, little else can be proven concretely except in the case of simple link and node models (such as point queues).\index{point queue}
The vast majority of theoretical results on variational inequalities and fixed point problems relies on regularity conditions, such as continuity or monotonicity of the mappings $\mc{N}$ or $\mc{B}$.
For realistic link models and node models, these mappings are not continuous.
For example, queue spillback\index{queue spillback} imposes discontinuities on the network loading mapping.
This means that we usually cannot prove that dynamic user equilibrium always exists or is unique, and indeed later in this chapter we present counterexamples where this need not be the case.
However, the variational inequality can be used to define gap measures that can measure the degree of compliance with the equilibrium principle.
In practice, whether the increased realism of a dynamic traffic model outweighs these disadvantages is case-dependent, and these considerations enter into how the appropriate modeling tools are chosen.

\section{Solving for Dynamic Equilibrium}
\label{sec:duealgos}

This section describes how dynamic user equilibrium solutions can be found.
All of the algorithms in this section follow the same general framework, based on iterating between network loading, finding time-dependent shortest paths, and a path updating procedure, and they share common termination criteria.
Three specific path updating procedures are discussed here: the convex combinations method, simplicial decomposition, and gradient projection.
All of these methods are analogues of similar algorithms for static traffic assignment, which were discussed in Chapters~\ref{chp:solutionalgorithms} and~\ref{chp:staticextensions}.
The presentation here is mostly self-contained, but a few references are made to these earlier chapters to avoid repeating a lengthy derivation or explanation.
These are all heuristics, unless we make stronger assumptions on the properties of the network loading.

This is demonstrated more fully in the next section, where we show that dynamic traffic assignment does not always share the neat solution properties of static assignment.
These methods are described only for the case of fixed departure times, but it is not difficult to extend any of them to the case of departure time choice with schedule delay.

\subsection{General framework}
\label{sec:duegeneral}

The general framework for almost all dynamic traffic assignment algorithms involves iterating three steps: network loading, as discussed in Chapter~\ref{chp:networkloading}, time-dependent shortest paths, as discussed in Chapter~\ref{chp:tdsp}, and updating the path flow matrix $H$, as discussed here.
The algorithm takes the following form:
\begin{enumerate}
\item Initialize the path flow matrix $H$ to a feasible value in $\mc{H}$. 
\item Perform the network loading using path flows $H$, obtaining path travel times $C = \mc{N}(H)$. 
\item Identify the time-dependent shortest paths. 
\item Check termination criteria; if sufficiently close to dynamic user equilibrium, stop.
\item Update the path flow matrix $H$ and return to step 2.
\end{enumerate}

A few steps in this algorithm warrant further explanation.
The choice of the initial solution is arbitrary.
One reasonable choice is to identify one or more low travel-time paths between each OD pair based on free-flow travel times, and to divide the demand in the OD matrix between them.
The effort and time involved in identifying the initial solution should be balanced against the benefit obtained by a more careful choice, as opposed to starting with a simpler choice and just running additional iterations.

The termination criteria in step 4 can take several forms.
The most theoretically sound termination criterion is a gap measure, comparing the path flows $H$ and travel times $C$ with the assumed behavioral rule for travelers.
Such measures are preferred to simpler criteria, such as stopping when the path flows do not change much from iteration to iteration --- if the latter occurs, there is no way to tell whether this is because we are near an equilibrium solution, or because the algorithm is unable to make any more progress and has gotten stuck near a non-equilibrium solution.
For the case of fixed demand, the behavioral rule is that travelers must take a least travel-time path between their origin and destination, for their departure time.
After performing the time-dependent shortest path step, we can calculate the travel time on such a path for each OD pair $(r,s)$ and departure time $t$; call this value $\tau^*_{rs,t}$.
If all travelers in the network were on such paths, then the total travel time experienced by all travelers in the network would be given by
\labeleqn{sptt}{SPTT = \sum_{(r,s) \in Z^2} \sum_t d^{rs}_t \tau^*_{rs,t}\,.}
This is often known as the \emph{shortest path travel time}.\index{shortest path travel time}
This is contrasted with the \emph{total system travel time},\index{total system travel time} which is the actual total travel time spent by all travelers on the network:
\labeleqn{dtatstt}{TSTT = \sum_{\pi \in \Pi} \sum_t h^t_\pi \tau^t_\pi = H \cdot C\,.}

Clearly $TSTT \geq SPTT$, and $TSTT = SPTT$ only if all travelers are on least-time paths, which corresponds to dynamic user equilibria.
Therefore, the gap between these values is a quantitative measure of how close the path flows $H$ are to satisfying the equilibrium principle.
Two common gap values involve normalizing the difference $TSTT - SPTT$ in two ways.
The \emph{relative gap}\index{relative gap} $\gamma$ normalizes this difference by SPTT,
\labeleqn{relativegap}{\gamma = \frac{TSTT - SPTT}{SPTT}\,,}
while the \emph{average excess cost}\index{average excess cost} normalizes this difference by the total number of travelers on the network,
\labeleqn{averageexcesscost}{AEC = \frac{TSTT - SPTT}{\sum_{(r,s) \in Z^2} \sum_t d^{rs}_t}\,.}
The algorithm can be terminated whenever either of these measures is sufficiently small.
An advantage of the average excess cost is that it is measured in time units, and has the intuitive interpretation as being the average difference between a traveler's actual travel time, and the shortest path travel time available to them.
Similar expressions can be derived for departure time choice, with schedule delay.

Finally, the last step (updating the path flow matrix) requires the most care, and this section presents three alternatives.
The convex combinations method is the simplest and most economical in terms of computer memory.
The simplicial decomposition method converges faster, by storing previous solutions and using them to find a better improvement direction for $H$.
The third method is gradient projection, which involves computation of an approximate derivative of path travel times with respect to path flow.
This derivative enables the use of Newton's method to identify the amount of flow to shift between paths.

\subsection{Convex combinations method}

\index{dynamic traffic assignment!algorithms!convex combinations|(}
The convex combinations method is the simplest way to update the path choice matrix.
In this method, after computing time-dependent shortest paths, one identifies a ``target matrix'' $H^*$\label{not:Hstar} which would give the path flows travelers would choose if the path travel times were held constant at their current values $C$.
For the case of fixed departure times, for each OD pair $(r,s)$ and departure time $t$, all of the demand $d^{rs}_t$ is loaded onto its time-dependent shortest path in $H^*$, and all other path flows are set to zero.
For the case of departure time choice, for each OD pair $(r,s)$ the departure time $t^*$ and corresponding path with minimal schedule delay is identified, and all $d^{rs}$ travelers associated with this OD pair are assigned to that path and departure time, with all other path flows set to zero.
This is often called an \emph{all-or-nothing assignment}, because it involves choosing a single  alternative to assign an entire class of travelers.

The basic form of the convex combinations method updates $H$ using the formula
\labeleqn{convexcombinations}{H \leftarrow \lambda H^* + (1 - \lambda) H\,,}
where $\lambda \in [0, 1]$ is a step size showing how far to move in the direction of the target.
If $\lambda = 0$, then the path flow matrix is unchanged, whereas if $\lambda = 1$, the path flow matrix is set equal to the target.
Intermediate values produce a solution which is a weighted average of the current solution $H$ and the target $H^*$.
There are a number of ways to choose the step sizes $\lambda$.
The simplest alternative is to use a decreasing sequence of values, choosing for the $i$-th iteration the step size $\lambda_i = 1/(i + 1)$, so the step sizes form the sequence 1/2, 1/3, 1/4, and so on.
It is common to use sequences such that $\sum \lambda_i = +\infty$ and $\sum \lambda_i^2 < \infty$, because we do not want the sequence of $H$ values to converge too quickly (or else it may stop short of the equilibrium), but if the step size drops too slowly, the algorithm is likely to ``overshoot'' the equilibrium and oscillate.

More sophisticated variations are possible as well.
Rather than using a fixed step size, $\lambda$ can be chosen to minimize, say, the relative gap or average excess cost.
While this likely decreases the number of iterations required to reach a small gap, the computation at each iteration is increased.
In particular, each $\lambda$ value tested requires an entire network loading step to determine the corresponding gap, which can be a significant time investment.
It is also possible to vary the $\lambda$ value for different path flows and departure times.
Some popular heuristics are to use larger $\lambda$ values for OD pairs which are further from equilibrium, or to vary $\lambda$ for different departure times --- since travel times for later departure times depend on the path choices for earlier departure times, it may not make sense to invest much effort in equilibrating later departure times until earlier ones have stabilized.
\index{dynamic traffic assignment!algorithms!convex combinations|)}

\subsection{Simplicial decomposition}

\index{dynamic traffic assignment!algorithms!simplicial decomposition|(}
The simplicial decomposition algorithm is an extension of the convex combinations method that remembers the target matrices $H^*$ from earlier iterations.
This requires more computer memory, but storing these previous target matrices leads to a more efficient update of the path flow vector.
In particular, choosing an all-or-nothing assignment as a target matrix essentially restricts the choice of search direction to the corner points of the feasible region.
As a result, when the equilibrium solution is reached, the convex combinations method will ``zig-zag,'' taking a number of short, oblique steps rather than a more direct step towards the equilibrium problem.
(Figure~\ref{fig:cc_sd}).
Simplicial decomposition is able to do so by combining multiple target points into the search direction.
In this algorithm, the set $\mc{H}^*$\label{not:fancyHstar} is used to store all of the target path flow matrices found thus far.

\stevefig{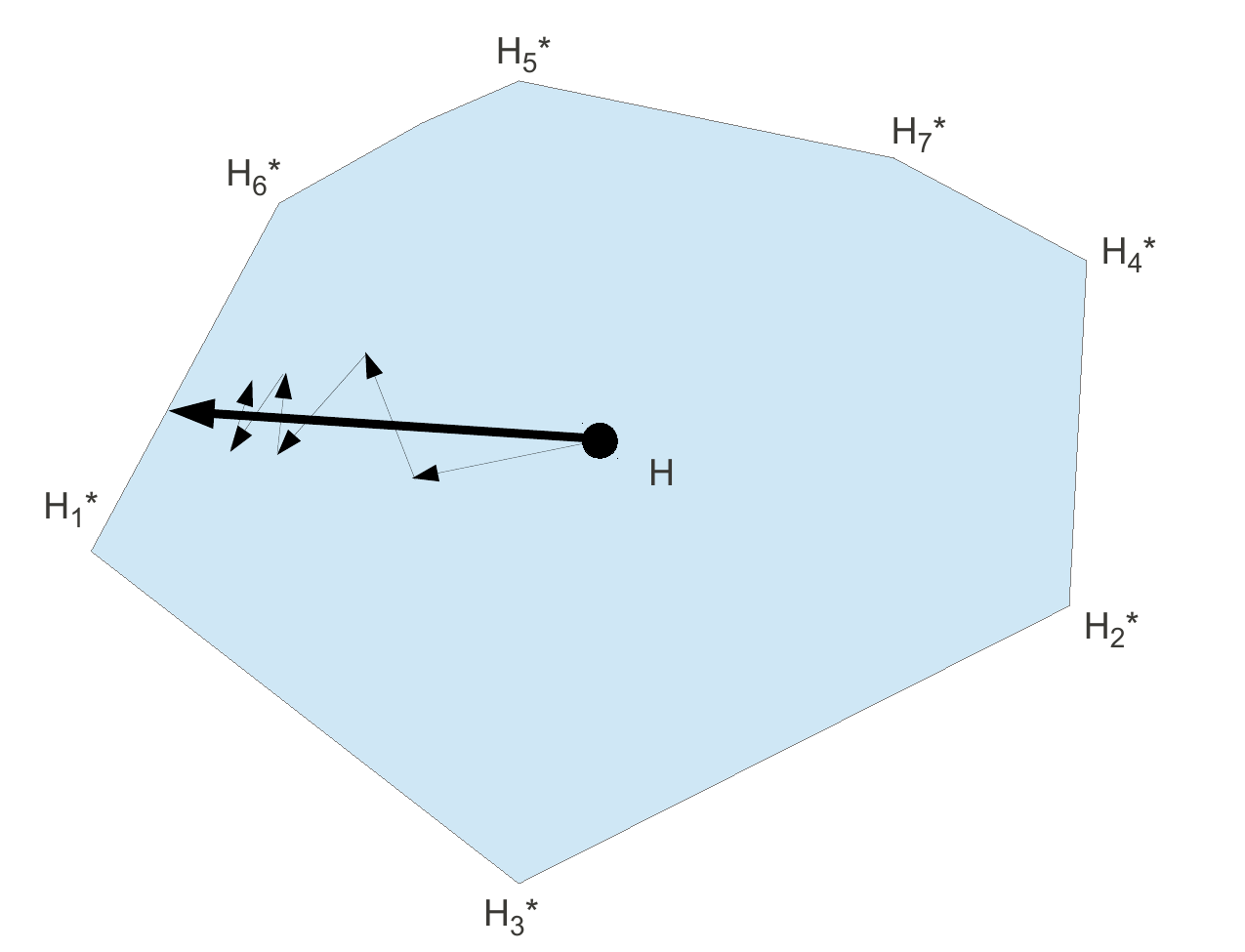}{Simplicial decomposition can move towards non-corner points.}{0.8\textwidth}

Many components of this algorithm are similar to the use of simplicial decomposition in Section~\ref{sec:staticsimplicial}, but we review them here.
Given a set $\mc{H}^* = \myc{H^*_1, H^*_2, \cdots, H^*_k}$, we say that a time-dependent path flow matrix $\hat{H}$ is a \emph{restricted equilibrium relative to $\mc{H}^*$} if it solves the variational inequality
\labeleqn{restrictedvisd}{\mc{N}(\hat{H}) \cdot (\hat{H} - H^*_i) \leq 0 \qquad \forall H^*_i \in \mc{H}^*\,.}
This is different from the variational inequality~\eqn{duevi} because the only possible choices for $H$ are the previously-generated target solutions, rather than \emph{any} of the feasible assignments.
Essentially, $H$ is a restricted equilibrium if none of the targets in $\mc{H}^*$ lead to improving directions in the sense that the total system travel time would be reduced by moving to some $H^*_i \in \mc{H}^*$ while fixing the travel times at their current values.

Simplicial decomposition alternates between adding a new target matrix to $\mc{H}^*$, and then finding a restricted equilibrium using the current matrices in $\mc{H}^*$.
Rather than finding an exact restricted equilibrium at each iteration (which is too computationally intensive), we take several ``inner iteration'' steps to move closer to a restricted equilibrium with the current set $\mc{H}^*$ before adding another target.
In each inner iteration, the current solution $H$ is adjusted to $H + \mu \Delta H$, where $\mu$ is a step size and $\Delta H$ is a direction which moves toward restricted equilibrium.
One choice for this direction is
\labeleqn{smithdta}{\Delta H = \frac{\sum_{H^*_i \in \mathcal{H}^*} \left[ \mc{N}(H) \cdot (H - H^*_i) \right]^+ (H^*_i - H)}{\sum_{H^*_i \in \mathcal{H}^*} \left[ \mc{N}(H) \cdot (H - H^*_i) \right]^+ }}
As in simplicial decomposition for static assignment (Section~\ref{sec:staticsimplicial}), this formula is simpler than it appears.
The direction is a weighted average of the directions $H^*_i - H$ (potential moves toward each target in $\mc{H}^*$), where the weight for each potential direction is the extent to which it improves upon the current solution: $\left[ \mc{N}(H) \cdot (H - H^*_i) \right]$ is the reduction in total system travel time obtained by moving from $H$ to $H^*_i$ while holding travel times constant.
If this term is negative, the weight is set to zero (the direction is not helpful), and the denominator simply serves as a normalizing factor by dividing by the sum of all weights.

The step size $\mu$ is chosen through trial-and-error.
One strategy is to iteratively test $\mu$ values in some sequence (say, $1, 1/2, 1/4, \ldots$) until we have found a solution acceptably closer to restricted equilibrium than $H$, as with the Armijo rule\index{Armijo rule} discussed in Appendix~\ref{sec:unconstrainedstepsize}.
``Acceptably closer'' can be calculated using the \emph{restricted average excess cost}\index{gap function!average excess cost!restricted}
\labeleqn{raecdta}{AEC' = \frac{\mc{N}(H) \cdot H - \min_{H^*_i \in \mathcal{H}^*} \myc{\mc{N}(H) \cdot H^*_i}}{\sum_{(r,s) \in Z^2} \sum_t d^{rs}_t}}
which is similar to the average excess cost, but instead of using the shortest path travel time uses the best available target vector in $\mc{H}^*$.

\emph{Unlike} the version of simplicial decomposition used for the static traffic assignment problem, there is no guarantee that a sufficiently small choice of $\mu$ will result in a reduction of restricted average excess cost.
However, in practice, this rule seems to work acceptably.

Putting all of this together, for simplicial decomposition, step 5 of the dynamic traffic assignment algorithm in Section~\ref{sec:duegeneral} involves performing all of these steps as a ``subproblem'':
\begin{enumerate}
\setcounter{enumi}{4}
\item \textbf{Subproblem:} Find an approximate restricted equilibrium $H$ using only the vectors in $\mathcal{H}^*$.
\begin{enumerate}[(a)]
  \item Find the improvement direction $\Delta H$ using equation~\eqn{smithdta}.
  \item Update $H \leftarrow H + \mu \Delta H$, with $\mu$ sufficiently small (to reduce $AEC'$).
  \item Perform network loading to update travel times.
  \item Return to step (a) of subproblem unless $AEC'$ is small enough.
\end{enumerate}
\end{enumerate}
Furthermore, the set $\mc{H}^*$ must be managed.
It is initialized to be empty, and whenever time-dependent shortest paths are found, a new all-or-nothing assignment (the one which would have been chosen as the sole target in the convex combinations method) is added to $\mc{H}^*$.
\index{dynamic traffic assignment!algorithms!simplicial decomposition|)}

\subsection{Gradient projection}

\index{dynamic traffic assignment!algorithms!gradient projection|(}
The gradient projection method updates the path flow matrix $H$ by using Newton's method.\index{Newton's method}
Consider first the simpler problem of trying to shift flows between two paths $\pi_1$ and $\pi_2$ (for the same departure time) to equalize their travel times.
Write $\tau_1(h_1, h_2)$ and $\tau_2(h_1, h_2)$ to denote the travel times on the two paths as a function of the flows $h_1$ and $h_2$ on the two paths.
If we were to shift $\Delta h$ vehicles from path 1 to path 2, the difference in travel times between the paths is given by
\labeleqn{newtonformula}{g(\Delta h) = \tau_1(h_1 - \Delta h, h_2 + \Delta h) - \tau_2(h_1 - \Delta h, h_2 + \Delta h)\,.}
We want to choose $\Delta h$ so that the difference $g(\Delta h)$ between the path travel times is equal to zero.

Using one step of Newton's method, an approximate value of $\Delta h$ is given by
\labeleqn{newtonstep}{\Delta h = -\frac{g(0)}{g'(0)}\,.}
The numerator $g(0)$ is simply the difference in travel times between the two paths before any flow is shifted.
To calculate the denominator, we need to know how the difference in travel times will change as flow is shifted from $\pi_1$ to $\pi_2$.
Computing the derivative of~\eqn{newtonformula} with the help of the chain rule, we have
\labeleqn{newtonder}{g'(0) = \myp{\pdr{t_1}{h_2} + \pdr{t_2}{h_1}} - \myp{\pdr{t_1}{h_1} + \pdr{t_2}{h_2}}\,.}
This formula is not easy to evaluate exactly, but we can make a few approximations.
The first term in~\eqn{newtonder} reflects the impact the flow on one path has on the \emph{other} path's travel time, while the second term reflects the impact the flow on one path has on \emph{its own} travel time.
Typically, we would expect the second effect will be larger in magnitude than the first effect (although exceptions do exist).
So, as a first step we can approximate the derivative as
\labeleqn{newtonderapprox}{g'(0) \approx -\myp{\pdr{t_1}{h_1} + \pdr{t_2}{h_2}}\,.}

The next task is to calculate the derivative of a path travel time with respect to flow along this path.
Unlike in static traffic assignment, there is no closed-form expression mapping flows to travel times, but rather a network loading procedure must be used.
For networks involving triangular fundamental diagrams (which can include the point queue model, cf.\ Section~\ref{sec:pqsqlwr}), the derivative of the travel time on a single link $(i,j)$ with respect to flow entering at time $t$ can be divided into two cases.
In the first case, suppose that the vehicle entering link $(i,j)$ at time $t$ exits at $t' = t + \tau_{ij}(t)$, and that the link is sending-constrained at that time (that is, all of the sending flow $S_{ij}(t')$ is able to move).
In this case, even if a marginal unit of flow is added at this time, the link will remain sending-constrained, and all of the flow will still be able to move.
No additional delay would accrue, so the derivative of the link travel time is
\labeleqn{uncongestedder}{\frac{d\tau_{ij}(t)}{dh} = 0\,.}

In the second case, suppose that at time $t'$ flows leaving link $(i,j)$ are constrained either by the capacity of the link or by the receiving flow of a downstream link.
In this case, not all of the flow is able to depart the link, and a queue has formed at the downstream end of $(i,j)$.
The clearance rate for this queue is given by $y_{ij}(t')$, so the incremental delay added by one more vehicle joining the queue is $1/y_{ij}(t')$, and
\labeleqn{congestedder}{\frac{\tau_{ij}(t)}{dh} = \frac{1}{y_{ij}(t')}\,.} 
(The denominator of this expression cannot be zero, since $t'$ is the time at which a vehicle is leaving the link.)

Therefore, we can calculate the derivative of an entire path's travel time inductively.
Assume that $\pi = [ r, i_1, i_2, \ldots, s]$, and that $t_i$ gives the travel time for arriving at each node $i$ in the path.
Then $dt_\pi/dh_\pi$ is obtained by summing expressions~\eqn{uncongestedder} and~\eqn{congestedder} for the uncongested and congested links in this path, taking care to use the correct time indices:
\labeleqn{pathder}{\frac{dt_\pi}{dh_\pi} \approx \sum_{(i,j) \in \pi} \frac{1}{y_{ij}(t_j)} [y_{ij}(t_j) < S_{ij}(t_j) ]\,,}
where the brackets are an indicator function equal to one if the statement in brackets is true, and zero if false.

We are almost ready to give the formula for~\eqn{newtonderapprox}, but we can make one more improvement to the approximation.
Two paths may share a certain number of links in common.
Define the \emph{divergence node}\index{node!divergence}\index{divergence node|see {node, divergence}} of two paths to be node where the paths diverge for the first time.
(Figure~\ref{fig:divergenode}).
Prior to the divergence node, the two paths include the same links, so there will be no effect of shifting flow from one path to the other.
Therefore, the sum in~\eqn{pathder} need only be taken for links beyond the divergence node.
Even if the paths rejoin at a later point, they may do so at different times, so we cannot say that there is no effect of shifting flow between common links downstream of a divergence node.
So, if $d(\pi_1, \pi_2)$ is the divergence node, then we have
\labeleqn{pathderfinal}{\frac{dt_\pi}{dh_\pi} \approx \sum_{(i,j) \in \pi : (i,j) > d(\pi_1, \pi_2)} \frac{1}{y_{ij}(t_j)} [y_{ij}(t_j) < S_{ij}(t_j) ]\,,}
where the ``$>$'' notation for links indicates links downstream of a node in a path.
This expression can then be used in~\eqn{newtonderapprox} and~\eqn{newtonstep} to find the approximate amount of flow which needs to be shifted from $\pi_1$ to $\pi_2$ to equalize their costs.

\stevefig{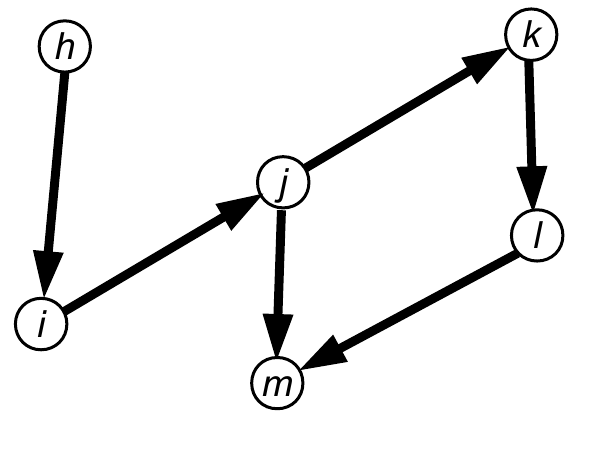}{Node $j$ is the divergence node of these two paths.}{0.6\textwidth}

The procedure for updating $H$ can now be described as follows:
\begin{enumerate}
\setcounter{enumi}{4}
\item \textbf{Subproblem:} Perform Newton updates between all paths and the target path found.
\begin{enumerate}[(a)]
  \item For each path $\pi$ with positive flow, let $\rho$ be the least-travel time path connecting the same OD pair and departure time.
  \item Calculate $\Delta h$ using equations~\eqn{pathderfinal},~\eqn{newtonderapprox}, and~\eqn{newtonstep}, using $\pi$ as the first path and $\rho$ as the second.
  \item If $\Delta h < h_\pi$, then update $h_\rho \leftarrow h_\rho + \Delta h$ and $h_\pi \leftarrow h_\pi - \Delta h$.
Otherwise, set $h_\rho \leftarrow h_\rho + h_\pi$ and $h_\pi \leftarrow 0$.
\end{enumerate}
\end{enumerate}
As stated, the path travel times and derivatives are not updated after each shift is performed.
Doing so would increase accuracy, but greatly increase computation time since a complete network loading would have to be performed.
\index{dynamic traffic assignment!algorithms!gradient projection|)}

\section{Properties of Dynamic Equilibria}
\label{sec:dueproperties}

\index{dynamic traffic assignment!equilibrium!properties|(}
The dynamic traffic assignment problem is much less well-behaved than the static traffic assignment problem,\index{static traffic assignment} where a unique equilibrium provably exists under natural conditions.
This is because the more realistic network loading procedures used in dynamic traffic assignment lack the regularity properties which make static assignment more convenient mathematically.
Recall that the dynamic user equilibrium path flows solve the variational inequality\index{dynamic traffic assignment!equilibrium!as variational inequality}
\labeleqn{duevi2}{\mathcal{N}(\hat{H}) \cdot (\hat{H} - H) \leq 0 \qquad \forall H \in \mc{H}\,.
}
Showing existence of a solution to this variational inequality would typically rely on showing that $\mathcal{N}$ is continuous.
However, when queue spillback is modeled, the travel times are discontinuous in the path flow variables.
Showing uniqueness of a solution to the variational inequality typically involves showing that $\mathcal{N}$ has some flavor of monotonicity, but common node models can be shown to violate monotonicity even when we restrict attention to networks involving only a single diverge and merge.
Thus, the failure of existence or uniqueness results is directly traced to the features that make dynamic traffic assignment a more realistic model.
As an analyst, you should be aware of this tradeoff.

This section catalogs a few examples of dynamic traffic assignment problems highlighting some unusual or counterintuitive results.
We begin with an example where the choice of node model results in no dynamic user equilibrium solution existing.
Two examples of equilibrium nonuniqueness are given, first where multiple user equilibrium solutions exist, and second where literally \emph{all} feasible path flow matrices are user equilibrium solutions, even though they vary widely in total system travel time.
We conclude with a dynamic equivalent to the Braess paradox which was developed by Daganzo, in which increasing the capacity on a link can make network-wide conditions arbitrarily worse, because of queue spillback.

\subsection{Existence: competition at intersections}
\label{sec:eggnetwork}

\index{dynamic traffic assignment!equilibrium!nonexistence|(}
The network in Figure~\ref{fig:eggnetwork} has two origin-destination pairs (A to B, and C to D).
Aside from these four centroids, there are two nodes representing intersections.
Turns are not allowed at these intersections, so $\Xi(1) = \myc{[A,1,2],[C,1,D]}$ and $\Xi(2) = \myc{[1,2,B],[C,2,D]}$.
As a result, each OD pair has only two routes available to it: $[A,B]$ and $[A,1,2,B]$ for A to B, and $[C,1,D]$ and $[C,2,D]$ for C to D.
For the sake of convenience, call these four routes ``top,'' ``bottom,'' ``left,'' and ``right,'' respectively.
The figure shows the \emph{travel times} with each link.
Two links in the network also carry a toll, expressed in the same units as travel time, and drivers choose routes to minimize the sum of travel time and toll.
Capacities and jam densities are large enough that congestion will never arise, and all links will remain at free flow.

\stevefig{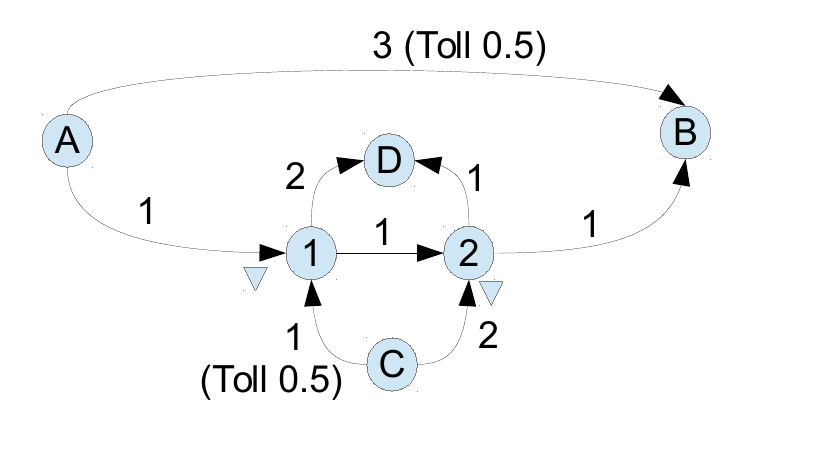}{Network for example with no equilibrium solution.
Link travel times shown.}{0.8\textwidth}

Nodes 1 and 2 have a distinctive node model, representing absolute priority of one approach over the other.
For Node 1, all flow on the movement $[A,1,2]$ must yield to flow on $[C,1,D]$.
This can be expressed with the following node model:
\begin{align}
   y_{C1D} &= \min \myc{S_{C1}, R_{1D}} \label{eqn:eggnode1} \\  
   y_{A12} &= \begin{cases}
               0                         & \mbox{ if } y_{C1D} > 0 \\
               \min \myc{S_{A1}, R_{12}} & \mbox{ otherwise } 
             \end{cases} \label{eqn:eggnode2}
\end{align}
Likewise, for Node 2, all flow on the movement $[C,2,D]$ must yield to flow on $[1,2,B]$.
That is, traffic from C to D has priority at node 1; and traffic from A to B has priority at node 2.
These relationships are indicated on the figure with triangles next to the approach that must yield.

Each OD pair has one unit of demand, that departs during the same time interval.
We now show that there is no assignment of demand to paths that satisfies the equilibrium principle, requiring any used path to have minimal travel time.
We begin first by showing that any equilibrium solution cannot split demand among multiple paths; each OD pair must place all of its flow on just one path or the other.
With the given demand, the cost on the bottom path is either 3 (if there is no flow on the left path from C to D) or 4 (if there is), in either case, the cost is different from that on the top path (3.5).
All demand will thus use either the top path or the bottom path, whichever is lower.
Similarly, for the OD pair from C to D, the cost on the right path is either 3 or 4, whereas that on the left is always 3.5.
Both paths cannot have equal travel time, so at most one path can be used.

Thus, there are only four solutions that can possibly be equilibria: all flow from A to B must either be on the top path or the bottom path; and all flow from C to D must either be on the left path or the right path.
First consider the case when all flow from A to B is on the top path.
The left path would then have a cost of 3.5, and the right path a cost of 3.
So all flow from C to D must be on the right path.
But then the travel times on the top and bottom paths are 3.5 and 3, respectively, so this solution cannot be an equilibrium (the travelers from A to B would prefer the bottom path, contradicting our assumption).

So now consider the other possibility, when all flow from A to B is on the bottom path.
The travel times on the left and right paths are now 3.5 and 4, so all travelers from C to D would pick the left path.
But then the travel times on the top and bottom paths are 3.5 and 4, respectively, so this solution is not an equilibrium either.

Therefore, there is \emph{no} assignment of vehicles to paths that satisfies the principle of user equilibrium.
This network is essentially the Ginger/Harold ``matching pennies'' game from Section~\ref{sec:notionofequilibrium} encoded into dynamic traffic assignment.
Game theorists\index{game theory} often resolve this type of situation by proposing a \emph{mixed-strategy equilibrium},\index{equilibrium!mixed-strategy} in which the players randomize their actions.
One can show that a mixed-strategy equilibrium exists in this network (Exercise~\ref{ex:eggmixed}), but this equilibrium concept is usually not applied to transportation planning.
This is for both computational reasons (when there are many players, as in practical planning networks, computing mixed-strategy equilibria is hard) and for modeling reasons (typically travelers do not randomize their paths with the intent of ``outsmarting'' other travelers).

Mathematically, the reason no equilibrium exists is because the node model specified by~\eqn{eggnode1} and~\eqn{eggnode2} is not continuous in the sending and receiving flows.

This is similar to static traffic assignment, where equilibrium existence cannot be guaranteed if the link performance functions are not continuous.
In dynamic traffic assignment, when building node models that capture absolute priority (as at yields, two-way stops, or signals with permitted turns), ensuring continuity is difficult.
In part, this is the motivation for using ratios of $\alpha$ values in Chapter~\ref{chp:networkloading} to reflect priorities at merges and general intersections.

But even if the link and node models are chosen to be continuous, dynamic equilibrium need not exist.
Repeating the arguments made in Chapter~\ref{sec:tapproperties} would also require the travel time calculations to be continuous in the link cumulative entrances and exits $N^\uparrow$ and $N^\downarrow$.
The use of a minimum to disambiguate multiple possible values of travel times in Equation~\eqn{inverseproper} can mean that small changes in a link's cumulative entries or exits could cause a large change in the travel time (imagine how the length of the dashed line in Figure~\ref{fig:noflow} would change if one of the corner points were shifted slightly).
Continuity could be assured only if the cumulative counts were \emph{strictly} increasing, that is, if vehicles were constantly flowing into and out of each link.
Therefore, although one can construct dynamic traffic assignment models which guarantee existence of equilibria, these require quite strong and  significant restrictions on the network loading and travel time calculation procedures.
\index{dynamic traffic assignment!equilibrium!nonexistence|)}

\subsection{Uniqueness: Nie's merge}
\label{sec:niemerge}

\index{dynamic traffic assignment!equilibrium!nonuniqueness|(}
The network in Figure~\ref{fig:niesmerge} consists of a single diverge and merge.
One can visualize this network as a stylized version of a freeway lane drop (Figure~\ref{fig:lanedropnie}), where the reduction from two lanes to one lane reduces the capacity from 80 vehicles per minute to 40 vehicles per minute, and where drivers either merge early (choosing the top lane at the diverge point) or late (waiting until the lane drop itself).
The inflow rate from upstream is 80 vehicles per minute.
Since the capacity of the downstream exit is only 40 vehicles per minute, the excess vehicles will form one or more queues in this network.

\stevefig{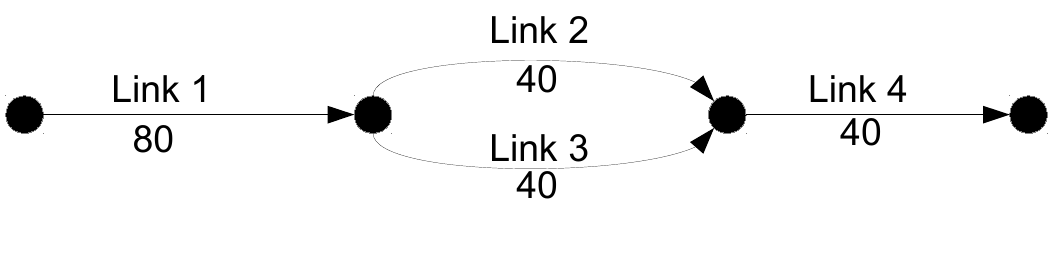}{Nie's merge network.
All capacity values in vehicles per minute.}{0.8\textwidth}
\stevefig{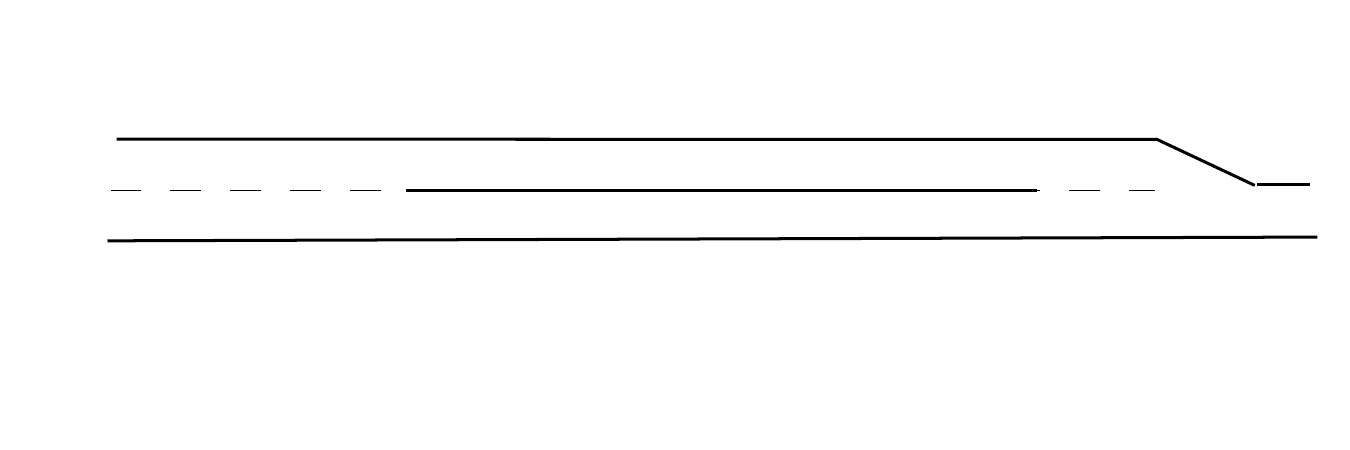}{A physical interpretation of Nie's merge.}{0.8\textwidth}

To simplify the analysis, assume that the time horizon is short enough that none of these queues will grow to encompass the entire length of the link.
With this assumption, queue spillback can be ignored, and we can focus on the issue of route choice.
Furthermore, under these assumptions, $S_1$ is always 80 vehicles per minute, and $R_2$, $R_3$, and $R_4$ are always 40 vehicles per minute.
We work with a timestep of one minute, and will express all flow quantities in vehicles per minute.\footnote{This is larger than what is typically used in practice, but simplifies the calculations and does not affect the conclusions of this example.}

In this network, there is only one choice of routes (the top or bottom link at the diverge).
We will restrict our attention to path flow matrices $H$ where the proportions of vehicles choosing the top and bottom routes are the same for all departure times, and show that there are multiple equilibrium solutions even when limiting our attention to such $H$ matrices.
There may be still more equilibria where the proportion of vehicles choosing each path varies over time.

Therefore, the ratio $h_{124} / h_{134}$ is constant for all time intervals, which means that the splitting fractions $p_{12}$ and $p_{13}$ at the diverge point are also constant, and equal to $h_{124}/40$ and $h_{134}/40$, respectively.

Using the diverge and merge models from Section~\ref{sec:diverge}, we can analyze the queue lengths and travel times on each link as a function of these splitting fractions.

It turns out that three distinct dynamic user equilibria exist:
\begin{description}
\item[Equilibrium I: ] If $p_{12} = 1$ and $p_{13} = 0$, then the diverge formula~\eqn{generaldiverge} gives the proportion of flow which can move as 
\labeleqn{nieeqmi}{\phi = \min \myc{ 1, \frac{40}{80}, \frac{40}{0} } = \frac{1}{2}\,.}
Therefore, the transition flows at each time step are $y_{12} = 40$ and $y_{13} = 0$.
Since $y_{12} < S_{12}$, a queue will form on the upstream link, and its sending flow will remain at $S_1 = 80$.
Therefore, once the first vehicles reach the merge, we will have $S_2 = 40$ and $S_3 = 0$.
Applying these proportions, together with $R_4 = 40$, the merge formula gives $y_{24} = 40$ and $y_{34} = 0$.
There will be no queue on link (2,4), so both link (2,4) and link (3,4) are at free-flow.
This solution is an equilibrium: the two  paths through the network only differ by the choice of link 2 or link 3, and both of these have the same travel time since they are both at free-flow.
In physical terms, this corresponds to all drivers choosing to merge early; the queue forms upstream of the point where everyone chooses to merge, and there is no congestion downstream.
\item[Equilibrium II: ] If $p_{12} = 0$ and $p_{13} = 1$, the solution is exactly symmetric to that of Equilibrium I.
A queue will again form on the upstream link, now because all drivers are waiting to take the bottom link.
In physical terms, this corresponds to all drivers choosing to merge early, but to the lane which is about to drop.
They then all merge back when the lane actually ends.
This solution does not seem especially plausible, but it does satisfy the equilibrium condition.
\item[Equilibrium III: ] If $p_{12} = p_{13} = 1/2$, drivers wish to split equally between the top and bottom links.
The proportion of flow which can move at the diverge is
\labeleqn{nieeqmiii}{\phi = \min \myc{ 1, \frac{40}{40}, \frac{40}{40} } = 1\,,}
so all vehicles can move and there is no queue at the diverge: $y_{12} = y_{13} = 40$.
Once these vehicles reach the merge, we will have $S_2 = S_3 = 40$ and $R_4 = 40$.
The merge formula~\eqn{mergecase2} then gives $y_{24} = y_{34} = 20$, so queues will form on both merge links.
However, since the inflow and outflow rates of links 2 and 3 are identical, the queues will have identical lengths, and so the travel times on these links will again be identical.
Therefore, this solution satisfies the principle of dynamic user equilibrium.
In physical terms, this is the case when no drivers change lanes until the lane actually ends, and queueing occurs at the merge point.
\end{description}

By examining the intermediate solutions, we can calculate the travel times on both paths for different values of the fraction $p_{12}$.
(We do not need to state the values of the other fraction, since $p_{13} = 1 - p_{12}$.)  This is shown in Figure~\ref{fig:niestability}, and the three equilibria correspond to the crossing points of the paths.
Note that all three of the equilibria share the same equilibrium travel time.
Although the queues in Equilibrium III are only half as long as those in Equilibria I and II, being split between two links, the outflow rates of these queues are also only half as great (20 veh/min instead of 40 veh/min).

\stevefig{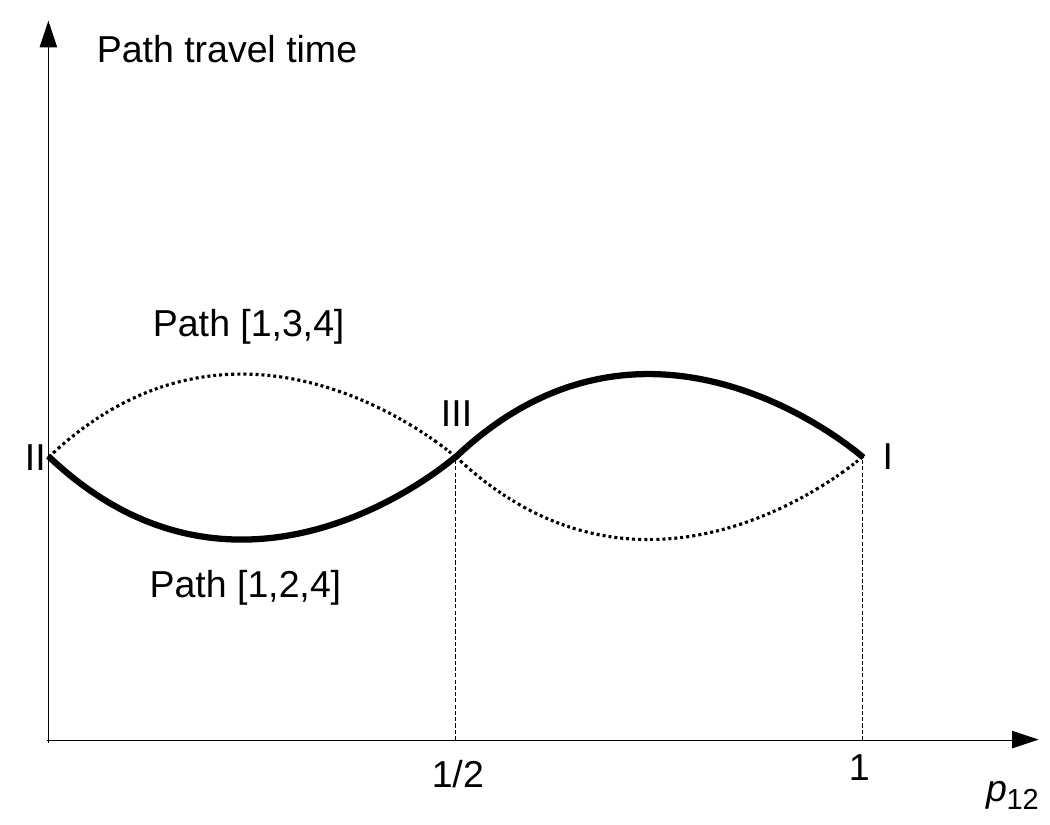}{Travel time on top and bottom paths as a function of splitting fraction, with three equilibrium solutions marked.}{0.8\textwidth}

This example shows that the dynamic user equilibrium solution is not unique, even in an extremely simple network.
This nonuniqueness has practical consequences.
The effect of a potential improvement to link 2 or 3 will depend crucially on the number of travelers on the link, which varies widely in all three equilibria.
One criterion for distinguishing which of these equilibria is more likely is \emph{stability},\index{equilibrium!stability} which explores what would happen if the equilibria were slightly perturbed.
As shown in Figure~\ref{fig:niestability}, if one begins at Equilibrium I and perturbs $p_{13}$ to a small value (reducing $p_{12}$ by the same amount), we see that the travel time on the top path increases, where as the travel time on the bottom path decreases.
It may seem odd that increasing flow on the bottom path decreases its travel time --- what is happening is that congestion at the diverge decreases (lowering the travel time of both paths), but congestion forms on the top link at the merge (increasing the travel time just of the top path).
Superimposing these effects produces the result in Figure~\ref{fig:niestability}.
As a result, travelers will switch from the (slower) top path to the (faster) bottom one, moving us even further away from Equilibrium I.
Therefore, Equilibrium I is not stable.
The same analysis holds for Equilibrium II.

Equilibrium III, on the other hand, is stable.
If a few travelers switch from the top to the bottom path, the travel time on the top path decreases and that on the bottom path increases.
Therefore, travelers will tend to switch back to the top path, restoring the equilibrium.
The same holds if travelers switch from the bottom path to the top path.
This gives us reason to believe that Equilibrium III is more likely to occur in practice than Equilibrium I or II.
This type of analysis is much more complicated in larger networks, and for the most part is completely unexplored.
Coming to a better understanding of the implications of nonuniqueness in large networks, as well as techniques for addressing this, is an important research question.

\subsection{Nie's merge redux}
\label{sec:bdmr}

Modifying the merge network from the previous section, we can produce an even more pathological result.
We now eliminate the lane drop entirely, but preserve the diverge/merge network, by setting the capacity of the downstream link to 80 vehicles per minute, the same as the upstream link.
This might represent a roadway section where lane changing is prohibited, perhaps in a work zone or with a high-occupancy/toll lane.
Now, at the merge point $S_2 + S_3 \leq R_4$ no matter what the flow pattern is on the network.
This means that the merge will always be freely-flowing, and there will be no queues at the merge.
Links 2, 3, and 4 will be at free-flow regardless of the path flows $H$.

However, queues can still occur at the \emph{diverge}, where the analysis is the same as before.
If $p_{12} = p_{13} = 1/2$, no queues will form at the diverge point, whereas for any other values of $p_{12}$ and $p_{13}$ formula~\eqn{nieeqmi} will predict queues on the upstream link, reaching their maximum length if $p_{12} = 1$ or $p_{13} = 1$.
However, \emph{all of these solutions satisfy the principle of equilibrium} because the only delay occurs on the upstream link, which is common to both paths.
There is zero relative gap or average excess cost, no matter what the values of $p$ are, since no traveler can reduce their travel time by choosing a different route.
Choosing a different route could influence the travel time of those further upstream, but under the assumption that drivers act only to minimize their own travel time, this influence on other travelers is of no concern.

Figure~\ref{fig:bdmrqueues} shows the queue lengths and total system travel time as $p$ varies.
The system-optimal solution is unique: if $p_{12} = p_{13} = 1/2$ there are no queues in the system and all vehicles experience free-flow travel times.
For any other values of $p_{12}$ and $p_{13}$, a queue will form at the diverge and some delay is experienced.
Yet \emph{all} possible $p$ values satisfy the equilibrium principle, since the delay is upstream of the diverge and drivers cannot choose another route to minimize their own travel time.
Furthermore, as the time horizon grows longer, the difference in total travel time between the worst of the equilibria (either $p_{12} = 0$ or $p_{12} = 1$) and the system optimum solution can grow arbitrarily large.\index{price of anarchy}
This is in contrast with the bounded ``price of anarchy'' which can often be found for static traffic assignment (cf.\ Section~\ref{sec:priceofanarchy}).

\stevefig{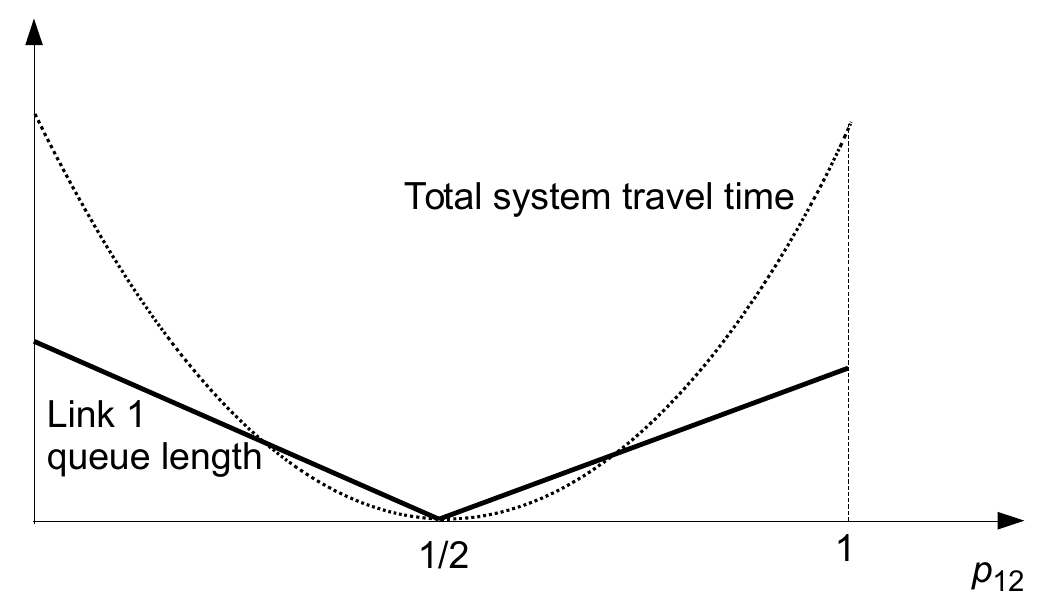}{Queue length and total travel time for different splitting proportions $p$ in the modified Nie's merge network. Note that all $p$ values represent equilibria.}{0.7\textwidth}

This effect is less of a practical concern than the original version of Nie's merge --- there are good behavioral reasons to doubt that a significant imbalance of travelers will choose one alternative over another identical one.
One argument is from entropy principles (Section~\ref{sec:maxentropy}): if travelers have the same behavior assumption, it is unlikely they would all choose one route over another with equal travel time.
Furthermore, the assumption of a triangular fundamental diagram implies that travel speeds remain at free-flow for all subcritical densities.
In practice the speed will drop slightly due to variations in preferred speeds and difficulties in overtaking at higher density, so drivers would likely prefer the route chosen by fewer travelers.

However, from the standpoint of \emph{modeling} the fact that all feasible solutions are equilibria poses significant challenges.
Literally any solution will have zero gap, and if an all-or-nothing assignment is chosen as the initial solution (as is sometimes done in implementations), dynamic traffic assignment software will report that a perfect equilibrium has been reached.
This example shows that initial solutions for dynamic traffic assignment should be carefully chosen, perhaps by spreading vehicles over multiple paths, or breaking ties stochastically in shortest path algorithms to avoid assigning all vehicles to the same path in the first all-or-nothing assignment.
\index{dynamic traffic assignment!equilibrium!nonuniqueness|)}

\subsection{Efficiency: Daganzo's paradox}
\label{sec:daganzoparadox}

\index{dynamic traffic assignment!equilibrium!inefficiency|(}
\index{queue spillback|(}
\index{node model!diverge|(}
The previous section showed that the worst user equilibrium solution can be arbitrarily worse than the system optimal solution, in terms of total travel time.
However, there was still one user equilibrium solution that was also system optimal (the case where $p_{12} = p_{13} = 1/2$).
This need not be in the case.
Here we present an example where the \emph{only} user equilibrium solution can be arbitrarily worse than system optimum.

Furthermore, it is ``paradoxical'' in the sense that increasing the capacity of the only congested link on the network can make the problem worse, and that reducing the capacity on this link can improve system conditions!
In this sense, it is a dynamic equivalent of the Braess paradox from Section~\ref{sec:motivatingexamples}.

Many have criticized the Braess paradox on the grounds that the link performance functions used in static assignment are unrealistic.
In the example shown below, queue spillback (a feature unique to dynamic network loading) is actually the critical factor in the paradox.

\stevefig{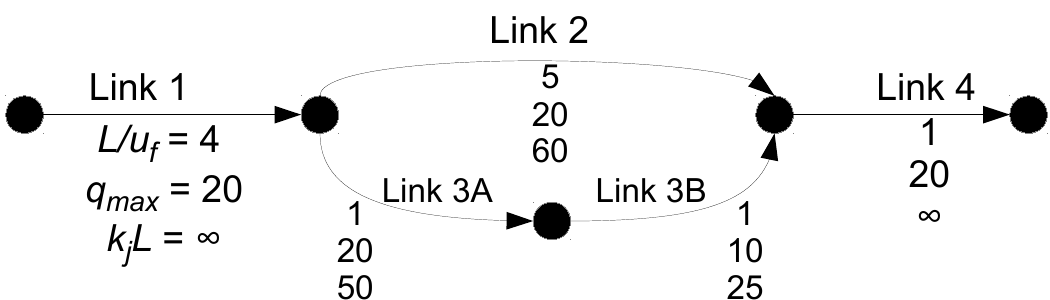}{Network for Daganzo's paradox.}{\textwidth}

See the network in Figure~\ref{fig:daganzonet}, where time is measured in units of time steps $\Delta t$.
Like the networks in the two previous sections, it consists of a single merge and diverge.
However, the free-flow travel times and capacities on the top and bottom links are now different: the top route is longer, but has a higher capacity, while the bottom route is shorter at free-flow, but has a bottleneck limiting the throughput on this route --- link 3B has only half the capacity of 3A.
We will use the spatial queue model of Section~\ref{sec:spatialqueue} to propagate traffic flow, although the same results would be obtained with an LWR-based model or anything else which captures queue spillback and delay.
The input demand is constant, at 20 vehicles per time step.

The capacity on the top route is high enough to accommodate all of the demand; if all of this demand were to be assigned to this route, the travel time would be 10 minutes per vehicle.
Assigning all vehicles to the top route is neither the user equilibrium nor the system optimum solution, but it does give an upper bound on the average vehicle travel time in the system optimal assignment --- it is possible to do better than this if we assign some vehicles to the bottom, shorter route.

To derive the user equilibrium solution, notice that initially all vehicles will choose the bottom path, since it has the lower travel time at free flow.
A queue will start forming on link 3A, since the output capacity is only 10 vehicles per time step (because of the series node model from Section~\ref{sec:linksinseries}, and the capacity of link 3B) and vehicles are entering at double this rate.
As the queue grows, the travel time on link 3A will increase as well.
With the spatial queue model and these inflow and outflow rates, you can show that the travel time for the $n$-th vehicle entering the link is
\labeleqn{daganzodelay}{1 + \frac{n}{q_{in}} \myp{\frac{q_{in}}{q_{out}} - 1} = 1 + \frac{n}{20}}
as long as the queue has not spilled back (see Exercise~\ref{chp:networkloading}.\ref{ex:pqtime}).

Based on equation~\eqn{daganzodelay}, when the 60th vehicle enters the link, it will spend four time units on link 3A, and thus its travel time across the entire bottom path would be the same as if it had chosen the top path.
At this point there are 40 vehicles on link 3A, as can be seen by drawing the $N^\uparrow$ and $N^\downarrow$ curves for this link.
From this point on, vehicles will split between the top and bottom paths to maintain equal travel times.

So far, so good; the first 60 vehicles that enter the network have a travel time of less than 10 minutes, and all the rest have a travel time of exactly 10 minutes.
Now see what happens if we increase the capacity on link 3B, with the intent of alleviating congestion by improving the bottleneck capacity.
If the capacity on 3B increases from 10 to 12, the story stays the same, except that it is the first 90 vehicles that have a travel time of less than 10 minutes.
We can see this by setting ~\eqn{daganzodelay} equal to four (the time needed to equalize travel times on the top and bottom paths), but with $q_{in} = 12$ instead of 10, and solving for $n$.
Network conditions indeed have improved.
But if the capacity increases still further, to 15, then equation~\eqn{daganzodelay} tells us that it is only after 180 vehicles have entered the system that travelers would start splitting between the top and bottom links.
By tracing the $N^\uparrow$ and $N^\downarrow$ curves, at this point in time 120 vehicles will have exited link 3A, meaning the queue length would be 60 vehicles.
\textbf{But the jam density of the bottom link only allows it to hold 50 vehicles.}  
The queue will thus spill upstream of the diverge node, and in this scenario, \emph{no} vehicles will opt to take the top path.
By the time a driver reaches the diverge point, the number of vehicles on the bottom link is 50, giving a travel time on 3A of $3\frac{1}{3}$ minutes, and a total travel time of $9 \frac{1}{3}$ minutes from origin to destination.
This is less than the travel time on the top path, so drivers prefer waiting in the queue to  taking the bypass route.

As a result, \emph{all} drivers will choose the bottom path, and the queue on the origin centroid connector will grow without bound, as will the travel times experienced by vehicles entering the network later and later.
By increasing the length of time vehicles enter the network at this rate, we can make the delays as large as we like.

We thus see that even with dynamic network loading, increasing the capacity on the only bottleneck link can make average travel times worse --- and in fact arbitrarily worse than the system optimum solution, where the average travel time is less than 10 minutes.
The reason for this phenomenon is the interaction between queue spillback and selfish route choice: in the latter scenario it would be better for the system for some drivers to choose the top route, even though it would worsen their individual travel time.

Of course, if the capacity on 3B were increased even further, all the way to 20 vehicles per time step, delays would drop again since there would be no bottleneck at all.
Exercise~\ref{ex:daganzoex} asks you to plot the average vehicle delay in this network as the capacity on link 3B varies from 0 to 25 vehicles per time step.
\index{dynamic traffic assignment!equilibrium!inefficiency|)}
\index{queue spillback|)}
\index{node model!diverge|)}

\subsection{Implications}

The purpose of these examples is to show that dynamic user equilibrium is complex.
Guaranteeing existence or uniqueness of dynamic equilibrium requires making strong assumptions on traffic flow propagation.
However, for some practical applications, using a more realistic traffic flow model is more important than mathematical properties of the resulting model.
Many people find comfort in the fact that we can at best solve for an equilibrium approximately, and thus dismiss the question of whether an equilibrium ``truly'' exists as akin to asking how many angels can dance on the head of a pin.

We emphasize that existence and uniqueness are not simply mathematical abstractions, and that they have significant implications for practice: if an equilibrium does not exist, should we really be ranking projects based on equilibrium analysis?
If multiple equilibria exist, what should we plan for?
Can we even find them all?
At the same time, we acknowledge that using static equilibrium to sidestep these difficulties is often unacceptable.
For many applications, the assumptions in link performance functions are simply too unrealistic.
Such is the nature of mathematical modeling in engineering practice, a topic taken up more fully in the next section.
\index{dynamic traffic assignment!equilibrium!properties|)}

\section{Practical Considerations}

We conclude this chapter with a discussion of what is required to implement a dynamic traffic assignment model in practice.
This discussion will take as implicit that you have decided that dynamic traffic assignment is the best tool for the problem at hand.
Recall the discussion in Sections~\ref{sec:introductionstatic} and~\ref{sec:introductiondynamic} about the relative strengths and weaknesses of static and dynamic traffic assignment.
As a general guideline, dynamic traffic assignment is best used when you need detailed representation of queues and other traffic phenomena, when your metrics explicitly track changes within the peak period, and when you have high confidence in the input data used to prepare the model.
If these are not true, you may want to consider a static traffic assignment (keeping in mind its own limitations).
In some applications you may want to treat route choice as an exogenous, given quantity rather than something to be determined by the model itself.
In such applications, a microsimulation model may be more appropriate.
Once you have decided that dynamic traffic assignment is the best tool for a specific problem, the main things needed to prepare a model are (1) the network topology, (2) ``supply-side'' data describing the physical infrastructure, and (3) ``demand-side'' data describing trip origins, destinations, and departure times.

Regarding network topology, we typically do not include literally every roadway segment in a region in our models.
The significant majority of roadway segments correspond to local roads or neighborhood streets with very little traffic, and including these in the model increases computational demands without providing much additional insight.
Commonly, dynamic traffic assignment networks include freeways, highways, and arterials (both major and minor), but the exact spatial resolution can vary by application.
The local roads and neighborhood streets which are not included in the model are represented instead by centroid connectors that load trips on and off of the network.
Centroid connectors play a much more significant role in affecting model results than is commonly understood\index{link!centroid connector!practical impact} --- especially for dynamic traffic assignment~\citep{james16_thesis, james16_trb}.
They should not be generated randomly, or by a rote rule such as creating three connectors for every centroid.
Rather, it matters to create connectors in a way that represents as best as possible the options available for travelers wanting to enter or leave a specific neighborhood.
Guidance on this topic is still emerging, and a topic of active research.

Regarding supply-side data, you need to decide which link and node models you want to use, and obtain the parameters you need to implement them.
The basic node and link models introduced in Chapter~\ref{chp:networkloading} are intentionally designed to make use of relatively easy-to-obtain data: link capacity\index{flow!capacity} (easily estimated as the product of the number of lanes and an assumed per-lane capacity\footnote{To one significant digit, per-lane capacity is approximately 2000 vehicles per hour.}, or with a professional reference such as the Highway Capacity Manual, or with field data) and free-flow time\index{speed!free-flow speed} (estimated as the link length divided by the speed limit, or observed free-flow speed data) are used in all link models.
Jam density\index{density!jam} can be estimated by assuming a spatial headway (front bumper to front bumper) between vehicles when stopped in queue, and multiplying the resulting density by the number of lanes.
With capacity, free-flow speed, and jam density, you can also implement any of the LWR-based link models based on a triangular fundamental diagram.
Using a trapezoidal fundamental diagram requires estimating a fourth parameter, typically the backward wave speed.
It is less clear how to do this without field data; estimates in the range of $1/6$ to $1/2$ the free-flow speed are common.

The fancier node models can require additional parameters.
In particular, modeling signalized intersections requires information on the signal timing\index{traffic signal timing} plan, to estimate the green time available to each turning movement.
Signalized intersections pose a particular challenge when used in future-year scenarios, since timing plans are periodically updated based on observed flows.
Unsignalized intersections may also be converted to signalized intersections if demand increases.
Unfortunately, the network loading is sensitive to having the correct node model in place, which is one of the reasons dynamic traffic assignment models must be used cautiously for long-range forecasts when the infrastructure can change significantly from present conditions.
One way to approach this issue is to update the node model between selected iterations of the dynamic user equilibrium process (in step 4 of the general equilibrium framework of Section~\ref{sec:duegeneral}).
Signal timings can be adjusted based on the link flows from the most recent network loading, following traffic engineering practice; a simple heuristic is to assume that phasing remains the same, but adjust green times to be proportional to critical lane group volume.
A ``warrant analysis'' can also be used to identify unsignalized intersections may be converted to signalized ones based on predicted volumes.
\cite{mannering20}, \cite{garber14}, and other transportation engineering texts provide information on these procedures.
While imperfect, these processes at least provide predictions of how traffic control schemes may change as link flows change in the future.

Regarding demand-side data, you need to produce a time-dependent OD matrix indicating the number of trips leaving each origin, for each destination, at each departure time.
Directly estimating a time-dependent OD matrix from observed link counts is very difficult: a time-dependent OD matrix has much higher dimensionality than the number of links in the network.
One approach is to estimate a \emph{static} OD matrix, as described in Sections~\ref{sec:tapdata} or~\ref{sec:odme}--\ref{sec:sensitivity_litreview}, and then ``profile'' this demand by assuming how it is distributed over time.
Demand profiles can be estimated from survey data or travel diaries, although they likely vary from one OD pair to another.
Another approach is to determine the demand profile endogenously, using some of the departure time choice models described in Section~\ref{chp:tdsp} and Section~\ref{sec:dueprinciple}.
While the algorithms described in this chapter are not guaranteed to converge to departure time equilibrium --- indeed, some researchers have even shown that departure time equilibrium in a single bottleneck is unstable for some common iterative processes \citep{iryo19} --- in large networks including departure time choice as part of the iteration in dynamic traffic assignment does not seem to cause problems \citep{levin15departuretime}.
Activity-based models, or emerging data sources, are newer techniques that have the potential to greatly improve estimation of time-dependent OD matrices by directly providing the information needed.

After a model is constructed, a calibration process is usually needed before the model outputs appear reasonable.
Dynamic traffic assignment models can be unforgiving of errors in the input data; omitting a zero from a link capacity, for instance, may cause gridlock due to queue spillback.
For this reason, point queue models are more tolerant of errors in the input data \citep{boyles19} precisely because they cannot capture queue spillback; this is a curse with accurate inputs, but a blessing with inaccurate ones.
Dynamic traffic assignment produces a large amount of data as output, including time-dependent link flows and densities, and time-dependent travel times between any nodes in the network.
This data can be compared to field data for calibration.
Given the quantity of data produced by the models, effective visualization tools can be very helpful in identifying potential problems with the input data.

The advice given in this section is necessarily generic, and individual projects may have particular considerations that warrant breaking the rules of thumb we give here.
We conclude by remarking once again that dynamic traffic assignment models are much more realistic than static models in terms of representing traffic flow and congestion, but there are tradeoffs in terms of computation, stability, and tolerance of errors in inputs.
No tool is right for every task. 
Rather, experienced practitioners know how to match the available tools to the job at hand.
By reading this book, we hope that you have gained the insight to understand the advantages and disadvantages of both static and dynamic traffic assignment models, and to make educated decisions about the right tool for a particular project.
Finally, research on dynamic traffic assignment is continuing to progress, and perhaps some of the challenges described above can be resolved with further attention and thought.
As researchers, we would be delighted for you to contribute to work in this field.

\section{Historical Notes and Further Reading}

Formulating dynamic traffic assignment models is more difficult than doing so for static traffic assignment.
For this reason, this chapter has eschewed highly detailed mathematical formulations of the equilibrium principle in favor of simply expressing the dynamic user equilibrium principle and proposing heuristics.
Readers looking for a more careful mathematical formulation of dynamic user equilibrium are referred to \cite{friesz23} and the references therein.

It is particularly difficult to formulate this problem if we seek a user equilibrium, rather than a system optimum, and if if there are multiple destinations.
For this reason, the first dynamic traffic assignment models~\citep{merchant78model,merchant78optimality} were restricted to a single destination and a system optimum was sought.
These authors used an ``exit function'' link model which is no longer common; using more recent link models, it is possible to find the system-optimum solution with a single destination by solving a linear program~\citep{ziliaskopoulos00}.
When there are multiple destinations, ensuring that flow to one destination does not overtake that headed to another imposes nonconvexity on the feasible flows~\citep{carey92}.

For this reason, variational inequality formulations are more common than mathematical optimization.
Examples of these include \cite{friesz93}, \cite{wie95}, and \cite{chen98}.
Optimal control\index{dynamic traffic assignment!equilibrium!as optimal control} approaches have also been proposed~\citep{friesz89,ran93}, as have formulations as a nonlinear complementarity problem~\citep{ban12}.\index{dynamic traffic assignment!equilibrium!as nonlinear complementarity problem}
Fixed point approaches are also common --- \cite{bargera_dta06} and \cite{bellei05} are just two examples --- but as with the static assignment problem, are more useful for specifying the problem then for solving it.

Using link or path flows as the main decision variable is the most intuitive choice, and therefore the most common in the literature.\index{dynamic traffic assignment!flow representation}
However, the use of splitting proportions is becoming more common~\citep{szeto12,nezamboylesdta,gentile16}.
At the end of Section~\ref{sec:flowtracking}, we remarked that the formulas~\eqn{downstreamdisaggregate} and~\eqn{upstreamdisaggregate} for disaggregating flow among paths was an approximation, because of the assumption that vehicles on different paths are uniformly distributed within the sending flow.
This issue is described briefly in \cite{daganzo95}, and in more detail in \cite{yperman_diss} and \cite{blumberg09}.
\cite{bargera22} present a method for disaggregating sending flows in a way which exactly respects the first-in, first-out principle.
For more on the theoretical basis for working with disaggregate $N^\pi$ values in the LWR model, see \cite{zhang02} and \cite{lebacque05}.

For the convex combinations and simplicial decomposition algorithms, refer to the references in Section~\ref{sec:algosreference}.
For gradient projection as specialized to dynamic network loading, see \cite{nezamboylesdta} and \cite{gentile16}.

The equilibrium existence counterexample in Section~\ref{sec:eggnetwork} is adapted from \cite{eggnetwork}; see \cite{iryo10} for additional examples.
The uniqueness counterexample in Section~\ref{sec:niemerge} is from \cite{nie10eqa}, and its special case in Section~\ref{sec:bdmr} is from \cite{boyles_bdmr}.
The efficiency counterexample in Section~\ref{sec:daganzoparadox} is from \cite{daganzo98}.
One consequence of these counterexamples is that queue spillback significantly complicates the finding and interpretation of dynamic user equilibria; see also the discussion in \cite{boyles19}.
\index{dynamic traffic assignment!equilibrium|)}
\index{dynamic traffic assignment|)}

\section{Exercises}
\label{sec:dynamiceqm_exercises}

\begin{enumerate}
\item \diff{37} \emph{(Equivalence of link-based and path-based flow representations)}.
Given splitting proportions $\alpha_{hij,s}^t$ for each destination $s$, time interval $t$, and turning movement $[h,i,j]$, show how ``equivalent'' path flow values $h_t^\pi$ can be found.
Then, if given path flows $h_t^\pi$, show how ``equivalent'' $\alpha_{hij,s}^t$ values can be found.
(``Equivalent'' means that the link cumulative counts $N^\uparrow$ and $N^\downarrow$ would be the same for all time steps after performing network loading, possibly with a small error due to time discretization that would shrink to zero as $\Delta t \rightarrow 0$.) \label{ex:pathlinkequivalence}
\item \diff{37} Consider the four-link network shown in Figure~\ref{fig:4cell}, and perform network loading using the cell transmission model.
(Each link is one cell long.) During the first time interval, 10 vehicles enter Link 1 on the top path; during the second time interval, 5 vehicles enter Link 1 on the top path and 5 on the bottom path; and during the third time interval, 10 vehicles enter Link 1 on the bottom path.
No other vehicles enter the network.

Interpolating as necessary, what time do the first and last vehicles on the top path exit cell 4?  the first and last vehicles on the bottom path?  What is the derivative of the travel time on the top path for a vehicle leaving at the start of interval 2? \label{ex:4cell}
\stevefig{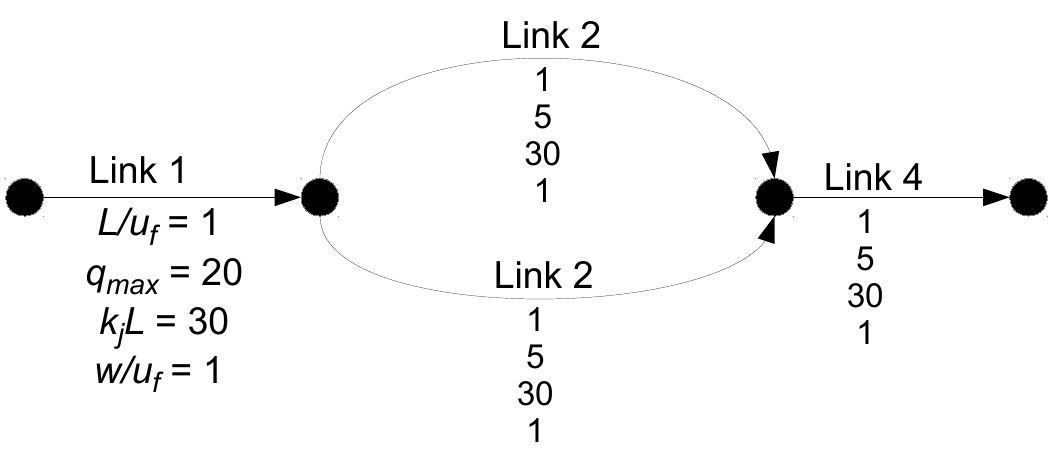}{Network for Exercise~\ref{ex:4cell}.}{\textwidth}
\item \diff{12} In the course of the convex combinations algorithm, assume that the $H$ and $H^*$ matrices are as below, and $\lambda = 1/3$.
What is the new $H$ matrix?
\labeleqn{convexex}{H = \vect{ 6 & 12 \\ 30 & 24} \qquad H^* = \vect{ 18 & 0 \\ 0 & 54}}
\item \diff{23} Consider a network with only one origin-destination pair connected by four paths, with three departure time intervals.
Departure times are fixed.
At some point in the simplicial decomposition algorithm, $\mathcal{H}^*$ contains the following three matrices:
\labeleqn{simplicialmatrices}{
\begin{bmatrix} 20 & 0 & 0 & 0 \\ 0 & 10 & 0 & 0\\ 0 & 0 & 30 & 0 \end{bmatrix} \begin{bmatrix} 0 & 20 & 0 & 0 \\ 0 & 10 & 0 & 0\\ 0 & 30 & 0 & 0 \end{bmatrix} \begin{bmatrix} 20 & 0 & 0 & 0 \\ 0 & 0 & 0 & 10\\ 0 & 0 & 30 & 0 \end{bmatrix}
}
and the current path flow and travel time matrices are:
\labeleqn{simplicialht}{
H = \begin{bmatrix} 14 & 6 & 0 & 0 \\ 0 & 8 & 0 & 2\\ 0 & 9 & 21 & 0 \end{bmatrix} 
\qquad 
T(H) = \begin{bmatrix} 20 & 20 & 24 & 27 \\ 30 & 34 & 37 & 40 \\ 44 & 35 & 36 & 40 \end{bmatrix} 
}
What are the unrestricted and restricted average excess costs of the current solution?  What is the matrix $\Delta H$?

\item \diff{35} Assume that there is a single OD pair, two paths, and three departure times; 15 vehicles depart during the first, 10 during the second, and 5 during the third.
Let $H_{\pi t}$ denote the number of vehicles departing on path $\pi$ at time $t$, and that the path travel times are related to the path flows by the following equations:\footnote{In practice these would come from performing network loading and calculating path travel times as described in this chapter; these functions are provided here for your convenience.} 
\begin{align*}
T_{11} &= 15  \\
T_{21} &= 13 + 3H_{11} + H_{21} \\
T_{12} &= 15 + H_{11}^2 \\
T_{22} &= 13 + H_{11} + 3H_{12} + H_{21}^2 + H_{22} \\
T_{13} &= 15 + H_{12}^2 + \frac{1}{2} H_{11}^2 \\
T_{23} &= 13 + H_{12} + 3H_{13} + H_{22}^2 + H_{23} \\
\end{align*}
Find the path flows obtained from three iterations of the convex combinations method with step sizes $\lambda_1 = 1/2$, $\lambda_2 = 1/4$, and $\lambda_3 = 1/6$ (so you should find a total of four $H^*$ matrices, counting the initial matrix, and take three weighted averages).
Your initial matrix should be the shortest paths with zero flow.
What is the resulting AEC?  \emph{Break any ties in favor of path 1.}
\label{ex:convexcombinations}
\item \diff{35} Using the same network and demand as in Exercise~\ref{ex:convexcombinations}, now assume that you are solving the same problem with simplicial decomposition, and at some stage $\mathcal{H}^*$ contains the following two matrices:
$$ \vect{15 & 0 \\ 0 & 10 \\ 5 & 0} \qquad \vect{15 & 0 \\ 10 & 0 \\ 0 & 5} $$
and that the current solution is 
$$ H = \vect{15 & 0 \\ 4 & 6 \\ 3 & 2} $$
\begin{enumerate}[(a)]
\item What is the average excess cost of the current solution?
\item What is the restricted average excess cost?
\item What is the search direction from $H$?
\item Give the updated $H$ matrix and new restricted average excess cost after taking a step with $\mu = 0.1$.
\item If we terminate the subproblem and return to the master algorithm, what matrix (if any) do we add to $\mathcal{H}^*$?
\end{enumerate}
\item \diff{35} Again using the network and demand from Exercise~\ref{ex:convexcombinations}, now start with the initial solution 
$$ H = \vect{10 & 5 \\ 5 & 5 \\ 0 & 5} $$
Using the gradient projection method, identify new $H$ and $C$ matrices.
Try both exact and quasi-Newton steps, and see which gives the greater reduction in $AEC$.
\item \diff{33} Modify the three solution algorithms in this chapter (convex combinations, simplicial decomposition, and gradient projection) to handle the case where departure times are not fixed.
\item \diff{46} Find the mixed-strategy equilibrium in the network in Figure~\ref{fig:eggnetwork}.
(For each of the two vehicles, indicate the probability it will choose each of the two paths available to it; these probabilities should be such that the \emph{expected} travel times are the same for both options.) \label{ex:eggmixed}
\item \diff{35} Consider a network with one origin, one destination, two possible departure time intervals (1 and 2), and two paths (A and B).
A total of 16 vehicles depart during time interval 1, and 4 depart during time interval 2.
Let $t^\pi_\tau$ refer to the travel time on the $\pi$-th path when departing at the $\tau$-th time interval, with $h^\pi_\tau$ defined similarly for path flows.
Suppose the travel times are related to the path flows as follows:
\begin{align*}
t^A_1 &= 0.5 h^A_1 + 5 h^A_2 + 6 \\
t^B_1 &= 0.5 h^B_1 + 3 h^B_2 + 10 \\
t^A_2 &= 0.3 h^A_1 + 0.6 h^A_2 + 0.8 \\
t^B_2 &= 0.2 h^B_1 + 0.4 h^B_2 + 2 \\
\end{align*}
and consider the following four path flow matrices:
$$ H_1 = \vect{0 & 16 \\ 4 & 0} \qquad H_2 = \vect{12 & 4 \\ 2 & 2} \qquad H_3 = \vect{8 & 8 \\ 2 & 2} \qquad H_4 = \vect{16 & 0 \\ 0 & 4}$$
(M. Netter, 1972)\footnote{If you are wondering how the travel time on a path can be influenced by drivers leaving at a later time, recall that FIFO violations can occur due to phenomena such as express lanes opening, allowing ``later'' vehicles to overtake ``earlier'' ones and delay them in the network.}
\begin{enumerate}[(a)]
 \item Which of the above path flow matrices satisfy the principle of user equilibrium?
 \item Which of the above path flow matrices represent \emph{efficient} equilibria?
 \item Perturb $H_3$ in the following way: adjust $h^A_2$ from 2 to 2.1 (so $h^B_2$ becomes 1.9), changing all of the travel times.
 Then, adjust $h^A_1$ and $h^B_1$ until $t^A_1 = t^B_1$ (restoring equilibrium for the first departure time).
 Then, adjust $h^A_2$ and $h^B_2$ until $t^A_2 = t^B_2$ (restoring equilibrium for the second departure time), and so on, alternating between the two departure times, until you reach a new equilibrium solution.
 Which equilibrium solution do you end up at?  Is $H_3$ stable?
\end{enumerate}
\item \diff{33} Consider an instance of the Daganzo paradox network (Figure~\ref{fig:daganzonet}), but with travel times $\tau_2 = 10$, $\tau_{3A} = 4$ and $\tau_{3B} = 1$, capacities of $c_{3B} = 10$ and $c_1 = c_2 = c_{3A} = c_4 = 50$ vehicles per time step, and an inflow rate of $Q = 50$ vehicles per time step.
Assume that the maximum number of vehicles which can fit on link 3A is $(k_j L)_{3A} = 300$ vehicles.
\begin{enumerate}[(a)]
 \item How many vehicles should be on links 3A at any time to maintain equilibrium?
 \item Does this queue fit on the link?
 \item The capacity on 3B is now increased to $c_{3B} = 25$.
 At the new equilibrium, have drivers' travel times increased, decreased, or stayed the same?
 \item The capacity is now increased to $c_{3B} = 50$.
 At the new equilibrium, have drivers' travel times increased, decreased, or stayed the same (relative to $c_{3B} = 25$)?
\end{enumerate}
\item \diff{55} In the network of Figure~\ref{fig:daganzonet}, find the equilibrium solution in terms of the capacity on link 3B.
Plot the average travel time as this capacity ranges between 0 and 25 vehicles per time step.
Assume a time horizon of 20 time steps.
(You may do this either by deriving a closed-form expression for the average travel time, or by using one of the algorithms in this chapter to solve for the approximate equilibrium solution.) \label{ex:daganzoex}
\end{enumerate}

\label{sec:due_exercises}

\appendix

\clearpage
\addcontentsline{toc}{part}{Appendices}

\chapter{Mathematical Concepts}
\label{chp:mathbackground}

This appendix reviews the mathematical background needed to understand the mathematical formulations in this book.
This involves understanding certain properties of matrices, sets, and functions.
The descriptions are relatively terse, because there are countless textbooks and other resources with additional explanations and examples.
It is also a very incomplete overview, only covering concepts which will be used later in the book and largely assuming that the reader has seen this material at some point in the past.
The section also collects a number of useful results related to matrices, sets, and functions.
The proofs of most of these are left as exercises.

\section{Indices and Summation}
\label{sec:indicessummation}

\index{subscripts and superscripts}
\index{index notation|(}
We commonly work with attributes of links or nodes, such as the total number of trips starting or ending at a node, or the travel time on a link.
In such cases, it makes sense to use the same variable letter for the number of trips produced at a node (say $P$), and to indicate the productions at a specific node with a subscript or superscript ($P_1$ or $P^3$.)  Unfortunately, the superscript notation can be ambiguous, since the superscript in $P^3$ can be interpreted both as an exponent, and as the index referring to a particular node.
In most cases it will be clear which meaning is intended, and parentheses can be used in cases which are unclear: $(P^3)^2$ is the square of the productions at node 3.
To avoid confusion, subscripts are generally preferred to superscripts, but superscripts can make the notation more compact when there are multiple indices.
An example is $\delta_{ij}^\pi$ from Section~\ref{sec:introductionstatic}, where the variable $\delta$ is indexed both by a link $(i,j)$ and a path $\pi$.
We have striven to be consistent with which indices appear as subscripts and which appear as superscripts.
However, in a few occasions, being absolutely rigid in this regard would create cluttered formulas.
In such cases, we have sacrificed rigor for readability, and hope that you are not thrown off by an index typically appearing as a subscript in the superscript position, and vice versa.

Subscript or superscript indices also allow formulas to be written more compactly.
A common example is the summation notation,\index{summation notation|(} such as
\labeleqn{sumexample}{\sum_{i=1}^5 P_i = P_1 + P_2 + P_3 + P_4 + P_5\,.}
The left-hand side of Equation~\eqn{sumexample} is used as shorthand for the right-hand side.
More formally, the left-hand side instructs us to choose all the values of $i$ between 1 and 5 (inclusive); collect the terms $P_i$ for all of these values of $i$; and to add them together.
A variant of this notation is
\labeleqn{sumexample2}{\sum_{i \in N} P_i\,,}
which expresses the sum of the productions from all nodes in the set $N$.
Here $i$ ranges over all elements in the set $N$, rather than between two numbers as in~\eqn{sumexample}.
We can also add conditions to this range by using a colon.
If we only wanted to sum over nodes whose products are less than, say, 500, we can write
\labeleqn{sumexample4}{\sum_{i \in N : P_i < 500} P_i\,.}
When it is exceptionally clear what values or set $i$ is ranging over, we can simply write
\labeleqn{sumexample3}{\sum_i P_i\,,}
but this ``abbreviated'' form should be used sparingly, and avoided if there is any ambiguity at all as to what values $i$ should take in the sum.

When there are multiple indices, a double summation can be used:
\labeleqn{doublesumexample}{\sum_{i = 1}^3 \sum_{j = 2}^3 x_{ij} = x_{12} + x_{22} + x_{32} + x_{13} + x_{23} + x_{33}\,.}
You can think of a double sum either as a ``nested'' sum:
\labeleqn{doublesumnest}{\sum_{i=1}^3 \sum_{j=2}^3 x_{ij} = \sum_{i=1}^3 \myp{\sum_{j=2}^3 x_{ij}} = \sum_{i=1}^3 \myp{x_{i2} + x_{i3}}\,,}
which expands to the same thing as the right-hand side of~\eqn{doublesumexample}, or as summing over all combinations of $i$ and $j$ such that $i$ is between 1 and 3, and $j$ is between 2 and 3.
Triple sums, quadruple sums, and so forth behave in exactly the same way, and are common when there are many indices.

The summation index is often called a \emph{dummy variable},\index{dummy variable} because the indices do not have an absolute meaning.
Rather, they are only important insofar as they point to the correct numbers to add up.
For instance, $\sum_{i=1}^5 P_i$ and $\sum_{i=0}^4 P_{i+1}$ are exactly the same, because both expressions have you add up $P_1$ through $P_5$.
Likewise, $\sum_{i=1}^5 P_i$ and $\sum_{j=1}^5 P_j$ are exactly the same, the fact that we are counting from 1 to 5 using the variable $j$ instead of $i$ is of no consequence.

Related to this, it is \emph{wrong} to refer to a summation index outside of the sum itself.
A formula such as $x_i + \sum_{i=1}^5 y_i$ is incorrect.
For the formula to make sense, $x_i$ needs to refer to one specific value of $i$.
But using $i$ as the index in the sum $\sum_{i=1}^5$ means that $i$ must range over all the values between 1 and 5.
Does $y_i$ refer to the index of summation (which ranges from 1 to 5), or to the one specific value of $i$ used outside the sum?  If you want to refer to a specific node, as well as to index over all nodes, you can add a prime to one of them, as in
\labeleqn{correctsumprime}{x_i + \sum_{i' = 1}^5 y_{i'}\,,}
or you can either use a different letter altogether, as in
\labeleqn{correctsumdifferentletter}{x_i + \sum_{j=1}^5 y_j\,.}
Both conventions are common.

In transportation network analysis, we frequently have to sum over all of the links which are in the forward or reverse star\index{reverse star} of a node (Section~\ref{sec:datastructures}).
If the link flows are denoted by a variable $x$, the total flow entering node $i$ is the sum of the link flows in its reverse star.
This can be written in several ways:
\labeleqn{reversestarsum}{\sum_{(h,i) \in \Gamma^{-1}(i)} x_{hi} \qquad  \qquad \sum_{(h,i) \in A} x_{hi}}
The notation on the right can be a bit confusing at first glance, since it looks like we are summing the flow of \emph{every} link (the whole set $A$), rather than just the links entering node $i$.
The critical point is that in this formula, $i$ refers to one specific node which was previously chosen and defined outside of the summation.
\emph{The only variable we are summing over in Equation~\eqn{reversestarsum} is $h$.}  In this light, $i$ is fixed, and the right-most sum is over all the values of $h$ such that $(h,i)$ is a valid link (that is, $(h,i) \in A$).
These are exactly the links which form the reverse star of $i$.

\index{summation notation!properties|(}
Almost all of the sums we will see in this book involve only a finite number of terms.
These sums are much easier to work with than infinite sums, and have the useful properties listed below.
(Some of these properties do not apply to sums involving infinitely many terms.)
\begin{enumerate}
\item You can factor constants out of a sum:
\labeleqn{sumdistributive}{\sum_i cx_i = c \sum_i x_i\,,}
no matter what values $i$ ranges over.
This follows from the distributive property for sums: $a(b + c) = ab + ac$.
(A ``constant'' here is any term which does not depend on $i$.) 
\item If what you are summing is itself a sum, you can split them up:
\labeleqn{sumsplit}{\sum_i (x_i + y_i) = \sum_i x_i + \sum_i y_i\,.}  This follows from the commutative property for sums: $a + b = b + a$, so you can rearrange the order in which you add up the terms.
\item You can exchange the order of a double summation:
\labeleqn{doublesumswap}{\sum_i \sum_j x_{ij} = \sum_j \sum_i x_{ij}}
no matter what values $i$ and $j$ range over.
This also follows from the commutative property.
\end{enumerate}
None of these properties are anything new; they simply formalize how we can operate with the $\sum$ notation using the basic properties of addition.
These properties can also be combined.
It is common to exchange the order of a double sum, and then factor out a term, to perform a simplification:
\labeleqn{simplifyswapsum}{\sum_i \sum_j c_j x_{ij} = \sum_j \sum_i c_j x_{ij} = \sum_j \left( c_j \sum_i x_{ij} \right) \,.}
The last step is permissible because $c_j$ is a ``constant'' relative to a sum over $i$.
We could not do this step at the beginning, because $c_j$ is not a constant relative to a sum over $j$.
This kind of manipulation is helpful if, say, $\sum_i x_{ij}$ has a convenient form or a previously-calculated value.
\index{summation notation!properties|)}
\index{index notation|)}
\index{summation notation|)}

\index{product notation}
Occasionally it is useful to refer to products over an index.
In analogy with the summation notation $\sum$ used for sums over an index, products over an index are written with the notation $\prod$.\label{not:prod}
As an example
\labeleqn{prodexample}{\prod_{i=1}^5 P_i = P_1 P_2 P_3 P_4 P_5}
Exercise~\ref{ex:productproperties} asks you to investigate which of the properties of summation notation carry over to the product notation.

\section{Vectors and Matrices}
\label{sec:matrices}

\index{scalar}
\index{vector|(}
\index{vector!dimension}
A \emph{scalar} is a single real number, such as 2, $-3$, $\sqrt{2}$, or $\pi$.
The notation $x \in \bbr$ means that $x$ is a real number (a scalar).
A \emph{vector} is a collection of scalars; the \emph{dimension} of a vector is the number of scalars it contains.
For instance, $\vect{3 & -5}$ is a two-dimensional vector, $\vect{6 & 0 & \pi}$ is a three-dimensional vector, and so forth.
It is possible to write vectors either horizontally, with its component scalars in a row, or vertically, with its component scalars in a column.
For the most part it doesn't matter whether vectors are written in a row or in a column; the exception is in formulas involving multiplication between vectors and matrices, as described below, where row and column vectors must be distinguished.
Vectors are usually denoted with boldfaced lower-case letters, like $\mb{x}$; the notation $\mb{x} \in \bbr^n$ indicates that $\mb{x}$ is an $n$-dimensional vector.\footnote{There are some exceptions; for instance, shortest-path labels are traditionally denoted with an upper-case $L$.}  (So, $\vect{3 & -5} \in \bbr^2$ and $\vect{6 & 0 & \pi} \in \bbr^3$.)  Individual components of vectors are often denoted with subscripts or superscripts, as in $x_1$ or $x^1$.
The \emph{zero vector}\index{vector!zero} is a vector with zeros for all of its components, and is denoted $\mb{0}$.

Two vectors of the same dimension can be added together by adding the corresponding components of each vector.
If $\mb{x} = \vect{1 & 2}$ and $\mb{y} = \vect{3 & 4}$, then $\mb{x + y} = \vect{4 & 6}$.
Multiplying a vector by a scalar means multiplying each component of the vector by the scalar, so $3\mb{x} = \vect{3 & 6}$ and $-\mb{y} = \vect{-3 & -4}$.
The \emph{dot product}\index{vector!dot product} of two vectors of the same dimension is defined as
\labeleqn{dotproduct}{\mb{x} \cdot \mb{y} = \sum_i x_i y_i \,,}
where the sum is taken over all vector components; with the example above, $\mb{x} \cdot \mb{y} = 1 \times 3 + 2 \times 4 = 11$.

The dot product of two vectors is a scalar.
The magnitude of a vector $\mb{x}$ is given by its $\emph{norm}$\index{vector!norm} $|\mb{x}| = \sqrt{\mb{x} \cdot \mb{x}}$.
This norm provides a measure of distance between two vectors; the distance between $\mb{x}$ and $\mb{y}$ is given by $|\mb{x} - \mb{y}|$.

The dot product can also be written
\labeleqn{dotproductangle}{\mb{x} \cdot \mb{y} = |\mb{x}| |\mb{y}| \cos \theta}
where $\theta$\label{not:thetaang} is the angle between the vectors $\mb{x}$ and $\mb{y}$ if they are both drawn from a common point.
In particular, $\mb{x}$ and $\mb{y}$ are perpendicular if $\mb{x} \cdot \mb{y} = 0$.

A collection of vectors $\mb{x_1}, \ldots, \mb{x_n}$ is \emph{linearly independent}\index{linear (in)dependence} if the only solution to the equation $a_1 \mb{x_1} + \cdots + a_n \mb{x_n} = 0$ is $a_1 = \cdots = a_n = 0$.
Otherwise, these vectors are \emph{linearly dependent}.
\index{vector|)}

\index{matrix|(}
A matrix is a rectangular array of scalars.
If a matrix has $m$ rows and $n$ columns, it is called an $m \times n$ matrix, and is an element of $\bbr^{m \times n}$.
A matrix is \emph{square}\index{matrix!square} if it has the same number of rows and columns.
In this book, matrices are denoted by boldface capital letters, such as $\mb{X}$ or $\mb{Y}$.
In the examples that follow, let $\mb{X}$, $\mb{Y}$, and $\mb{Z}$ be defined as follows:
\[ \mb{X} = \vect{1 & 2 \\ 3 & 4} \qquad \mb{Y} = \vect{0 & -1 \\ -3 & 5} \qquad \mb{Z} = \vect{-1 & 1 & -2 \\ 3 & -5 & 8} \,.\]
Elements of matrices are indexed by their row first and column second, so $x_{11} = 1$ and $x_{12} = 2$.

Addition and scalar multiplication of matrices works in the same way as with vectors, so
\[ \mb{X + Y} = \vect{1 + 0 & 2 - 1 \\ 3 - 3 & 4 + 5} = \vect{1 & 1 \\ 0 & 9} \]
and
\[ 5\mb{X} = \vect{5 & 10 \\ 15 & 20} \,.
\]
The \emph{transpose}\index{matrix!transpose} of a matrix $\mb{A}$, written $\mb{A}^T$, is obtained by interchanging the rows and columns, so that if $\mb{A}$ is an $m \times n$ matrix, $\mb{A}^T$ is an $n \times m$ matrix.
As examples we have
\[ \mb{X}^T = \vect{1 & 3 \\ 2 & 4} \qquad \mb{Z}^T = \vect{-1 & 3 \\ 1 & -5 \\ -2 & 8} \,.\]

\index{matrix!multiplication|(}
Matrix multiplication is somewhat less intuitive.
Two matrices can only be multiplied together if the number of columns in the first matrix is the same as the number of rows in the second.
(The reason for this will become clear when the operation is defined.)  For instance, you can multiply a $2 \times 3$ matrix by a $3 \times 3$ matrix, but you can't multiply a $2 \times 3$ matrix by another $2 \times 3$ matrix.
This immediately suggests that the order of matrix multiplication is important, since two matrices may have compatible dimensions in one order, but not in the other.
Multiplying an $m \times n$ matrix by a $n \times p$ matrix creates an $m \times p$ matrix, defined as follows.
Let $\mb{A} \in \bbr^{m \times n}$ and $\mb{B} \in \bbr^{n \times p}$.
Then the product $\mb{C} = \mb{AB} \in \bbr^{m \times p}$ has elements defined as
\labeleqn{matrixmult}{c_{ij} = \sum_{k=1}^n a_{ik} b_{kj}\,.}
If you imagine that each row of the first matrix is treated as a vector, and that each column of the second matrix is treated as a vector, then each component of the product matrix is the dot product of a row from the first matrix and a column from the second.
For this dot product to make sense, these row and column vectors have to have the same dimension, that is, the number of columns in the first matrix must equal the number of rows in the second.
Using the matrices defined above, we have
\[ \mb{XY} =  \vect{1 \times 0 + 2 \times -3 & 1 \times -1 + 2 \times 5 \\ 3 \times 0 + 4 \times -3 & 3 \times -1 + 4 \times  5} = \vect{-6 & 9 \\ -12 & 17} \,.\]
Observe that the element in the first row and first column of the product matrix is the dot product of the first row from $\mb{X}$ and the first column from $\mb{Y}$; the element in the first row and second column of the product is the dot product of the first row from $\mb{X}$ and the second column of $\mb{Y}$; and so forth.
You should be able to verify that
\[ \mb{XZ} = \vect{5 & -9 & 14 \\ 9 & -17 & 26 } \,.\]

A matrix and a vector can be multiplied together, if you interpret a row vector as an $1 \times n$ matrix or a column vector as a $m \times 1$ matrix.
A dot product of two vectors $\mb{x}$ and $\mb{y}$ of equal dimension can be written as a matrix multiplication: $\mb{x} \cdot \mb{y} = \mb{x}^T \mb{y}$ if they are both column vectors, or as $\mb{x} \mb{y}^T$ if they are both row vectors.
Here, the distinction between row and column vectors is important, because matrices can only be multiplied if their dimensions are compatible.
In matrix multiplications, the convention used in this text is that vectors such as $\mb{x}$  are column vectors, and a row vector is denoted by $\mb{x}^T$.
Again, this distinction is only relevant in matrix multiplication, and for other purposes row and column vectors can be treated interchangeably.

\emph{It is worth repeating that matrix multiplication is \emph{not} commutative, so that usually $\mb{XY} \neq \mb{YX}$ (and in fact both products may not even exist, depending on the dimensions of $\mb{X}$ and $\mb{Y}$)}, although there are some exceptions.
You may wonder why this definition of matrix multiplication is used instead of other, seemingly simpler, approaches.
One reason is that this definition actually ends up representing real-world calculations more frequently than other definitions.
For instance, it can be used to compactly write a set of equations, as is common in optimization problems (see Section~\ref{sec:optvector}).
\index{matrix!multiplication|)}

\index{matrix!symmetric}
\index{matrix!diagonal}
\index{matrix!identity}
A square matrix $\mb{A}$ is \emph{symmetric} if $a_{ij} = a_{ji}$ (in other words, for a symmetric matrix $\mb{A} = \mb{A}^T$), and it is \emph{diagonal} if $a_{ij} = 0$ unless $i = j$ (that is, all its elements are zero except on the diagonal from upper-left to lower-right).
A very special diagonal matrix is the \emph{identity matrix}, which has 1's along the diagonal and 0's everywhere else.
The notation $\mb{I}$ denotes an identify matrix of any size, so we can write
\[ \mb{I} = \vect{1 & 0 \\ 0 & 1} \qquad \mbox{or} \qquad \mb{I} = \vect{1 & 0 & 0 & 0 \\ 0 & 1 & 0 & 0 \\ 0 & 0 & 1 & 0 \\ 0 & 0 & 0 & 1}\,.
\]
In practice, the dimension is usually obvious because of the requirements of matrix multiplication.
The identity matrix has the unique property that $\mb{AI} = \mb{IA} = \mb{A}$ for any matrix $\mb{A}$.

\index{matrix!inverse|(}
A square matrix $\mb{A}$ is \emph{invertible} if there is another square matrix (call it $\mb{A^{-1}}$) such that $\mb{AA^{-1} = A^{-1}A = I}$.
Multiplying by an inverse matrix should be thought of as the equivalent of ``matrix division.''  Just as it is not possible to divide by all numbers (division by zero is undefined), not all matrices have inverses.
One can show that a square matrix is invertible if and only if the vectors forming its rows are linearly independent.
The matrices seen in this book will all be invertible.

Computing the matrix inverse is a bit tedious, and is rarely needed in transportation network analysis.
However one special case is worth mentioning: 
\begin{prp}
\label{prp:diagonalinverse}
A diagonal matrix $\mb{A}$ is invertible if and only if all its diagonal entries are nonzero; in this case, its inverse $\mb{A^{-1}}$ is also a diagonal matrix, whose diagonal entries are the reciprocal of the diagonal entries in $\mb{A}$.

\end{prp}
So, for instance, 
\[ \left( \vect{5 & 0 \\ 0 & -3} \right)^{-1} = \vect{1/5 & 0 \\ 0 & -1/3} \,.\]
\index{matrix!inverse|)}

\index{matrix!positive (semi)definite|(}
Finally, a symmetric matrix $\mb{A} \in \bbr^{n \times n}$ is \emph{positive definite} if the matrix product $\mb{x}^T\mb{Ax}$ is strictly positive for any nonzero vector $\mb{x} \in \bbr^n$, and \emph{positive semidefinite} if $\mb{x}^T\mb{Ax} \geq 0$ for any $\mb{x} \in \bbr^n$ whatsoever.
As examples, consider the matrices
\[ \mb{A} = \vect{10 & 1 \\ 1 & 5} \qquad \mb{B} = \vect{0 & 0 \\ 0 & 1} \qquad \mb{C} = \vect{1 & 10 \\ 10 & 5} \,. \]
The matrix $\mb{A}$ is positive definite, because for any vector $\mb{x} = \vect{x_1 & x_2}$, the matrix product is
\[ \mb{x}^T\mb{Ax} = \vect{x_1 & x_2} \vect{10 & 1 \\ 1 & 5} \vect{x_1 \\ x_2} = 10x_1^2 + 2x_1 x_2 + 5x_2^2 \,.\]
The expression on the right is always positive, because it can be rewritten as
\[ 9x_1^2 + 4x_2^2 + (x_1 + x_2)^2 \,.\]
None of those three terms can be negative, and since $\mb{x \neq 0}$, at least one of these terms is strictly positive.

The matrix $\mb{B}$ is not positive definite, since if $\mb{x} = \vect{1 & 0}$ then $\mb{x}^T\mb{B x} = 0$, which is not strictly positive.
However, it is positive semidefinite, since the matrix product is
\[ \mb{x}^T\mb{B x} = \vect{x_1 & x_2} \vect{0 & 0 \\ 0 & 1} \vect{x_1 \\ x_2} = x_2^2 \]
which is surely nonnegative.

The matrix $\mb{C}$ is neither positive definite nor positive semidefinite, since if $\mb{x} = \vect{1 & -1}$ then
\[ \mb{x}^T \mb{C x} = \vect{1 & -1} \vect{1 & 10 \\ 10 & 5} \vect{1 \\ -1} = -14 < 0 \,.\]
Notice that for a matrix to be positive definite or semidefinite, we have to check whether a condition holds for all possible nonzero vectors.
If the condition fails for even one vector, the matrix is not positive definite or semidefinite.

Checking positive definiteness or semidefiniteness can be tedious.
One test is to find all the eigenvalues of the matrix; if they are all positive, the matrix is positive definite, and if all are non-negative it is positive semidefinite.
However, diagonal matrices arise fairly often in transportation network analysis, and this case is easy.

\begin{prp}
\label{prp:diagonalpd}
A diagonal matrix is positive definite if and only if all its diagonal entries are strictly positive.
A diagonal matrix is positive semidefinite if and only if all its diagonal entries are nonnegative.
\end{prp}

There are also times where we will apply the concept of positive definiteness to non-symmetric matrices.
The idea is the same -- we want $\mb{x}^T\mb{Ax}$ to be strictly positive for all nonzero $\mb{x}$ -- but the eigenvalue test does not apply, and there are non-symmetric matrices which have strictly positive eigenvalues but do not satisfy $\mb{x}^T\mb{Ax} > 0$ for all nonzero $\mb{x}$.
However, we can form the \emph{symmetric part} of the matrix $\mb{A}$ by calculating $\frac{1}{2} (\mb{A} + \mb{A}^T)$.
It is easy to show that this is always a symmetric matrix, and that $\mb{x}^T\mb{Ax} > 0$ if and only if $\mb{x^T} \myp{ \frac{1}{2} (\mb{A} + \mb{A}^T)} \mb{x} > 0$.
So, when we refer to a non-symmetric matrix being positive definite, what we will mean is that its symmetric part $\frac{1}{2} (\mb{A} + \mb{A}^T)$ is positive definite.
\index{matrix!positive (semi)definite|)}

\index{matrix!determinant}
The \emph{determinant} of a square matrix is occasionally useful in transportation network analysis (but less so than in other applied mathematics fields).
For a $1 \times 1$ matrix, its determinant is simply the value of the single entry in the matrix.
For an $n \times n$ matrix, the determinant can be computed by the following procedure: select any row or column of the matrix; for each entry in this row or column, compute the determinant of the $(n - 1) \times (n - 1)$ matrix resulting from deleting both the row and column of this entry; and alternately add and subtract the resulting determinants.
For our purposes, the determinant can be used to concisely express other matrix properties. For instance, one can show that a square matrix is invertible if and only if its determinant is not zero, and that it is positive definite if and only if the determinants of all the square submatrices including the element in the first row and column are positive.
It also appears in the characterization of totally unimodular matrices, which are discussed in Section~\ref{sec:totallyunimodular} as a property of certain integer optimization problems which makes them easier to solve.
\index{matrix|)}

\section{Sets}
\label{sec:sets}

\index{set|(}
A \emph{set} is a collection of any type of objects, denoted by a plain capital letter, such as $X$ or $Y$.
In transportation network analysis, we work with sets of numbers, sets of nodes, sets of links, sets of paths, sets of origin-destination pairs, and so forth.
Sets can contain either a finite or infinite number of elements.
As examples, let's work with the sets
\[ X = \{ 1, 2, 4 \}  \qquad Y = \{1, 3, 5, 7 \} \,,\]
using curly braces to denote a set, and commas to list the elements.
Set membership is indicated with the notation $\in$, so $1 \in X$, $1 \in Y$, $2 \in X$, but $2 \notin Y$.
The \emph{union}\index{set!union} of two sets $X \cup Y$ is the set consisting of elements either in $X$ or in $Y$ (or both), so
\[ X \cup Y = \{ 1, 2, 3, 4, 5, 7 \}\,.
\]
The \emph{intersection}\index{set!intersection} of two sets $X \cap Y$\label{not:cap} is the set consisting only of elements in both $X$ and $Y$, so
\[ X \cap Y = \{ 1 \}\,.
\]
A set is a \emph{subset}\index{set!subset} of another set, denoted by $\subseteq$, if \emph{all} of its elements also belong to the second set.
With these sets $X \nsubseteq Y$.
Even though the element 1 is in both $X$ and $Y$, the elements 2 and 4 are in $X$ but not $Y$.
We do have $X \cap Y \subset X$ and $X \cap Y \subset Y$, however.  

In this book, sets will take one of three forms:
\begin{enumerate}
\item Sets that consist of a finite number of elements, which can be listed as with $X$ and $Y$ above.
This includes the set of nodes in a network, the set of acyclic paths in a network, and so forth.
For a finite set, the notation $|X|$ indicates the number of elements in $X$.
\item Sets which are intervals on the real line, such as $[-1, 1]$ or $(0, 2]$.
These intervals are sets containing all real numbers between their endpoints; a square bracket next to an endpoint means that the endpoint is included in the set, while parentheses mean that the endpoint is not included.
Intervals usually contain infinitely many elements.
\item Sets which contain all objects satisfying one or more conditions.
For instance, the set $\{ x \in \bbr : x^2 < 4 \}$ contains all real numbers whose square is less than four; in this case it can be written simply as the interval $(-2, 2)$.
A more complicated set is $\{ (x,y) \in \bbr^2 : x + y \leq 3, |x| \leq |y| \}$.
This set contains all two-dimensional vectors which (if $x$ and $y$ are the two components of the vector) satisfy both the conditions $x + y \leq 3$ and $|x| \leq |y|$.
It can also be thought of as the intersection of the sets  $\{ (x,y) \in \bbr^2 : x + y \leq 3 \}$ and  $\{ (x,y) \in \bbr^2 : |x| \leq |y| \}$.
If there are \emph{no} vectors which satisfy all of the conditions, the set is \emph{empty},\index{set!empty} denoted $\emptyset$.
\end{enumerate}
We use the common mathematical conventions that $\bbr$ is the set of all real numbers, and $\bbz$ is the set of all integers.  
If we want to restrict attention to non-negative real numbers or integers (i.e., positive or zero), the notations $\bbr_+$ and $\bbz_+$ are used.
Superscripts (e.g., $\bbr^3$ or $\bbz^5$) indicate vectors of a particular dimension whose elements belong to a particular set: $\bbr^3$ is the set of 3-dimensional vectors of real numbers, and $\bbz^5$ the set of 5-dimensional integer vectors.
To refer to a matrix of a particular size, we indicate both dimensions in the superscript; for instance, $\bbz^{3 \times 5}_+$ is the set of matrices with 3 rows and 5 columns, each of whose elements is a non-negative integer.

Sets of scalars or vectors can be described in other ways.
Given any vector $\mb{x} \in \bbr^n$, the \emph{ball of radius $r$}\index{ball} is the set
\labeleqn{balldefinition}{B_r(\mb{x}) = \{ \mb{y} \in \bbr^n : |\mb{x - y}| < r \} \,,}
that is, the set of all vectors whose distance to $\mb{x}$ is less than $r$, where $r$ is some positive number.
A ball in one dimension is an interval; a two-dimensional ball is a circle; a three-dimensional ball is a sphere; a four-dimensional ball a hypersphere, and so on.

Given some set of real numbers $X$, the vector $\mb{x}$ is a \emph{boundary point}\index{boundary point} of $X$ if \emph{every} ball $B_r(\mb{x})$ contains both elements in $X$ and elements not in $X$, no matter how small the radius $r$.
Notice that the boundary points of a set need not belong to the set: 2 is a boundary point of the interval $(-2, 2)$.
A set is \emph{closed}\index{set!closed} if it contains all of its boundary points.
A set $X$ is \emph{bounded}\index{set!bounded} if every element of $X$ is contained in a sufficiently large ball centered at the origin, that is, if there is some $r$ such that $\mb{x} \in B_r(\mb{0})$ for all $\mb{x} \in X$.
A set is \emph{compact}\index{set!compact} if it is both closed and bounded.\footnote{
The definitions given in this paragraph apply whenever the elements of $X$ are real numbers.
A subtler touch is needed for other kinds of sets; see \cite{munkres00} or another text on topology if you are interested.}

These facts will prove useful:
\begin{prp}
\label{prp:setfacts1} Let $f(x_1, x_2, \cdots, x_n)$ be any linear function, that is, $f(\mb{x}) = a_1 x_1 + a_2 x_2 + \cdots + a_n x_n$ for some constants $a_i$, and let $b$ be any scalar.
\begin{enumerate}[(a)]
\item The set $\{ \mb{x} \in \bbr^n : f(\mb{x}) = b \}$ is closed.
\item The set $\{ \mb{x} \in \bbr^n : f(\mb{x}) \leq b \}$ is closed.
\item The set $\{ \mb{x} \in \bbr^n : f(\mb{x}) \geq b \}$ is closed.
\end{enumerate}
\end{prp}
\begin{prp}
\label{prp:setfacts2} Let $X$ and $Y$ be any sets of scalars or vectors.
\begin{enumerate}[(a)]
\item If $X$ and $Y$ are closed, so are $X \cap Y$ and $X \cup Y$.
\item If $X$ and $Y$ are bounded, so are $X \cap Y$ and $X \cup Y$.
\item If $X$ and $Y$ are compact, so are $X \cap Y$ and $X \cup Y$.
 \end{enumerate}
\end{prp}
Combining Propositions~\ref{prp:setfacts1} and~\ref{prp:setfacts2}, we can see that any set defined solely by linear equality or weak inequality constraints (any number of them) is closed.

\index{set!convex|see {convex set}}
We next discuss what it means for a set to be \emph{convex}.\index{convex set}
This notion is very important, and is described in more detail than the other concepts, starting with an intuitive definition.

If $X$ is convex, geometrically this means that line segments connecting points of $X$ lie entirely within $X$.
For example, the set in Figure~\ref{fig:convexset} is convex, while those in Figure~\ref{fig:nonconvexset1} and Figure~\ref{fig:nonconvexset2} are not.
Intuitively, a convex set cannot have any ``holes'' punched into it, or ``bites'' taken out of it.

\genfig{convexset}{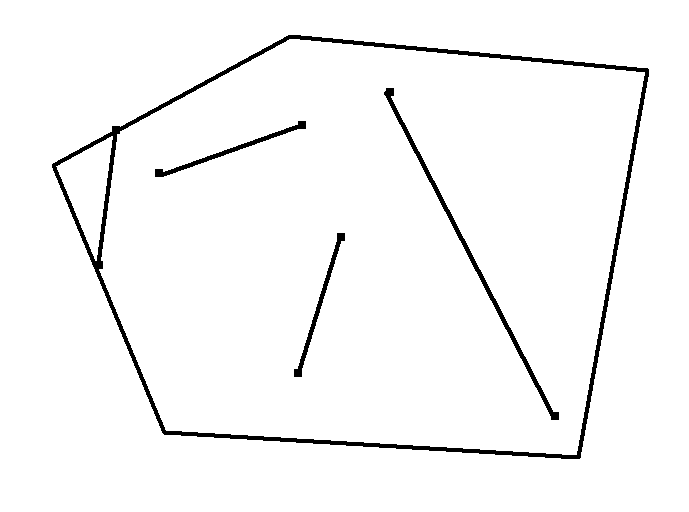}{A convex set.}{width=0.5\textwidth}
\genfig{nonconvexset1}{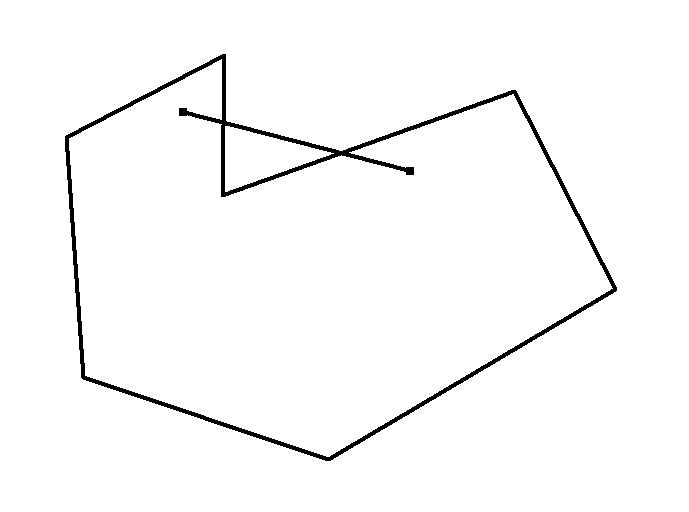}{A nonconvex set with a ``bite'' taken out of it.}{width=0.5\textwidth}
\genfig{nonconvexset2}{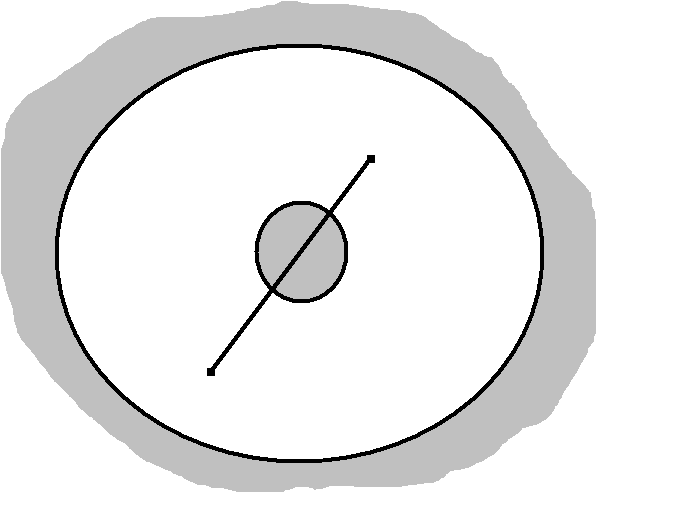}{A nonconvex set with a ``hole'' in it.}{width=0.5\textwidth}

Mathematically, we write this as follows:
\begin{dfn} A set $X \subseteq \bbr^n$ is convex if, for all $x_1, x_2 \in X$ and all $\lambda \in [0,1]$, the point $\lambda x_2 + (1 - \lambda) x_1 \in X$.
\end{dfn}
If this definition is not clear, notice that one way to express the line segment between any two points is $\lambda x_2 + (1 - \lambda) x_1$, and that you cover the entire line between $x_1$ and $x_2$ as $\lambda$ varies between 0 and 1, regardless of how close or far apart these two points are located.

\begin{exm}
Show that the one-dimensional set $X = \{ x : x \geq 0 \}$ is convex.
\end{exm}
\textbf{Solution}.
Pick any $x_1, x_2 \geq 0$ and any $\lambda \in [0, 1]$.
Because $x_1$, $x_2$, and $\lambda$ are all nonnegative, so are $\lambda x_2$ and $(1 - \lambda) x_1$, and therefore so is $\lambda x_2 + (1 - \lambda) x_1$.
Therefore $\lambda x_2 + (1 - \lambda) x_1$ belongs to $X$ as well.
$\blacksquare$

\begin{exm}
Show that the hyperplane $X = \myc{\mb{x} \in \bbr^n : \sum_{i = 1}^n a_i x_i - b = 0 }$ is convex.
\label{exm:hyperplane} \end{exm}
\textbf{Solution}.
This set is the same as $\{\mb{x} \in \bbr^n : \sum_{i = 1}^n a_i x_i  = b\}$ Pick any $\mb{x}, \mb{y} \in X$ and any $\lambda \in [0, 1]$.
Then
\begin{align*}
   \sum_{i=1}^n a_i (\lambda y_i + (1 - \lambda) x_i) &= \lambda \sum_{i=1}^n a_i y_i + (1 - \lambda) \sum_{i=1}^n a_i x_i \\
                                                               &= \lambda b + (1 - \lambda) b \\
                                                               &= b
\end{align*}
so $\lambda \mb{y} + (1 - \lambda) \mb{x} \in X$ as well.
$\blacksquare$

Sometimes, more complicated arguments are needed.

\begin{exm}
Show that the two-dimensional ball $B = \myc{ (x,y) : x^2 + y^2 \leq 1 }$ is convex.
\end{exm}
\textbf{Solution}.
Pick any vectors $\mb{a}, \mb{b} \in B$ and any $\lambda \in [0, 1]$.
We will write the components of these as $\mb{a} = \vect{a_x & a_y}$ and $\mb{b} = \vect{b_x & b_y}$.
The point $\lambda \mb{b} + (1 - \lambda) \mb{a}$ is the vector $\vect{ \lambda b_x + (1 - \lambda) a_x & \lambda b_y + (1 - \lambda) a_y}$.
To show that it is in $B$, we must show that the sum of the squares of these components is no greater than 1.
\begin{multline*}
(\lambda b_x + (1 - \lambda) a_x)^2 + (\lambda b_y + (1 - \lambda) a_y)^2 \\
   = \lambda^2(b_x^2 + b_y^2) + (1 - \lambda)^2 ( a_x^2 + a_y^2 ) + 2 \lambda (1 - \lambda) (a_x b_x + a_y b_y) \\
   \leq \lambda^2 + (1 - \lambda)^2 + 2 \lambda (1 - \lambda) (a_x b_x + a_y b_y)
\end{multline*}
because $\mb{a}, \mb{b} \in B$ (and therefore $a_x^2 + a_y^2 \leq 1$ and $b_x^2 + b_y^2 \leq 1$).
Notice that 
$a_x b_x + a_y b_y$ is simply the dot product of $\mb{a}$ and $\mb{b}$, which is equal to $|\mb{a}| |\mb{b}| \cos \theta$, where $\theta$ is the angle between the vectors $\mb{a}$ and $\mb{b}$.
Since $|\mb{a}| \leq 1$, $|\mb{b}| \leq 1$ (by definition of $B$), and since $\cos \theta \leq 1$ regardless of $\theta$, $a_x b_x + a_y b_y \leq 1$.
Therefore 
\begin{multline*}
\lambda^2 + (1 - \lambda)^2 + 2 \lambda (1 - \lambda) (a_x b_x + a_y b_y) \\
   \leq \lambda^2 + (1 - \lambda)^2 + 2 \lambda (1 - \lambda) 
   = (\lambda + (1 - \lambda))^2 = 1   
\end{multline*}
so the point $\lambda \mb{b} + (1 - \lambda) \mb{a}$ is in $B$ regardless of the values of $\mb{a}$, $\mb{b}$, or $\lambda$.
Thus $B$ is convex.
$\blacksquare$

\begin{exm}
Show that the set $C = \myc{ \vect{x & y} \in \bbr^2 : x^2 + y^2 > 1 }$ is not convex.
\end{exm}
\textbf{Solution}.
Let $\mb{b^1} = \vect{2 & 0}$, $\mb{b^2} = \vect{-2 & 0}$, $\lambda = 1/2$.
Then $\lambda \mb{b^2} + (1 - \lambda) \mb{b^1} = \vect{0 & 0} \notin C$ even though $\mb{b^1}, \mb{b^2} \in C$ and $\lambda \in [0, 1]$.
Therefore $C$ is not convex.
$\blacksquare$

Proving that a set is convex requires showing that something is true for \emph{all} possible values of $x_1, x_2 \in X$, and $\lambda \in [0, 1]$, whereas disproving convexity only requires you to pick one combination of these values where the definition fails.
Lastly, this result will prove useful in showing that a set is convex:
\begin{prp}
\label{prp:convexintersection}
 If $X$ and $Y$ are convex sets, so is $X \cap Y$.
\end{prp}

We conclude by discussing the projection\index{projection} operation; intuitively the idea is that we are given a point and a set, and want to find the point in the set closest to the given point.
Mathematically, if $X$ is a closed set of $n$-dimensional vectors, and $\mb{y} \in \bbr^n$ is any $n$-dimensional vector, the \emph{projection of $\mb{y}$ onto $X$}, written $\mr{proj}_X (\mb{y})$, is the vector $\mb{x}$ in $X$ which is ``closest'' to $\mb{y}$ in the sense that $|\mb{x} - \mb{y}|$ is minimal.
In Figure~\ref{fig:projection}, we have $\mr{proj}_X (a) = b$, $\mr{proj}_X (c) = d$, and $\mr{proj}_X (e) = e$.

\stevefig{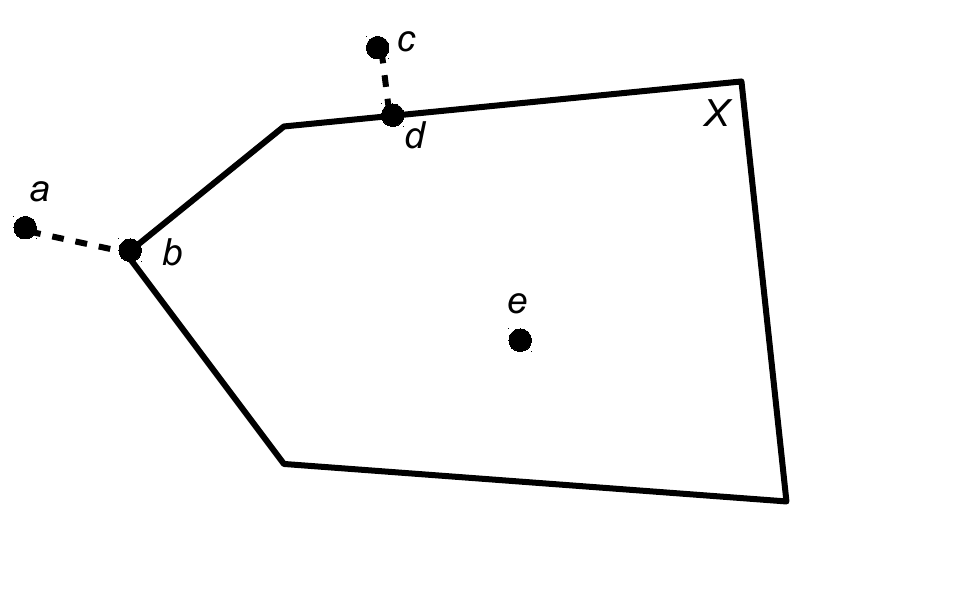}{Some examples of the projection operation.}{0.8\textwidth}

The projection operation is well-behaved when $X$ is a convex set.
If $X$ is convex, closed, and non-empty, the projection operation is always well-defined (there is exactly one point which is ``closest,'' so no need to worry about ties), and is continuous (as $\mb{y}$ varies, the projection $\mr{proj}_X (\mb{y})$ varies smoothly).
These facts will be useful when we develop algorithms that involve projection onto a convex set.
\index{set|)}

\section{Functions}
\label{sec:functions}

\index{function|(}
A function is a mapping between sets.
If the function $f$ maps set $X$ to set $Y$, then $f$ associates \emph{every} element of $X$ with \emph{some} single element of $Y$.
(Note that not every element of $Y$ needs to be associated with an element of $X$.)  The set $X$ is known as the \emph{domain} of $f$.
The set $X$ is called the \emph{domain}\index{function!domain} of $f$.
Examples include the familiar functions $f(x) = x^2$ and $g(x,y) = x^2 + y^2$.
The function $f$ maps $\bbr$ to $\bbr$, while $g$ maps $\bbr^2$ to $\bbr$.
An example of a function which maps $\bbr^2$ to $\bbr^2$ is the vector-valued function
\[ \mb{h}(x_1, x_2) = \vect{3x_1 + 2x_2 \\ -x_1 x_2} \]
The \emph{inverse}\index{function!inverse} of a function $f$, denoted $f^{-1}$, ``undoes'' the mapping $f$ in the sense that if $f(x) = y$, then $f^{-1}(y) = x$.
As an example, if $f(x) = x^3$, then $f^{-1}(x) = \sqrt[3]{x}$.
Not every function has an inverse, and inverse functions may need to be restricted to subsets of $X$ and $Y$.

The \emph{composition}\index{function!composition} of two functions $f$ and $g$, denoted $f \circ g$ involves substituting function $g$ into function $f$.
If $f(x) = x^3$ and $g(x) = 2x + 4$, then $f \circ g(x) = (2x + 4)^3$.
This can also be written as the function $f(g(x))$.

\index{function!continuous|see {continuous function}}
A function is \emph{continuous}\index{continuous function} if, for any point ${\hat{x}} \in X$, the limit $\lim_{{x} \rightarrow {\hat{x}}} f({x})$ exists and is equal to $f({\hat{x}})$.
Intuitively, continuous functions can be drawn without lifting your pen from the paper.
A function of a single variable is \emph{differentiable}\index{function!differentiable}\index{differentiable function} if, for any point $x \in X$, the limit\index{derivative}
\labeleqn{derivative}{\lim_{\Delta {x} \rightarrow \mb{0}} \frac{f({x} + \Delta {x}) - f({x})}{\Delta {x}} }
exists; that limit is then called the \emph{derivative} of $f$, and gives the slope of the tangent line at any point.
It can be shown that a differentiable function must be continuous.

\begin{prp}
\index{continuous function!properties}
\label{prp:continuousdifferentiable}
 Let $f$ and $g$ be continuous functions.
Then we have the following:
\begin{itemize}
\item The multiple $\alpha f$ is a continuous function for any scalar $\alpha$.
\item The sum $f + g$ is a continuous function.
\item The product $fg$ is a continuous function.
\item The composition $f \circ g$ is a continuous function.
\end{itemize}
Furthermore, all of these results hold if ``continuous'' is replaced by ``differentiable.''
\end{prp}

\begin{prp}
\label{prp:continuousprojection}
\index{projection!properties}
Let $X$ be a nonempty, convex set of $n$-dimensional vectors.
Then the projection function $\mr{proj}_X(\mb{x})$ is defined and continuous for all $\mb{x} \in \bbr^n$.
\end{prp}

Differentiability is more complicated when dealing with functions of multiple variables (that is, functions whose domain is $\bbr^2$, $\bbr^3$, and so forth).
The basic notion is the \emph{partial derivative},\index{partial derivative}\index{derivative!partial} in which all variables except one are assumed constant, and an ordinary derivative is taken with respect to the remaining variable.
For instance, if $f(x,y) = x^2 + 3xy + y^3$, the partial derivatives with respect to $x$ and $y$ are
 \[ \pdr{f}{x} = 2x + 3y \qquad \pdr{f}{y} = 3x + 3y^2 \,.
\]
Second partial derivatives are found in the same way, taking partial derivatives of $\pdr{f}{x}$ and $\pdr{f}{y}$.
If all of the second partial derivatives of a function are continuous at a point, the order of differentiation does not matter and 
\labeleqn{mixedpartials}{ \pdrc{f}{x}{y} = \pdrc{f}{y}{x} }
at that point.

All of these partial derivatives can be organized into different kinds of vectors and matrices.
The \emph{gradient}\index{gradient} of a function is the vector of all its partial derivatives.
For $f(x,y) = x^2 + 3xy + y^3$, this is
\[ \renewcommand*{\arraystretch}{1.5} \nabla f(x,y) = \vect{ \pdr{f}{x} \\ \pdr{f}{y} } = \vect{ 2x + 3y \\ 3x + 3y^2 }\,.
\]
A point $\mb{x}$ where $\nabla f(\mb{x}) = \mb{0}$ is called a \emph{stationary point}.
\index{stationary point}
(If $f$ is a function of a single variable, this is just a point where the derivative vanishes.)  The \emph{Hessian}\index{Hessian matrix} of a function is the matrix of all of its second partial derivatives.
For this same function, we have
\[ \renewcommand*{\arraystretch}{1.5} Hf(x,y) = \vect{ \pdrb{f}{x} & \pdrc{f}{x}{y} \\ \pdrc{f}{y}{x} & \pdrb{f}{y} } = \vect{ 2 & 3 \\ 3 & 6y} \,. \]
By equation~\eqn{mixedpartials}, if the Hessian matrix contains continuous functions, it is symmetric.

Now consider a vector-valued function of multiple variables $f: \bbr^n \rightarrow \bbr^m$.
(That is, $f$ has $n$ input variables and produces an $m$-dimensional vector as output.)
The \emph{Jacobian}\index{Jacobian matrix} of $f$ is a matrix of all of its \emph{first} partial derivatives.
For a concrete example, let 
\[ \mb{g}(x,y,z) = \vect{ g_1(x,y,z) \\ g_2(x,y,z) } = \vect{ x^2 + y^2 + z^2 \\ xyz }\,. \]
The Jacobian is the matrix
\[ \renewcommand*{\arraystretch}{1.5} J\mb{g}(x,y,z) = \vect{ \pdr{g_1}{x} & \pdr{g_1}{y} & \pdr{g_1}{z} \\ \pdr{g_2}{x} & \pdr{g_2}{y} & \pdr{g_2}{z} } = \vect{ 2x & 2y & 2z \\ yz & xz & xy } \,. \]
Be sure to note the difference between Hessians and Jacobians.
The Hessian is defined for scalar-valued functions, and contains all the second partial derivatives.
The Jacobian is defined for vector-valued functions, and contains all the first partial derivatives.

\index{convex function|(}
Finally, there is an important notion of function \emph{convexity}.\index{function!convex|see {convex function}}
Confusingly, this is a different idea than set convexity discussed in Section~\ref{sec:sets}, although there are some relationships and similar ideas between them.
Set and function convexity together play a pivotal role in optimization and network equilibrium problems, so function convexity is discussed at length here.
Geometrically, a convex function lies below its secant lines.
Remember that a secant line is the line segment joining two points on the function.
As we see in Figure~\ref{fig:convexsecant}, no matter what two points we pick, the function always lies below its secant line.
On the other hand, in Figure~\ref{fig:nonconvexsecant}, not every secant line lies above the function: some lie below it, and some lie both above and below it.
Even though we can draw some secant lines which are above the function, this isn't enough: \emph{every} possible secant must lie above the function.
For this concept to make sense, the domain $X$ of the function must be a convex set, an assumption which applies for the remainder of this section

\genfig{convexsecant}{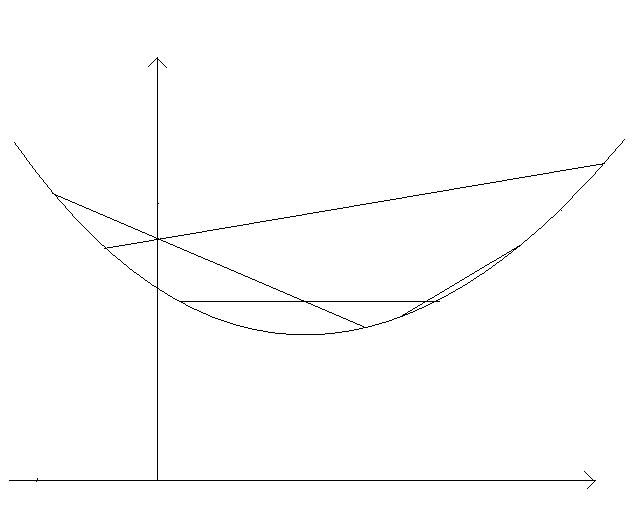}{A convex function lies below all of its secants.}{width=0.7\textwidth}
\genfig{nonconvexsecant}{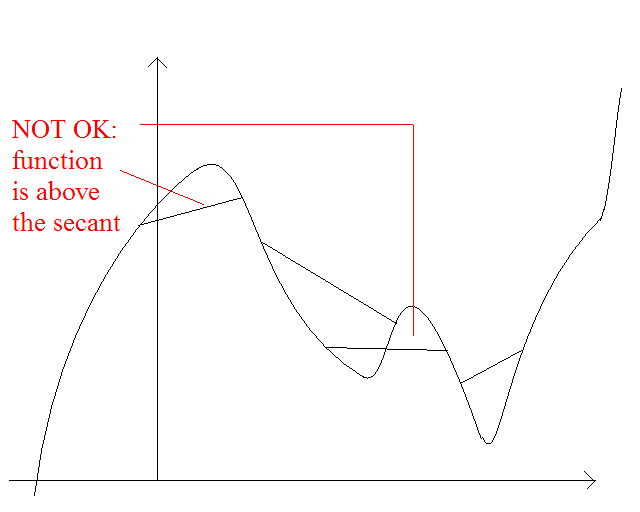}{A nonconvex function does not lie below all of its secants.}{width=0.7\textwidth}

The following definition makes this intuitive notion formal:
\begin{dfn} A function $f: X \rightarrow \bbr$ is \emph{convex} if, for every $x_1, x_2 \in X$ and every $\lambda \in [0, 1]$, 
\labeleqn{convexdefinition}{f((1 - \lambda) x_1 + \lambda x_2) \leq (1 - \lambda) f(x_1) + \lambda f(x_2)}
and \emph{strictly convex}\index{convex function!strict convexity} if 
\labeleqn{strictconvex}{f((1 - \lambda) x_1 + \lambda x_2) < (1 - \lambda) f(x_1) + \lambda f(x_2)}
for all distinct $x_1, x_2 \in X, \lambda \in (0, 1)$ \label{dfn:convex} \end{dfn}

Essentially, $x_1$ and $x_2$ are the two endpoints for the secant line.
Since this entire line segment must be above the function, we need to consider every point between $x_1$ and $x_2$.
This is what $\lambda$ does: as $\lambda$ varies between 0 and 1, the points $\lambda x_2 + (1 - \lambda) x_1$ cover every point between $x_1$ and $x_2$.
You can think of $\lambda x_2 + (1 - \lambda) x_1$ as a ``weighted average,'' where $\lambda$ is the weight put on $x_2$.
For $\lambda = 0$, all the weight is on $x_1$.
For $\lambda = 1$, all the weight is on $x_2$.
For $\lambda = 1/2$, equal weight is put on the two points, so the weighted average is the midpoint.
$\lambda = 1/3$ corresponds to the point a third of the way between $x_1$ and $x_2$.

So we need to say that, at all such intermediate $x$ values, the value of $f$ is lower than the $y$-coordinate of the secant.
The value of the function at this point is simply $f((1 - \lambda) x_1 + \lambda x_2)$.
Because the secant is a straight line, its $y$-coordinate can be seen as a weighted average of the $y$-coordinates of its endpoints, that is, $f(x_1)$ and $f(x_2)$.
This weighted average can be written as $(1 - \lambda) f(x_1) + \lambda f(x_2)$, so requiring the function to lie below the secant line is exactly the same as enforcing condition~\eqn{convexdefinition} for all possible secant lines: that is, for all $x_1, x_2 \in X$ and all $\lambda \in [0,1]$.

Figure~\ref{fig:convexdetail} explains this in more detail.
Along the horizontal axis, the secant endpoints $x_1$ and $x_2$ are shown, along with an intermediate point $\lambda x_2 + (1 - \lambda)x_1$.
The $y$-coordinates are also shown: at the endpoints, these are $f(x_1)$ and $f(x_2)$.
At the intermediate point, the $y$-coordinate of the function is $f(\lambda x_2 + (1 - \lambda)x_1)$, while the $y$-coordinate of the secant is $\lambda f(x_2) + (1 - \lambda) f(x_1)$.
Because the function is convex, the former can be no bigger than the latter.
Time spent studying this diagram is very well spent.
Make sure you understand what each of the four points marked on the diagram represents, and why the given mathematical expressions correctly describe these points.
Make sure you see what role $\lambda$ plays: as $\lambda$ increases from 0 to 1, the central vertical line moves from $x_1$ to $x_2$.
(What would happen if we picked $x_1$ and $x_2$ such that $x_1 > x_2$?  What if $x_1 = x_2$?)

\stevefig{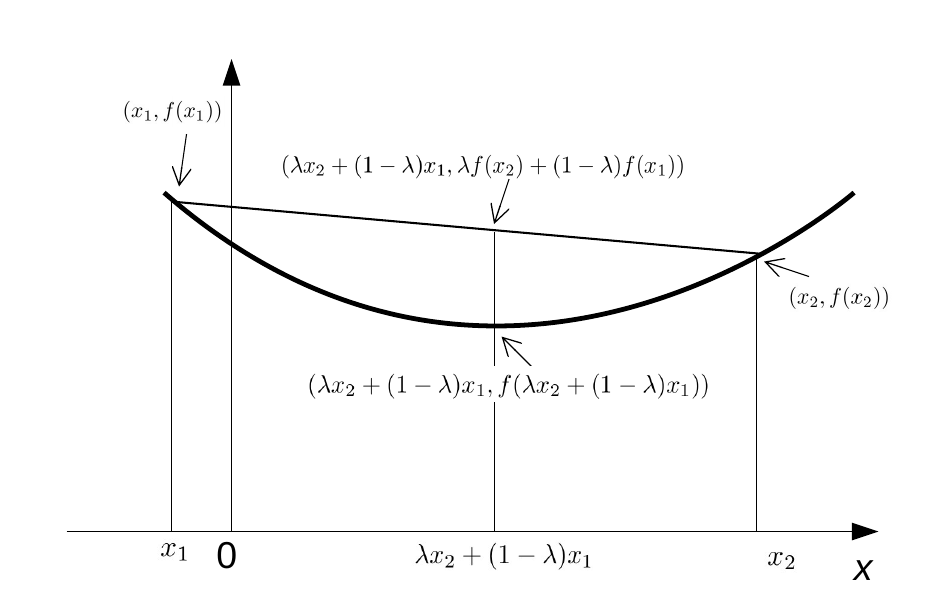}{All of the relevant points for the definition of convexity.}{\textwidth}

\begin{exm} Is the function $f(x) = |x|$, $x \in \bbr$ convex?  Is it strictly convex?  \end{exm}
\textbf{Solution}.
To see if $f$ is convex, we need to see if~\eqn{convexdefinition} is true; to see if it is strictly convex, we need to check~\eqn{strictconvex}.
Furthermore, these inequalities have to be true for \emph{every} $x_1, x_2 \in \bbr$, and \emph{every} $\lambda \in [0, 1]$.
It is not enough to simply pick a few values randomly and check the equations.
So, we have to work symbolically.
In this case,
\begin{align*}
f((1 - \lambda) x_1 + \lambda x_2)  &= \left|(1 - \lambda) x_1 + \lambda x_2\right| \\
                                    &\leq \left|(1 - \lambda) x_1 \right| + \left|\lambda x_2\right| \quad \mbox{by the triangle inequality} \\
                                    &= (1 - \lambda) |x_1| + \lambda |x_2| \quad \mbox{because $\lambda, 1 - \lambda \geq 0$} \\
                                    &= (1 - \lambda) f(x_1) + \lambda f(x_2)
\end{align*}
Therefore~\eqn{convexdefinition} is satisfied, so $f$ is convex.
To show that it is strictly convex, we would have to show that the inequality 
\[ \left|(1 - \lambda) x_1 + \lambda x_2\right| \leq \left|(1 - \lambda) x_1 \right| + \left|\lambda x_2\right| \]
can be replaced by a \emph{strict} inequality $<$.
However, we can't do this: for example, if $x_1 = 1$, $x_2 = 2$, $\lambda = 0.5$, the left side of the inequality ($|1/2 + 2/2| = 3/2$) is exactly equal to the right side ($|1/2| + |2/2| = 3/2$).
So $f$ is not strictly convex.
$\blacksquare$

Note that proving that $f(x)$ is convex requires a \emph{general} argument, where proving that $f(x)$ was not strictly convex only required a single counterexample.
This is because the definition of convexity is a ``for all'' or ``for every'' type of argument.
To prove convexity, you need an argument that allows for all possible values of $x_1$, $x_2$, and $\lambda$, whereas to disprove it you only need to give one set of values where the necessary condition doesn't hold.

\begin{exm} \label{exm:linearconvex} Show that every linear function $f(x) = ax + b$, $x \in \bbr$ is convex, but not strictly convex.\end{exm}
\textbf{Solution}.

\begin{align*}
f((1 - \lambda) x_1 + \lambda x_2)  &= a((1 - \lambda) x_1 + \lambda x_2) + b \\
                                    &= a((1 - \lambda) x_1 + \lambda x_2) + ((1 - \lambda) + \lambda) b \\
                                    &= (1 - \lambda) (a x_1 + b) + \lambda (ax_2 + b) \\
                                    &= (1 - \lambda) f(x_1) + \lambda f(x_2)
\end{align*}
So we see that inequality~\eqn{convexdefinition} is in fact satisfied as an \emph{equality}.
That's fine, so every linear function is convex.
However, this means we can't replace the inequality $\leq$ with the strict inequality $<$, so linear functions are not strictly convex. $\blacksquare$

Sometimes it takes a little bit more work, as in the following example:
\begin{exm}Show that $f(x) = x^2$, $x \in \bbr$ is strictly convex.\end{exm}
\textbf{Solution}.
Pick $x_1, x_2$ so that $x_1 \neq x_2$, and pick $\lambda \in (0,1)$.
\begin{align*}
f((1 - \lambda) x_1 + \lambda x_2)  &= ((1 - \lambda) x_1 + \lambda x_2)^2 \\
                                    &= (1 - \lambda)^2 x_1^2 + \lambda^2 x_2^2 + 2(1 - \lambda)\lambda x_1 x_2 
\end{align*}
What to do from here?  
Comparing term by term, because $\lambda \in (0, 1)$, we know that $\lambda^2 < \lambda$ and $(1 - \lambda)^2 < 1 - \lambda$.
Therefore $(1 - \lambda^2) x_1^2 < (1 - \lambda) x_1^2$ and $\lambda^2 x_2^2 < \lambda x_2^2$.
This takes care of two of the terms, all we have left is  
Since $x_1 \neq x_2$, $(x_1 - x_2)^2 > 0$.
Expanding, this means that $x_1^2 + x_2^2 > 2x_1 x_2$.
This means that 
\begin{align*}
(1 - \lambda)^2 x_1^2 + \lambda^2 x_2^2 + 2(1 - \lambda)\lambda x_1 x_2 &< (1 - \lambda)^2 x_1^2 + \lambda^2 x_2^2 + (1 - \lambda)(\lambda)(x_1^2 + x_2^2) \\
&= (1 - \lambda) x_1^2 + \lambda x_2^2 \mbox { (after some algebra) } \\
&= (1 - \lambda) f(x_1) + \lambda f(x_2)
\end{align*}
which proves strict convexity.
$\blacksquare$

This last example shows that proving convexity can require some effort, even for simple functions like $x^2$.
The good news is that there are often simpler conditions that we can check.
These conditions involve the first and second derivatives of a function.

\begin{prp} Let $f: X \rightarrow \bbr$ be a differentiable function, where $X$ is a subset of $\bbr$.
Then $f$ is convex if and only if
\[ f(x_2) \geq f(x_1) + f'(x_1)(x_2 - x_1) \]
for all $x_1, x_2 \in X$.
\label{prp:convex1}
\end{prp}

\begin{prp} Let $f: X \rightarrow \bbr$ be twice differentiable, where $X$ is a subset of $\bbr$, and let $f$ be twice differentiable on $X$.
Then $f$ is convex if and only if $f''(x) \geq 0$ for all $x \in X$ \label{prp:convex2}
\end{prp}

\index{convex function!strict convexity}
Equivalent conditions for strict convexity can be obtained in a natural way, changing $\geq$ to $>$ and requiring that $x_1$ and $x_2$ be distinct in Proposition~\ref{prp:convex1}.
Proposition~\ref{prp:convex2} also changes slightly in this case; $f''(x) > 0$ is sufficient for strict convexity but is not necessary.
\label{strictconvexity} You are asked to prove these statements in the exercises.
Essentially, Proposition~\ref{prp:convex1} says that $f$ lies \emph{above} its tangent lines (Figure~\ref{fig:convextangent}), while Proposition~\ref{prp:convex2} says that $f$ is always ``curving upward.'' (A convex function lies \emph{above} its tangents, but \emph{below} its secants.)

\genfig{convextangent}{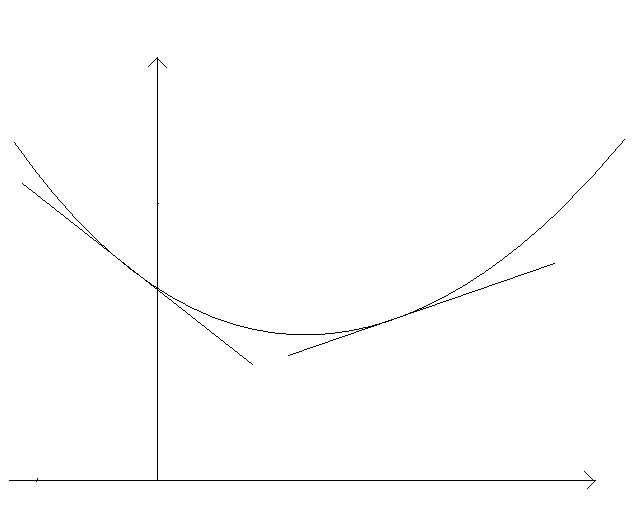}{A convex function lies above its tangents.}{width=0.5\textwidth}

These conditions are usually easier to verify than that of Definition~\ref{dfn:convex}.

\begin{exm} Show that $f(x) = x^2$ is strictly convex using Proposition~\ref{prp:convex1} \end{exm}
\textbf{Solution}.
Pick any $x_1, x_2 \in \bbr$ with $x_1 \neq x_2$.
We have $f'(x_1) = 2x_1$, so we need to show that
\[ x^2_2 > x_1^2 + 2x_1 (x_2 - x_1) \]
Expanding the right-hand side and rearranging terms, we see this is equivalent to
\[ x_1^2 - 2x_1x_2 + x_2^2 > 0 \]
or
\[ (x_1 - x_2)^2 > 0 \]
which is clearly true since $x_1 \neq x_2$.
Thus $f$ is strictly convex.
$\blacksquare$

\begin{exm} Show that $f(x) = x^2$ is strictly convex using Proposition~\ref{prp:convex2} \end{exm}
\textbf{Solution}.
$f''(x) = 2 > 0$ for all $x \in \bbr$, so $f$ is strictly convex. $\blacksquare$

Much simpler!  If $f$ is differentiable (or, better yet, twice differentiable) checking these conditions is almost always easier.

Furthermore, once we know that some functions are convex, we can use this to show that many other combinations of these functions must be convex as well.

\begin{prp} \label{prp:sumconvex}  If $f$ and $g$ are convex functions, and $\alpha$ and $\beta$ are positive real numbers, then $\alpha f + \beta g$ is convex as well.
\label{prp:linearcombo} \end{prp}
\begin{prp} If $f$ and $g$ are convex functions, then $f \circ g$ is convex as well.
\label{prp:composition} \end{prp}

Some common convex functions are $|x|$, $x^2$, $e^x$, and $ax + b$.
So, Proposition~\ref{prp:linearcombo} tells us that $3x^2 + 4|x|$ is convex.
It also tells us that any quadratic function $ax^2 + bx + c$ is convex as long as $a > 0$.
Proposition~\ref{prp:composition} says that the composition of two convex functions is convex as well.
For instance, $e^{x^2}$ is convex, and $x^4 = (x^2)^2$ is convex as well.

What about functions of more than one variable?  The ``shortcut'' conditions in Propositions~\ref{prp:convex1} and~\ref{prp:convex2} only apply if the domain of $f$ is one-dimensional.
It turns out that very similar conditions can be given for the multivariable case.
The multi-dimensional equivalent of the first derivative is the gradient, and the equivalent of the second derivative is the Hessian.

The equivalent conditions on convexity are
\begin{prp} Let $f: X \rightarrow \bbr$ be a function whose gradient exists everywhere on $X$.
Then $f$ is convex if and only if
\[ f(\mb{x_2}) \geq f(\mb{x_1}) + \nabla f(x_1) \cdot (\mb{x_2} - \mb{x_1}) \]
for all $\mb{x_1}, \mb{x_2} \in X$.
\label{prp:convex3}
\end{prp}

\begin{prp} Let $f: X \rightarrow \bbr$ be a function whose Hessian exists everywhere on $X$.
Then $f$ is convex if and only if $H(f)$ is positive semidefinite for all $\mb{x} \in X$ \label{prp:convex4}; and $f$ is strictly convex if $H(f)$ is positive definite for all $\mb{x} \in X$.
\end{prp}

Unfortunately, neither of these is as easy to check as the single-dimension equivalents.
In particular, it is rather tedious to check whether or not a matrix is positive semidefinite or not.
Fortunately, in many important cases in transportation network analysis, the Hessian is diagonal and Proposition~\ref{prp:diagonalpd} applies.

As a final note, one useful link between convex sets and convex functions is the following result:
\begin{prp} If $g: \bbr^n \rightarrow \bbr$ is a convex function, then the set $X = \{ x : g(x) \leq 0$ is a convex set.
\label{prp:levelset} \end{prp}
\index{convex function|)}
\index{function|)}

\section{Exercises}
\label{sec:mathematicalpreliminariesappendix_exercises}

\begin{enumerate}
\item \diff{22} For this exercise, let $x_{11} = 5$, $x_{12} = 6$, $x_{13} = 7$, $x_{21} = 4$, $x_{22} = 3$, and $x_{23} = 9$.
Evaluate each of these sums.
\label{ex:sumpractice}
\begin{enumerate}[(a)]
  \item $\sum_{i=1}^2 \sum_{j=1}^3 x_{ij}$
  \item $\sum_{i=1}^2 \sum_{j=2}^3 x_{ij}$
  \item $\sum_{j=1}^3 x_{1j}$
  \item $\sum_{j=1}^2 \sum_{i=1}^3 x_{ji}$
  \item $\sum_{j=1}^3 \sum_{k=0}^1 x_{k+1,j}$
  \item $\sum_{i=1}^2 \sum_{j=i}^3 x_{ij}$
\end{enumerate}
\item \diff{22} Repeat Exercise~\ref{ex:sumpractice}, but with products $\prod$ instead of sums $\sum$.
\item \diff{33} Section~\ref{sec:indicessummation} lists three properties of the summation notation $\sum$.
Formulate and prove analogous properties for the product notation $\prod$.\index{product notation!properties}
You can assume that the product involves only a finite number of factors.
\label{ex:productproperties}
\item \diff{22} Prove Proposition~\ref{prp:diagonalinverse}.
\item \diff{24} Prove Proposition~\ref{prp:diagonalpd}.
\item \diff{53} Prove that if a matrix is invertible, its inverse matrix is unique.
\item \diff{26} For each of the sets below, identify its boundary points and indicate whether or not it is closed, whether or not it is bounded, and whether or not it is convex.
\begin{enumerate}[(a)]
 \item $(5,10]$       
 \item $[4,6]$        
 \item $(0, \infty)$  
 \item $\{x: |x| > 5\}$
 \item $\{x: |x| \leq 5\}$
 \item $\{(x,y) : 0 \leq x \leq 4, -3 \leq y \leq 3\}$
 \item $\{(x,y) : 4x - y = 1\}$
\end{enumerate}
\item \diff{34} Prove Proposition~\ref{prp:setfacts1}.
\item \diff{34} Prove Proposition~\ref{prp:setfacts2}.
\item \diff{53} Identify the projections of the following points on the corresponding sets.
You may find it helpful to draw sketches.
\begin{enumerate}[(a)]
  \item The point $x = 3$ on the set $[0, 1]$.
  \item The point $x = \frac{1}{2}$ on the set $[0, 1]$.
  \item The point $(2,5)$ on the unit circle $x^2 + y^2 = 1$.
  \item The point $(2,5)$ on the line $x + y = 1$.
  \item The point $(2,5)$ on the line segment between $(0,1)$ and $(1,0)$.
  \item The point $(2,3)$ on the line segment between $(0,1)$ and $(1,0)$.
  \item The point $(1,2,3)$ on the sphere $x^2 + y^2 + z^2 = 2$.
\end{enumerate}
\item \diff{20} Prove Proposition~\ref{prp:convexintersection}.
\item \diff{23} Find the inverses of the following functions.
\begin{enumerate}[(a)]
  \item $f(x) = 3x + 4$
  \item $f(x) = 5e^x - 1$
  \item $f(x) = 1/x$
\end{enumerate}
\item \diff{57} Prove Proposition~\ref{prp:continuousdifferentiable} from first principles, using the definitions of continuity and differentiability in the text.
\item \diff{65} Prove Proposition~\ref{prp:continuousprojection}, using the formal definition of continuity.
\item \diff{22} Calculate the gradients of the following functions:
\begin{enumerate}[(a)]
  \item $f(x_1, x_2, x_3) = x_1^2 + 2x_2 x_3 + x_3^2$
  \item $f(x, y) = \frac{xy}{x^2 + y^2 + 1}$
  \item $f(x_1, x_2, d_1, d_2) = 3x_1 + 5x_2^2 + 3(d_1 + 4)^2 + 5(d_2 - 1)^3$
\end{enumerate}
\label{ex:gradients}
\item \diff{25} Calculate the Jacobian matrices of the following functions:
\begin{enumerate}[(a)]
  \item $f(x, y) = \vect{y \\ -x}$
  \item $f(x_1, x_2, x_3) = \vect{ x_1^2 + x_2 + x_3 \\ 3x_2 - 5x_3 }$
  \item $f(x, y) = \vect{3x + y \\ 3y + x \\ xy }$
\end{enumerate}
\item \diff{26} Calculate the Hessian matrices of functions from Exercise~\ref{ex:gradients}.
\item \diff{35}  Determine which of the following functions are convex, strictly convex, or neither.
Justify your answer rigorously.
\begin{enumerate}[(a)]
  \item $f(x) = \log x$ where $x \geq 1$.
  \item $f(x) = -x^3$ where $x \leq 0$.
  \item $f(x) = \sin x$ where $x \in \bbr$.
  \item $f(x) = e^{x^2 - 2x - 3}$ where $x \in \bbr$.
  \item $f(x) = \begin{cases} 0 & \mbox{if $x \leq 0$} \\ x & \mbox{if $x \geq 0$} \end{cases}$
\end{enumerate}
\item \diff{35} Determine which of the following sets is convex.
Justify your answer rigorously.
\begin{enumerate}[(a)]
  \item $X = \{ (x_1, x_2) \in \bbr^2 :  4x_1 - 3x_2 \leq 0 \}$
  \item $X = \{ x \in \bbr : e^x \leq 4 \}$
  \item $X = \{ (x_1, x_2, x_3) \in \bbr^3 : x_1 = 0 \}$
  \item $X = \{ (x_1, x_2, x_3) \in \bbr^3 : x_1 \neq 0 \}$
  \item $X = \{ 1 \}$
\end{enumerate}
\item \diff{31} Show that the integral of an increasing function is a convex function.
\item \diff{54} Show that a differentiable function $f$ of a single variable is convex if and only if $f(x) + f'(x)(y - x) \leq f(y)$ for all $x$ and $y$.
\label{ex:convex1d1d}
\item \diff{54} Show that a twice-differentiable function $f$ of a single variable is convex if and only if $f''$ is everywhere nonnegative.
\label{ex:convex2d1d}
\item \diff{68} Prove the claims of Exercises~\ref{ex:convex1d1d} and~\ref{ex:convex2d1d} for functions of multiple variables.
\item \diff{32} Page~\pageref{strictconvexity} states that ``$f''(x) > 0$ is sufficient for strict convexity but not necessary.'' Give an example of a strictly convex function where $f''(x) = 0$ at some point.
\end{enumerate}

\chapter{Optimization Concepts}
\label{chp:basicoptimization}

This appendix provides a brief introduction on formulating optimization problems.
We provide examples of linear, nonlinear, and integer formulations.
The focus in this appendix is on the basic terminology and formulation of these problems.
Specifically, we provide a detailed explanation of representing optimization formulations using index and matrix notations.
The following appendix goes into more detail on solution methods.

\section{Components of an Optimization Problem}

\index{optimization|(}
An optimization problem is a mathematical model that can help decision makers arrive at decisions which optimize a specific goal, given the constraints they face.
Every optimization problem has three components: an objective function, decision variables, and constraints.
When one talks about \emph{formulating}\index{optimization!formulation} an optimization problem, it means translating a ``real-world'' problem into the mathematical equations and variables which comprise these three components:

\begin{description}
\item [Objective function:]\index{optimization!objective function} A mathematical function which represent the goal of the decision maker, to be either maximized or minimized.
\item [Decision variables:]\index{optimization!decision variable} Variables representing the factors or choices that can be controlled by the decision maker, either directly or indirectly.
\item [Constraints:]\index{optimization!constraint} Equations or inequalities representing restrictions on the values the decision variables can take.
\end{description}

\index{optimization!objective function|(}
The objective function, often\footnote{But not always!  Often problem-specific notation will be used, such as $c$ for cost, and so forth.
The notation given in this section is what is traditionally used if we are referring to a generic optimization problem outside of a specific context.} denoted $f$ or $z$,\label{not:z} reflects a single quantity to be either maximized or minimized.
Examples in the transportation world include ``minimize congestion'', ``maximize safety'', ``maximize accessibility'', ``minimize cost'', ``maximize pavement quality'', ``minimize emissions,'' ``maximize revenue,'' and so forth.
You may object to the use of a \emph{single} objective function, since real-world problems typically involve many different and conflicting objectives from different stakeholders, and this objection is certainly valid.
There are several reasons why the scope of this book is restricted to single objectives.
From a historical perspective, if we don't know how to optimize a single objective function, then we have no hope of being able to optimize multiple objectives simultaneously, since the latter build on the former.
From a pedagogical perspective, the methods of multi-objective optimization are much more complex, and the basic concepts of optimization are best learned in a simpler context first.
From a mathematical perspective, the definition of ``simultaneously optimize'' is very tricky, since there is probably no plan which will, say, simultaneously minimize congestion and agency cost --- therefore multiobjective optimization is ``fuzzier,'' and this fuzziness can be confusing when first taught.
So, don't be alarmed by this restriction to a single objective, but do keep it in the back of your head if you use optimization beyond this book.
\index{optimization!objective function|)}

\index{optimization!decision variable|(}
The decision variables (often, but not always, denoted as the vector $\textbf{x}$) reflect aspects of the problem that you (or the decision maker) can control.
This can include both variables you can directly choose, as well as variables which you indirectly influence by the choice of other decision variables.
For example, if you are a private toll road operator trying to maximize your profit, you can directly choose the toll, so that is a decision variable.
You can't directly choose how many people drive on the road --- but because that's influenced by the toll you chose, the toll road volume should be a decision variable as well.
However, you should avoid including extraneous decision variables.
Every decision variable in your formulation should either directly influence the objective function, or influence another decision variable that affects the objective function.
\index{optimization!decision variable|)}

\index{optimization!constraint|(}
Constraints represent any kind of limitation on the values that the decision variables can take.
The most intuitive types of constraints are those which directly and obviously limit the choices you can make: you can't exceed a budget, you are required by law to provide a certain standard of maintenance, you are not allowed to change the toll by more than \$1 from its current value, and so forth.
One very frequent mistake is to omit ``obvious'' constraints (e.g., the toll can't be negative).
An optimization formulation must be complete and not leave out any constraint, no matter how obvious it may seem to you.
Real-world optimization problems are solved by computers, for which nothing is ``obvious.'' The second type of constraint is required to ensure consistency among the decision variables.
Following the toll road example from the previous paragraph, while you can influence both the toll (directly) and the roadway volume (indirectly), the roadway volume must be consistent with the toll you chose.
This relationship must be reflected by an appropriate constraint; say, a demand function $d(\tau)$ giving the demand $d$ for the tollway in terms of the toll $\tau$, which can be estimated in a variety of econometric ways.
Constraints are the primary way these linkages between decision variables can be captured.
\index{optimization!constraint|)}

A particular choice of decision variables $\textbf{x}$ is \emph{feasible}\index{optimization!feasibility} if it satisfies all of the constraints, and it is \emph{optimal}\index{optimization!optimality} if it maximizes or minimizes the objective function (whichever is appropriate for the given problem).
For some problems, we may need to distinguish local and global optima.
Section~\ref{sec:localglobal} discusses this in greater detail.
In short, we usually want to find a global optimum, but for some complicated optimization problems a local optimum may be the best we can do.

The following example explains these definitions:
\begin{exm}
A plant produces two types of concrete (Type 1 and Type 2).
The owner of the plant makes a profit of \$90 for each truckload of Type 1 concrete produced and a profit of \$120 for each truckload of Type 2 concrete produced.
Three materials are required to produce concrete: cement, aggregate, and water.
The plant requires 30 units of cement, 50 units of aggregate, and 60 units of water to produce one truckload of Type 1 concrete.
The plant requires 40 units of cement, 20 units of aggregate, and 90 units of water to produce one truckload of Type 2 concrete.
The plant is supplied 2000 units of cement, 2500 units of aggregate, and 4000 units of water daily.
How many truckloads of Type 1 and Type 2 cement should the plant produce to maximize daily profit?
\end{exm}

In any optimization problem, the first step is to identify the decision variables or the variables which can be controlled by the decision maker.
In this problem, the decision variables are the daily truckloads of Type 1 and Type 2 concrete produced by the plant.

Let $x$ represent the daily truckloads of Type 1 concrete produced by the plant and $y$ represent the daily truckloads of Type 2 concrete produced by the plant.

The second step is to identify the constraints or restrictions on the decision variables.
In this problem there are restrictions on the total volume of cement, aggregate, and water supplied to the plant daily which limits the total amount of Type 1 and Type 2 concrete which can be produced.

Each truckload of Type 1 concrete requires 30 units of cement.
Therefore, to produce $x$ truckloads of Type 1 concrete requires $30x$ units of cement.
Each truckload of Type 2 concrete requires 40 units of cement.
Therefore, to produce $y$ truckloads of Type 2 concrete requires $40y$ units of cement.
The plant is supplied 2000 units of cement daily.
Therefore, the constraint on the total cement consumed by the plant can be written as 
\begin{equation}
30x + 40y \leq 2000
\,.
\end{equation}
The plant is supplied 2500 units of aggregate and 4000 units of water daily.
Along similar lines, the constraint on the total volume of aggregate and water consumed by the plant can be written as
\begin{align}
50x + 20y &\leq 2500 \\
60x + 90y &\leq 4000
\,.
\end{align}
In addition, the daily truckloads of Type 1 and Type 2 concrete produced has to be greater than or equal to zero.
The non-negativity constraint can be represented as
\begin{equation}
x \geq 0, y \geq 0
\,.
\end{equation}

Each Type 1 concrete truckload generates a profit of \$90.
Therefore, $x$ truckloads of Type 1 concrete generates \$$90x$ profit.
Each Type 2 concrete truckload generates a profit of \$120.
Therefore, $y$ truckloads of Type 2 concrete generates a profit of \$$120y$.
The goal of the decision maker is to maximize this daily profit which can be written as
\begin{equation}
\max \,\, 90x + 120y
\,.
\end{equation}

The entire problem can be summarized as below.
The following set of equations represents an optimization formulation or a mathematical programming formulation for the concrete plant profit maximizing problem:
\optimizex{\max}{x,y}{90x + 120y}{
                    & 50x + 20y &\leq 2500 \\
                    & 60x + 90y &\leq 4000 \\
                    & x, y & \geq 0 
                    }Here $x$ and $y$ are written underneath $\max$ to indicate that these are the decision variables.
In cases where it is obvious what the decision variables are, we sometimes omit writing them below the $\max$ or $\min$ in the objective.
More complicated optimization problems might use letters to name things other than decision variables, and in those cases it is helpful to explicitly write down which variables the decision maker can affect.
The abbreviation ``s.t.'' stands for ``subject to'' or ``such that'' and indicates the constraints.

The above example has linear objective functions and constraints.
Optimization formulations can also have nonlinear functions as shown in the example below.
The objective in the above example is to maximize the profit.
When the decision maker is concerned with controlling costs rather than profit, the objective function often involves minimization.

\begin{exm}
A business has the option of setting up concrete plants at two locations.
The daily cost operating a plant at Location 1 per unit of production has two components: a fixed cost of 20 and an additional cost which increases by 0.5 for every unit of production.
Similarly, the daily cost of operating a plant at Location 2 per unit production has a fixed cost of 12 and an additional cost which increases by 1.2 for every unit of production.
The business has committed to supplying at least 12 units of concrete daily.
What is the total amount to be produced at each location to minimize cost and satisfy demand?
\end{exm}

There are two decision variables in this problem: the production in plant at Location 1 and production in plant at Location 2.
Let $x$ denote the units of production in plant at Location 1 and $y$ represent the units of production in plant at Location 2.

The business has agreed to supply at least 12 units of concrete daily.
Therefore, the sum of production in plants at Location 1 and Location 2 must be greater than equal to 12:
\begin{equation}
x + y \geq 12
\,.
\end{equation}
Also, the production at both plants cannot be negative which is represented as
\begin{equation}
x \geq 0 \qquad  y \geq 0
\,.
\end{equation}
The daily cost of operating a plant at Location 1 per unit of production is $20 + 0.5x$.
Since the plant produces $x$ units of concrete, the total cost of operating a plant at Location 1 is $(20+0.5x)x$.
Similarly, the total cost of operating a plant at Location 2 is $(12+1.2y)y$.
The goal of the business is to minimize the total cost of operation which is given as
\begin{equation}
\min \,\, (20+0.5x)x + (12+1.2y)y
\,.
\end{equation}
The optimization formulation can thus be summarized as 
\optimizex{\min}{x,y}{20x + 0.5x^2 + 12y + 1.2y^2}{
                     & x + y & \geq 12\\
                     & x, y & \geq 0
                      }Optimization models are widely used in resource allocation where a fixed resource has to be allocated among various tasks to maximize or minimize objectives.
A common variant of the resource allocation for civil engineers is the budget allocation problem.
In the budget allocation problem, a fixed budget has to be distributed among various projects.
The objective can vary depending on the nature of the decision maker.
For example, private entities might be interested in maximizing returns or profit whereas public entities might be interested in maximizing social welfare or equity.

\begin{exm}
A business is seeking to invest in two public-private partnership projects.
Project 1 yields an expected return of \$10 and a standard deviation of \$5 for each dollar invested.
Project 2 has an expected return of \$12 and a standard deviation of \$7 for each dollar invested.
For any investment portfolio, the business's utility is described as $U = \mu - \frac{1}{2}\sigma$, where $\mu$ represents the average return on investment and $\sigma$ the standard deviation of return.
The business has a total of \$1000 to invest.
Assume that the returns on investments on both projects are independent.
Determine the amount of money to be invested in both projects so that the business can maximize its utility?
\end{exm}

There are two decision variables in this problem: the amount of money invested in Project 1 ($x$), and the amount of money invested in Project 2 ($y$).
We know the total amount of money available is \$1000.
Therefore
\begin{equation}
x + y \leq 1000
\,.
\end{equation}
Also, the business cannot invest negative money in either project:
\begin{equation}
x \geq 0, y \geq 0
\,.
\end{equation}
Project 1 yields an average return of \$10 for each dollar invested.
Therefore, with an investment of \$$x$, the average return is \$$10x$.
Project 2 yields an average return of \$12 for each dollar invested.
Therefore, with an investment of \$$y$, the average return is \$$12y$.
The total average return on investment is given as $10x + 12y$.

The standard deviation of return on project 1 is \$5 for each dollar invested and for project 2 is \$7 for each dollar invested.
Therefore, the standard deviation of return on investments is given as $ \sqrt{25x^2 + 49y^2}$.

The business's utility is given as $U = \mu - \frac{1}{2}\sigma = (10x + 12y) -\frac{1}{2}\sqrt{25x^2 + 49y^2}$.
The optimization formulation can be summarized as follows:
\optimizex{\min}{x,y}{10x + 12y - 0.5\sqrt{25x^2 + 49y^2}}{
                    & x+y &\leq 1000 \\
                    & x, y & \geq 0
                    }

\section{Index Notation}

\index{index notation|(}
In the previous section, we used the symbols $x$ and $y$ to represent the decision variables.
Each symbol corresponds to a scalar.
In this section, we will develop an optimization formulation for a problem with more decision variables and constraints.
You will notice that using scalar based symbols will become cumbersome as the size of the formulation increases.
We will then introduce the set and index notations which will help represent the optimization formulation in a compact manner using more convenient notation.
First, let us consider the following problem.

\begin{exm}
A business has five factories located at Pittsburgh, Boston, Austin, Los Angeles, and Miami producing two types of window frames.
It costs \$10 and \$26 to produce one unit of window frame of Type 1 and Type 2 respectively.
The business has committed to supplying at least 300 units of window frame of Type 1 and at least 200 units of window frame of Type 2 each week.
Due to the number of workers employed at each location, the maximum number of window frames which can be produced at Pittsburgh, Boston, Austin, Los Angeles, and Miami are 100, 125, 100, 125, and 50 respectively.
The cost of producing one unit of window frame of Type 2 at Pittsburgh, Boston, Austin, Los Angeles, and Miami is \$40, \$40, \$15, \$20, and \$20 respectively.
Determine the number of window frames of Type 1 and Type 2 which the business has to produce at each location with the objective of minimizing cost?
\end{exm}

Let us first identify the decision variables.
In this problem, the business has to decide how many units of window frames of Type 1 and Type 2 to produce at the five locations: Pittsburgh, Boston, Austin, Los Angeles, and Miami.
Therefore, there are 10 decision variables.

Let $a$, $b$, $c$, $d$, and $e$ represent the number of units of window frames of Type 1 to be produced per week at Pittsburgh, Boston, Austin, Los Angeles, and Miami respectively.
Let $p$, $q$, $r$, $s$, and $t$ represent the number of units of window frames of Type 2 to be produced per week at Pittsburgh, Boston, Austin, Los Angeles, and Miami respectively.

The business has committed to supplying at least 300 units of window frames of Type 1 per week.
Therefore, the total amount of window frames of Type 1 produced per week must be greater than or equal to 300:
\begin{equation}
a + b + c + d + e \geq 300
\,.
\end{equation}
The business has committed to supplying at least 200 units of window frames of Type 2 per week.
Therefore, the total amount of window frames of Type 2 produced per week must be greater than or equal to 200:
\begin{equation}
p + q + r + s + t \geq 200
\,.
\end{equation}
The Pittsburgh factory can produce a maximum of 100 units of window frames per week.
Therefore, the sum of the total number of window frames of Type 1 produced per week at Pittsburgh and the total number of window frames of Type 2 produced per week at Pittsburgh must be less than or equal to 100:
\begin{equation}
a + p \leq 100
\,.
\end{equation}
The factory at Boston can produce a maximum of 125 units of window frames per week.
Therefore,
\begin{equation}
b + q \leq 125
\,.
\end{equation}
The factory at Austin can produce a maximum of 100 units of window frames per week.
Therefore,
\begin{equation}
c + r \leq 100
\,.
\end{equation}
Along similar lines, the factories at Los Angeles and Miami are constrained to produce a maximum of 125 and 50 window frames respectively:
\begin{align}
d + s &\leq 125 \\
e + t &\leq 50
\,.
\end{align}
The number of window frames of both types produced per week at all locations has to be greater than or equal to zero.
The non-negativity constraints are represented as 
\begin{equation}
a \geq 0, b \geq 0, c \geq 0, d \geq 0, e \geq 0
\,.
\end{equation}
\begin{equation}
p \geq 0, q \geq 0, r \geq 0, s \geq 0, t \geq 0
\,.
\end{equation}
The cost of producing $a$ units of window frames of Type 1 and $p$ units of window frames of Type 2 at Pittsburgh is $10a +40p$.
Along similar lines, the cost of producing specific number of window frames of both types at the other four locations can be determined.
The total production cost is the sum of production costs at each location which is given as $10a+10b+25c+30d+30e+40p+40q+15r+20s+20t$.
The objective is to minimize the total production costs.

The optimization formulation can thus be summarized as shown below
\optimizex{\min}{a,\ldots,e,p,\ldots,t}{10a+10b+25c+30d+30e +40p+40q+15r+20s+20t}{
&a + b + c + d + e &\geq 300 \\
&p + q + r + s + t &\geq 200 \\
&a + p &\leq 100 \\
&b + q &\leq 125 \\
&c + r &\leq 100 \\
&d + s &\leq 125 \\
&e + t &\leq 50 \\
&a,b,c,d,e &\geq 0\\
&p,q,r,s,t &\geq 0
}This formulation is correct, but a bit unwieldy.
To someone looking just at the formulation, it is hard to tell which variables refer to what.
Index notation, presented in the following subsections, can improve the presentation of formulations like this one.

\subsection{Single index notation}

This section introduces sets and indices, and explains how to use notation based on them to present formulations concisely.
First imagine that the business had factories at ten locations and was producing five types of window frames.
We would need 50 different symbols to denote the window frame types produced at all locations! The  formulation can be represented in a simpler way by using set and index notation.
This section first presents the single index notation.
A more compact multiple index notation is provided in the following section.

A set is a collection of similar objects.
Let $\mathcal{I}$\label{not:fancyI} denote the set of all locations.
For the above example, let
\begin{equation}
\mathcal{I} = \{ \text{Pittsburgh, Boston, Austin, Los Angeles, Miami} \} \,.
\end{equation} 
In the above example, we are given the names of locations of five factories.
In some cases, we may not know or we may not be interested in the names of the locations.
All we may know or care about is that there are five factories located in five different places.
In such case, the set $\mathcal{I}$ can be represented as
\begin{equation}
\mathcal{I} = \{ \text{1, 2, 3, 4, 5}
\} \,.
\end{equation}
The above notation can become cumbersome if the number of locations is high.
So a shorter form representation is 
\begin{equation}
\mathcal{I} = \{ \text{1, \ldots , 5}
\} \,.
\end{equation}
The above notation can be further generalized for any value $n$ denoting number of locations as
\begin{equation}
\mathcal{I} = \{ \text{1, \ldots , n} \,.
\}
\end{equation}
An index is used to refer to any element of the set.
The symbol $\in$ represents ``an element of''.
Therefore, $i \in \mathcal{I}$ denotes an index $i$ which is an element of the set $\mathcal{I}$.
Thus when $\mathcal{I} = \{ \text{Pittsburgh, Boston, Austin, Los Angeles, Miami} \}$, $i$ can refer to any of Pittsburgh, Boston, Austin, Los Angeles, or Miami.
Indices can be combined with symbols to represent decision variables and input parameters in a concise manner.

In the previous example, we used the symbols $a, b, c, d,$ and $e$ to represent the number of units of window frames  of Type 1 to be produced per week at Pittsburgh, Boston, Austin, Los Angeles, and Miami respectively.
Using the set-index notation, let $x_i$ denote the number of units of window frames per week of Type 1 to be produced at location $i \in \mathcal{I}$.
Thus $a$ corresponds to $x_{\text{Pittsburgh}}$, $b$ corresponds to $x_{\text{Boston}}$, and so on.

The business has to produce at least 300 units of window frames of Type 1 per week.
This constraint can be represented as
\begin{equation}
x_{\text{Pittsburgh}} + x_{\text{Boston}} + x_{\text{Austin}} + x_{\text{Los Angeles}} + x_{\text{Miami}} \geq 300
\,.
\end{equation}
The above constraint can be succinctly represented using the summation operator $\Sigma$
as
\begin{equation}
\sum_{i \in \mathcal{I}} x_i \geq 300 \,.
\end{equation}
If the set $\mathcal{I} = \{ \text{1,\ldots , 5} \}$, then the above constraint can also be represented using either of the following equations:
\begin{equation}
\sum_{i=1}^{5} x_i \geq 300 
\end{equation}
\begin{equation}
\sum_{1 \leq i \leq 5} x_i \geq 300 \,.
\end{equation}
Similarly, let $y_i$ denote the number of units of window frames per week of Type 2 to be produced at location $i \in \mathcal{I}$.
As before, the constraint on minimum number of units of Type 2 window frames produced can be written as any of the following:
\begin{equation}
\sum_{i \in \mathcal{I}} y_i \geq 200
\end{equation}
\begin{equation}
\sum_{i=1}^{5} y_i \geq 200
\end{equation}
\begin{equation}
\sum_{1 \leq i \leq 5} y_i \geq 200
\,.
\end{equation}
Now, look at the constraints which limits the number of window frames produced at each location.
Consider the case where the city names are written explicitly, so  $\mathcal{I} = \{ \text{Pittsburgh, Boston, Austin, Los Angeles, Miami} \}$.
The constraints are 
\begin{align}
x_{\text{Pittsburgh}} + y_{\text{Pittsburgh}} &\leq 100 \\
x_{\text{Boston}} +  y_{\text{Boston}} &\leq 125 \\
x_{\text{Austin}} + y_{\text{Austin}} &\leq 100 \\
x_{\text{Los Angeles}} + y_{\text{Los Angeles}} &\leq 125 \\
x_{\text{Miami}} + y_{\text{Miami}} &\leq 50
\,.
\end{align}
This will become very cumbersome to write as the number of locations and window frame types increase.
Let $u_i$ represent the maximum amount of window frame which can be produced at location $i \in \mathcal{I}$.
For example, $u_{\text{Pittsburgh}} = 100$, $u_{\text{Los Angeles}} = 125$.
The symbol $\forall$ represents ``for all''.
Therefore, $\forall i \in \mathcal{I}$ implies for all the elements in the set $\mathcal{I}$.
\emph{All} of these constraints can then be written as a single equation:
\begin{equation}
x_i + y_i \leq u_i \,\, \forall i \in \mathcal{I}
\,.
\end{equation}
The above equation denotes for each element $i \in \mathcal{I}$, the constraints $x_i + y_i \leq u_i$ holds.
Similarly the non-negativity constraints can be written as
\begin{align}
x_i \geq 0 \,\, & \forall i \in \mathcal{I}  \\
y_i \geq 0 \,\, & \forall i \in \mathcal{I}  
\,.
\end{align}
When $\mathcal{I} = \{ \text{1,\ldots , 5} \}$, the production limit at each location and non-negativity constraints can  also be represented using either of the two following set of equations:
\begin{align}
x_i + y_i \leq u_i \,\, & \forall i = 1, \ldots, 5 \\
x_i \geq 0 \,\, & \forall i = 1, \ldots, 5 \\
y_i \geq 0 \,\, & \forall i = 1, \ldots, 5  
\end{align}
or
\begin{align}
x_i + y_i \leq u_i \,\, & \forall  1 \leq i \leq 5 \\
x_i \geq 0 \,\, & \forall 1 \leq i \leq 5\\
y_i \geq 0 \,\, & \forall 1 \leq i \leq 5  
\,.
\end{align}
The objective function for the formulation was 
\begin{equation}
\min \,\, 10a+10b+25c+30d+30e+40p+40q+15r+20s+20t 
\,.
\end{equation}
Using the index-set notation, the objective function can be written
\begin{multline}
\min \,\, 10x_{\text{Pittsburgh}}+10x_{\text{Boston}}+25x_{\text{Austin}}+30x_{\text{Los Angeles}}+30x_{\text{Miami}} \\
+40y_{\text{Pittsburgh}}+40y_{\text{Boston}}+15y_{\text{Austin}}+20y_{\text{Los Angeles}}+20x_{\text{Miami}} \,.
\end{multline}

Let $v_i$ represent the cost of producing one window frame of Type 1 at location $i \in \mathcal{I}$ and $w_i$ represent the cost of producing one window frame of Type 2 at location $i \in \mathcal{I}$.
For example $v_{\text{Miami}} = 30, w_{\text{Miami}} = 20$.
The objective function can now be succinctly represented as
\begin{equation}
\min \,\, \sum_{i \in \mathcal{I}} (v_ix_i + w_iy_i)
\,.
\end{equation}
When $\mathcal{I} = \{ \text{1,\ldots , 5} \}$ the objective function can also be represented as
\begin{equation}
\min \,\, \sum_{i=1}^{5} (v_ix_i + w_iy_i)
\end{equation}
or
\begin{equation}
\min \,\, \sum_{1 \leq i \leq 5} (v_ix_i + w_iy_i)
\,.
\end{equation}
Therefore, the formulation can be rewritten as
\optimizex{\min}{x_i, y_i}{\sum_{i \in \mathcal{I}} (v_ix_i + w_iy_i)}{
& \sum_{i \in \mathcal{I}} x_i & \geq 300 \\
& \sum_{i \in \mathcal{I}} y_i & \geq 200 \\
& x_i + y_i & \leq u_i \,\, \forall i \in \mathcal{I} \\
& x_i &\geq 0 \,\,  \forall i \in \mathcal{I} \\
& y_i &\geq 0 \,\,  \forall i \in \mathcal{I} 
}or equivalently, with $\sum_{i=1}^5$ or $\sum_{1 \leq i \leq 5}$ on the summations.

In addition to making the formulation more compact, the set index notation also makes the formulation easier to understand, and easier to change (if there were more cities, all we would have to change is the definition of the set $\mc{I}$ or the number 5 to whatever the new number of cities is).
Notice also that now it is important to specify that $x_i$ and $y_i$ are the decision variables: $u_i$, $v_i$, and $w_i$ now represent \emph{given} problem data which we cannot change.

\subsection{Multiple index notation}

This section introduces a notation that further simplifies the example presented in the previous subsection .
Previously we had set
\[
\mathcal{I} = \{ \text{Pittsburgh, Boston, Austin, Los Angeles, Miami} \}
\,.
\]
or $\mathcal{I} = \{1, \ldots , 5\}$.
Two sets of decision variables $x_i,\, y_i$ were used to represent window frames of Type 1 and 2 produced at all locations $i \in \mathcal{I}$.
If the business was producing ten types of window frames, then even the single index notation becomes cumbersome, as we would need ten different subscripted symbols, $x_i, y_i, z_i, \ldots$ and so on.

This issue can be addressed using multiple indices.
Let us introduce the set $\mathcal{J}$\label{not:fancyJ} to model the two types of window frames.
Similar to the set of locations, $\mathcal{J}$ can be defined in two ways: $\mathcal{J} = \{ \text{Type 1, Type 2}\}$ or $\mathcal{J} = \{1,2\}$.
Let $j$ be an index referring to any element in $\mathcal{J}$, i.e., $j \in \mathcal{J}$.

Instead of using two symbols with subscripts for location, we will define the decision variable using one symbol with two subscripts.
Let $x_{ij}$ denote the number of type $j \in \mathcal{J}$ window frames produced per week at location $i \in \mathcal{I}$.

When the sets were $\mathcal{I} = \{ \text{Pittsburgh, Boston, Austin, Los Angeles, Miami} \}$ and $\mathcal{J} = \{ \text{Type 1, Type 2}\}$, in the single index notation two equations (one for Type 1, other for Type 2) were used to represent the total production must be greater than equal to demand:
\begin{align}
\sum_{i \in \mathcal{I}} x_i &\geq 300 \\
\sum_{i \in \mathcal{I}} y_i &\geq 200
\,.
\end{align}
Using the double index notation, the two equations can be summarized into one equation as shown below, by introducing $b_j$ to be the demand for window frames of type $j \in \mathcal{J}$.
For this example, $b_{\text{Type 1}} = 300$ and  $b_{\text{Type 2}} = 200$:
\begin{equation}
\sum_{i \in \mathcal{I}} x_{ij} \geq b_j \,\, \forall j \in \mathcal{J}
\,.
\end{equation}
Pay close attention to the two set element references, $\forall j \in \mathcal{J} $ on the right hand side and $i \in \mathcal{I}$ underneath the summation operator.
The $\forall j \in \mathcal{J} $ on the right hand side ensures that the equation $ \sum_{i \in \mathcal{I}} ()$ is repeated for each element $j \in \mathcal{J}$:
\begin{equation}
x_{\text{Pittsburgh,j}} + x_{\text{Boston,j}} + x_{\text{Austin,j}} + x_{\text{Los Angeles,j}} + x_{\text{Miami,j}} \geq b_j \,\, \forall j \in \mathcal{J}
\,.
\end{equation}
Since the set $\mathcal{J}$ has two elements, the equation is repeated twice representing the demand for Type 1 and Type 2 window frames.
In the single index notation $x_i + y_i \leq u_i \,\, \forall i \in \mathcal{I}$ is used to represent the constraint on maximum window frames which can be produced at each location.
In the double index notation, the left hand side can be made more compact using a summation operator:
\begin{equation}
\sum_{j \in \mathcal{J}} x_{ij} \leq u_i \,\, \forall i \in \mathcal{I}
\,.
\end{equation}
Pay close attention to the indices over which the summation ($j \in \mathcal{J}$) is happening and the elements on the right hand side ($i \in \mathcal{I}$).
The above expression repeats the following equation for each element or each location $i \in \mathcal{I}$:
\begin{equation}
x_{i,\text{Type 1}} +  x_{i,\text{Type 2}}  \leq u_i \,\, \forall i \in \mathcal{I}
\,.
\end{equation}
In the single index notation, two sets of equations were used to represent the non-negativity conditions - $ x_i \geq 0 \,\, \forall i \in \mathcal{I} $ and $ y_i \geq 0 \,\,  \forall i \in \mathcal{I}$.
This can be concisely represented using a single equation:
\begin{equation}
x_{ij} \geq 0 \,\, \forall i \in \mathcal{I}, j \in \mathcal{J} 
\,.
\end{equation}
You are enforcing $x_{ij} $ to be $\geq 0$ for each element or location $i \in \mathcal{I}$ as well as window frame type $ j \in \mathcal{J}$.
The order in which you reference the elements and sets on the right hand side after $\forall$ does not matter, i.e., the following two equations represent the same non-negativity constraints:
\begin{align}
x_{ij} \geq 0 \,\, &\forall i \in \mathcal{I}, j \in \mathcal{J} \\ 
x_{ij} \geq 0 \,\, &\forall j \in \mathcal{J}, i \in \mathcal{I}
\,.
\end{align}
With single index notation, the objective function was written $\min \,\, \sum_{i \in \mathcal{I}} (v_ix_i + w_iy_i) $.
In the double index notation the objective function can be represented as $ \min \,\, \sum_{i \in \mathcal{I}} (v_ix_{i, \text{Type 1}} + w_i x_{,:\text{Type 2}}) $.
We can make the representation even more compact and more intuitive by defining $v_{ij}$ as the cost of producing one window frame of type $j \in \mathcal{J}$ at location $i \in \mathcal{I}$.
The objective function can then be written
\begin{equation}
\min \,\,\sum_{i \in \mathcal{I}} \sum_{j \in \mathcal{J}} v_{ij}x_{ij}
\,.
\end{equation}
Note that the order in which we sum the objective function does not matter, i.e,
$\sum_{i \in \mathcal{I}} \sum_{j \in \mathcal{J}} v_{ij}x_{ij}  = \sum_{j \in \mathcal{J}} \sum_{i \in \mathcal{I}} v_{ij}x_{ij}$.
(See Section~\ref{sec:indicessummation}).

The formulation in double index notation is finally:
\optimizex{\min}{x_{ij}}{\sum_{i=1}^{5} \sum_{j=1}^{2} v_{ij}x_{ij}}
{
& \sum_{i \in \mathcal{I}} x_{ij} & \geq b_j \,\, \forall j \in \mathcal{J} \\
& \sum_{j \in \mathcal{J}} x_{ij} & \leq u_i \,\, \forall i \in \mathcal{I} \\
& x_{ij} & \geq 0 \,\, \forall i \in \mathcal{I}, j \in \mathcal{J} 
}

\subsection{Additional information on index notation}

This section provides a brief overview of some mathematical notation which is commonly seen in the optimization literature.

The set $\mathbb{R}$ contains all real numbers.
The set of all non-negative real numbers is commonly represented as $\mathbb{R_+}$.
Another way to represent the non-negativity constraint $x_{ij}  \geq 0 \,\, \forall i \in \mathcal{I}, j \in \mathcal{J}$ is:
\begin{equation}
x_{ij} \in \mathbb{R_+} \,\, \forall i \in \mathcal{I}, j \in \mathcal{J}
\,.
\end{equation}

In the above equation each decision variable $x_{ij}$ is restricted to lie in the set of non-negative real numbers.
If the decision variables can be positive or negative real numbers or zero, then the above equation can be modified as
\begin{equation}
x_{ij} \in \mathbb{R} \,\, \forall i \in \mathcal{I}, j \in \mathcal{J}
\,.
\end{equation}

The set $\mathbb{Z}$\label{not:bbz} contains all integers.\footnote{The use of the letter `Z' is from the German word ``Zahlen,'' meaning ``numbers.''} The set of all non-negative integers is commonly represented as $\mathbb{Z_+}$.\label{not:bbzp}
In certain type of optimization problems, called integer programs, the decision variables are restricted to be integers or non-negative integers which can be represented as follows:
\begin{align}
x_{ij} \in \mathbb{Z} & \,\, \forall i \in \mathcal{I}, j \in \mathcal{J} \\
x_{ij} \in \mathbb{Z_+} & \,\, \forall i \in \mathcal{I}, j \in \mathcal{J}
\,.
\end{align}

In some cases, there will be certain constraints defined on a subset of indices.
Let us assume in the above example, all the window frames produced on the east coast of the United States are first transported to a warehouse before delivery to retailers.
The warehouse can store only 150 window frames each week.
This restriction can be represented by the following constraint:
\begin{equation}
x_{\text{Boston,Type 1}} + x_{\text{Boston,Type 2}} + x_{\text{Miami,Type 1}} + x_{\text{Miami,Type 2}} \leq 150
\,.
\end{equation}
The above constraint can be more concisely represented by defining a subset $\mathcal{K(I)} = \{ i\in \mathcal{I}: i \text{ lies on East Coast}\}$.
In this expression, `$:$' means ``such that''.
The `$|$' symbol has the same meaning:.
$\mathcal{K(I)} = \{ i\in \mathcal{I}| i \text{ lies on East Coast}\}$.
Given the definition of the subset $\mathcal{K(I)}$, the constraint can be more succinctly presented as
\begin{equation}
\sum_{k \in \mathcal{K(I)}} \sum_{j \in \mathcal{J}} x_{kj} \leq 150
\,.
\end{equation}

The set notation can also be used to represent various mathematical expressions of decision variables in a clean manner.
Let $\mathcal{I} = \{1,\ldots,10\}$.
We want to represent the expression $x_2 + x_4 + x_6 + x_8 + x_{10}$ in a concise manner.
One way to do this would be to define $\mathcal{K(I)} = \{ i\in \mathcal{I}: i \text{ is even}\}$ and then write
\begin{equation}
x_2 + x_4 + x_6 + x_8 + x_{10} = \sum_{k \in \mathcal{K(I)}} x_k
\,.
\end{equation}
Another way to represent this would be to indicate the condition expression in the summation operator itself:
\begin{equation}
x_2 + x_4 + x_6 + x_8 + x_{10} = \sum_{ \{i \in \mathcal{I}: i \text{ is even}\} } x_i = \sum_{ \{i \in \mathcal{I}| i \text{ is even} \}} x_i 
\,.
\end{equation}
Along similar lines:
\begin{align}
x_3 + x_4 + x_5 &= \sum_{ \{i \in \mathcal{I}:  3 \leq i \leq 5 \} } x_i = \sum_{ \{ i \in \mathcal{I}| 3 \leq i \leq 5 \} } x_i \\
x_7 + x_8 + x_9 + x_{10} &= \sum_{ \{ i \in \mathcal{I}:  i \geq 7 \}} x_i = \sum_{ \{i \in \mathcal{I}| i \geq 7 \} } x_i \\
\end{align}
\index{index notation|)}

\section{Vector and Matrix Notations}
\label{sec:optvector}

\index{vector|(}
\index{matrix|(}
Vector and matrix notation, introduced in Section~\ref{sec:matrices}, is also very widely used to represent optimization formulations concisely.
This section explains these conventions.
To start, when we write equations or inequalities involving vectors, we mean that they apply to every component of those vectors.
For example, $\mb{x} \geq \mb{b}$, means that \emph{every} element in $\mb{x}$ is greater than or equal to the corresponding element in $\mb{b}$.
\begin{equation}
\mb{x} \geq \mb{b} \implies 
\begin{cases}
x_1 \geq b_1 \\
x_2 \geq b_2 \\
\vdots \\
x_n \geq b_n
\end{cases}
\,.
\end{equation}
Thus $\mb{x} \geq \mb{0}$ implies each element of $\mb{x}$ is greater than or equal to zero.
Similarly, 
\begin{equation}
\mb{Ax} = \mb{b} \implies 
\begin{cases}
a_{11}x_1 + a_{12}x_2 + \cdots + a_{1n}x_n  = b_1 \\
a_{21}x_1 + a_{22}x_2 + \cdots + a_{2n}x_n =  b_2 \\
\vdots \\
a_{m1}x_1 + a_{m2}x_2 + \cdots + a_{mn}x_n = b_n
\end{cases}
\,.
\end{equation}
Also,
\begin{equation}
\mb{b}^T\mb{x}= b_1x_1 + b_2x_2 + \ldots + b_nx_n = \sum_{i=1}^{n} b_ix_i
\end{equation}
is another way to express the dot product of two vectors $\mb{x}$ and $\mb{b}$.

Given the above information, we now reformulate the window frame optimization formulation using vectors and matrices.
Let $\mb{x} = (x_1,\ldots,x_{10} ) $ represent a column vector of decision variables, that is,
\begin{equation}
\mb{x} = \left[\begin{array}{c}
x_1 \\
x_2 \\
\ldots \\
x_{10}
\end{array} \right]
\,.
\end{equation}
In the vector $\mb{x}$, let $x_1,x_2,x_3,x_4,\text{ and } x_5$ represent the five decision variables corresponding to the amount of window frame of Type 1 produced at the five locations and $x_6,x_7,x_8,x_9,\text{ and } x_{10}$ represent the five decision variables corresponding to the amount of window frame of Type 2 produced at the five locations.
Let $\mb{c} = (c_1,\ldots,c_{10} ) $ represent a column vector of costs where $c_1, c_2, c_3, c_4,\text{ and } c_5$ represent the cost of producing one unit of Type 1 window frame at the five different locations and $c_6, c_7, c_8, c_9,\text{ and } c_{10}$ represent the cost of producing one unit of Type 2 window frame at the five different locations.

\begin{equation}
\mb{c} = \left[\begin{array}{c} c_1 \\ c_2 \\ c_3 \\ c_4 \\ c_5 \\ c_6 \\ c_7 \\ c_8 \\ c_9 \\ c_{10} \end{array} \right] = 
\left[\begin{array}{c} 10 \\ 10 \\ 25 \\ 30 \\ 30 \\ 40 \\ 40 \\ 15 \\ 20 \\ 20 \end{array} \right]
\,.
\end{equation}

The objective function in this case is
\begin{equation}
\min \,\, 10x_1 + 10x_2 + 25x_3 + 30x_4 + 30x_5 + 40x_6 + 40x_7 + 15x_8 + 20x_9 + 20x_{10}
\,.
\end{equation}
Let $\mb{c}^T$ represent the transpose of the vector $\mb{c}$, so $\mb{c}^T$ is a row vector.
The objective function can then be compactly represented as
\begin{equation}
\min \,\, \mb{c}^T \mb{x} = \sum_{i=1}^{10}c_ix_i
\,.
\end{equation}
The production of window frame of Type 1 must be greater than 300 and window frame of Type 2 must be higher than 200:
\begin{align}
x_1 + x_2 + x_3 + x_4 + x_5  \geq 300 \\
x_6 + x_7 + x_8 + x_9 + x_{10} \geq 200
\,.
\end{align}
Let us represent the above set of constraints in a compact manner.
Define a matrix $\mb{A}$ as
\begin{equation}
\mb{A} = \left[ \begin{array}{cccccccccc} 1&1&1&1&1&0&0&0&0&0 \\ 0&0&0&0&0&1&1&1&1&1     \end{array} \right]
\,.
\end{equation}
Let $\mb{b}$ be a column vector of demands:
\begin{equation}
\mb{b} = \left[ \begin{array}{c} b_1 \\ b_2 \end{array}\right] = \left[ \begin{array}{c} 300 \\ 200     \end{array} \right]
\,.
\end{equation}
The constraints can now be represented as
\begin{equation}
\mb{Ax} = \left[ \begin{array}{c} x_1 + x_2 + x_3 + x_4 + x_5 \\  x_6 + x_7 + x_8 + x_9 + x_{10}   \end{array} \right] \geq \left[ \begin{array}{c} 300 \\ 200     \end{array} \right] = \mb{b}
\,,
\end{equation}
or in a compact form, simply as
\begin{equation}
\mb{Ax} \geq \mb{b}
\,.
\end{equation}
The constraints on the number of window frames produced at each location is given as
\begin{align}
x_1 + x_6 &\leq 100 \\
x_2 + x_7 &\leq 125 \\
x_3 + x_8 &\leq 100 \\
x_4 + x_9 &\leq 125 \\
x_5 + x_{10} &\leq 50 
\,.
\end{align}
Define a column vector $\mb{u}$ and matrix $\mb{U}$ as
\begin{equation}
\mb{u} = \left[ \begin{array}{c} 100 \\ 125 \\ 100 \\ 125 \\ 50 \end{array}\right] 
\,;
\end{equation}
\begin{equation}
\mb{U} = \left[ \begin{array}{cccccccccc} 1&0&0&0&0&1&0&0&0&0 \\ 0&1&0&0&0&0&1&0&0&0
\\ 0&0&1&0&0&0&0&1&0&0 \\ 0&0&0&1&0&0&0&0&1&0 \\ 0&0&0&0&1&0&0&0&0&1 \end{array} \right]
\,.
\end{equation}
Thus the production at each location being lesser than the capacity can be written in a compact manner as
\begin{equation}
\mb{Ux} \leq \mb{u}
\,.
\end{equation}
The nonnegativity constraints can be represented as $\mb{x} \geq \mb{0}$.
Thus the final optimization formulation can be given as
\optimizex{\min}{\mb{x}}{\mb{c}^T\mb{x}}{
&\mb{Ax} &\geq \mb{b} \\
&\mb{Ux} &\leq \mb{u} \\
&\mb{x} &\geq \mb{0}
}

Another way to represent the non-negativity constraints is $\mb{x} \in \mathbb{R}_+^{10}$.
In the general case, when the decision variable vector $\mb{x}$ has $n$ elements, the nonnegativity constraints can be represented as $\mb{x} \in \mathbb{R}_+^n$.
If the decision variables are restricted to the set of positive integers, then $\mb{x} \in \mathbb{Z}_+^n$.
\index{vector|)}
\index{matrix|)}

\section{Examples of Basic Optimization Problems}
\label{sec:optimizationexamples}

This section provides several examples of optimization problems.
The first is a classic optimization problem known as the \emph{transportation problem}.

\begin{exm}
\label{exm:transportation}
\emph{(Transportation problem.)}
A timber company has three mills which produces wooden frames and five markets.
The three mills can produce 20, 40, and 30  units of wooden frames respectively on a daily basis.
The daily demands for the wooden frames at the five markets are 15, 30, 20, 15, and 10 respectively.
The cost to transport the wooden frames from the three mills to the five markets are shown below:
\begin{center}
\begin{tabular}{ |c|ccccc| } 
 \hline
 Mill Location/Market & 1 & 2 & 3 & 4 & 5\\ \hline
1 & 22 & 13 & 24 & 31 & 47 \\
2 & 35	& 24 & 11 &	27 & 38 \\
3 & 44	& 33 & 25 &	11 & 26 \\
 \hline
\end{tabular}
\end{center}
Formulate an optimization problem for satisfying demands at the markets while minimizing the transportation costs from the mills.
\end{exm}

Let $\mathcal{I}$ and $\mathcal{J}$ denote the set of mill and market locations, respectively.
The decision variables in this problem are $x_{ij}$, each of which represents the volume of wooden frames to be sent from location $i \in \mathcal{I}$ to location $j \in \mathcal{J}$.
For example, $x_{23}$ corresponds to the volume of wooden frames to be sent from Mill 2 to Market 3.

Let $u_i $ represent the amount of wooden frames which can be produced at mill $i \in \mathcal{I}$ and $d_j$ represent the amount of wooden frames needed at market $j \in \mathcal{J}$.
For example, $u_3 = 30$ and $d_4 =15$.

At each mill, the total volume of wooden frames transported must be less than or equal to the capacity of the mill.
For example, at Mill 2:
\begin{equation}
x_{21} + x_{22} + x_{23} + x_{24} + x_{25} \leq 40 
\,.
\end{equation}

This can be represented generally as
\begin{equation}
\sum_{j \in \mathcal{J}} x_{ij} \leq u_i \,\, \forall i  \in \mathcal{I} 
\,.
\end{equation}

Similarly, we must provide at least as many frames as are demanded.
The demand constraint can be represented as
\begin{equation}
\sum_{i \in \mathcal{I}} x_{ij} \geq d_j \,\, \forall j  \in \mathcal{J} 
\,.
\end{equation}

In addition, the volume of wooden frames transported between mill and market locations cannot be negative:
\begin{equation}
x_{ij} \geq 0 \,\, \forall i  \in \mathcal{I}, \forall j  \in \mathcal{J}  
\,.
\end{equation}

Let $c_{ij}$ represent the cost to transport a wooden frame from mill $i \in \mathcal{I}$ to location $j \in \mathcal{J}$.
For example $c_{24} = 27$.
The objective is to minimize the total transportation costs which is
\begin{equation}
\min \,\,\sum_{i \in \mathcal{I}} \sum_{j \in \mathcal{J}} c_{ij}x_{ij}
\end{equation}
The final formulation for the transportation problem can be summarized as
\optimizex{\min}{x_{ij}}{\sum_{i \in \mathcal{I}} \sum_{j \in \mathcal{J}} c_{ij}x_{ij}}
{
&\sum_{j \in \mathcal{J}} x_{ij} &\leq u_i \,\, \forall i  \in \mathcal{I} \\
&\sum_{i \in \mathcal{I}} x_{ij} &\geq d_i \,\, \forall j  \in \mathcal{J} \\
&x_{ij} &\geq 0 \,\, \forall i  \in \mathcal{I}, \forall j  \in \mathcal{J}  
}

This formulation is an example of a \emph{linear program}.
In a linear program, the objective functions and constraints are linear functions of the decision variables, and the solution can be any real number satisfying the constraints (it does not have to be an integer).
Notice that in the above example, we are not restricting the number of wooden frames to be transported between each pair of mill and market to be integer.
However, the solution to this specific linear program will always yield integer solutions as long as the input data is integer, although this is not true in general (see Section~\ref{sec:totallyunimodular}).

We now modify this problem by considering an additional factor.

\begin{exm}
The trucks used in transporting wooden frames from Mill 1 are very old, the frames may be damaged because of their poor suspension systems.
To prevent this, additional packing material is needed, and the amount depends on the destination market location.
The table below shows the amount of packing material per frame for each mill and market combination (Mills 2 and 3 have newer trucks that do not require special packaging.)  
\begin{center}
\begin{tabular}{ |c|ccccc| } 
 \hline
 Mill Location/Market & 1 & 2 & 3 & 4 & 5\\ \hline
1 & 3 & 7 & 3 & 1 & 0 \\
2 & 0	& 0 & 0 &	0 & 0 \\
3 & 0	& 0 & 0 &	0 & 0 \\
 \hline
\end{tabular}
\end{center}
If Mill 1 has 21 units of packing material available each day, formulate the problem of meeting demands while minimizing transportation costs.
\end{exm}

Let $r_{ij}$ represent the amount of packaging needed to transport a wooden frame from mill $i \in \mathcal{I}$ to market $j \in \mathcal{J}$.
For example, $r_{12} = 7$ and $r_{34} = 0$.
Let $R_i$ represent the total packing material available at mill $i$ each day, so $R_1 = 21$.
$R_2$ and $R_3$ can be assigned the value zero, because they do not need to have any packing material available.
The linear transportation problem formulation can be modified by adding an additional resource constraint as shown below:
\optimizex{\min}{x_{ij}}{\sum_{i \in \mathcal{I}} \sum_{j \in \mathcal{J}} c_{ij}x_{ij}}
{
& \sum_{j \in \mathcal{J}} x_{ij} &\leq u_i \,\, \forall i  \in \mathcal{I} \\
& \sum_{i \in \mathcal{I}} x_{ij} &\geq d_i \,\, \forall j  \in \mathcal{J} \\
& \sum_{j \in \mathcal{J}} r_{ij}x_{ij} &\leq R_i \,\, \forall i  \in \mathcal{I} \\
& x_{ij} &\geq 0 \,\, \forall i  \in \mathcal{I}, \forall j  \in \mathcal{J}  
}It turns out that adding the packing material constraint changes the optimization problem in such a way that the optimal solutions may not be integers.
If this condition is important, we must enforce it with an additional constraint, by replacing
\begin{equation}
x_{ij} \geq 0 \,\, \forall i  \in \mathcal{I}, \forall j  \in \mathcal{J} 
\end{equation}
with
\begin{equation}
x_{ij} \in \mathbb{Z_+} \,\, \forall i  \in \mathcal{I}, \forall j  \in \mathcal{J} 
\,.
\end{equation} 
The formulation now becomes an \emph{integer program}\index{optimization!integer}.
The methods needed to solve integer programs are different than the methods used to solve linear program, and integer programs are much harder to solve.
Sometimes, it may be adequate to solve the problem as a linear program and then convert its optimal solution to an integral one, say, by rounding --- if the values of the decision variables are in the hundreds or thousands, the effect of rounding is likely small.
However, for some integer programs this can lead to very poor solutions.
\index{integer program!see {optimization, integer}}

We now return to the original transportation problem formulation of Example~\ref{exm:transportation} without the packing material constraint.
In that example, the objective functions and constraints are all linear functions of the decision variables.
In many real world applications, it might not be possible to use linear functions to represent the objective function or constraints.
The next example modifies the transportation costs to reflect ``diseconomies of scale,'' where the unit cost of shipping increases with the quantity (perhaps the most efficient trucks are used first, but as more and more frames are shipped, you have to start using older and less fuel-efficient trucks).

\begin{exm} 
Assume now that the unit cost of transporting wooden frames between mill $i \in \mathcal{I}$ and market $j \in \mathcal{J}$ is  $c_{ij} + x_{ij}$.
(For example, if 10 wooden frames are being transported between mill $i \in \mathcal{I}$ and location $j \in \mathcal{J}$, then the cost of transporting each frame between the two locations is $c_{ij} + 10$, and the total transportation cost would be $(c_{ij} + 10) \times 10$).
Formulate an optimization problem to meet the demands at the markets while minimizing transportation costs.
\end{exm}

In this formulation, if $x_{ij}$ wooden frames are being transported then the total transportation cost between mill $i \in \mathcal{I}$ and market $j \in \mathcal{J}$ is $(c_{ij} + x_{ij}) \times x_{ij}$.
The constraints are the same, but the objective function now changes:
\optimizex{\min}{x_{ij}}{\sum_{i \in \mathcal{I}} \sum_{j \in \mathcal{J}} (c_{ij} + x_{ij})x_{ij}}
{
& \sum_{j \in \mathcal{J}} x_{ij} &\leq u_i \,\, \forall i  \in \mathcal{I} \\
& \sum_{i \in \mathcal{I}} x_{ij} &\geq d_i \,\, \forall j  \in \mathcal{J} \\
& x_{ij} &\geq 0 \,\, \forall i  \in \mathcal{I}, \forall j  \in \mathcal{J}  
}

In this formulation, the objective function is no longer linear in the decision variables.
The above formulation is an example of a \emph{nonlinear program}\index{optimization!nonlinear}.
In a nonlinear program, either the objective function or any of the constraints is nonlinear.
A nonlinear program would also be needed to represent ``economies of scale,'' where the unit shipping cost decreases with the amount sent.
\index{nonlinear program!see {optimization, nonlinear}}

If we wanted to enforce the decision variables to be integers, we would need to replace
\begin{equation}
x_{ij} \geq 0 \,\, \forall i  \in \mathcal{I}, \forall j  \in \mathcal{J} 
\end{equation}
with
\begin{equation}
x_{ij} \in \mathbb{Z_+} \,\, \forall i  \in \mathcal{I}, \forall j  \in \mathcal{J} 
\,.
\end{equation} 
in the final formulation.
This leads to a \emph{nonlinear integer program}, which again requires different solution methods.
Integer programs are also used in modeling selections or incorporating ``yes or no'' decisions.
In such problems often a decision variable which can only take the values 0 (for ``no'') or 1 (for ``yes'') is used to model the selection decisions.
This is called a \emph{binary variable}\index{binary variable}.

\begin{exm} 
The timber company is now interested in establishing factories to produce wooden doors.
The potential sites for establishing the factories are the same as the current locations of the timber mills.
There is a fixed cost (which has been amortized to a daily cost) of 300, 350, and 450 associated with establishing the factories at the three locations.
Factories built at each site would respectively have daily production capacities of 40, 80, and 60 doors.
The daily demands for the wooden doors at the five markets are 35, 15, 50, 20, and 40 respectively.
The unit cost of transporting wooden frames is assumed to be the same as the unit cost of transporting wooden frames which are given in the linear transportation example.
Formulate an optimization problem to meet the demands at the markets while minimizing the facility location \emph{and} transportation costs.
\end{exm}

As in Example~\ref{exm:transportation}, let $\mathcal{I}$ and $\mathcal{J}$ represent the set of potential factory locations and markets respectively.
There are two sets of decision variables.
The first set of decision variables $y_i$ takes the value 1 if a facility is located at $i \in \mathcal{I}$ and 0 otherwise.
For example the values $y_1 = 1$, $y_2=0$, and $y_3 =1$ would mean that factories are opened at sites 1 and 3.
The second set of decision variables $x_{ij}$ represents the volume of demand at market $j \in \mathcal{J}$ served by a facility at location $i \in \mathcal{I}$.

Let $d_j$ represent the demand for wooden doors at market $j \in \mathcal{J}$.
Therefore $d_1  = 35, d_2 = 15$ and so on.
Let $u_i$ represent the capacity of the factory at location $i \in \mathcal{I}$.
Note that $u_1 = 40, u_2 = 80,  u_3 = 60$.

If a factory is located at $i \in \mathcal{I}$, then the total volume of wooden doors supplied to all markets cannot exceed the production capacity of the factory.
If a factory is not located at $i \in \mathcal{I}$, then the total volume of total volume of wooden doors supplied to all markets must be zero.
This constraint can be represented as
\begin{equation}
\sum_{j \in \mathcal{J}} x_{ij} \leq u_iy_i \,\, \forall i  \in \mathcal{I}  \,.
\end{equation}
Note that when $y_i = 1$, the right hand side becomes $u_i$ (total production at an open factory cannot exceed its capacity).
When $y_i = 0$ the right hand side becomes zero (nothing can be produced at a factory which was never opened).
Similar to the transportation problem, the total volume supplied to each market location must meet the demand:
\begin{equation}
\sum_{i \in \mathcal{I}} x_{ij} \geq d_i \,\, \forall j  \in \mathcal{J} \,.
\end{equation}
In addition to the non-negativity constraints on the variable $x_{ij}$, there is also an additional binary restriction on the facility location decision variable $y_i$:
\begin{align}
x_{ij} \geq 0 &\,\, \forall i  \in \mathcal{I}, \forall j  \in \mathcal{J}   \\
y_i \in \{0,1\} &\,\, \forall i  \in \mathcal{I}
\,.
\end{align}
The objective function has two cost components: the cost of establishing the facility and the transportation costs.
The expression for the transportation costs is the same as in the linear transportation problem example.
Let $f_i$ represent the cost of establishing a facility at location $ i \in \mathcal{I}$.
The total facility location cost can be given as $ \sum_{i \in \mathcal{I}} f_iy_i$.
Therefore the objective function is
\begin{equation}
\min \,\, \sum_{i \in \mathcal{I}} f_iy_i + \sum_{i \in \mathcal{I}} \sum_{j \in \mathcal{J}} (c_{ij} + x_{ij})x_{ij}
\,.
\end{equation}
The final formulation is
\optimizex{\min}{x_{ij}, y_i}{\sum_{i \in \mathcal{I}} f_i y_i + \sum_{i \in \mathcal{I}} \sum_{j \in \mathcal{J}} (c_{ij} + x_{ij})x_{ij}}
{
& \sum_{j \in \mathcal{J}} x_{ij} &\leq u_iy_i \,\, \forall i  \in \mathcal{I} \\
& \sum_{i \in \mathcal{I}} x_{ij} &\geq d_i \,\, \forall j  \in \mathcal{J} \\
& x_{ij} &\geq 0 \,\, \forall i  \in \mathcal{I}, \forall j  \in \mathcal{J}  \\
& y_i &\in \{0,1\} \,\, \forall i  \in \mathcal{I}
}In this formulation, the facility location decision variable is restricted to be either 0 or 1 and cannot be any real number.
Therefore, the above formulation is an example of an integer program.
Binary variables are also common in project or portfolio selection problems, which appear often in civil engineering.
The next example shows such a problem.

\begin{exm}
A business is seeking to invest in two private-public partnership projects.
There is an option investing \$400, \$200, \$300, \$150, and \$250 in projects 1 through 5 respectively.
The expected and standard deviation of return on the five projects is shown in the table below.

\begin{center}
\begin{tabular}{ |c|c|c| } 
 \hline
 Project & Expected Value & Stdev \\ \hline
1 & 5000 & 300  \\
2 & 2000 & 150  \\
3 & 3500 & 125 \\
4 & 2000 & 225 \\
5 & 3750 & 75 \\ \hline
\end{tabular}
\end{center}
The returns on the five projects have a correlation coefficient of $+0.1$.
For any investment portfolio, the business balances risk and expected return with the utility function $U = \mu - \frac{1}{2}\sigma$ where $\mu$\label{not:mumean} is the average return on investment and $\sigma$\label{not:sigmarv} the standard deviation of return.
The business has a total of \$1000 to invest.
Formulate an optimization problem for determining the investment projects to maximize its utility.
\end{exm}

Let $\mathcal{I}$ represent the set of projects.
In this case, let $\mathcal{I} = \{1,2,3,4,5\}$.

The decision variables in this problems should reflect which projects are selected for investment.
This can be done using binary decision variables $x_i$ which take the value 1 if project $i \in \mathcal{I}$ is selected, and 0 if not.

The business has a total budget of \$1000.
Let $B$ represent the total budget and $b_i$ represent the amount which needs to be invested in each project $i \in \mathcal{I}$.
For example $b_4 = 150$ and $b_5 = 250$.
The budget constraint can be enforced as
\begin{equation}
\sum_{i \in \mathcal{I}} b_ix_i \leq B
\,.
\end{equation}
Now let us tackle the objective function.
Let $\mu_i$ represent the average return of investment in project $i \in \mathcal{I}$.
The average return on investment $\mu$ is given as
\begin{equation}
\mu = 5000x_1 + 2000x_2 + 3500x_3 + 2000x_4 + 3750x_5 =\sum_{i \in \mathcal{I}} \mu_ix_i 
\,.
\end{equation}
Given the standard deviation of returns $\sigma_i$ for each $i \in \mathcal{I}$ and the correlation coefficient $\rho_{ij} = 0.1$ for all $i \in \mathcal{I}$ and $j \in \mathcal{I}$.
The standard deviation of return is given as
\begin{multline}
\sigma = \sqrt{300^2x_1 + \ldots + 225^2x_4 + 75^2 x_5 +  0.1 \times 300 \times 150 + \ldots} \\ =\sqrt{ \sum_{i \in \mathcal{I}} \sigma_i^2x_i + \sum_{i \in \mathcal{I}} \sum_{j \in \mathcal{I}: j \neq i} \rho_{ij} \sigma_i \sigma_j x_i x_j }
\,.
\end{multline}

The optimization formulation is
\optimizex{\min}{x_i}{\sum_{i \in \mathcal{I}} \mu_ix_i  - \frac{1}{2} \sqrt{ \sum_{i \in \mathcal{I}} \sigma_i^2 x_i + \sum_{i \in \mathcal{I}} \sum_{j \in \mathcal{I}: j \neq i} \rho_{ij} \sigma_i \sigma_j x_i x_j }}
{
& \sum_{i \in \mathcal{I}} b_ix_i &\leq B \\
& x_i &= \{0,1\} \,\, \forall i  \in \mathcal{I}
}The above formulation is a nonlinear integer program.
The nonlinearity arises from two sources: both the square root and the $x_i x_j$ interaction terms in the objective function.
The integrality arises from the binary decision variables.

\section{More Examples of Optimization Formulations}
\label{sec:moreoptexamples}

The previous section was a tour of common types of optimization problems: linear programs, integer programs, and nonlinear programs.
This section provides additional examples of how optimization problems might be formulated in a transportation engineering setting.
As you read through these, think about how they relate to the concepts defined in the previous section (which ones are linear, integer, etc.).

\begin{exm}
\label{exm:transitfrequency}
\emph{(Transit frequency setting.)}  You are working for a public transit agency in a city, and must decide the frequency of service on each of the bus routes.
The bus routes are known and cannot change, but you can change how the city's bus fleet is allocated to each of these routes.
(The more buses assigned to a route, the higher the frequency of service.)  Knowing the ridership on each route, how should buses be allocated to routes to minimize the total waiting time?
\end{exm}
\solution{
To formulate this as an optimization problem, we need to identify an objective, decision variables, and constraints.
For this problem, as we do this we will be faced with other assumptions which must be made.
With this and other such problems, there may be more than one way to write down an optimization problem, what matters is that we \emph{clearly state all of the assumptions made}.
Another good guiding principle is to \emph{start with the simplest model which captures the important behavior}, which can then be refined by relaxing assumptions or replacing simple assumptions with more realistic ones.

We need some notation, so let $R$ be the set of bus routes.
We know the current ridership on each route $d_r$.
Here we will make our first assumption: that the demand on each route is inelastic and will not change based on the service frequency.
Like all assumptions, it is not entirely true but in some cases may be close enough to the truth to get useful results; if not, you should think about what you would need to replace this assumption, which is always a good exercise.
We are also given the current fleet size (which we will denote $N$), and must choose the number of buses associated with each route (call this $n_r$) --- thus the decision variables in this problem are the $n_r$ values.

The objective is to minimize the total waiting time, which is the sum of the waiting time for the passengers on each route.
How long must passengers wait for a bus?  If we assume that travelers arrive at a uniform rate, then the average waiting time will be half of the service headway.
How is the headway related to the number of buses on the route $n_r$?  Assuming that the buses are evenly dispersed throughout the time period we are modeling, and assuming that each bus is always in use, then the headway on route $r$ will be the time required to traverse this route ($T_r$) divided by the number of buses $n_r$ assigned to this route.
So, the average delay per passenger is half of the headway, or $T_r / (2 n_r)$, and the total passenger delay on this route is $(d_r T_r) / (2 n_r)$.
This leads us to the objective function
\labeleqn{busobjective}{D(\mb{n}) = \sum_{r \in R} \frac{d_r T_r}{2 n_r}\,,}
in which the total delay is calculated by summing the delay associated with each route.

What constraints do we have?  Surely we must run at least one bus on each route (or else we would essentially be canceling a route), and in reality as a matter of policy there may be some lower limit on the number of buses assigned to each route; for route $r$, call this lower bound $L_r$.
Likewise, there is some upper bound $U_r$ on the number of buses assigned to each route as well.
So, we can introduce the constraint $L_r \leq n_r \leq U_r$ for each route $r$.

Putting all of these together, we have the optimization problem
\optimize{
 \min_{\mb{n}} & D(\mb{n}) = \sum_{r \in R} \frac{d_r T_r}{2 n_r}  &  \\
 \mathrm{s.t.} & n_r \geq L_r & \forall r \in R \\
               & n_r \leq U_r & \forall r \in R \\
               & \sum_{r \in R} n_r \leq N                        \,.
}
}

\begin{exm}
\label{exm:maintenance}
\emph{(Scheduling maintenance.)}  You are responsible for scheduling routine maintenance on a set of transportation facilities (such as pavement sections or bridges.)  The state of these facilities can be described by a condition index which ranges from 0 to 100.
Each facility deteriorates at a known, constant rate (which may differ between facilities).
If you perform maintenance during a given year, its condition will improve by a known amount.
Given an annual budget for your agency, when and where you should perform maintenance to maximize the average condition of these facilities?  You have a 10 year planning horizon.
\end{exm}

\solution{
In contrast to the previous example, where the three components of the optimization problem were described independently, from here on problems will be formulated in a more organic way, describing a model built from the ground up.
(This is how optimization models are usually described in practice.)  After describing the model in this way, we will identify the objective function, decision variables, and constraints to write the optimization problem in the usual form.
We start by introducing notation based on the problem statement.

Let $F$ be the set of facilities, and let $c_f^t$ be the condition of facility $f$ at the \textbf{end}\footnote{This word is intentionally emphasized.
In this kind of problem it is very easy to get confused about what occurs at the start of period $t$, at the end of period $t$, during the middle of period $t$, etc.} of year $t$, where $t$ ranges from 1 to 10.
Let $d_f$ be the annual deterioration on facility $f$, and $i_f$ the amount by which the condition will improve if maintenance is performed.
So, if no maintenance is performed during year $t$, then
\labeleqn{nomaintenance}{c_f^t = c_f^{t-1} - d_f \qquad \forall f \in F, t \in \{ 1, 2, \ldots, 10\}\,.}
and if maintenance is performed we have
\labeleqn{maintenance}{c_f^t = c_f^{t-1} - d_f + i_f \qquad \forall f \in F, t \in \{ 1, 2, \ldots, 10\}\,.}
Both of these cases can be captured in one equation with the following trick: let $x_f^t$ equal one if maintenance is performed on facility $f$ during year $t$, and 0 if not.
Then 
\labeleqn{stateevolutionalmost}{c_f^t = c_f^{t-1} - d_f + x_f^t i_f \qquad \forall f \in F, t \in \{ 1, 2, \ldots, 10\}\,.}
Finally, the condition can never exceed 100 or fall below 0, so the full equation for the evolution of the state is
\labeleqn{stateevolution}{c_f^t = \begin{cases} 100 & \mbox{ if } c_f^{t-1} - d_f + x_f^t i_f > 100 \\ 0 & \mbox { if } c_f^{t-1} - d_f + x_f^t i_f < 0 \\ c_f^{t-1} - d_f + x_f^t i_f  & \mbox{otherwise}\,. \end{cases}}
for all $f \in F$ and $t \in \{ 1, 2, \ldots, 10\}$.
Of course, for this to be usable we need to know the current conditions of the facilities, $c_f^0$.

The annual budget can be represented this way: let $k_f$ be the cost of performing maintenance on facility $f$, and $B^t$ the budget available in year $t$.
Then
\labeleqn{budget}{\sum_{f \in F} k_f x_f^t \leq B^t \qquad \forall t \in \{1, \ldots 10\}\,.}  For the objective, we need the average condition across all facilities and all years; this is simply $\frac{1}{10 |F|} \sum_{f \in F} \sum_{t=1}^{10} c_f^t$.
The obvious decision variables are the maintenance variables $x_f^t$, but we also have to include $c_f^t$ because these are influenced by the maintenance variables.
As constraints, we need to include the state evolution equations~\eqn{stateevolution}, the budget constraints~\eqn{budget}, and, less obviously the requirement that $x_f^t$ be either 0 or 1.
Putting it all together, we have the optimization problem
   \optimize{
         \max_{\mb{x,c}} & \frac{1}{10 |F|} \sum_{f \in F} \sum_{t = 1}^{10} c_f^t  &  \\
         \mathrm{s.t.} & \sum_{f \in F} k_f x_f^t \leq B^t & \forall t \in \{1, \ldots 10\} \\
                       & c_f^t \mbox{ is given by~\eqn{stateevolution} } & \forall f \in F, t \in \{ 1, 2, \ldots, 10 \}    \\
                       & x_f^t \in \{ 0, 1 \} & \forall f \in F, t \in \{1, 2, \ldots, 10 \}\,.
   }
}

In this example, pay close attention to the use of formulas like~\eqn{budget}, which show up very frequently in optimization.
It is important to make sure that every ``index'' variable in the formula is accounted for in some way.
Equation~\eqn{budget} involves the variables $x_f^t$, but for which facilities $f$ and time periods $t$?  The facility index $f$ is summed over, while the time index $t$ is shown at right as $\forall t \in \{1, \ldots 10 \}$.
This means that a copy of~\eqn{budget} exists for each time period, and in each of these copies, the left-hand side involves a sum over all facilities at that time.
Therefore the one line~\eqn{budget} actually includes 10 constraints, one for each year.
It is common to forget to include things like $\forall t \in \{1, \ldots 10 \}$, or to try to use an index of summation outside of the sum as well (e.g., an expression like $B_f - \sum_{f \in F} k_f x_f^t$), which is meaningless.
Make sure that all of your indices are properly accounted for!

In this next example, the objective function is less obvious.

\begin{exm}
\label{exm:facilitylocation}
\emph{(A facility location problem.)}  In a city with a grid network, you need to decide where to locate three bus terminals.
Building the terminal at different locations costs a different amount of money.
Knowing the home locations of customers throughout the city who want to use the bus service, where should the terminals be located to minimize the construction cost and walking distance customers have to walk?  Assume that each customer will walk from their home location to the nearest terminal.
\end{exm}
\solution{This problem could easily become very complicated if we take into account the impact of terminal locations on bus routes, so let's focus on what the problem is asking for: simply locating terminals to minimize walking distance from customers' home locations.

\stevefig{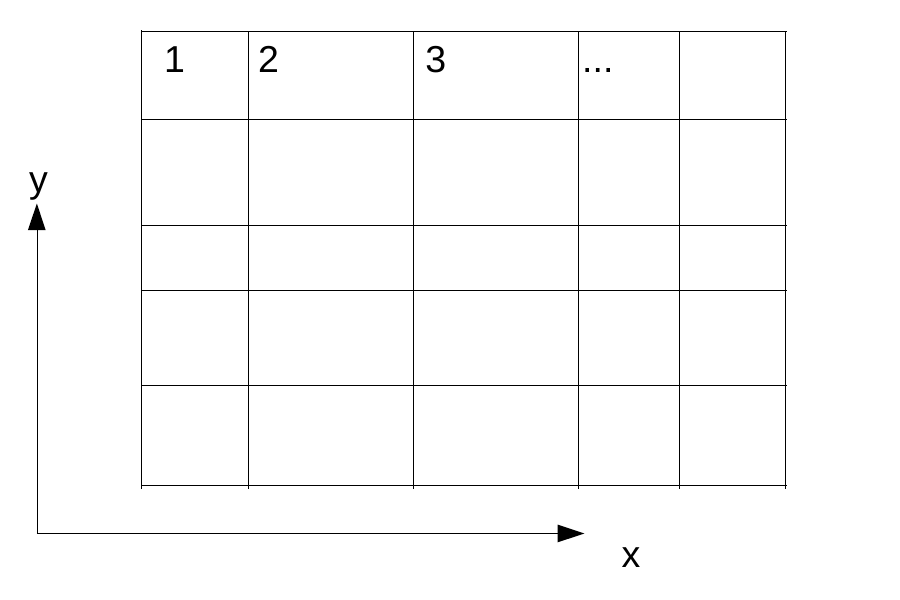}{Coordinates and labeling of intersections for 
Example~\ref{exm:facilitylocation}}{0.6\textwidth}

Number each of the intersections in the grid network (Figure~\ref{fig:gridnetwork}) from 1 to $I$, the total number of intersections.
Assume that the terminals will be located at these intersections, and let the variables $L_1$, $L_2$, and $L_3$ denote the numbers of the intersections where terminals will be built.
Let $C(i)$ be the cost of building a terminal at location $i$, so the total cost of construction is $C(L_1) + C(L_2) + C(L_3)$.
Let $P$ be the set of customers, and let $H_p$ denote the intersection that is the home location of customer $p$.

How can we calculate the walking distance between two intersections (say, $i$ and $j$)?  Figure~\ref{fig:gridnetwork} shows a coordinate system superimposed on the grid.
Let $x(i)$ and $y(i)$ be the coordinates of intersection $i$ in this system.
Then the walking distance between points $i$ and $j$ is
\labeleqn{manhattan}{d(i,j) = |x(i) - x(j)| + |y(i) - y(j)|\,.}
This is often called the \emph{Manhattan distance}\index{Manhattan distance} between two points, after one of the densest grid networks in the world.

So what is the walking distance $D(p)$ for customer $p$?  The distance from $p$ to the first terminal is $d(H_p, L_1)$, to the second terminal is $d(H_p, L_2)$, and to the third is $d(H_p, L_3)$.
The passenger will walk to whichever is closest, so $D(p,L_1,L_2,L_3) = \min \{ d(H_p, L_1), d(H_p, L_2), d(H_p, L_3) \}$ and the total walking distance is $\sum_{p \in P} D(p,L_1,L_2,L_3)$.

For this problem, the decision variables and constraints are straightforward: the only decision variables are $L_1$, $L_2$, and $L_3$ and the only constraint is that these need to be integers between 1 and $I$.
The tricky part is the objective function: we are instructed both to minimize total cost as well as total walking distance.
We have equations for each of these, but we can only have one objective function.
In these cases, it is common to form a \emph{convex combination} of the two objectives, introducing a weighting parameter $\lambda \in [0, 1]$.
That is, let
\labeleqn{objective}{f(L_1, L_2, L_3) = \lambda [C(L_1) + C(L_2) + C(L_3)] + (1 - \lambda) \mys{\sum_{p \in P} D(p,L_1,L_2,L_3)}\,. }

Look at what happens as $\lambda$ varies.
If $\lambda = 1$, then the objective function reduces to simply minimizing the cost of construction.
If $\lambda = 0$, the objective function is simply minimizing the total walking distance.
For a value in between 0 and 1, the objective function is a weighted combination of these two objectives, where $\lambda$ indicates how important the cost of construction is relative to the walking distance.

For concreteness, the optimization problem is
\optimize{
     \min_{L_1, L_2, L_3} & \lambda [C(L_1) + C(L_2) + C(L_3)] + (1 - \lambda) \mys{\sum_{p \in P} D(p,L_1,L_2,L_3)}  &  \\
     \mathrm{s.t.} & L_f \in \{ 1, 2, \ldots, I \} \qquad \qquad \qquad \qquad \forall f \in \{1, 2, 3 \} \,.&
}
}

In the final example in this section, finding a mathematical representation of a solution is more challenging.

\stevefig{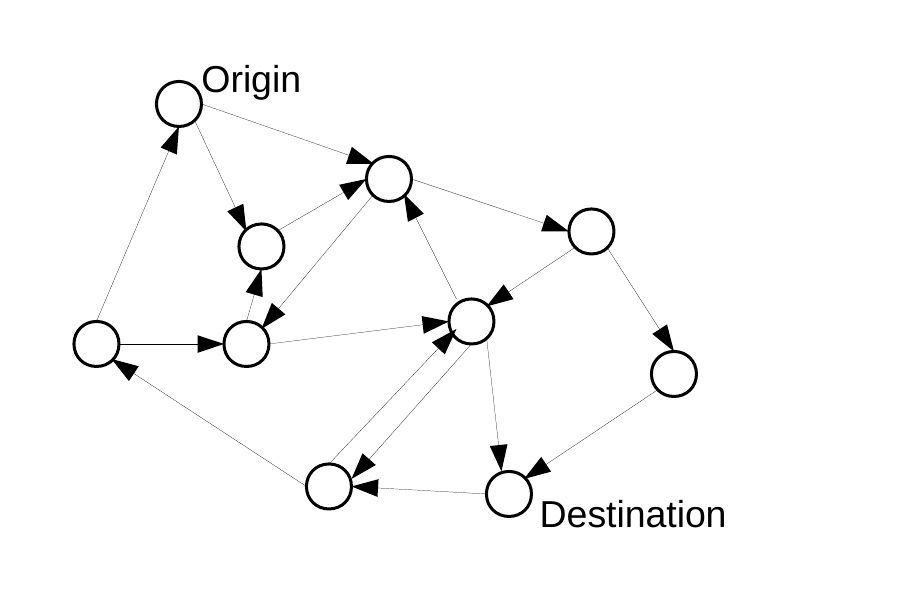}{Roadway network for Example~\ref{exm:shortestpath}}{0.6\textwidth}

\begin{exm}
\index{shortest path!optimization formulation|(}
\label{exm:shortestpath}
\emph{(Shortest path problem.)}  Figure~\ref{fig:spnetwork} shows a road network, with the origin and destination marked.
Given the travel time on each roadway link, what is the fastest route connecting the origin to the destination?
\end{exm}
\solution{
We've already presented some algorithms for solving this problem in Section~\ref{sec:shortestpath}, but here we show how it can be placed into the general framework of optimization problems.

Notation first: let the set of nodes be $N$, and let $r$ and $s$ represent the origin and destination nodes.
Let the set of links be $A$, and let $t_{ij}$ be the travel time on link $(i,j)$.
So far, so good, but how do we represent a route connecting two intersections?

Following Example~\ref{exm:maintenance}, introduce binary variables $x_{ij} \in \{ 0, 1 \}$, where $x_{ij} = 1$ if link $(i,j)$ is part of the route, and $x_{ij} = 0$ if link $(i,j)$ is not part of the route.
The travel time of a route is simply the sums of the travel times of the links in the route, which is $\sum_{(i,j) \in A} t_{ij} x_{ij}$.

We now have an objective function and decision variables, but what of the constraints?  Besides the trivial constraint $x_{ij} \in \{ 0, 1 \}$, we need constraints which require that the $x_{ij}$ values actually form a contiguous path which starts at the origin $r$ and ends at the destination $s$.
We do this by introducing a \emph{flow conservation constraint} at each intersection.
For node $i$, we will use $\Gamma(i)$ denotes the set of links which leave node $i$, and $\Gamma^{-1}(i)$ denotes the set of links which enter node $i$.
(See Section~\ref{sec:terminology} for more on this notation.)

Consider any contiguous path connecting intersection $r$ and $s$, and examine any node $i$.
One of four cases must hold: 
\begin{description}
\item[Case I: ] Node $i$ does not lie on the path at all.
Then all of the $x_{ij}$ values should be zero for $(i,j) \in \Gamma(i) \cup \Gamma^{-1}(i)$.
\item[Case II: ] Node $i$ lies on the path, but is not the origin or destination.
Then $x_{ij} = 1$ for exactly one $(i,j) \in \Gamma(i)$, and for exactly one $(i,j) \in \Gamma^{-1}(i)$.
\item[Case III: ] Node $i$ is the origin $r$.
Then all $x_{ij}$ values should be zero for $(i,j) \in \Gamma^{-1}(i)$, and $x_{ij} = 1$ for exactly one $(i,j) \in \Gamma(i)$.
\item[Case IV: ] Node $i$ is the destination $s$.
Then all $x_{ij}$ values should be zero for $(i,j) \in \Gamma(i)$, and $x_{ij} = 1$ for exactly one $(i,j) \in \Gamma^{-1}(i)$.
\end{description}
To combine these cases in an elegant way, look at the differences $\sum_{(i,j) \in \Gamma(i)} x_{ij} - \sum_{(h,i) \in \Gamma^{-1}(i)} x_{hi}$.
For cases I and II, this difference will be 0; for case III, the difference will be $+1$, and for case IV, the difference will be $-1$.

So, this leads us to the optimization formulation of the shortest path problem:
\labeleqn{fc}{
\begin{array}{rll} \ds
     \min_{\mb{x}} & \sum_{(i,j) \in A} t_{ij} x_{ij} \qquad&  \\
     \mathrm{s.t.} & \sum_{(i,j) \in \Gamma(i)} x_{ij} - \sum_{(h,i) \in \Gamma^{-1}(i)} x_{hi} =   \begin{cases}
                                                                          1 & \mbox{ if } i = r \\
                                                                          -1 & \mbox{ if } i = s \\
                                                                          0 & \mbox{ otherwise}
                                                                       \end{cases}
                                                                       & \forall i \in N  \\
                    & x_{ij} \in \{ 0, 1 \} & \forall (i,j) \in A\,.
\end{array}
}}

A careful reader will notice that if the four cases are satisfied for a solution, then the equations~\eqn{fc} are satisfied, but the reverse may not be true.
Can you see why, and is that a problem?
\index{shortest path!optimization formulation|)}

There is often more than one way to formulate a problem: for instance, we might choose to minimize congestion by spending money on capacity improvements, subject to a budget constraint.
Or, we might try to minimize the amount of money spent, subject to a maximum acceptable limit on congestion.
Choosing the ``correct'' formulation in this case may be based on which of the two constraints is harder to adjust (is the budget constraint more fixed, or the upper limit on congestion?).
Still, you may be troubled by this seeming imprecision.
This is one way in which modeling is ``as much art as science,'' which is not surprising --- since the human and political decision-making processes optimization tries to formalize, along with the inherent value judgments, (what truly is the objective?) are not as precise as they seem on the surface.
One hallmark of a mature practitioner of optimization is a flexibility with different formulations of the same underlying problem, and a willingness to engage in ``back-and-forth'' with the decision maker as together you identify the best formulation for a particular scenario.
In fact, in some cases it may not matter.
The theory of duality\index{optimization:duality} (which is beyond the scope of this book) shows that these alternate formulations sometimes lead to the same ultimate decision, which is comforting.

\section{General Properties of Optimization Problems}

After seeing some concrete examples of optimization formulations, we now define a general form for optimization problems.
This form is useful because it allows us to talk about properties of optimization problems, or methods for solving them, in a way that is independent of any specific context (shipping, investment, routing, etc.).
Any of the optimization problems discussed in this book can be represented in this form, perhaps after some conversions or transformation.

The general form is
\optimizex{\min}{\mb{x}}{f(\mb{x})}
{
&	 g_i(\mb{x}) \ge b_i& \forall i=1,2,..,m	\\
&	 h_i(\mb{x}) = 0	& \forall  i=1,2,..,l 	 
}If any variables are restricted to be integers, those constraints are listed in addition to the ones given above.

In the above problem $\mb{x} = (x_1,x_2,...,x_n)$ is a vector of decision variables.
The function $f(\mb{x})$ is the \emph{objective function} which represents the cost or utility associated with a decision $\mb{x}$.
The equations $g_i(\mb{x}) \ge b_i$ and $h_i(\mb{x}) = 0$ are \emph{constraint functions}\index{optimization!constraint function} which describe the constraints mathematically.
Let $X=\{\mb{x} \in \mathbb{R}^n: g_i(\mb{x}) \ge b_i\,\,\forall i=1,2,..m ; h_i(\mb{x}) = 0\,\,	 \forall  i=1,2,..l  \}$.
The set $X$ is called the \emph{feasible region}\index{optimization!feasible region} which represents the set of values which the decision variables can take or the set of values which satisfies all constraint functions.
Any value $\mb{x} \in X$ is a \emph{feasible solution}\index{optimization!feasible set} to the optimization problem.
The goal of the mathematical problem is to determine the \emph{optimal solution}\index{optimization!optimal solution} which is defined as the feasible solution with the smallest objective function value, i.e, to find $\mb{x^*} \in X$ such that $f(\mb{x^*}) \le f(\mb{x}) \,\, \forall \mb{x} \in X$.

Depending on the mathematical properties of the objective and constraint functions there are different categories of optimization problems.
Linear optimization problems are mathematical programs where the functions $ f(\mb{x})$,  $g_i(\mb{x})$, $h_i(\mb{x})$ are linear.
If any of the objective or constraint functions are nonlinear, then we have nonlinear optimization problems.
Integer programs are mathematical programs  where the decision variables have restricted to a set of  integers.
Depending on the nature of the mathematical program, there are different types of solution algorithms with varying degrees of efficiency.

\subsection{Local and global optima} 
\label{sec:localglobal}

\begin{figure}
  \includegraphics[width=\linewidth]{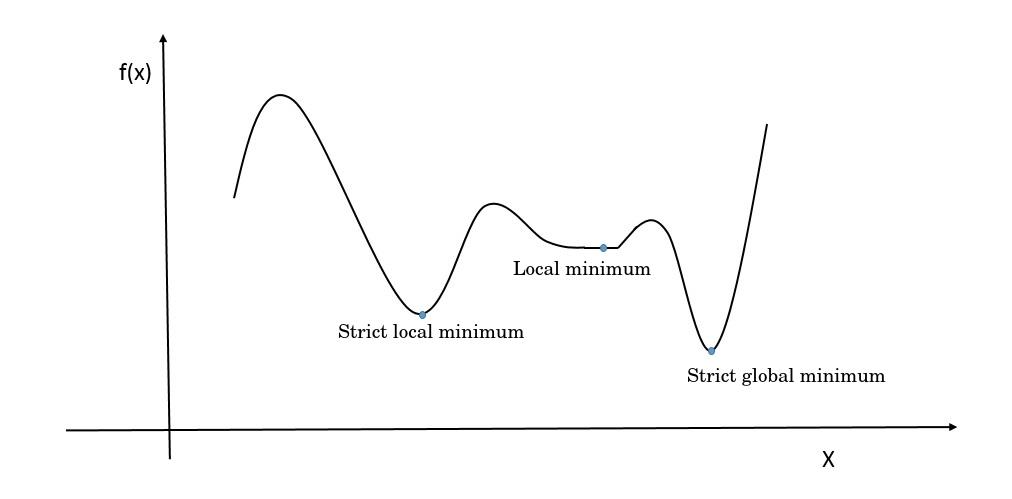}
  \caption{Local, strict local, and global minima.}
  \label{fig:localglobal}
\end{figure}

Objective functions can have both \emph{local} and \emph{global} maximum or minimum solutions (Figure~\ref{fig:localglobal}).
Let $f$ be a function defined on a region $X \subseteq \bbr^n$. 
Recall from Section~\ref{sec:sets} that a ball centered on $\mb{x}$ of radius $\epsilon$ is the set containing all points of $\bbr$ whose distance to $\mb{x}$ is less than $\epsilon$.

\begin{dfn}
A point $\mb{x^*} \in X$ is a \emph{local minimum}\index{optimization!local minimum} of a function $f(\mb{x})$ if there is a ball $B$ centered on $\mb{x}$ with positive radius such that
\begin{equation*}
f(\mb{x}) \geq f(\mb{x^*}) \,\, \forall \mb{x} \in B \cap X
\,.
\end{equation*}
\end{dfn}

\begin{dfn}
A point $\mb{x^*} \in X$ is a \emph{global minimum}\index{optimization!global minimum} of a function $f(\mb{x})$ if 
\begin{equation*}
f(\mb{x}) \geq f(\mb{x^*}) \,\, \forall \mb{x} \in X
\,.
\end{equation*}
\end{dfn}

Strict local and global minima are obtained by replacing the inequalities in the above definitions with strict inequalities:

\begin{dfn}
A point $\mb{x^*} \in X$ is a \emph{strict local minimum}\index{optimization!strict local minimum} of a function $f(\mb{x})$ if there is a ball $B$ centered on $\mb{x}$ with positive radius such that
\begin{equation*}
f(\mb{x}) > f(\mb{x^*}) \,\, \forall \mb{x} \in B \cap X, \mb{x} \neq \mb{x^*}
\,.
\end{equation*}
\end{dfn}

\begin{dfn}
A point $\mb{x^*} \in X$ is a \emph{strict global minimum}\index{optimization!strict global minimum} of a function $f(\mb{x})$ if 
\begin{equation*}
f(\mb{x}) > f(\mb{x^*}) \,\, \forall \mb{x} \in X, \mb{x} \neq \mb{x^*}
\,.
\end{equation*}
\end{dfn}

The corresponding terms for maximization problems are defined in analogous ways, by reversing the signs of the inequalities.

\subsection{Objective function transformations}
\label{sec:objectivetransform}

\index{optimization!objective function|(}
The following results are very useful in showing cases when two superficially different optimization problems may in fact be the same.
The first result shows that it is easy to convert any maximization problem to a minimization problem (or vice versa) by negating the objective function.
The second shows that constants may be freely added or subtracted to objective functions without changing the optimal solutions.
The third shows that multiplication or division by a positive constant does not change the optimal solutions.

\begin{prp} 
\label{prp:maxmin}
\index{optimization!maximization}
If the feasible set is $X$, the solution $\mb{\hat{x}}$ is a global maximum of $f$ if and only if $\mb{\hat{x}}$ is a global minimum of $-f$ for the same feasible set $X$.
\end{prp}
\begin{proof}
For the first part, assume that $\mb{\hat{x}}$ is the global maximum of $f$ on the set $X$.
Then $f(\mb{\hat{x}}) \geq f(\mb{x})$ for all $x \in X$.
Multiplying both sides by $-1$ reverses the sign of the inequality, giving $-f(\mb{\hat{x}}) \leq -f(\mb{x})$ for all $X$, so $\mb{\hat{x}}$ is a global minimum for $-f$.
For the second part, assume that $\mb{\hat{x}}$ is the global minimum of $-f$, so $-f(\mb{\hat{x}}) \leq -f(\mb{x})$ for all $x \in X$.
Again multiplying both sides by $-1$ gives $f(\mb{\hat{x}}) \geq f(\mb{x})$ for all $x \in X$, so $\mb{\hat{x}}$ is a global maximum of $f$.
\end{proof}

This result is convenient, because we do not need to develop two different theories for maximization and minimization problems.
Instead, we can focus just on one.
In this book, we develop results for minimization problems, a fairly common convention in engineering.\footnote{Economists often use maximization problems as their standard convention.}  Whenever you encounter a maximization problem, you can convert it to a minimization problem by negating the objective function, then using the minimization results and procedures.
(Of course, this choice is completely arbitrary; we could just as well have chosen to develop results only for maximization problems and in fact some other fields do just this.)

\begin{prp}
\label{prp:optimizationconstants}
Let $\mb{\hat{x}}$ be an optimal solution when the objective function is $f(\mb{x})$ and the feasible set is $X$.
Then for any constant $b$ not depending on $\mb{x}$, $\mb{\hat{x}}$ is also an optimal solution when the objective function is $f(\mb{x}) + b$ and the feasible set is $X$.

\end{prp}
\begin{proof}
By Proposition~\ref{prp:maxmin}, with no loss of generality we can assume that the optimization problem is a minimization problem and that $\mb{\hat{x}}$ is a global minimum.
Therefore $f(\mb{\hat{x}}) \leq f(\mb{x})$ for all $x \in X$.
We can add the constant $b$ to both sides of the equation, so $f(\mb{\hat{x}}) + b \leq f(\mb{x}) + b$ for all $x \in X$.
Therefore $\mb{\hat{x}}$ is also a global minimum (and thus optimal) when the objective function is $f(\mb{x}) + b$.

\end{proof}

\begin{prp}
\label{prp:constmultiple}
Let $\mb{\hat{x}}$ be an optimal solution when the objective function is $f(\mb{x})$ and the feasible set is $X$.
Then for any constant $c > 0$ not depending on $\mb{x}$, $\mb{\hat{x}}$ is also an optimal solution when the objective function is $cf(\mb{x})$ and the feasible set is $X$.
\end{prp}
\begin{proof}
See Exercise~\ref{ex:constmultiple}.
\end{proof}
\index{optimization!objective function|)}

\subsection{Existence of optima}

It is possible to formulate optimization problems that do not have solutions --- either because the constraints are contradictory and there is no feasible solution satisfying all of them, or because of the way the objective function and constraints are defined.
Optimization problems of the first type are called \emph{infeasible}\index{optimization!infeasible}.
For example, the problem ``minimize $x$ subject to $x \leq 0$'' has no optimal solution --- whatever $x$ you choose, you can find an $x$ which is even more negative.
Optimization problems of this type, where the objective can be made as negative as you like, are called \emph{unbounded}\index{optimization!unbounded}.
Some problems lack solutions even without falling into either of these categories; for example, ``minimize $x$ subject to $x > 0$'' has no optimal solution, for whatever $x$ you choose, $x/2$ is feasible but smaller.
This is why we tend to avoid strict inequalities in optimization problems.
But there are still problematic optimization formulations even without strict inequalities, such as ``minimize $1/x$ subject to $x \geq 1$.''

Practically speaking, if you are in a situation like the ones described above, you should take another look at your optimization problem: perhaps one constraint can be relaxed (maybe penalized in the objective, rather than strictly excluding solutions), or perhaps you missed a constraint that should have been there.
Systems with unbounded or truly infeasible solutions are rare.
This subsection provides a mathematical perspective on the topic, giving conditions on the constraints and objective function which can guarantee existence of a solution.

Keeping with our standard convention, we consider a problem of the form ``minimize $f(\mb{x})$ subject to $\mb{x} \in X$,'' where the set $X$ contains all solutions satisfying all constraints.
Weierstrass' theorem identifies gives sufficient conditions, under which an optimal solution exists. 

\begin{thm}
\label{thm:weierstrass}
\emph{(Weierstrass' theorem.)}
\index{Weierstrass' theorem}
\index{continuous function!applications}
\index{compact set!applications}
Let $f$ be a continuous, real-valued function defined on $X$, and let $X$ be a non-empty, compact set.
The optimization problem $\{ \min f(x): x \in X \}$ has a minimum solution.
\end{thm}

(For definitions of the terms in this theorem, see Appendix~\ref{chp:mathbackground}.)

\begin{exm}
\label{ex:1}
Consider the optimization problem  $\min (x-5)^2, x \in X = \{x: x \in \mathbb{R}, 3 \leq x \leq 7 \}$. Does a minimum exist? 
\end{exm}
\solution{
The function $(x-5)^2$ is continuous and real valued.
The feasible region is non-empty, closed, and bounded.
Therefore, the above optimization problem has a minimum.
By plotting the function, we can see that the minimum occurs at $x=5$.}

\begin{exm}
\label{ex:2}
Consider the optimization problem  $\min (x-5)^2, x \in X = \{x: x \in \mathbb{R}, 5 < x \leq 7 \}$. Does a minimum exist? 
\end{exm}
\solution{In the above variant, the feasible region is not closed.
Given any $y \in X$, we will be able to find another $z \in X$, such that $f(z) < f(y)$.
For example, if we pick $y =  5.001$, we can always find another $z = 5.0001$ whose objective function value is smaller.
Therefore, no minimum exists in this case. 
}

Note that Weierstrass' theorem provides sufficient and not necessary conditions (more on this in next subsection).
The continuity, closed, bounded assumption are not required for minima to exist.
We can have situations where minima exist without the above conditions being satisfied.

\begin{exm}
Consider the optimization problem  $\min (x-5)^2, x \in X = \{x: x \in \mathbb{R}, x > 3 \}$. Does a minimum exist? 
\end{exm}
\solution{In the above optimization problem, the feasible region is not closed or bounded. However, we do know that the minimum exists at $x=5$.}

\begin{exm}
Consider the objective function
\begin{equation*}
f(x) = \begin{cases}
          1 & $x = 0$ \\
          5 & \mbox{otherwise} 
      \end{cases}
\,.      
\end{equation*}
The feasible region is $X = \{x: x \in \mathbb{R}\}$. Does a minimum exist? 
\end{exm}
\solution{
In this example, the function is discontinuous, the feasible region is not closed or bounded. However, we still have a minimum at $x=0$.}

\subsection{Necessary and sufficient conditions} 
\label{sec:necessarysufficient}

\index{necessary condition|(}
\index{sufficient condition|(}
The previous subsection presented Weierstrass' theorem, which gave sufficient conditions for the existence of a minimum solution.
Some readers may be unfamiliar with the distinction between ``necessary'' and ``sufficient'' conditions, so we explain them here in the context of this theorem.
Consider two statements $A$ and $B$. 
Saying that $A$ is a \emph{necessary} condition for $B$ means that $B$ is true only if $A$ is true, written $B \implies A$ (read ``$B$ implies $A$'').
For example, every square is also a rectangle.
This means that a shape can be a square only if it is a rectangle; being a rectangle is necessary for being a square, so we can identify $A$ with ``the shape is a rectangle,'' and $B$ with ``the shape is a square.''
Notice that $A$ being true does not imply that $B$ is true; there are rectangles which are not squares.

In the above situation, we can also say that $B$ is a \emph{sufficient} condition for $A$, meaning that if $B$ is true, $A$ must also be true.
For example, ``being a square'' is sufficient for ``being a rectangle.''
In a statement like $B \implies A$, the arrow points from the sufficient condition to the necessary one.
When $B$ is a sufficient condition for $A$, $A$ can be true without $B$ being true.

Used alone, the terms ``necessary'' or ``sufficient'' mean that $A$ and $B$ are not exactly the same thing, but that one includes cases that the other does not have.
For example, there are some rectangles which are not squares.
If $A$ and $B$ \emph{do} mean exactly the same thing, we say that $A$ is necessary \emph{and} sufficient for $B$; $A$ is true if and only if (abbreviated ``iff'') $B$ is true.
This is written $A \iff B$, and means that both $A \implies B$ and $B \implies A$ are true.
An example are the statements ``$n$ is an even integer'' and ``$n = 2k$ for some integer $k$.''  
Each statement is true exactly when the other is true; in other words, they are \emph{equivalent} and we can replace one statement by the other whenever we like.

In Weierstrass' theorem, statement $A$ is ``the function $f$ is continuous and real-valued, and $X$ is non-empty, closed, and bounded,'' and statement $B$ is ``the function has a minimum value on $X$.''
Condition $A$ is sufficient, because any time $A$ is true, $B$ is true as well.
However, condition $A$ is not necessary, because there are examples where $A$ is false but $B$ is still true (Examples~\ref{ex:1} and~\ref{ex:2} above).
\index{necessary condition|)}
\index{sufficient condition|)}

\section{Exercises}
\label{sec:optimizationexercises}

\begin{enumerate}
\item \diff{11} Why is it generally a bad idea to use strict inequalities (`$>$', `$<$') in mathematical programs?  
\item \diff{12} Mathematically specify an optimization problem which has no optimal solution, and an optimization problem which has multiple optimal solutions.
\item \diff{24} For the following functions, identify all stationary points and global optima.
\begin{enumerate}[(a)]
 \item $f(x) = x^4 - 2x^3$
 \item $f(x_1, x_2) = 2x_1^2 + x_2^2 - x_1 x_2 - 7x_2$
 \item $f(x_1, x_2) = x_1^6 + x_2^6 - 3x_1^2 x_2^2$
 \item $f(x_1, x_2, x_3) = (x_1 - 4)^2 + (x_2 - 2)^2 + x_3^2 + x_1 x_2 + x_1 x_3 + x_2 x_3$
\end{enumerate}
\item \diff{42} Express the constraint ``$x \in \mathbb{Z}$'' (that is, $x$ must be an integer) in the form $g(x) = 0$ for some continuous function $g$.
(This shows that nonlinear programming can't be any easier than integer programming!) 
\item \diff{30} Can you ever improve the value of the objective function at the optimal solution by adding constraints to a problem?
If yes, give an example.
If no, explain why not.
\item \diff{21} Prove Proposition~\ref{prp:constmultiple}.
\label{ex:constmultiple}
\item \diff{27} Write out the objective and all of the constraints for the maintenance scheduling problem of Example~\ref{exm:maintenance}, using the bridge data shown in Table~\ref{tbl:maintenanceexample}, where the maintenance cost is expressed in millions of dollars.
Assume a two-year time horizon and an annual budget of \$5 million.
\label{ex:maintenanceexample}
\begin{table}
   \centering
   \caption{Data for Exercise~\ref{ex:maintenanceexample} \label{tbl:maintenanceexample}}
   \small
   \begin{tabular}{|c|cccc|}
      \hline
      Bridge  & Initial condition      & Deterioration & Repairability & Maintenance cost  \\
      $f$     & $c_f^0$                & $d_f$              & $i_f$           & $k_f$ \\
      \hline
      1        & 85                    & 5                  & 3               & 2 \\
      2        & 95                    & 10                 & 1               & 4 \\
      3        & 75                    & 1                  & 1               & 3 \\
      4        & 80                    & 6                  & 6               & 1 \\
      5        & 100                   & 3                  & 3               & 1 \\
      \hline
   \end{tabular}
\end{table}
\item \diff{44} Reformulate the transit frequency-setting problem of Example~\ref{exm:transitfrequency} so that the objective is to minimize the total amount of money spent, while achieving a pre-specified level of service (maximum delays).
\item \diff{53} You are responsible for allocating bridge maintenance funding for a state, and must develop optimization models to assist with this for the current year.
Political realities and a strict budget will require you to develop multiple formulations for the purposes of comparison.

Let $\mathcal{B}$ denote the set of bridges in the state.
There is an ``economic value'' associated with the condition of each bridge: the higher the condition, the higher the value to the state, because bridges in worse condition require more routine maintenance, impose higher costs on drivers who must drive slower or put up with potholes, and carry a higher risk of unforeseen failures requiring emergency maintenance.
To represent this, there is a value function $V_b(x_b)$ associated with each bridge $b \in \mathcal{B}$, giving the economic value of this bridge if $x_b$ is spent on maintenance this year.
You have a total budget of $X$ to spend statewide.
The state is also divided into $n$ districts, and each bridge belongs to one district.
Let $\mathcal{D}_i$ denote the set of bridges which belong to district $i$; different districts may contain different numbers of bridges.

In the following, be sure to define any additional notation you introduce.
There is more than one possible formulation that achieves the stated goals, so feel free to explain any parts of your formulations which may not be self-evident.
\begin{enumerate}[(a)]
   \item \textbf{Maximize average benefits: } Formulate an optimization model (objective function, decision variables, and constraints) to maximize the average economic value of bridges across the state.
   \item \textbf{Equitable allocation: } The previous model may recommend spending much more money in some districts than others, which is politically infeasible.
Formulate an optimization model based on maximizing the total economic value to the state, but ensuring that each district receives a relatively fair share of the total budget $X$.
   \item \textbf{Equitable benefits: } There is a difference between a fair allocation of money, and a fair allocation of \emph{benefits} (say, if the functions $V_b$ differ greatly between districts).
Defining \emph{benefit} as the difference between the economic value after investing in maintenance, and the ``do-nothing'' economic value, formulate an optimization model which aims to achieve a relatively fair distribution of benefits among districts.
\end{enumerate}
\item \diff{34}  Consider the mathematical program
\begin{align*}
\max \quad        & 3x_1 + 5x_2 + 6x_3  \\
\mbox{s.t } \quad &  x_1 + x_2 + x_3  = 1 \\
            & 2x_1 + 3x_2 +  2x_3  \leq 5 \\
            &  x_2 + x_3  \leq 1/2 \\
            & x_1, x_2, x_3 \geq  0 
\end{align*}
\begin{enumerate}
\item Are all of the constraints needed?
\item Solve this problem graphically.
\end{enumerate}
\item \diff{57}  In an effort to fight rising maintenance costs and congestion, a state transportation agency is considering a toll on a certain freeway segment during the peak period.
Imposing a toll accomplishes two objectives at once: it raises money for the state, and also reduces congestion by making some people switch to less congested routes, travel earlier or later, carpool, take the bus, and so forth.
Suppose that there are 10000 people who would want to drive on the freeway if there was no toll and no congestion, but the number of people who actually do is given by $$x = 10000e^{(15 -6\tau - t)/500}$$ where $\tau$ is the roadway toll (in dollars) and $t$ is the travel time (in minutes).
(That is, the higher the toll, or the higher the travel time, the fewer people drive.)  The travel time, in turn, is given by
$$t = 15\left[1 + 0.15\left(\frac{x}{c}\right)^4\right]$$
minutes, where $c = 8000$ veh/hr is the roadway capacity during rush hour.
Regulations prohibit tolls exceeding \$10 at any time.
Citizens are unhappy with both congestion and having to pay tolls.
After conducting multiple surveys, the agency has determined that citizen satisfaction can be quantified as 
$$s = 100 - t - \tau / 5$$
   \begin{enumerate}[(a)]
   \item Formulate two nonlinear programs based on this information, where the objectives are to either (a) maximize revenue or (b) maximize citizen satisfaction.
   \item Solve both of these problems, reporting the optimal  values for the decision variables and the objective function. 
   \item Comment on the two solutions.
(e.g., do they give a similar toll value, or very different ones?  How much revenue does the state give up in order to maximize satisfaction?)
   \item Name at least two assumptions that went into formulating this problem.
Do you think they are realistic?
        Pick one of them, and explain how the problem might be changed to eliminate that assumption and make it more realistic.
   \end{enumerate}
\item \diff{59} You are asked to design a traffic signal timing at the intersection of 8th \& Grand.
Assuming a simple two-phase cycle (where Grand Avenue moves in phase 1, and 8th Street in phase 2), no lost time when the signal changes, and ignoring turning movements, the total delay at the intersection can be written as
$$\frac{\lambda_1 (c - g_1)^2}{2c\left(1 - \frac{\lambda_1}{\mu_1}\right)} + \frac{\lambda_2 (c - g_2)^2}{2c\left(1 - \frac{\lambda_2}{\mu_2}\right)}$$
where $g_1$ and $g_2$ are the effective green time allotted to Grand Avenue and 8th Street, $c = g_1 + g_2$ is the cycle length, $\lambda_1$ and $\mu_1$ are the arrival rate and saturation flow for Grand Avenue, and $\lambda_2$ and $\mu_2$ are the arrival rate and saturation flow for 8th Street.

All signals downtown are given a sixty-second cycle length to foster good progression, and the arrival rate and saturation flow are 2200 veh/hr and 3600 veh/hr for Grand Avenue, respectively, and 300 veh/hr and 1900 veh/hr for 8th Street.
Furthermore, no queues can remain at the end of the green interval; this means that $\mu_i g_i$ must be at least as large as $\lambda_i c$ for each approach $i$.
   \begin{enumerate}[(a)]
   \item Why does the constraint $\mu_i g_i \geq \lambda_i c$ imply that no queues will remain after a green interval?
   \item Formulate a nonlinear program to minimize total delay.
   \item Simplify the nonlinear program so there is only one decision variable, and solve this nonlinear program using the bisection method of Section~\ref{sec:bisection} (terminate when $b - a \leq 1$).

   \item Write code to automate this process, and perform a sensitivity analysis by plotting the effective green time on Grand Ave as $\lambda_1$ varies from 500 veh/hr to 3000 veh/hr, in increments of 500 veh/hr.
Interpret your plot.
   \item Identify two assumptions in the above model.
Pick one of them, and describe how you would change your model to relax that assumption.
   \end{enumerate}
   
\item \diff{26} Find the global minima of the following functions in two ways: the bisection method, and using Newton's method to directly find a stationary point (if the Newton step leaves the feasible region, move to the boundary point closest to where Newton's method would go).
Run each method for five iterations, and see which is closer to the optimum, making the comparison based on the value of the objective function at the final points.
\begin{enumerate}[(a)]
 \item $f(x) = -\arctan x$, $x \in [0, 10]$ 
 \item $f(x) = x \sin (1 / (100x))$, $x \in [0.015, 0.04]$ 
 \item $f(x) = x^3$, $x \in [5,15]$
\end{enumerate}

\item \diff{51} Give a network and a feasible solution $\mathbf{x}$ to the mathematical formulation of the shortest path problem in Example~\ref{exm:shortestpath} where the links with $x_{ij} = 1$ do not form a contiguous path between $r$ and $s$, as alluded to at the end of the example.

\end{enumerate}

\chapter{Optimization Techniques}
\label{chp:fancyoptimization}

This appendix contains additional information on optimization problems, and specific techniques for solving them.
In turn, we present solution methods for one-dimensional optimization, linear optimization problems, unconstrained nonlinear optimization, constrained nonlinear optimization, and optimization problems with integrality constraints.
We close with an introduction to metaheuristics, general-purpose techniques which can be applied to any optimization problem, but without guarantees of finding an optimal solution.

\section{More Line Search Algorithms}
\label{sec:linesearch}

\index{optimization!line search|(}
\index{convex function!applications}
Section~\ref{sec:bisection} presented the bisection method as a way of solving optimization problems with a single decision variable$x$; a convex, differentiable objective function $f$; and a single, ``interval'' constraint of the form $a \leq x \leq b$.
(Actually, we can relax the assumption that $f$ is convex to $f$ being \emph{unimodal}\index{unimodal function}, that is, if every local minimum is also a global minimum.
Convexity is sufficient, but not necessary, for the function to be unimodal.)
Most optimization problems have more than one variable, but solution algorithms often rely on simpler ``subproblems'' which are often one-dimensional.
The Frank-Wolfe method from Section~\ref{sec:frankwolfe} is an example of this.

This section presents two other line search methods that can be applied to similar problems.
The \emph{golden section} method requires more iterations than bisection, but does not requiring taking the derivative of $f$.
It is useful if $f$ is not differentiable, or if the derivative is cumbersome to evaluate.
\emph{Newton's method}, on the other hand, requires the stronger assumption that $f$ is twice differentiable (and requires calculating both first and second derivatives), but usually requires fewer iterations than bisection.
 
\subsection{Golden section method}
\label{sec:goldensectionsearch}

\index{optimization!line search!golden section|(}
The golden section method does not use derivatives to narrow down where the minimum lies.
Instead, the logic is as follows: assume we know the value of the function at \emph{two} intermediate points (call them $c$ and $d$), as well as at the endpoints $a$ and $b$.
These points are in the order $a < c < d < b$.  
If $f$ is unimodal, then there are only two possibilities: either $f(a) \geq f(c) \geq f(d)$, or $f(c) \leq f(d) \leq f(b)$.
The reason is that a unimodal function must decrease before its minimum is reached, and increase afterwards.
So if $f(c) \geq f(d)$, then the minimum must be to the right of $c$, and the function must be decreasing everywhere to the left of $c$, so $f(a) \geq f(c) \geq f(d)$.
If $f(c) \leq f(d)$, then using a similar argument we know that $f(c) \leq f(d) \leq f(b)$.
In the first case, we know that the minimum is somewhere in the interval $[c, b]$, and in the second case we know it is somewhere in the interval $[a,d]$.
In either case, we have reduced the width of the interval where the minimum can lie, and we iterate by choosing two more intermediate points in the reduced interval, on till convergence.

This rule works no matter how the two intermediate points are chosen.
To make the method as efficient as possible, we want to re-use points already found if possible.
For instance, in the first case above, the interval is reduced to $[c, b]$.
The point $d$ is within this interval, and can serve as one of the two points we use at the next iteration.
We also want to reduce the interval as much as possible from one iteration to the next.
Since we don't know in advance which end we will trim from, $c$ should be as far from the left endpoint as $d$ is from the right, that is, $c - a = b - d$.
To be able to re-use the point $d$ in the next iteration, we want $d$ to be the same fraction of the distance from $c$ to $b$ as $c$ is from $a$ to $b$. 
(Figure~\ref{fig:goldensection}).
The same logic applies if we are in the second case and must eliminate the upper end.

Using these principles to solve for $c$ and $d$, we see that $c$ should be to the right of $a$ by $(3 - \sqrt{5})/2 \approx 38.2\%$ of the total length of the interval $[a,b]$, and that $d$ should be to the right of $a$ by $(\sqrt{5} - 1)/2 \approx 61.8\%$ of the length of $[a,b]$.
Both of our principles are satisfied: $c$ is 38.2\% to the right of $a$, and $d$ is 38.2\% to the left of $b$; and point $d$ lies 38.2\% of the way from $c$ to $b$, just as $c$ is 38.2\% of the way from $a$ to $b$.
This means that the lengths of the intervals $[c, b]$ and $[a, b]$ have the \emph{golden ratio} of 1 to 1.618.
(This ratio appears often in nature, art, design, and mathematics, being very closely related to the Fibonacci sequence.)

With this in mind, here are the steps of the algorithm:
\begin{description}
   \item[Step 0: Initialize.]  Let $\theta \leftarrow (3 - \sqrt{5})/2$,\label{not:thetags} set the iteration counter $k \leftarrow 0$, and initialize $a_0 \leftarrow a$, $b_0 \leftarrow b$, $c_0 \leftarrow a + \theta(b-a)$, $d_0 \leftarrow b - \theta(b-a)$.
   \item[Step 1: Determine which end to eliminate.]  If $f(c_k) \geq f(d_k)$, go to step 2; otherwise, go to step 3.
   \item[Step 2: Eliminate lower end.]  Update $a_{k+1} \leftarrow c_k$, $b_{k+1} \leftarrow b_k$, $c_{k+1} \leftarrow d_k$, and $d_{k+1} \leftarrow c_k + \theta(b_k - c_k)$.  Skip to step 4.   
   \item[Step 3: Eliminate upper end.]  Update $a_{k+1} \leftarrow a_k$, $b_{k+1} \leftarrow d_k$, $c_{k+1} \leftarrow a_k + \theta(c_k - a_k)$, and $d_{k+1} \leftarrow c_k$.  Proceed to step 4.
   \item[Step 4: Iterate.]  Increase the counter $k$ by 1 and check the termination criterion.
                            If $b_k - a_k < \epsilon$, then terminate; otherwise, return to step 1.
\end{description}

\genfig{goldensection}{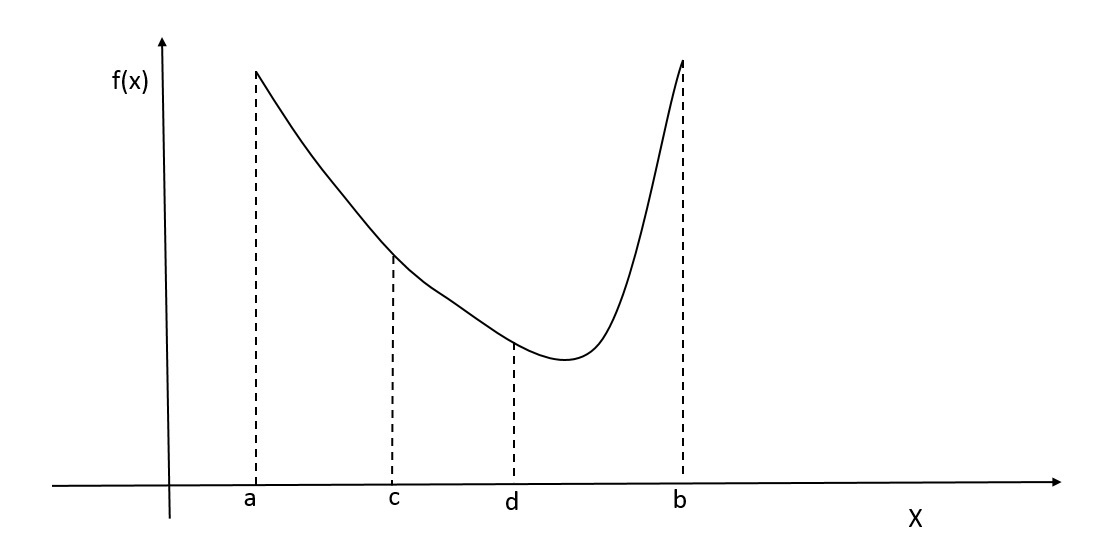}{Golden section method.}{width=\linewidth}

\begin{exm}
Find the minimum of the function $f(x) = (x-1)^2 + e^x$ in the interval $[0,2]$, within a tolerance of $\epsilon = 0.01$.
\end{exm}
\solution{
In the initialization phase, we set $k=0, a_0 = 0, b_0 = 2, \epsilon = 0.01$, and 
\begin{align*}
c_0 &= a_0 + \theta(b_0 -a_0) &= 0.764 \\
d_0 &= b_0 - \theta(b_0 -a_0) &= 1.236 \,.
\end{align*}
We now proceed to step 1.
Since $f(c_1) = 2.2025 < f(d_1) = 3.4975$, we decide to eliminate the upper end and perform step 3.
\begin{align*}
a_1 &= a_0 = 0 \\
b_1 &= d_0 = 1.236  \\
c_1 &= a_0 + \theta(c_0 - a_0) = 0.4722 \\
d_1 &= c_0 = 0.764  \,. 
\end{align*}

The interval is still wider than the tolerance $\epsilon$, so we return to the first step.
Now, $f(c_2) = 1.882 < f(d_2) = 2.2025$, so we again eliminate the upper end by performing step 3.
\begin{align*}
a_2 &= a_1 = 0  \\
b_2 &= d_1 =  0.764  \\
c_2 &= a_1 + \theta(c_1 - a_1) = 0.2918 \\
d_2 &= c_1 = 0.4721  \,.
\end{align*}
We proceed as follows until $b_k - a_k < \epsilon$, with the remaining steps shown in Table~\ref{tbl:goldensection}.
The optimal value is $\frac{a_{13} + b_{13}}{2} = 0.3150$.
}

We used this same example for bisection in Section~\ref{sec:bisection}.  
Golden section required more iterations than bisection, to find a solution with the same tolerance.
This is because bisection reduces the interval width by half at each iteration, whereas golden section reduces it by $\theta \approx 38\%$.  
On the other hand, we did not have to calculate any derivatives.
The value of the objective function $f$ was enough.

\begin {table}
\begin{center}
\caption{Demonstration of the golden section algorithm with $f(x) = (x-1)^2 + e^x$, $x \in [0,2]$.}
\begin{tabular}{| c | c c | c c | c c| }
\hline
 $k$ & $a_k$ & $b_k$ & $c_k$ & $d_k$ & $f(c_k)$ & $f(d_k)$  \\ 
\hline 
0  &   0   &   2   & 0.764 & 1.236 & 2.20254 & 3.49751 \\
1 & 0 & 1.236 & 0.472152 & 0.764 & 1.88206 & 2.20254 \\
2 & 0 & 0.764 & 0.291848 & 0.472152 & 1.84038 & 1.88206 \\
3 & 0 & 0.472152 & 0.180362 & 0.291848 & 1.86946 & 1.84038 \\
4 & 0.180362 & 0.472152 & 0.291848 & 0.360688 & 1.84038 & 1.84304 \\
5 & 0.180362 & 0.360688 & 0.249247 & 0.291848 & 1.84669 & 1.84038 \\
6 & 0.249247 & 0.360688 & 0.291848 & 0.318118 & 1.84038 & 1.83950 \\
7 & 0.291848 & 0.360688 & 0.318118 & 0.334391 & 1.83950 & 1.84012 \\
8 & 0.291848 & 0.334391 & 0.308100 & 0.318118 & 1.83956 & 1.83950 \\
9  & 0.308100 & 0.334391 & 0.318118 & 0.324348 & 1.83950 & 1.83963 \\
10 & 0.308100 & 0.324348 & 0.314306 & 0.318118 & 1.83948 & 1.83950 \\
11 & 0.308100 & 0.318118 & 0.311926 & 0.314306 & 1.83950 & 1.83948 \\
12 & 0.311926 & 0.318118 & 0.314306 & 0.315753 & 1.83948 & 1.83949  \\
\hline
\end{tabular}
\label{tbl:goldensection}
\end{center}
\end{table}

\begin{exm}
Find the minimum of the function $f(x) = 1+ e^{-x}\sin(-x)$ in the interval $[0,3]$ with a tolerance of 0.01.
\end{exm}
\solution{
See Table~\ref{tbl:goldensection_2} below.
The optimal value is $\frac{a_{13} + b_{13}}{2} = 0.7860$.
(Using calculus, the optimal solution is actually $\pi / 4 \approx 0.7854$.)
}
\index{optimization!line search!golden section|)}

\begin {table}
\begin{center}
\caption{Demonstration of the golden section algorithm with $f(x) = 1+ e^{-x}\sin(-x)$, $x \in [0,3]$.}
\begin{tabular}{| c | c c | c c | c c| }
\hline
 $k$ & $a_k$ & $b_k$ & $c_k$ & $d_k$ & $f(c_k)$ & $f(d_k)$  \\
\hline 
0  & 0 & 3 & 1.146 & 1.854 & 0.710349 & 0.849629 \\
1  & 0 & 1.854 & 0.708228 & 1.146 & 0.679624 & 0.710349 \\
2  & 0 & 1.146 & 0.437772 & 0.708228 & 0.726369 & 0.679624 \\
3  & 0.437772 & 1.146 & 0.708228 & 0.875457 & 0.679624 & 0.680064 \\
4  & 0.437772 & 0.875457 & 0.604968 & 0.708228 & 0.689418 & 0.679624 \\
5  & 0.604968 & 0.875457 & 0.708228 & 0.772130 & 0.679624 & 0.677660 \\
6  & 0.708228 & 0.875457 & 0.772130 & 0.811575 & 0.677660 & 0.677820 \\
7  & 0.708228 & 0.811575 & 0.747707 & 0.772130 & 0.678073 & 0.677660 \\
8  & 0.747707 & 0.811575 & 0.772130 & 0.787178 & 0.677660 & 0.677604 \\
9  & 0.772130 & 0.811575 & 0.787178 & 0.796507 & 0.677604 & 0.677643 \\
10 & 0.772130 & 0.796507 & 0.781442 & 0.787178 & 0.677608 & 0.677604 \\
11 & 0.781442 & 0.796507 & 0.787178 & 0.790752 & 0.677604 & 0.677612 \\
12 & 0.781442 & 0.790752 & 0.784999 & 0.787178 & 0.677603 & 0.677604 \\
\hline
\end{tabular}
\label{tbl:goldensection_2}
\end{center}
\end{table}

\subsection{Newton's method}
\label{sec:newtonlinesearch}

\index{Newton's method!for optimization|see {optimization, line search, Newton's method; optimization, nonlinear, Newton's method}}
\index{optimization!line search!Newton's method|(}
You may recall Newton's method from calculus: to find a zero of some function $g$, the method starts with an initial guess, which is iteratively updated by finding where the linear approximation of $g$ at that point has a zero.
Under the right conditions, Newton's method is almost miraculously fast, exhibiting quadratic convergence (roughly speaking, the number of correct digits doubles at each iteration).
There are a few downsides: it can only be applied if $g$ is differentiable, and Newton's method may not converge at all if you start too far away from the root, or if the shape of the function makes Newton's method ``overshoot'' and end up further away than it started.  
Nevertheless, Newton's method and its variants are very commonly used in optimization.
In this appendix, we will use it in several ways; and several of the algorithms mentioned in the main text make use of it as well.

Here we are specifically concerned with solving a one-dimensional optimization problem with an interval constraint.
If the objective function $f$ is convex, then it is enough to find a point where the derivative $f'$ vanishes.
So, we simply apply Newton's method to the \emph{derivative}, with $g \equiv f'$, to try to find $\hat{x}$ such that $f'(\hat{x}) = 0$.
Newton's method uses the derivative of $g$, which ends up being the \emph{second} derivative $f''$.
(If $f$ is not twice-differentiable, Newton's method cannot be applied.)
We have to make one minor modification to Newton's method: the line search cannot leave the feasible region $[a,b]$, so we truncate the search at these boundary points.
This can actually be helpful, since it prevents Newton's method from diverging.
There are still cases where Newton's method can fail; an example is given in Example~\ref{exm:newtonfail} below.

Unlike bisection or golden section, Newton's method is not an ``interval reduction'' method, where we gradually shrink the range of possible values where the optimum can lie.
So we need a different way to measure convergence.
It is common to stop when $f'$ is ``close enough'' to zero; we will let $\epsilon'$ denote this value. 

The steps of Newton's method for line search are:
\begin{description}
\item[Step 0: Initialize.]  Set the iteration counter $k \leftarrow 0$, and initialize $x_0$ to any point in $[a,b]$.  (If you have a good guess as to the minimum point, it can greatly speed things up.)
\item[Step 1: Check convergence.] If $|f'(x_k)| < \epsilon'$, then terminate.
\item[Step 2: Calculate recommended shift.]  Create a candidate point $\tilde{x} \leftarrow x_k - f'(x_k) / f''(x_k)$.  
\item[Step 3: Ensure feasibility.] Project candidate solution onto the feasible region by setting $x_{k+1} \leftarrow a$ if $\tilde{x} < a$, $x_{k+1} \leftarrow b$ if $\tilde{x} > b$, and $x_{k+1} \leftarrow \tilde{x}$ otherwise.
That is, set $x_{k+1} \leftarrow \mbox{proj}_{[a,b]} (\tilde{x})$.
\item[Step 4: Iterate.]  Increase $k$ by 1, and return to step 1.
\end{description}

\begin{exm}
Find the minimum of the function $f(x) = (x-1)^2 + e^x$ in the interval $[0,2]$ using Newton's method, with $\epsilon' = 0.0337$.
(In this and the next example, $\epsilon'$ is chosen to make the tolerance comparable to the $\epsilon$ value used for bisection and golden section; for this function, when $|f'(x)| < 0.0337$, $x$ is within 0.01 of its optimal value.)
\end{exm}
\solution{
Start by computing the formulas for the first and second derivative of $f$, since we will be using these often: $f'(x) = 2(x-1)+e^x$, and $f''(x) = 2+e^x$.

In the initialization phase, we set $k=0$.
For an initial guess, choose $x_0 = 1$.  (This makes for a fair comparison, since this is the starting point for bisection.)
The first and second derivatives are equal to $e$ and $e + 2$, respectively, so the new candidate point is $\tilde{x} = 1 - e/(e+2) = 0.4239$.
This lies within the boundary $[0,2]$, so we accept the candidate point as the next solution: $x_1 = 0.4239$.

At this new point, the first and second derivatives equal 0.3756 and 3.528, so the next candidate is $0.4239 - 0.3756/3.528 = 0.3174$.
We accept the candidate as the new point, so $x_2 = 0.3174$.  
At $x_2$, the derivative is $f'(x_2) = 0.0084 < \epsilon'$, so we terminate and report 0.3174 as the optimal solution.
}

Notice that Newton's method achieved in only two iterations the level of precision bisection reached in eight, and golden section reached in thirteen!
At this point, the solution given by Newton's method differs from the true optimum $x$ by roughly $2 \times 10^{-3}$.  
One more iteration of Newton's method would reduce this error to $1 \times 10^{-6}$, and yet another would reduce it to $3 \times 10^{-13}$.  
This is what we mean when we say its convergence rate is miraculous!

\begin{exm}
Find the minimum of the function $f(x) = 1+ e^{-x}\sin(-x)$ in the interval $[0,3]$ using Newton's method, with $\epsilon' = 0.00645$.
(Again, this choice of $\epsilon'$ ensures that when Newton's method terminates, $x$ is within 0.01 of its optimal value.)
\end{exm}
\solution{
Start by computing the formulas for the first and second derivative of $f$, since we will be using these often: $f'(x) = e^{-x}(\sin x - \cos x)$, and $f''(x) = 2e^{-x} \cos x$.

In the initialization phase, we set $k=0$.
For an initial guess, choose $x_0 = 1.5$.  (This makes for a fair comparison, since this is the starting point for bisection.)
The first and second derivatives are equal to $0.2068$ and $0.03157$, respectively, so the new candidate point is $\tilde{x} = 1.5 - 0.2068/0.03157 = -5.051$.
This is outside of the feasible interval $[0,3]$, so we project the candidate back onto the feasible region, choosing $x_1 = 0$.

At this new point, the first and second derivatives equal $-1$ and 2, so the next candidate is $0 - (-1)/2 = 0.5$.
We accept the candidate as the new point, so $x_2 = 0.5$.  
Another two iterations are needed: $x_3 = 0.5 - (-0.2415)/(1.065) = 0.7268$, and $x_4 = 0.7822$.
At $x_4$, the derivative is $|f'(x_4)| = 0.0021 < \epsilon'$, so we terminate and report 0.7822 as the optimal solution.
}

Newton's method is not infallible, however, as shown by the next example.
\begin{exm}
\label{exm:newtonfail}
Apply Newton's method to the function $f(x) = x \tan^{-1}x - \frac{1}{2} \log (1 + x^2)$ over $x \in [-2, 2]$, starting with the initial guess $x_0 = 1.3917452$. 
(This formula looks unwieldy, but plot the function!
It actually looks quite nice, and you can spot the minimum immediately by inspection.)
\end{exm}
\solution{
The first and second derivatives are $f'(x) = \tan^{-1}x$ and $f''(x) = 1 / (x^2 + 1)$.
Calculating $f'$ and $f''$ at the initial guess $x_0$, we find the next point to be $x_1 = -1.3914752$ --- exactly the opposite of where we started!
You can guess what happens next.  Calculating from $x_1$, we find the next point to be $x_2 = +1.3914752$, and so on \emph{ad infinitum}.
}

Again, this is a ``nice'' function --- it is even strictly convex!
Yet Newton's method fails to find the optimal solution.
You can verify that if we start with an initial guess closer to zero, Newton's method will converge, and if the initial guess is further away, Newton's method will diverge to whatever the endpoints are of the feasible interval.
\index{optimization!line search!Newton's method|)}
\index{optimization!line search|)}

\section{Linear Programming}
\label{sec:linearoptimization}

\index{linear programming|see {optimization, linear}}
\index{optimization!linear|(}
This section briefly outlines the properties of a linear optimization problem (commonly called a \emph{linear program}), terminology associated with solutions of linear programs, and two solution methods: a graphical method suitable for small problems, and the simplex method which is widely used for practical problems. 
Recall that a linear program is an optimization problem where the objective function and constraints are all linear in the decision variables; the constraints may either be equalities ($=$) or weak inequalities ($\geq$, $\leq$).

\subsection{Basic properties and graphical solution method}

This section will predominantly use the following example:
\optimizex{\min}{x,y}{10x + 26y}
{
& 11x + 3y &\geq 21 \\
&  6x + 20y &\geq 39  \\
&   x + y &\leq 9 \\
&   x, y & \geq 0
}
\begin{figure}
\centering
\begin{tikzpicture}
\path [fill=yellow!20!white] (1.5,1.5) -- (0,7) -- (0,9) -- (9,0) -- (6.5,0);
\draw (0,0) -- coordinate (x axis mid) (10,0);
\draw (0,0) -- coordinate (y axis mid) (0,10);
\foreach \x in {0,...,10}      \draw (\x,1pt) -- (\x,-3pt)  node[anchor=north] {\x};
\foreach \y in {0,...,10}      \draw (1pt,\y) -- (-3pt,\y) node[anchor=east] {\y}; 
\node[below=0.8cm] at (x axis mid) {$x$};
\node[rotate=90, left=0.8cm] at (y axis mid) {$y$};
\draw [blue, ultra thick, domain=0:21/11] plot (\x, { -11*\x/3 + 7});
\draw [blue, ultra thick, domain=0:6.5] plot (\x, { -3*\x/10 + 39/20});
\draw [blue, ultra thick, domain=0:9.0] plot (\x, { -\x + 9});
\draw [cyan, dashed, ultra thick, domain=0:9.6] plot (\x, { -10*\x/26 + 96/26});
\node [right] at (4.5,0.75) {$6x+20y \geq 39$};
\node [right] at (0.6,4.5) {$11x+3y \geq 21$};
\node [right] at (2.8,2.8) {$z=10x+26y$};
\node [left] at (2.9,6) {$x+y \leq 9$};
\draw [blue, ultra thick] (0,7)--(0,9);
\draw [blue, ultra thick] (6.5,0)--(9,0);
\node [below left] at (1.6,1.6) {(1.5,1.5)};
\draw [fill=black] (1.5,1.5)circle (3pt);
\end{tikzpicture}
\caption{An illustrative linear program. \label{fig:lp}}
\end{figure}

The feasible region for this problem is shown in Figure~\ref{fig:lp} as the area shaded in yellow.
The optimal solution to this linear program is the point $(x,y)$ in the feasible region which minimizes $10x + 26y$.
In this case, the optimal solution is $(1.5,1.5)$, with a corresponding objective function value of 54.
It is not immediately obvious how we found this point, since there is an infinite number of feasible solutions.

\index{optimization!linear!corner points}
This point happens to lie at a ``corner'' of the feasible region.
This is not a coincidence.
If a linear program has an optimal solution, and if its feasible region has a corner, then there is an optimal solution that lies at such a corner point.
As a result, we can usually confine our attention to the corner points, of which there are only a finite number.
In this particular case, the corner point feasible solutions are $ (0,9), (0,7), (1.5, 1.5), (6.5,0)$ and $(9,0)$.
Checking in each in turn, we see that the objective function is minimized at $(1.5, 1.5)$, so this is the optimal solution.

This provides us with a method for solving a linear program, if the linear program has an optimal solution and its feasible region has at least one corner.
\begin{itemize}
\item Plot all the constraints and identify the feasible region.
\item Determine the feasible corner points.
\item Evaluate the objective function at the feasible corner points.
\item The ``best'' feasible corner point is the optimal solution
\end{itemize}
For small problems this can be done graphically, as in Figure~\ref{fig:lp}.
For larger problems, it is difficult to identify all the corner points by inspection, and a different method is needed.

This intuitive description is enough for our purposes; for readers wanting a more technical definition of a ``corner,'' we can define it as a point that is not the midpoint of any line segment contained in the feasible region.
These points are also known as \emph{vertices}, or \emph{extreme points}\index{extreme point!see {optimization, linear, corner points}}.
A point along an edge of the feasible region can be drawn as the midpoint of a line segment drawn along this edge; a point in the interior of the feasible region is the midpoint of many different line segments.  
But a line segment drawn in the feasible region can only have a corner as one of its endpoints, never the midpoint.

Why do we need the caveats ``if a linear program has an optimal solution'' and ``if its feasible region has a corner''?
Consider the examples below.
Modify the constraints so the optimization problem has the following form:
\optimizex{\min}{x, y}{10x + 26y}{
&   x + y &\leq 9 \\
&   y &\leq 5 \\
&   y &\geq 0
}This produces the feasible region shown in Figure~\ref{fig:ulp}.
There is no restriction on how small $x$ can be, and as $x$ tends to $-\infty$, so does the objective function.
Therefore there is no optimal solution; the problem is said to be \emph{unbounded}\index{optimization!unbounded} since there is no limit to how much the objective function can be minimized.
(The same would hold true for a maximization problem if it is possible for the objective function to be made arbitrarily large.)

\begin{figure}
\centering
\resizebox{\textwidth}{!}{
\begin{tikzpicture}
\path [fill=yellow!20!white] (-4,0) -- (-4,5) -- (4,5) -- (9,0) -- (-4,0);
\draw (-5,0) -- coordinate (x axis mid) (10,0);
\draw (0,-1) -- coordinate (y axis mid) (0,10);
\foreach \x in {-4,...,-1}      \draw (\x,1pt) -- (\x,-3pt)  node[anchor=north] {\x};
\foreach \x in {1,...,10}      \draw (\x,1pt) -- (\x,-3pt)  node[anchor=north] {\x};
\foreach \y in {-1}    \draw (1pt,\y) -- (-3pt,\y) node[anchor=east] {\y}; 
\foreach \y in {1,...,10}   \draw (1pt,\y) -- (-3pt,\y) node[anchor=east] {\y}; 
\node[below=0.8cm] at (x axis mid) {$x$};
\node[rotate=90] at (-1,5.5) {$y$};
\draw [blue, ultra thick, domain=-1.0:9.0] plot (\x, { -\x + 9});
\draw [cyan, dashed, ultra thick, domain=-4.0:9.6] plot (\x, { -10*\x/26 + 96/26});
\draw [blue, ultra thick, domain=-4.0:9.0] plot (\x, {0});
\draw [blue, ultra thick, domain=-4.0:4.0] plot (\x, {5});
\draw [blue, ultra thick] (-4.0,0)--(9,0);
\draw [blue, ultra thick] (-4.0,5.0)--(4,5);
\node [left] at (4.0,2.0) {$z=10x+26y$};
\node [left] at (2.0,4.75) {$y \leq 5$};
\node [left] at (2.0,0.25) {$y \geq 0$};
\node [left] at (4.9,4) {$x+y \leq 9$};
\end{tikzpicture}
}
\caption{An unbounded linear program. \label{fig:ulp}}
\end{figure}

Or, consider this optimization problem:
\optimizex{\min}{x, y}{10x + 26y}{
&   x + y &\leq 9 \\
&   y &\leq 5 \\
&   x &\geq 9.5 \\
&   y &\geq 0
}There is no solution which obeys all the constraints (see Figure~\ref{fig:inflp}).
The constraints $x + y \leq 9$ and $y \geq 0$ mean that $x$ cannot be greater than 9; however there is also the constraint $x \geq 9.5$.
Such a problem is called \emph{infeasible}\index{optimization!infeasible}, because there is not even a feasible solution (let alone an optimal one).
As a rule of thumb, unbounded or infeasible problems often mean that you have missed something in your formulation.
In the real world, it is not possible to produce ``infinitely good'' solutions (surely there is some limitation; this would be a constraint you have missed), and there is usually some possible course of action, even if it is very unpleasant (your constraints are too restrictive).

\begin{figure}
\centering
\resizebox{\textwidth}{!}{
\begin{tikzpicture}
\draw (-5,0) -- coordinate (x axis mid) (10,0);
\draw (0,-1) -- coordinate (y axis mid) (0,10);
\foreach \x in {-4,...,-1}      \draw (\x,1pt) -- (\x,-3pt)  node[anchor=north] {\x};
\foreach \x in {1,...,10}      \draw (\x,1pt) -- (\x,-3pt)  node[anchor=north] {\x};
\foreach \y in {-1}    \draw (1pt,\y) -- (-3pt,\y) node[anchor=east] {\y}; 
\foreach \y in {1,...,10}    \draw (1pt,\y) -- (-3pt,\y) node[anchor=east] {\y}; 
\node[below=0.8cm] at (x axis mid) {$x$};
\node[rotate=90] at (-1,5.5) {$y$};
\draw [blue, ultra thick, domain=-1.0:9.0] plot (\x, { -\x + 9});
\draw [cyan, dashed, ultra thick, domain=-4.0:9.6] plot (\x, { -10*\x/26 + 96/26});
\draw [blue, ultra thick, domain=-4.0:9.0] plot (\x, {0});
\draw [blue, ultra thick, domain=-4.0:4.0] plot (\x, {5});
\draw [blue, ultra thick] (-4.0,0)--(9,0);
\draw [blue, ultra thick] (-4.0,5.0)--(4,5);
\draw [blue, ultra thick] (9.5,0)--(9.5,10);
\node [left] at (4.0,2.0) {$z=10x+26y$};
\node [left] at (2.0,4.75) {$y \leq 5$};
\node [left] at (2.0,0.25) {$y \geq 0$};
\node [left] at (4.9,4) {$x+y \leq 9$};
\node [left] at (9.5,5) {$x \geq 9.5$};
\end{tikzpicture}
}
\caption{An infeasible linear program. \label{fig:inflp}}
\end{figure}

Here is an optimization problem which has an optimal solution, but its feasible region has no corners (Figure~\ref{fig:nocornerlp}):
\optimizex{\min}{x, y}{x + y}{
&   x + y \geq &0
}
\begin{figure}
\centering
\begin{tikzpicture}
\path [fill=yellow!20!white] (-5,5) -- (5,5) -- (5,-5) -- (-5,5);
\draw (-5,0) -- coordinate (x axis mid) (5,0);
\draw (0,-5) -- coordinate (y axis mid) (0,5);
\foreach \x in {-5,...,-1}      \draw (\x,1pt) -- (\x,-3pt)  node[anchor=north] {\x};
\foreach \x in {1,...,5}      \draw (\x,1pt) -- (\x,-3pt)  node[anchor=north] {\x};
\foreach \y in {-5,...,-1}    \draw (1pt,\y) -- (-3pt,\y) node[anchor=east] {\y}; 
\foreach \y in {1,...,5}    \draw (1pt,\y) -- (-3pt,\y) node[anchor=east] {\y}; 
\node at (5, 0.5) {$x$};
\node[rotate=90] at (-0.5,5) {$y$};
\draw [blue, ultra thick, domain=-5.0:5.0] plot (\x, { -\x});
\draw [cyan, dashed, ultra thick, domain=-2.0:5.0] plot (\x, { -\x + 3});
\node [left] at (4.0,2.0) {$z=x+y$};
\node [left] at (-1.25,1.0) {$x + y \geq 0$};
\end{tikzpicture}
\caption{An linear program with no corner points. \label{fig:nocornerlp}}
\end{figure}
Any solution on the line $x + y = 0$ is feasible and optimal.
So this problem has infinitely many optimal solutions, but no corner points.
These examples show that the statement ``at least one corner point is optimal, so we only have to look at corner points'' needs a few technical qualifications (that there is in fact an optimal solution, and a corner point).
In practice these kinds of counterexamples are rare.

Also note that this statement does not say that \emph{only} corner points may be optimal, even with the qualifications that there are corner points and optimal solutions.
There very well may be a non-corner point which is optimal.
But in such cases there is \emph{also} an optimal corner point, so it is still acceptable to just check the corner points.
We will still find an optimal solution that way.
An example of such a problem is
\optimizex{\min}{x,y}{12x + 40y}{
& 11x + 3y &\geq 21 \\
&  6x + 20y &\geq 39  \\
&   x + y &\leq 9 \\
&   x, y &\geq 0
}(See Figure~\ref{fig:mult_lp}.)
Any feasible point on the line segment between $(1.5, 1.5)$ and $(6.5, 0)$ has the same, optimal objective function value of 78.
So a solution like $(3.5, 0.75)$ is optimal even though it is not a corner point.
However, the corner points $(1.5, 1.5)$ and $(6.5, 0)$ are both optimal as well, so we can still find \emph{an} optimal solution if we just check corners.
(For most purposes it is fine to identify a single optimal solution, and identifying all optima is unnecessary.)

\begin{figure}
\centering
\begin{tikzpicture}
\path [fill=yellow!20!white] (1.5,1.5) -- (0,7) -- (0,9) -- (9,0) -- (6.5,0);
\draw (0,0) -- coordinate (x axis mid) (10,0);
\draw (0,0) -- coordinate (y axis mid) (0,10);
\foreach \x in {0,...,10}      \draw (\x,1pt) -- (\x,-3pt)  node[anchor=north] {\x};
\foreach \y in {0,...,10}    \draw (1pt,\y) -- (-3pt,\y) node[anchor=east] {\y}; 
\node[below=0.8cm] at (x axis mid) {$x$};
\node[rotate=90, left=0.8cm] at (y axis mid) {$y$};
\draw [blue, ultra thick, domain=0:21/11] plot (\x, { -11*\x/3 + 7});
\draw [blue, ultra thick, domain=0:6.5] plot (\x, { -3*\x/10 + 39/20});
\draw [blue, ultra thick, domain=0:9.0] plot (\x, { -\x + 9});
\draw [cyan, dashed, ultra thick, domain=0:100/12] plot (\x, { -12*\x/40 + 100/40});
\node [right] at (4.5,0.75) {$6x+20y \geq 39$};
\node [right] at (0.6,4.5) {$11x+3y \geq 21$};
\node [right] at (3.0,2.0) {$z=12x+40y$};
\node [left] at (2.9,6) {$x+y \leq 9$};
\draw [blue, ultra thick] (0,7)--(0,9);
\draw [blue, ultra thick] (6.5,0)--(9,0);
\draw [fill=black] (1.5,1.5)circle (3pt);
\draw [fill=black] (6.5,0)circle (3pt);
\end{tikzpicture}
\caption{Multiple optimal solutions. \label{fig:mult_lp}}
\end{figure}

This graphical method is good for building intuition and illustrating the main ideas behind solving linear programs.
However, it becomes difficult as problems become larger.
A problem with three decision variables would require a 3D drawing; a problem with four decision variables would require drawing in 4 dimensions, which is beyond the scope of most of us mere mortals.
Practical-sized problems commonly involve thousands or even millions of decision variables, so a different technique is needed.
However, the idea of focusing on corner points can be translated to problems of this scale.
To do so requires a systematic procedure for generating and checking corner points.
This procedure is called the \emph{simplex algorithm}.
The mathematical tools for doing this rely on linear algebra, since corner points can be identified by solving systems of linear equations (in the examples in this section, each corner point was the intersection of two lines).
The remainder of this section is focused on presenting this method.

\subsection{Standard form of a linear program}

\index{optimization!linear!standard form|(}
The simplex algorithm is easier to explain if the linear program is converted to a standard form.
This standard form is not restrictive, and you can transform any linear program into an equivalent one which is in standard form.
A linear program is in standard form if:
\begin{itemize}
\item The objective function is to be minimized.
\item All decision variables have a non-negativity constraint (that is, $x_i \geq 0$ is a constraint for every decision variable $x_i$).
\item All other constraints are equality ($=$) constraints.
\end{itemize}

In the previous section, we considered this linear program:
\optimizex{\min}{x, y}{10x + 26y}{
& 11x + 3y &\geq 21\\
&  6x + 20y &\geq 39  \\
&   x + y &\leq 9 \\
&   x,y & \geq 0
}This is a minimization problem where all variables must be non-negative; but the remaining constraints are inequalities, not equalities.
We can fix this by introducing three new variables $s_1$, $s_2$, and $s_3$, and rewriting the optimization problem:
\optimizex{\min}{x, y, s_1, s_2, s_3}{10x + 26y}{
& 11x + 3y -s_1 &= 21\\
&  6x + 20y - s_2 &= 39  \\
&   x + y + s_3 &= 9 \\
&   x, y, s_1, s_2, s_3 &\geq 0
}This problem is equivalent to the original one in the sense that any solution to the standard-form problem can easily be translated to a solution in the original problem (just ignore the new variables); and a solution is feasible to one problem if and only if it is feasible to the other problem.
Given a solution $(x,y)$ to the original problem, let $s_1 = 11x + 3y - 21$, $s_2 = 6x + 20y - 30$, and $s_3 = 9 - (x+y)$.
If that solution is feasible, these new variables must all be non-negative, and the three equality constraints in the standard-form problem will be satisfied.
The reverse is true as well; for instance, since $11x + 3y = 21 + s_1$ in a feasible solution to the standard-form problem, and since $s_1 \geq 0$, we must have $11x + 3y \geq 21$ and the first constraint in the original problem is also satisfied.
The new variables introduced ($s_1$, $s_2$, and $s_3$) are often called \emph{slack}\index{optimization!linear!slack variable} variables, because they show how much the left-hand side of the constraint can be changed before it is violated.

This technique can be used to ensure that all constraints are equalities.
Similarly, non-standard problems violating the other requirements can also be converted into equivalent, standard-form linear programs.
If the objective function is maximization, rather than minimization, we can multiply the objective function by $-1$. 
(Minimizing $-f$ is the same as maximizing $f$; see Proposition~\ref{prp:maxmin} in Section~\ref{sec:objectivetransform}.)
If there is a decision variable which is ``free'' (it does not have to be non-negative), it can be replaced with two new decision variables.
For instance, if the variable $y$ can be either negative or positive, everywhere $y$ appears in the optimization problem we can replace it with $y^+ - y^-$, where $y^+$ and $y^-$ are two new decision variables with non-negativity constraints ($y^+ \geq 0$, $y^- \geq 0$).
This works because we represent any value (positive or negative) as a subtraction of two positive numbers.

It is common to write linear programs in matrix-vector form, since this is more convenient when there are a large number of decision variables and constraints.
In matrix notation, the linear program in standard form can be written as:
\optimizex{\min}{\mb{x}}{\mb{c} \cdot \mb{x}}{
& \mb{Ax} &= \mb{b} \\
& \mb{x} &\geq \mb{0}
}where $\mb{x}$ is the vector of decision variables, $\mb{c}$\label{not:clp} is a vector representing the coefficients in the objective function, $\mb{b}$\label{not:blp} is a vector representing the right-hand side of each equality constraint, and $\mb{A}$\label{not:Alp} is the matrix of coefficients in all of the equality constraints.
(Section~\ref{sec:matrices} provides a review of vector and matrix operations if you are unfamiliar with them.)
For the above example in standard form, these vectors and matrices are as follows:
\begin{align*}
\mb{x} &= \begin{bmatrix}
                x \\
                y \\
                s_1 \\
                s_2 \\
                s_3
                \end{bmatrix}
\end{align*}

\begin{align*}
    \mb{c} &= \begin{bmatrix}
                10 \\
                26 \\
                0 \\
                0 \\
                0
                \end{bmatrix}
\end{align*}

\begin{align*}
    \mb{A} &= \begin{bmatrix}
                11 &  3  & -1 & 0  & 0 \\
                6  &  20 & 0  & -1  & 0 \\
                1  &  1  & 0  & 0   & 1
                \end{bmatrix}
\end{align*}

\begin{align*}
    \mb{b} &= \begin{bmatrix}
                21 \\
                39 \\
                9 
                \end{bmatrix}
\,.
\end{align*}
\index{optimization!linear!standard form|)}

\subsection{Simplex algorithm}
\label{sec:simplex}

\index{optimization!linear!simplex method|(}
Usually, a linear program has an infinite number of feasible solutions.
However, as discussed above, it is enough to look at corner points of the feasible region, since at least one of them will be optimal (assuming that the problem has an optimal solution, and that there is a corner point, conditions we will take as given in this subsection).
The simplex algorithm is designed to systematically search these corner points.
Starting from one corner point, it moves to an adjacent corner point with an equal or smaller objective function value.
Since there are a finite number of corner points, and since linear programs have the property that you can reach an optimal corner from any other by a sequence of adjacent corners with nonincreasing objective function values, this procedure is guaranteed to find an optimal solution.
A corner point can be identified by solving systems of equations representing a subset of the constraints, so the simplex method is presented in terms of linear algebra.
This subsection provides the steps of the simplex method, along with some intuitive justification for them, but we do not provide rigorous proofs of correctness, nor do we handle all of the ``corner cases'' that can arise when applying the simplex method.
A reader interested in these details is referred to a standard text in linear programming, such as~\cite{bertsimas97}.

Assume we are given a standard-form linear program, where the decision variable $\mb{x}$ is a vector of dimension $n$,\label{not:nlp} the objective function coefficient $\mb{c}$ is a vector of dimension $n$, $\mb{b}$ is a vector of dimension $m$,\label{not:mlp} and the coefficient matrix $\mb{A}$ is a matrix of dimension $m \times n$.
We can assume there are at least as many variables as equality constraints ($n \geq m$); otherwise some of the constraints are redundant and can be removed.
(Indeed, if the number of independent constraints is equal to the number of variables, there is only one feasible solution!)
Let $\mb{A_i}$\label{not:Ai} be the vector corresponding to the $i$-th column of the matrix $\mb{A}$.

\begin{align*}
\mb{A_i} &= \begin{bmatrix}
        a_{i1} \\
        a_{i2} \\
        \vdots \\
        a_{im}
\end{bmatrix}
\end{align*}

\optimizex{\min}{\mb{x}}{\mb{c} \cdot \mb{x}}{
& \mb{Ax} &= \mb{b} \\
& \mb{x} &\geq 0
}

A \emph{basis}\index{optimization!linear!basis} $B$\label{not:Blp} refers to a set of $m$ linearly independent columns of $\mb{A}$.
These columns are represented by the square matrix $\mb{A_B}$, which is of size $m \times m$.
We will use $N$\label{not:NBlp} (for ``nonbasic'') to represent all of the other columns of $\mb{A}$.
By limiting ourselves to the decision variables $x_i$ corresponding to the columns in the basis, we now have exactly the same number of decision variables and constraints.
Fixing all of the other decision variables to zero, we can now solve the $m$ equality constraints in terms of these $m$ basic variables.
The resulting solution $\mb{x_B}$\label{not:xBlp} corresponds to an intersection point of the constraints.
If $\mb{x_B} \geq \mb{0}$, then this solution satisfies all of the constraints, and is therefore feasible.
The resulting value of $\mb{x}$, whose basic components solve the constraints ($\mb{x_B} = \mb{A_B}^{-1} \mb{b})$, and whose nonbasic components are zero ($\mb{x_N} = \mb{0}$) form a corner point solution.

Every corner point can be identified in this way.
So the question is to find some basis $B$ corresponding to an optimal solution.
Starting from some initial basis, the simplex method looks for an adjacent basis (one which is the same except for one column) with a lower value of the objective function.
If one exists, we update the basis and repeat the process to try to find a better one.
If no adjacent basis is better, we terminate with the optimal solution.
The simplex method can determine whether an adjacent basis is better by computing the \emph{reduced cost}\index{optimization!linear!reduced cost}, the rate of change in the objective function if a nonbasic column were to enter the basis in place of one of the existing columns.
Since the objective function is minimization, a negative reduced cost will represent a better solution.
Using matrix algebra, we can easily compute reduced costs, and new bases, using formulas shown below.

The steps of the simplex method are as follows:
\begin{enumerate}
\item \textbf{Initialization}:
    \begin{itemize}
    \item Determine a basis $\mb{A_B}$ consisting of $m$ linearly independent columns of $\mb{A}$.
    \item Calculate $\mb{x_B} \leftarrow \mb{A_B}^{-1} \mb{b}$.
    \item If $\mb{x_B} \geq \mb{0}$, continue; otherwise, this basis is infeasible.
          Choose a different basis and repeat this step.
          If you have already tried every basis, the linear program is infeasible.  
          Terminate the algorithm and report this.
    \item Set $\mb{x_N} \leftarrow 0$.
    \end{itemize}
\item \textbf{Compute reduced costs}:
    \begin{itemize}
    \item Compute the reduced cost of each nonbasic index $j \in N$ as $\bar{c}_j \leftarrow c_j - \mb{c_B}^T \mb{A_B}^{-1} \mb{A_j}$,\label{not:cj} where $\mb{c_B}$\label{not:cB} is the column vector of objective function coefficients for the basic variables.
    \item If $\bar{c}_j \geq 0$ for all $j \in N$, then the current solution is optimal.
          Terminate the algorithm and return the solution.
    \item Otherwise, choose some nonbasic variable $j \in N$ with $\bar{c}_j < 0$.
          The column representing this variable will enter the basis.
    \item Set $d_j \leftarrow 1$, and $d_i \leftarrow 0$ for all $i \in N - \{ j \}$.
    \end{itemize}
\item \textbf{Identify descent direction}:
    \begin{itemize}
    \item Compute the basic direction $\mb{d_B} \leftarrow -\mb{A_B}^{-1} \mb{A_j}$.
    \item If $d_i \geq 0$ for all $i \in B$, the linear program is unbounded.
          Terminate the algorithm and report this.
    \item Otherwise, $d_i < 0$ for at least one $i \in B$.  Compute the maximum step size
          \labeleqn{simplexstepsize}{
          \theta = \min_{i \in B: d_i < 0} \myp{ \frac{x_i}{-d_i}}
          \,.}
    \item Let $k$ be the index of a variable achieving this minimum, with $\theta = -x_i/d_i$.
          The column representing this variable will leave the basis (to be replaced with column $j$).
    \end{itemize}
\item \textbf{Compute new solution}:
    \begin{itemize}
    \item Set $x_i^{new} \leftarrow x_i + \theta d_i$ for all variables $i$ (both basic and nonbasic).
    \item Update the basis $\mb{A_B}$ by replacing $\mb{A_k}$ with $\mb{A_j}$.
    \item Update the index sets by removing $k$ from $B$, and adding $j$.  Likewise, remove $j$ from $N$ and add $k$. 
    \item Update the vectors $\mb{c_b},\mb{x_b}, \mb{x_N}$ appropriately based on $\mb{x_i^{new}}$.
    \item Return to step 2.
    \end{itemize}
\end{enumerate}

As an example, we show how the simplex method will operate on the following standard-form linear program (Figure~\ref{fig:simplexlp} plots the feasible region).

\optimizex{\min}{x,y,s_1, \ldots, s_4}{10x + 26y}{
& 11x + 3y -s_1 &= 21\\
&  6x + 20y - s_2 &= 39  \\
&   x + y + s_3 &= 9 \\
&   -\frac{1}{3}x + y + s_4 &= 6 \\
&   4x - 5y + s_5 &= 20 \\
&   x, y, s_1, s_2, s_3, s_4 &\geq 0
}

\begin{figure}
\centering
\begin{tikzpicture}
\path [fill=yellow!20!white] (1.5,1.5) -- (1/4,73/12) -- (9/4,27/4) -- (65/9, 16/9) -- (119/22,36/110) -- (1.5, 1.5);
\draw (0,0) -- coordinate (x axis mid) (10,0);
\draw (0,0) -- coordinate (y axis mid) (0,10);
\foreach \x in {0,...,10}      \draw (\x,1pt) -- (\x,-3pt)  node[anchor=north] {\x};
\foreach \y in {0,...,10}    \draw (1pt,\y) -- (-3pt,\y) node[anchor=east] {\y}; 
\node[below=0.8cm] at (x axis mid) {$x$};
\node[rotate=90, left=0.8cm] at (y axis mid) {$y$};
\draw [blue, ultra thick, domain=0:21/11] plot (\x, { -11*\x/3 + 7});
\draw [blue, ultra thick, domain=0:6.5] plot (\x, { -3*\x/10 + 39/20});
\draw [blue, ultra thick, domain=0:9.0] plot (\x, { -\x + 9});
\draw [blue, ultra thick, domain=0:5.0] plot (\x, { \x/3 + 6});
\draw [blue, ultra thick, domain=5.0:10.0] plot (\x, { 4*\x/5 - 4});
\draw [cyan, dashed, ultra thick, domain=0:9.6] plot (\x, { -10*\x/26 + 96/26});
\node [right] at (2.5,1.4) {$6x+20y \geq 39$};
\node [right] at (0.9,3.7) {$11x+3y \geq 21$};
\node [right] at (4.0,7.0) {$\frac{-1}{3}x + y \leq 6$};
\node [left] at (8.7,3.0) {$4x - 5y \leq 20$};
\node [right] at (2.8,2.8) {$z=10x+26y$};
\node [left] at (3.9,5) {$x+y \leq 9$};
\draw [blue, ultra thick] (0,7)--(0,9);
\draw [blue, ultra thick] (6.5,0)--(9,0);
\node [below left] at (1.6,1.6) {(1.5,1.5)};
\draw [fill=green] (1.5,1.5)circle (3pt);
\draw [fill=yellow] (1/4,73/12)circle (3pt);
\draw [fill=red] (9/4,27/4) circle (3pt); 
\draw [fill=orange] (65/9, 16/9) circle (3pt);
\draw [fill=black] (119/22,36/110) circle (3pt);
\end{tikzpicture}
\caption{Feasible region for simplex method example. \label{fig:simplexlp}}
\end{figure}

\underline{Initialization}: Set $B = \{1,2,3,4,7\}$ and $N = \{5,6\}$.
Then
\begin{equation*}
\begin{aligned}[c]
    \mb{b} &= \begin{bmatrix} 
        21 \\
        39  \\
        9  \\
        6 \\
        20 
        \end{bmatrix}
\end{aligned}
\begin{aligned}[c]
    \mb{A_B} &= \begin{bmatrix} 
        11 &  3 & -1 & 0 & 0 \\
        6  &  20 & 0 & -1 & 0 \\
        1  &  1 & 0  & 0  & 0 \\
        \frac{-1}{3} & 1 & 0 & 0 & 0 \\
        4 & -5 & 0 & 0 & 1 
        \end{bmatrix}
\end{aligned}
\begin{aligned}[c]
        \mb{c_B} &= \begin{bmatrix} 
        10 \\
        26 \\
        0  \\
        0 \\
        0 
        \end{bmatrix}
\end{aligned}
\end{equation*}

\begin{equation*}
 \begin{aligned}[c]
    \mb{x_B} &= \mb{A_B}^{-1} \mb{b}     = \begin{bmatrix} 
                                            x \\
                                            y  \\
                                            s_1  \\
                                            s_2 \\
                                            s_5 
                                            \end{bmatrix}
                                            &= \begin{bmatrix}
                                            2.25 \\
                                            6.75 \\
                                            24.00 \\
                                            109.5 \\
                                            44.75
                                            \end{bmatrix}
\end{aligned}
\qquad
\begin{aligned}[c]
    \mb{x_N} &= \begin{bmatrix} 
                        s_3 \\
                        s_4
                        \end{bmatrix}
                        &= \begin{bmatrix}
                        0 \\
                        0
                        \end{bmatrix}
\end{aligned}
\,.
\end{equation*}

The values of ${x_B}$ are non-negative, so this basis is a feasible corner point.
The objective function $z = \mb{c_B}^T\mb{x_B} = 198$.
The solution corresponds to the red dot in Figure \ref{fig:simplexlp}.

\underline{Compute reduced costs}: $N = \{5,6\}$.

\begin{align*}
\bar{c}_5 &= c_5 - \mb{c_B}^T \mb{A_B}^{-1} \mb{A_5} = -14 \\
\bar{c}_6 &= c_6 - \mb{c_B}^T \mb{A_B}^{-1} \mb{A_6} = -12
\,.
\end{align*}

As both reduced costs are negative, we can pick either of them.
Let us pick $j=5$.
Set $d_5 = 1, d_6 = 0$.

\underline{Identify descent direction}:

Now $B = \{1,2,3,4,7\}$.
Compute:

\begin{align*}
\mb{d_B} &= -\mb{A_B}^{-1} \mb{A_5} = \begin{bmatrix}
                                        -0.75 \\
                                        -0.25 \\
                                        -9.00 \\
                                        -9.50 \\
                                         1.75 
                                        \end{bmatrix}
\end{align*}

\begin{equation*}
\theta = \min \myc{ \frac{2.25}{0.75}, \frac{6.75}{0.25}, \frac{24.0}{9.0}, \frac{109.5}{9.5} } = \frac{24}{9}
\,.
\end{equation*}

The index $k = 3$.

\underline{Compute new solution}:

\begin{equation*}
\mb{x^{new}} = \mb{x} + \theta \mb{d} = \begin{bmatrix}
                                                    2.25 \\
                                                    6.75 \\
                                                    24 \\
                                                    109.5 \\
                                                    0 \\
                                                    0 \\
                                                    44.75 
                                                    \end{bmatrix} + \theta
                                                    \begin{bmatrix}
                                                    -0.75 \\
                                                    -0.25 \\
                                                    -9.00 \\
                                                    -9.50 \\
                                                    1.00 \\
                                                    0 \\
                                                    1.75 
                                                    \end{bmatrix}
                                                    = \begin{bmatrix}
                                                    0.25 \\
                                                    6.0833 \\
                                                    0 \\
                                                    84.1667 \\
                                                    2.6667 \\
                                                    0 \\
                                                    49.4167 
                                                    \end{bmatrix}
                                                    \,.
\end{equation*}

The new $B = \{1,2,5,4,7\}$ and $N = \{3,6\}$.

\begin{equation*}
 \begin{aligned}[c]
    \mb{x_B} &= \begin{bmatrix} 
                    x \\
                    y \\
                    s_3 \\
                    s_2 \\
                    s_5 
                    \end{bmatrix}
                    &= \begin{bmatrix}
                    0.25 \\
                    6.083 \\
                    2.6667 \\
                    84.1667 \\
                    49.4167
                    \end{bmatrix}
\end{aligned}
\qquad
\begin{aligned}[c]
    \mb{x_N} &= \begin{bmatrix} 
                    s_1 \\
                    s_4
                    \end{bmatrix}
                    &= \begin{bmatrix}
                    0 \\
                    0
                    \end{bmatrix}
\end{aligned}
\end{equation*}

\begin{equation*}
\begin{aligned}[c]
    \mb{A_B} &= \begin{bmatrix} 
                    11           & 3  & 0 &  0 & 0 \\
                    6            & 20 & 0 & -1 & 0 \\
                    1            & 1  & 1 &  0 & 0 \\
                    \frac{-1}{3} & 1  & 0 &  0 & 0 \\
                    4            & -5 & 0 &  0 & 1
                    \end{bmatrix}
\end{aligned}
\begin{aligned}[c]
    \mb{c_B} &= \begin{bmatrix} 
                    10 \\
                    26 \\
                    0  \\
                    0 \\
                    0 
                    \end{bmatrix}
\end{aligned}
\,.
\end{equation*}

The objective function value is $z = \mb{c_B}^T\mb{x_B} = 160.66588$.
The basic feasible solution corresponds to the yellow dot in Figure~\ref{fig:simplexlp}.

\underline{Compute reduced costs}: $N = \{3,6\}$.

\begin{align*}
\bar{c}_3 &= c_3 - \mb{c_B}^T \mb{A_B}^{-1} \mb{A_3} = 1.5556 \\
\bar{c}_6 &= c_6 - \mb{c_B}^T \mb{A_B}^{-1} \mb{A_6} = -21.333
\,.
\end{align*}

As the reduced cost $c_6 < 0$, pick $j=6$.
Set $d_6 = 1, d_3 = 0$.

\underline{Identify descent direction}:

Now $B = \{1,2,5,4,7\}$.
Compute:

\begin{align*}
\mb{d_B} &= -\mb{A_B}^{-1}\mb{A_6} = \begin{bmatrix}
                                        0.25 \\
                                        -0.9167 \\
                                        0.6667 \\
                                        -16.8333 \\
                                        -5.5833 
                                        \end{bmatrix}
\end{align*}

\begin{equation*}
\theta = \min \myc{ \frac{6.083}{0.967}, \frac{84.1667}{16.8333}, \frac{49.4167}{5.5833} } = \frac{84}{16.8333}
\,.
\end{equation*}

The index $k = 4$.

\underline{Compute new solution}:

\begin{equation*}
\mb{x^{new}} = \mb{x} + \theta \mb{d} = \begin{bmatrix}
                                                    0.25 \\
                                                    6.0833 \\
                                                    0 \\
                                                    84.1667 \\
                                                    2.6667 \\
                                                    0 \\
                                                    49.4167 
                                                    \end{bmatrix} + \theta
                                                    \begin{bmatrix}
                                                    0.25 \\
                                                    -0.9167 \\
                                                    0.0 \\
                                                    -16.8333 \\
                                                    0.6667 \\
                                                    1.0 \\
                                                    -5.5833 
                                                    \end{bmatrix}
                                                    = \begin{bmatrix}
                                                    1.5 \\
                                                    1.5 \\
                                                    0 \\
                                                    0 \\
                                                    6.0 \\
                                                    5.0 \\
                                                    21.5 
                                                    \end{bmatrix}
                                                    \,.
\end{equation*}

The new $B = \{1,2,5,6,7\}$ and $N = \{3,4\}$.

\begin{equation*}
 \begin{aligned}[c]
    \mb{x_B} &= \begin{bmatrix} 
                    x \\
                    y \\
                    s_3 \\
                    s_4 \\
                    s_5 
                    \end{bmatrix}
                    &= \begin{bmatrix}
                    1.5 \\
                    1.5 \\
                    6.0 \\
                    5.0 \\
                    21.5
                    \end{bmatrix}
\end{aligned}
\qquad
\begin{aligned}[c]
    \mb{x_N} &= \begin{bmatrix} 
                    s_1 \\
                    s_2
                    \end{bmatrix}
                    &= \begin{bmatrix}
                    0 \\
                    0
                    \end{bmatrix}
\end{aligned}
\end{equation*}

\begin{equation*}
\begin{aligned}[c]
    \mb{A_B} &= \begin{bmatrix} 
                    11           & 3  & 0 &  0 & 0 \\
                    6            & 20 & 0 &  0 & 0 \\
                    1            & 1  & 1 &  0 & 0 \\
                    \frac{-1}{3} & 1  & 0 &  1 & 0 \\
                    4            & -5 & 0 &  0 & 1 
                    \end{bmatrix}
\end{aligned}
\begin{aligned}[c]
    \mb{c_B} &= \begin{bmatrix} 
                    10 \\
                    26 \\
                    0  \\
                    0 \\
                    0 
                    \end{bmatrix}
\end{aligned}
\,.
\end{equation*}

The objective function $z = \mb{c_B}^T\mb{x_B} = 54.0$.
The basic feasible solution corresponds to the green dot in Figure~\ref{fig:simplexlp}.

\underline{Calculate the reduced cost}: $N = \{3,4\}$.

\begin{align*}
\bar{c}_3 &= c_3 - \mb{c_B}^T \mb{A_B}^{-1}\mb{A_3} = 0.2178 \\
\bar{c}_4 &= c_4 - \mb{c_B}^T \mb{A_B}^{-1}\mb{A_4} = 1.2673
\,.
\end{align*}

As both the reduced costs are positive, we have reached the optimal solution which is $x = 1.5, y = 1.5$ with an objective function value of 54.

The above procedure provides the basic framework for the simplex algorithm.
There are a few issues that can arise when applying this algorithm.
First, the issue of \emph{degeneracy}\index{optimization!linear!degeneracy}\index{degenerate solution}, which occurs when a basic variable has the value of zero. 
Usually it is the nonbasic variables that are equal to zero; but sometimes we need to set a basic variable to zero to solve the system of equations.
This can cause a few problems, and a poorly-designed implementation of the simplex algorithm may get stuck in an infinite loop in the presence of degeneracy, rotating among the same set of columns over and over again with no change in the objective function or the solution, just changing which zero variables count as basic and nonbasic.
Degeneracy can be addressed by a good tiebreaking rule that can prevent these infinite loops.
One such rule is \emph{Bland's rule}\index{Bland's rule}: among all the indices $j$ for which $\bar{c}_j$ is a negative reduced cost, choose the first of them (smallest $j$) to enter the basis; and among all the indices $k$ which achieve the minimum in~\eqn{simplexstepsize}, choose the first of them (smallest $k$) to leave the basis.

There are also a few steps which were not unambiguously specified in the algorithm above.
\begin{description}
\item [Identifying an initial basis:]
Choosing an initial basis which corresponds to a feasible solution can be difficult, especially if there are many columns.
One solution to this problem is the ``big-$M$\label{not:M} method,''\index{optimization!linear!big-$M$ method} which generates an auxiliary variable $y_i$\label{not:ylp} for each of the $m$ equality constraints.
To each of the constraints, we either add $y_i$ (if the right-hand side $b_i$ is positive) or subtract $y_i$ (if $b_i$ is negative).
In the objective function, we add $M y_i$ for each of these auxiliary variables, and we add a non-negativity constraint for each $y_i$.
By adding these variables, it is easy to get an initial solution by setting $y_i = |b_i|$ for each constraint.
If $M$ is a large enough constant, this solution will have a very poor objective function value, and \emph{any} feasible corner point solution to the original problem will have a lower objective.
So, we can start the simplex method at this solution, and as it moves to better and better solutions it will eventually remove all of the artificial variables from the basis.
The resulting optimal solution will be feasible to the original problem as well.
\item [Selecting a nonbasic variable to enter the basis:]
In the description of the simplex method, we left the choice $j$ of an entering variable arbitrary; we could choose any nonbasic variable with a negative reduced cost.
If there is more than one choice, we have to decide how to proceed, and this choice can affect the number of iterations required.
In the above example, in the first iteration, we had $\bar{c}_5 = -14$ and $\bar{c}_6 = -12$, either of them could have been chosen.
We chose the nonbasic variable corresponding to index $5$ ($s_3$) to enter the basis.
The simplex algorithm then traveled to the basic feasible solution corresponding to the yellow dot (Figure~\ref{fig:simplexlp}), and then reached the optimal solution in the next round.
Therefore, when we picked $s_3$, the simplex needed two more iterations to reach the optimal solution.
If we had instead chosen the nonbasic variable corresponding to index $6$ ($s_4$) to enter the basis, the algorithm would have traveled to the orange dot, then the black dot and finally reached the optimal solution.
Therefore, the choice of index $6$ would have resulted in one additional iteration.
Unfortunately, there is no known way to predict in advance how many iterations will be required.
There are two common methods which perform acceptably in practice: either pick the nonbasic variable with the most negative reduced cost; or among the nonbasic variables with negative reduced costs, pick the one with the smallest index.
\end{description}
\index{optimization!linear!simplex method|)}
\index{optimization!linear|)}

\section{Unconstrained Nonlinear Optimization}
\label{sec:unconstrainednonlinear}

\index{optimization!nonlinear|(}
This section discusses optimization problems of the form $\min f(\mb{x})$, where $\mb{x}$ is an $n$-dimensional real vector, with \emph{no} constraints whatsoever.
We begin by deriving optimality conditions, and then present methods for solving such problems.

\subsection{Optimality conditions}
\label{sec:unconstrainedconditions}

\index{optimization!nonlinear!optimality conditions|(}
In the main text, Section~\ref{sec:convexoptimization} gave optimality conditions for optimization problems with a convex objective function, linear equality constraints, and non-negativity constraints on the variables.
Not all optimization problems fall into these categories, and this section provides some information on these cases.
Recall that the purpose of optimality conditions is to provide mathematical equations or inequalities that can be used to check whether an optimal solution has been found or not.
In theory, it may be possible to obtain an optimal solution by solving such a system of equations and inequalities directly.
However, this is often very difficult.
Instead, the optimality conditions are mostly used in solution algorithms to know when an optimal solution has been found (or if we are close to optimality), and to provide guidance on how to improve a solution if it is not optimal.
Throughout this discussion, we make heavy reference to ``necessary'' and ``sufficient'' conditions, explained in Section~\ref{sec:necessarysufficient}.

This subsection deals with unconstrained nonlinear minimization problems of the form:
\begin{equation*}
\min \,\, \{f(\mb{x}):\mb{x} \in \bbr^n \}
\end{equation*}
The objective function may or may not be convex.
Some of the results require assumptions on differentiability, which we will state as needed.
These results will be stated without proof; readers wanting more explanation are referred to the books by~\cite{bertsekas_nlp} and~\cite{bazaara_nlp}.

This section presents several necessary and sufficient conditions (see Section~\ref{sec:necessarysufficient}) for optimal solutions of unconstrained nonlinear minimization problems.
Optimality conditions are important because they help identify if a given solution is optimal or not.
This can help in algorithm development to check if we can stop the algorithm or proceed further.
In certain specific cases, the optimality conditions can also help solve for the optimal solution or arrive at a set containing optimal solutions.

\begin{dfn}
If the function $f$ is differentiable, a \emph{stationary point}\index{stationary point} is a value $\mb{\hat{x}}$ such that $\nabla f(\hat{x}) = 0$.
\end{dfn}
In many cases, stationary point is either a local minimum or local maximum.
However, this is not always the case; for instance if $f(x) = x^3$, then $\hat{x} = 0$ is a stationary point, but the function has neither a minimum or maximum there.
A stationary point which is neither a local minimum or a local maximum is called a \emph{saddle point}\index{saddle point}.

\begin{thm}
\emph{(First-order necessary conditions for local minima.)}
If $f$ is differentiable at a local minimum point $\mb{\hat{x}}$, then $\mb{\hat{x}}$ is also a stationary point.
\end{thm}
This result is a necessary condition; it is ``first-order'' because it refers to the first derivative.
Therefore, if $\mb{\hat{x}}$ is a local minimum, then it is a stationary point.
But stationary points need not be local minima; they could also be local maxima or saddle points, for instance.

\begin{exm}
Determine the stationary points of the function $ f(x) = (x-5)^2$. Does a minimum exist? 
\end{exm}
\solution{
\begin{equation*}
\nabla f(\hat{x}) = 2(\hat{x}-5) = 0 \,\,\implies \hat{x} = 5.
\end{equation*}
Thus $\hat{x} = 5$ is the stationary point. By plotting the graph, one can easily determine that the stationary point corresponds to a local as well as global minimum.
}

\begin{exm}
Determine the stationary points of the function $f(x) = x^3 + x^2 -x + 1$. Does a minimum exist? 
\end{exm}
\solution{
\begin{equation*}
\nabla f(\hat{x}) = 3(\hat{x})^2 + 2\hat{x} - 1 = 0 \,\,\implies \hat{x} = 1/3, -1.
\end{equation*}
If you plot the graph, you will notice that $\hat{x} = 1/3$ corresponds to a local minimum and $\hat{x} = -1$ corresponds to a local maximum.
}

\begin{exm}
Determine the stationary points of the function $f(x_1,x_2) = 4(x_1 - 7)^2 + 5x_2^2(x_2-10)$.
\end{exm}
\solution{
\begin{equation*}
\frac{\partial f}{ \partial \hat{x}_1} = 8(\hat{x}_1 - 7) = 0 \implies \hat{x}_1 =7
\end{equation*}
\begin{equation*}
\frac{\partial f}{ \partial \hat{x}_2} = 3(\hat{x}_2)^2 - 20\hat{x}_2 = 0 \implies \hat{x}_2 = 0, 20/3.
\end{equation*}
The two stationary points are $(7,0)$ and $(7, 20/3)$.
However, we do not have enough information to determine if the stationary points are minima, maxima, or saddle points.
}

If the first order necessary condition does not provide clarity on whether a stationary point is a minimum or not, a second-order necessary condition can be used if the function is twice differentiable:
\begin{thm}
\emph{(Second-order necessary conditions for local minima.)}
If $f$ is twice differentiable at a point $\mb{\hat{x}}$ which is a local minimum, then (i) $\mb{\hat{x}}$ is a stationary point, and (ii) the Hessian matrix\index{Hessian matrix!applications} $Hf(\mb{\hat{x}})$ is positive semidefinite.
\end{thm}
The theorem can be extended as follows: if $\mb{\hat{x}}$ is a local minimum, then $\mb{\hat{x}}$ is stationary and the Hessian is negative semidefinite there; if $\mb{\hat{x}}$ is a saddle point, then $\mb{\hat{x}}$ is stationary and the Hessian is neither positive semidefinite nor negative semidefinite.
Again, these are \emph{necessary} conditions.
They must be true of every optimal solution, but not every point satisfying these conditions is optimal.

\begin{exm}
Consider the function $f(x) = (x-5)^5$. Identify the stationary point and check if the Hessian is positive semidefinite.
\end{exm}
\solution{
\begin{equation*}
\nabla f(\hat{x}) = 5(\hat{x} - 5)^4 = 0 \,\,\implies \hat{x} = 5.
\end{equation*}
The stationary point is $\hat{x} = 5$.
The Hessian at the stationary point is $f''(\hat{x}) = 20(\hat{x} - 5)^3$, which is nonnegative (actually zero) and thus positive semidefinite.
However if you plot the function, you will notice that $\hat{x} = 5$ does not correspond to a local minima, but instead a saddle point.
}

\begin{exm}
Consider the stationary point $(7,0)$ of the function $f(x_1,x_2) = 4(x_1 - 7)^2 + 5x_2^2(x_2-10)$. Is $(7,0)$ is a saddle point? 
\end{exm}
\solution{
At $(7,0)$ we have
\begin{align*}
\frac{\partial^2 f}{ \partial x_1^2} &= 8 \\
\frac{\partial^2 f}{ \partial x_1 \partial x_2} &= 0 \\
\frac{\partial^2 f}{ \partial x_2 \partial x_1} &= 0 \\
\frac{\partial^2 f}{ \partial x_2^2} &= 6x_2 - 20 
\,,
\end{align*}
so the Hessian is given by
\begin{equation*}
\begin{aligned}[c]
    H f &= \begin{bmatrix} 
                       8 &  0  \\                                  
                       0  &  -20 
           \end{bmatrix}
\end{aligned}
\,.
\end{equation*}
The Hessian is neither positive nor negative semidefinite.
Therefore, the stationary point $(7,0)$ is a saddle point.
}

A minor modification of the second-order necessary conditions gives a \emph{sufficient} condition on optimality.
\begin{thm}
\emph{(Second-order sufficiency conditions for a strict local minimum.)}
If the function $f$ is twice-differentiable at a point $\mb{\hat{x}}$, then $\mb{\hat{x}}$ is a strict local minimum if (i) $\mb{\hat{x}}$ is a stationary point and (ii) the Hessian matrix $H(\mb{\hat{x}})$ is positive definite.
\end{thm}
Likewise, if $H$ is negative definite, then we know $\mb{\hat{x}}$ is a strict local maximum.
As a sufficient condition, any point satisfying these conditions must be a local minimum.
It is not necessary; for instance if $f(x) = x^2$, then $\hat{x} = 0$ is a strict local minimum, but the Hessian is not positive definite there.

\begin{exm}
Verify if the stationary points of the function $f(x) = x^3 + x^2 -x + 1$ are local minima.
Does a minimum exist? 
\end{exm}
\solution{
The stationary points of the function are $\hat{x} = 1/3$ and $\hat{x} = -1$.
\begin{equation*}
H f(\hat{x}) = 6(\hat{x}) + 2 
\end{equation*}
The Hessian is positive definite at $\hat{x} = 1/3$ and negative definite at $\hat{x} = -1$.
Therefore, $\hat{x} = 1/3$ is a strict local minimum and $\hat{x} = -1$ is a strict local maximum.
}

\begin{exm}
Consider the stationary point $(7,20/3)$ of the function defined by $f(x_1,x_2) = 4(x_1 - 7)^2 + 5x_2^2(x_2-10)$.
Is $(7,20/3)$ is a strict local minimum? 
\end{exm}
\solution{
At $(7,20/3)$ we have
\begin{align*}
\frac{\partial^2 f}{ \partial x_1^2} &= 8 \\
\frac{\partial^2 f}{ \partial x_1 \partial x_2} &= 0 \\
\frac{\partial^2 f}{ \partial x_2 \partial x_1} &= 0 \\
\frac{\partial^2 f}{ \partial x_2^2} &= 6x_2 - 20
\,,
\end{align*}
so the Hessian is
\begin{equation*}
\begin{aligned}[c]
F f &= \begin{bmatrix} 
               8 &  0  \\                                  
               0  &  20 
               \end{bmatrix}
\end{aligned}
\,.
\end{equation*}
This matrix is positive definite, so $(7,20/3)$ is a strict local minimum.
}

Recall from Section~\ref{sec:convexoptimization} that if $f$ is a convex function, then any local minimum must also be a global minimum, and that if $f$ is strictly convex, any local minimum is the (unique) strict global minimum.
\index{optimization!nonlinear!optimality conditions|)}

\subsection{Solution framework}
\label{sec:unconstrainedsolution}

\index{optimization!nonlinear!algorithms|(}
Unconstrained optimization problems are commonly solved using the following algorithmic framework; the details of steps 1--3 are filled in later in this subsection.
\begin{description}
\item[Step 0: Initialize.]  Set the iteration counter $k = 0$, and initialize $\mb{x_0}$ by choosing any vector in $\bbr^n$.  (If you have a good guess as to the minimum point, it can greatly speed things up.)
\item[Step 1: Check convergence.]  Evaluate one of the termination criteria (Section~\ref{sec:unconstrainedconvergence}.)
If it is satisfied, return $\mb{x_k}$ as an approximately-optimal solution, and stop.
\item[Step 2:] Determine a descent direction $\mb{d_k}$.  (Section~\ref{sec:unconstraineddescentdir}; this step is where gradient descent and Newton's method differ.)
\item[Step 3:] Determine a step size $\mu_k$. (Section~\ref{sec:unconstrainedstepsize}.)
\item[Step 4: Iterate.]  Update the solution, $\mb{x_{k+1}} \leftarrow \mb{x_k} + \mu_k \mb{d_k}$, then increase $k$ by 1, and return to step 1.
\end{description}

The following subsections spell out choices for Steps 1, 2, and 3.
These choices are independent of each other, and you can combine any choice for one step with any choice for another.

\subsection{Convergence criteria}
\label{sec:unconstrainedconvergence}

\index{optimization!nonlinear!convergence criteria|(}
There are several stopping criteria available to decide when to terminate.
It is always best to select a stopping criterion to directly measure what we are trying to do.
In this case, we know that the minimum of $f$ must occur at a stationary point, which means $\nabla f = \mb{0}$.  
So, a natural stopping criteria is to terminate when $\nabla f$ is ``close enough'' to zero, that is, when $|| \nabla f || < \epsilon$.

One difficulty with this choice is that you need a good intuition for what values of $||\nabla f||$ correspond to ``good enough'' solutions for your particular application, and this is not always easy.
So, there are other convergence criteria that can also be used:
\begin{itemize}
\item You can stop when the solution stabilizes, $||\mb{x_{k}} - \mb{x_{k-1}}|| < \epsilon$.
\item You can stop when the objective function stabilizes, $|f(\mb{x_{k}}) - f(\mb{x_{k-1}})| < \epsilon$.
\item You can normalize the previous two inequalities to reflect relative stability (e.g., stop when the objective function decreases by less than 1\% between iterations).
\item You can stop after a pre-specified number of iterations.
\item You can stop after a certain amount of computation time has elapsed.
\end{itemize}
These methods are more intuitive to apply than $|| \nabla f || < \epsilon$.
The downside is that there is no guarantee that you are close to an optimal solution when they are satisfied.
This is particularly obvious for the last criteria; but even for the earlier ones, there is no way to tell the difference between a solution ``stabilizing'' for a good reason (you are close to the minimum) or for a bad reason (the algorithm is stuck somewhere suboptimal but can't make progress).

You can also use a combination of these methods; for instance, stopping when the objective changes by less than 1\%, or after one hour of run time (whichever comes first).
\index{optimization!nonlinear!convergence criteria|)}

\subsection{Descent direction}
\label{sec:unconstraineddescentdir}

There are two main ways to select the direction $\mb{d_k}$; the idea in both cases is that the gradient points in the direction where the function \emph{increases} as quickly as possible.
Therefore, the negative of the gradient is the direction where the function decreases as quickly as possible.
Since we are trying to minimize the objective, we will use the latter in some way.

In \emph{gradient descent}\index{optimization!nonlinear!gradient descent}, we simply use the gradient itself: $\mb{d_k} = -\nabla f(\mb{x_k})$.  
This method is simple and requires relatively little computation.

A more sophisticated choice is \emph{Newton's method}\index{optimization!nonlinear!Newton's method}, which is based on the following logic: we can take a quadratic approximation to the objective function at the current point.
We can identify the minimum point of that quadratic approximation more easily than the original function $f$ (which may be much more complicated), and then choose the search direction $\mb{d_k}$ to move towards that minimum point.
We will derive this direction, and then explain how it is related to Newton's method as you learned it in calculus, or as we used it in the previous section.

The quadratic approximation to $f$ at $\mb{x_k}$ is its second-order Taylor series, based on its gradient and Hessian:
\labeleqn{quadraticapprox}{
f(\mb{x}) \approx f(\mb{x_k}) + (\mb{x} - \mb{x_k})^T \nabla f(\mb{x_k}) + \frac{1}{2}(\mb{x} - \mb{x_k})^T H f(\mb{x_k}) (\mb{x} - \mb{x_k})^T
\,.
}
Let $\tilde{f}$ denote the right-hand side of equation~\eqn{quadraticapprox}.
The minimum point of the quadratic approximation is the point where $\nabla \tilde{f}$ vanishes, that is, where
\labeleqn{quadraticzero}{
\nabla \tilde{f}(\mb{x}) = \nabla f(\mb{x_k}) + Hf(\mb{x_k}) (\mb{x} - \mb{x_k}) = 0
\,.
}
Solving for $\mb{x}$, we obtain
\labeleqn{quadraticmin}{
\mb{x} = \mb{x_k} - (Hf(\mb{x_k}))^{-1} \nabla f(\mb{x_k})
\,.
}
So, the search direction is
\labeleqn{newtonsearch}{
\mb{d_k} = \mb{x} - \mb{x_k} = -(Hf(\mb{x_k}))^{-1} \nabla f(\mb{x_k})
\,.
}

Why is this method related to Newton's method?
Newton's method aims to find the zero of a function by finding the zero of its linear approximation.
In this case, we are trying to find the minimum point of $f$, which occurs where $\nabla f = 0$.
So, when we use Newton's method, we find where the linear approximation of $\nabla f$ vanishes.
But if the gradient of a function is linear, then the original function is quadratic, and the gradient will be zero at the minimum point of this quadratic.
Also recall the previous section, where in a one-dimensional setting we used Newton's method to try to find the point where $f'(x) = 0$, and ended up reducing $x$ by $f'(x) / f''(x)$.
The formula $\mb{d_k} = -(Hf(\mb{x_k}))^{-1} \nabla f(\mb{x_k})$ is analogous: the gradient is the generalization of the first derivative; the Hessian is the generalization of the second derivative; and multiplying by the inverse of a matrix is the generalization of dividing by a scalar (just as multiplying by $2^{-1}$ is the same as dividing by $2$).
Indeed, in the single-dimensional case, equation~\eqn{newtonsearch} simplifies to
\labeleqn{newtonsearch1d}{
d_k = -f'(x_k)/f''(x_k)
\,,
}
so the new point $x_{k+1} = x_k + d_k$ is given by exactly the same formula as in Step 2 of Newton's method for line search in Section~\ref{sec:newtonlinesearch}

Newton's method usually yields a ``better'' search direction, and fewer iterations are required to reach the minimum.
However, it can be computationally expensive to calculate all the elements in the Hessian matrix, and to calculate its inverse.
For this reason, there are a variety of ``quasi-Newton''\index{Newton's method!quasi} methods, which replace $H^{-1}$ in equation~\eqn{newtonsearch} by another matrix which is easier to calculate, but is approximately the same.

\subsection{Step size}
\label{sec:unconstrainedstepsize}

There are several procedures available to determine the step size $\mu_k$.
The simplest is to use a \emph{constant step size}\index{optimization!nonlinear!constant step size}, where at each iteration $\mu_k = \mu$, for some pre-determined constant $\mu$.
If the predetermined value is too big, the algorithm may overshoot the minimum, slowing convergence or even diverging.
If the step size is too small, the algorithm may take a long time to converge.
If you are determining the search direction using Newton's method, a choice of $\mu = 1$ may actually perform quite well, since this will move you exactly to the minimum point of the quadratic approximation.
This choice will probably not work well for gradient descent, since the ``units'' of the gradient are not the same as for the decision variables.  
So, if you use this approach, you should have a good intuition about the objective function.

A more sophisticated method is \emph{exact line search}\index{optimization!nonlinear!exact line search}, where you choose the value of $\mu_k$ so that $f(\mb{x_k} + \mu_k \mb{d_k})$ is minimized.
This is a one-dimensional optimization problem, which can be solved using the methods of Section~\ref{sec:linesearch}.\footnote{You will need to impose an upper bound on $\mu_k$ to do this, but in practice this is not usually very hard.}
While this choice yields the greatest improvement in the objective function at each iteration, it can be computationally expensive to compute this minimum.
(You would have to use the line search method at every iteration.)

A compromise between the two is a \emph{backtracking line search}\index{optimization!nonlinear!backtracking line search}, which also tries to minimize $f(\mb{x_k} + \mu_k \mb{d_k})$, but in a loose way.
Rather than hoping to find the exact minimum, we choose a value of $\mu_k$ which is ``good enough.''
Select a positive initial step size $\eta$, and two other constants $\beta$ and $\gamma$, both strictly between 0 and 1.
We then test the sequence of $\mu_k$ values $\mu_k = \eta, \beta \eta, \beta^2 \eta, \ldots$, stopping as soon as this inequality is satisfied:
\labeleqn{armijo}{
f(\mb{x_k} + \mu_k \mb{d_k}) - f(\mb{x_k}) <  \gamma \mu_k \nabla f(\mb{x_k})^T \mb{d_k}
}
The first value $\eta, \beta \eta, \beta^2 \eta, \ldots$ which satisfies this equation is chosen for $\mu_k$.
The intuition in formula~\eqn{armijo} is that the left-hand side shows how much $f$ changes if we take a step of size $\mu_k$.
We are trying to solve a minimization problem, so hopefully $f$ decreases, and the left-hand side is negative.
On the right-hand side, $\mu_k \nabla f(\mb{x_k})^T \mb{d_k}$ is how much we would expect $f$ to decrease based on its linear approximation.
Likewise, since $\mb{d_k}$ is a direction in which $f$ decreases, the right-hand side is also a negative number.
We stop at the first choice of $\mu_k$ for which the actual decrease in the objective function is at least a certain fraction $\gamma$ of what we would expect from the linear approximation; this is exactly what the condition~\eqn{armijo} checks.
Typical values of the constants in this method are $\eta = 1$, $\beta = 1/2$, and $\gamma = 1/10$, but you should experiment with different values for your specific problem.
This rule is often called the \emph{Armijo rule}\index{Armijo rule}\index{optimization!nonlinear!Armijo rule}.

\subsection{Examples}
\label{sec:unconstrainedexampls}

\begin{exm}
Apply the unconstrained optimization algorithm to the function $f(x) = x^2 - 10x + 20$.
Terminate when $|x_{k} - x_{k-1}| < 0.01$, use gradient descent for the direction, and use a constant step size of $\mu = 1$.
Repeat with a constant step size of $\mu = 0.1$.
\end{exm}
\solution{
From basic calculus we know that the minimum of the above function is at $\hat{x} = 5$.
For this example, we pick an initial value of $x_0 = 15$.  
At any point $x_k$, the descent direction is $d_k = - \nabla f(x_k) = 2x_k - 10$.
Then, using a constant step size of $\mu = 1$, we have the following: 

At the initial point, the descent direction is $d_0 = -20$, so $x_1 = 15 - 1 \times 20 = -5$.
Proceeding to the next iteration, we check the termination criterion.
Since $|x_1 - x_0| > 0.01 $, we continue.
The new descent direction is $d_1 = -20$. Therefore, $x_2 = -5  + 1 \times 20 = 15$.
Since $|x_2 - x_1| > 0.01 $, we proceed to the next iteration.
But we've returned to our initial point!
Notice that the solutions will continue to oscillate between $15$ and $-5$, due to the large step size.

Repeating using a smaller constant step size of $\mu = 0.1$ produces convergence: from the initial point, the descent direction is $d_1 = -20$ and $x_1 = 15  - 0.1 \times 20 = 13$.
We have $|x_1 - x_0| > 0.01 $, so we continue.
At this new point, the descent direction is $d_1 = -16$. Therefore, $x_2 = 13  - 0.1 \times 16 = 11.4$. 
Subsequent iterations are shown in Table~\ref{tbl:descentmethod_2}.
}

\begin {table}
\begin{center}
\caption{Gradient descent applied to $f(x) = x^2 - 10x + 20$, with a constant step size $\mu_k = 0.1$ \label{tbl:descentmethod_2}}
\begin{tabular}{| c | c |c|  }
\hline
$k$ &  $x_k$  & $f(x_k)$ \\
\hline
0  &  15  &  100 \\
1  &  13.0  &  64.0 \\
2  &  11.4  &  40.96 \\
3  &  10.12  &  26.2144 \\
4  &  9.096  &  16.7772 \\
5  &  8.2768  &  10.7374 \\
6  &  7.6214  &  6.8719 \\
7  &  7.0972  &  4.3980 \\
8  &  6.6777  &  2.8147 \\
9  &  6.3422  &  1.8014 \\
10  &  6.0737  &  1.1529 \\
11  &  5.8590  &  0.7379 \\
12  &  5.6872  &  0.4722 \\
13  &  5.5498  &  0.3022 \\
14  &  5.4398  &  0.1934 \\
15  &  5.3518  &  0.1238 \\
16  &  5.2815  &  0.0792 \\
17  &  5.2252  &  0.0507 \\
18  &  5.1801  &  0.0325 \\
19  &  5.1441  &  0.0208 \\
20  &  5.1153  &  0.0133 \\
21  &  5.0922  &  0.0085 \\
22  &  5.0738  &  0.0054 \\
23  &  5.0590  &  0.0035 \\
24  &  5.0472  &  0.0022 \\
25  &  5.0378  &  0.0014 \\
\hline
\end{tabular}
\end{center}
\end{table}

In this example, the constant step size plays a major role in how quickly we converge (if at all).
For reference, Table~\ref{tbl:descentmethod_3} shows how many iterations are required for different choices.
Similarly, the initial value chosen, and the parameters $\eta, \beta, \gamma$ of an inexact line search play an important role in convergence.

\begin {table}
\begin{center}
\caption{Number of iterations for gradient descent to minimize $f(x) = x^2 - 10x + 20$ with different step sizes \label{tbl:descentmethod_3}}
\begin{tabular}{| c | c |c | }
\hline
$\mu_k$ &  Iterations needed \\
\hline
0.1 & 27 \\
0.2 & 14 \\
0.3 & 9 \\
0.4 & 7 \\
0.5 & 3 \\
0.6 & 7 \\
0.7 & 10 \\
\hline
\end{tabular}
\end{center}
\end{table}

\begin{exm}
Apply the unconstrained optimization algorithm to the function $f(x,y) = (x-1)^4 + 5(y-2)^4 + xy$.
Terminate when $||\mb{x_k} - \mb{x_{k-1}}|| < 0.01$, and use gradient descent.
First use the algorithm with a constant step size of $\mu = 0.0025$, then solve again using backtracking line search with $\eta = 1$, $\beta = 0.1$, and $\gamma = 0.1$.
\end{exm}
\solution{
Notice that there are two decision variables; we will use the vector $\mb{x} = (x, y)$ to describe both decision variables together.
For a specific iteration, we will let their values be given by $\mb{x_k} = (x_k, y_k)$.

For this example, we pick an initial value of $(x_1,y_1) = (4,4)$.
At any point $k$, the descent direction is
\[
\mb{d_k} = - \nabla f(x_k,y_k) =  \vect{
                                -4(x_k-1)^3 - y_k \\
                                -20(y_k-2)^3 - x_k 
                                 }
\,.                              
\]

So from the initial point, the descent direction is
\[
\mb{d_0} = - \nabla f(x_0,y_0) =  \vect{ -4(x_0-1)^3 - y_0 \\ -20(y_0-2)^3 - x_0 }  = \vect{ -112 \\ -164 }
\]
and the new point is
\[
\mb{x_1} = \vect{ x_1 \\ y_1 }  = \vect{ 4 \\ 4 } +  0.0025 \vect{ -112 \\ -164 }  = \vect{ 3.72 \\ 3.59 }
\,. 
\]
For convergence, we check if $\sqrt{ (x_1 - x_0)^2 + (y_1 - y_0)^2} < 0.01$.
This is not true, so we increase $k$ to 1 and move to the next iteration.

The new descent direction is
\[
\mb{d_1} = - \nabla f(x_1,y_1) =  \vect{ -4(x_1-1)^3 - y_1 \\ -20(y_1-2)^3 - x_1 }  = \vect{ -84.04 \\ -84.11 }    
\]
and the new point is
\[
\mb{x_2} = \vect{ x_2 \\ y_2 }  = \vect{ 3.72 \\ 3.59 } +  0.0025 \vect{ -84.08 \\ -84.11 }  = \vect{ 3.51 \\ 3.38 } 
\,.
\]
Since $\sqrt{ (x_2 - x_1)^2 + (y_2 - y_1)^2} > 0.01$, we increase $k$ to 2 and continue.
The algorithm converges after 64 iterations, at the point $(1.66, 2.23)$.
The objective function has a value of 3.91.

Solving again using backtracking line search, we have the same initial points and directions: $\mb{x_0} = (4,4)$ and $\mb{d_0} = (-112,-164)$.
However, we have to do a little more work to determine the step size $\mu_0$.
We start by testing $\mu_0 = \eta = 1$.
With this choice, the new solution would be
\[
\vect{4 \\ 4} + 1 \vect{-112 \\ -164} = \vect{-108 \\ -160}
\,.
\]
The value of the objective function here is very large, approximately $3.58 \times 10^9$, so the left-hand side of~\eqn{armijo} is $3.58 \times 10^9 - 177 \approx 3.58 \times 10^9$.
For the right hand side, we calculate
\[
\gamma \mu_0 \nabla f(\mb{x_0})^T \mb{d_0} = (0.1)(1) \vect{112 & 164} \vect{-112 \\ -164} = -3944
\,.
\]
Inequality~\eqn{armijo} is clearly false (the left-hand side is very positive, the right-hand side is negative), so we try again with $\mu_0 = \beta \eta = 1/10$.
This choice corresponds to the solution
\[
\vect{4 \\ 4} + \frac{1}{10} \vect{-112 \\ -164} = \vect{-7.2 \\ -12.4}
\,.
\]
The objective function has a value of $2.19 \times 10^5$, so the left-hand side of~\eqn{armijo} is still approximately $2.19 \times 10^5$, while the right-hand side is
\[
\gamma \mu_0 \nabla f(\mb{x_0})^T \mb{d_0} = \frac{1}{100} \vect{112 & 164} \vect{-112 \\ -164} = -394.4
\,.
\]
which is again false.
Trying again with $\mu_0 = \beta^2 \eta = 1/100$, the new solution is $\mb{x} = (2.88, 2.36)$, and the left and right-hand sides of~\eqn{armijo} are now $-157.6$ and $-39.4$, respectively.
So we accept this choice: $\mu_0 = 1/100$, and $\mb{x_1} = (2.88, 2.36)$.

Proceeding similarly, you can verify that the next step sizes are $\mu_1 = 1/10$, $\mu_2 = 1/10$, $\mu_3 = 1$, and $\mu_4 = 0.1$, with the algorithm terminating after that step at the solution $\mb{x_5} = (0.236, 1.774)$, with objective value $0.772$.
Notice that each iteration of backtracking line search required more work, but in the end we only had to perform five iterations, rather than sixty-four.
}

\begin{exm}
Apply the unconstrained optimization method to the function $f(x,y) = (x-1)^4 + 5(y-2)^4 + xy$.
Terminate when $||\mb{x_k} - \mb{x_{k-1}}|| < 0.01$, and use Newton's method with a constant step size of $\mu = 1$.
\end{exm}
\solution{
We will use the same initial value $\mb{x_0} = (4, 4)$.
As before, the gradient at an arbitrary point is
\[
\nabla f(x_k,y_k) =  \vect{
                        4(x_k-1)^3 + y_k \\
                        20(y_k-2)^3 + x_k 
                     }
\,,
\]
and the Hessian is
\[
H f(x_k,y_k) =  \vect{
                        12(x_k-1)^2   &  1 \\
                        1             & 60(y_k -2)^2 
                     }
\,.                                         
\]

So, at the initial point we have
\[
\nabla f(4.0, 4.0) =   \vect{ 4(4-1)^3 + 4 \\  20(4-2)^3 + 4 } = \vect{ 112 \\ 164 } 
\]
and
\[
H f(4.0,4.0) =  \vect{
                        108   &  1 \\
                        1     & 240 
                     }
\,.                     
\]
The search direction from Newton's method is
\[
\mb{d_0} = (Hf(\mb{x_1}))^{-1} \nabla f(\mb{x_0}) = \vect{ -1.0307 \\  -0.6790 }
\,.
\]
Using the constant step size $\mu = 1$, the new point is
\[
\vect{ x_1 \\ y_1 }  = \vect{ 4 \\ 4 } +  1 \vect{ -1.0307 \\ -0.6790 }  = \vect{ 2.9692 \\ 3.3209 } 
\,.
\]
Since $\sqrt{ (x_1 - x_0)^2 + (y_1 - y_0)^2} > 0.01$, we increment $k$ to 1 and move to the next iteration.

At $(x_1,y_1) = (2.9692, 3.3209)$, the gradient and Hessian are
\[
\nabla f(2.9692, 3.3209) =   \vect{ 4(2.9692-1)^3 + 3.3209 \\  20(3.3209-2)^3 + 2.9692 } = \vect{ 38.8672 \\ 49.069 } 
\]
and
\[
H f(2.9692, 3.3209) =  \begin{bmatrix}
                        46.5353   &  1 \\
                        1     &  104.696 
                     \end{bmatrix}
\,,
\]
respectively.
So the search direction is
\[
\mb{d_1} = -(H f(\mb{x_1}))^{-1} \nabla f(\mb{x_1}) = \vect{ -0.7178 \\  -0.4618 }
\]
and
\[
\vect{ x_2 \\ y_2 }  = \vect{ 2.9692 \\ 3.3209 } +  1 \vect{ -0.7178 \\  -0.4618 }  = \vect{ 2.2513 \\ 2.8591 } 
\,,
\]
and so forth.
The algorithm terminates after sixteen iterations with the solution $\mb{x} = (0.2399, 1.756)$, at an objective function value of 0.7728.
}
\index{optimization!nonlinear!algorithms|)}

\section{Constrained Nonlinear Optimization}

The main text derived optimality conditions for a specific class of nonlinear optimization problems including the static traffic assignment problem (Section~\ref{sec:convexoptimization}).
Likewise, Sections~\ref{sec:convexmsa} and~\ref{sec:frankwolfe} introduced the method of successive averages and Frank-Wolfe algorithms for solving traffic assignment.
This section discusses how these techniques can be applied to other kinds of nonlinear optimization problems involving linear constraints.

In this section, we will consider optimization problems of the form
\optimizex{\min}{\mb{x}=(x_1, \ldots, x_n)}{f(\mb{x})}{
&  g_i(\mb{x}) \leq 0 & \forall  i=1,2,\ldots,m  \\
&  h_j(\mb{x}) = 0    & \forall  j=1,2,\ldots,\ell    
}In this formulation, $g_i$ and $h_j$ are functions representing the inequality and equality constraints, respectively.
We will concisely denote the feasible region by $X$.
That is, $X$ is the set
\[
X = \myc{\mb{x} \in \bbr^n: g_i(\mb{x}) \leq0\,\,\forall i=1,2,\ldots,m ; h_j(\mb{x}) = 0\,\,  \forall  j=1,2,\ldots,\ell}
\,.
\]
We first derive optimality conditions that apply to a general constrained nonlinear optimization problem, and then specialize them to the case where all of the constraints are linear and the objective function is convex.

\subsection{Optimality conditions}
\label{sec:constrainedoptimalityconditions}

\index{optimization!nonlinear!optimality conditions|(}
As with unconstrained optimization, there are first- and second-order necessary and sufficient conditions characterizing the local minima of constrained optimization problems.
The main distinction is that optimal solutions can now lie along the boundary of the feasible region, as well as in its interior.
It is possible to create examples where constraints behave irregularly.
To address this, we will need to define \emph{constraint qualifications}\index{constraint qualification}, which are essentially regularity conditions on the constraints.  
We begin with some definitions leading up to these.

\begin{dfn}
For any feasible solution $\mb{x} \in X$, an \emph{active} constraint\index{optimization!constraints!active} is an inequality constraint which is satisfied with equality.
That is, the $i$-th inequality constraint is active if $g_i(\mb{x}) = 0$.  
\end{dfn}

\begin{dfn}
For any feasible solution $\mb{x} \in X$, the \emph{active constraint set}\index{optimization!constraints!active set} is the set of all active constraints at $\mb{x}$: $A(\mb{x}) = \{i: g_i(\mb{x}) = 0\}$.\label{not:Ax} 
\end{dfn}

We now introduce several constraint qualifications, which we will discuss below.

\begin{dfn}
The \emph{linear independence constraint qualification}\index{constraint qualification!linear independence} holds at a feasible solution $\mb{x} \in X$ if the gradient vectors $\nabla g_i(\mb{x})$ for active constraints $i \in A(\mb{x})$, and for all equality constraints $\nabla h_j(\mb{x})$, are linearly independent.
\end{dfn}

\begin{dfn}
The \emph{Mangasarian-Fromovitz constraint qualification}\index{constraint qualification!Mangasarian-Fromovitz} holds at a feasible solution $\mb{x} \in $ if (i) the gradient vectors $\nabla h_j(\mb{\hat{x}})$ are linearly independent, and (ii) there exists a direction vector  $\mb{d} \in \bbr^n$ such that $\nabla g_i(\mb{\hat{x}})^T \mb{d} < 0$ for every $i \in A(\mb{x})$ and $\nabla h_j(\mb{\hat{x}})^T \mb{d} = 0$ for all $j = 1, \ldots, \ell$.
\end{dfn}

The linear independence constraint qualification is stronger, and you can show that whenever it holds, the Mangasarian-Fromovitz constraint qualification also holds.  
There are other constraint qualifications as well.
If the linear independence constraint qualification holds, then most other constraint qualification conditions hold.
In many nonlinear optimization formulations involving transportation networks, all of the constraints are linear, and a simpler constraint qualification holds:

\begin{dfn}
The \emph{linearity constraint qualification}\index{constraint qualification!linearity} holds if all of the functions $g_i$ and $h_j$ defining the constraints are affine functions.
\end{dfn}

Perhaps the most famous necessary optimality conditions for constrained optimization are the Karush-Kuhn-Tucker (KKT) conditions\index{Karush-Kuhn-Tucker (KKT) conditions|(}.
They make use of auxiliary variables $\bm{\iota}$\label{not:iota} and $\bm{\kappa}$, often called \emph{Lagrange multipliers}\index{Lagrange multipliers}.\footnote{Many authors use $\lambda$ and $\mu$ where we use $\kappa$ and $\iota$.
In this book, we have chosen to reserve $\lambda$ to mean a ``convex combinations'' step size between 0 and 1, and $\mu$ to mean a step size without such a limitation.
}

\begin{thm}
\emph{(First-order necessary conditions for a local minimum.)}
Let $f, g_i, h_j$ be continuously differentiable functions.
If (i) $\mb{\hat{x}}$ is a local minimum of $f$, and (ii) any of the above constraint qualifications hold at $\mb{\hat{x}}$, then there exist scalars $\iota_i$, $i = 1, 2, \ldots, m$ and $\kappa_j$, $j = 1, \ldots, \ell$, such that the following system holds:
\begin{align}
\nabla f(\mb{\hat{x}}) + \sum_{i=1}^{m} \iota_i \nabla g_i(\mb{\hat{x}}) + \sum_{j=1}^{\ell} \kappa_j \nabla h_j(\mb{\hat{x}}) &= \mb{0}  \label{eqn:NecKKT1}\\
h_j(\mb{\hat{x}}) &= 0 \,\, \forall j = 1,2,\ldots,\ell \label{eqn:NecKKT2} \\
g_i(\mb{\hat{x}}) &\leq 0 \,\, \forall i = 1,2,\ldots,m \label{eqn:NecKKT3} \\
\iota_i g_i(\mb{\hat{x}}) &= 0 \,\, \forall i = 1,2,\ldots,m \label{eqn:NecKKT4} \\
\iota_i &\geq 0 \,\, \forall i = 1,2,\ldots,m \label{eqn:NecKKT5}
\,.
\end{align}
\end{thm}

The above conditions can also be written in matrix form.
Let $\bm{\iota}$ and $\bm{\kappa}$ be the vectors whose components are $\iota_1,\iota_2,\ldots,\iota_m$ and $\kappa_1,\kappa_2,\ldots,\kappa_\ell$, and likewise let $\mb{g(x)}$ and $\mb{h(x)}$ be vectors with components $g_i(\mb{x})$ and $h_j(\mb{x})$.
Then we can rewrite the KKT conditions in terms of the Jacobians of $\mb{g}$ and $\mb{h}$ as
\begin{align*}
\nabla f(\mb{\hat{x}}) + J \mb{g}(\mb{\hat{x}}) \bm{\iota} + J \mb{h}(\mb{\hat{x}}) \bm{\kappa} &= \mb{0} \\
\mb{h}(\mb{\hat{x}}) &= \mb{0} \\
\mb{g}(\mb{\hat{x}}) &\leq \mb{0}  \\
\bm{\iota}^T \mb{g}(\mb{\hat{x}}) &= 0 \\
\bm{\iota} &\geq \mb{0} 
\,.
\end{align*}

Another  way of representing the first order necessary condition is using the Lagrangian function
\begin{align*}
\mc{L}(\mb{x},\bm{\iota},\bm{\kappa}) = f(\mb{\hat{x}}) + \sum_{i=1}^{m} \iota_i g_i(\mb{\hat{x}}) + \sum_{j=1}^{l} \kappa_j h_j(\mb{\hat{x}})
\,.
\end{align*}
We can rewrite the KKT functions in a simpler way using the Lagrangian.
Taken as a vector, the partial derivatives of $\mc{L}$ with respect to $\mb{x}$ (written as $\nabla_\mb{x} \mc{L}$) form the left-hand side of the first KKT condition~\eqn{NecKKT1}.
Similarly, the partial derivatives of $\mc{L}$ with respect to $\bm{\lambda}$ form the left-hand side of the second KKT condition~\eqn{NecKKT2}.
So we can write the KKT conditions as:
\begin{align*}
\nabla_\mb{x} \mc{L}(\mb{\hat{x}}, \bm{\iota}, \bm{\kappa}) &= \mb{0} \\
\nabla_{\bm{\lambda}} \mc{L}(\mb{\hat{x}}, \bm{\iota}, \bm{\kappa}) &= \mb{0} \\
\,.
\end{align*}

These are necessary conditions, meaning that they must be satisfied at any optimal solution.
Without additional restrictions, they are not sufficient, meaning that there may be non-optimal points which also satisfy those conditions.
Still, we can identify \emph{all} the points which satisfy the KKT conditions to generate a set of ``candidate solutions;'' if an optimal solution exists it must be one of them.
We give some examples of how to do this in the following subsection.
For large-scale problems this is not a practical approach (the following subsection describes some methods that can be used), but it may be reasonable for problems with only a few decision variables and constraints, or where there is a special structure which further simplifies these conditions.

With additional conditions on the objective function and constraints, we can give stronger results based on the KKT conditions:
For instance, under the linear independence constraint qualification, there are \emph{unique} vectors of Lagrange multipliers $\bm{\iota}$ and $\bm{\kappa}$ satisfying the KKT conditions at a local minimum.
The KKT conditions can also be \emph{sufficient}, under additional restrictions.
For example, by imposing convexity conditions on $f$ and the $g_i$, and linearity on the $h_j$, we have this result:
\begin{thm}
\emph{(First-order sufficient conditions for global minima.)}
Let $f$ and all $g_i$ be continuously differentiable convex functions, and let all $h_j$ be linear functions.  Then if there are $\mb{\hat{x}} \in \bbr^n$, $\bm{\iota} \in \bbr^m$, and $\bm{\kappa} \in \bbr^\ell$ satisfying the following conditions, $\mb{x*}$ is a global minimum of $f$ subject to $\mb{x} \in X$.
\begin{align}
\nabla f(\mb{\hat{x}}) + \sum_{i=1}^{m} \iota_i \nabla g_i(\mb{\hat{x}}) + \sum_{j=1}^{\ell} \kappa_j \nabla h_j(\mb{\hat{x}}) &= \mb{0} \\
h_j(\mb{\hat{x}}) &= 0 \,\, \forall j = 1,2,\ldots,\ell  \\
g_i(\mb{\hat{x}}) &\leq 0 \,\, \forall i = 1,2,\ldots,m  \\
\iota_i g_i(\mb{\hat{x}}) &= 0 \,\, \forall i = 1,2,\ldots,m  \\
\iota_i &\geq 0 \,\, \forall i = 1,2,\ldots,m 
\end{align}
\end{thm}
In this result, the conditions on $g_i$ and $h_j$ imply that the feasible region $X$ is a convex set; since $f$ is also a convex function, we have a convex optimization problem, and therefore any local minimum is also global.
This is what allowed us to convert the ``local'' necessary condition into a global sufficient condition.

As with unconstrained problems, we can also formulate a second-order sufficient condition, based on second partial derivatives.

\begin{thm}
\emph{(Second-order sufficient conditions for a local minimum.)}
If the objective $f$ and all $g_i$ and $h_j$ defining the constraints are twice continuously differentiable; the Hessian of $\mc{L}$ with respect to $\mb{x}$ is positive definite at $\mb{\hat{x}}$, $\bm{\iota}$, and $\bm{\kappa}$; and $\iota_i > 0$ for any $i \notin A(\mb{\hat{x}})$; then any solution satisfying the KKT conditions is a strict local minimum of $f$.
\end{thm}
This condition can be relaxed slightly.
Rather than requiring that the matrix $H_{\mb{xx}}$ consisting of all second partial derivatives $\partial^2 \mc{L} / \partial x_i \partial x_j$ be positive definite (requiring $\mb{d}^T H_\mb{xx}(\mb{\hat{x}}) \mb{d} > 0$ for \emph{all} nonzero $\mb{d} \in \bbr^n$), it is enough if $\mb{d}^T H_\mb{xx}(\mb{\hat{x}}) \mb{d} > 0$ holds for any nonzero $\mb{d} \in \bbr^n$ such that $\nabla h_j(\mb{\hat{x}})^T \mb{d} = 0$ for all $j$, and $\nabla g_i(\mb{\hat{x}})^T \mb{d} = 0$ for all active $i \in A(\mb{\hat{x}})$.

\subsection{Examples of applying KKT conditions}

This section shows how the KKT conditions can be directly solved to identify candidate solutions, and how to then identify an optimal solution among these candidates.
Again, this approach is not practical for large-scale problems, but is good for small problems, for building intuition about how the conditions work, and on occasion in problems with a special structure (for instance, if the conditions can be solved in closed form).

\begin{exm}
Solve the following optimization problem:
\begin{align*}
  \min  \quad &  x_2  & \\
  \mathrm{s.t.} \quad &  (x_1 - 1)^2 + x_2^2 \leq 1 \\
                \quad &  x_1 \geq 2
\end{align*}
Are the KKT conditions satisfied here?
\end{exm}

\solution{
In this formulation, there is only a single feasible point $\mb{x} = (2,0)$, so this solution must be optimal (both locally and globally).
To examine the KKT conditions at this point, we start by writing the two constraints can be written as
\begin{align*}
g_1(x_1, x_2) = (x_1 - 1)^2 + x_2^2 - 1 &\leq 0
g_2(x_1, x_2) = -x_1 + 2 &\leq 0
\,,
\end{align*}
and so their gradients are
\[
\nabla g_1(x_1, x_2) = \vect{2(x_1 - 1) \\ 2x_2}
\qquad \qquad
\nabla g_2(x_1, x_2) = \vect{-1 \\ 0}
\,.
\]
The gradient of the objective function is
\[
\nabla f(x_1, x_2) = \vect{ 0 \\ 1}
\,.
\]

From the first-order necessary conditions, the minimum should satisfy
\[
\nabla f(\mb{\hat{x}}) + \sum_{i=1}^{m} \iota_i \nabla g_i(\mb{\hat{x}}) + \sum_{j=1}^{\ell} \kappa_j \nabla h_j(\mb{\hat{x}}) = 0 
\]
For this particular example, this equation reduces to
\[
\nabla f(\mb{\hat{x}}) + \iota_1 \nabla g_1(\mb{\hat{x}}) + \iota_2 \nabla g_2(\mb{\hat{x}}) = 0 \\
\]
or
\[
\vect{ 0 \\ 1 }  + \iota_1 \vect{ 2(x_1^*-1) \\ 2x_2^* }    + \iota_2 \vect{ -1 \\ 0 } = 0
\,.
\]

At $(x_1^*, x_2^*) = (2,0)$, we have 
\[
\vect{ 0 \\ 1 }  + \iota_1 \vect{ 2 \\ 0 }    + \iota_2 \vect{ -1 \\ 0 } = 0
\,.
\]
But there are no values of $\iota_1$ and $\iota_2$ for which the above equation is satisfied.
So it is impossible to satisfy the KKT conditions at this point, even though it is optimal.
How can this be?
The answer is that the constraint qualifications are not satisfied at this point.
The first constraint is not linear, so the linearity constraint qualification fails.
The gradients of the two constraints are $(2,0)$ and $(-1,0)$, which are linearly dependent, so the linear independence constraint qualification also fails.
Similarly, you can show that the Mangasarian-Fromovitz constraint qualification fails at this point.
}

This example highlights the importance of constraint qualifications.
However, in traffic assignment, a majority of the formulations will have linear constraints, in which case the linear constraint qualification holds.
In such cases, you do not have to worry further.

Note that convexity alone does not imply constraint qualification.
In the above example both $g_1$ and $g_2$ are convex functions, so both the objective function and constraints are convex functions, so this is a convex optimization problem.

\begin{exm}
Solve the following optimization problem, where $a$, $b$, and $c$ are positive real numbers.
\optimizex{\min}{}{x_1 + x_2 + x_3}{
&  \frac{x_1^2}{a} + \frac{x_2^2}{b} + \frac{x_3^2}{c} &= 1 
}
\end{exm}
\solution{
We have $f(x_1, x_2, x_3) = x_1 + x_2 + x_3$ and a single equality constraint $
h(x_1, x_2, x_3) = \frac{x_1^2}{a} + \frac{x_2^2}{b} + \frac{x_3^2}{c} - 1$.
Calculating gradients, we have
\[
\nabla f(x_1, x_2, x_3) = \vect{ 1 \\ 1 \\ 1 }  
\qquad \qquad
\nabla h(x_1, x_2, x_3) = \vect{ \frac{2x_1}{a} \\ \frac{2x_2}{b} \\ \frac{2x_3}{c} }  
\,,
\]
and the condition
\[
\nabla f(x_1^*, x_2^*, x_3^*) + \kappa \nabla h(x_1^*, x_2^* x_3^*) = 0
\]
becomes
\[
\vect{ 1 \\ 1 \\ 1 }  + \kappa \vect{ \frac{2x_1^*}{a} \\ \frac{2x_2^*}{b} \\ \frac{2x_3^*}{c} } = 0
\,.
\]
Solving, we get $x_1^* = -\frac{a}{2v}$, $x_2^* = -\frac{b}{2v}$, and $x_3^* = -\frac{c}{2v}$.

We also have the condition $h(x_1^*, x_2^*, x_3^*) = 0$.
For this problem, we compute as follows:
\begin{align*}
\nabla h(x_1^*, x_2^*, x_3^*) = 0
&\implies   \frac{(x_1^*)^2}{a} + \frac{(x_2^*)^2}{b} + \frac{(x_3^*)^2}{c} = 1 \\
& \implies \frac{a}{4v^2} + \frac{b}{4v^2} + \frac{c}{4v^2} = 1 \\
& \implies v = \frac{\pm \sqrt{a+b+c}}{2}
\,.
\end{align*}

Therefore, we have two possible solutions for $(x_1^*, x_2^*, x_3^*)$:
\[ 
\myp{\frac{a}{\sqrt{a+b+c}}, \frac{b}{\sqrt{a+b+c}}, \frac{c}{\sqrt{a+b+c}}}
\]
and
\[
\myp{-\frac{a}{\sqrt{a+b+c}}, -\frac{b}{\sqrt{a+b+c}}, -\frac{c}{\sqrt{a+b+c}}}
\,.
\]
Checking both, the objective function is minimized at the second point, so this is optimal.
}

\begin{exm}
Write down the KKT conditions for the following optimization problem.
\optimizex{\min}{}{x_1^2 + x_2^2}{
& 11x_1 + 3x_2 &\geq 21 \\
& 6x_1 + 20x_2 &\geq 39 \\
& x_1 + x_2 &\leq 9 \\
& x_1 &\geq 0 \\
& x_2 &\geq 0
}
\end{exm}
\solution{
Since all constraints are linear, the linear constraint qualification is satisfied and the KKT conditions are indeed necessary.
We can rewrite the optimization problem in the form needed for the KKT conditions:
\optimizex{\min}{}{x_1^2 + x_2^2}{
& g_1(x_1,x_2) = -11x_1 - 3x_2 + 21  &\leq 0 \\
& g_2(x_1,x_2) = -6x_1 - 20x_2 + 39  &\leq 0  \\
& g_3(x_1,x_2) = x_1 + x_2 - 9       &\leq 0  \\
& g_4(x_1,x_2) = -x_1                &\leq 0  \\
& g_5(x_1,x_2) = -x_2                &\leq 0  
}

We calculate gradients as follows:
\[
\nabla f(x_1, x_2) = \vect{ 2x_1 \\ 2x_2 }  
\qquad 
\nabla g_1(x_1, x_2) = \vect{ -11 \\ -3 }  
\qquad 
\nabla g_2(x_1, x_2) = \vect{ -6 \\ -20 }  
\]
\[
\nabla g_3(x_1, x_2) = \vect{ 1 \\ 1 }  
\qquad 
\nabla g_4(x_1, x_2) = \vect{ -1 \\ 0 }  
\qquad 
\nabla g_5(x_1, x_2) = \vect{ 0 \\ -1 }  
\,.
\]

With these gradients, the first KKT condition~\eqn{NecKKT1} reduces to the two equations
\begin{align*}
2x_1^* - 11\iota_1 - 6\iota_2 + \iota_3 - \iota_4  &= 0 \\
2x_2^* - 3\iota_1 - 20\iota_2 + \iota_3       - \iota_5 &= 0
\,.
\end{align*}
The second KKT condition~\eqn{NecKKT2} is ignored since there are no equality constraints.
The third KKT condition~\eqn{NecKKT3} ensures feasibility:
\begin{align*}
g_1(x_1^*,x_2^*) \leq 0 &\implies -11x_1^* - 3x_2^* + 21 \leq 0 \\
g_2(x_1^*,x_2^*) \leq 0 &\implies -6x_1^* - 20x_2^* + 39 \leq 0  \\
g_3(x_1^*,x_2^*) \leq 0 &\implies  x_1^* + x_2^* - 9 \leq 0  \\
g_4(x_1^*,x_2^*) \leq 0 &\implies -x_1^* \leq 0  \\
g_5(x_1^*,x_2^*) \leq 0 &\implies -x_2^* \leq 0
\,.
\end{align*}
The fourth ensures complementarity~\eqn{NecKKT4}, that the $\iota_i$ values must be zero unless the constraint is active:
\begin{align*}
\mu_1g_1(x_1^*,x_2^*) = 0 &\implies  \iota_1(-11x_1^* - 3x_2^* + 21) = 0 \\
\mu_2g_2(x_1^*,x_2^*) = 0 &\implies  \iota_2(-6x_1^* - 20x_2^* + 39) = 0 \\
\mu_3g_3(x_1^*,x_2^*) = 0 &\implies  \iota_3(x_1^* + x_2^* - 9) = 0  \\
\mu_4g_4(x_1^*,x_2^*) = 0 &\implies  \iota_4(-x_1^*) = 0 \\
\mu_5g_5(x_1^*,x_2^*) = 0 &\implies  \iota_5(-x_2^*) = 0
\,.
\end{align*}
Finally, we require that the Lagrange multipliers all be non-negative:
\[
\iota_1, \iota_2, \iota_3, \iota_4, \iota_5 \geq 0
\,.
\]
}

If desired, the optimal solution(s) can be found by finding values of $\mb{x}$ and $\bm{\iota}$ satisfying all these conditions.
One way to do this is to start by identifying all 32 combinations of active/inactive constraints.
For each, set $\iota_i = 0$ for the inactive constraints, and solve equations~\eqn{NecKKT1} and~\eqn{NecKKT4} for the remaining variables as a system of equations.
If the solution satisfies the feasibility conditions~\eqn{NecKKT3} and has non-negative Lagrange multipliers, it satisfies all KKT conditions, and is a candidate solution.
Among these, the solutions with the lowest objective function value are the optima.

We illustrate this procedure in the following example, which has fewer constraints:
\begin{exm}
Solve the following optimization problem.
\optimizex{\min}{}{x_1^2 + 2x_2^2}{
& x_1 + 2x_2 &\leq 6 \\
& x_2 &\geq x_1 + 2 
}
\end{exm}
\solution{
Converting to the form required by the KKT conditions, we have $f(x_1, x_2) = x_1^2 + 2x_2^2$, and constraints $g_1(x_1,x_2) = x_1 + 2x_2 - 6$ and $g_2(x_1,x_2) = x_1 - x_2 + 2$, with gradients
\[
\nabla f(x_1, x_2) = \vect{ 2x_1 \\ 4x_2 }  
\qquad
\nabla g_1(x_1, x_2) = \vect{ 1 \\ 2 }  
\qquad
\nabla g_2(x_1, x_2) = \vect{ 1 \\ -1 }  
\,.
\]
The KKT conditions for this problem are
\begin{align*}
2x_1^* + u_1 + \iota_2 &= 0 \\
4x_2^* + 2u_1 - \iota_2 &= 0 \\
\\
x_1^* + 2x_2^* - 6 &\leq 0 \\
x_1^* - x_2^* + 2 &\leq 0 \\
\\
\iota_1(x_1^* + 2x_2^* - 6) &= 0 \\
\iota_2(x_1^* - x_2^* + 2) &= 0 \\
\\
\iota_1, \iota_2 &\geq 0
\,.
\end{align*}

The functions defining the constraints are linear, so the KKT conditions are indeed necessary.
We now find candidate solutions by considering all combinations of active and inactive constraints.

\underline{Case I: Both constraints inactive}:

In this case $\iota_1 = \iota_2 = 0$, so the equality conditions reduce to
\begin{align*}
2x_1^* &= 0 \\
4x_2^* &= 0 \,.
\end{align*}
The solution is $(x_1^*, x_2^*) = (0,0)$.
However this solution violates the constraint $g_2(x_1, x_2) = x_1^* - x_2^* + 2 \leq 0$.
So, there is no optimum solution where both constraints are inactive.

\underline{Case II: Only the second constraint is active}:

In this case $\iota_1 = 0$.
The equality conditions reduce to
\begin{align*}
2x_1^* + \iota_2 &= 0 \\
4x_2^* - \iota_2 &= 0 \\
\iota_2(x_1^* - x_2^* + 2) &= 0
\,.
\end{align*}
The first two equations give $x_1^* = -\iota_2/2$ and $x_2^* = \iota_2/4$. 
Because the second constraint is active, we know $x_1^* - x_2^* + 2 = 0$.
Therefore, the third equation will be satisfied automatically.
Substituting $x_1^* = -\iota_2/2$ and $x_2^* = \mu_2/4$ into the active constraint, we have $\iota_2 = 8/3$, and therefore $x_1^* = 4/3$ and $x_2^* = 2/3$.

This solution satisfies both constraints, and the Lagrange multipliers are non-negative, so it is a candidate for optimality.

\underline{Case III: Only the first constraint is active}:

In this case $\iota_2 = 0$, and the equality conditions are:
\begin{align*}
2x_1^* + \iota_1 &= 0 \\
4x_2^* + 2\iota_1 &= 0 \\
\iota_1(x_1^* + 2x_2^* - 6) &= 0 \,.
\end{align*}
Proceeding in the same way, the first two equations require $x_1^* = -\iota_1/2$ and $x_2^* = -\iota_2/2$. 
Likewise, because we assume the first constraint is active, we can replace the third equation with $x_1^* + 2x_2^* = 6$. 
But with the values of $x_1^*$ and $x_2^*$ for this case, this simplifies to $-\iota_1 = 6$.
This violates the requirement $\iota_1 \geq 0$, and therefore there cannot be an optimal solution corresponding to this case.

\underline{Case IV: Both constraints are active}:

In this case, we know both constraints are satisfied with equality:
\begin{align*}
x_1^* + 2x_2^*  &= 6 \\
x_1^* - x_2^*  &= -2  
\,.
\end{align*}
The only solution to these equations is $x_1^* = 2/3$ and $x_2^* = 8/3$.
It remains to see whether there are Lagrange multipliers satisfying the rest of the KKT conditions.
Substituting into the first two conditions, we have
\begin{align*}
\iota_1 + \iota_2 &= -2x_1^* = -\frac{4}{3} \\
2\iota_1 - \iota_2 &= -4x_2^* =  -\frac{32}{3}
\,.
\end{align*}
Solving this system, we find $\iota_1 < 0$, violating the non-negativity condition, and establishing that the solution where both constraints are active cannot be optimal.

We have exhausted all combinations of active and inactive constraints, and are left with a single candidate solution (from Case II).
Therefore, this one is the optimum, and the solution to the problem is $x_1^* = -\frac{4}{3}$ and $x_2^* = \frac{2}{3}$.
If there were multiple candidate solutions, we would evaluate the objective function at each to see which one(s) are minima.
}
\index{Karush-Kuhn-Tucker (KKT) conditions|)}
\index{optimization!nonlinear!optimality conditions|)}

\subsection{Algorithms for linear constraints}

\index{optimization!nonlinear!algorithms|(}
The remainder of this section will describe four algorithms for solving nonlinear optimization problems with linear constraints.
The first two algorithms are the method of successive averages (MSA) and Frank-Wolfe (FW) algorithm.
Both of these are called convex combination methods.
The third and fourth methods, gradient projection (GP) and manifold suboptimization (MS), modify the direction of steepest descent to ensure feasibility.
The algorithms are presented in increasing order of complexity and efficiency.
All of these algorithms will converge to a local optimum.
If the objective is convex, this local optimum will also be globally optimal.

We focus here on nonlinear optimization problems with linear constraints, as the traffic assignment problem have this structure.
These four algorithms have all been applied to traffic assignment as well.
In Chapter~\ref{chp:solutionalgorithms}, these algorithms are presented specifically in this context.
Here we treat them more generally, as they can be used for other nonlinear optimization problems as well.

There are many ways to modify these algorithms, and the best implementation often depends on the specific problem structure and the computing environment being used.
For example, all four methods have the common structure of ``find an improving direction'' and ``choose a step size in that direction.''
This section mainly emphasizes the former, since the algorithms differ more in this regard.
However, there are alternative choices for the latter as well.
Our presentations of FW, PG, GP, and MS in this subsection choose the step size which minimizes the objective function along the chosen direction.
This indeed results in the greatest improvement possible at each iteration.
However, this can be computationally demanding, and it might be more efficient to take a faster but less exact step.
For instance, the algorithms in Chapter~\ref{chp:solutionalgorithms} often use a single iteration of Newton's method in place of finding an exactly-optimal step size, since this tradeoff is usually worthwhile in traffic assignment.
If you are comparing the presentations of the algorithms in that chapter which are specialized to traffic assignment to the more general presentations here, such differences should be noted, but they are matters of implementation and not essential algorithmic distinctions.

The standard form of the nonlinear optimization problem with linear constraints is as follows.
The decision variables are $\mb{x} = \vect{x_1 & x_2 & \ldots & x_n}$.
The objective function $f$ is a convex function of $\mb{x}$.
The feasible region $X$ is defined by linear equality and inequality constraints, as in a linear program (Section~\ref{sec:linearoptimization}).
With linear constraints such as these, the feasible region is convex and closed.
We further assume that it is bounded, as is the case in traffic assignment and many other applications.

\subsection{Convex combination methods}

\index{optimization!nonlinear!method of successive averages|(}
\index{optimization!nonlinear!Frank-Wolfe|(}
All convex combinations methods can be described by the following framework (see Algorithm \ref{alg:feasible}).
They have three main steps: (i) determining a feasible descent direction, (ii) determining the step size, and (iii) updating the solution.
The convergence criterion for feasible direction methods can be chosen in a similar manner as for unconstrained optimization algorithms. 

\begin{algorithm}
\SetAlgoLined
\textbf{Initialize}: \\
$\mb{x}$: initial value \\
$\epsilon$: tolerance \\
\While{$ ConvergenceCriterion > \epsilon$ }{
    \textbf{Determine feasible descent direction}: $\mb{d}$ \\
    \textbf{Determine step size}: $\mu_k \geq 0 $ such that $f(\mb{x} + \lambda \mb{d}) < f(\mb{x})$ and $\mb{x} + \mu \mb{d} \in X$\\
    \textbf{Update}: $\mb{x} \leftarrow \mb{x} + \mu \mb{d}$ \\
}
\Return $\mb{x}$
\caption{Feasible directions method.}
\label{alg:feasible}
\end{algorithm}

When the current solution is $\mb{x}$, the feasible descent direction $\mb{d}$ is obtained by taking $\mb{d} = \mb{x^*} - \mb{x}$, where $\mb{x^*}$ is a feasible solution obtained by solving an approximation of the original nonlinear problem.
Both MSA and FW obtain the feasible direction by replacing the objective function with its first-order (linear) approximation at $\mb{x}$.
Because the constraints are linear, the resulting first-order approximation will be the linear program shown below, which can be solved using the simplex algorithm.
(In some applications, the first-order approximation can be solved more efficiently by exploiting a particular structure of the problem.
For example, in traffic assignment, one can show that an all-or-nothing assignment to shortest paths is optimal to this subproblem, which can be found using one of the algorithms in Section~\ref{sec:shortestpath}.

\optimizex{\min}{\mb{x^*}}{\nabla f(\mb{x})^T \mb{x^*}}{
& x \in X
}

The difference between MSA and FW is how the step size is chosen.
In MSA, at each iteration $k$, the step size is commonly chosen to be $\lambda = 1/(k+1)$, although any sequence of step sizes satisfying $\sum \lambda_k = \infty$ and $\sum \lambda_k^2 < \infty$ will work.
In FW, the step size is obtained by solving an auxiliary optimization problem in $\lambda$ alone.
This optimization problem has only a single variable, so it can be solved using any of the line search methods from Sections~\ref{sec:bisection} or~\ref{sec:linesearch}.
\optimizex{\min}{\lambda}{f(\mb{x} + \lambda \mb{d})}{
& 0 \leq \lambda \leq 1
}

FW picks the step size in a more intelligent manner than MSA and typically converges faster.
For this reason, MSA is rarely used in solving nonlinear optimization with linear constraints.
However, in certain complex traffic assignment problems, evaluating the objective function $f$ can be computationally expensive.
In such cases, a line search method in FW can require large amounts of time, and taking faster MSA steps outweighs the greater precision that each FW step would have.

\begin{exm}
Determine the optimal solution of the following optimization problem using MSA and FW.
\optimizex{\min}{}{2(x_1 - 1)^2 + (x_2 - 2)^2}{
& x_1 + 4x_2 - 2 &\leq 0   \\
&-x_1 + x_2      &\leq 0  \\
& x_1, x_2       &\geq 0
}
\end{exm}
\solution{
The KKT conditions for this problem can be solved to yield the optimal solution $x_1^* = 0.7878, x_2^* = 0.3030$.
The MSA and FW algorithms will converge to this solution over successive iterations; the advantage of these methods is that solving the KKT conditions for a large problem can be very difficult, whereas MSA and FW scale better with problem size.
For this objective function, the gradient at any point $(x_1, x_2)$ is
\begin{align*}
\nabla f(x_1,x_2) =  \vect{
                        4(x_1-1) \\
                        2(x_2-2) 
                     }
\,.
\end{align*}

\textbf{\underline{MSA}:}
Assume the initial solution is $x_1 = 0, x_2 = 0$, and the step size is $\lambda = 1/(k+1)$, where $k$ is the iteration number.
The search direction is obtained by solving the following linear program.
\optimizex{\min}{}{-4x^*_1 - 4x^*_2}{
& x^*_1 + 4x^*_2 - 2 &\leq 0   \\
&-x^*_1 + x^*_2 &\leq 0  \\
& x^*_1, x^*_2 &\geq 0
}
The optimal solution to this problem is $x^*_1 = 2, x^*_2 = 0$.
The search direction is $\mb{d} = \mb{x^*} - \mb{x} = \vect{2 & 0}$ and $\lambda = 1/2$, so the new solution is
\begin{align*}
x_1 &= 0 + \frac{1}{2}(2) = 1 \\
x_2 &= 0 + \frac{1}{2}(0) = 0
\,.
\end{align*}

At the start of the second iteration we have $x_1 = 1, x_2 = 0$.
Therefore, the search direction is obtained by solving the following linear program.
\optimizex{\min}{}{-4x^*_2}{
& x^*_1 + 4x^*_2 - 2 &\leq 0   \\
&-x^*_1 + x^*_2 &\leq 0  \\
& x^*_1, x^*_2 &\geq 0
}
The optimal solution to this linear program is $x^*_1 = 0.4, x^*_2 = 0.4$.
In iteration 2, the step size is $\lambda = \frac{1}{3}$.
Therefore, the updated solution is
\begin{align*}
 x_1 &= 1 + \frac{1}{3}(0.4-1) = 0.8 \\
 x_2 &= 0 + \frac{1}{3}(0.4-0) = 0.1333
 \,.
\end{align*}
Table \ref{tbl:MSA_appendix} shows the progress of MSA over further iterations.
The third column in this table is calculated as the sum of square deviation from the exact optimal solution, $\hat{x}_1 = 0.7878, \hat{x}_2 = 0.3030$.

\begin {table}
\begin{center}
\caption{\label{tbl:MSA_appendix} Convergence of the method of successive averages.}
\begin{tabular}{|c|c|c|c|c|}
\hline
$k$ & $x_1$  & $x_2$  & Distance to optimal  \\
\hline 
1   & 0      & 0      & 0.7124 \\
2   & 1      & 0      & 0.1368 \\
3   & 0.8    & 0.1333 & 0.0289 \\
4   & 0.7    & 0.2    &  0.01831 \\
5   & 0.96   & 0.16   & 0.0501 \\
6   & 0.866  & 0.2    & 0.0168 \\
7   & 0.800  & 0.2285 & 0.005 \\
20  & 0.780  & 0.28   & 0.0005 \\
53  & 0.7849 & 0.2943 & 8.33793 $\times 10^{-5}$ \\
\hline
\end{tabular}
\label{tbl:MSA}
\end{center}
\end{table}

\textbf{\underline{FW}:}
Again assume the initial solution to $x_1 = 0, x_2 = 0$.
The search direction is obtained by solving the same linear program as in MSA.
\optimizex{\min}{}{-4x^*_1 - 4x^*_2}{
& x^*_1 + 4x^*_2 - 2 &\leq 0   \\
&-x^*_1 + x^*_2 &\leq 0  \\
& x^*_1, x^*_2 &\geq 0
}
The optimal solution to the LP is again $x^*_1 = 2, x^*_2 = 0$, so the direction is the same: $\mb{d} = \vect{2 & 0}$.
In FW, the step size is obtained by solving the following minimization problem:
\optimizex{\min}{\lambda}{2 \myp{0 + \lambda(2-0) - 1}^2 + \myp{0 + \lambda(0-0) - 2}^2}{
& 0 \leq \lambda \leq 1  
}
The optimal solution is $\lambda = 0.5$, and the solution is updated as 
\begin{align*}
x_1 &= 1 + 0.5 \times (2-0) = 1 \\
x_2 &= 0 + 0.5 \times (0-0) = 0
\,.
\end{align*}
To this point, FW is proceeding in the same way as MSA.

At the start of the second iteration we have $x_1 = 2, x_2 = 0$.
Therefore, the search direction is obtained by solving the following linear program.
\optimizex{\min}{}{-4x^*_2}{
& x^*_1 + 4x^*_2 - 2 &\leq 0   \\
&-x^*_1 + x^*_2 &\leq 0  \\
& x^*_1, x^*_2 &\geq 0
}
The optimal solution to the linear program is $x^*_1 = 0.4, x^*_2 = 0.4$, as before.
However, at this point FW chooses the step size differently than MSA.
Solving the optimization problem
\optimizex{\min}{\lambda}{2 \myp{ 1 + \lambda(0.4-1) - 1 })^2 + \myp{0 + \lambda(0.4-0) - 2 }^2}{
& 0 \leq \lambda \leq 1
}
gives the step size $\lambda = 0.9091$.
Therefore, we obtain the new solution
\begin{align*}
x_1 &= 1 + 0.9091 \times (0.4 - 1) = 0.4545 \\
x_2 &= 0 + 0.9091 \times (0.4 - 0) = 0.3634
\,.
\end{align*}

\begin {table}
\begin{center}
\caption{\label{tbl:FW} Convergence of the Frank-Wolfe algorithm.}
\begin{tabular}{|c|c|c|c|c|}
\hline
$k$ & $x_1$  & $x_2$  & Distance to optimal \\
\hline 
1   & 0      & 0      & 0.7124 \\
2   & 1      & 0      & 0.1368 \\
3   & 0.4545 & 0.3636 & 0.1147 \\
4   & 0.7979 & 0.2828 & 0.0005 \\
20  & 0.7919 & 0.2947 & 8.6 $\times 10^{-5}$ \\
\hline
\end{tabular}
\end{center}
\end{table}

Table~\ref{tbl:FW} shows the convergence of FW over successive iterations.
Comparing with Table~\ref{tbl:MSA}, we see that FW converges in fewer iterations than MSA.
FW achieves an error of 0.0005 in 4 iterations; MSA required 20 iterations to achieve this level of precision. 
The tradeoff is that FW required more effort at each iteration to determine the step size $\lambda$.
In problems such as this, the tradeoff is clearly in favor of FW.
}
\index{optimization!nonlinear!method of successive averages|)}
\index{optimization!nonlinear!Frank-Wolfe|)}

\subsection{Gradient projection}

\index{optimization!nonlinear!gradient projection|(}
The gradient projection method operates by taking a step in the direction of steepest descent (ignoring feasibility), and \emph{then} jumping to the nearest feasible point if the new solution is infeasible.
Recall that $\mathrm{proj}_X(\mb{x})$ means ``the point in $X$ closest to $\mb{x}$.''
(If $\mb{x} \in X$ already, then $\mathrm{proj}_X(\mb{x}) = \mb{x}$.)
If the set $X$ is convex, then the projection operation is uniquely defined, and a continuous function of $\mb{x}$.
In general, projection cannot be easily computed.
When the constraints are linear (as we assume in this section), it is considerably easier.
For some problems, it can be exceptionally easy.
For instance, in traffic assignment, we can reformulate all constraints to be simple non-negativity constraints of the form $h^\pi \geq 0$, and ``projection'' simply means ``any negative path flow should be set to zero.''

The steps of gradient projection are shown in Algorithm~\ref{alg:GP}.
In this algorithm, $\theta$ is a pre-defined upper bound on the step size $\mu$.
\begin{algorithm}
  \SetAlgoLined
  \textbf{Initialize}: \\
  $\mb{x}$: initial value \\
  $\epsilon$: tolerance \\
  \While{$ConvergenceCriterion > \epsilon$ }{
      $\mb{d} \leftarrow -\nabla f(\mb{x})$ \\
      $\mu_k \leftarrow \arg \min_\mu \{ f(\mathrm{proj}_X(\mb{x} + \mu \mb{d})), 0 \leq \mu \leq \theta \} $ \\
      $\mb{x} \leftarrow  \mb{x} + \mu \mb{d}$ \\
  }
  \Return $\mb{x}$
\caption{\label{alg:GP} Gradient projection method.}
\end{algorithm}

\begin{exm}
Re-solve the previous optimization problem using the gradient projection method, with $\theta = 1$.
\optimizex{\min}{}{2(x_1 - 1)^2 + (x_2 - 2)^2}{
& x_1 + 4x_2 - 2 &\leq 0   \\
&-x_1 + x_2      &\leq 0  \\
& x_1, x_2       &\geq 0
}
\end{exm}
\solution{
We start from the same initial solution $x_1 = 0, x_2 = 0$. 
Here the search direction is $\mb{d} = -\nabla f(\mb{x}) = \vect{4 & 4}$.
To find the step size, we solve the optimization problem
\optimizex{\min}{\mu}{f(\mathrm{proj}_X(\mb{x} + \mu \mb{d}))}{
& 0 \leq \mu \leq 1 
}

The projection operation can be more easily seen on a plot of the feasible region, as in Figure~\ref{fig:projectionexample}.
For a given value of $\mu$, $\mb{x} + \mu \mb{d}$ is the point $(4 \mu, 4 \mu)$.
If $\mu \leq 1/10$, this coincides with the line $x_1 = x_2$, which is part of the feasible region.
As a result, projection does not change the point and $\mathrm{proj}_X (4 \mu, 4 \mu) = (4 \mu, 4 \mu)$.
Once $\mu > 1/10$, the point $(4 \mu, 4 \mu)$ violates the constraint $x_1 + 4 x_2 \leq 2$, and is infeasible.
For these points, we have to project them onto the nearest feasible point.
As an example, if $\mu = 1/4$, the point $(4 \mu, 4 \mu) = (1,1)$ is infeasible.
The closest feasible point is $(14/17, 5/17)$; this can be seen geometrically in Figure~\ref{fig:projectionexample} by drawing a line through $(1,1)$ perpendicular to the line corresponding to that constraint.
The expression mapping a value of $\mu$ in this range to the closest feasible point can be found algebraically, and the formula is shown below.
If $\mu > \frac{2}{3}$, the point on the line $x_1 + 4 x_2$ lies below the $x_2$-axis, which is also infeasible.
For $\mu$ values this large, the closest feasible point is simply $(2,0)$.
Therefore we have
\[
\mathrm{proj}_X(\mb{x_1} + \mu \mb{d_1}) =
\begin{cases}
(4 \mu, 4 \mu)                        & \mu \in [0, 1/10]   \\
\frac{1}{17} (48\mu + 2, -12 \mu + 8) & \mu \in [1/10, 2/3] \\
(2,0)                                       & \mu \geq 2/3
\end{cases}
\,.
\]
You can check that this formula is continuous.

\begin{figure}
\centering
\begin{tikzpicture}[x=5cm,y=5cm]
\path [fill=yellow!20!white] (0,0) -- (0.4,0.4) -- (2,0) -- (0,0);
\draw (0,0) -- coordinate (x1 axis mid) (2,0);
\draw (0,0) -- coordinate (x2 axis mid) (0,2);
\foreach \x in {0,...,4}      \draw (\x / 2,1pt) -- (\x / 2,-3pt);
\foreach \y in {0,...,4}       \draw (1pt,\y / 2) -- (-3pt,\y / 2);
\foreach \x [evaluate=\x as \xval using (1/2)*\x] in {0,1,...,4} 
  \draw (\xval, -3pt) node[anchor=north] {\pgfmathroundto{\xval}\pgfmathresult};
\foreach \y [evaluate=\y as \yval using (1/2)*\y] in {0,1,...,4} 
  \draw (-3pt, \yval) node[anchor=east] {\pgfmathroundto{\yval}\pgfmathresult};
\node[below=0.8cm] at (x1 axis mid) {$x_1$};
\node[left=0.8cm] at (x2 axis mid) {$x_2$};
\draw [blue, ultra thick, domain=0:2] plot (\x, \x);
\draw [blue, ultra thick, domain=0:2] plot (\x, {(2 - \x)/4});
\draw [blue, ultra thick, domain=0:2] plot (\x, 0);
\draw [blue, ultra thick] plot (0, 0) -- (0, 2);
\draw [red, ultra thick] plot (0, 0) -- (0.4,0.4) -- (2,0);
\draw [red] plot (1,1) -- (14/17, 5/17);
\node[rotate=45] at (1.5,1.65) {$-x_1 + x_2 \leq 0$};
\node[rotate=-14] at (1.5,0.25) {$x_1 + 4 x_2 \leq 2$};
\node[right=0.25cm] at (1,1) {$\mu = 1/4 \rightarrow (1,1)$};
\node[below=0.25cm,text width=2cm] at (14/17,5/17) {$\mathrm{proj}_X(1,1) = (14/17, 5/17)$};
\end{tikzpicture}
\caption{Feasible region and projection. \label{fig:projectionexample}}
\end{figure}
}

Evaluating the objective function in terms of $\mu$, we thus have
\begin{multline*}
f(\mathrm{proj}_X(\mb{x} + \mu \mb{d})) = \\
\begin{cases}
2(4\mu - 1)^2 + (4\mu - 2)^2           & \mu \in [0, 1/10]   \\
2 \myp{ \frac{1}{17} (48 \mu + 2) - 1}^2
+ \myp{ \frac{1}{17} (-12 \mu + 8) - 2}^2 & \mu \in [1/10, 2/3] \\
6                                            & \mu \geq 2/3
\end{cases}
\,.
\end{multline*}
Again, you can check that this function is continuous and convex.
Using a line search method, the minimum is found for $\mu = 0.236$, which corresponds to the point $(0.7878, 0.3030)$.
As shown above, this is actually the optimal solution, found in a single iteration!
Typically, GP only converges to the optimal solution in the limit, but as this example shows it can be extremely efficient in terms of iterations.
The drawback is that each iteration requires a significant amount of work.
But in problems where the projection can be calculated very efficiently (as in traffic assignment), GP can be an excellent algorithmic choice.
\index{optimization!nonlinear!gradient projection|)}

\subsection{Manifold suboptimization}

\index{optimization!nonlinear!manifold suboptimization|(}
Manifold suboptimization works by modifying the direction of steepest descent (related to the gradient of the objective) to ensure that it remains within the feasible set.
This method also involves projection, but unlike gradient projection, the projection is implicit through matrix multiplications which ensure that whichever constraints are exactly satisfied at the current solution remain satisfied.
There are several names that are used to describe this method, which unfortunately can cause confusion.
Some authors refer to this method as the \emph{projected gradient} method (distinguishing it from gradient projection, because the feasibility is used to identify the improving direction, rather than moving in a direction and then checking feasibility), and others even call it \emph{gradient projection} (for historical reasons).
We believe that calling it \emph{manifold suboptimization} avoids this confusion, and by doing so follow the convention of~\cite{bertsekas_nlp}.

Intuitively, the idea is to distinguish between ``active'' constraints\index{optimization!constraints!active} (where the left-hand side is exactly equal to the right-hand side) and ``inactive'' ones (where the left-hand side is strictly less than the right-hand side).  
The gradient is modified so that any active constraints continue to hold with equality; inactive constraints do not impose any such restriction, since their left-hand sides can vary some without violating feasibility.
Equality ($=$) constraints are always active.
Inequality ($\leq$) constraints may or may not be active.
As the algorithm proceeds, which constraints are active and which are inactive may vary, and some care must be taken to manage this set properly and to ensure feasibility of the step size (the left-hand side of an inactive constraint can vary some, but if it changes too much it may become violated).
In this subsection, we present the steps of the algorithm.
Readers wishing a derivation of these steps, or a proof that these steps indeed maintain feasibility and lead to an optimal solution, are referred to the books by~\cite{bazaara_nlp} and~\cite{bertsekas_nlp}.

Recall that throughout this section, we are assuming a feasible region defined by linear equality and inequality constraints.
For manifold suboptimization, we need an explicit description of this feasible set, so we will  write the feasible region $X$ as
\begin{equation*}
X = \{x: \mb{A_1 x} \leq \mb{b_1}, \mb{A_2 x} = \mb{b_2} \}
\,,
\end{equation*}
where $\mb{A_1}$ and $\mb{b_1}$ are the coefficient matrix and right-hand side vector for the inequality constraints, with $\mb{A_2}$ and $\mb{b_2}$ the same for the equality constraints.

For any feasible solution $\mb{x}$, let $\tilde{A}_1$ be the matrix containing the coefficients of the set of active inequality constraints, so
\begin{equation*}
\mb{\tilde{A}_1 x} = \mb{b_1}
\,.
\end{equation*}
Define $\tilde{A}$ as the matrix containing the coefficients of \emph{all} active constraints:
\begin{equation*}
\mb{\tilde{A}} =
\begin{bmatrix}
\mb{\tilde{A}_1} \\
\mb{A_2}
\end{bmatrix} 
\,.
\end{equation*}
It may be that the rows of the matrix $\mb{\tilde{A}}$ are not linearly independent; this indicates that one or more of the active constraints is redundant (they are implied by some of the others).
These redundant constraints may be removed.

Algorithm \ref{alg:MS} provides the steps for the manifold suboptimization method.
Note that the convergence criterion is different from the MSA and FW methods.
In this algorithm, $\bar{\theta}$ is a parameter giving a ``default'' maximum step size (which may be further restricted to ensure feasibility).

\begin{algorithm}
  \SetAlgoLined
  \textbf{Initialize}: \\
  $\mb{x}$: initial value \\
  $OptimalityFlag \leftarrow FALSE$ \\
  \While{$OptimalityFlag = FALSE$ }{
      $\mb{P} \leftarrow \mb{I} - \mb{\tilde{A}}^T(\mb{\tilde{A}}\mb{\tilde{A}}^T)^{-1}\mb{\tilde{A}}$\\
      $\mb{d} \leftarrow -\mb{P} \nabla f(\mb{x})$ \\
      \If{$\mb{d} \neq \mb{0}$}{
        $\theta \leftarrow \max \{\bar{\theta}, \mb{x} + \mu \mb{d} \in X \} $ \\
        $\mu \leftarrow \arg \min_\mu \{ f(\mb{x} + \mu \mb{d}), 0 \leq \mu \leq \theta \} $ \\
        $\mb{x} \leftarrow  \mb{x} + \mu \mb{d}$ \\
      } 

      \If{$\mb{d} = \mb{0}$}{
        $\mb{q} = -(\mb{\tilde{A}}\mb{\tilde{A}}^T)^{-1}\mb{\tilde{A}} \nabla f(\mb{x})$ \\
        \If{$q_j \geq 0 \,\, \forall j \mbox{ corresponding to } \mb{\tilde{A}_1}$}{
          $OptimalityFlag \leftarrow TRUE$\\
        }
        \Else{
          Delete row from $\mb{\tilde{A}}$ corresponding to most negative $q_j$ \\
        }
      } 
  }
  \Return $\mb{x}$
\caption{\label{alg:MS} Manifold suboptimization algorithm.}
\end{algorithm}

\begin{exm}
Re-solve the previous optimization problem using the manifold suboptimization method, using $\bar{\theta} = 1$.
\optimizex{\min}{}{2(x_1 - 1)^2 + (x_2 - 2)^2}{
& x_1 + 4x_2 - 2 &\leq 0   \\
&-x_1 + x_2      &\leq 0  \\
& x_1, x_2       &\geq 0
}
\end{exm}
\solution{
We start from the same initial solution $x_1 = 0, x_2 = 0$.  
At this point, the active constraints are $-x_1 + x_2 \leq 0$, $x_1 \geq 0$, and $x_2 \geq 0$.
We only use $-x_1 + x_2 \leq 0$ and $x_2 \geq 0$ as they together imply $x_1 \geq 0$.
(Mathematically, the three vectors corresponding to the coefficients in these constraints are linearly dependent: $\vect{-1 & 1}$, $\vect{0 & 1}$, and $\vect{1 & 0}$.)

\[
\mb{\tilde{A}} = \begin{bmatrix} -1 & 1 \\ 0 & -1 \end{bmatrix} 
\]

\[
\mb{P} = \mb{I} - \mb{\tilde{A}}^T(\mb{\tilde{A}}\mb{\tilde{A}}^T)^{-1}\mb{\tilde{A}} = \begin{bmatrix} 0 & 0 \\ 0 & 0 \end{bmatrix} 
\]

\[
\mb{d} = -\mb{P} \nabla f(\mb{x}) = \begin{bmatrix} 0 \\ 0 \end{bmatrix}
\,.
\]

Since $\mb{d}=\mb{0}$, we calculate $\mb{q}$:
\[
\mb{q} = -(\mb{\tilde{A}}\mb{\tilde{A}}^T)^{-1}\mb{\tilde{A}} \nabla f(\mb{x_1}) = \begin{bmatrix} -4 \\ -8 \end{bmatrix} 
\,.
\]

Since the second row is the most negative, we delete the second row of the matrix $\mb{\tilde{A}}$.
\[
\mb{\tilde{A}} = \begin{bmatrix} -1 & 1 \end{bmatrix} 
\]

\[
\mb{P} = \mb{I} - \mb{\tilde{A}}^T(\mb{\tilde{A}}\mb{\tilde{A}}^T)^{-1}\mb{\tilde{A}} = \begin{bmatrix} 0.5 & 0.5 \\ 0.5 & 0.5 \end{bmatrix} 
\]

\[
\mb{d_1} = -\mb{P} \nabla f(\mb{x_1}) = \begin{bmatrix} 4 \\ 4 \end{bmatrix}
\,.
\]

The new solution can be obtained as:

\begin{align*}
x_1 &= 0 + 4 \mu \\
x_2 &= 0 + 4 \mu
\,.
\end{align*}

We need to find the step size $\mu$ which minimizes  $2(0 + 4 \mu - 1)^2 + (0 + 4 \mu- 2)^2$ while retaining feasibility.
Note that
\[
x_1 + 4x_2 - 2 \leq0 \implies \mu \leq 0.1
\,,
\]
which is less than $\bar{\theta} = 1$, so $\theta = 0.1$.
Therefore, we obtain the step size by solving the optimization problem
\optimizex{\min}{\mu}{2(0 +  4 \mu - 1)^2 + (0 + 4 \mu - 2)^2}{
& 0 \leq \mu \leq 0.1
}The solution is $\mu = 0.1$, so the new solution to the original problem is
\begin{align*}
x_1 &= 0 + 4 \mu = 0.4\\
x_2 &= 0 + 4 \mu = 0.4
\,.
\end{align*}
If we do one more iteration of MS, we get $x_1 =  0.7878, x_2 =  0.3030$, which is optimal.
MS converges to the optimal solution much faster than MSA or FW (and about as fast as GP), because its steps are closer to the steepest descent direction.
Section~\ref{sec:pathbasedalgorithms} has more discussion of these reasons in the specific context of traffic assignment.
}
\index{optimization!nonlinear!manifold suboptimization|)}
\index{optimization!nonlinear!algorithms|)}
\index{optimization!nonlinear|)}

\section{Integer Programming}
\label{sec:integerprogramming}

\index{optimization!integer|(}
All of the optimization formulations presented above assume that the decision variables are continuous in nature.
In some cases, this is not a tenable assumption.
For example, when determining the location to establish a warehouse, a fractional solution splitting the warehouse location between two nodes is not meaningful.
Integer programs are a special class of optimization problems, where some or all decision variables must take integer values.
Integer programs have wide applications:
\begin{description}
\item[Vehicle routing:] Determine routes for a specified number of trucks, starting from a depot location to service a given number of customers with demands, subject to truck capacity constraints at minimal cost.
\item[Facility location:] Determine where to locate facilities in a network, satisfying customer demand, with an objective such as minimizing transportation costs, distances, or maximum distance to facilities.
\item[Transit network design:] Determine the optimal routes and schedules of a bus fleet to satisfy customer demands at minimal costs. 
\end{description}
As always, there is no such thing as a free lunch; for reasons discussed below, integer problems are significantly harder to solve, and often it is impractical to find a provably optimal solution.

There are several types of integer programs.
An integer linear program (ILP)\index{optimization!integer!linear} is a linear program where \emph{all} the variables are restricted to be integers.
A mixed integer linear program (MILP)\index{optimization!integer!mixed} is a linear program where only some of the decision variables are assumed to be integer. 
A MILP takes the form\label{not:Eilp}\label{not:dilp}
\optimizex{\min}{\mb{x}, \mb{y}}{\mb{c} \cdot \mb{x} + \mb{d} \cdot \mb{y}}
{
&    \mb{Ax} + \mb{Ey}  \leq \mb{b} \\
&     \mb{x} \geq 0, \mb{y} \geq 0  \\
&     \mb{y} \in  \bbz^p
}where $\mb{x}, \mb{c} \in \bbr^n$, $\mb{y}, \mb{d} \in \bbr^p$, $\mb{A} \in \bbr^{m \times n}$, $\mb{E} \in \bbr^{m \times p}$, and $\mb{b} \in \bbr^m$ with $m$, $n$, and $p$\label{not:pilp} respectively denoting the number of constraints, continuous decision variables, and integer decision variables.
In an ILP, $n = 0$ and the terms involving $\mb{x}$ in this formulation disappear.
A special case of an ILP is the binary integer linear program\index{optimization!integer!binary} (BILP) where the decision variables are restricted to be either 0 or 1, i.e. $\mb{y}\in \{0,1\}^p $.

Integer programs are tougher to solve than linear programs as optimal solutions need not occur at the corner points of a convex polyhedron.
Indeed, the feasible region itself is not a convex set; $y_1 = 1$ and $y_2 = 2$ may be feasible solutions, but their average $\frac{1}{2} y_1 + \frac{1}{2}y_2 = 1.5$ is not.
Solving the corresponding linear program and then rounding the solution may not produce a good solution either.
(One notable exception is discussed in Section~\ref{sec:totallyunimodular}.)

Consider the following ILP with two decision variables $y_1$ and $y_2$:
\optimizex{\min}{y_1, y_2}{10y_1 + 26y_2}{
& 11y_1 + 3y_2  \geq 21  \\
&  6y_1 + 20y_2 \geq 39  \\
&   y_1,    y_2 \geq 0   \\
&   y_1,    y_2 \in \bbz
}

\begin{figure}
\centering
\begin{tikzpicture}
\path [fill=yellow!20!white] (1.5,1.5) -- (0,7) -- (0,10) -- (10,10) -- (10,0) -- (6.5,0);
\draw (0,0) -- coordinate (x axis mid) (10,0);
\draw (0,0) -- coordinate (y axis mid) (0,10);
\foreach \x in {0,...,10}      \draw (\x,1pt) -- (\x,-3pt)     node[anchor=north] {\x};
\foreach \y in {0,...,10}       \draw (1pt,\y) -- (-3pt,\y) node[anchor=east] {\y}; 
\node[below=0.8cm] at (x axis mid) {$y_1$};
\node[left=0.8cm] at (y axis mid) {$y_2$};
\draw [blue, ultra thick, domain=0:21/11] plot (\x, { -11*\x/3 + 7});
\draw [blue, ultra thick, domain=0:6.5] plot (\x, { -3*\x/10 + 39/20});
\draw [cyan, dashed, ultra thick, domain=0:9.6] plot (\x, { -10*\x/26 + 96/26});
\node at (6.5,0.5) {$6y_1+20y_2 \geq 39$};
\node at (2,5) {$11y_1+3y_2 \geq 21$};
\node [right] at (6,1.5) {$z=10y_1+26y_2$};
\draw [blue, ultra thick] (0,7)--(0,10);
\draw [blue, ultra thick] (6.5,0)--(10,0);
\node [below left] at (1.6,1.6) {(1.5,1.5)};
\draw [fill=black] (1.5,1.5)circle (3pt);
\draw [fill=red] (2,2)circle (3pt);
\draw [fill=green] (4,1)circle (3pt);
\draw [fill=red] (5,1)circle (3pt);
\draw [fill=red] (6,1)circle (3pt);
\draw [fill=red] (3,2)circle (3pt);
\draw [fill=red] (4,2)circle (3pt);
\draw [fill=red] (5,2)circle (3pt);
\draw [fill=red] (6,2)circle (3pt);
\draw [fill=red] (2,3)circle (3pt);
\draw [fill=red] (3,3)circle (3pt);
\draw [fill=red] (4,3)circle (3pt);
\draw [fill=red] (5,3)circle (3pt);
\draw [fill=red] (6,3)circle (3pt);
\draw [fill=red] (1,4)circle (3pt);
\draw [fill=red] (2,4)circle (3pt);
\draw [fill=red] (3,4)circle (3pt);
\draw [fill=red] (4,4)circle (3pt);
\draw [fill=red] (5,4)circle (3pt);
\draw [fill=red] (6,4)circle (3pt);
\draw [<-] [thick] (1.6,1.6) -- (4,6);
\node at (4,6.3) {Optimal LP Solution};
\draw [<-] [thick] (4.1,1.1) -- (5,5);
\node [right] at (4,5.3) {Optimal IP Solution};
\end{tikzpicture}
\caption{Difference between LP and IP solution. \label{fig:round}}
\end{figure}

Without the integer restrictions on the variables, the optimal solution of the linear program would be $y_1 = 1.5$, $y_2 = 1.5$, with an optimal objective function value of 54 (see Figure~\ref{fig:round}).
When the LP solutions are rounded up ($y_1=2,y_2=2$) or rounded down ($y_1=1,y_2=1$) the resulting points are sub-optimal or infeasible for the integer program.
The optimal solution for the integer program is in fact $y_1=4,y_2=1$, with an optimal objective function value of 66.

Still, while ignoring the integrality constraints and solving the resulting linear program may not lead to an optimal solution, this \emph{linear program relaxation}\index{optimization!integer!linear program relaxation} is still very useful, because it provides a lower bound on the optimal value of the objective function.  
In the above example, the optimal solution without integrality constraints had an objective function value of 54.
Restricting the feasible set more by adding the integer constraints cannot improve this.
Therefore, even if we don't know the optimal integer solution of the optimization problem, we know its objective function value cannot be lower than 54.

Furthermore, any feasible integer solution provides an upper bound on the optimal value for the objective.
For instance, the rounded solution $y_1 = y_2 = 2$ is feasible, and has an objective function value of 72.  
Even if we don't know whether this solution is optimal or not, this solution tells us that the optimal objective function value can't be any higher than 72.
We therefore have a corresponding pair of upper and lower bounds.
Solving the LP relaxation, and calculating the objective function at the feasible solution $(2,2)$ tell us that the optimal objective function value (often denoted $z^*$) must lie between 54 and 72.

A common framework for solving integer programs is to work to bring these bounds closer together.
Upper bounds can be tightened by finding better feasible solutions (with lower objective function values that are closer to the optimal solution).
Lower bounds can be tightened by solving ``partial relaxations'' where only some of the integrality constraints are enforced; these will have higher objective function values than a full relaxation, which are also closer to the optimal solution (but from the other side).
Ultimately, we determine a decreasing sequence of upper bounds $z_{UB}^1 \geq z_{UB}^2 \geq \cdots z^*$ and an increasing sequence of lower bounds $ z_{LB}^1 \leq z_{LB}^2 \ldots \cdots \leq z^*$.
When the difference between the bounds is small enough (say, within $\epsilon)$, we know that the feasible solution corresponding to the best-known upper bound is within $\epsilon$ of being optimal to the original integer program.

The success of such a framework depends on the strength of the bounds.
The tighter the upper and lower bounds are, the greater the guarantee we can provide for the solution we find.
Effectively solving these problems often requires intelligently exploiting the problem structure, and expert domain knowledge, in order to provide good feasible solutions (for upper bounds), and good partial relaxations (for lower bounds).
The branch and bound algorithm described in Section~\ref{sec:branchandbound} is one framework for doing this.
Before presenting this method, we will provide another concrete example of an integer program for facility location, and how it relates to the selection of good bounds.

Consider the case of locating a certain number of facilities in a region (say, warehouses, or fire stations).
Let $I$ denote the set of locations that demand what the facility provides, and $J$ the set of potential locations where facilities can be built.
Let $c_{ij}$ denote the cost of meeting the demand at location $i\in I$ from facility $j\in J$.
Let $f_j$ represent the fixed cost of locating a facility at site $j \in J$.
There are two sets of decision variables corresponding to location and assignment decisions.
The location-related decision variable $y_j$ takes the value 1 if a facility is located at site $j$ and 0 otherwise.
The decision variable $x_{ij}$ denotes the fraction of demand for center $i$ met by facility $j$.
The $x_{ij}$ variables can be continuous (its demand can be provided by multiple facilities), whereas the $y_j$ values must be integer (a facility cannot be ``half-built'').
This is called the uncapacitated facility location problem (UFLP), and can be formulated as follows:

\begin{align}
\textbf{O-UFLP}: &   &  \nonumber\\
\min_\mb{x,y} \quad & \sum_{i\in I} \sum_{j \in J} c_{ij}x_{ij} + \sum_{j\in J} f_{j}y_{j}  & \label{eqn:uflobj}\\
\mathrm{s.t.}    \quad & \sum_{j \in J} x_{ij} = 1 & \forall i \in I  \label{eqn:ufl1}\\
                      & \sum_{i \in I} x_{ij} \leq |I|y_j & \forall j \in J \label{eqn:ufl2}\\
                      & y_j \in \{0,1\} & \forall j \in J \label{eqn:ufl3}\\
                      & x_{ij} \geq 0 & \forall i \in I , \forall j \in J \label{eqn:ufl4} \,.
\end{align}

The objective of the UFLP is to locate facilities and assign demand points to facilities so that the overall facility location costs and transportation costs are minimized, as shown in equation~\eqn{uflobj}.
Constraint~\eqn{ufl1} ensures all demands are met.
Constraint~\eqn{ufl2} ensures that the demand points are assigned to open facilities only. 
(If $y_j = 0$, no facility is located at site $j$, so no demand can be served from there; if $y_j = 1$, then potentially all demand could be served from there.)

Now consider a new optimization model consisting of the same objective function~\eqn{uflobj}, along with constraints~\eqn{ufl1},~\eqn{ufl3}, and~\eqn{ufl4}, but with constraint~\eqn{ufl2} replaced by: $x_{ij} \leq y_j  \forall i \in I , \forall j \in J$.

\begin{align}
\textbf{S-UFLP}: &   &  \nonumber\\
\min_\mb{x,y} \quad & \sum_{i\in I} \sum_{j \in J} c_{ij}x_{ij} + \sum_{j\in J} f_jy_j  & \label{eqn:suflobj}\\
\mathrm{s.t.}    \quad & \sum_{j \in J} x_{ij} = 1 & \forall i \in I  \label{eqn:sufl1}\\
                      & x_{ij} \leq y_j  & \forall i \in I , \forall j \in J \label{eqn:sufl2}\\
                      & y_j \in \{0,1\} & \forall j \in J \label{eqn:sufl3}\\
                      & x_{ij} \geq 0 & \forall i \in I , \forall j \in J \label{eqn:sufl4} \,.
\end{align}

Let $P$ and $Q$ denote the feasible regions of the linear relaxations for models O-UFLP and S-UFLP, respectively.
We can show that $Q \subseteq P$, that is, any solution feasible to S-UFLP is also feasible to O-UFLP.
Any point which satisfies $x_{ij} \leq y_j$ must also satisfy $x_{ij} \leq |I|y_j$, as can be seen by summing up $x_{ij}$ along all values of $i$.
Now, when restricted to integer values of $\mb{y}$, the set of feasible solutions is the same in both cases.
However, when considering their LP relaxations, there are solutions which are feasible to O-UFLP which are not feasible to S-UFLP.
For example, if $y_1 = 1/2$, a solution of the form $x_{11} = 1$ and $x_{i1} = 0$ for $i = \myc{2, \ldots, n}$ satisfies $x_{i1} \leq |I|y_1$ but not $x_{11} \leq y_1$.
As a result of this, the optimal solution to a relaxation based on S-UFLP will have a greater (or equal) objective function value than the same relaxation based on O-UFLP, producing a tighter (and more useful) lower bound on the original integer program.
We therefore say that S-UFLP is a \emph{stronger} (or \emph{tighter}) formulation for UFLP.
This discussion can be generalized in the following result:
\begin{prp}
 Consider two feasible regions $P$ and $Q$ for the linear relaxation of the same minimization integer program $z = \min \{\mb{c^T x }: \mb{x}\in X \cap \bbz^n \}$
 with $Q \subseteq P$. If $z_{LP}^P = \min \{\mb{c^Tx}:\mb{x}\in P \}$ and  $z_{LP}^Q = \min \{\mb{c^T x}: \mb{x} \in Q \}$, then $z_{LP}^P \leq z_{LP}^Q$.
\end{prp}
 
\begin{figure}
\centering
\begin{tikzpicture}
\draw [ultra thick] [<->] (0,10) -- (0,0) -- (12,0);
\draw [fill=blue] (2,2)circle (3pt);
\draw [fill=blue] (4,2)circle (3pt);
\draw [fill=blue] (6,2)circle (3pt);
\draw [fill=blue] (8,2)circle (3pt);
\draw [fill=blue] (10,2)circle (3pt);

\draw [fill=blue] (2,4)circle (3pt);
\draw [fill=red] (4,4)circle (3pt);
\draw [fill=red] (6,4)circle (3pt);
\draw [fill=red] (8,4)circle (3pt);
\draw [fill=blue] (10,4)circle (3pt);

\draw [fill=blue] (2,6)circle (3pt);
\draw [fill=red] (4,6)circle (3pt);
\draw [fill=red] (6,6)circle (3pt);
\draw [fill=red] (8,6)circle (3pt);
\draw [fill=blue] (10,6)circle (3pt);

\draw [fill=blue] (2,8)circle (3pt);
\draw [fill=blue] (4,8)circle (3pt);
\draw [fill=blue] (6,8)circle (3pt);
\draw [fill=blue] (8,8)circle (3pt);
\draw [fill=blue] (10,8)circle (3pt);

\draw [blue, ultra thick] (3.6,3.6) -- (8.9,3.6) -- (8.2,6.6) -- (3.6,6.6) -- (3.6,3.6);
\node at (5,3.9) {\textbf{Q}};
\draw [blue, dashed, ultra thick] (2.5,2.5) -- (9.2,2.4) -- (9.4,7.8) -- (3.4,7.4) -- (2.5,2.5);
\node at (5,2.8) {\textbf{P}};

\end{tikzpicture}
\caption{Stronger and weaker IP formulations. \label{fig:form}}
\end{figure}
 
The intuition is shown in Figure~\ref{fig:form}.
Note that  $Q \cap \bbz^n = P \cap \bbz^n $.
However since $Q \subseteq P$, the corresponding LP solution will be larger.
Stronger formulations leading to tighter linear programming bounds will result in faster solutions, as long as the time need to solve the tighter linear programming relaxation is not significantly higher.

\subsection{Branch and bound algorithm}
\label{sec:branchandbound}

\index{optimization!integer!branch-and-bound algorithm|(}
Branch and bound is one of the most popular algorithms used to solve integer programs, and is the basis for many commercial integer programming solvers.
It is essentially an intelligent enumeration algorithm, based on the principle of ``divide and conquer.''
The original integer program is replaced by a series of easier-to-solve relaxations, whose solutions are used to further guide the search toward the optimal solution of the integer program.
The branch and bound algorithm is illustrated using the same integer programming example (see figure ~\ref{fig:bb}) used before:
\optimizex{\min}{y_1, y_2}{10y_1 + 26y_2}{
& 11y_1 + 3y_2  \geq 21  \\
&  6y_1 + 20y_2 \geq 39  \\
&   y_1,    y_2 \geq 0   \\
&   y_1,    y_2 \in \bbz
}We will first show how the algorithm works on this specific example, and then provide the general method below.

\begin{figure}
\centering
\resizebox{\textwidth}{!}{
\begin{tikzpicture}
[reg/.style={rectangle,draw=green!50,fill=green!20,thick,
    inner sep=3pt,minimum size=12mm},
fath/.style={rectangle,draw=red!50,fill=red!20,thick,
    inner sep=3pt,minimum size=12mm},
pre/.style={<-,shorten <=1pt,>=stealth',semithick},
post/.style={->,thick}]
\node[reg,align = center] (n1) {$z_{LP1}^* = 54$ \\ $y_1=1.5,y_2=1.5$ };
\node[reg,align = center, below right = 1.00cm and -0.2cm of n1] (n2)  {$z_{LP2}^* = 65.64$ \\ $y_1=1.36,y_2=2$ };
\node[reg,align = center, below left = 1.00cm and -0.2cm of n1] (n3)  {$z_{LP3}^* = 57.67$ \\ $y_1=3.17,y_2=1$ };
\path (n1) edge [post] node[auto,pos=0.75] {$y_2\geq 2$} (n2);
\path (n1) edge [post] node[auto,swap,pos=0.75] {$y_2\leq 1$} (n3);

\node[fath,align = center, below right = 1.00cm and 0.25cm of n2] (n4)  {$z_{LP4}^* = 72$ \\ $y_1=2,y_2=2$ };
\node[fath,align = center, below = 1.00cm  of n2] (n5)  {$z_{LP5}^* = 96.67$ \\ $y_1=1,y_2=3.33$ };
\path (n2) edge [post] node[auto,pos=0.75] {$y_1\geq 2$} (n4);
\path (n2) edge [post] node[auto,swap] {$y_1\leq 1$} (n5);

\node[reg,align = center, below = 1.00cm of n3] (n6)  {$z_{LP6}^* = 59.5$ \\ $y_1=4,y_2=0.75$};
\node[fath,align = center, below left = 1.00cm and 0.2cm of n3] (n7)  {LP7 \\ Infeasible};

\path (n3) edge [post] node[auto] {$y_1 \geq 4$} (n6);
\path (n3) edge [post] node[auto,swap,pos=0.75] {$y_1 \leq 3$} (n7);

\node[fath,align = center, below right= 1.00cm and 0.25 cm of n6] (n8)  {$z_{LP8}^* = 66$ \\ $y_1=4,y_2=1$};
\node[reg,align = center, below = 1.00cm of n6] (n9)  {$z_{LP9}^* = 65$ \\ $y_1=6.5,y_2=0$};

\path (n6) edge [post] node[auto,pos=0.75] {$y_2\geq 1$} (n8);
\path (n6) edge [post] node[auto,swap] {$y_2\leq 0$} (n9);

\node[fath,align = center, below right= 1.00cm and 0.25 cm of n9] (n10)  {$z_{LP10}^* = 70$ \\ $y_1=7,y_2=0$};
\node[fath,align = center, below = 1.00cm of n9] (n11)  {LP11 \\ Infeasible};

\path (n9) edge [post] node[auto,pos=0.8] {$y_1\geq 7$} (n10);
\path (n9) edge [post] node[auto,swap] {$y_1\leq 6$} (n11);

\end{tikzpicture}
}
\caption{Branch and bound example.  \label{fig:bb}}
\end{figure}

The branch and bound algorithm starts by initializing the upper bound $z_{UB} = \infty$ and lower bound $z_{LB} = -\infty$, since we do not yet know any information about the optimal solution.
The linear relaxation of the above integer program is solved.
Let $X_{LP1}$ denote the feasible region of the LP relaxation.
Thus $X_{LP1}=\{\mb{y} \in \bbr^2_+:11y_1 + 3y_2 \geq 21,6y_1 + 20y_2 \geq 39 \}$.
The optimal value of the linear programming objective function is $z_{LP1}^*=54$, and the solution is $y_1=1.5,y_2=1.5$.
The best lower bound is now updated to to 54, $z_{LB} = 54$.

Now, the optimal solution either has $y_2$ as an integer greater than or equal to 2, or as an integer less than or equal to 1.
We consider each option in turn.
(We could have also branched on $y_1$; it would be instructive for you to solve the problem that way.)
We generate two linear programs by adding the constraints $y_2 \geq 2$ to one, and and $y_2 \leq 1$ to the other.
The feasible region of the second linear program $X_{LP2} = X_{LP1} \cap \{y_2 \geq 2\}$ and the feasible region of the third linear program $X_{LP3} = X_{LP1} \cap \{y_2 \leq 1\}$.
Solving the two linear programs we get $z_{LP2}^*=65.64$ with the solution $y_1=1.36,y_2=2$ and $z_{LP3}^*=57.67$ with the solution $y_1=3.17,y_2=1$.

Let us first examine LP2.
Since $y_2$ is already an integer in its optimal solution, we don't need to branch again on $y_2$.
However, $y_1$ is not an integer, and in the optimal solution it must be either greater than or equal to 2, or less than or equal to 1.
So the constraints $y_1 \geq 2$ and $y_2 \leq 1$ are added to the feasible regions of LP2 to generate new linear programs LP4 and LP5, with the respective feasible regions $X_{LP4} = X_{LP2} \cap \{y_1 \geq 2\}$ and $X_{LP5} = X_{LP2} \cap \{y_1 \leq 1\}$. 

Solving LP4 we get the optimal objective value to be 74 with the solution $y_1=2,y_2=2$.
The solution to this linear programming relaxation is an integer solution which is feasible for the original integer program.
Therefore $z_{LP5}^*$ is as an upper bound to the optimal integer solution.
Moreover since $z_{LP5}^* < z_{UB}$ we update $z_{UB}=74$; this is the first upper bound we have the optimal objective function value..
Since this is an integer solution, we do not search further in this direction.
This is called ``fathoming by optimality.''

Solving LP5 we get the optimal objective value to be 96.67 with the solution $y_1=1,y_2=3.33$.
Since the LP objective function is greater than the current best upper bound 74, there is no point in further exploring in this direction also and the search is stopped here.
This is called ``fathoming by bound.''
Notice that this cuts off a considerable portion of the feasible region.
We know we do not have to explore it any further, since the best possible objective function values there are at least 96.67, whereas we have already found a solution $(2,2)$ with an objective of 74.

Let us now look at LP3.
Since $y_2$ has integer solution we don't need to branch on $y_2$.
But with respect to $y_1$, the optimal solution is either an integer at least 4, or at most 3.
So the constraints $y_1 \geq 4$ and $y_1 \leq 3$ is added to the feasible region of LP3 to get new linear programs LP6 and LP7, with respective feasible regions $X_{LP6} = X_{LP3} \cap \{y_1 \geq 4\}$ and $X_{LP7} = X_{LP3} \cap \{y_2 \leq 3\}$.

LP7 turns out to be infeasible, so the search is stopped along this direction.
This is called ``fathoming by infeasibility.''
LP6, on the other hand, has an optimal objective function value of 59.5 with $y_1=4,y_2=0.75$.
We branch again on $y_2$, adding the constraints $y_2 \geq 1$ and $y_2 \leq 0$ to the feasible region of LP6 to obtain new linear programs LP8 and LP9 with feasible regions $X_{LP8} = X_{LP6} \cap \{y_2 \geq 1\}$ and $X_{LP9} = X_{LP3} \cap \{y_2 \leq 0\}$. 

LP8 has an integer solution $y_1=4,y_2=1$, with an objective function value of 66.
This solution is feasible for the original integer program.
Moreover since $z_{LP8}^* < z_{UB}$, we have found a better feasible solution, and the upper bound is updated to 66: $z_{UB}=66$.
We fathom the search in this direction by integrality. 

Looking at the solution of LP9, we generate two linear programs, LP10 and LP11, by adding the constraints $y_1 \geq 7$ and $y_2 \leq 6$ to the feasible region of LP9.
The feasible region of the new linear programs are $X_{LP10} = X_{LP9} \cap \{y_1 \geq 7\}$ and $X_{LP11} = X_{LP9} \cap \{y_1 \leq 6\}$.
LP11 is infeasible, and therefore we fathom by infeasibility.
LP10 has an integer solution, but the objective function value is greater than the current upper bound.
Therefore we fathom by bound as well as integrality.
There are no more directions left to search.
Therefore the optimal solution corresponds to the current best upper bound of 66. 

The general steps of branch and bound are shown in Algorithm~\ref{alg:branchandbound}.
\index{optimization!integer!branch-and-bound algorithm|)}

\begin{algorithm}
    \textbf{Initialization}: Set $z_{UB} \leftarrow -\infty$ and $z_{LB} \leftarrow \infty$. Add the LP relaxation of the original integer program to a queue $\mb{Q}$.\;
    \While{ $\mb{Q}$ is not empty}
    {
        Select an element CurrLP from $\mb{Q}$ and remove it\;
        Attempt to solve CurrLP\;
        \eIf{CurrLP is infeasible}
        {
            Continue while loop (fathom by infeasibility)\;
        }
        {
           Let $\mb{x}$ be the optimal solution to CurrLP, and $z$ the objective value\;
           \eIf{$\mb{x}$ satisfies all of the original integrality constraints}
           {
               $z_{UB} \leftarrow \min \myc{z_{UB}, z}$\;
               If $z_{UB}$ changed, mark $\mb{x}$ as the best known solution\;
               Continue while loop (fathom by optimality)\;
           }
           {
                   \eIf{$z \geq z_{UB}$}
                   {
                       Continue while loop (fathom by bound)\;
                   }
                   {
                           $z_{LB} \leftarrow \max \myc{z_{LB}, z}$\;
                           Choose a fractional variable to branch on\;
                           Generate two new LPs by adding $\geq$ and $\leq$ constraints to the feasible region of CurrLP\;
                           Add the two new LPs to $\mb{Q}$\;
                }
           }
      }
    }
\caption{Branch and bound algorithm.\label{alg:branchandbound}}
\end{algorithm}

\subsection{Total unimodularity}
\label{sec:totallyunimodular}

In general, we cannot exactly solve an integer program only by solving its linear programming relaxation (although we can use these LP relaxations in a method such as branch and bound).
There is one notable exception, which shows up frequently enough in network optimization problems to be worth knowing.

\begin{dfn}
A matrix $\mb{A}$ is \emph{totally unimodular}\index{totally unimodular matrix}\index{matrix!totally unimodular} if the determinant of every square submatrix of $\mb{A}$ has the value $-1$, $0$, or $1$. 
\end{dfn}

Consider an integer linear program $z = \min \{\mb{c^T} \mb{x}:\mb{x}\in X \cap \bbz^n \}$ where $X=\{\mb{x} \in \bbr^n:\mb{Ax} \leq \mb{b}, \mb{x} \geq 0 \}$, and the corresponding linear programming relaxation $z_{LP} =  \min \{\mb{c^T} \mb{x}:\mb{x}\in X \}$.
If the matrix $\mb{A}$ is totally unimodular, then the LP relaxation has integral solutions. The vertex points of $X$ are integer points. 
\emph{That is, if the coefficients representing the constraints form a totally unimodular matrix, we can ignore the integer constraints altogether.
The resulting linear program will have an integer optimal solution which is feasible to the original integer program.}
In such a case, we can simply apply a linear programming method to find the optimal solution.

\index{shortest path|(}
The concept of total unimodularity is explained using a modification to the standard shortest path problem from Section~\ref{sec:shortestpath}.
In this modification, each link is associated with \emph{two} attributes: a cost for traversing it, and the amount of some other resource consumed when traveling it (time, battery charge, etc.)
Consider the network shown in figure~\ref{fig:network}.
The network has four nodes and five arcs. Let $c_{12}$ and $r_{12}$ denote the cost and the amount of resource consumed by traversing arc  $(1,2)$.
The goal of the shortest path problem is to determine the minimum cost path from $1$ to $4$.

\begin{figure}
\centering
\begin{tikzpicture}
  [
  place/.style={circle,draw=black,fill=blue!20,thick, inner sep=0pt,minimum size=8mm,font=\sffamily\small}, 
  transition/.style={rectangle,draw=black!50,fill=black!20,thick, inner sep=0pt,minimum size=4mm},
  bend angle=45,node distance=4cm,
  pre/.style={<-,shorten <=1pt,>=stealth',semithick},
  post/.style={->,shorten >=1pt,>=stealth',semithick},
  cost/.style={sloped, font=\sffamily\small},
  nscost/.style={font=\sffamily\small}
  ]

  \node[place] (A) {$1$};
  \node[place] (B) [below right of=A] {$2$};
  \node[place] (C) [above right of=A] {$3$};
  \node[place] (D) [above right of=B] {$4$};

  \path
  (A) edge[post] node[cost, below] {$(c_{12}, r_{12})$} (B)
      edge[post] node [cost, above] {$(c_{13}, r_{13})$} (C)

  (B) edge[post] node [nscost,, right] {$(c_{23}, r_{23})$} (C)
       edge[post] node [cost, below] {$(c_{24}, r_{24})$} (D)

  (C) edge[post] node [cost,above] {$(c_{34}, r_{34})$} (D);    
  \end{tikzpicture}
\caption{Network for shortest path example  \label{fig:network}}
\end{figure}

Let $x_{12}$,$x_{13}$,$x_{23}$,$x_{24}$, and $x_{34}$ denote the decision variables.
These are equal to 1 if the corresponding arc lies on the shortest path, and 0 otherwise.
The shortest path formulation for the above network is:
\optimizex{\min}{}
{c_{12}x_{12} + c_{13}x_{13} + c_{23}x_{23} + c_{24}x_{24} + c_{34}x_{34}}
{
&    x_{12} + x_{13} &= 1          \\
&  -x_{12} + x_{23}  + x_{24} &= 0 \\
&  -x_{13} - x_{23}  + x_{34} &= 0 \\
&   -x_{24}  - x_{34} &= -1 \\
&  x_{12},  x_{13}, x_{23}, x_{24}, x_{34} & \in \{0,1\}
}The constraint matrix is
\begin{equation*}
\mb{A} = 
\begin{bmatrix} 
1 & 1 & 0 & 0 & 0 \\
-1 & 0 & 1 & 1 & 0 \\
0 & -1 & -1 & 0 & 1 \\
0 & 0 & 0 & -1 & -1 
\end{bmatrix}
\,.
\end{equation*}
All the square submatrices of $A$ have determinant $0$, $1$, or $-1$.
Therefore we can solve the shortest path problem by solving the corresponding linear programming relaxation. 
(The shortest path problem \emph{could} be solved using branch and bound; but it is  better to solve it as a linear program, and even better to solve it using a specialized method like those in Section~\ref{sec:shortestpath}.)

\index{shortest path!resource-constrained|(}
The resource-constrained shortest path problem is a variant where the shortest path is sought, but with an upper bound on the total amount of resources consumed.
An example might be finding the fastest route in a network with toll roads, where the total monetary cost does not exceed a given budget.
We represent this with an additional constraint in the shortest path formulation:
\begin{equation*}
r_{12}x_{12} + r_{13}x_{13} + r_{23}x_{23} + r_{24}x_{24} + r_{34}x_{34} \leq R
\,.
\end{equation*}
In this case the constraint matrix becomes
\begin{equation*}
\mb{A} = 
\begin{bmatrix} 
1 & 1 & 0 & 0 & 0 \\
-1 & 0 & 1 & 1 & 0 \\
0 & -1 & -1 & 0 & 1 \\
0 & 0 & 0 & -1 & -1 \\
r_{12} & r_{13} & r_{23} & r_{24} & r_{34} 
\end{bmatrix}
\,.
\end{equation*}
Depending on the values of the resources consumed, the constraint matrix need not be totally unimodular.
In this case the linear programming relaxation does not always have an integer solution.
For this reason, the resource-constrained shortest path problem is much harder to solve than the traditional shortest path problem.
\index{shortest path|)}
\index{shortest path!resource-constrained|)}
\index{optimization!integer|)}

\section{Metaheuristics}
\label{sec:metaheuristics}

\index{optimization!metaheuristics|(}
\index{optimization!heuristics|see {optimization, metaheuristics}}
In optimization, there is often a tradeoff between how widely applicable a solution method is, and how efficient or effective it is at solving specific problems.
For any one specific problem, a tailor-made solution process is likely much faster than a generally-applicable method, but generating such a method requires more effort and specialized knowledge, and is less ``rewarding'' in the sense that the method can only be more narrowly applied.
Some of the most general techniques are \emph{heuristics}, which are not guaranteed to find the global optimum solution, but tend to work reasonably well in practice.
In practice, they only tend to be applied for very large or very complicated problems which cannot be solved exactly in a reasonable amount of time with our current knowledge.

Many engineers are initially uncomfortable with the idea of a heuristic.
After all, the goal of an optimization problem is to find an optimal solution, so why should we settle for something which is only approximately ``optimal,'' often without any guarantees of how approximate the solution is?
First, for very complicated problems a good heuristic can often return a reasonably good solution in much less time than it would take to find the exact, global optimal solution.
For many practical problems, the cost improvement from a reasonably good solution to an exactly optimal one is not worth the extra expense (both time and computational hardware) needed, particularly if the heuristic gets you within the margin of error based on the input data.

Heuristics are also very, very common in psychology and nature.
If I give someone a map and ask them to find the shortest-distance route between two points in a city by hand, they will almost certainly not formulate an mathematical model and solve it to provable optimality.
Instead, they use mental heuristics (rules of thumb based on experience) and can find paths which are actually quite good.
Many of the heuristics are inspired by things seen in nature.

An example is how ant colonies find food.
When a lone wandering ant encounters a food source, it returns to the colony and lays down a chemical pheromone.
Other ants who stumble across this pheromone begin to follow it to the food source, and lay down more pheromones, and so forth.
Over time, more and more ants will travel to the food source, taking it back to the colony, until it is exhausted at which point the pheromones will evaporate.
Is this method the optimal way to gather food?  Perhaps not, but it performs well enough for ants to have survived for millions of years!

Another example is the process of evolution through natural selection.
The human body, and many other organisms, function remarkably well in their habitats, even if their biology is not exactly ``optimal.''\footnote{For instance, in humans the retina is ``backwards,'' creating a blind spot; in giraffes, the laryngeal nerve takes an exceptionally long and roundabout path; and the descent of the testes makes men more vulnerable to hernias later in life.}  One of the most common heuristics in use today, and one described below, is based on applying principles of natural selection and mutation to a ``population'' of candidate solutions to an optimization problem, using an evolutionary process to identify better and better solutions over time.
This volume describes two heuristics: \emph{simulated annealing}, and \emph{genetic algorithms}, both of which can be applied to many different optimization problems.

Simulated annealing was proposed by~\cite{kirkpatrick83}, building on an earlier algorithm of~\cite{metropolis53}.
Genetic algorithms for optimization were first proposed by \cite{rechenberg73} and \cite{schwefel77}.
Other metaheuristics not described in the text include GRASP~\citep{feo95}, tabu search~\citep{glover89,glover90}, ant swarms~\citep{dorigo04}, and bee swarms~\citep{pham05}, among many others.
Note that there is some recent criticism of some of the more fanciful analogies among novel metaheuristics~\citep{sorensen15,camachovillalon22,camachovillalon23}; we agree that it is best to focus on how these heuristics balance finding good solutions without being trapped in local optima, rather than on the details of exactly how a particular real-world process is encoded.

\subsection{Simulated annealing}
\label{sec:simulatedannealing}

\index{optimization!metaheuristics!simulated annealing|(}
Simulated annealing is a simple heuristic that makes an analogy to metallurgy, but this analogy is best understood after the heuristic is described.

As a better starting point for understanding simulated annealing, consider the following ``local search'' heuristic.
(For now, descriptions will be a bit vague; more precise definitions will follow soon.)  Local search proceeds as follows:
\begin{enumerate}
   \item Choose some initial feasible solution $\mb{x} \in X$, and calculate the value of the objective function $f(\mb{x})$.
   \item Generate a new feasible solution $\mb{x'} \in X$ which is close to the current solution $\mb{x}$.
   \item Calculate $f(\mb{x'})$.
   \item If $f(\mb{x'}) \leq f(\mb{x})$, the new solution is better than the old solution.
So update the current solution by setting $\mb{x}$ equal to $\mb{x'}$.

   \item Otherwise, the old solution is better, so leave $\mb{x}$ unchanged.
   \item Return to step 2 and repeat until we are unable to make further progress.
\end{enumerate}
One way to visualize local search is a hiker trying to find the lowest elevation point in a mountain range.
In this analogy, the park boundaries are the feasible set, the elevation of any point is the objective function, and the location of the hiker is the decision variable.
In local search, starting from his or her initial position, the hiker looks around and finds a nearby point which is lower than their current point.
If such a point can be found, they move in that direction and repeat the process.
If every neighboring point is higher, then they stop and conclude that they have found the lowest point in the park. 

It is not hard to see why this strategy can fail; if there are multiple local optima, the hiker can easily get stuck in a point which is not the lowest (Figure~\ref{fig:localoptima}).
Simulated annealing attempts to overcome this deficiency of local search by allowing a provision for occasionally moving in an uphill direction, in hopes of finding an even lower point later on.
Of course, we don't always want to move in an uphill direction and have a preference for downhill directions, but this preference cannot be absolute if there is any hope of escaping local minima.

\stevefig{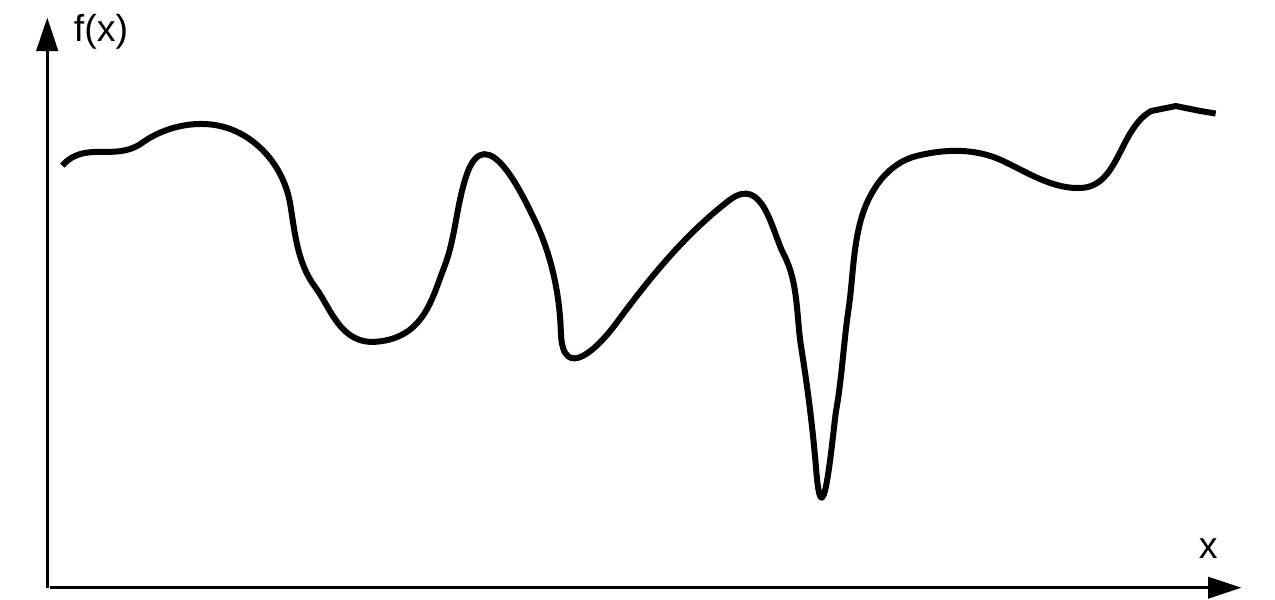}{Moving downhill from the current location may not lead to the global optimum.}{0.7\textwidth}

Simulated annealing accomplishes this by making the decision to move or not probabilistic, introducing a \emph{temperature} parameter $T$ which controls how likely you are to accept an ``uphill'' move.
When the temperature is high, the probability of moving uphill is large, but when the temperature is low, the probability of moving uphill becomes small.
In simulated annealing, the temperature is controlled with a \emph{cooling schedule}.
Initially, the temperature is kept high to encourage a broad exploration of the feasible region.
As the temperature decreases slowly, the solution is drawn more and more to lower areas.
The cooling schedule can be defined by an initial temperature $T_0$, a final temperature $T_f$, the number of search iterations $n$ at a given temperature level, and a scaling factor $k \in (0, 1)$ which is applied every $n$ iterations to reduce the temperature.

Finally, because the search moves both uphill and downhill, there is no guarantee that the final point of the search is the best point found so far.
So, it is worthwhile to keep track of the best solution $\mb{\hat{x}}$ encountered during the algorithm.
(This is analogous to the hiker keeping a record of the lowest point observed, and returning to that point when done searching.)

So, the simulated annealing algorithm can be stated as follows:
\begin{enumerate}
   \item Choose some initial feasible solution $\mb{x} \in X$, and calculate the value of the objective function $f(\mb{x})$.
   \item Initialize the best solution to the initial one $\mb{\hat{x}} \leftarrow \mb{x}$.
   \item Set the temperature to the initial temperature: $T \leftarrow T_0$.
   \item Repeat the following steps $n$ times:
   \begin{enumerate}[(a)]
      \item Generate a new feasible solution $\mb{x'}$ which neighbors $\mb{x}$.
      \item If $f(\mb{x'}) < f(\mb{\hat{x}})$, it is the best solution found so far, so update $\mb{\hat{x}} \leftarrow \mb{x'}$.
      \item If $f(\mb{x'}) \leq f(\mb{x})$, it is a better solution than the current one, so update $\mb{x} \leftarrow \mb{x'}$.
      \item Otherwise, update $\mb{x} \leftarrow \mb{x'}$ with probability $\exp(-[f(\mb{x'}) - f(\mb{x})] / T)$.
   \end{enumerate}
   \item If $T > T_f$, then reduce the temperature ($T \leftarrow kT$) and return to step 4.
   \item Report the best solution found $\mb{\hat{x}}$.
\end{enumerate}

The key step is step 4d.
Notice how the probability of ``moving uphill'' depends on two factors: the temperature, and how much the objective function will increase.
The algorithm is more likely to accept an uphill move if it is only slightly uphill, or if the temperature is high.
The exponential function captures these effects while keeping the probability between 0 and 1.
A few points require explanation.

\textbf{How should the cooling schedule be chosen?}
Unfortunately, it is hard to give general guidance here.
Heuristics often have to be ``tuned'' for a particular problem: some problems do better with higher temperatures and slower cooling ($k$ values closer to 1, $n$ larger), others work fine with faster cooling.
When you use simulated annealing, you should try different variations of the cooling schedule to identify one that works well for your specific problem.

\textbf{How should an initial solution be chosen?}
It is often helpful if the initial solution is relatively close to the optimal solution.
For instance, if the optimization problem concerns business operations, the current operational plan can be used as the initial solution for further optimization.
However, it's easy to go overboard with this.
You don't want to spent so long coming up with a good initial solution that it would have been faster to simply run simulated annealing for longer starting from a worse solution.
The ideal is to think of a good, quick rule of thumb for generating a reasonable initial solution; failing that, you can always choose the initial solution randomly. 

\textbf{How do I define a neighboring solution for step 4a?}
Again, this is problem-specific, and one of the decisions that must be made when applying simulated annealing.
A good neighborhood definition should involve points which are ``close'' to the current solution in some way, but ensure that feasible solutions are connected in the sense that any two feasible solutions can be reached by a chain of neighboring solutions.
For the examples in Section~\ref{sec:optimizationexamples}, some possibilities (emphasis on \emph{possibilities}, there are other choices) are:
\begin{description}
   \item[Transit frequency setting problem: ] The decision variables $\mb{n}$ are the number of buses assigned to each route.
Given a current solution $\mb{n}$, a neighboring solution is one where exactly one bus has been assigned to a different route.

   \item[Scheduling maintenance: ] The decision variables are $\mb{x}$, indicating where and when maintenance is performed, and $\mb{c}$, indicating the condition of the facilities.
In this problem, the constraints completely define $\mb{c}$ in terms of $\mb{x}$, so for simulated annealing it is enough to simply choose a feasible solution $\mb{x}$, then calculate $\mb{c}$ using the state evolution equations.
In this problem, there is no reason to perform less maintenance than the budget allows each year.
If the initial solution is chosen in this way, then a neighboring solution might be one where one maintenance activity on a facility is reassigned to another facility that same year.
If the resulting solution is infeasible (because the resulting $\mb{c}$ values fall below the minimum threshold), then another solution should be generated.
   \item[Facility location problem: ] The decision variables are the intersections that the three terminals are located at, $L_1$, $L_2$, and $L_3$.
Given current values for these, in a neighboring solution two of the three terminals are at the same location, but one of the three has been assigned to a different location.
   \item[Shortest path problem: ] The decision variables specify a path between the origin and destination.
Given a current path, a neighboring path is one which differs in only intersection.
(Can you think of a way to express this mathematically?)
\end{description}

\textbf{How do I perform a step with a given probability?}
Most programming languages have a way to generate a uniform random variable between 0 and 1.
If we want to perform a step with probability $p$, generate a random number from the continuous uniform $(0,1)$ distribution.
If this number is less than $p$, perform the step; otherwise, do not.

\subsubsection{Example with facility location}

This section demonstrates simulated annealing, on an instance of the facility location problem from Section~\ref{sec:optimizationexamples}.
In this problem instance, the road grid consists of ten north-south and ten east-west streets.
The cost of locating a terminal at each of the locations is shown in Figure~\ref{fig:facilitycost} (these numbers were randomly generated).
There are 30 customers, randomly located throughout the grid, as shown in Figure~\ref{fig:customerlocations}.
In these figures, the coordinate system $(x,y)$ is used to represent points, where $x$ represents the number of blocks to the right of the upper-left corner, and $y$ represents the number of blocks below the upper-left corner.

\genfig{facilitycost}{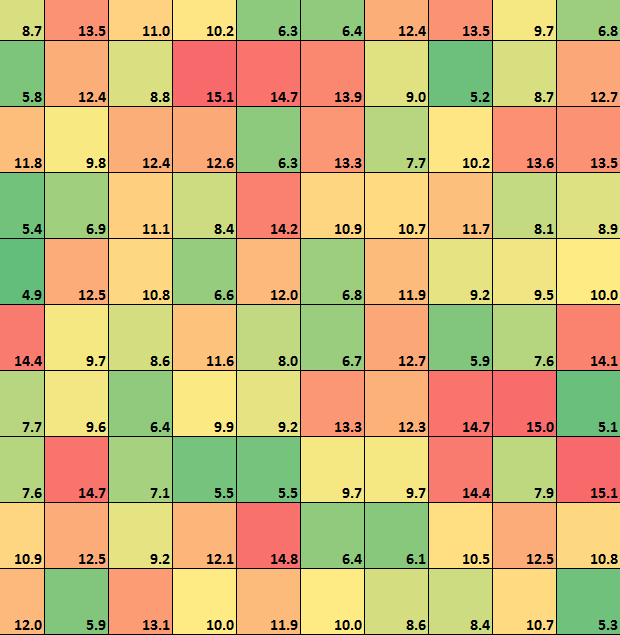}{Cost of locating facilities at each intersection.}{width=0.8\textwidth}
\genfig{customerlocations}{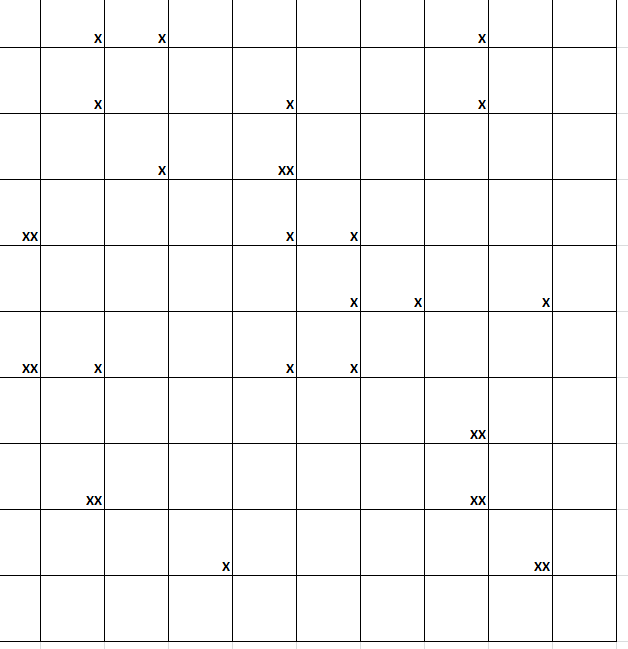}{Locations of customers.}{width=0.8\textwidth}

The cooling schedule is as follows: the initial temperature is $T_0 = 1000$, the final temperature is $T_f = 100$, $k = 0.75$, and the temperature is reduced every 8 iterations.
(These values were chosen through trial-and-error, and work reasonably well for this problem.)  The first several iterations are shown below.

The algorithm begins with an initial solution, generated randomly.
Suppose this solution locates the three facilities at coordinates (5,1), (8,6), and (6,1), respectively; the total cost associated with this solution is 158.0 cost units.
A neighboring solution is generated by picking one of the facilities (randomly) and assigning it to another location (randomly).
Suppose that the third facility is reassigned to location (1,0), so the new candidate solution locates the facilities at (5,1), (8,6), and (1,0).
This solution has a cost of 136.4 units, which is lower than the current solution.
So, simulated annealing replaces the current solution with the new one.

In the next iteration, suppose that the first facility is randomly chosen to be reassigned to (8,2), so the candidate solution is (8,2), (8,6), and (1,0).
This solution has a cost of 149.1, which is higher than the current cost of 136.4.
The algorithm will not always move to such a solution, but will only do so with probability calculated as in Step 4c: 
\[ p = \exp(-[f(\mb{x'}) - f(\mb{x})] / T) = \exp(-[149.1 - 136.4] / 100) = 0.881 \,.\]
This probability is high, because (being one of the early iterations) the temperature is set high.
If the temperature were lower, the probability of accepting this move would be lower as well.

Supposing that the move is accepted, the algorithm replaces the current solution with the candidate and continues as before.
If, on the other hand, the move is rejected, the algorithm generates another candidate solution based on the same current solution (5,1), (8,6), and (1,0).
The algorithm continues in the same way, reducing the temperature by 25\% every 8 iterations. 

The progress of the algorithm until termination is shown in Figure~\ref{fig:saprogress}.
The solid line shows the cost of the current solution, while the dashed line tracks the cost of the best solution found so far.
A few observations are worth making.
First, in the early iterations the cost is highly variable, but towards termination the cost becomes more stable.
This is due to the reduction in temperature which occurs over successive iterations.
When the temperature is high, nearly any move will be accepted so one expects large fluctuations in the cost.
When the temperature is low, the algorithm is less likely to accept cost-increasing moves, so fewer fluctuations are seen.
Also notice that the best solution was found shortly after iteration 700.
The final solution is not the best, although it is close.

\genfig{saprogress}{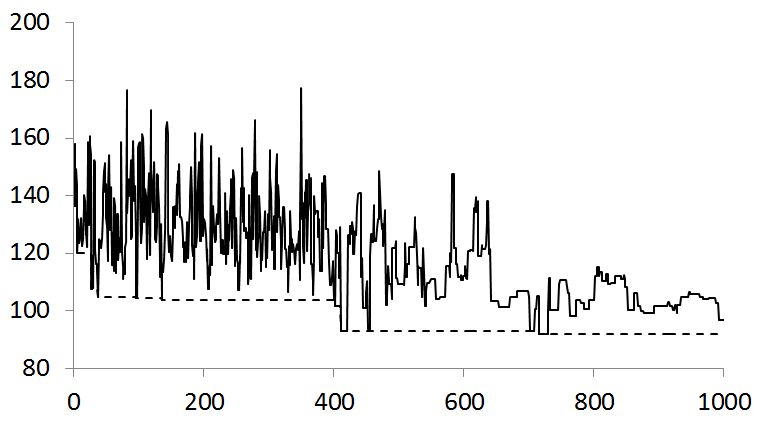}{Progress of simulated annealing.
Solid line shows current solution cost, dashed line best cost so far.}{width=0.8\textwidth}

The best facility locations found by simulated annealing are shown in Figure~\ref{fig:optimalsolution}, with a total cost of 92.1.
\index{optimization!metaheuristics!simulated annealing|)}

\genfig{optimalsolution}{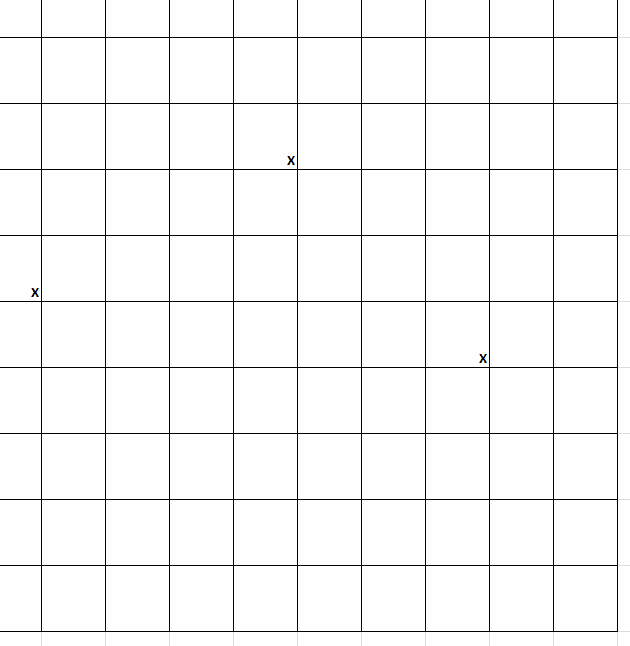}{Locations of facilities found by simulated annealing (cost 92.1).}{width=0.7\textwidth}

\subsection{Genetic algorithms}
\label{sec:geneticalgorithm}

\index{optimization!metaheuristics!genetic algorithm|(}
Another heuristic is the genetic algorithm.
In contrast to simulated annealing, where there is a single search happening through the feasible region (recall the analogy with the hiker), genetic algorithms maintain a larger \emph{population} of solutions which reflect a diverse set of points within the feasible region.
Genetic algorithms work by attempting to improve the population from one iteration to the next.

The process of improvement is intended to mimic the processes of natural selection, reproduction, and mutation which are observed in biology.
For the sake of our purposes, \emph{natural selection} means identifying solutions in the population which have good (low) values of the objective function.
The hope is that there are certain aspects or patterns in these solutions which make them good, which can be maintained in future generations.
Given an initial population as the first generation, subsequent generations are created by choosing good solutions from the previous generation, and ``breeding'' them with each other in a process called \emph{crossover} which mimics sexual reproduction: new ``child'' solutions are created by mixing characteristics from two ``parents.''  Lastly, there is a random \emph{mutation} element, where solutions are changed externally with a small probability.
If all goes well, after multiple generations the population will tend towards better and better solutions to the optimization problem.

At a high level, the algorithm is thus very straightforward and presented below.
However, at a lower level, there are a number of careful choices which need to be made about how to implement each step for a particular problem.
Thus, each of the steps is described in more detail below.
The high-level version of genetic algorithms is as follows:
\begin{enumerate}
   \item Generate an initial population of $N$ feasible solutions (generation 0), set generation counter $g \leftarrow 0$.
   \item Create generation $g + 1$ in the following way, repeating each step $N$ times:
   \begin{enumerate}[(a)]
      \item Choose two ``parent'' solutions from generation $g$.
      \item Combine the two parent solutions to create a new solution.
      \item With probability $p$, mutate the new solution.
   \end{enumerate}
   \item Increase $g$ by 1 and return to step 2 unless done.
\end{enumerate}
Although not listed explicitly, it is a good idea to keep track at every stage of the best solution $\mb{\hat{x}}$ found so far.
Just like in simulated annealing, there is no guarantee that the best solution will be found in the last generation.
So, whenever a new solution is formed, calculate its objective function value, and record the solution if it is better than any found so far.

Here are explanations of each step in more detail.
It is very important to realize that \emph{there are many ways to do each of these steps} and that what is presented below is one specific example intended to give the general flavor while sparing unnecessary complications at this stage.

\textbf{Generate an initial population of feasible solutions: } In contrast with simulated annealing, where we had to generate an initial feasible solution, with genetic algorithms we must generate a larger population of initial feasible solutions.
While it is still desirable to have these initial feasible solutions be reasonably good, it is also very important to have some diversity in the population.
Genetic algorithms work by combining characteristics of different solutions together.
If there is little diversity in the initial population the difference between successive generations will be small and progress will be slow.

\textbf{Selection of parent solutions: } Parent solutions should be chosen in a way that better solutions (lower objective function values) are more likely to be chosen.
This is intended to mimic natural selection, where organisms better adapted for a particular environment are more likely to reproduce.
One way to do this is through \emph{tournament selection}, where a number of solutions from the old generation are selected randomly, and the one with the best objective function value is chosen as the first parent.
Repeating the ``tournament'' again, randomly select another subset of solutions from the old generation, and choose the best one as the second parent.
The number of entrants in the tournament is a parameter that you must choose.

\textbf{Combining parent solutions: }  This is perhaps the trickiest part of genetic algorithms: how can we combine two feasible solutions to generate a new feasible solution which retains aspects of both parents?
The exact process will differ from problem to problem.
Here are some ideas, based on the example problems from Section~\ref{sec:optimizationexamples}:
\begin{description}
   \item[Transit frequency setting problem: ] The decision variables are the number of buses on each route.
Another way of ``encoding'' this decision is to make a list assigning each bus in the fleet to a corresponding route.
(Clearly, given such a list, we can construct the $n_r$ values by counting how many times a route appears.)  Then, to generate a new list from two parent lists, we can assign each bus to either its route in one parent, or its route in the other parent, choosing randomly for each bus.

   \item[Scheduling maintenance: ] As described above, optimal solutions must always exhaust the budget each year.
So, along the lines of the transit frequency setting problem, make a list of the number of maintenance actions which will be performed in each year, noting the facility assigned to each maintenance action in the budget.
To create a new solution, for each maintenance action in the budget choose one of the facilities assigned to this action from the two parents, selecting randomly.
Once this list is obtained, it is straightforward to calculate $\mb{x}$ and $\mb{c}$.
   \item[Facility location problem: ] For each facility in the child solution, make its location the same as the location of that facility in one of the two parents (chosen randomly).
\end{description}

\textbf{Mutating solutions: } Mutation is intended to provide additional diversity to the populations, avoiding stagnation and local optimal solutions.
Mutation can be accomplished in the same way that neighbors are generated in simulated annealing.
The probability of mutation should not be too high --- as in nature, most mutations are harmful --- but enough to escape from local optima.
The mutation probability $p$ can be selected by trial and error, determining what works well for your problem. 

\subsubsection{Example with facility location}

This section demonstrates genetic algorithms on the same facility location problem used to demonstrate simulated annealing.
The genetic algorithm is implemented with a population size of 100, over ten generations.
Tournaments of size 3 are used to identify parent solutions, and the mutation probability is 0.05.
The initial population is generated by locating all the facilities are completely at random.

To form the next generation, 100 new solutions need to be created, each based on combining two solutions from the previous generation.
These two parents are chosen using the tournament selection rule, as shown in Figure~\ref{fig:newchromosome}.
The two winners of the tournament locate the three facilities at (2,3), (1,7), (8,7); and (1,4), (7,6), (5,0) respectively.
These combine in the following way:
\begin{enumerate}
   \item The first facility is located either at (2,3) or (1,4); randomly choose one of them, say, (2,3).
   \item The second facility is located either at (1,7) or (7,6); randomly choose one of them, say, (1,7).
   \item The third facility is located either at (8,7) or (5,0); randomly choose one of them, say, (5,0).
\end{enumerate}
This gives a new solution (2,3), (1,7), (5,0) in the next generation, as shown in Figure~\ref{fig:newchromosome}.
With 5\% probability, this solution will be ``mutated.''  If this solution is selected for mutation, one of the three facility is randomly reassigned to another location.
For instance, the facility at (1,7) may be reassigned to (0,4).
This process is repeated until all 100 solutions in the next generation have been created.

\stevefig{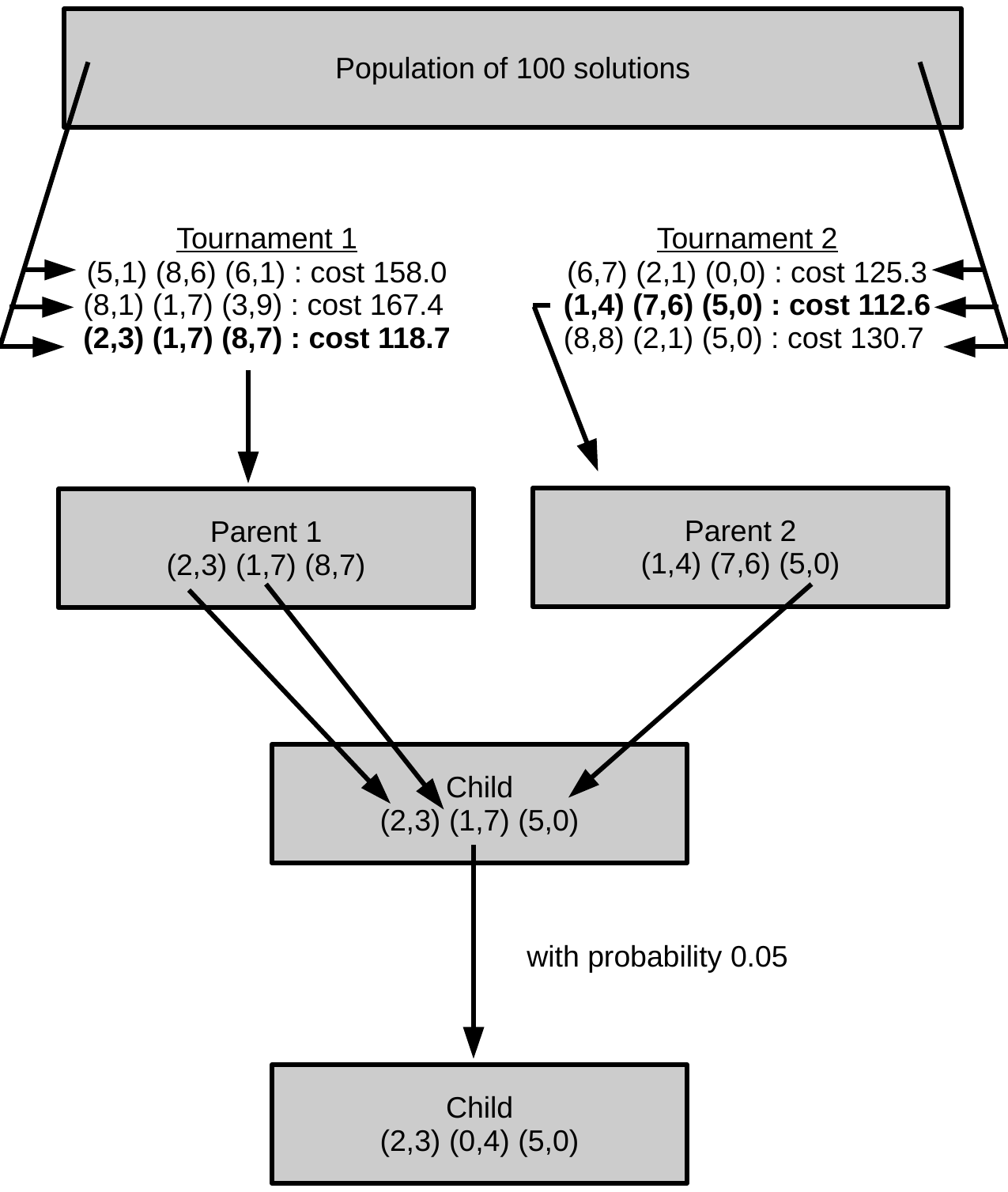}{Generation of a new solution: reproduction and mutation.}{0.8\textwidth}

Figure~\ref{fig:gaprogress} shows the progress of the algorithm over ten generations.
The solid line shows the average cost in each generation.
The dashed line shows the cost of the best solution found so far, and the crosses show the cost of each of the solutions comprising the generations.
Since lower-cost alternatives are more likely to win the tournaments and be selected as parents, the average cost of each generation decreases.
Figure~\ref{fig:geneticalgorithm} shows the locations of the terminals in the best solution found.
\index{optimization!metaheuristics!genetic algorithm|)}
\index{optimization!metaheuristics|)}
\index{optimization|)}

\genfig{gaprogress}{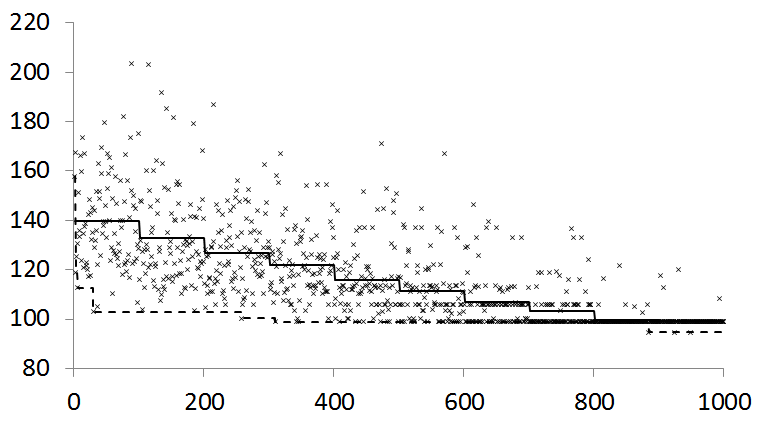}{Progress of genetic algorithm.
Solid line shows average generation cost, dashed line best cost so far, crosses cost of individual solutions.}{width=0.8\textwidth}

\genfig{geneticalgorithm}{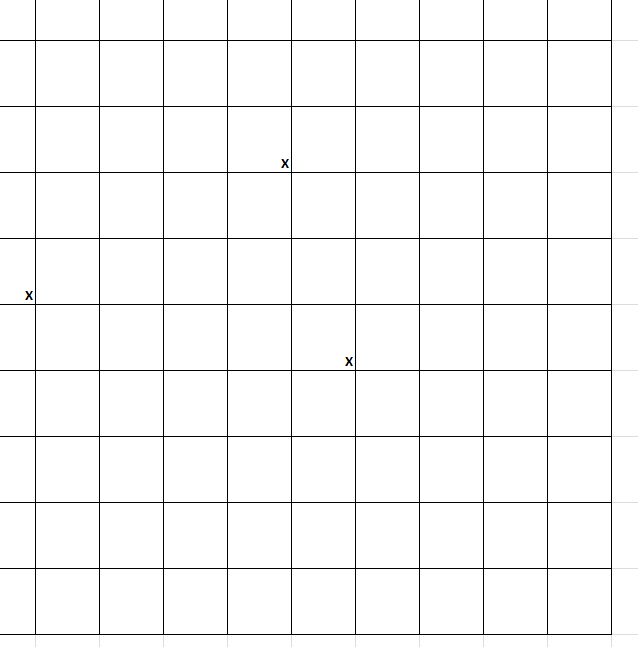}{Locations of facilities found by genetic algorithm (cost 94.7).}{width=0.7\textwidth}

\chapter{Algorithms and Complexity}
\label{chp:algorithmcomplexity}

Appendices~\ref{chp:basicoptimization} and~\ref{chp:fancyoptimization} introduced the concepts of formulating optimization problem,  provided examples of different types and categories of formulations, and presented algorithms for solving them.
This appendix focuses more on what specifically an algorithm is, and how complexity analysis can be used to compare the performance of alternative algorithms for solving a problem.

\section{Algorithms}
\label{sec:algorithms}

\index{algorithm|(}
The word \emph{algorithm} appears often in applied mathematics and computer science to denote a concrete procedure for accomplishing a task, specified in enough detail that it can be implemented by a computer (often in a programming language like C or Python).
We adopt~\citeapos{knuth1} criteria for operationalizing this definition: an algorithm must be \emph{finite} (always stopping in a finite number of steps), \emph{definite} (each step is unambiguous), \emph{effective} (each operation is simple enough to be done in a finite amount of time), and map \emph{inputs} to at least one \emph{output}.
An algorithm can have any finite number of inputs, even zero: an algorithm with no inputs would compute the same thing every time it is run, but that still might be useful (e.g., to answer a difficult question).
However, it must have at least one output; an algorithm must report \emph{something} when it terminates.

For example, a step such as ``add integers $x$ and $y$'' is definite and effective.
A step such as ``multiply $\pi$ by 2'' is definite, but not effective since $\pi$ is an irrational number with an infinite number of digits (there is no way to multiply \emph{all} the digits by 2).   
Such a step can be made effective by specifying ``round $\pi$ to 15 digits and multiply it by 2.''
A step such as ``pick an integer between 1 and 10'' is effective, but it is not definite (\emph{how} exactly should I pick such an integer?)
One way to make it definite is to specify ``randomly pick an integer between 1 and 10, with each possibility equally likely.''

Another example of a ``definite but not effective'' step is ``add 1 to $x$ if there are one million consecutive zeros in the decimal expansion of $\pi$.''
It is definite: either there are that many consecutive zeros or there are not.
However, we don't already know whether that is true, and there is no simple way to find out.
We could spend a lot of time computing digits of $\pi$ and maybe we will find this many zeros; but maybe we won't, and we don't know if it's because there aren't any, or if we just haven't computed far enough.
With our current knowledge of mathematics, there is no way to answer that question in a finite amount of time, so this step is not effective.

Many algorithms contain three main components:
\begin{description}
\item [Initialization:] These are the first steps performed.
The initialization steps prepare for the main algorithm, assembling the inputs, and providing initial or default values to variables which will be used later on.
\item [Body:] The steps in the body form the main core of the algorithm.
Many algorithms are \emph{iterative}, and the steps in the body are repeated multiple times.
\item [Stopping criterion:] An algorithm must be finite, and not run forever.
However, we often don't know how many iterations are required in advance.
Therefore, we periodically check a stopping criterion; once it is satisfied, we can stop.
There may also be a few final processing and ``clean up'' steps needed at the end of the algorithm.
\end{description}

You were probably taught an algorithm for adding two positive integers when you were a child (although that word was probably not used).
One way to do this is to work from right to left, adding the digits one at a time.
If the sum of the digits is 9 or less, then that is the corresponding digit in the sum; if it is 10 or more, subtract 10 to get the digit in the sum, and ``carry'' by adding 1 to the \emph{next} digit to the left.
For example, to add $158$ and $371$, you would start by adding $8 + 1 = 9$ for the rightmost digit; $5 + 7 = 12$, so enter $2$ for the next digit and carry a 1 for the next digit: $1 + 3 + 1 = 5$, so the sum is $529$.

This can be specified as an algorithm in the following way.
For simplicity, we assume that the two integers we are adding ($a$ and $b$) have the same number of digits $n$, and they are specified as $a = a_n a_{n-1} \ldots a_2 a_1$ and $b = b_n b_{n-1} \ldots b_2 b_1$.  
(In the previous example, $a = 158$, and $a_3 = 1$, $a_2 = 5$, and $a_3 = 8$; and likewise for $b = 371$.)
We also assume that we know how to add two single digits, perhaps with an ``addition table'' telling us the value of adding $0 + 0$, $0 + 1$, and so on, up through $9 + 9$.
The steps below describe the addition algorithm, calculating the digits in the sum $s$.

Note the use of the $\leftarrow$ operator to mean \emph{assignment}: calculate the expression on the right-hand side, and store the result in the variable on the left-hand side.
For instance, $c \leftarrow 0$ sets the variable \ttt{c} to zero.
Many programming languages use the $=$ sign to denote this (you would type \ttt{c = 0}).
In this text, we will reserve $=$ to refer to \emph{mathematical equality}, the statement that the left-hand side and right-hand side are currently equal to each other.
So a statement like $x = x + 1$ can never be true (no number is equal to one plus itself).
However, a statement like $x \leftarrow x + 1$ is perfectly fine; it means add one to the current value in $x$, and then store that new value in $x$.
Some programming languages use context to determine whether $=$ is meant in the sense of assignment or equality, while others use different operators for the two (for instance, C and Python use \ttt{=} for assignment, and \ttt{==} for equality).
To avoid confusion, we will use the different characters $\leftarrow$ for assignment, and $=$ for equality.

\begin{enumerate}
\item Initialization: set the carry digit $c \leftarrow 0$, and the current digit $i \leftarrow 1$.
\item Compute $s_i \leftarrow c + a_i + b_i$.
\item If $s_i \geq 10$, set $c \leftarrow 1$ and $s_i \leftarrow s_i - 10$; otherwise, set $c \leftarrow 0$.
\item If $i < n$, set $i \leftarrow i + 1$ and return to step 2.
\item If $c = 1$, set $s_{n+1} \leftarrow 1$.
\item Return the integer $s$ and terminate.
\end{enumerate}

In this algorithm, step 1 is the initialization, and steps 2--3 are the body.
Step 4 contains the stopping criterion, expressed in a ``negative'' way (``repeat as long as $i < n$'' rather than ``stop repeating if $i \geq n$.'')
Step 5 contains a final step before we terminate (if there is a carry digit after adding the $n$-th place, we need to add a 1 at the left-hand side of the sum, which will have one more digit than the original numbers.)

In general, when reading an algorithm, \emph{you should work out a small example on a separate piece of paper to see how the steps work.}  
Trying to mentally think through all the steps adds unnecessary difficulty, especially for complex algorithms.
It is much easier to see \emph{how} it works on a small example or two, and then study the specific steps to see how they fit together.\footnote{This is similar to studying a play or movie.
The best way to do this is to first \emph{watch} a performance, and then turn to the script to study it in more detail.
Starting with the script is much harder, and not really how a movie is meant to be experienced.
So too with algorithms.
}
So, let's see how this algorithm works for the numbers $a = 158$ and $b = 371$. 
These are three-digit numbers, so $n = 3$.
After the initialization (step 1), $c = 0$ and $i = 1$, meaning we are working with the first (right-most) digit and there is nothing to carry so far.
Step 2 adds the two right-most digits: $s_1 \leftarrow c + a_1 + b_1 = 9$, giving the right-most digit in the sum.  
This is less than 10, so the carry digit is set to zero (it is already zero, so nothing changes).
The current digit is $i = 1$, which is less than 3; so in step 4 we increase it and return to step 2 with $i = 2$.
Repeating with the second digit, we compute $s_2 \leftarrow c + a_2 + b_2 = 12$.
This is greater than 10, so in step 3, set the carry digit $c$ to 1, and reduce the sum digit by 10, so $s_2 = 2$, the second digit from the right in the sum.
In step 4, the current digit is still less than the number of digits, so we increase $i$ and return to step 2 with $i = 3$.
With the third digit, we compute $s_3 \leftarrow c + a_3 + b_3 = 5$; since this is less than 10, we set the carry digit to zero.
Now $i = n$, so we meet the stopping criterion and move to step 5.
Since the carry digit is zero, we are done and return the sum $s = 529$ (collecting the individual digits computed above).
Step 5 would have handled the case when there is a carry digit in the end, meaning the sum has more digits than the addends; for instance, the sum of $999$ and $111$ is $1110$, and this step would have added the `1' at the very left.

This algorithm satisfies all the criteria.
It is finite, since we repeat the body once for every digit, and any integer has a finite number of digits.
It is definite, since every step is spelled out precisely without any room for ambiguity.
It is effective, since every step can be done in a finite amount of time.
It takes two inputs (the numbers $a$ and $b$), and produces one output (the sum $s$).

An analogy is often made between algorithms, and cooking recipes.
For example, cooking pasta can be expressed by the following recipe:
\begin{enumerate}
\item Bring a large pot of water to a rolling boil.
\item Add the pasta, and stir until the water returns to a boil.
\item Wait 9--11 minutes (or until pasta has desired consistency).
\item Drain pasta in a colander.
\end{enumerate}
Which of the criteria of an algorithm are satisfied by a recipe like this?
Which of them are violated, and are they violated in a ``weak'' way (where you could satisfy our definition by adding some more explanation to the step) or a ``strong'' way (a fundamental violation that can't be easily fixed)?

\section{Pseudocode}
\label{sec:pseudocode}

\index{pseudocode|(}
The previous section gave the formal definition of an algorithm.
Sometimes, specifying an algorithm formally --- either writing out all the steps in full detail in English (as in the addition example above) or with source code in a programming language --- may be unnecessarily verbose or confusing.
\emph{Pseudocode} is a less formal way of specifying the steps of an algorithm.
Pseudocode can be thought of as a midway point between intuitive English descriptions (which express the general idea, but are usually ambiguous) and computer code (unambiguous, but harder to read and understand).
The goal of pseudocode is to express the main ideas of an algorithm while avoiding unnecessary pedantry.
It communicates the algorithm in enough detail that a programmer can write actual code in a specific language, but without the ``fiddly bits'' inherent in specific languages (declaring variables, calling library functions, details of data structure implementations, and so forth).
Algorithm~\ref{alg:max} is an example of pseudocode, and illustrates some common conventions.

\begin{algorithm}[h]
\SetAlgoLined
\textbf{Input: } $A = \{ A[1], A[2], \ldots, A[n] \}$: Set, list, or array containing the $n$ numbers \\
\textbf{Initialize:} \\
$max\_value \leftarrow -\infty$ \\
\For{$i\leftarrow 1$ \KwTo $n$}{
    \If {$A[i] > max\_value$}{
        $ max\_value \leftarrow A[i]$
    }
}
\caption{Finding the maximum of $n$ numbers. \label{alg:max}}
\end{algorithm}

This pseudocode is for an algorithm that finds the maximum of $n$ numbers.
As an example, let us find the maximum of five numbers: $7$, $3$, $5$, $9$, and $1$.
The five numbers are assigned to the list or set or array $A$, so that $A[1] = 7, A[2] = 3, \ldots, A[5] = 1$.\footnote{Note that some programming languages index arrays from 0, rather than 1; using pseudocode lets us abstract from this.
The assumption is that a programmer working in such a language would make the necessary changes.}
The variable $max\_value$ is initialized to a large negative number.

The \emph{for loop} is a control structure which repeatedly executes a set of instructions a fixed number of times, as specified by an iterator variable.
In this case, $i$ is the iterator variable which takes values from 1 to $n$.
The variable $i$ takes the values 1, 2, 3, 4, and 5 in turn; for each of these values, the set of instructions indented under the for loop is executed.
Inside the for loop is an \emph{if} block.
This is another control structure which executes the indented instructions when (and only when) a certain condition is satisfied.

Each \emph{iteration} corresponds to one time that the instructions in the for loop are executed.
We now illustrate a few iterations.
At first, the iterator $i$ takes the value 1.
The if statement then checks whether $A[1] = 7 > max\_value = -\infty$.
Since $7 > -\infty$, the instructions under the if statement are executed, and $max\_value$ now takes the value 7.
Next, the iterator variable $i$ takes the value 2.
The if statement checks whether $A[1] = 3 > max\_value = 7$.
This is not true, so the instructions under the if statement are not executed, and $max\_value$ remains at 7.
After five iterations, we will find that $max\_value$ is equal to 9, which is indeed the maximum of the given five numbers.
The stopping criterion for this algorithm is when the for loop has been executed for all the numbers in the array.

\index{Newton's method!algorithm|(}
As another example, the pseudocode in Algorithm~\ref{alg:new} represents the algorithm to find the root of a differentiable function $f(x)$ using Newton's method.
This pseudocode indicates some conventions commonly found in optimization algorithms.
Like the other algorithms seen so far, Newton's method repeats a ``body'' of steps iteratively until the stopping criterion is met.
The variable $k$ keeps track of the number of iterations, and $x_k$ is the solution at iteration $k$.

\begin{algorithm}[H]
\SetAlgoLined
\textbf{Initialize:} \\
$x_0$: initial guess\\
$\epsilon$: tolerance \\
$x_1 \leftarrow x_0 - \frac{f(x_0)}{f'(x_0)}$\\
$k \leftarrow 1$ \\
\While{$|x_{k} - x_{k-1}| > \epsilon$}{
    $x_{k+1} \leftarrow x_k - \frac{f(x_k)}{f'(x_k)}$ \\
    $k \leftarrow k+1$
}
\caption{Finding the root of a function $f(x)$. \label{alg:new}}
\end{algorithm}

In this implementation, the \emph{while loop} control structure is used to control the number of repetitions.
In the while loop, a set of instructions is repeated as long a condition is met.
The condition in the while loop corresponds to the stopping criterion, expressed in a ``negative'' way (repeat as long as the stopping criterion is false).
In Newton's method, the algorithm stops when the change in solution across consecutive iterations $|x_{k} - x_{k-1}|$ is less than a pre-specified tolerance level $\epsilon$.
The tolerance level is chosen depending on the level of accuracy desired.

We illustrate a few iterations of this algorithm for the function $f(x) = x^2 - 3$, with the initial guess $x_0 = 1$.
The tolerance $\epsilon$ is set at 0.001, which represents the desired level of accuracy in the output.
$x_1$ is now calculated using the formula to be equal to 2 and iteration counter $k$ is set to 1.
The condition in the while loop $|x_1 - x_0|$ is equal to 1 which is greater than $epsilon$.
Therefore, we now calculate $x_3$ to be equal to 1.75 and increase the value of $k$ to 2.
This process is repeated until the difference between $x_k$ and $x_{k-1}$ is less than 0.001.
The stopping criterion we use ($|x_{k} - x_{k-1}| \leq \epsilon$) is one of several widely used stopping criteria.
Another popular criterion is to stop the iterations when the iteration counter hits a pre-specified maximum number $N$.

An attentive reader will notice that to be a valid algorithm, it must be finite, meaning that the convergence criterion must eventually be true.
If the stopping criterion is a maximum number of iterations $N$, this will clearly happen at some point.
For the stopping criterion $|x_{k} - x_{k-1}| \leq \epsilon$, it is not obvious that it will eventually be satisfied (and indeed, for some functions $f$ and initial guesses it may never be).
\emph{Analysis of algorithms}\index{algorithm!analysis} is the sub-field of computer science which studies algorithms to identify conditions for finiteness (will it eventually stop?), proofs of correctness (are the outputs always the ones we want?), and the time and computer memory they need before stopping.
For the most part, this book will not dwell on these details except where relevant for applying transportation network models.

In this case, the advantage of the $|x_{k} - x_{k-1}| \leq \epsilon$ stopping criterion is that it is directly related to the precision in the solution; after a fixed number of iterations $N$ we may or may not be close to the actual root of the function; or we may have found the solution many iterations ago but kept on going unnecessarily until we hit $N$ iterations.
However, the $k = N$ stopping criterion has the advantage of guaranteeing that it will eventually be satisfied.
We can combine the two stopping criteria to gain the advantages of both, by changing the condition in the while loop, so it reads ``while $|x_{k} - x_{k-1}| > \epsilon$ and $k < N$.''
This lets us stop as soon as we have found the solution, and also guarantees the algorithm will not run forever.
\index{Newton's method!algorithm|)}
\index{pseudocode|)}

\section{Algorithmic Efficiency}

\index{algorithm!efficiency|(}
Many different algorithms can be devised to solve a particular problem.
For example, there are several methods available in the literature to sort a list of  numbers.
How can we choose the ``best'' one for a particular application?
Algorithm efficiency can be characterized in several ways: run time, memory use, and so on.
In this section, we focusing solely on characterizing efficiency in terms of run time; similar methods can be applied to characterize memory use and other resources.

One way to do that would be to implement each algorithm in a programming language, and measure how long each one takes.
A difficulty with this approach is that this kind of performance comparison depends on the operating system, programming language, and the hardware specification of the computer used to test the algorithms.
You would also need to decide what kinds of tests to run: how large should the list of numbers be, how many tests are needed, and so forth.
This section presents an alternative approach based on \emph{worst-case asymptotic analysis}, identifying the critical bottlenecks in the algorithm as the problem becomes harder and harder.
By ``worst-case,'' we mean that our intention is to find a conservative performance guarantee; for particular inputs the algorithm may perform better than our analysis shows, but it will never perform worse.
An alternative is to examine ``average-case'' performance under a typical range of inputs, but this is trickier because what is ``typical'' for one person's application of an algorithm may be completely different than another's.  
Worst-case performance lets us focus on the algorithm itself, and not one individual use case or implementation.

A small example illustrates what we mean by worst-case analysis.
Imagine you are the owner of a casino, and you have developed a new game in which a gambler wagers a certain amount of money.
A 6-sided die is then rolled.
If the value is 1 or 2, the gambler loses 50\% of the amount wagered; if the outcome is 3 or 4, the entire wager is lost; and if the outcome is 5 or 6, the gambler wins an amount of money equal to the wager (so their initial stake is doubled).
The process of determining the amount won can be represented by Algorithm~\ref{alg:gamble}.

\begin{algorithm}[h]
\SetAlgoLined
\textbf{Initialize:} \\
Initialize $money\_bet$ to the wager: \\
Obtain $outcome$ by rolling a 6-sided die. \\
\If {outcome=1  or outcome=2}{
    $money\_returned \leftarrow money\_bet / 2$
}
\ElseIf {outcome=3  or outcome=4}{
    $money\_returned \leftarrow money\_bet$
}
\ElseIf {outcome=5  or outcome=6}{
    $money\_returned \leftarrow money\_bet \times 2$
}
\Return $money\_returned$
\caption{Gambling algorithm for demonstrating worst-case analysis.}
\label{alg:gamble}
\end{algorithm}

We want to measure how long it will take for this algorithm to run.
We can either measure this in terms of ``clock time'' (imagine using a stopwatch to time how long it takes from start to finish), or by counting the number of operations that need to be performed.
Algorithmic analysis usually uses the latter approach, because it keeps the focus on the algorithm rather than a specific hardware.
The exact same algorithm will run much faster on a supercomputer than on your phone.
As computers continue to advance, it is more useful to say ``this algorithm requires X steps'' than ``it takes Y seconds to run.''

This requires defining what we mean by an ``operation.''
Some computations take much more time than others; computing the square root of a number, or the cosine of a number, takes significantly longer than changing the sign of a number from positive to negative, for instance.
We will define a basic operation to be (i) a simple arithmetic computation (add, subtract, multiply, divide); (ii) testing whether one value is equal to, greater than, or less than another; and (iii) storing a value in a variable using the assignment operator $\leftarrow$.\footnote{
By ``value'' in this section, we mean either an integer or a floating-point number rounded to a finite number of digits.}

In Algorithm~\ref{alg:gamble}, there are two basic operations associated with initialization.
Depending on the value of $outcome$, it takes a different number of basic operations to determine which if block is executed.
If $outcome = 1$, we only require one comparison to see that the money returned should be half the bet (we don't have to check whether $outcome = 2$, because the first if statement has an ``or'' in it; as soon as we know $outcome = 1$ we can jump to the statement in the block.)
If $outcome = 2$, we need two comparisons: the first checks whether outcome is 1 (false), the second checks whether outcome is 2 (true), and only then we execute the statement.
If $outcome = 3$, we need three comparisons; two failed comparisons testing whether outcome is 1 or 2, and a third successful one.
In the worst case, $outcome = 6$, and we need six comparisons to determine which block of code to execute.
Within the if blocks, there are either zero arithmetic operations (if the outcome is 3 or 4) or one arithmetic operation (all other outcomes) multiplying or dividing the wager by 2.
There is always one assignment operation in each block.

Putting this together, how many basic operations does the algorithm need?
In the best case, $outcome = 1$, and a total of five basic operations are needed (two initialization, one comparison, one arithmetic, and one assignment).
In the worst case, $outcome = 6$, and ten operations are needed (two initialization, six comparison, one arithmetic, one assignment).
The average case can be treated by representing the operation count as a random variable.
Its expected value can be calculated by identifying all possible operation counts under every outcome, multiplying by the probability of that outcome, and adding.
The resulting number is 7.2, representing the long-run average of the number of operations needed.

Which one of these is most useful?
The average case might seem to be best, but its calculation was more involved, and it would have to be recomputed if a different kind of die were used where all six outcomes were not equally likely.
By contrast, the best case and worst case were much simpler to compute, and they would stay the same even if the die were to be weighted more to some outcomes than others.
Between these two, the worst case is more useful for planning purposes.
The best case is optimistic to a fault; all it tells you is that you'll need at least five basic operations (but nothing about how bad it could get).
This is akin to a shady investment pitch which only tells you how much you could earn, and nothing about how much you stand to lose.
The worst case, by contrast, is a more solid basis for planning and allocating resources.
You will never need more than ten operations (and often will need less), so if you plan for this you will never exceed your resource budget.

In this case, the performance of the algorithm requires a constant number of steps.
What if it is repeated $n$ times (either the same gambler playing the game multiple times, or several gamblers playing)?
The worst case is that ten operations are needed each time; in this case $10n$ basic operations are needed to run the game.
Again, this is a conservative estimate.
The only time $10n$ operations are needed is if a 6 is rolled all $n$ times the game is played, which is unlikely.
However, we can absolutely guarantee that you will never need more than this.

As another example, return to Algorithm~\ref{alg:max} for finding the largest of $n$ numbers.
In the for loop, for each value of $i$ we perform at most two basic operations: one comparison operation (testing whether $A[i] > max\_value$) and one assignment operation ($ max\_value \leftarrow A[i]$).
Since $i$ takes values from 1 to $n$, these two operations will be repeated $n$ times, giving a total of $2n$ operations.
There is also one basic operation associated with initialization, for a total of $2n + 1$ basic operations in the worst case.
(The best case has fewer than this; for example, if the $A[1]$ is the largest of the $n$ numbers, then only $n + 2$ operations are needed, because $max\_value$ is only assigned once.)

To further simplify the worst-case analysis, we introduce an \emph{asymptotic notation}\index{asymptotic notation} which highlights the bottleneck as the problem grows harder.
For finding the maximum of $n$ numbers, the one operation needed for initialization becomes insignificant if $n$ is large.
If $n$ is equal to a million, say, what really matters are the (up to) two million operations needed in the body of the algorithm, not the one needed for initialization.
So we would like to focus on the $2n$ part and not the $+1$ part.
Even here, we can simplify further.
The most important part is that it grows with $n$ (and not $n^2$, say): doubling $n$ will also double the number of steps needed.
The factor of $2$ is less important.

To see why, imagine that someone proposes an alternative algorithm for finding the maximum of $n$ numbers, and worst case analysis shows that it may require up to $n^2 + 4n + 5$ operations.
As $n$ grows large, the $n^2$ part will dominate the remaining terms.
Furthermore, if $n$ were to double, $n^2$ would increase by a factor of four.
This is evidence that Algorithm~\ref{alg:max} is a better choice.
As $n$ gets larger and larger, the number of operations needed by Algorithm~\ref{alg:max} increases in direct proportion, whereas that needed by the hypothetical alternative grows with the square of $n$.
The fact that the first algorithm is $2n$, as opposed to $n$ is \emph{much less important} than that it is linear in $n$ and not quadratic.

Another reason to focus on the $n$ and not the $2$ is to abstract certain implementation details.  
In Algorithm~\ref{alg:max}, do we count an initialization step associated with the first line, setting the values in the array $A$?
Maybe it is passed directly as an address in computer memory, to values that already exist; in this case there is no new operation.
On the other hand, if we have to assign them from other values in memory, we need an additional $n$ assignment operations, for a total of $3n + 1$.
Or maybe we have to read them from a file, requiring $n$ operations to read them from the drive and $n$ assignment operations to store them in $A$, for a total of $4n + 1$.
In all of these cases, the most important part is that the computational requirements scale proportionally to $n$, since this will tell you how the run time will change if you give it a larger or smaller problem instance.
Whether the coefficient is 2, 3, or 4 is of secondary importance.

We will introduce an asymptotic notation that lets us simplify in this way; we want a concise way to express that as $n$ grows large, the number of steps required by Algorithm~\ref{alg:max} grows in direct proportion to $n$, without fiddling around with the constants in the full expression $2n + 1$.
The notation will say that the number of steps required by Algorithm~\ref{alg:max} is $O(n)$ (read ``big-oh of $n$'')\index{big-$O$}.\label{not:bigO}
Technically, we say that a function $f(n)$ is $O(g(n))$ if there are constants $A$ and $n_0$ such that $|f(n)| \leq A g(n)$ whenever $n > n_0$.  
Practically, we can identify such a function $g(n)$ by selecting the dominant term in $f(n)$ as $n \rightarrow \infty$, and removing any constant coefficient it has.
For example, $n^2 + 4n + 5$ is $O(n^2)$, and $4 \cdot 2^n + 6 n^3 - 8 \log n$ is $O(2^n)$.
The initial version of the gambling algorithm (Algorithm~\ref{alg:gamble}) required at most 10 steps; we can say this is $O(1)$.

Some common asymptotic notations that arise are:
\begin{itemize}
\item $O(1)$, for algorithms that run in \emph{constant time} (independent of input size);
\item $O(\log n)$, sometimes called \emph{logarithmic time} or \emph{sublinear time} (the run time increases more slowly than the problem size);
\item $O(n)$, or \emph{linear time} (run time scales in direct proportion to the problem size);
\item $O(n^2)$, or \emph{quadratic time};
\item $O(n^3)$, $O(n^4)$, and other \emph{polynomial time} algorithms of the form $O(n^c)$ for some constant $c$.
\item $O(2^n)$, or \emph{exponential time};
\item $O(n!)$, or \emph{factorial time}.
\end{itemize}
These are listed in increasing order of difficulty.
Generally, algorithms up to polynomial time are considered to be reasonably efficient.
Algorithms with exponential time complexity are difficult to solve for large instances in reasonable amount of time.

Some people think that this classification is not very useful, because computers continue to get faster over time.
However, the difference between complexity classes is absolutely critical, and will remain so for the foreseeable future.
Imagine you have two algorithms for a network optimization problem with $n$ nodes; one runs in quadratic time, $O(n^2)$, and the other runs in exponential time, $O(2^n)$.
Also imagine that using either algorithm, on your current computer you can just barely solve a problem with 1000 nodes in an acceptable amount of time.
After some time, you upgrade your computer so that it is twice as powerful as the one you had before.
Roughly speaking, the quadratic-time algorithm on this better machine can handle a network which is larger by a factor of $\sqrt{2}$, so roughly 1400 nodes. 
The exponential-time algorithm, on the other hand, can only handle a network with one more node, since $2^{1001}$ is twice as large as $2^{1000}$!
Faster hardware will not overcome the disadvantages of an exponential-time algorithm, at least for large problems.

However, just because an algorithm is exponential doesn't mean it is useless.
For one, the $O(\cdot)$ notation is asymptotic, and tells you what will happen as $n$ tends to infinity.
An exponential time algorithm may be perfectly fine working with the problem sizes you are dealing with; it just means that it does not scale well, and for large enough problems you will run into difficulties.
It is also worst case, and there are algorithms which may have exponential worst case complexity, yet behave more like polynomial time in practice.
The simplex method for solving linear programs (Section~\ref{sec:simplex}), and the branch and bound method for solving mixed-integer linear programs (Section~\ref{sec:branchandbound}) are examples of algorithms which are exponential in the worst case, yet offer better performance in practice, especially when implemented well and tailored to a specific application.

For one more example of complexity calculation, Algorithm~\ref{alg:bub} presents the ``bubble sort'' algorithm for sorting an array of $n$ numbers into increasing order.
The algorithm has two loops: an outer loop using the iterator $i$, and an inner loop with the iterator $j$.
When $i=1$, there are $4(n-1)$ operations associated with the $j$ loop.
When $i=2$, there are $4(n-2)$ operations, and so on.
Therefore in total there are $4(n-1) + 4(n-2) + \ldots + 1 = 2n(n-1)$ operations.
The bubble sort is therefore $O(n^2)$ in complexity, and has quadratic time.
This implies that if the input size doubles, the running time quadruples.
There are faster sorting methods with $O(n \log n)$ complexity; if you have to sort a large array of numbers this will perform much better than bubble sort.
However, bubble sort is much easier to implement, and may perform acceptably if you only have to sort a small array.
\index{algorithm|)}

\begin{algorithm}[h]
\SetAlgoLined
\textbf{Initialize:} \\
$A$: Set, list, or array containing the n numbers \\
\For{$i\leftarrow 1$ \KwTo $n-1$}{
    \For{$j\leftarrow 1$ \KwTo $n-i$}{
        \If {$A[j] > A[j+1]$}{
            $ temp \leftarrow A[j]$ \\
            $ A[j] \leftarrow A[j+1]$ \\
            $ A[j+1] \leftarrow temp$ \\
        }
    }
}
\caption{Bubble sort}
\label{alg:bub}
\end{algorithm}

\section{Complexity Classes}

\index{complexity class|(}
The previous section introduced asymptotic notation, and how it is used to characterize the worst-case run time of an algorithm, measured by counting the number of basic operations needed.
This is useful for comparing two algorithms for the same problem.
As the input size of a problem becomes larger and larger, any $O(1)$ algorithm will eventually run faster than any $O(n)$ algorithm, which will run faster than any $O(n^2)$ algorithm, which will run faster than any $O(2^n)$ algorithm, and so forth.
But can we say anything about how hard a \emph{problem} is?
In other words, I might prefer a $O(n^2)$ algorithm to an $O(2^n)$ one. 
But is $O(n^2)$ the best I can do?
Or with some more ingenuity and work, could we come up with something even better?

Complexity classes aim to understand hardness in terms of the underlying problems.
As explained in more detail below, one such class of problems is called $P$\index{P (complexity class)}; these are problems for which a polynomial-time algorithm exists to solve it.
A larger class of problems is called $NP$\index{NP (complexity class)}; for some of these problems, we do not know of any polynomial-time algorithm that will always solve them.
Some $NP$ problems are called $NP$-complete\index{NP-complete}, which means that any polynomial-time algorithm solving them would \emph{also} solve any other problem in $NP$.
Such an algorithm would be extremely valuable, and decades of work has gone into trying to find one.
This search has thus far proved fruitless.

This is the core of the famous $P \stackrel{?}{=} NP$ question, arguably the most important open problem in theoretical computer science, and one carrying a \$1 million bounty from the Clay Mathematics Institute.
Answering this question either way (providing such an algorithm, or definitively proving that none exists) would be a tremendous accomplishment.
Indeed, even though this core question has as yet gone unanswered, the search for its answer has yielded many other important findings in computer science.
\cite{aaronson16} provides an excellent survey of the history of this question, attempts to answer it, and the value this search has provided

Many students learning about complexity classes initially find them an esoteric theoretical conceit.
However, there is great practical significance!
If your supervisor asks you to solve a particular problem, they may not be happy if the best you can come up with is an exponential time algorithm, since it will scale badly for large-scale problems.
However, if you can show this problem is $NP$-complete, then you can point out that it's not just you who can't find a better algorithm --- the brightest minds in computer science have also been unable to find a better algorithm for this problem!
In other words, it's not you, it's the problem that's inherently difficult.
If you can suggest a simplification of the problem that brings it into $P$, then this directly suggests a faster way to solve it.
\index{algorithm!efficiency|)}

The next section provides a brief overview of complexity classes, emphasizing applications to transportation network optimization problems.

\subsection{Decision problems vs.\ optimization}
\label{sec:decision}

In computer science, complexity classes are usually defined in terms of \emph{decision problems}\index{decision problem}, where the answer is either ``yes'' or ``no.''
This book primarily discusses \emph{optimization problems}\index{optimization!and decision problems}, where the answer involves a solution in terms of decision variables and an objective function.
However, optimization problems can be transformed into decision problems.
Here are a few examples of how to do this:

\index{shortest path}
\paragraph{Shortest path problem:} 
Consider a directed network $G=(N,A)$, where $N$ is the set of nodes, and $A$ the set of directed links, with non-negative integer costs $c_{ij}$ for each $(i,j) \in A$.
Let $r$ and $s$ denote the origin and destination node.

\underline{Optimization problem:} 
Find the path $\pi$ between $r$ and $s$ minimizing its cost $c^\pi = \sum_{(i,j) \in \pi} c_{ij}$.

\underline{Decision problem:} For a given value $C$, is there a path $\pi$ between $r$ and $s$ whose cost $c^\pi$ does not exceed $C$?

\index{traveling salesperson problem}
\paragraph{Traveling salesperson:} 
Consider a complete network $G=(N,A)$, where $N$ is the set of nodes and $A$ the set of links, with integer costs $c_{ij}$.

\underline{Optimization problem:} 
Determine a minimum cost tour which visits each node exactly once.  (In other words, find a permutation of nodes $( p_1, p_2,\ldots, p_n )$ which minimizes $\sum_{i=1}^{n-1}c_{p_ip_{i+1}} + c_{p_np_1}$.)

\underline{Decision problem:} For a given value $C$, is there a tour visiting each node exactly once whose cost does not exceed $C$?  (In other words, is there a permutation of nodes  such that $\sum_{i=1}^{n-1}c_{p_ip_{i+1}} + c_{p_np_1} \leq C$?)

\index{knapsack problem}
\paragraph{Knapsack problem:} 
Consider a set of objects $I$ with positive integer sizes $a_i$ and profits $p_i$ for each $i \in I$, and a positive integer $U$ denoting the capacity.

\underline{Optimization problem:} Find a subset of objects $J_I \subset I$ maximizing total profit $\sum_{j \in J_I} p_j $ such that $\sum_{j \in J_I} a_j \leq U$.

\underline{Decision problem:} For a given profit level $P$, is there a set of objects $J_I \subset I$ such that $\sum_{j \in J_I} p_j \geq P$ and $\sum_{j \in J_I} a_j \leq U$.

\index{shortest path!resource constrained}
\paragraph{Resource-constrained shortest path problem:} 
Consider a directed network $G=(N,A)$, where $N$ is the set of nodes, and $A$ the set of directed links, with positive integer costs $c_{ij}$ and resources $r_{ij}$ for each $(i,j) \in A$.
Let $r$ and $s$ denote the origin and destination node, and $R$ the resource budget.

\underline{Optimization problem:} Find the path $\pi$ between origin node $s$ and destination node $t$ with minimum cost $c_\pi = \sum_{(i,j) \in \pi} c_{ij}$, subject to the resource usage constraint $\sum_{(i,j) \in \pi} r_{ij} \leq R$.

\underline{Decision version:} For a given value $C$, is there a path $\pi$ between $s$ and $t$ such that $c_\pi \leq C$ and $\sum_{(i,j) \in \pi} r_{ij} \leq R$?

\index{longest path}
\textbf{Longest acyclic path problem:} 
Consider a directed network $G=(N,A)$, where $N$ is the set of nodes, and $A$ the set of directed links, with positive integer costs $c_{ij}$ for each $(i,j) \in A$.
Let $r$ and $s$ denote the origin and destination node.

\underline{Optimization problem:} 
Find the acyclic path $\pi$ between $r$ and $s$ with maximal cost $c^\pi = \sum_{(i,j) \in \pi} c_{ij}$.

\underline{Decision problem:} For a given value $C$, is there an acyclic path $\pi$ between $r$ and $s$ whose cost $c^\pi$ is at least $C$?

From the standpoint of solution, a polynomial time algorithm for the optimization problem also leads to a polynomial time algorithm for the decision problem, and vice versa.
For example, if we have a polynomial time algorithm for the optimization version of the shortest path problem, we can use it to solve the decision problem: just find the shortest path using the optimization version; if its cost is less than or equal to $C$, then the answer to the decision problem is ``yes.''
Otherwise it is ``no.''
Or, if we have a polynomial time algorithm for the decision version of the problem and want to find the exact cost of the shortest path (not just whether it is less than $C$ or not), we can repeatedly solve the decision version problem for different versions of $C$ using a bisection-like approach (cf.\ Section~\ref{sec:bisection}).
This requires solving the decision problem multiple times, but a polynomial number of times; substituting one polynomial into another still results in a polynomial, so the overall method is still polynomial time (even if the polynomial has a higher degree).

\subsection{The classes $P$, $NP$, $NP$-hard, and $NP$-complete}

This subsection introduces four classes of problems that arise often in optimization and elsewhere in computer science: $P$, $NP$, $NP$-hard, and $NP$-complete.
They are formally defined in terms of decision problems like the ones in the previous subsection, where the answer to a posed question is either ``yes'' or ``no.''
These four classes are not mutually exclusive, so a given problem may exist in more than one of these classes at the same time --- for example, the standard shortest path problem is in both $P$ and $NP$.
They are also not exhaustive; there are decision problems which do not lie in any of these four classes.

We start by defining the classes $P$ (``polynomial'') and $NP$ (``nondeterministic polynomial,'' and \emph{not} ``non-polynomial'' as is often guessed by people seeing the acronym for the first time).
The definitions given here are not fully rigorous, since we should also state exactly what our model of computation is, what we count as a basic operation, the exact format of the problem input, and so on.
Readers interested in these details should consult~\cite{sipser06} or another standard text in the theory of computation.
The definitions below will suffice if you are using common hardware, programming languages, and data formats.

\index{P (complexity class)}
\begin{dfn} A decision problem belongs to the class $P$ if, given an instance of size $n$, the answer can be obtained using an algorithm whose worst-case running time is polynomial, that is, $O(n^k)$ for some constant $k$.
\end{dfn}

\index{NP (complexity class)}
\begin{dfn} A decision problem belongs to the class $NP$ if, given an instance of size $n$ whose answer is ``yes,'' there is a polynomial-time algorithm to confirm that the answer is in fact ``yes.''
\end{dfn}

Worded intuitively, a problem in $P$ is easy to solve; a problem in $NP$ is easy to \emph{check}, in the sense that someone claiming the answer is ``yes'' can provide evidence that you can verify in polynomial time.  
A few examples are shown below; keep in mind that all of these refer to the \emph{decision} versions of these problems described in the previous subsection, not the optimization versions.

The decision version of the shortest path problem with non-negative link costs is in $P$, since Dijkstra's algorithm (and others) can find the shortest path in polynomial time.
So we simply compute the shortest path; if its cost is less than or equal to $C$, the answer to the decision problem is ``yes'' (there is a path of cost at most $C$), otherwise the answer is ``no'' (no such path exists).

The decision version of the traveling salesperson problem is in $NP$.
If the answer is ``yes'' (there is a tour visiting each node exactly once whose cost is less than or equal to $C$), someone can convince you of this by providing the tour itself as evidence.
You can check in polynomial time whether this tour is valid (it visits every node once and only once) and whether its total cost does not exceed $C$.
We have said nothing about how someone might find such a tour in the first place; simply that given such a tour as supposed evidence that the answer is ``yes,'' you can check this evidence in polynomial time.
Likewise, we have not said anything about ``no'' instances of the problem, in which all tours have a cost greater than $C$.
Indeed, there does not seem to be an easy way to demonstrate that a given instance of the traveling salesperson problem is ``no.''
Just providing a tour whose cost is more than $C$ does not prove that there isn't some other tour with a lower cost.
But this doesn't matter; the class $NP$ is only concerned with ``yes'' instances of decision problems.\footnote{The complexity class co-$NP$ is the analogous class for ``no'' instances.}

Similarly, the decision versions for the knapsack problem, longest path problem, and resource-constrained shortest path problems are in $NP$, since in any case where the answer is ``yes'' you can provide a feasible knapsack assignment or path as easily-checked evidence that the answer really is ``yes.''

Every decision problem in $P$ is also in $NP$.
Given a ``yes'' instance of the problem, I can check whether the answer is indeed ``yes'' by simply answering the question myself using a polynomial-time algorithm (which exists because the problem is in $P$).
So the decision version of the shortest path problem is in both $P$ and $NP$.

We do not know if the reverse is true, that every problem in $NP$ is also in $P$.
For instance, we know that the traveling salesperson problem is in $NP$, but we do not know if it is in $P$ or not.
We do not currently know of any polynomial time algorithm for solving it, but this does not mean that none exists.
Perhaps there is such an algorithm, but if so nobody has been clever enough to find it yet.
Most computer scientists believe that no such algorithm exists, but nobody has been able to prove that either.
Answering the question either way would be a tremendous advance in computer science.

The reason is that algorithms used to solve one kind of problem can often be used to solve other kinds of problems as well.
Many important problems (even problems with no obvious ``network'' in them) can be reduced to the traveling salesperson problem in the sense that an algorithm for the latter can also be used to solve the former, by constructing a network, choosing costs in the right way, and then translating the resulting tour back to the original problem.
So, finding a good way to solve the traveling salesperson problem would directly provide good ways to solve many other problems as well.
This idea is formalized in the following definition:

\begin{dfn}
A decision problem $A$ can be polynomially reduced to decision problem $B$ if (a) there is a procedure that can transform an instance of $A$ to an instance of $B$ in polynomial time; and (b) the answer for an instance of $A$ is ``yes'' if and only if the answer for the corresponding instance of ``B'' under this transformation is also ``yes.''
\end{dfn}

The implication is that if we know a polynomial time algorithm for $B$, then we can use it to solve $A$.
The number of steps needed to do this is the sum of the number of steps needed to construct an instance of problem $B$ from the given instance of $A$, and the number of steps needed to answer that instance of $B$.
If there is a polynomial reduction from $A$ to $B$, and if $B$ is solvable in polynomial time, then the total number of steps required is also polynomial.
Intuitively, the problem $A$ is no harder than $B$, since any method for $B$ can also be used to solve $A$.
(The problem $A$ could be easier than $B$, in the sense that a specialized method for $A$ might be a faster approach than translating the problem to $B$ and solving that; we just know it can't be harder.)

We give a few examples of these kinds of transformations for decision problems in the next subsection.
The main reason for introducing polynomial reductions is to define two additional complexity classes: $NP$-hard, and $NP$-complete.

\index{NP-hard}
\begin{dfn}
A decision problem $X$ belongs to the class $NP$-hard if \emph{every} problem in the class $NP$ has a polynomial reduction to $X$.
\end{dfn}

In other words, any problem in class $NP$-hard is at least as difficult to solve as any problem in $NP$, since any problem in class $NP$ can be converted to an $NP$-hard problem.
$NP$-hard problems are not necessarily in $NP$; they might be strictly more difficult in the sense that we might not be able to verify a ``yes'' answer to an $NP$-hard problem in polynomial time.

\index{NP-complete}
\begin{dfn}
A decision problem belongs to the class $NP$-complete if it is both in classes $NP$ and $NP$-hard.
\end{dfn}

The problems in $NP$-complete can be thought of as the most difficult problems in $NP$, because any other $NP$ problem can be converted to an $NP$-complete problem which will have the same answer.
The consequence of this is that if someone can develop a polynomial time algorithm for \emph{any} problem in $NP$-complete, then we can answer all problems in $NP$ in polynomial time as well; therefore the classes $P$ and $NP$ will be the same: $P = NP$.
On the other hand, proving that no polynomial time algorithm exists to answer an $NP$-complete problem would show that there are some problems in $NP$ which are not in $P$, so $P \neq NP$.
This makes the class $NP$-complete very important.
The decision version of the traveling salesperson problem is an example of an $NP$-complete problem.
The proof of this statement is difficult and beyond the scope of this book.

We summarize the discussion in this subsection in Figure~\ref{fig:hardness}, which illustrates the relationships between the four complexity classes discussed.
The top panel shows how the classes overlap if $P \neq NP$; the bottom panel shows what it would look like if $P = NP$.
Again, determining whether $P$ and $NP$ are in fact equal is perhaps the most important open question in computer science.
Most experts believe that $P \neq NP$ is the likelier answer to this question, but this is not a universally-held opinion, and in any case belief is not proof (there are plenty of true mathematical results whose discovery came as a surprise to expert mathematicians).

\stevefig{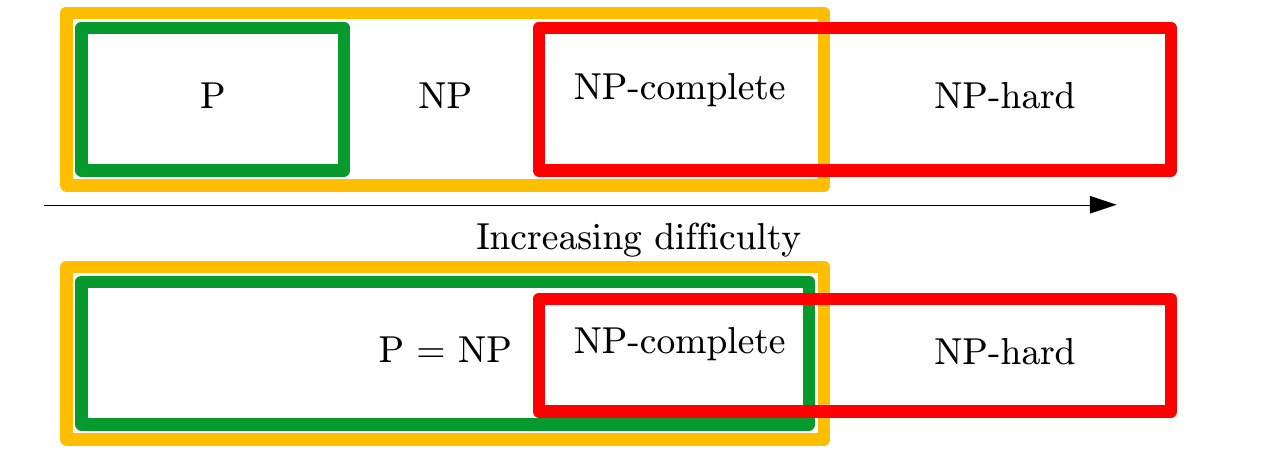}{P vs.\ NP vs.\ NP-complete vs.\ NP-hard under the assumptions $P \neq NP$ and $P = NP$.}{\linewidth}

As a final word, we comment on some differences between decision problems (the kind used to formally define the complexity classes in this section) and optimization problems (the kind arising more when doing transportation network analysis in practice).
For instance, the decision version of the traveling salesperson problem is in $NP$ (a given tour can be verified to be feasible and have cost less than $C$ in polynomial time), but the optimization version is not (there is no easy way to prove that a given tour is in fact minimal across all feasible tours).  
Sometimes, researchers in operations research use this terminology a bit loosely, and refer to a problem as being ``$NP$-complete'' (or similarly) without explicitly saying whether they are referring to the decision or optimization versions of these problems.
This is unfortunate, and if you are conducting research involving the complexity of transportation network problems we recommend that you become familiar with the distinction between the two.

\subsection{Examples demonstrating $NP$-completeness}

One of the primary uses of complexity classes is to establish that a given problem is \emph{intrinsically} hard.
In particular, if you can show a problem is $NP$-complete, then expecting to develop a polynomial time algorithm is more than simply an ambitious goal; it would be a groundbreaking advance in computer science.
For practical application, we are therefore content with heuristics that provide reasonably good solutions in a reasonable amount of time, rather than exact methods, since we do not know of any that scale well with problem size (at least in the worst case).
This subsection provides some examples of how you can show a decision problem is $NP$-complete, using the idea of polynomial reductions introduced above.

The typical method for showing a decision problem $X$ is $NP$-complete has two steps:
\begin{itemize}
\item Showing that $X$ is in NP.
That is, given an instance of $X$ where the answer is ``yes,'' it is possible to demonstrate this in polynomial time.
\item Providing a polynomial reduction from \emph{any other} $NP$-complete problem to $X$.  
This shows that $X$ is at least as difficult as any problem in $NP$ (and therefore $NP$-hard); together with the proof that $X$ is also in $NP$, this establishes $NP$-completeness.
\end{itemize}
The second of these is more challenging and can require more creativity than the first.
Notice that this method requires that we already know some problems are $NP$-complete.
Currently, hundreds of decision problems are known to be $NP$-complete, and you can use any of them as the starting point for the second step in a proof of $NP$-completeness.
(The very first proof of $NP$-completeness for a decision problem was much harder, since the above scheme could not be followed.  Instead, one had to give a general argument that any $NP$ problem could be reduced to that problem!)
\cite{garey79} is a good reference listing a variety of problems known to be $NP$-complete.

For example, the following decision problems are known to be $NP$-complete:
\begin{description}
\item[Hamiltonian cycle problem: ] Given a directed network $G = (N, A)$, is there a cycle which visits every node exactly once?
\item[Partition problem: ] Given a set $S$ of positive integers (possibly with multiple copies of each number), can the set be partitioned into two subsets $S_1$ and $S_2$ such that the sum of the integers in $S_1$ equals the sum of the integers in $S_2$?
\end{description}
Using these, we can prove that the traveling salesperson problem, knapsack problem, and resource-constrained shortest path decision problems introduced above are all $NP$-complete.
Let us look at the NP-completeness proof of the decision problems discussed above.

\begin{exm}
Show that the decision version of the traveling salesperson problem (given in Section~\ref{sec:decision}) is $NP$-complete.
\end{exm}
\solution{
The first step is to show that this problem is in $NP$; we have already shown in the previous section that we can check that a tour is valid, and that its total cost is at most $C$, in linear time, so we can verify a ``yes'' instance in polynomial (in fact, linear) time.

The second step involves a polynomial reduction from a known $NP$-complete problem.
In this case, we select the Hamiltonian cycle problem, and produce a reduction that would allow us to answer the Hamiltonian cycle problem by solving a suitable traveling salesperson problem.
(The order of the reduction is critical; showing that we could answer the traveling salesperson problem by solving a suitable Hamiltonian cycle problem would not be helpful.
That would show that the Hamiltonian cycle problem is at least as hard as the traveling salesperson problem; but we already know that the Hamiltonian cycle problem is $NP$-complete.
Our goal is to say something about the traveling salesperson problem, so we have to proceed in the other direction and show that we can use it to answer a Hamiltonian cycle question.)

So, assume we are given an instance of the Hamiltonian cycle problem with a network $G = N(A, A)$, and let $n = |N|$ denote the number of nodes in the network.
A Hamiltonian cycle visits every node exactly once before returning to the starting node, so such a cycle will contain $n$ links.
Construct a new network $G' = (N', A')$ in the following manner: the node sets are identical ($N' = N$), but $G'$ is a complete graph, meaning there is a link $(i,j)$ between \emph{every} pair of nodes.
The costs in $G'$ are chosen in the following way: $c_{ij} = 1$ if $(i,j) \in A$, and $c_{ij} = 2$ if $(i,j) \not\in A$.
That is, if the link exists in the original network, it has a cost of 1; if the link did not, it has a cost of 2.
This transformation can be accomplished in polynomial time, since the number of links in $G'$ is $n(n-1)$.

Now, consider the decision version of the traveling salesperson problem in $G'$, with $C = n$.
If the answer to this problem is ``yes,'' then there is a tour visiting every node and involving exactly $n$ links, all with a cost of 1.
All of these links therefore exist in the original network $G$, and correspond to a Hamiltonian cycle, so the answer to the Hamiltonian cycle problem is also ``yes.''

Or, if the answer to the traveling salesperson problem is ``no,'' then every tour in $G'$ has at least one link with a cost of 2.
No such links exist in the original network $G$, so there is no Hamiltonian cycle either: the answer to that problem is also ``no.''

This completes the reduction; we have already shown that it can be accomplished in polynomial time, and that the traveling salesperson problem is in $NP$.
Therefore it is $NP$-complete.
}

Essentially, we showed that any procedure which can solve the traveling salesperson problem could \emph{also} solve the Hamiltonian path problem, and therefore it cannot be any easier.

\begin{exm}
\label{exm:knapsack}
Show that the decision version of the knapsack problem (given in Section~\ref{sec:decision}) is $NP$-complete.
\index{knapsack problem!NP-completeness}
\end{exm}
\solution{
We start by showing that a ``yes'' instance can be verified in polynomial time to establish that the problem is in $NP$.
If the answer is ``yes,'' one can prove this by furnishing a subset $J_I$ whose total profit is at least $P$, and whose total size is no greater than $U$.
The total profit and total size can be calculated in linear time by simply summing the profit and size of each object in $J_I$.

The remainder of the argument establishes a polynomial reduction from an $NP$-complete problem.
For this example, we will reduce from the partition decision problem, with a given set of positive integers $S$.
From this, we will construct an instance of the knapsack decision problem whose answer (``yes'' or ``no'') will always be the same as that for the given partition problem.
Let $S = \myc{b_1, b_2, \ldots, b_n }$ (possibly with repetitions).
Create an instance of the knapsack problem where the set of objects is the set of integers $I = \myc{1, 2, \ldots, |S|}$.  
The size of each object, and its profit, are both set equal to the corresponding integer in $S$, so $a_i = p_i = b_i$.
The available capacity is half of the sum of the integers in $S$, so $U = \frac{1}{2} \sum_{i \in S} b_i$.
The minimum profit $P$ is set to this same amount as well: $P = \frac{1}{2} \sum_{i \in S} b_i$.
This reduction can be done in linear time.

Assume that the answer to the knapsack problem is ``yes.'' 
Then there is a subset of objects $J_I$ such that $\sum_{j \in J_I} a_j \leq U$ and $\sum_{j \in J_I} p_j \geq P$.
But $U = P = \frac{1}{2} \sum_{i \in S} b_i$, so this implies $\sum_{J_I} a_j = \frac{1}{2} \sum_{i \in S} b_i$.
Furthermore, each $a_j = b_j$, so $\sum_{j \in J_I} b_j = \frac{1}{2} \sum_{i \in S} b_i$.
That is, the integers represented by $J_I$ form exactly half of the total sum of the integers in $S$.
Therefore, the remaining integers not in $J_I$ must also form exactly half the total sum of the integers in $S$, so the sets $J_I$ and $I - J_I$ have equal sums, and the answer to the partition problem is ``yes.''

Now assume that the answer to the knapsack problem is ``no.''
Then every subset of objects $J_I$ either has $\sum_{j \in J_I} a_j > U$ (total size exceeds capacity) or $\sum_{j \in J_I} p_j < P$ (total profit falls short of target).
Since $U = P = \frac{1}{2} \sum_{i \in S} b_i$ and $a_i = p_i = b_i$, this means that no subset $J_I$ has the property that its integers add up to exactly half of the total of the integers in $S$.
This means that the answer to the partition problem must also be ``no,'' since if there were sets $S_1$ and $S_2$ whose sums were equal, they would both have to equal $\frac{1}{2} \sum_S b_i$, and $J_I = S_1$ would be a solution to the knapsack problem.

We can therefore answer the partition problem (known to be $NP$-complete) by solving the knapsack problem, so the knapsack problem must be at least as hard.
Since we have shown that the knapsack problem is in $NP$, it is therefore $NP$-complete as well.
}

\begin{exm}
\label{exm:rcsp}
Show that the decision version of the resource-constrained shortest path problem (given in Section~\ref{sec:decision}) is $NP$-complete.
\index{shortest path!resource-constrained!NP-completeness}
\end{exm}
\solution{
This problem is in $NP$, since a ``yes'' instance can be verified in linear time once an appropriate path is provided as evidence; simply check that the path is valid, and compute two sums over its links to verify that its total cost is at most $C$, and that its total resource consumption is at most $R$.

To show $NP$-completeness, we reduce the partition decision problem to the resource-constrained shortest path problem.
Given the set $S = \myc{b_1, b_2, \ldots, b_n}$, we construct a network $G = (N, A)$ with node and link sets defined as follows.
Let $U = \sum_{b_i \in S} b_i$, and $V = 2nU$.
We create $n + 1$ nodes; and between each consecutive pair of nodes $i$ and $i + 1$, we create two parallel links.
In each pair of parallel links, the top link has a cost of $V$, and a resource consumption of $0$.
The bottom link has a cost of $V - b_i$, and a resource consumption of $b_i$.
(See Figure~\ref{fig:rcsp}.)
Finally, for this instance of the resource-constrained shortest path decision problem, set the cost target $C = nV - (U/2)$, the resource limit to $U/2$, and select nodes 1 and $n + 1$ as the origin and destination.
The number of links and nodes is linear in the size of the original set $S$, so this transformation can be done in polynomial time.

\begin{figure}
\includegraphics[width=\linewidth]{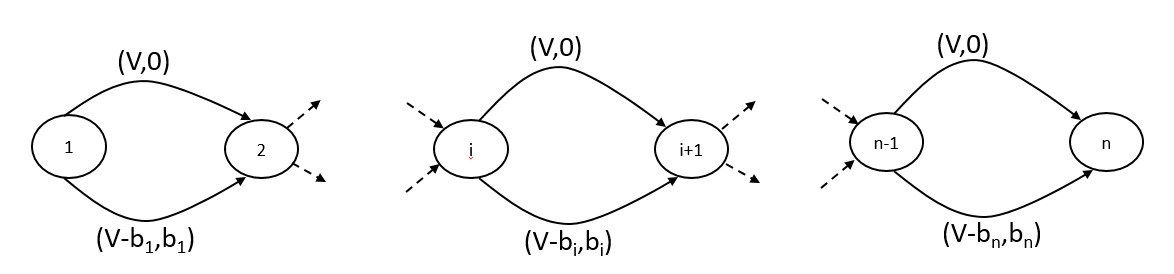}
\caption{Network transformation to reduce the partition problem to resource-constrained shortest paths.}
\label{fig:rcsp}
\end{figure}

Assume that the resulting resource-constrained shortest path problem is a ``yes'' instance.
Let $\pi$ be a path whose cost is at most $C$, and whose resource consumption is at most $R$.
This path consists of $n$ links, each of which is either the top or bottom link in the $n$ pairs of parallel links in the network.
Let $S_1$ contain the values of $S$ corresponding to ``bottom'' links, and $S_2$ the values corresponding to ``top'' links.
Then the total cost of this path is $nV - \sum_{b_i \in S_1} b_i$, and its total resource consumption is $\sum_{b_i \in S_1} b_i$.
Since this is a ``yes'' instance, the cost target is achieved and $nV - \sum_{b_i \in S_1} b_i \leq C = nV - (U/2)$, so $\sum_{b_i \in S_1} \geq U/2$.  
Likewise, the resource constraint is satisfied, so $\sum_{b_i \in S_1} b_i \leq R = U/2$.
Therefore $\sum_{b_i \in S_1}$ is exactly $U/2$, that is, the sum of the entries in set $S_1$ is exactly half the total of the sum of the entries in $S$.
We can thus conclude that the sum of the values in sets $S_1$ and $S_2$ are equal, and this is also a ``yes'' instance of the partition problem.

Or, if this is a ``no'' instance of the resource-constrained shortest path problem, then every path either exceeds the cost target or the resource limitation.
Following the same construction as in the previous paragraph, this means that for every set $S_1$ we have either $\sum_{b_i \in S_1} b_i < U/2$ or $\sum_{b_i \in S_1} b_i > U/2$.
Therefore, there is no subset of the values in $S$ whose sum is exactly half of the total sum of the values in $S$; so it is not possible to partition $S$ into two sets $S_1$ and $S_2$ whose sums are exactly equal.
Therefore this is a ``no'' instance of the partition problem as well.
}

You may have noticed that this proof was fairly similar to that used to show that the knapsack problem is $NP$-complete.
In fact, once Example~\ref{exm:knapsack} was done, and we established that the knapsack problem was $NP$-complete, we could have done Example~\ref{exm:rcsp} in a faster way by reducing the knapsack problem to the resource-constrained shortest path problem (rather than starting from the partition problem anew).  
Verifying this would be a very instructive exercise.

\begin{exm}
Show that the decision version of the longest acyclic path problem (given in Section~\ref{sec:decision}) is $NP$-complete.
\index{longest path!NP-completeness}
\end{exm}
\solution{
For a ``yes'' instance, verifying that a given path is feasible and has cost at least $C$ can be done in linear time, so this problem is in $NP$.

To show $NP$-completeness, we reduce the Hamiltonian cycle problem to the longest acyclic path problem.
The transformation proceeds as follows: let $G = (N, A)$ be the given network for the Hamiltonian cycle problem (where $n$ is the number of nodes).
Create a network $G'$ with $n + 1$ nodes in the following way: select any node $i$ in $N$, and split it into two nodes $i^+$ and $i^-$.  
Any link $(i,j)$ in the forward star of node $i$ appears in $G'$ as a link from $i^+$ to $j$.
Any link $(h,i)$ in the reverse star of node $i$ appears in $G'$ as a link from $h$ to $i^-$.
All other links in $G$ appear in $G'$ as they are, connecting the same nodes.
Select $i^+$ as the origin, and $i^-$ as the destination, and finally select the cost threshold $C = n$ to create an instance of the longest acyclic path problem in $G'$.

If this is a ``yes'' instance for the longest acyclic path problem, then there is a path $\pi$ whose cost is at least $n$.
Since the path is acyclic, it cannot revisit any node; and since there are $n + 1$ nodes in $G'$ the path $\pi$ must pass through all of them.
This acyclic path can be translated to a Hamiltonian cycle on $G$ by replacing $i^+$ and $i^-$ in $\pi$ by the single node $i$.

If this is a ``no'' instance for the longest acyclic path problem, then every path $\pi$ in $G'$ has cost less than $n$, meaning it does not pass through every node.
In this case, the original Hamiltonian cycle problem must also be a ``no'' instance, since any Hamiltonian cycle can be translated to a path in $G'$ from $i^+$ to $i^-$ passing through every node.
}

Two final comments on this subject.
First, pay attention to the problem specifics.
Minor changes can have a significant impact on the complexity of a problem.
For example, the Hamiltonian cycle problem (determining whether there is a path visiting each \emph{node} exactly once) is $NP$-complete; but the Eulerian cycle problem (is there a path visiting each \emph{link} exactly once) can actually be solved in linear time.
The shortest path problem (with positive link costs) can be solved in polynomial time; the shortest path problem with negative link costs and the longest (acyclic) path problem are both $NP$-complete.

Second, just because a problem is $NP$-hard or $NP$-complete does not mean the problems are impossible to solve, or completely impractical.
Such problems are often solved either using heuristics (which are inexact, but may be good enough), exact methods which are usually fast enough in practice (even if in the worst case the number of steps is exponential), or for small problem instances where an exponential time algorithm may still be acceptable.
As some concrete examples, integer optimization problems (Section~\ref{sec:integerprogramming}) are NP-hard.
But several algorithms have been developed that often solve real world instances efficiently.
Knapsack problems can be solved by pseudo-polynomial algorithms, a complexity class where the number of operations is bounded by the \emph{value} of the input, and not its absolute size.
These problems are also commonly solvable in practice.
\index{complexity class|)}

\chapter{Table of Notation}
\label{chp:notation}

This appendix lists all of the mathematical notation used in this book.
We have attempted to use distinct notation for all of the major concepts used repeatedly in the book, but sometimes the same symbol is re-used either due to overwhelming convention, or when we believe the mnemonic advantage of a particular notation outweighs any possible confusion.
Minor concepts appearing only once or twice ``collide'' more frequently with the major notation, but only in situations where we feel there is no ambiguity.
In such cases, the minor usage of the notation is \emph{italicized} in the table below.

The table is sorted alphabetically, with capital letters preceding lowercase, and Roman letters preceding Greek.
Assorted mathematical symbols are shown at the end.
The third column shows the page where the notation is first used.
This is not necessarily where the concept is \emph{defined}; in particular, a number of mathematical concepts are explained most fully in one of the appendices.
In such cases, the index can help identify the location of definitions.

\tablehead{\textbf{Symbol} & \textbf{Definition} & \textbf{First use} \\ \hline}
\tabletail{\hline \multicolumn{3}{r}{\emph{Italics indicate less important usages of conflicting notation.}}\\}
\begin{supertabular}{lll}
$\mb{0}$ & Zero vector & \pageref{not:zero} \\
$A$ & Set of links & \pageref{not:A} \\
$A(\mb{x})$ & Set of active constraints at $\mb{x}$ & \pageref{not:Ax} \\
$\mb{A}$ & Constraint coefficient matrix in a linear program & \pageref{not:Alp} \\
$\mb{A_i}$ & The $i$-th column of constraint matrix $\mb{A}$ & \pageref{not:Ai} \\
$A^r$ & Set of links used by trips from origin $r$ & \pageref{not:Ar} \\
$\hat{A}^{rs}$ & Set of allowable links from $r$ to $s$ & \pageref{not:Ahatrs} \\
$\hat{A}^r$ & \emph{Set of unused shortest path links from $r$} & \pageref{not:Ahatr} \\
$A_s$ & \emph{Attractions to destination $s$} & \pageref{not:As} \\
$AEC$ & Average excess cost & \pageref{not:aec} \\
$B$ & Set of basic columns in a linear program & \pageref{not:Blp} \\
$\mathcal{B}$ & Bush (rooted at $r$) & \pageref{not:Bush} \\
$\mathcal{B}$ & Behavior rule in dynamic assignment & \pageref{not:fancyB} \\
$B_r(\mb{x})$ & Ball of radius $r$ centered at $\mb{x}$ & \pageref{not:B} \\
$C$ & Matrix of time-dependent path travel times/costs & \pageref{not:Ctimes} \\
$C$ & \emph{Generic curve in $(x,t)$ space} & \pageref{not:Ccurve} \\
$C$ & \emph{Signal cycle length} & \pageref{not:C} \\
$\mb{C}$ & Vector of path costs in terms of path flows & \pageref{not:Ch} \\
$C(t)$ & \emph{Generalized cost when reaching destination at $t$} & \pageref{not:Cvickrey} \\
$CS$ & Consumer surplus & \pageref{not:CS} \\
$D$ & \emph{Largest value in OD matrix} & \pageref{not:Dmax} \\
$D_r(t)$ & Total demand starting at $r$ in interval $t$ & \pageref{not:Drt} \\
$D^{rs}$ & Demand function for OD pair $(r,s)$ & \pageref{not:Drs} \\
$D_i$; $D_{ij}$ & Cost derivative label for node $i$/link $(i,j)$ & \pageref{not:Dij} \\
$E[\cdot]$ & Expected value & \pageref{not:E} \\
$\mb{E}$ & Coefficient matrix for integer variables & \pageref{not:Eilp} \\
$F$ & Generic multifunction & \pageref{not:Fmulti} \\
$G$ & The network itself & \pageref{not:G} \\
$G_i$; $G_{hij}$ & \emph{Green time for signal approach/turning movement} & \pageref{not:Gi} \\
$H$ & Set of (static) feasible assignments & \pageref{not:H} \\
$H$ & Hessian matrix & \pageref{not:Hessian} \\
$H$ & Matrix of time-dependent path flows & \pageref{not:Hdta} \\
$H^*$ & All-or-nothing ``target'' path flow matrix & \pageref{not:Hstar} \\
$\hat{H}$ & Dynamic user equilibrium path flow matrix & \pageref{not:Hhat} \\
$\mb{H}$ & Vector of path flows in terms of path costs & \pageref{not:Hc} \\
$\mc{H}$ & Set of (dynamic) feasible assignments & \pageref{not:H} \\
$\mc{H}^*$ & Set of previously-generated $\mc{H}$ matrices & \pageref{not:fancyHstar} \\
$\mb{I}$ & Identity matrix & \pageref{not:Identity} \\
$I$ & Set of alternatives & \pageref{not:I} \\
$\mc{I}$ & Index set & \pageref{not:fancyI} \\
$J$ & Jacobian matrix & \pageref{not:Jacobian} \\
$J$ & \emph{Parameter in Ak\c{c}elik delay function} & \pageref{not:J} \\
$\mc{J}$ & Index set & \pageref{not:fancyJ} \\
$K$ & Generic compact and convex set & \pageref{not:K} \\
$K_{ij}$ & \emph{Per-dollar capacity improvement on link $(i,j)$} & \pageref{not:Kij} \\
$K_{rs}$ & \emph{Proportionality constant in entropy maximization} & \pageref{not:Krs} \\
$L$ & Link length & \pageref{not:L} \\
$\mc{L}$ & Lagrangian function & \pageref{not:fancyL} \\
$L_i^r$ & Shortest path label from origin $r$ to node $i$ & \pageref{not:Lir} \\
$L_i$; $L_{ij}$ & Shortest path label for node $i$/link $(i,j)$ & \pageref{not:Lij} \\
$L_{ij}$; $\mb{L}$ & Link likelihood in Dial's method; matrix of all & \pageref{not:Lijlik} \\
$M$ & Large positive constant & \pageref{not:M} \\
$M_i$; $M_{ij}$ & Average cost label for node $i$/link $(i,j)$ & \pageref{not:Mij} \\
$MEC$ & Maximum excess cost & \pageref{not:mec} \\
$N$ & Set of nodes & \pageref{not:N} \\
$N$ & Set of nonbasic columns in a linear program & \pageref{not:NBlp} \\
$N$ & \emph{Number of travelers in Vickrey's bottleneck} & \pageref{not:Nvickrey} \\
$\mc{N}$ & Network loading mapping & \pageref{not:fancyN} \\
$N^\uparrow$ & Cumulative entries to a link & \pageref{not:Nup} \\
$N^\uparrow_{hi,\pi}(t)$ & Cumulative entries to $(h,i)$ on path $\pi$ at time $t$ & \pageref{not:Huphipit} \\
$N^\downarrow$ & Cumulative exits from a link & \pageref{not:Ndown} \\
$N^\downarrow_{hi,\pi}(t)$ & Cumulative exits from $(h,i)$ on path $\pi$ at time $t$ & \pageref{not:Huphipit} \\
$N(x,t)$ & Cumulative vehicles passing point $x$ by time $t$ & \pageref{not:Nxt} \\
$N_{hij}$ & \emph{Number of lanes for turning movement $[h,i,j]$} & \pageref{not:Nhij} \\
$N_{rs}$ & \emph{Intrazonal trips between $r$ and $s$} & \pageref{not:Nrs} \\
$O$ & Order of magnitude upper bound & \pageref{not:bigO} \\
$P_{ij}$ & Link constant for Markov property & \pageref{not:Pij} \\
$P_r$ & Productions from origin $r$ & \pageref{not:Pr} \\
$Q(k)$ & Fundamental (flow-density) diagram & \pageref{not:Q} \\
$Q(t)$ & \emph{Queue delay experienced if reaching destination at $t$}  & \pageref{not:Qvickrey} \\
$\mb{Q}$ & \emph{Matrix of $x^2$ coefficients in quadratic} & \pageref{not:Qbold} \\
$R$ & Receiving flow for a link & \pageref{not:R} \\
$\tilde{R}$ & Remaining receiving flow & \pageref{not:Rtwiddle} \\
$\bbr$ & Set of all real numbers & \pageref{not:bbr} \\
$\bbr_+$ & Set of non-negative real numbers & \pageref{not:bbrp} \\
$\bbr^n$ & Set of $n$-dimensional real vectors & \pageref{not:bbrn} \\
$S$ & Sending flow for a link & \pageref{not:S} \\
$\mc{S}$ & Generalized cost mapping in dynamic assignment & \pageref{not:fancyS} \\
$\tilde{S}$ & Remaining sending flow & \pageref{not:Stwiddle} \\
$\hat{S}$ & Delayed sending flow & \pageref{not:Shat} \\
$S_{hi,\pi}$ & Amount of $(h,i)$ sending flow from path $\pi$ & \pageref{not:Shipit} \\
$SEL$ & Scan eligible list & \pageref{not:SEL} \\
$\cdot^T$ & Matrix transpose & \pageref{not:transpose} \\
$T$ & Length of analysis period & \pageref{not:T} \\
$T$ & \emph{Travel time in terms of flow} & \pageref{not:Tminor} \\
$T^\uparrow_{ij}(n)$ & Time when vehicle $n$ entered link $(i,j)$ & \pageref{not:Tupijn} \\
$T^\downarrow_{ij}(n)$ & Time when vehicle $n$ exited link $(i,j)$ & \pageref{not:Tdownijn} \\
$TMF$ & Total misplaced flow & \pageref{not:TMF} \\
$TSC$ & Total system cost & \pageref{not:TSC} \\
$TSTT$ & Total system travel time & \pageref{not:TSTT} \\
$U_{ij}$ & Longest path label for node $i$/link $(i,j)$ & \pageref{not:Uij} \\
$U_i$; $U^\pi$ & Total utility associated with alternative $i$ (path $\pi$) & \pageref{not:Upi} \\
$V_{ab}$; $\mb{V}$ & Logit term for segments from $a$ to $b$; matrix of all & \pageref{not:Vab} \\
$V_i$; $V_m$ & Observed utility for alternative $i$/mode $m$ & \pageref{not:Vm} \\
$\mr{VI}(K,\mb{f})$ & Variational inequality on set $K$ with function $\mb{F}$ & \pageref{not:vikf} \\
$W_i$; $W_{ij}$ & Node and link weights in Dial's method & \pageref{not:Wi} \\
$X$ & Set of feasible link assignments & \pageref{not:X} \\
$\mc{X}$ & Set of generated target flow vectors & \pageref{not:fancyX} \\
$X(\mb{x}, \mc{X})$ & Set of feasible link assignments combining $\mb{x}$ and $\mc{X}$ & \pageref{not:XxX} \\
$X'$ & Restricted set of feasible link assignments & \pageref{not:Xprime} \\
$X$ & \emph{Flow in terms of travel time} & \pageref{not:Xminor} \\
$X_i$ & \emph{Degree of saturation for approach $i$} & \pageref{not:Ximinor} \\
$Y_i$ & Set of variables satisfying inequality constraint $i$ & \pageref{not:Yi} \\
$Z$ & Set of zones & \pageref{not:Z} \\
$\mathbb{Z}$ & Set of integers & \pageref{not:bbz} \\
$\mathbb{Z}_+$ & Set of non-negative integers & \pageref{not:bbz} \\
$Z^2$ & Set of all OD pairs & \pageref{not:Z2} \\
$Z^2(\zeta)$ & Set of OD pairs using pair of alternate segments $\sigma$ & \pageref{not:Z2zeta} \\
$Z_j$ & \emph{Set of variables satisfying equality constraint $j$} & \pageref{not:Zj} \\
$Z_\zeta$ & Relevant origins for pair of alternate segments $\zeta$ & \pageref{not:Zzeta} \\
$a$ & \emph{Lower bound for one-dimensional optimization} & \pageref{not:a} \\
$a_{ij}$ & Coefficient of $x_i$ in the $j$-th equality constraint & \pageref{not:aij} \\
$a_{ij}$ & \emph{Coefficient of $x_{ij}$ in linear delay function} & \pageref{not:aij2} \\
$b$ & \emph{Upper bound for one-dimensional optimization} & \pageref{not:b} \\
$\mb{b}$ & Vector of right-hand sides of constraints & \pageref{not:blp} \\
$\mb{b}$ & \emph{Vector of $x$ coefficients in quadratic} & \pageref{not:bbold} \\
$b_{ij}$ & Constant term in linear delay function & \pageref{not:bij} \\
$b_j$ & Right-hand side of equality constraint $j$ & \pageref{not:bj} \\
$\mb{c}$ & Vector of objective coefficients in a linear program & \pageref{not:clp} \\
$\mb{c_B}$ & Vector of basis objective coefficients & \pageref{not:cB} \\
$c^\pi$; $\mb{c}$ & Travel time on path $\pi$; vector of all & \pageref{not:cpi} \\
$\hat{c}^\pi$ & Effective travel time on path $\pi$ & \pageref{not:chatpi} \\
$c^\sigma$ & Travel time on segment $\sigma$ & \pageref{not:csigma} \\
$c_{ij}$ & Generalized cost of link $(i,j)$ & \pageref{not:cij} \\
$c_{ij}(t)$ & Generalized for link $(i,j)$ if entering at time $t$ & \pageref{not:cijt} \\
$c^0_{ij}$ & Constant for link $(i,j)$ for allowabilty & \pageref{not:c0ij} \\
$\bar{c}_j$ & Reduced cost for variable $x_j$ in a linear program & \pageref{not:cj} \\
$\mb{d}$ & Direction vector & \pageref{not:dbold} \\
$\mb{d}$ & Objective coefficients for integer variables & \pageref{not:dilp} \\
$d_{hij}$ & Delay associated with turning movement $[h,i,j]$ & \pageref{not:dhij} \\
$d^{rs}$ & Demand from origin $r$ to destination $s$ & \pageref{not:drs} \\
$d^{rs}_t$ & Demand from $r$ to $s$ departing at time $t$ & \pageref{not:drst} \\
$\overline{d}^{rs}$ & Upper bound on demand from $r$ to $s$ & \pageref{not:dbarrs} \\
$\hat{d}^{rs}$; $\mb{\hat{d}}$ & Demand from $r$ to $s$ at elastic demand equilibrium& \pageref{not:hatdrs} \\
$f$ ($\mb{f}$) & Generic (vector-valued) objective or other function  & \pageref{not:f} \\
$f: X \rightarrow Y$ & Function with domain $X$ taking values in $Y$ & \pageref{not:fXY} \\
$f: X \rightrightarrows Y$ & Multifunction from domain $X$ to subsets of $Y$ & \pageref{not:fXmY} \\
$f'$ & Derivative of single-variable function $f$ & \pageref{not:fprime} \\
$f'$ & Directional derivative of multi-variable function $f$ & \pageref{not:ftotder} \\
$\hat{f}$ & Optimal or equilibrium value of a function $f$ & \pageref{not:fhat} \\
$\bar{f}$, $\ubar{f}$ & Upper and lower bounds on $\hat{f}$ & \pageref{not:fbar} \\
$f(t)$ & Schedule delay cost if arriving at $t$ & \pageref{not:fdelay} \\
$g$ & Difference in path travel times & \pageref{not:g} \\
$g_a$; $g_r$ & Acceptance and rejection gap & \pageref{not:ga} \\
$g_i^s$ & Estimate of cost from node $i$ to destination $s$ & \pageref{not:gis} \\
$g(\mb{x})$ & \emph{Generic inequality constraint} & \pageref{not:gcon} \\
$h^\pi$; $\mb{h}$ & Flow on path $\pi$; vector of all path flows & \pageref{not:hpi} \\
$h^\pi_t$ & Flow on path $\pi$ departing at time $t$ & \pageref{not:hpit} \\
$\mb{\hat{h}}$ & Equilibrium path flow vector & \pageref{not:hhat} \\
$h(\mb{x})$ & \emph{Generic equality constraint} & \pageref{not:hcon} \\
$(h,i)$ & Index for a link entering node $i$ & \pageref{not:hi} \\
$[h,i,j]$ & Index for a turning movement from $(h,i)$ to $(i,j)$ & \pageref{not:hij} \\
$i$ & Index for a node & \pageref{not:i} \\
$i$ & \emph{Index for an inequality constraint} & \pageref{not:icon} \\
$i : t$ & Time-expanded node; physical node $i$ at time $t$ & \pageref{not:it} \\
$(i,j)$ & Index for a link & \pageref{not:ij} \\
$\inf$ & Infimum & \pageref{not:inf} \\
$j$ & Index for a node & \pageref{not:j}\\
$j$ & \emph{Index for an equality constraint} & \pageref{not:jcon} \\
$k$ & Iteration counter & \pageref{not:k} \\
$k$ & Density at a point & \pageref{not:kdensity} \\
$k_c$ & Critical density & \pageref{not:kc} \\
$k_j$ & Jam density & \pageref{not:kj} \\
$(k,\ell$) & Index for a link & \pageref{not:kl} \\
$\ell_i^s$ & Length of shortest path from $i$ to $s$ & \pageref{not:lis} \\
$\log$ & \textbf{Natural} logarithm \footnotesize{(We do not use  base-10 logarithms.)} & \pageref{not:log} \\
$\operatorname{med}$ & Median & \pageref{not:med} \\
$m$ & Number of links & \pageref{not:mlinks} \\
$m$ & Number of constraints in a linear program & \pageref{not:mlp} \\
$m$ & \emph{Index for a mode of travel} & \pageref{not:m} \\
$n$ & Number of nodes & \pageref{not:n} \\
$n$ & Number of decision variables in a linear program & \pageref{not:nlp} \\
$n$ & Index for a vehicle & \pageref{not:nveh} \\
$n$ & \emph{Number of some other object} & \pageref{not:nn} \\
$n(x,t)$ & Number of vehicles in cell $x$ at time $t$ & \pageref{not:nxt} \\
$o(i)$ & Topological order of node $i$ & \pageref{not:oi} \\
$p$ & Probability of some event & \pageref{not:p} \\
$p$ & Number of integer variables & \pageref{not:pilp} \\
$p_{hij}$ & Fraction of flow on $(h,i)$ turning onto $(i,j)$ & \pageref{not:phij} \\
$p_i$ & Probability of choosing alternative $i$ & \pageref{not:p_i} \\
$\mr{proj}_K(\mb{x})$ & Projection of vector $\mb{x}$ onto set $K$ & \pageref{not:proj} \\
$q$ & Flow at a point & \pageref{not:q} \\
$q_i^r; \mb{q}^r$ & Backnode label for node $i$ from origin $r$; vector of all & \pageref{not:qir} \\
$q_{max}$ & Capacity (in traffic flow theoretic sense) & \pageref{not:qmax} \\
$q_{max}^\uparrow$ & Capacity at upstream link end & \pageref{not:qarrow} \\
$q_{max}^\downarrow$ & Capacity at downstream link end & \pageref{not:qarrow} \\
$q_{max}^{hij}$ & Oriented capacity for turning movement $[h,i,j]$ & \pageref{not:qmaxhij} \\
$r$ & Index for an origin & \pageref{not:r}\\
$\hat{r}$ & Origin for which demand changes & \pageref{not:rhat} \\
$(r,s)$ & Index for an OD pair & \pageref{not:rs} \\
$s$ & Index for an destination & \pageref{not:s}\\
$s$ & \emph{Bottleneck capacity in Vickrey's bottleneck} & \pageref{not:svickrey} \\
$s_i$; $s_{hij}$ & \emph{Saturation flow for approach/turning movement} &  \pageref{not:si} \\
$\hat{s}$ & Destination for which demand changes & \pageref{not:shat} \\
$t$ & Travel time & \pageref{not:t}\\
$t$ & Time index & \pageref{not:tindex} \\
$t^*$ & Preferred arrival time & \pageref{not:tstar} \\
$t_0$ & Departure time & \pageref{not:t0} \\
$t_c$ & Critical gap & \pageref{not:tc} \\
$t_f$ & Follow-up gap & \pageref{not:tf} \\
$t_S$; $t_E$ & Start and end times of congestion in Vickrey model & \pageref{not:tS} \\
$t_{ij}$; $\mb{t}$ & Travel time on link $(i,j)$; vector of all link times & \pageref{not:tij}\\
$t'_{ij,x}$; $t'_{ij,y}$ & Partial derivatives of $t_{ij}$ with respect to $x$ and $y$ & \pageref{not:tijx} \\
$\hat{t}_{ij}$ & Marginal cost on link $(i,j)$ & \pageref{not:thatij} \\
$\tilde{t}_{ij}$ & Transformed cost for link $(i,j)$ & \pageref{not:ttwiddleij} \\
$t^0_{ij}$ & Free-flow travel time on link $(i,j)$ & \pageref{not:t0ij}\\
$\mb{t}^\rightarrow$ & Vector of times on Gartner links & \pageref{not:tarrow} \\
$u$ & Vehicle speed at a point & \pageref{not:u} \\
$u_{AB}$ & Speed of shockwave separating $A$ and $B$ & \pageref{not:uAB} \\
$u_f$ & Free-flow speed & \pageref{not:uf} \\
$u_{ij}$ & Capacity on link $(i,j)$ & \pageref{not:uij}\\
$\mr{vect}$ & Vector formed by specified components & \pageref{not:vect} \\
$x$ & Spatial coordinate & \pageref{not:xcoord} \\
$x, \mb{x}$ & Generic decision or other variable; scalar and vector & \pageref{not:x} \\
$\hat{x}, \mb{\hat{x}}$ & Equilibrium or optimum value of variable(s) & \pageref{not:xhat} \\
$\mb{x_B}$, $\mb{x_N}$ & Basic/nonbasic variables in a linear program & \pageref{not:xBlp} \\
$x_{ij}$; $\mb{x}$ & Flow on link $(i,j)$; vector of all link flows & \pageref{not:xij}\\
$x^*_{ij}$; $\mb{x^*}$ & Target flow on link $(i,j)$; vector of all & \pageref{not:xstar} \\
$\mb{x'}$ & Feasible link flows for restricted set & \pageref{not:xprime} \\
$\mb{\hat{x}'}$ & Restricted equilibrium link flows & \pageref{not:xprimehat} \\
$x_{ij}^r$; $\mb{x}^r$ & Flow on link $(i,j)$ from origin $r$; vector of all & \pageref{not:xijr} \\
$\hat{x}_{ij}^r$; $\mb{\hat{x}}^r$ & Equilibrium flow on $(i,j)$ from origin $r$; vector of all & \pageref{not:xhatijr} \\
$\mb{x}^\rightarrow$ & Vector of link flows on Gartner links & \pageref{not:xarrow} \\
$\mb{\tilde{x}}$ & Candidate equilibrium solution & \pageref{not:xtwiddle} \\
$y$ & Auxiliary decision variable & \pageref{not:ylp} \\
$y$ & \emph{Generic vertical coordinate} & \pageref{not:y} \\
$y(x,t)$ & Number of vehicles entering cell $x$ at time $t$ & \pageref{not:yxt} \\
$y_{hij}(t)$ & Flow on turning movement $[h,i,j]$ during interval $t$ & \pageref{not:yhij} \\
$y_{hi}(t)$ & Total flow leaving $(h,i)$ during interval $t$ & \pageref{not:yhi} \\
$y_{ij}$ & \emph{Link improvement parameter in network design} & \pageref{not:yndp} \\
$z$ & Generic notation for objective function & \pageref{not:z} \\
$\Gamma(i)$ & Forward star for node $i$ & \pageref{not:Gammai} \\
$\Gamma^{-1}(i)$ & Reverse star for node $i$ & \pageref{not:Gamma-1i} \\
$\mb{\Delta}$ & Link-path adjacency matrix & \pageref{not:Delta} \\
$\Delta_{ri}$ & Shorthand for right-hand side of~\eqn{ddbushflowc} and~\eqn{ddbushflowcorigin}& \pageref{not:Deltari}  \\
$\Delta H$ & Set of feasible path flow shifts & \pageref{not:DeltaH} \\
$\Delta N$ & Increment of vehicle flow & \pageref{not:DeltaN} \\
$\Delta h$ & Flow shifted between paths & \pageref{not:Deltah} \\
$\Delta t$ & Time step & \pageref{not:Deltat} \\
$\Delta x$ & Spatial discretization & \pageref{not:Deltax} \\
$\overline{\Delta h}^r$ & Maximum flow shift from origin $r$ & \pageref{not:Deltahath} \\
$\Theta$ & Weighting factor among objective terms & \pageref{not:Theta} \\
$\Lambda$ & \emph{Duration of analysis period} & \pageref{not:Lambda} \\
$\Lambda^r_i$ & Derivative of $L^r_i$ in sensitivity analysis & \pageref{not:Lambdari} \\
$\Xi(i)$ & Set of turning movements for node $i$ & \pageref{not:Xii} \\
$\Pi$ & Set of all acyclic paths & \pageref{not:Pi} \\
$\Pi^{rs}$ & Set of acyclic paths from $r$ to $s$ & \pageref{not:Pirs} \\
$\Pi_{[h,i,j]}$ & Set of paths including turn movement $[h,i,j]$ & \pageref{not:Pihij} \\
$\hat{\Pi}^{rs}$ & Set of used paths from $r$ to $s$ & \pageref{not:Pihatrs} \\
$\Sigma_{ab}$ & Set of segments from $a$ to $b$ & \pageref{not:Sigmaab} \\
$\Phi$ & Standard normal cumulative distribution function & \pageref{not:Phi} \\
$\Omega$ & Set of active turning movements & \pageref{not:Omega} \\
$\alpha$ & Shape parameter in link performance function & \pageref{not:alpha} \\
$\alpha$ & \emph{Weighting factor in conjugate Frank-Wolfe} & \pageref{not:alphacfw} \\
$\alpha_r$ & \emph{Trip balancing parameter for origin $r$} & \pageref{not:alphar} \\ 
$\alpha_{hi}$ & Fraction of node $i$ flow entering from $h$ & \pageref{not:alphahi} \\
$\alpha_{hij}$ & Turning movement flow increment & \pageref{not:alphahij} \\
$\alpha_{hij,s}^t$ & Fraction of flow headed for $s$ leaving $(h,i)$ on $(i,j)$ at $t$ & \pageref{not:alphahijst} \\
$\beta$ & Shape parameter in link performance function & \pageref{not:beta} \\
$\beta$ & Penalty factor for early arrival & \pageref{not:betasd} \\
$\beta_{ij}$ & Link $(i,j)$ Lagrange multiplier in max.\ entropy & \pageref{not:betaij} \\
$\beta_{ij}$ & Link priority parameter for merge or conflict & \pageref{not:betaconflict} \\
$\beta_s$ & \emph{Trip balancing parameter for destination $s$} & \pageref{not:betas} \\ 
$\gamma$ & Relative gap & \pageref{not:gammagap} \\
$\gamma$ & Penalty factor for late arrival & \pageref{not:gammasd} \\
$\gamma$ & \emph{Destination choice parameter} & \pageref{not:gamma} \\
$\gamma$ & \emph{Shorthand for $(2\beta - 1)/(2 \beta - 2)$ in conical function} & \pageref{not:gammavdf} \\
$\gamma_S$ & Smith gap & \pageref{not:gammaS} \\
$\gamma_{rs}$ & OD pair $(r, s)$ Lagrange multiplier in max.\ entropy & \pageref{not:gammars} \\
$\delta_{ij}^\pi$ & Number of times link $(i,j)$ is in path $\pi$ (typ.\ 0 or 1) & \pageref{not:deltaij} \\
$\epsilon$ & Small positive number & \pageref{not:epsilon} \\
$\epsilon$; $\epsilon_i$ & Unobserved utility (for alternative $i$) & \pageref{not:epsilonutility} \\
$\zeta$ & Shorthand for $f$ in terms of step size & \pageref{not:zeta} \\
$\zeta$ & Pair of alternate segments & \pageref{not:zetapas} \\
$\theta$ & Largest possible step size & \pageref{not:thetastep} \\
$\theta$ & Gumbel distribution parameter & \pageref{not:theta} \\
$\theta$ & \emph{Generic angle} & \pageref{not:thetaang} \\
$\theta$ & \emph{Golden section constant} $(3 - \sqrt{5})/2$  & \pageref{not:thetags} \\
$\iota$ & Karush-Kuhn-Tucker multiplier & \pageref{not:iota} \\
$\kappa$ & Lagrange multiplier & \pageref{not:kappa} \\
$\kappa$ & Generic shortest-path time for an OD pair & \pageref{not:kappagen} \\
$\kappa_{rs}$ & Shortest-path time/cost from $r$ to $s$ & \pageref{not:kappars} \\
$\lambda$ & Weighting factor in a convex combination & \pageref{not:lambda} \\
$\mu$ & Step size when moving in a direction & \pageref{not:mu} \\
$\mu$ & Mean of a random variable & \pageref{not:mumean} \\
$\bar{\mu}$ & Greatest possible step size & \pageref{not:mubar} \\
$\xi^r_{ij}$ & Derivative of $x^r_{ij}$ in sensitivity analysis & \pageref{not:xirij} \\
$\pi$ & Index for a path & \pageref{not:pi} \\
$\pi^c$ & Companion path & \pageref{not:pic} \\
$\hat{\pi}_{rs}$ & Basic path for OD pair $(r,s)$ & \pageref{not:pirshat} \\
$\rho^\pi$ & Auxiliary cost for path $\pi$ & \pageref{not:rhopi} \\
$\sigma$ & Path segment & \pageref{not:sigmaseg} \\
$\sigma$ & Standard deviation of a random variable & \pageref{not:sigmarv} \\
$\sigma^\zeta_1$; $\sigma^\zeta_2$ & Segments in pair $\zeta$ & \pageref{not:sigmapas} \\
$\sigma_L$; $\sigma_U$ & Shortest and longest path segments & \pageref{not:sigma} \\
$\tau_{ij}(t)$ & Travel time on link $(i,j)$ if entering at time $t$ & \pageref{not:tauijt} \\
$\phi$ & Fraction of diverge flow that can move & \pageref{not:phidiverge} \\
$\phi$ & \emph{Least-squares function} & \pageref{not:phi} \\
$\phi_{rs}$ & Friction function between $r$ and $s$ & \pageref{not:phirs} \\
$\omega$ & Dummy variable of integration & \pageref{not:omega} \\
$\cap$ & Set intersection & \pageref{not:cap} \\
$\cup$ & Set union & \pageref{not:union} \\
$\emptyset$ & Empty set & \pageref{not:emptyset} \\
! & Factorial & \pageref{not:fact} \\
$\gg$ & Much greater than & \pageref{not:gg} \\
$\leftarrow$ & Assignment (setting a variable, see Appendix~\ref{sec:algorithms}) & \pageref{not:leftarrow} \\
$\nabla f$ & Gradient of multivariate function $f$ & \pageref{not:nabla} \\
$\oplus$ & Segment concatenation & \pageref{not:oplus} \\
$\perp$ & Orthogonality (zero dot product) & \pageref{not:perp} \\
$\prod$ & Product & \pageref{not:prod} \\
$\subseteq$ & Subset with possible equality & \pageref{not:subseteq} \\
$\supseteq$ & Superset with possible equality & \pageref{not:superseteq} \\
$[ \cdot ] $ & Indicator function, 1 if $\cdot$ is true and 0 if false & \pageref{not:indicator} \\
$[\cdot]^+$ & Positive part $\max \{ \cdot, 0 \}$ & \pageref{not:bracketplus} \\
\end{supertabular}

\bibliography{references}

\printindex

\end{document}